\documentclass[10pt,a4paper]{article}
\usepackage[utf8]{inputenc}
\usepackage[american]{babel}
\usepackage{amsmath,amssymb,amsthm,mathtools}
\usepackage{enumitem}
\usepackage[dvipsnames]{xcolor}
\usepackage{graphicx}
\usepackage{array}
\usepackage[left=2cm,right=2cm,top=2cm,bottom=2cm]{geometry}
\usepackage{booktabs,longtable}
\usepackage{pgfplots}
\pgfplotsset{compat=1.12}
\usetikzlibrary{calc,quotes,angles,arrows.meta,positioning,decorations.pathreplacing}
\usepgfplotslibrary{fillbetween} \usepackage{hyperref}
\hypersetup{pdftitle={Leafwise thermalization near vacuum for the inhomogeneous Landau equation with
very soft and Coulomb potentials}}

\numberwithin{equation}{section} 
 \setlist[enumerate]{itemsep=3pt} \setlist[itemize]{itemsep=3pt}

\newtheorem{theo}{Theorem}[section]
\newtheorem{defin}{Definition}[section]
\newtheorem{lem}{Lemma}[section]
\newtheorem{prop}{Proposition}[section]
\newtheorem{remark}{Remark}[section]
\newtheorem{cor}{Corollary}[section]
\newtheorem{hypothesis}{Inputs}[section]
\newcounter{rq}\renewcommand{\therq}{R\arabic{rq}}

\newcommand{\R}{\mathbb{R}}
\newcommand{\C}{\mathbb{C}}
\newcommand{\eps}{\varepsilon}
\newcommand{\dd}{\,\mathrm{d}}
\newcommand{\jap}[1]{\langle #1\rangle}
\newcommand{\ab}{\bar a}
\newcommand{\bb}{\bar b}
\newcommand{\cb}{\bar c}
\newcommand{\Qf}{Q_{\mathrm{fib}}}
\newcommand{\nn}[1]{|\!|\!|#1|\!|\!|}
\newcommand{\ip}[2]{\langle\!\langle #1,#2\rangle\!\rangle}
\newcommand{\mom}{\boldsymbol{\ell}}
\newcommand{\tn}{\widetilde\nabla}
\newcommand{\Dc}{\mathfrak D}
\newcommand{\Nn}{\mathfrak N}
\newcommand{\cQ}{\mathcal Q}
\newcommand{\rad}{\mathfrak r}

\title{Leafwise thermalization near vacuum for the inhomogeneous Landau\\ equation with very soft
and Coulomb potentials} \author{Xingyu Li} \date{}
\begin{document}

\maketitle

\begin{abstract}
We study the inhomogeneous Landau equation with very soft and Coulomb potentials,
$\gamma\in[-3,-2]$, for data near vacuum. In this range the collisions along a free
trajectory have a divergent time integral, and the dynamics is organized along \emph{leaves}, the
velocity distributions in the variable $z=x-tv$: on
each leaf the solution follows a homogeneous Landau flow with a leaf-dependent collision rate, first
far from equilibrium and then near a Maxwellian whose parameters approximately solve a compressible
Euler system with the free-expansion prefactor $t^{-2}$.

For data of size $\eps$, analytic in velocity and Gaussian in space, we prove global existence of
classical solutions, relaxation of every leaf (not uniform in $x/t$), and convergence of the leaf
parameters at rate $t^{-1}$. The solution scatters in $L^1$ to the free transport of a thermalized
profile that generally remains at a fixed relative distance from the data. The $L^1$ distance to the
scattering state decays like $(1+\eps\log t)^{-(p-3)/p}$ at $\gamma=-2$ and $(1+\eps
t^{|\gamma|-2})^{-(p-3)/p}$ for $\gamma<-2$, where $\jap v^{-p}$ is the velocity decay of the data.
These rates are sharp when the spatial profiles stay away from Maxwellians. The relaxation of each
leaf is governed by its collision clock: at $\gamma=-2$ the leaf relaxes like a power of $t$
proportional to its collision rate, and for $\gamma<-2$ like a stretched exponential of its clock. We also prove the
near-Maxwellian and homogeneous results we need beyond the literature, including at $\gamma=-2$.
\end{abstract}

\noindent{\small\textit{2020 Mathematics Subject Classification.} 35Q20, 35Q84, 82C40, 35B40.}

\setcounter{tocdepth}{2}
\tableofcontents

\section{Introduction}\label{sec:intro}

\subsection{The Landau equation near vacuum}

The Landau equation, introduced by Landau in 1936 for a plasma of charged particles \cite{Lan36},
describes a dilute gas in which collisions are dominated by small deflections. It reads
\begin{equation}\label{eq:landau}
\partial_tf+v\cdot\nabla_xf=Q(f,f),\qquad f(0,x,v)=f_0(x,v)\ge0,
\end{equation}
for the density $f=f(t,x,v)$ of particles at position $x\in\R^3$ with velocity $v\in\R^3$. The
collision operator acts only in $v$:
\begin{equation}\label{eq:Q}
Q(f,h)=\nabla_v\cdot\int_{\R^3}a(v-w)\big(f(w)\nabla_vh(v)-h(v)\nabla_wf(w)\big)\dd w,\qquad
a(u)=|u|^{\gamma+2}\Big(I-\frac{u\otimes u}{|u|^2}\Big).
\end{equation}
The matrix $a(u)$ is the projection onto $u^\perp$, multiplied by $|u|^{\gamma+2}$. The parameter
$\gamma$ describes the interaction between particles: hard potentials $\gamma\ge0$, moderately soft
potentials $-2<\gamma<0$, very soft potentials $-3<\gamma<-2$, and the Coulomb case $\gamma=-3$,
which is the physically most important one \cite{LP81,Vil02}. The Landau equation is also the
grazing limit of the Boltzmann equation \cite{DLD92,Des92,AV04}. We consider $\gamma\in[-3,-2]$.

Most of the mathematical theory of \eqref{eq:landau} concerns two regimes. Near a global Maxwellian
there is a well-developed perturbation theory, including for very soft and Coulomb potentials
\cite{Guo02,SG08,CTW,CM,KGH20,DLSS21}. For the spatially homogeneous equation, in which $f$ does not
depend on $x$, the theory of weak and smooth solutions, entropy dissipation and convergence to
equilibrium is extensive \cite{Vil98,DV00a,DV00b,D15,CDH,GIMV19,GGIV22,DHJ24,GS}. Recently, Guillen
and Silvestre proved that the homogeneous Landau--Coulomb equation does not blow up \cite{GS}. The
inhomogeneous equation with large data and soft potentials is far less understood, and we refer
to \cite{CSS,HST20,HS20,SS23} for local theory, conditional regularity and continuation criteria.
Golding, Henderson and Silvestre \cite{GHS26} recently proved a continuation criterion in a weighted
$L^\infty$ norm, without bounds on the hydrodynamic quantities, for the inhomogeneous Landau and
Boltzmann equations. For soft potentials it essentially rules out the hydrodynamic blow-up mechanism
(see also \cite{GH25} for the Coulomb case). For very hard potentials, by contrast, the
inhomogeneous equation develops finite-time singularities, through a self-similar implosion of the
compressible Euler equations \cite{BCGJVY26}. For the Landau--Fermi--Dirac equation with Coulomb
potential, global smooth solutions with large data are constructed in \cite{GH26}, using the
$L^\infty$ bound that comes from the Pauli exclusion principle.

\begin{figure}[ht]
\centering
\begin{tikzpicture}[scale=0.9]
\begin{scope}
\draw[->] (-2.2,0)--(2.4,0) node[right]{$x$};
\draw[->] (0,-1.6)--(0,1.8) node[above]{$v$};
\fill[blue!15] (-0.8,-1)--(0.8,-1)--(0.8,1)--(-0.8,1)--cycle;
\draw[blue!70!black,thick] (-0.8,-1)--(0.8,-1)--(0.8,1)--(-0.8,1)--cycle;
\node at (0,-2.0) {$t=0$};
\end{scope}
\begin{scope}[xshift=6.5cm]
\draw[->] (-4.2,0)--(4.4,0) node[right]{$x$};
\draw[->] (0,-1.6)--(0,1.8) node[above]{$v$};
\fill[blue!15] (-3.8,-1)--(-2.2,-1)--(3.8,1)--(2.2,1)--cycle;
\draw[blue!70!black,thick] (-3.8,-1)--(-2.2,-1)--(3.8,1)--(2.2,1)--cycle;
\draw[red!70!black,very thick] (1.5,0.233)--(1.5,0.767);
\draw[red!70!black,dashed] (1.5,-1.4)--(1.5,1.5);
\node[red!70!black,above] at (1.5,1.5) {$x$ fixed};
\draw[decorate,decoration={brace,amplitude=3pt}] (1.6,0.767)--(1.6,0.233)
node[midway,right=4pt,font=\small]{$\simeq 2/t$};
\fill (1.5,0.5) circle (1.3pt);
\draw[->,font=\small] (0.6,-0.7) node[below]{$v=x/t$} -- (1.45,0.45);
\node at (0,-2.0) {$t=3$};
\end{scope}
\end{tikzpicture}
\caption{Free transport $f(t,x,v)=f_0(x-tv,v)$ shears the support of $f_0$ (left) along the lines
$x-tv=\mathrm{const}$ (right). At a fixed point $x$ only velocities in an interval of length
$\simeq2/t$ around $x/t$ are present (dispersion). The vertical segment is what we call the
\emph{leaf} through $x$.}\label{fig:phase}
\end{figure}
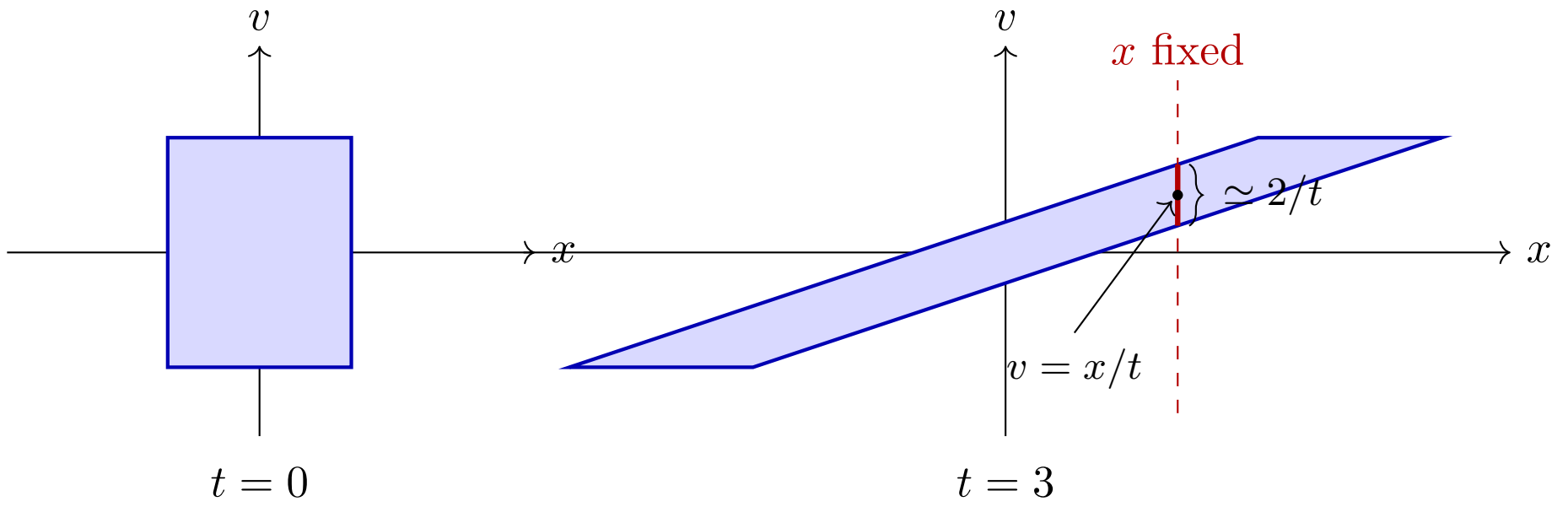

\emph{Near vacuum} is the regime in which $f_0$ is small and localized in $x$. It is the natural
regime for a gas expanding into empty space, and it has a long history for the Boltzmann equation
\cite{IS84,BT85,Ham85,Tos86,Ark11,Alo09}. The underlying mechanism is dispersion: under free
transport $f(t,x,v)=f_0(x-tv,v)$, the density $\int f\dd v$ decays like $t^{-3}$, because the
particles at a given point $x$ come from a region of velocities of volume $\simeq t^{-3}$
(Figure~\ref{fig:phase}).
For the Vlasov--Poisson system the same mechanism gives global small solutions
\cite{BD85,HRV11,Smu16,IPWW22}. For the Landau equation the difficulty is that $Q$ involves
derivatives in $v$, and these grow linearly in $t$ along free transport. Luk \cite{Luk} overcame
this for moderately soft potentials $\gamma\in(-2,0)$, using the vector fields $t\nabla_x+\nabla_v$
that commute with free transport. Chaturvedi then extended the result to hard potentials
\cite{Ch20} and to the Boltzmann equation with moderately soft potentials \cite{Ch21}.

\subsection{The collision clock and the threshold
\texorpdfstring{$\gamma=-2$}{gamma=-2}}\label{ssec:critical}

\begin{figure}[ht]
\centering
\begin{tikzpicture}
\begin{axis}[width=0.62\textwidth,height=0.38\textwidth,xmin=1,xmax=60,ymin=0,ymax=7,
 xlabel={$t$},ylabel={$\tau(t)=\int_1^t s^{-(3+\gamma)}\dd s$},legend pos=north west,legend
style={font=\small},
 tick label style={font=\small}]
\addplot[domain=1:60,samples=200,thick,blue] {(1-x^(-0.5))/0.5};
\addlegendentry{$\gamma=-1.5$}
\addplot[domain=1:60,samples=200,thick,black] {ln(x)};
\addlegendentry{$\gamma=-2$}
\addplot[domain=1:60,samples=200,thick,red!70!black] {(x^(0.5)-1)/0.5};
\addlegendentry{$\gamma=-2.5$}
\addplot[domain=1:60,samples=200,thick,dashed,orange!80!black] {x-1};
\addlegendentry{$\gamma=-3$}
\addplot[domain=1:60,dotted,blue] {2};
\end{axis}
\end{tikzpicture}
\caption{The collision clock $\tau(t)$ for several values of $\gamma$. For $\gamma>-2$ it converges
(dotted line: the limit $2$ for $\gamma=-1.5$), and collisions are a perturbation of free transport.
For $\gamma\le-2$ it diverges, logarithmically at $\gamma=-2$ and like $t^{|\gamma|-2}$ for
$\gamma<-2$.}\label{fig:clock}
\end{figure}
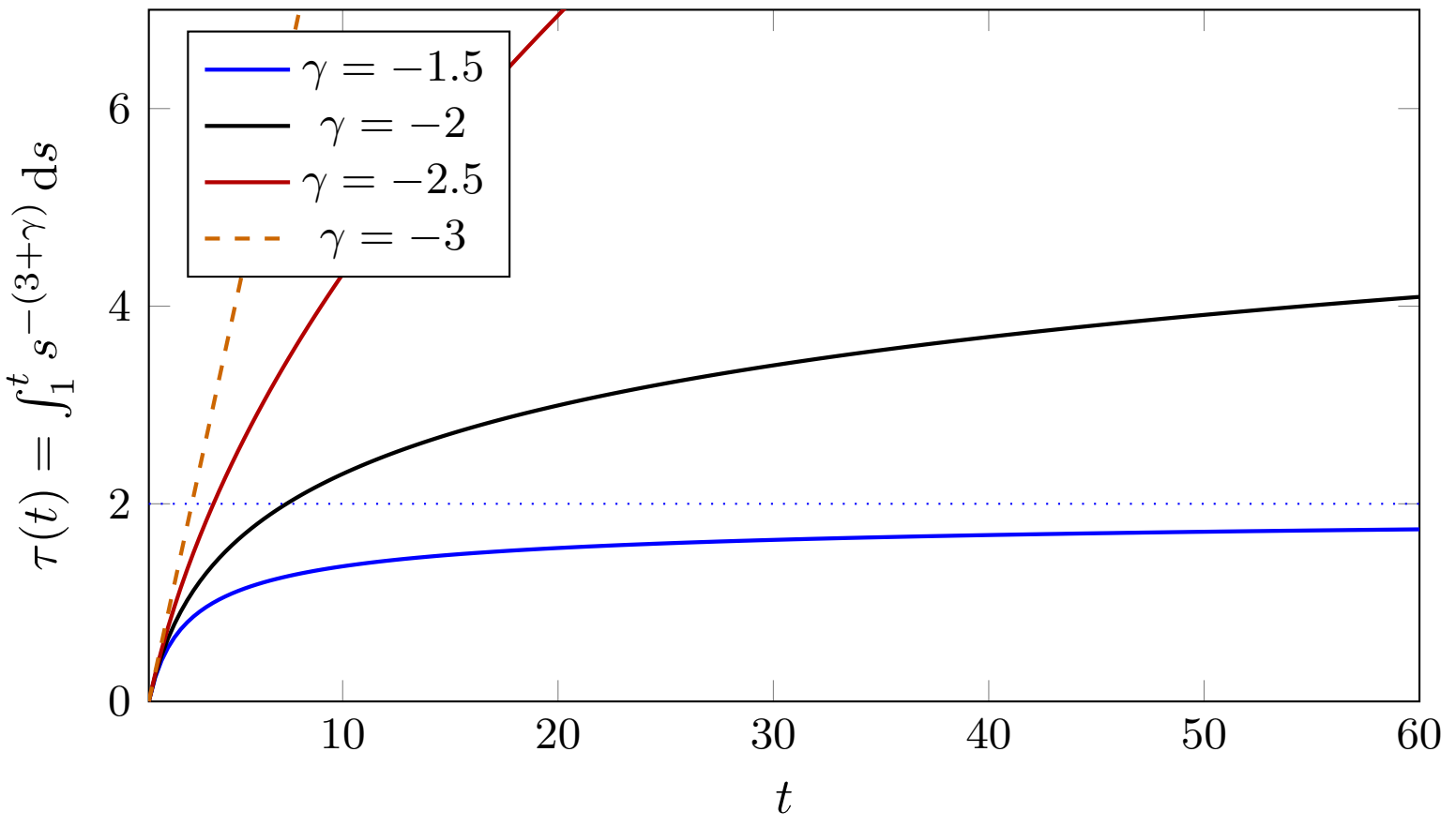

Consider a particle moving freely along $x=y+tv$. The other particles that it meets at time $t$ have
velocities within $O(t^{-1})$ of $x/t$, and their density is $O(t^{-3})$. The relative velocities
are therefore $O(t^{-1})$, and the kernel $a$, which is homogeneous of degree $\gamma+2$,
contributes $t^{-(\gamma+2)}$. Together with the derivative $\nabla_v$ of $Q$, which gives a factor
$t$ on functions that vary on the scale $t^{-1}$ in $v$, we find that the collision operator acts
along free transport with the strength
\begin{equation}\label{eq:kappa}
\kappa(t)=t^{-(3+\gamma)} .
\end{equation}
This heuristic is in fact exact, by the leaf identity below.
The total effect of the collisions up to time $t$ is measured by the \emph{collision clock}
\begin{equation}\label{eq:tau}
\tau(t)=\int_1^t\kappa(s)\dd s=\begin{cases}\dfrac{1-t^{-(\gamma+2)}}{\gamma+2}\ \longrightarrow\
\dfrac1{|\gamma+2|},&\gamma>-2,\\[2mm]\log t,&\gamma=-2,\\[1mm]\dfrac{t^{|\gamma|-2}-1}{|\gamma|-2},&\gamma<-2.\end{cases}
\end{equation}
For $\gamma>-2$ the clock stops: each particle undergoes only a bounded, and for small data small,
amount of collisional deflection over infinite time, and the solution scatters to a solution of free
transport. This regime is treated in \cite{Luk,Ch20,Ch21}. For $\gamma\le-2$ the clock runs forever,
as Figure~\ref{fig:clock} shows. However small the data, after an effective collision time of order
one the collisions have changed the local velocity distributions by an amount of order one. The
solution can then no longer be a perturbation of free transport, and the problem is nonlinear. As
far as we know, no near-vacuum result was available for $\gamma\le-2$.

\subsection{Leaves}\label{ssec:leaves}

For $\gamma\le-2$ we follow \emph{leaves} rather than individual particles.

\emph{Profile variables.} We write $f$ in the variables that are constant along free transport,
\begin{equation}\label{eq:profile}
g(t,y,v):=f(t,y+tv,v),\qquad y=x-tv .
\end{equation}
Then $\partial_tg=Q(f,f)(t,y+tv,v)$: the transport term disappears, and the derivatives transform as
$\nabla_x\mapsto\nabla_y$ and $t\nabla_x+\nabla_v\mapsto\nabla_v$. The vector fields
$t\nabla_x+\nabla_v$ are the ones used by Luk \cite{Luk}, and before by Bardos and Degond
\cite{BD85} and in the vector-field method for Vlasov systems \cite{Smu16,FJS17,Bi20,Wa22}.

\begin{figure}[ht]
\centering
\begin{tikzpicture}[scale=1.0]
\draw[->] (-3.3,0)--(3.6,0) node[right]{$x$};
\draw[->] (0,0)--(0,4.2) node[above]{$t$};
\foreach \s in {-0.8,-0.4,0,0.4,0.8}{
 \draw[gray!60,thin] (0,0)--({3.9*\s},3.9);}
\draw[red!70!black,very thick] (0,0)--(1.56,3.9) node[above]{ray $x=t\xi$};
\foreach \a/\b in {-0.5/0.55,0.3/0.25,0.6/0.18,-0.2/0.45,0.1/0.36}{
 \draw[blue!70!black,thin] (\a,0)--({\a+3.9*\b},3.9);}
\fill[red!70!black] (1.2,3) circle (1.6pt);
\node[font=\small,right] at (1.25,2.85) {$(t,x)$};
\draw[<->,font=\small] (-1.2,0.2)--(-0.5,0.2) node[midway,below]{$z$};
\end{tikzpicture}
\caption{Space-time picture. Rays $x=t\xi$ (gray) and free trajectories $x=y+tv$ (blue). The
trajectories that pass through the point $(t,x)$ start from different $y$, and the leaf at $(t,x)$
records their distribution in $z=y$. As $t\to\infty$ their velocities $v=(x-z)/t$ approach
$\xi=x/t$.}\label{fig:rays}
\end{figure}
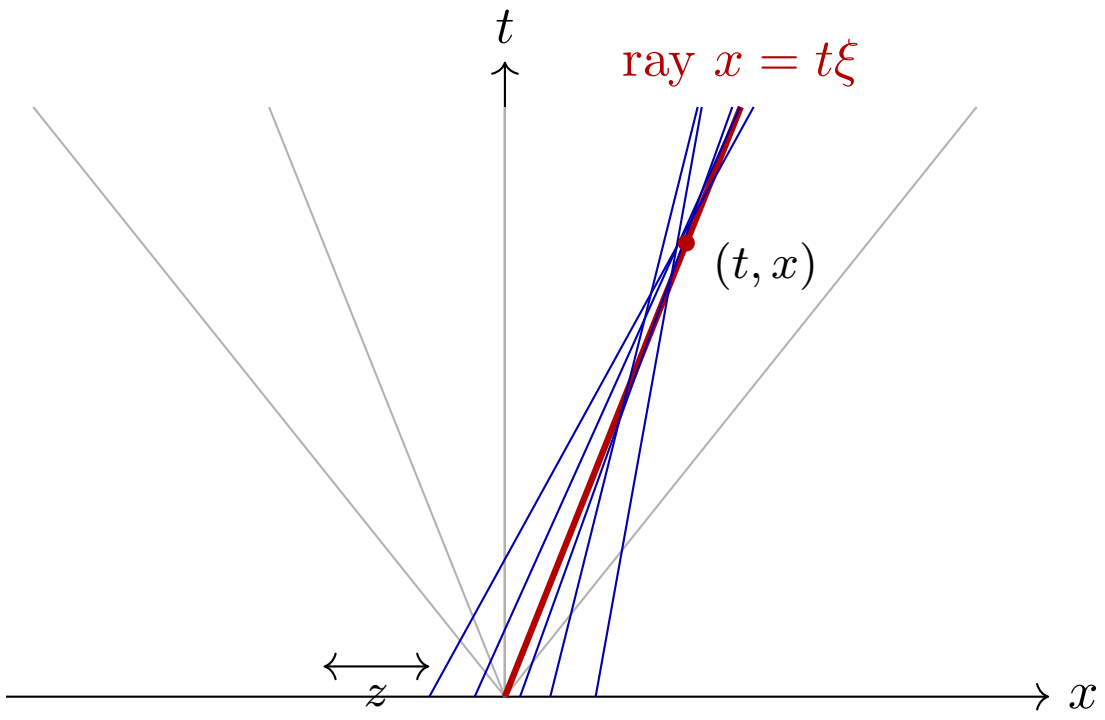

\emph{Leaves.} Fix a time $t$ and a point $x$. The particles at $x$ have velocities $v$, and the
corresponding profile variable is $z=x-tv$. The \emph{leaf} through $x$ is the velocity distribution
at $x$, written in the variable $z$:
\begin{equation}\label{eq:leafdef}
g^{(x)}(z):=g\big(t,z,(x-z)/t\big)=f\big(t,x,(x-z)/t\big),\qquad z\in\R^3 .
\end{equation}
In phase space it is the vertical segment of Figure~\ref{fig:phase}. In the variables $(y,v)$ it is
the affine $3$-plane $\{y+tv=x\}$. It is convenient to label leaves by the \emph{asymptotic
velocity} $\xi=x/t$, so that the leaf $\xi$ moves along the ray $x=t\xi$ in space-time, as in
Figure~\ref{fig:rays}.

\emph{The leaf identity.} On each leaf, the inhomogeneous Landau equation is exactly a homogeneous
Landau equation in the variable $z$, with the collision clock
\eqref{eq:kappa}:
\begin{equation}\label{eq:leafid0}
\partial_tg(t,y,v)=\kappa(t)\,Q\big(g^{(x)},g^{(x)}\big)(y),\qquad x=y+tv .
\end{equation}
We prove it in Proposition~\ref{prop:leaf} by an elementary change of variables in the integral
\eqref{eq:Q}.
The difficulty is that the leaf through a given point $(y,v)$ changes with time:
the left side of \eqref{eq:leafid0} is a derivative at fixed $(y,v)$, not at fixed $(x,z)$. The
difference is a transport along the leaves of size $t^{-2}z\cdot\nabla_\xi$, which we call the
\emph{leafwise transport}. Replacing the leaf through $(y,v)$ by the fiber $\{v=\mathrm{const}\}$,
we obtain the \emph{fibered} equation $\partial_tG=\kappa\,Q(G(\cdot,v),G(\cdot,v))$, a family of
homogeneous Landau equations parametrized by $v$, and we show in Section~\ref{sec:leafgeom} that
the error has relative size $t^{-1}$.

\emph{Example: free transport.} Take product data
$f_0(x,v)=\eps\rho_0(v)F_0(x)$, a special case of our main theorem, and ignore the collisions. Then
$g(t,y,v)=f_0(y,v)$, and by \eqref{eq:leafdef} the leaf through $x=t\xi$ is
\begin{equation}
g^{(x)}(z)=\eps\rho_0\big(\xi-z/t\big)F_0(z)\ \longrightarrow\
\eps\rho_0(\xi)\,F_0(z)\qquad(t\to\infty).
\end{equation}
Thus every leaf carries a copy of the \emph{spatial} profile $F_0$, with the mass $\eps\rho_0(\xi)$.
This example contains the mechanism that we use throughout the paper: free transport turns
the spatial profile of the data into the velocity distribution seen at each point, and the
collisions then act on $F_0$ as if it were a velocity distribution, at a rate proportional to
$\eps\rho_0(\xi)$. Unless $F_0$ is a Maxwellian, every leaf therefore starts, after normalization,
near the distance $d_0=\|F_0-M_0\|_{L^2}$ from equilibrium, which is not assumed
small, and has to relax through the nonlinear homogeneous Landau dynamics. If the spatial profile
depends on the velocity, $f_0=\eps\rho_0(v)F_0(x,v)$, the leaf $\xi$ carries the profile
$F_0(\cdot,\xi)$. Our main theorem covers such data at $t=0$ when $F_0(\cdot,v)$ is a Gaussian with
a rate independent of $v$, times the factor $e^{-\Phi_1(x)}$ of (D$_{\rm iso}$) and a bounded
analytic factor, as in Definition~\ref{def:data}. Theorem~\ref{thm:t1} covers more general profiles
when they are posed, up to a small error, at a positive time $t_1$. Since $\rho_0(\xi)\to0$ as
$|\xi|\to\infty$, leaves with large $|\xi|$ carry little mass, collide slowly, and relax late.

The leaf picture combines several classical ideas. Free transport sends a particle with velocity $v$
along the ray $x/t\to v$. For collisionless systems this leads to the description of the large-time
dynamics in terms of the asymptotic velocity $\xi=x/t$ and of the profile \eqref{eq:profile}. We
refer to the scattering theory for Vlasov--Poisson \cite{BD85,HRV11,IPWW22,FOPW23} and to the
vector-field method \cite{Smu16,FJS17,Bi20}. In this picture the local velocity distribution at $x$
concentrates at $\xi=x/t$, with width $t^{-1}$, which is the leaf.

Leaves are also related to homoenergetic affine flows. For the Boltzmann equation, Truesdell
\cite{Tru56} and Galkin \cite{Gal58} introduced solutions of the form $f(t,x,v)=g(t,v-L(t)x)$ with
$L(t)=(I+tL_0)^{-1}L_0$. For them the inhomogeneous equation reduces exactly to a homogeneous
equation with an additional linear drift. Their existence was proved by Cercignani \cite{Cer89}, and
their self-similar long-time behavior was studied by Bobylev, Caraffini and Spiga \cite{BCS96},
James, Nota and Vel\'azquez \cite{JNV19a,JNV19b,JNV20} and Kepka \cite{Kep23}. Free expansion into
vacuum corresponds to the isotropic case $L(t)=t^{-1}I$, and then $v-L(t)x=-(x-tv)/t=-z/t$: the leaf
variable. Our leaves are local versions of homoenergetic flows: on each ray the gas looks like an
isotropically expanding homoenergetic gas, with parameters that vary slowly from ray to ray.

Hydrodynamics enters once a leaf has relaxed: it is close to a Maxwellian in $z$. The
parameters of these Maxwellians are the local density, mean velocity and temperature of the gas in
the moving frame of free expansion. Their evolution is governed by the conservation laws, as in the
derivation of fluid dynamics from kinetic theory \cite{Caf80,BGL91,SR09}. Because the local
temperature is $O(t^{-2})$ and the spatial scale is $t$, the resulting compressible Euler system,
that of Theorem~\ref{thm:Meuler}, carries a factor $t^{-2}$. This system is the kinetic counterpart
of the expanding, affine solutions of the compressible Euler equations \cite{Ser97,Gra98,Sid17},
whose velocity field is close to $x/t$. The link between the Landau equation and compressible Euler
dynamics also appears in the study of singularities: for very hard potentials the hydrodynamic
fields of the singular solutions of \cite{BCGJVY26} follow an imploding Euler solution, while for
soft potentials such implosions are essentially excluded \cite{GH25,GHS26}. Here, in contrast, the
Euler dynamics is expanding.

On each leaf the evolution is a homogeneous Landau flow run with the clock
$\lambda(\xi)\tau(t)$, where $\lambda(\xi)$ is the collision rate of the leaf. We take its
relaxation to equilibrium, and the regularity along the way, from the homogeneous theory
\cite{D15,CDH,GS} and the near-equilibrium theory \cite{CM,Mou06}.

\subsection{Main result}

We give an informal statement here, and the precise ones follow in the next section. The data are
of the form $f_0(x,v)=\eps\rho_0(v)F_0(x,v)$, where $\rho_0>0$ is real-analytic in a strip and
decays like $\jap v^{-p}$, and $F_0(\cdot,v)$ is a perturbed Gaussian. The perturbation may depend
on $v$, analytically and boundedly, as made precise in Definition~\ref{def:data}. Under free
transport every leaf acquires an asymptotic mass $\rho_\infty(\xi)=\eps\rho_0(\xi)$ and asymptotic
temperature $\theta_\infty(\xi)$ equal to the
variance of $F_0(\cdot,\xi)$, so that the leaf $\xi$ collides at the rate
$\lambda(\xi)=\rho_\infty(\xi)\theta_\infty(\xi)^{\gamma/2}$, and its effective time is
\begin{equation}\label{eq:sxi}
s_\xi(t)=\lambda(\xi)\big(\tau(t)-\tau(T_0)\big).
\end{equation}
We call the data \emph{non-thermal} if the profiles $F_0(\cdot,\xi)$ stay at a positive distance
$d_0>0$ from Maxwellians, in the sense of~\eqref{eq:d0}. For product data $f_0=\eps\rho_0(v)F_0(x)$
this means that $F_0$ is not a Maxwellian.

\begin{theo}[Informal version of Theorems~\ref{thm:main}, \ref{thm:scatter}, \ref{thm:rate2}
and~\ref{thm:rate3}] Let $\gamma\in[-3,-2]$ and let $f_0=\eps\rho_0F_0$ be as above, with $\eps$
small. Then:
\begin{enumerate}[label=(\roman*)]
\item there is a global classical nonnegative solution, which is real-analytic in the asymptotic
velocity, uniformly in time;
\item every leaf relaxes: once $s_\xi(t)$ exceeds a fixed threshold $S_*$, the leaf is close to a
Maxwellian in Gaussian norm;
\item the parameters of these Maxwellians follow a compressible Euler system with the prefactor
$t^{-2}$, and converge at the rate $t^{-1}$, uniformly in $\xi$;
\item the relaxation is not uniform: the time $t_\xi$ at which the leaf $\xi$ relaxes tends to
infinity as $|\xi|\to\infty$;
\item every leaf relaxes completely, and the solution scatters in $L^1$ to the free transport of a
profile whose fibers are Maxwellians. If the constant $d_1$ of Theorem~\ref{thm:scatter}(iv) is
positive, in particular for non-thermal product data, the relative distance of this profile from the
data is bounded below independently of $\eps$;
\item the convergence has a rate: $\|g(t)-g_{\rm sc}\|_{L^1}\le
C\|f_0\|_{L^1}\big(1+\eps(\tau(t)-\tau(T_0))\big)^{-(p-3)/p}$, with the clock $\tau$
of~\eqref{eq:tau}, that is, $(1+\eps\log(t/T_0))^{-(p-3)/p}$ at $\gamma=-2$ and a rate of order
$(\eps t^{|\gamma|-2})^{-(p-3)/p}$ for $\gamma<-2$. If $d_1>0$, in particular for non-thermal
product data, this rate is sharp;
\item after it has entered, the leaf $\xi$ relaxes like $t^{-c\lambda(\xi)}$ at $\gamma=-2$, and
like $\exp(-c\,s_\xi(t)^{a_\gamma})$, $a_\gamma=2/(14+|\gamma|)$, up to an error term made precise
in Theorem~\ref{thm:rate3}, for $\gamma<-2$.
\end{enumerate}
\end{theo}

\begin{figure}[ht]
\centering
\begin{tikzpicture}
\begin{axis}[width=0.6\textwidth,height=0.42\textwidth,xmin=0,xmax=5,ymin=0,ymax=60,
 xlabel={$|\xi|$},ylabel={$t$},tick label style={font=\small},axis on top]
\addplot[name path=entry,domain=0:5,samples=100,thick,red!70!black] {1+ (1+x^2)^(1.0)*1.5};
\addplot[name path=band,domain=0:5,samples=100,thick,dashed,red!70!black] {1+ (1+x^2)^(1.0)*3};
\path[name path=top] (axis cs:0,60)--(axis cs:5,60);
\path[name path=bot] (axis cs:0,0)--(axis cs:5,0);
\addplot[blue!12] fill between[of=entry and bot];
\addplot[orange!20] fill between[of=band and entry];
\addplot[green!12] fill between[of=top and band];
\node[font=\small] at (axis cs:3.9,8) {crude regime};
\node[font=\small] at (axis cs:1.2,40) {relaxed regime};
\node[font=\small,rotate=48] at (axis cs:3.55,40) {band};
\node[font=\small,red!70!black] at (axis cs:4.2,52) {$t=t_\xi$};
\end{axis}
\end{tikzpicture}
\caption{Schematic picture of the three regimes in the $(|\xi|,t)$ plane, for $\gamma<-2$. The
entry time $t_\xi$, defined by $s_\xi(t_\xi)=S_*$, increases with $|\xi|$ because the collision rate
$\lambda(\xi)$ decreases. At every time some leaves are still in the crude regime. At $\gamma=-2$
the entry times for the model clock are $t_\xi=T_0e^{S_*/\lambda(\xi)}$.}\label{fig:entry}
\end{figure}
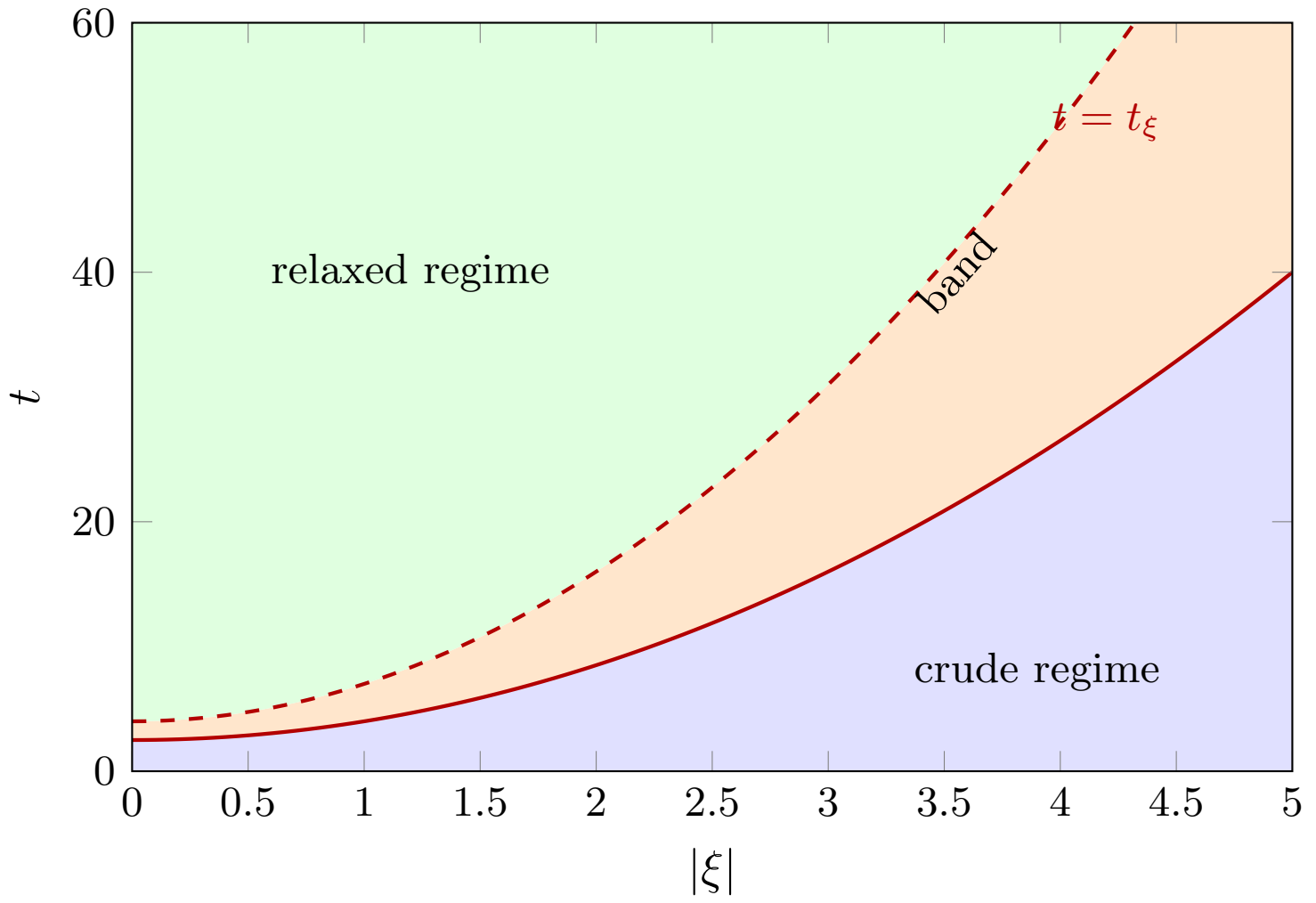

The non-uniformity in (iv) is forced by finite mass: since $\rho_0(\xi)\to0$ as $|\xi|\to\infty$,
the leaves at large $\xi$ collide more and more slowly, and their effective time $s_\xi(t)$ is small
at every time. By Lemma~\ref{lem:Wtail} and Theorem~\ref{thm:main}(iv), for non-thermal data these
leaves are still far from equilibrium. Figure~\ref{fig:entry} shows the resulting picture in the
$(|\xi|,t)$ plane. Each leaf goes through three regimes. In the \emph{crude regime} $s_\xi<S_*$ it
follows the homogeneous Landau flow far from equilibrium. In a \emph{transition band} around
$s_\xi=S_*$ it enters a neighborhood of a Maxwellian. Finally, in the \emph{relaxed regime}
$s_\xi>S_*$ it stays near the Maxwellian with the moving parameters.

\subsection{Strategy of the proof}\label{ssec:strategy}

\begin{figure}[ht]
\centering
\begin{tikzpicture}[node distance=0.55cm and 0.5cm,
 box/.style={draw,rounded corners,align=center,font=\small,minimum height=0.9cm,text
width=3.5cm,fill=#1},box/.default=white]
\node[box=gray!10] (data) {Data $\eps\rho_0F_0$\\ analytic in $v$, Gaussian in $x$};
\node[box=blue!8,below=of data] (layer) {Initial layer $[0,t_1]$\\ adapted boxes
(\S\ref{sec:layer})};
\node[box=blue!8,below=of layer] (disp) {Dispersed layer $[t_1,T_0]$\\ fixed-frame cell energies
(\S\ref{ssec:disp})};
\node[box=orange!12,below left=0.8cm and -0.6cm of disp] (crude) {Crude regime $s_\xi<S_*$\\
fibered comparison (\S\ref{sec:fib}),\\ homogeneous flow (\S\ref{sec:homog})};
\node[box=green!10,below right=0.8cm and -0.6cm of disp] (relax) {Relaxed regime $s_\xi>S_*$\\
Gaussian analytic leaf energies (\S\S\ref{sec:weights},\ref{sec:cells})};
\node[box=yellow!15,below=2.4cm of disp] (band) {Band and bootstrap\\ one radius, cells
(\S\ref{sec:boot})};
\node[box=green!10,right=0.6cm of relax] (euler) {Leaf Euler system\\ Cauchy--Kovalevskaya
(\S\ref{ssec:macro})};
\node[box=gray!10,below=of band] (main) {Theorem~\ref{thm:main}};
\draw[-{Latex}] (data)--(layer);
\draw[-{Latex}] (layer)--(disp);
\draw[-{Latex}] (disp)--(crude);
\draw[-{Latex}] (disp)--(relax);
\draw[-{Latex}] (crude)--(band);
\draw[-{Latex}] (relax)--(band);
\draw[-{Latex}] (euler)--(relax);
\draw[-{Latex}] (band)--(main);
\end{tikzpicture}
\caption{The structure of the proof.}\label{fig:flow}
\end{figure}
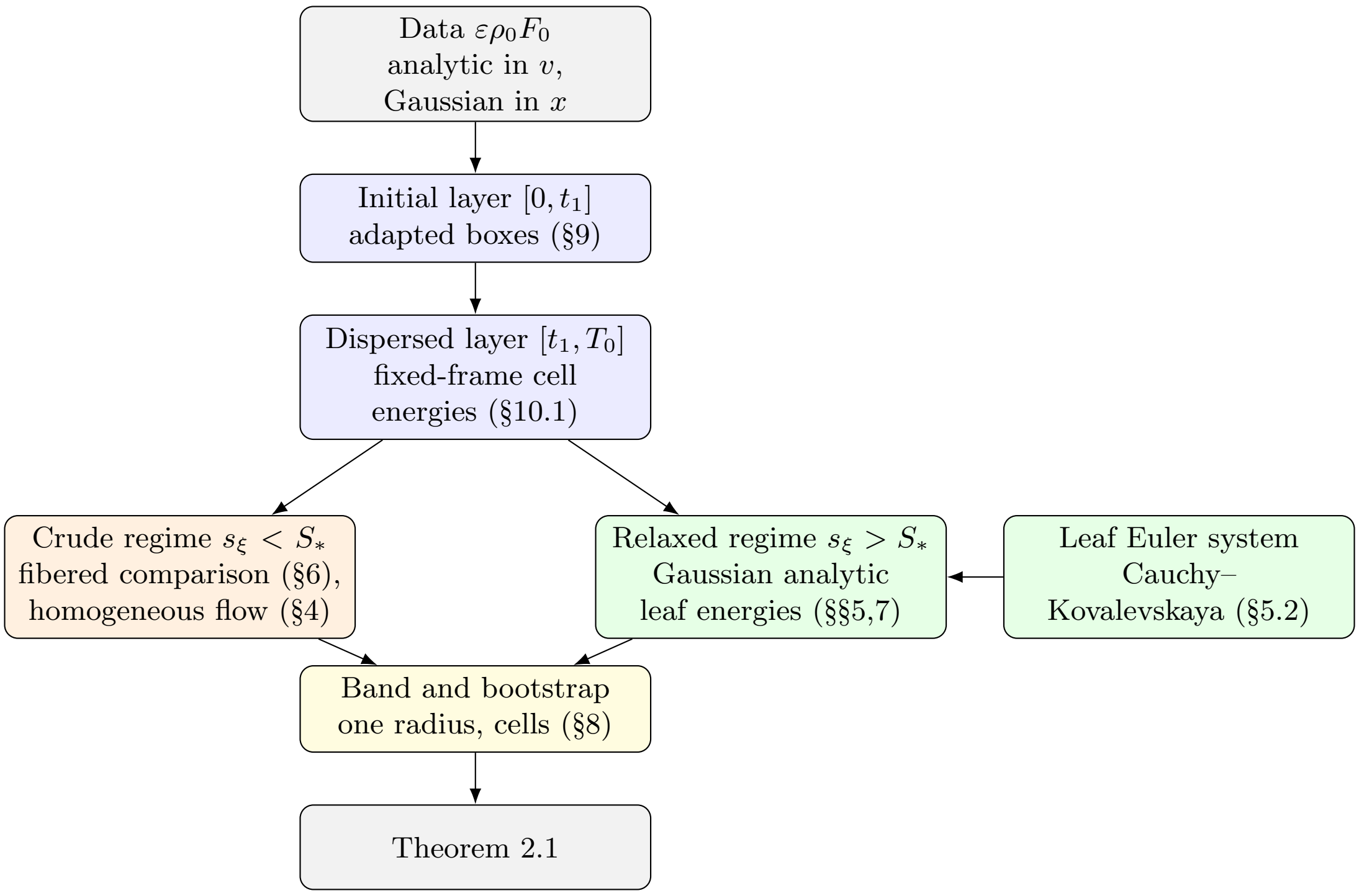

The logical structure of the proof is summarized in Figure~\ref{fig:flow}.

\subsubsection{Three time intervals}\label{sssec:three}

We split the time axis at two times $t_1<T_0$ that depend only on the data class:
\begin{equation}\label{eq:three}
[0,\infty)=\underbrace{[0,t_1]}_{\text{initial layer}}\cup\underbrace{[t_1,T_0]}_{\text{dispersed
layer}}\cup[T_0,\infty).
\end{equation}
On the initial layer the solution has not yet dispersed: the leaves are not yet thin, and we
propagate the solution directly from the data. On the dispersed layer the leaves have formed, but
the effective collision time has advanced only by
\begin{equation}\label{eq:disp-eff}
\eps\big(\tau(T_0)-\tau(t_1)\big)=O\big(\eps T_0^{\kappa_\gamma}\big),\qquad\kappa_\gamma=|\gamma|-2,
\end{equation}
where $T_0^{\kappa_\gamma}$ is read as $1+\log T_0$ at $\gamma=-2$. This is small for
$\eps\le\eps_0(T_0)$, and the solution is still close to free transport. After $T_0$ the collisions
are no longer small, and each leaf goes through the three regimes of Figure~\ref{fig:entry}. Most of
the work concerns $t\ge T_0$, and we describe this part first.

\subsubsection{Reduction to homogeneous dynamics}\label{sssec:st1}

By the leaf identity \eqref{eq:leafid0}, each leaf evolves by a homogeneous Landau equation in $z$,
run with the clock $\kappa(t)$, up to the leafwise transport $t^{-2}z\cdot\nabla_\xi$. The leaf
through $(y,v)$ differs from the fiber $\{v=\mathrm{const}\}$ by velocities of size $|z|/t$. By
Proposition~\ref{prop:Rest}, the fibered equation therefore differs from the leaf identity by a
remainder of relative size $O(t^{-1})$. This remainder stays small over a bounded effective time
$s_\xi$ when $T_0$ is large.

On the crude regime $s_\xi<S_*$ we therefore compare, in Section~\ref{sec:fib}, the leaf with the
\emph{fibered} flow, a family of homogeneous Landau flows indexed by $v$. Section~\ref{sec:homog}
collects what we need from the homogeneous theory of \cite{GS,CDH,D15}: the flow stays smooth, it
obeys Gaussian bounds on bounded time intervals, and it enters a neighborhood of a Maxwellian after
an effective time $S_*$ that depends only on the data class. For $\gamma\in(-3,-2]$ this needs a
rate of convergence to equilibrium that is uniform on the class of data. We prove this rate in
Appendix~\ref{app:hom}.

\subsubsection{Derivative loss of nonlocal leaf energies}\label{sssec:st2}

After entry a leaf is close to a Maxwellian. To keep it there we use an energy estimate for the
linearized operator $\mathcal L$, measured on each leaf in the weighted
norms of the near-equilibrium theory \cite{CM}. In the Gaussian weights that we use, which are
weaker than $\mu^{-1/2}$, only a part of $\mathcal L$ is dissipative. The remaining bounded part is
absorbed through the time decay of $e^{s\mathcal L}$, by adding to the energy the term
\begin{equation}\label{eq:nonlocal-term}
\int_0^\infty\|e^{s\mathcal L}\varphi\|^2_{L^2}\dd s
\end{equation}
as in Proposition~\ref{prop:triG}(i). For $\gamma<-2$ the decay of $e^{s\mathcal L}$ is not
exponential, since then $\mathcal L$ has no spectral gap. This term is used
for every $\gamma\in[-3,-2]$, and it makes the norm \emph{nonlocal} in $z$.

Consider now the leafwise transport. Let $H(t,\xi,z)$ be a leaf quantity, $D_t$ the derivative along
the ray $x=t\xi$ at fixed $(\xi,z)$, and $E_n$ the energy of $\partial^n_\xi H$. Differentiating
$D_tH=t^{-2}z\cdot\nabla_\xi H+\dots$ $n$ times, we get
\begin{equation}\label{eq:st-loss-eq}
D_t\partial^n_\xi H=t^{-2}z\cdot\nabla_\xi\partial^n_\xi H+\dots ,
\end{equation}
so the equation for the $n$-th derivative contains the $(n+1)$-st. For an energy that is
\emph{local} in $z$, such as $\iint w^2H^2\dd z\dd\xi$, this term is a transport in $\xi$ at fixed
$z$, and by Lemma~\ref{lem:localnoloss} it integrates by parts without loss. For a nonlocal norm it
does not: by Proposition~\ref{prop:commloss}, the commutator of the norm with the multiplication by
$z$ leaves one derivative $\nabla_\xi$, which on leaves is one derivative $\partial_v$.
Hence the natural estimate has the form
\begin{equation}\label{eq:st-loss}
\frac d{dt}E_n\le a(t)E_n+C_0t^{-2}\sqrt{E_nE_{n+1}} ,
\end{equation}
with $a(t)$ integrable, and in finite regularity it does not close.

Local energies whose weights are exact inverse Maxwellians on every leaf would avoid the nonlocal
term, but they are not available in general, for two reasons, explained in
Subsection~\ref{ssec:obstructions}.
First, a static weight that is the inverse of a Maxwellian in
$y$ on every fiber has Gaussian leaf restrictions only in the class of Gaussian solutions of free
transport, which contains the traveling Maxwellians of \cite{CL}. Second, a weight that follows the
leaf parameters loses, wherever the leaf temperature varies with $\xi$, a power $|z|^3$ of weight at
the rate $t^{-2}$. Our data are never Gaussian solutions of free transport, since $\rho_0$ decays
only polynomially, and for them we do not have such an energy.

\subsubsection{Analytic energies}\label{sssec:st3}

A function of $\xi$ is real-analytic with radius $\rad$ if its derivatives grow at most like
$n!\,\rad^{-n}$. We measure this by the analytic energy
\begin{equation}\label{eq:st-analytic}
\mathbb A(t)=\sum_{n\ge0}\Big(\frac{\rad(t)^n}{n!}\Big)^2E_n(t),
\end{equation}
with a radius $\rad(t)$ that decreases in time. Differentiating \eqref{eq:st-analytic} produces,
besides the terms of \eqref{eq:st-loss}, the term
\begin{equation}\label{eq:st-gain}
2\sum_nn\,\frac{\rad'}{\rad}\Big(\frac{\rad^n}{n!}\Big)^2E_n ,
\end{equation}
which is negative and grows with $n$. In the same weights, the loss term
$C_0t^{-2}\sqrt{E_nE_{n+1}}$ is of size $C_0t^{-2}\frac{n+1}\rad$ times a geometric mean of the
$n$-th and $(n+1)$-st terms of $\mathbb A$. Hence \eqref{eq:st-gain} absorbs the loss if
\begin{equation}\label{eq:st-radius}
|\rad'(t)|\ge\tfrac34C_0t^{-2},
\end{equation}
and Lemma~\ref{lem:analytic} shows that, since $t^{-2}$ is integrable, the radius decreases only by
$\frac34C_0/T_0$ over infinite time, as in the abstract Cauchy--Kovalevskaya theorem
\cite{Ovs65,Ni72,Nis77,Caf90} and in the analytic and Gevrey proofs of Landau damping
\cite{MV11,BMM16}, where a loss of one derivative at an integrable rate is also compensated by a
shrinking radius.

Analyticity is not produced by the collisions: it is assumed for the data, in
$v$, and propagated. Moreover, since $\nabla_\xi$ on leaves is $\partial_v$, analyticity in the
asymptotic velocity $\xi$ is analyticity of the profile $g(t,y,\cdot)$ in $v$. This is the form in
which it appears in Theorem~\ref{thm:main}(i).

\subsubsection{Gaussian weights}\label{sssec:st4}

The leaves collide at different rates $\lambda(\xi)$. Therefore a derivative $\partial_\xi$ of the
collision term produces, besides the linearized operator applied to the derivative, a term
$(\partial_\xi\lambda)Q(\psi,\psi)$, where $\psi$ is the leaf written in the frame of its
parameters. Its linear part couples the mass of the perturbation, which is measured only in the
energy norm, to the velocity tail of its derivative. We call this the \emph{mass-tail coupling}.
In polynomial weights the coupling is not bounded by the dissipation norm, whatever the power. In
Gaussian weights $e^{q|z|^2}$, however, the drift of the linearized operator against the weight
gains a full factor $\jap z^2$ in the zero-order part of the dissipation norm, and the coupling is
bounded by the dissipation. For the weights $e^{q\jap z^s}$ with $s<2$ it is not. We prove these
facts in Subsection~\ref{sssec:obs2}. For this reason the data are taken Gaussian
in space: by the example of Subsection~\ref{ssec:leaves}, the spatial profile becomes the velocity
distribution on the leaves.

\subsubsection{The macroscopic part and the leaf Euler system}\label{sssec:st5}

The linearized operator controls only the part of the perturbation that is orthogonal to its kernel,
which is spanned by the Maxwellian, its momentum and its energy. The macroscopic part, that is, the
leaf density, mean velocity and temperature, is not damped. We therefore linearize around the
\emph{local leaf Maxwellian}, the Maxwellian with the same leaf moments as the solution. By
Lemma~\ref{lem:Mfree}, the microscopic part then has vanishing moments, and the macroscopic part is
carried by the leaf parameters.

By Theorem~\ref{thm:momtransport}, the leaf moments obey an exact transport law, which for
Maxwellian leaves becomes, by Theorem~\ref{thm:Meuler}, a compressible Euler system with the
prefactor $t^{-2}$. In the time
\begin{equation}\label{eq:st-sigma}
\sigma=T_0^{-1}-t^{-1}\in[0,T_0^{-1})
\end{equation}
the prefactor disappears, and the whole evolution lasts the finite time $T_0^{-1}$. In
Lemma~\ref{lem:Mref} we solve the system by the abstract Cauchy--Kovalevskaya theorem in strips of
analyticity. The leaf parameters therefore converge, at the rate $t^{-1}$, to a limit that is not a
global equilibrium: it depends on $\xi$ and is determined by the data.

\subsubsection{Uniformly local energies}\label{sssec:st6}

The arguments of Subsections~\ref{sssec:st2}--\ref{sssec:st5} cannot be applied to each leaf
separately, because the leafwise transport moves information from leaf to leaf: the leaf $\xi$ at
time $t$ depends on the leaves $\xi+z\sigma(t)$, at distances up to $|z|/T_0$, which is unbounded on
the support of a Gaussian in $z$. Proposition~\ref{prop:ACtr} shows that, already for the
collisionless model, the analytic energy of a single leaf is not controlled by its initial value,
whatever the radius.

We therefore use energies that are \emph{uniformly local} in $\xi$. We integrate the leaf energies
over cells
\begin{equation}\label{eq:st-cells}
\xi\in\xi_c+O(\rad_0),
\end{equation}
of the size of the analyticity radius, with one radius for all leaves. On a cell the transport
integrates by parts onto the cutoff and produces only a factor $|z|/\rad_0$. By
Lemma~\ref{lem:CUtr}, the slowly shrinking Gaussian weight absorbs this factor, up to a Gronwall
term of size $t^{-2}$. Pointwise bounds in $\xi$ are recovered by the Sobolev inequality of
Lemma~\ref{lem:CUsob} on cubes of side $\rad_0$, which loses two derivatives but no power of
$\rad_0$, and these two derivatives are again absorbed by the decrease of the radius.

\subsubsection{The entry front and the bootstrap}\label{sssec:st7}

The leaves enter the relaxed regime at different times. There is no time after which all leaves have
$s_\xi\ge S_*$, and, for non-thermal data, none after which all leaves are near equilibrium (see
Subsection~\ref{sssec:obs3}). Hence the crude comparison of Subsection~\ref{sssec:st1} and the
near-equilibrium energies of Subsections~\ref{sssec:st2}--\ref{sssec:st6} must be used at the same
time, on different leaves.

We organize the bootstrap by the leaf clock $s_\xi$. The nonlocal part \eqref{eq:nonlocal-term} of
the energy is switched on in a band
\begin{equation}\label{eq:st-band}
S_b\le s_\xi\le S_b+1,\qquad S_*<S_b\le S_*+1,
\end{equation}
just after entry. By Proposition~\ref{prop:Bband}, in this band a damping weight in the effective
time absorbs the deficit of the partially switched energy and the terms that come from the moving
front. The bootstrap assumptions, stated in Theorem~\ref{thm:CUboot}, are:
\begin{itemize}
\item[(U1)] a bound on the deviation of the leaf parameters from the reference Euler flow;
\item[(U2)] bounds on the cell energies;
\item[(U3)] a smaller bound on the energies of the cells that have reached the band;
\item[(U4)] Gaussian bounds on the leaves that are crude or have recently entered, which follow
from (U1)--(U3).
\end{itemize}
We choose the constants in the order of Appendix~\ref{app:const}. Since the entry level
$\delta_1$ is fixed before the entry time $S_*$ and no information on neighboring leaves is used at
entry, this order is consistent.

\subsubsection{The initial and the dispersed layers}\label{sssec:st8}

On $[0,t_1]$ the solution is propagated directly from the data, through the relative density
$G=g/f_0$. Since the density is small, its equation has small coefficients, and it is nearly
degenerate. We use energies on boxes adapted to the dispersion, of widths
\begin{equation}\label{eq:st-boxes}
\text{radial}\ \simeq\big(\jap y+(\log\jap v)^{1/2}\big)^{-1},\qquad\text{transverse}\
\simeq\big(1+(\log\jap v)^{1/2}\big)^{-1},
\end{equation}
together with a null structure of the Landau kernel and a velocity weight. We do this in
Section~\ref{sec:layer}. The endpoint $\gamma=-2$ needs in addition a transit-time argument, which
we give in Subsection~\ref{ssec:endpoint}.

On $[t_1,T_0]$ the collisions act for the small effective time \eqref{eq:disp-eff}. Energies in a
fixed frame suffice there, in cells of a fixed size $\rad_1$, the radius of analyticity on this
interval, without damping or loss of radius. In Subsection~\ref{ssec:disp} they provide, at $T_0$,
the data required by the bootstrap.

\subsubsection{Time regularity of the energies}\label{sssec:st9}

The energy inequalities above are differential inequalities for energies of the solution itself,
which is not smooth enough for them to be differentiated directly. At the top $y$-order the
collision term is a distribution of order two along the leaves, and the sums over the levels $n$ can
be differentiated term by term only if $\sum_nn\,w_n^2E_n$ is integrable in time. We use all
energy inequalities in integrated form. The necessary time regularity comes from a viscous
approximation, which keeps a fixed fraction of the collision dissipation, and we prove it in
Subsection~\ref{sssec:CUac}. We treat the homogeneous trajectories near $\sigma=0$ in the same
way in Subsection~\ref{sssec:RSshort}, and identify them with the solution of \cite{GS} by a
uniqueness argument.

\subsection{Further comments and related work}

\emph{Vacuum stability.} For the Boltzmann equation near vacuum, global existence for small data
goes back to Illner and Shinbrot \cite{IS84}, using the Kaniel--Shinbrot iteration \cite{KS78} (see
also \cite{BT85,Ham85,Tos86,Glas96,Ark11,Alo09}). For the Landau equation, \cite{Luk,Ch20} treat
$\gamma>-2$ and hard potentials, and \cite{Ch21} the Boltzmann equation with moderately soft
potentials. The weakly collisional Vlasov--Poisson--Landau system is treated in \cite{CLN23}.
Chaturvedi and Luk \cite{CL} prove linear stability of traveling Maxwellians for very soft
potentials. These are Gaussian solutions of free transport, with a leaf temperature that does not
depend on $\xi$, as the proof of Proposition~\ref{prop:travelling} shows. Their restriction to every
leaf is a Maxwellian, so near them, by Lemma~\ref{lem:localnoloss}, the static weight $
\mathcal M^{-1}$ yields a local top-order energy in finite regularity, with neither derivative loss
nor weight loss.

\emph{Comparison with \cite{CL}.} The work \cite{CL} is the closest to ours. It concerns the
linearization of the Landau equation with Coulomb potential around a traveling Maxwellian: the
microscopic part of the linear solution decays at an inverse polynomial rate, while the macroscopic
part stays bounded but in general does not decay. The second fact matches the picture of the present
paper, in which the leaf parameters do not relax to a fixed state but converge, at the rate
$t^{-1}$, to the end state of the Euler system with prefactor $t^{-2}$. Compared with \cite{CL}, our
result is nonlinear, concerns data near vacuum rather than near a traveling Maxwellian (by (D3) the
asymptotic leaf density $\eps\rho_0(\xi)$ is bounded above and below by multiples of $\jap\xi^{-p}$,
so $\log\rho_0$ is not a polynomial of degree $\le2$, and by that proposition no
static weight that is a Maxwellian in $y$ on every fiber has Gaussian leaf restrictions for these
data), and covers the whole range $\gamma\in[-3,-2]$, including the endpoint $\gamma=-2$. For this
we need analyticity in $v$: away from the traveling Maxwellians we know no loss-free top-order
energy with Gaussian structure, for the two reasons explained above, and we absorb
the derivative loss of the nonlocal energies by a radius of analyticity that shrinks at an
integrable rate. (For the homogeneous Landau equation, analytic and Gevrey regularity in $v$ is
classical \cite{CLX10,MPX13}.) Whether finite regularity in $v$ suffices near vacuum for
$\gamma\le-2$ is left open.

\emph{The regime $\gamma\le-2$.} Below the threshold $\gamma=-2$, near-vacuum solutions thermalize
leaf by leaf. Their leaf parameters converge, at the rate $t^{-1}$, to a limit determined by the
data through the Euler system with prefactor $t^{-2}$, which has finite total duration. The
solutions scatter, but the scattering state is thermalized and in general not close to the data,
by Theorem~\ref{thm:scatter}(iv).

The $L^1$ distance to the scattering state is bounded by
$C\|f_0\|_{L^1}\big(1+\eps(\tau(t)-\tau(T_0))\big)^{-(p-3)/p}$, where $p$ is the decay exponent of
$\rho_0$. This bound is logarithmic in $t$ at $\gamma=-2$ and of order $(\eps t)^{-(p-3)/p}$ in the
Coulomb case. The bound counts the leaves that have not yet relaxed, and it is sharp if $d_1>0$, by
Theorems~\ref{thm:rate2} and~\ref{thm:rate3}. We give a heuristic discussion in
Remark~\ref{rem:rate}.
Convergence in stronger norms is left open.

\emph{Inputs from the literature.} We list the results that we quote as Inputs~\ref{hyp:E} in
Subsection~\ref{ssec:results}. We need a few of them slightly beyond the
setting in which they are stated: the near-Maxwellian estimates at $\gamma=-2$, the
entropy--entropy-dissipation inequality of \cite{CDH} for $-4<\gamma\le-2$, the Gaussian bounds
(GS4)--(GS5), and the behavior near $\sigma=0$ of the solution of \cite{GS}, which we derive from
\cite{HST20}. We prove these extensions in Appendix~\ref{app:inputs}.

\subsection{Structure of the paper}

Section~\ref{sec:setting} fixes the notation and the data class, lists the results from
the literature and states the main theorem. Section~\ref{sec:leafgeom} develops the geometry of
leaves (the leaf identity, the remainder, and the exact transport law of leaf moments) and explains
the two obstructions that lead to analyticity in $v$ and to Gaussian weights.
Section~\ref{sec:homog} collects the properties of homogeneous trajectories.
Section~\ref{sec:weights} sets up the Gaussian leaf energies and derives and solves the leaf Euler
system. Section~\ref{sec:fib} treats the crude regime by a fibered comparison.
Section~\ref{sec:cells} introduces the uniformly local analytic energies, and Section~\ref{sec:boot}
runs the bootstrap. Section~\ref{sec:layer} treats the initial layer and, in
Subsection~\ref{ssec:endpoint}, the endpoint $\gamma=-2$ in all parts of the proof.
Section~\ref{sec:proofmain} treats the dispersed layer, assembles the proof, and proves the
scattering theorem and the rates.

We collect the proofs of the main sections that are long or technical in a final subsection of
the section that uses them: Subsection~\ref{ssec:tri} for the trilinear estimates in Gaussian
weights, Subsection~\ref{ssec:fib} for the tools for the fibered comparison, including the
short-time construction of homogeneous trajectories, Subsections~\ref{ssec:cellsproofs}
and~\ref{sssec:CUac} for the uniformly local energies, including the local existence theory and the
time regularity of the cell energies, and, for the initial layer,
Subsections~\ref{ssec:tools}--\ref{ssec:drift} for the tools on adapted boxes, the energy estimates
and the drift term at the top $y$-levels. On a first reading these subsections can be
skipped. The appendices contain the rest. Appendix~\ref{app:inputs} contains the complements to the
literature: two complements to the near-Maxwellian theory and its extension to $\gamma=-2$, and the
proofs of the results on the homogeneous equation (the uniform $L^1$ rate, the entry theorem, and
the Gaussian bounds). Appendix~\ref{app:auxest} collects auxiliary results, and the symbolic check
of the leaf Euler system is given in an ancillary file. Appendix~\ref{app:regimes} contains the
frame terms, the band balance and the
crude leaf inequality, the viscous approximation, and the order in which the constants are chosen.

\subsection*{Acknowledgements}
The author is supported by Universit\'e de Toulouse III Paul Sabatier and PRIN 2022 ``Turbulent
effects vs Stability in Equations from Oceanography'' (TESEO), project number 2022HSSYPN.

\subsection*{Use of AI tools}
In preparing this work the author used the AI assistant Claude (Anthropic) to check computations 
(including a symbolic verification of the leaf Euler system,
provided as an ancillary file), to search for errors in the proofs, and to edit the text. All
results and proofs were checked by the author, who takes full responsibility for the content of the
paper.

\section{Setting, data class and main results}\label{sec:setting}

\emph{Notation.} Throughout, $c$ and $C$ denote positive constants that may change from line to
line, and we state their dependence where it matters. Constants introduced in a proof are local to
that proof unless they are referred to elsewhere by name (for instance $C_{\rm far}$ or
$c_{\rm ret}$), and a few letters (for instance $B$, $p$, $\sigma$ and $\vartheta$) are reused with
a local meaning in independent parts of the paper, where that meaning is fixed when the letter is
introduced. Collision rates are denoted by $\lambda$ ($\lambda_p$, $\lambda_\xi$, $\lambda_{\rm
ref}$, $\lambda_v$), and analyticity radii by $\rad$ ($\rad(t)$, $\rad_0$, $\rad_1$, $\rad_k$,
$\rad_{\rm post}$). In a few self-contained proofs $\lambda$ also denotes a local exponent or
constant. We choose the constants in a fixed order, which is listed in Table~\ref{tab:const} of
Appendix~\ref{app:const}.

\subsection{The collision operator}\label{ssec:collop}
Fix $\gamma\in[-3,-2]$. With $a$ as in \eqref{eq:Q}, we put
\begin{equation}
b_i(u)=\partial_ja_{ij}(u)=-2|u|^\gamma u_i,\qquad
c(u)=\partial_ib_i(u)=\begin{cases}-2(\gamma+3)|u|^{\gamma},&\gamma>-3,\\-8\pi\delta_0,&\gamma=-3.\end{cases}
\end{equation}
For a function $F$ on $\R^3$ we write $\ab[F]=a*F$, $\bb[F]=b*F$ and $\cb[F]=c*F$. The bilinear
operator and its flux are
\begin{equation}
Q(F,G)=\nabla\cdot B(F,G),\qquad B(F,G)=\ab[F]\nabla G-\bb[F]G,\label{eq:flux}
\end{equation}
so that also $Q(F,G)=\ab[F]:\nabla^2G-\cb[F]G$. Since $a$ is even, we have
\begin{equation}\label{eq:Bflux}
B_i(F,G)(y)=\int a_{ij}(y-z)\big[F(z)\partial_jG(y)-G(y)\partial_jF(z)\big]\dd z .
\end{equation}
The collision invariants are $1$, $z_i$ and $|z|^2$: $\int Q(F,F)\psi=0$ for $\psi$ in their span.
The Maxwellians
\begin{equation}
M_p(z)=\rho(2\pi\theta)^{-3/2}e^{-|z-y_0|^2/(2\theta)},\qquad p=(\rho,y_0,\theta),
\end{equation}
satisfy $Q(M_p,M_p)=0$. We write $\mu=M_{(1,0,1)}$.

\emph{Frames.} For $p=(\rho,y_0,\theta)$ we define
$T_pF(z)=\rho\,\theta^{-3/2}F(\theta^{-1/2}(z-y_0))$, so that $T_p\mu=M_p$. By the homogeneity of
$a$, we have
\begin{equation}\label{eq:frame}
Q(T_pF,T_pF)=\lambda_p\,T_p\,Q(F,F),\qquad\lambda_p=\rho\,\theta^{\gamma/2}.
\end{equation}
Therefore a homogeneous Landau trajectory normalized by its moments runs with the clock multiplied
by $\lambda_p$, the \emph{collision rate} of the parameters $p$.

\subsection{Profile, leaves, clocks}\label{ssec:leaf}

The profile is $g(t,y,v)=f(t,y+tv,v)$, and for a function $h(t,y,v)$ the leaf through $x$ is
$h^{(x)}(z)=h(t,z,(x-z)/t)$ (see \eqref{eq:profile}, \eqref{eq:leafdef}). We write $\xi=x/t$ and use
$(t,\xi,z)$ as coordinates on phase space, with $y=z$ and $v=\xi-z/t$. We denote by $D_t$ the
derivative along rays at fixed $(\xi,z)$. The operators $\partial_y$, $\partial_v$ act on profiles,
and $D=\partial_v-t\partial_y$ is the image of $\partial_v$ under this change of variables.

The \emph{leaf moments} and \emph{leaf parameters} are
\begin{equation}
m_h(t,x)=\int h^{(x)}(z)\,\mom(z)\dd z,\qquad\mom(z)=(1,z,|z|^2),\qquad
m_g(t,t\xi)=\rho_g\big(1,y_{0,g},|y_{0,g}|^2+3\theta_g\big)(t,\xi),
\end{equation}
and $p_g(t,\xi)=(\rho_g,y_{0,g},\theta_g)$. We also use the log-variables
$\varpi=(\ell,y_0,\vartheta)$ with $\ell=\log\rho$, $\vartheta=\log\theta$, and
$\varpi'=(\nabla_\xi\ell,y_0,\vartheta)$. For complex values we write
$|\varpi'|:=|\nabla_\xi\ell|+|y_0|+|\vartheta|$, with the Euclidean norms on $\C^3$, $\C^3$ and
$\C$.

The clocks are $\kappa(t)=t^{-(3+\gamma)}$, $\tau(t)=\int_1^t\kappa$ as in
\eqref{eq:kappa}--\eqref{eq:tau}, and the effective time $s_\xi=\lambda_{\rm
ref}(\xi)(\tau(t)-\tau(T_0))$ of the leaf $\xi$, defined in Subsection~\ref{ssec:cells} and modeled
by~\eqref{eq:sxi}. We write $\kappa_\gamma:=|\gamma|-2$ for $\gamma\in(-3,-2)$ and $\kappa_{-3}:=1$,
so that $\int_{t_1}^{T_0}\kappa\le CT_0^{\kappa_\gamma}$. At $\gamma=-2$, we read
$T_0^{\kappa_\gamma}$ as $1+\log T_0$.

\emph{Weighted spaces.} For the explicit Gaussian weights $e^{c|z|^2}$, $L^2(e^{c|z|^2})$ has the
norm $(\int|f|^2e^{c|z|^2})^{1/2}$. For polynomial weights and for the leaf weights $m$ we follow
\cite{CM}: $\|f\|_{L^2(\jap v^k)}=\|\jap v^kf\|_{L^2}$ and $\|f\|_{L^2(m)}=\|mf\|_{L^2}$, and
likewise for $X(m)$ and $Y(m)$.

\subsection{The data class}\label{ssec:data}

For $R>0$ we write $\mathbb S_R=\{v\in\mathbb C^3:|\mathrm{Im}\,v|<R\}$.

\begin{defin}[The data class, and the class of reference profiles]\label{def:data}
Let $\Lambda,\Lambda_\Phi\in(0,1]$, $L\ge1$, $p\ge4$, $C_F,\beta_F,A_F,C_\Phi,C_{\rm pv}>0$ and
$N_y\ge20$; put $N_x=\max\{20,N_y\}$. Consider functions $f_0(x,v)=\eps\rho_0(v)F_0(x,v)$, $\eps>0$,
and the conditions:
\begin{enumerate}[label=(D\arabic*)]
\item $F_0\ge0$, $\int F_0(x,v)\dd x=1$ for every $v\in\R^3$, $F_0\le C_Fe^{-\beta_F|x|^2}$, and
$\sum_{|\alpha|\le N_x}\sup_{x,v}e^{\beta_F|x|^2/2}|\nabla^\alpha_xF_0|\le A_F$;
\item \emph{(D$_{\rm iso}$)} $F_0(x,v)=c_0(v)\exp\big(-\beta_F|x|^2-\Phi_1(x)-\Phi_{\rm
v}(x,v)\big)$, where $c_0(v)>0$ normalizes the mass and:
\begin{itemize}
\item $\Phi_1$ is real-analytic, and $\nabla\Phi_1$ extends holomorphically to the strip
$\{x\in\mathbb C^3:|\mathrm{Im}\,x|<\Lambda_\Phi\}$, with $|\nabla\Phi_1|\le C_\Phi$ there. By the
Cauchy estimates, $\|\nabla^j\Phi_1\|_{L^\infty(\R^3)}\le C(j)\,C_\Phi\Lambda_\Phi^{1-j}$ for
$j\ge1$;
\item $\Phi_{\rm v}$ is real on $\R^3\times\R^3$ and extends holomorphically to $\{x\in
\mathbb C^3:|\mathrm{Im}\,x|<\Lambda_\Phi\}\times\mathbb S_\Lambda$, with $|\Phi_{\rm v}|\le C_{\rm pv}$
there, and $|\mathrm{Im}\,\Phi_{\rm v}(x,v)|\le1$ for $x\in\R^3$, $v\in\mathbb S_\Lambda$.
\end{itemize}
The Gaussian rate $\beta_F$ does not depend on $v$, and the velocity enters $F_0$ only through the
bounded perturbation $\Phi_{\rm v}$. For $\Phi_{\rm v}=0$, $F_0$ does not depend on $v$;
\item $\rho_0>0$ extends holomorphically to $\mathbb S_\Lambda$, and
\begin{equation}
A_F^{-1}\rho_0(v)\le|\rho_0(v+i\eta)|\le A_F\rho_0(v),\qquad\jap
v|\nabla\log\rho_0|+|\nabla^2\log\rho_0|\le L,\qquad A_F^{-1}\jap v^{-p}\le\rho_0\le A_F\jap v^{-p};
\end{equation}
\item $F_0>0$, and for every $x\in\R^3$, $F_0(x,\cdot)$ extends holomorphically to
$\mathbb S_\Lambda$, the extension and its $x$-derivatives of order $\le N_x$ are continuous on
$\R^3\times \mathbb S_\Lambda$, and for $v\in\R^3$, $|\eta|<\Lambda$ and $|\alpha|\le N_x$,
\begin{equation}\begin{gathered}
|F_0(x,v+i\eta)|\le A_FF_0(x,v),\qquad e^{\beta_F|x|^2/2}|\nabla^\alpha_xF_0(x,v+i\eta)|\le A_F,\\
\jap v\,|\nabla_v\log F_0(x,v)|+|\nabla_v^2\log F_0(x,v)|\le L .
\end{gathered}\end{equation}
\end{enumerate}
The class $\Dc^\sharp$ of \emph{reference profiles} consists of the $f_0$ that satisfy (D1), (D3)
and (D4). The \emph{data class} $\Dc\subset\Dc^\sharp$ consists of the $f_0$ that satisfy (D1),
(D$_{\rm iso}$), (D3) and (D4). For $\Phi_{\rm v}=0$ these are products,
$f_0(x,v)=\eps\rho_0(v)F_0(x)$, and for them (D4) follows from (D1) and (D$_{\rm iso}$).
\end{defin}

The main theorem, Theorem~\ref{thm:main}, concerns data in $\Dc$, posed at $t=0$. The larger class
$\Dc^\sharp$ enters in Theorem~\ref{thm:t1}, for data posed at a positive time: there the spatial
profile $F_0(\cdot,v)$ need not be of the form (D$_{\rm iso}$). For instance, its Gaussian rate and
its center may depend on the velocity. From Section~\ref{sec:leafgeom} on, except in the initial
layer of Section~\ref{sec:layer}, we write all statements for $f_0\in\Dc^\sharp$, and the
constants said to depend on $\Dc$ depend on the parameters of the definition.

\emph{Examples.} A typical product element of $\Dc$ is $\rho_0=\jap v^{-p}$ with $\Lambda=\frac12$,
and $\Phi_1$ constant, so that $F_0$ is the centered isotropic Gaussian
$(\beta_F/\pi)^{3/2}e^{-\beta_F|x|^2}$, or $\Phi_1=c\jap x$ with
$c>0$, $\Lambda_\Phi=\frac12$ and $C_\Phi=c$ (on the strip, $|x|^2\le|1+x\cdot x|$, so $|\nabla\jap
x|\le1$).

An element of $\Dc$ that is not a product is
$F_0(x,v)=c_0(v)\exp\big(-\beta_F|x|^2-\frac12\psi(x)b(v)\big)$, with $\psi(x)=x_1(1+x\cdot
x)^{-1}$, $b(v)=1+(1+v\cdot v)^{-1}$ and $\Lambda=\frac14$, where $c_0(v)$ normalizes the mass. On
$\mathbb S_{1/4}$ we have $|1+v\cdot v|\ge\frac{15}{16}$, so $b$ is holomorphic there with
$|b|\le\frac{31}{15}$ and $\jap{\mathrm{Re}\,v}|\nabla b|\le2$, and $|\psi|\le\frac12$ on $\R^3$.
Hence $|\mathrm{Im}(\psi b)|\le\frac8{15}$ on $\R^3\times\mathbb S_{1/4}$, the integral
$c_0(v)^{-1}=\int e^{-\beta_F|x|^2-\frac12\psi(x)b(v)}\dd x$ is holomorphic with real part bounded
below, and (D1) and (D4) hold with constants depending on $\beta_F$. The profile also satisfies
(D$_{\rm iso}$), with $\Phi_1=0$, $\Phi_{\rm v}=\frac12\psi(x)b(v)$, $\Lambda_\Phi=\frac12$ and
$C_{\rm pv}=\frac{31}{30}$. Indeed, for $|\mathrm{Im}\,x|<\frac12$ we have $\mathrm{Re}(1+x\cdot
x)\ge\frac34+|\mathrm{Re}\,x|^2$, so $\psi$ is holomorphic there with
$|\psi|\le(|\mathrm{Re}\,x|+\frac12)/(\frac34+|\mathrm{Re}\,x|^2)\le1$. Moreover,
$|\mathrm{Im}\,\Phi_{\rm v}|\le\frac4{15}$ for real $x$.

For this profile the shape, the mean and the temperature of $F_0(\cdot,v)$ depend on $v$, and
$d_0>0$ and $d_1>0$ in \eqref{eq:d0} and in Theorem~\ref{thm:scatter}(iv). Indeed, on $\R^3$ the
function $b$ takes its values in $[1,2]$, and $c\,e^{-\beta_F|x|^2-s\psi(x)}$ is a Maxwellian for no
$s\ne0$.

With $b(v)=(1+v\cdot v)^{-1}$ instead, the profiles $F_0(\cdot,v)$ tend to the centered Gaussian as
$|v|\to\infty$, so that $d_0=d_1=0$, although no profile is a Maxwellian.

Under free transport the leaf parameters of $f_0\in\Dc^\sharp$ converge to
\begin{equation}
\rho_\infty(\xi)=\eps\rho_0(\xi),\qquad y_{0,\infty}(\xi)=\int xF_0(x,\xi)\dd
x,\qquad\theta_\infty(\xi)=\tfrac13\int|x-y_{0,\infty}(\xi)|^2F_0(x,\xi)\dd x,
\end{equation}
and the collision rate of the leaf $\xi$ is
$\lambda(\xi)=\rho_\infty(\xi)\theta_\infty(\xi)^{\gamma/2}$. For product data, $y_{0,\infty}$ and
$\theta_\infty$ are constants. In general they are real-analytic functions of $\xi$ with
$|y_{0,\infty}|\le C$, $c\le\theta_\infty\le C$ and, by (D4) and $\nabla_vF_0=F_0\nabla_v\log F_0$,
$\jap\xi\big(|\nabla y_{0,\infty}|+|\nabla\theta_\infty|\big)\le CL$ and
$|\nabla^2y_{0,\infty}|+|\nabla^2\theta_\infty|\le CL^2$, with $c,C$ depending only on $C_F$ and
$\beta_F$. For the temperature no constant is lost: since $\int(x-y_{0,\infty})F_0\dd x=0$, we have
$\nabla\theta_\infty=\frac13\int|x-y_{0,\infty}|^2F_0\,\nabla_v\log F_0\dd x$, hence
$\jap\xi|\nabla\log\theta_\infty(\xi)|\le\sup_x\jap\xi|\nabla_v\log F_0(x,\xi)|\le L$. By (D3) and
these bounds, we have $|\nabla_\xi\log\lambda|\le CL$, and $\lambda(\xi)\to0$ as $|\xi|\to\infty$.
The leafwise data norm
\begin{equation}
B:=\sup_{v\in\R^3}\sup_{|\eta|<\Lambda}\sum_{|\alpha|\le N_y}\int
e^{\beta_F|y|^2/2}\frac{|\partial^\alpha_yf_0(y,v+i\eta)|^2}{(\eps\rho_0(v))^2}\dd y\le
C(A_F,\beta_F)
\end{equation}
is a data constant. Here we use (D3) and (D4) for the bound.

For $\xi\in\R^3$ we denote by $M_0^\xi:=M_{(1,y_{0,\infty}(\xi),\theta_\infty(\xi))}$ the
Maxwellian with the mass, mean and temperature of $F_0(\cdot,\xi)$, and we set
\begin{equation}\label{eq:d0}
d_0(\xi):=\|F_0(\cdot,\xi)-M_0^\xi\|_{L^2(\R^3)},\qquad d_0:=\inf_{\xi\in\R^3}d_0(\xi) .
\end{equation}
For product data these do not depend on $\xi$: $M_0^\xi=M_0$ and $d_0(\xi)=d_0=\|F_0-M_0\|_{L^2}$. A
Maxwellian of unit mass is determined by its mean and temperature. Thus $d_0(\xi)=0$ if and only if
$F_0(\cdot,\xi)$ is a Maxwellian. We say that \emph{$F_0$ is a Maxwellian} if $F_0(\cdot,\xi)$ is a
Maxwellian for every $\xi$. For product data in $\Dc$, the gradient $\nabla\Phi_1$ is bounded, so
this happens if and only if $\Phi_1$ is affine. In this case $F_0$ is an isotropic Gaussian of
temperature $(2\beta_F)^{-1}$, and (D1) forces its center to be $0$. Within $\Dc$ a Maxwellian $F_0$
is always of product form: if every $F_0(\cdot,v)$ is a Maxwellian, then $\Phi_1+\Phi_{\rm
v}(\cdot,v)$ is a quadratic polynomial, which is affine with a linear part independent of $v$, since
$\nabla\Phi_1$ and $\Phi_{\rm v}$ are bounded. Hence $F_0$ does not depend on $v$. A Gaussian with
another covariance does not belong to $\Dc$, but it belongs to $\Dc^\sharp$, with another $\beta_F$,
and Theorem~\ref{thm:t1} covers it as a profile posed at the time $t_1$, not at $t=0$. For
$\Phi_1=c\jap x$, $c>0$, we have $d_0>0$. The quantity $d_0$ appears only in the last clause of
Theorem~\ref{thm:main}(iv). Note that it is positive if the constant $d_1$ of
Theorem~\ref{thm:scatter}(iv) is, since $d_1\le C(\Dc)\,d_0(\xi)^{1/2}$ for every $\xi$, as follows
from $\|h\|_{L^1}\le C\|\jap x^2h\|_{L^2}\le C\|h\|^{1/2}_{L^2}\|\jap x^4h\|^{1/2}_{L^2}$ with
$h=F_0(\cdot,\xi)-M_0^\xi$. For product data, $d_0>0$ if and only if $F_0$ is not a Maxwellian, and
likewise for $d_1$. As in the introduction, the data are called non-thermal if $d_0>0$. For product
data this means that $F_0$ is not a Maxwellian, or equivalently that $d_1>0$. For $f_0\in\Dc^\sharp$
it can happen that no profile $F_0(\cdot,\xi)$ is a Maxwellian and still $d_0=d_1=0$ (see the
examples).

\emph{Role of the form {\rm(D$_{\rm iso}$)}.} The initial layer, treated in
Section~\ref{sec:layer}, uses the form (D$_{\rm iso}$) in two ways. The first is a null structure.
The only unbounded part of $\zeta=(\partial_v-t\nabla_y)\log f_0$
is $2\beta_Ft\,y$, and the corresponding part of $\zeta(y+tw,v-w)-\zeta(y,v)$ is $2\beta_Ft^2w$,
which lies in the kernel of $a(w)$. Lemmas~\ref{lem:GMwt} and~\ref{lem:ULaform} rest on this, and it
needs a Gaussian rate that does not depend on $v$. The second is the comparison of values at distant
velocities. The pairing of radial losses and gains, and the
moment bounds along the collision manifold, in Lemmas~\ref{lem:RBperp}, \ref{lem:RBdir}
and~\ref{lem:RBrad}, compare the values of $f_0$ at distant velocities, with constants that must be
uniform as $t\to0$. For $\rho_0$ we use the two-sided bound of (D3), and for the perturbation we use
$|\Phi_{\rm v}|\le C_{\rm pv}$.
A Gaussian rate that depends on $v$, a profile that is translated by an amount depending on $v$, or
a perturbation $e^{-\Phi_{\rm v}}$ that is not bounded above and below are excluded by this method.
Such profiles are covered at a positive time by Theorem~\ref{thm:t1}.

The Gaussian weight in space, the analyticity in $v$, and the special form (D$_{\rm iso}$) are used
in different places. As shown in Subsection~\ref{ssec:obstructions}, Gaussian decay makes the
mass-tail coupling dissipative. The transport loss is controlled by analyticity in $v$. The form
(D$_{\rm iso}$) matters only in the initial layer. We have not tried to optimize
these assumptions.

\subsection{Results from the literature and main results}\label{ssec:results}

\subsubsection{Results from the literature}

The proof uses published results on the linearized Landau operator near the Maxwellian $\mu$, on the
spatially homogeneous Landau equation, on entropy methods, and some general tools from the theory of
parabolic equations. We list them here in the form in which we use them. Where we need a variant,
we prove it in Appendix~\ref{app:CM} or~\ref{app:hom}.

\emph{Notation for the near-Maxwellian theory.} Let $\mu=M_{(1,0,1)}$, $
\mathcal Lf:=Q(\mu,f)+Q(f,\mu)$, and $E_\mu:=L^2(\mu^{-1/2})$, with $\langle f,g\rangle_{E_\mu}=\int
fg\mu^{-1}$. On $E_\mu$ the operator $\mathcal L$ is self-adjoint and nonpositive, with
$\ker\mathcal L=\mathrm{span}\{\mu,v_1\mu,v_2\mu,v_3\mu,|v|^2\mu\}$. We denote by $\Pi_0$ the
orthogonal projection in $E_\mu$ onto $\ker\mathcal L$ (equivalently, the projection along the
functions with vanishing moments $\int f$, $\int vf$, $\int|v|^2f$), and we put $\Pi:=I-\Pi_0$. We
write $\ab[\varphi]=a*\varphi$, $\bb[\varphi]=b*\varphi$, $\cb[\varphi]=c*\varphi$, and
$J_\beta(v):=\int|v-w|^\beta\mu(w)\dd w$. For $v\ne0$, $P_v$ is the orthogonal projection onto $\R
v$, and $\tn:=P_v\nabla+\jap v(I-P_v)\nabla$ is the anisotropic gradient. Following
\cite[\S2.2]{CM}, $\mathcal L=\mathcal A+\mathcal B$ with
\begin{equation}
\mathcal Af:=\ab[f]:\nabla^2\mu-\cb[f]\mu+M\chi_Rf,\qquad
\mathcal Bf:=\ab[\mu]:\nabla^2f-\cb[\mu]f-M\chi_Rf,
\end{equation}
where $\chi_R=\chi(\cdot/R)$ for a fixed cutoff $\chi\in C^\infty_c$, $
\mathbf 1_{B_1}\le\chi\le\mathbf 1_{B_2}$, and $M,R>0$ are large constants chosen in each statement. The
weights are $m=\jap v^k$ (polynomial, $\sigma=0$) and $m=e^{\kappa\jap v^2}$ with $0<\kappa<\frac12$
(Gaussian, $\sigma=2$). The dissipation norm of \cite[(3.1)]{CM} is
\begin{equation}
\|f\|^2_{Y(m)}:=\big\|\jap v^{(\gamma+\sigma)/2}mf\big\|^2_{L^2}+\big\|\jap
v^{\gamma/2}\tn(mf)\big\|^2_{L^2},
\end{equation}
and we put $X_k:=L^2(\jap v^k)$, $Y_k:=Y(\jap v^k)$, and $Z_k$ for the dual of $Y_k$ with respect to
the pairing of $X_k$. For a weight $m$ we define, as in \cite[(2.10)]{CM},
\begin{equation}
\tilde\zeta_m:=\ab[\mu]\frac{\nabla m}m\cdot\frac{\nabla m}m+\bb[\mu]\cdot\frac{\nabla
m}m-\frac12\cb[\mu] .
\end{equation}

\emph{Notation for the homogeneous equation.} $F=F(\sigma,v)$ denotes a solution of $\partial_\sigma
F=Q(F,F)$, $i(F):=\int|\nabla F|^2/F$ its Fisher information, $\|F\|_{L^p_k}:=\|\jap v^kF\|_{L^p}$,
and $F$ is \emph{normalized} if $\int F=1$, $\int Fv=0$ and $\int F|v|^2=3$. We write $
\mathcal H(F)=\int F\log F$, $\mathcal H(F|\mu)=\int F\log(F/\mu)$, $M_l(F):=\int\jap v^lF$, and $
\mathcal D(F):=-\int Q(F,F)\log F$ for the entropy dissipation.

\begin{hypothesis}[Results from the literature]\label{hyp:E}
Let $\gamma\in[-3,-2]$.
\begin{enumerate}[label=(E\arabic*)]
\item \emph{Near-Maxwellian theory.} For $\gamma\in[-3,-2)$ these are results of \cite{CM}. At
$\gamma=-2$, which \cite{CM} does not cover, (a) and (b) hold by published results whose range
includes $\gamma=-2$ and by Lemma~\ref{lem:ell1}, and (c), (d) and (f), in the weights in which we
use them, are proved (Lemmas~\ref{lem:CMpoly} and~\ref{lem:CMbil}, and items~(2)--(5) of the proof
of Proposition~\ref{prop:CMm2}) in Appendix~\ref{app:CM}, except~(e),
which is replaced there by the spectral gap of Lemma~\ref{lem:CMgap}.
\begin{enumerate}[label=(\alph*)]
\item \emph{Coercivity} \cite[(2.6)]{CM}, where it summarizes several works: \cite[Thm.~1.3]{Mou06}
proves it with $\nabla$ in place of $\tn$ in the second term, for $\gamma\ge-3$, and the anisotropic
second term is \cite[Lemma~5 and Cor.~1]{Guo02}, also for $\gamma\ge-3$. We use only the first
term. There is $c_\gamma>0$ with
\begin{equation}
\langle\mathcal Lf,f\rangle_{E_\mu}\le-c_\gamma\Big(\big\|\jap v^{(\gamma+2)/2}\Pi
f\big\|^2_{E_\mu}+\big\|\jap v^{\gamma/2}\tn\big(\mu^{-1/2}\Pi f\big)\big\|^2_{L^2}\Big)
\end{equation}
for all $f$ in the domain of $\mathcal L$ in $E_\mu$.
\item \emph{The coefficients} \cite[Lemma~2.1]{CM}, with~(ii) from \cite[Lemma~3]{Guo02}. For~(i),
which \cite{CM} take from \cite[Props.~2.3, 2.4]{DeL97}, a paper on hard potentials, we use instead
the symmetry argument and Lemma~\ref{lem:ell1} (item~(1) of the proof of
Proposition~\ref{prop:CMm2}), for every $\gamma\in[-3,-2]$. Part~(iv) rests in \cite{CM} on
\cite[Lemma~2.4(d)]{CTW}, which is stated there for $\gamma\ge-2$. However, its proof, from
(i)--(iii) and \cite[Lemma~2]{Guo02} (which assumes only $\theta>-3$), does not use this
restriction. The statements are: (i) $\ab[\mu](v)$ has the simple eigenvalue $\ell_1(v)>0$ with
eigenvector $v$ and the double eigenvalue $\ell_2(v)>0$ on $v^\perp$, and $\ell_1(v)\sim2\jap
v^\gamma$, $\ell_2(v)\sim\jap v^{\gamma+2}$ as $|v|\to\infty$. (ii) $|\partial^\beta\ab[\mu](v)|\le
C_\beta\jap v^{\gamma+2-|\beta|}$, and $\ab[\mu](v)\xi\cdot\xi\ge K\big(\jap v^\gamma|P_v\xi|^2+\jap
v^{\gamma+2}|(I-P_v)\xi|^2\big)$ with $K>0$. (iii)
$\mathrm{tr}\,\ab[\mu]=\ell_1+2\ell_2=2J_{\gamma+2}$ and $\bb[\mu](v)=-\ell_1(v)v$. (iv) For
$|v|>1$, $|\partial^\beta\ell_1|\le C_\beta\jap v^{\gamma-|\beta|}$ and $|\partial^\beta\ell_2|\le
C_\beta\jap v^{\gamma+2-|\beta|}$.
\item \emph{Dissipativity of $\mathcal B$} \cite[Lemmas~2.2, 2.3]{CM}. For $m=\jap v^k$ with
$k>\frac{\gamma+3}2$, $\limsup_{|v|\to\infty}\tilde\zeta_m\jap v^{-\gamma}\le2(\frac{\gamma+3}2-k)$.
For $m=e^{\kappa\jap v^2}$, $0<\kappa<\frac12$, $\limsup_{|v|\to\infty}\tilde\zeta_m\jap
v^{-2-\gamma}\le4\kappa(2\kappa-1)$. For each such $m$ and $f$ smooth and rapidly decaying,
\begin{equation}
\int(\mathcal Bf)fm^2=-\int\ab[\mu]\nabla(mf)\cdot\nabla(mf)+\int(\tilde\zeta_m-M\chi_R)f^2m^2 ,
\end{equation}
and consequently:
\begin{itemize}
\item[(S0)] for $k>\frac{\gamma+3}2$ there are $M,R$ and $K^{\mathcal B}_k>0$ with $\langle
\mathcal Bh,h\rangle_{X_k}\le-K^{\mathcal B}_k\|h\|^2_{Y_k}$.
\end{itemize}
\item \emph{Decay and regularization of $e^{t\mathcal B}$} \cite[Lemmas~3.1, 3.2]{CM}, polynomial
weights:
\begin{itemize}
\item[(S1)] for $\frac{\gamma+3}2<k<k_1$ and $0\le b<\frac{k_1-k}{|\gamma|}$, $\|e^{t
\mathcal B}\|_{X_{k_1}\to X_k}\le C\jap t^{-b}$;
\item[(S2)] for $\frac32<k<k_1$ and $0\le b<\frac{k_1-k}{|\gamma|}$, $\|e^{t
\mathcal B}\|_{Z_{k_1}\to X_{k-|\gamma|/2}}\le C(t\wedge1)^{-1/2}\jap t^{-b}$.
\end{itemize}
(In \cite{CM} the rate is $\jap t^{-(k_1-k_\sharp)/|\gamma|}$ with any $k_\sharp\in(k,k_1)$, hence
the strict inequality on $b$. In \cite[(3.6)]{CM} the source space $H^{-1}_*(\jap v^{k_1})$ is
$Z_{k_1}$ and the target space $L^2(\jap v^k\jap v^{\gamma/2})$ is $X_{k-|\gamma|/2}$.)
\item \emph{Decay in a small space} \cite[Prop.~3.3]{CM}, for $\gamma\in[-3,-2)$:
\begin{itemize}
\item[(S3)] for $\vartheta\in(\frac12,1)$,
\begin{equation}\label{eq:Theta-s}
\Theta_{\rm s}(t):=\|e^{t\mathcal L}\Pi\|_{L^2(\mu^{-\vartheta})\to E_\mu}\le Ce^{-\lambda
t^{2/|\gamma|}}
\end{equation}
for some $\lambda>0$. In particular, $\Theta_{\rm s}$ is integrable and $\Theta_{\rm s}(t)\le C\jap
t^{-1}$, and we use only these two properties.
\end{itemize}
\item \emph{The collision operator in weighted spaces} \cite[Lemmas~4.1--4.3, Cor.~4.4]{CM}: (i)
for $-3<\alpha<0$ and $\theta>3$, $\int|v-w|^\alpha\jap w^{-\theta}\dd w\le C\jap v^\alpha$. (ii)
For $\frac3{5+\gamma}<p\le\infty$ and $\theta>2+3(1-\frac1p)$, $|(a*f)(v)v\cdot
v|+|(a*f)(v)v|+|(a*f)(v)|\le C\jap v^{\gamma+2}\|f\|_{L^p(\jap v^\theta)}$, and for
$\frac3{4+\gamma}<p\le\infty$ and $\theta'>3(1-\frac1p)$, $|(b*f)(v)|\le C\jap
v^{\gamma+1}\|f\|_{L^p(\jap v^{\theta'})}$. (iii) For every weight $m$ as above with $m\ge1$,
$\theta>\frac72$ and $\theta'>\frac94$,
\begin{equation}
\big|\langle Q(f,g),h\rangle_{L^2(m)}\big|\le C\Big(\|f\|_{L^2(\jap
v^\theta)}\|g\|_{Y(m)}+\|f\|_{H^1(\jap v^{\theta'})}\|g\|_{L^2(m)}\Big)\|h\|_{Y(m)} ;
\end{equation}
(iv) for $m=\jap v^k$, $k>\max\{\frac72,\frac94+\frac{|\gamma|}2\}$ (the second condition is needed
for the embedding of \cite[Cor.~4.4]{CM}, and we use only $k\ge7$), with $X=L^2(m)$, $Y=Y(m)$ and
$Z$ the dual of $Y$ for the pairing of $X$, $|\langle Q(f,g),h\rangle_X|\le
C(\|f\|_X\|g\|_Y+\|f\|_Y\|g\|_X)\|h\|_Y$ and $\|Q(f,g)\|_Z\le C(\|f\|_X\|g\|_Y+\|f\|_Y\|g\|_X)$.
\end{enumerate}
Instead of \cite[Lemma~2.5]{CM} we use its variant on the spaces $L^2(\jap v^s)$,
Lemma~\ref{lem:RVA}, whose part~(i) also provides the bound that the proof of (S3) in \cite{CM}
takes from \cite[Lemma~2.5]{CM}. What we need from \cite[Thm.~3.5, Prop.~3.6,
Cor.~3.7, Prop.~5.1]{CM} is derived from (a)--(f) in Proposition~\ref{prop:CMthr} and in the proof
of Theorem~\ref{thm:enter}. For part~(b) we also use \cite[Lemmas~2, 3]{Guo02} and the proof of
\cite[Lemma~2.4(d)]{CTW}, for every $\gamma\in[-3,-2]$.
\item \emph{Homogeneous Landau theory} \cite{GS}, partly based on \cite{D15,Si17,CSS}. Let
$F_0\ge0$ have finite mass, energy and entropy.
\begin{itemize}
\item[(GS1)] \cite[Thms.~1.1, 1.2]{GS}: if $F_0\le C_{\rm dat}e^{-\beta|v|^2}$, there is a unique
global classical solution, positive and Schwartz for $\sigma>0$, and $i(F(\sigma))$ is
nonincreasing.

We also use the local theory that \cite{GS} quote as \cite[Thm.~2.1]{GS}, from \cite{HST20}, for
$\gamma\in[-3,0)$: for $F_0\in L^\infty_{k_0}$, $k_0=\max\{5,15/(5+\gamma)\}$, there are $T>0$ and a
classical solution on $[0,T]$ (in \cite[Thm.~1.2]{HST20} the hypothesis is $F_0\in L^\infty_k$ for
some $k>k_0$, which holds for all data we use), strictly positive and $C^\infty$ for $\sigma>0$,
with bounds on all derivatives for $\sigma>0$, which takes the data continuously if $F_0$ is
continuous, and which can be continued as long as $\|F(\sigma)\|_{L^\infty_{k_0}}$ stays bounded.

The result of \cite{HST20} is stated for the inhomogeneous equation, and \cite{GS} use it for data
independent of $x$. We use two of its properties. First, the solution is bounded in $L^\infty_{k}$
on $[0,T]$, $k>k_0$, with a bound and a time $T$ that depend only on $\gamma$,
$\|F_0\|_{L^\infty_k}$ and the lower bound parameters of the data \cite[Thm.~1.2(i)]{HST20}. Second,
if $F_0$ is continuous, $F(\sigma)\to F_0$ as $\sigma\to0^+$, uniformly on compact sets, and
uniformly if $F_0$ is uniformly continuous in $x$ \cite[Prop.~3.1]{HST20}. Therefore, for data
independent of $x$, the convergence is uniform in $v$. Together with the continuity of $F_0$ and of
$F$ for $\sigma>0$, this implies that $F$ is continuous on $[0,T]\times\R^3$. The continuation
criterion is \cite[Thm.~1.3]{HST20}, in terms of a quantity bounded by the $L^\infty_k$ norm, $k>5$.
The a priori estimate is \cite[Thm.~1.1]{GS}:
along classical solutions, $i(F)$ is nonincreasing, for $\gamma\in[-3,1]$.

From these, (GS2), (GS3) and our Propositions~\ref{prop:GS4} and~\ref{prop:GS5},
we derive in Lemma~\ref{lem:GS1q} that, for normalized continuous $F_0\le C_{\rm
dat}e^{-\beta|v|^2}$, the solution of (GS1) is continuous on $[0,\infty)\times\R^3$ with $F(0)=F_0$,
and is bounded by a Maxwellian on every $[0,T]$. We use this only qualitatively, in
Lemma~\ref{lem:RSuniq}, to identify the solution near $\sigma=0$ with the one constructed in
Lemma~\ref{lem:RSshortF}. We prove the quantitative bounds near $\sigma=0$ there.
\item[(GS2)] \cite[Thm.~2.4]{GS}, from \cite[Prop.~2]{D15} (see also \cite[Prop.~4]{D15} and (CDH3)
below): for $k\ge0$, $\|F(\sigma)\|_{L^1_k}\le C(S)$ for $\sigma\le S$, with $C(S)$ depending only
on $S$, $\|F_0\|_{L^1_k}$ and the mass, energy and entropy of $F_0$. (\cite[Prop.~2]{D15} states it
for $\gamma=-3$ and $k>2$, with the momentum also among the data, which mass and energy bound,
and \cite{GS} note that the proof applies for $\gamma\in[-3,0]$.)
\item[(GS3)] \cite[Thm.~2.2]{GS}, from \cite{Si17} (cited in \cite{GS} as Thm.~3.7): if
$\|F(\sigma)\|_{L^p_k}\le W$ on $[\sigma_1,S]$, with $p>\frac3{5+\gamma}$, $p\ge1$ and $k\ge k_{\rm
GS}(\gamma,p)$, then $\|F(\sigma)\|_{L^\infty}\le C(1+(\sigma-\sigma_1)^{-3/(2p)})$ there, with $C$
depending only on $W$, the mass, the energy and the entropy (not on $S$). The exponent $\frac3{2p}$
is that of \cite{GS}, while \cite[Thm.~3.8]{Si17} yields the factor $(\sigma-\sigma_1)^{-3/2}$ for
$\sigma-\sigma_1\le T_*$, which suffices, since (GS3) is used only where $\sigma-\sigma_1$ is
bounded below (Lemma~\ref{lem:Abarrier}(i), Appendix~\ref{sssec:homRM}, Lemma~\ref{lem:GS1q}). For
$\gamma\in[-3,-2)$ and $\frac3{5+\gamma}<p<\frac3{3+\gamma}$ ($p<\infty$ at $\gamma=-3$),
\cite[Thm.~3.8]{Si17} (\cite{GS} refer to \cite[Thm.~3.7]{Si17}, which covers $\gamma\in[-2,0]$)
requires $\int\jap v^\kappa F^p\le W_0$ with
$\kappa>\max\{p(5+\gamma)-3,\,2+\frac92(1-\frac1p)-3(1+\frac\gamma2)\}$, so $k_{\rm
GS}(\gamma,p)=\kappa/p$. For $p=\frac52$ this is $k>\frac{17}{25}+\frac35|\gamma|$. At $\gamma=-2$,
\cite[Thm.~3.7]{Si17} provides the bound from the mass, energy and entropy alone.
\item[(GS4)] \cite[Thm.~2.3]{GS}, from \cite{CSS}: $L^\infty$ and $L^1_k$ bounds ($k$ large) on
$[0,S]$ and $F_0\le C_{\rm dat}e^{-\beta|v|^2}$ with $\beta$ small imply $F\le C_{\rm
Gauss}e^{-\beta|v|^2}$ on $[0,S]$. \cite{GS} give a sketch, and for $\gamma\in[-3,-2]$ we prove this
in Proposition~\ref{prop:GS4}, with $\beta\le\beta_0(\gamma)$, moments of order $k_4(\gamma)$, and
$C_{\rm Gauss}$ independent of $S$.
\item[(GS5)] \cite[Prop.~B.2]{GS}: a Gaussian bound on $[0,T]$ implies Gaussian bounds for $\nabla
F$ and $\nabla^2F$ on $[T/2,T]$, with a slightly smaller exponent. \cite{GS} give a sketch, and for
$\gamma\in[-3,-2]$ we prove this in Proposition~\ref{prop:GS5}.
\end{itemize}
\item \emph{Entropy methods.} Let $F$ be smooth, positive and normalized. We put $P_{ij}:=\int
Fv_iv_j$ and, for $i\ne j$,
\begin{equation}
R_{ij}(v,w):=(v_i-w_i)\Big(\frac{\partial_jF}F(v)-\frac{\partial_jF}F(w)\Big)-(v_j-w_j)\Big(\frac{\partial_iF}F(v)-\frac{\partial_iF}F(w)\Big).
\end{equation}
\begin{itemize}
\item[(CDH1)] \cite[Prop.~3]{CDH}: for $i\ne j$,
\begin{equation}
\frac{\partial_iF}F(v)=\frac{v_jP_{ij}+v_iP_{ii}+\int R_{ij}(v,w)F(w)\,[w_iP_{ij}-w_jP_{ii}]\dd
w}{P_{ij}^2-P_{ii}P_{jj}},\end{equation}\begin{equation}
\frac{\partial_jF}F(v)=\frac{v_iP_{ij}+v_jP_{jj}+\int R_{ij}(v,w)F(w)\,[w_iP_{jj}-w_jP_{ij}]\dd
w}{P_{ij}^2-P_{ii}P_{jj}},
\end{equation}
and $v_i\frac{\partial_jF}F(v)-v_j\frac{\partial_iF}F(v)=\int R_{ij}(v,w)F(w)\dd w$.
\item[(CDH2)] \cite[Prop.~5, (32)]{CDH}: if $\mathcal H(F)\le\bar H$, then $\inf_{i\ne
j}(P_{ii}P_{jj}-P_{ij}^2)\ge c(\bar H)>0$. The proof in \cite{CDH} uses a bound on $\sup_{|A|\le
q}\int_AF$ in terms of the entropy \cite[(35)]{CDH}. We use it in the form $\sup_{|A|\le
q}\int_AF\le Mq+(\bar H^++C_*)/\log M$, where $C_*$ is an absolute bound of $\int F\log^-F$ over
normalized $F$. The term $C_*$ accounts for the part of $\int_AF\log F$ where $F<1$, and this form
also covers $\bar H\le0$, which can occur in our setting since $\mathcal H(\mu)<0$.
\item[(CDH3)] \cite[Lemma~8]{CDH}: for $-4<\gamma\le-2$, $l>2$ and an H-solution with normalized
data $F_0\in L^1_l$, $M_l(F(\sigma))\le
C\,M_l(F_0)+C\,l^{(l-6)\frac{|\gamma+1|}{\gamma+4}-\gamma}\,\sigma$, with $C=C(\gamma,
\mathcal H(F_0))$.
\item[(CKP)] \cite{Csi67,Kul67}: for probability densities $F,G$, $\|F-G\|_{L^1}^2\le2\int
F\log(F/G)$.
\item[(LSI)] \cite[Cor.~5.7.2]{BGL}: for a probability measure $e^{-W}\dd v$ on $\R^3$ with $W$
smooth and $\nabla^2W\ge I$, $\int f^2\log f^2\,e^{-W}-\big(\int f^2e^{-W}\big)\log\big(\int
f^2e^{-W}\big)\le2\int|\nabla f|^2e^{-W}$. We use its consequence on convex sets, (EL4) in
Appendix~\ref{app:hom}, obtained by approximation.
\end{itemize}
In Lemma~\ref{lem:EL1} we prove from (CDH1), (CDH2) and the argument of \cite[pp.~395--396]{CDH}
the entropy--entropy-dissipation inequality \cite[(19)]{CDH} for $-4<\gamma\le-2$, and
\cite[Thm.~1]{CDH} is its case $\gamma=-3$.
\item \emph{General tools.}
\begin{itemize}
\item[(CK)] \emph{Abstract Cauchy--Kovalevskaya theorem} \cite[Theorem]{Nis77}, which contains
\cite[Thm.~1.1]{Ni72}. Let $(B_\rho)_{0<\rho\le\rho_0}$ be a scale of Banach spaces, $B_\rho\subset
B_{\rho'}$ and $\|\cdot\|_{\rho'}\le\|\cdot\|_\rho$ for $\rho'<\rho$. Suppose that for some
$R,\eta>0$ and all $0<\rho'<\rho<\rho_0$, $(u,t)\mapsto F(u,t)$ is continuous from
$\{\|u\|_\rho<R\}\times\{|t|<\eta\}$ into $B_{\rho'}$, that $\|F(u,t)-F(v,t)\|_{\rho'}\le
C\|u-v\|_\rho/(\rho-\rho')$ for $\|u\|_\rho,\|v\|_\rho<R$ and $|t|<\eta$, with $C$ independent of
$t,u,v,\rho,\rho'$, and that $F(0,t)$ is continuous in $t$ with $\|F(0,t)\|_\rho\le
K/(\rho_0-\rho)$. Then there is $a>0$ and a unique $u$, continuously differentiable with values in
$B_\rho$ and with $\|u(t)\|_\rho<R$ for $|t|<a(\rho_0-\rho)$, $0<\rho<\rho_0$, such that $u'=F(u,t)$
and $u(0)=0$.
\item[(PS)] \emph{Schauder theory for linear parabolic equations}
$\partial_tu=a_{ij}\partial_{ij}u+b_i\partial_iu+cu+f$ with $a_{ij}$ uniformly elliptic. (i)
\cite[Ch.~IV, Thm.~5.2]{LSU}: on a bounded cylinder $\Omega\times(0,T)$ with $\partial\Omega\in
H^{l+2}$, $l>0$ not an integer, coefficients and $f$ in $H^{l,l/2}$, and initial and boundary data
in $H^{l+2}$ and $H^{l+2,l/2+1}$ satisfying the compatibility conditions of order $[l/2]+1$, the
first boundary value problem has a unique solution in $H^{l+2,l/2+1}$, bounded by the norms of the
data. (ii) \cite[Ch.~IV, Thm.~10.1]{LSU}: local estimates in $H^{l+2,l/2+1}$ on subcylinders, in
terms of $f$ in $H^{l,l/2}$, the data on the part of the boundary that is included, and $\max|u|$.
(iii) \cite[Thm.~1]{Kn80}: interior Schauder estimates in $x$ for coefficients that are H\"older
continuous in $x$, uniformly in $t$, and only measurable in $t$. (iv) \cite[Thm.~4.1]{Lor00}, with
the polynomial weight of degree $0$ ($p\equiv2$, equivalent to $p\equiv1$): the estimate of
$\sup_t\|u(t)\|_{C^{2+\theta}(\R^n)}$ on strips $[0,T]\times\R^n$, up to $t=0$, for coefficients
that are bounded and measurable and $\theta$-H\"older continuous in $x$ uniformly in $t$
\cite[(H1)--(H3)]{Lor00}. (v) \cite[Chs.~IV, V]{Lie96}: Schauder estimates and existence on bounded
domains.
\item[(KS)] \emph{Krylov--Safonov} \cite[Thm.~7.37, Cor.~7.41, Cor.~7.45]{Lie96}: H\"older
estimates, interior and at the parabolic boundary, for nondivergence parabolic equations with
bounded measurable uniformly elliptic coefficients.
\item[(C)] \emph{Compactness and fixed points}: the Aubin--Lions--Simon lemma \cite[Cor.~4]{Sim87}
(if $X\subset B\subset Y$ with $X\subset B$ compact, a family bounded in $L^\infty(0,T;X)$ with time
derivatives bounded in $L^r(0,T;Y)$, $r>1$, is relatively compact in $C([0,T];B)$), the
Lions--Magenes lemma \cite[Ch.~XVIII, \S1, Thms.~1--2]{DL} (in a Gelfand triple $V\subset H\subset
V'$, $u\in L^2(0,T;V)$ with $u'\in L^2(0,T;V')$ lies in $C([0,T];H)$, and $\frac
d{dt}\|u\|_H^2=2\langle u',u\rangle$), and Schauder's fixed point theorem.
\end{itemize}
\end{enumerate}
\end{hypothesis}

The near-Maxwellian theory of \cite{CM} treats $\gamma\in[-3,-2)$. At $\gamma=-2$ we prove the
estimates of \cite{CM} that we use in Appendix~\ref{app:CM}, following the proofs of \cite{CM}, with
the decay in a small space \cite[Prop.~3.3]{CM} replaced by the spectral gap of the linearized
operator (the weight $\jap v^{(\gamma+2)/2}$ in \cite[(2.6)]{CM} equals $1$). In that appendix
we also prove the bound on the regularizing part $\mathcal A$ on the spaces $L^2(\jap v^s)$, which
for small $s$ are not contained in $L^1$. The weight thresholds given by
our argument may not be optimal, and we do not compare them with those of \cite{CM}.

\subsubsection{The main theorem}

\begin{theo}[Main theorem]\label{thm:main}
Let $\gamma\in[-3,-2]$, $N_y\ge20$, and $f_0=\eps\rho_0F_0\in\Dc$. There are constants
\begin{equation}
q_0,\ \delta_1,\ \delta_s,\ S_*,\ \rad_{\rm post},\ T_0,\ C>0,
\end{equation}
depending only on $\gamma$ and the class $\Dc$, and $\eps_0>0$, such that for $\eps\le\eps_0$ the
following hold.
\begin{enumerate}[label=(\roman*)]
\item \emph{Global existence and analyticity.} There is a global, classical, nonnegative solution
of \eqref{eq:landau}. For $t\ge T_0$ its profile $g(t,y,\cdot)$ is real-analytic, with radius at
least $\rad_{\rm post}/2$ in $v$, uniformly in $(t,y)$.
\item \emph{Every leaf relaxes.} For every $\xi$, $s_\xi(t)\to\infty$ as $t\to\infty$. If $t\ge
T_0$ and $s_{x/t}(t)\ge S_*$, write $p=p_g(t,x/t)=(\rho,y_0,\theta)$. Then
\begin{equation}
\Big\|\frac{f(t,x,(x-\cdot)/t)-M_{p}}{\rho}\Big\|_{L^2(e^{q_0\theta_-|\cdot-y_0|^2/(8\theta)})}\le
C\delta_1 ,
\end{equation}
where $\theta_->0$ is the lower bound of the fiber temperatures at $T_0$ from (H1).
\item \emph{Leaf hydrodynamics.} Let $\varpi'_*$ be the solution of the leaf Euler system of
Theorem~\ref{thm:Meuler} without micro fluxes, with data $\varpi'_g(T_0)$ (Lemma~\ref{lem:Mref}).
Then $\|\varpi'_g(t)-\varpi'_*(t)\|_{L^\infty_\xi}\le\delta_s$ for all $t\ge T_0$, and there is
$\varpi'_\infty$ with $\|\varpi'_g(t)-\varpi'_\infty\|_{L^\infty_\xi}\le C/t$.
\item \emph{Non-uniformity.} The entry time $t_\xi$, defined by $s_\xi(t_\xi)=S_*$, tends to
$\infty$ as $|\xi|\to\infty$. For every $t$ there are leaves with $s_\xi(t)<S_*$. Suppose in
addition that the data are non-thermal, that is, $d_0>0$ in~\eqref{eq:d0}. There is $C_{\rm nM}$,
depending only on $\gamma$ and $\Dc$, such that, if $T_0\ge C_{\rm nM}/d_0$ and $C_{\rm
nM}\eps\,(1+\tau(T_0))\le d_0$, then at every $t\ge T_0$ there are $x$ at which the left side
of~(ii) is $\ge d_0/2$. In particular, for no $t\ge T_0$ does a bound of the form~(ii) with right
side $<d_0/2$ hold uniformly in $x$. If $F_0$ is a Maxwellian, then at
$t=T_0$ the bound of~(ii) holds uniformly in $x$, with right side $C(L/T_0+\eps
T_0^{\kappa_\gamma})$ (Remark~\ref{rem:nonunif}). For such data the non-uniformity concerns the
clock $s_\xi$, and not necessarily the distance from equilibrium, which, at times $t>T_0$, is not
shown to be small uniformly in $x$.
\end{enumerate}
Here $s_\xi$ is the leaf clock of Subsection~\ref{ssec:cells}, and its rate $\lambda_{\rm ref}(\xi)$
is the collision rate of the fiber $g(T_0,\cdot,\xi)$, and $\lambda_{\rm ref}/\lambda\in[\frac12,2]$
for $\eps\le\eps_0(T_0)$ (by the closeness clause of (LWP)(a3), Section~\ref{sec:proofmain}). At
$\gamma=-2$ we have $s_\xi(t)=\lambda_{\rm ref}(\xi)\log(t/T_0)$, and the entry times are
$t_\xi=T_0e^{S_*/\lambda_{\rm ref}(\xi)}$.
\end{theo}

All requirements on $T_0$ in Theorem~\ref{thm:main} are lower bounds. Hence we can enlarge $T_0$
to satisfy the conditions of the last clause of~(iv). For this clause $T_0$ and $\eps_0$ then depend
also on $d_0$, and (i)--(iii) hold for the new $T_0$, with the clock $s_\xi$ started there.

We use the smallness of $\eps$ to make the effective collision time on $[0,T_0]$ small. We do not
assume that $f_0$ is small relative to a Maxwellian. The near-equilibrium regime is entered leaf by
leaf, after an effective time $S_*$ that depends only on $\Dc$.

We split the proof at a time $t_1=t_1(\Dc)$, before which the solution is close to free transport in
a strong sense. The part of the proof after $t_1$ does not use (D$_{\rm iso}$), and it leads to a
statement of independent interest.

\begin{theo}[Data posed at $t_1$]\label{thm:t1}
Let $\gamma\in[-3,-2]$ and $N_y\ge20$. Pose the Cauchy problem at time $t_1=t_1(\Dc)$ with data
$f(t_1)$ whose profile $g(t_1)$ is close to a reference profile
$f_0=\eps\rho_0(v)F_0(x,v)\in\Dc^\sharp$ (that is, $f_0$ satisfies (D1), (D3) and (D4), but it need
not be a product), in the following sense, for some constant $C_{t_1}$:
\begin{itemize}
\item $|g(t_1)-f_0|\le C_{t_1}\eps f_0$ pointwise;
\item (G1)--(G3) hold with constants within a factor $1+C_{t_1}\eps$ of those of
Lemma~\ref{lem:LWfree}(i);
\item relative to $\eps$, $g(t_1)$ is within $C_{t_1}\eps$ of $f_0$ in the raw leaf norm $
\mathfrak B_\xi$ of Subsection~\ref{sssec:FPtr}, with the weight $e^{2q_d|z|^2}$, at the radius $\rad_1$
and the $y$-weight $\delta_{y,\ell}$, that is, \eqref{eq:dispdata} with $C_{\rm lay}\eps^2t_1^2$
replaced by $C_{t_1}^2\eps^2$ (data parameters, Appendix~\ref{app:const}).
\end{itemize}
The exact free data are a special case. Here $t_1$ and the data parameters
$(\rad_1,\delta_{y,\ell})$ are those of Appendix~\ref{app:const}. Since (D$_{\rm iso}$) is not
assumed, they are part of the hypothesis when Theorem~\ref{thm:t1} is used alone.

Then the conclusions of Theorem~\ref{thm:main} hold, with $\rho_\infty$, $y_{0,\infty}(\xi)$,
$\theta_\infty(\xi)$, $\lambda(\xi)$ and $d_0$ those of $f_0$, and with $\eps_0$ depending also on
$C_{t_1}$ (through the $\eps$-conditions of Proposition~\ref{prop:disp} and
Corollary~\ref{cor:LWPa}, where $C_{t_1}^2$ takes the place of $C_{\rm lay}t_1^2$, and the condition
$C_{t_1}\eps\le1$). The constant $C_w$ of Proposition~\ref{prop:disp} then depends also on
$C_{t_1}$, and so do the constants that are built from it: $C_{\rm nM}$ in
Theorem~\ref{thm:main}(iv), and $C_{\rm sh}$, $c_{\rm sh}$ in Theorems~\ref{thm:rate2}(iii)
and~\ref{thm:rate3}(iii).
\end{theo}

Theorem~\ref{thm:main} is Theorem~\ref{thm:t1} together with the initial layer, which brings data
posed at $t=0$ to the form required at $t_1$: Theorem~\ref{thm:layer} for $\gamma<-2$ and
Theorem~\ref{thm:EPlayer} at $\gamma=-2$. The second one holds for every $\gamma\in[-3,-2]$, with
other energies and the bound $O(\eps t^{3/4})$. Its proof follows that of the first.

Theorem~\ref{thm:main}(ii) bounds the distance of a relaxed leaf from its Maxwellian by a fixed
small constant. In fact this distance tends to zero, and with it the solution converges, in the
profile variables, to a limit whose fibers are Maxwellians.

\begin{theo}[Complete relaxation, and scattering to a thermalized state]\label{thm:scatter}
Under the assumptions of Theorem~\ref{thm:main} or of Theorem~\ref{thm:t1}, let $g$ be the profile
of the solution constructed there, and let $\mathcal N_\xi(t)$, $\hat p_g(t,\xi)$ and $d_\xi(t)$ be
the normalized leaf, its normalized parameters and its distance from the Maxwellian with the same
moments, as in Lemma~\ref{lem:Wtail}(ii).
\begin{enumerate}[label=(\roman*)]
\item \emph{Complete relaxation.} For every $\xi\in\R^3$, $d_\xi(t)\to0$ as $t\to\infty$, and the
same holds for the distance in the weighted norm of Theorem~\ref{thm:main}(ii).
\item \emph{Asymptotic parameters.} There is a continuous function $p_{\rm sc}=(\rho_{\rm
sc},y_{\rm sc},\theta_{\rm sc})$ on $\R^3$ with
\begin{equation}
\sup_{\xi\in\R^3}\Big(\Big|\log\frac{\rho_g(t,\xi)}{\rho_{\rm sc}(\xi)}\Big|+|y_{0,g}(t,\xi)-y_{\rm
sc}(\xi)|+|\theta_g(t,\xi)-\theta_{\rm sc}(\xi)|\Big)\le\frac Ct\qquad(t\ge T_0),
\end{equation}
with a constant $C$ that depends only on $\gamma$, $\Dc$ and $T_0$, and not on $\eps$. Moreover,
$\rho_{\rm sc}$ is of class $C^1$, $\rho_{\rm sc}>0$, $|\nabla\log\rho_{\rm sc}|\le C$,
$\frac12\theta_-\le\theta_{\rm sc}\le\theta_++1$ and $|y_{\rm sc}|\le L+1$, where $\theta_\pm$ are
the bounds of the fiber temperatures at $T_0$ from (H1) (Subsection~\ref{ssec:fibstanding}). In
particular $\mathcal N_\xi(t)\to M_{(1,y_{\rm sc}(\xi),\theta_{\rm sc}(\xi))}$ in $L^2(\R^3)$ for
every $\xi$.
\item \emph{Scattering.} Let $g_{\rm sc}(y,v):=M_{p_{\rm sc}(v)}(y)$. Then $g_{\rm sc}\in
L^1(\R^3_y\times\R^3_v)$, and
\begin{equation}
\big\|f(t,x,v)-g_{\rm sc}(x-tv,v)\big\|_{L^1(\R^3_x\times\R^3_v)}=\|g(t)-g_{\rm
sc}\|_{L^1(\R^3_y\times\R^3_v)}\longrightarrow0\qquad(t\to\infty).
\end{equation}
\item \emph{The scattering state is not close to the data.} Let
\begin{equation}
d_1(v):=\inf_{c,M}\|cM-F_0(\cdot,v)\|_{L^1(\R^3)},\qquad d_1:=\inf_{v\in\R^3}d_1(v),
\end{equation}
where the first infimum is over $c\ge0$ and the Maxwellians $M$ of unit mass. Then $d_1(\cdot)$ is
continuous (it is $1$-Lipschitz with respect to the $L^1$-distance of the profiles, and $v\mapsto
F_0(\cdot,v)$ is continuous into $L^1$ by (D1) and (D4)), $d_1(v)>0$ if and only if $F_0(\cdot,v)$
is not a Maxwellian, and
\begin{equation}
\|g_{\rm sc}-f_0\|_{L^1(\R^6)}\ge\eps\int\rho_0(v)\,d_1(v)\dd v\ge d_1\,\|f_0\|_{L^1(\R^6)} .
\end{equation}
The constants $d_1(v)$ and $d_1$ depend only on $F_0$, not on $\eps$. The middle term is positive as
soon as $F_0(\cdot,v)$ is not a Maxwellian for $v$ in a set of positive measure. In this case
$\|g_{\rm sc}-f_0\|_{L^1}\ge c\,\|f_0\|_{L^1}$ with $c=\int\rho_0d_1(v)\dd v/\int\rho_0>0$, which
does not depend on $\eps$, also when $d_1=0$. For product data $d_1(v)=d_1$, and $d_1>0$ if and only
if $F_0$ is not a Maxwellian. Under the assumptions of Theorem~\ref{thm:t1}, $f_0$ is the reference
profile. Since $|g(t_1)-f_0|\le C_{t_1}\eps f_0$, the same lower bounds hold with $g(t_1)$ in place
of $f_0$, up to $C_{t_1}\eps\|f_0\|_{L^1}$.
\end{enumerate}
\end{theo}

Thus the solution scatters, but the profile $g_{\rm sc}$ of the limiting free solution is a
Maxwellian in $y$ on every fiber, and by~(iv), if $d_1>0$, the map from the data to the scattering
state is not close to the identity, however small the data are. We give rates for the convergence in
(i) and (iii) in Theorems~\ref{thm:rate2} and~\ref{thm:rate3} below. If $d_0>0$ and $T_0$ and
$\eps$ satisfy the conditions of the last clause
of Theorem~\ref{thm:main}(iv), the convergence in~(i) is not uniform in $\xi$: as the proof of
this part shows, at every $t\ge T_0$ there are leaves with $d_\xi(t)\ge d_0/2$. If
$d_0=0$, in particular if $F_0$ is a Maxwellian, we do not know whether it is uniform (see
Remark~\ref{rem:nonunif}). The free solution $g_{\rm sc}(x-tv,v)$ is, for large $t$ and at each $x$,
close to a Maxwellian in $v$ with density $t^{-3}\rho_{\rm sc}(x/t)$, mean velocity $(x-y_{\rm
sc}(x/t))/t$ and temperature $t^{-2}\theta_{\rm sc}(x/t)$. This is the hydrodynamic state of free
expansion that the leaf Euler system reaches at the end of its finite duration.

At the endpoint $\gamma=-2$ the dissipation controls the energy, and the convergence in
Theorem~\ref{thm:scatter} comes with a rate. Recall that then $s_\xi(t)=\lambda_{\rm
ref}(\xi)\log(t/T_0)$.

\begin{theo}[Rates at $\gamma=-2$]\label{thm:rate2}
Let $\gamma=-2$, under the assumptions of Theorem~\ref{thm:main} or of Theorem~\ref{thm:t1}. Let
$S_b\in(S_*,S_*+1]$ be the constant of Subsection~\ref{ssec:cells}, and for every $\xi$ let
\begin{equation}
t^b_\xi:=T_0\,e^{(S_b+1)/\lambda_{\rm ref}(\xi)},\qquad a_\xi:=\min\big\{c_{\rm rate}\lambda_{\rm
ref}(\xi),\tfrac12\big\},
\end{equation}
so that $s_\xi(t^b_\xi)=S_b+1$. Here $c_{\rm rate}>0$ depends only on $\Dc$. There is a constant
$C$, which depends only on $\Dc$ and $T_0$, and not on $\eps$, such that:
\begin{enumerate}[label=(\roman*)]
\item \emph{Relaxation of the leaves.} For every $\xi$ and every $t\ge t^b_\xi$, with $x=t\xi$ and
$(\rho,y_0,\theta)=p_g(t,\xi)$,
\begin{equation}
\Big\|\frac{f(t,x,(x-\cdot)/t)-M_{p_g(t,\xi)}}{\rho}\Big\|_{L^2(e^{q_0\theta_-|\cdot-y_0|^2/(8\theta)})}\le
C\Big(\frac{t^b_\xi}t\Big)^{a_\xi/2} .
\end{equation}
If $c_{\rm rate}\lambda_{\rm ref}(\xi)\le\frac12$, the right side is
$C\exp\big(-\frac{c_{\rm rate}}2(s_\xi(t)-S_b-1)\big)$: the leaf relaxes exponentially in its own
clock.
\item \emph{Scattering.} With $g_{\rm sc}$ as in Theorem~\ref{thm:scatter} and $p\ge4$ the exponent
of (D3), for all $t\ge T_0$,
\begin{equation}
\|g(t)-g_{\rm sc}\|_{L^1(\R^6)}\le C\,\|f_0\|_{L^1(\R^6)}\,\big(1+\eps\log(t/T_0)\big)^{-(p-3)/p} .
\end{equation}
\item \emph{The rate in~(ii) is sharp.} Suppose in addition that the constant $d_1$ of
Theorem~\ref{thm:scatter}(iv) is positive (this holds, in particular, for non-thermal product data).
Then $d_1\in(0,1]$. There are constants $C_{\rm sh}$ and $c_{\rm sh}>0$, which depend only on $\Dc$,
and not on $T_0$, $\eps$ or $d_1$, such that, if $T_0\ge C_{\rm sh}/d_1$ and $C_{\rm
sh}\,\eps\,(1+\log T_0)\le d_1$, then for all $t\ge T_0$,
\begin{equation}
\|g(t)-g_{\rm sc}\|_{L^1(\R^6)}\ge c_{\rm sh}\,d_1^{\,2-3/p}\,\|f_0\|_{L^1(\R^6)}\,\big(1+\eps\log(t/T_0)\big)^{-(p-3)/p}
.
\end{equation}
As in Theorem~\ref{thm:main}(iv), we can enlarge $T_0$ to satisfy the first condition, and then
$T_0$, $\eps_0$ and the constant $C$ of~(i) and~(ii) depend also on $d_1$.
\end{enumerate}
\end{theo}

The rate in~(i) is algebraic in $t$, with an exponent proportional to the collision rate of the
leaf, which is of order $\eps\jap\xi^{-p}$. The rate in~(ii) is logarithmic: at time $t$ only the
leaves with $\eps\jap\xi^{-p}\log(t/T_0)\gtrsim1$ have relaxed, and, once $\eps\log(t/T_0)\ge1$, the
others carry a fraction of the mass of order $(\eps\log(t/T_0))^{-(p-3)/p}$. By~(iii), for $T_0$ and
$\eps$ as there, this rate cannot be improved if $d_1>0$: a leaf with $s_\xi(t)\le c\,d_1$ has
changed little since the time $T_0$, when it was close to $\rho_g F_0(\cdot,\xi)$, so its
$L^1$-distance from the corresponding leaf of $g_{\rm sc}$ is still at least $\frac14d_1$ times its
mass. If $d_1=0$, in particular if $F_0$ is a Maxwellian, we do not know whether the convergence is
faster.

For $\gamma<-2$ the clock of a leaf is $s_\xi(t)=\lambda_{\rm ref}(\xi)\big(\tau(t)-\tau(T_0)\big)$,
with $\tau(t)=(t^{|\gamma|-2}-1)/(|\gamma|-2)$, and the dissipation norm is weaker than the energy
norm at large $|z|$. The convergence still has a rate.

\begin{theo}[Rates for $\gamma<-2$]\label{thm:rate3}
Let $\gamma\in[-3,-2)$, under the assumptions of Theorem~\ref{thm:main} or of Theorem~\ref{thm:t1}.
Let $S_b\in(S_*,S_*+1]$ be the constant of Subsection~\ref{ssec:cells}, and for every $\xi$ let
$t^b_\xi$ be the time with $s_\xi(t^b_\xi)=S_b+1$. Put
\begin{equation}
\varrho_\xi(t):=\frac{t^{1+\gamma}}{\lambda_{\rm
ref}(\xi)},\qquad a_\gamma:=\frac2{14+|\gamma|},\qquad
k_\gamma:=12+|\gamma| .
\end{equation}
There are constants $C$ and $c>0$, which depend only on $\gamma$, $\Dc$ and $T_0$, and not on
$\eps$, such that:
\begin{enumerate}[label=(\roman*)]
\item \emph{Relaxation of the leaves.} For every $\xi$ and every $t\ge t^b_\xi$, with $x=t\xi$ and
$(\rho,y_0,\theta)=p_g(t,\xi)$,
\begin{multline}
\Big\|\frac{f(t,x,(x-\cdot)/t)-M_{p_g(t,\xi)}}{\rho}\Big\|_{L^2(e^{q_0\theta_-|\cdot-y_0|^2/(8\theta)})}\\
\le
C\exp\Big(-c\,\big(s_\xi(t)-S_b-1\big)^{a_\gamma}\Big)+C\,\varrho_\xi(t)\,\Big(\log\Big(e+\frac1{\varrho_\xi(t)}\Big)\Big)^{k_\gamma/2}
.
\end{multline}
Here $\varrho_\xi(t)\le\big((|\gamma|-2)\,t\big)^{-1}$ for $t\ge t^b_\xi$.
\item \emph{Scattering.} With $g_{\rm sc}$ as in Theorem~\ref{thm:scatter} and $p\ge4$ the exponent
of (D3), for all $t\ge T_0$,
\begin{equation}
\|g(t)-g_{\rm sc}\|_{L^1(\R^6)}\le C\,\|f_0\|_{L^1(\R^6)}\,\Big(1+\eps\big(\tau(t)-\tau(T_0)\big)\Big)^{-(p-3)/p}
.
\end{equation}
\item \emph{The rate in~(ii) is sharp.} Suppose in addition that the constant $d_1$ of
Theorem~\ref{thm:scatter}(iv) is positive (this holds, in particular, for non-thermal product data).
Then $d_1\in(0,1]$. There are constants $C_{\rm sh}$ and $c_{\rm sh}>0$, which depend only on
$\gamma$ and $\Dc$, and not on $T_0$, $\eps$ or $d_1$, such that, if $T_0\ge C_{\rm sh}/d_1$ and
$C_{\rm sh}\,\eps\,T_0^{|\gamma|-2}\le d_1$, then for all $t\ge T_0$,
\begin{equation}
\|g(t)-g_{\rm sc}\|_{L^1(\R^6)}\ge c_{\rm sh}\,d_1^{\,2-3/p}\,\|f_0\|_{L^1(\R^6)}\,\Big(1+\eps\big(\tau(t)-\tau(T_0)\big)\Big)^{-(p-3)/p}
.
\end{equation}
As in Theorem~\ref{thm:rate2}(iii), we can enlarge $T_0$ to satisfy the first condition, and then
$T_0$, $\eps_0$ and the constants of~(i) and~(ii) depend also on $d_1$.
\end{enumerate}
\end{theo}

Since $\tau(t)-\tau(T_0)=(t^{|\gamma|-2}-T_0^{|\gamma|-2})/(|\gamma|-2)$, the rate in~(ii) is
algebraic. For large $t$ it is of order $(\eps t^{|\gamma|-2})^{-(p-3)/p}$, and in the Coulomb case
$\gamma=-3$ of order $(\eps t)^{-(p-3)/p}$. With the clock in place of the logarithm, (ii) and~(iii)
are the same statements as in Theorem~\ref{thm:rate2}, and the mechanism is the same. At time $t$
the leaves with $\eps\jap\xi^{-p}(\tau(t)-\tau(T_0))\lesssim1$ have not relaxed, and they carry a
fraction of the mass of order $(\eps(\tau(t)-\tau(T_0)))^{-(p-3)/p}$.

In~(i) the first term is the decay of the deviation that the leaf has when it leaves the band. It is
a stretched exponential in the clock of the leaf, because the dissipation controls the energy only
on balls, with a constant that grows like a power of the radius, while the Gaussian weight controls
the rest. We do not expect the exponent $a_\gamma$ to be optimal. The local part of the energy alone
would lead to $2/|\gamma|$, the exponent of the linearized decay (S3), and the difference comes from
the nonlocal part of the norm, which we bound in a polynomial weight.

The second term is the response to the leafwise transport, which forces the micro part at the rate
$t^{-2}$ while the collisions damp it at the rate $\kappa\lambda_\xi$. The ratio of these rates is
$\varrho_\xi(t)$, up to a constant, and for large $t$ it dominates. Remark~\ref{rem:rate} explains
this heuristically.

In the proof of~(i) we use a second level-$0$ energy, in the weight $\hat m_\xi^{1/2}$, and the two
additional requirements~\ref{rq:Rdc} and~\ref{rq:Rrate} on the constants.

\section{Geometry of leaves}\label{sec:leafgeom}

In this section we derive the equation of a single leaf and explain why energies of a single leaf in
the norms of \cite{CM} do not close.

\subsection{The leaf identity}\label{ssec:leafid}

Recall from \eqref{eq:leafdef} that the leaf through $x$ at time $t$ is the function
$g^{(x)}(z)=g(t,z,(x-z)/t)$, the restriction of the profile $g$ to the set of points $(y,v)$ with
$y+tv=x$. The leaf identity says that the collision operator acts on each leaf as a homogeneous
Landau operator.

\begin{prop}[Leaf identity]\label{prop:leaf}
Let $f$ be a smooth solution of \eqref{eq:landau}, and let $g(t,y,v)=f(t,y+tv,v)$ be its profile.
For all $t>0$,
\begin{equation}\label{eq:leafid}
\partial_tg(t,y,v)=\kappa(t)\,\cQ_t[g](y,v),\qquad\cQ_t[g](y,v):=Q\big(g^{(x)},g^{(x)}\big)(y),\quad
x=y+tv .
\end{equation}
\end{prop}

\begin{proof}
\emph{Step 1: change of variables.} Fix $t$ and $x$, and let $F(w):=f(t,x,w)$. Then
$F(w)=g^{(x)}(x-tw)$. We substitute $z=x-tw$ in the flux \eqref{eq:Bflux}, written in the variable
$w$. Then we have
\begin{equation}\label{eq:leafid-sub}
\dd w=t^{-3}\dd z,\qquad v-w=\frac{z-y}t\ \ (y=x-tv),\qquad a(v-w)=t^{-(2+\gamma)}a(y-z),
\end{equation}
where the last identity is the homogeneity of $a$. Moreover $\nabla F(w)=-t\,(\nabla
g^{(x)})(x-tw)$.

\emph{Step 2: the flux and the operator.} Inserting \eqref{eq:leafid-sub} in the flux, we obtain
\begin{equation}\label{eq:leafid-flux}
B(F,F)(v)=-t^{-(4+\gamma)}B\big(g^{(x)},g^{(x)}\big)(x-tv).
\end{equation}
Since $\nabla_v=-t\nabla_y$ on functions of $y=x-tv$, it follows from \eqref{eq:leafid-flux} that
$Q(F,F)(v)=t^{-(3+\gamma)}Q(g^{(x)},g^{(x)})(y)$.

\emph{Step 3.} By the definition of the profile, we have
$\partial_tg(t,y,v)=(\partial_tf+v\cdot\nabla_xf)(t,y+tv,v)=Q(f,f)(t,y+tv,v)$. By Step~2 this is
$t^{-(3+\gamma)}Q(g^{(x)},g^{(x)})(y)$, and $t^{-(3+\gamma)}=\kappa(t)$ by \eqref{eq:kappa}.
\end{proof}

Therefore $\cQ_t$ is the homogeneous Landau operator in the variable $y$, acting on the restriction
of $g$ to a leaf. The identity \eqref{eq:leafid} holds for every $\gamma$ and every $t>0$, and it
needs no smallness. It explains the clock $\kappa$ of Subsection~\ref{ssec:critical}: in
$\kappa(t)=t^{-3}\cdot t^{-(\gamma+2)}$, the factor $t^{-3}$ is the dispersion of the density, and
the factor $t^{-(\gamma+2)}$ comes from the homogeneity of the kernel.

\subsection{Calculus on leaves}\label{ssec:leafcalc}

We now record how derivatives in the profile variables $(y,v)$ act on leaves. Along a ray $x=t\xi$
we write $D_t$ for the derivative in $t$ at fixed $(\xi,z)$.

\begin{lem}[Leaf calculus]\label{lem:leafgeom}
Let $h=h(t,y,v)$ be smooth, and write $H(t,x,z)=h^{(x)}(z)$. Then:
\begin{enumerate}[label=(\roman*)]
\item $(\cQ_t[g])^{(x)}=Q(g^{(x)},g^{(x)})$;
\item $(\partial_vh)^{(x)}=t\nabla_xH=\nabla_\xi H$ and $(\partial_yh)^{(x)}=\nabla_zH+\nabla_xH$;
\item at fixed $(x,z)$, $\partial_tH=(\partial_th)^{(x)}-t^{-2}(x-z)\cdot(\partial_vh)^{(x)}$;
\item along rays, at fixed $(\xi,z)$,
\begin{equation}\label{eq:leaftransport}
D_tH=(\partial_th)^{(x)}+t^{-2}z\cdot\nabla_\xi H .
\end{equation}
\end{enumerate}
\end{lem}

\begin{proof}
(i) is \eqref{eq:leafid}, evaluated on the leaf. (ii) and (iii) follow from the chain rule applied
to $H(t,x,z)=h(t,z,(x-z)/t)$. For (iv), we differentiate $H(t,t\xi,z)$ in $t$ and get
\begin{equation}
D_tH=\partial_tH+\xi\cdot\nabla_xH=(\partial_th)^{(x)}-t^{-2}(x-z)\cdot(\partial_vh)^{(x)}+t^{-2}x\cdot(\partial_vh)^{(x)},
\end{equation}
by (ii), (iii) and $\xi=x/t$. Since $(\partial_vh)^{(x)}=\nabla_\xi H$, we obtain
\eqref{eq:leaftransport}.
\end{proof}

The term $t^{-2}z\cdot\nabla_\xi$ in~(iv) is the \emph{leafwise transport}. It
appears because the leaf through a fixed point $(y,v)$ moves: the particles that are at $x$ at time
$t$ are spread, at time $t+\dd t$, over leaves whose asymptotic velocities differ by $t^{-2}z\,\dd
t$. The transport has the integrable rate $t^{-2}$, but it loses one derivative in $\xi$, that is,
one derivative in $v$. This loss is the origin of the first obstruction of
Subsection~\ref{ssec:obstructions}.

\subsection{The fibered equation and the remainder}\label{ssec:remainder}

In this subsection we replace the leaf through $(y,v)$ by the fiber $\{v=\mathrm{const}\}$, and
estimate the error. We write
\begin{equation}\label{eq:fibered}
\partial_tg=\kappa\big(\Qf[g]+R\big),\qquad\Qf[g](y,v)=Q\big(g(t,\cdot,v),g(t,\cdot,v)\big)(y).
\end{equation}
The fibered operator $\Qf$ acts on each fiber separately, and the remainder $R$ measures the
difference between the leaf and the fiber through $(y,v)$. To interpolate between them, we put, for
$\vartheta\in[0,1]$,
\begin{equation}
G_\vartheta(z)=g\big(t,z,v-\vartheta\tfrac{z-y}{t}\big),\qquad
H_\vartheta(z)=(z-y)\cdot(\partial_vg)\big(t,z,v-\vartheta\tfrac{z-y}{t}\big),\qquad\tn_\vartheta=\partial_y-\tfrac\vartheta
t\partial_v ,
\end{equation}
so that $G_1=g^{(x)}$, $G_0=g(t,\cdot,v)$, $\partial_\vartheta G_\vartheta=-t^{-1}H_\vartheta$ and
$H_\vartheta(y)=0$.

\begin{prop}[The remainder]\label{prop:Rformula}
We have
\begin{equation}
R=-\frac1t\int_0^1\Big[2\,\ab[G_\vartheta]_{ij}(y)\,(\tn_\vartheta)_i\partial_{v_j}g+\ab[H_\vartheta]_{ij}(y)\,(\tn_\vartheta)_i(\tn_\vartheta)_jg-\cb[H_\vartheta](y)\,g\Big]\dd\vartheta
,
\end{equation}
with all unbarred quantities evaluated at $(t,y,v)$. For $\gamma=-3$ the last term vanishes.
\end{prop}

\begin{proof}
Since $G_1=g^{(x)}$ and $G_0=g(t,\cdot,v)$, and $\partial_\vartheta G_\vartheta=-t^{-1}H_\vartheta$,
we have
\begin{equation}\label{eq:R-interp}
R=\int_0^1\partial_\vartheta Q(G_\vartheta,G_\vartheta)(y)\dd\vartheta=-\frac1t\int_0^1\big[Q(H_\vartheta,G_\vartheta)+Q(G_\vartheta,H_\vartheta)\big](y)\dd\vartheta
.
\end{equation}
We write $Q(F,G)=\ab[F]:\nabla^2G-\cb[F]G$. In $Q(G_\vartheta,H_\vartheta)(y)$ the term
$\cb[G_\vartheta]H_\vartheta(y)$ vanishes, because $H_\vartheta(y)=0$. By the chain rule along
$z\mapsto(z,v-\vartheta(z-y)/t)$, we have
\begin{equation}
\nabla^2G_\vartheta(y)=\tn_\vartheta\otimes\tn_\vartheta
g,\qquad\partial_{ij}H_\vartheta(y)=(\tn_\vartheta)_j\partial_{v_i}g+(\tn_\vartheta)_i\partial_{v_j}g
.
\end{equation}
Inserting these in \eqref{eq:R-interp}, and using that $\ab$ is symmetric, we obtain the formula.
For $\gamma=-3$, $\cb[H_\vartheta](y)=-8\pi H_\vartheta(y)=0$.
\end{proof}

For $N,M\ge0$ we set $\Nn_k(t)=\sum_{|\alpha|+|\beta|\le k}\sup_{y,v}\jap y^N\jap
v^{M+2}|\partial_y^\alpha\partial_v^\beta g|$.

\begin{prop}[Size of the remainder]\label{prop:Rest}
Let $t\ge1$ and $N\ge6$. Then $\jap y^{N-(3+\gamma)}\jap v^{M+2}|R(t,y,v)|\le Ct^{-1}\Nn_2(t)^2$,
uniformly in $\gamma\in[-3,-2]$. The derivatives $\partial^\alpha_y\partial^\beta_vR$ satisfy the
same bound with $\Nn_{k+2}^2$, $k=|\alpha|+|\beta|$.
\end{prop}

\begin{proof}
\emph{Step 1: the coefficients.} The bound $|\partial^\alpha g(z,w)|\le\Nn_{|\alpha|}\jap z^{-N}$
holds at every velocity $w$, so it needs no weight in the sheared velocity $v-\vartheta(z-y)/t$. For
$p\in[-2,1]$ and $N>3+\max(p,0)$, we have
\begin{equation}\label{eq:Rest-conv}
\int|y-z|^p\jap z^{-N}\dd z\le C\jap y^p .
\end{equation}
Applying \eqref{eq:Rest-conv} with $p=2+\gamma$, $3+\gamma$ and $1+\gamma$ to the kernels of $\ab$
and $\cb$, we get
\begin{equation}\label{eq:Rest-coef}
|\ab[G_\vartheta](y)|\lesssim\Nn_0\jap y^{2+\gamma},\qquad|\ab[H_\vartheta](y)|\lesssim\Nn_1\jap
y^{3+\gamma},\qquad|\cb[H_\vartheta](y)|\lesssim(3+\gamma)\Nn_1\jap y^{1+\gamma}.
\end{equation}

\emph{Step 2: the other factors, and conclusion.} In Proposition~\ref{prop:Rformula} the factors
that multiply the coefficients are bounded by $(1+t^{-1})^2\Nn_2\jap y^{-N}\jap v^{-M-2}$. Combining
this with \eqref{eq:Rest-coef} and $t\ge1$, we obtain the bound on $R$.

\emph{Step 3: derivatives.} The fields $\partial_x$ and $t\nabla_x+\nabla_v$ are derivations of
$Q(f,f)$. Thus $\partial^\alpha_y\partial^\beta_vR$ is a Leibniz sum of polarized remainders, to
which Steps~1 and~2 apply.
\end{proof}

For comparison, for $\gamma\in(-2,0)$ the same computation yields
$\partial_tg=\kappa\Qf[g]+O(\eps^2t^{-(4+\gamma)})$, hence $g-g_\infty=O(\eps^2t^{-(2+\gamma)})$.
For $\gamma\in(-2,-1]$ this is the rate of \cite[Thm.~1.3]{Luk}. For $\gamma\le-2$ the factor
$\kappa$ is not integrable, and the fibered approximation is useful only on bounded effective time
intervals. We use it in this way in Section~\ref{sec:fib}.

\subsection{The exact transport law of leaf moments}\label{ssec:moments}

The collision operator conserves mass, momentum and energy on each leaf. Hence only the leafwise
transport changes the leaf moments.
Theorem~\ref{thm:momtransport} below makes this exact. Recall the moment vector
$\mom(z)=(1,z,|z|^2)$ and the leaf moments $m_h(t,x)=\int h^{(x)}\mom\dd z$.

In this subsection $g$ is the profile of a classical solution of \eqref{eq:landau} on an open
interval $I\subset(0,\infty)$, and all derivatives written below are assumed to exist and to be
continuous. For a function $h$ of $(t,y,v)$ we write $H(t,\xi,z):=h^{(t\xi)}(z)$, and
$G(t,\xi,z):=g^{(t\xi)}(z)$. For an integer $a\ge0$ we say that $h$ satisfies $(\mathcal D_a)$ if
for all multi-indices $j,k$ with $|j|\le a$ and every compact set $K\subset I\times\R^3$
\begin{equation}
\sup_{(t,\xi)\in K}\,\sup_{z}\,\jap z^7\big|\partial_z^j\partial_\xi^kH(t,\xi,z)\big|<\infty .
\end{equation}
By Lemma~\ref{lem:leafgeom}(ii), we have
$(\partial^\alpha_y\partial^\beta_vh)^{(t\xi)}=(\nabla_z+t^{-1}\nabla_\xi)^\alpha\nabla_\xi^\beta
H$, so if $h$ satisfies $(\mathcal D_a)$ and $|\alpha|\le a$, then
$\partial^\alpha_y\partial^\beta_vh$ satisfies $(\mathcal D_{a-|\alpha|})$. For a polynomial $q$ of
degree $\le3$, $\jap z^{-7}|q|\le C\jap z^{-4}$ is integrable on $\R^3$ and $o(|z|^{-2})$ at
infinity. This justifies the integrations by parts in $z$ against such $q$, and the differentiations
in $\xi$ under the integral sign, below. Only finitely many $z$-derivatives enter, and
in Remark~\ref{rem:Dsol} we verify the condition for the solutions constructed in this paper.

\begin{lem}[Leaf moments of derivatives]\label{lem:commmom}
Let $\alpha,\beta$ be multi-indices, and let $h$ satisfy $(\mathcal D_{|\alpha|})$. Let $q$ be a
polynomial in $z$ of degree $\le3$. Then, for each $t\in I$, as functions of $\xi$ (with $x=t\xi$),
\begin{equation}
\int(\partial^\alpha_y\partial^\beta_vh)^{(x)}\,q\dd
z=\sum_{\alpha'\le\alpha}\binom\alpha{\alpha'}(-1)^{|\alpha'|}\,t^{-|\alpha-\alpha'|}\,\partial_\xi^{\alpha-\alpha'+\beta}\int
h^{(x)}\,\partial^{\alpha'}q\dd z .
\end{equation}
\end{lem}

\begin{proof}
We write $H(\xi,z)=h^{(t\xi)}(z)$ at fixed $t$. By Lemma~\ref{lem:leafgeom}(ii), applied to $h$ and
to its derivatives, we have $(\partial_vk)^{(x)}=\nabla_\xi K$ and
$(\partial_yk)^{(x)}=(\nabla_z+t^{-1}\nabla_\xi)K$ for every smooth $k$, with
$K(\xi,z)=k^{(t\xi)}(z)$. Since $\nabla_z$ and $\nabla_\xi$ commute, we get
\begin{equation}
(\partial^\alpha_y\partial^\beta_vh)^{(x)}=(\nabla_z+t^{-1}\nabla_\xi)^\alpha\nabla_\xi^\beta
H=\sum_{\alpha'\le\alpha}\binom\alpha{\alpha'}t^{-|\alpha-\alpha'|}\,\partial_z^{\alpha'}\partial_\xi^{\alpha-\alpha'+\beta}H
.
\end{equation}
We integrate against $q$, move $\partial_z^{\alpha'}$ onto $q$ by parts, and take
$\partial_\xi^{\alpha-\alpha'+\beta}$ out of the integral. Both steps are justified by $(
\mathcal D_{|\alpha|})$, since $|\alpha'|\le|\alpha|$.
\end{proof}

\begin{theo}[Transport of leaf moments]\label{thm:momtransport}
Let $g$ satisfy $(\mathcal D_a)$ for some $a\ge2$. For all multi-indices $\alpha,\beta$ with
$|\alpha|\le a-2$, the function $h=\partial^\alpha_y\partial^\beta_vg$ satisfies, along each ray
$x=t\xi$, $t\in I$,
\begin{equation}\label{eq:momlawh}
D_t\,m_h(t,t\xi)=t^{-2}\,\nabla_\xi\cdot\int z\otimes\mom(z)\,h^{(t\xi)}(z)\dd z .
\end{equation}
In particular,
\begin{equation}\label{eq:momlaw}
D_t\,m_g(t,t\xi)=t^{-2}\,\nabla_\xi\cdot\int z\otimes\mom(z)\,g^{(t\xi)}(z)\dd z .
\end{equation}
There are no lower-order terms.
\end{theo}

\begin{proof}
\emph{Step 1: the transport identity.} By \eqref{eq:leaftransport}, applied to $h$, we have
\begin{equation}\label{eq:mom-step1}
D_tm_h=\int(\partial_th)^{(x)}\mom\dd z+t^{-2}\int z\cdot\nabla_\xi h^{(x)}\,\mom\dd
z=\int(\partial_th)^{(x)}\mom\dd z+t^{-2}\nabla_\xi\cdot\int z\otimes\mom\,h^{(x)}\dd z .
\end{equation}
We only need to show that the first integral on the right vanishes.

\emph{Step 2: the collision term vanishes.} By \eqref{eq:leafid}, we have
$\partial_th=\kappa(t)\,\partial^\alpha_y\partial^\beta_v\cQ_t[g]$. Lemma~\ref{lem:commmom}, applied
with $\cQ_t[g]$ in place of $h$ (we check its hypothesis in Step~3), expresses
$\int(\partial^\alpha_y\partial^\beta_v\cQ_t[g])^{(x)}\mom\dd z$ through $\xi$-derivatives of the
integrals
\begin{equation}\label{eq:mom-inv}
\int(\cQ_t[g])^{(x)}\psi\dd z=\int Q(g^{(x)},g^{(x)})\psi\dd
z,\qquad\psi\in\mathrm{span}\{1,z_i,|z|^2\},
\end{equation}
where the equality is Lemma~\ref{lem:leafgeom}(i), and $\psi$ runs over the derivatives of the
components of $\mom$. Each integral in \eqref{eq:mom-inv} vanishes for every $\xi$, since $\psi$ is
a collision invariant. Hence the first integral in \eqref{eq:mom-step1} vanishes, and
\eqref{eq:momlawh} follows.

\emph{Step 3: the hypothesis of the lemma.} We check that $\cQ_t[g]$ satisfies
$(\mathcal D_{a-2})$. Then $(\mathcal D_{|\alpha|})$ holds, since $|\alpha|\le a-2$. At fixed
$(t,\xi)$ the leaf restrictions of $\cQ_t[g]$ are $Q(G,G)=\ab[G]:\nabla_z^2G-\cb[G]G$. The kernel of
$\ab$ is bounded by $|z|^{\gamma+2}$, with $\gamma+2\in[-1,0]$, and
$\cb[G]=-2(\gamma+3)|\cdot|^\gamma*G$ for $\gamma>-3$, $\cb[G]=-8\pi G$ for $\gamma=-3$. Since
$\partial^j_z\partial^k_\xi G=O(\jap z^{-7})$ for $|j|\le a$, locally uniformly in $(t,\xi)$, the
coefficients $\ab[G]$, $\cb[G]$ and their derivatives $\partial^j_z\partial^k_\xi$, $|j|\le a$, are
bounded, locally uniformly. By the Leibniz rule, we have $\partial^j_z\partial^k_\xi Q(G,G)=O(\jap
z^{-7})$ for $|j|\le a-2$.

\emph{Step 4: differentiation under the integral of Step~1.} The operator $D_t$
commutes with $\nabla_z$ and $\nabla_\xi$, and, as in Steps~1 and~2, we have
\begin{equation}
D_tG=\kappa(t)Q(G,G)+t^{-2}z\cdot\nabla_\xi G .
\end{equation}
Therefore $D_th^{(t\xi)}$ is a sum of terms $t^{-i}\partial^j_z\partial^k_\xi D_tG$ and
$t^{-i-1}\partial^j_z\partial^k_\xi G$ with $|j|\le|\alpha|$. By Step~3 these are $O(\jap z^{-6})$,
locally uniformly in $(t,\xi)$, and $\jap z^{-6}|\mom|\le C\jap z^{-4}$ is integrable. Hence $m_h$
is differentiable along the rays, with $D_tm_h=\int D_th^{(t\xi)}\mom\dd z$.
\end{proof}

The identity \eqref{eq:momlaw} is the system of local conservation laws, written in the variables of
the leaves. It is a compressible Euler-type system for the five fields $m_g$, with fluxes given by
third moments and with the prefactor $t^{-2}$. The prefactor reflects free expansion: $\nabla_x$
acts on the scale $t$, and the local temperature is $O(t^{-2})$. Since $\int^\infty t^{-2}\dd
t<\infty$, the hydrodynamic evolution has finite total duration. Theorem~\ref{thm:Meuler} closes
this system up to small errors.

\begin{remark}[Verification of $(\mathcal D_a)$]\label{rem:Dsol}
Let $g$ be a classical solution in the cell norms on $[T_0,T')$, as in Corollary~\ref{cor:CUmain}.
The solution of Theorem~\ref{thm:main} is one, with $T'=\infty$. Then $g$ satisfies $(\mathcal D_a)$
on $I=(T_0,T')$ for every $a\le N_y-2$.

\emph{The micro part.} We fix $T<T'$ and write $g=\mathcal M[g]+\hat r$. On $[T_0,T]$ the cell
energies are bounded and (F3) holds, so the static parts (i) and (ii) of Lemma~\ref{lem:BTsob} hold
with a bound $B=B(T)$ in place of $B_\xi$. The proof of~(ii) applies to $|a'|\le N_y-2$ and to every
$b$, with a constant that depends on $|b|$, since~(i) holds for every $|\beta|$.

By Lemma~\ref{lem:leafgeom}(ii), we have $\nabla_\xi H=(\partial_vh)^{(x)}$ and
$\nabla_zH=(\partial_yh)^{(x)}-t^{-1}(\partial_vh)^{(x)}$. Thus $\partial^j_z\partial^k_\xi(\hat
r^{(t\xi)})$ is a combination of restrictions $t^{-i}(\partial^{a'}_y\partial^{b}_v\hat r)^{(t\xi)}$
with $|a'|\le|j|$. Each of them is bounded by $C(T,j,k)e^{-(q/4)|\hat z|^2}$, which is $O(\jap
z^{-7})$ uniformly in $\xi$, since $\theta_g\le\theta_++1$ and $|y_{0,g}|\le L+1$ by (F3), that is,
by (U1).

\emph{The Maxwellian part.} The restriction $\mathcal M[g]^{(t\xi)}=M_{p_g(t,\xi)}$ has Gaussian
derivatives. Indeed, $\rho_g$ is continuous, hence bounded on compact sets, and the derivatives of
every order of $\nabla_\xi\ell_g$, $y_{0,g}$ and $\vartheta_g$ are bounded on $[T_0,T]$, by
Lemma~\ref{lem:Mref} and the bound $\sup_Q\mathfrak m_{Q^*}\le\delta_U^2$ in the definition, that
is, (U1), through the Cauchy estimates of Lemma~\ref{lem:TWcalc}(i). (F3) states these bounds for
$k\le J$, and the same argument applies to every $k$, with constants that depend on $k$.

To avoid circularity we use only these qualitative local bounds, and not the mass part of (F3),
which is derived from Theorem~\ref{thm:Meuler} and hence from Theorem~\ref{thm:momtransport}.
\end{remark}

\subsection{Two obstructions}\label{ssec:obstructions}

The two obstructions below explain why the solution has to be controlled in norms that are analytic
in $v$ and Gaussian in $z$. Readers interested only in the proof may skip to
Lemma~\ref{lem:analytic} and Proposition~\ref{prop:Gthreshold}, which are the only results of this
subsection used later.

\subsubsection{Derivative loss through the leafwise transport}\label{sssec:obs1}

Energies in the profile variables lose nothing under transport, since the profile equation has no
transport term. Indeed, for any weight $w>0$ and any solution of $\partial_th=\kappa(
\mathcal L_t[g]h+S)$, we have
\begin{equation}\label{eq:prof-energy}
\frac{d}{dt}\iint w^2h^2\dd y\dd v=\iint(\partial_tw^2)h^2\dd y\dd v+2\kappa\iint w^2h(
\mathcal L_t[g]h+S)\dd y\dd v .
\end{equation}
The coercivity of the linearized operator, however, is a property of each leaf separately. The
natural near-equilibrium norms of \cite{CM} are \emph{nonlocal} in $z$: they contain terms of the
form $\int_0^\infty\|S(\tau)\varphi\|^2\dd\tau$, where $S$ is a semigroup.
The next proposition shows that such norms do not commute with the leafwise transport.

\begin{prop}[Loss of nonlocal leaf energies]\label{prop:commloss}
Let $\mathcal N(\varphi)=\iint K(z,z')\varphi(z)\varphi(z')\dd z\dd z'$ be a symmetric quadratic
form, independent of $(t,\xi)$, and let $E(t)=\int\mathcal N\big(H(t,\xi,\cdot)\big)\dd\xi$ with
$H(t,\xi,z)=h^{(t\xi)}(z)$. The contribution of the leafwise transport in \eqref{eq:leaftransport}
to $\frac{d}{dt}E$ is
\begin{equation}\label{eq:commloss}
t^{-2}\int\!\!\iint K(z,z')\,(z'-z)\cdot H(\xi,z)\,\nabla_\xi H(\xi,z')\dd z\dd z'\dd\xi .
\end{equation}
It vanishes for every $h$ if $K(z,z')=k(z)\delta(z-z')$, that is, if $\mathcal N(\varphi)=\int
k\varphi^2$. Conversely, if it vanishes for every $h$, then $K$ is supported on the diagonal.
\end{prop}

Locality of $\mathcal N$ alone is not enough: for $\mathcal N(\varphi)=\int|\nabla\varphi|^2$ the
contribution is $2t^{-2}\iint\nabla_zH\cdot\nabla_\xi H$, which does not vanish in general, because
of the commutator $[\nabla_z,z\cdot\nabla_\xi]=\nabla_\xi$. In general the contribution involves one
derivative $\nabla_\xi H=(\partial_vh)^{(x)}$.

\begin{proof}
By the symmetry of $K$, the contribution is
\begin{equation}\label{eq:commloss-1}
I:=2t^{-2}\int\!\!\iint K(z,z')\,H(\xi,z)\,z'\cdot\nabla_\xi H(\xi,z')\dd z\dd z'\dd\xi .
\end{equation}
Integrating by parts in $\xi$, then exchanging $z$ and $z'$ and using again the symmetry of $K$, we
get
\begin{equation}\label{eq:commloss-2}
I=-2t^{-2}\int\!\!\iint K(z,z')\,z\cdot\nabla_\xi H(\xi,z')\,H(\xi,z)\dd z\dd z'\dd\xi .
\end{equation}
The average of \eqref{eq:commloss-1} and \eqref{eq:commloss-2} is \eqref{eq:commloss}.
\end{proof}

Thus the loss is one derivative, at the integrable rate $t^{-2}$. Its coefficient does not depend on
the order of differentiation. By Lemma~\ref{lem:localnoloss}, local energies avoid the loss. Among
them, the one with the best structure has on each leaf the weight $M^{-1}$, where $M$ is the leaf
Maxwellian. Indeed, the linearized operator around $M$ is symmetric in $L^2(M^{-1})$, so there is no
macro--micro cross term. The next two results describe the restrictions on such weights: the first
treats static weights, and the second weights that follow the leaves. Recall that the restriction of
$h(y,v)$ to the leaf through $x$ is $h^{(x)}(z)=h(z,(x-z)/t)$, and that $M_p$,
$p=(\rho,y_0,\theta)$, is the Maxwellian with these
parameters. We call a function of $z$ \emph{Gaussian} if its logarithm is a
polynomial of degree $\le2$ whose quadratic part is negative definite. By a \emph{Gaussian solution
of free transport} we mean a solution $f(t,x,v)=F(x-tv,v)$ with $\log F$ a polynomial of degree
$\le2$ in $(y,v)$. Its quadratic part in $(y,v)$ need not be definite (for $\rho\equiv1$,
$y_0\equiv0$, $\theta\equiv1$ in the next proposition it is $-\frac12|y|^2$).

\begin{prop}[Static weights with Gaussian leaf restrictions]\label{prop:travelling}
Let $\rho,\theta>0$ and $y_0$ be continuous functions on $\R^3$, with $\rho$ bounded, and let
$\mathcal M(y,v):=M_{(\rho(v),y_0(v),\theta(v))}(y)$. Suppose that at two distinct times $t_1,t_2>0$
the leaf restrictions $\mathcal M^{(x)}$ are Gaussian for every $x\in\R^3$. Then $\theta$ is
constant, $y_0$ is affine and $\log\rho$ is a polynomial of degree $\le2$. In particular
$\log\mathcal M$ is a polynomial of degree $\le2$ in $(y,v)$, that is, $\mathcal M(x-tv,v)$ is a
Gaussian solution of free transport. Every global Maxwellian of finite mass, in particular every
traveling Maxwellian of \cite{CL}, is of this form, since its logarithm is a combination of $1$,
$v$, $|v|^2$, $y$, $y\cdot v$, $|y|^2$ and the components of $y\times v$, with a negative
coefficient of $|y|^2$. Conversely, if $\theta$ is constant, $y_0$ affine and $\log\rho$ a
polynomial of degree $\le2$, then $\mathcal M^{(x)}$ is Gaussian for every $x$ and every $t$ large
enough.
\end{prop}

\begin{proof}
Fix $t\in\{t_1,t_2\}$. Since $z=x-tv$ is an affine bijection between $z$ and $v$, the hypothesis
means that
\begin{equation}\label{eq:trav-F}
F_x(v):=\log\mathcal M(x-tv,v)=-\tfrac q2|x|^2+q\,x\cdot\big(tv+y_0(v)\big)+c(v)-\tfrac
q2\big|tv+y_0(v)\big|^2
\end{equation}
is a polynomial of degree $\le2$ in $v$ for every $x$, where $q=1/\theta$ and
$c=\log\rho-\frac32\log(2\pi\theta)$.

\emph{Step 1: the coefficients are polynomials.} As a function of $x$, $F_x(v)$ is a polynomial of
degree $\le2$ with coefficients $q$, $q(tv+y_0)$ and $c-\frac q2|tv+y_0|^2$. Taking finite
differences in $x$, we see that each of these coefficients is a polynomial of degree $\le2$ in $v$,
at $t=t_1$ and at $t=t_2$.

\emph{Step 2: $\theta$ is constant.} The difference of $q(tv+y_0)$ at $t_1$ and $t_2$ is
$(t_1-t_2)qv$, so $qv_i$ is a polynomial of degree $\le2$ for each $i$. Since $q$ is a polynomial,
of degree $d$ say, we have $d+1\le2$. An affine function that is positive on $\R^3$ is constant, so
$q$ is constant.

\emph{Step 3: $y_0$ is affine.} Now $q(tv+y_0)$ shows that $y_0$ is a polynomial of degree $\le2$.
The difference of $c-\frac q2|tv+y_0|^2$ at $t_1$ and $t_2$ is $-\frac
q2\big((t_1^2-t_2^2)|v|^2+2(t_1-t_2)v\cdot y_0\big)$, so $v\cdot y_0(v)$ has degree $\le2$, and then
$\log\rho=\frac q2|y_0|^2+P_2$ with a polynomial $P_2$ of degree $\le2$. If $y_0$ had a nonzero
quadratic part $P$, then along a ray $v=s\omega$ with $P(\omega)\ne0$ we would have
$|y_0(s\omega)|^2\ge cs^4$ for $s$ large, and $\rho$ would be unbounded. Hence $y_0$ is affine, and
$\log\rho$ is a polynomial of degree $\le2$.

\emph{Step 4: the converse.} We write $y_0(v)=Av+b$. On the leaf through $x$, $v=(x-z)/t$, so
$\log\mathcal M^{(x)}$ is a polynomial of degree $\le2$ in $z$ with quadratic part $-\frac
q2|(I+t^{-1}A)z|^2+t^{-2}R(z)$, where $R$ is the quadratic part of $\log\rho$. It is negative
definite for $t$ large.
\end{proof}

\begin{lem}[Leaf-adapted Maxwellian weights]\label{lem:leafadapted}
Let $\bar p(\xi)=(\bar\rho,\bar y_0,\bar\theta)(\xi)$ be $C^1$, with $\bar\rho,\bar\theta>0$, and
let $w^2(t,y,v):=M_{\bar p(\xi)}(y)^{-1}$ with $\xi=v+y/t$. On the leaf through $x$, $\xi=x/t$ is
constant, so the restriction of $w^2$ is exactly $M_{\bar p(x/t)}^{-1}$. At fixed $(y,v)$,
\begin{equation}\label{eq:dtlogw}
\partial_t\log w^2=t^{-2}\Big[y\cdot\nabla\log\bar\rho-\frac32\,\frac{y\cdot\nabla\bar\theta}{\bar\theta}+\frac{\big((y\cdot\nabla)\bar
y_0\big)\cdot(y-\bar y_0)}{\bar\theta}+\frac{|y-\bar
y_0|^2\,(y\cdot\nabla\bar\theta)}{2\bar\theta^2}\Big],
\end{equation}
with all coefficients evaluated at $\xi$.
\end{lem}

\begin{proof}
At fixed $(y,v)$, we have $\partial_t\xi=-y/t^2$. It remains to differentiate $-\log M_{\bar
p(\xi)}(y)=-\log\bar\rho+\frac32\log(2\pi\bar\theta)+|y-\bar y_0|^2/(2\bar\theta)$ in $\xi$.
\end{proof}

We compare the three kinds of weights. Consider first the leaf-adapted weights. By
Lemma~\ref{lem:localnoloss}, the energy $\iint w^2h^2$ with
the weight of Lemma~\ref{lem:leafadapted} loses no derivative, and on each leaf its weight is
exactly the inverse of a Maxwellian. However, by \eqref{eq:prof-energy} its time derivative contains
$\iint(\partial_tw^2)h^2$. Wherever $\nabla\bar\theta\ne0$, the last term of \eqref{eq:dtlogw} is
cubic in $y$, and it has no sign. Next, a Gaussian weight whose exponent decreases in time, as in
Section~\ref{sec:weights}, produces only a good term $-Ct^{-2}|y|^2$. But it is no longer exactly
the inverse of the leaf Maxwellian. Finally, the static weights $w^2=\mathcal M^{-1}$, with
$\mathcal M$ a Maxwellian in $y$ on every fiber, have $\partial_tw^2=0$. By
Proposition~\ref{prop:travelling}, they have Gaussian leaf restrictions only in the class of
Gaussian solutions of free transport. This class contains the traveling Maxwellians, which are the
setting of \cite{CL}. The proposition says nothing about static weights of other forms. For example,
$\exp(-|y\times v|^2-|y|^2-|v|^2)$ has Gaussian, though not Maxwellian, leaf restrictions at every
time, since $(x-tv)\times v=x\times v$, but its logarithm has degree $4$. We do not know a loss-free
top-order energy with Gaussian structure outside the class of Gaussian solutions of free
transport. In
general we absorb the derivative loss by analyticity, in the following elementary way.

\begin{lem}[Absorption of an integrable one-derivative loss]\label{lem:analytic}
Let $E_n(t)\ge0$, $n\ge0$, be absolutely continuous on $[T_1,\infty)$ and satisfy
\begin{equation}\label{eq:abs-hyp}
\frac{d}{dt}E_n\le a(t)E_n+C_0\,b(t)\sqrt{E_nE_{n+1}},
\end{equation}
with $a,b\ge0$ integrable and $C_0$ independent of $n$. Let $\rad_0>0$ with
$\frac34C_0\int_{T_1}^\infty b\le\rad_0/2$, and set
\begin{equation}\label{eq:abs-lambda}
\rad(t)=\rad_0-\tfrac34C_0\int_{T_1}^tb,\qquad
\mathbb A(t)=\sum_{n\ge0}\Big(\frac{\rad(t)^n}{n!}\Big)^2E_n(t).
\end{equation}
Assume also that $\mathbb A$ is bounded on $[T_1,T]$ for every $T$. Then $\rad\ge\rad_0/2$ and
\begin{equation}\label{eq:abs-concl}
\mathbb A(t)\le\mathbb A(T_1)\exp\big(\int_{T_1}^t(a+C_0b/\rad_0)\big) .
\end{equation}

\end{lem}

\begin{proof}
\emph{Step 1: the truncated energy.} The bound $\rad\ge\rad_0/2$ follows from
\eqref{eq:abs-lambda} and the choice of $\rad_0$. Let $a_n=\rad^n\sqrt{E_n}/n!$ and $
\mathbb A_{\le N}:=\sum_{n\le N}a_n^2$, which is absolutely continuous. By \eqref{eq:abs-hyp}, we
have
\begin{equation}\label{eq:abs-1}
\frac{d}{dt}\mathbb A_{\le N}\le\sum_{n\le N}\Big[2n\frac{\rad'}{\rad}a_n^2+a(t)a_n^2+C_0b\,\frac{n+1}\rad
a_na_{n+1}\Big].
\end{equation}

\emph{Step 2: the loss is absorbed.} We have
\begin{equation}\label{eq:abs-2}
\sum_{n<N}(n+1)a_na_{n+1}\le\tfrac12a_0^2+\tfrac32\sum_{1\le m\le N}ma_m^2 .
\end{equation}
Since $\rad'=-\frac34C_0b$, the coefficient of $ma_m^2$, $1\le m\le N$, in \eqref{eq:abs-1} does
not exceed $(-\frac32C_0b+\frac32C_0b)/\rad=0$. Hence we get
\begin{equation}\label{eq:abs-3}
\frac{d}{dt}\mathbb A_{\le N}\le\Big(a+\frac{C_0b}{2\rad}\Big)\mathbb A_{\le
N}+C_0b\,\frac{N+1}\rad\,a_Na_{N+1}.
\end{equation}

\emph{Step 3: the boundary term.} The integral over $[T_1,T]$ of the last term of \eqref{eq:abs-3}
is bounded by $\frac{2C_0}{\rad_0}(N+1)x_N$, where
\begin{equation}\label{eq:abs-xN}
x_N:=\int_{T_1}^Tb\,(a_N^2+a_{N+1}^2) .
\end{equation}
Since $\sum_Nx_N\le2\int_{T_1}^Tb\,\mathbb A<\infty$, there is a sequence $N_j\to\infty$ with
$(N_j+1)x_{N_j}\to0$. Otherwise $x_N\ge c/(N+1)$ for some $c>0$ and all large $N$, and the series
would diverge.

\emph{Step 4.} Applying Gronwall's inequality to the absolutely continuous functions $
\mathbb A_{\le N_j}$, with $\rad\ge\rad_0/2$, we obtain for $t\in[T_1,T]$
\begin{equation}\label{eq:abs-4}
\mathbb A_{\le N_j}(t)\le\Big(\mathbb A_{\le N_j}(T_1)+\frac{2C_0}{\rad_0}(N_j+1)x_{N_j}\Big)\exp\Big(\int_{T_1}^t\big(a+C_0b/\rad_0\big)\Big).
\end{equation}
Letting $j\to\infty$, we obtain the claim, because $\mathbb A_{\le N}$ increases to $\mathbb A$. We
use no bound on $\sum_nn\,a_n^2$. This completes the proof.
\end{proof}

With $b(t)=t^{-2}$ the total loss of radius is $\frac34C_0/T_1$, which is small for $T_1$ large. The
abstract Cauchy--Kovalevskaya theorem \cite{Ovs65,Ni72,Nis77,Caf90} rests on the same mechanism, in
the form used for Landau damping \cite{MV11} and for boundary layers \cite{SC98}. In the actual
proof, in Section~\ref{sec:cells}, we use the weights $w_n=\rad^n(n+1)^5/n!$ rather than
$\rad^n/n!$, so that we can estimate products.

\subsubsection{The mass--tail coupling}\label{sssec:obs2}

Here we show that the coupling between the mass of a leaf and its tail is
controlled by the dissipation exactly for Gaussian weights. The collision rate
$\lambda(\xi)=\rho_\infty\theta_\infty^{\gamma/2}(\xi)$ varies from leaf to leaf. In the frame
$T_{p_g}$ of the leaf, the leaf solves $\partial_\tau\psi=\lambda(\xi)Q(\psi,\psi)$ by
\eqref{eq:frame}, and a derivative $\partial_\xi$ of this equation produces the term
$(\partial_\xi\lambda)Q(\psi,\psi)$. Linearized around $\mu$, its part that is not controlled by the
dissipation is the \emph{mass-tail form}
\begin{equation}\label{eq:masstail}
\mathfrak b_m(g,h):=\int\bb[\mu]\cdot\nabla(m^2h)\,g\dd z ,
\end{equation}
in which the first factor carries mass and is measured in the energy norm, not in the dissipation
norm. The dissipation norm of $\mathcal L=Q(\mu,\cdot)+Q(\cdot,\mu)$ in the weight $m$ is the norm
of $H^1_*(m)$ in \cite[(3.1)]{CM},
\begin{equation}\label{eq:Ynorm}
\|f\|^2_{Y(m)}:=\big\|\jap z^{(\gamma+\sigma)/2}mf\big\|^2_{L^2}+\big\|\jap
z^{\gamma/2}\tn(mf)\big\|^2_{L^2},\qquad\tn:=P_z\nabla+\jap z(I-P_z)\nabla,
\end{equation}
where $P_z$ is the projection onto $\R z$. Here $\sigma=0$ for $m=\jap z^k$, $\sigma=s$ for
$m=e^{q\jap z^s}$ with $s\in(0,2)$, and $\sigma=2$ for the Gaussian weights $m=e^{q|z|^2+\omega}$.
For Gaussian weights $Y(m)$ is the norm $Y_2(m)$ of Subsection~\ref{ssec:tri}.

We use the eigenvalue $\ell_1(z):=\ab[\mu](z)\hat z\cdot\hat z$. By \cite[Lemma~2.1(c)]{CM}, we have
$\bb[\mu](z)=-\ell_1(z)z$. Since $|(z-w)\times\hat z|=|w\times\hat z|$, we get
\begin{equation}\label{eq:ell1}
\ell_1(z)=\int|z-w|^\gamma|w_\perp|^2\mu(w)\dd w ,
\end{equation}
where $w_\perp$ is the component of $w$ orthogonal to $\hat z$. The integral converges, since
$|w_\perp|\le|z-w|$ and $\gamma+2>-3$. For radial functions we write $u(r)$, $r=|z|$. For any
function, $\partial_r:=\hat z\cdot\nabla$ is the radial derivative.

\begin{lem}[The radial eigenvalue]\label{lem:ell1}
We have $\ell_1\le C$ on $\R^3$, $\ell_1(z)\le C\jap z^\gamma$, and for $|z|=r\ge1$,
\begin{equation}
\ell_1(z)=2r^\gamma\big(1+O(r^{-2})\big),\qquad\big|\partial_r\ell_1(z)\big|\le Cr^{\gamma-1} .
\end{equation}
\end{lem}

\begin{proof}
We fix a unit vector $e$, write $z=re$, and split the integral at $|w|=r/2$. On $\{|w|\ge r/2\}$,
$|w_\perp|^2|z-w|^\gamma\le|z-w|^{\gamma+2}$ and $|w_\perp|^2|z-w|^{\gamma-1}\le|z-w|^{\gamma+1}$,
both locally integrable. Since $\mu(w)\le e^{-r^2/16}\mu(w)^{1/2}(2\pi)^{-3/4}$ there, this part of
$\ell_1$ and of $\partial_r\ell_1=\gamma\int|z-w|^{\gamma-2}(r-w\cdot e)|w_\perp|^2\mu\dd w$ is
$O(e^{-r^2/16})$. On $\{|w|<r/2\}$, $|z-w|\ge r/2$, so the second part of $\partial_r\ell_1$ is
bounded by $Cr^{\gamma-1}\int|w|^2\mu$. For $\ell_1$ there, we expand
$|z-w|^\gamma=r^\gamma\big(1-\gamma\,w\cdot e/r+O(|w|^2/r^2)\big)$, uniformly for $|w|<r/2$. The
term linear in $w\cdot e$ integrates to zero against $|w_\perp|^2\mu$ over the ball, which is
invariant under $w\cdot e\mapsto-w\cdot e$. Since $\int|w_\perp|^2\mu=2$, we obtain the expansion.
For the global bounds, we have $\ell_1\le\mathrm{tr}\,\ab[\mu]\le2\int|z-w|^{\gamma+2}\mu(w)\dd w\le
C$, as $\gamma+2\in[-1,0]$. Together with the expansion, this implies $\ell_1\le C\jap z^\gamma$.
\end{proof}

\begin{prop}[The mass tail is not dissipative in polynomial weights]\label{prop:sharp}
Let $m=\jap z^k$, $k\ge0$, let $\eta\in C^\infty_c(-1,1)$ be nonzero, and put
$g_R(z)=R^{-k}\eta(|z|-R)$, $h_R(z)=R^{-k}\eta'(|z|-R)$. As $R\to\infty$,
\begin{equation}\label{eq:sharp}
\mathfrak b_m(g_R,h_R)=8\pi R^{\gamma+3}\|\eta'\|^2_{L^2}\big(1+O(R^{-1})\big),\qquad\|g_R\|_{Y(m)}\|h_R\|_{Y(m)}=O(R^{\gamma+2}).
\end{equation}
Consequently $\mathfrak b_m$ is not bounded on $Y(m)\times Y(m)$. Moreover, for
$\varphi_R=g_R+ih_R$,
\begin{equation}\label{eq:sharp-im}
|\mathrm{Im}\langle\mathcal L\varphi_R,\varphi_R\rangle_{L^2(m)}|\ge
cR\,|\mathrm{Re}\langle\mathcal L\varphi_R,\varphi_R\rangle_{L^2(m)}|
\end{equation}
for $R$ large. Hence the numerical range of $\mathcal L$ in $L^2(m)$ is contained in no sector
$\{|\mathrm{Im}\,\zeta|\le C|\mathrm{Re}\,\zeta|\}$ with vertex at $0$.
\end{prop}

For $|z|$ large, $\bb[\mu](z)=-2|z|^\gamma z(1+O(|z|^{-2}))$ points
inward with size $|z|^{\gamma+1}$. The dissipation in the radial direction carries only the factor
$|z|^\gamma$, because $a(z)\hat z\cdot\hat z=0$. A radial bump at distance $R$ is therefore subject
to a coupling of size $R^{\gamma+3}$ against a dissipation of size $R^{\gamma+2}$.

\begin{proof}
All functions are radial and supported in $\{R-1<r<R+1\}$, $r=|z|$, and $R\ge2$. For radial $u$,
$\bb[\mu]\cdot\nabla u=-\ell_1r\,\partial_ru$, and $\ab[\mu]\nabla u\cdot\nabla
v=\ell_1\partial_ru\,\partial_rv$, because $\hat z$ is an eigenvector of $\ab[\mu](z)$. Hence, for
radial $u,v$, we have
\begin{equation}\label{eq:bradial}
\mathfrak b_m(u,v)=-4\pi\int_0^\infty\ell_1r^3\,\partial_r(m^2v)\,u\dd r .
\end{equation}
We put $F(r):=\ell_1(r)r^3m(r)^2R^{-2k}$. On the support, it follows from Lemma~\ref{lem:ell1}
that $F=2R^{\gamma+3}(1+O(R^{-1}))$ and $F'=O(R^{\gamma+2})$. Moreover,
$\ell_1r^3R^{-2k}\partial_r(m^2)=O(R^{\gamma+2})$.

\emph{Step 1: the form.} By \eqref{eq:bradial}, with
$\partial_r(m^2h_R)=R^{-k}(m^2\eta''+\partial_r(m^2)\eta')$, we have
\begin{equation}
\begin{aligned}
\mathfrak b_m(g_R,h_R)&=-4\pi\int F\eta''\eta\dd r+O(R^{\gamma+2})=4\pi\int F\eta'^2\dd r+4\pi\int
F'\eta'\eta\dd r+O(R^{\gamma+2})\\
&=8\pi R^{\gamma+3}\|\eta'\|^2_{L^2}\big(1+O(R^{-1})\big).
\end{aligned}
\end{equation}
Here and below $\eta$ and its derivatives are evaluated at $r-R$.

\emph{Step 2: the norms.} For radial $u$, $\tn(mu)=\partial_r(mu)\hat z$, so
$\|u\|^2_{Y(m)}=4\pi\int r^2\jap r^\gamma\big((mu)^2+(\partial_r(mu))^2\big)\dd r$. Since
$mg_R=\eta\,(1+O(R^{-1}))$ and $\partial_r(mg_R)=\eta'(1+O(R^{-1}))+O(R^{-1})\eta$, we get
$\|g_R\|^2_{Y(m)}\le CR^{\gamma+2}(\|\eta\|^2+\|\eta'\|^2)$. In the same way $\|h_R\|^2_{Y(m)}\le
CR^{\gamma+2}(\|\eta'\|^2+\|\eta''\|^2)$.

\emph{Step 3: the numerical range.} We use the sesquilinear product $\langle
u,v\rangle_{L^2(m)}=\int u\bar v\,m^2$, and we note that the operator $\mathcal L$ is real. For real
$u,v$ in $L^2(m)$, we have $\langle Q(\mu,u),v\rangle_{L^2(m)}=-\int\ab[\mu]\nabla
u\cdot\nabla(m^2v)+ \mathfrak b_m(u,v)$, by \eqref{eq:flux}. Moreover $\langle
Q(u,\mu),v\rangle_{L^2(m)}=O(e^{-R^2/8})$, because $Q(u,\mu)=\ab[u]:\nabla^2\mu-\cb[u]\mu$, and
$\ab[u],\cb[u]$ grow at most polynomially in $R$, while $\mu,\nabla^2\mu$ are $O(e^{-R^2/4})$ on the
support of $v$. For the diffusion part, we have $\int\ab[\mu]\nabla
u\cdot\nabla(m^2v)=4\pi\int\ell_1r^2\partial_ru\,\partial_r(m^2v)\dd r$. With
$u,v\in\{g_R,h_R\}$ this is $O(R^{\gamma+2})$, by the bounds on the radial eigenvalue used
above. For the form, using
\eqref{eq:bradial} and one integration by parts, we get
\begin{equation}
\mathfrak b_m(g_R,g_R)=O(R^{\gamma+2}),\qquad\mathfrak b_m(h_R,h_R)=O(R^{\gamma+2}),\qquad
\mathfrak b_m(h_R,g_R)=-8\pi R^{\gamma+3}\|\eta'\|^2_{L^2}\big(1+O(R^{-1})\big).
\end{equation}
Indeed $\int F\eta\eta'=-\frac12\int F'\eta^2$, and $-4\pi\int F\eta'\cdot\eta'=-8\pi
R^{\gamma+3}\|\eta'\|^2(1+O(R^{-1}))$ up to $O(R^{\gamma+2})$. Hence we have
\begin{equation}
\mathrm{Re}\langle\mathcal L\varphi_R,\varphi_R\rangle=\langle
\mathcal Lg_R,g_R\rangle+\langle\mathcal Lh_R,h_R\rangle=O(R^{\gamma+2}),
\end{equation}
and
\begin{equation}
\mathrm{Im}\langle\mathcal L\varphi_R,\varphi_R\rangle=\langle
\mathcal Lh_R,g_R\rangle-\langle\mathcal Lg_R,h_R\rangle=-16\pi
R^{\gamma+3}\|\eta'\|^2_{L^2}\big(1+O(R^{-1})\big).
\end{equation}
\end{proof}

\begin{prop}[The Gaussian threshold]\label{prop:Gthreshold}
\leavevmode
\begin{enumerate}[label=(\roman*)]
\item If $m=e^{q|z|^2+\omega}$ with $q>0$ and $|\nabla\omega|\le C_1$, then
\begin{equation}\label{eq:Gthr-up}
|\mathfrak b_m(u,\varphi)|\le C(1+q+C_1)\|u\|_{Y(m)}\|\varphi\|_{Y(m)} .
\end{equation}

\item If $m=e^{q\jap z^s}$ with $q>0$ and $s\in(0,2)$, then there are $R_0,c>0$ and, for $R\ge
R_0$, radial $u_R,\varphi_R$ supported in $\{R-1<|z|<R+1\}$ with
\begin{equation}\label{eq:Gthr-low}
|\mathfrak b_m(u_R,\varphi_R)|\ge cR^{1-s/2}\|u_R\|_{Y(m)}\|\varphi_R\|_{Y(m)} .
\end{equation}

\end{enumerate}
Thus, among the weights $\jap z^k$, $e^{q\jap z^s}$ with $s\in(0,2)$,
and the Gaussian weights of~(i), $\mathfrak b_m$ is bounded on $Y(m)\times Y(m)$ exactly for the
Gaussian weights.
\end{prop}

For a Gaussian weight, the gradient of $m^2$ in the mass-tail form is of the same order as the
zero-order term that the weight adds to the dissipation norm, so the coupling is controlled by the
dissipation. For a subgaussian weight it is not, and radial bumps as in
Proposition~\ref{prop:sharp}, now oscillating, give the lower bound.

\begin{proof}
(i) Here $\sigma=2$. By $\bb[\mu]=-\ell_1z$, we have
$\bb[\mu]\cdot\nabla(m^2\varphi)=-\ell_1|z|\,m\big(\partial_r(m\varphi)+(\partial_r\log
m)\,m\varphi\big)$, and $|\partial_r\log m|\le2q|z|+C_1$. Lemma~\ref{lem:ell1} implies that
$\ell_1|z|\le C\jap z^{\gamma+1}$. Since $|\partial_rF|=|P_z\nabla F|\le|\tn F|$, we obtain
\begin{equation}
|\mathfrak b_m(u,\varphi)|\le C\!\int\!\jap z^{\gamma/2}|\tn(m\varphi)|\,\jap
z^{(\gamma+2)/2}|mu|\dd z+C(q+C_1)\!\int\!\jap z^{(\gamma+2)/2}|m\varphi|\,\jap
z^{(\gamma+2)/2}|mu|\dd z,
\end{equation}
and the claim follows from the Cauchy--Schwarz inequality.

(ii) Here $\sigma=s$. Let $\eta\in C^\infty_c(-1,1)$ be nonzero, $\nu:=R^{s/2}$, and define
$u_R,\varphi_R$ by $mu_R=\eta(r-R)\cos\nu r$ and $m\varphi_R=\eta(r-R)\sin\nu r$. We write
$G(r):=\ell_1(r)r^3$. On the support, we have $G=2R^{\gamma+3}(1+O(R^{-1}))$ and
$G'=O(R^{\gamma+2})$ by the same lemma, and $\partial_r\log m=qsr\jap r^{s-2}=O(R^{s-1})$. By
\eqref{eq:bradial} and $\partial_r(m^2\varphi)=m\big(\partial_r(m\varphi)+(\partial_r\log
m)m\varphi\big)$, we get
\begin{equation}
\mathfrak b_m(u_R,\varphi_R)=-4\pi\!\int G\big(\eta'\sin\nu r+\nu\eta\cos\nu r+(\partial_r\log
m)\eta\sin\nu r\big)\eta\cos\nu r\dd r .
\end{equation}
The term with $\eta'$ is $O(R^{\gamma+3})$, and the term with $\partial_r\log m$ is
$O(R^{\gamma+2+s})$. In the main term we write $\cos^2\nu r=\frac12+\frac12\cos2\nu r$. One
integration by parts shows that $\int G\eta^2\cos2\nu r\dd r=O(R^{\gamma+3}\nu^{-1})$, since
$(G\eta^2)'=O(R^{\gamma+3})$. Hence we have
\begin{equation}
\mathfrak b_m(u_R,\varphi_R)=-4\pi\nu R^{\gamma+3}\|\eta\|^2_{L^2}\big(1+O(R^{-1})+O(\nu^{-1})+O(R^{s/2-1})\big),
\end{equation}
because $R^{\gamma+3}=o(\nu R^{\gamma+3})$ and $R^{\gamma+2+s}=\nu R^{\gamma+3}R^{s/2-1}$. Therefore
$|\mathfrak b_m(u_R,\varphi_R)|\ge2\pi\nu R^{\gamma+3}\|\eta\|^2$ for $R\ge R_0$. For the norms, we
have $\tn(mu_R)=\partial_r(mu_R)\hat z$ with $\partial_r(mu_R)=\eta'\cos\nu r-\nu\eta\sin\nu r$, so
\begin{equation}
\begin{aligned}
\|u_R\|^2_{Y(m)}&=4\pi\!\int r^2\big(\jap r^{\gamma+s}(mu_R)^2+\jap
r^\gamma(\partial_r(mu_R))^2\big)\dd r\\
&\le CR^{\gamma+2}\big(R^s\|\eta\|^2+\|\eta'\|^2+\nu^2\|\eta\|^2\big)\le
CR^{\gamma+2+s}\big(\|\eta\|^2+\|\eta'\|^2\big),
\end{aligned}
\end{equation}
and the same for $\varphi_R$. The ratio is at least $c\,\nu R^{\gamma+3}/R^{\gamma+2+s}=cR^{1-s/2}$.
\end{proof}

\subsubsection{Non-uniform relaxation}\label{sssec:obs3}

Here we show that the leaves far out in $\xi$ enter the relaxed regime arbitrarily late. Part~(i) of
the next lemma shows that their clocks run arbitrarily slowly. Part~(ii) shows that a leaf
cannot move far from its state at $T_0$ before its clock has run.

\begin{lem}[Unrelaxed tail leaves]\label{lem:Wtail}
\leavevmode
\begin{enumerate}[label=(\roman*)]
\item Let $\rho\in L^1(\R^3)$, $\rho>0$, with $|\nabla\log\rho|\le L$. Then $\rho(\xi)\to0$ as
$|\xi|\to\infty$. Hence $\lambda(\xi)\to0$. Applied to $\rho=\rho_\xi$, the fiber masses, it implies
that $\lambda_{\rm ref}\le\theta_-^{\gamma/2}\rho_\xi\to0$. Thus, for every $t$, the effective time
$s_\xi(t)$ tends to $0$ as $|\xi|\to\infty$, and the entry time $t_\xi$ tends to $\infty$.
\item For $t\ge T_0$ let $\mathcal N_\xi(t)(z):=g(t,z,\xi-z/t)/\rho_g(t,\xi)$ be the normalized
leaf through $x=t\xi$, let $\hat p_g(t,\xi):=(1,y_{0,g},\theta_g)(t,\xi)$, and let
\begin{equation}\label{eq:Wtail-d}
d_\xi(t):=\big\|\mathcal N_\xi(t)-M_{\hat p_g(t,\xi)}\big\|_{L^2(\R^3)}
\end{equation}
be its distance from the Maxwellian with the same moments. As long as the bounds~(U1)--(U2) of
Theorem~\ref{thm:CUboot} hold on $[T_0,t]$, with (F3) of Subsection~\ref{ssec:inputs}, the data at
$T_0$ satisfy Corollary~\ref{cor:LWPa}(a3), and $\rad_0T_0a^{1/2}\ge1$ (a lower bound on $T_0$), with
$a:=\frac14\min\{q_0,\beta/\theta_+\}$,
\begin{equation}\label{eq:Wtail-bound}
\big\|\mathcal N_\xi(t)-\mathcal N_\xi(T_0)\big\|_{L^2(e^{a|z|^2/2})}+|d_\xi(t)-d_\xi(T_0)|\le
C_{\rm tail}\big(s_\xi(t)+T_0^{-1}\big),
\end{equation}
where $C_{\rm tail}$ depends only on the constants fixed before $T_0$ (Appendix~\ref{app:const}),
and not on $T_0$ or $\eps$.
\end{enumerate}
Consequently, if $d_\xi(T_0)\ge2\delta$ on a sequence $|\xi_k|\to\infty$ and $C_{\rm
tail}T_0^{-1}\le\delta/2$, then at every $t\ge T_0$ there are leaves with $d_\xi(t)\ge\delta$.
\end{lem}

\begin{proof}
(i) If $\rho(\xi_k)\ge c$ along a sequence $|\xi_k|\to\infty$, then $\rho\ge ce^{-L}$ on the unit
balls around the $\xi_k$, infinitely many of which are disjoint, which contradicts $\rho\in L^1$.
The fiber masses satisfy $|\nabla\log\rho_\xi|\le L$ by (H1), that is, Corollary~\ref{cor:LWPa}(a3),
and $\int\rho_\xi\dd\xi<\infty$ since the total mass is finite.

(ii) \emph{Step 1: the integral identity.} By Proposition~\ref{prop:leaf}, we have
$\partial_tg=\kappa\cQ_t[g]$ at fixed $(y,v)$. Since $\xi-z/t=(\xi-z/s)+z(s^{-1}-t^{-1})$, we get
\begin{equation}\label{eq:Wtail-id}
g(t,z,\xi-z/t)-g(T_0,z,\xi-z/T_0)=\int_{T_0}^t\kappa(s)\,\cQ_s[g](z,\xi-z/t)\dd
s+\big(g(T_0,z,\xi-z/t)-g(T_0,z,\xi-z/T_0)\big).
\end{equation}

\emph{Step 2: a supremum over nearby leaves.} We fix $a:=\frac14\min\{q_0,\beta/\theta_+\}$. At time
$s$ the point $(z,\xi-z/t)$ lies on the leaf $\xi'(z)=\xi+z(s^{-1}-t^{-1})$, and
$|\xi'(z)-\xi|\le|z|/T_0$: at each $z$ a different leaf is sampled. Let $u(s,z,\xi')$ be the
normalized leaf $\xi'$ at time $s$, or one of its $z$-derivatives of order $\le2$. For fixed $z$,
Sobolev's inequality in $\xi'$ on cubes of side $\rad_0$ implies
\begin{equation}\label{eq:Wtail-sob}
\sup_{|\xi'-\xi|\le|z|/T_0}|u(s,z,\xi')|^2\le C\sum_{|\varsigma|\le2}\rad_0^{2|\varsigma|-3}\int_{|\xi''-\xi|\le|z|/T_0+2\rad_0}\big|\partial^\varsigma_\xi
u(s,z,\xi'')\big|^2\dd\xi'' .
\end{equation}
We multiply by $e^{a|z|^2}$ and integrate in $z$. Exchanging the integrals, the $z$-integral at
fixed $\xi''$ runs over $|z|\ge T_0(|\xi''-\xi|-2\rad_0)_+$, where $e^{a|z|^2}\le
e^{2a|z|^2}e^{-aT_0^2(|\xi''-\xi|-2\rad_0)_+^2}$. Moreover, $\rad_0^{-3}\int
e^{-aT_0^2(|\xi''-\xi|-2\rad_0)_+^2}\dd\xi''\le C$ for $\rad_0T_0a^{1/2}\ge1$. Under~(U2),
Lemma~\ref{lem:TWcalc}(i) bounds $\|\partial^\varsigma_\xi u(s,\cdot,\xi'')\|_{L^2(e^{2a|z|^2})}$,
$|\varsigma|\le2$, uniformly in $s$ and $\xi''$. Here $e^{2a|z|^2}\le C\hat m_\xi^2$ by
Lemma~\ref{lem:TG0}(i), since $2a\le q_0/2$. Thus $z\mapsto\sup_{\xi'}|u(s,z,\xi')|$ is bounded in
$L^2(e^{a|z|^2})$, uniformly in $s$.

\emph{Step 3: the two terms of \eqref{eq:Wtail-id}.} The leaf collision operator is
$\ab[h]:\nabla^2h-\cb[h]h$, with $|\ab[h_{\xi'}]|+|\cb[h_{\xi'}]|\le
C\sup_{\xi'}\|h_{\xi'}\|_{L^1\cap L^\infty}$ (for $\gamma=-3$, $\cb[h]=-8\pi h$). Here
$L':=\frac54L$, as in Subsection~\ref{ssec:cells}. By (F3),
$\rho_g(s,\xi')/\rho_{\xi'}\in[\frac12,2]$ and $\lambda(p_g(s,\xi'))/\lambda_{\rm
ref}(\xi')\in[c_\lambda,c_\lambda^{-1}]$. Moreover, $|\nabla\log\rho_\xi|\le L\le2L'$ and
$|\nabla\log\lambda_{\rm ref}|\le2L'$. The term $\cQ_s[g]$ is quadratic in the leaf: on the leaf
$\xi'$ it is $\rho_g(s,\xi')^2Q(u,u)$, with $u$ the normalized leaf. Since the leaf temperatures lie
in $[\frac12\theta_-,\theta_++1]$, $\rho_g\le C\lambda(p_g)$, and so we have
$\kappa\rho_g(s,\xi')^2/\rho_g(T_0,\xi)\le C\kappa\lambda_{\rm ref}(\xi)e^{4L'|z|/T_0}$. Thus the
first term, divided by $\rho_g(T_0,\xi)$, is bounded in $L^2(e^{a|z|^2/2})$ by
$C\int_{T_0}^t\kappa\lambda_{\rm ref}(\xi)\dd s=Cs_\xi(t)$. Here the factors $e^{4L'|z|/T_0}$ are
controlled by the Gaussian weight. The second term is bounded by
$|z|T_0^{-1}\sup|\partial_vg(T_0)|$, hence, by (H1) and $a<\beta/\theta_+$, by $CT_0^{-1}$ in the
same norm. Hence the left side of the identity, divided by $\rho_g(T_0,\xi)$, is bounded in
$L^2(e^{a|z|^2/2})$ by $C(s_\xi(t)+T_0^{-1})$, and, integrating in $z$, so is
$|\rho_g(t,\xi)-\rho_g(T_0,\xi)|/\rho_g(T_0,\xi)$.

\emph{Step 4: normalization.} We write
\begin{equation}\label{eq:Wtail-norm}
\mathcal N_\xi(t)-\mathcal N_\xi(T_0)=\frac{g(t,z,\xi-z/t)-g(T_0,z,\xi-z/T_0)}{\rho_g(t,\xi)}+
\mathcal N_\xi(T_0)\,\frac{\rho_g(T_0,\xi)-\rho_g(t,\xi)}{\rho_g(t,\xi)},
\end{equation}
Here $\rho_g(T_0,\xi)\le4\rho_g(t,\xi)$ by (F3), and $\mathcal N_\xi(T_0)$ is bounded in
$L^2(e^{a|z|^2})$, as shown above. This proves the bound on $\mathcal N_\xi(t)-\mathcal N_\xi(T_0)$
in \eqref{eq:Wtail-bound}. The mean and the temperature of the normalized leaf are moments of it, so
they change by at most a constant multiple of this bound, with a constant that depends only on $a$
and $L$. Moreover, the map $(y_0,\theta)\mapsto M_{(1,y_0,\theta)}$ is Lipschitz into $L^2$ for
$\theta\in[\frac12\theta_-,\theta_++1]$. By the triangle inequality, we obtain the bound on
$d_\xi(t)-d_\xi(T_0)$.

For the last statement, we note that $s_{\xi_k}(t)\to0$ by~(i), so $d_{\xi_k}(t)\ge2\delta-C_{\rm
tail}(s_{\xi_k}(t)+T_0^{-1})\ge\delta$ for $k$ large.
\end{proof}

\begin{remark}[Maxwellian spatial profile]\label{rem:nonunif}
Let $d_0(\xi)$ and $d_0$ be as in~\eqref{eq:d0}. By the proof of Theorem~\ref{thm:main}(iv), at
$t=T_0$ the normalized leaf $\xi$ is within $C(L/T_0+\eps T_0^{\kappa_\gamma})$ of $F_0(\cdot,\xi)$
in $L^2(e^{a'|z|^2})$, $a'=q_d/2$, and its mean and temperature are within the same distance of
$(y_{0,\infty}(\xi),\theta_\infty(\xi))$. Therefore
$|d_\xi(T_0)-d_0(\xi)|\le C(L/T_0+\eps T_0^{\kappa_\gamma})$ for every $\xi$.

If $d_0>0$, the last clause of Theorem~\ref{thm:main}(iv) follows from this and
Lemma~\ref{lem:Wtail}.

If $F_0$ is a Maxwellian, then $d_0=0$, and at $T_0$ every leaf is within $C(L/T_0+\eps
T_0^{\kappa_\gamma})$ of a Maxwellian, uniformly in $x$, also in the Gaussian norm of
Theorem~\ref{thm:main}(ii), whose weight is $\le Ce^{a'|z|^2}$. Thus for a Maxwellian $F_0$ no lower
bound of the kind in Theorem~\ref{thm:main}(iv) holds at $T_0$ beyond $C(L/T_0+\eps
T_0^{\kappa_\gamma})$, and the hypothesis $d_0>0$ cannot be replaced by a lower bound that does not
depend on $F_0$. For such data the non-uniformity in part~(iv) concerns the clock
$s_\xi$, and not necessarily the distance from equilibrium.

For data that are not products there is an intermediate case: $d_0=0$ although $F_0$ is not a
Maxwellian, for instance when $F_0(\cdot,\xi)$ tends to a Maxwellian as $|\xi|\to\infty$. Only the
leaves with $|\xi|$ large enter the argument, so the last clause of part~(iv) holds, with another
constant $C_{\rm nM}$, with $d_0$ replaced by $\limsup_{|\xi|\to\infty}d_0(\xi)$, when this is
positive, by the last statement of that lemma. When it vanishes, we do not know whether such a lower
bound holds.
\end{remark}

Therefore the bootstrap cannot be started from a fixed time after which every leaf is near
equilibrium. In Section~\ref{sec:boot} we organize it by the leaf clock $s_\xi$, with a front that
moves out to $|\xi|=\infty$.

\section{Homogeneous trajectories}\label{sec:homog}

In this section we collect the properties of homogeneous Landau trajectories that we use on the
leaves, in particular a uniform $L^1$ rate and the entry theorem, whose proofs are in
Appendix~\ref{app:hom}.

Throughout this section, $F(\sigma)$ solves $\partial_\sigma F=Q(F,F)$ on $\R^3$.

\subsection{The class of initial leaves}\label{ssec:homclass}

All bounds of this section are uniform on the following class of normalized initial leaves. Let
$\mathfrak K(C_{\rm dat},\beta)$ be the set of $H_0\ge0$ with unit mass, zero mean, unit temperature
and $H_0\le C_{\rm dat}e^{-\beta|z|^2}$. For $N\ge1$ and a fixed large $k_N$ let
\begin{equation}\label{eq:KN}
\mathfrak K^N(C_{\rm dat},\beta,A_N)=\Big\{H_0\in\mathfrak K(C_{\rm dat},\beta):\ \sum_{|\alpha|\le2N+4}\sup_z\jap
z^{k_N}|\nabla^\alpha H_0|\le A_N\Big\}.
\end{equation}
For data in $\Dc$, and for the data of Theorem~\ref{thm:t1} with a reference profile in
$\Dc^\sharp$, the normalized fibers at $T_0$ belong to $\mathfrak K^7(C_{\rm dat},\beta,A_7)$ with
constants depending only on $\Dc$, provided $\eps\le\eps_0(T_0)$, by Corollary~\ref{cor:LWPa}(a3).
Indeed, under free transport the leaf through $x=t\xi$ is
\begin{equation}\label{eq:freeleaf}
z\mapsto\eps\rho_0(\xi-z/t)F_0(z,\xi-z/t),
\end{equation}
which is a small perturbation of $\eps\rho_0(\xi)F_0(z,\xi)$ for $t$ large, by (D3) and (D4).

\subsection{Inputs from the homogeneous theory}\label{ssec:hominputs}

Here we list what we use from the homogeneous theory. We use the results (GS1)--(GS5) of \cite{GS},
stated in (E2) of Inputs~\ref{hyp:E}, where $i(F)$ is the Fisher information. For (GS4) and (GS5),
which are \cite[Thm.~2.3, Prop.~B.2]{GS}, we write out proofs for $\gamma\in[-3,-2]$ with the
constants that we need, in Subsection~\ref{sssec:homGauss}. From these, with all constants tracked,
we obtain a Maxwellian bound on every bounded interval,
\begin{equation}\label{eq:homMaxw}
F(\sigma)\le C_M(S)e^{-\beta_M|v|^2}\qquad(\sigma\le S),
\end{equation}
and Gaussian bounds on $N$ derivatives along the trajectory, uniformly for $F_0\in\mathfrak K^N$.
This is Proposition~\ref{prop:RS}, which we prove in Appendix~\ref{app:hom}. A Maxwellian bound
uniform for all times is not available in general \cite[Remark~11.1]{GS}. However,
Lemma~\ref{lem:Abarrier} provides one after the trajectory has entered a neighborhood of $\mu$.

\subsection{Entry into the near-Maxwellian regime}\label{ssec:entry}

By Theorem~\ref{thm:enter}, every trajectory that starts in $\mathfrak K^N$ enters a fixed
neighborhood of $\mu$ after a time that depends only on the class. The main tool is the uniform
$L^1$ rate of Proposition~\ref{prop:Lrate}. We say that $F$ is a \emph{$\delta$-admissible
background on $[s_*,\infty)$ at level $k_1$} if, with $X_k=L^2(\jap v^k)$ and $Y_k=H^1_*(\jap v^k)$
the dissipation space of \cite{CM},
\begin{equation}\label{eq:AB}
\text{(AB1)}\ \ \sup_{\sigma\ge s_*}\|F(\sigma)-\mu\|_{X_{k_1}}\le\delta,\qquad\text{(AB2)}\ \
\int_{s_*}^\infty\|F(\sigma)-\mu\|^2_{Y_{k_1}}\dd\sigma\le\delta^2 .
\end{equation}

\begin{prop}[Uniform $L^1$ rate]\label{prop:Lrate}
Let $\gamma\in[-3,-2]$ and $l>2+\frac73|\gamma|$, and put
$\vartheta_l:=\frac{3(l-2)}{4|\gamma|}-\frac74>0$. Let $F(0)\ge0$ have unit mass, zero mean and unit
temperature, and let $F$ be the solution of (GS1). Assume:
\begin{itemize}
\item $\mathcal H(F(0))\le\bar H$, $\int\jap v^lF(0)\le A_l$ and $F(0)\le C_{\rm
dat}e^{-\beta|v|^2}$;
\item the Fisher information along $F$ satisfies $\sup_{\sigma>0}i(F(\sigma))\le I_0$;
\item $F$ is an H-solution with the entropy identity from $\sigma=0$ (by Lemma~\ref{lem:Rqual} this
holds for $F(0)\in\mathfrak K^N$, $N\ge2$).
\end{itemize}
Then, for all $\sigma\ge0$,
\begin{equation}\label{eq:Lrate}
\mathcal H\big(F(\sigma)\,|\,\mu\big)\le C_\vartheta(1+\sigma)^{-\vartheta_l},\qquad\|F(\sigma)-\mu\|_{L^1}\le(2C_\vartheta)^{1/2}(1+\sigma)^{-\vartheta_l/2},
\end{equation}
with $C_\vartheta=C_\vartheta(\gamma,l,\bar H,I_0,A_l)$. On $\mathfrak K^N(C_{\rm dat},\beta,A_N)$,
$N\ge2$, the quantities $\bar H$, $I_0$ and $A_l$ are bounded in terms of $(C_{\rm
dat},\beta,A_N,l)$ (for $I_0$ by Lemma~\ref{lem:FPgl}(ii), Lemma~\ref{lem:RSshortF}(iii) and (GS1),
see the proof), so the rate is uniform there. Since $\vartheta_l\to\infty$ as $l\to\infty$, we have
$\|F(\sigma)-\mu\|_{L^1}\le C_m(1+\sigma)^{-m}$ for every $m$.
\end{prop}

We give the proof in Appendix~\ref{app:hom}. It uses four inputs. The first one is the
entropy--entropy-dissipation inequality
\begin{equation}\label{eq:EEDcdh}
I_\gamma(F|\mu)\le KM_{2-\gamma}(F)\mathcal D(F)
\end{equation}
of \cite[Thm.~1, (19)]{CDH}, which we prove for $-4<\gamma\le-2$ in Lemma~\ref{lem:EL1}. The others
are the linear growth of the moments \cite[Lemma~8]{CDH}, the monotonicity of the Fisher
information (GS1), and the logarithmic Sobolev inequality on balls, which follows from
\cite[Cor.~5.7.2]{BGL} by approximation.
In the proof we localize the entropy to a ball of radius $R$, where $\jap v^\gamma\ge\jap R^\gamma$.
We control the entropy outside the ball by the moments and the Fisher information, and optimize in
$R$. The resulting differential inequality for $\mathcal H(F|\mu)$ yields the rate, and the
Csisz\'ar--Kullback--Pinsker inequality implies the $L^1$ bound. For $\gamma\in(-3,-2]$ no rate that
is uniform on $\mathfrak K^N$ seems to be available in the literature.
For $\gamma=-3$ the proposition also gives an alternative proof of a polynomial rate,
uniform on $\mathfrak K^N$ (compare \cite[Thm.~2]{CDH}).

\begin{lem}[Qualitative regularity]\label{lem:Rqual}
Let $N\ge2$, $H_0\in\mathfrak K^N(C_{\rm dat},\beta,A_N)$, and let $H$ be the solution of (GS1).
Then for every $k$ and every $b>0$:
\begin{itemize}
\item $H-\mu\in C([0,b];X_k)\cap L^2(0,b;Y_k)$;
\item $\sigma\mapsto\frac12\|H(\sigma)-\mu\|^2_{X_k}$ is Lipschitz on $[0,b]$ with derivative
$\langle Q(H,H),H-\mu\rangle_{X_k}$.
\end{itemize}
Moreover, $H$ is an H-solution in the sense of \cite[Def.~1]{CDH}, with the kernel
$|z|^{\gamma+2}\Pi(z)$ in place of the Coulomb kernel, with the entropy identity $
\mathcal H(H(\sigma))+\int_0^\sigma\mathcal D(H)=\mathcal H(H_0)$.
\end{lem}

The proof is also given there. By the Maxwellian bound \eqref{eq:homMaxw} and the
propagation of weighted derivative bounds on bounded intervals of Subsection~\ref{ssec:hominputs},
we can bound $H-\mu$ in weighted $C^2$, uniformly on $[0,b]$, with Gaussian decay. This implies the
continuity in $X_k$ and the bound in $L^2(0,b;Y_k)$. Together with Glaeser's inequality, it also
justifies differentiating $\mathcal H(H(\sigma))$ under the integral. The lemma provides the
qualitative input of the entry theorem below, namely the H-solution property, the
entropy identity, and $H-\mu\in C([0,b];X_{k_1})\cap L^2(0,b;Y_{k_1})$.

\begin{theo}[Entry]\label{thm:enter}
Let $\gamma\in[-3,-2]$, $k_1\ge7$ and $\delta>0$. There is
$T_*=T_*(\gamma,C_{\rm dat},\beta,A_N,k_1,\delta)$ such that every trajectory with $F(0)\in
\mathfrak K^N(C_{\rm dat},\beta,A_N)$, $N\ge2$, is a $\delta$-admissible background on $[s_*,\infty)$ for some
$s_*\le T_*$.
\end{theo}

We give the proof in Subsection~\ref{sssec:homenter}, for all
$\gamma\in[-3,-2]$ at once, in three steps. First, Proposition~\ref{prop:Lrate} yields convergence
in $L^1$ at a polynomial rate, uniformly on $\mathfrak K^N$, and together with the moment bounds
\cite[Lemma~8]{CDH} we get convergence in weighted $L^1$. Next, the Fisher information, which is
nonincreasing by (GS1), bounds $F$ in $L^3$ uniformly in time, and by interpolation we obtain
convergence in $X_{k_1}$ at every time. Finally, an inequality of the form of \cite[Prop.~5.1]{CM},
derived in Step~3 of that proof and justified along the trajectory by Lemma~\ref{lem:Rqual},
implies (AB1) and (AB2) from a time $s_*\le T_*$ at which $F-\mu$ is small enough in $X_{k_1}$.
At $\gamma=-2$ the inputs from \cite{CM} are those checked in Appendix~\ref{app:CM}.

\subsection{Uniform bounds after trapping}\label{ssec:trapping}

The first lemma below provides bounds that hold for all times after entry, namely an $L^\infty$
bound, bounds on the coefficients $\ab[F]$ and $\cb[F]$, and a Gaussian bound uniform in time. With
them, we show in the second lemma that the convergence becomes Gaussian and uniform on the class.

\begin{lem}[Uniform bounds after trapping]\label{lem:Abarrier}
Let $F(0)\in\mathfrak K^N$, $N\ge2$, and let $s_*\le T_*$ be the entry time of
Theorem~\ref{thm:enter} for a fixed small $\delta_0$, with $k_1\ge k_\flat(\gamma)$, where
$k_\flat(-3):=11$ and
\begin{equation}\label{eq:kflat}
k_\flat(\gamma):=\max\{7,\lfloor\frac{3|\gamma|}{2(3-|\gamma|)}\rfloor+1\}
\end{equation}
for $\gamma\in(-3,-2]$ (so $k_\flat(-2)=7$, $k_\flat(-\frac52)=8$, and $k_\flat(\gamma)\to\infty$ as
$\gamma\to-3^+$). There are constants $K_\infty,C_a,c_G,C_c,\beta',C_\flat$, depending only on
$(\gamma,C_{\rm dat},\beta,A_N,k_1)$, such that:
\begin{enumerate}[label=(\roman*)]
\item $\|F(\sigma)\|_{L^\infty}\le K_\infty$ for $\sigma\ge s_*+1$;
\item for $\sigma\ge s_*+1$, $|v|\ge2$ and $\zeta\perp v$,
\begin{equation}\label{eq:Abar-a}
\ab[F](v)\hat v\cdot\hat v\le C_a\jap v^\gamma,\qquad\ab[F](v)\zeta\cdot\zeta\ge c_G\jap
v^{\gamma+2}|\zeta|^2,
\end{equation}
and $-\cb[F](v)\le C_c\jap v^{\gamma}$ for $\gamma>-3$, $F(\sigma,v)\le C_c\jap v^{-2}$ for
$\gamma=-3$, $\sigma\ge s_*+2$;
\item with $\beta'=\min\{\beta_M/2,\frac18,c_G/(2C_a)\}$, $F(\sigma,v)\le C_\flat e^{-\beta'|v|^2}$
for all $\sigma\ge0$.
\end{enumerate}
\end{lem}

\begin{proof}
\emph{Step 1: the $L^\infty$ bound (i).} We apply (GS3) on $[s_*,S]$ for every $S$, together with
the bound
\begin{equation}\label{eq:Abar-L52}
\|\jap v^{2k_1/5}F\|_{L^{5/2}}\le\|\jap v^{k_1}F\|_{L^2}^{2/5}\|F\|_{L^3}^{3/5},
\end{equation}
which follows from the uniform bound on the Fisher information (GS1) and the trapping. Note that
its weight $\frac25k_1\ge\frac{14}5$ exceeds $k_{\rm GS}(\gamma,\frac52)\le\frac{62}{25}$, and
$\frac3{5+\gamma}<\frac52<\frac3{3+\gamma}$.

\emph{Step 2: the coefficients (ii).} We use the identity $\ab[F](v)\hat v\cdot\hat
v=\int|v-w|^\gamma|w_\perp|^2F(w)\dd w$, where $w_\perp$ is the component of $w$ orthogonal to $\hat
v$, split into $|w|\le|v|/2$, $|w|>|v|/2$ and $|v-w|<1$, with the moment bounds that hold after
trapping. For $\gamma>-3$ we bound the part of $-\cb[F](v)$ from $|v-w|<1$ by H\"older's
inequality with an exponent $q'\ge2$ such that $\frac1{q'}<1-\frac{|\gamma|}3$ and
$\frac{2k_1}{q'}\ge|\gamma|$, which exists since $k_1\ge
k_\flat(\gamma)>\frac{3|\gamma|}{2(3-|\gamma|)}$, and by $\|F\|_{L^{q'}(B_1(v))}\le
K_\infty^{1-2/q'}\|F\|_{L^2(B_1(v))}^{2/q'}\le C\jap v^{-2k_1/q'}$, which follows from (AB1). This
is $\le C\jap v^\gamma$. The threshold is sharp under these constraints, since a spike of height
$K_\infty$ on a ball of radius $\jap v^{-2k_1/3}$ is compatible with them. The lower bound in
\eqref{eq:Abar-a} is the quantitative ellipticity estimate of Lemma~\ref{lem:RGellip}. For
$\gamma=-3$, the bound on $F$ is a Moser iteration with explicit constants \cite{Mos64}, which we
carry out in Proposition~\ref{prop:ABiv}. It requires $k_1\ge11$ for the moment $M_9$.

\emph{Step 3: the Gaussian bound (iii), a supersolution.} Let $\Phi=C_\flat e^{-\beta'|v|^2}$. Then
we have $\nabla^2\Phi=\Phi(4\beta'^2v\otimes v-2\beta'I)$, and by (ii), for $|v|\ge2$, we get
\begin{equation}\label{eq:Abar-super}
\ab[F]:\nabla^2\Phi\le\Phi\big(4\beta'^2|v|^2\ab[F]\hat v\cdot\hat
v-2\beta'\,\mathrm{tr}_\perp\ab[F]\big)\le\Phi\jap
v^{\gamma+2}(4\beta'^2C_a-4\beta'c_G)\le-2\beta'c_G\Phi\jap v^{\gamma+2}.
\end{equation}
By (ii), the reaction term $-\cb[F]\Phi$ is smaller than $\beta'c_G\jap v^{\gamma+2}\Phi$ for
$|v|\ge R_*$, where we take $R_*\ge2$ so that~\eqref{eq:Abar-super} applies. Hence we have
$\ab[F]:\nabla^2\Phi-\cb[F]\Phi<0$ on $\{|v|>R_*\}$ for $\sigma\ge s_*+2$.

\emph{Step 4: comparison.} We choose $C_\flat>\max\{C_M(T_*+2),K_\infty e^{\beta'R_*^2}\}$
(strictly, so that the contact cannot occur first on $\{|v|=R_*\}$ or at $\sigma=s_*+2$). On $[0,s_*+2]$, we
have $F\le C_Me^{-\beta_M|v|^2}\le\Phi$. Afterwards, we have $F\le K_\infty<\Phi$ on $|v|\le R_*$,
and $F<\Phi$ for $|v|$ large on every bounded interval. At a first contact point $(\sigma_c,v_c)$
with $|v_c|>R_*$ we would have $F=\Phi$, $\nabla^2(F-\Phi)\le0$ and $\partial_\sigma(F-\Phi)\ge0$,
and hence
\begin{equation}
\partial_\sigma F=\ab[F]:\nabla^2F-\cb[F]F\le\ab[F]:\nabla^2\Phi-\cb[F]\Phi<0=\partial_\sigma\Phi ,
\end{equation}
which is a contradiction. This completes the proof.
\end{proof}

\begin{lem}[Uniform Gaussian relaxation]\label{lem:Aentry}
For every $\delta'>0$ there is $S_\flat(\delta')$, depending only on
$(\delta',\gamma,C_{\rm dat},\beta,A_N,k_1)$, such that for $F(0)\in\mathfrak K^N$, $N\ge2$,
$\sigma\ge S_\flat(\delta')$ and $q\le\beta'/4$,
\begin{equation}
\|F(\sigma)-\mu\|_{L^2(e^{2q|v|^2})}\le\delta' .
\end{equation}
\end{lem}

\begin{proof}
We split at $R\ge1$. Inside, we have $\int_{|v|\le R}e^{2q|v|^2}|F-\mu|^2\le
e^{2qR^2}(K_\infty+1)\|F-\mu\|_{L^1}$. Outside, by Lemma~\ref{lem:Abarrier}(iii), we have
$\int_{|v|>R}e^{2q|v|^2}|F-\mu|^2\le C(C_\flat^2+1)e^{-\beta'R^2}$. We choose $R$ to make the second
term small, and then $S_\flat\ge T_*+1$ (for the bound $\|F-\mu\|_{L^\infty}\le K_\infty+1$ of
part~(i) of that lemma) such that the $L^1$ norm, which tends to zero uniformly on $
\mathfrak K^N$ by the proof of Theorem~\ref{thm:enter}, makes the first term small.
\end{proof}

This lemma makes the entry level $\delta_1$ independent of the entry time. Indeed, we choose
$\delta_1$ first, and the entry threshold $S_*$ afterwards.

\section{Gaussian leaf energies and the leaf Euler system}\label{sec:weights}

\subsection{Gaussian leaf energies}\label{ssec:gle}

On the leaf $\xi$ we normalize by the frame of the leaf parameters, and we measure the micro part in
Gaussian weights. For a leaf $\xi$ at time $t$, let $p=p_g(t,\xi)$ and
\begin{equation}\label{eq:gle-psi}
\psi_\xi(\hat z)=\big(T_p^{-1}(g-\mathcal M[g])^{(t\xi)}\big)(\hat z),\qquad\hat m_\xi(\hat
z)=e^{q_\xi|\hat z-\hat z_0|^2+\tilde\omega(\hat z)},
\end{equation}
where $\mathcal M[g]$ is the local leaf Maxwellian of Definition~\ref{def:Mloc}, and
\begin{equation}\label{eq:gle-q}
q(t)=q_0-C_2(T_0^{-1}-t^{-1}),\qquad q_\xi=q(t)\theta_g,\qquad\hat z_0=-\theta_g^{-1/2}y_{0,g} .
\end{equation}
Therefore $q(t)$ decreases slowly. The function $\hat m_\xi$ is the normalized leaf weight of
$A_{\rm ref}(v)e^{q(t)|y|^2}$ in the frame $p$. The
factor $e^{\tilde\omega}$ is slowly varying, and it comes from the amplitude
\begin{equation}\label{eq:gle-Aref}
A_{\rm ref}(\xi):=A_{p_{\rm ref}(\xi)},\qquad A_p:=\rho^{-1}\theta^{3/4}\ \text{ for
}p=(\rho,y_0,\theta),
\end{equation}
where
\begin{equation}\label{eq:gle-pref}
p_{\rm ref}(\xi):=p_g(T_0,\xi)=\big(e^{\ell_{\rm ref}(\xi)},y_{0,g}(T_0,\xi),\theta_g(T_0,\xi)\big)
\end{equation}
are the leaf parameters at $T_0$, and $\ell_{\rm ref}$ and $\varpi'_{\rm ref}=\varpi'_g(T_0,\cdot)$
are as in Lemma~\ref{lem:Mref}. The amplitude $A_{\rm ref}$ does not depend on $t$, and, by
Corollary~\ref{cor:LWPa}(a4), $|\nabla\log A_{\rm ref}|\le L$ on $\R^3$. Lemma~\ref{lem:Gleafweight}
computes the static model of this weight, namely $A_\infty(v)e^{q|y|^2}$ in the frame of the limit
parameters of free transport. Lemma~\ref{lem:TG0}(i) carries the computation over to the frames
$p_g(t,\xi)$ and to the factor $A_{\rm ref}$.

The \emph{weight class} $\mathcal W(q_\pm,\eps_1)$ consists of the weights
\begin{equation}\label{eq:gle-W}
m=e^{q|\hat z|^2+\omega},\qquad q\in[q_-,q_+],\quad|\nabla\omega|\le\eps_1,\quad|\omega(0)|\le\omega_0
.
\end{equation}
It is stable under the changes of frame along the evolution.

We need the following two properties. Let $\mathcal L=Q(\mu,\cdot)+Q(\cdot,\mu)$.

\begin{prop}[Gaussian hypodissipativity and trilinear bounds]\label{prop:triG}
Let $q_{\rm tri}:=\frac18$. There are $C_\gamma$ and $c_W,C_W>0$ such that for $m\in
\mathcal W(q_\pm,\eps_1)$ with $q_+\le q_{\rm tri}$ and $C_\gamma\eps_1^2\le q_-$:
\begin{enumerate}[label=(\roman*)]
\item \emph{(hypodissipativity)} there is a norm $\|\cdot\|_{\tilde X(m)}$, equivalent to
$X(m)=L^2(m)$ with constants $C_W$, such that
\begin{equation}\label{eq:triG-hypo}
\langle\mathcal L\varphi,\varphi\rangle_{\tilde X(m)}\le-c_W\|\varphi\|^2_{Y(m)}
\end{equation}
for moment-free $\varphi$ (in the energies of Section~\ref{sec:cells} it is realized by the form
$\eta\|\cdot\|^2_{X}+\chi\nn{\cdot}^2_\natural$), where $\|\cdot\|_{Y(m)}$ is the anisotropic
dissipation norm of \cite{CM};
\item \emph{(trilinear bound)}
\begin{equation}\label{eq:triG-tri}
|\langle Q(f,g),h\rangle_{L^2(m)}|\le C\|f\|_{X(m)}\|g\|_{Y(m)}\|h\|_{Y(m)}
\end{equation}
plus symmetric terms, without the mass-tail defect;
\item \emph{(arbitrary backgrounds)} for every background $G$ in the class of
Theorem~\ref{thm:RGF}(iv), the linearization around $G$ is hypodissipative in $L^2(m)$ up to a
bounded term, uniformly in the background.
\end{enumerate}
\end{prop}

\begin{proof}[Comments on the proof]
Part~(i) combines the weak dissipativity of $\mathcal B$, \cite[Lemma~2.3]{CM}, with the hybrid norm
of \cite[(3.12), Prop.~3.6]{CM}, which here take the form of Proposition~\ref{prop:CMthr} and (CM1)
of Appendix~\ref{app:band}. Part~(ii) is \cite[Lemma~4.3]{CM}. We redo both for the class
$\mathcal W$. The only new point is the mass-tail term, which we bound by
Proposition~\ref{prop:Gthreshold}(i). Where \cite{CM} use \cite[Lemma~2.5]{CM}, we use the version
of Lemma~\ref{lem:RVA}. Part~(iii) is a direct computation in the Gaussian weight, using the
ellipticity of $\ab[G]$. The complete proofs are given in Proposition~\ref{prop:trilinG} and
Corollary~\ref{cor:W1G}, and Lemmas~\ref{lem:GW} and~\ref{lem:Whypo}.
\end{proof}

The \emph{analytic leaf energy} at radius $\rad$ is
\begin{equation}\label{eq:Aenergy}
\hat{\mathbb A}_\xi(t)=\sum_{n\ge0}w_n(\rad)^2\sum_{|\beta|=n,\ |\alpha|\le
N_y}\delta_y^{2|\alpha|}\,\big\|T_p^{-1}\big(\partial_y^\alpha\partial_v^\beta(g-
\mathcal M[g])\big)^{(t\xi)}\big\|^2_{X(\hat m_\xi)},\qquad w_n(\rad)=\frac{\rad^n(n+1)^5}{n!},
\end{equation}
normalized by the leaf mass. We write $\hat{\mathbb A}_{\ge1,\xi}$ for the part with
$n+|\alpha|\ge1$. By Lemma~\ref{lem:algebra}, the weights $(n+1)^5$ make the analytic norms an
algebra, and by Lemma~\ref{lem:TWcalc}(iv) and Lemma~\ref{lem:TWabs}, the shifted products are
absorbed by the radius.

\subsection{The macroscopic part: a leaf Euler system}\label{ssec:macro}

On each leaf we split the solution into a Maxwellian with the leaf moments and a moment-free micro
part. Then the leaf parameters solve a compressible Euler system, slowed down by $t^{-2}$, in
which the micro part contributes a stress and a heat flux and nothing else.

\subsubsection{The local leaf Maxwellian}

\begin{defin}[Local leaf Maxwellian]\label{def:Mloc}
Let $\mathcal M[g](t,y,v):=M_{p_g(t,\xi)}(y)$ with $\xi=v+y/t$, and $\hat r:=g-\mathcal M[g]$.
\end{defin}

Note that on the leaf through $x$, $\xi=x/t$ is constant, so $\mathcal M[g]^{(x)}=M_{p_g(t,x/t)}$
is an exact Maxwellian in $z$ with the leaf moments of $g$.

\begin{lem}[The micro part is moment-free]\label{lem:Mfree}
For all $\alpha,\beta$, $t$ and $x$, the leaf moments of $\partial^\alpha_y\partial^\beta_v\hat r$
vanish.
\end{lem}

\begin{proof}
By construction, we have $\int\hat r^{(t\xi)}q\dd z=0$ for all $q\in\mathrm{span}\{1,z_i,|z|^2\}$
and all $\xi$. By Lemma~\ref{lem:leafgeom}(ii), we have $(\partial_vh)^{(x)}=\nabla_\xi h^{(x)}$ and
$(\partial_yh)^{(x)}=\nabla_zh^{(x)}+t^{-1}\nabla_\xi h^{(x)}$. Thus the leaf moments of
$\partial^\alpha_y\partial^\beta_v\hat r$ are linear combinations of $\xi$-derivatives of moments of
$\hat r$ against polynomials in the same span, which is closed under differentiation in $z$. Hence
all of them vanish.
\end{proof}

\begin{lem}[Equation for $\hat r$]\label{lem:Meq}
Let $\cQ_t(h_1,h_2)(y,v):=Q\big(h_1^{(x)},h_2^{(x)}\big)(y)$ with $x=y+tv$, the polarization of
$\cQ_t[g]$, and $\mathcal L_t[g]h=\cQ_t(
\mathcal M[g],h)+\cQ_t(h,\mathcal M[g])$. Then
\begin{equation}\label{eq:Meq}
\partial_t\hat r=\kappa\big(\mathcal L_t[g]\hat r+\cQ_t(\hat r,\hat r)\big)-\partial_t
\mathcal M[g],\qquad\partial_t\mathcal M[g]=\sum_k(\partial_{p_k}M)_{p_g(t,\xi)}(y)\big[D_tp_{g,k}-t^{-2}y\cdot\nabla_\xi
p_{g,k}\big](t,\xi).
\end{equation}
On each leaf the functions $\partial_{p_k}M_{p_g}$ span the kernel of the linearized operator.
Therefore the micro part of the forcing $-\partial_t\mathcal M[g]$ is
\begin{equation}\label{eq:Meq-force}
t^{-2}(I-\mathcal P_t)\Big[\sum_k(z\cdot\nabla_\xi p_{g,k})\,\partial_{p_k}M_{p_g}\Big],
\end{equation}
which is of size $O(t^{-2})$. Here $\mathcal P_t$ is the projection onto this kernel on each leaf
(in the frame, the projection $\Pi_0$ onto $\ker\mathcal L$).
\end{lem}

\begin{proof}
Since $\cQ_t(\mathcal M,\mathcal M)=0$ leafwise, we have $\cQ_t(g,g)=\mathcal L_t[g]\hat
r+\cQ_t(\hat r,\hat r)$. At fixed $(y,v)$, $\xi=v+y/t$ satisfies $\partial_t\xi=-y/t^2$, which
implies the formula for $\partial_t\mathcal M$ in \eqref{eq:Meq}. The coefficients $D_tp_g$ are
constant on each leaf, so their contribution lies in the kernel.
\end{proof}

For this reason we linearize around the local leaf Maxwellian. Around the fibered flow the forcing
would be $O(t^{-1})$, which is not integrable for the analytic energies, while around the local
leaf Maxwellian it is $O(t^{-2})$, by \eqref{eq:Meq-force}.

\subsubsection{The leaf Euler system}

For the normalized micro part $\psi_\xi=T^{-1}_{p_g(t,\xi)}\hat r^{(t\xi)}$ we define the traceless
stress and the heat flux
\begin{equation}\label{eq:stress}
s_{ij}[\psi]=\int\psi\,\hat z_i\hat z_j\dd\hat z,\qquad h_i[\psi]=\int\psi\,\hat z_i|\hat
z|^2\dd\hat z .
\end{equation}

\begin{theo}[Leaf Euler system]\label{thm:Meuler}
Along rays, with $\nabla=\nabla_\xi$, $s=s[\psi_\xi]$ and $h=h[\psi_\xi]$,
\begin{equation}\label{eq:Euler}\begin{aligned}
D_t\ell&=t^{-2}\big(\nabla\cdot y_0+y_0\cdot\nabla\ell\big),\\
D_ty_{0,j}&=t^{-2}\Big((y_0\cdot\nabla)y_{0,j}+e^{\vartheta}\partial_j(\vartheta+\ell)+e^{\vartheta}\big(\partial_is_{ij}+s_{ij}\partial_i(\ell+\vartheta)\big)\Big),\\
D_t\vartheta&=t^{-2}\Big(y_0\cdot\nabla\vartheta+\tfrac23\nabla\cdot
y_0+\tfrac23s_{ij}\partial_iy_{0,j}+\tfrac13e^{\vartheta/2}\big(\nabla\cdot
h+h\cdot\nabla(\ell+\tfrac32\vartheta)\big)\Big).
\end{aligned}\end{equation}
\end{theo}

\begin{proof}
\emph{Step 1: the conservative form.} By Theorem~\ref{thm:momtransport}, we have
\begin{equation}\label{eq:Euler-cons}
D_tm_g=t^{-2}\,\nabla\cdot\int z\otimes\mom\,g^{(x)}\dd z .
\end{equation}
We split $g^{(x)}=M_{p_g}+\hat r^{(x)}$.

\emph{Step 2: the fluxes.} The Maxwellian fluxes are $\rho y_{0,i}$,
$\rho(y_{0,i}y_{0,j}+\theta\delta_{ij})$ and $\rho y_{0,i}(|y_0|^2+5\theta)$. By
Lemma~\ref{lem:Mfree} and the change of variables $z=y_0+\theta^{1/2}\hat z$, we see that the
micro fluxes are $0$, $\rho\theta s_{ij}$ and $\rho(2\theta y_{0,j}s_{ij}+\theta^{3/2}h_i)$.

\emph{Step 3: the log-variables.} Solving \eqref{eq:Euler-cons} for $(D_t\ell,D_ty_0,D_t\vartheta)$,
we obtain \eqref{eq:Euler}. The computation is elementary, and a short symbolic check is included as
an ancillary file with the arXiv version.
\end{proof}

Without the micro fluxes, this is the compressible Euler system for a perfect monatomic
gas. It is written in the log-variables, in the moving frame of free expansion, and it is slowed
down by $t^{-2}$. In the \emph{Euler time}
\begin{equation}\label{eq:Eulertime}
\sigma=T_0^{-1}-t^{-1}
\end{equation}
the system is autonomous, and \emph{the whole evolution lasts the finite Euler time $1/T_0$}. The
system depends on the density through $\nabla\ell$ alone. Hence the amplitude $\eps$ does not
appear, and the positivity of $\rho$ and $\theta$ is automatic.

\begin{lem}[The reference Euler flow]\label{lem:Mref}
Suppose that $\varpi'_{\rm ref}=\varpi'_g(T_0,\cdot)$ extends holomorphically to $\mathbb S_{R_1}$
with $|\varpi'_{\rm ref}|\le L_1$ there (this is Corollary~\ref{cor:LWPa}(a4)). Put $\ell_{\rm
ref}:=\ell_g(T_0,\cdot)$. On $\mathbb S_{R_1}$,
\begin{equation}\label{eq:ellref}
\ell_{\rm ref}=\log\eps+\log\rho_0+\operatorname{Log}(\rho_g/(\eps\rho_0))(T_0),
\end{equation}
where $\log\rho_0$ is the logarithm on the convex strip $\mathbb S_\Lambda\supset\mathbb S_{R_1}$
that is real on $\R^3$ ((D3)), and $\operatorname{Log}$ is the principal branch
($\mathrm{Re}\,\rho_g/(\eps\rho_0)\ge\frac12$ by the proof of Corollary~\ref{cor:LWPa}(a4)).
Therefore, $\ell_{\rm ref}$ is the holomorphic extension of the real function $\log\rho_g(T_0)$, and
$\nabla\ell_{\rm ref}$ is the first component of $\varpi'_{\rm ref}$. There is $T_{\rm E}(L_1,R_1)$
such that for $T_0\ge T_{\rm E}$ the system of Theorem~\ref{thm:Meuler} with $s=h=0$, written for
$u=(\ell-\ell_{\rm ref},y_0,\vartheta)$ as in the proof, with $u(0)=(0,y_{0,g},\vartheta_g)(T_0)$,
has a unique solution $u=(u_1,y_{0*},\vartheta_*)$ on $[T_0,\infty)$, holomorphic and bounded on
$\mathbb S_{R_1/2}$ together with $\nabla u$. Put $\varpi'_*:=(\nabla\ell_{\rm ref}+\nabla
u_1,y_{0*},\vartheta_*)$. Then
\begin{equation}\label{eq:Mref-bd}
\sup_{\mathbb S_{R_1/2}}|\varpi'_*(t)-\varpi'_{\rm
ref}|\le\frac{C(L_1,R_1)}{T_0},\qquad\sup_{
\mathbb S_{R_1/2}}|\varpi'_*(t)-\varpi'_*(\infty)|\le\frac{C(L_1,R_1)}t .
\end{equation}
We take $T_{\rm E}$ also with $C(L_1,R_1)/T_{\rm E}\le\frac1{10}$, so that $|\varpi'_*|\le L_1+1$ on
$\mathbb S_{R_1/2}$.

\emph{Normalization of $\ell_*$.} Put $\ell_*:=\ell_{\rm ref}+u_1$, $\rho_*:=e^{\ell_*}$ and
$\varpi_*:=(\ell_*,y_{0*},\vartheta_*)$. Then $\nabla\ell_*$ is the first component of $\varpi'_*$
(by definition), $\varpi_*(T_0)=\varpi_g(T_0)$, and $\sup_{\mathbb S_{R_1/2}}|\ell_*(t)-\ell_{\rm
ref}|\le C(L_1,R_1)/T_0$. Since $\mathcal F$ has real coefficients and $\nabla\ell_{\rm ref}$,
$u(0)$ are real on $\R^3$, uniqueness applied to $\overline{u(\sigma,\bar\xi)}$ shows that $u$ is
real on $\R^3$. Hence so are $\varpi'_*$ and $\ell_*$, and $\rho_*>0$ there. Any solution
$(q,y_0,\vartheta)$ of the $\xi$-differentiated system for $\varpi'$ with the data $\varpi'_{\rm
ref}$, holomorphic and bounded on $\mathbb S_{R_1/2}$ together with $\nabla y_0$, satisfies
$q=\nabla(\ell_{\rm ref}+v_1)$ with $v_1:=\int_0^\sigma(\nabla\cdot y_0+y_0\cdot q)$, and
$(v_1,y_0,\vartheta)$ solves the system for $u$. Uniqueness holds in this class, and in the class of
the statement, with no bound on the size of the solutions. Indeed, let two solutions of the system
for $u$, continuously differentiable in $\sigma$ with values in the bounded holomorphic functions on
$ \mathbb S_{R_1/2}$, be bounded by $R'$ there and agree at some $\sigma_1$. The uniqueness part of
\cite[Theorem]{Nis77}, applied with $\rho_0=R_1/2$ and the ball of radius $3R'$ to their differences
from the common value at $\sigma_1$ (the Cauchy estimates yield its hypothesis on $F(0,t)$),
provides a time $a(R',L_1,R_1)>0$ such that they agree on $\mathbb S_{R_1/4}$ for
$\sigma\in[\sigma_1,\sigma_1+a R_1/4]$, and hence on $\mathbb S_{R_1/2}$ by the identity theorem.
Finitely many such steps cover $[0,1/T_0]$, so $\varpi'_*$ is the unique such solution.
\end{lem}

\begin{proof}
\emph{Step 1: the system as an abstract Cauchy problem.} In the Euler time \eqref{eq:Eulertime} the
system involves $\ell$ through $\nabla\ell$ alone, and it reads
\begin{equation}\label{eq:Mref-F}
\partial_\sigma u=\mathcal F(u,\nabla u;\nabla\ell_{\rm ref}),\qquad u=(\ell-\ell_{\rm
ref},y_0,\vartheta),\qquad u(0)=(0,y_{0,g},\vartheta_g)(T_0).
\end{equation}
Here $\mathcal F$ is entire in $u$, it is affine in $\nabla u$, with coefficients entire in $u$
(this is the leaf Euler system without the micro fluxes), and it is affine in the
coefficient $\nabla\ell_{\rm ref}$. The affinity in $\nabla u$ yields the Lipschitz condition of
\cite{Ni72,Nis77} with a single factor $(s-s')^{-1}$. On the scale of spaces $X_s$ of bounded
holomorphic functions on $\mathbb S_s$, we have
\begin{equation}\label{eq:Mref-scale}
\|\partial_\xi u\|_{X_{s'}}\le\frac{C\|u\|_{X_s}}{s-s'} .
\end{equation}

\emph{Step 2: existence on $[0,1/T_0]$.} The abstract Cauchy--Kovalevskaya theorem, in the form of
\cite[Theorem]{Nis77}, which contains \cite[Thm.~1.1]{Ni72}, applied to $u-u(0)$ with $\rho_0=R_1$
(its hypothesis on $F(0,t)$ holds since $\|\nabla u(0)\|_{X_s}\le CL_1/(R_1-s)$ by the Cauchy
estimates), yields a unique solution for $\sigma\in[0,c(R_1-s)/C(L_1,R_1))$ in $X_s$. With
$s=\frac78R_1$, this interval contains $[0,1/T_0]$ when $T_0\ge T_{\rm E}$.

\emph{Step 3: the bounds.} By the Cauchy estimates on discs of radius $R_1/8$, we have $|\nabla
u|+|\partial_\sigma u|\le C(L_1,R_1)$ on $\mathbb S_{3R_1/4}$ and $|\partial_\sigma\nabla u|\le
C(L_1,R_1)$ on $\mathbb S_{R_1/2}$. The bounds \eqref{eq:Mref-bd} follow, since
$\sigma\in[0,T_0^{-1}]$ and $\sigma_\infty-\sigma=1/t$. This completes the proof.
\end{proof}

Then we control the full leaf parameters by an analytic hierarchy for the deviation
\begin{equation}\label{eq:dvarpi}
\delta\varpi:=\varpi_g-\varpi_*,\qquad\varpi=(\ell,y_0,\vartheta),
\end{equation}
in which the micro fluxes $s,h$ act as sources, as in Proposition~\ref{prop:UFforce}(f) and in
(M1)--(M3) of its proof. Differentiating in $\xi$ loses one derivative per order, and we absorb
the loss as in Lemma~\ref{lem:analytic}, with $b(t)=t^{-2}$. We note three facts about
\eqref{eq:dvarpi}. Since $\ell_*$ is normalized as in Lemma~\ref{lem:Mref}, we have
$\delta\varpi(T_0)=0$. The function $\ell_*$ itself is unbounded. However, it enters only through
$\nabla\ell_*$ and through $\delta\ell=\ell_g-\ell_*$, namely in $\rho_*/\rho_g=e^{-\delta\ell}$
(see Lemma~\ref{lem:UFref}), in $u_*$ and $U_*$ of Lemma~\ref{lem:TWcoef} and \eqref{eq:TWdtheta}
(where $u_{*,1}(\xi)=\ell_*(\xi)-\ell_*(\xi_j)-\delta\ell(\xi_j)$), and, at level $0$, in
$\mathfrak m_Q$ and $\mathfrak M_c$. Finally, under (U1), we have
$|\delta\varpi|\le C'_{\rm Sob}\delta_U\le\frac14$, by Lemma~\ref{lem:TWcalc}(i).

\subsection{Trilinear estimates in Gaussian weights}\label{ssec:tri}

In this subsection we prove the bilinear and trilinear estimates of \cite{CM} in the Gaussian leaf
weights, which contain Proposition~\ref{prop:triG}. Proposition~\ref{prop:TLsrc} provides the source
estimate on a single leaf that is used in Proposition~\ref{prop:Bband}.

Compared with the polynomial weights of \cite{CM}, the estimates differ as follows. First, the
$\eps$-defect $\eps^{|\gamma|/2}\|\cdot\|_X$ of the polynomial estimates disappears, since the
zero-order part $\jap{\hat z}^{(\gamma+2)/2}m$ of $Y_2(m)$ controls every term generated by
$\nabla\log m$. Moreover, the mass tail needs no separate functional. It is bounded by
$Y_2\times Y_2$, so $\langle Q(f,g),h\rangle_{L^2(m)}$ is bounded by one energy-type and two
dissipative norms in every configuration. Finally, on a single leaf all factors live at the same
$\xi$. Then the summations need only Lemma~\ref{lem:algebra}, with no Sobolev embedding in $\xi$, so
no extra levels are used.

\subsubsection{The weight class}

For $0<q_-\le q_+\le1/8$, $\eps_1\ge0$ and $\omega_0\ge0$, let $\mathcal W(q_\pm,\eps_1)$ be the set
of weights on $\R^3_{\hat z}$ of the form
\begin{equation}\label{eq:Wclass}
m=e^{\Phi},\qquad\Phi=q|\hat z|^2+\omega,\qquad q\in[q_-,q_+],\qquad|\nabla\omega|\le\eps_1,\qquad|\omega(0)|\le\omega_0
.
\end{equation}
The parameter $\omega_0$ is fixed by Lemma~\ref{lem:TG0}(i) and is not displayed, and all
constants of this subsection may depend on it. By the normalization of $\omega(0)$, we have the
two-sided bound
\begin{equation}\label{eq:Wclass-2side}
e^{-\omega_0}e^{q_-|\hat z|^2-\eps_1|\hat z|}\le m\le e^{\omega_0}e^{q_+|\hat z|^2+\eps_1|\hat z|},
\end{equation}
We need this bound for the estimates that are not homogeneous of degree $2$ in $m$. We need no
bound on $\nabla^2\omega$, and no smallness of $\eps_1$, except in the coercivity of
Lemma~\ref{lem:GW}(i). Recall the norms
\begin{equation}\label{eq:Y2norm}
\|u\|^2_{Y_2(m)}=\|\jap{\hat z}^{(\gamma+2)/2}mu\|^2_{L^2}+\|\jap{\hat
z}^{\gamma/2}\tn(mu)\|^2_{L^2},\qquad\|u\|_\flat=\|\jap{\hat z}^4u\|_{H^1}.
\end{equation}

\begin{lem}[Membership and comparisons]\label{lem:TG0} 
\begin{enumerate}[label=(\roman*)]
\item Let $q\theta_+\le1/8$ and $t\ge1$. The normalized leaf weights of Lemma~\ref{lem:Gleafweight}
belong to the class $\mathcal W(q\theta_-,q\theta_+,\eps_1(t))$, with
\begin{equation}
\eps_1(t):=2\theta_+^{1/2}L(q+t^{-1})\le2\theta_+^{1/2}L(q+1),\qquad\omega_0:=(q+2)(L+1)^2 .
\end{equation}
The same holds for the weights $\hat m_\xi$, in the frames $p_g(t,\xi)$ and with the factor
$A_{\rm ref}$, with $(L,\theta_-,\theta_+)$
replaced by $(L+1,\theta_-/2,\theta_++1)$, the bounds that (F3) gives for $p_g$, and with $\omega_0$
replaced by a constant $\omega_0'=C(L,\theta_\pm)$. For these weights also $|\nabla^2\omega|\le
C(L_1,R_1,\theta_\pm)t^{-2}$ and $e^{\pm\omega}\le C\jap{\hat z}^pe^{Cq\jap{\hat z}}$, with $p$ the
exponent of (D3). Thus (F2) and Lemma~\ref{lem:GW} apply to $\hat m_\xi$ with these constants (only
the bound on $\nabla^2\omega$ depends on $L_1,R_1$, which are fixed before $T_1$, and it is used
only in the alternative argument for Lemma~\ref{lem:GW}(ii), while $T_1$ is fixed by the bound on
$\nabla\omega$).
\item For $m\in\mathcal W$ and $\theta,k\ge0$,
\begin{equation}
\|u\|_{L^2(\jap{\hat z}^\theta)}+\|u\|_{X_k}\le C\|u\|_{X(m)},\qquad\|u\|_{H^1(\jap{\hat
z}^\theta)}+\|u\|_\flat+\|u\|_{Y_k}\le C\|u\|_{Y_2(m)} .
\end{equation}
\item If $P$ is a polynomial of degree $\le3$ with coefficients $c$, then
$\|P\mu\|_{X(m)}+\|P\mu\|_{Y_2(m)}+\|P\mu\|_\flat\le C|c|$.
\end{enumerate}
The constants depend only on $(\gamma,q_\pm,\eps_1,\omega_0,\theta,k)$.
\end{lem}

\begin{proof}
(i) \emph{The static weights.} By Lemma~\ref{lem:Gleafweight}, we have
$|\nabla\omega_x|\le\eps_1(t)$. At $\hat z=0$ we have $z=y_{0\infty}$, so
$\omega_x(0)=q|y_{0\infty}|^2+\log\big(A_\infty(\xi-y_{0\infty}/t)/A_\infty(\xi)\big)$. Since
$|y_{0,\infty}|\le L/2$ by Subsection~\ref{ssec:disp}, and $|\nabla\log A_\infty|\le2L$, we have
$0\le q|y_{0\infty}|^2\le qL^2/4$, and the logarithm does not exceed $L^2/t$ in absolute value.
Hence $|\omega_x(0)|\le\omega_0$.

\emph{The weights $\hat m_\xi$.} In the frames $p_g$, we have $\hat m_\xi=e^{q|z|^2}A_{\rm
ref}(\xi-z/t)/A_{\rm ref}(\xi)$ with $z=y_{0,g}+\theta_g^{1/2}\hat z$. The same computation,
with the center $y_{0,g}$, the bounds $|y_{0,g}|\le L+1$, $\theta_g\in[\frac12\theta_-,\theta_++1]$
of (F3) and the bound $|\nabla\log A_{\rm ref}|\le L$ of Corollary~\ref{cor:LWPa}(a4), yields
$|\nabla\omega|\le2(\theta_++1)^{1/2}(L+1)(q+t^{-1})$ and a constant $\omega_0'=C(L,\theta_\pm)$.

\emph{Second derivatives and pointwise bounds.} We have
$|\nabla^2\omega|=t^{-2}\theta_g|(\nabla^2\log A_{\rm ref})(\xi-z/t)|\le
C(L_1,R_1,\theta_\pm)t^{-2}$, since $\nabla\log A_{\rm ref}=-\nabla\ell_g+\frac34\nabla\vartheta_g$
at $T_0$ and the Cauchy estimates on $\mathbb S_{R_1}$ bound $\nabla^2\ell_g$ and
$\nabla^2\vartheta_g$ on $\R^3$ by $C(L_1,R_1)$, by the same corollary. The bound
on $e^{\pm\omega}$ follows as for the static weights, since
$\rho_g(T_0)/(\eps\rho_0)\in[\frac12,2]$ and $\theta_g(T_0)\in[\theta_-,\frac23\theta_+]$ on $\R^3$,
as in the proof of that corollary.

(ii) We have $m\ge e^{-\omega_0}e^{q_-|\hat z|^2-\eps_1|\hat z|}\ge c_N\jap{\hat z}^N$ for every
$N$, with $c_N=c_N(q_-,\eps_1,\omega_0)$. For the gradient,
\begin{equation}
m\nabla u=\nabla(mu)-(\nabla\Phi)mu,\qquad|\nabla\Phi|\le2q_+|\hat z|+\eps_1,
\end{equation}
which implies that $\jap{\hat z}^\theta|\nabla u|\le C_N\jap{\hat z}^{-N}\big(|\nabla(mu)|+\jap{\hat
z}|mu|\big)$. Since $|\nabla v|\le|\tn v|$, this is bounded by the two parts of $Y_2(m)$.

(iii) We have $m\,P\mu\le e^{\omega_0}e^{(q_+-1/2)|\hat z|^2+\eps_1|\hat z|}|P|$, and $q_+\le1/8$.
\end{proof}

\subsubsection{A bilinear estimate in Gaussian weights}

\begin{lem}[Bilinear estimate in Gaussian weights]\label{lem:QmG}
Let $m\in\mathcal W(q_\pm,\eps_1)$, $\theta>7/2$, $\theta'>9/4$. Then
\begin{equation}
|\langle Q(f,g),h\rangle_{L^2(m)}|\le C\Big(\|f\|_{L^2(\jap{\hat
z}^\theta)}\|g\|_{Y_2(m)}+\|f\|_{H^1(\jap{\hat z}^{\theta'})}\|g\|_{X(m)}\Big)\|h\|_{Y_2(m)},
\end{equation}
with $C=C(\gamma,q_+,\eps_1,\theta,\theta')$ and no $\eps$-defect.
\end{lem}

\begin{proof}
\emph{Step 1: decomposition.} Put $G=mg$ and $H=mh$. Integrating by parts as in
\cite[Lemma~4.3]{CM}, and using $\partial_i(m^{-1}G)=m^{-1}(\partial_iG-\partial_i\Phi\,G)$ and
$\partial_j(mH)=m(\partial_jH+\partial_j\Phi\,H)$, we have
\begin{equation}\label{eq:QmG-dec}
\langle Q(f,g),h\rangle_{L^2(m)}=A_1+A_2+A_3+A_4+B_1+B_2,
\end{equation}
where
\begin{equation}\begin{aligned}
A_1&=-\!\int(a*f)\nabla G\cdot\nabla H, & A_2&=-\!\int(a*f)\nabla G\cdot\nabla\Phi\,H, &
A_3&=\int(a*f)\nabla\Phi\cdot\nabla H\,G,\\
A_4&=\int(a*f)\nabla\Phi\cdot\nabla\Phi\,GH, & B_1&=\int(b*f)\cdot\nabla H\,G, &
B_2&=\int(b*f)\cdot\nabla\Phi\,GH .
\end{aligned}\end{equation}
Only first derivatives of $\Phi$ occur.

\emph{Step 2: coefficient bounds.} By \cite[Lemma~4.2(i)]{CM} with $p=2$, with the unit vector
$e:=\hat z/|\hat z|$, the quantities $(a*f)e\cdot e$, $(a*f)e$ and $a*f$ are bounded by $C\jap{\hat
z}^{\gamma}$, $C\jap{\hat z}^{\gamma+1}$ and $C\jap{\hat z}^{\gamma+2}$ times $\|f\|_{L^2(\jap{\hat
z}^\theta)}$. Hence we have
\begin{equation}\label{eq:QmG-coef}
|(a*f)\xi\cdot\zeta|\le C\|f\|_{L^2(\jap{\hat z}^\theta)}\jap{\hat
z}^\gamma|\tilde\xi||\tilde\zeta|,\qquad\tilde\xi=P_{\hat z}\xi+\jap{\hat z}(I-P_{\hat z})\xi .
\end{equation}
For $\nabla\Phi=2q\hat z+\nabla\omega$ we have $|\widetilde{\nabla\Phi}|\le2q|\hat
z|+\eps_1\jap{\hat z}\le C\jap{\hat z}$. By \cite[Lemma~4.2(ii)]{CM} with $p=4$ and Sobolev's
inequality, we get $|b*f|\le C\jap{\hat z}^{\gamma+1}\|f\|_{H^1(\jap{\hat z}^{\theta'})}$.

\emph{Step 3: the six terms.} We write $\|f\|$ for $\|f\|_{L^2(\jap{\hat z}^\theta)}$.
\begin{itemize}
\item $|A_1|\le C\|f\|\int\jap{\hat z}^\gamma|\tn G||\tn H|\le C\|f\|\,\|g\|_{Y_2}\|h\|_{Y_2}$.
\item $|A_2|+|A_3|\le C\|f\|\int\jap{\hat z}^{\gamma+1}\big(|\tn G||H|+|G||\tn H|\big)$, and we
split $\jap{\hat z}^{\gamma+1}=\jap{\hat z}^{\gamma/2}\jap{\hat z}^{(\gamma+2)/2}$.
\item $|A_4|\le C\|f\|\int\jap{\hat z}^{\gamma+2}|G||H|$, which is the product of the zero-order
parts of $Y_2$.
\item $|B_1|\le C\|f\|_{H^1}\int\jap{\hat z}^{\gamma+1}|G||\nabla H|\le C\|f\|_{H^1}\|\jap{\hat
z}^{(\gamma+2)/2}G\|\,\|h\|_{Y_2}\le C\|f\|_{H^1}\|g\|_X\|h\|_{Y_2}$, using $\jap{\hat
z}^{(\gamma+2)/2}\le1$ since $\gamma\le-2$.
\item $|B_2|\le C\|f\|_{H^1}\int\jap{\hat z}^{\gamma+2}|G||H|\le C\|f\|_{H^1}\|g\|_X\|h\|_{Y_2}$,
by the same argument.
\end{itemize}
\end{proof}

The analogue of this lemma in polynomial weights, which is based on \cite[Lemma~4.3]{CM},
needs an $\eps$-defect, because there $Y$ has only the zero-order weight $\jap{\hat z}^{\gamma/2}$,
while $\nabla\log m$ carries the tangential size $\eps\jap{\hat z}$. In $Y_2(m)$ the zero-order
weight is $\jap{\hat z}^{(\gamma+2)/2}$, which matches $|\widetilde{\nabla\Phi}|\lesssim\jap{\hat
z}$ exactly. For the centered weights $m=\jap v^k$ of Appendix~\ref{app:CM}, $\nabla\log m$ is
radial and of size at most $k\jap v^{-1}$, and Lemma~\ref{lem:CMbil} has no defect.

\subsubsection{The trilinear proposition}

\begin{prop}[Trilinear bounds in Gaussian weights]\label{prop:trilinG}
Let $m\in\mathcal W(q_\pm,\eps_1)$, and write $X=X(m)$ and $Y_2=Y_2(m)$. For all $f,g,h$:
\begin{enumerate}[label=(\roman*)]
\item
\begin{equation}\label{eq:NtriG-i}
\big|\langle Q(f,g),h\rangle_{L^2(m)}-\nu(f)\,\mathfrak b_m(g,h)\big|\le
C\|f\|_\flat\|g\|_{Y_2}\|h\|_{Y_2};
\end{equation}

\item $|\mathfrak b_m(g,h)|\le C\|g\|_{Y_2}\|h\|_{Y_2}$;
\item hence $|\langle Q(f,g),h\rangle_{L^2(m)}|\le C\|f\|_\flat\|g\|_{Y_2}\|h\|_{Y_2}$, and also
\begin{equation}
|\langle Q(f,g),h\rangle_{L^2(m)}|\le C\big(\|f\|_X\|g\|_{Y_2}+\|f\|_{Y_2}\|g\|_X\big)\|h\|_{Y_2};
\end{equation}
\item
\begin{equation}\label{eq:NtriG-nat}
|\ip h{Q(f,g)+Q(g,f)}_\natural|\le C\big(\|f\|_{X}\|g\|_{Y_2}+\|f\|_{Y_2}\|g\|_{X}\big)\|h\|_{Y_2};
\end{equation}

\item if $P\in\ker\mathcal L$, i.e.\ $P$ is a Gaussian polynomial of the form $(c_0+c_1\cdot\hat
z+c_2|\hat z|^2)\mu$, then
\begin{equation}
|\langle Q(f,P)+Q(P,f),h\rangle_{L^2(m)}|\le C|c|\big(\|f\|_X+\|f\|_{Y_2}\big)\|h\|_{Y_2} .
\end{equation}
\end{enumerate}
\end{prop}

\begin{proof}
(i) We write
\begin{equation}\label{eq:triG-i}
\langle Q(f,g),h\rangle_{L^2(m)}=-\!\int(a*f)\nabla
g\cdot\nabla(m^2h)+\int(b*f)\,g\cdot\nabla(m^2h).
\end{equation}
The first integral is $A_1+A_2+A_3+A_4$ of Lemma~\ref{lem:QmG}, which is bounded by
$C\|f\|_{L^2(\jap{\hat z}^4)}\|g\|_{Y_2}\|h\|_{Y_2}$. In the second one, we insert
\begin{equation}\label{eq:triG-split}
b*f=\nu(f)\bb[\mu]+b*F,\qquad F:=f-\nu(f)\mu,\qquad\nu(f)=\int f,\qquad\int F=0 .
\end{equation}
By \eqref{eq:masstail}, the first part is $\nu(f)\mathfrak b_m(g,h)$. For the second part we use
\begin{equation}\label{eq:triG-bF}
|b*F(\hat z)|\le C\jap{\hat z}^\gamma\|F\|_\flat .
\end{equation}
Indeed, by \cite[Lemma~4.2(ii)]{CM} with $p=4$ and Sobolev's inequality, we have
\begin{equation}\label{eq:triG-b1}
\int|\hat z-w|^{\gamma+1}|u(w)|\dd w\le C\jap{\hat z}^{\gamma+1}\|u\|_{H^1(\jap
w^{5/2})},\qquad\|u\|_{H^1(\jap w^{5/2})}:=\|\jap w^{5/2}u\|_{H^1},
\end{equation}
which is \eqref{eq:triG-bF} for $|\hat z|\le2$. For $|\hat z|>2$, since $\int F=0$, we have
\begin{equation}
b*F(\hat z)=\int_{|w|<|\hat z|/2}\big(b(\hat z-w)-b(\hat z)\big)F(w)\dd w+\int_{|w|\ge|\hat
z|/2}\big(b(\hat z-w)-b(\hat z)\big)F(w)\dd w .
\end{equation}
On $\{|w|<|\hat z|/2\}$, we have $|b(\hat z-w)-b(\hat z)|\le C|w||\hat z|^\gamma$, since $b$ is
homogeneous of degree $\gamma+1$. On the complement, we have $|b(\hat z)||F(w)|\le C|\hat
z|^\gamma|w||F(w)|$ and $|b(\hat z-w)||F(w)|\le4|\hat z|^{-1}|\hat z-w|^{\gamma+1}|w||F(w)|$, and
we can apply~\eqref{eq:triG-b1} with $u=|w|F$. Finally, we have $\int|w||F|+\|\jap
w^{5/2}|w|F\|_{H^1}\le C\|F\|_\flat\le C\|f\|_\flat$. This proves the claimed bound, with
$\|f\|_\flat$ on the right. Moreover, we have $\nabla(m^2h)=m\big(\nabla H+\nabla\Phi\,H\big)$ with
$|\nabla\Phi|\le C\jap{\hat z}$. Hence this part is bounded by
\begin{equation}
C\|f\|_\flat\!\int\jap{\hat z}^\gamma|G|\big(|\nabla H|+\jap{\hat z}|H|\big)\le
C\|f\|_\flat\|g\|_{Y_2}\|h\|_{Y_2},
\end{equation}
because $\jap{\hat z}^{\gamma/2}\le\jap{\hat z}^{(\gamma+2)/2}$ and $\jap{\hat
z}^{(\gamma+1)/2}\le\jap{\hat z}^{(\gamma+2)/2}$.

(ii) This is Proposition~\ref{prop:Gthreshold}(i) with $C_1=\eps_1$.

(iii) The first form follows from (i), (ii) and $|\nu(f)|\le C\|f\|_\flat$. The second form
follows from Lemma~\ref{lem:QmG} with Lemma~\ref{lem:TG0}(ii), applied once with $(f,g)$ and once
with the roles of the norms exchanged through $\|f\|_{L^2(\jap{\hat z}^\theta)}\le C\|f\|_X$ and
$\|f\|_{H^1(\jap{\hat z}^{\theta'})}\le C\|f\|_{Y_2}$.

(iv) This is (CM1)--(CM2) of Appendix~\ref{app:band} at level $7$, together with the comparison
used in~(iii).

(v) We apply (iii), second form, to $Q(f,P)$ and to $Q(P,f)$, and use Lemma~\ref{lem:TG0}(iii) for
the norms of $P$.
\end{proof}

\subsubsection{Coercivity and the linear terms in Gaussian weights}

\begin{cor}[Coercivity and linear estimates in Gaussian weights]\label{cor:W1G}
Let $m\in\mathcal W(q_\pm,\eps_1)$, and suppose the conclusion of Lemma~\ref{lem:GW}(i) holds for
$m$. Then for all $f$ with $\|f\|_X\le1$ and all $\varphi$:
\begin{itemize}
\item[(W1$_G$)]
\begin{equation}\label{eq:W1G}
\int m^2\varphi\,L_{\mu+f}\varphi\le-c_W\|\varphi\|^2_{Y_2}+C_W\|\varphi\|^2_{L^2}+C_W\big(\|f\|_X\|\varphi\|_{Y_2}+\|f\|_{Y_2}\|\varphi\|_X\big)\|\varphi\|_{Y_2};
\end{equation}

\item[(W2$_G$)] for $u\in\ker\mathcal L$,
\begin{equation}\label{eq:W1G-mom}
\Big|\int m^2\varphi\big(Q(f,u)+Q(u,f)\big)\Big|\le
C_W\big(\|f\|_X+\|f\|_{Y_2}\big)|\mom(u)|\,\|\varphi\|_{Y_2} .
\end{equation}

\end{itemize}
They hold with $Y_2$ in place of $Y$ and, unlike their analogues in polynomial weights,
\emph{without} an $\eps$-defect. In particular, they imply Lemma~\ref{lem:GW}(ii).
\end{cor}

\begin{proof}
We write $L_{\mu+f}=\mathcal B+(\mathcal A_0+M\chi_R)+L'_f$, with $\mathcal A_0=Q(\cdot,\mu)$.
\begin{itemize}
\item By Lemma~\ref{lem:GW}(i), we have $\int
m^2\varphi\,\mathcal B\varphi\le-c\|\varphi\|^2_{Y_2}$.
\item The term $\mathcal A_0\varphi$ equals $\ab[\varphi]:\nabla^2\mu-\cb[\varphi]\mu$ and is
Gaussian-localized. We have $|\ab[\varphi]|\le C\jap{\hat z}^{\gamma+2}\|\varphi\|_{L^2(\jap{\hat
z}^4)}$, and $m^2|\nabla^2\mu|\le Ce^{-|\hat z|^2/8}$ since $q_+\le1/8$. The term $\cb[\varphi]\mu$
is $-8\pi\varphi\mu$ for $\gamma=-3$, and for $\gamma>-3$ it is bounded by the
Hardy--Littlewood--Sobolev inequality. Hence $|\int m^2\varphi\,
\mathcal A_0\varphi|+M\int\chi_Rm^2\varphi^2\le C\|\varphi\|^2_{L^2(\jap{\hat z}^\theta)}$.
\item For every $\delta>0$, we have $\|\varphi\|^2_{L^2(\jap{\hat
z}^\theta)}\le\delta\|\varphi\|^2_{Y_2}+C_\delta\|\varphi\|^2_{L^2}$, because $\jap{\hat
z}^{2\theta}\le\delta\jap{\hat z}^{\gamma+2}m^2$ outside a ball $B_{R_\delta}$.
\item $L'_f\varphi=Q(f,\varphi)+Q(\varphi,f)$ is bounded by Proposition~\ref{prop:trilinG}(iii),
second form.
\end{itemize}
Choosing $\delta$ in the third item with $C\delta\le c/2$, where $c$ is the constant of the first
item, we obtain (W1$_G$) with $c_W:=c/2$, half the coercivity constant. From now on this is the
constant $c_W$ of (W1$_G$), in particular in Lemma~\ref{lem:GW}(ii). Moreover, (W2$_G$) is part~(v)
of the same proposition, since $u=(c_0+c_1\cdot\hat
z+c_2|\hat z|^2)\mu$ and its moment vector $\int u\,\mom$ determine each other through the
invertible Gram matrix of $(1,\hat z,|\hat z|^2)$ in $L^2(\mu)$, so that $|c|\le C|\mom(u)|$.
\end{proof}

\subsubsection{The leafwise source estimate}

To bound the source terms of the analytic hierarchy on a single leaf, we work in the leafwise
formulation of Section~\ref{sec:boot}, around the local leaf Maxwellian.
We fix a leaf $\xi$ and write $\varphi_{\alpha\beta}=T^{-1}_{p_g(t,\xi)}(\hat
r_{\alpha\beta})^{(t\xi)}$. By Lemma~\ref{lem:Mfree}, these functions are moment-free. The Leibniz
sources are
\begin{equation}\label{eq:Leibsrc}
\Sigma_{\alpha\beta}=\sum_{0<(\alpha',\beta')\le(\alpha,\beta)}\tbinom\alpha{\alpha'}\tbinom\beta{\beta'}\Big[Q(
\mathfrak g_{\alpha'\beta'},\varphi_{\alpha''\beta''})+Q(\varphi_{\alpha''\beta''},
\mathfrak g_{\alpha'\beta'})\Big]+N_{\alpha\beta},
\end{equation}
with $\mathfrak g_{\alpha'\beta'}=T^{-1}(\partial^{\alpha'}_y\partial^{\beta'}_v
\mathcal M[g])^{(t\xi)}$ and $N_{\alpha\beta}$ the Leibniz sum of $Q(\varphi_{\cdot},\varphi_{\cdot})$. All
of these are functions on the leaf $\xi$.

We define the leaf analytic energy and dissipation by
\begin{equation}\label{eq:leafAD}
\mathbb A_\xi=\sum_nw_n^2\!\sum_{I_n}\delta_y^{2|\alpha|}\,e_\xi[\varphi_{\alpha\beta}],\qquad
\mathbb D_\xi=\sum_nw_n^2\!\sum_{I_n}\delta_y^{2|\alpha|}\,\lambda_\xi\|\varphi_{\alpha\beta}\|^2_{Y_2}.
\end{equation}

\begin{prop}[Leafwise sources]\label{prop:TLsrc}
Let $\rad\in[\rad_0/2,\rad_0]$ with $\rad_0\le r_{\rm C}/8$, $r_{\rm C}=c_{\rm
C}R_1/(L_1+1)$ the Cauchy radius of Lemma~\ref{lem:UFref}. On every leaf, the leaf weight lies in
$\mathcal W(q_\pm,\eps_1)$ (Lemma~\ref{lem:TG0}(i)). Assume (F3), requirement~\ref{rq:TWa} and
$\|\delta\varpi\|_{\mathcal A}\le2\delta_U$. Then
\begin{equation}\label{eq:TLsrc}
\sum_nw_n^2\!\sum_{I_n}\delta_y^{2|\alpha|}\,2\lambda_\xi\big|\mathcal S^{\rm
src}_\xi(\varphi_{\alpha\beta},\Sigma_{\alpha\beta})\big|\le
C\Big(\frac{\rad_0}{R_1}+\delta_y+\|\delta\varpi\|_{\mathcal A}+\mathbb A_\xi^{1/2}+
\mathbb A_\xi\Big)\mathbb D_\xi .
\end{equation}
Here
\begin{equation}
\|\delta\varpi\|_{\mathcal A}:=\Big(\sum_{n\ge0}w_n(\rad)^2\sum_{|\beta|=n}\big|\partial^\beta_\xi\delta\varpi(\xi)\big|^2\Big)^{1/2}
\end{equation}
is the analytic norm, at $\xi$, of the deviation of the leaf parameters from the reference Euler
flow. It is the norm $\|wA\|_{\ell^2}$ of Lemma~\ref{lem:algebra} with
$a_\beta=|\partial^\beta_\xi\delta\varpi(\xi)|$, and it enters here only as a quantity. In the crude
leaf inequality the bound $\|\delta\varpi\|_{\mathcal A}\le2\delta_U$
is a hypothesis. However, in the cell inequality we treat the $\delta\varpi$-part instead through
the remark after Proposition~\ref{prop:CUtri} and (H$_d$). There is no
mass-tail functional and no macroscopic term.
\end{prop}

\begin{proof}
The local part of $\mathcal S^{\rm src}_\xi$ is $\eta\langle\varphi,\Sigma\rangle_{L^2(\hat m)}$,
and the nonlocal part is $\ip\varphi\Sigma_\natural$. We bound both by
Proposition~\ref{prop:trilinG}, namely by (iii) for the local part and by (iv) for the nonlocal
part. Each bound is the product of one energy-type norm and two dissipative norms, all on the leaf
$\xi$.

\emph{Step 1: the linear couplings.} By part~(iii) of that proposition, first form, we have
\begin{equation}
|\langle Q(\mathfrak g,\varphi_{\alpha''\beta''}),\varphi_{\alpha\beta}\rangle|\le C\|
\mathfrak g\|_\flat\|\varphi_{\alpha''\beta''}\|_{Y_2}\|\varphi_{\alpha\beta}\|_{Y_2}.
\end{equation}
By Lemma~\ref{lem:QmG} with $f=\varphi_{\alpha''\beta''}$, and with $\|f\|_{L^2(\jap{\hat
z}^\theta)}+\|f\|_{H^1(\jap{\hat z}^{\theta'})}\le C\|f\|_{Y_2}$ from Lemma~\ref{lem:TG0}(ii), we
get
\begin{equation}
|\langle Q(\varphi_{\alpha''\beta''},\mathfrak g),\varphi_{\alpha\beta}\rangle|\le
C\|\varphi_{\alpha''\beta''}\|_{Y_2}\big(\|\mathfrak g\|_{Y_2}+\|
\mathfrak g\|_X\big)\|\varphi_{\alpha\beta}\|_{Y_2}.
\end{equation}
Thus both couplings are $\le C\|
\mathfrak g\|\,\|\varphi_{\alpha''\beta''}\|_{Y_2}\|\varphi_{\alpha\beta}\|_{Y_2}$, with
\begin{equation}\label{eq:frak-g-norm}
\|\mathfrak g\|:=\|\mathfrak g\|_\flat+\|\mathfrak g\|_{Y_2}+\|\mathfrak g\|_X .
\end{equation}
The mass tail is included in this bound.

\emph{Step 2: the coefficients $\mathfrak g_{\alpha'\beta'}$.} The coefficients $
\mathfrak g_{\alpha'\beta'}$ are Gaussian polynomials whose coefficients are polynomials in
$\partial^{\le|\beta'|}_\xi\varpi_g$. We split $\varpi_g=\varpi_*+\delta\varpi$. The Gaussian
coefficients are holomorphic in $\xi$ only on the polydiscs of radius $\simeq r_{\rm C}$, where
$\mathrm{Re}\,\theta^{-1}>0$. By Lemma~\ref{lem:Mref} and the part of the proof of
Lemma~\ref{lem:UFref} that does not use (U1), their $\varpi_*$-parts have Cauchy bounds at the
radius $r_{\rm C}$. Hence we have
\begin{equation}
\|\mathfrak g_{\alpha'\beta'}\|\le C_{\alpha'}\,\beta'!\,r_{\rm
C}^{-|\beta'|}\big(1+\|\delta\varpi\|_{\mathcal A}\big)
\end{equation}
for the $\varpi_*$-part. The $\delta\varpi$-part is $\mathfrak g_{\alpha'\beta'}(u_g)-
\mathfrak g_{\alpha'\beta'}(u_*)$. It is not bilinear, since $\mathfrak g_{\alpha'\beta'}$ is an
analytic, non-polynomial function of the parameters $u$ (through $e^\ell$, $e^{-\vartheta}$ and the
Gaussian). We write it as $\int_0^1D\mathfrak g_{\alpha'\beta'}(u_*+\tau(u_g-u_*))[u_g-u_*]\dd\tau$,
as in the proof of Lemma~\ref{lem:TWcoef}. The norm $\|\cdot\|_{\mathcal A}$ at the single point
$\xi$ is an algebra, by Lemma~\ref{lem:algebra} with $p=5$, and the proof of
Lemma~\ref{lem:TWcalc}(vi) applies
to it without change, with the point $\xi$ in place of the cube and $u_g(\xi)$ as the center $w_0$.
Its smallness condition has the form $C(\rad_0/r_{\rm C}+\|\delta\varpi\|_{\mathcal A})\le
r_p/(16dC'_{\rm alg})$, since the levels $\ge1$ of $u_*$ are bounded by $C\rad_0/r_{\rm C}$. Here
$C'_{\rm alg}$ is the algebra constant of the pointwise norm, which need not be $C_{\rm alg}$. The
condition holds by requirement~\ref{rq:TWa}, whose constant $C$ we take large enough to cover this
case and the factor $2$ in $\|\delta\varpi\|_{\mathcal A}\le2\delta_U$. As for the third bound of
Lemma~\ref{lem:TWcoef}, we obtain $\|\mathfrak g_{\alpha'\beta'}(u_g)-
\mathfrak g_{\alpha'\beta'}(u_*)\|_{\mathcal A}\le C\|\delta\varpi\|_{\mathcal A}$, with $C$ depending only on
constants chosen before $\delta_1$, the norm $\|\cdot\|$ above being used on each coefficient.

\emph{Step 3: summation of the linear couplings.} We sum these couplings by
the algebra lemma with the factor $\|\delta\varpi\|_{\mathcal A}$. Since
$w_n/w_{n-j}=\rad^j\frac{(n-j)!}{n!}\big(\frac{n+1}{n-j+1}\big)^5$ and $n+1\le(j+1)(n-j+1)$, we have
\begin{equation}
\frac{w_n\binom n j\,j!\,r_{\rm C}^{-j}}{w_{n-j}}\le\Big(\frac{\rad}{r_{\rm C}}\Big)^j(j+1)^5,
\end{equation}
which implies that the $\beta'\ne0$ couplings sum to $C(\rad_0/r_{\rm C})\mathbb D_\xi\le
C(L_1)(\rad_0/R_1)\mathbb D_\xi$. The couplings with $\beta'=0$, $\alpha'\ne0$ have $O(1)$
coefficients. They go from $y$-order $|\alpha''|<|\alpha|$ to $|\alpha|$, so they carry
$\delta_y^{|\alpha'|}\le\delta_y$ and sum to $C\delta_y\mathbb D_\xi$. Since the leaf factor
$\lambda_\xi$ appears once and both $\varphi$'s are measured in $Y_2$, each term is $\le
C\|
\mathfrak g\|\,\lambda_\xi^{1/2}\|\varphi_{\alpha''\beta''}\|_{Y_2}\cdot\lambda_\xi^{1/2}\|\varphi_{\alpha\beta}\|_{Y_2}$.

\emph{Step 4: the quadratic terms.} For $Q(\varphi_1,\varphi_2)$ with $\varphi_1$ of level $k$ and
$\varphi_2$ of level $n-k$, by part~(iii) of the trilinear proposition, second form, we have the
bound
\begin{equation}
C\big(\|\varphi_1\|_X\|\varphi_2\|_{Y_2}+\|\varphi_1\|_{Y_2}\|\varphi_2\|_X\big)\|\varphi_{\alpha\beta}\|_{Y_2}.
\end{equation}
We take $a_\beta=\|\varphi_1\|_X$ or $\lambda_\xi^{1/2}\|\varphi_1\|_{Y_2}$, and similarly for $b$,
$d$. Then the algebra lemma with $p=5$ implies
\begin{equation}
\sum_nw_n^2\sum\cdots\le C\Big(\sum w^2E\Big)^{1/2}\Big(\sum
w^2\lambda_\xi\|\cdot\|^2_{Y_2}\Big)^{1/2}\Big(\sum w^2\lambda_\xi\|\cdot\|^2_{Y_2}\Big)^{1/2}=C\,
\mathbb A_\xi^{1/2}\,\mathbb D_\xi .
\end{equation}
We take no supremum over $\xi$, since all factors are on the same leaf. The nonlocal parts are
treated in the same way with (iv). The proof does not use the term $C\mathbb A_\xi\mathbb D_\xi$
of the statement. This completes the proof.
\end{proof}

\section{The crude regime: fibered comparison}\label{sec:fib}

In this section we control a leaf before it enters the near-Maxwellian regime, that is, while its
clock $s_\xi$ is bounded, by comparison with the fibered flow. The main result is
Theorem~\ref{thm:RGF}.

We make the comparison on a bounded effective time interval $\{s_\xi\le S^\sharp\}$,
$S^\sharp=4S$, for a level $S$. In the bootstrap we have $S=S_\sharp$, as in Lemma~\ref{lem:BTU4}.

\emph{The fibered flow.} Let $G$ be the solution of
\begin{equation}\label{eq:fibflow}
\partial_tG=\kappa\,Q\big(G(t,\cdot,v),G(t,\cdot,v)\big),\qquad G(T_0)=g(T_0).
\end{equation}
By \eqref{eq:frame}, on the fiber $v$ it is a homogeneous trajectory in the frame $P_v$ of the fiber
parameters, that is,
\begin{equation}\label{eq:fibflow-frame}
G(t,\cdot,v)=T_{P_v}F_v(s_v(t)),\qquad s_v=\lambda_v\big(\tau(t)-\tau(T_0)\big).
\end{equation}
By Section~\ref{sec:homog}, the trajectories $F_v$ are uniformly smooth and Gaussian
on bounded intervals.

\emph{The remainder in Gaussian weights.} We have to redo Proposition~\ref{prop:Rest} in Gaussian
weights, and without derivatives in $v$ of the fibered flow, which are not uniformly controlled.
There are three points that require care.

The first one is a derivative loss in the remainder. By Proposition~\ref{prop:Rformula},
$\partial^\alpha_yR$ contains $|\alpha|+2$ derivatives of $g$, one of them in $v$. The comparison
energy controls only $|\alpha|+|\beta|\le N$, so it cannot bound these local factors. We take the
two extra derivatives from a crude bound. Since $R$ carries $t^{-1}$, the crude constant affects
only the lower bound on $T_0$.

The second one is a cycle in the Gaussian rates. A crude bound enters $R$ through local factors
measured in $L^2(e^{2q|y|^2})$. Consequently it must decay at a Gaussian rate strictly larger than
the rate $q$ of the comparison, and hence larger than the rate of its conclusion. The bootstrap
supplies crude bounds only at its own rate $q_0$, and it needs the conclusion of
Theorem~\ref{thm:RGF} at that same rate, through (G4), the Gaussian $C^2$ bound on the background
that we need for the second half $Q(\varphi,\mu+f)$ of the linear term $L_{\mu+f}\varphi$ (note
that Lemma~\ref{lem:Whypo} treats only $Q(\mu+f,\varphi)$). We break the cycle with the Maxwellian
barrier for the real solution of Lemma~\ref{lem:RGbar}. It requires no Gaussian bound from the
bootstrap, just the moment and $L^\infty$ bounds (G1), (G3) on the leaves, and it yields a
\emph{structural} Gaussian rate. Then, by interpolation against a crude bound with no Gaussian
weight, we obtain in Lemma~\ref{lem:RGloc} the local factors at a fixed fraction of that rate.

The third point is that we need no $v$-derivatives of the fibered flow. We run the comparison fiber
by fiber, with $y$-derivatives only. On a leaf, $v$-derivatives enter only through
$\nabla_z=\partial_y-t^{-1}\partial_v$, with the factor $t^{-1}$, and they are covered by the crude
bound. Hence we need no bounds on $\partial^\beta_vG$, and no derivatives $D^k\Phi_\sigma$,
$k\ge2$, of the homogeneous flow.

Throughout this section, $\gamma\in[-3,-2]$, and constants are uniform in $\gamma$ on this range
unless a dependence is stated. In Subsection~\ref{ssec:fib} we provide Gaussian Sobolev bounds for
homogeneous trajectories on short intervals, and bounds on the local factors, on the remainder and
on the fibered flow.

\subsection{A Maxwellian barrier for the real solution}\label{ssec:barrier}

In the lemma below we bound the solution by a Gaussian in $y$ at a structural
rate, using, besides a Gaussian bound at $T_0$, only the moment and $L^\infty$ bounds on the leaves.

\emph{Idea of the proof.} The upper barrier is a sum of two Gaussians in $y$, with rates chosen so
that $\kappa Q(g^{(x)},\cdot)$ applied to them is controlled by their time derivative. At a first
touching point, we compare $\partial_tg$ with the barrier by the leaf identity. The domain is
non-compact, and we handle this by a penalization in $\xi$.

The leaf restriction of $g$ through $x$ at time $t$ is $g^{(x)}(z)=g(t,z,x/t-z/t)$, and $\xi=x/t$ is
constant along it. By Proposition~\ref{prop:leaf}, we have
\begin{equation}\label{eq:bar-leafid}
\partial_tg(t,y,v)=\kappa(t)Q(G_1,G_1)(y)\quad\text{exactly},\qquad G_1=g^{(y+tv)}.
\end{equation}
There is no remainder in this form.

\begin{lem}[A Maxwellian barrier for the real solution]\label{lem:RGbar}
Let $g\ge0$ be a classical solution on $[T_0,T]$, $T_0\ge2$. Let $\rho:\R^3\to(0,\rho_+]$ be $C^1$
with $|\nabla\log\rho|\le L$, $L\ge1$. Assume:
\begin{itemize}
\item[(a)] for all $t\in[T_0,T]$ and $x$, with $\xi=x/t$,
\begin{equation}\label{eq:bar-a}
M_-\rho(\xi)\le\int g^{(x)}\le M_+\rho(\xi),\qquad\int|z|^2g^{(x)}\le E_+\rho(\xi),\qquad\sup_z\jap
z^8g^{(x)}(z)\le K\rho(\xi);
\end{equation}
\item[(b)] $g(T_0,y,v)\le C_D\,\rho(v+y/T_0)\,e^{-r_D|y|^2}$;
\item[(c)] (qualitative) $\sup_{[T_0,T]\times\R^6}e^{q_1|y|^2}g(t,y,v)/\rho(v+y/t)<\infty$ for some
$q_1>0$.
\end{itemize}
Let $r:=\min\{r_D,\,r_{\rm b}\}$ with $r_{\rm b}:=M_-/(48E_+)$, and let $T_0\ge2/r$. Then, for
$t\in[T_0,T]$,
\begin{equation}\label{eq:bar}
g(t,y,v)\le C_{\rm bar}\,\rho(v+y/t)\,e^{-(r/2)|y|^2},\qquad C_{\rm
bar}:=e^{(L+1)^2/(4T_0)}\max\{C_D,\,2Ke^{rR_B^2}\},
\end{equation}
where $R_B=R_B(\gamma,r,M_\pm,E_+,K)$ is given in Step~1.
\end{lem}

\begin{proof}
\emph{Step 1: Gaussian supersolutions on one leaf.} Let $G\ge0$ on $\R^3$ satisfy the three bounds
of (a) with $\rho(\xi)$ replaced by a number $\rho_0>0$. Let $\kappa_G:=-\cb[G]\ge0$, so that
$\kappa_G=8\pi G$ for $\gamma=-3$. For $s>0$ let $b_s(z)=e^{-s|z|^2}$ and
\begin{equation}
\mathcal L_Gb:=\ab[G]:\nabla^2b+\kappa_Gb .
\end{equation}
\emph{Claim:} if $0<s\le r_{\rm b}$ and $|z|\ge R(s)$, then
\begin{equation}\label{eq:bar-claim}
\mathcal L_Gb_s(z)\le-\tfrac13sM_-\rho_0|z|^{\gamma+2}b_s(z).
\end{equation}

We have $\nabla^2b_s=(4s^2z\otimes z-2sI)b_s$, so $\ab:\nabla^2b_s=(4s^2|z|^2\,\ab\hat z\cdot\hat
z-2s\,\mathrm{tr}\,\ab)b_s$. As in the proof of Lemma~\ref{lem:Whypo}, we have
\begin{equation}
\ab[G](z)\hat z\cdot\hat z=\int|z-w|^\gamma|w_\perp|^2G(w)\dd
w,\qquad\mathrm{tr}\,\ab[G](z)=2\int|z-w|^{\gamma+2}G(w)\dd w .
\end{equation}
Let $|z|\ge1$. We split at $|w|=|z|/2$.
\begin{itemize}
\item On $\{|w|\le|z|/2\}$: $|z-w|^\gamma\le2^{|\gamma|}|z|^\gamma\le8|z|^\gamma$, and
$|z-w|^{\gamma+2}\ge(\frac32|z|)^{\gamma+2}\ge\frac23|z|^{\gamma+2}$. By Chebyshev's inequality,
$\int_{|w|\le|z|/2}G\ge(M_--4E_+/|z|^2)\rho_0\ge\frac12M_-\rho_0$ if $|z|^2\ge8E_+/M_-$.
\item On $\{|w|>|z|/2\}$: $|w_\perp|\le|z-w|$, $\jap w^{-8}\le2^8\jap z^{-8}$, and
$\int_{|w|>|z|/2}\jap w^{-8}\dd w\le C\jap z^{-5}$. Hence $\int_{|w|>|z|/2}|z-w|^{\gamma+2}G\le
C_{\rm abs}K\rho_0\jap z^{-5}$, with $C_{\rm abs}$ absolute.
\end{itemize}
Therefore, we have $\ab\hat z\cdot\hat z\le\rho_0(8E_+|z|^\gamma+C_{\rm abs}K\jap z^{-5})$ and
$\mathrm{tr}\,\ab\ge\frac23M_-\rho_0|z|^{\gamma+2}$. Now we come to $\kappa_G$. If $\gamma=-3$, then
$\kappa_G(z)=8\pi G(z)\le8\pi K\rho_0\jap z^{-8}$. If $\gamma>-3$, the same split, with
$2(\gamma+3)\int_{|u|<1}|u|^\gamma\dd u=8\pi$, yields
$\kappa_G(z)\le16M_+\rho_0|z|^\gamma+C_{\rm abs}K\rho_0\jap z^{-5}$. Altogether, we obtain
\begin{equation}
\frac{\mathcal L_Gb_s}{\rho_0b_s}\le\big(32s^2E_+-\tfrac43sM_-\big)|z|^{\gamma+2}+4C_{\rm
abs}Ks^2|z|^{-3}+16M_+|z|^{\gamma}+(C_{\rm abs}+8\pi)K|z|^{-5}.
\end{equation}
For $s\le r_{\rm b}$ the bracket is $\le-\frac23sM_-$. Since $|z|^{\gamma+2}\ge|z|^{-1}$, the other
three terms are $\le\frac13sM_-|z|^{\gamma+2}$ as soon as we have
\begin{equation}
|z|\ge R(s):=\max\Big\{1,\Big(\frac{8E_+}{M_-}\Big)^{\frac12},\Big(\frac{48C_{\rm
abs}Ks}{M_-}\Big)^{\frac12},\Big(\frac{192M_+}{sM_-}\Big)^{\frac12},\Big(\frac{6(C_{\rm
abs}+8\pi)K}{sM_-}\Big)^{\frac14}\Big\}.
\end{equation}
This proves the claim. On $\{|z|\le R\}$, the same bounds imply $|
\mathcal L_Gb_s|\le\Lambda(s,R)\rho_0b_s$, with $\Lambda$ depending on $(s,R,K,M_+)$.

Let $R_B:=\sup_{s\in[r/2,r]}R(s)$, which is finite and explicit.

\emph{Step 2: the barrier family.} For $\delta\in(0,1]$ let
\begin{equation}\label{eq:bar-B}
B_\delta(t,y,v):=C_B\phi(t)\rho(\xi)\big(1+\delta\jap\xi^2\big)e^{-r(t)|y|^2}+\delta\phi_1(t)\rho(\xi)e^{-q'(t)|y|^2},\qquad\xi=v+y/t,
\end{equation}
where:
\begin{itemize}
\item $C_B:=\max\{C_D,2Ke^{rR_B^2}\}$;
\item $r(t):=r-(T_0^{-1}-t^{-1})\in[r/2,r]$, since $T_0\ge2/r$;
\item $q':=\min\{q_1/2,r/2\}$ and $q'(t):=q'\big(1-(T_0^{-1}-t^{-1})\big)\in[q'/2,q']$;
\item $\phi(t):=e^{(L+1)^2(T_0^{-1}-t^{-1})/4}$ and $\phi_1(t):=e^{(L+1)^2(T_0^{-1}-t^{-1})/(4q')}$.
\end{itemize}
\begin{itemize}
\item[(2a)] $\partial_tB_\delta\ge0$ at fixed $(y,v)$. Since $\partial_t\xi=-y/t^2$ and
$|2\delta\xi|/(1+\delta\jap\xi^2)\le1$, the logarithmic derivative of the first term is at least
\begin{equation}
t^{-2}\big(\tfrac14(L+1)^2-(L+1)|y|+|y|^2\big)=t^{-2}\big(|y|-\tfrac12(L+1)\big)^2\ge0 .
\end{equation}
That of the second is at least $t^{-2}\big((L+1)^2/(4q')-L|y|+q'|y|^2\big)\ge0$.
\item[(2b)] The restriction of $B_\delta$ to the leaf through $x$ is $A_1b_{r(t)}+A_2b_{q'(t)}$,
where $A_1=C_B\phi\rho(\xi)(1+\delta\jap\xi^2)$ and $A_2=\delta\phi_1\rho(\xi)$ are constant on the
leaf.
\item[(2c)] $g(T_0)<B_\delta(T_0)$ everywhere, by (b), $r\le r_D$ and $C_B\ge C_D$. Moreover
$g/B_\delta\to0$ as $|(y,v)|\to\infty$, uniformly in $t$. For $|y|\to\infty$, use (c) and $q'\le
q_1/2$. For $|y|\le R$ and $|v|\to\infty$, $|\xi|\to\infty$, $g\le K\rho(\xi)$ by (a), and
$B_\delta\ge C_B\delta\rho(\xi)\jap\xi^2e^{-rR^2}$.
\end{itemize}

\emph{Step 3: no touching.} We fix $\delta$, and suppose that $g\le B_\delta$ fails somewhere on
$[T_0,T]$. Let $\bar t$ be the supremum of the times $t$ with $g\le B_\delta$ on
$[T_0,t]\times\R^6$. By (2c) and continuity, we have $\bar t>T_0$, and there is $(y_1,v_1)$ with
$(g-B_\delta)(\bar t,y_1,v_1)=0=\max(g-B_\delta)(\bar t,\cdot)$. Since $g-B_\delta\le0$ before $\bar
t$, we have
\begin{equation}
\partial_t(g-B_\delta)(\bar t,y_1,v_1)\ge0 .
\end{equation}
Let $x_1=y_1+\bar t v_1$ and $G_1=g^{(x_1)}$. Then $G_1\le A_1b_{r(\bar t)}+A_2b_{q'(\bar t)}$ on
$\R^3$, with equality at $y_1$. Hence $\nabla^2G_1(y_1)\le\nabla^2B^{(x_1)}_\delta(y_1)$, and since
$\ab[G_1]\ge0$, we get
\begin{equation}
\partial_tg(\bar t,y_1,v_1)=\kappa Q(G_1,G_1)(y_1)\le\kappa\big(A_1\mathcal L_{G_1}b_{r(\bar
t)}+A_2\mathcal L_{G_1}b_{q'(\bar t)}\big)(y_1).
\end{equation}
Note that $|y_1|\ge R_B$, since for $|y|<R_B$, (a) implies $g\le
K\rho(\xi)<C_B\rho(\xi)e^{-rR_B^2}\le A_1b_{r(t)}(y)$. By Step~1 with $\rho_0=\rho(\xi_1)$, the
first term satisfies
\begin{equation}
A_1\mathcal L_{G_1}b_{r(\bar t)}(y_1)\le-\tfrac16rM_-\rho(\xi_1)|y_1|^{\gamma+2}A_1b_{r(\bar
t)}(y_1)<0 .
\end{equation}
Let $R_{\rm b}:=\max\{R_B,\sup_{s\in[q'/2,q']}R(s)\}$.
\begin{itemize}
\item If $|y_1|\ge R_{\rm b}$, the second term is $\le0$ too.
\item If $R_B\le|y_1|\le R_{\rm b}$, the first term is $\le-c_\flat C_B\rho(\xi_1)^2$, with
$c_\flat:=\frac16rM_-R_{\rm b}^{\gamma+2}e^{-rR_{\rm b}^2}$, and the second is
$\le\Lambda_1\delta\phi_1(T)\rho(\xi_1)^2$, with $\Lambda_1:=\sup_{s\in[q'/2,q']}\Lambda(s,R_{\rm
b})$.
\end{itemize}
For $\delta\le\delta_{1,\rm b}:=c_\flat C_B/(2\Lambda_1\phi_1(T))$ the sum is negative in both
cases. Thus we have $\partial_tg(\bar t,y_1,v_1)<0\le\partial_tB_\delta(\bar t,y_1,v_1)$, which is a
contradiction. Hence $g\le B_\delta$ on $[T_0,T]$ for every $\delta\le\delta_{1,\rm b}$. Letting
$\delta\to0$ and using $\phi\le e^{(L+1)^2/(4T_0)}$ and $r(t)\ge r/2$, we obtain the claimed bound.
This completes the proof.
\end{proof}

Note that the barrier uses neither (G2) nor any Gaussian bound from the bootstrap. Its rate $r$
depends only on the data rate and on $M_-/E_+$, which are structural. The qualitative hypothesis (c)
only makes the maximum attained, and no constant depends on $q_1$.

\subsection{Standing assumptions}\label{ssec:fibstanding}

Here we list the hypotheses under which we prove the fibered comparison. From now on, in this
section and in Subsection~\ref{ssec:fib}, $P_v=(\rho_v,y_{0,v},\theta_v)$ are the mass, mean and
temperature of $g(T_0,\cdot,v)$, and we assume:
\begin{itemize}
\item[(H1)] $v\mapsto P_v$ is $C^1$, $\rho_v\in(0,\rho_+]$, $|y_{0,v}|\le L$,
$\theta_v\in[\theta_-,\theta_+]$, and $|\nabla\log\rho_v|+|\nabla\log\theta_v|\le L$, with $L\ge1$.
\end{itemize}
We define the fiber rate and the fiber clock by
\begin{equation}\label{eq:fiberclock}
\lambda_v:=\lambda(P_v)=\rho_v\theta_v^{\gamma/2},\qquad
s_v(t):=\lambda_v\big(\tau(t)-\tau(T_0)\big).
\end{equation}
Then we have $|\nabla\log\lambda_v|\le(1+\frac{|\gamma|}2)L\le\frac52L$ and $\rho_v\le
c_\theta\lambda_v$, with $c_\theta:=\max\{1,\theta_+\}^{3/2}$. We apply Lemma~\ref{lem:RGbar} with
$\rho(\xi):=\rho_\xi$, and we write $r$, $C_{\rm bar}$ for its constants. Finally, we assume
\begin{itemize}
\item[(H3)] (crude bound) for some $J\ge1$ and all $t\in[T_0,T]$,
\begin{equation}\label{eq:H3}
\sum_{|a|+|b|\le J}\ \sup_{y,v}\ \rho_v^{-1}\big|\partial^a_y\partial^b_vg(t,y,v)\big|\le K_{\rm
cr} .
\end{equation}
\end{itemize}
No Gaussian weight is required in (H3).

For the fibered flow $G$, with $G(T_0)=g(T_0)$, we assume:
\begin{itemize}
\item[(H2)] (KD$^{n+2}$) $F_{v,0}:=T_{P_v}^{-1}g(T_0,\cdot,v)\in
\mathfrak K^{n+2}(C_{\rm dat},\beta,A_{n+2})$ for all $v$, with constants independent of $v$.
\end{itemize}
In the notation of Lemma~\ref{lem:sym}, $T_{P}F(y)=\rho\theta^{-3/2}F\big((y-y_0)/\theta^{1/2}\big)$
and $\lambda(P)=\rho\theta^{\gamma/2}$.

\subsection{The fibered comparison}\label{ssec:fibcomp}

In the next proposition we show that, on each fiber and as long as its clock is bounded,
the solution stays within $O(T_0^{-1})$ of the fibered flow, in Gaussian weights.

\begin{prop}[Fibered comparison]\label{prop:RGcompF}
Assume (H1)--(H3) and the hypotheses (a), (c) of Lemma~\ref{lem:RGbar} with $\rho(\xi)=\rho_\xi$,
and let the fibered trajectories satisfy the conclusions of Proposition~\ref{prop:RS} on
$[0,S^\sharp]$. Then hypothesis (b) of that lemma holds by
(H1)--(H2), with $r_D=\beta/(4\theta_+)$ and a constant $C(C_{\rm dat},\beta,L,\theta_\pm)\le C_D$.
Indeed, we have
\begin{equation}\label{eq:RGcomp-b}
g(T_0,y,v)\le C_{\rm dat}\theta_-^{-3/2}\rho_ve^{-\beta|y-y_{0,v}|^2/\theta_+}
\end{equation}
and $\rho_v\le\rho_\xi e^{L|y|}$. Let $T_0\ge\max\{2/r,\,(L+1)^2/4\}$ (the first is used in
that lemma, the second for $e^{(L+1)^2/(4T_0)}\le e$), $n\ge1$, $J\ge n+3$, and
\begin{equation}\label{eq:qnat}
0<q\le q_\natural:=\min\Big\{\frac1{8(L^2+3\theta_+)},\ \frac{q_S}{32\theta_+},\
\frac{r_{n+2}}4\Big\}.
\end{equation}
Put $\hat w:=(g-G)/\rho_v$, $m:=e^{q|y|^2}$, and
\begin{equation}\label{eq:compE}
\hat{\mathfrak E}_v(t):=\sum_{|\alpha|\le n}\epsilon_*^{2|\alpha|}\big\|m\,\partial^\alpha_y\hat
w(t,\cdot,v)\big\|^2_{L^2_y},
\end{equation}
where $\epsilon_*\in(0,1]$ is fixed in the proof. There are constants $C_1$, $\delta_{0,\rm f}$,
$C_{\rm comp}$ and $C_{S,\rm f}$, depending on $(n,q,S^\sharp,L,\theta_\pm,C_{\rm
dat},\beta,A_{n+2})$ and on $(M_\pm,E_+,J)$ through $r_\natural$, such that the following holds. If
\begin{equation}\label{eq:compast}
(\ast)\qquad C_{\rm comp}\,C_{\rm loc}^2\,S^\sharp e^{C_1S^\sharp}\,T_0^{-1}\le\delta_{0,\rm f}/2,
\end{equation}
then for all $t\in[T_0,T]$ and all $v$ with $s_v(t)\le S^\sharp$:
\begin{equation}\label{eq:comp}
\hat{\mathfrak E}_v(t)^{1/2}\le\delta_{c,\rm f}(t,v):=C_{\rm comp}\,C_{\rm
loc}^2\,s_v(t)e^{C_1s_v(t)}\,T_0^{-1},
\end{equation}
and, for $|\alpha|\le n-2$,
\begin{equation}
\sup_ye^{(q/2)|y|^2}\big|\partial^\alpha_y(g-G)(t,y,v)\big|\le C_{S,\rm f}\,\rho_v\,\delta_{c,\rm
f}(t,v).
\end{equation}
\end{prop}

\begin{proof}
We fix $v$, and let $t_v$ be the last time in $[T_0,T]$ with $s_v(t_v)\le S^\sharp$.

\emph{Step 1: the equation.} Fiberwise, we have $\partial_tg=\kappa(Q(g,g)+R)$ and
$\partial_tG=\kappa Q(G,G)$. With $g=G+w$ and $\hat G=G/\rho_v$, by bilinearity we get
\begin{equation}\label{eq:comp-eq}
\partial_t\hat w=\kappa\Big[\rho_v\big(Q(\hat G,\hat w)+Q(\hat w,\hat G)+Q(\hat w,\hat w)\big)+\hat
R\Big],\qquad\hat R:=R/\rho_v .
\end{equation}
The weight $m$ is static, so by Lemma~\ref{lem:localnoloss} there is no transport term. Since
$\partial_y$ is a derivation of $Q$, for $|\alpha|=k\le n$ and $\hat h_\alpha:=\partial^\alpha_y\hat
w$, we have
\begin{equation}\label{eq:comp-energy-id}
\frac12\frac{d}{dt}\|m\hat h_\alpha\|^2=\kappa\rho_v\sum_{a\le\alpha}\binom\alpha
a\Big[T_2(a)+T_3(a)+T_4(a)\Big]+\kappa\langle\partial^\alpha\hat R,\hat h_\alpha\rangle_{L^2(m)},
\end{equation}
where $T_2(a)=\langle Q(\partial^a\hat G,\partial^{\alpha-a}\hat w),\hat h_\alpha\rangle_{L^2(m)}$,
and $T_3$, $T_4$ are the same with $Q(\partial^a\hat w,\partial^{\alpha-a}\hat G)$ and
$Q(\partial^a\hat w,\partial^{\alpha-a}\hat w)$. By Lemmas~\ref{lem:RGloc} and~\ref{lem:RGGbounds},
$\partial^\alpha_yg$ and $\partial^\alpha_yG$ decay at Gaussian rates $>q$, namely the rate
$r_\natural$ for $g$, and the rate $q_S/(8\theta_+)\ge4q$ for $G$. These bounds
hold, by the first lemma with $j\le n+2\le J-1$, for the derivatives of $g$ of total order
$\le n+2$ and, by part~(ii) of the second, for the $y$-derivatives of $G$ of order $\le n+2$,
and hence, through the equations of $g$ and of $G$, for $\partial_t\partial^\alpha_y\hat w$ with
$|\alpha|\le n$. Hence $\hat{\mathfrak E}_v$ is finite and continuously differentiable on
$[T_0,t_v]$, the identity~\eqref{eq:comp-energy-id} holds classically, the integrations by parts in
$y$ that are used in Step~2 have no boundary terms, and $\hat{\mathfrak E}_v(T_0)=0$.

\emph{Step 2: the terms.} We write $E_\alpha=\|m\hat h_\alpha\|^2$, $D_\alpha=\|\hat
h_\alpha\|^2_{Y_2(m)}$, $E_\le=\sum_{|\alpha|\le n}E_\alpha$ and $D_\le=\sum_{|\alpha|\le
n}D_\alpha$.
\begin{itemize}
\item[(T1)] $T_2(0)\le-c_{W,\rm f}D_\alpha+C_{W,\rm f}E_\alpha$, by
Lemma~\ref{lem:RGGbounds}(iii) and Lemma~\ref{lem:RFtools}(vi).
\item[(T2)] For $a\ne0$, by part~(iii) of the same lemma with $f=\partial^a\hat G$, we have
$|T_2(a)|\le CK_G\big(D_{\alpha-a}^{1/2}+E_{\alpha-a}^{1/2}\big)D_\alpha^{1/2}$. Here the bound
$\|\partial^a\hat G\|_{H^1(\jap y^3)}\le CK_G$ uses order $\le n+1$ and $\jap y^4\le Ce^{2q|y|^2}$.
\item[(T3)] By part~(ii) of that lemma with $h=\partial^{\alpha-a}\hat G$, we have $|T_3(a)|\le
CK_GE_a^{1/2}E_\alpha^{1/2}$. Here $h$ has order $\le n+2$, and $e^{2q\jap y}\le e^{3q}e^{2q|y|^2}$.
\item[(T4)] By part~(iii) with $f=\partial^a\hat w$, $h=\partial^{\alpha-a}\hat w$, we
obtain
\begin{equation}
|T_4(a)|\le C\big(\|\partial^a\hat w\|_{L^2(\jap y^4)}D_{\alpha-a}^{1/2}+\|\partial^a\hat
w\|_{H^1(\jap y^3)}E_{\alpha-a}^{1/2}\big)D_\alpha^{1/2}.
\end{equation}
If $|a|\le n-1$, then $\|\partial^a\hat w\|_{H^1(\jap y^3)}\le CE_\le^{1/2}$. If $|a|=n$, then
$a=\alpha$, and $\|\hat h_\alpha\|_{H^1(\jap y^3)}\le CD_\alpha^{1/2}$ by
part~(iv), while $E_0\le E_\le$. In all cases we have $|T_4(a)|\le
CE_\le^{1/2}(D_\le+E_\le)$.
\end{itemize}

\emph{Step 3: summation.} Let $\hat{\mathfrak D}_v:=\sum_\alpha\epsilon_*^{2|\alpha|}D_\alpha$. Let
$N_n$ be the number of pairs $(\alpha,a)$ with $|\alpha|\le n$, $a\le\alpha$.
\begin{itemize}
\item For (T2), by Young's inequality and
$\epsilon_*^{2|\alpha|}\le\epsilon_*^2\epsilon_*^{2|\alpha-a|}$ for $a\ne0$, we have
\begin{equation}
\sum\epsilon_*^{2|\alpha|}|T_2(a)|\le\frac{c_{W,\rm f}}4\hat{
\mathfrak D}_v+C''K_G^2\epsilon_*^2\big(\hat{\mathfrak D}_v+\hat{\mathfrak E}_v\big).
\end{equation}
We choose $\epsilon_*$ with $C''K_G^2\epsilon_*^2\le c_{W,\rm f}/4$.
\item For (T3), $\epsilon_*^{2|\alpha|}\le\epsilon_*^{|a|}\epsilon_*^{|\alpha|}$ implies
$\sum\epsilon_*^{2|\alpha|}|T_3(a)|\le CK_G\hat{\mathfrak E}_v$.
\item For (T4), it follows from $E_\le\le\epsilon_*^{-2n}\hat{\mathfrak E}_v$ and
$D_\le\le\epsilon_*^{-2n}\hat{\mathfrak D}_v$ that $\sum\epsilon_*^{2|\alpha|}|T_4(a)|\le
C\epsilon_*^{-3n}\hat{\mathfrak E}_v^{1/2}(\hat{\mathfrak D}_v+\hat{\mathfrak E}_v)$.
\end{itemize}
Let $\delta_{0,\rm f}:=c_{W,\rm f}\epsilon_*^{3n}/(4C)$. As long as $\hat{
\mathfrak E}_v\le\delta_{0,\rm f}^2$, we have
\begin{equation}\label{eq:comp-diff}
\frac12\frac{d}{dt}\hat{\mathfrak E}_v+\frac{c_{W,\rm f}}4\kappa\rho_v\hat{\mathfrak D}_v\le
C_{\sharp,\rm f}\kappa\rho_v\hat{\mathfrak E}_v+\kappa\hat{\mathfrak E}_v^{1/2}\|\hat
R\|_n,\qquad\|\hat R\|_n:=\Big(\sum_{|\alpha|\le n}\|m\partial^\alpha\hat R\|^2\Big)^{1/2},
\end{equation}
with $C_{\sharp,\rm f}:=C_{W,\rm f}+CK_G+c_{W,\rm f}$. The constants $c_{W,\rm f},C_{W,\rm f},K_G$,
hence $\epsilon_*,\delta_{0,\rm f},C_{\sharp,\rm f}$, come from these two lemmas.

\emph{Step 4: Gronwall.} By Lemma~\ref{lem:RGremF}, we have $\|\hat R\|_n\le C_RN_n^{1/2}C_{\rm
loc}^2\rho_vt^{-1}$. Let $u_\eta:=(\hat{\mathfrak E}_v+\eta^2)^{1/2}$. Then we have $u_\eta'\le
C_{\sharp,\rm f}\kappa\rho_vu_\eta+\kappa\|\hat R\|_n$. Using $\rho_v\le c_\theta\lambda_v$, $\dot
s_v=\kappa\lambda_v$ and $t\ge T_0$, and letting $\eta\to0$, we obtain
\begin{equation}
\hat{\mathfrak E}_v(t)^{1/2}\le\int_{T_0}^te^{C_1(s_v(t)-s_v(t'))}\,C_RN_n^{1/2}c_\theta C_{\rm
loc}^2\,\frac{\kappa(t')\lambda_v}{T_0}\dd t'\le C_{\rm comp}C_{\rm
loc}^2\,\frac{s_v(t)e^{C_1s_v(t)}}{T_0},
\end{equation}
with $C_1:=C_{\sharp,\rm f} c_\theta$ and $C_{\rm comp}:=C_RN_n^{1/2}c_\theta$.

\emph{Step 5: continuity.} The set of $t\in[T_0,t_v]$ with $\hat{\mathfrak E}_v\le\delta_{0,\rm
f}^2$ on $[T_0,t]$ is closed and contains $T_0$. On it, $(\ast)$ implies $\hat{
\mathfrak E}_v^{1/2}\le\delta_{0,\rm f}/2$, so the set is also open in $[T_0,t_v]$. Hence it is all of
$[T_0,t_v]$. The argument is fiber by fiber, since the remainder is bounded through
the lemma of Step~4, which uses no information on other fibers beyond (H3) and the barrier.

\emph{Step 6: pointwise bounds.} For $|\alpha|\le n-2$, by Sobolev's inequality, we have
\begin{equation}
\|e^{(q/2)|y|^2}\partial^\alpha\hat w\|_{L^\infty}\le C\sum_{|\beta|\le2}\|\jap
y^2e^{(q/2)|y|^2}\partial^{\alpha+\beta}\hat w\|\le C\epsilon_*^{-n}\hat{\mathfrak E}_v^{1/2}.
\end{equation}
This completes the proof.
\end{proof}

With $S^\sharp=4S$, the factor $e^{C_1S^\sharp}$ of the proposition is $e^{4C_1S}$,
and $C_{\rm loc}^2$ takes the place of $(1+K_{\rm glob})^2$.

\subsection{From fibers to leaves}\label{ssec:fib2leaf}

In the next theorem we transfer the fiberwise comparison to the leaves. A leaf
at time $t$ samples the fibers $v=\xi-z/t$. On the region $Z_t$ defined below, these fibers have
clocks not exceeding $\frac43s_\xi$, so that Proposition~\ref{prop:RGcompF} applies to them, while
outside $Z_t$ the bounds of Lemma~\ref{lem:RGloc} and the Gaussian decay suffice. For a leaf through
$x=t\xi$ we write $v(z):=\xi-z/t$ and $\hat g^{(x)}:=g^{(x)}/\rho_\xi$. Let $c_{\rm
fib}:=\frac25\log\frac43$ and $Z_t:=\{|z|\le c_{\rm fib}t/L\}$.

\begin{theo}[Gaussian regularity of the leaves before entry]\label{thm:RGF}
Assume the hypotheses of Proposition~\ref{prop:RGcompF} with $n\ge3$, $S^\sharp=4S$ for a level
$S>0$, and $(\ast)$. Put
\begin{equation}\label{eq:epsT}
\delta^\sharp_{c,\rm f}:=C_{\rm comp}C_{\rm loc}^2S^\sharp
e^{C_1S^\sharp}T_0^{-1},\qquad\varepsilon_T(t):=\delta^\sharp_{c,\rm f}+C_{\rm loc}t^{-1}+C_{\rm
loc}e^{-r_\natural c_{\rm fib}^2t^2/(8L^2)} .
\end{equation}
Let $t\in[T_0,T]$ and let $\xi$ satisfy $s_\xi(t)\le3S$. Then for $j\le n-2$ and all $z$:
\begin{enumerate}[label=(\roman*)]
\item \emph{(closeness)}
\begin{equation}\label{eq:RGF-i}
\big|\nabla^j_z\hat g^{(x)}(z)-\mathbf 1_{Z_t}(z)\,\rho_\xi^{-1}(\partial^j_yG)(t,z,v(z))\big|\le
C_T\,\varepsilon_T(t)\,e^{-(q/4)|z|^2};
\end{equation}

\item \emph{(Gaussian bounds)}
\begin{equation}\label{eq:RGF-ii}
\sum_{j\le n-2}\big(\|e^{(q/4)|z|^2}\nabla^j_z\hat
g^{(x)}\|_{L^\infty}+\|e^{(q/8)|z|^2}\nabla^j_z\hat g^{(x)}\|_{L^2}\big)\le C_T(K_G+\varepsilon_T);
\end{equation}

\item \emph{(moments)}
\begin{equation}\label{eq:RGF-iii}
\big|\int\hat g^{(x)}-1\big|+\big|\int z\hat g^{(x)}-y_{0,\xi}\big|+\big|\int|z|^2\hat
g^{(x)}-|y_{0,\xi}|^2-3\theta_\xi\big|\le C_M\varepsilon_T(t);
\end{equation}

\item \emph{((G1)--(G3))} if moreover $C_M\varepsilon_T(T_0)\le\varepsilon_\natural(L,\theta_-)$
and $C(q)C_T\varepsilon_T(T_0)\le1$, then $\hat g^{(x)}$ satisfies (G1) with mass in
$[\frac12,\frac32]$, center $|y_0^x|\le L+1$ and temperature
$\theta^x\in[\frac12\theta_-,\theta_++1]$; (G3) with constant $2K_3+1$, with $K_3$ of
Lemma~\ref{lem:RGGbounds}; and (G2) with $c_G$ from Lemma~\ref{lem:RGellip} with $M_-=\frac12$,
$E_{\rm mom}=2(L^2+3\theta_++1)$, $K=2K_G+1$.
\end{enumerate}
The constants $C_T,C_M$ depend on $(n,q,L,r_\natural)$ and on $C_{S,\rm f}$ (through
$C_T\ge2C_{S,\rm f}$ in~(i) and the term $C_{S,\rm f}\delta^\sharp_{c,\rm f}$ in~(iii)), hence on
$K_G$ and $S^\sharp$. They are fixed after these and before $K$, which depends on $C_T$ through
$K_4$ (Lemma~\ref{lem:BTU4}).
\end{theo}

\begin{proof}
\emph{Step 1: leaf derivatives.} Since $\frac{d}{dz}g(t,z,\xi-z/t)=(\partial_y-t^{-1}\partial_v)g$,
we have
\begin{equation}\label{eq:RGF-Ej}
\nabla^j_zg^{(x)}(z)=(\partial^j_yg)(t,z,v(z))+\mathcal E_j(z),
\end{equation}
where $\mathcal E_j$ is a sum of at most $2^j$ terms $t^{-i}\partial^{j-i}_y\partial^i_vg$ with
$i\ge1$. By Lemma~\ref{lem:RGloc} at order $j\le n+2\le J-1$, we have
$|\mathcal E_j(z)|\le2^jt^{-1}C_{\rm loc}\rho_{v(z)}e^{-r_\natural|z|^2}$.

\emph{Step 2: fibers sampled by the leaf.} We have $|v(z)-\xi|=|z|/t$ and
$|\nabla\log\lambda|\le\frac52L$. For $z\in Z_t$ this implies $\lambda_{v(z)}\le
e^{5c_{\rm fib}/2}\lambda_\xi=\frac43\lambda_\xi$, so $s_{v(z)}(t)\le\frac43\cdot3S=S^\sharp$.
Moreover, we have $\rho_{v(z)}\le e^{c_{\rm fib}}\rho_\xi\le2\rho_\xi$.

\emph{Step 3: (i) on $Z_t$.} We write $\partial^j_yg=\partial^j_yG+\partial^j_y(g-G)$ at the fiber
$v(z)$. By Proposition~\ref{prop:RGcompF}, we have $|\partial^j_y(g-G)|\le C_{S,\rm
f}\rho_{v(z)}\delta^\sharp_{c,\rm f} e^{-(q/2)|z|^2}$. Together with the bound on $\mathcal E_j$ and
$r_\natural\ge4q$, we deduce that the difference in (i) is $\le\big(2C_{S,\rm f}\delta^\sharp_{c,\rm
f}+2^{j+1}C_{\rm loc}t^{-1}\big)e^{-(q/4)|z|^2}$.

\emph{Step 4: (i) off $Z_t$.} Here every term of $\nabla^j_zg^{(x)}$ is bounded by
the same lemma, and $\rho_{v(z)}\le\rho_\xi e^{L|z|}$. Therefore, we have $|\nabla^j_z\hat
g^{(x)}|\le2^jC_{\rm loc}e^{L|z|-r_\natural|z|^2}$. We write $r_\natural|z|^2\ge\frac
q4|z|^2+\frac12r_\natural|z|^2+\frac7{16}r_\natural|z|^2$, which holds since $q\le r_\natural/4$.
Then we have $e^{L|z|-r_\natural|z|^2/2}\le e^{L^2/(2r_\natural)}$ and, on $|z|>c_{\rm fib}t/L$,
$e^{-7r_\natural|z|^2/16}\le e^{-r_\natural c_{\rm fib}^2t^2/(8L^2)}$.

\emph{Step 5: (ii).} On $Z_t$, by Lemma~\ref{lem:RGGbounds}(ii), we have
$\rho_\xi^{-1}|\partial^j_yG(t,z,v(z))|\le2|\partial^j_y\hat G(t,z,v(z))|\le2K_Ge^{-2q|z|^2}$. We
combine this with (i), and use $\|e^{(q/8)|z|^2}e^{-(q/4)|z|^2}\|_{L^2}=C(q)$.

\emph{Step 6: (iii).} We compare $\hat g^{(x)}(z)$ with $\hat g(t,z,\xi):=g(t,z,\xi)/\rho_\xi$. By
the mean value theorem in $v$ and the same lemma, we have
\begin{equation}
|\hat g^{(x)}(z)-\hat g(t,z,\xi)|\le\frac{|z|}t\sup_{v'\in[v(z),\xi]}\frac{|\partial_vg(t,z,v')|}{\rho_\xi}\le\frac{C_{\rm
loc}}t|z|e^{L|z|-r_\natural|z|^2}.
\end{equation}
Since $s_\xi(t)\le3S\le S^\sharp$, the comparison at the fiber $\xi$ yields $|\hat
g(t,z,\xi)-\hat G(t,z,\xi)|\le C_{S,\rm f}\delta^\sharp_{c,\rm f} e^{-(q/2)|z|^2}$. By
Lemma~\ref{lem:RGGbounds}(i), the moments of $\hat G(t,\cdot,\xi)$ are
$(1,y_{0,\xi},|y_{0,\xi}|^2+3\theta_\xi)$. Then (iii) follows by integrating the two differences
against $1,|z|,|z|^2$.

\emph{Step 7: (iv).} By (iii), the mass $M_x$, center $y_0^x=M_x^{-1}\int z\hat g^{(x)}$ and
temperature $\theta^x=(3M_x)^{-1}\int|z-y_0^x|^2\hat g^{(x)}$ differ from $(1,y_{0,\xi},\theta_\xi)$
by at most $C(L)C_M\varepsilon_T$. We choose $\varepsilon_\natural(L,\theta_-)$ so that this does
not exceed $\frac12$, $1$ and $\min\{\frac12\theta_-,1\}$ in the three components. For (G3), on
$Z_t$ we have $\rho_\xi^{-1}G(t,z,v(z))\le2\hat G(t,z,v(z))\le2K_3\jap z^{-8}$, and with (i) for
$j=0$, we get $\sup\jap z^8\hat g^{(x)}\le2K_3+C(q)C_T\varepsilon_T$. For (G2), we apply
Lemma~\ref{lem:RGellip} with $M_-=\frac12$, $\int|z|^2\hat g^{(x)}/M_x\le2(L^2+3\theta_++1)$, and
$\|\hat g^{(x)}\|_{L^\infty}\le2K_G+1$, by (i) with $j=0$, the bound of Step~5 and
$C(q)C_T\varepsilon_T\le1$, as in the proof of (ii). This completes the proof.
\end{proof}

At the effective time $S_*$, this theorem and Lemma~\ref{lem:Aentry} together imply
Corollary~\ref{cor:Aleafentry}, which says that the normalized micro part on the leaf is smaller
than $\delta_1$ in $X(\hat m_\xi)$, uniformly in $\xi$. In other words, the leaf enters the relaxed
regime.

\subsection{Tools for the fibered comparison}\label{ssec:fib}
The tools that we need for the fibered comparison are short-time smoothness of homogeneous
trajectories in Gaussian form, the local factors and the remainder in fiber-normalized Gaussian
norms, bounds on the fibered flow, and the rates and constants of the comparison.

For the short-time smoothness of homogeneous trajectories in Lemma~\ref{lem:RSshortF} we use a
ladder of Gaussian weights, which avoids a polynomial loss in the intermediate terms. We discuss
the derivative count in Remark~\ref{rem:RSshortF}.

\subsubsection{Notation and tools}

For a weight $m$ on $\R^3$, $\|u\|_{X(m)}=\|mu\|_{L^2}$ and
\begin{equation}
\|u\|^2_{Y_2(m)}=\|\jap y^{(\gamma+2)/2}mu\|^2_{L^2}+\|\jap
y^{\gamma/2}\tn(mu)\|^2_{L^2},\qquad\tn=P_y\nabla+\jap y(I-P_y)\nabla,
\end{equation}
as before. Recall that $Q(f,h)=\ab[f]:\nabla^2h-\cb[f]h$, with
$|a(u)|\le|u|^{\gamma+2}$, and $\cb[f]=-8\pi f$ for $\gamma=-3$, $\cb[f]=c*f$ with
$c(u)=-2(\gamma+3)|u|^\gamma$ for $\gamma>-3$.

\begin{lem}[Tools]\label{lem:RFtools}
Let $0<q\le1/8$, $q'\in[q/2,q]$ and $m=e^{q'|y|^2}$. The constants below depend only on $q$.
\begin{enumerate}[label=(\roman*)]
\item $\sup_y|\ab[f](y)|\le C\|f\|_{X(m)}$.
\item \emph{Non-divergence bound.} For all $f,h,\varphi$,
\begin{equation}
|\langle Q(f,h),\varphi\rangle_{L^2(m)}|\le C\|f\|_{X(m)}\Big(\|m\nabla^2h\|_{L^2}+\|mh\|_{L^2}+\|e^{2q'\jap
y}h\|_{L^\infty}\Big)\|\varphi\|_{X(m)} .
\end{equation}
\item \emph{Divergence bound} (Lemma~\ref{lem:QmG} with $\omega=0$). With $\theta=4$, $\theta'=3$,
\begin{equation}
|\langle Q(f,h),\varphi\rangle_{L^2(m)}|\le C\Big(\|f\|_{L^2(\jap
y^\theta)}\|h\|_{Y_2(m)}+\|f\|_{H^1(\jap y^{\theta'})}\|h\|_{X(m)}\Big)\|\varphi\|_{Y_2(m)} .
\end{equation}
\item $\|u\|_{L^2(\jap y^4)}\le C\|u\|_{X(m)}$, $\|u\|_{H^1(\jap y^3)}\le C\|u\|_{Y_2(m)}$, and
$\|u\|_{Y_2(m)}\le\|m\nabla u\|_{L^2}+2\|\jap y mu\|_{L^2}$.
\item
\begin{equation}\label{eq:RF-sob}
\|e^{2q'\jap y}u\|_{L^\infty}\le C\sum_{|\beta|\le2}\|m\partial^\beta u\|_{L^2} .
\end{equation}

\item If $G\ge0$ satisfies (G1)--(G3) of Lemma~\ref{lem:Whypo} and $q\le(8(L^2+3\theta_+))^{-1}$,
the conclusion of that lemma holds for every $q'\in[q/2,q]$ with the same constants,
which we denote by $c_{W,\rm f},C_{W,\rm f}$ (they depend on $K_3,K_G$, and are not the
near-equilibrium constants $c_W,C_W$).
\end{enumerate}
\end{lem}

\begin{proof}
(i) We split $\ab[f](y)=\int a(y-z)f(z)\dd z$ at $|y-z|=1$. For the near part, we have
$\int_{|u|<1}|u|^{2(\gamma+2)}\dd u=4\pi/(2\gamma+7)\le4\pi$, so this part is
$\le(4\pi)^{1/2}\|f\|_{L^2}$. For the far part, we have $|u|^{\gamma+2}\le1$, so this part is
$\le\|f\|_{L^1}\le\|m^{-1}\|_{L^2}\|mf\|_{L^2}$, and $\|m^{-1}\|^2_{L^2}\le(\pi/q)^{3/2}$.

(ii) The $\ab$-part is bounded by $\sup|\ab[f]|\,\|m\nabla^2h\|_{L^2}\|m\varphi\|_{L^2}$, and we
use (i). For $\gamma=-3$ the $\cb$-part is
$8\pi\int|f||h||\varphi|m^2\le8\pi\|h\|_{L^\infty}\|mf\|\|m\varphi\|$. For $\gamma>-3$ we split
$c=c_<+c_>$ at $|u|=1$. Then we have $\|c_<\|_{L^1}=8\pi$ and
$\|c_>\|^2_{L^2}=4(\gamma+3)^2\cdot4\pi/(-2\gamma-3)\le16\pi$, since $-2\gamma-3\ge1$.
\begin{itemize}
\item Far part. We have $|c_>*f|\le\|c_>\|_{L^2}\|f\|_{L^2}$, so it contributes
$\le4\sqrt\pi\|mf\|\,\|mh\|\,\|m\varphi\|$.
\item Near part. For $|y-z|<1$ we have $|y|^2\le|z|^2+2|z|+1$ and $|z|\le|y|+1$, hence $m(y)\le
e^{3q'}e^{2q'|y|}m(z)$. Therefore, we have
\begin{equation}
\iint|c_<(y-z)||f(z)||h(y)||\varphi(y)|m(y)^2\le e^{3q'}\|e^{2q'|\cdot|}h\|_{L^\infty}\big\||c_<|*(m|f|)\big\|_{L^2}\|m\varphi\|_{L^2},
\end{equation}
and by Young's inequality, this is $\le8\pi
e^{3q'}\|e^{2q'\jap\cdot}h\|_{L^\infty}\|mf\|\|m\varphi\|$.
\end{itemize}
(iii) Here $m\in\mathcal W(q/2,q,0)$, and $\theta>\frac72$, $\theta'>\frac94$.

(iv) The first two bounds are Lemma~\ref{lem:TG0}(ii). For the third one, we use $|\tn
w|\le\jap y|\nabla w|$, $\jap y^{1+\gamma/2}\le1$ and $\jap y^{(\gamma+2)/2}\le1$, and then
$\nabla(mu)=m\nabla u+2q'y\,mu$ and $2q'\le\frac14$.

(v) We apply the Sobolev embedding $H^2(\R^3)\subset L^\infty$ to $e^{2q'\jap y}u$, and use
$|\partial^\beta e^{2q'\jap y}|\le Ce^{2q'\jap y}$ for $|\beta|\le2$ and $e^{2q'\jap y}\le
e^{2q'}m$, which follows from $2\jap y\le|y|^2+2$.

(vi) In the proof of Lemma~\ref{lem:Whypo}, $q'$ enters at three places, namely the bound
$q'\int|w|^2G\le M/8$, which holds for $q'\le q$; the radius $R_0$ beyond which
$-\frac72q'M|z|^{\gamma+2}$ dominates $CK|z|^{\gamma+1}$, which is $\le CK/(qM)$ for $q'\ge q/2$;
and the bound $\tilde\zeta\le C(1+q)^2$ on $\{|z|\le R_0\}$. All three are uniform on $[q/2,q]$.
\end{proof}

Note that part~(ii) rests on the weighted $L^2$ bound for $\ab$ of part~(i), compare
\cite[Lemma~4.2(i)]{CM}, and it handles $\cb$ without the Hardy--Littlewood--Sobolev inequality and
without a polynomial loss. We use it again for the term $a=\alpha$ in the proof of
Lemma~\ref{lem:RSmed}.

\subsubsection{Short-time smoothness of homogeneous trajectories}\label{sssec:RSshort}

Let $F$ be the solution of $\partial_\sigma F=Q(F,F)$ given by (GS1), with $F(0)=F_0\in
\mathfrak K^N(C_{\rm dat},\beta,A_N)$. We fix an integer $m$ with $4\le m\le N+2$, and
$0<q\le\min\{\beta/16,1/24\}$.
We define the \emph{ladder}
\begin{equation}
q_k:=q\Big(1-\frac k{2m}\Big),\qquad\mathfrak m_k:=e^{q_k|v|^2}\qquad(0\le k\le m),
\end{equation}
so that $q_0=q$, $q_m=q/2$, and $q_j-q_k\ge q/(2m)$ for $j<k$. Let
\begin{equation}
E^\natural_m(\sigma):=\sum_{|\alpha|\le m}\|\mathfrak m_{|\alpha|}\partial^\alpha
F(\sigma)\|^2_{L^2},\qquad D^\natural_m(\sigma):=\sum_{|\alpha|\le m}\|\partial^\alpha
F(\sigma)\|^2_{Y_2(\mathfrak m_{|\alpha|})} .
\end{equation}
For a \emph{base} $\lambda>0$ we write $q_k(\lambda):=\lambda(1-\frac k{2m})$, $
\mathfrak m_k(\lambda):=e^{q_k(\lambda)|v|^2}$, and $E^{(\lambda)}_m$, $D^{(\lambda)}_m$ for the same sums
with $\mathfrak m_k(\lambda)$ in place of $\mathfrak m_k$, so that $E^\natural_m=E^{(q)}_m$. Since
$\mathfrak m_k(\lambda)$ increases with $\lambda$, so do $E^{(\lambda)}_m$ and $D^{(\lambda)}_m$.
We prepare the proof by three lemmas, namely a viscous approximation, a uniqueness statement, and
an energy identity. In them $\Gamma_c(v):=e^{-c|v|^2}$, and the coefficients are
$\ab[h](v)=\int a(w)h(v-w)\dd w$, and similarly $\bb[h]$, $\cb[h]$, so that
$Q(h,g)=\nabla\cdot(\ab[h]\nabla
g-\bb[h]g)=\ab[h]:\nabla^2g-\cb[h]g$. Here $\partial_ja_{ij}=b_i$ and $\nabla\cdot b=c$ in $
\mathcal D'(\R^3)$, with $c=-8\pi\delta_0$ for $\gamma=-3$, so that $\nabla\cdot\ab[h]=\bb[h]$ and
$\nabla\cdot\bb[h]=\cb[h]$ whenever $h\in L^1\cap L^\infty$.

\begin{lem}[The viscous homogeneous problem]\label{lem:RSvisc}
Let $\nu\in(0,1]$ and $c_0\in(0,\frac18]$ with $\nu c_0\le\frac1{16}$, and let $F_*\ge0$ be
$C^\infty$, with $F_*\le C_*\Gamma_{c_0}$ and $|\partial^{\mathbf k} F_*|\le C_{
\mathbf k}\Gamma_{c_0}$ for every $\mathbf k$. Put $c_\nu(\sigma):=c_0/(1+4\nu c_0\sigma)$.
\begin{enumerate}[label=(\roman*)]
\item There are $\tau_\nu=\tau_\nu(C_*,c_0,\nu)\in(0,1]$ and a solution $F$ of
\begin{equation}
\partial_\sigma F=Q(F,F)+\nu\Delta F,\qquad F(0)=F_*,
\end{equation}
on $[0,\tau_\nu]$ with $0\le F\le2C_*\Gamma_{c_\nu(\sigma)}$, which is $C^\infty$ on
$[0,\tau_\nu]\times\R^3$, and $|\partial^{\mathbf k} F|\le
C(\mathbf k,\delta')\Gamma_{c_\nu(\sigma)-\delta'}$ on $[0,\tau_\nu]$ for every $\mathbf k$ and
$\delta'\in(0,\frac{c_0}8]$, with constants that depend also on $\nu$ and on the $C_{\mathbf k}$.
\item If such a solution exists on $[0,T)$, $T\le1$, with
$A_0:=\sup_{[0,T)}\sup_vF/\Gamma_{c_\nu(\sigma)}<\infty$, then it extends, with the same properties,
to $[0,T+s]$ for some $s>0$. Moreover $F\le C_*e^{K_{\rm b}\sigma}\Gamma_{c_\nu(\sigma)}$, where
$K_{\rm b}=C(1+M_F)$ and
\begin{equation}\label{eq:MF}
M_F:=\sup_{[0,T)}\sup_v\big(|\ab[F]|+|\cb[F]|+\ab[F](v)v\cdot v\big) .
\end{equation}

\item $\int F$ and $\int vF$ are constant, and $\int|v|^2F(\sigma)=\int|v|^2F_*+6\nu\sigma\int F_*$.
\end{enumerate}
\end{lem}

\begin{proof}
The proofs of Lemmas~\ref{lem:LWlin} and~\ref{lem:LWregex}(i)--(iv), in which the uniqueness~(ii)
enters the continuation~(iv), apply with $\R^6$, $D$, $\eta_v(\Delta_y+\Delta_v)$, $t_*$ and the
barrier $\eps\rho_0(v)e^{-c|y|^2}$ replaced by $\R^3_v$, $\nabla_v$, $\nu\Delta$, $0$ and
$\Gamma_c(v)$. We only need to discuss the points where the argument uses the structure.
\begin{itemize}
\item \emph{Coefficients.} If $0\le h\le A\Gamma_{c_h}$ with $c_h\ge c_0/2$, then
$|\ab[h]|+|\bb[h]|+|\cb[h]|\le C(c_0)A$, because the kernels $|w|^{\gamma+2}$, $|w|^{\gamma+1}$ and,
for $\gamma>-3$, $|w|^\gamma$ are locally integrable and $\Gamma_{c_h}$ is integrable. Moreover,
$-\cb[h]\ge0$, and, since $a(u)u=0$, we have $\ab[h](v)v\cdot v=\int a(v-w)w\cdot w\,h(w)\dd w\le
C(c_0)A$. This null structure is used in parts~(i) and~(iii) of the first lemma. If $h$ is H\"older
continuous in $(\sigma,v)$, so are $\ab[h]$ and $\cb[h]$, by the same splitting of the kernels.
Here there is no transport, so there is no restriction at $\gamma=-2$.
\item \emph{Barrier.} With $c=c_\nu(\sigma)$, we have $\partial_\sigma\Gamma_c=4\nu
c^2|v|^2\Gamma_c$ and
$\nu\Delta\Gamma_c=\nu(4c^2|v|^2-6c)\Gamma_c\le4\nu c^2|v|^2\Gamma_c$, so the viscous term is
compensated exactly. Since $\ab[h]\ge0$, we get $\ab[h]:\nabla^2\Gamma_c=\ab[h]:(4c^2v\otimes
v-2cI)\Gamma_c\le4c^2\ab[h]v\cdot v\,\Gamma_c$, and $-\cb[h]\Gamma_c\le M_h\Gamma_c$. Thus
$C_*e^{K_{\rm b}\sigma}\Gamma_{c_\nu}$ is a supersolution for $K_{\rm b}=C(1+M_h)$, and the
comparison principle of part~(ii) of that lemma (with the function $\rho^{-2}\sup|
\mathfrak w|(|v|^2+K_0\sigma)$) yields the bound in~(ii) and the bound $F\le2C_*\Gamma_{c_\nu}$
of~(i), for $\tau_\nu\le(\log2)/K_{\rm b}$.
\item \emph{Regularity up to $\sigma=0$.} Since the data are $C^\infty$, with Gaussian bounds on
all derivatives, the induction of the proof of Lemma~\ref{lem:LWregex}(i) can be carried up to
$\sigma=0$ at every order, with \cite[Thm.~4.1]{Lor00} in the local form of the paragraph before
Lemma~\ref{lem:LWlin}, as in the proof of Lemma~\ref{lem:CUtstar}(i). In $\R^3_v$ there is no
restriction on the order. This implies the derivative bounds of~(i), and~(ii) is the continuation
argument of parts~(iii), (iv) of the second lemma.
\item (iii) We integrate the equation against $1$, $v$ and $|v|^2$ over a ball of radius $\varrho$
and let $\varrho\to\infty$, using the Gaussian bounds on $F$ and its derivatives. Recall that $Q$
conserves these moments, $\int\Delta F=\int v\Delta F=0$, and $\int|v|^2\Delta F=6\int F$.
\end{itemize}
\end{proof}

\begin{lem}[Uniqueness against a regular solution]\label{lem:RSuniq}
Let $0\le s<s'$, and let $F_1$, $F_2$ be solutions of $\partial_\sigma F=Q(F,F)$ on $[s,s']$ that
are continuous on $[s,s']\times\R^3$, $C^1$ in $\sigma$ and $C^2$ in $v$ on $(s,s']\times\R^3$ (with
$\partial_\sigma F_i$ and $\nabla^2F_i$ continuous there), with $F_1(s)=F_2(s)$. Assume that
$F_1\ge0$, that $F_1+|F_2|\le C'e^{-\beta'|v|^2}$ on $[s,s']$ for some $C',\beta'>0$, and that
$\sup_{[s,s']}\|\jap v^2\nabla^2F_2\|_{L^2}<\infty$. Then $F_1=F_2$ on $[s,s']$.
\end{lem}

\begin{proof}
Let $d:=F_1-F_2$, $\Psi:=\jap v^4$ and $Y(\sigma):=\int\Psi d(\sigma)^2$. By the Gaussian bounds,
$Y$ is continuous on $[s,s']$ (dominated convergence), and $Y(s)=0$. By bilinearity, we have
$\partial_\sigma d=Q(F_1,d)+Q(d,F_2)$. The coefficients of $F_1$ are bounded, that is,
$|\ab[F_1]|+|\bb[F_1]|+|\cb[F_1]|\le C$, as in Lemma~\ref{lem:RSvisc}, and $\ab[F_1]\ge0$ is
Lipschitz, with $\nabla\cdot\ab[F_1]=\bb[F_1]$ and $\nabla\cdot\bb[F_1]=\cb[F_1]$.

Let $\zeta_0\in C^\infty_c(\R^3)$ be radial and nonincreasing in $|v|$, with $\zeta_0=1$ on the unit
ball, put $\zeta:=\zeta_0(\cdot/R)^2$, and let $s<s_1<\sigma\le s'$. On $[s_1,\sigma]$ the function
$d$ is $C^1$ in $\sigma$ and $C^2$ in $v$, and $\zeta\Psi d$ has compact support, so we have
\begin{equation}
\int\zeta\Psi d^2\Big|_{s_1}^{\sigma}=2\int_{s_1}^\sigma\!\!\int\zeta\Psi
d\big(Q(F_1,d)+Q(d,F_2)\big) .
\end{equation}
In the first term we write $Q(F_1,d)=\nabla\cdot(\ab[F_1]\nabla d-\bb[F_1]d)$ and integrate by
parts. With $\phi:=\zeta\Psi$, we obtain
\begin{equation}
\int\phi\,d\,Q(F_1,d)=-\int\phi\,\ab[F_1]\nabla d\cdot\nabla d-\int d\,\ab[F_1]\nabla
d\cdot\nabla\phi-\frac12\int\big(\phi\,\cb[F_1]+\bb[F_1]\cdot\nabla\phi\big)d^2+\int\bb[F_1]\cdot\nabla\phi\,d^2
,
\end{equation}
where the term with $\bb[F_1]\phi\,d\cdot\nabla d=\frac12\bb[F_1]\phi\cdot\nabla(d^2)$ was
integrated by parts once more. By the Cauchy--Schwarz inequality for the nonnegative form
$\ab[F_1]$, the first two terms are bounded by
$\frac14\int\ab[F_1]\nabla\phi\cdot\nabla\phi\,d^2/\phi$
on $\{\phi>0\}$. Since $\nabla\phi=\Psi\nabla\zeta+\zeta\nabla\Psi$, $|\nabla\zeta|^2\le
CR^{-2}\zeta$ and $|\nabla\Psi|\le C\Psi$, we have $|\nabla\phi|^2/\phi+|\nabla\phi|\le C\Psi$.
Hence the first term is bounded by $CY$, with $C$ independent of $R$.

For the second term, we have $Q(d,F_2)=\ab[d]:\nabla^2F_2-\cb[d]F_2$. By the proof of
Lemma~\ref{lem:RFtools}(i) with $q=\frac18$, which provides a bound in terms of $\|d\|_{L^2}$ and
$\|d\|_{L^1}$, we have $|\ab[d]|\le
C(\|d\|_{L^2}+\|d\|_{L^1})\le CY^{1/2}$, since $\|d\|_{L^1}\le\|\jap v^{-2}\|_{L^2}\|\jap
v^2d\|_{L^2}$. Hence we have $\int\phi|d||\ab[d]||\nabla^2F_2|\le CY^{1/2}\cdot Y^{1/2}\|\jap
v^2\nabla^2F_2\|_{L^2}$. For $\gamma=-3$, we have $\cb[d]=-8\pi d$, and $\int\phi d^2|F_2|\le CY$.
For $\gamma>-3$, we split $c=c_<+c_>$ at $|w|=1$, as in the proof of part~(ii) of that lemma. Then
$\|c_<*d\|_{L^2}\le8\pi\|d\|_{L^2}$ and $|c_>*d|\le C\|d\|_{L^1}$, and with $\|\Psi
F_2\|_{L^\infty}+\|\Psi^{1/2}F_2\|_{L^2}\le C$, the term is bounded by $CY$.

Hence we have $\int\zeta\Psi d^2(\sigma)\le\int\zeta\Psi d^2(s_1)+C\int_{s_1}^\sigma Y$. Letting
$R\to\infty$, by monotone convergence, and then $s_1\downarrow s$, we obtain $Y(\sigma)\le
C\int_s^\sigma Y$. By Gronwall's inequality, we deduce that $Y\equiv0$. This completes the proof.
\end{proof}

\begin{lem}[The energy identity]\label{lem:RSident}
Let $m\ge4$, $\lambda\in(0,\frac1{24}]$ and $0\le s<s'$, and let $F\ge0$ solve $\partial_\sigma
F=Q(F,F)$ on $[s,s']$, with $F\in C([s,s'];H^{m-1}_{\rm loc})$ and
\begin{equation}
\sup_{[s,s']}E^{(\lambda)}_m+\int_s^{s'}D^{(\lambda)}_m<\infty .
\end{equation}
Then, for $|\alpha|=k\le m$, the function $\sigma\mapsto\|\mathfrak m_k(\lambda)\partial^\alpha
F(\sigma)\|^2_{L^2}$ is absolutely continuous on $[s,s']$, and for a.e.\ $\sigma$ its derivative is
\begin{equation}\label{eq:RSident}
2\sum_{a\le\alpha}\binom\alpha a\langle Q(\partial^aF,\partial^{\alpha-a}F),\partial^\alpha
F\rangle_{L^2(\mathfrak m_k(\lambda))},
\end{equation}
each pairing being integrable in $\sigma$. The pairings are taken in the divergence form
$-\int\big(\ab[f]\nabla g-\bb[f]g\big)\cdot\nabla\big(\mathfrak m_k^2h\big)$, which is the form of
Lemma~\ref{lem:QmG}, and which coincides with the non-divergence form of Lemma~\ref{lem:RFtools}(ii)
when $g$ has two derivatives in $L^2(\mathfrak m_k^2)$.
\end{lem}

\begin{proof}
We write $\mathfrak m:=\mathfrak m_k(\lambda)$, $u:=\partial^\alpha F$ and $I:=[s,s']$. Since
$\partial^\alpha$ is a derivation of $Q$, we have in $\mathcal D'((s,s')\times\R^3)$
\begin{equation}
\partial_\sigma u=\nabla\cdot G_1,\qquad G_1:=\sum_{a\le\alpha}\tbinom\alpha
a\big(\ab[\partial^aF]\nabla\partial^{\alpha-a}F-\bb[\partial^aF]\partial^{\alpha-a}F\big).
\end{equation}
\emph{Integrability.} Since $q_j(\lambda)\ge\lambda/2$ and $\mathfrak m_j(\lambda)\ge c\jap v^4$,
the energy bounds $\partial^\beta F$ in $L^2(\jap v^8)$ for $|\beta|\le m$, hence in $L^1\cap L^2$,
and in $L^\infty$ for $|\beta|\le m-2$, by Sobolev's inequality. Therefore, by
Lemma~\ref{lem:RFtools}(i), we have $|\ab[\partial^aF]|\le C$ for every $a\le\alpha$. If $|a|\le
m-2$, then $|\bb[\partial^aF]|\le
C$, and $\partial^{\alpha-a}F\in L^2$. If $|a|\ge m-1$, then $|\alpha-a|\le1$,
$\partial^{\alpha-a}F$ is bounded by a Gaussian, by Sobolev's inequality with the weights $
\mathfrak m_j(\lambda)$, $j\le3$, and $\bb[\partial^aF]\in L^2+L^\infty$, since $b$ is the sum of a function
in $L^1$ and a bounded one. Finally, $\nabla\partial^{\alpha-a}F\in L^2$ for $a\ne0$, and $\nabla
u\in L^2(I;L^2_{\rm loc})$ by the bound on $\int D^{(\lambda)}_m$. Thus $G_1\in L^2(I;L^2_{\rm
loc})$, $\nabla u\in L^2(I;L^2_{\rm loc})$, and $u\in L^\infty(I;L^2_{\rm loc})\cap C(I;H^{-1}_{\rm
loc})$.

\emph{The identity.} Lemma~\ref{lem:CUident}(i) holds, with the same proof, in $\R^3$ with $\nabla$
in place of $D$ and $G_0=0$. Applying it with $\Psi:=\zeta\mathfrak m^2$, where
$\zeta=\zeta_0(\cdot/R)^2$ is as in the proof of Lemma~\ref{lem:RSuniq}, we have
\begin{equation}
\int\zeta\mathfrak m^2u^2\Big|_{\sigma_1}^{\sigma_2}=-2\int_{\sigma_1}^{\sigma_2}\!\!\int\nabla(\zeta
\mathfrak m^2u)\cdot G_1 .
\end{equation}
Let $R\to\infty$. The left side converges by monotone convergence. On the right, we write
$\nabla(\zeta\mathfrak m^2u)=\zeta\nabla(\mathfrak m^2u)+\mathfrak m^2u\nabla\zeta$. The integrand
$\nabla(\mathfrak m^2u)\cdot G_1$ is integrable over $I\times\R^3$.

For the terms with $|a|\le m-1$ (with $f=\partial^aF$, $g=\partial^{\alpha-a}F$, $h=u$), the
decomposition $A_1+\dots+B_2$ in the proof of Lemma~\ref{lem:QmG} is pointwise, and its coefficient
bounds, $|\ab[f]\xi\cdot\zeta|\le C\|f\|_{L^2(\jap v^4)}\jap v^\gamma|\tilde\xi||\tilde\zeta|$ with
$\tilde\xi=P_v\xi+\jap v(I-P_v)\xi$, and $|\bb[f]|\le C\|f\|_{H^1(\jap v^3)}\jap v^{\gamma+1}$,
bound it by products of the integrands of $\|g\|_{Y_2(\mathfrak m)}$ or $\|g\|_{X(\mathfrak m)}$
with those of $\|u\|_{Y_2(\mathfrak m)}$. These products are integrable over $I\times\R^3$, since
$\|u\|_{Y_2(\mathfrak m)}\in L^2(I)$ and, for $a\ne0$, $\|g\|_{Y_2(\mathfrak m)}$ is bounded on $I$
as in~(L2).

For $a=\alpha$ (so $k=m$), we have $|\ab[u]|\le C(\|u\|_{L^2}+\|u\|_{L^1})$ as in the uniqueness
proof, $\bb[u]\in L^\infty(I;L^2+L^\infty)$, and $F+|\nabla F|\le Ce^{-r|v|^2}$
with $r\in(q_m(\lambda),q_3(\lambda))$, by Sobolev's inequality with the weights $
\mathfrak m_j(\lambda)$, $j\le3$ (here $m\ge4$). With $|\nabla(\mathfrak m^2u)|\le
\mathfrak m(|\nabla(\mathfrak mu)|+2\lambda|v|\mathfrak m|u|)$ this term is bounded by
$C\Theta\big(\jap v^{\gamma/2}|\tn(\mathfrak mu)|+\jap v^{(\gamma+2)/2}\mathfrak m|u|\big)$, with
$\Theta:=\jap v^{-\gamma/2}\mathfrak m(|\nabla F|+|\bb[u]|F)\in L^\infty(I;L^2)$. Thus the first
part converges by dominated convergence.

The same coefficient bounds imply $|\ab[f]\xi|\le C\jap v^{\gamma+1}|\tilde\xi|$, and $\jap
v^{\gamma+1}\le\jap v^{\gamma/2}$. It follows that $\mathfrak m|G_1|\in L^2(I\times\R^3)$ and
$\mathfrak m|u|\in L^2(I\times\R^3)$, while $|\nabla\zeta|\le C/R$. Hence the second part tends to
$0$.

The limit identity is the integrated form of the claim, with each pairing integrable in $\sigma$.
The two forms of the pairings coincide by one integration by parts, with the cutoff $\zeta$, when
$g$ has two derivatives in $L^2(\mathfrak m^2)$.
\end{proof}

\begin{lem}[Short-time Gaussian Sobolev bounds, full form]\label{lem:RSshortF}
\begin{enumerate}[label=(\roman*)]
\item The initial energy is bounded: $E^\natural_m(0)\le C(m,C_{\rm dat},\beta,A_N,q)$.
\item There are $\sigma_0=\sigma_0(m,q,E^\natural_m(0))>0$ and $c=c(m,q,E^\natural_m(0))>0$ such
that
\begin{equation}
E^\natural_m(\sigma)+c\int_0^\sigma D^\natural_m\le2E^\natural_m(0)\qquad\text{on }[0,\sigma_0].
\end{equation}
\item On $[0,\sigma_0]$,
\begin{equation}\label{eq:RSshort-pt}
\sum_{|\alpha|\le m-2}\|e^{(q/4)|v|^2}\partial^\alpha F\|_{L^\infty}\le C\,E^\natural_m(0)^{1/2} .
\end{equation}

\item \emph{(Restart.)} Let $\sigma'\ge0$ and $q^+\in(q,2q]$. Suppose that $F(\sigma')\le
C_*e^{-c_*|v|^2}$ with $c_*\ge8q^+$, that $E^{(q^+)}_m(\sigma')<\infty$, and that $F\le
C'e^{-\beta'|v|^2}$ on $[\sigma',\sigma'+1]$ for some $C',\beta'>0$. Then (ii) and (iii) hold on
$[\sigma',\sigma'+\sigma_0]$, with $E^\natural_m(\sigma')$ in place of $E^\natural_m(0)$ and
$\sigma_0=\sigma_0(m,q,E^\natural_m(\sigma'))$. Moreover $F\in
C([\sigma',\sigma'+\sigma_0];H^{m-1}_{\rm loc})$, $F$ is $C^1$ in $\sigma$ and $C^2$ in $v$ on
$[\sigma',\sigma'+\sigma_0]\times\R^3$, and $\int_{\sigma'}^{\sigma'+\sigma_0}D^\natural_m<\infty$.
\end{enumerate}
The function $\sigma_0(m,q,\cdot)$ is nonincreasing. Moreover, $\sigma_0$ and $c$ are bounded below, and the
constant $C$ of~(iii) is bounded above, uniformly for $q$ in compact subsets of $(0,\frac1{24}]$ and
bounded values of $E^\natural_m(\sigma')$.
\end{lem}

\begin{proof}
\emph{Step 0: the data.} Lemma~\ref{lem:RGinterp} with $2N+4$ polynomially weighted derivatives,
i.e.\ with $N+2$ in place of $N$, and with the Maxwellian bound $F_0\le C_{\rm dat}e^{-\beta|v|^2}$,
yields $|\nabla^jF_0|\le Ce^{-\beta|v|^2/4}$ for $j\le N+2$. Since $q\le\beta/16$, we have $\|
\mathfrak m_k\nabla^jF_0\|_{L^2}\le C\|e^{-3\beta|v|^2/16}\|_{L^2}$. This is (i). In the same way,
we have $E^{(2q)}_m(0)\le C(m,C_{\rm dat},\beta,A_N,q)$, since $2q\le\beta/8$. Therefore, the
hypotheses of~(iv) hold at $\sigma'=0$, with $q^+=2q$, $c_*=\beta\ge16q$, and with the Gaussian
bound on $[0,1]$ from Lemma~\ref{lem:GS1q}, and (ii), (iii) are the case $\sigma'=0$ of~(iv). Thus
we only need to prove~(iv). Steps~1--4 are the a priori estimate, written for $F$ on
$[\sigma',\sigma'+\sigma]$, with $\sigma'=0$ to simplify the notation. We apply them in Step~5 to
viscous approximations, for which the identity of Step~2 holds classically. In Step~6 we pass to
the limit, and in Step~7 we identify the limit with $F$.

\emph{Step 1: four ladder facts.} We write $E=E^\natural_m$. For $|\alpha|\le m$ and $k\le m$, we
have the following.
\begin{itemize}
\item[(L1)] if $|\alpha|\le k$, then $\|\partial^\alpha F\|_{X(\mathfrak m_k)}\le E^{1/2}$, since
$q_k\le q_{|\alpha|}$;
\item[(L2)] if $|\alpha|<k$, then $\|\partial^\alpha F\|_{Y_2(\mathfrak m_k)}\le C_LE^{1/2}$.
Indeed, by Lemma~\ref{lem:RFtools}(iv), $\|u\|_{Y_2(\mathfrak m_k)}\le\|\mathfrak m_k\nabla
u\|+2\|\jap v\mathfrak m_ku\|$. With $u=\partial^\alpha F$, the first term is $\le\|
\mathfrak m_{|\alpha|+1}\nabla\partial^\alpha F\|$ because $|\alpha|+1\le k$. The second one is
$\le2\sup_v\big(\jap ve^{-(q_{|\alpha|}-q_k)|v|^2}\big)\|\mathfrak m_{|\alpha|}\partial^\alpha F\|$,
and $q_{|\alpha|}-q_k\ge q/(2m)$;
\item[(L3)] $\|\partial^\alpha F\|_{L^2(\jap v^4)}\le CE^{1/2}$, and $\|\partial^\alpha
F\|_{H^1(\jap v^3)}\le CE^{1/2}$ if $|\alpha|\le m-1$, since $\jap v^4\le Ce^{(q/2)|v|^2}\le
C\mathfrak m_j$ for all $j$;
\item[(L4)] $\|\mathfrak m_m\nabla^2F\|+\|\mathfrak m_mF\|+\|e^{2q_m\jap v}F\|_{L^\infty}+\|\jap
v^8F\|_{L^\infty}\le CE^{1/2}$. The first two follow from $q_m\le q_2$. The third is
part~(v) of the same lemma with $q'=q_m$. The fourth follows from the third.
\end{itemize}

\emph{Step 2: the identity.} Since $\partial^\alpha$ is a derivation of $Q$, for $|\alpha|=k$, we
have
\begin{equation}
\frac12\frac{d}{d\sigma}\|\mathfrak m_k\partial^\alpha F\|^2_{L^2}=\sum_{a\le\alpha}\binom\alpha
a\big\langle Q(\partial^aF,\partial^{\alpha-a}F),\partial^\alpha F\big\rangle_{L^2(\mathfrak m_k)} .
\end{equation}
In Steps~3--4 we assume this identity and the finiteness of $E$, which hold for the approximations
of Step~5.

\emph{Step 3: the terms.} Assume that $E\le B:=4E^\natural_m(0)+1$ on $[0,\sigma]$. We fix
$|\alpha|=k$.
\begin{itemize}
\item[(A)] $a=0$. Mass, momentum and energy are conserved, so $F$ satisfies (G1) of
Lemma~\ref{lem:Whypo} with $M=1$, $y_0=0$, $\theta=1$, and $q\le1/24=(8\cdot3)^{-1}$. (G3) holds
with $K\le CB^{1/2}$ by (L4). (G2) holds by Lemma~\ref{lem:RGellip} with $M_-=1$, $E_{\rm mom}=3$
and $\|F\|_{L^\infty}\le K$. Then by Lemma~\ref{lem:RFtools}(vi), this term is $\le-c_{W,\rm
f}\|\partial^\alpha F\|^2_{Y_2(\mathfrak m_k)}+C_{W,\rm f}\|\partial^\alpha F\|^2_{X(
\mathfrak m_k)}$, with $c_{W,\rm f},C_{W,\rm f}$ depending on $(q,B)$.
\item[(B)] $a\ne0$, and either $k<m$ or $a\ne\alpha$. We apply part~(iii) of the tools lemma with
$f=\partial^aF$, $h=\partial^{\alpha-a}F$, $\varphi=\partial^\alpha F$ and weight $\mathfrak m_k$,
and we have
\begin{itemize}
\item $\|f\|_{L^2(\jap v^4)}\le CE^{1/2}$ by (L3);
\item $\|h\|_{Y_2(\mathfrak m_k)}\le C_LE^{1/2}$ by (L2), since $|\alpha-a|<k$;
\item $\|f\|_{H^1(\jap v^3)}\le CE^{1/2}$ by (L3), since $|a|\le m-1$ in both cases;
\item $\|h\|_{X(\mathfrak m_k)}\le E^{1/2}$ by (L1).
\end{itemize}
Thus this term is $\le CE\,\|\partial^\alpha F\|_{Y_2(\mathfrak m_k)}$.
\item[(C)] $k=m$ and $a=\alpha$. By part~(ii) of the tools lemma with weight $\mathfrak m_m$ and
(L4), this term is bounded by
\begin{equation}
CE^{1/2}\|\partial^\alpha F\|^2_{X(\mathfrak m_m)}\le CE^{3/2}.
\end{equation}
\end{itemize}

\emph{Step 4: closing.} There are at most $C_m$ pairs $(\alpha,a)$. By Young's inequality, we have
$CE\sum_\alpha\|\partial^\alpha F\|_{Y_2}\le\frac{c_{W,\rm f}}2D^\natural_m+CE^2$. Hence, as long as
$E\le B$, we get
\begin{equation}
\frac12\frac{d}{d\sigma}E+\frac{c_{W,\rm f}}2D^\natural_m\le C_{W,\rm f}E+CE^2+CE^{3/2}\le
C_BE,\qquad C_B:=C_{W,\rm f}+1+C(B+B^{1/2}).
\end{equation}
The term $1$ is needed only in Step~5. Let $\sigma_0:=\log(3/2)/(2C_B)$. Then we have
$E(\sigma)\le e^{2C_B\sigma}E(0)\le\frac32E(0)<B$ on $[0,\sigma_0]$, and by continuity the
assumption $E\le B$ holds there. Integrating, we obtain
\begin{equation}
E(\sigma)+c_{W,\rm f}\int_0^\sigma D^\natural_m\le E(0)+2C_B\int_0^\sigma E\le
E(0)\big(1+\tfrac32\log\tfrac32\big)\le2E(0).
\end{equation}
The last bound is the a priori form of~(ii), and we prove (ii) itself in Steps~5--7. For (iii), by
Sobolev's inequality, we have
\begin{equation}
\|e^{(q/4)|v|^2}\partial^\alpha F\|_{L^\infty}\le C\sum_{|\beta|\le2}\|\jap
v^2e^{(q/4)|v|^2}\partial^{\alpha+\beta}F\|_{L^2},
\end{equation}
and $\jap v^2e^{(q/4)|v|^2}\le C\mathfrak m_m$.

\emph{Step 5: viscous approximations.} Let $F_*:=F(\sigma')$, normalized by the conservation laws,
and we shift time so that $\sigma'=0$. Let $\varphi\in C^\infty_c(B_1)$ be radial, $\varphi\ge0$,
$\int\varphi=1$, $\varphi_\epsilon:=\epsilon^{-3}\varphi(\cdot/\epsilon)$ for $\epsilon\in(0,1]$,
$M_\varphi:=\int|w|^2\varphi$, $\ell_\epsilon:=(1+\epsilon^2M_\varphi/3)^{1/2}\in[1,2]$, and
\begin{equation}
F^\epsilon_*(v):=\ell_\epsilon^3\,(\varphi_\epsilon*F_*)(\ell_\epsilon v) .
\end{equation}
Then $F^\epsilon_*\ge0$ is $C^\infty$ and normalized. Indeed, by the radial symmetry of $\varphi$,
the mass and the momentum are those of $F_*$, and
$\int|v|^2F^\epsilon_*=\ell_\epsilon^{-2}(3+\epsilon^2M_\varphi)=3$. Since $|\ell_\epsilon
v-z|\ge|v|-1$ for $|z|\le1$, we have $F^\epsilon_*\le8e^{c_*}C_*\Gamma_{c_*/2}$, and every
derivative of $F^\epsilon_*$ is bounded by $C(\epsilon,\kappa)\Gamma_{c_*/2}$. \emph{Convergence of
the data.} Let $q^+_k:=q_k(q^+)$, so that $q^+_k-q_k\ge(q^+-q)/2>0$. Since $e^{q_k|v|^2}\le
Ce^{q^+_k(|v|-1)_+^2}\le Ce^{q^+_k|\ell_\epsilon v-z|^2}$ for $|z|\le1$, the maps
$u\mapsto\ell_\epsilon^{3+|\alpha|}(\varphi_\epsilon*u)(\ell_\epsilon\,\cdot)$ are bounded from
$L^2(e^{2q^+_k|v|^2})$ to $L^2(e^{2q_k|v|^2})$, uniformly in $\epsilon$, and they converge to the
identity on $C_c$, hence strongly on $L^2(e^{2q^+_k|v|^2})$. With $\partial^\alpha
F^\epsilon_*=\ell_\epsilon^{3+|\alpha|}(\varphi_\epsilon*\partial^\alpha F_*)(\ell_\epsilon\,\cdot)$
and $E^{(q^+)}_m(F_*)<\infty$, we deduce that $E^\natural_m(F^\epsilon_*)\to E^\natural_m(F_*)$.

Let $c_0:=\min\{c_*/2,\frac18\}$, $\nu:=\epsilon$ (so $\nu c_0\le\frac1{16}$ for
$\epsilon\le\frac12$), and let $F_\nu$ be the solution of Lemma~\ref{lem:RSvisc} with the data
$F^\epsilon_*$. By part~(iii) of that lemma, $F_\nu$ has mass $1$ and momentum $0$, and its
temperature is $\theta(\sigma)=1+2\nu\sigma\in[1,2]$ for $2\nu\sigma\le1$. We use the ladder
\begin{equation}
q^\nu_k(\sigma):=\frac{q_k}{1+2\nu\sigma},\qquad\mathfrak m^\nu_k:=e^{q^\nu_k(\sigma)|v|^2},
\end{equation}
and the energies $E_\nu$, $D_\nu$ of $F_\nu$ with these weights. The derivatives of $F_\nu$ have
Gaussian bounds at the rate $c_\nu-\frac{c_0}8\ge\frac{27}{40}c_0>q$ up to $\sigma=0$ (since
$c_*\ge8q^+>8q$ and $q\le\frac1{24}$), so $E_\nu$ is continuous on the interval of existence, $C^1$
there, with $E_\nu(0)=E^\natural_m(F^\epsilon_*)$, and the identity of Step~2 holds with two
additional terms, namely the viscous term and the time derivative of the weight. For
$u=\partial^\alpha F_\nu$ and $q=q^\nu_k(\sigma)$, $\mathfrak m=\mathfrak m^\nu_k$, we have
\begin{equation}
\nu\langle\mathfrak m^2u,\Delta u\rangle=-\nu\|\mathfrak m\nabla
u\|^2+\nu\!\int(6q+8q^2|v|^2)\mathfrak m^2u^2,\qquad\frac12\partial_\sigma
\mathfrak m^2=\partial_\sigma q\,|v|^2\mathfrak m^2=-\frac{2\nu q}{1+2\nu\sigma}|v|^2\mathfrak m^2 .
\end{equation}
Since $4q(1+2\nu\sigma)=4q_k\le1$, the terms with $|v|^2$ add up to a nonpositive quantity, and the
rest is $\le6\nu q\|\mathfrak mu\|^2\le\|\mathfrak mu\|^2$.

Steps~1, 3 and~4 apply to $F_\nu$, for $2\nu\sigma\le1$, with the following changes. The rungs
$q^\nu_k(\sigma)$ lie in $[q/4,q]$ and their gaps are $\ge q/(4m)$, so (L1)--(L4) hold with the
constants changed by absolute factors. In~(A), (G1) holds with $M=1$, $y_0=0$ and
$\theta=\theta(\sigma)\le2$, (G2) holds by the same ellipticity lemma with $E_{\rm mom}=6$, and (G3)
holds as before.

In the proof of Lemma~\ref{lem:Whypo} the exponent appears in the three places listed in the proof
of part~(vi) of the tools lemma, namely the bound $q'\int|w|^2G\le M/8$, which holds since
$q^\nu_k(\sigma)\cdot3\theta(\sigma)=3q_k\le\frac18$, the radius $R_0\le CK/(qM)$, and the bound
$\tilde\zeta\le C(1+q)^2$. The last two are uniform for $q'\in[q/4,q]$. Therefore,
(A) holds with constants depending only on $(q,B)$, which we denote again by $c_{W,\rm f}$ and
$C_{W,\rm f}$.

The viscous remainder is controlled by the term $1$ in $C_B$. Finally, as long as $E_\nu\le B$, the
coefficients of $F_\nu$ satisfy $M_{F_\nu}\le C(B)$, by (L4) and $\int F_\nu=1$, so
Lemma~\ref{lem:RSvisc}(ii) implies $F_\nu\le8e^{c_*}C_*e^{K_{\rm b}(B)\sigma}\Gamma_{c_\nu}$, and
$F_\nu$ can be continued.

We take $B:=4E^\natural_m(F_*)+2$, which is $\ge4E_\nu(0)+1$ for $\epsilon$ small, and
$\sigma_0:=\min\{\log(3/2)/(2C_B),\frac12\}$. Then $\nu\sigma_0\le\frac12$ for $\nu\le1$. By Step~4
and the continuation, $F_\nu$ exists on $[0,\sigma_0]$, and we have
\begin{equation}
E_\nu(\sigma)+c_{W,\rm f}\int_0^\sigma D_\nu\le2E_\nu(0)\qquad\text{on }[0,\sigma_0],
\end{equation}
with $\sigma_0$ and $c_{W,\rm f}$ independent of $\nu$ and $\epsilon$. Moreover, $\sigma_0$ is
nonincreasing in $E^\natural_m(F_*)$.

\emph{Step 6: the limit.} By (L3)--(L4), $F_\nu$ is bounded in $L^\infty(0,\sigma_0;H^m(
\mathcal C))$ for every compact $\mathcal C$, uniformly in $\nu$, and so is $\partial_\sigma
F_\nu=Q(F_\nu,F_\nu)+\nu\Delta F_\nu$ in $L^\infty(0,\sigma_0;H^{m-2}(\mathcal C))$. Indeed, we
write $Q(F,F)=\ab[F]:\nabla^2F-\cb[F]F$, with $\partial^\kappa\ab[F]=\ab[\partial^\kappa F]$ bounded
for $|\mathbf k|\le m$ by part~(i) of the tools lemma, and $\partial^{
\mathbf k}\cb[F]=\cb[\partial^{\mathbf k} F]$ bounded for $|\mathbf k|\le m-2$, since $\partial^{\mathbf k}
F\in L^1\cap L^\infty$ (for $\gamma=-3$, $\cb[F]F=-8\pi F^2$).

By the Aubin--Lions--Simon lemma \cite[Cor.~4]{Sim87}, a subsequence $\nu\to0$ converges in
$C([0,\sigma_0];H^{m-1}(\mathcal C))$ for every compact $\mathcal C$, and by interpolation in
$C([0,\sigma_0];H^{m-s}(\mathcal C))$ for $s\in(0,\frac12)$, hence in $C([0,\sigma_0];C^2(
\mathcal C))$, since $m-s>\frac72$.

Let $\tilde F$ be the limit. By the uniform Gaussian bound $\|e^{2q^\nu_m\jap
v}F_\nu\|_{L^\infty}\le CB^{1/2}$ of (L4), the coefficients converge locally uniformly, which
implies that $\tilde F(\sigma)=F_*+\int_0^\sigma Q(\tilde F,\tilde F)$ with a continuous integrand.
Hence $\tilde F\ge0$
is $C^1$ in $\sigma$ and $C^2$ in $v$ on $[0,\sigma_0]\times\R^3$, and solves $\partial_\sigma\tilde
F=Q(\tilde F,\tilde F)$ with $\tilde F(0)=F_*$.

The weights $\mathfrak m^\nu_k$ increase to $\mathfrak m_k$ as $\nu\to0$. By weak lower
semicontinuity with converging weights, as in Lemma~\ref{lem:LCd}(i), on compact sets and with an
exhaustion, and since $\partial^\alpha F_\nu$ and $\nabla\partial^\alpha F_\nu$ converge weakly in
$L^2$ on compact subsets of $(0,\sigma_0)\times\R^3$, the bound of Step~5 passes to the limit, and
we get $E^\natural_m(\tilde F(\sigma))+c_{W,\rm f}\int_0^\sigma D^\natural_m(\tilde
F)\le2E^\natural_m(F_*)$ on $[0,\sigma_0]$. Moreover, $\tilde F\in C([0,\sigma_0];H^{m-1}_{\rm
loc})$, and~(iii) holds for $\tilde F$ as in Step~4.

\emph{Step 7: identification.} We apply Lemma~\ref{lem:RSuniq} on $[0,\sigma_0]$ with $F_1:=F$ and
$F_2:=\tilde F$. Its hypotheses hold. Indeed, the solution $F\ge0$ is classical for $\sigma>0$ by
(GS1), and continuous up to $\sigma=0$, that is, up to $\sigma'$ before the shift. For $\sigma'=0$
this follows from Lemma~\ref{lem:GS1q}, and for $\sigma'>0$ because $F$ is classical there.
Moreover,
$F\le C'e^{-\beta'|v|^2}$ on $[0,1]$ by the hypothesis of~(iv). Finally, $\tilde F$ is classical,
$|\tilde F|\le Ce^{-(q/4)|v|^2}$ by~(iii), and $\|\jap v^2\nabla^2\tilde F\|_{L^2}\le C\|
\mathfrak m_2\nabla^2\tilde F\|_{L^2}\le CE^{1/2}$. Therefore, we have $F=\tilde F$ on
$[0,\sigma_0]$, which implies~(iv) with $c=c_{W,\rm f}$. The constants of Steps~3--5 depend on $q$
only through the ladder and through the tools lemma and Lemma~\ref{lem:Whypo}, continuously and with
$q$ bounded away from $0$. This yields the uniformity stated after~(iv), and completes the proof.
\end{proof}

\begin{remark}[Consequences, and a derivative count]\label{rem:RSshortF}
\begin{enumerate}[label=(\arabic*)]
\item The lemma provides Gaussian bounds, with the weight $e^{(q/2)|v|^2}$, at every
order $m\le N+2$, with the bound $2E^\natural_m(0)$. The Maxwellian bound (RM) of
Proposition~\ref{prop:RS} on bounded intervals follows from bounds for $|\alpha|\le2$ in
$L^2\cap L^\infty$ with some Gaussian weight, that is, from part~(iii) with $m=4$, as shown in
Appendix~\ref{app:hom}.
\item \emph{Derivative count.} Lemma~\ref{lem:RGinterp} turns $2N+4$ polynomially weighted
derivatives into $N+2$ Gaussian ones. Thus the lemma holds for $m\le N+2$. Lemma~\ref{lem:RSmed}
starts from it, so it holds for $m\le N+2$. Then property (RS) yields Gaussian bounds on $N$
derivatives, not $2N$, for data in $\mathfrak K^N$. Lemma~\ref{lem:RGGbounds} below needs Gaussian
bounds on $n+2$ derivatives of the fibered
trajectories, and by (RS) they hold for data in $\mathfrak K^{n+2}$. In the
assembly we have $n=5$, and we assume (KD$^7$).
\end{enumerate}
\end{remark}

\subsubsection{Local factors and the remainder}

We use the standing assumptions (H1)--(H3), and the conclusion
of Lemma~\ref{lem:RGbar} with $\rho(\xi)=\rho_\xi$, whose rate and constant we denote by $r$ and
$C_{\rm bar}$.

\begin{lem}[Local factors]\label{lem:RGloc}
Under (H1), (H3) and the conclusion of Lemma~\ref{lem:RGbar}, for $j:=|a|+|b|\le J-1$,
\begin{equation}
|\partial^a_y\partial^b_vg(t,y,v)|\le C_{\rm loc}\,\rho_v\,e^{-r_j|y|^2},\qquad r_j:=\Big(1-\frac
jJ\Big)\frac r8,
\end{equation}
with
\begin{equation}\label{eq:Cloc}
C_{\rm loc}:=2C_Je^{2L}\big(C_{\rm bar}e^{L^2/r+r}+K_{\rm cr}\big) .
\end{equation}

\end{lem}

\begin{proof}
Since $|\nabla\log\rho|\le L$ and $|\xi-v|=|y|/t\le|y|$, we have $\rho_\xi\le\rho_ve^{L|y|}$. With
$L|y|\le\frac r4|y|^2+L^2/r$, by the barrier we get $g\le C_{\rm
bar}e^{L^2/r}\rho_ve^{-(r/4)|y|^2}$.
On the ball $B_2(y,v)\subset\R^6$ we have $\rho_{v'}\le e^{2L}\rho_v$ and
$(|y'|)^2\ge(|y|-2)_+^2\ge\frac12|y|^2-4$. Hence we have
\begin{equation}
\|g\|_{L^\infty(B_2(y,v))}\le C'e^{2L}\rho_ve^{-(r/8)|y|^2},\quad C':=C_{\rm
bar}e^{L^2/r+r},\qquad\|g\|_{C^J(B_2(y,v))}\le e^{2L}\rho_vK_{\rm cr}.
\end{equation}
The Gagliardo--Nirenberg inequality on balls of $\R^6$, $\|\nabla^ju\|_{L^\infty(B_1)}\le
C_J\big(\|u\|^{1-j/J}_{L^\infty(B_2)}\|u\|^{j/J}_{C^J(B_2)}+\|u\|_{L^\infty(B_2)}\big)$, together
with $X^{1-j/J}Y^{j/J}\le X+Y$ applied after extracting the Gaussian factor, implies the claim.
\end{proof}

\begin{lem}[The remainder in fiber-normalized Gaussian norms]\label{lem:RGremF}
Let $n\ge0$, $J\ge n+3$, $r_\natural:=r_{n+2}$, and $0<q\le r_\natural/4$. Under the hypotheses of
Lemma~\ref{lem:RGloc}, for every fiber $v$, $t\in[T_0,T]$ and $|\alpha|\le n$,
\begin{equation}
\rho_v^{-1}\big\|e^{q|y|^2}\partial^\alpha_yR(t,\cdot,v)\big\|_{L^2_y}\le C_R\,C_{\rm
loc}^2\,\rho_v\,t^{-1},\qquad C_R=C(n,L,r_\natural).
\end{equation}
\end{lem}

\begin{proof}
\emph{Polarization.} For two profiles $g^1,g^2$ we put
$G^i_\vartheta(z)=g^i(t,z,v-\vartheta(z-y)/t)$
and $R(g^1,g^2):=Q(G^1_1,G^2_1)(y)-Q(G^1_0,G^2_0)(y)$. Writing $z=y+u$ and using translation
invariance, we have
\begin{equation}
Q(G^1_\vartheta,G^2_\vartheta)(y)=Q\big(g^1(t,y+\cdot,v-\vartheta\cdot/t),\,g^2(t,y+\cdot,v-\vartheta\cdot/t)\big)(0).
\end{equation}
Consequently $\partial_y$ acts through $g^1,g^2$, and
$\partial^\alpha_yR=\sum_{a\le\alpha}\binom\alpha aR(\partial^ag,\partial^{\alpha-a}g)$. By the
computation of Proposition~\ref{prop:Rformula}, polarized, we obtain
\begin{equation}
R(g^1,g^2)=-\frac1t\int_0^1\Big[\ab[H^1_\vartheta](y):(\tn_\vartheta\otimes\tn_\vartheta)g^2-\cb[H^1_\vartheta](y)\,g^2+2\,\ab[G^1_\vartheta](y):(\tn_\vartheta\otimes\partial_v)g^2\Big]\dd\vartheta,
\end{equation}
with $H^1_\vartheta(z)=(z-y)\cdot(\partial_vg^1)(t,z,v-\vartheta(z-y)/t)$, local factors at
$(t,y,v)$, and $\cb[H^1_\vartheta](y)=0$ for $\gamma=-3$.

\emph{Coefficients.} Here $g^1=\partial^ag$ with $|a|\le n$, so $g^1$ and $\partial_vg^1$ have order
$\le n+1\le J-1$. By Lemma~\ref{lem:RGloc} and $\rho_{v'}\le\rho_ve^{L|v'-v|}$ with
$|v'-v|\le|z-y|$, we have
\begin{equation}
|\ab[G^1_\vartheta](y)|+|\ab[H^1_\vartheta](y)|+|\cb[H^1_\vartheta](y)|\le3C_{\rm
loc}\rho_v\sup_{p\in[-2,1]}I_p(y),\quad I_p(y):=\int|y-z|^pe^{L|z-y|}e^{-r_\natural|z|^2}\dd z .
\end{equation}
For $\cb$ we used $|c(u)||u|\le2(\gamma+3)|u|^{\gamma+1}$, and the factor $2(\gamma+3)$ enters only
the constant. We split at $|y-z|=1$ and use $e^{L|z-y|}\le e^{L|y|}e^{L|z|}$. The near part is
$\le4\pi e^{L|y|+L^2/(4r_\natural)}$, since $\int_{|u|<1}|u|^p\dd u\le4\pi$ for $p\ge-2$. The far
part is $\le\jap ye^{L|y|}\int\jap ze^{L|z|-r_\natural|z|^2}\dd z$, since $|y-z|^p\le\jap y\jap z$
for $p\le1$. Hence we have $I_p(y)\le C_Ie^{(L+1)|y|}$.

\emph{Local factors.} With $g^2=\partial^{\alpha-a}g$, they have order $\le n+2\le J-1$. Since
$|\tn_\vartheta u|\le|\partial_yu|+|\partial_vu|$ componentwise for $t\ge1$, by the same lemma they
are $\le4C_{\rm loc}\rho_ve^{-r_\natural|y|^2}$.

\emph{Conclusion.} Combining these bounds, we have $\rho_v^{-1}|\partial^\alpha_yR|\le
C_n\,t^{-1}C_{\rm loc}^2\rho_ve^{(L+1)|y|-r_\natural|y|^2}$. Since $q\le r_\natural/4$,
$\|e^{q|y|^2}e^{(L+1)|y|-r_\natural|y|^2}\|_{L^2_y}\le\|e^{(L+1)|y|-\frac34r_\natural|y|^2}\|_{L^2}=C(L,r_\natural)$.
\end{proof}

\subsubsection{The fibered flow}

Let $G$ be the fibered flow, with $G(T_0)=g(T_0)$, and assume (H1) and (H2).
\begin{lem}[Bounds on the fibered flow]\label{lem:RGGbounds}
Assume (H1), (H2), and let the fibered trajectories satisfy the conclusions of
Proposition~\ref{prop:RS} on $[0,S^\sharp]$. Let $0<q\le q_S/(32\theta_+)$,
with $q_S$ of (RS). For every $v$ and $t$ with $s_v(t)\le S^\sharp$, the normalized fiber $\hat
G:=G(t,\cdot,v)/\rho_v$ satisfies:
\begin{enumerate}[label=(\roman*)]
\item $\hat G\ge0$, $\int\hat G=1$, $\int y\hat G=y_{0,v}$, $\int|y-y_{0,v}|^2\hat G=3\theta_v$;
\item
\begin{equation}\label{eq:KG}
\sum_{|\alpha|\le n+2}\big\|e^{2q|y|^2}\partial^\alpha_y\hat G\big\|_{L^2\cap L^\infty}\le K_G,
\end{equation}
and $\sup_y\jap y^8\hat G\le K_3$;
\item (G1)--(G3) of Lemma~\ref{lem:Whypo} with $M=1$, the constants $(L,\theta_+,K_3)$, and $c_G$
from Lemma~\ref{lem:RGellip} with $M_-=1$, $E_{\rm mom}=L^2+3\theta_+$, $K=K_G$.
\end{enumerate}
The constants $K_G,K_3$ depend on $(n,S^\sharp,C_{\rm dat},\beta,A_{n+2},L,\theta_\pm)$.
\end{lem}

\begin{proof}
By Lemma~\ref{lem:sym}(i) and uniqueness in (GS1), $G(t,\cdot,v)=T_{P_v}F_v(s_v(t))$, where $F_v$ is
the homogeneous trajectory from $F_{v,0}$. Indeed, we have
$\partial_t\big(T_{P_v}F_v(s_v)\big)=\kappa\lambda_vT_{P_v}Q(F_v,F_v)=\kappa Q(G,G)$. Since
$F_{v,0}$ has the moments of $\mu$ and they are conserved, (i) follows. By (RS) on $[0,S^\sharp]$,
for data in $\mathfrak K^{n+2}$, we have $|\nabla^jF_v(\sigma,z)|\le Ce^{-q_S|z|^2/4}$ for $j\le
n+2$. Then we get
\begin{equation}
\begin{aligned}
|\partial^\alpha_y\hat G(y)|&=\theta_v^{-3/2-|\alpha|/2}\big|(\nabla^\alpha
F_v)\big(s,(y-y_{0,v})/\theta_v^{1/2}\big)\big|\\
&\le C\theta_-^{-(n+5)/2}e^{-q_S|y-y_{0,v}|^2/(4\theta_+)}\le C'e^{-q_S|y|^2/(8\theta_+)},
\end{aligned}
\end{equation}
by $|y-y_0|^2\ge\frac12|y|^2-L^2$. Since $2q\le q_S/(16\theta_+)$, (ii) follows. Part (iii) follows
from (i), (ii) and the quantitative ellipticity estimate.
\end{proof}

\subsubsection{Rates, constants and inputs from the bootstrap}

\begin{remark}[The rates and the constants of the comparison]\label{rem:RGfull}
\begin{enumerate}[label=(\alph*)]
\item \emph{Form of the conclusion.} Theorem~\ref{thm:RGF} states closeness against
$\rho_\xi^{-1}\partial^j_yG$ sampled along the leaf, not against $\nabla^j_zG^{(x)}$, since the
latter would need $\partial_vG$. In Corollary~\ref{cor:Aleafentry} we have $j=0$, where the two
agree, so its first step holds with $\varepsilon_T(t)$ in place of $Ce^{4C_1S_\sharp}C_{\rm
loc}^2S_\sharp T_0^{-1}$. The bound (G4) follows from part~(ii) of the theorem with $j\le2$, so
$n\ge4$, and we take $n=5$.
\item \emph{Rates.} The comparison runs at $q\le q_\natural$, and the leaf bounds have rate $q/4$
(sup) and $q/8$ ($L^2$) in $z$. In the frame, we have $|\hat z|^2\le4(|z|^2+(L+1)^2)/\theta_-$,
since $\theta_g\ge\frac12\theta_-$ by (F3), so $e^{q_0|\hat z|^2}\le Ce^{(4q_0/\theta_-)|z|^2}$. To
leave room for the factors $\jap{\hat z}$ and $e^\omega$, we choose $q_0$ such that
\begin{equation}
q_0\le\theta_-\,q_\natural/64 .
\end{equation}
Further, we need $q_0\le q_{\rm GW}$, the exponent of Lemma~\ref{lem:GW} in the frames $p_g$ (with
$L+1$, $\frac12\theta_-$, $\theta_++1$), for Lemma~\ref{lem:CUbuf} and Corollary~\ref{cor:W1G}.
Moreover, we need $q_0\max\{1,\theta_+\}\le\beta'/4$ and $q_0\theta_+\le c_\beta/4$, with
$c_\beta:=\beta'\theta_-/(16\theta_+)$, for Corollary~\ref{cor:Aleafentry}, and
$q_0\theta_+\le\frac1{64}$ for Lemma~\ref{lem:UFref}. Since $q_\natural$ depends only on
$L,\theta_\pm,\beta_M$, on the data rate $\beta$, on $M_-/E_+$ and on $(n,J)$, this requirement
involves only structural constants.
\item \emph{Counts in the assembly.} We take $n=5$, $J=10$, the count of Lemma~\ref{lem:BTsob},
and (KD$^7$), as in (LWP)(a3) and Remark~\ref{rem:RSshortF}(2). Then $J\ge n+3$, and the theorem
yields leaf bounds of order $n-2=3\ge2$.
\end{enumerate}
\end{remark}

\begin{remark}[Bounds supplied by the bootstrap]\label{rem:RGfullinputs}
For the fibered comparison and its transfer to the leaves, we use the following bounds from the
bootstrap of Theorem~\ref{thm:CUboot} on $[T_0,T]$.
\begin{enumerate}[label=(\arabic*)]
\item \emph{(G1) and (G3) on every leaf}, normalized by $\rho_\xi$, as in
Lemma~\ref{lem:RGbar}(a). (G1) holds on every leaf by Lemma~\ref{lem:BTG12}, from (U1) alone, which
fixes $M_\pm$ and $E_+=2E_{\rm mom}$. (G3) is Lemma~\ref{lem:BTsob}(iv) on every leaf, with $K^{\rm
cr}_{\rm rel}$, and with $K_{\rm rel}$ on the relaxed ones.
\item \emph{The crude bound (H3)} on the whole region. It follows from (U1)--(U3), with the
constant $K_{\rm cr}$, by the leaf Sobolev and Cauchy estimates of part~(iii) of the same
lemma. Note that (H3) has no Gaussian weight. This is needed because the bootstrap provides Gaussian
bounds only at the rate $q_0$.
\item \emph{Qualitative Gaussian decay}, hypothesis~(c) of the barrier lemma. It follows from the
Gaussian leaf energies by part~(v) of the leaf Sobolev lemma.
\end{enumerate}
\end{remark}

\section{Uniformly local analytic energies}\label{sec:cells}

In this section we construct energies for which the analytic hierarchy closes despite the leafwise
transport.

\subsection{Energies on a single leaf}\label{ssec:oneleaf}

We cannot close the analytic hierarchy with the energy of a single leaf, as the next proposition
shows.

\begin{prop}[The leafwise transport on a single leaf]\label{prop:ACtr}
Let $q>0$, $X_q:=L^2(e^{2q|z|^2}\dd z)$, and $w_n=w_n(\rad)=\rad^n(n+1)^p/n!$ with $\rad>0$
and $p\ge0$.
\begin{enumerate}[label=(\roman*)]
\item There is no constant $C$ such that, for all $\varphi,\varphi'\in X_q$,
\begin{equation}\label{eq:ACtr-i}
|\langle\varphi,z_1\varphi'\rangle_{X_q}|\le C\|\varphi\|_{X_q}\|\varphi'\|_{X_q}.
\end{equation}
\item For every sequence $E_n\ge0$,
\begin{equation}\label{eq:ACtr-ii}
\sum_{n\ge0}w_n^2E_{n+1}\ge2^{-2p}\rad^{-2}\sum_{m\ge1}m^2w_m^2E_m .
\end{equation}
In particular, for every $c>0$ there is a sequence $(E_n)$, nonzero at one $n$ only, with
$\sum_nw_n^2E_{n+1}>c\sum_nn\,w_n^2E_n$.
\item \emph{(The collisionless model.)} Let $T_0>0$, and let $H(t,\xi,z)$ solve
$D_tH=t^{-2}z\cdot\nabla_\xi H$ for $t\ge T_0$, that is, \eqref{eq:leaftransport} without the
collision term. Let
\begin{equation}\label{eq:ACtr-Am}
\mathbb A^{\rm m}(t):=\sum_n\Big(\frac{\rad_0^n}{n!}\Big)^2\sum_{|\beta|=n}\|\partial^\beta_\xi
H(t,0,\cdot)\|^2_{X_q}
\end{equation}
be the analytic energy of the leaf $\xi=0$. For all $\rad_0>0$, $t>T_0$ and $M>0$ there is a
solution, entire in $\xi$, with $\mathbb A^{\rm m}(T_0)\le1$ and $\|H(t,0,\cdot)\|_{X_q}\ge M$.
\end{enumerate}
\end{prop}

\begin{proof}
(i) Let $\psi\ge0$ be a nonzero smooth function supported in the unit ball, let $R\ge2$, and put
$\varphi=\varphi'=e^{-q|z|^2}\psi(z-Re_1)$. Then we have $\|\varphi\|^2_{X_q}=\|\psi\|^2_{L^2}$ and
$\langle\varphi,z_1\varphi\rangle_{X_q}=\int z_1\psi(z-Re_1)^2\dd z\ge(R-1)\|\psi\|^2_{L^2}$.

(ii) For $m\ge1$, $w_{m-1}/w_m=\frac m\rad\big(\frac m{m+1}\big)^p\ge2^{-p}\frac m\rad$, and
$\sum_{n\ge0}w_n^2E_{n+1}=\sum_{m\ge1}w_{m-1}^2E_m$. For the second claim we take $E_n=\delta_{n,m}$
with $m>4^p\rad^2c$.

(iii) Let $\sigma(t)=T_0^{-1}-t^{-1}$, so that $\sigma'=t^{-2}$ and $\sigma(T_0)=0$. For data $H_0$,
the solution is
\begin{equation}\label{eq:ACtr-sol}
H(t,\xi,z)=H_0(\xi+\sigma(t)z,z).
\end{equation}
We fix $t>T_0$ and $M>0$, and put $K:=M/\rad_0$. We choose $\chi\ge0$ smooth with $\|\chi\|_{X_q}=1$
and support in $\{z_1\ge2\rad_0/\sigma(t)\}$, and put $H_0(\xi,z)=e^{-\rad_0K}e^{K\xi_1}\chi(z)$.
Then we have $\partial^\beta_\xi H(T_0,0,\cdot)=0$ unless $\beta=ne_1$, and
$\partial^n_{\xi_1}H(T_0,0,\cdot)=e^{-\rad_0K}K^n\chi$, which implies that
\begin{equation}
\mathbb A^{\rm m}(T_0)=e^{-2\rad_0K}\sum_n\frac{(\rad_0K)^{2n}}{(n!)^2}\le
e^{-2\rad_0K}\Big(\sum_n\frac{(\rad_0K)^n}{n!}\Big)^2=1 .
\end{equation}
On the other hand, $H(t,0,z)=e^{-\rad_0K}e^{K\sigma(t)z_1}\chi(z)\ge e^{\rad_0K}\chi(z)$ on the
support of $\chi$, and hence $\|H(t,0,\cdot)\|_{X_q}\ge e^{\rad_0K}=e^M\ge M$.
\end{proof}

\begin{remark}[The leafwise transport between leaves]\label{rem:ACtr}
The energy of a single leaf, $\hat{\mathbb A}_\xi$, evolves under the leafwise transport
$t^{-2}z\cdot\nabla_\xi$, which involves the neighboring leaves. On a single leaf this term pairs
the order-$n$ quantity with the order-$(n+1)$ quantity through the unbounded factor $z$. By
Proposition~\ref{prop:ACtr}(i), in the model norm $X_q$ this pairing is not bounded by the leaf
energies, so it does not have the form $C_0b\sqrt{E_nE_{n+1}}$ required by Lemma~\ref{lem:analytic}.
If the Gaussian weight controls $|\hat z|$ (the term (T8) of Subsection~\ref{sssec:bandproof}), what
remains is of the form $Ct^{-2}(E_n+E_{n+1})$ at the weight $w_n^2$, which by \eqref{eq:ACtr-ii}
loses two levels, while the radius term of that lemma supplies only
$2(|\rad'|/\rad)\sum_nn\,w_n^2E_n$. Part~(iii) of the proposition shows the reason. By
\eqref{eq:ACtr-sol}, the leaf energy at $(t,\xi)$ depends on the leaves $\xi+z\,\sigma(t)$, at
unbounded distance. Here $\sigma$ is the Euler time. Even without collisions, this energy is not
controlled by its value at $T_0$, whatever the radius at time $t$. The data in (iii) are not bounded
uniformly in $\xi$, so a bootstrap with bounds at all leaves excludes them. However, such a
bootstrap must then use the neighboring leaves. By Lemma~\ref{lem:CUtr}, in energies that are
integrated in $\xi$ the transport integrates by parts and loses no level.
\end{remark}

The proposition and the remark show that we cannot close the analytic hierarchy one leaf at a time.
On the other hand, in $L^2_\xi$, by Lemma~\ref{lem:localnoloss}, the transport integrates out. We
therefore use energies that are uniformly local in $\xi$.

We can do this for the following two reasons.
\begin{enumerate}[label=(\alph*)]
\item \emph{No loss of radius at entry.} This is Lemma~\ref{lem:CUnodrop}. Here
$\hat{\mathbb A}_{\ge1,\xi}$ and $\mathbb D_{\ge1,\xi}$ denote the parts of the leaf energy and of
the dissipation with $(\alpha,\beta)\ne(0,0)$, not only those with $|\beta|\ge1$. Therefore
$\hat{\mathbb A}_{\ge1}$ includes the $y$-derivatives, and every quadratic term has both factors in
this part, since $W^2_{\alpha\beta}=W_{\alpha_10}W_{\alpha_2\beta}W_{\alpha\beta}$ exactly for
$\alpha=\alpha_1+\alpha_2$. The crude regime already makes it smaller than any prescribed constant
chosen before the radius, and in particular smaller than $\delta_1^2$. Thus the analyticity radius
$\rad(t)$ and the weight $\delta_y$ can be the \emph{same for all leaves}. Only the damping
$\omega(s_\xi)$ and the nonlocal coefficient $\chi(s_\xi)$ depend on the leaf, and both have bounded
logarithmic gradient in $\xi$.
\item \emph{Cells at the analyticity scale.} We integrate the leaf densities against
$e^{-2\jap{(\xi-\xi_c)/\rad_0}}$, normalized by $\rad_0^{-3}$. By Lemma~\ref{lem:CUtr}, the local
transport then integrates by parts in $\xi$ and produces only a factor $|z|$, which is controlled by
the shrinking Gaussian weight. By Lemma~\ref{lem:CUsob}, the Sobolev inequality in $\xi$ on cubes of
side $\rad_0$ loses no power of $\rad_0$, only two levels with polynomial factors. By
Proposition~\ref{prop:CUtri}, the trilinear sums close with one factor taken on a cube and two
factors integrated.
\end{enumerate}
The result is Theorem~\ref{thm:CUboot}, stated and proved in Section~\ref{sec:boot}, which avoids
this obstruction. Subsection~\ref{ssec:cellsproofs} contains the local existence in cell norms,
Proposition~\ref{prop:CULW}. The constants are summarized in Remark~\ref{rem:CUstatus}.

\subsection{Regimes, damping, cells and cubes}\label{ssec:cells}

The objects of the uniformly local energy method are the leaf clock and the four regimes of a leaf,
the damping, the cell energies \eqref{eq:cellE}, the graded weights of Remark~\ref{rem:UFgrade} and
the cube energies \eqref{eq:cubeE}.

\paragraph{The leaf clock.} For every $\xi$ let
\begin{equation}\label{eq:leafclock-ref}
\lambda_{\rm ref}(\xi):=\lambda(P_\xi)=\rho_\xi\theta_\xi^{\gamma/2},\qquad s_\xi(t):=\lambda_{\rm
ref}(\xi)\big(\tau(t)-\tau(T_0)\big),
\end{equation}
where $\lambda_{\rm ref}(\xi)$ is the collision rate of the fiber parameters at $T_0$. Thus
$s_\xi=s_v|_{v=\xi}$ is the fiber clock, and Theorem~\ref{thm:RGF} and
Corollary~\ref{cor:Aleafentry} apply to it as stated. By (H1),
$|\nabla\log\lambda_{\rm ref}|\le\max\{1,\frac{|\gamma|}2\}L\le2L'$. We write
$\lambda_\xi:=\lambda(p_g(t,\xi))$ for the rate of the current leaf parameters. Under (U1), by (F3),
we have $\lambda_\xi/\lambda_{\rm ref}(\xi)\in[c_\lambda,c_\lambda^{-1}]$ with a constant
$c_\lambda$ which we take $\le\frac12\big(\theta_-/(2\theta_++2)\big)^{3/2}$. In the fibered
comparison, $\lambda_\xi$ denotes $\lambda_v$ at $v=\xi$, that is, $\lambda_{\rm ref}(\xi)$. Note
that $\lambda_{\rm ref}(\xi)\ne\lambda(p_{\rm ref}(\xi))$ in general: $p_{\rm ref}$,
$\varpi'_{\rm ref}$, $A_{\rm ref}$ refer to the leaf at $T_0$.

\paragraph{Entry levels.} We take
\begin{equation}\label{eq:Sstar}
S_*:=S_\flat(\delta'')+2,\qquad S_\sharp:=S_*+1,\qquad\delta'':=\min\big\{1,\big(\delta_1/(6C_eC_{\rm
pt})\big)^{5/2}\big\},
\end{equation}
with the constant $C_{\rm pt}$ of the proof of Corollary~\ref{cor:Aleafentry}, which does not depend
on $\delta_1$.
The corollary then applies to every leaf with $s_\xi\in[S_*-\frac12,3S_\sharp]$:
the fibers that such a leaf samples on $|z|\le t/(15LS_\sharp)$ have $s_v\ge
S_\flat(\delta'')+\frac12$, and the rest of the leaf contributes $\le CLS_\sharp/T_0$, as in its
proof.

Let $d_*>0$, fixed in requirement~\ref{rq:CUc} below, and put
\begin{equation}\label{eq:Sb}
S_b:=S_*e^{2L'd_*},\qquad L':=\tfrac54L .
\end{equation}
Here $|\nabla\log\lambda_{\rm ref}|\le2L'$ by (H1), as noted for the leaf clock above.

A leaf is in one of four regimes:
\begin{itemize}
\item \emph{crude} if $s_\xi<S_*$;
\item \emph{pre-band} if $S_*\le s_\xi<S_b$;
\item in the \emph{band} if $S_b\le s_\xi\le S_b+1$;
\item \emph{post-band} if $s_\xi>S_b+1$.
\end{itemize}
The nonlocal coefficient is $\chi(s)$, nondecreasing, with $0\le\chi'\le2$, $\chi=0$ for $s\le S_b$
and $\chi=1$ for $s\ge S_b+1$. It is constructed in Appendix~\ref{app:band}, and made $C^\infty$ in
Subsection~\ref{ssec:cellsproofs}.

\paragraph{Damping.} Let $H(0)=0$ and
\begin{equation}\label{eq:damping}
H'=C_1^U\ \text{on }[0,S_*),\qquad H'=C_1^B\ \text{on }[S_*,S_b+1],\qquad H'=0\ \text{on
}(S_b+1,\infty),
\end{equation}
and $\omega_\xi(t):=e^{-2H(s_\xi(t))}$. Hence $H$ is Lipschitz with constant $C_1^U$.

In the pre-band regime the background is $\delta_c$-close: by Corollary~\ref{cor:Aleafentry},
$\|f_\xi\|_X\le\delta_1/C_e$, and the requirement $\delta_1/C_e\le\delta_c/2$ implies the closeness.
Hence Proposition~\ref{prop:Bband} applies there with $\chi=0$. Its balance (deficit and front term)
needs only $H'\ge C_1^B$ and closeness.

\paragraph{Cells.} For $\xi_c\in\R^3$ put
\begin{equation}\label{eq:chic}
\chi_c(\xi):=e^{-\jap{(\xi-\xi_c)/\rad_0}},\qquad\text{so
that}\qquad|\nabla\chi_c|\le\rad_0^{-1}\chi_c,\qquad\rad_0^{-3}\int\chi_c^2=c_\chi,
\end{equation}
with $c_\chi$ an absolute constant. For each $(\alpha,\beta)$, we write
$h_{\alpha\beta}:=\hat r_{\alpha\beta}$, and let $\ell_\xi[h]$ be the local density in $z$-variables
as in Proposition~\ref{prop:FM}:
\begin{equation}\label{eq:elldens}
\ell_\xi[h]:=\Big(\frac{A_{p_g(t,\xi)}}{A_{\rm ref}(\xi)}\Big)^2\int m_x(z)^2|h^{(x)}(z)|^2\dd
z=\big\|T^{-1}_{p_g}h^{(x)}\big\|^2_{X(\hat m_\xi)},\qquad m_x(z)=W\big(z,\xi-z/t\big),
\end{equation}
with $W=A_{\rm ref}(v)e^{q(t)|y|^2}$. The leaf density and dissipation are
\begin{equation}\label{eq:leafed}
e_\xi[\varphi]=\eta\,\ell_\xi[h]+\chi(s_\xi)\nn\varphi^2_\natural,\qquad
d_\xi[\varphi]=\lambda_\xi\|\varphi\|^2_{Y_2(\hat m_\xi)},\qquad\varphi=T^{-1}_{p_g}h^{(x)}.
\end{equation}
Put
\begin{equation}\label{eq:En}
E_n(\xi):=\sum_{|\beta|=n}\sum_{|\alpha|\le N_y}\delta_y^{2|\alpha|}e_\xi[\varphi_{\alpha\beta}],
\end{equation}
and define $D_n(\xi)$ in the same way with $d_\xi$. The \emph{cell energy} and \emph{cell
dissipation} are
\begin{equation}\label{eq:cellE}
E^c_n:=\rad_0^{-3}\!\int\chi_c^2\,\omega_\xi\,E_n(\xi)\dd\xi,\qquad
\mathfrak A_c:=\sum_nw_n^2E^c_n,\qquad\mathfrak D_c:=\sum_nw_n^2\rad_0^{-3}\!\int\chi_c^2\,\omega_\xi\,D_n(\xi)\dd\xi
.
\end{equation}
The parts with $(\alpha,\beta)\ne0$ are $\mathfrak A_{c,\ge1}$ and $\mathfrak D_{c,\ge1}$.

\begin{lem}[Frame generator and the leaf equation]\label{lem:leafeq}
\leavevmode
\begin{enumerate}[label=(\roman*)]
\item Write a frame as $p=(p_\alpha,p_\beta,p_a)$ with $e^{p_\alpha}=\rho\theta^{-3/2}$,
$e^{p_\beta}=\theta^{1/2}$, $p_a=y_0$, so that $T_pF(z)=e^{p_\alpha}F(e^{-p_\beta}(z-p_a))$
(Lemma~\ref{lem:sym}). If $p(t)$ is a $C^1$ curve of frames, then
$\frac{d}{dt}(T_{p(t)}F)=T_{p(t)}\big(\mathcal G_{p(t)}[\dot p]F\big)$, where, for a vector $\dot
p=(\dot p_\alpha,\dot p_\beta,\dot p_a)$,
\begin{equation}\label{eq:framegen}
\mathcal G_p[\dot p]:=\dot p_\alpha-\dot p_\beta\,\hat z\cdot\nabla-e^{-p_\beta}\dot p_a\cdot\nabla
,
\end{equation}
with $\hat z\cdot\nabla$ and $\nabla$ acting on $F$ in its own (normalized) variable $\hat z$.
Consequently, if $\varphi=T^{-1}_{p(t)}H(t)$, then
$\frac{d}{dt}\varphi=T^{-1}_p\frac{d}{dt}H-\mathcal G_p[\dot p]\varphi$. Along rays, in leaf
parameters, $D_tp_\alpha=D_t\ell-\frac32D_t\vartheta$, $D_tp_\beta=\frac12D_t\vartheta$ and
$e^{-p_\beta}D_tp_a=e^{-\vartheta/2}D_ty_0$. The operator $\mathcal G_p[\dot p]$ maps moment-free
functions to moment-free functions.
\item Let $h=h_{\alpha\beta}=\partial^\alpha_y\partial^\beta_v\hat r$, $p=p_g(t,\xi)$ and
$\varphi=\varphi_{\alpha\beta}=T^{-1}_ph^{(x)}$, $x=t\xi$. Then, along rays (with $D_t$ of
Lemma~\ref{lem:leafgeom}(iv)),
\begin{equation}\label{eq:leafeq-phi}
D_t\varphi=T^{-1}_p(\partial_th)^{(x)}+t^{-2}\,T^{-1}_p\big[z\cdot(\partial_vh)^{(x)}\big]-
\mathcal G_p[D_tp]\varphi,
\end{equation}
where
\begin{equation}\label{eq:leafeq-h}
\partial_th=\partial^\alpha_y\partial^\beta_v\big[\kappa\big(\mathcal L_t[g]\hat r+\cQ_t(\hat
r,\hat r)\big)\big]-\partial^\alpha_y\partial^\beta_v\partial_t\mathcal M[g]
\end{equation}
by Lemma~\ref{lem:Meq}. Hence the transport term is the full, unprojected
$t^{-2}T^{-1}_p[z\cdot(\partial_vh)^{(x)}]$, and the forcing is the full
$-T^{-1}_p(\partial^\alpha_y\partial^\beta_v\partial_t\mathcal M[g])^{(x)}$. The functions $\varphi$
and $D_t\varphi$ are moment-free, and so are the collision terms and $\mathcal G_p[D_tp]\varphi$.
\end{enumerate}
\end{lem}

\begin{proof}
(i) We put $w=e^{-p_\beta}(z-p_a)$, the normalized variable. Then we have
$\frac{d}{dt}\big[e^{p_\alpha}F(w)\big]=e^{p_\alpha}\big[\dot p_\alpha F-\dot p_\beta\,w\cdot\nabla
F-e^{-p_\beta}\dot p_a\cdot\nabla F\big](w)$, which is the first claim. Differentiating
$T_pT^{-1}_pH=H$, we obtain the formula for $\frac{d}{dt}\varphi$, and the leaf parameters follow
from $p_\alpha=\ell-\frac32\vartheta$, $p_\beta=\frac12\vartheta$, $p_a=y_0$. For
$\psi\in\mathrm{span}\{1,\hat z_i,|\hat z|^2\}$ we have
$\int\psi\,\hat z\cdot\nabla\varphi=-\int\varphi\,(3\psi+\hat z\cdot\nabla\psi)$ and
$\int\psi\,\partial_i\varphi=-\int\varphi\,\partial_i\psi$, and both $3\psi+\hat z\cdot\nabla\psi$
and $\partial_i\psi$ lie in the same span.

(ii) By Lemma~\ref{lem:leafgeom}(ii) and~(iv),
$D_t[h^{(x)}]=(\partial_th)^{(x)}+t^{-2}z\cdot(\partial_vh)^{(x)}$. Applying (i) with $H=h^{(x)}$
along the ray, we obtain the identity. By Lemma~\ref{lem:Mfree} the leaf moments of $h$ vanish for
every $t$ and $\xi$, so $\varphi$ and $D_t\varphi$ are moment-free. The collision terms
$\kappa\,\cQ_t(\cdot,\cdot)$ have vanishing leaf moments for every $\xi$. Thus, by the same
argument, so do their derivatives $\partial^\alpha_y\partial^\beta_v$. Finally,
$\mathcal G_p[D_tp]\varphi$ is moment-free by~(i).
\end{proof}

\begin{remark}[A radius for each $y$-order]\label{rem:UFgrade}
By the leaf calculus, $\partial_y=\partial_z+t^{-1}\partial_\xi$ on a leaf. Therefore
$\partial_y^\alpha$ applied to a coefficient built from $p_g(t,\xi)$ produces up to $|\alpha|$
additional $\xi$-derivatives of $p_g$, each with a factor $t^{-1}$. This happens in two places. In
the coefficient terms (C) of Proposition~\ref{prop:ACpf}, and in the couplings of
Proposition~\ref{prop:TLsrc} with $\beta'=0$, $\alpha'\ne0$, the coefficient
$\mathfrak g_{\alpha'\beta'}$ contains $t^{-k}\partial^{|\beta'|+k}_\xi\delta\varpi$ with
$k\le|\alpha'|$, and these proofs treat it as a coefficient of level $|\beta'|$. In the forcing
(T12), the leaf restriction contains $t^{-k}\partial^{|\beta|+k}_\xi\Theta$, by \eqref{eq:UFleaf}
below, and its $\delta\varpi$-part and its micro moments have $\xi$-level up to $|\beta|+k+1$. The
$\varpi_*$-parts are controlled by the Cauchy estimates at a radius $r_{\rm C}\gg\rad_0$. For
$\delta\varpi$ and the micro moments, the weight ratio
$w_{|\beta|}/w_{|\beta|+k+1}\simeq((|\beta|+k+1)/\rad)^{k+1}$ is unbounded in $|\beta|$. For $k\ge1$
the loss is not of the form required by Lemma~\ref{lem:analytic}.

\emph{The graded weights.} Let $\epsilon_y:=1/(4N_y+4)$ and
\begin{equation}\label{eq:graded}
\rad_k(t):=\rad(t)-k\,\epsilon_y\rad_0,\qquad
W_{\alpha\beta}:=w_{|\beta|}\big(\rad_{|\alpha|}(t)\big)\,\delta_y^{|\alpha|} .
\end{equation}
\begin{enumerate}[label=(\alph*)]
\item $\rad_k\ge\rad_0/4$ for $k\le N_y$, and all $\rad_k$ decrease at the rate
$\rad'$. Thus the absorption lemma applies to each $y$-order separately, with the same $C_0$.
\item \emph{Monotonicity.} If $\alpha=\alpha'+\alpha''$, then
$\rad_{|\alpha|}\le\min\{\rad_{|\alpha'|},\rad_{|\alpha''|}\}$, so
$\rad_{|\alpha|}^n\le\rad_{|\alpha'|}^k\rad_{|\alpha''|}^{n-k}$. Hence every splitting of
the weights used in Lemma~\ref{lem:algebra}, Propositions~\ref{prop:ACpf}(N), \ref{prop:TLsrc}
and~\ref{prop:CUtri} holds without change. The same is true for the couplings from lower to higher
$y$-order.
\item \emph{Gap.} For $|\alpha|\ge1$, $j\ge0$ and $m\ge1$,
\begin{equation}\label{eq:UFgap}
\frac{w_j(\rad_{|\alpha|})}{w_{j+m}(\rad)}\le\frac{(j+m)!}{j!}\frac{\rad_{|\alpha|}^j}{\rad^{j+m}}\le(j+m)^m(1-\epsilon_y)^j\rad^{-m}\le
C(m,\epsilon_y)\rad_0^{-m},
\end{equation}
because $\rad_{|\alpha|}/\rad\le1-\epsilon_y\rad_0/\rad\le1-\epsilon_y$.
\item Lemma~\ref{lem:CUsob} and Lemma~\ref{lem:CUcomp} hold with $\rad_0/4$ in place of $\rad_0/2$,
with additional absolute factors $4^{|\gamma|}$. Lemma~\ref{lem:CUnodrop} is unchanged: for each
$y$-order the radius is the same on all leaves.
\end{enumerate}
With \eqref{eq:UFgap}, a coefficient $t^{-k}\partial^{|\beta'|+k}\delta\varpi$ with $k\ge1$ (so
$|\alpha|\ge|\alpha'|\ge1$) is bounded at weight $w_{|\beta'|}(\rad_{|\alpha|})$ by
$C(\rad_0t)^{-k}$ times its weight $w_{|\beta'|+k}(\rad)$. Since $\rad_0T_0\ge1$ by
requirement~\ref{rq:UFd}, the bounds of the coefficient terms (C) and of the leafwise sources hold
with $\delta_U$ replaced by $C\delta_U$. We treat the forcing in Proposition~\ref{prop:UFforce}.
From now on the cell and cube energies are defined with the graded weights. For $k\le N_y$ and
$n\ge0$ we put
\begin{equation}
E^c_{k,n}:=\rad_0^{-3}\!\int\chi_c^2\,\omega_\xi\sum_{|\alpha|=k,\,|\beta|=n}\delta_y^{2|\alpha|}e_\xi[\varphi_{\alpha\beta}]\dd\xi,
\end{equation}
\begin{equation}
\mathfrak A_c:=\sum_{k\le N_y}\sum_nw_n(\rad_k)^2E^c_{k,n}=\sum_{(\alpha,\beta)}W^2_{\alpha\beta}\,\rad_0^{-3}\!\int\chi_c^2\omega_\xi
e_\xi[\varphi_{\alpha\beta}]\dd\xi,
\end{equation}
and similarly for $\mathfrak D_c$, $\mathfrak a_Q$, $\mathfrak d_Q$ and the parts $\ge1$. Every
formula below written with $\sum_nw_n^2E^c_n$ (the radius term, the one-level losses) is understood
$y$-order by $y$-order, with $w_n(\rad_k)$ and $E^c_{k,n}$ in place of $w_n(\rad)$ and $E^c_n$.
By~(a), all $\rad_k$ decrease at the same rate, and thus the absorption lemma applies to each
$y$-order.

\begin{lem}[No coupling from a higher to a lower $y$-order]\label{lem:nodown}
In the hierarchy for $\varphi_{\alpha\beta}$, every term involves the components
$\varphi_{\alpha'\beta'}$ only with $\alpha'\le\alpha$ componentwise, coefficients whose $y$-order
is $\le|\alpha|$, and derivatives of $\varphi_{\alpha\beta}$ in its own variable $\hat z$.
\end{lem}

\begin{proof}
We go through the terms of the equation for $\partial^\alpha_y\partial^\beta_v\hat r$, restricted to
a leaf.
\begin{itemize}
\item \emph{Collision terms.} By Leibniz's rule, $\partial^\alpha_y\partial^\beta_vQ(g,g)$ is a sum
of $Q(\partial^{\alpha'}_y\partial^{\beta'}_vg,\partial^{\alpha''}_y\partial^{\beta''}_vg)$ with
$\alpha'+\alpha''=\alpha$. The Maxwellian factors produce coefficients of $y$-order $\le|\alpha|$.
\item \emph{Transport.} The local transport $t^{-2}z\cdot\nabla_\xi$ acts on $\varphi_{\alpha\beta}$
and raises only $\beta$ (since $\partial_\xi=\partial_v$). The unprojected transport
$t^{-2}T_p^{-1}[z\cdot(\partial_vh)^{(x)}]$ does the same.
\item \emph{Frame motion.} The terms $\mathcal G_p[Dp]\varphi_{\alpha\beta}$ are first-order
operators in $\hat z$ acting on $\varphi_{\alpha\beta}$ itself, by Lemma~\ref{lem:sym}(iii), with
coefficients built from $Dp$, and they do not change $(\alpha,\beta)$.
\item \emph{Coefficient terms.} On a leaf, $\partial_y=\partial_z+t^{-1}\partial_\xi$. Thus
$\partial^\alpha_y$ of a coefficient built from $p_g(t,\xi)$ produces
$t^{-k}\partial^{|\beta'|+k}_\xi$ with $k\le|\alpha'|\le|\alpha|$: a shift into $\xi$-levels, that
is, into $v$-levels, not into $y$-levels.
\item \emph{Forcing.} By the identity \eqref{eq:UFleaf}, the forcing for $(\alpha,\beta)$ involves
$\partial^{\alpha-\alpha'}_z\partial^{\beta+\alpha'}_\xi\Theta$ with $\alpha'\le\alpha$, again a
shift into $\xi$-levels.
\item \emph{Other terms.} The kernel term $\Pi_0[(y_0+\theta^{1/2}\hat
z)\cdot\varphi_{\alpha,\beta+e}]$, the nonlocal transport and coercivity, and the weight terms are
diagonal in $\alpha$ or raise only $\beta$.
\end{itemize}
In none of these terms does a component or a coefficient of $y$-order $>|\alpha|$ appear.
\end{proof}

By~(b) and this lemma, every splitting of the weights in the trilinear and coefficient
terms satisfies $W_{\alpha\beta}\le W_{\alpha'\beta'}W_{\alpha''\beta''}\cdot(\text{binomial
factors})$, as in the ungraded case.
\end{remark}

\paragraph{Cubes.} Let $\xi_j\in\rad_0\mathbb Z^3$ and
\begin{equation}\label{eq:cubes}
Q_j:=\xi_j+\big[-\tfrac{\rad_0}2,\tfrac{\rad_0}2\big)^3,\qquad
Q_j^*:=\xi_j+\big[-\tfrac{3\rad_0}2,\tfrac{3\rad_0}2\big)^3 .
\end{equation}
The \emph{undamped cube energy} and the \emph{damped cube dissipation} are
\begin{equation}\label{eq:cubeE}
\mathfrak a_Q:=\sum_{(\alpha,\beta)}W_{\alpha\beta}^2\rad_0^{-3}\!\int_Qe_\xi[\varphi_{\alpha\beta}]\dd\xi,\qquad
\mathfrak d_Q:=\sum_{(\alpha,\beta)}W_{\alpha\beta}^2\rad_0^{-3}\!\int_Q\omega_\xi
d_\xi[\varphi_{\alpha\beta}]\dd\xi ,
\end{equation}
with parts $\mathfrak a_{Q,\ge1}$ and $\mathfrak d_{Q,\ge1}$. The macroscopic cube energy is
\begin{equation}\label{eq:cubem}\mathfrak m_Q:=\sum_nw_n^2\rad_0^{-3}\!\int_Q\sum_{|\beta|=n}|\partial^\beta_\xi\delta\varpi|^2\dd\xi
.\end{equation}

\paragraph{Requirements on $\rad_0$ and $d_*$.} These are imposed when $\rad_0$ is chosen, unless
noted, and requirement~\ref{rq:CUd} is a lower bound on $T_0$. Let $C_\varpi:=C(L_1,R_1)+1$ bound
$|\varpi'_g|+|\nabla(y_{0,g},\vartheta_g)|$. It enters Lemma~\ref{lem:CUtr} and the frame ratio in
Lemma~\ref{lem:CUcomp} through $|\nabla_\xi\log A_{p_g}|$, which involves only $\nabla\ell$,
$\nabla y_0$ and $\nabla\vartheta$. By Lemma~\ref{lem:Mref}, (U1) below and Lemma~\ref{lem:CUsob},
this holds as soon as $C\delta_U\le\rad_0$, which is a condition on $T_0$ and $\eps$, imposed when
these are chosen.
\begin{itemize}
\item[] \refstepcounter{rq}\label{rq:CUa}\textbf{\therq.} $\epsilon_0:=12\sqrt3\,L'\rad_0$
satisfies $8C_1^US_\sharp\epsilon_0\le1$.
\item[] \refstepcounter{rq}\label{rq:CUb}\textbf{\therq.}
$\rad_0\big(4C_1^UL'(S_\sharp+2)+(4+C_d)C_\varpi+8L\big)\le1$, with $C_d=C_d(\gamma,\theta_-)$
the frame-derivative constant of Lemma~\ref{lem:CUcomp}(iii).
\item[] \refstepcounter{rq}\label{rq:CUc}\textbf{\therq.}
\begin{equation}\label{eq:dstar}
d_*:=\rad_0\big(2C_1^US_\sharp+\log(C_{\rm far}C_U/\delta_1^2)\big)
\end{equation}
satisfies $8C_1^US_*L'd_*\le1$ and $S_b\le S_\sharp$. Here $C_{\rm far}\ge1$ is the absolute
constant of Step~4 of Theorem~\ref{thm:CUboot}. Thus $d_*>0$, since $C_U\ge1\ge\delta_1^2$.
\item[] \refstepcounter{rq}\label{rq:CUd}\textbf{\therq.} $T_0\ge C_{\rm tr}/(\rad_0^2C_2)$, with
$C_{\rm tr}$ from Lemma~\ref{lem:CUtr}. The Gronwall rate $C_{\rm tr}/(\rad_0^2C_2)$ of the
transport remainder is part of $a_1$, whose integral is bounded in requirement~\ref{rq:UFd}.
\end{itemize}
Each of these is an upper bound on $\rad_0$ in terms of constants chosen before $\rad_0$ (among
them $\delta_1$), or a lower bound on $T_0$.

\subsection{The transport term in a cell}\label{ssec:trint}

For the local part of the cell energy, the transport loses no level.

\emph{Idea of the proof.} For a local energy, we can integrate the transport
$t^{-2}z\cdot\nabla_\xi$ by parts in $\xi$ at fixed $z$. It then falls on the cutoff $\chi_c^2$,
which produces a factor $\rad_0^{-1}$, and on $\omega_\xi$ and the frame, which produce bounded
factors. The remaining factor $|z|$ is controlled by the decrease of the Gaussian rate $q(t)$. We
treat the nonlocal part of the energy leaf by leaf.

\begin{lem}[Local transport in a cell]\label{lem:CUtr}
Let
\begin{equation}\label{eq:Lc}
L^c:=\sum_nw_n^2\rad_0^{-3}\!\int\chi_c^2\omega_\xi\sum_{|\beta|=n,\alpha}\delta_y^{2|\alpha|}\eta\,\ell_\xi[h_{\alpha\beta}]\dd\xi
\end{equation}
be the local part of $\mathfrak A_c$. In $\frac{d}{dt}L^c$, consider two contributions:
\begin{itemize}
\item the local transport, that is, the term $t^{-2}z\cdot\nabla_\xi[h^{(x)}]$ of $D[h^{(x)}]$;
\item the Gaussian weight (T8).
\end{itemize}
Their sum is bounded by
\begin{equation}\label{eq:CUtr}
-C_2t^{-2}\sum_nw_n^2\rad_0^{-3}\!\int\chi_c^2\omega_\xi\sum\delta_y^{2|\alpha|}\eta\Big(\frac{A_p}{A_{\rm
ref}}\Big)^2\!\!\int m_x^2|z|^2|h^{(x)}_{\alpha\beta}|^2\dd z\dd\xi+\frac{C_{\rm
tr}}{\rad_0^2C_2}t^{-2}L^c ,
\end{equation}
with $C_{\rm tr}$ absolute. There is no loss of level.
\end{lem}

\begin{proof}
\emph{Step 1: integration by parts in $\xi$.} By Lemma~\ref{lem:leafgeom}(iv), for any profile
function $h$,
\begin{equation}\label{eq:CUtr-D}
D\big[h^{(x)}\big](z)=(\partial_th)^{(x)}(z)+t^{-2}z\cdot\nabla_\xi\big[h^{(x)}\big](z),
\end{equation}
and $\nabla_\xi[h^{(x)}_{\alpha\beta}]=(\partial_vh_{\alpha\beta})^{(x)}$ is the next level. For
\emph{fixed} $(\alpha,\beta)$, the transport differentiates the \emph{same} function in $\xi$. We
put $\Omega:=\chi_c^2\omega_\xi(A_p/A_{\rm ref})^2m_x^2$, a function of $(\xi,z)$. The transport
contributes
\begin{equation}\label{eq:CUtr-ibp}
2t^{-2}\eta\,\rad_0^{-3}\!\iint\Omega\,h\,z\cdot\nabla_\xi h\dd
z\dd\xi=-t^{-2}\eta\,\rad_0^{-3}\!\iint(\nabla_\xi\Omega)\cdot z\,|h|^2\dd z\dd\xi .
\end{equation}
There is no boundary term in $\xi$: $\chi_c$ decays exponentially, and the cube energies are
bounded, by the bootstrap. For the components of $y$-order $N_y$, the integrability in $z$ that this
integration by parts uses is not known at a fixed time. In the profile variables we need no
integration by parts: the two contributions together are the term with $\mathfrak w_c$ of
Lemma~\ref{lem:CUcompid}, a part of the time derivative of the weight, and the bound below is the
pointwise bound on $\mathfrak w_c$ stated there.

\emph{Step 2: the logarithmic gradient of $\Omega$.} We have $|\nabla_\xi\log\Omega|\le\Gamma$, with
\begin{equation}\label{eq:CUtr-Gamma}
\Gamma:=\frac2{\rad_0}+4C_1^UL'(S_\sharp+2)+2C_\varpi+8L\le\frac3{\rad_0}
\end{equation}
by requirement~\ref{rq:CUb}. Indeed, we have $|\nabla\log\chi_c^2|\le2/\rad_0$ and
$|\nabla_\xi\log\omega_\xi|=2H'|\nabla s_\xi|\le4C_1^UL's_\xi$, and $H'=0$ for $s_\xi>S_b+1$.
Moreover, the frame ratio $(A_{p_g}/A_{\rm ref}(\xi))^2$ contributes $\le2C_\varpi+4L$, and $m_x^2$
contributes $\le4L$, by Lemma~\ref{lem:CUcomp}(iii). \emph{Step 3: absorption by the Gaussian
weight.} By (T8), $D(m_x^2)/m_x^2\le t^{-2}(-2C_2|z|^2+4L|z|)$. Using
\begin{equation}\label{eq:CUtr-young}
(\Gamma+4L)|z|\le C_2|z|^2+\frac{(\Gamma+4L)^2}{4C_2},\qquad(\Gamma+4L)^2\le\frac{16}{\rad_0^2},
\end{equation}
in \eqref{eq:CUtr-ibp}, we obtain \eqref{eq:CUtr}.
\end{proof}

We treat the nonlocal part of the transport on each leaf, as in (T10). There it is bounded by
$Ct^{-2}\|\varphi\|_{X_7}\|\varphi'\|_{X_8}\le C't^{-2}e_\xi[\varphi]^{1/2}e_\xi[\varphi']^{1/2}$.
After integration with $\chi_c^2\omega_\xi$ and Cauchy--Schwarz, this becomes
$C_0t^{-2}\sqrt{E^c_nE^c_{n+1}}$, which is the form required by Lemma~\ref{lem:analytic}. The kernel
part of the transport term of Lemma~\ref{lem:leafeq}(ii) is treated in the same way. In
$z$-variables the local part uses the full forcing and the unprojected transport, so no kernel term
arises there.

\subsection{Sobolev inequality and trilinear sums on cubes}\label{ssec:cubesob}

We state three tools for the uniformly local energies, which are proved in
Subsection~\ref{ssec:cellsproofs}. By Lemma~\ref{lem:CUsob}, on cubes of side $\rad_0$, Sobolev's
inequality in $\xi$ loses no radius: derivatives in $\xi$ are measured in the analytic weights, and
a derivative on a cube of side $\rad_0$ produces the factor $\rad_0^{-1}$ that the weight
$w_{n+1}/w_n\simeq\rad/(n+1)$ provides. Proposition~\ref{prop:CUtri} then estimates the trilinear
sources of the hierarchy in uniformly local form, with one factor in $\sup_c$ and two in the cell,
as in the theory of uniformly local Sobolev spaces \cite{Kat75}. Finally, Lemma~\ref{lem:CUcomp}
shows that the leaf energies on neighboring cubes are comparable.

\subsection{The cell inequality}\label{ssec:cellineq}

The cell inequality, Proposition~\ref{prop:CUcell}, is the differential inequality for the cell
energy $\mathfrak A_c$ on which the bootstrap rests. We obtain it by integrating the leaf
inequalities of the different regimes over a cell.

We first fix the coercivity constants. Let
\begin{equation}\label{eq:cmin}
c_{\rm cl}:=\frac{c_W}{4C_W},\quad c_{\rm buf}:=\frac{\eta c_W}2=\frac{c_W}{8C_W}=\frac{c_{\rm
cl}}2,\quad c_{\min}:=\min\{\eta c_h(K),\ c_{\rm buf}\},\quad c'_{\min}:=c_{\rm cl},\quad
c''_{\min}:=c_{\rm buf}.
\end{equation}
Hence we have $c_{\min}\le c''_{\min}\le c'_{\min}$. These are the coercivity constants of the crude
leaves in Proposition~\ref{prop:ACpf}, of the buffer leaves in Lemma~\ref{lem:CUbuf}, and, by
Remark~\ref{rem:Bbandrange}, of the pre-band, band and post-band leaves in
Proposition~\ref{prop:Bband}. Each of them is at least $c_{\min}$. Therefore every leaf has
coercivity at least $c_{\min}$, including the buffer leaves below, whose coercivity is
$c_{\rm buf}\ge c_{\min}$, and every leaf with $s_\xi\ge S_*$ has coercivity at least $c'_{\min}$.
The crude leaves with $s_\xi\in[S_*-\frac12,S_*)$, which we call \emph{buffer leaves}, have
coercivity at least $c_{\rm buf}$ by the next lemma. The constants $c'_{\min}$, $c_{\rm buf}$ and
$c''_{\min}$ are chosen before $\delta_1$ and do not depend on $K$. The coercivity constants
$c(\xi)$ used below are those that remain in these two propositions after their own absorptions: the
linear terms there provide $2\eta c_h(K)$ and $2\eta c_W$, and the absorptions made there use at
most half of these, so the further absorptions below stay within the available coercivity.

\begin{lem}[Coercivity on the buffer leaves]\label{lem:CUbuf}
Let $\xi$ be a leaf with normalized background $G=\mu+f_\xi$ in its frame, such that $f_\xi$
satisfies (G4) with constant $K$ (Remark~\ref{rem:BTnorm}) and $\|f_\xi\|_{X(\hat
m_\xi)}\le\delta_c$. Then for every $\varphi$,
\begin{equation}\label{eq:CUbuf}
\int\hat m_\xi^2\,\varphi\,\big(Q(G,\varphi)+Q(\varphi,G)\big)\le-\frac{c_W}2\|\varphi\|^2_{Y_2(\hat
m_\xi)}+C_{\rm buf}(1+K^2)\|\varphi\|^2_{X(\hat m_\xi)},
\end{equation}
with $C_{\rm buf}=C(q_\pm,C_W,c_W)$. Consequently, in the crude leaf inequality of
Proposition~\ref{prop:ACpf} for such a leaf, the linear term of every component
$\varphi_{\alpha\beta}$ of the hierarchy can be estimated with the coercivity $c_{\rm buf}=\eta
c_W/2$ in place of $\eta c_h(K)$, the remainder $\eta C_{\rm
buf}(1+K^2)\kappa\lambda_\xi\|\varphi_{\alpha\beta}\|^2_X$ being absorbed by the damping on crude
leaves, provided
\begin{equation}
c_\lambda C_1^U\ge C_{\rm buf}(1+K^2).
\end{equation}
At $(\alpha,\beta)=(0,0)$ the collision term is $L_Gf_\xi-Q(f_\xi,f_\xi)=Q(G,f_\xi)+Q(f_\xi,\mu)$.
Moreover, the quadratic part satisfies
\begin{equation}\label{eq:CUbuf-cubic}
|\langle Q(f_\xi,f_\xi),f_\xi\rangle_{X(\hat m_\xi)}|\le CK\|f_\xi\|^2_{X(\hat m_\xi)}
\end{equation}
by (L2) of the proof of that proposition, and it is controlled in the same way, with $C_{\rm buf}$
enlarged by a structural constant.
\end{lem}

\begin{proof}
\emph{Step 1: the inequality \eqref{eq:CUbuf}.} We apply Corollary~\ref{cor:W1G}, (W1$_G$), with
$f=f_\xi$ (so $\|f\|_X\le\delta_c\le1$) and $m=\hat m_\xi$. The conclusion of Lemma~\ref{lem:GW}(i)
holds for the leaf weight in the frame $p_g$. This weight is off-centered, and by
Lemma~\ref{lem:TG0}(i) it has $|\nabla\omega|\le\eps_1(t)=2(\theta_++1)^{1/2}(L+1)(q+t^{-1})$, which
is admissible for $q\le q_0$ and $t\ge T_1$ by the hypodissipativity lemma and the discussion of the
Gaussian hypodissipativity in Appendix~\ref{app:band}. Then we have
\begin{equation}\label{eq:CUbuf-1}
\int\hat m^2\varphi L_{\mu+f}\varphi\le-c_W\|\varphi\|^2_{Y_2}+C_W\|\varphi\|^2_{L^2}+C_W\|f\|_X\|\varphi\|^2_{Y_2}+C_W\|f\|_{Y_2}\|\varphi\|_X\|\varphi\|_{Y_2}
.
\end{equation}
By requirement~\ref{rq:Wdc}, $C_W\delta_c\le c_W/4$, and thus the third term is bounded by
$\frac{c_W}4\|\varphi\|^2_{Y_2}$. For the fourth, $\|f\|_{Y_2}$ is bounded by (G4): in
$\|f\|^2_{Y_2}=\|\jap{\hat z}^{(\gamma+2)/2}\hat mf\|^2+\|\jap{\hat z}^{\gamma/2}\tn(\hat mf)\|^2$
we have $|\tn(\hat mf)|\le\jap{\hat z}\big(\hat m|\nabla f|+|\nabla\hat m||f|\big)$ and
$|\nabla\hat m|\le C(q_\pm)\jap{\hat z}\hat m$, while $\jap{\hat z}^{\gamma/2+1}\le1$ since
$\gamma\le-2$. Hence $\|f\|_{Y_2}\le C(q_\pm)\sum_{j\le1}\|\jap{\hat z}\hat m\nabla^jf\|_{L^2}\le
C(q_\pm)K$. By Young's inequality, we have
$C_W\|f\|_{Y_2}\|\varphi\|_X\|\varphi\|_{Y_2}\le\frac{c_W}4\|\varphi\|^2_{Y_2}+C_W^2C(q_\pm)^2c_W^{-1}K^2\|\varphi\|^2_X$.
Finally, we have $\|\varphi\|_{L^2}\le C_m\|\varphi\|_X$ by Lemma~\ref{lem:TG0}(ii). This proves the
inequality, with $C_{\rm buf}$ depending also on $(\eps_1,C_m)$.

\emph{Step 2: the consequence.} Recall that in the convention of Proposition~\ref{prop:ACpf} the
time derivative of $\eta\|\varphi\|_X^2$ produces
$2\eta\kappa\lambda_\xi\langle\varphi,L_G\varphi\rangle$, so the coercive part is
$2\eta\cdot\frac{c_W}2=2c_{\rm buf}$ and the remainder is
$2\eta C_{\rm buf}(1+K^2)\kappa\lambda_\xi\|\varphi\|^2_X$. The damping on crude leaves is
$-2H'(s_\xi)\dot s_\xi\,\omega_\xi e_\xi$ with $H'=C_1^U$ and
$\dot s_\xi=\kappa\lambda_{\rm ref}(\xi)\ge c_\lambda\kappa\lambda_\xi$, and
$e_\xi\ge\eta\|\varphi\|^2_X$ for each component, since $\chi=0$ on crude leaves. It supplies
$2C_1^Uc_\lambda\kappa\lambda_\xi\eta\|\varphi\|_X^2$, so the remainder is absorbed when
$c_\lambda C_1^U\ge C_{\rm buf}(1+K^2)$. The other terms of that proposition are unchanged, and
their smallness conditions are taken with $\min\{\eta c_h(K),c_{\rm buf}\}$ in place of
$\eta c_h(K)$.
\end{proof}

\begin{prop}[Cell inequality]\label{prop:CUcell}
Suppose that on $[T_0,T]$ the following hold.
\begin{itemize}
\item[(H$_a$)] For every cube with $Q^*\subset\{s<S_b+1\}$, $\mathfrak a_{Q^*,\ge1}\le c_\sharp$;
for every other cube, $\mathfrak a_{Q^*,\ge1}\le\delta_b^2$. Moreover $C_Nc_\sharp^{1/2}\le
c_{\min}/8$ and $C_N\delta_b\le c'_{\min}/8$.
\item[(H$_b$)] For every cube $Q$ that contains a leaf with $s_\xi\ge S_*$: $
\mathfrak a_{Q^*}\le\delta_b^2$ and $\sup_{Q^*}\|f_\xi\|_X\le\delta_c$, with $C_N\delta_b\le c''_{\min}/4$ and
$\delta_b^2\le\delta_s$.
\item[(H$_c$)] Every leaf with $s_\xi\le3S_\sharp$ satisfies (G1)--(G4), with the constants of
Lemma~\ref{lem:BTG12} in (G1)--(G2) and the constant $K$ in (G3)--(G4).
\item[(H$_d$)] $\sup_Q\mathfrak m_{Q^*}\le(2\delta_U)^2$.
\item[(H$_e$)] $\sup_Q\mathfrak a_{Q^*}\le B_{\rm cr}$, so that $\mathfrak p_\xi\le\mathfrak p_*$ on
every leaf (Corollary~\ref{cor:FMuse}(d), and under (U2) this is the proof of
Lemma~\ref{lem:BTsob}(i), while in Subsection~\ref{ssec:cellsproofs} it is
$B_{\rm loc}\le B_{\rm cr}$).
\end{itemize}
Assume also the requirements of Proposition~\ref{prop:ACpf} and of Appendix~\ref{app:const} for
Propositions~\ref{prop:Bband} and~\ref{prop:TLsrc}. Then, for every cell $c$,
\begin{multline}\label{eq:CUcell}
\frac{d}{dt}\mathfrak A_c+\frac{c_{\min}}4\kappa\,
\mathfrak D_c\le\sum_n\Big[2n\frac{\rad'}{\rad}w_n^2E^c_n+C_0t^{-2}w_n^2\sqrt{E^c_nE^c_{n+1}}\Big]\\
+t^{-2}a_1\,\mathfrak A_c+\mathfrak f^c+\mathfrak q_c ,
\end{multline}
where
\begin{equation}\label{eq:a1rc}
a_1:=C+4\mathfrak p_*+\frac{C_{\rm tr}}{\rad_0^2C_2},\qquad
\mathfrak q_c:=C(K)\,\rad_0^{-3}\!\int_{\{s_\xi<S_*\}}\chi_c^2\,\kappa\lambda_\xi\,\omega_\xi\,\dd\xi\cdot\Big(\frac{\rad_0}{R_1}+\delta_y+\delta_U\Big)^2,
\end{equation}
$\mathfrak p_*$ is the bound of $\mathfrak p$ from the same corollary, and $\mathfrak f^c$ is the
forcing (T12), estimated in Proposition~\ref{prop:UFforce}. For the parts with $(\alpha,\beta)\ne0$
the inequality does not follow in the same way: the background terms then produce a rate, which we
treat in Proposition~\ref{prop:UFge1}. For the solution on $[T_0,T]$, $T<T_{\max}$, under the
assumptions stated at the beginning of Subsection~\ref{sssec:CUac}, the inequality is understood,
and proved, for the energies at the radii $\rad_k-\iota$, $\iota\in(0,\epsilon_y\rad_0]$, which are
absolutely continuous, for a.e.\ $t$ (Corollary~\ref{cor:CUac}), and this is how it is used in
Theorem~\ref{thm:CUboot}. The same applies to the two propositions just named, to the first only in
its parts (f), (g).
\end{prop}

\begin{proof}
\emph{Step 1: the leaf inequalities, regime by regime.} We integrate the leaf inequalities against
$\rad_0^{-3}\chi_c^2\omega_\xi\dd\xi$, regime by regime. On each leaf, every term except the local
transport has already been estimated.
\begin{itemize}
\item \emph{Crude leaves:} Proposition~\ref{prop:ACpf}. Its quadratic term (N) is replaced by
Proposition~\ref{prop:CUtri} and (H$_a$). Its coefficient terms (C) are pointwise in $\xi$ for the
$\varpi_*$-part. For the $\delta\varpi$-part, we use the remark after that proposition and (H$_d$).
The background part produces $\mathfrak q_c$.
\item \emph{Buffer leaves} ($s_\xi\in[S_*-\frac12,S_*)$) \emph{that lie in $Q^*$ for a cube $Q$
meeting $\{s\ge S_*\}$:} as for crude leaves, but with the linear term estimated by
Lemma~\ref{lem:CUbuf}. Its hypotheses hold by (H$_c$) and by the closeness in (H$_b$), which is
taken on $Q^*$. The coercivity is at least $c_{\rm buf}$.
\item \emph{Pre-band and band leaves:} Proposition~\ref{prop:Bband}, with $\chi=0$ in the pre-band
regime. In its proof:
\begin{itemize}
\item (T4) is the background term of the cubewise trilinear bound, bounded with (H$_b$);
\item the sources (T5) are those of Proposition~\ref{prop:TLsrc}: linear couplings with pointwise
coefficients, and quadratic terms as for crude leaves;
\item (T6), (T7), (T8), (T10), (T12) are per leaf, and (T11) is Proposition~\ref{prop:FM}, with
$\mathfrak p_\xi\le\mathfrak p_*$ by (H$_e$).
\end{itemize}
\item \emph{Post-band leaves:} the same, with $\chi=1$ and $H'=0$.
\item \emph{Transport.} The local part is Lemma~\ref{lem:CUtr}, which uses part of (T8). The
nonlocal part and the kernel term are of the form $C_0t^{-2}\sqrt{E^c_nE^c_{n+1}}$, as shown after
that lemma.
\item \emph{Radius.} The radius term comes from $\frac{d}{dt}w_n(\rad)^2=2n\frac{\rad'}{\rad}w_n^2$,
as in Lemma~\ref{lem:analytic}. Since $\rad$ is the same for all leaves, we can take it out of the
$\xi$-integral.
\item \emph{Damping.} $D\omega_\xi=-2H'(s_\xi)\dot s_\xi\omega_\xi$ on each leaf, for a.e.\ $t$. We
use this damping in each regime. Since $H$ is piecewise linear, $\omega_\xi$ is Lipschitz in $t$,
uniformly in $\xi$, and the identity holds except at the at most two times per leaf with
$s_\xi\in\{S_*,S_b+1\}$. By Corollary~\ref{cor:CUac}, the cell energies at the radii $\rad_k-\iota$,
with a margin $\iota>0$, are absolutely continuous, and the differential inequalities below hold for
a.e.\ $t$. Where $\chi<1$, both one-sided derivatives of $H$ are $\ge C_1^B$, since
$C_1^U\ge C_1^B$.
\end{itemize}
\emph{Step 2: the terms to be absorbed.} The coercive parts are leafwise: they add up to
$\rad_0^{-3}\int\chi_c^2\,c(\xi)\,\kappa\,\omega_\xi D(\xi)\dd\xi$, with $c(\xi)\ge c_{\min}$ on
every leaf, $c(\xi)\ge c'_{\min}$ on the leaves with $s_\xi\ge S_*$, and $c(\xi)\ge c_{\rm buf}$ on
the buffer leaves of the previous item. We show that the quadratic terms, the background terms, and
the coefficient and $\delta\varpi$-terms each remove at most $\frac14c(\xi)$, leaf by leaf, the last
ones under the smallness conditions of the crude leaf inequality, as before.

\emph{Step 3: the quadratic terms.} We use the cubewise form of Proposition~\ref{prop:CUtri}. We
split the cubes into the set $J_{\rm cr}$ of those with $Q_j^*\subset\{s<S_b+1\}$ and the set
$J_{\rm rel}$ of the others. If $j\in J_{\rm rel}$, some leaf of $Q_j^*$ has $s\ge S_b+1$, and by
Lemma~\ref{lem:CUcomp}(i) we deduce that every leaf of $Q_j^*$ has
$s\ge(1-\epsilon_0)(S_b+1)\ge S_b+\frac12>S_*$, since
$\epsilon_0(S_b+1)\le2\epsilon_0S_\sharp\le\frac12$ by requirement~\ref{rq:CUa}. Thus
$U_{\rm rel}:=\bigcup_{j\in J_{\rm rel}}Q_j^*$ contains no crude leaf. Since
$\chi_c(\xi_j)^2\le e^{12}\chi_c(\xi)^2$ on $Q_j^*$ and every point lies in at most $27$ cubes
$Q_j^*$, for any set $J$ of cubes we have
\begin{equation}\label{eq:CUcell-cubes}
\sum_{j\in J}\chi_c(\xi_j)^2\,\mathfrak d_{Q_j^*}\le27e^{12}\,\rad_0^{-3}\!\int_{U_J}\chi_c^2\,\omega_\xi
D(\xi)\dd\xi,\qquad U_J:=\bigcup_{j\in J}Q_j^*,
\end{equation}
and the factor $27e^{12}$ is part of $C_N$. By (H$_a$), the quadratic terms are therefore bounded by
\begin{equation}\label{eq:CUcell-quad}
C_N\delta_b\,\kappa\,\rad_0^{-3}\!\int_{U_{\rm rel}}\chi_c^2\omega_\xi
D\dd\xi+C_Nc_\sharp^{1/2}\,\kappa\,\rad_0^{-3}\!\int\chi_c^2\omega_\xi D\dd\xi .
\end{equation}
On a crude leaf, including a buffer leaf, only the second integral contributes, with coefficient
$\le c_{\min}/8\le c(\xi)/8$, since $c(\xi)\in\{\eta c_h(K),c_{\rm buf}\}$ is at least $c_{\min}$.
In particular the quadratic term (N) of the crude leaf inequality on a buffer leaf, whose
condition reads $C_N\eta^{1/2}c_\sharp^{1/2}\le c_{\rm buf}/4$, is covered by $C_Nc_\sharp^{1/2}\le
c_{\min}/8$, since $\eta\le1$. On a leaf with $s_\xi\ge S_*$ the coefficient is $\le
c'_{\min}/8+c_{\min}/8\le c(\xi)/4$.

\emph{Step 4: the background terms.} The background terms occur only on leaves with $s_\xi\ge S_*$,
and we use the background variant of the cubewise trilinear bound on the cubes $Q$ that meet $\{s\ge
S_*\}$. In its first configuration the supremum is $\sup_Q\|f_\xi\|_X\le\delta_c$, and the two
dissipations are those of the leaves of $Q$ with $s\ge S_*$. In its second configuration the
supremum of $(\omega\lambda_\xi)^{1/2}\|f_\xi\|_{Y_2}$ over $Q$ is bounded, by Lemma~\ref{lem:CUsob}, by
the dissipation of the levels $(0,\gamma)$, $|\gamma|\le2$, on $Q^*$. As in Step~3,
every leaf of $Q^*$ has $s\ge(1-\epsilon_0)S_*\ge S_*-\frac12$, since $\epsilon_0S_*\le\frac12$ by
the same requirement. Thus these dissipations are supported on leaves that are either in $\{s\ge
S_*\}$ or buffer leaves, where $c(\xi)\ge c''_{\min}$. With (H$_b$), $C_N\delta_b\le c''_{\min}/4$,
the background terms therefore remove at most $\frac14c''_{\min}\le\frac14c(\xi)$ on each leaf, as
explained in Remark~\ref{rem:CUbg}. The first configuration is leafwise: it is the term (T4)
above, controlled by requirement~\ref{rq:Wdc}.

\emph{Step 5.} Hence the dissipation left is at least $\frac14c(\xi)\ge\frac14c_{\min}$ on every
leaf, and \eqref{eq:CUcell} follows.
\end{proof}

\begin{remark}[The background terms and the constant $c_{\min}$]\label{rem:CUbg}
The splitting of (H$_a$) by the type of cube removes $c_h(K)$ from the condition on $\delta_b$ in
the quadratic terms. In the background terms, the supremum of the background in $Y_2$ over a cube
$Q$ that meets $\{s\ge S_*\}$ is controlled by the dissipation on $Q^*$, which may contain crude
leaves with $s\in[S_*-\frac12,S_*)$. With the crude coercivity $\eta c_h(K)$ on these leaves, the
condition on $\delta_b$ would involve $c_h(K)$, and $K$ depends on $\delta_1$ through $S^\sharp$.
Lemma~\ref{lem:CUbuf} avoids this. By Corollary~\ref{cor:Aleafentry}, the buffer leaves are
$\delta_c$-close to equilibrium, so their linear term is hypodissipative with the constant $c_W/2$
of the near-equilibrium theory. There the background appears in two ways: through $\|f_\xi\|_X$,
which is small pointwise, and through $\|f_\xi\|_{Y_2}\le C(q_\pm)K$, which multiplies a term of
lower order controlled by the crude damping. In particular the buffer leaves need no background
estimate on a larger set of cubes, and the argument does not propagate outwards. All conditions on
$\delta_b$ involve only constants chosen before $\delta_1$.
\end{remark}

\subsection{Proofs for the uniformly local energies}\label{ssec:cellsproofs}

Here we prove the results stated above and the tools for Section~\ref{sec:boot}: one
radius for all leaves, comparability on cubes, the Sobolev inequality on cubes, the trilinear sums
in uniformly local form, and local existence in cell norms.

Throughout, the notation is that of this section and of Section~\ref{sec:boot}. In particular
$w_n=w_n(\rad)=\rad^n(n+1)^5/n!$, as in \eqref{eq:Aenergy}, and
$W_{\alpha\beta}=w_{|\beta|}(\rad_{|\alpha|})\delta_y^{|\alpha|}$ are the graded weights.

\subsubsection{One radius for all leaves}

\begin{lem}[No radius drop]\label{lem:CUnodrop}
Recall that $\hat{\mathbb A}_{\ge1,\xi}$ is the part with $(\alpha,\beta)\ne(0,0)$
(Subsection~\ref{ssec:oneleaf}(a)), and choose $c_\sharp\le\delta_1^2/4$. Suppose that
at the entry time $t_\xi$ of a leaf:
\begin{itemize}
\item $\|f_\xi\|_{X(\hat m)}\le\delta_1/C_e$, as given by Corollary~\ref{cor:Aleafentry};
\item $\hat{\mathbb A}_{\ge1,\xi}(t_\xi)\le c_\sharp$, at the radius $\rad(t_\xi)$ and with the
weight $\delta_y$.
\end{itemize}
Then $\hat{\mathbb A}_\xi(t_\xi)\le\frac54\delta_1^2$, at the same radius and weight, provided
$C_e^2\ge\eta+C_\natural$. Thus we may take
\begin{equation}
r_\xi(t):=\rad(t)=\rad_0-\int_{T_0}^t|\rad'|,\qquad\delta_{y,\xi}:=\delta_y\qquad\text{for
all }\xi .
\end{equation}
Then $\rad_{\rm post}:=\rad_0/2$, and the radius need not be reduced at the entry time of a
leaf.
\end{lem}

\begin{proof}
We have $\hat{\mathbb A}_\xi=e_\xi[\varphi_{00}]+\hat{\mathbb A}_{\ge1,\xi}$, and
$e_\xi[\varphi_{00}]=\eta\|f_\xi\|^2_X+\chi\nn{f_\xi}^2_\natural\le(\eta+C_\natural)\|f_\xi\|^2_X\le\delta_1^2$.
Moreover, $c_\sharp\le\delta_1^2/4$. Hence we have
$\hat{\mathbb A}_\xi(t_\xi)\le\delta_1^2+\frac14\delta_1^2=\frac54\delta_1^2$.

The choice of $c_\sharp$ is compatible with requirement~\ref{rq:Wcsharp}, which bounds $c_\sharp$
above by a constant depending on $K$. Here it is also bounded by $\delta_1^2/4$, and $\delta_1$ is
fixed before $K$. In the crude regime of the bootstrap, we obtain the bound
$\hat{\mathbb A}_{\ge1}\le c_\sharp$ by making $\rad_0$, $\delta_y$, $T_0^{-1}$ and $\eps$ small,
and these are chosen after $c_\sharp$.
\end{proof}

\subsubsection{Comparability on cubes}

We use the cells, cubes and requirements defined above.

\begin{lem}[Comparability on cubes]\label{lem:CUcomp}
Under requirements~\ref{rq:CUa}--\ref{rq:CUb}, for $\xi,\xi'\in Q_j^*$:
\begin{enumerate}[label=(\roman*)]
\item $|s_\xi/s_{\xi'}-1|\le\epsilon_0$ and $e^{-1/2}\le\omega_\xi/\omega_{\xi'}\le e^{1/2}$;
\item $e^{-6}\le\chi_c(\xi)/\chi_c(\xi')\le e^6$;
\item the quadratic forms $\ell_\xi$ and $d_\xi$ (as forms in $h^{(x)}$, a function of $z$) are
comparable within the factor $e^{6}$, and so are the leaf rates $\lambda(p_g(t,\xi))$.
\end{enumerate}
Consequently, $\mathfrak a_{Q_j}\le e^{5}\,\omega_{\xi_j}^{-1}\,\mathfrak A_{\xi_j}$, $
\mathfrak d_{Q_j}\le e^5\mathfrak D_{\xi_j}$, and $\sum_j\chi_c(\xi_j)^2
\mathfrak d_{Q_j^*}\le27e^{12}\mathfrak D_c$.
\end{lem}

\begin{proof}
(i) We have $s_\xi=\lambda_{\rm ref}(\xi)(\tau(t)-\tau(T_0))$ and
$|\nabla\log\lambda_{\rm ref}|\le2L'$, with $|\xi-\xi'|\le3\sqrt3\rad_0$, which implies that
$|\log(s_\xi/s_{\xi'})|\le6\sqrt3L'\rad_0=\epsilon_0/2$. Since $\epsilon_0\le1$ by
requirement~\ref{rq:CUa} and $e^x-1\le2x$ for $x\in[0,1]$, we obtain
$|s_\xi/s_{\xi'}-1|\le\epsilon_0$, and hence $|s-s'|\le\epsilon_0\max(s,s')$. Moreover,
$|H(s)-H(s')|\le C_1^U|s-s'|\le C_1^U\epsilon_0\max(s,s')$. If $\max(s,s')\le2S_\sharp$ this is
$\le\frac14$ by the same requirement. If $\min(s,s')>S_b+1$, then $H(s)=H(s')$. In the remaining
case we have $s,s'\le(1+\epsilon_0)(S_b+1)\le2S_\sharp$.

(ii) This follows from $|\jap a-\jap b|\le|a-b|\le3\sqrt3$.

(iii) For $\ell_\xi$ the weight is $(A_{p_g}/A_{\rm ref}(\xi))^2m_x^2$, with
$m_x(z)=A_{\rm ref}(\xi-z/t)e^{q|z|^2}$ and $A_p=\rho^{-1}\theta^{3/4}$, so $A_{\rm ref}$ enters
twice. Since $|\nabla\log A_{p_g}|\le|\nabla\ell_g|+\frac34|\nabla\vartheta_g|\le C_\varpi$ and, by
Corollary~\ref{cor:LWPa}(a4), $|\nabla\log A_{\rm ref}|\le L\le2L$, we have
\begin{equation}
|\nabla_\xi\log\big((A_{p_g}/A_{\rm ref}(\xi))^2m_x^2\big)|\le2|\nabla\log A_{p_g}|+2|\nabla\log
A_{\rm ref}(\xi)|+2|\nabla\log A_{\rm ref}(\xi-z/t)|\le2C_\varpi+8L .
\end{equation}
For $d_\xi$, we write
\begin{equation}
\|\varphi\|^2_{Y_2(\hat m)}=\int\jap{\hat z}^{\gamma+2}|v|^2+\jap{\hat z}^\gamma\big(|\nabla
v|^2+|\hat z|^2|\nabla v|^2-(\hat z\cdot\nabla v)^2\big),\qquad v=\hat m\varphi ,
\end{equation}
which is the definition of $\tn$ written with polynomial coefficients. Under a change of the frame
center by $a$ and of the scale, a direct computation shows that the $\xi$-derivative of the form, at
fixed weight, is bounded by $C_d|\nabla\varpi_g|\le C_dC_\varpi$ times the form, with
$C_d=C_d(\gamma,\theta_-)$. The weight $\hat m$ contributes $2C_\varpi+8L$, as for $\ell_\xi$. For
the center, we use $|\hat z|\,|(a\cdot\hat r)|\tau|^2-\rho(a\cdot\tau)|\le|a|(\rho^2+\jap{\hat
z}^2|\tau|^2)$, where $\nabla v=\rho\hat r+\tau$ and $\tau\perp\hat r$. Hence the logarithmic
$\xi$-derivatives of both forms are $\le(4+C_d)C_\varpi+8L$. Integrating over a segment of length
$\le3\sqrt3\rad_0$ and using requirement~\ref{rq:CUb}, we obtain the factor $e^{3\sqrt3}\le e^6$.
The rate $\lambda(p_g)=\rho_g\theta_g^{\gamma/2}$ has
$|\nabla\log\lambda(p_g)|\le|\nabla\ell_g|+\frac32|\nabla\vartheta_g|\le\frac32C_\varpi$, and is
treated in the same way.

For the consequences, note that $\chi_{\xi_j}\ge e^{-2}$ on $Q_j$ and that $\omega$ varies by at
most $e^{1/2}$ on $Q_j$, which yields the factor $e^{4}e^{1/2}\le e^5$. Each point lies in at most
$27$ cubes $Q_j^*$.
\end{proof}

\subsubsection{Sobolev at the analyticity scale}

We put $\tilde w_n:=w_n/(n+2)^2$.

\begin{lem}[Cube Sobolev without loss of radius]\label{lem:CUsob}
There is an absolute $C_S$ such that, for every cube $Q_j$ and every $(\alpha,\beta)$ with
$|\beta|=n$:
\begin{equation}\begin{aligned}
\sup_{\xi\in Q_j}\ell_\xi[h_{\alpha\beta}]&\le C_S\sum_{|\gamma|\le2}\rad_0^{2|\gamma|-3}\!\int_{Q_j^*}\ell_{\xi'}[h_{\alpha,\beta+\gamma}]\dd\xi',\\
\sup_{\xi\in Q_j}\omega_\xi d_\xi[\varphi_{\alpha\beta}]&\le
C_S\sum_{|\gamma|\le2}\rad_0^{2|\gamma|-3}\!\int_{Q_j^*}\omega_{\xi'}d_{\xi'}[\varphi_{\alpha,\beta+\gamma}]\dd\xi'
.
\end{aligned}\end{equation}
Consequently, with the weights $\tilde w$,
\begin{equation}
\sum_n\tilde w_n^2\sum_{|\beta|=n,\alpha}\delta_y^{2|\alpha|}\sup_{Q_j}\eta\ell_\xi[h_{\alpha\beta}]\le
C_S'\,\mathfrak a_{Q_j^*},
\end{equation}
and similarly for $\omega_\xi d_\xi$ with $\mathfrak d_{Q_j^*}$, and for the parts with
$(\alpha,\beta)\ne0$.
\end{lem}

\begin{proof}
We fix the quadratic form $q_j:=\ell_{\xi_j}$ at the center of the cube. By
Lemma~\ref{lem:CUcomp}(iii), $e^{-6}q_j\le\ell_\xi\le e^6q_j$ on $Q_j^*$. The map
$\xi\mapsto u(\xi):=h^{(x)}_{\alpha\beta}$, with $x=t\xi$, takes values in the Hilbert space of
$q_j$, and $\partial^\gamma_\xi u=h^{(x)}_{\alpha,\beta+\gamma}$. By the scaled Sobolev embedding
$H^2(Q^*)\subset L^\infty(Q)$ in $\R^3$, for Hilbert-space-valued functions, we have
\begin{equation}
\sup_Q\|u\|^2\le C\sum_{|\gamma|\le2}\rad_0^{2|\gamma|-3}\int_{Q^*}\|\partial^\gamma u\|^2 .
\end{equation}
Then we return to $\ell_{\xi'}$ with the factor $e^6$. We treat the dissipation in the same way,
with the form $\omega_{\xi_j}d_{\xi_j}$ and parts~(i), (iii) of that lemma. We use no derivative of
a form, only the comparability.

For the weighted sum, we have $w_n\le(2/\rad_0)^j(n+2)^jw_{n+j}$ for $j\le2$, because
$\rad\ge\rad_0/2$. Hence we have $\rad_0^{|\gamma|}\tilde w_n\le4w_{n+|\gamma|}$. Each multi-index
of order $n+|\gamma|$ arises from at most $10$ pairs $(\beta,\gamma)$. The nonlocal part of $e_\xi$
is $\le C_\natural\ell_\xi$, so the local part suffices for upper bounds.
\end{proof}

The scale is essential. On cubes of side $1$ the same argument produces a factor $\rad_0^{-4}$,
which the smallness of the higher levels, of order $\rad_0^2$, cannot compensate. On cubes of
side $\rad_0$ each $\xi$-derivative produces a factor $\rad_0^{-1}$, which is the loss of one
level in the analytic norm.

\subsubsection{Trilinear sums in uniformly local form}

\begin{prop}[Uniformly local trilinear bound]\label{prop:CUtri}
For every leaf let $\tau_\xi(f,g;h)$ be a trilinear form with
\begin{equation}
|\tau_\xi(f,g;h)|\le C_\tau\big(\|f\|_X\|g\|_{Y_2}+\|f\|_{Y_2}\|g\|_X\big)\|h\|_{Y_2},
\end{equation}
in the norms $X(\hat m_\xi)$ and $Y_2(\hat m_\xi)$. Put
\begin{equation}
\mathcal N_\xi:=\sum_{(\alpha,\beta)}W_{\alpha\beta}^2\!\!\sum_{0<(\alpha',\beta')<(\alpha,\beta)}\!\!\tbinom\alpha{\alpha'}\tbinom\beta{\beta'}\big|\tau_\xi(\varphi_{\alpha'\beta'},\varphi_{\alpha''\beta''};\varphi_{\alpha\beta})\big|
.
\end{equation}
Then, with $C_N$ depending only on $(N_y,C_\tau,\eta)$,
\begin{equation}
\rad_0^{-3}\!\int\chi_c^2\,\omega_\xi\,\kappa\lambda_\xi\,\mathcal N_\xi\dd\xi\le
C_N\,\kappa\sum_j\chi_c(\xi_j)^2\,\mathfrak a_{Q_j^*,\ge1}^{1/2}\,\mathfrak d_{Q_j^*}\le
C_N\,\kappa\Big(\sup_j\mathfrak a_{Q_j^*,\ge1}\Big)^{1/2}\mathfrak D_c .
\end{equation}
(The absolute factor $27e^{12}$ of Lemma~\ref{lem:CUcomp} is included in $C_N$.) If the sum is
restricted to leaves in a set $\Sigma$, the supremum may be restricted to cubes $Q_j$ that meet
$\Sigma$.

The same holds for the \emph{background terms}, in which $(\alpha',\beta')=0$ and
$\varphi_{00}=f_\xi$, with three changes:
\begin{itemize}
\item the term $\|f\|_X\|g\|_{Y_2}\|h\|_{Y_2}$ is bounded by $\sup_{Q_j}\|f_\xi\|_X\cdot
\mathfrak d$;
\item the term $\|f\|_{Y_2}\|g\|_X\|h\|_{Y_2}$ is bounded with $\mathfrak a_{Q_j}$, the
\emph{total} undamped cube energy, in place of $\mathfrak a_{Q_j^*,\ge1}$;
\item the factor $W^2_{\alpha\beta}$ is split as $W_{\alpha\beta}\cdot W_{\alpha\beta}$, without
binomials.
\end{itemize}
\end{prop}

\begin{proof}
We fix a cube $Q=Q_j$. For a splitting with $|\beta'|=k$ and $|\beta''|=n-k$, call the factor of
lower $v$-level the \emph{low} factor. The bound on $\tau$ is symmetric in $(f,g)$, so we may assume
$k\le n/2$. We distribute $\omega_\xi\lambda_\xi=(\omega_\xi\lambda_\xi)^{1/2}(\omega_\xi\lambda_\xi)^{1/2}$
over the two factors measured in $Y_2$. We distinguish two configurations.
\begin{itemize}
\item \emph{Low factor in $X$.} We take $\sup_Q$ of its undamped $X$-norm. By Lemma~\ref{lem:CUsob},
the $\tilde w$-weighted sum of these suprema is $\le C\mathfrak a_{Q^*,\ge1}^{1/2}$. The two other
factors are damped dissipations, integrated over $Q$.
\item \emph{Low factor in $Y_2$.} We take $\sup_Q$ of $(\omega\lambda_\xi)^{1/2}\|\cdot\|_{Y_2}$. By
the same lemma, the $\tilde w$-weighted sum is $\le C\mathfrak d_{Q^*}^{1/2}$. The high factor is
measured in the undamped $X$-norm and integrated over $Q$, so its contribution is
$\le\mathfrak a_{Q,\ge1}^{1/2}$. The third factor is a damped dissipation.
\end{itemize}
In both cases,
\begin{equation}
\rad_0^{-3}\!\int_Q\big(\sup_Qa_{\beta'}\big)\,b_{\beta''}(\xi)\,d_\beta(\xi)\dd\xi\le
A_{\beta'}B_{\beta''}D_\beta,
\end{equation}
where $A_{\beta'}=\sup_Qa_{\beta'}$, $B_{\beta''}=(\rad_0^{-3}\int_Qb_{\beta''}^2)^{1/2}$ and
$D_\beta=(\rad_0^{-3}\int_Qd_\beta^2)^{1/2}$. Then we apply Lemma~\ref{lem:algebra} to $(A,B,D)$.

The low factor carries $\tilde w_k$ instead of $w_k$, so the extra factor $(k+2)^2$ appears. With
$p=5$, $\varrho_{n,k}(k+2)^2\le4(k+1)^{-2}$, and $\sum_k(k+1)^{-2}<\infty$. By Young's inequality
$\ell^1*\ell^2\cdot\ell^2$, with Cauchy--Schwarz in $k$, we obtain
\begin{equation}
\le C\Big(\sum_k\tilde w_k^2A_k^2\Big)^{1/2}\|wB\|_{\ell^2}\|D\|_{\ell^2}
\end{equation}
The $y$-binomials are bounded by $2^{N_y}$, and
$\delta_y^{2|\alpha|}=\delta_y^{|\alpha'|}\delta_y^{|\alpha''|}\delta_y^{|\alpha|}$ splits exactly.
Both factors have $(\alpha,\beta)\ne0$, so only the parts $\ge1$ enter. We multiply by
$\chi_c^2\le e^{12}\chi_c(\xi_j)^2$ on $Q_j$ and sum over $j$. For the last inequality we use
Lemma~\ref{lem:CUcomp}.

For the background terms, $f$ is the low factor (level $0$). In the first configuration we have
$\sup_Q\|f\|_X$ directly. In the second, the cube Sobolev lemma bounds
$\sup_Q(\omega\lambda_\xi)^{1/2}\|f\|_{Y_2}$ by the damped dissipation of the levels $(0,\gamma)$
with $|\gamma|\le2$, with weights $\rad_0^{2|\gamma|}\le Cw_{|\gamma|}^2$. The factor
$g=\varphi_{\alpha\beta}$ ranges over all levels, so the total $\mathfrak a_Q$ enters.
\end{proof}

The $\delta\varpi$-parts of the coefficient terms are differences
$\mathfrak g(u_g)- \mathfrak g(u_*)$ of analytic functions of the parameters. They are not bilinear
in $(\delta\varpi,\varphi)$. The proof of the third bound of Lemma~\ref{lem:TWcoef}, which combines
Taylor's formula with Lemma~\ref{lem:TWcalc}(vi), in the local theory in the truncated form of
Remark~\ref{rem:TWtrunc}, shows that their cube norms $S_Q$, $E_Q$ are bounded by $C_{\rm Max}$
times those of $\delta\varpi$, and $E_Q(\delta\varpi)\le C\mathfrak m_Q^{1/2}$. The smallness
condition of~(vi) is requirement~\ref{rq:TWa}. With this, we sum them in the same way, the
difference playing the role of $\delta\varpi$. We take the supremum on the factor of lower level,
that is, on $\partial^{\le k}\delta\varpi$ via the cube Sobolev inequality for $\mathfrak m_{Q^*}$
if $\delta\varpi$ is the lower factor, and on $\varphi$ via the cube Sobolev lemma otherwise. Thus
we obtain $C(\sup_Q\mathfrak m_{Q^*})^{1/2}$ times cell dissipation and energy terms.

\subsubsection{Local existence in the uniformly local space}

\begin{lem}[Viscous products in the analytic weights]\label{lem:CUvisc}
Let $w_n=w_n(\rad)=\rad^n(n+1)^5/n!$, $\tilde w_n:=w_n/(n+2)^2$, and
$A_n(f):=(\sum_{|\beta|=n}|\partial^\beta f|^2)^{1/2}$.
\begin{enumerate}[label=(\roman*)]
\item $A_{n+1}(f)^2\le\sum_{|\beta|=n}|\nabla\partial^\beta f|^2\le3A_{n+1}(f)^2$ and $A_k(\nabla
f)\le\sqrt3A_{k+1}(f)$.
\item $(n+1)w_n\le64\rad w_{n-1}$ for $n\ge1$, and $w_N\ge(\frac23)^5\frac{N+1}\rad w_{N+1}$
for $N\ge1$.
\item For $0\le k\le n/2$, $w_n\binom nk(k+1)\le\dfrac{C_4(n+1)}{\rad^2(k+2)}\,\tilde
w_{k+1}\,w_{n-k+1}$, with $C_4=16$.
\end{enumerate}
Let $\langle f,g\rangle_c:=\int\Psi_cf\cdot g\dd\xi$ with a weight $\Psi_c\ge0$, $\|f\|_c:=\langle
f,f\rangle_c^{1/2}$, and let $u$ be a function of $\xi$ with
\begin{equation}
a_n:=w_n\|A_n(u)\|_c,\qquad g_n:=w_n\Big(\sum_{|\beta|=n}\|\nabla\partial^\beta u\|_c^2\Big)^{1/2} .
\end{equation}
For functions $F,G$ of $\xi$ (vector-valued, or $H$-valued with a bounded bilinear product) put
\begin{equation}
\hat E^{\le N}_c(G):=\Big(\sum_{k\le N}w_{k+1}^2\|A_k(G)\|_c^2\Big)^{1/2},\qquad\hat S^{\le
k}(F):=\sup_Q\Big(\sum_{j\le k}\tilde w_{j+1}^2\sup_{\xi\in Q}A_j(F)(\xi)^2\Big)^{1/2},
\end{equation}
where $Q$ runs over the cubes $Q_j$. These are the norms of Lemma~\ref{lem:TWcalc} shifted by one
level, which are the natural norms of a first derivative. The supremum over $Q$ is taken outside the
sum over the levels. Then, for every $N\ge1$,
\begin{equation}
\sum_{n=1}^Nw_n^2\sum_{|\beta|=n}\big|\langle\partial^\beta
u,\partial^\beta(FG)\rangle_c\big|\le\frac{C_{\rm vp}}\rad\Big(\sum_{n\le
N-1}g_n^2\Big)^{1/2}\Big(\hat S^{\le N/2}(F)\,\hat E^{\le N}_c(G)+\hat E^{\le N}_c(F)\,\hat S^{\le
N/2}(G)\Big),
\end{equation}
with $C_{\rm vp}$ absolute (it is not the constant $C_V$ of Lemma~\ref{lem:LCvisc}). If $F$ is of
order zero (for instance $F=e^{-\vartheta}$), the same holds with $\tilde w_k$, $w_k$ in place of
$\tilde w_{k+1}$, $w_{k+1}$ in the norms of $F$, and with $C_{\rm vp}$ in place of $C_{\rm
vp}/\rad$. The term $n=0$ is $\le a_0\sup|F|\,\|G\|_c$.
\end{lem}

\begin{proof}
Part (i) follows by counting the multi-indices. (ii) follows from
$w_n/w_{n-1}=\frac\rad n(\frac{n+1}n)^5$. For (iii), inserting the weights, we see that the left
side divided by $\frac{n+1}{\rad^2(k+2)}\tilde w_{k+1}w_{n-k+1}$ equals
$\frac{(k+1)^2(k+3)^2}{(k+2)^4}\cdot\frac{(n+1)^4(n-k+1)}{(n-k+2)^5}\le\Big(\frac{n+1}{n-k+2}\Big)^4\le16$,
since $n+1\le2(n-k+2)$ for $k\le n/2$.

By Cauchy--Schwarz in $\beta$, we see that the inner sum is $\le\langle A_n(u),A_n(FG)\rangle_c$. By
the pointwise Leibniz bound~\eqref{eq:TWpt}, we have
$A_n(FG)\le\sum_k\binom nk(\min(k,n-k)+1)A_k(F)A_{n-k}(G)$. For $k\le n/2$, taking the supremum on
$F$ and using (iii), we obtain
\begin{equation}
w_nA_n(FG)\le\frac{C_4(n+1)}{\rad^2}\sum_{k\le n/2}\frac{\sigma_k}{k+2}\,w_{n-k+1}A_{n-k}(G)+(\text{the
symmetric sum}),\qquad\sigma_k:=\tilde w_{k+1}\sup A_k(F).
\end{equation}
Therefore every Leibniz term loses one level: two factors with $k+1$ and $n-k+1$ derivatives on $u$
are paired with a test function with $n$. We move the loss onto the test factor. By (ii) and (i), we
have $(n+1)a_n\le64\rad w_{n-1}\|A_n(u)\|_c\le64\rad g_{n-1}$. Hence we see that the $n$-th term is
$\le\frac{64C_4}\rad g_{n-1}\sum_k\frac{\sigma_k}{k+2}y_{n-k}$ with $y_m:=w_{m+1}\|A_m(G)\|_c$. We
sum over $n$ by Cauchy--Schwarz and Young's inequality $\ell^1*\ell^2\to\ell^2$, with
$\sum_k\frac{\sigma_k}{k+2}\le\zeta(2)^{1/2}\hat S^{\le N/2}(F)$. The symmetric sum is the same with
$F$ and $G$ exchanged. For an $F$ of order zero the ratio in (iii) is $\le C(n+1)/(\rad(k+1)^2)$ for
$k\le n/2$, and $\le C/(\rad(n-k+2))$ for $k>n/2$, a term without loss. The levels used are
$\le n/2$ on the factor in supremum, and $\le n+1$ on $u$ in the cell norm.

Up to the integration against $\Psi_c$ the argument is pointwise in $\xi$. We apply it on each cube
$Q$, with the weight $\Psi_c\mathbf 1_Q$ in place of $\Psi_c$ and with $\sup_{\xi\in Q}$ in
$\sigma_k$. In this way we obtain the bound with $g_{n,Q}$, $\hat E^{\le N}_{c,Q}$ and
$\hat S_Q^{\le N/2}$, the quantity under $\sup_Q$ in the definition of $\hat S$. Since
$\|\cdot\|_c^2=\sum_Q\|\cdot\|_{c,Q}^2$, by the Cauchy--Schwarz inequality over the cubes we have
\begin{equation}
\sum_Q\Big(\sum_ng_{n,Q}^2\Big)^{1/2}\hat S^{\le N/2}_Q(F)\,\hat E^{\le N}_{c,Q}(G)\le\hat S^{\le
N/2}(F)\Big(\sum_ng_n^2\Big)^{1/2}\hat E^{\le N}_c(G),
\end{equation}
and the same for the symmetric term. The constants do not depend on $N$.
\end{proof}

\begin{lem}[Gronwall for a family]\label{lem:LCgron}
Let $U^c\ge0$, $c\in\R^3$, be locally Lipschitz on $(t_*,T]$ with $\limsup_{s\downarrow
t_*}U^c(s)\le U^c(t_*)$ for each $c$, let $S:=\sup_cU^c$ be bounded on $[t_*,T]$ and lower
semicontinuous on $(t_*,T]$, and let $\frac{d}{dt}U^c\le AU^c+BS+\sigma$ for a.e.\ $t\in(t_*,T]$ and
every $c$, with constants $A,B,\sigma\ge0$ (the constant source $\sigma$ is $0$ in most
applications). Then
\begin{enumerate}[label=(\roman*)]
\item
\begin{equation}\label{eq:LCgron}
U^c(t)\le e^{A(t-t_*)}U^c(t_*)+\int_{t_*}^te^{A(t-s)}\big(BS(s)+\sigma\big)\dd s
\end{equation}
for every $c$. This holds also when $\frac{d}{dt}U^c\le AU^c+BS'+\sigma$ with any bounded measurable
$S'\ge0$ in place of $S$, with $S'$ in the integral;
\item $S(t)\le e^{(A+B)(t-t_*)}\big(S(t_*)+\sigma(t-t_*)\big)$.
\end{enumerate}
\end{lem}

\begin{proof}
(i) For $t_*<s<t$, we integrate
$\frac{d}{d\tau}\big(e^{-A\tau}U^c(\tau)\big)\le e^{-A\tau}\big(BS(\tau)+\sigma\big)$ over $[s,t]$.
Note that $S$ is measurable, being lower semicontinuous, and bounded. Then we let $s\downarrow t_*$
and use $\limsup_{s\downarrow t_*}U^c(s)\le U^c(t_*)$. (ii) We take the supremum over $c$ in~(i),
with $U^c(t_*)\le S(t_*)$. Then we see that $\phi(t):=e^{-A(t-t_*)}S(t)$ satisfies
$\phi(t)\le S(t_*)+\sigma(t-t_*)+B\int_{t_*}^t\phi$, since
$\int_{t_*}^te^{-A(s-t_*)}\sigma\dd s\le\sigma(t-t_*)$. The first two terms are nondecreasing in
$t$, and thus Gronwall's inequality implies $\phi(t)\le e^{B(t-t_*)}\big(S(t_*)+\sigma(t-t_*)\big)$.
\end{proof}

\begin{lem}[Pointwise bounds from leaf energies]\label{lem:LCdata}
Let $t_*\ge1$, $q_*>0$, $\delta\in(0,q_*]$, $\underline\rad\in(0,1]$ and $M_1\ge0$, and let
$g_*\ge0$ satisfy, for all $\xi\in\R^3$, all $\alpha$ with $|\alpha|\le N_y$ and all $\beta$,
\begin{equation}
\big\|e^{q_*|z|^2}\,\partial^\alpha_y\partial^\beta_vg_*(z,\xi-z/t_*)\big\|_{L^2_z}\le
M_1\,\eps\rho_0(\xi)\,|\beta|!\,\underline\rad^{-|\beta|}.\tag{H}
\end{equation}
Then, for $|\alpha|\le N_y-2$ and all $\beta$,
\begin{equation}
|\partial^\alpha_y\partial^\beta_vg_*(y,v)|\le C_*\,|\beta|!\,(\underline\rad/2)^{-|\beta|}\,\Gamma_{q_*-\delta/4}(y,v),\qquad\Gamma_c(y,v)=\eps\rho_0(v)e^{-c|y|^2},
\end{equation}
with $C_*=C(M_1,q_*,\delta,\underline\rad,L)$.
\end{lem}

\begin{proof}
We fix $\xi$, $\alpha$ with $|\alpha|\le N_y-2$, and $\beta$, and put
$u(z):=e^{(q_*-\delta/8)|z|^2}\partial^\alpha_y\partial^\beta_vg_*(z,\xi-z/t_*)$. Since
$\partial_z[k(z,\xi-z/t_*)]=\big((\partial_y-t_*^{-1}\partial_v)k\big)(z,\xi-z/t_*)$, each
$\partial^\varsigma_zu$, $|\varsigma|\le2$, is a sum of at most $C$ terms
$P(z)e^{(q_*-\delta/8)|z|^2}\partial^{\alpha+\varsigma_1}_y\partial^{\beta+\varsigma_2}_vg_*(z,\xi-z/t_*)$,
with $|\varsigma_1|+|\varsigma_2|\le2$ and $P$ a polynomial of degree $\le2$ with bounded
coefficients ($t_*\ge1$). Since $\jap z^2e^{(q_*-\delta/8)|z|^2}\le C_\delta e^{q_*|z|^2}$, by (H)
we have $\|\partial^\varsigma_zu\|_{L^2}\le
CM_1\eps\rho_0(\xi)(|\beta|+2)!\,\underline\rad^{-|\beta|-2}$, and
$(|\beta|+2)!\le6\cdot2^{|\beta|}|\beta|!$. By Sobolev's inequality $H^2(\R^3)\subset L^\infty$, we
obtain
\begin{equation}
|\partial^\alpha_y\partial^\beta_vg_*(z,\xi-z/t_*)|\le
CM_1\underline\rad^{-2}\,\eps\rho_0(\xi)\,|\beta|!\,(\underline\rad/2)^{-|\beta|}e^{-(q_*-\delta/8)|z|^2}.
\end{equation}
Given $(y,v)$, we take $z=y$ and $\xi=v+y/t_*$. By (D3),
$\rho_0(\xi)\le\rho_0(v)e^{L|y|}\le C(L,\delta)\rho_0(v)e^{\delta|y|^2/8}$, which proves the claim.
\end{proof}

\begin{lem}[Uniqueness]\label{lem:LWuniq}
Let $0<t_*<T$ and $c>0$, and let $g_1\ge0$ and $g_2$ be classical solutions on $[t_*,T]$ of
$\partial_tg=D\cdot(\ab[g]Dg-\bb[g]g)$, $D=\partial_v-t\nabla_y$, with $g_1(t_*)=g_2(t_*)$ and
$|\partial^{\mathbf k}_{y,v}g_i|\le C\Gamma_c$ on $[t_*,T]$ for $|\mathbf k|\le2$. Then $g_1=g_2$.
Here $\ab[h](y,v)=\int a(w)h(y+tw,v-w)\dd w$, and similarly $\bb$, $\cb$. With this convention the
equation is $\partial_tg=\kappa\,\cQ_t[g]$ of Proposition~\ref{prop:leaf}, and the factor $\kappa$
of the leaf normalization does not appear.
\end{lem}

\begin{proof}
Let $d:=g_1-g_2$. Then we have
\begin{equation}
\partial_td=D\cdot\big(\ab[g_1]Dd-\bb[g_1]d\big)+D\cdot\big(\ab[d]Dg_2-\bb[d]g_2\big).
\end{equation}
We fix $q\in(0,c)$ and $m\in(\frac32,p]$, and let
$\Psi:=e^{2q|y|^2}\jap v^{-2m}(\eps\rho_0(v))^{-2}$. By the Gaussian bounds,
$\mathcal E(t):=\iint d^2\Psi$ is finite and $C^1$, and the integrations by parts below have no
boundary terms. Here $D\log\Psi=V-4qt\,y$ with $|V|\le C(m,L)$.

\emph{The first term.} We apply the identity of Lemma~\ref{lem:LWenergy} with $h=g_1$, $u=d$, and
Cauchy--Schwarz for the form $\ab[g_1]\ge0$ (we only use $g_1\ge0$). We bound the coefficients as
follows. First, $|\ab[g_1]|+|\bb[g_1]|+|\cb[g_1]|\le C$ (the singularities of $a,b,c$ are integrable
for $\gamma>-3$, and $\cb[g]=-8\pi g$ for $\gamma=-3$). Next,
$\ab[g_1]y\cdot y=\int a(w)y''\cdot y''\,g_1(y'',v-w)\dd w\le C$ with $y''=y+tw$, by
Lemma~\ref{lem:Bnull}. Finally, $|\bb[g_1]\cdot y|\le C(1+t)$, since
$b(w)\cdot y=b(w)\cdot y''+2t|w|^{\gamma+2}$. Hence we have
$\ab[g_1]D\Psi\cdot D\Psi/\Psi+|\bb[g_1]\cdot D\Psi|+|\cb[g_1]|\Psi\le C(1+T)^2\Psi$, and the
contribution of the first term is bounded by $C(1+T)^2\mathcal E$.

\emph{The second term.} Since $\partial_w[\psi(y+tw,v-w)]=-(D\psi)(y+tw,v-w)$ for every function
$\psi$, an integration by parts in $w$ shows that $D\cdot\ab[d]=\bb[d]$ and $D\cdot\bb[d]=\cb[d]$.
Thus $D\cdot(\ab[d]Dg_2-\bb[d]g_2)=\ab[d]:D^2g_2-\cb[d]g_2$, which contains no derivative of $d$.
With $|D^2g_2|+|g_2|\le C(1+T)^2\Gamma_c$, the second term is bounded by
\begin{equation}
C\iint\Psi|d|\,\Gamma_c\int k_\gamma(w)\,|d(y+tw,v-w)|\dd w\,\dd y\dd v,
\end{equation}
with $k_\gamma(w)=|w|^{\gamma+2}+2(\gamma+3)|w|^\gamma$ for $\gamma>-3$ and $k_\gamma(w)=|w|^{-1}$
for $\gamma=-3$, plus, for $\gamma=-3$, the local term $8\pi\iint\Psi d^2g_2\le C\mathcal E$. We
write it as $\iint u\,(\mathcal Ku)$ with $u=\Psi^{1/2}|d|$ and the kernel
$\mathcal K=\Psi^{1/2}(y,v)\Gamma_c(y,v)k_\gamma(w)\Psi^{-1/2}(y+tw,v-w)$. Then
\begin{equation}
\mathcal K\le\eps\,k_\gamma(w)\,\frac{\rho_0(v-w)\jap{v-w}^m}{\jap
v^m}\,e^{-(c-q)|y|^2}e^{-q|y+tw|^2}\le C\,k_\gamma(w)\,e^{-(c-q)|y|^2}e^{-q|y+tw|^2},
\end{equation}
since $\eps\le1$ and $\rho_0\jap\cdot^m\le A_F$ for $m\le p$. For fixed $(y,v)$,
$\int\mathcal K\dd w\le C\int k_\gamma(w)e^{-q|y+tw|^2}\dd w\le C(t_*,q)$. For fixed $(y',v')$, with
$(y,v)=(y'-tw,v'+w)$ (a change of variables with Jacobian $1$),
$\int\mathcal K\dd w\le C\int k_\gamma(w)e^{-(c-q)|y'-tw|^2}\dd w\le C(t_*,c-q)$. By Schur's test,
the second term is bounded by $C\mathcal E$.

Therefore we have $\mathcal E'\le C\mathcal E$ on $[t_*,T]$ with $\mathcal E(t_*)=0$, which implies
that $\mathcal E\equiv0$.
\end{proof}

\begin{prop}[Local existence in cell norms]\label{prop:CULW}
Let $t_*\ge T_0$. Let $\rad_k(t_*)=\rad(t_*)-k\epsilon_y\rad_0$, $k\le N_y$, be graded
radii at $t_*$, with $\rad(t_*)\le\rad_0$ and
$\rad_k(t_*)\ge\underline\rad:=\rad_0/2$, and let $q_*$ be the Gaussian exponent at $t_*$.
Let $H_\epsilon\le H$ be the one-sided mollification of $H$ defined in~(2) of the proof, at the
scale $\epsilon=1/(40C_1^U)$, and let $\Omega_\xi\in\{1,\Omega^\epsilon_\xi\}$, where
$\Omega^\epsilon_\xi:=e^{-2H_\epsilon(s_\xi(t_*))}$ is frozen at $t_*$. Write $\mathfrak B^\Omega_c$
and $\mathfrak B^\Omega_{c,\ge1}$ for the joint energies $\mathfrak A_c+K_{\mathfrak M}
\mathfrak M_c$ and $\mathfrak A_{c,\ge1}+K_{\mathfrak M}\mathfrak M_c$ of Subsection~\ref{ssec:UF}, with
$K_{\mathfrak M}=T_0$, with $\omega_\xi$ replaced by $\Omega_\xi$, and at the radii
\begin{equation}
\tilde\rad_k(t):=\rad_k(t_*)-\int_{t_*}^tr_{\rm loc},\qquad r_{\rm
loc}(t):=\Big(\tfrac34C_0'+\big(C_1^F+C_1^{\rm mac}\big)T_0^{1/2}\Big)t^{-2},
\end{equation}
with the constants $C_0'$, $C_1^F$, $C_1^{\rm mac}$ of Remark~\ref{rem:UFradius} (those of the
bootstrap). Write $\mathfrak B^{(1)}$ when $\Omega\equiv1$, and put
$M_\natural:=e^{2H(S_b+1)}(c_\chi C_U+\mu_\sharp)$ and $m_\sharp:=c_\sharp/(54e^5)$. Assume:
\begin{itemize}
\item[(H1$'$)] $\sup_c\mathfrak B^{(1)}_c(t_*)\le M$ with $m_\sharp\le M\le M_\natural$, and
$\sup_c\mathfrak B^{(1)}_{c,\ge1}(t_*)\le m_\sharp$;
\item[(H2$'$)] (G1)--(G3) hold, in the normalizations of Remark~\ref{rem:BTnorm}, on every leaf
with a constant $K_{\rm all}\ge K$, and with the constant $K/2$ on every leaf with
$s_\xi(t_*)\le3S_\sharp$; $\sup_Q\mathfrak m_{Q^*}(t_*)\le\delta_U^2$; and
$\|f_\xi\|_X\le\delta_c/2$ on every leaf with $s_\xi(t_*)\ge S_b-\frac14$.
\end{itemize}
Put $v_*:=\sup_{t\ge T_0,\,\xi}\dot s_\xi(t)<\infty$. There is $\tau=\tau(M,K_{\rm
all},q_*)\in(0,\min\{1,1/(4v_*)\}]$, which depends also on the constants up to $T_0$ but not on
$t_*$, with the following properties.
\begin{enumerate}[label=(\alph*)]
\item There is a classical solution $g\ge0$ on $[t_*,t_*+\tau]$, unique among the classical
solutions with $|\partial^{\mathbf k}_{y,v}g|\le C\Gamma_c$ for $|\mathbf k|\le2$ and some $c>0$
(Lemma~\ref{lem:LWuniq}). On this interval $\sup_c\mathfrak B^{(1)}_c\le2M$, $\sup_c
\mathfrak B^{(1)}_{c,\ge1}\le2m_\sharp$, and
\begin{equation}
\sup_Q\mathfrak m_{Q^*}\le\frac{c_\sharp}{eT_0}\le\delta_U^2 ;
\end{equation}
(G1)--(G3) hold with $2K_{\rm all}$, and with $K$ on the leaves with $s_\xi(t_*)\le3S_\sharp$.
Moreover, $\|f_\xi\|_X\le\delta_c$ on the leaves with $s_\xi(t_*)\ge S_b-\frac14$. No bound of the
type (G4) is assumed or propagated.
\item There are $A,B$, depending on $(M,K_{\rm all},q_*)$ and on the constants chosen up to $T_0$,
but not on $t_*$, such that for both choices of $\Omega$, every cell $c$, $X\in\{
\mathfrak B^\Omega,\mathfrak B^\Omega_{\ge1}\}$ and $t\in[t_*,t_*+\tau]$,
\begin{equation}
X_c(t)\le e^{A(t-t_*)}X_c(t_*)+B(t-t_*)e^{(A+B)(t-t_*)}\sup_{c'}
\mathfrak B^\Omega_{c'}(t_*)+\sigma_X(t-t_*)e^{(A+B)(t-t_*)} ,
\end{equation}
where $\sigma_X:=(1+B)\sigma_F$ for $X=\mathfrak B^\Omega$ and $\sigma_X:=\sigma_F^{\ge1}+B\sigma_F$
for $X=\mathfrak B^\Omega_{\ge1}$, with the constant sources $\sigma_F$, $\sigma_F^{\ge1}$ of
Proposition~\ref{prop:UFforce}(a), (d), defined in the proof of Corollary~\ref{cor:CUtstar}. They
depend on $M$ and on the constants chosen up to $T_0$, but not on $t_*$. The last term cannot be
dropped: the reference forcing makes $\mathfrak A_c$ grow even from data with vanishing cell
energies, such as the local Maxwellians of the reference flow. For $\Omega=\Omega^\epsilon$ this is
proved in the same continuity argument as~(a), and it uses the undamped truncated thresholds of that
argument, as in the discussion of cubes in point~(2) of the proof.
\item For every $(\alpha,\beta)$ with $|\alpha|\le N_y$ and every compact $
\mathcal C\subset\R^3_y\times\R^3_v$, $\partial^\alpha_y\partial^\beta_vg\in C([t_*,t_*+\tau];H^{-1}(
\mathcal C))$, and it is bounded in $L^2(\mathcal C)$.
\item For every cell $c$,
\begin{equation}
\int_{t_*}^{t_*+\tau}\kappa\,\mathfrak D^{(1)}_c\dd t\le\frac{3M}{c_{\rm ret}},\qquad
\mathfrak D^{(1)}_c:=\sum_{(\alpha,\beta)}W^2_{\alpha\beta}\rad_0^{-3}\!\int\chi_c^2\,d_\xi[\varphi_{\alpha\beta}]\dd\xi
,
\end{equation}
where $\mathfrak D^{(1)}_c$ is the undamped cell dissipation, at the radii $\tilde\rad_k$, and
$c_{\rm ret}>0$ is the constant of the paragraph after Lemma~\ref{lem:LCtrunc}.
\end{enumerate}
\end{prop}

\begin{proof}
In the proof we use the viscous approximation of Appendix~\ref{app:visc}: Lemmas~\ref{lem:LWlin},
\ref{lem:LWregex}, \ref{lem:LCvisc}, \ref{lem:LCtrunc}, \ref{lem:LCgron} and~\ref{lem:LCd}, which
treat the Gaussian weight under the viscosity, the viscosity on the nonlocal part with the
$\xi$-dependent weights and frames, the nonlocal collision term on $\{\chi>0\}$, and the limit for
all derivative families. We take the cutoff $\chi$ to be $C^\infty$, with $|\chi'|+|\chi''|\le C$.
The argument consists of the following six points.

\emph{(1) The joint family, the radius and the constants.} Besides the transport and kernel
couplings $C_0t^{-2}\sqrt{E^c_nE^c_{n+1}}$, treated in the lemma on the viscous terms, the cell
hierarchy contains the top part of the forcing, Proposition~\ref{prop:UFforce}(c), the macroscopic
hierarchy
of part~(f) of that proposition, with $\upsilon=K_{\mathfrak M}^{1/2}$ as in Step~7 of its
proof, and the spread losses of the shifted Leibniz terms, with constants $C_1^F$ and $C_1^{\rm
mac}$.

As in Step~3 of Theorem~\ref{thm:CUboot}, we run the estimate for the joint family with
$K_{\mathfrak M}=T_0$. All that the forcing proposition needs from the solution is listed in
Remark~\ref{rem:UFinputs}. In the local theory they are supplied as follows:
\begin{itemize}
\item the cube bound on $\sup_Q\mathfrak a_{Q^*}$ by the first bullet of the list below;
\item (U1), $C\delta_U\le\rad_0$, requirement~\ref{rq:TWa} and the frame bounds (F3) by the
second bullet of that list;
\item $t\ge T_0\ge T_{\rm E}(L_1,R_1)$ and $\rad_0t\ge1$ by the third bullet of that list;
\item the comparability of Lemma~\ref{lem:CUcomp}(i)--(iii) for $\Omega^\epsilon$, by the
discussion of cubes in point~(2);
\item $\omega\le1$ is replaced by $\Omega^\epsilon\le1$, since $H_\epsilon\ge0$ (point~(2)), and
$D\omega\le0$ by $D\Omega^\epsilon=0$, since $\Omega^\epsilon$ is frozen at $t_*$ (point~(2));
\item the structural conditions $\rad_0\le r_{\rm C}/2$ and
requirements~\ref{rq:CUa}--\ref{rq:CUb} hold as in the bootstrap;
\item the graded radii $\tilde\rad_k$ lie in $[\rad_0/4,\rad_0]$, since
$\rad_0/2\le\rad_k(t_*)\le\rad_0$ by hypothesis and $\int_{t_*}^\infty r_{\rm
loc}\le\int_{T_0}^\infty|\rad'|\le\rad_0/4$ by Remark~\ref{rem:UFradius}.
\end{itemize}

\begin{itemize}
\item \emph{The cube bound.} By comparability on cubes, which gives $\chi_{\xi_i}^2\ge e^{-4}$ on
$Q_i$, $\sup_Q\mathfrak a_{Q^*}\le B_{\rm loc}:=27e^4\cdot2M\le54e^4M_\natural\le B_{\rm cr}$, since
$c_\chi C_U+\mu_\sharp\le e\,c_\chi C_U$. By the same remark and the paragraph after
Proposition~\ref{prop:UFforce}, every solution-dependent constant of that proposition is a
nondecreasing function of this bound, evaluated in the bootstrap at $B_{\rm cr}$. Hence in the local
theory these constants do not exceed their bootstrap values, and $r_{\rm loc}$, defined with the
bootstrap values, dominates at every $y$-order the local one-level losses, at the rate $\frac34C_0'$
of Lemma~\ref{lem:analytic}, and the spread losses of Lemma~\ref{lem:TWabs}. This is why
$M\le M_\natural$ is assumed.
\item \emph{The macroscopic field at the level $\delta_U$.} By the truncated thresholds of~(3) (and,
for the limit solution, by~(a)),
$K_{\mathfrak M}\mathfrak M^{(1)}_c\le \mathfrak B^{(1)}_{c,\ge1}\le2m_\sharp$, and thus, as for the
cube bound, we have
\begin{equation}\label{eq:CULW-mQ}
\sup_Q\mathfrak m_{Q^*}\le27e^4\sup_c\mathfrak M^{(1)}_c\le\frac{54e^4m_\sharp}{T_0}=\frac{c_\sharp}{eT_0}\le\delta_U^2
,
\end{equation}
as in the proof of Lemma~\ref{lem:CUH2}(i). Therefore (U1) holds on the interval, and with it the
smallness conditions of the forcing proposition on $\delta\varpi$: the bound on $C_\varpi$
before requirement~\ref{rq:CUa}, $|\nabla\delta\varpi|\le C\delta_U/\rad_0\le1$, the composition
condition of Lemma~\ref{lem:TWcalc}(vi) in requirement~\ref{rq:TWa}, the frame bounds (F3), and the
absorption of the $\delta\varpi$-parts of the coefficient terms. All hold with the bootstrap
constants. Inside the continuity argument of~(3) this bound holds as long as the truncated joint
threshold holds.
\item \emph{Times.} The forcing proposition needs $\rad_0t\ge1$, the reference flow from $T_0$ and
$K_{\mathfrak M}=T_0$. Hence $t_*\ge T_0$. The rate requirements~\ref{rq:UFd} and~\ref{rq:TWb} on
$T_0$ serve only the Gronwall factor over $[T_0,\infty)$, and on the short interval they are
replaced by $\tau$.
\end{itemize}

\emph{(2) The frozen weight.} We extend $H$ by $H(s)=0$ for $s<0$. Since $H(0)=0$, $H$ stays
nondecreasing and $C_1^U$-Lipschitz. Let $0\le\varrho\in C^\infty_c((0,1))$ with $\int\varrho=1$,
and put $H_\epsilon(s):=\int_0^1H(s-\epsilon r)\varrho(r)\dd r$ with $\epsilon:=1/(40C_1^U)$. Since
$H$ is nondecreasing and $C_1^U$-Lipschitz, we have
\begin{equation}
H-\tfrac1{40}\le H_\epsilon\le H,\qquad0\le H_\epsilon'\le C_1^U,\qquad H_\epsilon'=0\ \text{on}\
[S_b+1+\epsilon,\infty),\qquad|H_\epsilon''|\le C_1^U\|\varrho'\|_{L^1}/\epsilon .
\end{equation}
Hence we have $\omega_\xi(t)\le e^{-2H(s_\xi(t_*))}\le\Omega^\epsilon_\xi\le
e^{1/20}\omega_\xi(t_*)$ for $t\ge t_*$. The function $H$ itself is only Lipschitz and concave, and
$\partial^2_\xi e^{-2H(s_\xi)}$ has a nonnegative singular part on $\{s_\xi=S_b+1\}$, where
$\chi=1$. This would produce in Lemma~\ref{lem:LCvisc}(ii) a trace of the nonlocal density on that
hypersurface, which the energies do not control. For this reason $H$ is mollified.
\begin{itemize}
\item \emph{Gradient and Hessian.} Since $|\nabla s_\xi|\le2L's_\xi$, as in Lemma~\ref{lem:CUtr},
and $H_\epsilon'=0$ after $S_b+1+\epsilon$,
$|\nabla_\xi\log\Omega^\epsilon|\le4C_1^UL'\min\{s_\xi(t_*),S_b+1+\epsilon\}\le4C_1^UL'(S_\sharp+2)$,
uniformly in $\xi$ and $t_*$. With this bound, requirement~\ref{rq:CUb} implies $\Gamma\le3/\rad_0$
in that lemma. Using also $|\nabla^2s_\xi|\le CL'^2s_\xi$, we obtain
$|\nabla^2\log\Omega^\epsilon|\le C_\Omega:=C(C_1^U)^2L'^2(S_\sharp+2)^2+CC_1^UL'^2(S_\sharp+2)$.
\item \emph{Local part of the viscosity.} With
$\tilde\Omega:=\chi_c^2\Omega(A_p/A_{\rm ref})^2m_x^2$, Lemma~\ref{lem:LCvisc}(i) applies to
$\Omega=\Omega^\epsilon$, by the bound on $\nabla\log\Omega^\epsilon$ just stated, and we have
\begin{equation}
2\eta_v\!\iint\tilde\Omega\,h(\Delta_y+\Delta_v)h\le-\eta_v\!\iint\tilde\Omega\big(|\nabla_\xi
h|^2+|(\nabla_z+t^{-1}\nabla_\xi)h|^2\big)+\eta_v\!\iint\tilde\Omega\big(32q^2|z|^2+C\rad_0^{-2}\big)|h|^2
.
\end{equation}
\item \emph{Nonlocal part of the viscosity.} $\Psi_\Omega:=\rad_0^{-3}\chi_c^2\Omega\chi(s_\xi)$
is $C^2$, with $|\partial^2\Psi_\Omega|\le
C\rad_0^{-3}\chi_c^2\Omega(\rad_0^{-2}+C_\Omega+|\chi'|+|\chi''|)$. Thus
part~(ii) of the same lemma applies, with $\rad_0^2C_\Omega\le C$ by the same requirement, and
with $\rad_0C^{(2)}_\varpi\le C$, since $|\nabla^2_\xi p_g|\le
C(L_1,R_1)+C\delta_U\rad_0^{-2}$ by Lemma~\ref{lem:Mref}, Lemma~\ref{lem:TWcalc}(i) and (U1), and
$C\delta_U\le\rad_0$.
\item \emph{Cubes.} Part~(i) of the comparability lemma holds for $\Omega^\epsilon$: $H_\epsilon$ is
$C_1^U$-Lipschitz and constant after $S_b+1+\epsilon$, and in the mixed case
$(1+\epsilon_0)(S_b+1+\epsilon)\le2S_\sharp$. In Proposition~\ref{prop:CUtri} the sole role of the
damping is in this comparability. The supremum factor is undamped. Its bound
$\sup_Q \mathfrak a_{Q^*,\ge1}\le27e^4\cdot2m_\sharp=c_\sharp/e$ comes from the run with
$\Omega\equiv1$, and hence from the undamped truncated run of the same continuity argument,
Corollary~\ref{cor:CUtstar}, and for this reason the run with $\Omega^\epsilon$ is carried out
together with it.
\item \emph{Terms absorbed by the damping in the bootstrap.} Since $D\Omega=0$, the terms that the
damping absorbs in the bootstrap are Gronwall terms here: (L1)--(L2) of Proposition~\ref{prop:ACpf},
the remainder of Lemma~\ref{lem:CUbuf}, and the deficit $(\frac12-\chi)_+\|\varphi\|_{L^2}^2$ and
the front term of Proposition~\ref{prop:Bband} on the leaves with $\chi<1$. All are $\le C(K_{\rm
all},M)\kappa\lambda_\xi e_\xi$, and $\kappa\lambda_\xi\le C$ for $t\ge T_0$.
\end{itemize}

\emph{(3) The continuity argument for the viscous approximation.} We work with the \emph{graded}
truncation
\begin{equation}
\mathcal T_N:=\big\{(\alpha,\beta):|\alpha|\le N_y,\ |\alpha|+|\beta|\le N\big\},\qquad
U^c_N:=\sum_{(\alpha,\beta)\in\mathcal T_N}W^2_{\alpha\beta}\rad_0^{-3}\!\int\chi_c^2\Omega\,e_\xi[\varphi_{\alpha\beta}]\dd\xi+K_{
\mathfrak M}\sum_{n\le N}w_n^2\mathfrak M^c_n,
\end{equation}
and similarly $U^c_{N,\ge1}$ (the parts with $(\alpha,\beta)\ne0$). We use the graded truncation
$\mathcal T_N$ rather than a truncation in $|\beta|$ alone, since the forcing
of part~(d) of the forcing proposition, the coefficient terms and the viscous micro forcing raise
the $\xi$-level of $\delta\varpi$ and of the micro moments by up to $N_y+1$. With $\mathcal T_N$,
every term of the hierarchy of a component $(\alpha,\beta)\in\mathcal T_N$ involves only components
in $\mathcal T_N$, macroscopic levels $\le N$, and the boundary quantities
\begin{itemize}
\item[(a)] $\varphi_{\alpha,\beta+e_i}$ with $|\alpha|+|\beta|=N$ (transport and kernel part, at
each $y$-order);
\item[(b)] $\delta\varpi$ and the micro moments $\tilde M_k$ at the $\xi$-level $N+1$ (forcing,
($\Theta$1) and ($\Theta$5) in the proof of that proposition, and the viscous micro
forcing; their level is $|\beta|+|\alpha'|+1\le N+1$);
\item[(c)] $\mathfrak M^c_{N+1}$ and $E^c_{0,N+1}$, by part~(f) of the same proposition.
\end{itemize}
The coefficient terms reach only the $\xi$-level $|\beta|+|\alpha|\le N$, and Lemma~\ref{lem:nodown}
excludes couplings to higher $y$-orders. The boundary terms are absorbed at the truncation level
below. Since $\mathcal T_N\uparrow\{(\alpha,\beta):|\alpha|\le N_y\}$, monotone convergence as
$N\to\infty$ recovers the full sums. We work with the hypotheses
\begin{equation}
\sup_cU^c_N\le2M,\qquad\sup_cU^c_{N,\ge1}\le2m_\sharp
\end{equation}
of the continuity argument. Both are needed: the cube bound $B_{\rm loc}$ of~(1), which fixes the
solution-dependent constants of the forcing proposition and hence the domination by $r_{\rm loc}$,
uses the first, and the Gronwall constants themselves depend on it. Inside the continuity argument
the smallness of $\delta\varpi$ comes from the second, through the truncated version
of~\eqref{eq:CULW-mQ} in~(1): $\sup_Q\mathfrak m^{\le N}_{Q^*}\le c_\sharp/(eT_0)\le\delta_U^2$, and
thus $S^{\le N-2}_Q(\delta\varpi)\le C\delta_U$ by Remark~\ref{rem:TWtrunc}, where $N-2\ge N/2+1$
for $N\ge6$.

\emph{No mollification of the data; what the data provide.} The data $g_*:=g(t_*)$ are not
regularized. Their regularity is finite in $y$: the energies contain the $y$-orders
$|\alpha|\le N_y$ only, so Lemma~\ref{lem:BTsob}(ii), whose proof uses the $y$-order $|a|+2$, covers
$|a|\le N_y-2$ (and every order in $v$), and it carries the pointwise exponent $q/4$, which does not
control energies at $q_\delta$. The right inputs are the following. We use Lemma~\ref{lem:LCdata},
whose inputs are (G1), (F3) and the energies, but not (G4). Its hypothesis~(H) holds with
$q_*-\delta/8$ in place of $q_*$ and at a radius $\rad_D:=c_D\underline\rad$,
$c_D=c(L_1,R_1,\theta_\pm)\in(0,1]$, as we now check.

\emph{The micro part.} By Lemmas~\ref{lem:CUsob} and~\ref{lem:CUcomp}, for $\xi$ in a cube $Q$, we
have $\sum_n\tilde w_n^2\sum_{|\beta|=n,\alpha}\delta_y^{2|\alpha|}\eta\,\ell_\xi[h_{\alpha\beta}]\le
C'_S\,\mathfrak a_{Q^*}\le C'_S\cdot27e^5M$, with $\tilde w_n\ge\frac14\underline\rad^n/n!$. By the
definition of $\ell_\xi$, $|\nabla\log A_{\rm ref}|\le L$ by Corollary~\ref{cor:LWPa}(a4),
$\rho_g\le C\eps\rho_0$ and (F3), we deduce that
$\|e^{(q_*-\delta/8)|z|^2}h_{\alpha\beta}^{(x)}\|_{L^2_z}$ is bounded by
$C\eps\rho_0(\xi)|\beta|!\,\underline\rad^{-|\beta|}$.

\emph{The Maxwellian part.} For $\mathcal M[g_*]$, the factorial bounds on
$\partial^k_\xi\varpi_g(t_*,\cdot)$, from Lemma~\ref{lem:Mref} for $\varpi_*$ and from
$S_Q(\delta\varpi)\le C\delta_U$ for $\delta\varpi$, and the composition estimate of part~(vi) of
the cube calculus lemma for the entire map $\varpi\mapsto M_{p(\varpi)}(z)$ (as in
Lemma~\ref{lem:TWcoef}) yield the same bound at the radius $\rad_D$. The Gaussian
$e^{-|z-y_0|^2/(2\theta)}$, with $1/(2\theta)\ge1/(2\theta_++2)\ge4q_0>q_*$, since
$q_0\le q_{\rm GW}\le1/(8(\theta_++1))$ by Remark~\ref{rem:RGfull}(b), controls the factor
$e^{C\rad_D|z|}$ from complex $y_0$.

\emph{Conclusion.} Applying Lemma~\ref{lem:LCdata} with $(q_*-\delta/8,\delta/2,\rad_D)$ in place of
$(q_*,\delta,\underline\rad)$, we obtain
\begin{equation}
g_*\ge0,\qquad|\partial^\alpha_y\partial^\beta_vg_*(y,v)|\le
C_*\,|\beta|!\,(\rad_D/2)^{-|\beta|}\,\Gamma_{q_*-\delta/4}(y,v)\qquad(|\alpha|\le N_y-3,\
\beta\ \text{arbitrary}),\tag{D}
\end{equation}
with $C_*=C_*(M,K_{\rm all},\underline\rad,q_*,\delta,\delta_y,L_1,R_1,\theta_\pm)$, where
$\Gamma_c$ is the Gaussian barrier of the viscous approximation. Moreover, by the cube Sobolev
lemma, at the full energy exponent, we have
\begin{equation}
\sup_\xi\ell_\xi[h_{\alpha\beta}](t_*)\le C(N)\,M\qquad(|\alpha|\le N_y,\ |\beta|\le N+2).\tag{D$'$}
\end{equation}
The exponents form the chain $q_*\to q_*-\frac\delta4=c_0$ (data barrier of Lemma~\ref{lem:LWregex},
with $c_0=q_*-\frac\delta4$) $\to c_\eta(t)-\frac\delta8\ge q_\delta(t)+\frac\delta2$ (viscous
derivative bounds, for $\eta_v\le\delta/(64q_*^2)$), which leaves the margin $\delta/2$ over the
energy exponent $q_\delta$.

With these inputs, Lemma~\ref{lem:CUtstar} below shows that the truncated energies are continuous in
$t$ at $t_*$, and so is $\sup_c$ of them. The continuity argument itself is
Corollary~\ref{cor:CUtstar}.

We only need the upper bound $\limsup_{t\downarrow t_*}\sup_cU^c_N(t)\le\sup_cU^c_N(t_*)$.

The uniformity in $c$ has two sources. In $\xi$ it comes from the energies, not from the smoothness
of $g_*$. Indeed, at fixed $N$ every term of $U^c_N$ has level $n\le N$, and the levels $n+1$, $n+2$
are bounded in every cell, $E^c_{k,n+j}\le w_{n+j}(\rad_k)^{-2}\mathfrak A_c\le C(N)M$ for $j\le2$,
so the analytic tail plays no role. In $z$, at the $y$-orders $N_y-5,\dots,N_y$, the family
$\{h_{\alpha\beta}(\cdot,\xi)\}_\xi$ is only bounded in $L^2_z$ at the energy weight, and the
viscous flow need not be uniformly strongly continuous at $t_*$. The data are not regularized, since
a mollification would change the leaf parameters nonlinearly ($\delta\varpi=\Phi(m_g)-\varpi_*$). In
Lemma~\ref{lem:CUtstar}, mollified data enter only qualitatively, cell by cell.

We write the macroscopic hierarchy, part~(f) of the forcing proposition and (M1)--(M3) in its
proof, in $\varpi$ rather than in the conserved variables $m$, because the margins of
Lemma~\ref{lem:CUright} are in $\delta\varpi$, with leading coefficient $e^{A(t-t_*)}$, and the
conversion by
part~(vi) of the cube calculus lemma produces a factor $C^2>1$ in front of the data term.

\emph{The macroscopic part at the viscous level.} For
$\partial_tg=\kappa \mathcal Q_t[g]+\eta_v(\Delta_y+\Delta_v)g$, the identities
$(\Delta_yh)^{(x)}=\Delta_zh^{(x)}+2t^{-1}\nabla_z\cdot\nabla_\xi h^{(x)}+t^{-2}\Delta_\xi h^{(x)}$
and $(\Delta_vh)^{(x)}=\Delta_\xi h^{(x)}$ imply, for the leaf moments
$m=\int g^{(x)}\boldsymbol\ell\dd z$ with $\boldsymbol\ell=(1,z,|z|^2)$,
\begin{equation}
D_tm=t^{-2}\nabla_\xi\cdot\!\int g^{(x)}z\otimes\boldsymbol\ell\dd
z+\eta_v\Big[(1+t^{-2})\Delta_\xi m-2t^{-1}\nabla_\xi\cdot\!\int g^{(x)}\nabla_z\boldsymbol\ell\dd
z+\int g^{(x)}\Delta_z\boldsymbol\ell\dd z\Big],
\end{equation}
with $\nabla_z\boldsymbol\ell=(0,I,2z)$ and $\Delta_z\boldsymbol\ell=(0,0,6)$. The leaf parameters
are $\varpi=\Phi(m)$, with $\Phi$ holomorphic for $\rho,\theta>0$, and thus the viscous part of the
equation for $\varpi$ is
\begin{equation}
\eta_v(1+t^{-2})\Delta_\xi\varpi-\eta_v(1+t^{-2})\,D^2\Phi(m)\big[\nabla_\xi m,\nabla_\xi
m\big]+\eta_vD\Phi(m)\Big[-2t^{-1}\nabla_\xi\cdot(\rho_gI,2\rho_gy_{0,g})+(0,0,6\rho_g)\Big],
\end{equation}
with $\nabla_\xi m=D\Phi^{-1}(\varpi)\nabla_\xi\varpi$. With $\varpi=(\ell,y_0,\vartheta)$,
$\ell=\log\rho_g$, $y_0=y_{0,g}$, $\vartheta=\log\theta_g$, and $a:=1+t^{-2}$, this is the explicit
system
\begin{equation}\begin{aligned}
V_\ell&=a\eta_v\big(\Delta\ell+|\nabla\ell|^2\big),\qquad V_{y_0}=a\eta_v\big(\Delta
y_0+2(\nabla\ell\cdot\nabla)y_0\big)-2\eta_vt^{-1}\nabla\ell,\\
V_\vartheta&=a\eta_v\Big(\Delta\vartheta+|\nabla\vartheta|^2+2\nabla\ell\cdot\nabla\vartheta+\tfrac23e^{-\vartheta}|\nabla
y_0|^2\Big)-\tfrac43\eta_vt^{-1}e^{-\vartheta}\nabla\cdot y_0+2\eta_ve^{-\vartheta}.
\end{aligned}\end{equation}
Consequently $\ell$ enters only through $\nabla\ell$, the only non-polynomial coefficient is the
entire function $e^{-\vartheta}$, and no composition estimate, hence no smallness condition, is
needed for these terms. Every non-Laplacian term is of the form
$a\eta_v\,c(\varpi)\,\partial_i\varpi_\mu\,\partial_j\varpi_\nu$,
$\eta_vt^{-1}c(\varpi)\partial_i\varpi_\mu$ or $2\eta_ve^{-\vartheta}$, with
$c\in\{1,2,\frac23e^{-\vartheta},-2,-\frac43e^{-\vartheta}\}$. At a level $|\beta|=n\le N$,
integrated against $\Psi_c:=\rad_0^{-3}\chi_c^2\Omega$ and $w_n^2$, with $a_n$ and $g_n$ of
Lemma~\ref{lem:CUvisc} for $u=\delta\varpi$:
\begin{itemize}
\item one integration by parts and Young's inequality, as in~(2), yield the good term
$-\eta_v\rad_0^{-3}\int\chi_c^2\Omega w_n^2|\nabla_\xi\partial^\beta\delta\varpi|^2$ up to
$C\eta_v\rad_0^{-2}$ times the energy at level $n$;
\item \emph{the nonlinear terms are not bounded by $C\eta_v$ times the energy}: in the weights $w_n$
every Leibniz term of a product of two first derivatives loses exactly one level (already the terms
with one derivative on one factor carry the factor $\binom n1=n$), and a literal application of
part~(iii) of the cube calculus lemma to $\nabla\varpi\otimes\nabla\varpi$ would yield a bound cubic
in the good terms. Instead, we apply the viscous product lemma to each full product, with
$F=c(\varpi)\partial_i\varpi_\mu$ and $G=\partial_j\varpi_\nu$, and to the first-order terms and the
heating with an $F$ of order zero. This includes the parts with $n+1$ derivatives on one factor, so
no split is needed. We put
\begin{equation}
\Sigma:=1+\sup_QS^{\le N/2+1}_Q\big(\varpi-\varpi(\xi_Q)\big)+\sup_QS^{\le
N/2+1}_Q\big(e^{-\vartheta}\big)\le C_\Sigma=C(L_1,R_1,\theta_\pm),
\end{equation}
by the truncated Lemma~\ref{lem:TWcalc}(i) (Remark after it: the supremum factor carries at most
$N/2+1$ derivatives, so at most $N/2+3\le N$ levels are used, for $N\ge6$), comparability on cubes
and the Cauchy estimates for $\varpi_*$. For the bound by a structural constant we use the
continuity hypothesis: the composition with $e^{-\vartheta}$ is exponential,
$S_Q(e^{-\vartheta})\le e^{-\vartheta(\xi_Q)}\exp\big(C_{\rm
alg}S_Q(\vartheta-\vartheta(\xi_Q))\big)$, and $S^{\le N-2}_Q(\delta\varpi)\le C\delta_U\le1$ by the
truncated threshold, while $S_Q(\vartheta_*-\vartheta_*(\xi_Q))\le C\rad_0/r_{\rm C}$. Hence
$N_\Sigma$ below is structural. We apply the viscous product lemma only to products that contain a
factor of $\delta\varpi$; the pure reference products form $V(\varpi_*)$ and produce the source term
below. The cell factors are bounded by
$\hat E^{\le N}_c(\nabla\varpi)\le\sqrt3\big((\mathfrak M^{c,\le
N})^{1/2}+w_{N+1}\|A_{N+1}(\delta\varpi)\|_c+C_{\rm ref}\big)$ and, by parts~(i)--(ii) of the same
lemma, $w_{N+1}\|A_{N+1}(\delta\varpi)\|_c\le\frac{8\rad}{N+1}g_N$. By Young's inequality, for
$N\ge N_\Sigma:=CC_\Sigma$, a structural constant, we obtain
\begin{equation}
\frac{d}{dt}\mathfrak M^{c,\le N}\Big|_{\rm visc}\le-\frac38\,\eta_v\sum_{n\le
N}g_n^2+C\eta_v\rad_0^{-2}\Sigma^2\,\mathfrak M^{c,\le N}+C\eta_v\rad_0^{-2}C_{\rm
ref}^2\sup_{c'}\mathfrak M^{c',\le N}+C_{\rm ref}\,\eta_v\big(\mathfrak M^{c,\le N}\big)^{1/2};
\end{equation}
the fractions are $\frac14$ for the integration by parts and $\frac14+\frac18$ for the viscous
product lemma. The last term collects the sources of the reference part $\varpi_*$, which only
involves $\nabla\ell_*$, $y_{0*}$, $\vartheta_*$ and their derivatives, bounded on
$\mathbb S_{R_1/4}$ by the bounds of the reference flow and the Cauchy estimates, and
$C_{\rm ref}=C(L_1,R_1,\theta_\pm)$. The third term comes from the mixed products
$\partial\delta\varpi\cdot\partial\varpi_*$, whose factor
$\hat S(\partial\delta\varpi)\le C\sup_{c'}(\mathfrak M^{c',\le N})^{1/2}$ is a supremum over cells,
and it is of the form $B\,S$ of Lemma~\ref{lem:LCgron}. Of the good term
$\frac38K_{\mathfrak M}\eta_v\sum g_n^2$ of the joint family, the micro forcing below uses
$\frac1{16}+\frac1{16}$, the viscous micro forcing at all $y$-orders $\ge1$ together uses
$\frac1{16}$, and $\frac3{16}K_{\mathfrak M}\eta_vg_N^2$ remains for the truncation level. We obtain
the part used by the viscous micro forcing by Young's inequality with the fraction
$\frac1{16}(\#\{\alpha:1\le|\alpha|\le N_y\})^{-1}$ for each $\alpha$, after summing over $\beta$
and over $\alpha'\le\alpha$ at each level by Cauchy--Schwarz against the components
$w_{m-1}\|\nabla\partial^{\beta+\alpha'}\delta\varpi\|_c$ of $g_{m-1}$. This part also covers the
Leibniz parts of the viscous micro forcing at $y$-orders $\ge1$, at the boundary level $N$ and at
the interior levels alike. All constants are independent of $N$, $\eta_v$ and $\delta$;
\item \emph{the viscous forcing of the micro equations.} The micro part $\hat r=g-\mathcal M[g]$
satisfies the micro equation with the additional forcing
\begin{equation}\label{eq:Fvisc}
\mathfrak F_{\rm visc}:=\eta_v(\Delta_y+\Delta_v)\mathcal M[g]-DM(\varpi)[V(\varpi)] .
\end{equation}
The terms $DM[a\eta_v\Delta\varpi]$ cancel, and thus
\begin{equation}
\mathfrak F_{\rm visc}^{(x)}=\eta_v\Big(aD^2M[\nabla_\xi\varpi,\nabla_\xi\varpi]+2t^{-1}\nabla_z\cdot
DM[\nabla_\xi\varpi]+\Delta_zM\Big)-DM\big[V(\varpi)-a\eta_v\Delta\varpi\big],
\end{equation}
which contains no second derivative of $\varpi$ and is moment-free on every leaf. Its top part at
the micro level $(0,\beta)$, $|\beta|=n$, is $\eta_v\,c(\varpi)(z)\,\nabla\partial^\beta\varpi$, a
crossing into the macroscopic level $n+1$ with coefficient $\eta_v$. We bound it by
\begin{equation}\label{eq:visc-top}
C\eta_vw_n(E^c_{0,n})^{1/2}g_n\le\frac1{16}K_{\mathfrak M}\eta_vg_n^2+C\eta_vK_{
\mathfrak M}^{-1}w_n^2E^c_{0,n},
\end{equation}
where we use the macroscopic good term, which carries $K_{\mathfrak M}\ge1$. At $n=N$ it uses $g_N$,
which is present in the truncated family, so it needs no condition on $N$.

We treat the Leibniz parts by the viscous product lemma with the micro test factor: by its
part~(ii), $(n+1)w_n\|\hat r_{0n}\|\le64\rad w_{n-1}(L^c_n)^{1/2}$, which is absorbed by a quarter
of the micro good term $\eta_v\eta w_{n-1}^2L^c_n$ of~(2), up to
$C\eta_v\eta^{-1}\rad^{-2}\Sigma'^2$ times the energy and, through the level $N+1$ of $\varpi$,
$C\eta^{-1}\Sigma'^2(N+1)^{-2}\eta_vg_N^2\le\frac1{16}K_{\mathfrak M}\eta_vg_N^2$ for
$N\ge C\Sigma'\eta^{-1/2}$.

Here $\Sigma':=\Sigma+C_{\rm Max}$: the Leibniz parts also carry the $S$-norms of the $
\mathcal H_j$-valued coefficients $DM_p$, $D^2M_p$, $\partial_zDM_p$, which are $\le C_{\rm Max}$ by
Lemma~\ref{lem:TWcoef} and Remark~\ref{rem:TWtrunc}. The coefficients are $\mathcal H_j$-valued and
holomorphic in the leaf parameters, as in that lemma.

The $y$-derivatives $\partial^\alpha_y=(\partial_z+t^{-1}\partial_\xi)^\alpha$ fall on $\varpi$ with
the factor $t^{-|\alpha'|}$. Inside $\mathcal T_N$ they reach at most the $\xi$-level $N+1$, and we
treat the boundary level $N+1$ at the truncation level below, item~(b).

The parts of $\mathfrak F_{\rm visc}$ that contain only the reference field (the top part with
$\nabla\partial^\beta\varpi_*$ in place of $\nabla\partial^\beta\delta\varpi$, and the Leibniz
products of derivatives of $\varpi_*$ alone) are not bounded through $g_n$, which is defined with
$u=\delta\varpi$. They form the micro reference source. Paired with the micro components, they are
bounded by $C'_{\rm ref}\eta_v(U^c_N)^{1/2}$, with $C'_{\rm ref}=C(L_1,R_1,\theta_\pm,N_y)$ from the
Cauchy estimates of Lemma~\ref{lem:UFref} (at $y$-orders $\ge1$ with the gap~\eqref{eq:UFgap}).

Its zeroth-order part cancels: $\Delta_zM_p=DM_p[(0,0,2e^{-\vartheta})]$, which is the heating term
$2\eta_ve^{-\vartheta}$ of $V_\vartheta$, so the heating produces no micro source. This source
enters $\sigma_{\rm ref}$ of point~(4)(b).
\end{itemize}
\emph{The truncation level.} At $n=N$ the joint hierarchy contains the crossings
$C_0^Ft^{-2}w_N^2\sqrt{E^c_{0,N}}\sqrt{\mathfrak M^c_{N+1}}$, from part~(c) of the forcing
proposition, and
$K_{\mathfrak M}(C_0+C_0^F)t^{-2}w_N^2\sqrt{\mathfrak M^c_N}(\sqrt{
\mathfrak M^c_{N+1}}+C\sqrt{E^c_{0,N+1}})$, from part~(f). In the joint variables
$E^{(0)}_n=E^c_{0,n}+K_{\mathfrak M}\mathfrak M^c_n$:
\begin{itemize}
\item the macro--macro crossing is homogeneous in $K_{\mathfrak M}$;
\item $CK_{\mathfrak M}(C_0+C_0^F)w_N^2\sqrt{
\mathfrak M_N}\sqrt{E_{0,N+1}}=C(C_0+C_0^F)K_{\mathfrak M}^{1/2}w_N^2\sqrt{K_{\mathfrak M}
\mathfrak M_N}\sqrt{E_{0,N+1}}$ with $C(C_0+C_0^F)K_{\mathfrak M}^{1/2}\le C_0'$ for the values of $C_0$,
$C_0^F$ before their enlargement by the factor $C+2$ in Step~7 of the proof of
that proposition (after it, only the bound on the sum of the crossings in the next
sentence is used);
\item $C_0^Fw_N^2\sqrt{E_{0,N}}\sqrt{\mathfrak M_{N+1}}=C_0^FK_{
\mathfrak M}^{-1/2}w_N^2\sqrt{E_{0,N}}\sqrt{K_{\mathfrak M}\mathfrak M_{N+1}}$.
\end{itemize}
By part~(g) of the same proposition, established in Step~7 of its proof, these crossings, together
with the transport and kernel crossings at $y$-order~$0$, add up to at most
$C_0't^{-2}w_N^2\sqrt{E^{(0)}_N}\sqrt{E^{(0)}_{N+1}}$.

We split this sum once (not each crossing separately: the one-level part $\frac34C_0'$ of
$r_{\rm loc}$ has no margin), by $w_N\le\frac{N+1}\rad w_{N+1}$, as
$\frac{N+1}{2\rad}C_0't^{-2}(\tilde a_N^2+\tilde b_{N+1}^2)$ with $\tilde a_N^2=w_N^2E^{(0)}_N$ and
$\tilde b_{N+1}^2=w_{N+1}^2(E_{0,N+1}+K_{\mathfrak M}\mathfrak M_{N+1})$. The radius absorbs
$\tilde a_N^2$, as at every level $n<N$.

For $\tilde b_{N+1}^2$, by Lemma~\ref{lem:LCtrunc}(i), the micro viscosity supplies the good term
$-\eta_v\eta w_N^2L^c_{N+1}$. Moreover $E_{0,N+1}\le(\eta+C)L^c_{N+1}$, since the energy contains
the local part with the coefficient $\eta$ and the nonlocal part $\le C\ell_\xi$. The remaining
macroscopic good term is $\frac3{16}K_{
\mathfrak M}\eta_vg_N^2\ge\frac3{16}(\frac23)^{10}K_{
\mathfrak M}\eta_v\big((N+1)/\rad\big)^2w_{N+1}^2\mathfrak M^c_{N+1}$, by parts~(i)--(ii) of
the viscous product lemma
(the combinatorial factor of (i) goes in the favorable direction).

Therefore the joint good term is $\ge c\,\eta_v\eta(\eta+C)^{-1}((N+1)/\rad)^2\tilde b_{N+1}^2$, and
half of it dominates the crossings for $N+1\ge N'_\eta:=C\rad_0C_0'/(\eta\,\eta_vt_*^2)$, as in the
lemma on the viscous terms (note $C_0^F\le C_0'$). In this way we treat the boundary quantities (a)
at $y$-order $0$ and (c).

\emph{The boundary at $y$-orders $k\ge1$.} (a) is the same split, done for each $y$-order $k$ with
$\rad_k$ in place of $\rad$ and with the micro good term
$\eta_v\eta\,w_{N-k}(\rad_k)^2L^c_{k,N-k+1}$. It needs $N-k+1\ge N'_\eta$, which holds for
$N\ge N'_\eta+N_y$. For (b), we take $|\alpha|=k\ge1$, $|\beta|=N-k$ and $|\alpha'|=k$. By the
gap~\eqref{eq:UFgap} with $j=N-k$, $m=k+1$, and $\rad\ge\rad_0/4$, we have
\begin{equation}
W_{\alpha\beta}\,t^{-k}\le\delta_y^k\,(N+1)^{k+1}(1-\epsilon_y)^{N-k}\rad^{-k-1}t^{-k}\,w_{N+1}(\rad)\le
C_k\,\rad_0^{-1}(\rad_0t)^{-k}\,w_{N+1}(\rad),
\end{equation}
since $(N+1)^{k+1}(1-\epsilon_y)^{N-k}\le C(k,\epsilon_y)$. Thus the crossing is $\le
C_kC_0^F\rad_0^{-1}t^{-2}\tilde a_{\alpha\beta}\tilde b_{N+1}\le Ct^{-2}\rad_0^{-1}(\tilde
a_{\alpha\beta}^2+\tilde b_{N+1}^2)$, with $\tilde
a_{\alpha\beta}:=W_{\alpha\beta}e_{\alpha\beta}^{1/2}$.

The first part is a Gronwall term with the integrable rate $Ct^{-2}\rad_0^{-1}$, independent of
$N$. The second is absorbed by the joint good term $\ge
c\,\eta_v\eta(\eta+C)^{-1}((N+1)/\rad)^2\tilde b_{N+1}^2$ as soon as $N+1\ge
N''_\eta:=\big(C\rad_0/(\eta\,\eta_vt_*^2)\big)^{1/2}$. This does not follow from $N+1\ge
N'_\eta$ in general, since $C_0'\rad_0/(\eta\,\eta_vt_*^2)$ need not be $\ge1$.

Unlike at $y$-order $0$, this crossing carries no factor $(N+1)/\rad$, because the gap supplies
$(1-\epsilon_y)^{N-k}$. The moments $\tilde M_k$ at level $N+1$ are bounded by
$\|\varphi_{0,N+1}\|_X$ by Lemma~\ref{lem:TWmom}(i), hence by the micro good term at $(0,N)$.

The viscous micro forcing is the same with the extra factor $\eta_v$, and it is absorbed by the
part of $\frac1{16}K_{\mathfrak M}\eta_vg_N^2$ reserved for it. The interior viscous crossings at
$y$-orders $\ge1$, at the macroscopic levels $m=|\beta|+|\alpha'|+1\le N$, are
$\le C\eta_v\rad_0^{-1}\tilde a_{\alpha\beta}g_{m-1}$ and are absorbed by the same part of
$g_{m-1}$, leaving Gronwall terms $C\eta_v\rad_0^{-2}K_{\mathfrak M}^{-1}\tilde a^2_{\alpha\beta}$.

Half of the joint good term is used for the crossings at $y$-order $0$ (this is the half in the
definition of $N'_\eta$) and the other half for those of item~(b) at $y$-orders $\ge1$ (this is the
constant in $N''_\eta$).

With
\begin{equation}\label{eq:Nstar}
N\ge N_*:=\max\{N_\eta,N'_\eta,N''_\eta,N_\Sigma,C\Sigma'\eta^{-1/2},6\}+2N_y,
\end{equation}
the count of the truncation remark holds at every $y$-order, since $(N-N_y)/2+3\le N-N_y$. The
crossing of the viscous micro forcing was treated above without condition on $N$. The spread losses
stay within $\mathcal T_N$.

The continuity argument at fixed $N\ge N_*$ then runs with thresholds and $\tau$ independent of $N$,
$\eta_v$ and $\delta$. We let $N\to\infty$ by monotone convergence, and then $\eta_v,\delta\to0$ by
the weak lower semicontinuity of Lemma~\ref{lem:LCd}. It applies to the macroscopic part as well.
The viscous family converges in $C^k(\mathcal C)$ for $k\le N_y-4$ and every compact $\mathcal C$,
and $\partial^\alpha_y\partial^\beta_vg_{\eta_v}$ converges weakly in $L^2(\mathcal C)$ for
$|\alpha|\le N_y$ and all $\beta$ (the derivative families are bounded in $L^2(\mathcal C)$ and
converge in $H^{-1}(\mathcal C)$), and the moments, which are integrals over all $z$, converge
because of the tail bound, uniform in $\eta_v$ and in $\xi\in Q$,
\begin{equation}
\int_{|z|>R}\big|(\partial^\beta_vg_{\eta_v})^{(x)}\big|\jap z^2\dd z\le\Big(\sup_{\xi\in
Q}\ell_\xi[h_{0\beta}]\Big)^{1/2}\Big(\int_{|z|>R}m_x^{-2}\jap z^4\dd
z\Big)^{1/2}\longrightarrow0\qquad(R\to\infty),
\end{equation}
where the supremum is bounded by the cube Sobolev lemma and the truncated cell bounds. Moreover,
$\rho_g$ is bounded below on compact sets, since $|\delta\ell|\le C\delta_U$. Hence
$\partial^\beta\varpi_{\eta_v}\to\partial^\beta\varpi$ locally uniformly, each $\mathfrak M^c_n$ is
lower semicontinuous by Fatou's lemma in $\xi$, and $\mathfrak M_c$ by Fatou's lemma over $n$. We
use no continuity of the full sums, and no margin in the radius.

\emph{(4) The conclusions.} (a) is Corollary~\ref{cor:CUtstar}, passed to the limit as in~(3), with
(1)--(3) and (6). The coercivity is that of the paragraph after Lemma~\ref{lem:LCtrunc}: on leaves
with $s_\xi<S_b$ on the interval, which have $s_\xi(t_*)\le3S_\sharp$, it uses (G1)--(G3) with $K$,
and (5) below in place of (G4). On leaves with $s_\xi\ge S_b$ it is that of the closeness, given by
Proposition~\ref{prop:Bband} and Remark~\ref{rem:Bbandrange}, which uses $\|f_\xi\|_X\le\delta_c$
and not the constant of (G3). Since $\tau\le1/(4v_*)$, a leaf with $s_\xi(t_*)<S_b-\frac14$ does not
reach $S_b$ on the interval. Thus the Gronwall constants, hence $\tau$, $A$ and $B$, depend on
$K_{\rm all}$, but not on $t_*$: the rates are bounded by $C(K_{\rm all},M)\kappa\lambda_\xi\le C$,
by $Ct^{-2}\le CT_0^{-2}$, by $C_{\rm tr}/(\rad_0^2C_2t^2)$ and by $C\eta_v\rad_0^{-2}$. (b) is
Lemma~\ref{lem:LCgron}(i), together with~(ii), applied to the truncated sums of the viscous
approximation and passed to the limit as in~(3). For the viscous approximation, (b) holds with
$\sigma_F$ and $\sigma_F^{\ge1}$ replaced by $\sigma_F+\sigma_{\rm ref}$ and
$\sigma_F^{\ge1}+\sigma_{\rm ref}$, where
\begin{equation}
\sigma_{\rm ref}:=\big(C_{\rm ref}T_0^{1/2}+C'_{\rm ref}\big)\eta_v(2M)^{1/2}
\end{equation}
collects the viscous reference sources of~(3): the macroscopic one, \begin{equation}C_{\rm
ref}\eta_vK_{\mathfrak M}^{1/2}(K_{\mathfrak M}\mathfrak M^{c,\le N})^{1/2}\le C_{\rm
ref}\eta_vT_0^{1/2}(U^c_N)^{1/2},\end{equation}
and the micro one, $\le C'_{\rm ref}\eta_v(U^c_N)^{1/2}$. Under the thresholds we have
$(U^c_N)^{1/2}\le(2M)^{1/2}$. This term vanishes as $\eta_v\to0$. To keep the thresholds inside the
continuity argument we add $\sigma_{\rm ref}\tau\le m_\sharp/10$ to the definition of $\eta_0$.

(c) is Lemma~\ref{lem:LCd}(ii), whose hypotheses hold by the paragraphs after it.

(d): by the same corollary with $\Omega\equiv1$, we have
\begin{equation}\label{eq:CULW-diss}
c_{\rm ret}\int_{t_*}^{t_*+\tau}\kappa\mathfrak D^c_N\dd t\le3M
\end{equation}
for the truncated viscous dissipation, uniformly in $N$, $\eta_v$ and $\delta$. Then we let
$N\to\infty$, by monotone convergence.

For the limit $(\eta_v,\delta)\to0$, we fix $(\alpha,\beta)$. The quantity
$\int_{t_*}^{t_*+\tau}\kappa\rad_0^{-3}\int\chi_c^2d_\xi[\varphi_{\alpha\beta}]\dd\xi\dd t$ is the
square of a weighted $L^2$ norm, over $(t_*,t_*+\tau)\times\R^3_\xi\times\R^3_z$, of
$h_{\alpha\beta}^{(x)}$ and of $\nabla_zh_{\alpha\beta}^{(x)}$, namely
$\lambda_\xi\|\varphi\|^2_{Y_2(\hat m_\xi)}$ with $\varphi=T^{-1}_{p_g}h^{(x)}_{\alpha\beta}$,
written in the variable $z$. Its weights are continuous functions of $z$, of the leaf parameters
$p_g(t,\xi)$ and of the exponent $q_\delta$, and they increase as $\delta\downarrow0$. Along the
sequence, $h_{\alpha\beta}[g_\eta]$ converges weakly in $L^2$ on the compact subsets of
$(t_*,t_*+\tau)\times\R^6$ (by the compactness lemma at every $t$, with the uniform bounds, and the
locally uniform convergence of the Maxwellian part). Its $z$-derivatives are bounded in $L^2$ on
these sets, by~\eqref{eq:CULW-diss}, so they converge weakly as well. The leaf parameters converge
locally uniformly.

By part~(i) of the same lemma on compact sets, in the variables $(t,y,v)$, by an exhaustion, and by
Fatou's lemma in the sum over $(\alpha,\beta)$, we can pass to the limit in the bound, as for
$\int\kappa\mathsf D^{(k)}_c$ in the proof of Proposition~\ref{prop:disp}.

By Lemma~\ref{lem:LWregex}(i)--(ii), the viscous solution is unique at fixed $\eta_v$. The limit
$\eta_v,\delta\to0$ is a single limit (the data are not regularized), its initial value is $g(t_*)$
by the bound in $C([t_*,t_*+\tau];H^{-1}(\mathcal C))$ of~(c), and Lemma~\ref{lem:LWuniq} identifies
the limit with the solution.

The hypotheses of the uniqueness lemma hold. Indeed, Lemma~\ref{lem:LCdata}, applied at every
$t\in[t_*,t_*+\tau]$ with $t$ in place of $t_*$, and with~(H) obtained from the cell bounds of~(a)
as in point~(3), now with the radius $\rad_0/4$ in place of $\underline\rad$ and the exponent
$q(t)\ge q_*-C_2\tau/T_0^2$, so that the rate is $q(t)-\delta/4\ge q_*/2$ since $\delta\le q_*/4$,
implies that $|\partial^{\mathbf k}_{y,v}g(t)|\le C\Gamma_{q_*/2}$ for $|\mathbf k|\le2$, uniformly
on the interval. The same holds, with constants depending on the bounds, for every classical
solution in the cell norms in the sense of the proof of Corollary~\ref{cor:CUmain}, on every compact
subinterval, in particular for the solution on $[T_0,T_{\max})$ in Lemma~\ref{lem:CUright} and in
the proof of that corollary. We apply the uniqueness lemma with the smaller of the two rates.

\emph{Definition of $\tau$.} $\tau$ is the largest number in $(0,\min\{1,1/(4v_*)\}]$ with
$\tau\le K/(2L_3)$ and $\tau\le3\delta_c^2/(4L_{\rm cl})$ (point~(6)), $e^{2C\tau}\le\frac32$,
$2C\tau M\le m_\sharp/6$, $\sigma_F\tau\le M/10$ and $\sigma_F^{\ge1}\tau\le m_\sharp/10$, where
$C=C(K_{\rm all},M)$ is the Gronwall constant and $\sigma_F$, $\sigma_F^{\ge1}$ are the constant
sources of the continuity argument. We compute these sources with $T_0^{-2}\ge t^{-2}$ in place of
$t^{-2}$, so $\tau$ does not depend on $t_*$, $N$, $\eta_v$ or $\delta$.

\emph{(5) Avoiding (G4).} In the local theory no smallness is required of the Gronwall constants, so
we can replace every use of (G4) by a bound in terms of the energies. We put
\begin{equation}
G_\xi(t):=\sum_{|a|+|b|\le2}\ell_\xi[\hat r_{ab}](t)^{1/2}\le C\rad_0^{-2}\delta_y^{-2}B_{\rm
loc}^{1/2}=:G_*(M),
\end{equation}
by Lemma~\ref{lem:CUsob} and the static Lemma~\ref{lem:BTsob}(i) with $B=B_{\rm loc}$, applied with
$\rad_0/4$ in place of $\rad_{\rm post}$: the local radii satisfy only $\tilde\rad_k\ge\rad_0/4$
(point~(1)), and by Remark~\ref{rem:UFgrade}(d) this changes only absolute factors. The same reading
of $\rad_{\rm post}$ applies to the uses of its part~(ii) below and in point~(6). Since
$\nabla_{\hat z}=\theta^{1/2}\nabla_z$ and $\partial_z=\partial_y-t^{-1}\partial_v$ (by (F1) and
(F3)), $\sum_{j\le2}\|\hat m_\xi\nabla^j_{\hat z}f_\xi\|_{L^2}\le C(\theta_\pm)G_\xi$. By the static
part~(ii) with $a=b=0$, we have
$\|e^{(q/4)|\hat z|^2}f_\xi\|_{L^\infty}\le C\rad_{\rm post}^{-2}\delta_y^{-2}B_{\rm loc}^{1/2}\le
CG_*(M)$. Polynomial weights are dominated by the leaf weight:
$\|\jap{\hat z}^\theta u\|_{L^2}+\|u\|_{L^1}\le C(q_-,\theta)\|u\|_{X(\hat m_\xi)}$. (G4) enters the
local theory in three places.
\begin{itemize}
\item[($\alpha$)] \emph{The term $Q(\varphi,f_\xi)$ of the crude linear term}
Proposition~\ref{prop:ACpf}. In non-divergence form,
$Q(\varphi,f)=(a*\varphi):\nabla^2f-(c*\varphi)f$ with $c=\nabla\cdot b$, so
$c*\varphi=-2(\gamma+3)|\cdot|^\gamma*\varphi$ for $\gamma>-3$ and $c*\varphi=-8\pi\varphi$ for
$\gamma=-3$. We write $m=\hat m_\xi$. For the first term,
$|a*\varphi|\le C\jap{\hat z}^{\gamma+2}\|\varphi\|_{L^2(\jap{\hat z}^\theta)}\le
C\|\varphi\|_{X(m)}$, since $\gamma+2\le0$, and thus
\begin{equation}
\Big|\int m^2\varphi\,(a*\varphi):\nabla^2f_\xi\Big|\le\|a*\varphi\|_{L^\infty}\|m\varphi\|_{L^2}\|m\nabla^2f_\xi\|_{L^2}\le
C\,G_\xi\,\|\varphi\|^2_{X(m)} .
\end{equation}
For the second term with $\gamma=-3$, we have
$|\int m^2\varphi^2f_\xi|\le\|f_\xi\|_{L^\infty}\|m\varphi\|^2_{L^2}$. For $\gamma>-3$, we split
$|\cdot|^\gamma=k_0+k_\infty$ with $k_0=|\cdot|^\gamma\mathbf 1_{\{|\cdot|\le1\}}\in L^1$ and
$|k_\infty|\le1$. Since $m\in\mathcal W(q_\pm,\eps_1)$, for $|\hat z-\hat w|\le1$ we have
$m(\hat z)\le m(\hat w)e^{C(1+|\hat z|)}$, with $C=C(q_+,\eps_1)$. Hence
\begin{equation}
\Big|\int m^2\varphi f_\xi\,(k_0*\varphi)\Big|\le\big\|e^{C(1+|\hat
z|)}f_\xi\big\|_{L^\infty}\int|m\varphi|\,\big(k_0*|m\varphi|\big)\le C\big\|e^{(q/4)|\hat
z|^2}f_\xi\big\|_{L^\infty}\|k_0\|_{L^1}\|m\varphi\|^2_{L^2},
\end{equation}
and $|\int m^2\varphi f_\xi(k_\infty*\varphi)|\le\|\varphi\|_{L^1}\|m\varphi\|_{L^2}\|mf_\xi\|_{L^2}$.
Altogether, we obtain
\begin{equation}
\big|\langle Q(\varphi,f_\xi),\varphi\rangle_{L^2(\hat m_\xi^2)}\big|\le C\Big(\|\hat
m_\xi\nabla^2f_\xi\|_{L^2}+\|\hat m_\xi f_\xi\|_{L^2}+\big\|e^{(q/4)|\hat
z|^2}f_\xi\big\|_{L^\infty}\Big)\|\varphi\|^2_{X(\hat m_\xi)}\le C\,G_*(M)\,\|\varphi\|^2_{X(\hat
m_\xi)} .
\end{equation}
This is a Gronwall term with rate $C\kappa\lambda_\xi G_*(M)$, without a $Y_2$ factor, so the coercivity
$c_h(K)$ and the threshold $m_\sharp$ are unchanged.
\item[($\beta$)] \emph{Lemma~\ref{lem:CUbuf}}, where (G4) bounds $\|f_\xi\|_{Y_2}$. The radial part
of $\nabla\hat m_\xi$ is $2q_\xi\hat z\hat m_\xi$, up to $O(\eps_1)\hat m_\xi$, with
$q_\xi=q(t)\theta_g(t,\xi)$, and its tangential part is $O(\eps_1)\hat m_\xi$, so, with
$\gamma\le-2$, $\|f_\xi\|_{Y_2}\le C(q_\pm)\sum_{j\le1}\|\hat m_\xi\nabla^jf_\xi\|_{L^2}\le
C\,G_*(M)$. The remainder $C_{\rm buf}(1+G_*(M)^2)\kappa\lambda_\xi\|\varphi\|_X^2$ is a Gronwall
term, and no damping is needed to control it on a bounded interval.
\item[($\gamma$)] \emph{The level-$(0,0)$ background dissipation} in the inequality for the parts
$\ge1$, in the paragraph after Lemma~\ref{lem:LCtrunc}:
$\kappa\lambda_\xi\|f_\xi\|^2_{Y_2}\le C\kappa\lambda_\xi G_*(M)^2$, by~($\beta$). On close leaves,
by the divergence-form bound of Lemma~\ref{lem:QmG}, we have
$|\langle Q(f_\xi,\varphi),\varphi\rangle|\le C(\|f_\xi\|_{L^2(\jap{\hat
z}^4)}\|\varphi\|_{Y_2}+\|f_\xi\|_{H^1(\jap{\hat
z}^3)}\|\varphi\|_X)\|\varphi\|_{Y_2}\le(C\delta_c+\frac
c8)\|\varphi\|^2_{Y_2}+Cc^{-1}G_*(M)^2\|\varphi\|^2_X$, which is used in the same way.
\end{itemize}
The constant $G_*(M)$ enters only $A$, $B$ and $\tau$.

\emph{(6) (G3) on the interval.} By the leaf equation for $\hat g^{(x)}=\rho_\xi^{-1}g^{(x)}$ along
rays,
\begin{equation}
D_t\hat g^{(x)}=\kappa\rho_\xi\,Q\big(\hat g^{(x)},\hat g^{(x)}\big)+t^{-2}z\cdot\nabla_\xi\hat
g^{(x)}+t^{-2}\big(z\cdot\nabla_\xi\log\rho_\xi\big)\hat g^{(x)},
\end{equation}
we have $\jap z^8|D_t\hat g^{(x)}|\le\kappa\rho_\xi\jap z^8\big(|\ab[\hat g^{(x)}]||\nabla^2_z\hat
g^{(x)}|+|\cb[\hat g^{(x)}]||\hat g^{(x)}|\big)+t^{-2}\jap z^9\big(|\nabla_\xi\hat g^{(x)}|+L|\hat
g^{(x)}|\big)$. We write $\hat g^{(x)}=\rho_\xi^{-1}M_{p_g}+\rho_\xi^{-1}\hat r^{(x)}$. The
Maxwellian part is explicit, with $|\partial^k_\xi p_g|\le C_k$ by (F3), which holds on the interval
by~(1). For the micro part, by the static Lemma~\ref{lem:BTsob}(ii) with $|a|+|b|\le2$ (its proof
uses levels $\le4$) and $B=B_{\rm loc}$, we have
$e^{(q/4)|\hat z|^2}|T^{-1}_p(\partial^a_y\partial^b_v\hat r)^{(x)}|\le C\rad_{\rm
post}^{-4}\delta_y^{-4}B_{\rm loc}^{1/2}$, and
$e^{-(q/4)|\hat z|^2}\jap z^{10}\le C(q,L,\theta_\pm)$. With $\kappa\rho_\xi\le C$ and the standard
bounds on $\ab,\cb$, we obtain
\begin{equation}
\big|N^{(3)}_\xi(t)-N^{(3)}_\xi(t_*)\big|\le L_3(t-t_*),\qquad
L_3=C(q,L,\theta_\pm)\big(1+\rad_{\rm post}^{-4}\delta_y^{-4}B_{\rm loc}^{1/2}\big)^2,
\end{equation}
uniformly in $\xi$. For the viscous approximant $g_\eta$, $D_t\hat g^{(x)}$ contains in addition
$\eta_v\rho_\xi^{-1}((\Delta_y+\Delta_v)g_\eta)^{(x)}$. By the same static bound with $|a|+|b|\le4$
(levels $\le6\le N_y$), its $\jap z^8$-weighted supremum is
$\le C\eta_v(1+\rad_{\rm post}^{-6}\delta_y^{-6}B_{\rm loc}^{1/2})$, so $L_3$ is replaced by
$L_3+C(1+\rad_{\rm post}^{-6}\delta_y^{-6}B_{\rm loc}^{1/2})$, uniformly in $\eta_v\le1$, with
$B_{\rm loc}$ given by the truncated thresholds. With $\tau\le K/(2L_3)$, we obtain (G3) with
$K/2+L_3\tau\le K$ on the leaves with $s_\xi(t_*)\le3S_\sharp$ and with
$K_{\rm all}+K/2\le2K_{\rm all}$ on the others, on the interval, for $g_\eta$ inside the continuity
argument and for the limit.

\emph{Closeness on the interval.} From the leaf equation for $f_\xi$, Lemma~\ref{lem:leafeq}(ii)
with $\alpha=\beta=0$, we deduce that
\begin{equation}
\Big|\frac{d}{dt}\|f_\xi\|^2_{X(\hat m_\xi)}\Big|\le L_{\rm
cl}:=C(q_\pm,\theta_\pm,L_1,R_1)\big(1+F^{(0)}+\sup|
\mathfrak e|\big)\big(1+G_*(M)\big)^3\qquad\text{on }[t_*,T'],
\end{equation}
uniformly in $\xi$ and in $\eta_v\le1$. The ingredients are:
\begin{itemize}
\item the hypodissipativity $\langle L_\mu f,f\rangle_X\le C\|f\|^2_{Y_2}$;
\item $|\langle Q(f,f),f\rangle_X|\le C\|f\|_X\|f\|^2_{Y_2}$;
\item $\|f\|_{Y_2}\le CG_*$ (point~(5)($\beta$));
\item the transport, $t^{-2}\|\jap{\hat z}\hat mT_p^{-1}(\partial_v\hat r)^{(x)}\|\,\|f\|_X\le
Ct^{-2}G_*^2$, since $\|\hat z\hat mu\|\le(2q_\xi)^{-1}\|\hat m\nabla u\|+C\|\hat mu\|$, with
$q_\xi=q(t)\theta_g(t,\xi)\ge\frac12q(t)\theta_-$;
\item the frame motion, $|\langle\mathcal G_p[Dp]f,f\rangle_X|\le Ct^{-2}|\mathfrak e|G_*^2$;
\item the forcing $t^{-2}F^{(0)}\|f\|_X$;
\item the kernel part of the level-$0$ forcing, $\le Ct^{-2}\|\varphi_{0e}\|_X$ by
Proposition~\ref{prop:UFforce}(e), which applies in the local theory through
Remark~\ref{rem:UFinputs};
\item the time derivative of the weight $\hat m_\xi$, whose $q'(t)$-part has a good sign;
\item for $g_\eta$ only, the term $C\eta_v(G_*^2+1)$.
\end{itemize}
Here $\sup|\mathfrak e|\le C$ by (F3) and the micro moments, and $F^{(0)}$ is the level-$0$ forcing
constant of the same part~(e).

With $\tau\le3\delta_c^2/(4L_{\rm cl})$, we have
$\|f_\xi\|^2_X\le\frac{\delta_c^2}4+L_{\rm cl}\tau\le\delta_c^2$ on the leaves with
$s_\xi(t_*)\ge S_b-\frac14$. The conditions $\tau\le K/(2L_3)$ and
$\tau\le3\delta_c^2/(4L_{\rm cl})$ are part of the definition of $\tau$.

A pointwise bound is available here, unlike for (G4), because the weight $\jap z^8$ of (G3) is
polynomial and is dominated by $e^{(q/4)|\hat z|^2}$, while the weight of (G4) is the full Gaussian
$\hat m_\xi$, and the energies have no margin in the exponent.
\end{proof}

\begin{lem}[Continuity at $t_*$ of the truncated viscous energies]\label{lem:CUtstar}
In the setting of point~(3) of the proof of Proposition~\ref{prop:CULW}, let $\delta\in(0,q_*/4]$,
$\eta_v\le\eta_0$, and fix $N\ge N_*$, with
the graded truncation $\mathcal T_N$. Let $g_\eta$ be the viscous solution with the data
$g_*=g(t_*)$, and let $X^c\in\{U^c_N,U^c_{N,\ge1}\}$. There are $h\in(0,\min\{\tau,\tau_\eta\}]$,
with $\tau_\eta$ the existence time of Lemma~\ref{lem:LWregex}(i), and $A_{\rm cr},B'_{\rm
cr}<\infty$, depending on $(\eta_v,N,\delta,C_*,M,K_{\rm all},t_*)$ but not on $c$ or $t$, such
that:
\begin{enumerate}[label=(\roman*)]
\item \emph{(low orders up to $t_*$)} for $|\alpha|\le N_y-4$, every $\beta$ and $t\in[t_*,t_*+h]$,
\begin{equation}\label{eq:CUtstar-pt}
|\partial^\alpha_y\partial^\beta_vg_\eta(t)|\le C_1(\alpha,\beta)\Gamma_{c_\eta(t)-\delta/8},
\end{equation}
and the same holds for $\partial_t\partial^\alpha_y\partial^\beta_vg_\eta$ if $|\alpha|\le N_y-6$;
\item \emph{(crude inequality, cell by cell)} on $(t_*,t_*+h]$, classically, $\frac{d}{dt}X^c\le
A_{\rm cr}X^c+B'_{\rm cr}$;
\item \emph{(upper bound from $t_*$, cell by cell)} for every $c$ and $t\in[t_*,t_*+h]$,
\begin{equation}
X^c(t)\le e^{A_{\rm cr}(t-t_*)}X^c(t_*)+B'_{\rm cr}(t-t_*)e^{A_{\rm cr}(t-t_*)} ;
\end{equation}
\item each $X^c$ is continuous on $[t_*,T_\eta)$, where $[t_*,T_\eta)$ is the interval of existence
of $g_\eta$, and so is $S_X(t):=\sup_cX^c(t)$.
\end{enumerate}
\end{lem}

\begin{proof}
\emph{(i) Schauder up to the initial time.} Let $Z=B\times[t_*,t_*+\tau_\eta]$, with $B$ a unit ball
in $\R^6$, where $g_\eta$ exists by Lemma~\ref{lem:LWregex}(i). Since $\tau_\eta\le1$, the constants
of the local form are those for intervals of length~$1$.

For $\partial_tu=a^{ij}\partial_{ij}u+cu+f$ on $2Z$, with $a\ge\eta_vI$, the coefficients and $f$
H\"older continuous in the space variables, uniformly in $t$, and $u(t_*)=u_*$, the estimate up to
the initial time with coefficients H\"older in space \cite[Thm.~4.1]{Lor00}, in the local form of
Appendix~\ref{app:visc}, bounds $D^2u$ in $C^\alpha_x$ and $\partial_tu$ in $L^\infty$ on $Z$ by
$\sup_{2Z}|u|$, $\|u_*\|_{C^{2+\alpha}(2B)}$ and the spatial H\"older norm of $f$ on $2Z$, with a
constant depending on $\eta_v$, $\alpha$ and the spatial H\"older norms of the coefficients. There
are no compatibility conditions, since only the bottom of $2Z$ carries data.

Here $a$ is the matrix, in the variables $(y,v)\in\R^6$, of the second-order part
$\ab[g_\eta]:D^2+\eta_v(\Delta_y+\Delta_v)$, with $D=\partial_v-t\nabla_y$, and thus $a\ge\eta_vI$.

The coefficients $\partial^\kappa\ab[g_\eta]=\ab[\partial^\kappa g_\eta]$ and
$\partial^\kappa\cb[g_\eta]$ are H\"older continuous in $(y,v)$, uniformly in $t$ and in the center
of $Z$, as soon as $\partial^\kappa g_\eta$ has a Gaussian bound and a spatial H\"older bound
(Lemma~\ref{lem:LWlin}(i) and the induction in the proof of Lemma~\ref{lem:LWregex}(i)). This is
that induction, with every step taken up to $t_*$ by \cite{Lor00}, using the data
$\partial^\kappa g_*\in C^{2+\alpha}(2B)$, which (D) provides for $|\mathbf k|+3$ derivatives within
the stated $y$-orders. We use no H\"older continuity in $t$, uniform in $(y,v)$. At $\gamma=-2$ it
fails, by Lemma~\ref{lem:LWlin}(i).

For the weights, we have $\sup_{2Z}\Gamma_c\le C_\delta\Gamma_{c-\delta/8}$ at the center of $Z$.
Moreover, $\|g_*\|_{C^{k+2+\alpha}(2B)}\le
C_\delta\mathcal N_{k+3}(g_*;c)\Gamma_{c-\delta/8}(\text{center})$. Thus the bounds are Gaussian and
the constants do not depend on the center, hence not on $c$. By (D), we have $k+3\le N_y-3$ in $y$.
Therefore $|\alpha|\le N_y-4$, and $|\alpha|\le N_y-6$ for $\partial_t$, since the equation loses
two derivatives.

In $v$ the induction is anisotropic: we fix $\alpha$ and raise $\beta$. The commutator of
$\partial_v$ with $\ab[g]:D^2$ is $\ab[\partial_vg]:D^2$, a forcing at the same $y$-order, and thus
all $\beta$ are reached.

Consequently the frames $p_g(t,\xi)$, $\rho_g/\rho_\xi$ and $\partial^k_\xi p_g$ for every $k$, are
Lipschitz on $[t_*,t_*+h]$ uniformly in $\xi$ (they are $z$-integrals of $\partial^\beta_vg_\eta$
with uniform Gaussian tails), and for $h$ small (F3) holds with $\frac14,4$ in place of $\frac12,2$,
as do the closeness on $\{s_\xi(t_*)\ge S_b-\frac14\}$ and (G1)--(G3).

\emph{(ii)} We redo the computation of Lemma~\ref{lem:LCtrunc}, of the paragraph after it and of
point~(3), but we bound every factor that is not a component of $X^c$ by~(i) instead of by the
continuity hypothesis: the frames and leaf masses; the coefficients $\ab$, $\cb$ and the Maxwellian
coefficients; the macroscopic fields $\partial^\beta\delta\varpi$, $|\beta|\le N+1$, in supremum;
and the low factor of each Leibniz split of a trilinear term, chosen as the factor with $y$-order
$\le N_y/2\le N_y-6$ and $v$-order $\le N$. What remains is a quadratic form in the components of
$X^c$, in the levels $N+1$ of the top crossings, and in one $z$-derivative, with coefficients
bounded by $C(N,C_1)$. The levels $N+1$ are absorbed by the viscous good terms as at the truncation
level in point~(3) of that proof, and the $z$-derivatives by $-\eta_v\iint\tilde\Omega|L_1h|^2$, by
$-c_h\kappa\lambda\|\cdot\|^2_{Y_2}$ (Young's inequality, with constants depending on $N$) and by
the shrinking exponent for the weight loss $32\eta_vq^2|z|^2$. No $\sup_{c'}$ remains, since all
quantities that are nonlocal in $\xi$ are of low order. For $U^c_{N,\ge1}$ the coupling to level
$(0,0)$ is a source bounded by~(i).

\emph{(iii)} Let $\phi_\epsilon$ be a standard mollifier in $(y,v)$ and
$g^\epsilon_*:=\phi_\epsilon*g_*\ge0$. Since $\phi_\epsilon*\Gamma_c\le C\Gamma_{c-\delta/16}$ for
$\epsilon\le\epsilon(\delta)$, (D) holds for $g^\epsilon_*$ at $q_*-\frac5{16}\delta$, uniformly in
$\epsilon$, and $g^\epsilon_*$ has Gaussian bounds on all derivatives, with constants depending on
$\epsilon$. With $c_0'=q_*-\frac5{16}\delta$ the margin is $c_\eta-\frac\delta8\ge
q_\delta+\frac\delta2$ for $\eta_v\le\delta/(64q_*^2)$.

(a) The viscous solutions $g^\epsilon_\eta$ exist on $[t_*,t_*+h]$ with $h$ independent of
$\epsilon$, by Lemma~\ref{lem:LWregex}(i), and $X^c[g^\epsilon_\eta]\in C^1([t_*,t_*+h])$, since for
smooth data~(i) holds at all orders. Part~(ii) holds for $g^\epsilon_\eta$ with the same
$A_{\rm cr},B'_{\rm cr}$, which depend only on the low-order bounds of~(i), uniform in $\epsilon$.
Integrating from $t_*$, we obtain
$X^c[g^\epsilon_\eta](t)\le e^{A_{\rm cr}(t-t_*)}X^c[g^\epsilon_*]+B'_{\rm cr}(t-t_*)e^{A_{\rm
cr}(t-t_*)}$.

(b) \emph{Stability.} $d:=g^\epsilon_\eta-g_\eta$ solves
$\partial_td=\ab[g^\epsilon_\eta]:D^2d-\cb[g^\epsilon_\eta]d+\eta_v\Delta
d+\ab[d]:D^2g_\eta-\cb[d]g_\eta$. By~(i), $|D^2g_\eta|\le C\Gamma_{c_\eta-\delta/8}$ up to $t_*$. By
the maximum principle with the barrier $\Gamma_{c'}e^{K(t-t_*)}$ at a lower rate $c'$, as in
part~(ii) of the same lemma, we obtain $\sup_{[t_*,t_*+h]}\sup|d/\Gamma_{c'}|\le C\epsilon$.

By the interior estimate \cite{Kn80}, through the induction of part~(iii) of that lemma, uniform in
$\epsilon$ on $[t_*+s,t_*+h]$, we have $\partial^\kappa g^\epsilon_\eta\to\partial^\kappa g_\eta$
locally uniformly for $t>t_*$, with a uniform Gaussian majorant. Thus
$X^c[g^\epsilon_\eta](t)\to X^c[g_\eta](t)$ for every $t\in(t_*,t_*+h]$ and every $c$.

(c) \emph{The data, cell by cell.} We fix $c$. For $|\alpha|\le N_y$, $|\beta|\le N+1$,
$\partial^\alpha_y\partial^\beta_vg^\epsilon_*=\phi_\epsilon*\partial^\alpha_y\partial^\beta_vg_*$
converges in $L^2(W_c)$, $W_c:=\chi_c(v+y/t_*)^2A_{\rm ref}(v)^2e^{2(q_*-\delta)|y|^2}$, since
$W_c(y,v)\le C_\delta W'_c(y',v')$ for $|(y-y',v-v')|\le1$, because
$A_{\rm ref}(v)\le e^LA_{\rm ref}(v')$ by Corollary~\ref{cor:LWPa}(a4), with $W'_c$ the same weight
at the exponent $q_*$, so that $\phi_\epsilon*$ is bounded from $L^2(W'_c)$ to $L^2(W_c)$ uniformly
in $\epsilon$, and the data lie in $L^2(W'_c)$ by $\mathfrak A_c(t_*)\le M$.

The moments and their $\xi$-derivatives converge locally uniformly with uniform bounds, by (D). Thus
the frames, $\mathcal M[g^\epsilon_*]$ and $\delta\varpi[g^\epsilon_*]$ converge, and the nonlocal
part converges as in Lemma~\ref{lem:LCd}(iii). Hence we get $X^c[g^\epsilon_*]\to X^c[g_*]$.

In this limit we need no uniformity in $c$: $c$ is fixed, and in particular the non-convexity of the
energies in $g$ plays no role. Letting $\epsilon\to0$ in~(a) at fixed $c$ and $t>t_*$, we
obtain~(iii).

\emph{(iv)} By~(iii), we have $\limsup_{t\downarrow t_*}X^c(t)\le X^c(t_*)$. By~(i),
$g_\eta(t)\to g_*$ locally uniformly with derivatives of order $\le N_y-4$. The components of
$X^c(t)$ are bounded in $L^2(W_c)$ by~(iii), so they converge weakly, and weak lower semicontinuity
with converging frames, as in Lemma~\ref{lem:LCd}(iii), implies that
$\liminf_{t\downarrow t_*}X^c(t)\ge X^c(t_*)$. For $t\ge t_*+s>t_*$, Lemma~\ref{lem:LWregex}(iii)
with the margin $\delta/2$ makes every $X^c$ Lipschitz on $[t_*+s,T]$, uniformly in $c$, for every
$T<T_\eta$. Hence $S_X$ is continuous on $(t_*,T_\eta)$. At $t_*$ it is lower semicontinuous, as a
supremum of continuous functions, and
$S_X(t)\le e^{A_{\rm cr}(t-t_*)}S_X(t_*)+B'_{\rm cr}(t-t_*)e^{A_{\rm cr}(t-t_*)}\to S_X(t_*)$
by~(iii).
\end{proof}

\begin{cor}[The continuity argument of point~(3)]\label{cor:CUtstar}
Let $T'$ be the supremum of the $t\in[t_*,t_*+\tau]$ such that $\sup_cU^c_N\le2M$ and
$\sup_cU^c_{N,\ge1}\le2m_\sharp$ on $[t_*,t]$. Then $T'=t_*+\tau$, and on $[t_*,t_*+\tau]$ the
bounds of Proposition~\ref{prop:CULW}(a), (b) hold for the truncated viscous energies, with the data
$g(t_*)$ unchanged. Moreover $c_{\rm ret}\int_{t_*}^{t_*+\tau}\kappa\,\mathfrak D^c_N\dd t\le3M$ for
every cell, for both choices of $\Omega$.
\end{cor}

\begin{proof}
At $t_*$ the values are $\le M$ and $\le m_\sharp$, and by Lemma~\ref{lem:CUtstar}(iv) the
hypotheses hold on $[t_*,t_*+h']$ for some $h'>0$, which implies that $T'>t_*$. On $(t_*,T')$ the
sharp inequalities of Lemma~\ref{lem:LCtrunc}, of the paragraph after it and of point~(3) hold, in
the form
\begin{equation}
\frac{d}{dt}U^c_N+c_{\rm ret}\,\kappa\,\mathfrak D^c_N\le
CU^c_N+C\sup_{c'}U^{c'}_N+\sigma,\qquad\frac{d}{dt}U^c_{N,\ge1}\le
CU^c_{N,\ge1}+C\sup_{c'}U^{c'}_N+\sigma^{\ge1},
\end{equation}
with the dissipation $c_{\rm ret}\kappa\mathfrak D^c_N$ that is left, as in the paragraph after that
lemma, with $C=C(K_{\rm all},M)$, $\sigma:=\sigma_F+\sigma_{\rm ref}$ and
$\sigma^{\ge1}:=\sigma_F^{\ge1}+\sigma_{\rm ref}$. Here $C$, $\sigma_F$ and $\sigma_F^{\ge1}$ are
independent of $N$, $\eta_v$, $\delta$ and $t_*$. Moreover, $\sigma_{\rm ref}$ is the viscous
reference source of point~(4)(b), which is $O(\eta_v)$, with $\sigma_{\rm ref}\tau\le m_\sharp/10$,
and
\begin{equation}\begin{aligned}
\sigma_F&:=CT_0^{-2}\big(F_*+\rad_0^{-1}(C_{\rm mom}+P_1^{1/2})P_1^{1/2}\big)(2M)^{1/2},\\
\sigma_F^{\ge1}&:=CT_0^{-2}\Big(F_*\big(\tfrac{\rad_0}{r_{\rm
C}}+\delta_y\big)+\rad_0^{-1}(C_{\rm mom}+P_1^{1/2})P_1^{1/2}\Big)(2m_\sharp)^{1/2}
\end{aligned}\end{equation}
are the constant sources of Proposition~\ref{prop:UFforce}, which applies in the local theory
without change by Remark~\ref{rem:UFinputs}: the reference forcing~(a), together with the first term
of $\mathfrak b_c$ for the parts $\ge1$, and the constant source of the shifted Leibniz terms
in~(d). There we have $\mathfrak A_c^{1/2}\le(U^c_N)^{1/2}\le(2M)^{1/2}$ and
$\mathfrak B_{c,\ge1}^{1/2}\le(U^c_{N,\ge1})^{1/2}\le(2m_\sharp)^{1/2}$ under the hypotheses of the
continuity argument, and $t^{-2}\le T_0^{-2}$. The second term of $\mathfrak b_c$ is
$\le C(U^c_{N,\ge1}+U^c_N)$. These sources do not vanish in general: for data equal to the local
Maxwellians of the reference flow all cell energies vanish at $t_*$, but the micro part of
$\partial_t\hat r$ at level $(0,0)$ does not vanish where $\nabla\varpi_*\ne0$. We use
Lemma~\ref{lem:LCgron}, whose hypotheses hold by~(iii) and~(iv) of the preceding lemma. By its
part~(ii), with $A=B=C$ and the source $\sigma$, we have $\sup_cU^c_N\le e^{2C\tau}(M+\sigma\tau)$.
By its part~(i), with $S'=\sup_{c'}U^{c'}_N\le2M$ and the source $\sigma^{\ge1}$, we have
$\sup_cU^c_{N,\ge1}\le e^{C\tau}(m_\sharp+2C\tau M+\sigma^{\ge1}\tau)$. With $e^{2C\tau}\le\frac32$,
$2C\tau M\le m_\sharp/6$, $\sigma_F\tau\le M/10$ and $\sigma_F^{\ge1}\tau\le m_\sharp/10$
(definition of $\tau$, point~(4)), $\sigma_{\rm ref}\tau\le m_\sharp/10$ (definition of $\eta_0$,
point~(4)(b)) and $m_\sharp\le M$ (hypothesis (H1$'$)), these are
\begin{equation}
\le\tfrac32\big(1+\tfrac1{10}+\tfrac1{10}\big)M<2M,\qquad\le\big(\tfrac32\big)^{1/2}\big(1+\tfrac16+\tfrac1{10}+\tfrac1{10}\big)m_\sharp<2m_\sharp
.
\end{equation}
The same lemma implies the cell-by-cell bound of~(b). With $d:=t-t_*\le\tau\le1$, by (i) and (ii) we
have, for $X=U_N$, $X^c(t)\le e^{Cd}X^c(t_*)+Cde^{2Cd}S(t_*)+\sigma de^{Cd}+C\sigma d^2e^{2Cd}$, and
for $X=U_{N,\ge1}$ the same with $\sigma^{\ge1}de^{Cd}$ in place of $\sigma de^{Cd}$, which is the
form stated in~(b), with $\sigma_{\rm ref}$ added to $\sigma_F$ and to $\sigma_F^{\ge1}$. In the
limit $\eta_v\to0$ it disappears. \emph{The dissipation.} We integrate the first inequality over
$[s,t_*+\tau]$, $s>t_*$, drop $U^c_N(t_*+\tau)\ge0$, and let $s\downarrow t_*$, using~(iii) of the
preceding lemma. Under the thresholds, we obtain
\begin{equation}
c_{\rm ret}\int_{t_*}^{t_*+\tau}\kappa\,\mathfrak D^c_N\dd t\le
U^c_N(t_*)+\tau\big(4CM+\sigma\big)\le3M ,
\end{equation}
since $U^c_N(t_*)\le M$, $4C\tau\le2\log\frac32\le1$ and $\sigma\tau\le M/5$. \emph{Existence of
$g_\eta$.} Lemma~\ref{lem:LWregex}(i) provides $g_\eta$ on $[t_*,t_*+\tau_\eta]$ only, with
$\tau_\eta=\tau_\eta(C_*,c_0,t_*)$. On $[t_*,\min\{T',T_\eta\})$ the truncated thresholds and
Lemma~\ref{lem:LCdata} at each time (levels $\le2$ in $(y,v)$ suffice for $\mathcal N_0$) bound
$\ab[g_\eta]$, $\cb[g_\eta]$ and the constant of the null structure independently of the barrier
constant. Since $g_\eta$ is the bounded solution of the linear problem with $h=g_\eta$,
Lemma~\ref{lem:LWlin}(iii) implies $g_\eta\le C_*e^{K_{\rm b}(t-t_*)}\Gamma_{c_\eta(t)}$, with
$K_{\rm b}=K_{\rm b}(M,t_*)$, up to $\min\{T',T_\eta\}$, and Lemma~\ref{lem:LWregex}(iv) yields
$T_\eta>T'$. By the continuity of $S_X$ in~(iv) of the preceding lemma, the bounds extend beyond
$T'$ unless $T'=t_*+\tau$. Hence $T'=t_*+\tau$ and $T_\eta>t_*+\tau$. Here $\tau\le1$, as the linear
estimates and the bound $c_\eta(t)\ge c_0-4\eta_vc_0^2(t-t_*)$ are stated on intervals of length
$\le1$.
\end{proof}

\begin{lem}[The bootstrap gives the hypotheses of the local theory]\label{lem:CUH2}
Assume the requirements of Theorem~\ref{thm:CUboot}.
\begin{enumerate}[label=(\roman*)]
\item If (U2) holds at a time $t$, then (U1) holds at $t$.
\item If (U1)--(U4) hold on $[T_0,T]$, then (H1$'$)--(H2$'$) of Proposition~\ref{prop:CULW} hold at
every $t_*\in[T_0,T]$, with
\begin{equation}
M:=e^{2H(S_b+1)}\big(c_\chi C_U+\mu_\sharp\big),\qquad
m_\sharp:=e^{2H(S_b+1)}\mu_\sharp=\frac{c_\sharp}{54e^5},\qquad K_{\rm all}:=\max\{K,K_{\rm rel}\}.
\end{equation}
These constants do not depend on $T$ or $t_*$. More precisely, (H1$'$)--(H2$'$) hold at a time
$t_*\in[T_0,T]$ as soon as (U1)--(U4) hold on $[T_0,t_*]$ and (U1)--(U4) hold at $t_*$ with the
halved constants that the theorem supplies. The bootstrap radii
satisfy the radius hypothesis of the proposition: they are graded as in
Remark~\ref{rem:UFgrade}, $\rad(t_*)=\rad_0-\int_{T_0}^{t_*}|\rad'|\le\rad_0$, and
$\rad_k(t_*)\ge\rad_0/2$ by \eqref{eq:radrate} and requirement~\ref{rq:TWb}, as in
Remark~\ref{rem:UFradius}.
\end{enumerate}
\end{lem}

\begin{proof}
(i) By (U2), we have $K_{\mathfrak M}\mathfrak M_c\le\mathfrak B_{c,\ge1}\le2\mu_\sharp$, with
$K_{\mathfrak M}=T_0$. By Lemma~\ref{lem:CUcomp}, as in Remark~\ref{rem:UFradius}, we obtain
\begin{equation}
\sup_Q\mathfrak m_{Q^*}\le27e^5e^{2H(S_b+1)}\sup_c
\mathfrak M_c\le27e^5e^{2H(S_b+1)}\frac{2\mu_\sharp}{T_0}=\frac{c_\sharp}{T_0}\le\frac{c_\chi
C_U}{T_0}\le\delta_U^2 ,
\end{equation}
since $c_\sharp\le c_\chi$, $C_U\ge1$, and $\delta_U\ge(c_\chi C_U/T_0)^{1/2}$ because the constant
in the definition of $\delta_U$ is $\ge1$. This is (U1), with a margin.

(ii) By Theorem~\ref{thm:CUboot}, applied on $[T_0,t_*]$, (U1)--(U4) hold at $t_*$ with halved
constants. To avoid circularity, we use that theorem below only through these halved bounds at
$t_*$; Lemmas~\ref{lem:BTG12} and~\ref{lem:BTsob}, and Step~1 of the proof of the theorem at $t_*$,
with Step~2 and Corollary~\ref{cor:Aleafentry} on $[T_0,t_*]$, use only the hypotheses (U1)--(U4) on
$[T_0,t_*]$. In particular, we have $\mathfrak A_c\le c_\chi C_U$ and
$\mathfrak B_{c,\ge1}\le\mu_\sharp$, and thus
$\mathfrak B_c=\mathfrak A_c+K_{\mathfrak M}\mathfrak M_c\le\mathfrak A_c+\mathfrak B_{c,\ge1}\le
c_\chi C_U+\mu_\sharp$. Since $H$ is constant on $[S_b+1,\infty)$, $\omega_\xi\ge e^{-2H(S_b+1)}$
for every $\xi$ and $t$, and thus the undamped energies are bounded by $e^{2H(S_b+1)}$ times the
damped ones. This proves (H1$'$) with $M=M_\natural$ and $m_\sharp$ as stated, and
$C_N(54e^5m_\sharp)^{1/2}=C_Nc_\sharp^{1/2}\le c_{\min}/8$ by the choice of $c_\sharp$.

For (H2$'$):
\begin{itemize}
\item \emph{(G1)--(G3).} (G1) holds on every leaf by Lemma~\ref{lem:BTG12}, from (U1). On
$\{s_\xi\le3S_\sharp\}$, (G3) with $K/2$ is the halved (U4) of the theorem, and (G2) follows by the
same lemma. On the other leaves, $s_\xi>3S_\sharp\ge S_b+1$, so $B_\xi=\delta_b^2$ in the leaf
Sobolev lemma, and part~(iv) of that lemma implies (G3) with $K_{\rm rel}$, and then (G2) follows as
before. (G4) is not part of (H2$'$).
\item \emph{$\delta\varpi$.} By the halved (U1) of the theorem, $\sup_Q
\mathfrak m_{Q^*}\le\delta_U^2$.
\item \emph{Closeness.} A leaf with $s_\xi\ge S_b-\frac14$ has $s_\xi\ge S_*-\frac14$, since
$S_b\ge S_*$. Step~1 of the proof of the theorem yields $\|f_\xi\|_X\le\delta_c/2$ on such leaves:
by the entry corollary if $s_\xi\le3S_\sharp$, and by the second case of Step~1 and
Lemma~\ref{lem:CUsob} otherwise, with the requirements $\delta_1/C_e\le\delta_c/2$ and
$C\delta_b\le\delta_c/2$.
\end{itemize}

\end{proof}

\begin{lem}[Lower semicontinuity]\label{lem:CUsemi}
Let $T_n\in[T_0,T_{\max})$ increase to $T<T_{\max}$, and assume that (U1)--(U4) hold on $[T_0,T_n]$
for every $n$. Then the functions $\Phi_2,\Phi_3$ of the proof of Corollary~\ref{cor:CUmain} satisfy
$\Phi_i(T)\le\liminf_n\Phi_i(T_n)$.
\end{lem}

\begin{proof}
We choose $n_0$ with $T-T_{n_0}<\tau_\flat/2$, with $\tau_\flat$ of the proof of
Corollary~\ref{cor:CUmain}. Note that this time is a lower bound for the time $\tau$ of
Proposition~\ref{prop:CULW} with the constants of Lemma~\ref{lem:CUH2}(ii) and any
$q_*\in[q_0/2,q_0]$, so it does not depend on $n_0$. By the same lemma, the proposition applies at
$t_*=T_{n_0}$, and by its part~(c), for every $(\alpha,\beta)$ with $|\alpha|\le N_y$ and every
compact $\mathcal C$, we have $\partial^\alpha_y\partial^\beta_vg(T_n)\to\partial^\alpha_y\partial^\beta_vg(T)$
in $H^{-1}(\mathcal C)$, with a bound in $L^2(\mathcal C)$. Hence
$\partial^\alpha_y\partial^\beta_vg(T_n)\rightharpoonup\partial^\alpha_y\partial^\beta_vg(T)$ weakly
in $L^2(\mathcal C)$ for $|\alpha|\le N_y$ and all $\beta$, and, by Sobolev's inequality in $\R^6$,
$g(T_n)\to g(T)$ in $C^k(\mathcal C)$ for $k\le N_y-4$. We only use these orders. Hence the leaf
parameters, the frames $p_g(t,\xi)$, $\delta\varpi$ and its derivatives converge locally uniformly
in $\xi$, and so do the weights $\chi(s_\xi(t))$, $\omega_\xi(t)$, $W_{\alpha\beta}(t)$ (through
$\rad_k(t)$), $e^{2q(t)|y|^2}$, and $\psi(s_{\xi_c}(t))$. Moreover, the components of the cell
energies, at a fixed level $(\alpha,\beta)$, converge weakly in $L^2$ on compact sets, with weights
that converge locally uniformly. By weak lower semicontinuity with converging weights, as in
Lemma~\ref{lem:LCd}(i) for the local part and~(iii), by convexity, for the nonlocal part, and
Fatou's lemma in the sum over levels, we obtain $\mathfrak A_c(T)\le\liminf_n\mathfrak A_c(T_n)$,
and similarly for $\mathfrak A_{c,\ge1}$. For $\mathfrak M_c$ and $\mathfrak B_{c,\ge1}$ the leaf
moments are integrals over all $z$. They converge, together with their $\xi$-derivatives, by the
tail bound of point~(3) of the proof of the local proposition, applied to the family $g(T_n)$, and
thus $\delta\varpi(T_n)\to\delta\varpi(T)$ with all derivatives, locally uniformly, and Fatou's
lemma applies. The nonlocal part is a convex continuous functional of the local part for a fixed
frame, and the change of frame contributes $o(1)$ on bounded sets, as in~(iii) of the same lemma. A
product of a nonnegative lower semicontinuous function and a nonnegative continuous function is
lower semicontinuous, and so is a supremum of lower semicontinuous functions.
\end{proof}

\begin{lem}[Upper bound to the right]\label{lem:CUright}
Let $T\in[T_0,T_{\max})$, and assume that (U1)--(U4) hold on $[T_0,T]$. Then there is $h>0$ with
$T+h<T_{\max}$ such that $\Phi_2,\Phi_3\le1$ on $[T,T+h]$. Consequently (U1)--(U4) hold on
$[T_0,T+h]$.
\end{lem}

\begin{proof}
By Theorem~\ref{thm:CUboot}, we have $\Phi_2,\Phi_3\le\frac12$ at $T$. By Lemma~\ref{lem:CUH2}(ii),
Proposition~\ref{prop:CULW} applies at $t_*=T$. Let $\tau$, $A$, $B$ be its constants, take
$\Omega=\Omega^\epsilon$, and let $0<h\le\tau$. In particular $T+h<T_{\max}$.

\emph{Comparison.} For $t\ge T$, $\omega_\xi(t)\le e^{-2H(s_\xi(T))}\le\Omega^\epsilon_\xi$, because
$H$ is nondecreasing, $s_\xi$ is increasing and $H_\epsilon\le H$. Moreover,
$\rad_k(t)\le\tilde\rad_k(t)$, because $|\rad'|\ge r_{\rm loc}$ by Remark~\ref{rem:UFradius}, and
$w_n(\rad)$ is increasing in $\rad$. All other weights are the same functions of $(t,\xi,z)$ in both
energies. Hence we have $\mathfrak B_c(t)\le\mathfrak B^\Omega_c(t)$ and
$\mathfrak B_{c,\ge1}(t)\le\mathfrak B^\Omega_{c,\ge1}(t)$. At $t=T$, we have
$\Omega^\epsilon_\xi\le e^{1/20}\omega_\xi(T)$, and thus $X^\Omega_c(T)\le e^{1/20}X_c(T)$ for every
$X$. We put $\beta_*:=c_\chi C_U+\mu_\sharp$, so that $\sup_c\mathfrak B_c(T)\le\beta_*$ as in the
proof of that lemma, and $\theta(h):=e^{Ah}+Bhe^{(A+B)h}$. Since
$\mu_\sharp\le c_\sharp/(54e^5)\le c_\chi C_U/(54e^5)$, we have $\beta_*\le\frac{11}{10}c_\chi C_U$.

\emph{$\Phi_2$.} Let $\sigma_*$ be the larger of the two source constants $\sigma_X$ of
Proposition~\ref{prop:CULW}(b). By this bound, we have
$\mathfrak A_c(t)\le \mathfrak B^\Omega_c(t)\le\theta(h)e^{1/20}\beta_*+\sigma_*he^{(A+B)h}\le2c_\chi
C_U$ if $\theta(h)\le\frac{19}{11}e^{-1/20}$ and $\sigma_*he^{(A+B)h}\le\frac1{10}c_\chi C_U$. Both
conditions hold for $h$ small, since $\theta(0)=1$. Moreover, we have
\begin{equation}
\mathfrak B_{c,\ge1}(t)\le e^{Ah+1/20}\mu_\sharp+Bhe^{(A+B)h+1/20}\beta_*+\sigma_*he^{(A+B)h}\le2\mu_\sharp
\end{equation}
if $e^{Ah+1/20}\le\frac32$ and $Bhe^{(A+B)h+1/20}\beta_*+\sigma_*he^{(A+B)h}\le\frac12\mu_\sharp$.

\emph{$\Phi_3$.} Let $L_\psi$ be the Lipschitz constant of $\psi$. Then we have
$\psi(s_{\xi_c}(t))\le\psi(s_{\xi_c}(T))+L_\psi v_*h$. At $T$, we have
$\psi(s_{\xi_c}(T)) \mathfrak B_c(T)\le\psi\mathfrak A_c(T)+\mathfrak B_{c,\ge1}(T)\le C'_{\rm
cl}\omega(S_b)\delta_1^2+\mu_\sharp\le(C'_{\rm cl}+1)\omega(S_b)\delta_1^2$, since
$\mu_\sharp\le e^{-2H(S_b+1)}\delta_1^2\le\omega(S_b)\delta_1^2$. Hence
\begin{multline}
\psi(s_{\xi_c}(t))\mathfrak A_c(t)\le e^{Ah+1/20}(C'_{\rm cl}+1)\omega(S_b)\delta_1^2+\big(L_\psi
v_*h\,\theta(h)+Bhe^{(A+B)h}\big)e^{1/20}\beta_*\\
+(1+L_\psi v_*)\sigma_*he^{(A+B)h}\le2C'_{\rm cl}\omega(S_b)\delta_1^2
\end{multline}
for $h$ small, because $C'_{\rm cl}\ge6$ and $e^{1/20}\cdot\frac76<2$.

\emph{(U1) and (U4).} By the two previous paragraphs, (U2) and (U3) hold on $[T_0,T+h]$. Then (U1)
holds on $[T_0,T+h]$ by Lemma~\ref{lem:CUH2}(i), and (U4) by Lemma~\ref{lem:BTU4}, applied on
$[T_0,T+h]$. This uses no time regularity of the (G3)--(G4) norms.

Choosing $h$ as the smallest of the bounds for $\Phi_2$ and $\Phi_3$, we obtain the claim. Openness
needs only some $h>0$, and we use no uniformity of $h$ in $T$.
\end{proof}

\subsubsection{Time regularity of the cell energies of the solution}\label{sssec:CUac}

The cell inequalities of this section and of Section~\ref{sec:boot} are differential inequalities
for the cell energies of the solution itself, not of an approximation. In this subsection we state
in which sense these energies can be differentiated in time. There are two difficulties. First, the
components $h_{\alpha\beta}$ of $y$-order $|\alpha|=N_y$ have no further $y$-derivative in the
energy space. Moreover, the collision term of their equation is a distribution of order two along
the leaves, and the pairing with it is defined only if the dissipation of the component is finite.
The local theory provides it, in Proposition~\ref{prop:CULW}(d). Second, the sums over the levels
$|\beta|$ can be differentiated term by term only if the radius term $\sum_nn\,w_n^2E^c_n$ is
integrable in time, which is not known at the radii $\rad_k(t)$. Here we give up a margin $\iota>0$
in the radius, which is removed at the end, in Step~3 of the proof of Theorem~\ref{thm:CUboot}.

Throughout this subsection, $T\in[T_0,T_{\max})$, $g$ is the classical solution in the cell norms of
the proof of Corollary~\ref{cor:CUmain}, (U1)--(U4) hold on $[T_0,T]$, and (U1)--(U4) hold with the
halved constants of that theorem at each of the times $t_j:=T_0+j\tau_\flat<T$, $j\ge0$, with
$\tau_\flat$ of the same proof. To avoid circularity, we do not use that theorem in this subsection.
The second assumption holds at $t_0=T_0$ by Step~0 of that proof, and the proof obtains it at the
other times $t_j$ by induction on $j$, in Step~2 (time regularity and radii with a margin). For
$|\alpha|\le N_y$, $\beta$ arbitrary and a cell $c$, we put
\begin{equation}
E^c_{\alpha\beta}(t):=\rad_0^{-3}\!\int\chi_c^2\,\omega_\xi\,e_\xi[\varphi_{\alpha\beta}]\dd\xi,\qquad
D^c_{\alpha\beta}(t):=\rad_0^{-3}\!\int\chi_c^2\,d_\xi[\varphi_{\alpha\beta}]\dd\xi ,
\end{equation}
so that $\mathfrak A_c=\sum W^2_{\alpha\beta}E^c_{\alpha\beta}$, and let
\begin{equation}\label{eq:D1c}
\mathfrak D^{(1)}_c:=\sum W^2_{\alpha\beta}D^c_{\alpha\beta}\ge\mathfrak D_c
\end{equation}
be the undamped cell dissipation.

\begin{lem}[The dissipation is integrable]\label{lem:CUdiss}
For every cell $c$,
\begin{equation}
\int_{T_0}^T\kappa\,\mathfrak D^{(1)}_c\dd t\le\frac{3M}{c_{\rm
ret}}\Big(\frac{T-T_0}{\tau_\flat}+1\Big),
\end{equation}
with $M$ of Lemma~\ref{lem:CUH2}(ii), $c_{\rm ret}$ of Proposition~\ref{prop:CULW}(d) and
$\tau_\flat$ of the proof of Corollary~\ref{cor:CUmain}. In particular $\kappa D^c_{\alpha\beta}\in
L^1(T_0,T)$ for every $(\alpha,\beta)$, and $\varphi_{\alpha\beta}(t,\xi,\cdot)\in Y_2(\hat m_\xi)$
for a.e.\ $(t,\xi)$.
\end{lem}

\begin{proof}
Let $t_j:=T_0+j\tau_\flat$ for the integers $j\ge0$ with $t_j<T$. The last clause of
Lemma~\ref{lem:CUH2}(ii) rests on the halved bounds (U1)--(U4) at $t_*$ and on Step~1 of the proof
of Theorem~\ref{thm:CUboot} at $t_*$. That step is not static, since it involves the comparison on
$[T_0,t_*]$ of Corollary~\ref{cor:Aleafentry}, but it rests on the hypotheses (U1)--(U4) on
$[T_0,t_*]$ alone. With the standing assumption on the times $t_j$, that clause therefore gives the
hypotheses (H1$'$)--(H2$'$) of the local proposition at $t_*=t_j$, with constants that do not depend
on $j$, and with $q_*=q(t_j)\in[q_0/2,q_0]$. Thus the time $\tau$ of that proposition is at least
$\tau_\flat$. By the uniqueness of Lemma~\ref{lem:LWuniq}, used as in the proof of
Corollary~\ref{cor:CUmain}, $g$ coincides on $I_j:=[t_j,\min\{t_j+\tau_\flat,T\}]$, where both are
defined, with the solution of the proposition started at $t_j$. The intervals $I_j$ cover $[T_0,T]$.
Therefore, by part~(d) of the proposition, we have
$\int_{I_j}\kappa \mathfrak D^{(1)}_c\le3M/c_{\rm ret}$, at the radii
$\tilde\rad_k(t)=\rad_k(t_j)-\int_{t_j}^tr_{\rm loc}$. Since $|\rad'|\ge r_{\rm loc}$ by
Remark~\ref{rem:UFradius}, we have $\rad_k(t)\le\tilde\rad_k(t)$, and the dissipation is
nondecreasing in the radii. Summing over $j$, we obtain the bound. The last statement follows by
Fubini's theorem, since the weights are positive.
\end{proof}

\begin{lem}[Two energy identities]\label{lem:CUident}
Let $I=[a,b]$ be a compact interval.
\begin{enumerate}[label=(\roman*)]
\item \emph{Profile variables.} Let $D=\partial_v-t\nabla_y$, and let $u\in L^\infty(I;L^2_{\rm
loc}(\R^6))\cap C(I;H^{-1}_{\rm loc}(\R^6))$ satisfy $\partial_tu=D\cdot G_1+G_0$ in $
\mathcal D'((a,b)\times\R^6)$, with $Du,G_1\in L^2(I;L^2_{\rm loc})$ and $G_0\in L^1(I;L^2_{\rm loc})$. Then,
for every Lipschitz function $\Psi$ on $I\times\R^6$ with support in $I\times\mathcal C$ for a
compact set $\mathcal C$, and for all $a\le t'\le t\le b$,
\begin{equation}
\iint\Psi\,u^2\dd y\dd v\,\Big|_{t'}^{t}=\int_{t'}^{t}\!\iint\Big(\partial_s\Psi\,u^2-2D(\Psi
u)\cdot G_1+2\Psi\,u\,G_0\Big)\dd y\dd v\dd s .
\end{equation}
\item \emph{The nonlocal pairing.} Let $\varphi\in L^\infty(I;X_7)\cap L^2(I;Y_7)$ have a weak
derivative $\varphi'=R_1+R_0$ with $R_1\in L^2(I;Z_7)$ and $R_0\in L^2(I;X_7)$. Then $\varphi\in
C(I;X_7)$, the function $t\mapsto\nn{\varphi(t)}^2_\natural$ is absolutely continuous on $I$, and
its derivative is $2\ip{R_1+R_0}{\varphi}_\natural$ for a.e.\ $t$.
\end{enumerate}
\end{lem}

\begin{proof}
(i) Let $\phi_\epsilon$ be a mollifier in $(y,v)$ and $u^\epsilon:=\phi_\epsilon*u$, and similarly
$G_1^\epsilon$, $G_0^\epsilon$. The coefficients of $D$ do not depend on $(y,v)$, so
$\partial_tu^\epsilon=D\cdot G_1^\epsilon+G_0^\epsilon$ and $Du^\epsilon=(Du)^\epsilon$. For every
$k$ and every compact set, $u^\epsilon\in C(I;C^k)$, because $u\in C(I;H^{-1}_{\rm loc})$, and the
right side lies in $L^1(I;C^k)$. Hence $u^\epsilon$ is absolutely continuous in $t$ with values in
$C^k$ on compact sets, the function $t\mapsto\iint\Psi|u^\epsilon|^2$ is absolutely continuous, and
after an integration by parts its derivative is
$\iint(\partial_t\Psi|u^\epsilon|^2-2D(\Psi u^\epsilon)\cdot G_1^\epsilon+2\Psi u^\epsilon
G_0^\epsilon)$. We integrate over $[t',t]$ and let $\epsilon\to0$. For every $s\in I$,
$u(s)\in L^2_{\rm loc}$, with the bound of $L^\infty(I;L^2_{\rm loc})$. Indeed, for a compact set
$\mathcal C'$ we choose $s_n\to s$ with
$\|u(s_n)\|_{L^2( \mathcal C')}\le\|u\|_{L^\infty(I;L^2(\mathcal C'))}$, and a subsequence converges
weakly in $L^2(\mathcal C')$, and its limit is $u(s)$, by the continuity in $H^{-1}(\mathcal C')$.
Hence $u^\epsilon(s)\to u(s)$ in $L^2(\mathcal C)$. Moreover, $u^\epsilon\to u$,
$Du^\epsilon\to Du$, $G_1^\epsilon\to G_1$ in $L^2(I\times\mathcal C)$, $G_0^\epsilon\to G_0$ in
$L^1(I;L^2(\mathcal C))$.

(ii) Let $V:=X_7\cap Y_7$. Then $V\subset X_7\subset V'$ is a Gelfand triple, and $V'=X_7+Z_7$,
where $Z_7$ is the dual of $Y_7$ for the $X_7$-pairing, as in Appendix~\ref{app:band}. Thus
$\varphi\in L^2(I;V)$ and $\varphi'\in L^2(I;V')$, and by the Lions--Magenes lemma \cite[Ch.~XVIII,
\S1, Thm.~1]{DL}, we have $\varphi\in C(I;X_7)$. Its proof, by regularization in time, provides
$\varphi_\epsilon\in C^1(I;V)$ with $\varphi_\epsilon\to\varphi$ in $L^2(I;V)$ and in $C(I;X_7)$,
and $\varphi_\epsilon'\to\varphi'$ in $L^2(I;V')$. The form $\ip\cdot\cdot_\natural$ is symmetric,
bounded on $X_7\times X_7$ by Lemma~\ref{lem:lowweight} and on $Z_7\times Y_7$ by (CM1), hence on
$V'\times V$. Thus we have
$\frac{d}{dt}\nn{\varphi_\epsilon}^2_\natural=2\ip{\varphi_\epsilon'}{\varphi_\epsilon}_\natural$.
Then we integrate over $[t',t]$ and let $\epsilon\to0$.
\end{proof}

\begin{lem}[The energy identity of one component]\label{lem:CUcompid}
Let $|\alpha|\le N_y$, let $\beta$ be arbitrary, and let $c$ be a cell. Put $\xi=v+y/t$ and
\begin{equation}
\Theta_c(t,y,v):=\chi_c(\xi)^2\,\omega_\xi(t)\Big(\frac{A_{p_g(t,\xi)}}{A_{\rm
ref}(\xi)}\Big)^2A_{\rm ref}(v)^2e^{2q(t)|y|^2},
\end{equation}
so that $\eta\rad_0^{-3}\iint\Theta_c|h_{\alpha\beta}|^2\dd y\dd v$ is the local part of
$E^c_{\alpha\beta}$. Then $E^c_{\alpha\beta}$ is absolutely continuous on $[T_0,T]$. For a.e.\ $t$
its derivative is the integral against $\rad_0^{-3}\chi_c^2\dd\xi$ of the terms that the formal
differentiation of $\omega_\xi e_\xi[\varphi_{\alpha\beta}]$ along the rays produces (the terms
listed in Subsection~\ref{sssec:bandproof}, before they are estimated, and they are the same on all
leaves, with $\chi(s_\xi)=0$ for $s_\xi\le S_b$), read as follows.
\begin{itemize}
\item The pairings of $\varphi_{\alpha\beta}$ with the collision terms, in the local and in the
nonlocal part, are the weak forms of Lemma~\ref{lem:QmG}, of Proposition~\ref{prop:trilinG} and of
(CM2). By Lemma~\ref{lem:CUdiss} they are defined for a.e.\ $(t,\xi)$.
\item The local transport (T9) and the Gaussian weight term (T8) are replaced, together, by
\begin{equation}
\eta\rad_0^{-3}\!\iint\Theta_c\,\mathfrak w_c\,|h_{\alpha\beta}|^2\dd y\dd v,\qquad
\mathfrak w_c:=-t^{-2}\,y\cdot\nabla_\xi\log\Big(\chi_c^2\,\omega_\xi\big(A_{p_g}/A_{\rm
ref}\big)^2\Big)(t,\xi)-2C_2t^{-2}|y|^2 .
\end{equation}
The terms (T8) and (T9) take this form after the integration by parts in $\xi$ of
Lemma~\ref{lem:CUtr}, and
\begin{equation}
\mathfrak w_c\le t^{-2}\Big(\frac{C_{\rm tr}}{\rad_0^2C_2}-C_2|y|^2\Big).
\end{equation}
\end{itemize}
Each of these terms is integrable on $[T_0,T]$.
\end{lem}

\begin{proof}
We write $u:=h_{\alpha\beta}$ and $I:=[T_0,T]$. Sums over $(\alpha',\beta')$ run over
$\alpha'\le\alpha$, $\beta'\le\beta$, with $(\alpha'',\beta''):=(\alpha-\alpha',\beta-\beta')$.

\emph{The equation in divergence form.} On a leaf,
$Q(F,G)=\nabla_z\cdot\big((a*F)\nabla_zG-(b*F)G\big)$, and $\nabla_z[k^{(x)}]=-t^{-1}(Dk)^{(x)}$ for
a profile function $k$. Hence we have $\cQ_t(h_1,h_2)=D\cdot\mathcal J_t(h_1,h_2)$, where
$\mathcal J_t$ is the bilinear flux whose restriction to the leaf through $x$ is
$-t^{-1}\big((a*h_1^{(x)})\nabla_zh_2^{(x)}-(b*h_1^{(x)})h_2^{(x)}\big)$. The derivatives
$\partial_y$, $\partial_v$ commute with $D$, and they act on $\mathcal J_t(h_1,h_2)$ by the Leibniz
rule, since $\ab[h](y,v)=\int a(w)h(y+tw,v-w)\dd w$, and similarly for $\bb$. By
Lemma~\ref{lem:Meq}, in $\mathcal D'((T_0,T)\times\R^6)$, we have
\begin{equation}
\partial_tu=D\cdot G_1+G_0,\qquad G_1:=\kappa\sum\tbinom\alpha{\alpha'}\tbinom\beta{\beta'}\Big(
\mathcal J_t(\mathcal M',\hat r'')+\mathcal J_t(\hat r',\mathcal M'')+\mathcal J_t(\hat r',\hat
r'')\Big),\qquad G_0:=-\partial^\alpha_y\partial^\beta_v\partial_t\mathcal M[g],
\end{equation}
where $\hat r':=\partial^{\alpha'}_y\partial^{\beta'}_v\hat r$, $
\mathcal M':=\partial^{\alpha'}_y\partial^{\beta'}_v\mathcal M[g]$, and similarly with two primes.

\emph{The hypotheses of Lemma~\ref{lem:CUident}(i).} By (U2), Lemma~\ref{lem:CUsob} and
Proposition~\ref{prop:CULW}(c), applied on the intervals used for the dissipation bound, we have
$u\in L^\infty(I;L^2_{\rm loc})\cap C(I;H^{-1}_{\rm loc})$, and the Maxwellian part is continuous.
On the leaves, $Du$ is $-t\nabla_zh^{(x)}_{\alpha\beta}$. By the integrability of the dissipation
and Lemma~\ref{lem:TG0}(ii), $\nabla_zh^{(x)}_{\alpha'\beta'}\in L^2(I\times Q\times\R^3_z)$ for
every cube $Q$ and every $(\alpha',\beta')$, and $\dd\xi\dd z=\dd y\dd v$. Thus
$Du\in L^2(I;L^2_{\rm loc})$. In $G_1$, every factor $\nabla_zh_2^{(x)}$ is of this kind or is the
gradient of a Maxwellian part, which is bounded by a Gaussian. The coefficients satisfy
$|a*k^{(x)}|\le C\|k^{(x)}\|_{L^2(\jap z^4)}$, which is bounded uniformly in $(t,\xi)$ for
$k=\hat r'$ and $k=\mathcal M'$, by the Sobolev bounds used above and (F3). If $|\alpha'|\le N_y-1$,
then also $|b*k^{(x)}|\le C\|k^{(x)}\|_{H^1(\jap z^{5/2})}$ is bounded uniformly. If
$|\alpha'|=N_y$, then $\alpha''=0$, the factor $h_2$ is bounded by a Gaussian by
Lemma~\ref{lem:LCdata}, and $b*k^{(x)}\in L^2_{\rm loc}(\dd z)$ with a uniform bound, because $b$ is
the sum of an integrable function and a bounded one. Hence $G_1\in L^2(I;L^2_{\rm loc})$. Finally,
$G_0$ is bounded on compact sets, by (F3).

\emph{The local part.} We take $\zeta\in C^\infty_c(\R^3)$ radial and nonincreasing in $|\cdot|$,
with $0\le\zeta\le1$ and $\zeta=1$ on the unit ball (so that $\zeta(\cdot/R)$ increases with $R$),
and we apply the identity~(i) with $\Psi:=\zeta(y/R)^2\zeta(v/R)^2\Theta_c$, which is Lipschitz on
compact sets: $\omega_\xi$ is Lipschitz, and $p_g$ is $C^1$. Since $\xi$ is constant on the leaves,
$D$ annihilates the functions of $(t,\xi)$, and we have
\begin{equation}
\partial_t\Theta_c=\Theta_c\big(\mathfrak w_{0}+\mathfrak w_c\big),\qquad
\mathfrak w_0:=-2H'(s_\xi)\dot s_\xi+2\partial_t\log A_{p_g}(t,\xi),\qquad D\log\Theta_c=2\nabla\log A_{\rm
ref}(v)-4q\,t\,y,
\end{equation}
for a.e.\ $t$. Indeed, at fixed $(y,v)$, $\partial_t\xi=-y/t^2$ and $q'=-C_2t^{-2}$. The first term
of $\mathfrak w_0$ is the damping, and the second is the local part of the frame term (T11), and
$|\mathfrak w_0|\le C(\kappa\lambda_\xi+t^{-2})$. We obtain the bound on $\mathfrak w_c$ from
$|\nabla_\xi\log(\chi_c^2\omega_\xi(A_{p_g}/A_{\rm ref})^2)|\le\Gamma\le3/\rad_0$ as in the proof of
Lemma~\ref{lem:CUtr}, and $\Gamma|y|\le C_2|y|^2+\Gamma^2/(4C_2)$. Now we let $R\to\infty$ in this
identity.
\begin{itemize}
\item The boundary terms increase to $\iint\Theta_cu^2$ at $t$ and at $t'$, which are finite.
\item The term with $\mathfrak w_0$ and the positive part of $\mathfrak w_c$, which is bounded,
converge by dominated convergence, and the negative part of $\mathfrak w_c$ by monotone convergence.
\item In $-2D(\Psi u)\cdot G_1$, the part in which $D$ falls on the cutoff is
$O(R^{-1})\Theta_c|u||G_1|$. In the rest, we change variables to $(\xi,z)$. For each term of $G_1$
the integral over $z$ is one of the weak forms of the first item of the statement, multiplied by
$\chi_c^2\omega_\xi$ and with the factor $\zeta^2$ under the integral: for instance the terms
$\mathcal J_t(\mathcal M[g],u)+\mathcal J_t(\hat r,u)$ give, up to the factor
$\eta\rad_0^{-3}\chi_c^2\omega_\xi$, the pairing
$2\kappa\lambda_\xi\langle\varphi,Q(\mu+f_\xi,\varphi)\rangle_{X(\hat m_\xi)}$ in the form
$A_1+\dots+A_4+B_1+B_2$ of Lemma~\ref{lem:QmG}. The proofs of this lemma and of
Proposition~\ref{prop:trilinG} bound the integrands pointwise, so each integrand is integrable over
$I\times\R^3_\xi\times\R^3_z$, with the bounds of these statements, by the integrability of the
dissipation and the uniform bounds on the leaf energies. Therefore these terms converge by dominated
convergence, and the part with the derivative of the cutoff tends to $0$.
\item The term with $G_0$ converges by dominated convergence, and it is the local part of the
forcing (T12).
\end{itemize}
Since all other terms have finite limits, the integral of $\Theta_c\mathfrak w_c^-u^2$ over
$[t',t]\times\R^6$ is finite. This proves the identity for the local part, with every term
integrable.

\emph{The nonlocal part.} By the integrability of the dissipation and Fubini's theorem, for a.e.\
$\xi$ the functions $\varphi_{\alpha'\beta'}(\cdot,\xi,\cdot)$ lie in $L^2(I;Y_2(\hat m_\xi))$,
hence in $L^2(I;Y_7)$, and, by Lemma~\ref{lem:TG0}(ii), in $L^\infty(I;X_7)$. We fix such a $\xi$
and put $\varphi(t):=\varphi_{\alpha\beta}(t,\xi,\cdot)$. The identity \eqref{eq:leaftransport}
holds in $\mathcal D'((T_0,T)\times\R^3_\xi\times\R^3_z)$ for $h=h_{\alpha\beta}$, with
$(\partial_th)^{(x)}=-t\,\nabla_z\cdot G_1^{(x)}+G_0^{(x)}$.

Testing it with products of functions of $\xi$, of $t$ and of $z$, we see that for a.e.\ $\xi$ the
function $t\mapsto h^{(x)}_{\alpha\beta}$, $x=t\xi$, has a weak derivative, given by the right side
of this identity. Since the frame $p_g(t,\xi)$ is $C^1$ in $t$, $\varphi$ has the weak
derivative given by Lemma~\ref{lem:leafeq}(ii).

We control the terms of this derivative as follows.
\begin{itemize}
\item The collision part $\kappa\lambda_\xi(L_{\mu+f_\xi}\varphi+\Sigma_\varphi)$ lies in
$L^2(I;Z_7)$, by (CM2) and the first sentence of \emph{The nonlocal part}.
\item The frame term $\mathcal G_p[D_tp]\varphi$ lies in $L^\infty(I;Z_7)$, as in the estimate of
the nonlocal part of (T11) in Proposition~\ref{prop:FM}.
\item The transport and the forcing lie in $L^\infty(I;X_7)$, as in (T10) and (T12).
\end{itemize}

By the identity~(ii), $t\mapsto\nn{\varphi(t)}^2_\natural$ is absolutely continuous, and
its derivative is the sum of the nonlocal parts of the terms of Subsection~\ref{sssec:bandproof}.

The factor $\chi_c^2\omega_\xi\chi(s_\xi)$ is Lipschitz in $t$, and its derivative produces the
damping and the front term (T7). The bounds of that subsection for these terms, with the
integrability of the dissipation, show that they are integrable over $I\times\R^3_\xi$ against
$\chi_c^2\dd\xi$. Thus the integral over $\xi$ is absolutely continuous in $t$, and its derivative
is the integral of the derivatives.
\end{proof}

\begin{cor}[Cell energies with a margin in the radius]\label{cor:CUac}
Let $\iota\in(0,\epsilon_y\rad_0]$. Let $\mathfrak A^\iota_c$, $\mathfrak A^\iota_{c,\ge1}$,
$\mathfrak M^\iota_c$, $\mathfrak B^\iota_c$, $\mathfrak B^\iota_{c,\ge1}$ and $\mathfrak D^\iota_c$
be the cell quantities of Subsections~\ref{ssec:cells} and~\ref{ssec:UF} at the radii
$\rad_k(t)-\iota$ in place of $\rad_k(t)$, that is, with
\begin{equation}\label{eq:Wiota}
W^\iota_{\alpha\beta}:=\delta_y^{|\alpha|}w_{|\beta|}(\rad_{|\alpha|}(t)-\iota)
\end{equation}
in place of $W_{\alpha\beta}$ and $w_n(\rad(t)-\iota)$ in place of $w_n(\rad(t))$.
\begin{enumerate}[label=(\roman*)]
\item The radii $\rad_k-\iota$ are graded as in Remark~\ref{rem:UFgrade}, they lie in
$[\rad_0/2,\rad_0]$, and they decrease at the rate $\rad'$ of Remark~\ref{rem:UFradius}. Each of the
quantities above does not exceed the corresponding quantity at the radii $\rad_k$, and increases to
it as $\iota\downarrow0$.
\item For every $t$,
\begin{equation}\label{eq:CUac-sum}
\sum_{(\alpha,\beta)}(1+|\beta|)(W^\iota_{\alpha\beta})^2E^c_{\alpha\beta}\le C_\iota\,\mathfrak A_c
\end{equation}
with $C_\iota:=\sup_n(1+n)(1-\iota/\rad_0)^{2n}<\infty$, and the same holds for $\mathfrak M_c$.
\item The energies $\mathfrak A^\iota_c$, $\mathfrak A^\iota_{c,\ge1}$, $\mathfrak M^\iota_c$,
$\mathfrak B^\iota_c$, $\mathfrak B^\iota_{c,\ge1}$ are absolutely continuous on $[T_0,T]$, and for
a.e.\ $t$ their derivatives are the sums of the derivatives of their terms, which are given by
Lemma~\ref{lem:CUcompid} and, for $\mathfrak M^c_n$, by the leaf Euler system as in
Proposition~\ref{prop:UFforce}(f). The series of the absolute values of these derivatives converge
in $L^1(T_0,T)$.
\end{enumerate}
Consequently Lemma~\ref{lem:CUtr}, Propositions~\ref{prop:CUcell}, \ref{prop:UFforce}
and~\ref{prop:UFge1} and Lemma~\ref{lem:TWabs} hold for the energies at the radii $\rad_k-\iota$,
as inequalities for a.e.\ $t\in[T_0,T]$ between derivatives of absolutely continuous functions,
under their hypotheses at the radii $\rad_k$ and with the same constants.
\end{cor}

\begin{proof}
(i) By \eqref{eq:radrate} and requirement~\ref{rq:TWb}, we have $\rad(t)\ge\frac34\rad_0$, and thus
$\rad_{N_y}(t)=\rad(t)-N_y\epsilon_y\rad_0\ge\rad_0\big(\frac34-\frac{N_y}{4N_y+4}\big)=\rad_0\big(\frac12+\epsilon_y\big)$,
and $\rad_k-\iota\ge\rad_0/2$. The differences $\rad_k-\rad_{k+1}=\epsilon_y\rad_0$ and the rate are
unchanged, and $(\rad_k-\iota)/(\rad-\iota)\le1-\epsilon_y$ for $k\ge1$, as in~\eqref{eq:UFgap},
because $\rad-\iota\le\rad_0$. The weights $w_n(\rad)$ are increasing and continuous in $\rad$, so
the last statement holds by monotone convergence.

(ii) This follows from $w_n(\rad_k-\iota)/w_n(\rad_k)=(1-\iota/\rad_k)^n\le(1-\iota/\rad_0)^n$.

(iii) Each $\mathfrak M^c_n$ is Lipschitz on $[T_0,T]$: $\partial^\beta_\xi\delta\varpi$ is
continuously differentiable along the rays, because the moments of $\partial^\beta_vg$ are, by the
Gaussian bounds of Lemma~\ref{lem:LCdata} on
$\partial_t\partial^\beta_vg=\partial^\beta_v(\kappa\cQ_t[g])$, and $\omega_\xi$ is Lipschitz. By
Lemma~\ref{lem:CUcompid}, each $(W^\iota_{\alpha\beta})^2E^c_{\alpha\beta}$ is absolutely
continuous, with the derivative
\begin{equation}\label{eq:CUac-deriv}
2|\beta|\frac{\rad'}{\rad_{|\alpha|}-\iota}(W^\iota_{\alpha\beta})^2E^c_{\alpha\beta}+(W^\iota_{\alpha\beta})^2\frac{d}{dt}E^c_{\alpha\beta}
.
\end{equation}
A series of absolutely continuous functions that converges at one point, and whose derivatives
converge absolutely in $L^1$, converges to an absolutely continuous function, and can be
differentiated term by term.

Thus we only need to bound the sum of the $L^1$ norms. The first term is bounded by
$4|\rad'|\rad_0^{-1}C_\iota\mathfrak A_c$, by~(ii). For the second, the estimates of this section,
of Section~\ref{sec:boot} and of Appendix~\ref{app:band} are estimates of the terms of the energy
identity of one component at a fixed time. They are obtained by the Cauchy--Schwarz inequality, so
they bound the sums of the absolute values, and they hold at the radii $\rad_k-\iota$ by~(i). They
yield the following bounds, with a constant $C_T$ that may depend on $T$ and on the constants of the
bootstrap:
\begin{itemize}
\item the damping and the front term are bounded by $C\kappa\mathfrak A_c$;
\item the collision terms (the linear terms, the sources and the coefficient terms) are bounded by
$C_T\kappa\mathfrak D^{(1)}_c+C_T\kappa\mathfrak A_c$, where the norms $\|f_\xi\|_X+\|f_\xi\|_{Y_2}$
of the background are bounded as in point~(5)($\beta$) of the proof of Proposition~\ref{prop:CULW};
\item the frame terms are bounded by $Ct^{-2}(\mathfrak A_c+\mathfrak D^{(1)}_c)$;
\item the nonlocal transport, the kernel part and the forcing are bounded by
$C_Tt^{-2}(1+C_\iota\sup_{c'}\mathfrak B_{c'})$, where the one-level losses, which carry a factor
$(1+|\beta|)/\rad_0$, and the spread losses are covered by~(ii);
\item the positive part of $\mathfrak w_c$ contributes at most $Ct^{-2}\rad_0^{-2}\mathfrak A_c$.
\end{itemize}
All these bounds are integrable on $[T_0,T]$, by (U2) and the integrability of the dissipation.

We still have to bound the Gaussian terms
\begin{equation}\label{eq:Gab}
G_{\alpha\beta}:=\eta\rad_0^{-3}\iint\Theta_c\mathfrak w_c^-|h_{\alpha\beta}|^2\ge0 .
\end{equation}
Integrating the derivative~\eqref{eq:CUac-deriv} over $[T_0,T]$, and using $\rad'\le0$ and
$E^c_{\alpha\beta}\ge0$, we obtain
\begin{equation}
\int_{T_0}^T(W^\iota_{\alpha\beta})^2G_{\alpha\beta}\dd
t\le(W^\iota_{\alpha\beta})^2E^c_{\alpha\beta}(T_0)+\int_{T_0}^T(W^\iota_{\alpha\beta})^2\big|\tfrac{d}{dt}E^c_{\alpha\beta}+G_{\alpha\beta}\big|\dd
t ,
\end{equation}
and the sum over $(\alpha,\beta)$ of the right side is finite by the bounds just listed. The
macroscopic energies are treated in the same way, with Proposition~\ref{prop:UFforce}(f) and~(ii).

\emph{The consequence.} The proofs of the statements listed estimate, at a fixed time, the sums over
$(\alpha,\beta)$ of the terms of the component identity. The bound of Lemma~\ref{lem:CUtr} is the
bound on $\mathfrak w_c$. The radii enter them through the properties of~(i) and through the
weights. Their hypotheses, (H$_a$)--(H$_e$), (H$_a'$) and (U1)--(U4), either do not involve the
radii or are upper bounds for cell and cube energies, which hold a fortiori at the radii
$\rad_k-\iota$. By~(iii) the left sides are derivatives of absolutely continuous functions, for
a.e.\ $t$.
\end{proof}

\section{The bootstrap}\label{sec:boot}

In this section we propagate the cell energies of Section~\ref{sec:cells} on $[T_0,\infty)$. This is
the content of the bootstrap theorem, Theorem~\ref{thm:CUboot}.

\subsection{The leaf clock and the band}\label{ssec:band}

\begin{figure}[ht]
\centering
\begin{tikzpicture}
\begin{axis}[width=0.6\textwidth,height=0.38\textwidth,xmin=1,xmax=40,ymin=0,ymax=6,
 xlabel={$t$},ylabel={$s_\xi(t)$},legend pos=outer north east,legend style={font=\small},tick label
style={font=\small}]
\addplot[domain=1:40,samples=200,thick,blue] {1.0*(x^(0.5)-1)/0.5};
\addlegendentry{$\lambda(\xi)=1$}
\addplot[domain=1:40,samples=200,thick,red!70!black] {0.4*(x^(0.5)-1)/0.5};
\addlegendentry{$\lambda(\xi)=0.4$}
\addplot[domain=1:40,samples=200,thick,orange!80!black] {0.15*(x^(0.5)-1)/0.5};
\addlegendentry{$\lambda(\xi)=0.15$}
\addplot[domain=1:40,dashed,black] {2};
\addplot[domain=1:40,dotted,black] {3};
\node[font=\small] at (axis cs:36,2.25) {$S_*$};
\node[font=\small] at (axis cs:36,3.25) {$S_b$};
\end{axis}
\end{tikzpicture}
\caption{The leaf clocks $s_\xi(t)=\lambda_{\rm ref}(\xi)(\tau(t)-\tau(T_0))$ for $\gamma=-2.5$ and
three collision rates. A leaf enters the relaxed regime when its clock reaches $S_*$, and the band
when it reaches $S_b$. Leaves with small $\lambda_{\rm ref}(\xi)$ enter late.}\label{fig:sclock}
\end{figure}
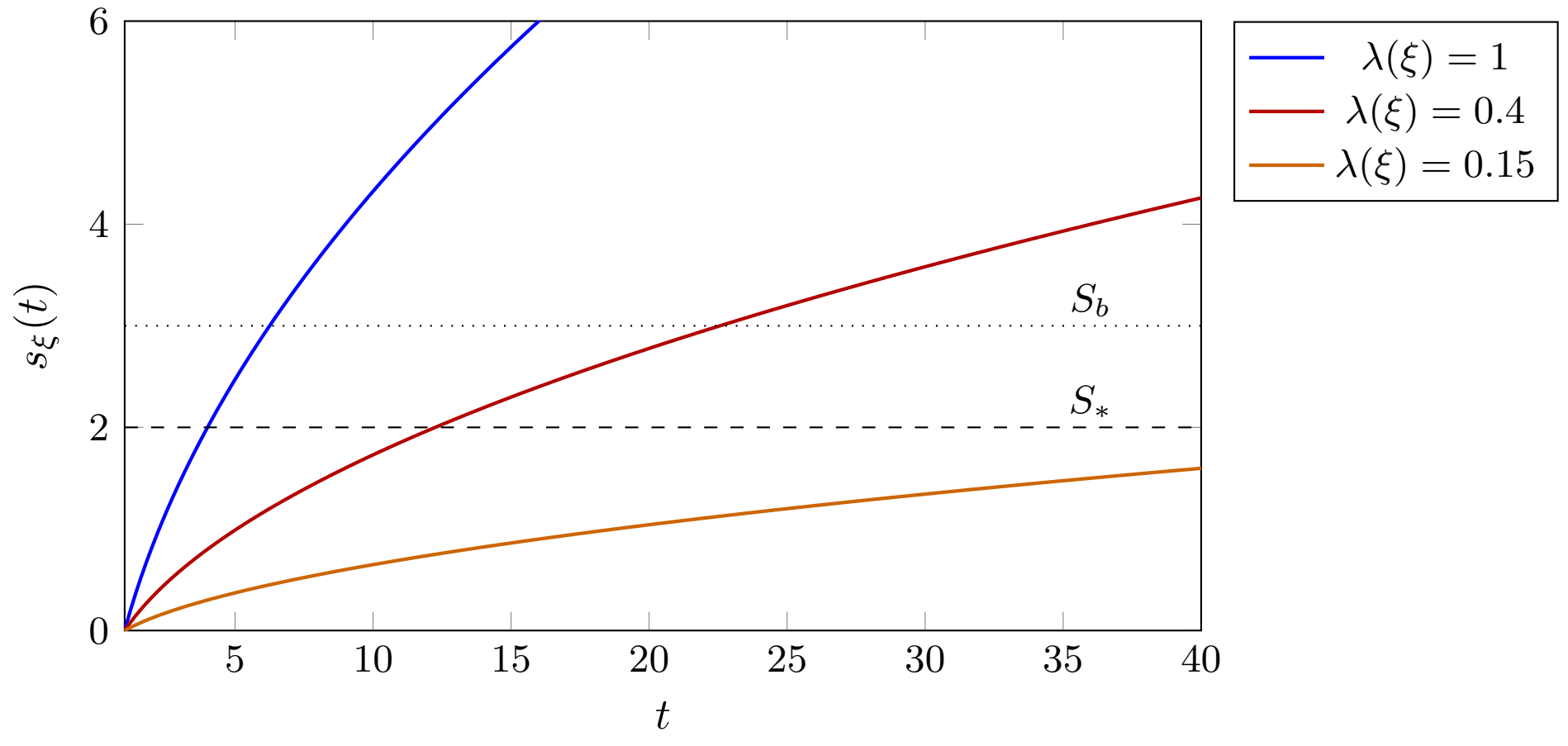

We run the bootstrap on $[T_0,T]$ and organize it by the leaf clock $s_\xi(t)$, as in
Figure~\ref{fig:sclock}. Let $S_*$ be the entry threshold and $[S_b,S_b+1]$ the band of
Subsection~\ref{ssec:cells}.

On the leaf $\xi$ the energy is weighted by the damping factor $\omega_\xi(t)=e^{-2H(s_\xi(t))}$ of
\eqref{eq:damping}, which depends on the effective time. The function $H$ increases up to the end of
the band and is constant after it. Differentiating the weight in time, we obtain the damping term
\begin{equation}\label{eq:dampterm}
-2H'(s_\xi)\,\lambda_{\rm ref}(\xi)\,\kappa\,\omega_\xi e_\xi .
\end{equation}
This term balances, term by term, the growth of the energy in the band, where the leaf is close to,
but not yet inside, the regime of Proposition~\ref{prop:triG}.

\begin{prop}[The balance in the band]\label{prop:Bband}
Suppose that on the leaf $\xi$:
\begin{itemize}
\item the normalized micro part is smaller than $\delta_c$ in $X(\hat m_\xi)$;
\item the damping satisfies $H'\ge C_1^B$ wherever $\chi<1$;
\item $C_{\rm r}t^{-2}\le\frac{\eta c_W}4\kappa\lambda_\xi$ wherever $\chi>0$, where $C_{\rm
r}t^{-2}\|\varphi\|^2_{Y_2(\hat m_\xi)}$ bounds the nonlocal part of (T11) in
Subsection~\ref{sssec:bandproof} ($C_{\rm r}=C_{11}\mathfrak p_*$, Corollary~\ref{cor:FMuse}(a));
\item $\mathfrak p_\xi\le\mathfrak p_*$ (this holds under (U1)--(U2) and the standing requirement
$C\delta_U\le\rad_0$, by Corollary~\ref{cor:FMuse}(d) with $B=B_{\rm cr}$).
\end{itemize}
Then, for each component $\varphi$ of the hierarchy,
\begin{multline}\label{eq:Bband}
D_t(\omega e_\xi[\varphi])\le-\tfrac{c_W}{4C_W}\kappa\omega d_\xi[\varphi]+C\kappa\lambda_\xi\omega
b_\xi e_\xi[\varphi]+2\kappa\lambda_\xi\omega|\mathcal S^{\rm src}_\xi(\varphi)|\\+Ct^{-2}\omega
e_\xi[\varphi]+C_0t^{-2}\omega\,e_\xi[\varphi]^{1/2}e_\xi[\varphi^+]^{1/2}+\omega\,\mathcal T^{\rm
loc}_\xi(\varphi)+Ct^{-2}\omega e_\xi[\varphi]^{1/2}F_\xi,
\end{multline}
where $b_\xi$ is the square of the size of the micro part, $\mathcal S^{\rm src}_\xi$ are the
sources, $\varphi^+$ is the component of the next order, $F_\xi$ is the forcing of
Lemma~\ref{lem:Meq}, and $\mathcal T^{\rm loc}_\xi(\varphi)$ is the local transport (T9),
$2t^{-2}\eta\langle\varphi,(y_0+\theta^{1/2}\hat
z)\cdot\varphi^+\rangle_{X(\hat m_\xi)}$, together with the Gaussian weight term (T8); it is not
estimated on the leaf. Here $e_\xi[\varphi^+]:=\sum_je_\xi[\varphi^{(\beta+e_j)}]$. In cells we
treat it by Lemma~\ref{lem:CUtr}.
\end{prop}

The proof is in Subsection~\ref{sssec:bandproof}. We add three comments on the statement. First,
the left side of \eqref{eq:Bband} stands for the sum of the pairings that the formal
differentiation of $\omega e_\xi[\varphi]$ along the ray produces. The energy of a single leaf need
not be differentiable along its ray (see the discussion of the meaning of the inequality in the
proof of Lemma~\ref{lem:UFlevel0}). Therefore we apply the inequality only after integration in
$\xi$, in Proposition~\ref{prop:CUcell} and, for the viscous approximation, in
Subsection~\ref{ssec:cellsproofs}, or, at level $0$, in the integrated form~\eqref{eq:L0int}. The
time regularity of the cell energies is established in Subsection~\ref{sssec:CUac}. The same
applies to Proposition~\ref{prop:ACpf}. Second, for $\alpha\ne0$ the additional coefficient terms of
Remark~\ref{rem:UFgrade} are part of $\mathcal S^{\rm src}_\xi$. Finally, the local transport
cannot be bounded leaf by leaf in a form that closes in the analytic sum. Indeed, by
Proposition~\ref{prop:ACtr}(ii) and Remark~\ref{rem:ACtr}, such a bound loses two levels. Thus we
keep it as $\mathcal T^{\rm loc}_\xi$.

\begin{remark}[Hypotheses of Proposition~\ref{prop:Bband}]\label{rem:Bbandrange}
The proof of \eqref{eq:Bband} in Subsection~\ref{sssec:bandproof} uses the following hypotheses of
the proposition:
\begin{enumerate}[label=(\alph*)]
\item the background is $\delta_c$-close to equilibrium;
\item the damping satisfies $H'\ge C_1^B$ wherever $\chi<1$, and there it absorbs the deficit
$(\frac12-\chi)_+\|\varphi\|_{L^2}^2$ and the front term;
\item $t^{-2}\ll\kappa\lambda_\xi$ wherever $\chi>0$, which we use to absorb the nonlocal part of
(T11) by (T1), under requirement~\ref{rq:Wrate};
\item the frame bound $\mathfrak p_\xi\le\mathfrak p_*$, used in (T11).
\end{enumerate}
Hypothesis~(c) holds for $s\ge S_b$ once
\begin{equation}\label{eq:Bband-T0}
T_0\ge\frac{4C_{\rm r}}{\eta c_Wc_\lambda S_*\min\{1,|\gamma|-2\}}\qquad\text{for }\gamma\in[-3,-2)
\end{equation}
which is that requirement, or Lemma~\ref{lem:EPreq} at $\gamma=-2$. Indeed,
$\kappa\lambda_\xi\ge c_\lambda\dot s_\xi\ge c_\lambda s_\xi\kappa/\tau$ and
$\kappa/\tau\ge\min\{1,|\gamma|-2\}/t$. The first inequality follows from (F3), which implies
$\lambda_\xi\ge c_\lambda\lambda_{\rm ref}(\xi)$ under (U1).
The term $\frac{\eta c_W}4$ comes from the coercivity. In the proof we have
$2\eta c_W-\frac12\eta c_W-\frac14\eta c_W-\eta c_W=\frac14\eta c_W$, where (T1) provides $2\eta
c_W$, (T4) uses $\frac12\eta c_W$ by requirement~\ref{rq:Wdc}, Young's inequality in (T4) uses
$\frac14\eta c_W$, and the conclusion uses $\eta c_W$. Hence the
same inequality, with the coercivity $c_{\rm cl}:=c_W/(4C_W)$, holds on every $\delta_c$-close leaf
with $\mathfrak p_\xi\le\mathfrak p_*$, $\chi=1$ and any $H'\ge0$ (post-band leaves), and on every
$\delta_c$-close leaf with $\mathfrak p_\xi\le\mathfrak p_*$, $\chi=0$ and $H'\ge C_1^B$ (for
instance buffer leaves with the crude damping $H'=C_1^U\ge C_1^B$, or with the auxiliary damping of
Lemma~\ref{lem:UFlevel0}). In particular it holds on $[S_*,S_b+1]$ and after the band.
\end{remark}

\subsection{The bootstrap statement}\label{ssec:bootstate}

We carry out the bootstrap in uniformly local form in the bootstrap theorem below. Its structure is
as follows. On the crude region, Theorem~\ref{thm:RGF} yields the Gaussian bounds, and the cell
inequality of Proposition~\ref{prop:CUcell}, built from the crude leaf inequality of
Proposition~\ref{prop:ACpf} on each leaf and from Lemma~\ref{lem:CUtr} for the local transport,
yields the analytic bounds with the damping $e^{-2C_1^Us}$. This is Step~3 of the proof of the
theorem. At entry, Lemma~\ref{lem:Aentry} and the same theorem give, in
Corollary~\ref{cor:Aleafentry}, the smallness $\delta_1$, using only information on the leaf itself.
In the pre-band regime and in the band, Proposition~\ref{prop:Bband} with the damping weight
closes the analytic bound (U3) at the entry level. On the relaxed region, the cell inequality, with
the parts of order $\ge1$ controlled by a level-$0$ estimate and an integral Gronwall argument over
the cubes, closes (U2)--(U3), and the macroscopic hierarchy closes (U1).
Subsection~\ref{ssec:inputs} passes from the cell energies to the inputs of the fibered comparison.
The theorem states these bounds for the cell energies, and Corollary~\ref{cor:CUmain}
reads them back pointwise in $\xi$. Continuation follows from the local existence theorem in the
cell norms, Proposition~\ref{prop:CULW}, which uses the viscous approximation of
Appendix~\ref{app:visc}. The constants are chosen in the order of Appendix~\ref{app:const}.

\subsection{The bootstrap in uniformly local form}\label{ssec:bootproof}

The bootstrap theorem below yields global existence by the continuation argument of
Corollary~\ref{cor:CUmain}.

\begin{theo}[Uniformly local bootstrap]\label{thm:CUboot}
Fix the constants in the order of Appendix~\ref{app:const}, with the following specifications:
\begin{itemize}
\item $c_\sharp$ as in Lemma~\ref{lem:CUnodrop} and (H$_a$);
\item $S_*$ as in the paragraph on entry levels;
\item $d_*$ and $S_b$ together with $\rad_0$;
\item requirements~\ref{rq:CUa}--\ref{rq:CUd} and~\ref{rq:UFa}--\ref{rq:TWb} in addition;
\item the graded weights of Remark~\ref{rem:UFgrade}, the rate of the radius of
Remark~\ref{rem:UFradius}, and $\delta_U$ as defined there.
\end{itemize}
Put
\begin{equation}\label{eq:deltab}
\delta_b^2:=\tfrac54C_be^{2C_1^B}\delta_1^2,\qquad C_b:=54e^5C_{\rm cl}'+28,
\end{equation}
with $C_{\rm cl}'$ structural. Require:
\begin{itemize}
\item $C_N\delta_b\le\min\{c''_{\min}/4,\,c'_{\min}/8\}$ and $\delta_b^2\le\delta_s$, which are
conditions on $\delta_1$;
\item for the continuation argument (Lemma~\ref{lem:CUH2}), $\delta_1/C_e\le\delta_c/2$ and
$C\delta_b\le\delta_c/2$, and $c_\sharp\le c_\chi$. The constant $C$ in the
definition of $\delta_U$ is taken $\ge1$;
\item independently of the continuation, $C'_{\rm Sob}\delta_U\le\rad_0\min\{\delta_s,1\}/4$
with $C'_{\rm Sob}:=8C_{\rm Sob}$ (Lemma~\ref{lem:TWcalc}(i)), split as in Remark~\ref{rem:UFradius}
into conditions on $T_0$ and on $\eps$. This condition implies Theorem~\ref{thm:main}(iii) and the
condition $C\delta_U\le\rad_0$ for $C_\varpi$.
\end{itemize}
Let
\begin{equation}\label{eq:psiramp}
s_\psi:=\max\{S_*,S_b-\tfrac14\},\qquad\psi(s):=\min\Big\{1,\max\Big\{0,\frac{s-s_\psi}{S_b-s_\psi}\Big\}\Big\}.
\end{equation}
Suppose that on $[T_0,T]$, $T<T_{\max}$:
\begin{itemize}
\item[(U1)] $\sup_Q\mathfrak m_{Q^*}\le(2\delta_U)^2$;
\item[(U2)] for every cell, $\mathfrak A_c\le2c_\chi C_U$ and $\mathfrak B_{c,\ge1}\le2\mu_\sharp$,
where $\mathfrak B_{c,\ge1}=\mathfrak A_{c,\ge1}+K_{\mathfrak M}\mathfrak M_c$ is the joint energy
of Subsection~\ref{ssec:UF} and $\mu_\sharp:=e^{-2H(S_b+1)}c_\sharp/(54e^5)$;
\item[(U3)] for every cell,
\begin{equation}\label{eq:U3}
\psi\big(s_{\xi_c}(t)\big)\mathfrak A_c(t)\le2C_{\rm cl}'\,\omega(S_b)\delta_1^2,
\end{equation}
with $\psi$ as in \eqref{eq:psiramp}. Therefore $\psi$ is Lipschitz with
$L_\psi=(S_b-s_\psi)^{-1}<\infty$ (since $d_*>0$, $S_b>S_*$), $\psi=0$ on $(-\infty,s_\psi]$, in
particular on the crude regime, and $\psi=1$ on $[S_b,\infty)$. In particular, $
\mathfrak A_c(t)\le2C_{\rm cl}'\,\omega(S_b)\delta_1^2$ for every cell with $s_{\xi_c}(t)\ge S_b$;
\item[(U4)] (G3) and (G4) with constant $K$ on $\{s_\xi\le3S_\sharp\}$, in the normalizations of
Remark~\ref{rem:BTnorm}.
\end{itemize}
Then the same hold with constants halved (with $K/2$ in (G3) and (G4)). By Lemma~\ref{lem:BTU4},
(U4) on $[T_0,T]$ is itself a consequence of (U1)--(U3) on $[T_0,T]$, and we list it among the
hypotheses only for convenience. At $t=T_0$ they hold in the halved form. Moreover, (U2) implies
$\mathfrak a_{Q^*,\ge1}\le c_\sharp$ for every cube, hence (H$_a$) of Proposition~\ref{prop:CUcell}
(since $c_\sharp\le\delta_1^2/4\le\delta_b^2$) and (H$_a'$) of Proposition~\ref{prop:UFge1}.
Together with (U3) and the conclusions of Step~2 of the proof, it implies (H$_b$).
\end{theo}

\begin{proof}
\emph{Step 0: $t=T_0$.} Every leaf has $s_\xi=0$ and $\omega_\xi=1$. By (LWP)(a1) of
Section~\ref{sec:proofmain}, $\mathfrak A_c(T_0)\le c_\chi C_DB$ and $\mathfrak A_{c,\ge1}(T_0)\le
c_\chi C_D\big((\rad_0/\Lambda)^2+\delta_y^2\big)B$. Also $\delta\varpi(T_0)=0$, since
$\varpi_*(T_0)=\varpi_g(T_0)$ by Lemma~\ref{lem:Mref}, with the normalization of $\ell_*$ stated
there, and Corollary~\ref{cor:LWPa}(a4). Thus $\mathfrak m_Q(T_0)=\mathfrak M_c(T_0)=0$, and
$\mathfrak B_{c,\ge1}(T_0)=\mathfrak A_{c,\ge1}(T_0)\le\mu_\flat/32\le\mu_\sharp/32$ by
requirement~\ref{rq:UFc}, since $\mu_\flat\le\mu_\sharp$ by Step~3.

Now we consider (U4) at $T_0$. By the bounds just proved, (U1) and (U2) hold at $T_0$, and (U3) at
$T_0$ is trivial, since $\psi(s_{\xi_c}(T_0))=\psi(0)=0$. Applying Lemma~\ref{lem:BTU4} with
$T=T_0$, we obtain (G3) and (G4) with $K/2$ on every leaf. Its proof uses Theorem~\ref{thm:RGF} at
$s_v=0$, whose hypotheses at $T_0$ are (H1)--(H2) for $g(T_0)$ by clause (a3) of (LWP)(a) in
the same section, and (H3), (a), (c) of Lemma~\ref{lem:RGbar} at $t=T_0$, from
Lemma~\ref{lem:BTsob}. No separate bound on the free leaves is needed.

\emph{Step 1: (H$_a$)--(H$_b$).} In this step we use (U2), (U3), and, for the closeness in (H$_b$),
Corollary~\ref{cor:Aleafentry}, which needs the hypotheses of the same theorem. These are
supplied by Step~2 below, which rests on (U1)--(U4) and on the second case of the present step. That
case rests on (U3) and Lemma~\ref{lem:CUcomp} alone, so Step~2 may be read first. By (U2) and
the same lemma, we have for every cube
\begin{equation}\label{eq:boot-Ha}
\mathfrak a_{Q^*,\ge1}\le27e^5e^{2H(S_b+1)}\sup_c
\mathfrak A_{c,\ge1}\le27e^5e^{2H(S_b+1)}\cdot2\mu_\sharp=c_\sharp ;
\end{equation}
this is a static consequence of (U2). This proves (H$_a$) and (H$_a'$).

For (H$_b$), let $Q$ contain a leaf with $s_\xi\ge S_*$. By Lemma~\ref{lem:CUcomp}(i), every leaf of
$Q^*$ has $s\ge S_*(1-\epsilon_0)\ge S_*-\frac12$. We consider two cases.
\begin{itemize}
\item If $s<S_b+1$ on $Q^*$, then $s\le3S_\sharp$ there, and by the entry corollary we have
$\|f_{\xi'}\|_X\le\delta_1/C_e$ on $Q^*$. Hence the level-$0$ part of $\mathfrak a_{Q^*}$ is
$\le27\delta_1^2$. Since $\mathfrak a_{Q^*,\ge1}\le c_\sharp\le\delta_1^2/4$ now holds for
every cube, we obtain $\mathfrak a_{Q^*}\le28\delta_1^2\le\delta_b^2$.
\item Otherwise some leaf of $Q^*$ has $s\ge S_b+1$. Then every leaf of $Q^*$ has $s\ge S_b$, and
so do the centers of the $27$ cubes that make up $Q^*$. By (U3), the same lemma and
$H(s)-H(S_b)\le C_1^B$, we have
\begin{equation}\label{eq:boot-Hb}
\mathfrak a_{Q^*}\le27\cdot e^5e^{2C_1^B}\omega(S_b)^{-1}\cdot2C_{\rm
cl}'\,\omega(S_b)\delta_1^2\le\delta_b^2 .
\end{equation}
\end{itemize}
The pointwise bound on $\|f_\xi\|_X$ follows from Lemma~\ref{lem:CUsob} at level $0$, together with
the requirement $C\delta_b\le\delta_c/2$, with $C$ the constant of that lemma
at level~$0$. It holds on $Q^*$, as (H$_b$) requires. Indeed, a leaf of $Q^*$ has
$s\ge S_*-\frac12$. If $s\le3S_\sharp$, the corollary implies
$\|f_\xi\|_X\le\delta_1/C_e\le\delta_c/2$. If $s>3S_\sharp\ge S_b+1$, the cube $Q'$ of this leaf has
$\mathfrak a_{Q'^*}\le\delta_b^2$ by the second case above, and by the same Sobolev lemma we get
$\|f_\xi\|_X\le C\delta_b\le\delta_c/2$.

\emph{Step 2: the inputs of the fibered comparison, and (U4).} By Lemma~\ref{lem:BTsob} below, whose
hypotheses are (U1)--(U3),
the crude bound (H3), the barrier hypotheses on the relaxed leaves and the
qualitative Gaussian decay hold on $[T_0,T]$. Then, by the lemma applied in Step~0,
Theorem~\ref{thm:RGF} implies (U4) with $K/2$.

\emph{Time regularity and radii with a margin.} In Steps~3--5 and in Proposition~\ref{prop:UFU2} we
integrate the cell inequalities in time. For the energies at the radii $\rad_k(t)$ we do not
know that the sums over the levels can be differentiated term by term. Therefore we run these steps
for the energies at the radii $\rad_k(t)-\iota$, with a fixed $\iota\in(0,\epsilon_y\rad_0]$.

Corollary~\ref{cor:CUac} needs, besides (U1)--(U4) on $[T_0,T]$, the halved bounds at the times
$t_j=T_0+j\tau_\flat<T$ of the beginning of Subsection~\ref{sssec:CUac}, and we obtain these by
induction on $j$. They hold at $t_0=T_0$ by Step~0. If they hold at $t_0,\dots,t_j$, we put
$T_{j+1}:=\min\{t_{j+1},T\}$. Then the corollary applies with $T_{j+1}$ in place of $T$,
since (U1)--(U4) hold on $[T_0,T_{j+1}]$, and Steps~1--5 (Steps~1--2 above, Steps~3--5 below), run
on $[T_0,T_{j+1}]$, yield the halved bounds on $[T_0,T_{j+1}]$, in particular at $t_{j+1}$ if
$t_{j+1}<T$. After finitely many steps we obtain the conclusion on $[T_0,T]$. In the rest of this
paragraph and in Steps~3--5, $T$ stands for $T_{j+1}$.

By the corollary these energies are absolutely continuous on $[T_0,T]$, and
Propositions~\ref{prop:CUcell}, \ref{prop:UFforce}, \ref{prop:UFge1} and Lemma~\ref{lem:TWabs} hold
for them for a.e.\ $t$, with the same constants and under the hypotheses at the radii $\rad_k$.
Their values at $T_0$ and at $t_c$, and the bootstrap bounds that are used for the sources with a
square root and in Lemma~\ref{lem:UFfam}, are not larger than at the radii $\rad_k$, while the
level-$0$ quantities $\rho_j$, $\mathfrak s_c$ and $\mathfrak q_c$ do not depend on the radius.
Therefore Steps~3--5 yield the bounds (U1)--(U3) with half the constants for the energies at the
radii $\rad_k-\iota$, for every such $\iota$. Letting $\iota\downarrow0$ and using monotone
convergence in the sums over the levels as in~(i) of the corollary, we get them at the radii
$\rad_k$. In Steps~3--5 below, and in the proof of Proposition~\ref{prop:UFU2}, $\mathfrak A_c$, $
\mathfrak B_c$, $\mathfrak B_{c,\ge1}$, $\mathfrak M_c$ and $\mathfrak D_c$ stand for the quantities at the
radii $\rad_k-\iota$.

\emph{Step 3: all cells.} We apply Proposition~\ref{prop:CUcell} to the joint energy $
\mathfrak B_c=\mathfrak A_c+K_{\mathfrak M}\mathfrak M_c$, and Lemma~\ref{lem:TWabs}, which contains
Lemma~\ref{lem:analytic}, for each $y$-order, with the radius decreasing at the rate of
Remark~\ref{rem:UFradius}. The one-level losses and the spread losses are absorbed by the radius
term.

\emph{(3a) The Gronwall rates.} The source of Proposition~\ref{prop:UFforce}(g) is $\le
Ct^{-2}T_0^{1/2}\mathfrak B_c$, a Gronwall term. The total Gronwall rate is the rate $a(t)$ of the
proof of Proposition~\ref{prop:UFU2}, including the rate $C_0'\rad_0^{-1}t^{-2}$ of the
absorption of the one-level losses by the radius, with $a_1$ containing $C_{\rm
tr}/(\rad_0^2C_2)$ by requirement~\ref{rq:CUd}, and with the term $CT_0^{1/2}t^{-2}$ from~(g).
By requirements~\ref{rq:UFd} and~\ref{rq:TWb}, the integral of the total rate is $\le\log2$. Thus the Gronwall factor
is $\le2$, once, for all these rates together.

\emph{(3b) The square-root sources.} We write out the forcing by
parts~(a)--(d) of the forcing proposition, and the macroscopic hierarchy by
parts~(f)--(g) of the same proposition, applied to the joint energy $\mathfrak B_c$. We use the
bootstrap bound $\mathfrak B_c\le2c_\chi C_U+2\mu_\sharp\le3c_\chi C_U$.
\begin{itemize}
\item The reference forcing~(a) is a square-root source $\le Ct^{-2}F_*\mathfrak B_c^{1/2}$. Its
time integral is $\le CT_0^{-1}F_*(c_\chi C_U)^{1/2}\le c_\chi/8$, by requirement~\ref{rq:UFd}.
\item The constant source of the shifted Leibniz terms of part~(d) is, for
the full energy, a square-root source $\le Ct^{-2}\rad_0^{-1}(C_{\rm
mom}+P_1^{1/2})P_1^{1/2}\mathfrak A_c^{1/2}$. With the Gronwall factor, its contribution is $\le
CT_0^{-1}\rad_0^{-1}(C_{\rm mom}+P_1^{1/2})P_1^{1/2}(c_\chi C_U)^{1/2}\le c_\chi/8$, by
requirement~\ref{rq:TWb}.
\end{itemize}

\emph{(3c) The crude remainder.} We have
\begin{equation}\label{eq:boot-rc}
\int\mathfrak q_c\dd t\le C(K)c_\chi c_\lambda^{-1}S_*\Big(\frac{\rad_0}{R_1}+\delta_y+\delta_U\Big)^2\le\frac{\mu_\flat}{32},
\end{equation}
because each crude leaf spends effective time $\le S_*$ in the crude regime, by
requirement~\ref{rq:crude}.

\emph{(3d) Conclusion.} Hence we have
\begin{equation}\label{eq:boot-U2}
\mathfrak A_c(t)\le2\mathfrak B_c(T_0)+C(K)\Big(S_\sharp\Big(\frac{\rad_0}{R_1}+\delta_y+\delta_U\Big)^2+T_0^{-1}\Big)+\frac{c_\chi}4\le
c_\chi C_U ,
\end{equation}
where $C(K)T_0^{-1}\le c_\sharp/4$ by requirement~\ref{rq:WcrudeT}. We check this term by term.
\begin{itemize}
\item $2\mathfrak B_c(T_0)\le2c_\chi C_DB+\mu_\sharp/8$ by Step~0 (the part beyond $
\mathfrak A_c(T_0)$, $K_{\mathfrak M}\mathfrak M_c(T_0)$, vanishes, and $\mu_\sharp/8$ is a margin).
\item The crude remainder, bounded with $S_*$ in \eqref{eq:boot-rc} and doubled by the Gronwall
factor, is $\le\mu_\flat/16$.
\item $C(K)T_0^{-1}\le c_\sharp/4\le c_\chi/4$, since $c_\sharp\le c_\chi$.
\item The last term is $c_\chi/4$.
\end{itemize}
Since $c_\chi C_U-2c_\chi C_DB=c_\chi$ and $\mu_\sharp/8+\mu_\flat/16+c_\chi/4+c_\chi/4<c_\chi$, so
that $2c_\chi C_DB+\mu_\sharp/8+\mu_\flat/16+c_\chi/4+c_\chi/4\le c_\chi C_U$, this is the first
part of (U2) with half the constant. For the parts $\ge1$, Proposition~\ref{prop:UFU2} implies
$\mathfrak B_{c,\ge1}\le\mu_\sharp$, which is the second part of (U2) with half the constant. In
both, we apply the requirements of Appendix~\ref{app:const} with the threshold
$\mu_\flat=e^{-2C_1^US_\sharp-4C_1^B}c_\sharp/(54e^5)\le\mu_\sharp$, which involves only constants
chosen before $\rad_0$. For the crude remainder this is the requirement used in~(3c).

\emph{Step 4: entered cells, i.e.\ (U3).} Let $t_c$ be the time with $s_{\xi_c}(t_c)=s_\psi$. For
$t<t_c$, $\psi(s_{\xi_c}(t))=0$, so we only need to bound $\mathfrak A_c(t)$ for $t\ge t_c$. We put
$s_0:=s_\psi\in[S_*,S_b]$, and $a=2L'd_*+\epsilon_0/2$ as in the proof of Lemma~\ref{lem:UFnear}, so
that $C_1^US_*a\le\frac5{16}$.

The center leaf itself has $\omega(s_0)\le e^{2C_1^B(S_b-S_*)}\omega(S_b)\le2\,\omega(S_b)$,
since $2C_1^B(S_b-S_*)\le2C_1^US_*\cdot2L'd_*e^{2L'd_*}\le\frac12e^{1/8}$, but below we need the
bound for the near-region leaves. The ramp of $\psi$ stays out of the crude regime: otherwise these
leaves would carry $\omega_\xi/\omega(S_b)\ge e^{2C_1^U(S_*-s_0)}$, and $C'_{\rm cl}$ could not be
structural.

\emph{(4a) The value at $t_c$.} We split $\mathfrak A_c(t_c)$ into the near region $\{|\xi-\xi_c|\le
d_*\}$ and the far region.
\begin{itemize}
\item \emph{Near region.} By that lemma, every cube $Q$ that meets it has
$S_*-\frac12\le S_*e^{-a}\le s<S_b+1\le3S_\sharp$ on $Q^*$ at $t_c$. Consequently the argument of
the first case of Step~1 applies to $Q$, although $Q$ need not contain a leaf with $s\ge S_*$.
From the entry corollary on $Q^*$ and $\mathfrak a_{Q^*,\ge1}\le c_\sharp$, we obtain an undamped
cube energy $\le28\delta_1^2$. Note that no bound (U3) of another cell is used. The damping there is
$\le\omega(S_*e^{-a})=e^{2C_1^U(S_*-S_*e^{-a})+2C_1^B(S_b-S_*)}\omega(S_b)\le
e^{2C_1^US_*a+\frac12e^{1/8}}\omega(S_b)\le e^{3/2}\omega(S_b)$, since $2C_1^US_*a\le\frac58$. Hence
the near region contributes $\le Cc_\chi\,\omega(S_b)\delta_1^2$.
\item \emph{Far region.} By (U2) and the comparability on cubes, the far region contributes
\begin{equation}\label{eq:boot-far}
\le\sum_{|\xi_j-\xi_c|>d_*}e^{-2\jap{(\xi_j-\xi_c)/\rad_0}}e^{6}\mathfrak A_{\xi_j}\le C_{\rm
far}C_Ue^{-d_*/\rad_0}\le\omega(S_b)\delta_1^2 ,
\end{equation}
by the choice of $d_*$ in requirement~\ref{rq:CUc}, since $\omega(S_b)\ge e^{-2C_1^US_b}\ge
e^{-2C_1^US_\sharp}$. Here we used
\begin{equation}\label{eq:boot-farsum}\sum_{|\xi_j-\xi_c|>d_*}e^{-2\jap{(\xi_j-\xi_c)/\rad_0}}\le
C(d_*/\rad_0)^2e^{-2d_*/\rad_0}.\end{equation}
\end{itemize}
\emph{(4b) The evolution after $t_c$.} For $t\ge t_c$, we apply the two estimates of Step~3
from $t_c$. There are three contributions besides the initial
value. The additive term $C(K)T_0^{-1}$ of Step~3 comes from Gronwall rates, namely the frame terms
and the radius on crude leaves, which belong to $a_1$ and not to the sources. In any case, by
the requirement used in Step~3 we have $C(K)T_0^{-1}\le
e^{-2C_1^US_\sharp}c_\sharp/4\le\omega(S_b)\delta_1^2/16$, since $\omega(S_b)\ge
e^{-2C_1^US_\sharp}$ and $c_\sharp\le\delta_1^2/4$, so it could also be put into the absolute
constant $C$ below.
\begin{itemize}
\item Leaves of the cell still in the crude regime contribute $\int_{t_c}^\infty
\mathfrak q_c\le\omega(S_b)\delta_1^2$, by requirement~\ref{rq:crude} and $\omega(S_b)\ge
e^{-2C_1^US_\sharp-4C_1^B}$.
\item The forcing, as in~(3b), contributes $(CT_0^{-1}F_*)^2$ and
the constant sources of the shifted Leibniz terms, $\le\big(CT_0^{-1}\rad_0^{-1}(C_{\rm
mom}+P_1^{1/2})P_1^{1/2}\big)^2$. Both are $\le\frac13\omega(S_b)\delta_1^2$, the first by the
$Y_{\mathfrak b}$-item of requirement~\ref{rq:UFd}, since
$\mu_\flat/32\le\frac13e^{-2C_1^US_\sharp-4C_1^B}\delta_1^2$, and the second by the second entry of
the minimum in requirement~\ref{rq:TWb}. Both are read with $\omega(S_b)\ge
e^{-2C_1^US_\sharp-4C_1^B}$.
\item We run Step~4 for the joint energy $\mathfrak B_c$, because of the top part of~(c) in
the same proposition. At $t_c$, $K_{\mathfrak M}\mathfrak M_c(t_c)\le
\mathfrak B_{c,\ge1}(t_c)\le2\mu_\sharp\le\omega(S_b)\delta_1^2/(108e^5)$ by (U2), since
$c_\sharp\le\delta_1^2/4$ and $\omega(S_b)\ge e^{-2H(S_b+1)}$. The source of~(g) contributes $\le
CT_0^{-1}c_\chi C_U\le\frac13\omega(S_b)\delta_1^2$ by the same requirement.
\end{itemize}
\emph{(4c) Conclusion.} Enlarging the absolute constant $C$ so that it controls the small
contribution of $K_{\mathfrak M}\mathfrak M_c$, we get
\begin{equation}\label{eq:boot-U3}
\mathfrak A_c(t)\le2\big(Cc_\chi+3\big)\omega(S_b)\delta_1^2 .
\end{equation}
With $C_{\rm cl}':=2(Cc_\chi+3)$, which is structural, this is (U3) with half the
constant, since $\psi\le1$.

\emph{Step 5: (U1).} By Remark~\ref{rem:UFradius}, the joint energy controls the macroscopic cell
energies, and we have $\mathfrak M_c\le C(c_\chi C_U+F_*^2)T_0^{-1}+4\mathfrak M_c(T_0)$. By
Lemma~\ref{lem:CUcomp}, we deduce that $\sup_Q\mathfrak m_{Q^*}\le\delta_U^2$ with $\delta_U$ as
defined there. This is (U1) with half the constant (in the square root). This completes the proof.
\end{proof}

\begin{cor}[Theorem~\ref{thm:main}, with the uniformly local bootstrap]\label{cor:CUmain}
Under the hypotheses of Theorem~\ref{thm:main}, and with the constants of Theorem~\ref{thm:CUboot},
its conclusions hold, with these readings:
\begin{itemize}
\item in (i) the radius is $\rad_0/4$ and the bound is $\tilde w$-weighted, so the factorial
bound on $\partial^\beta_v$ holds up to a factor $(|\beta|+2)^2$, which does not affect
real-analyticity;
\item in (ii) the leafwise bound holds pointwise in $\xi$: by Corollary~\ref{cor:Aleafentry} on the
leaves with $s_\xi\in[S_*,3S_\sharp]$, and by (U3) and Lemma~\ref{lem:CUsob} on the others, as in
Step~1 of the proof of the bootstrap theorem, with $\delta_1^2$ replaced by $C\delta_b^2$.
\end{itemize}
\end{cor}

\begin{proof}
We run a continuation argument with Theorem~\ref{thm:CUboot}.

\emph{Step 1: the continuation quantities.} Put
\begin{equation}\label{eq:Phi23}
\Phi_2:=\max\Big\{\frac{\sup_c\mathfrak A_c}{2c_\chi C_U},\frac{\sup_c
\mathfrak B_{c,\ge1}}{2\mu_\sharp}\Big\},\qquad\Phi_3:=\frac{\sup_c\psi(s_{\xi_c})\mathfrak A_c}{2C'_{\rm
cl}\omega(S_b)\delta_1^2}.
\end{equation}
Then $\Phi_2\le1$ and $\Phi_3\le1$ at a time $t$ are (U2) and (U3) at $t$. The other two bootstrap
bounds are not part of the continuation. Indeed, by Lemma~\ref{lem:CUH2}(i), (U1) at $t$ follows
from (U2) at $t$, and by Lemma~\ref{lem:BTU4}, (U4) on $[T_0,T]$ follows from (U1)--(U3) on
$[T_0,T]$, without any time regularity of the (G3)--(G4) norms.

\emph{Step 2: local existence after $T_0$.} A \emph{classical solution in the cell norms} on
$[T_0,T')$ is a classical solution $g\ge0$ such that, for every $T<T'$,
$\sup_{[T_0,T]}\sup_c\mathfrak B^{(1)}_c<\infty$, $\sup_{[T_0,T]}\sup_Q
\mathfrak m_{Q^*}\le\delta_U^2$, and (G1) and (F3) hold on $[T_0,T]$, with constants that may depend on $T$.
For such a solution, Lemma~\ref{lem:LCdata}, applied at every $t\in[T_0,T]$ with~(H) obtained as in
point~(3) of the proof of Proposition~\ref{prop:CULW}, implies $|\partial^{\mathbf k}_{y,v}g|\le
C_T\Gamma_{q_0/4}$ on $[T_0,T]$ for $|\mathbf k|\le2$. By part~(a) of that proposition, the
solution it gives is a classical solution in the cell norms on its
interval. By Lemma~\ref{lem:LWuniq} two classical solutions in the cell norms with the same value at
some time coincide on their common interval.

By Step~0 of the proof of the theorem, (U1)--(U4) hold at $T_0$ in the halved form. Then
Lemma~\ref{lem:CUH2}(ii), with $T=T_0$, implies (H1$'$)--(H2$'$) at $t_*=T_0$. Indeed, $
\mathfrak B^{(1)}_c(T_0)=\mathfrak A_c(T_0)\le c_\chi C_DB\le M$ and $
\mathfrak B^{(1)}_{c,\ge1}(T_0)\le\mu_\flat/32\le m_\sharp$, since $\omega_\xi(T_0)=1$ and $
\mathfrak M_c(T_0)=\mathfrak m_Q(T_0)=0$, and the closeness clause of (H2$'$) is void, since
$s_\xi(T_0)=0<S_b-\frac14$. Applying part~(a) of the local proposition at $t_*=T_0$ to the value
$g(T_0)$ of the solution of Proposition~\ref{prop:disp}, we obtain a classical solution in the cell
norms on $[T_0,T_0+\tau]$. As at $t_1$ in the proof of that proposition, paragraph
\emph{The limit}, both pieces are classical up to $T_0$, with Gaussian bounds on their derivatives
of order $\le2$ as above, so their one-sided time derivatives at $T_0$ are the
right side of the equation at $g(T_0)$, and together they form one classical solution on
$[0,T_0+\tau]$. In particular, the maximal time $T_{\max}$ below satisfies $T_{\max}>T_0$.

\emph{Step 3: continuation.} Let $T_{\max}$ be the maximal time of existence of the classical
solution in the cell norms, and let $\mathcal T$ be the set of $T\in[T_0,T_{\max})$ with
$\Phi_2,\Phi_3\le1$ on $[T_0,T]$. Then (U1)--(U4) hold on $[T_0,T]$ for every $T\in\mathcal T$.
\begin{itemize}
\item \emph{$\mathcal T$ is nonempty.} It contains $T_0$ by Step~0 of the proof of the theorem.
\item \emph{$\mathcal T$ is closed in $[T_0,T_{\max})$.} Let $T_n\in\mathcal T$ increase to
$T<T_{\max}$. Then $\Phi_i\le1$ on $[T_0,T)$, and $\Phi_i(T)\le\liminf_n\Phi_i(T_n)\le1$ by the
lower semicontinuity of Lemma~\ref{lem:CUsemi}.
\item \emph{$\mathcal T$ is open in $[T_0,T_{\max})$.} If $T\in\mathcal T$, by
Lemma~\ref{lem:CUright} there is $h>0$ with $T+h<T_{\max}$ and $\Phi_2,\Phi_3\le1$ on $[T,T+h]$, so
$T+h\in\mathcal T$.
\end{itemize}
Therefore $\mathcal T=[T_0,T_{\max})$.

\emph{Step 4: $T_{\max}=\infty$.} Suppose $T_{\max}<\infty$. By Lemma~\ref{lem:CUH2}(ii), with
$\tau_\flat$ as defined below, the hypotheses (H1$'$)--(H2$'$) of the local proposition hold
at $t_*:=\max\{T_0,T_{\max}-\tau_\flat/2\}$ with the constants $M$, $m_\sharp$, $K_{\rm all}$ of
that lemma, which do not depend on $t_*$, and with $\underline\rad=\rad_0/2$
and $q_*\ge q_0/2$. Put
\begin{equation}\label{eq:tauflat}
\tau_\flat:=\inf_{q_*\in[q_0/2,q_0]}\tau(M,K_{\rm all},q_*)>0;
\end{equation}
here every $q_*$-dependent constant of the proposition (through $C_*$, $\delta\le q_*/4$,
$\eta_0$, $C(q,L)$, $L_3$, $G_*$ and the weight loss $32q^2|z|^2$) is bounded on this range. Thus
part~(a) of the proposition extends the solution to $[T_0,t_*+\tau_\flat]$, with
$t_*+\tau_\flat>T_{\max}$, which is a contradiction. Hence $T_{\max}=\infty$. The argument uses no
continuity of the bootstrap quantities in $t$, and no estimate on an interval beyond $T$ other than
the local proposition, whose hypotheses are taken at the time $T$ itself. The conclusions
(iii) and (iv) then follow as in Section~\ref{sec:proofmain}.
\end{proof}

\subsection{The forcing and the parts of order \texorpdfstring{$\ge1$}{>=1} in
cells}\label{ssec:UF}

Proposition~\ref{prop:UFU2} closes the bounds for the parts of order $\ge1$ of the cell energies,
and Proposition~\ref{prop:UFforce} estimates the forcing. The cell inequality of
Proposition~\ref{prop:CUcell} is stated for the full cell energy. For the parts with
$(\alpha,\beta)\ne0$ the background terms on the leaves with $s_\xi\ge S_*$ behave differently. In
the second configuration of the background variant of Proposition~\ref{prop:CUtri}, the supremum
falls on the dissipation of the background $f_\xi$ itself, at level $(0,0)$. In the inequality for
the full energy this dissipation is on the left side, and it is absorbed. In the inequality for the
parts $\ge1$ it is not, and by Young's inequality we are left with a term of the form $\big(\sup
\mathfrak a_{\ge1}\big)\mathfrak D^{(0,0)}_c$, whose time integral is not small. A leafwise Gronwall argument
does not transfer to cells either: the rate $\kappa\lambda_\xi\|f_\xi\|^2_{Y_2}$ depends on the
leaf, and $\int\sup_\xi\kappa\lambda_\xi\|f_\xi\|^2_{Y_2}\dd t$ is not finite in general, since the
front of entering leaves moves out forever.

We treat this term as a rate that multiplies the parts $\ge1$ themselves. By
Lemma~\ref{lem:UFlevel0}, a leafwise estimate at level $0$ alone, its time integral on each cube is
$O(\delta_1^2)$. Since no hierarchy is involved, the obstruction of
Remark~\ref{rem:ACtr} does not arise. In Proposition~\ref{prop:UFU2}, an integral Gronwall argument
over the lattice of cubes then closes with \emph{one} threshold $c_\sharp$ for all cells.
The forcing proposition writes out the forcing (T12) in the analytic sum, which requires the
graded weights of Remark~\ref{rem:UFgrade}, and the macroscopic hierarchy jointly with the micro
energies at $y$-order~$0$.

We write
\begin{equation}
\Sigma(t):=\{\xi:s_\xi(t)\ge S_*\},\qquad c_{\rm cl}:=\frac{c_W}{4C_W},
\end{equation}
so that $c_{\rm cl}=c'_{\min}$ is the coercivity constant of Proposition~\ref{prop:Bband}, which
holds on $\Sigma$ and on the buffer leaves by Remark~\ref{rem:Bbandrange}, and $c_{\min}\le c_{\rm
buf}=c_{\rm cl}/2$. We assume $C_1^U\ge\max\{C_1^B,1\}$, as required in
Appendix~\ref{app:const}.

\subsubsection{The forcing in cells}\label{sssec:UFforce}

We use the following terminology.
\begin{itemize}
\item $\Pi_0$ is the projection onto the kernel of the linearized operator along the moment-free
functions, $\Pi_0h=\sum_k\mathrm{mom}_k(h)e_k$, where $e_k\in\ker\mathcal L$ are the Gaussian
polynomials dual to the moments $1,\hat z_i,|\hat z|^2$. It commutes with $\mathcal L$ and
$e^{\tau\mathcal L}$, and $\|\Pi_0h\|_{X(\hat m)}\le C|\mathrm{mom}(h)|$.
\item $\hat F_{\alpha\beta}=(I-\Pi_0)T_p^{-1}(\partial^\alpha_y\partial^\beta_v\partial_t
\mathcal M[g])^{(x)}$ is the micro forcing.
\item A \emph{shifted Leibniz term} is a term of the Leibniz expansion in which $\xi$-derivatives
fall on a coefficient while the factor of highest level is not differentiated to the top.
\item A \emph{spread loss} is a term bounded by the right side of Lemma~\ref{lem:TWcalc}(iv) times
$t^{-2}$. By Lemma~\ref{lem:TWabs} below, such terms are controlled by a decreasing radius, which
adds a Gronwall rate proportional to their constant.
\item The norms $E_Q,S_Q$ of the same lemma are the cube energy and the cube supremum of
the analytic sums.
\item The constant $C_{\rm mom}:=C(1+B_{\rm cr}^{1/2})$, with $B_{\rm cr}$ of
Lemma~\ref{lem:BTsob}, bounds the moments of order $\le3$ of the normalized leaves through the
energies, by part~(ii) of that lemma and Lemma~\ref{lem:TWmom} below. It is fixed before
$T_0$.
\item The constant of the top part is $C_0^F\le C(1+C_{\rm mom}+P_1^{1/2})$, where $C$ depends on
$(L_1,R_1,\theta_\pm,q_0,N_y)$ through the Maxwellian coefficients of Lemma~\ref{lem:TWcoef} and
($\Theta$1) in the proof below. It is also fixed before $T_0$.
\item Every solution-dependent constant of Proposition~\ref{prop:UFforce} ($P_1$, $C_{\rm mom}$,
$C_0^F$, $C_1^F$, $C_1^{\rm mac}$, $\mathfrak p_*$, $a_F$, $F_0$) is a nondecreasing function of a
bound $B$ on $\sup_Q\mathfrak a_{Q^*}$, evaluated at $B=B_{\rm cr}$.
\item $E^{(0)}_n:=E^c_{0,n}+K_{\mathfrak M}\mathfrak M^c_n$.
\end{itemize}
The forcing term is the term (T12) in the numbering of Subsection~\ref{sssec:bandproof}.

\paragraph{The forcing on a leaf.} Let $\Phi(t,\xi,z):=M_{p_g(t,\xi)}(z)$. Then $
\mathcal M[g](t,y,v)=\Phi(t,v+y/t,y)$, and by Lemma~\ref{lem:Meq}, we have
\begin{equation}
\partial_t\mathcal M[g](t,y,v)=t^{-2}\,\Theta\big(t,v+y/t,y\big),\qquad\Theta:=\sum_k(\partial_{p_k}M)_{p_g(t,\xi)}(z)\big[
\mathfrak e_k(t,\xi)-z\cdot\nabla_\xi p_{g,k}(t,\xi)\big],
\end{equation}
with $\mathfrak e:=t^2Dp_g$. Since $(\partial_{p_k}M)_{p_g}
\mathfrak e_k=(\partial_{\varpi_k}M)_{p_g}\,t^2D\varpi_{g,k}$ (note $\mathfrak e_1=\rho_g\,t^2D\ell_g$), only
$t^2D\varpi_g$ enters $\Theta$. By Theorem~\ref{thm:Meuler}, $t^2D\varpi_g$ is an analytic function
of the first derivatives of $\varpi_g$ (with $\ell_g$ entering only through $\nabla\ell_g$) and of
the micro moments $s,h,\nabla s,\nabla h$. It is linear in the first derivatives and in $(\nabla
s,\nabla h)$. For a profile function $U(t,\xi,z)$ and $u(t,y,v)=U(t,v+y/t,y)$, we have
$\partial_vu\leftrightarrow\partial_\xi U$ and
$\partial_yu\leftrightarrow(\partial_z+t^{-1}\partial_\xi)U$, and $u^{(x)}=U(t,x/t,\cdot)$. Hence we
obtain
\begin{equation}\label{eq:UFleaf}
t^2\big(\partial^\alpha_y\partial^\beta_v\partial_t
\mathcal M[g]\big)^{(x)}=\sum_{\alpha'\le\alpha}\tbinom\alpha{\alpha'}t^{-|\alpha'|}\,\partial_z^{\alpha-\alpha'}\partial_\xi^{\beta+\alpha'}\Theta(t,\xi,\cdot)
.
\end{equation}
Let $\Theta_*$ be $\Theta$ with $(p_g,\mathfrak e)$ replaced by $(p_*,\mathfrak e_*)$, where $p_*$
is the reference Euler flow of Lemma~\ref{lem:Mref} and $\mathfrak e_*:=t^2Dp_*$. We put
$\delta\Theta:=\Theta-\Theta_*$.

\begin{lem}[Cauchy estimates for the reference forcing]\label{lem:UFref}
There are $r_{\rm C}=c_{\rm C}R_1/(L_1+1)$, with $c_{\rm C}$ absolute, and $C_{{\rm ref},a}$
depending on $(L,L_1,R_1,\theta_\pm,q_0)$, such that, for every $t\ge T_0$ at which the bound (U1)
holds (only the last step of the proof uses it, and only at the time $t$), real $\xi_0$,
multi-indices $a$ with $|a|\le N_y$ and $b$,
\begin{equation}
\Big(\frac{A_{p_g}}{A_{\rm ref}}\Big)(\xi_0)\,\big\|m_{x_0}\,\partial^a_z\partial^b_\xi\Theta_*(t,\xi_0,\cdot)\big\|_{L^2(\dd
z)}\le C_{{\rm ref},a}\,b!\,r_{\rm C}^{-|b|},\qquad x_0=t\xi_0 .
\end{equation}
\end{lem}

\begin{proof}
By Lemma~\ref{lem:Mref}, $\varpi'_*(t,\cdot)$ is holomorphic and bounded on $\mathbb S_{R_1/2}$,
uniformly in $t$, so $t^2D\varpi_*=\partial_\sigma\varpi_*$ is holomorphic and bounded on
$\mathbb S_{R_1/4}$ by the Cauchy estimates, since by Theorem~\ref{thm:Meuler} it is a polynomial
in $\varpi'_*$, $\nabla\varpi'_*$ and $e^{\vartheta_*}$. In $\Theta_*$ we have
$(\partial_{p_k}M)_{p_*} \mathfrak e_{*k}=(\partial_{\varpi_k}M)_{p_*}\partial_\sigma\varpi_{*k}$.

We take $c_{\rm C}:=1/120$. The polydisc $\Delta(\xi_0,\sqrt3r_{\rm C})$ lies in the ball of radius
$3r_{\rm C}\le R_1/40$ about the real point $\xi_0$, inside $\mathbb S_{R_1/4}$. On $
\mathbb S_{R_1/4}$, by the Cauchy estimate on discs of radius $R_1/4$ and $|\varpi'_*|\le L_1+1$ on
$ \mathbb S_{R_1/2}$, we have $|\partial_ey_{0*}|+|\partial_e\vartheta_*|\le4(L_1+1)/R_1$
for every unit $e\in\C^3$.

Since $\varpi'_*(\xi_0)$ is real, and the segment $[\xi_0,\xi']$ lies in
that ball, $|\vartheta_*(\xi')-\vartheta_*(\xi_0)|+|y_{0*}(\xi')-y_{0*}(\xi_0)|\le3r_{\rm
C}\cdot4(L_1+1)/R_1=12c_{\rm C}=\frac1{10}$ on the polydisc. In particular
$|\mathrm{Im}\,\vartheta_*|+|\mathrm{Im}\,y_{0*}|\le\frac1{10}$ there. The point $\mathrm{Re}\,\xi'$
lies in the polydisc, so also $|y_{0*}(\mathrm{Re}\,\xi')-y_{0*}(\xi_0)|\le\frac1{10}$.

On $\R^3$, by Corollary~\ref{cor:LWPa}(a4) and Lemma~\ref{lem:Mref}, we have
$\theta_-\le\theta_g(T_0)\le\frac23\theta_+$ and $|\vartheta_*-\vartheta_g(T_0)|\le
C(L_1,R_1)/T_0\le\frac1{10}$, which implies
$e^{-1/10}\theta_-\le\theta_*(\xi_0)\le\frac23e^{1/10}\theta_+$.
Hence, on the polydisc, since $\mathrm{Re}\,e^{w}\ge e^{-|w|}$ for $w\in\C$, $|w|\le1$ (applied to
$w=\vartheta_*(\xi_0)-\vartheta_*(\xi')$), we have
\begin{equation}
\mathrm{Re}\,\theta_*^{-1}=e^{-\mathrm{Re}\,\vartheta_*}\cos(\mathrm{Im}\,\vartheta_*)\ge
e^{-1/10}\theta_*(\xi_0)^{-1}\ge\frac6{5\theta_+},\qquad|\mathrm{Im}\,\theta_*^{-1}|\le\tan\tfrac1{10}\,\mathrm{Re}\,\theta_*^{-1},\qquad|\theta_*^{-1}|\le\frac{e^{1/5}}{\theta_-}
.
\end{equation}
We write $y_{0*}(\xi')=y_{\rm r}+iy_{\rm i}$, with $y_{\rm r},y_{\rm i}\in\R^3$, and $\bar
y:=y_{0*}(\mathrm{Re}\,\xi')$, so $|y_{\rm i}|\le\frac1{10}$ and $|y_{\rm r}-\bar y|\le\frac15$. Let
$\theta_*^{-1}=\alpha+i\beta$, so $\alpha\ge6/(5\theta_+)$, $|\beta|\le\alpha/9$ and $\alpha\le
e^{1/5}/\theta_-$. For real $z$, with the bilinear square and $w:=z-y_{\rm r}$, we have
\begin{equation}
\mathrm{Re}\big[\theta_*^{-1}(z-y_{0*})\cdot(z-y_{0*})\big]=\alpha(|w|^2-|y_{\rm
i}|^2)+2\beta\,w\cdot y_{\rm i}\ge\alpha\big(\tfrac9{10}|w|^2-\tfrac1{20}\big),
\end{equation}
by $2|\beta||y_{\rm i}||w|\le\alpha|w|/45\le\alpha(|w|^2+1)/90$. Since
$|w|^2\ge\frac12|z-\bar y|^2-|y_{\rm r}-\bar y|^2\ge\frac12|z-\bar y|^2-\frac1{25}$, half of it is
$\ge\frac9{40}\alpha|z-\bar y|^2-C\alpha\ge|z-\bar y|^2/(4\theta_+)-C/\theta_-$. Moreover
$|\rho_*(\xi')|\le e^{3r_{\rm C}(L_1+1)}\rho_*(\mathrm{Re}\,\xi')\le
e^{R_1/40}\rho_*(\mathrm{Re}\,\xi')$, by $|\nabla\ell_*|\le L_1+1$ on $\mathbb S_{R_1/2}$, and
$\Theta_*$ is $M_{p_*}$ times a polynomial of degree $\le3$ in $z$ whose coefficients are bounded by
$C(L_1,R_1,\theta_-)$ on the polydisc. Thus, for real $z$, we obtain
\begin{equation}
|\partial^a_z\Theta_*(t,\xi',z)|\le C_{{\rm ref},a}\,\rho_*(\mathrm{Re}\,\xi')\jap
z^{|a|+4}e^{-|z-y_{0*}(\mathrm{Re}\,\xi')|^2/(4\theta_+)} .
\end{equation}
For $\xi'$ in the polydisc, $|\mathrm{Re}\,\xi'-\xi_0|\le3r_{\rm C}\le R_1/40$, so
$\rho_*(\mathrm{Re}\,\xi')\le e^{R_1/40}\rho_*(\xi_0)$ and
$|y_{0*}(\mathrm{Re}\,\xi')-y_{0*}(\xi_0)|\le\frac1{10}$. The first bound holds because
$|\nabla\ell_*|\le L_1+1$ on $\mathbb S_{R_1/2}$ and $3r_{\rm
C}(L_1+1)=R_1/40$. The second holds because $|\nabla y_{0*}|\le4(L_1+1)/R_1$ on $\mathbb S_{R_1/4}$
(shown above) and $12r_{\rm C}(L_1+1)/R_1=\frac1{10}$. The weight is $m_{x_0}(z)=A_{\rm
ref}(\xi_0-z/t)e^{q|z|^2}$, with $A_{\rm ref}(\xi_0-z/t)\le A_{\rm ref}(\xi_0)e^{2L|z|/t}$, since
$|\nabla\log A_{\rm ref}|\le L$ by~(a4) again, and $q\theta_+\le\frac1{64}$ by the
choice of $q_0$ in Appendix~\ref{app:const}. Hence the $L^2(m_{x_0}^2\dd z)$-norm of
$\partial^a_z\Theta_*(t,\xi',\cdot)$ is $\le C_{{\rm ref},a}A_{\rm ref}(\xi_0)\rho_*(\xi_0)$,
uniformly on the polydisc. The map $\xi'\mapsto\partial^a_z\Theta_*(t,\xi',\cdot)$ is holomorphic
with values in this Hilbert space, by Morera's theorem. Applying the Cauchy estimates for
Hilbert-space-valued holomorphic functions, we get the bound with the factor $A_{\rm ref}\rho_*$.
Finally $(A_{p_g}/A_{\rm ref})A_{\rm ref}\rho_*=\theta_g^{3/4}\rho_*/\rho_g\le C$, by (U1). The
preceding bounds do not use (U1) and hold for every $t\ge T_0$, since
$\mathrm{Re}\,\theta_*^{-1}\ge6/(5\theta_+)$ on $\Delta(\xi_0,\sqrt3r_{\rm C})$, and the
coefficients of the Gaussian polynomials in $\varpi_*$ have Cauchy bounds at the radius $r_{\rm C}$.
\end{proof}

We write $E^c_{0,n}$ for the part of $E^c_n$ with $\alpha=0$, and we define the \emph{damped
macroscopic cell energies}
\begin{equation}
\mathfrak M^c_n:=\rad_0^{-3}\!\int\chi_c^2\,\omega_\xi\!\sum_{|\beta|=n}|\partial^\beta_\xi\delta\varpi(t,\xi)|^2\dd\xi,\qquad
\mathfrak M_c:=\sum_nw_n(\rad)^2\mathfrak M^c_n ,
\end{equation}
with $\delta\varpi=\varpi_g-\varpi_*$ and $\varpi=(\ell,y_0,\vartheta)$, as in
Subsection~\ref{ssec:macro}, and, for a constant $K_{\mathfrak M}\ge1$ fixed below, the joint
energies $\mathfrak B_c:=\mathfrak A_c+K_{\mathfrak M}\mathfrak M_c$ and $
\mathfrak B_{c,\ge1}:=\mathfrak A_{c,\ge1}+K_{\mathfrak M}\mathfrak M_c$. Let $\mathfrak A^{00}_c$ be the
level-$(0,0)$ part of $\mathfrak A_c$.

\begin{prop}[The forcing (T12) in cells]\label{prop:UFforce}
Assume (U1), (U2), (U3), (U4) on $[T_0,T]$, requirements~\ref{rq:CUa}, \ref{rq:CUb}
and~\ref{rq:TWa}, $\rad_0\le r_{\rm C}/2$ and $\rad_0T_0\ge1$. Let $\mathfrak f^c$ be the
contribution of (T12), local and nonlocal, to $\frac{d}{dt}\mathfrak A_c$, and $
\mathfrak f^c_{\ge1}$ its part with $(\alpha,\beta)\ne0$. Then:
\begin{enumerate}[label=(\alph*)]
\item \emph{Reference part.} There is $F_*=F_*(L_1,R_1,N_y,\theta_\pm,q_0)$ such that the part of
$\mathfrak f^c$ from $\Theta_*$ is $\le Ct^{-2}\mathfrak A_c^{1/2}F_*$. Its part with
$(\alpha,\beta)\ne0$ is
\begin{equation}\label{eq:UF-a}
\le Ct^{-2}\mathfrak A_{c,\ge1}^{1/2}F_*(\rad_0/r_{\rm C}+\delta_y) .
\end{equation}

\item \emph{Kernel part.} The part from $\Pi_0$ of the full forcing is $\le
C_0t^{-2}\sum_nw_n^2\sqrt{E^c_nE^c_{n+1}}$, separately for each $y$-order.
\item \emph{Top part.} The part from $\delta\Theta$ in which all $\xi$-derivatives fall on a first
derivative of $\delta\varpi$ or on $(\nabla s,\nabla h)$, at $y$-order $0$, is
\begin{equation}
\le C_0^Ft^{-2}\sum_nw_n^2\sqrt{E^c_{0,n}}\Big(\sqrt{E^c_{0,n+1}}+\sqrt{\mathfrak M^c_{n+1}}\Big) .
\end{equation}
The top part at $y$-orders $\ge1$ belongs to~(d), by ($\Theta$5) in the proof.
\item \emph{Remainder.} The rest is $\le t^{-2}a_F\,\mathfrak B_c$ in $\frac{d}{dt}
\mathfrak A_c$. In $\frac{d}{dt}\mathfrak A_{c,\ge1}$ it is $\le t^{-2}a_F\mathfrak B_{c,\ge1}+\mathfrak b_c$,
with
\begin{equation}
\mathfrak b_c:=Ct^{-2}\mathfrak B_{c,\ge1}^{1/2}\Big(F_*\big(\tfrac{\rad_0}{r_{\rm
C}}+\delta_y\big)+C\big(\mathfrak A^{00}_c\big)^{1/2}\Big).
\end{equation}
Here $a_F=C\big(1+\mathfrak p_*+P_1^{1/2}\big)^3\big(1+\rad_0^{-1}\big)$, and the remainder does
not contain the shifted Leibniz terms, in which a factor of lower level carries $\xi$-derivatives of
the coefficients. These are spread losses with constant $C_1^F:=C(C_{\rm mom}+P_1^{1/2})$, absorbed
by the radius, which adds the Gronwall rate $4C_1^FT_0^{1/2}\rad_0^{-1}t^{-2}$
(Lemma~\ref{lem:TWabs} and ($\Theta$2) in the proof). In addition,
they produce Gronwall terms with rate $Ct^{-2}\rad_0^{-1}(C_{\rm mom}+P_1^{1/2})$ and a constant
source
\begin{equation}\label{eq:UF-d}
Ct^{-2}\rad_0^{-1}(C_{\rm mom}+P_1^{1/2})P_1^{1/2}\mathfrak A_c^{1/2}
\end{equation}
($\mathfrak B_{c,\ge1}^{1/2}$ in place of $\mathfrak A_c^{1/2}$ for the parts $\ge1$), by
($\Theta$2) and ($\Theta$5) in the proof below.
\item \emph{Level $0$.} On each leaf, the micro forcing at level $(0,0)$ satisfies
\begin{equation}\label{eq:UF-e}
t^2\|\hat F_{00,\xi}\|_{X(\hat m)}\le F_0:=C(L_1,R_1)(1+\mathfrak p_*),
\end{equation}
and its kernel part is $\le Ct^{-2}\|\varphi_{0e}\|_X$.
\item \emph{The macroscopic hierarchy.} For each $n$,
\begin{equation}
\frac{d}{dt}\mathfrak M^c_n\le(C_0+C_0^F)t^{-2}\sqrt{\mathfrak M^c_n}\Big(\sqrt{
\mathfrak M^c_{n+1}}+C\sqrt{E^c_{0,n+1}}\Big)+t^{-2}\sqrt{\mathfrak M^c_n}\,\mathfrak S^c_n+\mathfrak L^c_n,
\end{equation}
where the spread losses satisfy, for every $\upsilon\ge1$,
\begin{equation}\sum_nw_n^2\mathfrak L^c_n\le C_1^{\rm
mac}\rad^{-1}t^{-2}\sum_n(n+1)\big((1+\upsilon)\,b_n^2+\upsilon^{-1}e_n^2\big),\end{equation}
where $b_n=w_n(\mathfrak M^c_n)^{1/2}$, $e_n=w_n(E^c_{0,n})^{1/2}$ $(n\ge1)$ and $e_0:=0$, and the
rest satisfies
\begin{equation}\label{eq:UF-f}
\sum_nw_n^2\sqrt{\mathfrak M^c_n}\,\mathfrak S^c_n\le C'_\varpi\mathfrak M_c+C
\mathfrak M_c^{1/2}\mathfrak A_c^{1/2},
\end{equation}
with $C'_\varpi=C(L_1,R_1,\theta_\pm,q_0,N_y)$. (The one-level constant contains $C_0^F$, because
the coefficients of $\partial\varpi$ in Theorem~\ref{thm:Meuler} contain the micro moments.)
\item \emph{Joint absorption.} Let $K_{\mathfrak M}:=T_0$. For the family $E^c_{0,n}+K_{
\mathfrak M}\mathfrak M^c_n$, all one-level losses of (b), (c), (f), of the transport and of the nonlocal
terms are of the form $C_0't^{-2}\sqrt{E^{(0)}_nE^{(0)}_{n+1}}$, with
$C_0':=2C_0+C_0^F+(C_0+C_0^F)T_0^{1/2}$, where $C_0$ and $C_0^F$ are $C+2$ times the constants
obtained in the proofs of (b), (c) and~(f), $C$ being the constant inside the bracket of~(f) (Step~7
of the proof), and (b), (c) and~(f) hold with these enlarged constants. The sources of (f)
contribute
\begin{equation}\label{eq:UF-g}
\le Ct^{-2}T_0^{1/2}\big(K_{\mathfrak M}\mathfrak M_c\big)^{1/2}\mathfrak A_c^{1/2}
\end{equation}
to $\frac{d}{dt}\mathfrak B_c$.
\end{enumerate}
\end{prop}

\begin{remark}[Proposition~\ref{prop:UFforce} in the local theory]\label{rem:UFinputs}
The proof of the proposition uses the solution only through (U1), a bound $B$ on
$\sup_Q\mathfrak a_{Q^*}$ (from (U2)), and the structural conditions
$\rad_0\le r_{\rm C}/2$, $\rad_0t\ge1$, requirements~\ref{rq:CUa}, \ref{rq:CUb}
and~\ref{rq:TWa}, $\omega\le1$, $D\omega\le0$ and the comparability of
Lemma~\ref{lem:CUcomp}(i)--(iii) on cubes (which uses that $H$ is $C_1^U$-Lipschitz), the frame
bounds (F3) of Subsection~\ref{ssec:inputs}, $C\delta_U\le\rad_0$, graded radii in
$[\rad_0/4,\rad_0]$, $t\ge T_0\ge T_{\rm E}(L_1,R_1)$ of Lemma~\ref{lem:Mref}, and $
\mathfrak p_*$, the bound of Corollary~\ref{cor:FMuse}(d), valid on all leaves. Its solution-dependent
constants are nondecreasing functions of $B$. Hence it applies without
change in the local theory of Subsection~\ref{ssec:cellsproofs}, with $B=B_{\rm loc}\le B_{\rm cr}$,
as in point~(1) of the proof of Proposition~\ref{prop:CULW}.
\end{remark}

We give the proof below, after three lemmas. It uses the graded weights of
Remark~\ref{rem:UFgrade}, the Cauchy estimates of Lemma~\ref{lem:UFref}, and the calculus on cubes
of Lemma~\ref{lem:TWcalc}. Here $P_1:=C\rad_0^{-2}e^{2H(S_b+1)}c_\chi C_U$ bounds $\sup_\xi
e_\xi[\varphi_{0e}]$, as shown in the proof of Lemma~\ref{lem:UFlevel0}, $F_0$ is the level-$0$
forcing of part~(e), and $C_{\rm mom}$ bounds the moments of the leaves. All three are fixed before
$T_0$.

\subsubsection{Proof of Proposition~\ref{prop:UFforce}}\label{sssec:UFforceproof}

For a family $E_n\ge0$ write $a_n:=w_n(\rad)E_n^{1/2}$.

\begin{lem}[Absorption of spread losses]\label{lem:TWabs}
Let $E_n(t)\ge0$ be absolutely continuous, and suppose that
\begin{equation}
\frac{d}{dt}\sum_nw_n^2E_n\le\sum_n2n\frac{\rad'}{\rad}a_n^2+a(t)\sum_na_n^2+b(t)\,\frac{C_1}{\rad}\sum_n(n+1)\,a_n^2,
\end{equation}
with $a,b\ge0$ integrable. If $\rad'\le-C_1b$ and $\rad\ge\rad_0/4$, then
\begin{equation}
\frac{d}{dt}\sum_nw_n^2E_n\le\Big(a+\frac{4C_1b}{\rad_0}\Big)\sum_nw_n^2E_n .
\end{equation}
The hypothesis holds, in particular, if the loss is
$b(t)\frac{C_1}\rad\sum_n\sum_{j=0}^{n+1}c_j(n+1)a_na_{n+1-j}$ with $c_j\ge0$,
$\sum_jc_j(j+1)\le1$, and also if such a bound holds on each cube with cube-dependent $c^Q_j$ and is
then summed with the cell weight (see the proof). The same holds for a joint family, with $a_n$
replaced by $(a_n^2+\tilde a_n^2)^{1/2}$ on both sides, and for each $y$-order separately, with the
radius $\rad_{|\alpha|}$ of Remark~\ref{rem:UFgrade}, since all $\rad_k$ decrease at the same
rate. Lemma~\ref{lem:analytic} is the case $c_0=1$, $c_j=0$ for $j\ge1$, up to the value of $C_1$.
If the differential inequality is available only for the truncated sums $\sum_{n\le N}$, as in
that lemma, the terms with $j\ge1$ involve only levels $\le N$, the term $j=0$ at
$n=N$ is the boundary term $b\frac{C_1}\rad(N+1)a_Na_{N+1}$, and the conclusion follows as in
Steps~3--4 of its proof, provided $\sum_nw_n^2E_n$ is bounded on bounded
time intervals. In the cube-summed case the lemma holds with $C_1$ replaced by $C_{\rm cube}C_1$,
both in the condition $\rad'\le-C_1b$ and in the conclusion, where $C_{\rm cube}$ comes from
Lemma~\ref{lem:CUcomp} and the overlap of the cubes $Q^*$. In Proposition~\ref{prop:UFforce} this
factor is included in $C_1^F$ and $C_1^{\rm mac}$.
\end{lem}

\begin{proof}
For the particular form, we use $2a_na_i\le a_n^2+a_i^2$. For the first half we get
$\frac12\sum_jc_j\sum_n(n+1)a_n^2\le\frac12\sum_n(n+1)a_n^2$. In the second half we put $i=n+1-j$,
so that $n+1=i+j\le(i+1)(j+1)$, and we obtain
$\frac12\sum_ic_j(j+1)(i+1)a_i^2\le\frac12\sum_i(i+1)a_i^2$.
On a cube, with $c^Q_j$, $x^Q_n$, $y^Q_i$ in place of $c_j$, $a_n$, $a_i$, the same argument yields
$\frac12\sum_n(n+1)(x^Q_n)^2+\frac12\sum_i(i+1)(y^Q_i)^2$; multiplied by the cube factor and by
$\chi_c(\xi_Q)^2\omega_{\xi_Q}$ and summed over $Q$ with Lemma~\ref{lem:CUcomp}, this is $\le
C\rad^{-1}\sum_n(n+1)a_n^2$ in the cell variables. In general the loss is
\begin{equation}
\le\frac{C_1b}{\rad}\sum_n(n+1)a_n^2=\frac{C_1b}{\rad}\sum_nna_n^2+\frac{C_1b}{\rad}\sum_na_n^2
.
\end{equation}
The radius term is $\le-2C_1b\sum_nna_n^2/\rad$, which absorbs the first part. The second is $\le
C_1b\rad^{-1}\sum_na_n^2\le4C_1b\rad_0^{-1}\sum_na_n^2$, since $\rad\ge\rad_0/4$.
\end{proof}

\paragraph{The coefficient spaces.} (Here $\mathcal H_j$ is a space of functions of $z$; $H$
remains the damping exponent.) On a cube $Q=Q_j$ we use the frame of its center,
$p_j:=p_g(t,\xi_j)$, and the space $\mathcal H_j$ of functions of $z$ with the norm
\begin{equation}
\|u\|_{\mathcal H_j}:=\Big(\frac{A_{p_j}}{A_{\rm ref}(\xi_j)}\Big)\big\|\jap
z^{N_0}\,m_{t\xi_j}(z)\,e^{q_0|z|^2}\big(|u|+|\nabla_zu|\big)\big\|_{L^2(\dd z)},
\end{equation}
with $N_0$ structural. By Lemma~\ref{lem:CUcomp}(iii) and Lemma~\ref{lem:TG0}(ii), for every $\xi\in
Q^*$ the norms $\|T^{-1}_{p_g(\xi)}\cdot\|_\flat$, $\|T^{-1}_{p_g(\xi)}\cdot\|_{X(\hat m_\xi)}$ and
$\|T^{-1}_{p_g(\xi)}\cdot\|_{Y_2(\hat m_\xi)}$ are $\le C\|\cdot\|_{\mathcal H_j}$. We write the
leaf parameters in log-variables relative to the center,
\begin{equation}
u_g(\xi):=\big(\ell_g(\xi)-\ell_g(\xi_j),\,y_{0,g}(\xi),\,\vartheta_g(\xi)\big),\qquad
w_0:=u_g(\xi_j),
\end{equation}
and
\begin{equation}\label{eq:ustar}
u_*(\xi):=\big(\ell_*(\xi)-\ell_g(\xi_j),\,y_{0*}(\xi),\,\vartheta_*(\xi)\big)
\end{equation}
for the reference flow. Then $\rho_je^{u_{*,1}}=\rho_*$, so that $\Phi_a(u_*)=\partial^a_zM_{p_*}$
and $\Gamma_\mu(U_*)$ below are the reference coefficients, and $u_g-u_*=\delta\varpi$.

\begin{lem}[Maxwellian coefficients]\label{lem:TWcoef}
For $|a|\le N_y+2$ put $\Phi_a(\xi):=\partial^a_zM_{p_g(t,\xi)}\in \mathcal H_j$. Assume (U1), the
frame bounds (F3) of Subsection~\ref{ssec:inputs} and requirement~\ref{rq:TWa}. Then on every cube
\begin{equation}
S_Q(\Phi_a)+E_Q(\Phi_a)\le C_{\rm Max},\qquad S^{\ge1}_Q(\Phi_a)+E^{\ge1}_Q(\Phi_a)\le C_{\rm
Max}\Big(\rad_0+\frac{\rad_0}{r_{\rm C}}+\delta_U\Big),
\end{equation}
and
\begin{equation}\label{eq:TWcoef-diff}
S_Q(\Phi_a(u_g)-\Phi_a(u_*))+E_Q(\Phi_a(u_g)-\Phi_a(u_*))\le C_{\rm Max}\delta_U .
\end{equation}
The constant $C_{\rm Max}$ depends only on $(L_1,R_1,\theta_\pm,q_0,N_y)$. The same
holds for $z^b\Phi_a$ with $|b|\le1$, and for $T^{-1}_{p_g}$-normalized Gaussian polynomials with
coefficients holomorphic in $u$.
\end{lem}

\begin{proof}
In the variables $u$, $M_p(z)=\rho_je^{u_1}(2\pi e^{u_3})^{-3/2}\exp(-|z-u_2|^2/(2e^{u_3}))$. Hence
$\Phi_a=\rho_j\mathcal N_a(u)$, where
\begin{equation}\label{eq:Na}
\mathcal N_a(u)=\partial^a_z\big[e^{u_1}(2\pi e^{u_3})^{-3/2}e^{-|z-u_2|^2/(2e^{u_3})}\big]
\end{equation}
is holomorphic on the polydisc $\Delta(w_0,r_p)$, $r_p=\frac1{10}$, with values in
$\rho_j^{-1}\mathcal H_j$. Its norm there is bounded.

We write $\theta_j:=\theta_g(\xi_j)\le\theta_++1$ (by (F3)), $\bar y:=\mathrm{Re}\,u_2$, $\tilde
y:=\mathrm{Im}\,u_2$, $|\tilde y|\le\sqrt3r_p=\sqrt3/10$ and $|\bar
y|\le|y_{0,g}(\xi_j)|+\sqrt3r_p\le L+1+\sqrt3/10\le L+2$, by (F3).

On $\Delta(w_0,\frac1{10})$, we have $\mathrm{Re}\,e^{-u_3}\ge\frac9{10}\theta_j^{-1}$ and
$|\mathrm{Im}\,e^{-u_3}|\le\frac19\theta_j^{-1}$, which implies
$\mathrm{Re}[e^{-u_3}(z-u_2)\cdot(z-u_2)]\ge\theta_j^{-1}(\frac45|z-\bar y|^2-C)$, and $|
\mathcal N_a(u)(z)|+|\nabla_z\mathcal N_a(u)(z)|\le
C\jap{z-\bar y}^{|a|+1}e^{-\frac25|z-\bar y|^2/\theta_j}$.

The $\mathcal H_j$-weight is $A_{p_j}\jap z^{N_0}\big(A_{\rm ref}(\xi_j-z/t)/A_{\rm
ref}(\xi_j)\big)e^{(q+q_0)|z|^2}\le A_{p_j}\jap z^{N_0}e^{2L|z|/t}e^{(q+q_0)|z|^2}$, since
$|\nabla\log A_{\rm ref}|\le L$ by Corollary~\ref{cor:LWPa}(a4). With $q\le q_0$ and
$q_0\max\{1,\theta_+\}\le\frac1{64}$ (choice of $q_0$, Appendix~\ref{app:const}), we get
$4(q+q_0)\le1/(4\theta_j)$. Since $|\bar y|\le L+2$, we have $\jap z\le C\jap{z-\bar y}$ and
$e^{4L|z|/t}\le
Ce^{4L|z-\bar y|}$ for $t\ge1$. Thus the integrand of $\|\mathcal N_a(u)\|^2_{\mathcal H_j}$ is $\le
C\jap{z-\bar y}^{2N_0+2|a|+2}e^{4L|z-\bar y|}e^{-|z-\bar y|^2/(2\theta_j)}$, which is integrable
uniformly on the polydisc.

The factor $\rho_jA_{p_j}=\theta_j^{3/4}$ is bounded, by (F3), as in Lemma~\ref{lem:UFref}. We apply
Lemma~\ref{lem:TWcalc}(vi). The quantity $S_Q(u_g-w_0)$ has two pieces.
\begin{itemize}
\item \emph{The reference.} $u_*-u_*(\xi_j)$ has $S_Q\le C(\rad_0+\rad_0/r_{\rm C})$, by
Lemma~\ref{lem:Mref} and~(v) of the calculus lemma: its level $0$ is $\sup_Q|u_*-u_*(\xi_j)|\le
C\rad_0$, and its levels $\ge1$ are $\le C\rad_0/r_{\rm C}$. Moreover,
$u_*(\xi_j)-w_0=-\delta\varpi(\xi_j)$ is $\le C\delta_U$.
\item \emph{The deviation.} $S_Q(\delta\varpi)\le C_{\rm Sob}E_{Q^*}(\delta\varpi)\le C
\mathfrak m_{Q^*}^{1/2}\le C\delta_U$, by~(i) of the calculus lemma and (U1).
\end{itemize}
By requirement~\ref{rq:TWa}, we have
\begin{equation}\label{eq:TWa-bound}
S_Q(u_g-w_0)\le r_p/(16dC_{\rm alg}),
\end{equation}
and the same bound for $u_*+\tau(u_g-u_*)-w_0$, $\tau\in[0,1]$, which has the same two pieces. Then
(vi) implies the first bound. For the second, $\Phi_a(w_0)$ is constant, so its levels $\ge1$
vanish, and (vi) bounds the rest by $C\,S_Q(u_g-w_0)$ or $C\,E_Q(u_g-w_0)$, which are $\le
C(\rad_0+\rad_0/r_{\rm C}+\delta_U)$, the level-$0$ size included. For the third, we write
$\Phi_a(u_g)-\Phi_a(u_*)=\int_0^1D\Phi_a(u_*+\tau(u_g-u_*))[u_g-u_*]\dd\tau$ and use (vi) for
$D\Phi_a$, (iii), and $E_Q(\delta\varpi)\le C\mathfrak m_Q^{1/2}\le C\delta_U$. Here we apply (vi)
on the smaller polydisc $\Delta(w_0,r_p/2)$, where $\|D\Phi_a\|\le2B/r_p$ by the Cauchy estimates.
The condition of~(vi), $S_Q(u_*+\tau(u_g-u_*)-w_0)\le(r_p/2)/(8dC_{\rm alg})$ is~\eqref{eq:TWa-bound} with
$u_*+\tau(u_g-u_*)$ in place of $u_g$, which is why the requirement has $16$ in place of
$8$.
\end{proof}

\begin{lem}[Micro moments]\label{lem:TWmom}
Let
\begin{equation}\label{eq:tildeMk}
\tilde M_k(\xi):=\rho_j^{-1}\int\hat r^{(t\xi)}(z)\,z^{\otimes k}\dd z,
\end{equation}
$k=2,3$. Then:
\begin{enumerate}[label=(\roman*)]
\item $\partial^\beta\tilde M_k=\rho_j^{-1}\int(\partial^\beta_v\hat r)^{(t\xi)}z^{\otimes k}\dd
z$, and $|\partial^\beta\tilde M_k(\xi)|\le C\|\varphi_{0\beta}(\xi)\|_{X(\hat m_\xi)}$ for $\xi\in
Q^*$;
\item the stress and the heat flux of Theorem~\ref{thm:Meuler} are $s=\Psi_s(u_g)\tilde M_2$ and
$h=\Psi_h(u_g)(\tilde M_2,\tilde M_3)$, with $\Psi_s,\Psi_h$ entire and matrix-valued;
\item on every cube, $E_Q(\tilde M_k)\le C\eta^{-1/2}\mathfrak a_{Q^*}^{1/2}$ and
\begin{equation}\label{eq:TWmom-S}
S_Q(\tilde M_k)\le C\eta^{-1/2}(\sup_Q\mathfrak a_{Q^*})^{1/2}\le CB_{\rm
cr}^{1/2}=C\rad_0P_1^{1/2};
\end{equation}
in particular $S_Q(\tilde M_k)\le CP_1^{1/2}$ and $\sup_Q|\tilde M_k|\le4S_Q(\tilde M_k)\le C_{\rm
mom}$.
\end{enumerate}
\end{lem}

\begin{proof}
(i) Since $\partial_v=\partial_\xi$ on leaves by Lemma~\ref{lem:leafgeom}, $\partial^\beta$ commutes
with the $z$-integral (differentiate under the integral sign at fixed $z$). In the frame $p_g(\xi)$,
$\int(\partial^\beta\hat r)^{(x)}z^{\otimes k}=\rho\int\varphi_{0\beta}(y_0+\theta^{1/2}\hat
z)^{\otimes k}\dd\hat z$, and moments of order $\le3$ are bounded by the $X$-norm, as in the proof
of Proposition~\ref{prop:FM}(iii). Finally, $\rho(\xi)/\rho_j\le e^{C\rad_0}$ on $Q^*$. (ii) By
Lemma~\ref{lem:Mfree}, the moments of $\hat r^{(x)}$ against $1,z_i,|z|^2$ vanish. Hence
$s=(\rho\theta)^{-1}\int\hat r(z-y_0)^{\otimes2}=(\rho\theta)^{-1}\int\hat rz^{\otimes2}$, and $h$
is $\rho^{-1}\theta^{-3/2}$ times a linear combination of $\int\hat rz^{\otimes3}$ and
$y_0\otimes\int\hat rz^{\otimes2}$, and we write $\rho^{-1}=\rho_j^{-1}e^{-u_1}$. (iii) By
Lemma~\ref{lem:TWcalc}(i), (i) and the definition of the cube energy, we have $S_Q(\tilde M_k)^2\le
C_{\rm Sob}^2E_{Q^*}(\tilde M_k)^2\le
C\sum_nw_n^2\rad_0^{-3}\int_{Q^*}\sum_{|\beta|=n}\|\varphi_{0\beta}\|_X^2\le C\eta^{-1}
\mathfrak a_{Q^*}$, and $\mathfrak a_{Q^*}\le B_{\rm cr}$ by (U2) and Lemma~\ref{lem:CUcomp}. A bound on the
level $1$ alone, such as $P_1$, would not suffice, since $S_Q$ involves all levels.
\end{proof}

We take $C_{\rm mom}\ge C_{\rm Max}$ and $C_{\rm mom}\ge\sup_Q|\tilde M_k|$, which is possible by
part~(iii) of the lemma, enlarging the constant in its definition. This does not change the
order of the constants.

\paragraph{The structure of $\Theta$.} Let $U:=(\ell,y_0,\vartheta,\tilde M_2,\tilde M_3)$ on a
cube, with $\ell$ entering through $\nabla\ell$ and through $u_1=\ell-\ell_g(\xi_j)$. By
Theorem~\ref{thm:Meuler} and Lemma~\ref{lem:TWmom}(ii), the Euler right-hand sides are finite sums
$\sum_\mu c_\mu(U)\,\partial_{i_\mu}U_{k_\mu}$ with $c_\mu$ entire. In the log-variables,
$T^{-1}_p(\partial_{p_k}M)_pDp_k$ is a Gaussian polynomial of degree $\le2$ times
$(D\ell,Dy_0,D\vartheta)$, and $z\cdot\nabla_\xi p_{g,k}$ is a first derivative of $U$ times a
holomorphic function of $U$. Hence we have
\begin{equation}\label{eq:TWtheta}
\Theta=\sum_{\mu=1}^{\mu_0}\Gamma_\mu(U)\,\partial_{i_\mu}U_{k_\mu},
\end{equation}
where each $\Gamma_\mu$ is an $\mathcal H_j$-valued holomorphic function of the type of
Lemma~\ref{lem:TWcoef}, and $\mu_0$ is absolute. Each $\Gamma_\mu$ is \emph{affine} in $(\tilde
M_2,\tilde M_3)$, with coefficients holomorphic in $u_g$. Indeed, the Euler system is affine in
$(s,h,\nabla s,\nabla h)$, and $s,h$ are linear in $\tilde M_k$. This is important since $S_Q(\tilde
M)\le
CP_1^{1/2}$ is not small. We apply the composition estimate Lemma~\ref{lem:TWcalc}(vi) only in the
variables $u_g$, and the micro moments enter through the products (ii)--(iv) of that lemma. The
reference has $U_*=(\ell_*,y_{0*},\vartheta_*,0,0)$. With $\delta U:=U-U_*$, we have
\begin{equation}\label{eq:TWdtheta}
\Theta-\Theta_*=\sum_\mu\Gamma_\mu(U)\,\partial_{i_\mu}\delta
U_{k_\mu}+\sum_\mu\big(\Gamma_\mu(U)-\Gamma_\mu(U_*)\big)\,\partial_{i_\mu}U_{*k_\mu}.
\end{equation}

\begin{proof}[Proof of Proposition~\ref{prop:UFforce}]
\emph{Step 1: the kernel part, (b).} By Lemma~\ref{lem:leafeq}(ii), the equation for
$\varphi=\varphi_{\alpha\beta}$ has the full forcing
$-T^{-1}_p(\partial^\alpha_y\partial^\beta_v\partial_t\mathcal M)^{(x)}$ and the unprojected
transport $t^{-2}T^{-1}_p[z\cdot(\partial_vh)^{(x)}]$, and by the same lemma the functions
$\varphi$, $D\varphi$, the collision terms and $\mathcal G_p[Dp]\varphi$ are moment-free. Therefore
the kernel part of the forcing equals the kernel part of the transport:
\begin{equation}
\Pi_0T^{-1}_p\big(\partial^\alpha_y\partial^\beta_v\partial_t
\mathcal M\big)^{(x)}=t^{-2}\,\Pi_0\big[(y_0+\theta^{1/2}\hat z)\cdot\varphi_{\alpha,\beta+e}\big].
\end{equation}
Its $X(\hat m)$-norm is $\le Ct^{-2}\|\varphi_{\alpha,\beta+e}\|_X$, because $\Pi_0$ is bounded by
the moments. Summing with Cauchy--Schwarz, we obtain (b), separately for each $y$-order. The local
part of the energy uses the full forcing, as Lemma~\ref{lem:CUtr} requires the unprojected
transport, while the nonlocal part uses the projected transport and the micro forcing
$\hat F_{\alpha\beta}$. Since the nonlocal form is nonnegative and
$\langle\!\langle\varphi\rangle\!\rangle^2_\natural\le
C_\natural\|\varphi\|^2_X$ with the constant $C_\natural$ from (E1), by Cauchy--Schwarz we
get $|\langle\varphi,\hat F\rangle_\natural|\le C_\natural\|\varphi\|_X\|\hat F\|_X$. Thus in both
parts it remains to bound $\|\hat F_{\alpha\beta}\|_X$, with $\|(I-\Pi_0)\cdot\|_X\le C\|\cdot\|_X$.

\emph{Step 2: the reference part, (a).} By~\eqref{eq:UFleaf} and Lemma~\ref{lem:UFref}, we have
\begin{equation}
t^2W_{\alpha\beta}\|\hat F^*_{\alpha\beta}\|_X\le C\sum_{\alpha'\le\alpha}t^{-|\alpha'|}\,\delta_y^{|\alpha|}\,w_{|\beta|}(\rad_0)\,(|\beta|+|\alpha'|)!\,r_{\rm
C}^{-|\beta|-|\alpha'|}.
\end{equation}
Since $(|\beta|+|\alpha'|)!\le|\beta|!(|\beta|+N_y)^{N_y}$ and $tr_{\rm C}\ge1$, this is $\le
C(N_y)\delta_y^{|\alpha|}(\rad_0/r_{\rm C})^{|\beta|}(|\beta|+1)^{5+N_y}$. With $\rad_0\le
r_{\rm C}/2$ the sum of the squares over $(\alpha,\beta)$ is finite, which defines $F_*$. Over
$(\alpha,\beta)\ne0$ it is $\le C(\rad_0/r_{\rm C}+\delta_y)^2$. We bound the pairing in a cell
by Cauchy--Schwarz with the weight $\rad_0^{-3}\chi_c^2\omega$, using $\omega\le1$.

\emph{Step 3: the terms of $\delta\Theta$.} We apply $\partial^\beta$ to~\eqref{eq:TWdtheta} at
$y$-order $0$, and to~\eqref{eq:UFleaf} for $y$-order $\ge1$. Every resulting term is of exactly one
of the following types. Their $X$-norms are bounded by the $\mathcal H_j$-norms, and the energies of
$\delta U$ are those of $\delta\varpi$ (macroscopic) and of $\tilde M_k$, which are bounded by the
micro energies at the same $\xi$, by Lemma~\ref{lem:TWmom}(i). Both carry the damping $\omega_\xi$
when paired with $\varphi_{0\beta}$ in the cell, since the pairing is on the same leaf.
\begin{enumerate}[label=($\Theta$\arabic*)]
\item \emph{Top:} $\Gamma_\mu(U)\,\partial^\beta\partial_i\delta U_k$, at $y$-order $0$. Its
$X$-norm is $\le C\max_\mu\sup_Q\|\Gamma_\mu(U)\|_{
\mathcal H_j}\big(|\partial^\beta\nabla\delta\varpi|+\sum_i\|\varphi_{0,\beta+e_i}\|_X\big)$, and each
multi-index of order $n+1$ arises from at most three of order $n$. Pairing with $\varphi_{0\beta}$,
integrating against $\rad_0^{-3}\chi_c^2\omega$ and using Cauchy--Schwarz, we obtain (c), with
$C_0^F=C\max_\mu\sup_Q\|\Gamma_\mu(U)\|_{\mathcal H_j}\le C(C_{\rm Max}+P_1^{1/2})\le C(1+C_{\rm
mom}+P_1^{1/2})$, the term $P_1^{1/2}$ coming from the affine dependence on $\tilde M_k$ and
part~(iii) of the micro-moment lemma.
\item \emph{Shifted Leibniz:} $\binom\beta{\beta'}\partial^{\beta'}\Gamma_\mu(U)\,\partial^{\beta-\beta'}\partial_i\delta
U_k$, $\beta'\ne0$. By Lemma~\ref{lem:TWcalc}(iv) with $u=\Gamma_\mu(U)$ and $v=\delta U_k$, these
are spread losses of the form of Lemma~\ref{lem:TWabs}, cube by cube, with
$c^Q_j=\sigma_j/(\zeta(2)^{1/2}S^{\ge1}_Q(u))$ and with the constant $C_1^F=C(C_{\rm
mom}+P_1^{1/2})$. Since the $c^Q_j$ depend on $Q$, we carry out the arithmetic step of that lemma
on each cube before the cubes are summed. Indeed, the supremum factor is
$S^{\ge1}_Q(\Gamma_\mu(U))\le C(C_{\rm Max}+S_Q(\tilde M))\le C(C_{\rm Max}+P_1^{1/2})$, by
part~(vi) of the calculus lemma and Lemmas~\ref{lem:TWcoef}--\ref{lem:TWmom}.

This factor is not small, because the micro moments are only bounded. Thus ($\Theta$2) is a spread
loss, absorbed by the radius at the rate $C_1^Ft^{-2}$, and not a Gronwall
term.

The second term of this bound is $\le Ct^{-2}\rad_0^{-1}S_Q(\delta
U)E^{\ge1}_Q(\Gamma_\mu(U))\mathfrak A_c^{1/2}$, and $E^{\ge1}_Q(\Gamma_\mu(U))\le C(C_{\rm
Max}+E_{Q^*}(\delta U))$. Indeed $\Gamma_\mu=\Gamma^0_\mu(u_g)+\Gamma^1_\mu(u_g)\cdot\tilde M$ is
affine in $\tilde M$, and the part $E_Q(\Gamma^1_\mu(u_g))\,S_Q(\tilde M)$ of the tame product of
Lemma~\ref{lem:TWcalc}(iii) uses $S_Q(\tilde M)\le C_{\rm Sob}E_{Q^*}(\tilde M)$, from~(i) of the
same lemma. Hence the second term is a Gronwall term with rate
$Ct^{-2}\rad_0^{-1}(C_{\rm mom}+P_1^{1/2})$ plus a constant source $Ct^{-2}\rad_0^{-1}(C_{\rm
mom}+P_1^{1/2})P_1^{1/2}\mathfrak A_c^{1/2}$, which is treated like the reference forcing.
\item \emph{Reference coefficient:} $\partial^\beta\big[(\Gamma_\mu(U)-\Gamma_\mu(U_*))\partial_iU_{*k}\big]$.
The factor $\partial_iU_*$ has Cauchy bounds at the radius $r_{\rm C}$, so by
Lemma~\ref{lem:TWcalc}(v) the shift does not change the bound. By~(iii) and~(vi) of the same lemma
and Lemma~\ref{lem:TWcoef}, this is a Gronwall term, $\le Ct^{-2}C_{\rm Max}\big(E_Q(\delta
U)+S_Q(\delta U)\big)\mathfrak A_c^{1/2}$.
\item \emph{Reference forcing:} $\partial^\beta\Theta_*$, which is Step~2.
\item \emph{Higher $y$-order: all terms with $|\alpha|\ge1$.} These are the terms
of~\eqref{eq:UFleaf} for every $\alpha'\le\alpha$, \emph{including} $\alpha'=0$. In each of
them we apply the same decomposition to $\partial^{\alpha-\alpha'}_z\Theta$, with $\Gamma_\mu$
replaced by $\partial^{\alpha-\alpha'}_z\Gamma_\mu$, which is still of the same type,
since $|\alpha-\alpha'|\le N_y$. Their top factors have $\xi$-level
$\le|\beta|+|\alpha'|+1$ and are $\delta\varpi$-derivatives or micro moments of $y$-order $0$,
measured with the weight of $y$-order $0$, while $\varphi_{\alpha\beta}$ carries
$W_{\alpha\beta}=w_{|\beta|}(\rad_{|\alpha|})\delta_y^{|\alpha|}$. Put $m:=|\alpha'|+1\ge1$. For
the unshifted factors, by the gap~\eqref{eq:UFgap} with $j=|\beta|$, we have
\begin{equation}
w_{|\beta|}(\rad_{|\alpha|})\,\delta_y^{|\alpha|}\,t^{-|\alpha'|}\le
C(N_y,\epsilon_y)\,\rad_0^{-1}(\rad_0t)^{-|\alpha'|}\,w_{|\beta|+|\alpha'|+1}(\rad).
\end{equation}
For the shifted factors, the bound~(iv) of the calculus lemma at the level $n=|\beta|+|\alpha'|$
produces in
addition $(n+1)/\rad$, and since $\rad_{|\alpha|}/\rad\le1-\epsilon_y$, we get
\begin{equation}
w_{|\beta|}(\rad_{|\alpha|})\,t^{-|\alpha'|}\,\frac{n+1}\rad\le(n+1)^{|\alpha'|+1}(1-\epsilon_y)^{|\beta|}\rad^{-|\alpha'|-1}t^{-|\alpha'|}\,w_n(\rad)\le
C(N_y,\epsilon_y)\,\rad_0^{-1}(\rad_0t)^{-|\alpha'|}\,w_n(\rad).
\end{equation}
Therefore all terms of the types ($\Theta$1)--($\Theta$3) at $y$-order $\ge1$, the top terms
included, are Gronwall terms with rate $Ct^{-2}\rad_0^{-1}(C_0^F+C_1^F)\le t^{-2}a_F$, plus
constant sources of the type of ($\Theta$2). They do not couple different families as spread losses:
the gap makes them Gronwall terms within the cell energy.
\item \emph{Kernel part:} the projection $\Pi_0$ of the full forcing, which is Step~1.
\end{enumerate}
\emph{Rates.} The list ($\Theta$1)--($\Theta$6) is exhaustive. The Gronwall terms are the second
part of ($\Theta$2), ($\Theta$3) and ($\Theta$5), with rates
\begin{equation}
\le Ct^{-2}\big(C_{\rm Max}+C_{\rm mom}+P_1^{1/2}\big)\big(1+\rad_0^{-1}\big)\le t^{-2}a_F,
\end{equation}
since $C_{\rm mom}\le C(1+B_{\rm cr}^{1/2})\le C(1+\rad_0P_1^{1/2})$. Their constant sources are
those of requirement~\ref{rq:TWb}. With ($\Theta$2) and ($\Theta$5) this is~(d). The factor
$\mathfrak p_*$ in $a_F$ is only an upper bound, since $\mathfrak e$ is expressed through
Theorem~\ref{thm:Meuler}.

\emph{Step 4: the parts $\ge1$.} For $(\alpha,\beta)\ne0$ the reference part is~(a). In
($\Theta$1)--($\Theta$5) every factor of highest level is in the part $\ge1$ or is macroscopic, with
one exception: if all derivatives fall on $\Gamma_\mu$, in ($\Theta$2) and ($\Theta$3), the
undifferentiated factor contains the level-$0$ moments $\tilde M_k$, that is, $s[\psi]$, $h[\psi]$
at level $0$. These multiply a Gaussian factor of level $|\beta|$ with the weight $W_{\alpha\beta}$,
whose reference part sums to $C(\rad_0/r_{\rm C}+\delta_y)\le C$. They yield the second term of
$\mathfrak b_c$.

\emph{Step 5: level $0$, (e).} At $(\alpha,\beta)=0$, $\hat F_{00}=t^{-2}(I-\Pi_0)T^{-1}_p\Theta$,
and the $\mathfrak e$-term is in the kernel. What remains is
$-t^{-2}(I-\Pi_0)\big[\sum_k(z\cdot\nabla p_{g,k})(\partial_{p_k}M)_{p_g}\big]$, a Gaussian
polynomial with coefficients $\nabla\varpi_g$. By Lemma~\ref{lem:Mref}, we have $|\nabla\varpi_*|\le
C(L_1,R_1)$. By the cube Sobolev inequality and (U1), we also have $|\nabla\delta\varpi|\le
C\delta_U/\rad_0\le1$. This implies $F_0=C(L_1,R_1)$, and the factor $1+\mathfrak p_*$ in~(e) is an
upper bound only. The kernel part is Step~1.

\emph{Step 6: the macroscopic hierarchy, (f).} We apply $\partial^\beta$ to
\begin{equation}
D\delta\varpi=t^{-2}\Big[\sum_\mu c_\mu(U)\partial_iU-\sum_\mu c_\mu(U_*)\partial_iU_*\Big],
\end{equation}
which is Theorem~\ref{thm:Meuler} for $\varpi_g$ and $\varpi_*$, with $s,h$ written by
part~(ii) of the micro-moment lemma, and Theorem~\ref{thm:momtransport} with $\alpha=0$, applied to
$\partial^\beta_vg$ (by Lemma~\ref{lem:commmom}, this is $\partial^\beta_\xi$ of~\eqref{eq:momlaw}).
Every term is of one of the types ($\Theta$1)--($\Theta$3), with $\Gamma_\mu$ replaced by the scalar
coefficients $c_\mu$:
\begin{itemize}
\item[(M1)] $c_\mu(U)\partial^\beta\partial_i\delta U$: the one-level losses of~(f), into the
macroscopic energy at level $n+1$ if $\delta U_k$ is a component of $\delta\varpi$, or into the
micro energy at level $n+1$ if it is $\tilde M_k$, by~(i) of the micro-moment lemma;
\item[(M2)] shifted Leibniz terms: spread losses with the constant $C_1^{\rm mac}=C(C_{\rm
mom}+P_1^{1/2})$, by the bound~(iv) of the calculus lemma. Its second term belongs to
$\mathfrak S^c_n$. For $\delta U_k=\tilde M_k$ it is $\le C\mathfrak M_c^{1/2}\mathfrak A_c^{1/2}$,
by part~(i) of the same lemma, Lemma~\ref{lem:TWcoef} and~(iii) of the micro-moment lemma. For
$\delta U_k$ a
component of $\delta\varpi$ it is $\le C(1+r_{\rm C}^{-1})C_{\rm Max}\mathfrak M_c+C
\mathfrak M_c^{1/2}\mathfrak A_c^{1/2}$, by the same Sobolev bound, (U1) and
$C\delta_U\le\rad_0$, and
the first part is included in $C'_\varpi\mathfrak M_c$;
\item[(M3)] $\partial^\beta[(c_\mu(U)-c_\mu(U_*))\partial_iU_*]$: Gronwall terms, $\le
Ct^{-2}(E_Q(\delta U)+S_Q(\delta U))\mathfrak M_c^{1/2}$. They contain the linearization at the
reference, for example $(\delta y_0\cdot\nabla)y_{0*}$, whose rate is part of $C'_\varpi$, not
small.
\end{itemize}
Since $D$ acts at fixed $\xi$ and $D\omega\le0$, we have
$D(\omega|\partial^\beta\delta\varpi|^2)\le2\omega|\partial^\beta\delta\varpi||D\partial^\beta\delta\varpi|$.
We integrate against $\rad_0^{-3}\chi_c^2$, which does not depend on $t$, and apply Cauchy--Schwarz
with the weight $\chi_c^2\omega$ on both factors. Then the micro factor is the damped cell energy at
the same $\xi$ as the macroscopic factor. In $\mathfrak S^c_n$ the coefficients from $\varpi_*$ give
$C'_\varpi$, the bilinear $\delta\varpi$-terms take the supremum on the factor of lower order by the
cube Sobolev inequality and (U1), and the micro moments below the top give $C
\mathfrak M_c^{1/2}\mathfrak A_c^{1/2}$. \emph{The weighted form of the spread losses.} By
the bound~(iv) of the calculus lemma,
on each cube $Q$ the terms (M2) at level $n$, summed over $n$ with the
weights $w_n^2$, are
\begin{equation}
\le\frac C\rad\,S^{\ge1}_Q(c_\mu(U))\sum_n\sum_{j\ge1}c^Q_j\,(n+1)\,x^Q_n\big(y^{{\rm
mac},Q}_{n+1-j}+y^{{\rm mic},Q}_{n+1-j}\big),\qquad\sum_jc^Q_j(j+1)\le1,
\end{equation}
where $x^Q_n$ and $y^{{\rm mac},Q}_m$ are macroscopic cube energies at the levels $n$ and $m$, and
$y^{{\rm mic},Q}_m$ are the micro ones that enter through $\tilde M$, as in~(M1).
In that bound only $1\le j\le n/2$ occurs, so $m=n+1-j\ge1$. Thus the level-$(0,0)$
micro energy never enters a spread loss, here or in ($\Theta$2). This is needed for the families of
the parts $\ge1$. For any $\upsilon>0$, by the arithmetic of the proof of the absorption lemma with
$2xy\le\upsilon x^2+\upsilon^{-1}y^2$ in place of $2xy\le x^2+y^2$, we get
\begin{equation}
\sum_n\sum_jc^Q_j(n+1)x^Q_ny^Q_{n+1-j}\le\frac\upsilon2\sum_n(n+1)(x^Q_n)^2+\frac1{2\upsilon}\sum_m(m+1)(y^Q_m)^2
.
\end{equation}
We take $\upsilon=1$ for $y^{{\rm mac},Q}$ and the given $\upsilon\ge1$ for $y^{{\rm mic},Q}$. Then
we multiply by $\chi_c(\xi_Q)^2\omega_{\xi_Q}$ and sum over $Q$ with Lemma~\ref{lem:CUcomp}. The
factor $S^{\ge1}_Q(c_\mu(U))\le C(C_{\rm mom}+P_1^{1/2})$ and the cube constant go into $C_1^{\rm
mac}$. This yields the bound on $\sum_nw_n^2\mathfrak L^c_n$ in~(f). This proves~(f), with (M2)
counted among the spread losses of Remark~\ref{rem:UFradius}.

\emph{Step 7: joint absorption, (g).} In the family $E^{(0)}_n:=E^c_{0,n}+K_{\mathfrak M}
\mathfrak M^c_n$, the cross terms are
\begin{equation}
C_0^F\sqrt{E^c_{0,n}}\sqrt{\mathfrak M^c_{n+1}}=C_0^FK_{
\mathfrak M}^{-1/2}\sqrt{E^c_{0,n}}\sqrt{K_{\mathfrak M}
\mathfrak M^c_{n+1}},\end{equation}\begin{equation}(C_0+C_0^F)CK_{\mathfrak M}\sqrt{
\mathfrak M^c_n}\sqrt{E^c_{0,n+1}}=(C_0+C_0^F)CK_{\mathfrak M}^{1/2}\sqrt{K_{\mathfrak M}
\mathfrak M^c_n}\sqrt{E^c_{0,n+1}},
\end{equation}
from (c) and from (f) times $K_{\mathfrak M}$. With $K_{\mathfrak M}=T_0\ge1$, all one-level losses
of (b), (c), (f) times $K_{\mathfrak M}$, of the transport and of the nonlocal terms are bounded in
total by $\big(2C_0+C_0^F+(C+2)(C_0+C_0^F)T_0^{1/2}\big)\sqrt{E^{(0)}_nE^{(0)}_{n+1}}$, which is
$\le C_0'\sqrt{E^{(0)}_nE^{(0)}_{n+1}}$ after enlarging $C_0$ and $C_0^F$ by the structural factor
$C+2$. The $C$ of~(f) contains $\eta^{-1/2}$ and the frame bounds. Since both constants are defined
only up to such factors, (c) and (f) remain true. The source of (f), times $K_{\mathfrak M}$, is
$CK_{\mathfrak M}\mathfrak M_c^{1/2}\mathfrak A_c^{1/2}=CK_{\mathfrak M}^{1/2}(K_{
\mathfrak M}\mathfrak M_c)^{1/2}\mathfrak A_c^{1/2}$. This proves~(g). The spread losses of the joint family
are those of ($\Theta$2), with the constant $C_1^F$, and those of (M2) times $K_{\mathfrak M}$. For
the latter we take $\upsilon=K_{\mathfrak M}^{1/2}$ in~(f). With $a_n^2:=w_n^2E^{(0)}_n\ge
e_n^2+K_{\mathfrak M}b_n^2$ (for the parts $\ge1$, $a_0^2=K_{\mathfrak M}b_0^2$, and $e_0=0$ in~(f))
and $K_{\mathfrak M}\ge1$, we obtain
\begin{equation}
K_{\mathfrak M}\sum_nw_n^2\mathfrak L^c_n\le C_1^{\rm
mac}\rad^{-1}t^{-2}\sum_n(n+1)\big((K_{\mathfrak M}+K_{\mathfrak M}^{3/2})b_n^2+K_{
\mathfrak M}^{1/2}e_n^2\big)\le2C_1^{\rm mac}T_0^{1/2}\rad^{-1}t^{-2}\sum_n(n+1)a_n^2 .
\end{equation}
Hence the joint family satisfies the hypothesis of the absorption lemma. We put the factor $2$
into $C_1^{\rm mac}$, and the factor $C_{\rm cube}$ of its cube-summed case,
which enters ($\Theta$2), into $C_1^F$ (for (M2) it is already in $C_1^{\rm
mac}$, Step~6). Thus the spread losses are absorbed with the rate $(C_1^F+C_1^{\rm
mac})T_0^{1/2}t^{-2}$, where $T_0^{1/2}=K_{\mathfrak M}^{1/2}$ comes from the choice of $\upsilon$.
This is the rate of that remark. With $\upsilon=1$ the micro part would carry the full
factor $K_{\mathfrak M}=T_0$, and the radius would lose a quantity of order $C_1^{\rm mac}$, which
is not small.

\emph{Inputs.} In Steps 1--7 we use (U1) (through $\delta U$, $|\nabla\delta\varpi|\le
C\delta_U/\rad_0$ and requirement~\ref{rq:TWa}); the cube bound on $\sup_Q\mathfrak a_{Q^*}$
(through~(iii) of the micro-moment lemma and Lemma~\ref{lem:CUsob}); the conditions $\rad_0\le
r_{\rm C}/2$, $\rad_0t\ge1$, $\omega\le1$, $D\omega\le0$; the comparability of
Lemma~\ref{lem:CUcomp}(i)--(iii), which uses requirements~\ref{rq:CUa}--\ref{rq:CUb} and that $H$ is
$C_1^U$-Lipschitz; the frame bounds (F3); the condition $C\delta_U\le\rad_0$; graded radii in
$[\rad_0/4,\rad_0]$ as in the absorption lemma; and $t\ge T_0\ge T_{\rm E}(L_1,R_1)$ for
Lemma~\ref{lem:Mref}. This is Remark~\ref{rem:UFinputs}. In the local theory each of these is
supplied by Proposition~\ref{prop:CULW}, points~(1) and~(2). This completes the proof.
\end{proof}

\begin{remark}[The decrease of the radius and the size of $\delta\varpi$]\label{rem:UFradius}
With $K_{\mathfrak M}=T_0$, the one-level part of the rate is $C_0't^{-2}$ (Lemma~\ref{lem:analytic}
needs $\frac34C_0'$, and the remaining $\frac14C_0'$ is the margin $|\rad'|-r_{\rm loc}$ below), and
its total loss is $C_0'/T_0\le C(C_0+C_0^F)T_0^{-1/2}$.

With the spread losses of Proposition~\ref{prop:UFforce}(d) and of the macroscopic hierarchy, whose
constants are $C_1^F$ and $C_1^{\rm mac}:=C(C_{\rm mom}+P_1^{1/2})$, the radius decreases at the
rate
\begin{equation}\label{eq:radrate}
\rad'=-\big(C_0'+(C_1^F+C_1^{\rm mac})T_0^{1/2}\big)t^{-2}
\end{equation}
as in Step~7 of the proof of that proposition. Here the constant of the joint family
is included in $C_1^F$ and $C_1^{\rm mac}$, so that Lemma~\ref{lem:TWabs} yields the Gronwall rate
$4(C_1^F+C_1^{\rm mac})T_0^{1/2}\rad_0^{-1}t^{-2}$ with no further constant. The total loss is
$\le C(C_0+C_0^F+C_1^F+C_1^{\rm mac})T_0^{-1/2}\le\rad_0/4$ by requirement~\ref{rq:TWb}.

In Lemma~\ref{lem:CUnodrop}, $|\rad'|$ is the full rate defined in this remark. The rate
$|\rad'|$ is at least the rate $r_{\rm loc}=\big(\frac34C_0'+(C_1^F+C_1^{\rm
mac})T_0^{1/2}\big)t^{-2}$ of the local theory of Proposition~\ref{prop:CULW}, with the constants
$C_0'$, $C_1^F$, $C_1^{\rm mac}$ of this remark, so that $|\rad'|-r_{\rm
loc}=\frac14C_0't^{-2}\ge0$. The local theory uses the same proposition, hence
the same absorbed constants, so $r_{\rm loc}$ still dominates the local spread losses.

The joint energy controls $\mathfrak M_c$ only with the factor $K_{\mathfrak M}^{-1}=T_0^{-1}$. The
source of part~(g) of the forcing proposition has time integral $\le CT_0^{-1/2}(\sup
\mathfrak A_c)^{1/2}$ against $\mathfrak B_c^{1/2}$. Therefore Step~3 of Theorem~\ref{thm:CUboot}, run for
$\mathfrak B_c$, yields
\begin{equation}
K_{\mathfrak M}\mathfrak M_c(t)\le\mathfrak B_c(t)\le C\big(c_\chi C_U+F_*^2\big)+4T_0
\mathfrak M_c(T_0) .
\end{equation}
Hence $\mathfrak M_c\le C(c_\chi C_U+F_*^2)T_0^{-1}+4\mathfrak M_c(T_0)$. By Lemma~\ref{lem:CUcomp},
we have $\mathfrak m_{Q^*}\le27e^5e^{2H(S_b+1)}\sup_c\mathfrak M_c$. This yields Step~5 of the proof
of that theorem, with
\begin{equation}
\delta_U:=Ce^{H(S_b+1)}\Big(\big(c_\chi C_U+F_*^2\big)^{1/2}T_0^{-1/2}+\eps T_0^{\kappa_\gamma}\Big)
\end{equation}
which is the $\delta_U$ of (U1) in the theorem. Every requirement on $\delta_U$ is a
smallness condition. Since $\delta_U$ is the sum of a part $Ce^{H(S_b+1)}(c_\chi
C_U+F_*^2)^{1/2}T_0^{-1/2}$, which does not depend on $\eps$, and a part $Ce^{H(S_b+1)}\eps
T_0^{\kappa_\gamma}$, each such requirement is split: the first part is imposed as a lower bound on
$T_0$, in terms of constants chosen before $T_0$, and the second as an upper bound on $\eps$,
which is chosen after $T_0$.
\end{remark}

\subsubsection{The level-\texorpdfstring{$0$}{0} dissipation of close leaves}\label{sssec:UFlevel0}

For each $\xi$ let $t^-_\xi$ be the time at which $s_\xi(t^-_\xi)=S_*-\frac12$, and $t^-_\xi=\infty$
if there is none. Since $s_\xi(T_0)=0$ and $s_\xi$ is increasing, $t^-_\xi>T_0$.

\begin{lem}[Level-$0$ dissipation of close leaves]\label{lem:UFlevel0}
Assume (U1), (U2), (U3), (U4) on $[T_0,T]$, the requirements of Theorem~\ref{thm:CUboot}, and
requirement~\ref{rq:UFd}. Then for every $\xi$ and $t\le T$,
\begin{equation}
\int_{t^-_\xi}^{t}\kappa\,d_\xi[f_\xi]\dd t'\le C_{L0}\,\delta_1^2,\qquad
C_{L0}:=\frac{8e^{4C_1^B}}{c_{\rm cl}} .
\end{equation}
The constant $C_{L0}$ is structural.
\end{lem}

\begin{proof}
\emph{Closeness.} For $t\ge t^-_\xi$ we have $s_\xi\ge S_*-\frac12$.
\begin{itemize}
\item If $s_\xi\le3S_\sharp$, then by Corollary~\ref{cor:Aleafentry} we have
$\|f_\xi\|_X\le\delta_1/C_e\le\delta_c$.
\item If $s_\xi>3S_\sharp$, then $s_\xi>S_b+1$, and by Step~1 of the proof of the bootstrap theorem
we have $\|f_\xi\|_X\le\delta_c$.
\end{itemize}

\emph{An auxiliary damping.} We put
\begin{equation}
\tilde H(s):=C_1^B\min\big\{(s-S_*+\tfrac12)_+,\;S_b-S_*+\tfrac32\big\},\qquad\tilde\omega_\xi:=e^{-2\tilde
H(s_\xi)} .
\end{equation}
Then $\tilde H'=C_1^B$ on $[S_*-\frac12,S_b+1]$ and $\tilde H'=0$ after. By
requirement~\ref{rq:CUc}, we have
$S_b-S_*=S_*(e^{2L'd_*}-1)\le4L'S_*d_*\le\frac1{2C_1^U}\le\frac12$. Hence
$\tilde\omega\ge e^{-4C_1^B}$.

\emph{The level-$0$ inequality.} We run the proof of Proposition~\ref{prop:Bband} for $\beta=0$,
with
$\tilde\omega$ in place of $\omega$. We use the numbering (T1)--(T12) of its terms in
Subsection~\ref{sssec:bandproof}. The terms are as follows.
\begin{itemize}
\item (T1)--(T4). At level $0$ the collision part is $L_\mu f+\mathcal Q(f,f)$. By
Proposition~\ref{prop:trilinG}, $|\langle
\mathcal Q(f,f),f\rangle|\le2C\|f\|_X\|f\|^2_{Y_2}\le2C\delta_c\|f\|^2_{Y_2}$. By
requirement~\ref{rq:Wdc}, this term can be put into the coercive part, and no rate remains.
\item (T5). There are no lower levels, so there are no sources.
\item (T6)--(T7). On $[S_*-\frac12,S_b+1]$ the rate $C_1^B$ absorbs the deficit
$(\frac12-\chi)_+\|f\|^2_{L^2}$ and the front term, as in the proposition. Only
closeness is used here. After $S_b+1$ we have $\chi=1$, so there is neither deficit nor front.
\item (T8)--(T10). In $z$-variables the local transport is $2t^{-2}\eta(A_p/A_{\rm ref})^2\int
m_x^2\,h\,z\cdot\partial_\xi h\dd z$, with $h=h_{00}$ and $\partial_\xi h=h_{0,e}$. Young's
inequality and (T8) with $C_2\ge C$ (a condition on the rate $C_2$ of $q(t)$) leave
$Ct^{-2}(e_\xi[f]+e_\xi[\varphi_{0e}])$. The nonlocal transport yields the same bound.
\item (T11). By Corollary~\ref{cor:FMuse}(a) and~(d), $\le4t^{-2}\mathfrak p_*e_\xi[f]$, plus a
nonlocal part absorbed by (T1) under requirement~\ref{rq:Wrate}, or Lemma~\ref{lem:EPreq} at
$\gamma=-2$.
\item (T12). By Proposition~\ref{prop:UFforce}(e) below, $\le
Ct^{-2}e_\xi[f]^{1/2}\big(F_0+Ce_\xi[\varphi_{0e}]^{1/2}\big)$.
\end{itemize}
Hence we obtain
\begin{equation}
D\big(\tilde\omega e_\xi[f]\big)+c_{\rm cl}\kappa\tilde\omega\,d_\xi[f]\le t^{-2}(C+4
\mathfrak p_*)\tilde\omega e_\xi[f]+Ct^{-2}\big(F_0^2+P_1\big),
\end{equation}
where $P_1$ bounds $\sup_\xi e_\xi[\varphi_{0e}]$, and the first $C$ is not larger than the constant
$C$ of $a_1$.

\emph{The meaning of the inequality.} The energy of a single leaf need not be differentiable along
its ray: the transport term contains the factor $z$, and the weight at time $t$ leaves no reserve
for it. The inequality is meant in integrated form: for every positive $\psi\in C^1([t^-_\xi,T])$
and all $t^-_\xi\le t'\le t\le T$,
\begin{equation}\label{eq:L0int}
\Big[\psi\,\tilde\omega\,e_\xi[f]\Big]_{t'}^{t}+\int_{t'}^{t}\psi\,c_{\rm
cl}\kappa\,\tilde\omega\,d_\xi[f]\dd\tau\le\int_{t'}^{t}\Big(\psi'\,\tilde\omega\,e_\xi[f]+\frac{\psi}{\tau^{2}}\big((C+4
\mathfrak p_*)\tilde\omega\,e_\xi[f]+C(F_0^2+P_1)\big)\Big)\dd\tau .
\end{equation}
To prove \eqref{eq:L0int}, we write the local part of the energy in $z$-variables, $\eta(A_p/A_{\rm
ref})^2\int m_x^2h^2\dd z$, and replace $m_x^2$ by $\zeta(z/R)^2m_x^2$, where $\zeta\in
C^\infty_c(\R^3)$, $0\le\zeta\le1$ and $\zeta=1$ on the unit ball. The truncated energy is Lipschitz
along the ray. Indeed, $g$ is a classical solution, so the truncated local part, an integral over a
compact set, is continuously differentiable. By the pointwise Gaussian bounds of
Lemma~\ref{lem:BTsob}(ii), $t\mapsto f_\xi(t)$ is Lipschitz with values in $X_8$, on which, by
Lemma~\ref{lem:lowweight}, $\nn{\cdot}^2_\natural$ is a bounded quadratic form. Moreover, $\chi$ is
$C^1$, while $\tilde\omega$ is Lipschitz. The derivative, which exists for a.e.\ $t$, is the sum of
the pairings from which the terms listed above arise, the local ones with the factor $\zeta(z/R)^2$
under the integral. We take the collision pairing in nondivergence form, namely
$2\eta\kappa(A_p/A_{\rm ref})^2\int\zeta(z/R)^2m_x^2\,h\,Q(g^{(x)},g^{(x)})\dd z$ with
$Q(g^{(x)},g^{(x)})=\ab[g^{(x)}]:\nabla^2_zg^{(x)}-\cb[g^{(x)}]g^{(x)}$, so that no derivative falls
on the cutoff. In $z$-variables the local part depends on the frame through the factor $(A_p/A_{\rm
ref})^2$ alone, by Proposition~\ref{prop:FM}(i).
The weight term (T8) and the local transport (T9) are
bounded pointwise in $z$. Indeed, by Young's inequality and the condition $C_2\ge C$ of the item
(T8)--(T10), we have
\begin{equation}\label{eq:T8T9-cut}
D(m_x^2)\,h^2+2t^{-2}m_x^2\,h\,z\cdot\partial_\xi h\le t^{-2}m_x^2\big(Ch^2+|\partial_\xi h|^2\big),
\end{equation}
and hence their contribution is $\le Ct^{-2}(e_\xi[f]+e_\xi[\varphi_{0e}])$ for every $R$. In all
other local terms the factor $\zeta(z/R)^2$ multiplies a function of $z$ that is integrable, with a
bound that is uniform on $[t',t]$. Indeed, $m_xh$, $m_x\nabla_zh$, $m_x\nabla^2_zh$ and
$m_x\partial_\xi h$ are in $L^2_z$ with such bounds, by Lemma~\ref{lem:CUsob}, (U2) and the
fact~(F1) of Subsection~\ref{ssec:inputs}, the coefficients of the collision operator are bounded,
and the forcing has finite $X$-norm, as in the item (T12) above. The nonlocal part of the energy
has polynomial weights and needs no cutoff. Now we multiply the truncated energy by
$\psi\tilde\omega$, which is Lipschitz, write the difference of the values at $t$ and $t'$ as the
integral of the derivative, and let $R\to\infty$. The truncated energies increase to $e_\xi[f]$,
which is bounded on $[t',t]$. Thus the boundary terms converge, and so do the integrals in which the
truncated energy is multiplied by $(\psi\tilde\omega)'$, by dominated convergence (this factor may
change sign). The other terms, apart from (T8) and (T9), converge by dominated convergence, in $z$
and then in time, to the pairings without cutoff, which are estimated as above. For (T8) and (T9) no
convergence is needed, since we use their bound~\eqref{eq:T8T9-cut}, which is uniform in $R$, before
taking the limit. This proves \eqref{eq:L0int}. We use only this integrated form, below and in the
proof of Theorem~\ref{thm:rate2}.

\emph{The bound $P_1$.} By the cube Sobolev inequality with $(\alpha,\beta)=(0,e)$, using
$w_{1+|\gamma|}\ge c\rad_0^{1+|\gamma|}$, then Lemma~\ref{lem:CUcomp} and (U2), we have
\begin{equation}
e_\xi[\varphi_{0e}]\le C\rad_0^{-2}\sup_Q\mathfrak a_{Q^*}\le C\rad_0^{-2}e^{2H(S_b+1)}c_\chi
C_U=:P_1 .
\end{equation}

\emph{Integration.} We have $\tilde\omega(t^-_\xi)=1$ and
$e_\xi[f](t^-_\xi)=\eta\|f\|^2_X\le\delta_1^2$, since $\chi=0$ there and $C_e^2\ge\eta$. We use
the integrated inequality on $[t^-_\xi,t]$ with $\psi(\tau):=\exp\big(-(C+4
\mathfrak p_*)\int_{t^-_\xi}^\tau\tau'^{-2}\dd\tau'\big)$, for which the terms with $e_\xi[f]$ on the right
cancel. The Gronwall factor $1/\psi$ is $\le2$ by the Gronwall item of requirement~\ref{rq:UFd},
whose rate contains $a_1\ge C+4\mathfrak p_*$. Hence we have
\begin{equation}
\int_{t^-_\xi}^tc_{\rm cl}\kappa\tilde\omega\,d_\xi[f]\le2\delta_1^2+CT_0^{-1}(F_0^2+P_1)\le4\delta_1^2
\end{equation}
by the same requirement. Dividing by $c_{\rm cl}\min\tilde\omega\ge c_{\rm cl}e^{-4C_1^B}$, we
obtain the lemma.
\end{proof}

At level $0$ the next level enters through the pointwise bound $P_1$ alone, with the integrable
factor $t^{-2}$, so the two-level loss of Remark~\ref{rem:ACtr}, which comes from the summation over
all levels, does not arise.

\subsubsection{The cell inequality for the parts \texorpdfstring{$\ge1$}{>=1}}\label{sssec:UFge1}

For a cube $Q_j$ we put
\begin{equation}
\rho_j(t):=\begin{cases}\kappa\,\rad_0^{-3}\displaystyle\int_{Q_j^*}d_\xi[f_\xi]\dd\xi&\text{if
$Q_j$ meets $\Sigma(t)$},\\[1ex]0&\text{otherwise.}\end{cases}
\end{equation}
This quantity is the undamped level-$0$ dissipation on $Q_j^*$. We also put
\begin{equation}
\mathfrak s_c:=C_{\mathfrak s}\Big(\frac{\rad_0}{R_1}+\delta_y\Big)^2\kappa\,\rad_0^{-3}\!\int_{\Sigma}\chi_c^2\,\omega_\xi\,d_\xi[f_\xi]\dd\xi,\qquad
\mathfrak s'_c:=C_{\mathfrak s'}\,\delta_U^2\sum_j\chi_c(\xi_j)^2\rho_j(t),
\end{equation}
with $C_{\mathfrak s}$, $C_{\mathfrak s'}$ structural. The first comes from the $\varpi_*$-part of
the coefficients, the second from their $\delta\varpi$-part.

\begin{prop}[The cell inequality for the parts $\ge1$]\label{prop:UFge1}
Assume the hypotheses of Proposition~\ref{prop:CUcell}, with (H$_a$) replaced by:
\begin{itemize}
\item[(H$_a'$)] every cube satisfies $\mathfrak a_{Q^*,\ge1}\le c_\sharp$; moreover
$C_Nc_\sharp^{1/2}\le c_{\min}/4$ and requirement~\ref{rq:UFb} holds.
\end{itemize}
Then for every cell $c$,
\begin{multline}
\frac{d}{dt}\mathfrak A_{c,\ge1}+\frac{c_{\min}}8\kappa\,\mathfrak D_{c,\ge1}\le
\mathcal T^c_{\ge1}+t^{-2}a_1\,\mathfrak A_{c,\ge1}+\mathfrak f^c_{\ge1}+\mathfrak q_c+
\mathfrak s_c+\mathfrak s'_c\\
+C_\rho\sum_j\chi_c(\xi_j)^2\rho_j(t)\,\mathfrak A_{\xi_j,\ge1} .
\end{multline}
Here:
\begin{itemize}
\item $\mathcal T^c_{\ge1}$ is the radius term plus the one-level losses
$C_0t^{-2}\sum_nw_n^2\sqrt{E^c_{n,\ge1}E^c_{n+1,\ge1}}$, as in that proposition;
\item $\mathfrak f^c_{\ge1}$ is the forcing, estimated in Proposition~\ref{prop:UFforce};
\item $\mathfrak q_c$ is as in the full cell inequality;
\item $\mathfrak A_{\xi_j,\ge1}$ is the energy of the cell centered at $\xi_j$;
\item $C_\rho$ is structural.
\end{itemize}
\end{prop}

\begin{proof}
The terms are those of the proof of Proposition~\ref{prop:CUcell}, restricted to
$(\alpha,\beta)\ne0$. Only the terms in which the level $(0,0)$ enters change. In all others both
factors are in the part $\ge1$, and the proof applies without change with (H$_a'$).

On crude leaves the level $(0,0)$ enters in two places: the linear terms (L1)--(L2) of
Proposition~\ref{prop:ACpf}, which are pointwise on the leaf and absorbed by the damping, and
the coefficient terms with the background, which give $\mathfrak q_c$. On buffer leaves the
linear term is estimated by Lemma~\ref{lem:CUbuf}, and the background enters there through
$\|f_\xi\|_X\le\delta_c$ and through $\|f_\xi\|_{Y_2}\le C(q_\pm)K$, both pointwise on the leaf, and
the remainder is controlled by the crude damping. On $\Sigma$ there are two
kinds of terms.

\emph{(i) Background terms.} On a leaf $\xi\in\Sigma$ and for $(\alpha,\beta)\ne0$, the background
part of (T4) is bounded by
\begin{equation}
2C_\tau\kappa\lambda_\xi W^2_{\alpha\beta}\big(\|f\|_X\|\varphi\|_{Y_2}+\|f\|_{Y_2}\|\varphi\|_X\big)\|\varphi\|_{Y_2},\qquad\varphi=\varphi_{\alpha\beta},
\end{equation}
with $C_\tau$ structural, by Proposition~\ref{prop:trilinG}. By (H$_b$), we have
$\|f\|_X\le\delta_c$, and $2C_\tau\delta_c\le c_{\rm cl}/16$ by requirement~\ref{rq:Wdc}. Applying
Young's inequality to the second part, we get
\begin{equation}
\le\frac{c_{\rm cl}}8\kappa\lambda_\xi W^2_{\alpha\beta}\|\varphi\|^2_{Y_2}+\frac{16C_\tau^2}{c_{\rm
cl}\eta}\,\kappa\,d_\xi[f]\,W^2_{\alpha\beta}e_\xi[\varphi] .
\end{equation}
We sum over $(\alpha,\beta)\ne0$ and write
\begin{equation}\label{eq:Age1}
A_{\ge1}(\xi):=\sum_{(\alpha,\beta)\ne0}W^2_{\alpha\beta}e_\xi[\varphi_{\alpha\beta}] .
\end{equation}
The first part is $\frac{c_{\rm cl}}8\kappa\mathbb D_{\ge1,\xi}$, which uses one eighth of the
coercivity on $\Sigma$. We integrate the second part against $\rad_0^{-3}\chi_c^2\omega_\xi$ over
$\Sigma$, cube by cube.
\begin{itemize}
\item On $Q_j$, $\chi_c^2\le e^{\sqrt3}\chi_c(\xi_j)^2$.
\item By Lemma~\ref{lem:CUsob} for the level $(0,0)$, in its undamped form,
\begin{equation}
\sup_{Q_j}\kappa d_\xi[f]\le C_S\kappa\sum_{|\gamma|\le2}\rad_0^{2|\gamma|-3}\!\int_{Q_j^*}d_{\xi'}[\varphi_{0\gamma}]\dd\xi'
.
\end{equation}
The term $\gamma=0$ is $C_S\rho_j$ when $Q_j$ meets $\Sigma$. For $|\gamma|\ge1$, we use
$\rad_0^{2|\gamma|}\le Cw^2_{|\gamma|}$ and $1\le e^{1/2}\omega_{\xi'}/\omega_{\xi_j}$ on
$Q_j^*$, by Lemma~\ref{lem:CUcomp}(i). Thus these terms are $\le C\omega_{\xi_j}^{-1}\kappa\,
\mathfrak d_{Q_j^*,\ge1}$.
\item $\rad_0^{-3}\int_{Q_j}\omega_\xi A_{\ge1}\le e^{1/2}\omega_{\xi_j}\mathfrak a_{Q_j,\ge1}$,
and $\omega_{\xi_j}\mathfrak a_{Q_j,\ge1}\le e^5\mathfrak A_{\xi_j,\ge1}$ by the same lemma.
\end{itemize}
Hence the second part is bounded by
\begin{equation}
C_\rho\sum_j\chi_c(\xi_j)^2\rho_j\,\mathfrak A_{\xi_j,\ge1}+C'\kappa\sum_j\chi_c(\xi_j)^2\,
\mathfrak a_{Q_j,\ge1}\,\mathfrak d_{Q_j^*,\ge1} .
\end{equation}
By (H$_a'$) and the same lemma, we deduce that the last sum is $\le C'c_\sharp\,27e^{12}\kappa
\mathfrak D_{c,\ge1}\le\frac{c_{\min}}{32}\kappa\mathfrak D_{c,\ge1}$, by requirement~\ref{rq:UFb}.

\emph{(ii) Coefficient terms with the background.} These are $\mathcal Q(
\mathfrak g_{\alpha\beta},f)+\mathcal Q(f,\mathfrak g_{\alpha\beta})$, the linear couplings of
Proposition~\ref{prop:TLsrc} with $(\alpha'',\beta'')=0$. On a leaf of $\Sigma$, the bound used
there is $C\|\mathfrak g_{\alpha\beta}\|\,\|f\|_{Y_2}\|\varphi\|_{Y_2}$. The coefficient $
\mathfrak g_{\alpha\beta}$ carries all the derivatives $\partial^\alpha_y\partial^\beta_v$, and we split it as
$\mathfrak g=\mathfrak g^*+\mathfrak g^\delta$, where $\mathfrak g^*$ is built from the reference
flow $\varpi'_*$ and $\mathfrak g^\delta$ is linear in the derivatives of $\delta\varpi$ (with
coefficients that are bounded, by the cube Sobolev inequality and (U1), when they involve
$\delta\varpi$ at lower order). The shifts $t^{-k}\partial^{|\beta|+k}_\xi$ produced by the
$y$-derivatives are controlled by Remark~\ref{rem:UFgrade}(c).
\begin{itemize}
\item \emph{The $\varpi_*$-part.} By the Cauchy estimates for $\varpi'_*$ at radius $r_{\rm
C}\gg\rad_0$, $W_{\alpha\beta}\|\mathfrak g^*_{\alpha\beta}\|\le\varkappa_{\alpha\beta}$ with
$\sum_{(\alpha,\beta)\ne0}\varkappa^2_{\alpha\beta}\le C(\rad_0/R_1+\delta_y)^2$, pointwise in
$\xi$. By Cauchy--Schwarz,
\begin{equation}
\le C\kappa\lambda_\xi\Big(\sum\varkappa^2\Big)^{1/2}\|f\|_{Y_2}\,
\mathbb D_{\ge1,\xi}^{1/2}\le\frac{c_{\rm cl}}{16}\kappa
\mathbb D_{\ge1,\xi}+C\Big(\frac{\rad_0}{R_1}+\delta_y\Big)^2\kappa\,d_\xi[f] .
\end{equation}
Integrated against $\rad_0^{-3}\chi_c^2\omega_\xi$ over $\Sigma$, the second part is $
\mathfrak s_c$.
\item \emph{The $\delta\varpi$-part.} Here $\delta\varpi$ is the factor of highest level, and (U1)
controls it only in cube norms. Following the rule of the remark after Proposition~\ref{prop:CUtri},
we take the supremum on the factor of lower level, which is $f_\xi$, and integrate $\delta\varpi$
over the cube. On a leaf $\xi\in\Sigma\cap Q_j$, by Young's inequality we have
\begin{equation}
C\kappa\lambda_\xi\sum_{(\alpha,\beta)\ne0}W_{\alpha\beta}\|
\mathfrak g^\delta_{\alpha\beta}\|\,\|f\|_{Y_2}W_{\alpha\beta}\|\varphi_{\alpha\beta}\|_{Y_2}\le\frac{c_{\rm
cl}}{32}\kappa\mathbb D_{\ge1,\xi}+C\Big(\sup_{Q_j}\kappa\,d_{\xi'}[f]\Big)\sum_{(\alpha,\beta)\ne0}W^2_{\alpha\beta}\|
\mathfrak g^\delta_{\alpha\beta}(\xi)\|^2 .
\end{equation}
As in~(i), $\sup_{Q_j}\kappa\,d_\xi[f]\le C_S\rho_j+C\omega_{\xi_j}^{-1}\kappa\,
\mathfrak d_{Q_j^*,\ge1}$, and $\rad_0^{-3}\int_{Q_j}\omega_\xi\sum W^2_{\alpha\beta}\|
\mathfrak g^\delta_{\alpha\beta}\|^2\le C\omega_{\xi_j}\mathfrak m_{Q_j^*}\le4C\omega_{\xi_j}\delta_U^2$, by
(U1) and the comparability on cubes. Multiplying by $\chi_c^2\le e^{\sqrt3}\chi_c(\xi_j)^2$ and
summing over $j$, we see that the second part is bounded by
\begin{equation}
C_{\mathfrak s'}\delta_U^2\sum_j\chi_c(\xi_j)^2\rho_j+C''\delta_U^2\,27e^{12}\,\kappa\,
\mathfrak D_{c,\ge1}=\mathfrak s'_c+C''\delta_U^2\,27e^{12}\,\kappa\,\mathfrak D_{c,\ge1}.
\end{equation}
The last term is $\le\frac{c_{\min}}{32}\kappa\mathfrak D_{c,\ge1}$ by requirement~\ref{rq:UFe}.
\end{itemize}

\emph{Collecting.} The terms (i) and (ii) \emph{replace} the background terms and the background
part of the coefficient terms in the proof of the full cell inequality, while the other losses there
(at most $\frac14c(\xi)$ for the quadratic terms and at most $\frac14c(\xi)$ for the remaining
coefficient terms) are unchanged. On a leaf of $\Sigma$, where $c(\xi)=c_{\rm cl}$, (i) and (ii)
remove $\big(\frac18+\frac1{16}+\frac1{32}\big)c_{\rm cl}=\frac7{32}c(\xi)$ of the leaf coefficient,
and
the cell dissipation loses
$\frac{c_{\min}}{32}+\frac{c_{\min}}{32}=\frac{c_{\min}}{16}\le\frac1{16}c(\xi)$ (from~(i) and from
the $\delta\varpi$-part of~(ii)). The total loss does not exceed
$\big(\frac14+\frac14+\frac7{32}+\frac2{32}\big)c(\xi)=\frac{25}{32}c(\xi)$, so at least
$\frac7{32}c(\xi)\ge\frac18c_{\min}$ remains. On the other leaves only the cell loss
$\frac1{16}c_{\min}$ is added. Hence at least $\frac{c_{\min}}8\kappa\mathfrak D_{c,\ge1}$ remains.
\end{proof}

\subsubsection{Closing the parts \texorpdfstring{$\ge1$}{>=1}}\label{sssec:UFclose}

We put
\begin{equation}
C_{\rm lat}:=\sup_{c\in\R^3}\sum_j\chi_c(\xi_j)^2,
\end{equation}
which is an absolute constant, because $\chi_c(\xi_j)^2=e^{-2\jap{k+a}}$ with $k\in\mathbb Z^3$ and
$a=(\xi_{j_0}-\xi_c)/\rad_0$. We also put
\begin{equation}
\mu_\sharp:=\frac{e^{-2H(S_b+1)}c_\sharp}{54e^5},\qquad\mu_\flat:=\frac{e^{-2C_1^US_\sharp-4C_1^B}c_\sharp}{54e^5}\le\mu_\sharp
,
\end{equation}
where the inequality uses $H(S_b+1)\le C_1^US_*+C_1^B(S_b+1-S_*)\le C_1^US_\sharp+2C_1^B$. The
constant $\mu_\flat$ involves only constants chosen before $\rad_0$.

\begin{lem}[Integral Gronwall over the lattice]\label{lem:UFfam}
For $c\in\R^3$ and $t\in[T_0,T]$ let $\mathfrak U_c(t)\ge0$, with $M:=\sup_{c,t}
\mathfrak U_c(t)<\infty$. Let $\rho_j\ge0$ with $\int_{T_0}^T\rho_j\le P$ for all $j$. Suppose that, for some
$B\ge0$, $\theta\in[0,\frac14]$ and $K_\rho\ge0$, for all $c$ and $t$,
\begin{equation}
\mathfrak U_c(t)\le B+\theta M+K_\rho\sum_j\chi_c(\xi_j)^2\int_{T_0}^t\rho_j(t')\,
\mathfrak U_{\xi_j}(t')\dd t' .
\end{equation}
If $K_\rho C_{\rm lat}P\le\frac14$, then $M\le2B$.
\end{lem}

\begin{proof}
The last term is $\le K_\rho C_{\rm lat}PM\le\frac14M$. Taking the supremum over $c$ and $t$, we
obtain $M\le B+\frac12M$, which proves the lemma.
\end{proof}

Note that we use no continuity in $c$ or $t$, and no supremum over leaves appears inside a time
integral. The rate $\rho_j$ is integrated in time on each cube first, and only then summed against
the cell weight.

\begin{prop}[The parts $\ge1$]\label{prop:UFU2}
Assume on $[T_0,T]$ the bootstrap hypotheses of Theorem~\ref{thm:CUboot} and the conclusions
(H$_b$), (H$_c$) of Steps~1--2 of its proof, together with requirements~\ref{rq:CUa}--\ref{rq:CUd}
and~\ref{rq:UFa}--\ref{rq:TWb}; that is, (U1), (U3), (U4) and
\begin{itemize}
\item[(U2)] for every cell, $\mathfrak A_c\le2c_\chi C_U$ and $\mathfrak B_{c,\ge1}\le2\mu_\sharp$,
\end{itemize}
where $\mathfrak B_{c,\ge1}:=\mathfrak A_{c,\ge1}+K_{\mathfrak M}\mathfrak M_c$ is the joint energy
of Proposition~\ref{prop:UFforce}. Then $\mathfrak B_{c,\ge1}\le\mu_\sharp$ for every cell on
$[T_0,T]$. Consequently:
\begin{enumerate}[label=(\roman*)]
\item every cube satisfies $\mathfrak a_{Q^*,\ge1}\le c_\sharp/2$; in particular (H$_a'$) holds;
\item every cell satisfies $U^c_{\ge1}\le m_\sharp:=c_\sharp/(54e^5)$, where $U^c_{\ge1}$ is the
undamped cell energy of the parts $\ge1$.
\end{enumerate}
\end{prop}

\begin{proof}
\emph{The evolution of the joint family.} The joint energy is $\mathfrak B_{c,\ge1}=
\mathfrak A_{c,\ge1}+K_{\mathfrak M}\mathfrak M_c$, with $K_{\mathfrak M}=T_0$.
\begin{itemize}
\item For $\mathfrak A_{c,\ge1}$ we use Proposition~\ref{prop:UFge1}, with the forcing $
\mathfrak f^c_{\ge1}$ written out by Proposition~\ref{prop:UFforce}(a)--(d), namely the reference
part, the source $\mathfrak b_c$ of~(d), the one-level losses of~(b), (c), the Gronwall term
$t^{-2}a_F \mathfrak B_{c,\ge1}$ of~(d), and the shifted Leibniz terms. By ($\Theta$2) and
($\Theta$5) in the proof of that proposition, the latter are spread losses with constant $C_1^F$,
Gronwall terms with rate $Ct^{-2}\rad_0^{-1}(C_{\rm mom}+P_1^{1/2})$, and a constant source
\begin{equation}
\mathfrak b'_c:=Ct^{-2}\rad_0^{-1}\big(C_{\rm mom}+P_1^{1/2}\big)P_1^{1/2}\,
\mathfrak B_{c,\ge1}^{1/2}.
\end{equation}
\item For $K_{\mathfrak M}\mathfrak M_c$ we use part~(f) of the same proposition, with
$\upsilon=K_{\mathfrak M}^{1/2}$, multiplied by $K_{\mathfrak M}$. It yields one-level losses, which
join those of the micro energies at $y$-order~$0$ by part~(g); spread
losses with constant $C_1^{\rm mac}T_0^{1/2}$ in the joint variables, as in (M2) and Step~7 of the
proof; a Gronwall term with rate $C'_\varpi t^{-2}$; and the
source of~(g),
\begin{equation}
\mathfrak b^{\rm mac}_c:=Ct^{-2}T_0^{1/2}\big(K_{\mathfrak M}\mathfrak M_c\big)^{1/2}
\mathfrak A_c^{1/2}\le Ct^{-2}T_0^{1/2}\big(2c_\chi C_U\big)^{1/2}\mathfrak B_{c,\ge1}^{1/2},
\end{equation}
by (U2). Here the \emph{full} energy $\mathfrak A_c$ appears, including the level $(0,0)$, because
the micro moments of all levels below the top enter the macroscopic hierarchy.
\end{itemize}
Applying Lemma~\ref{lem:TWabs}, which contains Lemma~\ref{lem:analytic}, to each $y$-order and to
the joint family at $y$-order $0$, we see that the one-level losses and the spread losses are
absorbed by the radius, which decreases at the rate of Remark~\ref{rem:UFradius}. The absorption of
the one-level losses by the radius itself produces a Gronwall rate. In the analytic lemma it is
$C_0b/\rad_0$, and for the joint family at $y$-order $0$ the one-level constant is
$C_0'=2C_0+C_0^F+(C_0+C_0^F)T_0^{1/2}$ by part~(g) of the forcing proposition, so the rate is
bounded by $C_0'\rad_0^{-1}t^{-2}$ on every $y$-order. The absorption of the spread losses produces
the additional Gronwall rate $4(C_1^F+C_1^{\rm mac})T_0^{1/2}\rad_0^{-1}t^{-2}$.
Therefore the total Gronwall rate is
\begin{equation}
a(t):=t^{-2}\Big(a_1+a_F+C'_\varpi+C\rad_0^{-1}\big(C_{\rm
mom}+P_1^{1/2}\big)+C_0'\rad_0^{-1}\Big)+4(C_1^F+C_1^{\rm mac})T_0^{1/2}\rad_0^{-1}t^{-2},
\end{equation}
and $\int_{T_0}^\infty a\le\frac12\log2+\frac12\log2=\log2$, by requirements~\ref{rq:UFd}
and~\ref{rq:TWb}. The Gronwall factor is therefore $\le2$, and integrating, we obtain
\begin{equation}
\mathfrak B_{c,\ge1}(t)\le2\mathfrak B_{c,\ge1}(T_0)+2\!\int_{T_0}^t\big(\mathfrak b_c+
\mathfrak b'_c+\mathfrak b^{\rm mac}_c+\mathfrak q_c+\mathfrak s_c+
\mathfrak s'_c\big)+2C_\rho\sum_j\chi_c(\xi_j)^2\!\int_{T_0}^t\rho_j\,\mathfrak A_{\xi_j,\ge1}.
\end{equation}
We bound the terms in turn. Let $M:=\sup_{c,t}\mathfrak B_{c,\ge1}\le2\mu_\sharp$, and recall
$\mu_\flat\le\mu_\sharp$.
\begin{itemize}
\item \emph{Initial value.} By Step~0 of the proof of the bootstrap theorem, we have
$\mathfrak M_c(T_0)=0$ and
$\mathfrak A_{c,\ge1}(T_0)\le\mu_\flat/32$ by requirement~\ref{rq:UFc}, so $2
\mathfrak B_{c,\ge1}(T_0)\le\mu_\flat/16\le\frac18\mu_\sharp$.
\item \emph{Crude remainder.} As in Step~3 of the same proof, $2\int
\mathfrak q_c\le\frac1{16}\mu_\flat\le\frac18\mu_\sharp$, by requirement~\ref{rq:crude}.
\item \emph{Coefficient remainders.} On $\Sigma(t')$ every leaf has $t'\ge t^-_\xi$, and
$\omega\le1$. By Fubini's theorem and Lemma~\ref{lem:UFlevel0}, we get $2\int\mathfrak s_c\le2C_{
\mathfrak s}(\rad_0/R_1+\delta_y)^2c_\chi C_{L0}\delta_1^2\le\frac1{16}\mu_\sharp$, by
the same requirement. With the bound on $\int\rho_j$ below, we also get $2\int\mathfrak s'_c\le2C_{
\mathfrak s'}\delta_U^2C_{\rm lat}P\le\frac1{16}\mu_\sharp$, by requirement~\ref{rq:UFe}.
\item \emph{Sources with a square root.} Each of $\mathfrak b_c$, $\mathfrak b'_c$, $
\mathfrak b^{\rm mac}_c$ is of the form $Yt^{-2}\mathfrak B_{c,\ge1}^{1/2}$ plus, for $\mathfrak b_c$, a part
with $(\mathfrak A^{00}_c)^{1/2}\le(2c_\chi C_U)^{1/2}$ in place of a constant. Hence
$2\int\le2YT_0^{-1}M^{1/2}\le\frac1{16}M+16Y^2T_0^{-2}$. With
\begin{equation}
Y_{\mathfrak b}=C\big(F_*+(c_\chi C_U)^{1/2}\big),\qquad Y_{
\mathfrak b'}=C\rad_0^{-1}\big(C_{\rm mom}+P_1^{1/2}\big)P_1^{1/2},\qquad Y_{\rm mac}=CT_0^{1/2}(c_\chi
C_U)^{1/2},
\end{equation}
the same two requirements imply $16Y^2T_0^{-2}\le\frac1{32}\mu_\sharp$ for each of
the three. For the macroscopic source, $16Y_{\rm mac}^2T_0^{-2}=CT_0^{-1}c_\chi C_U$, which decays
only like $T_0^{-1}$, because of the weight $K_{\mathfrak M}=T_0$.
\item \emph{The rate term.} If $Q_j$ meets $\Sigma(t')$, then by Lemma~\ref{lem:CUcomp}(i) every
leaf of $Q_j^*$ has $s\ge S_*(1-\epsilon_0)\ge S_*-\frac12$, so $t'\ge t^-_{\xi'}$ for all $\xi'\in
Q_j^*$. Hence, by Fubini's theorem and the dissipation bound of close leaves, we have
\begin{equation}
\int_{T_0}^T\rho_j\le\rad_0^{-3}\!\int_{Q_j^*}\int_{t^-_{\xi'}}^T\kappa\,d_{\xi'}[f_{\xi'}]\dd
t\dd\xi'\le27C_{L0}\delta_1^2=:P .
\end{equation}
\end{itemize}
Moreover, $\mathfrak A_{\xi_j,\ge1}\le\mathfrak B_{\xi_j,\ge1}$. Adding up, we obtain
\begin{equation}
\mathfrak B_{c,\ge1}(t)\le\Big(\tfrac18+\tfrac18+\tfrac1{16}+\tfrac1{16}+\tfrac3{32}\Big)\mu_\sharp+\tfrac3{16}M+2C_\rho\sum_j\chi_c(\xi_j)^2\!\int_{T_0}^t\rho_j\,
\mathfrak B_{\xi_j,\ge1}\le\tfrac{15}{32}\mu_\sharp+\tfrac3{16}M+\cdots .
\end{equation}
Lemma~\ref{lem:UFfam} applies with $\mathfrak U_c=\mathfrak B_{c,\ge1}$,
$B=\frac{15}{32}\mu_\sharp$, $\theta=\frac3{16}\le\frac14$, $K_\rho=2C_\rho$ and $P$ as above. Its
condition $2C_\rho C_{\rm lat}\cdot27C_{L0}\delta_1^2\le\frac14$ holds with a margin, since
requirement~\ref{rq:UFa} implies $2C_\rho C_{\rm lat}\cdot27C_{L0}\delta_1^2\le\frac18$. Indeed,
we have $M\le\frac{15}{32}\mu_\sharp+\frac3{16}M+\frac18M$, so
$M\le\frac{16}{11}\cdot\frac{15}{32}\mu_\sharp=\frac{15}{22}\mu_\sharp\le\mu_\sharp$.

For (i), by the comparability on cubes we have $\mathfrak a_{Q^*,\ge1}\le27e^5e^{2H(S_b+1)}\sup_c
\mathfrak A_{c,\ge1}\le c_\sharp/2$. For (ii), we have $\omega\ge e^{-2H(S_b+1)}$ everywhere, so
$U^c_{\ge1}\le e^{2H(S_b+1)}\mathfrak A_{c,\ge1}\le m_\sharp$.
\end{proof}

The threshold $c_\sharp$ is the one of the crude regime, and it is now required in every cell. This
is possible because the rate $\rho_j$ multiplies the parts $\ge1$ themselves. Indeed, the size of
the parts $\ge1$ is set by the initial data, the crude remainders and the coefficient remainders,
all of which are made small by $\rad_0$, $\delta_y$, $T_0$ and $\eps$. The level-$0$ size
$\delta_1^2$ enters through the factor $C_{\rm lat}P$ alone, and requirement~\ref{rq:UFa} makes it
$\le\frac18$.

\subsubsection{The near region at entry}\label{sssec:UFnear}

\begin{lem}[The near region lies below the band's end]\label{lem:UFnear}
Assume requirements~\ref{rq:CUa} and~\ref{rq:CUc}, and $C_1^U\ge1$. Let $c$ be a cell,
$s_0\in[S_*,S_b]$, and $t_c$ the time with $s_{\xi_c}(t_c)=s_0$. If a cube $Q_j$ meets
$\{|\xi-\xi_c|\le d_*\}$, then
\begin{equation}
S_*-\tfrac12\le s_{\xi'}(t_c)<S_b+1\qquad\text{for all }\xi'\in Q_j^* .
\end{equation}
\end{lem}

\begin{proof}
Since $|\nabla\log\lambda_{\rm ref}|\le2L'$, we have $|\log s_\xi-\log s_{\xi'}|\le2L'|\xi-\xi'|$. A
point $\xi'\in Q_j^*$ lies within $d_*+3\sqrt3\rad_0$ of $\xi_c$, and
$2L'\cdot3\sqrt3\rad_0=\epsilon_0/2$. We put $a:=2L'd_*+\epsilon_0/2$. Then we have $s_0e^{-a}\le
s_{\xi'}(t_c)\le s_0e^{a}\le S_be^a$.

\emph{Upper bound.} By requirement~\ref{rq:CUc}, we have $2L'd_*\le1/(4C_1^US_*)\le\frac18$, since
$S_*\ge2$. Hence $S_b\cdot2L'd_*=S_*e^{2L'd_*}\cdot2L'd_*\le\frac{e^{1/8}}{4C_1^U}$. By
requirement~\ref{rq:CUa}, $S_b\epsilon_0\le S_\sharp\epsilon_0\le\frac1{8C_1^U}$, and
$\epsilon_0\le\frac1{8S_\sharp}\le\frac1{24}$. Thus $S_ba\le\frac{e^{1/8}}4+\frac1{16}<\frac12$, and
$a\le\frac18+\frac1{48}$ implies $e^a<2$. Therefore
\begin{equation}
S_b(e^a-1)\le S_b\,a\,e^a<1 .
\end{equation}
\emph{Lower bound.} We have $s_0e^{-a}\ge S_*e^{-a}\ge S_*(1-a)\ge S_*-\frac5{16}$, since $S_*a\le
C_1^US_*a\le\frac14+\frac1{16}$, by the same two requirements, in the forms $8C_1^US_*L'd_*\le1$
and $8C_1^US_\sharp\epsilon_0\le1$, and $C_1^U\ge1$.
\end{proof}

Consequently in Step~4 of Theorem~\ref{thm:CUboot}, every cube $Q$ that meets the near region has
$S_*-\frac12\le s<S_b+1\le3S_\sharp$ on $Q^*$, and the argument of the first case of Step~1 applies
to it, although $Q$ need not contain a leaf with $s\ge S_*$. Its undamped cube energy is
$\le27\delta_1^2+\mathfrak a_{Q^*,\ge1}\le28\delta_1^2$, by Corollary~\ref{cor:Aleafentry} and
Proposition~\ref{prop:UFU2}(i). The constant $C'_{\rm cl}=2(Cc_\chi+3)$ of Step~4 is therefore
structural, and no bound (U3) of another cell is used. Bounding the near region through the bound
(U3) of neighboring cells would instead multiply the constant of (U3) by an absolute factor larger
than $1$.

\subsubsection{Requirements}\label{sssec:UFreq}

These are part of Appendix~\ref{app:const}.
\begin{itemize}
\item[] \refstepcounter{rq}\label{rq:UFa}\textbf{\therq.} $432\,C_\rho C_{\rm
lat}C_{L0}\delta_1^2\le1$. All three constants are structural. (The factor $432$ instead of the
critical value $216$ leaves a margin in Lemma~\ref{lem:UFfam}.)
\item[] \refstepcounter{rq}\label{rq:UFb}\textbf{\therq.} $C'\,27e^{12}c_\sharp\le
c_{\min}/32$, with $C'$ of Proposition~\ref{prop:UFge1}, in addition to requirement~\ref{rq:Wcsharp}
and $c_\sharp\le\delta_1^2/4$.
\item[] \refstepcounter{rq}\label{rq:UFc}\textbf{\therq.} $\rad_0\le r_{\rm C}/8$
(Proposition~\ref{prop:TLsrc}; it contains $\rad_0\le r_{\rm C}/2$); $\epsilon_y=1/(4N_y+4)$;
$32c_\chi C_DB\big((\rad_0/\Lambda)^2+\delta_y^2\big)\le\mu_\flat$; and $32C_{
\mathfrak s}(\rad_0/R_1+\delta_y)^2c_\chi C_{L0}\delta_1^2\le\mu_\flat$.
\item[] \refstepcounter{rq}\label{rq:UFd}\textbf{\therq.} $\rad_0T_0\ge1$;
$C(C_0+C_0^F)T_0^{-1/2}\le\rad_0/4$; $CT_0^{-1}(F_0^2+P_1)\le2\delta_1^2$; $\int_{T_0}^\infty
t^{-2}\big(a_1+a_F+C'_\varpi+C\rad_0^{-1}(C_{\rm
mom}+P_1^{1/2})+C_0'\rad_0^{-1}+CT_0^{1/2}\big)\le\frac12\log2$, where
$C_0'=2C_0+C_0^F+(C_0+C_0^F)T_0^{1/2}$, so that the term $C_0'\rad_0^{-1}$ contributes
$(2C_0+C_0^F)\rad_0^{-1}T_0^{-1}+(C_0+C_0^F)\rad_0^{-1}T_0^{-1/2}$ (it is the rate of the
absorption of the one-level losses by the radius); the last term is the rate of the source of
Proposition~\ref{prop:UFforce}(g) for the full energy, used in Step~3 of Theorem~\ref{thm:CUboot};
$16Y_{\mathfrak b}^2T_0^{-2}=CT_0^{-2}\big(F_*^2+c_\chi C_U\big)\le\mu_\flat/32$;
$CT_0^{-1}F_*(c_\chi C_U)^{1/2}\le c_\chi/8$ (the reference forcing in Step~3 of
Theorem~\ref{thm:CUboot}); and $16Y_{\rm mac}^2T_0^{-2}=CT_0^{-1}c_\chi C_U\le\mu_\flat/32$, and
$CT_0^{-1}c_\chi C_U\le\frac13e^{-2C_1^US_\sharp-4C_1^B}\delta_1^2$ (Step~4 of
Theorem~\ref{thm:CUboot}), with $Y_{\mathfrak b},Y_{\rm mac}$ of the proof of
Proposition~\ref{prop:UFU2}. All are lower bounds on $T_0$ in terms of earlier constants.
\item[] \refstepcounter{rq}\label{rq:UFe}\textbf{\therq.} With $\delta_U$ of
Remark~\ref{rem:UFradius}: $32C_{\mathfrak s'}\delta_U^2C_{\rm
lat}\,27C_{L0}\delta_1^2\le\mu_\flat$; $C''\delta_U^2\,27e^{12}\le c_{\min}/32$, with $C''$ of the
proof of Proposition~\ref{prop:UFge1}; and $32c_\chi T_0(C\eps T_0^{\kappa_\gamma})^2\le\mu_\flat$.
(Since $\delta\varpi(T_0)=0$ by Step~0 of Theorem~\ref{thm:CUboot}, the last clause is not used.
It is harmless, since a larger $\delta_U$ only weakens (U1).)

\item[] \refstepcounter{rq}\label{rq:crude}\textbf{\therq.} $6C(K)c_\chi
c_\lambda^{-1}S_\sharp\big((\rad_0/R_1)^2+\delta_y^2+\delta_U^2\big)\le\mu_\flat/16$, with $C(K)$
the constant of the crude remainder $\mathfrak q_c$ in Proposition~\ref{prop:CUcell}. All constants
involved, other than $\rad_0$, $\delta_y$ and $\delta_U$, are fixed before $\rad_0$.
The $\rad_0$- and $\delta_y$-parts are upper bounds on these, and the $\delta_U$-part is split into
conditions on $T_0$
and on $\eps$ as in Remark~\ref{rem:UFradius}. Since $\mu_\flat/16\le e^{-2C_1^US_\sharp-4C_1^B}\delta_1^2$ (by
$c_\sharp\le\delta_1^2/4$), it also implies $\int\mathfrak q_c\le\omega(S_b)\delta_1^2$ in Step~4 of
Theorem~\ref{thm:CUboot}.
\item[] \refstepcounter{rq}\label{rq:TWa}\textbf{\therq.}
$C\big(\rad_0+\rad_0/r_{\rm C}+\delta_U\big)\le r_p/(16dC_{\rm alg})$, with
$r_p=\frac1{10}$ and $d=5$ the number of components of $u_g$.
\item[] \refstepcounter{rq}\label{rq:TWb}\textbf{\therq.}
$C\big(C_0+C_0^F+C_1^F+C_1^{\rm mac}\big)T_0^{-1/2}\le\rad_0/4$; $4(C_1^F+C_1^{\rm
mac})T_0^{-1/2}\rad_0^{-1}\le\frac12\log2$; and, for the constant sources of the types
($\Theta$2) and ($\Theta$5) in the proof of Proposition~\ref{prop:UFforce},
\begin{equation}
\big(CT_0^{-1}\rad_0^{-1}(C_{\rm mom}+P_1^{1/2})P_1^{1/2}\big)^2\le\min\big\{\mu_\flat/512,\;\tfrac13e^{-2C_1^US_\sharp-4C_1^B}\delta_1^2\big\},
\end{equation}that is, $16Y_{\mathfrak b'}^2T_0^{-2}\le\mu_\flat/32$; and
$CT_0^{-1}\rad_0^{-1}(C_{\rm mom}+P_1^{1/2})P_1^{1/2}(c_\chi C_U)^{1/2}\le c_\chi/8$ (Step~3 of
Theorem~\ref{thm:CUboot}). Here $C_1^F$ and $C_1^{\rm mac}$ involve
$P_1=C\rad_0^{-2}e^{2H(S_b+1)}c_\chi C_U$, so this is a lower bound on $T_0$ in terms of
constants chosen before $T_0$.
\end{itemize}
Requirement~\ref{rq:TWa} is an upper bound on $\rad_0$ and a smallness condition on
$\delta_U$. Requirement~\ref{rq:TWb} replaces the radius part of requirement~\ref{rq:UFd}.
Requirements~\ref{rq:UFa} and~\ref{rq:UFb} are upper bounds on $\delta_1$ and on
$c_\sharp$. Requirement~\ref{rq:UFc} bounds $\rad_0$ and $\delta_y$ from above
in terms of constants chosen before $\rad_0$, and requirement~\ref{rq:UFd} bounds $T_0$ from below.
Requirement~\ref{rq:UFe} holds for $T_0$ large and then $\eps$ small.

The three sources with a square root in the proof of Proposition~\ref{prop:UFU2} are $
\mathfrak b_c$, which comes from the reference part and the level-$0$ moments in the forcing
(Proposition~\ref{prop:UFforce}(d)); the constant source $\mathfrak b'_c$ of the shifted Leibniz
terms; and the source $\mathfrak b^{\rm mac}_c$ of the macroscopic hierarchy in part~(g) of the
same proposition. The last one is the reason for the condition $CT_0^{-1}c_\chi
C_U\le\mu_\flat/32$ in requirement~\ref{rq:UFd}: because of the weight $K_{\mathfrak M}=T_0$, its
integral is of order $T_0^{-1/2}$, not $T_0^{-1}$.

\subsection{Leaf Sobolev and Cauchy estimates and the inputs of the fibered
comparison}\label{ssec:inputs}

The fibered comparison of Section~\ref{sec:fib} needs three inputs from the bootstrap on $[T_0,T]$,
listed also in Remark~\ref{rem:RGfullinputs}. These are the bounds (G1) and (G3) on every leaf,
normalized by $\rho_\xi$, as in Lemma~\ref{lem:RGbar}(a); the crude bound (H3) of
Subsection~\ref{ssec:fibstanding}, without Gaussian weight, on the whole region; and the qualitative
Gaussian decay of~(c) in the same lemma. On the unrelaxed leaves, (G3)--(G4) are the bootstrap
hypothesis (U4), which must in turn be improved. By Lemma~\ref{lem:BTG12}, (G1)--(G2) follow from
(U1). We derive the inputs from the cell energies in Lemma~\ref{lem:BTsob}, and we improve (U4) in
Lemma~\ref{lem:BTU4}. Throughout, $J:=10$ and $n:=5$, so that $J\ge n+3$ and $N_y\ge J+2$.

We use the following elementary facts.
\begin{itemize}
\item[(F1)] \emph{Leaf derivatives.} In the leaf variables $(z,\xi)$, with $y=z$ and $v=\xi-z/t$,
we have $\partial_v=\partial_\xi$ and $\partial_y=\partial_z+t^{-1}\partial_\xi$. Hence
$\nabla_z\big(h^{(x)}\big)=\big((\partial_y-t^{-1}\partial_v)h\big)^{(x)}$, and $\nabla^k_z$ of the
leaf restriction of $\hat r_{\alpha\beta}=\partial^\alpha_y\partial^\beta_v\hat r$ is a sum of at
most $2^k$ restrictions of $t^{-i}\hat r_{\alpha+a,\beta+b}$, $|a|+|b|=k$, $|b|=i$.
\item[(F2)] \emph{Weighted Sobolev inequality in $z$.} Let $m=e^{q|\hat z|^2+\omega}$ be a leaf
weight of the class $\mathcal W(q_\pm,\eps_1)$, with $\omega$ as in Lemma~\ref{lem:Gleafweight} or,
in the frames $p_g$, as in Lemma~\ref{lem:TG0}(i). Applying $H^2(\R^3)\subset L^\infty(\R^3)$ to
$e^{(q/4)|\hat z|^2}\varphi$, we get terms $|\hat z|^ke^{(q/4)|\hat z|^2}|\nabla^j\varphi|=|\hat
z|^ke^{-(3q/4)|\hat z|^2}e^{-\omega}\,m|\nabla^j\varphi|$, $k+j\le2$. Since $e^{\pm\omega}\le
C\jap{\hat z}^{p}e^{Cq\jap{\hat z}}$, with $p$ the exponent of (D3), the factor $e^{-(3q/4)|\hat
z|^2}$ controls the polynomial and the $e^{\pm\omega}$ losses, and we get
\begin{equation}
\sup_{\hat z}e^{(q/4)|\hat z|^2}|\varphi(\hat z)|\le C(q,L,\eps_1)\sum_{k\le2}\|m\nabla^k_{\hat
z}\varphi\|_{L^2}.
\end{equation}
\item[(F3)] \emph{Frames.} By (U1), Lemma~\ref{lem:Mref} and Lemma~\ref{lem:TWcalc}(i) at level
$0$, $|p_g(t,\xi)-p_{\rm ref}(\xi)|\le C(L_1,R_1)(\delta_U+T_0^{-1})$, where $p_{\rm ref}(\xi)$ are
the leaf parameters at $T_0$. For the mass this needs one more step, since
$\varpi'=(\nabla_\xi\ell,y_0,\vartheta)$ does not contain $\ell=\log\rho_g$. We integrate
$D_t\ell=t^{-2}(\nabla\cdot y_0+y_0\cdot\nabla\ell)$, the first equation of
Theorem~\ref{thm:Meuler}, along the rays, and obtain
$|\ell_g(t,\xi)-\ell_g(T_0,\xi)|\le\sqrt3\,C_\varpi(1+C_\varpi)/T_0$, since $|\nabla\cdot
y_0|\le\sqrt3|\nabla y_0|$. The same applies to $\partial^k_\xi\rho_g$. The leaf through $T_0\xi$
samples the fibers $v=\xi-z/T_0$, and $|\nabla\log\rho_v|+|\nabla\log\theta_v|\le L$ by (H1). With
the bound on $\partial_vg(T_0)$ of (LWP)(a3), the leaf parameters at $T_0$ and the fiber parameters
$P_\xi$ differ by $O(L/T_0)$, with a constant $C(L,C'_D,\beta,\theta_\pm)$. For
$\delta_U$ small and $T_0$ large we get
\begin{equation}
\tfrac12\rho_\xi\le\rho_g(t,\xi)\le2\rho_\xi,\qquad|y_{0,g}(t,\xi)|\le
L+1,\qquad\tfrac12\theta_-\le\theta_g(t,\xi)\le\theta_++1 .
\end{equation}
Moreover $|\partial^k_\xi p_g(t,\xi)|\le C_k(L_1,R_1,\rad_{\rm post})$ for $k\le J$, by the
Cauchy estimates for $\varpi'_*$ on $\mathbb S_{R_1/2}$ and for $\delta\varpi$ at radius
$\rad_{\rm post}/2$, by Lemma~\ref{lem:CUsob} and (U1).
\end{itemize}

\begin{remark}[Normalizations]\label{rem:BTnorm}
In (U4), Lemma~\ref{lem:BTU4} and Lemma~\ref{lem:RGbar}(a), the bounds (G1)--(G3) refer to the
normalized leaf $\hat g^{(x)}=\rho_\xi^{-1}g^{(x)}$ in the variable $z$, and (G4) refers to the
normalized micro part $f_\xi$ in the frame variable $\hat z$. Explicitly, (G3) with constant $K$ is
$N^{(3)}_\xi:=\sup_z\jap z^8\hat g^{(x)}(z)\le K$, and (G4) with constant $K$ is
\begin{equation}
N^{(4)}_\xi:=\max_{j\le2}\Big(\big\|\jap{\hat z}\hat m_\xi\nabla^j_{\hat
z}f_\xi\big\|_{L^2}+\big\|\jap{\hat z}\hat m_\xi\nabla^j_{\hat z}f_\xi\big\|_{L^\infty}\Big)\le K,
\end{equation}
where
\begin{equation}\label{eq:fxi}
f_\xi:=T^{-1}_{p_g(t,\xi)}\big(g^{(x)}-M_{p_g(t,\xi)}\big)=T^{-1}_{\hat p}\big(\hat g^{(x)}-M_{\hat
p}\big),
\end{equation}
$\hat p=(\rho_g/\rho_\xi,y_{0,g},\theta_g)$ (the mass is contained in $T_p$), and where $\hat m_\xi$
is the leaf weight of Section~\ref{sec:weights}, given by Lemma~\ref{lem:Gleafweight} and
Lemma~\ref{lem:TG0}(i), at the energy exponent $q(t)$. The weight of (G4) is thus the full Gaussian
weight of the energies, not a polynomial weight. The bound (G4) enters
Lemma~\ref{lem:CUbuf} and the crude leaf inequality, always as a conclusion of
Lemma~\ref{lem:BTU4}.

The crude leaf inequality Proposition~\ref{prop:ACpf} and Lemma~\ref{lem:Whypo} use (G1)--(G3) for
$G_\xi=T_p^{-1}g^{(x)}$ in the frame, $p=p_g(t,\xi)$. By (F3), $T_p^{-1}$ changes masses, centers
and temperatures by bounded factors, and $\jap{\hat z}^8\le C(L,\theta_\pm)\jap z^8$. Therefore
$G_\xi$ satisfies (G1) with mass $1$, center $0$ and temperature $1$, and (G3) with the constant
$\bar K:=C_{\rm fr}(L,\theta_\pm)K$, where $C_{\rm fr}\ge1+C(L,\theta_\pm)$, so that $\bar K\ge
K+C(L,\theta_\pm)$ for $K\ge1$, as the frame step needs.

Applying Lemma~\ref{lem:RGellip} in the frame (mass $1$, second moment $3$, $L^\infty$ bound
$\bar K$), we obtain (G2) for $G_\xi$ with a constant $c_G(\bar K)$. We write
$c_h(K):=c_h(\bar K,c_G(\bar K))$ for the coercivity constant of Lemma~\ref{lem:Whypo} with these
constants. In~(L1) of the proof of Proposition~\ref{prop:ACpf} we apply the lemma with
$16\jap{L_0}^8\bar K$ in place of $\bar K$ in (G3), and $c_h(K)$ denotes the constant obtained
there.

Since (G3) with a constant $K'\le K$ implies (G3) with $K$, and (G2) is obtained from
Lemma~\ref{lem:RGellip} with the $L^\infty$ bound $\bar K$, the coercivity $c_h(K)$ is available on
every leaf on which (G3) holds with a constant $\le K$.
\end{remark}

\begin{lem}[(G1) and (G2) from the macroscopic bound]\label{lem:BTG12}
Assume (U1) on $[T_0,T]$ and the smallness of $\delta_U+T_0^{-1}$ in (F3). Then for every
$t\in[T_0,T]$ and every $\xi$, the normalized leaf $\hat g^{(x)}$ satisfies (G1) with
\begin{equation}
M\in[\tfrac12,2],\qquad|y_0|\le L+1,\qquad\theta\in[\tfrac12\theta_-,\theta_++1].
\end{equation}
If moreover $\sup_z\jap z^8\hat g^{(x)}(z)\le K$, then $\hat g^{(x)}$ satisfies (G2) with
$c_G=c(\gamma,\frac12,E_{\rm mom},K)$ of Lemma~\ref{lem:RGellip}, $E_{\rm
mom}:=2\big((L+1)^2+3(\theta_++1)\big)$.
\end{lem}

\begin{proof}
By Lemma~\ref{lem:Mfree}, the micro part has no leaf moments, so the mass, center and temperature of
$g^{(x)}$ are those of $M_{p_g(t,\xi)}$, namely $(\rho_g,y_{0,g},\theta_g)$. By (F3), we have
$\rho_g/\rho_\xi\in[\frac12,2]$, $|y_{0,g}|\le L+1$ and $\theta_g\in[\frac12\theta_-,\theta_++1]$.
This is (G1), and it implies $\int|z|^2\hat g^{(x)}\le M E_{\rm mom}$. Applying
Lemma~\ref{lem:RGellip} with $M_-=\frac12$, the second moment $E_{\rm mom}$ and the $L^\infty$
bound $K$, we obtain (G2).
\end{proof}

Thus (G1) and (G2) are not bootstrapped: they follow from (U1) on every leaf, and (G2) from (G3) in
addition. The constant $c_G$ here is not the structural constant $c_G$ of
Lemma~\ref{lem:Abarrier}.

\begin{lem}[Leaf Sobolev and Cauchy estimates]\label{lem:BTsob}
Assume (U1)--(U3) of Theorem~\ref{thm:CUboot} on $[T_0,T]$, and the requirements of that theorem.
Put $B_{\rm cr}:=54e^5e^{2H(S_b+1)}c_\chi C_U$. Let $t\in[T_0,T]$, $x=t\xi$, and let $Q$ be the cube
that contains $\xi$. Put $B_\xi:=B_{\rm cr}$ if $Q^*\subset\{s<S_b+1\}$, and $B_\xi:=\delta_b^2$
otherwise.
\begin{enumerate}[label=(\roman*)]
\item For $|\alpha|\le N_y$ and $|\beta|=k$,
\begin{equation}
\ell_\xi[\hat r_{\alpha\beta}]\le C\,(k!)^2\,(\rad_{\rm
post}/2)^{-2k}\,\delta_y^{-2|\alpha|}\,B_\xi .
\end{equation}
\item For $|a|+|b|\le J$ and all $\hat z$, with $p=p_g(t,\xi)$,
\begin{equation}
e^{(q/4)|\hat z|^2}\big|T_p^{-1}\big(\partial^a_y\partial^b_v\hat r\big)^{(x)}(\hat z)\big|\le
C_J\,\rad_{\rm post}^{-J-2}\,\delta_y^{-J-2}\,B_\xi^{1/2},
\end{equation}
with $C_J=C_J(q,L,\theta_\pm)$.
\item The crude bound (H3) holds with $J=10$ on $[T_0,T]$, with
\begin{equation}\label{eq:Kcr}
K_{\rm cr}:=C_J'\big(1+\rad_{\rm post}^{-J-2}\delta_y^{-J-2}B_{\rm cr}^{1/2}\big),
\end{equation}
where $C'_J=C'_J(q,L_1,\theta_\pm,R_1,\rad_{\rm post})$, through the bounds on $\partial^k_\xi
p_g$ in (F3). This is admissible, since $K_{\rm cr}$ is chosen after $\rad_0$ and $\delta_y$.
\item Every leaf satisfies, normalized by $\rho_\xi$, the bounds of Lemma~\ref{lem:RGbar}(a) with
\begin{equation}
M_-=\tfrac12,\qquad M_+=2,\qquad E_+=2E_{\rm mom},\qquad K=K_M+C_2'\,\rad_{\rm
post}^{-2}\,\delta_y^{-2}\,B_\xi^{1/2},
\end{equation}
where
\begin{equation}\label{eq:KM}
K_M:=\sup_{\xi,z}\jap z^8\rho_\xi^{-1}M_{p_g(t,\xi)}(z)\le C(L,\theta_\pm)
\end{equation}
and $C_2'=C_2'(q,L,\theta_\pm)$. Hence
\begin{equation}\label{eq:Krel}
K\le K_{\rm rel}:=K_M+C_2'\rad_{\rm post}^{-2}\delta_y^{-2}\delta_b
\end{equation}
on every leaf with $B_\xi=\delta_b^2$, in particular on every leaf with $s_\xi\ge S_b+1$, and
\begin{equation}\label{eq:Kcrrel}
K\le K^{\rm cr}_{\rm rel}:=K_M+C_2'\rad_{\rm post}^{-2}\delta_y^{-2}B_{\rm cr}^{1/2}
\end{equation}
on every leaf.
\item $g$ satisfies the qualitative Gaussian decay of~(c) in the same lemma.
\end{enumerate}
Parts (i), (ii) and (iv) are static: at a time $t$ they use only a bound $\mathfrak a_{Q^*}(t)\le B$
on the cube of the leaf, graded radii $\ge\rad_{\rm post}$, and the frame bounds (F3) at $t$, and
then they hold with $B$ in place of $B_\xi$.
\end{lem}

\begin{proof}
(i) \emph{Pointwise bounds in $\xi$, and Cauchy estimates.} If $Q^*\subset\{s<S_b+1\}$, then by
Lemma~\ref{lem:CUcomp} and (U2) we have $\mathfrak a_{Q^*}\le27e^5\omega^{-1}\cdot2c_\chi C_U\le
B_{\rm cr}$, since $\omega^{-1}\le e^{2H(S_b+1)}$. Otherwise some leaf of $Q^*$ has $s\ge S_b+1$,
and the second case of Step~1 of the proof of Theorem~\ref{thm:CUboot}, which rests on (U3) and the
same lemma alone, implies $\mathfrak a_{Q^*}\le\delta_b^2$. In particular, we make no use here of
the comparison with the fibered flow. By Lemma~\ref{lem:CUsob}, we have for each $(\alpha,\beta)$
with $|\beta|=k$ and each $\xi\in Q$
\begin{equation}
\eta\,\tilde w_k(\rad)^2\,\delta_y^{2|\alpha|}\,\ell_\xi[\hat r_{\alpha\beta}]\le
C_S'\,\mathfrak a_{Q^*}\le C_S'B_\xi .
\end{equation}
Here $\tilde w_k(\rad)^{-1}=k!(k+2)^2/(\rad^k(k+1)^5)\le4\,k!\,\rad^{-k}$, and every graded
radius satisfies $\rad_{|\alpha|}(t)\ge\rad_{\rm post}=\rad_0/2$, since the total loss of
radius is $\le\rad_0/4$ by requirement~\ref{rq:TWb}, and $N_y\epsilon_y\le\frac14$ by
Remark~\ref{rem:UFgrade}(a). This is~(i). It is a Cauchy estimate: the analytic weights in $v$ turn
the energy bound into a factorial bound on each $v$-derivative.

(ii) \emph{Leaf Sobolev.} We fix $|a|+|b|\le J$. By the definition of $\ell_\xi$ in
Subsection~\ref{ssec:cells}, $\ell_\xi[h]=\|T_p^{-1}h^{(x)}\|^2_{X(\hat m_\xi)}$, and $\hat m_\xi$
lies in $\mathcal W(q_\pm,\eps_1)$ by Lemma~\ref{lem:TG0}(i). We apply (F2) to
$\varphi=T_p^{-1}(\hat r_{ab})^{(x)}$. Since $\nabla_{\hat z}=\theta^{1/2}\nabla_z$ in the frame,
and by (F1) and (F3), $\|\hat m_\xi\nabla^k_{\hat z}\varphi\|_{L^2}$, $k\le2$, is bounded by
$C(\theta_\pm)$ times a finite sum of terms $t^{-i}\ell_\xi[\hat r_{a+a',b+b'}]^{1/2}$ with
$|a'|+|b'|=k$. These have $|a+a'|\le J+2\le N_y$ and $|b+b'|\le J+2$. By~(i) each is $\le
C_J\rad_{\rm post}^{-J-2}\delta_y^{-J-2}B_\xi^{1/2}$.

(iii) \emph{The crude bound.} We write $g=\mathcal M[g]+\hat r$. For the micro part, we undo the
normalization. Since $T_p\psi(z)=\rho\theta^{-3/2}\psi((z-y_0)/\theta^{1/2})$, by~(ii) and (F3) we
have
\begin{equation}
|\partial^a_y\partial^b_v\hat r(t,y,v)|\le C_J\,\rho_g(t,\xi)\,e^{-(q/4)|\hat z|^2}\rad_{\rm
post}^{-J-2}\delta_y^{-J-2}B_{\rm cr}^{1/2},\qquad\xi=v+y/t,\ z=y,
\end{equation}
since $B_\xi\le B_{\rm cr}$ (we may assume $\delta_b^2\le B_{\rm cr}$). For the Maxwellian part,
$\mathcal M[g](t,y,v)=M_{p_g(t,v+y/t)}(y)$. By the chain rule, $\partial^a_y\partial^b_v
\mathcal M[g]$ is a sum of products of derivatives of $(p,y)\mapsto M_p(y)$ of total order $\le J$ with
derivatives $t^{-i}\partial^k_\xi p_g$, $k\le J$. By (F3) these are bounded by $C_J\rho_g(t,\xi)\jap
y^{2J}e^{-|y-y_{0,g}|^2/(2\theta_g)}$. Finally, $\rho_g(t,\xi)\le2\rho_\xi$ by (F3), and
$\rho_\xi\le\rho_ve^{L|y|/t}\le\rho_ve^{L|y|}$ since $|\xi-v|=|y|/t$ and $|\nabla\log\rho|\le L$.
The Gaussian factors in $y$ control $\jap y^{2J}e^{L|y|}$, with a constant depending on
$(q,L,\theta_\pm)$. Taking the supremum over $(y,v)$, we obtain (H3).

(iv) By Lemma~\ref{lem:Mfree}, the micro part is moment-free on every leaf, so the leaf moments of
$g^{(x)}$ are those of $M_{p_g(t,\xi)}$: mass $\rho_g$, center $y_{0,g}$ and temperature $\theta_g$.
By (F3), we have $\frac12\rho_\xi\le\int g^{(x)}\le2\rho_\xi$ and
$\int|z|^2g^{(x)}\le2\rho_\xi\big((L+1)^2+3(\theta_++1)\big)$. For the pointwise bound, we
apply~(ii) with $a=b=0$. The Sobolev step involves the levels $|a'|+|b'|\le2$, so the factor is
$C\rad_{\rm post}^{-2}\delta_y^{-2}$, and with the normalization we get
\begin{equation}
\jap z^8\rho_\xi^{-1}g^{(x)}(z)\le K_M+C\rad_{\rm
post}^{-2}\delta_y^{-2}\,B_\xi^{1/2}\,\sup_z\jap z^8e^{-(q/4)|\hat z|^2},
\end{equation}
which is the stated bound after the choice of $C_2'$. No smallness is required: $K_{\rm rel}$ and
$K^{\rm cr}_{\rm rel}$ are chosen after $\rad_0$ in Appendix~\ref{app:const}, and they enter only
the barrier constant $C_{\rm bar}$, not the barrier rate.

(v) By~(ii) at $a=b=0$, on each leaf $|\hat r^{(x)}(z)|\le C\rho_g(t,\xi)e^{-(q/4)|\hat z|^2}$ with
a constant that may depend on $T$, while the Maxwellian part is Gaussian with the parameters of
(F3). Since $|\hat z|^2\ge|z|^2/(2\theta_g)-C$, we deduce that $g^{(x)}(z)\le C(T)\rho_\xi
e^{-q_1|z|^2}$ for some $q_1>0$, which is the claim.
\end{proof}

\begin{lem}[(U4) from (U1)--(U3)]\label{lem:BTU4}
Assume (U1)--(U3) of Theorem~\ref{thm:CUboot} on $[T_0,T]$, the requirements of that theorem, and
those on $T_0$ in Remark~\ref{rem:RGconst}, with the level $S=S_\sharp$, so that
$S^\sharp=4S_\sharp$. Let $K$ be chosen as
\begin{equation}
K:=\max\big\{2(2K_3+1),\ 2K_4\big\},\qquad K_4:=C(L,\theta_\pm,q_\natural)\,C_T(K_G+1)+C(L,\theta_\pm),
\end{equation}
with $K_3,K_G$ from Lemma~\ref{lem:RGGbounds} and $C_T$ from Theorem~\ref{thm:RGF}, and apply
Lemma~\ref{lem:RGbar} with the $L^\infty$ constant $K_{\rm bar}:=\max\{K,K^{\rm cr}_{\rm rel}\}$.
Then every leaf with $s_\xi(t)\le3S_\sharp$, $t\in[T_0,T]$, satisfies (G3) and (G4)
with constant $K/2$, in the normalizations of Remark~\ref{rem:BTnorm}. In particular (U4) on
$[T_0,T]$, with $K/2$, is a consequence of (U1)--(U3) on $[T_0,T]$. (G1) holds on every leaf by
Lemma~\ref{lem:BTG12}, and (G2) on every leaf where (G3) holds, in particular on the crude leaves,
where it is used.
\end{lem}

\begin{proof}
We check the hypotheses of Proposition~\ref{prop:RGcompF} and Theorem~\ref{thm:RGF}, with $n=5$ and
$J=10$.
\begin{itemize}
\item (H1) and (H2) concern the data at $T_0$. They are clause (a3) of (LWP)(a)
in Section~\ref{sec:proofmain}, which holds for $\eps\le\eps_0(T_0)$ with the constants $L$,
$\theta_\pm$, $(C_{\rm dat},\beta,A_7)$, which are fixed first. Here Subsection~\ref{ssec:hominputs}
concerns the leaves, not the fibers.
\item (H3) with $J=10$ is Lemma~\ref{lem:BTsob}(iii).
\item Hypothesis~(a) of Lemma~\ref{lem:RGbar}, normalized by $\rho_\xi$: on the leaves with
$s_\xi\le3S_\sharp$ the moments are those of Lemma~\ref{lem:BTG12}, with $M_-=\frac12$, $M_+=2$ and
$E_+=2E_{\rm mom}$, and the $L^\infty$ bound is part~(iv) of the leaf Sobolev lemma with
$B_\xi\le B_{\rm cr}$, that is, $K^{\rm cr}_{\rm rel}$. Note that (U4) is not used.
On the others, $s_\xi>3S_\sharp\ge S_b+1$
(since $S_b\le S_\sharp$ and $S_\sharp\ge1$), and it is the same bound with $K_{\rm
rel}$. The constants $M_\pm,E_+$ are therefore structural, and so is the barrier rate $r$. The
barrier constant $C_{\rm bar}$ depends on $K_{\rm bar}$, which is fixed before $T_0$.
\item Hypothesis~(b) holds by (H1)--(H2), as stated in
that proposition, and hypothesis~(c) is part~(v) of the leaf Sobolev lemma.
\item The conclusions of Proposition~\ref{prop:RS} on $[0,S^\sharp]$ follow from (E2), for data in
$\mathfrak K^{n+2}=\mathfrak K^7$ , by Remark~\ref{rem:RSshortF}.
\item The smallness conditions $(\ast)$, $C_M\varepsilon_T(T_0)\le\varepsilon_\natural$ and
$C(q)C_T\varepsilon_T(T_0)\le1$ are the requirements on $T_0$ of Remark~\ref{rem:RGconst}.
\end{itemize}
The leaves with $s_\xi(t)\le3S_\sharp$ sample only fibers with $s_v\le4S_\sharp=S^\sharp$ on $Z_t$,
as in the proof of the same theorem. Hence part~(iv) of that theorem applies, and $\hat
g^{(x)}$ satisfies (G3) with constant $2K_3+1\le K/2$.

For (G4), by part~(ii) of the same theorem with $j\le n-2=3$, we have
$\|e^{(q/4)|z|^2}\nabla^j_z\hat g^{(x)}\|_{L^\infty}+\|e^{(q/8)|z|^2}\nabla^j_z\hat
g^{(x)}\|_{L^2}\le C_T(K_G+1)$. The normalized
micro part is $f_\xi=T^{-1}_{\hat p}(\hat g^{(x)}-M_{\hat p})$, as in Remark~\ref{rem:BTnorm}, where
$T^{-1}_{\hat p}$ is the change of variables $z=y_{0,g}+\theta_g^{1/2}\hat z$ times
$(\rho_\xi/\rho_g)\theta_g^{3/2}$, with $\rho_\xi/\rho_g\in[\frac12,2]$ by (F3), and the Maxwellian
part has Gaussian bounds with constant $C(L,\theta_\pm)$. In the frame, $|\hat
z|^2\le4(|z|^2+(L+1)^2)/\theta_-$, since $\theta_g\ge\frac12\theta_-$ by (F3), so $\jap{\hat z}\hat
m_\xi\le C\,e^{(4q_0/\theta_-)|z|^2}\jap z\le C\,e^{(q/16)|z|^2}\jap z\le C\,e^{(q/8)|z|^2}$, by the
choice $q_0\le\theta_-q_\natural/64$ of Remark~\ref{rem:RGfull}(b), with the comparison run at
$q=q_\natural$. The change of variables gives the factor $(\rho_\xi/\rho_g)\theta_g^{3/2+j/2}\le
C(L,\theta_\pm)$, and the weight comparison a factor
$C(q_\natural)$, which also controls $e^{\omega}\le
C\jap{\hat z}^{p}e^{Cq\jap{\hat z}}$ from Lemma~\ref{lem:TG0}(i), with $p$ of (D3). Hence
$\|\jap{\hat z}\hat m_\xi\nabla^jf_\xi\|_{L^2\cap L^\infty}\le
C(L,\theta_\pm,q_\natural)C_T(K_G+1)+C(L,\theta_\pm)=K_4\le K/2$ for $j\le2$.

The constant $K$ is chosen as follows. The barrier rate $r$ depends only on the data rate and on
$M_-/E_+$, by the barrier lemma,
so $q_\natural$, $r_\natural$ and $C_T$ do not depend on $K$, and
neither do $K_3$, $K_G$ and $K_M$. The constants $K$, $K_{\rm rel}$ and $K^{\rm cr}_{\rm rel}$ enter
only $C_{\rm bar}$ and $C_{\rm loc}$, hence $\varepsilon_T$, that is, the lower bound on
$T_0$.
\end{proof}

With these two lemmas, Step~2 of the proof of Theorem~\ref{thm:CUboot} is complete.

\section{The initial layer and the endpoint \texorpdfstring{$\gamma=-2$}{gamma=-2}}\label{sec:layer}

In this section we bring data posed at $t=0$ to the form required at the time $t_1$ by
Theorem~\ref{thm:t1}. This is done in Theorem~\ref{thm:layer} for $\gamma\in[-3,-2)$, and in
Theorem~\ref{thm:EPlayer} for $\gamma\in[-3,-2]$.

On $[0,t_1]$ the collisions act with strength $O(\eps)$, and
the solution stays close to free transport. The difficulty is to propagate analyticity in $v$ and
Gaussian decay in $y$, starting from data that are only close to vacuum, without the leaf structure,
which is not yet established.

\subsection{The equation for the relative density}\label{ssec:reldens}

In the initial layer we use the density relative to free transport, whose value at $t=0$ is $1$.
Write $g=f_0G$, where $f_0(y,v)=\eps\rho_0(v)F_0(y,v)$ is the profile of free transport, and
$D=\partial_v-t\nabla_y$. Then
\begin{equation}\label{eq:Geq}
\partial_tG=D\cdot(A\,DG)+(2A\zeta-\bb)\cdot DG+S_0G,\qquad A=\ab[f_0G],
\end{equation}
where $\zeta=D\log f_0$ and $S_0$ is an explicit zero-order coefficient. The coefficients are taken
at $(y,v)$ and the convolutions are along $w\mapsto(y+tw,v-w)$. The coefficient $A$ is of size
$\eps$ but degenerate in two ways: in the radial direction $\hat y$ it is only of size $\eps\jap
y^{-2}$, and its derivatives in $y$ lose a factor $\jap y$ radially.

\emph{The radius in $v$.} We put $R:=\Lambda/2$. In the initial layer, $R$ denotes this radius in
all
bounds of the form $C^{|\beta|}\beta!R^{-|\beta|}$ and in all conditions on $\rad/R$. For
$v\in\R^3$ the polydisc $\Delta(v,R):=\{v+\zeta:\zeta\in\mathbb C^3,\ |\zeta_j|<R\}$ lies in
$\mathbb S_\Lambda$, since $|\mathrm{Im}\,\zeta|\le|\zeta|<\sqrt3R<\Lambda$. Write $\zeta=a+i\eta$
with $a,\eta\in\R^3$. By the upper half of (D3), we have $|\rho_0(v+a+i\eta)|\le A_F\rho_0(v+a)$. Using this and
$|\nabla\log\rho_0|\le L$, we have $|\rho_0(v+\zeta)|\le A_Fe^{L|a|}\rho_0(v)\le
A_Fe^{L\Lambda}\rho_0(v)$ on $\Delta(v,R)$, as in Lemma~\ref{lem:dispfree}. Applying the Cauchy
estimates on the closed polydiscs of radius $r<R$ about $v$ and letting $r\uparrow R$, we obtain
\begin{equation}\label{eq:rho0Cauchy}
|\partial^\beta\rho_0(v)|\le A_Fe^{L\Lambda}\,\beta!\,R^{-|\beta|}\rho_0(v)\qquad(v\in\R^3,\
\beta\in\mathbb N^3).
\end{equation}

\emph{The Gaussian part and the perturbation.} By (D$_{\rm iso}$) we write
\begin{equation}\label{eq:FGU}
F_0(y,v)=F_{\rm G}(y)\,\mathcal U(y,v),\qquad F_{\rm G}(y):=c_{\rm
G}e^{-\beta_F|y|^2-\Phi_1(y)},\qquad\mathcal U(y,v):=\frac{c_0(v)}{c_{\rm G}}e^{-\Phi_{\rm v}(y,v)},
\end{equation}
where $c_{\rm G}>0$ normalizes $F_{\rm G}$. For product data $\mathcal U=1$.

\begin{lem}[The perturbation]\label{lem:Ufac}
Put $K_{\mathcal U}:=e^{2C_{\rm pv}}$ and $C_{\mathcal U}:=2C_{\rm pv}+\log(1/\cos1)+\frac\pi2$.
Then:
\begin{enumerate}[label=(\roman*)]
\item $K_{\mathcal U}^{-1}\le\mathcal U\le K_{\mathcal U}$ on $\R^3\times\R^3$, and hence
$K_{\mathcal U}^{-1}F_{\rm G}\le F_0\le K_{\mathcal U}F_{\rm G}$;
\item $\log\mathcal U$ extends holomorphically to $\{|\mathrm{Im}\,y|<\Lambda_\Phi\}\times
\mathbb S_\Lambda$, with $|\log\mathcal U|\le C_{\mathcal U}$ there;
\item for real $(y,v)$ and multi-indices $a,b$ with $|a|+|b|\ge1$, with $\partial_y^a$ taken in any
orthonormal frame of $\R^3_y$ (in particular $|\nabla_y\log\mathcal U|\le2C_{
\mathcal U}/\Lambda_\Phi$, in the frame of $\nabla_y\log\mathcal U$),
\begin{equation}
|\partial_y^a\partial_v^b\log\mathcal U(y,v)|\le C_{
\mathcal U}\,a!\,b!\,(\Lambda_\Phi/2)^{-|a|}R^{-|b|},\qquad|\partial_y^a\partial_v^b\,\mathcal U(y,v)|\le
K_{\mathcal U}e^{C_{\mathcal U}}\,a!\,b!\,(\Lambda_\Phi/2)^{-|a|}R^{-|b|}\,\mathcal U(y,v);
\end{equation}
\item $\jap v|\nabla_v\log\mathcal U|\le L$.
\end{enumerate}
\end{lem}

\begin{proof}
Put $J(v):=c_{\rm G}/c_0(v)=\int F_{\rm G}(x)e^{-\Phi_{\rm v}(x,v)}\dd x$, by (D1). For real $v$,
$e^{-C_{\rm pv}}\le J\le e^{C_{\rm pv}}$, and $\mathcal U=J^{-1}e^{-\Phi_{\rm v}}$. Together
these imply~(i).
For $v\in\mathbb S_\Lambda$ and real $x$, $\mathrm{Re}\,e^{-\Phi_{\rm v}}=e^{-\mathrm{Re}\,\Phi_{\rm
v}}\cos(\mathrm{Im}\,\Phi_{\rm v})\ge e^{-C_{\rm pv}}\cos1$. Hence $J$ is holomorphic on $
\mathbb S_\Lambda$, with $\mathrm{Re}\,J\ge e^{-C_{\rm pv}}\cos1$ and $|J|\le e^{C_{\rm pv}}$, and its
principal logarithm satisfies $|\log J|\le C_{\rm pv}+\log(1/\cos1)+\frac\pi2$. Hence $\log
\mathcal U=-\log J-\Phi_{\rm v}$ is holomorphic on $\{|\mathrm{Im}\,y|<\Lambda_\Phi\}\times\mathbb S_\Lambda$
and bounded by $C_{\mathcal U}$. This proves~(ii). The polydiscs of radii $\Lambda_\Phi/2$ in $y$
and
$R=\Lambda/2$ in $v$ about a real point, in any orthonormal frame of $\R^3_y$, lie in this domain,
since $|\mathrm{Im}\,y|<\sqrt3\Lambda_\Phi/2$ and $|\mathrm{Im}\,v|<\sqrt3R$ there, as in
Subsection~\ref{ssec:reldens}. Then the Cauchy estimates yield the first bound of~(iii), and the
second follows in the same way from $|\mathcal U|\le e^{C_{\mathcal U}}$ and~(i). Finally
$\nabla_v\log\mathcal U=\nabla_v\log F_0$, so~(iv) is part of (D4).
\end{proof}

\emph{Convention.} In this section and in Subsections~\ref{ssec:tools}--\ref{ssec:drift}, the
profile $F_0$ often appears only as a majorant or a minorant, evaluated at an integration point
$(y'',v'')$. We then write $F_0(y'')$ for $F_0(y'',v'')$.
 By Lemma~\ref{lem:Ufac}(i), every such
bound holds with $F_{\rm G}(y'')$ in place of $F_0(y'')$, up to the factor $K_{\mathcal U}$, which
enters the constants. In particular $F_{\rm G}\le K_{\mathcal U}C_Fe^{-\beta_F|y|^2}$, by (D1).

\begin{lem}[Null structure]\label{lem:Bnull}
We have $a(w)y=a(w)(y+tw)$ for all $w,y$. Consequently $A\hat y\cdot\hat y\le C\eps\jap y^{-2}$ for
$|y|\ge1$.
\end{lem}

\begin{proof}
The identity follows from $a(w)w=0$. For the consequence, $a(w)\hat y\cdot\hat
y=|y|^{-2}a(w)y''\cdot y''$ with $y''=y+tw$, and $|y''|^2$ is controlled by the Gaussian factor
$F_0(y'')$.
\end{proof}

Therefore a radial derivative, which produces a factor $t\jap y$ in \eqref{eq:Geq}, is paired with a
radial gain $\jap y^{-1}$. This suggests that we use energies on boxes adapted to the dispersion.

\subsection{Adapted boxes}\label{ssec:boxes}

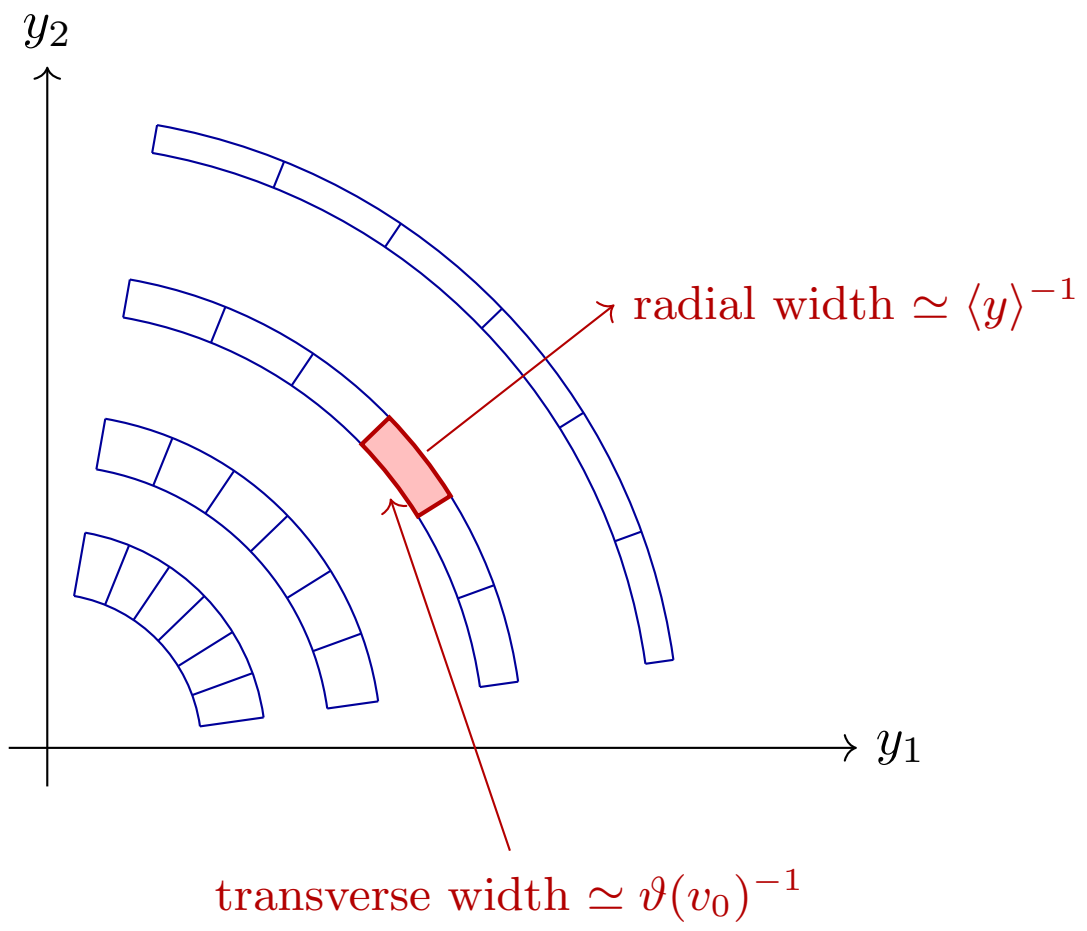
\begin{figure}[ht]
\centering
\begin{tikzpicture}[scale=0.85]
\draw[->] (-0.3,0)--(6.3,0) node[right]{$y_1$};
\draw[->] (0,-0.3)--(0,5.3) node[above]{$y_2$};
\foreach \r/\w in {1.2/0.5,2.2/0.4,3.4/0.3,4.7/0.22}{
 \foreach \a in {8,20,32,44,56,68,80}{
  \draw[blue!60!black,thin] ({\r*cos(\a)},{\r*sin(\a)}) -- ({(\r+\w)*cos(\a)},{(\r+\w)*sin(\a)});
 }
 \draw[blue!60!black,thin] (8:\r) arc (8:80:\r);
 \draw[blue!60!black,thin] (8:\r+\w) arc (8:80:\r+\w);
}
\fill[red!25] (32:3.4) arc (32:44:3.4) -- (44:3.7) arc (44:32:3.7) -- cycle;
\draw[red!70!black,thick] (32:3.4) arc (32:44:3.4) -- (44:3.7) arc (44:32:3.7) -- cycle;
\draw[->,red!70!black] (38:3.75)--(38:5.6) node[right,font=\small]{radial width $\simeq\jap
y^{-1}$};
\draw[->,red!70!black] (3.6,-0.8) node[below,font=\small,align=center]{transverse width
$\simeq\vartheta(v_0)^{-1}$} -- (36:3.3);
\end{tikzpicture}
\caption{Boxes adapted to the dispersion (schematic). The radial width shrinks like $\jap y^{-1}$,
so that a radial derivative produces a factor $\jap y$, which is paired with the radial gain of the
null
structure (Lemma~\ref{lem:Bnull}). The transverse width $\vartheta(v_0)^{-1}$ comes from the
transverse moments of the Gaussian $F_0$ along the collision manifold
(Lemma~\ref{lem:RBperp}).}\label{fig:boxes}
\end{figure}

We estimate the initial layer by energies on boxes adapted to the dispersion. Let
$\vartheta(v)=1+\sqrt{\log\jap v}$ and $h_r(y,v)=\jap y+\sqrt{\log\jap v}$. For each velocity ball
$B(v_0,1)$, cover $\R^3_y$ by boxes $Q$ of radial width $c/h_r$ and transverse width $c/\vartheta$,
with frames adapted to the direction $\hat y$, as in Figure~\ref{fig:boxes}. On each box the
energies are
\begin{equation}
a_{\mathbf k}=W_{\mathbf k}\|u_{\mathbf k}\varphi_Q\|_{L^2},\qquad
u_{\mathbf k}=\partial^\alpha_Q\partial^\beta_vG,\qquad
W_{\mathbf k}=\delta_y^{|\alpha|}h_r^{-a_r}\vartheta^{-a_\perp}\tau_\alpha\jap{v_0}^\sigma
w_\beta(\rad(t)),
\end{equation}
where $a_r,a_\perp$ count the radial and transverse derivatives, $\tau_\alpha=(h_r\vartheta^2)^{-1}$
on the two top $y$-levels, $\sigma=|\gamma+2|/2$, and $w_\beta$ are the analytic weights. On the
three top $y$-levels the weight is further multiplied by $(|\beta|+1)^{-1/2}$, and on the levels
$\ge N_y-4$ by $\vartheta(v_0)^{-1}$.

The coefficient $h$ in a commutator is evaluated at the points $(y+tw,v-w)$, which lie on a
translate of the box $Q$, a \emph{slab}, inside a box $P''$ of a different shape. The following
lemmas control the weights of $P''$ relative to those of $Q$.

\begin{lem}[Transverse moments]\label{lem:RBperp}
Assume (D$_{\rm iso}$). Let $e$ be a unit vector with $|e\cdot y|\le2$, and $g=f_0G$ with
$\frac12\le G\le2$. Then for all $m\ge0$ and $t>0$,
\begin{equation}
\int a(w)\xi\cdot\xi\,|e\cdot(y+tw)|^m g(y+tw,v-w)\dd w\le
C^m\big(\vartheta(v)^m+m^{m/2}\big)A\xi\cdot\xi ,
\end{equation}
with $C$ independent of $t,y,v,e,\xi$.
\end{lem}

\begin{proof}[Idea of the proof]
Slice $w$ along $e$. On each slice the quadratic factor $|w\times\xi|^2$ is controlled everywhere by
its average over a core interval of length $\simeq1/t$, and the tail is controlled by the Gaussian
factor of $F_0$ in the direction $e$, which (D$_{\rm iso}$) makes explicit. The only loss that grows
with $v$ comes from the ratio of the values of $\rho_0$ (the ratio of the values of $\mathcal U$ is
bounded by $K_{\mathcal U}^2$), which is $\le C\jap v^p$ and yields the factor $\vartheta(v)$. The
complete proof is given in Lemma~\ref{lem:RBdir} below, which has the same statement with
$\Theta_\perp=\vartheta$.
\end{proof}

\begin{lem}[Slabs]\label{lem:BKslab}
This is Lemma~\ref{lem:slab}. If $S$ is a slab of the shape of $Q$ that meets $P''$, then
$\int_S|u|^2\le CL\frac{|S|}{|P''|}\max_P\int_{\tilde P}|u|^2$, the maximum over the boxes $P$ of
the family of $P''$ that meet $S$, where $L$ is the ratio of the widths of $P''$ and $S$. If $u$ has
two more $y$-derivatives in the hierarchy, the factor $L$ can be dropped.
\end{lem}

The slab factor $L$ appears only at the two top $y$-levels, where the top weight $\tau$ of both the
target and the coefficient compensates it, up to moments that, by Lemma~\ref{lem:FPslab}, are
bounded uniformly in $t$. The velocity weight $\jap{v_0}^\sigma$ serves this purpose: without it,
the transverse weights of the coefficient box give powers of $\log\jap{v}$ that are not bounded
uniformly.

\subsection{The closed estimate}\label{ssec:closed}

The main result on the initial layer is Theorem~\ref{thm:layer}, and we outline its proof after the
statement.

\begin{theo}[The initial layer]\label{thm:layer}
Let $\gamma\in[-3,-2)$, $f_0\in\Dc$, $N_y\ge20$. There is $\eps_0'(\Dc)>0$ such that for
$\eps\le\eps_0'$ the solution exists on $[0,t_1]$, the box energies satisfy $
\mathbb E(t)^{1/2}\le2C_4\eps t$, and at $t_1$ the solution satisfies the hypotheses of
Theorem~\ref{thm:t1}.
\end{theo}

\begin{proof}[Structure of the proof]
\emph{The a priori estimate.} We make an energy estimate on each box for the equation of $u_{
\mathbf k}$ obtained by differentiating \eqref{eq:Geq}. We classify the commutator terms by which
factor is low or high in $v$ and in $y$. Each of the four configurations is estimated by a trilinear
bound
adapted to it, either \emph{relative} (compared with $A$) or \emph{absolute} (integrated against
$F_0\rho_0$).

We write the drift term $\bb$ either directly, with the kernel $b$, or in the form
$\bb[F]=\ab\langle DF\rangle$, depending on the configuration, so that a Cauchy--Kovalevskaya loss
occurs only in two corner cases at the top $y$-levels. There, half of the loss is absorbed by the
weakening $(|\beta|+1)^{-1/2}$ of the top weights, and the other half, as in Lemma~\ref{lem:FPhalf},
by the decrease of the radius on the target's own box.

Factors that are not paired with the weight of the target are put on the side of the estimate where
the target is an energy. By Lemma~\ref{lem:FPYside}, they are bounded by the trace of $\ab$, which
is an absolute bound. The remaining crossings into the top levels are controlled by
Lemma~\ref{lem:RBperp}
and its radial version.

\emph{Construction and transfer.} We construct the solution by a vanishing-viscosity approximation
of frozen linear problems, with a viscosity adapted to the boxes and analytic in $v$, and an
iteration that contracts in $L^\infty$. The transfer to the leafwise norms of Theorem~\ref{thm:t1}
at $t_1$ uses the velocity weight to make the bounds uniform in $\xi$.

The complete proof occupies the rest of this section. We collect the tools in
Subsection~\ref{ssec:tools}. The energies on adapted boxes and the closed estimate,
Proposition~\ref{prop:FPclose}, are in Subsection~\ref{ssec:layer}, and the $\bb$-part of the drift
at the top $y$-levels, Proposition~\ref{prop:FPfix}, is treated in Subsection~\ref{ssec:drift}.
\end{proof}

\subsection{Tools for the initial layer}\label{ssec:tools}

We collect the tools used in the initial layer: comparison of coefficients, localization, the weight
terms, the adapted covering and its frames, crude bounds at all finite orders, moments along the
collision manifold, the $\ab$-form of $\bb$, the trilinear estimates and the algebra of the analytic
weights on boxes, the Sobolev and null-structure estimates, the neighbor sums and the slab
estimates, and the Gaussian hypodissipativity estimates.

\begin{lem}[Energy identity without ellipticity]\label{lem:LWenergy}
Let $h\ge0$ and let $\Psi\ge0$ be a smooth weight on $\R^6_{y,v}$. Then for real $u$,
\begin{multline}
\iint D\cdot\big(\ab[h]Du-\bb[h]u\big)\,u\Psi=-\iint\ab[h]Du\cdot Du\,\Psi-\iint(\ab[h]Du\cdot
D\Psi)\,u\\
+\frac12\iint(\bb[h]\cdot D\Psi)u^2-\frac12\iint\cb[h]u^2\Psi .
\end{multline}
Consequently, if $\Psi>0$,
\begin{equation}
\iint D\cdot\big(\ab[h]Du-\bb[h]u\big)\,u\Psi\le-\frac12\iint\ab[h]Du\cdot
Du\,\Psi+\iint\Big(\frac{\ab[h]D\Psi\cdot D\Psi}{2\Psi}+\frac{|\bb[h]\cdot
D\Psi|}{2}+\frac{|\cb[h]|\Psi}{2}\Big)u^2 .
\end{equation}
\end{lem}

\begin{proof}
The vector field $D$ is constant on $\R^6$, so it integrates by parts without boundary terms. We
write $\iint\bb u\cdot Du\,\Psi=\frac12\iint\bb\cdot D(u^2)\Psi$ and integrate by parts once more,
using $D\cdot\bb=\cb$. For the bound, note that $\ab[h]\ge0$ as a matrix because $h\ge0$. Then the
Cauchy--Schwarz inequality for this form implies $|\ab Du\cdot D\Psi|\,|u|\le\frac12\ab Du\cdot
Du\,\Psi+\frac12(\ab
D\Psi\cdot D\Psi/\Psi)u^2$.
\end{proof}

\begin{lem}[Coefficient comparison]\label{lem:R3comp}
Let $0<\beta_-\le\beta_+\le2\beta_-$, and suppose
\begin{equation}
c_1\rho_0(v')e^{-\beta_+|y'|^2}\le g(y',v')\le c_2\rho_0(v')e^{-\beta_-|y'|^2}.
\end{equation}
Put $A:=\ab[g]$ and
\begin{equation}
\Theta(y,v):=\jap{R_\Theta(y,v)}+1,\qquad R_\Theta(y,v)^2:=8|y|^2+\frac{2p}{\beta_-}\log\jap
v+C_\Theta,
\end{equation}
where $p$ is the exponent of (D3) and the constant $C_\Theta\ge0$ is fixed in the proof. Then for
$t\le t_1$, all $(y,v)$ and all $\zeta$,
\begin{equation}
\int a(w)\zeta\cdot\zeta\,\jap{y+tw}\,g(y+tw,v-w)\dd w\le\Theta(y,v)\,A\zeta\cdot\zeta,\qquad
c\,h_r\le\Theta\le C\,h_r,
\end{equation}
with $h_r=\jap y+\sqrt{\log\jap v}$ as in Subsection~\ref{ssec:boxes}. Here $C_\Theta$ depends only
on an upper bound for $c_2/c_1$ and on $\beta_\pm$, $t_1$, $\gamma$, $\Dc$, and $c,C$ depend only on
$p/\beta_-$ and $C_\Theta$. None of them depends on $\eps$ when $c_1,c_2$ are proportional to
$\eps$.
Under (D$_{\rm iso}$), $g=f_0G$ with $\frac12\le G\le2$ satisfies the hypothesis with
$\beta_-=\frac23\beta_F$, $\beta_+=\frac43\beta_F$ and $c_2/c_1=4K_{
\mathcal U}^2e^{3C_\Phi^2/(2\beta_F)}$, and so does $f_0$. In
Subsections~\ref{ssec:tools}--\ref{ssec:drift}, $\Theta$ is this one function, with $C_\Theta$
chosen for this ratio. Consequently, if $|h|\le M\jap
y g$ pointwise, then
\begin{equation}
|\ab[h]\zeta\cdot\zeta'|\le M\Theta\,(A\zeta\cdot\zeta)^{1/2}(A\zeta'\cdot\zeta')^{1/2}.
\end{equation}
\end{lem}

\begin{proof}
Let $k(w)=a(w)\zeta\cdot\zeta\ge0$, $y'=y+tw$, and $J=\int k\,g(y',v-w)\dd w=A\zeta\cdot\zeta$.

\emph{Lower bound for $J$.} For $|w|\le1$ we have $\rho_0(v-w)\ge e^{-L}\rho_0(v)$ and
$|y'|^2\le2|y|^2+2t_1^2$. Moreover $\int_{B_1}a(w)\zeta\cdot\zeta\dd w=c_\gamma|\zeta|^2$ with
$c_\gamma>0$, since $\gamma+2>-3$. Hence we have $J\ge c\,c_1\rho_0(v)e^{-2\beta_+|y|^2}|\zeta|^2$,
with
$c=c(\gamma,L,\beta_+,t_1)>0$ (the factor $e^{-L}e^{-2\beta_+t_1^2}c_\gamma$).

\emph{The far region.} Since $\jap{y'}e^{-\beta_-|y'|^2}\le Ce^{-\beta_-\varrho^2/2}$ for
$|y'|>\varrho$, and $\int|w|^{\gamma+2}\rho_0(v-w)\dd w\le C$, we have
\begin{equation}
I_\varrho:=\int_{|y'|>\varrho}k\jap{y'}g(y',v-w)\dd w\le C\,c_2|\zeta|^2e^{-\beta_-\varrho^2/2}.
\end{equation}
Using $\rho_0(v)\ge A_F^{-1}\jap v^{-p}$, we deduce that
$I_\varrho/J\le C(c_2/c_1)e^{-\beta_-\varrho^2/2+2\beta_+|y|^2}\jap v^p$. Here the factor $c_2/c_1$
comes from the upper bound for $g$ in
$I_\varrho$ and the lower bound in $J$. For $\varrho=R_\Theta(y,v)$ the exponent is
$(2\beta_+-4\beta_-)|y|^2-\beta_-C_\Theta/2\le-\beta_-C_\Theta/2$, since $\beta_+\le2\beta_-$, and
we fix $C_\Theta$ so that $C(c_2/c_1)e^{-\beta_-C_\Theta/2}\le1$. Then $I_{R_\Theta}\le J$ and $\int
k\jap{y'}g\le\jap{R_\Theta}J+I_{R_\Theta}\le\Theta J$. Finally, $\max\{1,R_\Theta\}\le\Theta\le
R_\Theta+2$ and $\max\{\sqrt8|y|,(2p/\beta_-)^{1/2}(\log\jap v)^{1/2}\}\le
R_\Theta\le\sqrt8|y|+(2p/\beta_-)^{1/2}(\log\jap v)^{1/2}+C_\Theta^{1/2}$ imply $c\,h_r\le\Theta\le
C\,h_r$.

\emph{Under (D$_{\rm iso}$).} Since $|\nabla\Phi_1|\le C_\Phi$, $|\Phi_1(y)-\Phi_1(0)|\le
C_\Phi|y|\le\frac{\beta_F}3|y|^2+\frac{3C_\Phi^2}{4\beta_F}$. Thus \begin{equation}c_{\rm
G}e^{-\Phi_1(0)}e^{-3C_\Phi^2/(4\beta_F)}e^{-\frac43\beta_F|y|^2}\le F_{\rm G}(y)\le c_{\rm
G}e^{-\Phi_1(0)}e^{3C_\Phi^2/(4\beta_F)}e^{-\frac23\beta_F|y|^2}.\end{equation} Multiplying by
$\eps\rho_0(v)\,\mathcal U(y,v)\,G$, and using $K_{\mathcal U}^{-1}\le\mathcal U\le K_{\mathcal U}$
(Lemma~\ref{lem:Ufac}(i)) and $\frac12\le G\le2$, we obtain the claim.

\emph{The consequence.} Since $a(w)\ge0$,
$|a(w)\zeta\cdot\zeta'|\le(a(w)\zeta\cdot\zeta)^{1/2}(a(w)\zeta'\cdot\zeta')^{1/2}$. With $|h|\le
M\jap{y'}g$ at the integration point and the Cauchy--Schwarz inequality in $w$, we get
\begin{equation}
|\ab[h]\zeta\cdot\zeta'|\le M\Big(\int a(w)\zeta\cdot\zeta\jap{y'}g\Big)^{1/2}\Big(\int
a(w)\zeta'\cdot\zeta'\jap{y'}g\Big)^{1/2}\le M\Theta\,(A\zeta\cdot\zeta)^{1/2}(A\zeta'\cdot\zeta')^{1/2}.
\end{equation}
\end{proof}

\emph{Unit balls.} Lemmas~\ref{lem:ULloc} and~\ref{lem:ULtri} are stated on unit balls of
$\R^6_{y,v}$, and we use them on adapted boxes in Subsections~\ref{ssec:layer}
and~\ref{ssec:drift}, as explained in Remark~\ref{rem:ULbox}. We fix $\varphi\in C^\infty_c(B_2)$
with $0\le\varphi\le1$,
$\varphi=1$ on $B_1$ and $|\nabla\varphi|\le2$, where $B_r$ is the ball of radius $r$ about $0$ in
$\R^6$, and put $B_r(z_0):=z_0+B_r$. For $z_0=(y_0,v_0)\in\R^6$ and a function $u$ on $\R^6$ put
\begin{equation}
\chi_{z_0}:=\varphi(\cdot-z_0)^2,\qquad\Theta_{z_0}:=\jap{y_0}+\sqrt{\log\jap{v_0}}=h_r(y_0,v_0),\qquad\|u\|_{\rm
ul}:=\sup_{z'\in\R^6}\|u\|_{L^2(B_4(z'))}.
\end{equation}
Since $t\le t_1$, $|D\varphi(\cdot-z_0)|\le C(t_1)\mathbf 1_{B_2(z_0)}$. For $z\in B_4(z_0)$ we have
$h_r(z)\le Ch_r(z_0)$ and $h_r(z_0)\le Ch_r(z)$. Hence, by Lemma~\ref{lem:R3comp}, we have
$C^{-1}\Theta_{z_0}\le\Theta\le C\Theta_{z_0}$ on $B_4(z_0)$.

\begin{lem}[Localization and lower-order terms]\label{lem:ULloc}
For every $z_0\in\R^6$ and every real $u$,
\begin{equation}
\iint\big[D\cdot(ADu)+(2A\zeta-\bb)\cdot Du+S_0u\big]u\,\chi_{z_0}\le-\frac12\iint A\,Du\cdot
Du\,\chi_{z_0}+C\eps\iint_{B_2(z_0)}u^2 .
\end{equation}
\end{lem}

\begin{proof}
We apply Lemma~\ref{lem:LWenergy} with $\Psi=\chi_{z_0}$. Then we have $A D\chi\cdot
D\chi/\chi=4AD\varphi\cdot D\varphi\le C|A|\mathbf 1_{B_2}$, and $|A|,|\bb|\le C\eps$ uniformly. For
$\gamma=-3$ we also have $|\cb|=8\pi g\le C\eps$.

For the drift, we integrate by parts and obtain
\begin{equation}
\iint(2A\zeta\cdot Du)u\chi=-\iint\big(D\cdot(A\zeta)\chi+A\zeta\cdot D\chi\big)u^2 .
\end{equation}
Here $|D\cdot(A\zeta)|,|A\zeta|\le C\eps$ by Lemma~\ref{lem:GMwt}, which uses the null structure of
Lemma~\ref{lem:Bnull} and does not depend on the weight. We treat the term with $\bb$ in the same
way, using $D\cdot\bb=\cb$. Finally $|S_0|\le C\eps$ by part~(i) of the weight lemma.
\end{proof}

\begin{lem}[Localization costs $O(\eps)$]\label{lem:RBloc}
Let $\gamma\in[-3,-2)$. For a cutoff $\chi_Q$ adapted to $Q$, $|D\sqrt{\chi_Q}|^2_A\le
C\eps\,\mathbf 1_{\tilde Q}$. (At $\gamma=-2$ the bound is $C\eps\ell_Q\mathbf 1_{\tilde Q}$, see
Subsection~\ref{ssec:endpoint}, item~(2).)
 Here $|\xi|_A^2:=A\xi\cdot\xi$ and $\tilde Q$ is a fixed
dilate of $Q$.
\end{lem}

\begin{proof}
We have $D\sqrt\chi=\partial_v\sqrt\chi-t\big(O(\Theta_r)\hat y+O(\Theta_\perp)e\big)$. The radial
part contributes
\begin{equation}
t^2\Theta_r^2A\hat y\cdot\hat y\le Ct^2(\jap y^2+\log\jap v)\,\eps\jap y^{-2}\cdot\min(1,\jap
v^{\gamma+2}\jap y^2),
\end{equation}
using $A\hat y\cdot\hat y\le C\eps|y|^{-2}$ (null structure) and $|A|\le C\eps\jap v^{\gamma+2}$.
The transverse part contributes
$t^2\Theta_\perp^2|A|\le C\eps\log\jap v\,\jap v^{\gamma+2}$. Both are
$\le C\eps$, since $\gamma+2<0$. The $v$-part is $O(\eps)$ as before.

\end{proof}

\begin{lem}[Weight terms via the null structure]\label{lem:GMwt}
Write $\Psi=\Psi_{v_0}$. Then
\begin{equation}
\frac{D\Psi}{\Psi}=\zeta+2D\log\chi_{v_0}-\frac{2Pt\,y}{\jap y^2}.
\end{equation}
Moreover:
\begin{enumerate}[label=(\roman*)]
\item $A\zeta\cdot\zeta\le C\eps$ and $|A\zeta|\le C\eps$;
\item $|\bb[g]\cdot\zeta|\le C\eps$, using $b(w)\cdot y=b(w)\cdot(y+tw)+2t|w|^{\gamma+2}$;
\item $|D\cdot(A\zeta)|\le C\eps$.
\end{enumerate}
Consequently the weight terms of Lemma~\ref{lem:LWenergy} with weight $\Psi$ are bounded by
$\frac14\iint A\,Du\cdot Du\,\Psi+C\eps\iint|u|^2\Psi$, plus neighboring-ball terms from
$D\log\chi_{v_0}$.
\end{lem}

\begin{proof}
(i) Under the bound $G\le2$ we have $g\le2f_0$. By (D$_{\rm iso}$) and~\eqref{eq:FGU}, $\zeta=D\log
f_0=\zeta_1+2\beta_Ft\,y$ with
\begin{equation}\label{eq:zeta1}
\zeta_1(y,v)=\nabla\log\rho_0(v)+t\nabla\Phi_1(y)+D\log\mathcal U(y,v),
\end{equation}
which implies $|\zeta_1|\le L_\zeta:=2L+t_1\big(C_\Phi+2C_{\mathcal U}/\Lambda_\Phi\big)$, by
Lemma~\ref{lem:Ufac}(iii)--(iv). Since $|a(w)|\le|w|^{\gamma+2}$ and
$\int|w|^{\gamma+2}\rho_0(v-w)\dd w\le C$, we have $|A|\le C\eps$. By Lemma~\ref{lem:Bnull}, we
have $a(w)y=a(w)y''$ with $y''=y+tw$. Thus $|Ay|\le\int|w|^{\gamma+2}|y''|g(y'',v-w)\dd w\le C\eps$
and $Ay\cdot y\le C\eps$, since the Gaussian factor $F_0(y'')$ controls $|y''|^2$. Hence we get
$|A\zeta|\le|A||\zeta_1|+2\beta_Ft_1|Ay|\le C\eps$ and
$A\zeta\cdot\zeta\le2A\zeta_1\cdot\zeta_1+8\beta_F^2t^2Ay\cdot y\le C\eps$.

(ii) Since $b(w)=-2|w|^\gamma w$ and $y=y'-tw$, we have $b(w)\cdot y=b(w)\cdot y'+2t|w|^{\gamma+2}$.
Hence
\begin{equation}
|\bb\cdot y|\le|\bb[y'g]|+2t\int|w|^{\gamma+2}g\le C\eps .
\end{equation}
The bounded parts of $\zeta$ are immediate.

(iii) We have $D\cdot(A\zeta)=\bb\cdot\zeta+A:D\zeta$, and $D\zeta$ is bounded, by the same bounds
as in~(i) for the part $D\log\mathcal U$ of $\zeta_1$.

The consequence follows from Lemma~\ref{lem:LWenergy} and $\ab D\Psi\cdot D\Psi/\Psi\le
C\,(A\zeta\cdot\zeta+\mathrm{tr}A)\Psi$.
\end{proof}

\begin{lem}[Adapted covering]\label{lem:BKcov}
Fix $v_0\in\mathbb Z^3$ and write $\vartheta=\vartheta(v_0)$, $h_r(y)=h_r(y,v_0)$. There is a family
$\{Q\}$ of boxes in $\R^3_y$, with centers $y_Q$ and orthonormal frames $(\hat e_Q,e^1_Q,e^2_Q)$,
such that:
\begin{enumerate}[label=(\roman*)]
\item each $Q$ is either \emph{outer}, with $|y_Q|\ge2\vartheta$, or \emph{inner}, with
$|y_Q|<3\vartheta$; an outer $Q$ has width $\ell_r(Q)=c_{\rm box}/h_r(y_Q)$ along $\hat e_Q$ and
width $\ell_\perp(Q)=c_{\rm box}/\vartheta$ along $e^{1,2}_Q$, and an inner $Q$ is a cube of side
$\ell_r(Q)=\ell_\perp(Q):=c_{\rm box}/(5\vartheta)\in[c_{\rm box}/(5h_r(y_Q)),c_{\rm
box}/h_r(y_Q)]$; all widths are $\le c_{\rm box}$;
\item if $Q$ is outer, then $\hat e_Q=\hat y_Q$; if $Q$ is inner, $(\hat
e_Q,e^1_Q,e^2_Q):=(e_3,e_1,e_2)$ for a fixed orthonormal basis $(e_1,e_2,e_3)$;
\item the dilates $\tilde Q:=4Q$ have overlap at most $N_0$, and $\bigcup Q=\R^3$; we write $Q'\sim
Q$ (neighbors) if $|v_0(Q)-v_0(Q')|\le2$ and $\tilde Q\cap\tilde Q'\ne\emptyset$, and we take $N_0$
so large that, in addition, every $Q$ has at most $N_0$ neighbors;
\item for $y\in\tilde Q$: $h_r(y)\le2h_r(y_Q)$ and $h_r(y_Q)\le2h_r(y)$; if $Q$ is outer,
$|e^i_Q\cdot y|\le2c_{\rm box}/\vartheta\le2$ for $i=1,2$, $\hat e_Q\cdot y\ge\frac78|y_Q|$, and the
angle between $\hat y$ and $\hat e_Q$ is at most $4c_{\rm box}/(\vartheta|y_Q|)\le C/(\vartheta\jap
y)$; if $Q$ is inner, $|y|\le4\vartheta$;
\item if $\tilde Q\cap\tilde Q'\ne\emptyset$, where $Q'$ may belong to the family of any
$v_0'\in\mathbb Z^3$ with $|v_0'-v_0|\le2$, then $|y_Q-y_{Q'}|\le1$, $\ell_r(Q)\simeq\ell_r(Q')$,
and the frames satisfy $|\hat e_Q-\hat e_{Q'}|\le C\,(\vartheta\jap{y_Q})^{-1}$ when $Q$ and $Q'$
are outer.
\end{enumerate}
All constants depend only on $c_{\rm box}$, which is a small absolute constant, $c_{\rm
box}\le\frac18$ (used in Lemma~\ref{lem:EPtransit}(ii)).
\end{lem}

\begin{proof}
\emph{Inner region} $|y|<3\vartheta$. Here $\vartheta\le h_r\le4\vartheta$. We use the cubes of the
lattice $(c_{\rm box}/(5\vartheta))\mathbb Z^3$ whose centers lie in this region, with the fixed
frame. Their side $c_{\rm box}/(5\vartheta)$ is comparable to $c_{\rm box}/h_r(y_Q)$. They cover
$\{|y|<2\vartheta\}$, and the outer boxes below cover $\{|y|\ge2\vartheta\}$.

\emph{Outer region} $|y|\ge2\vartheta$. We use polar coordinates $y=r\omega$, and write $h_r(r)=\jap
r+\sqrt{\log\jap{v_0}}$. First we define the shells by $r_{k+1}=r_k+c_{\rm box}/(2h_r(r_k))$,
starting from $r_0=2\vartheta$. Since $|\partial_rh_r|\le1$ and $h_r\ge1$, consecutive shells have
comparable thickness. Then, on the shell $[r_k,r_{k+1}]$, we cover $\mathbb S^2$ by caps of angular
radius $c_{\rm box}/(4\vartheta r_k)$, with bounded overlap. This is possible since $c_{\rm
box}/(4\vartheta r_k)\le c_{\rm box}$. For a cap with center $\omega_Q$, let $Q$ be the box with
center $y_Q:=\frac12(r_k+r_{k+1})\omega_Q$, frame $\hat e_Q=\omega_Q$, $e^{1,2}_Q$ tangent at
$\omega_Q$, and the widths of~(i). It contains the product of the shell and the cap. Indeed,
transversally we have
$r_{k+1}\sin(c_{\rm box}/(4\vartheta r_k))\le\frac{c_{\rm box}}{4\vartheta}\big(1+\frac{c_{\rm
box}}{2h_r(r_k)r_k}\big)\le\frac{c_{\rm box}}{2\vartheta}$, and radially the product lies within
$\frac{c_{\rm box}}{4h_r(r_k)}+\frac{c_{\rm box}^2}{32\vartheta^2r_k}\le\frac{c_{\rm
box}}{2(1+c_{\rm box}/4)h_r(r_k)}\le\frac{c_{\rm box}}{2h_r(y_Q)}$
of $|y_Q|$, since $h_r(y_Q)\le(1+\frac{c_{\rm box}}4)h_r(r_k)$, $h_r(r_k)\le\frac32r_k$ and $c_{\rm
box}\le\frac18$.

Within $\tilde Q$, since $e^i_Q\cdot y_Q=0$, we have $|e^i_Q\cdot
y|=|e^i_Q\cdot(y-y_Q)|\le2c_{\rm box}/\vartheta\le2$. The widths are $\le c_{\rm box}$, so
$\operatorname{diam}\tilde Q\le4\sqrt3c_{\rm box}\le1$. This, together with $|\nabla h_r|\le1$ and
$h_r\ge1$, implies the comparability of $h_r$ in~(iv), in both regions. In the inner region
$|y|\le3\vartheta+1\le4\vartheta$ on $\tilde Q$. For~(v), let $y\in\tilde Q\cap\tilde Q'$. Then
$|y_Q-y_{Q'}|\le\frac12(\operatorname{diam}\tilde Q+\operatorname{diam}\tilde Q')\le1$, so
$h_r(y_Q)\simeq h_r(y_{Q'})$, and $\ell_r(Q)\simeq\ell_r(Q')$ if both are outer. If one is inner,
then $\vartheta\le h_r\le9\vartheta$ on $\tilde Q\cup\tilde Q'$ (also for $|v_0-v_0'|\le2$, since
$\vartheta(v_0')\le\vartheta(v_0)+\sqrt{\log3}$), so both $\ell_r$ are comparable to
$c_{\rm box}/\vartheta$. If both are outer, $|e^i_Q\cdot y|\le2c_{\rm box}/\vartheta$ and $\hat
e_Q\cdot y\ge|y_Q|-\frac14\ge\frac78|y_Q|$, so the angle between $\hat y$ and $\hat e_Q$ is
$\le4c_{\rm box}/(\vartheta|y_Q|)$. This holds for every outer $Q$, and, since $\jap y\le2|y_Q|$
on $\tilde Q$, it proves the angle bound of~(iv). The same holds for $Q'$, and
$|y_{Q'}|\ge|y_Q|-1\ge\frac12|y_Q|$, which gives $|\hat e_Q-\hat e_{Q'}|\le
C(\vartheta\jap{y_Q})^{-1}$.

\emph{The two regions} overlap in $2\vartheta\le|y|<3\vartheta$, where both families have sizes
$\simeq c_{\rm box}/\vartheta$. Bounded overlap: if $y\in\tilde Q'$ for an outer $Q'$, then
$\big||y|-|y_{Q'}|\big|\le2\ell_r(Q')+8c_{\rm box}^2/(\vartheta^2|y|)\le C\ell_r(Q')$ and the angle
between $\hat y$ and $\hat e_{Q'}$ is $\le4c_{\rm box}/(\vartheta|y_{Q'}|)$. Since
$h_r(y_{Q'})\simeq h_r(y)$ and consecutive shells have comparable thickness, this leaves $O(1)$
shells and, on each, $O(1)$ caps (the caps have bounded overlap). Moreover, there are no more than
$5^3$ inner cubes with $y\in\tilde Q'$. Thus the bounded overlap~(iii) holds. For neighbors, there
are $33$ points
$v_0'\in\mathbb Z^3$ with $|v_0'-v_0|\le2$, and for each of them the boxes of the family of $v_0'$
that meet $\tilde Q$ have widths comparable to those of $Q$, since
$|\vartheta(v_0')-\vartheta(v_0)|\le\sqrt{\log3}$. Their number is bounded by a volume count, and
$N_0$ can be taken so large that every $Q$ has at most $N_0$ neighbors.
\end{proof}

\begin{lem}[Cutoffs and frames]\label{lem:BKframe}
\leavevmode
\begin{enumerate}[label=(\roman*)]
\item There are $\varphi_Q\in C^\infty_c(\tilde Q\times B(v_0,2))$ of the form
$\varphi_Q(y,v)=\psi_Q(y)\chi(v-v_0)$, with one fixed $\chi\in C^\infty_c(B(0,2))$, $0\le\chi\le1$,
$\chi=1$ on $B(0,1)$, with $\varphi_Q=1$ on $2Q\times B(v_0,1)$ and
\begin{equation}
|\partial_{\hat e_Q}\varphi_Q|\le Ch_r,\qquad|\partial_{e_Q}\varphi_Q|\le
C\vartheta,\qquad|\nabla_v\varphi_Q|\le C .
\end{equation}
\item For neighbors $Q,Q'$, let $R_{QQ'}$ be the rotation taking the frame of $Q'$ to that of $Q$.
A derivative $\partial^{\alpha}_{Q'}$ of order $a\le N_y$ is a combination of the derivatives
$\partial^{\alpha'}_Q$ with $|\alpha'|=a$. The coefficients satisfy, in the weights
$h_r^{-a_r}\vartheta^{-a_\perp}$ of Subsection~\ref{ssec:boxes},
\begin{equation}
\max\Big\{\frac{h_r^{-a'_r}\vartheta^{-a'_\perp}}{h_r^{-a_r}\vartheta^{-a_\perp}},\
\frac{h_r^{-a_r}\vartheta^{-a_\perp}}{h_r^{-a'_r}\vartheta^{-a'_\perp}}\Big\}\,|\mathrm{coeff}|\le
C^{a}.
\end{equation}
\end{enumerate}
\end{lem}

\begin{proof}
(i) We take the products of one-dimensional cutoffs in the frame coordinates of $Q$, at the scales
$\ell_r$ and $\ell_\perp$.

(ii) A transverse vector $e'$ of $Q'$ has radial component $|e'\cdot\hat e_Q|\le
C/(\vartheta\jap{y_Q})$ in the frame of $Q$, by Lemma~\ref{lem:BKcov}(v). Re-expressing it, we have
\begin{equation}\label{eq:BK-cost}
\frac{|e'\cdot\hat e_Q|\,h_r}{\vartheta}\le\frac{C(\jap{y_Q}+\vartheta)}{\vartheta^2\jap{y_Q}}\le C
.
\end{equation}
The radial vector of $Q'$ has a transverse component of the same size, which yields
$|\cdot|\,\vartheta/h_r\le C$. The ratios in the other direction are bounded in the same way.
Indeed,
a cross coefficient is multiplied by $|\cdot|\,\vartheta/h_r\le|\cdot|\le1$ or
$|\cdot|\,h_r/\vartheta\le C$, by
$\vartheta\le h_r$ and~\eqref{eq:BK-cost}. If $Q$ and $Q'$ are both inner, the frames coincide. If
one of them is inner, then $|y_Q|,|y_{Q'}|\le3\vartheta+2$, so $\vartheta\le h_r\le9\vartheta$ on
$\tilde Q\cup\tilde Q'$, for either choice of $v_0$ in $h_r$ and $\vartheta$ (since
$\vartheta(v_0')\le\vartheta(v_0)+\sqrt{\log3}$ for $|v_0-v_0'|\le2$). A derivative of order $a$ is
a sum of at most $3^a$ products of rotation coefficients, each $\le1$ in modulus, and each
re-expressed derivative contributes a factor not exceeding $\max(h_r/\vartheta,\vartheta/h_r)\le9$.
Hence the bound holds with $C=27$.
\end{proof}

\begin{prop}[All finite orders]\label{prop:Bn}
(This proposition is needed only for the variant in Lemma~\ref{lem:FPg}, which the main line of the
proof does not use.) Let $G^j$ be the iterates of Lemma~\ref{lem:FPiter} and $G$ their limit. There
is $C$, depending only on $\Dc$ and $t_1$, and for every $n$ there is $\eps_n(t_1)\in(0,\eps_0]$,
nonincreasing in $n$, such that for $\eps\le\eps_n(t_1)$, all $j$ and all $t\in[0,t_1]$, every
$H\in\{G^j,G\}$ is of class $C^n$ in $(y,v)$ and
\begin{equation}
|\nabla^k_{y,v}(H-1)|\le C^{k+1}\,k!\,\eps t\qquad(0\le k\le n).
\end{equation}
\end{prop}

\begin{proof}
We write $z=(y,v)$, $\partial_i$ ($1\le i\le6$) for the derivatives in $z$,
$u_I:=\partial_{i_1}\cdots\partial_{i_k}u$ for $I\in\{1,\dots,6\}^k$, $|\nabla^ku|^2:=\sum_Iu_I^2$,
and $|X|^2_{\mathcal M}:=\mathcal MX\cdot X$. Here $c_0,C$ depend only on $\Dc$ and $t_1$. Fix $n$.
We show $|\nabla^k(H-1)|\le\Gamma\varkappa^kk!\,\eps t$ for $k\le n+1$, with $\Gamma,\varkappa$
fixed below, by induction on $j$ ($G^0=1$). For $G$ it follows by uniform convergence and the
Arzel\`a--Ascoli theorem.

\emph{(a) Explicit factors.} The function $\Phi_1$ extends holomorphically to
$\{|\mathrm{Im}\,y|<\Lambda_\Phi\}$, with $|\Phi_1(y+a+i\eta)-\Phi_1(y)|\le C_\Phi(|a|+|\eta|)$ for
real $a,\eta$. By the Cauchy estimates on polydiscs of radius $\min\{\Lambda_\Phi,1\}/4$, we have
$|\partial_y^aF_{\rm G}|\le C^{|a|+1}a!\,e^{\beta_F|y|}F_{\rm G}$. Using~\eqref{eq:rho0Cauchy},
Lemma~\ref{lem:Ufac}(iii) and the Leibniz rule, we get $|\partial^\nu f_0|\le
C^{|\nu|+1}\nu!\,e^{\beta_F|y|}f_0$ for every $\nu$. The Gaussian absorbs
$e^{\beta_F|y''|}\jap{y''}^2$, so, as in
Lemma~\ref{lem:FPdecay}(i), $\int|k(w)|\jap{y''}^2|\partial^\nu f_0(z)|\dd w\le
C^{|\nu|+1}\nu!\,\eps$ for $k\in\{a,b,c\}$ (a local term for $c$ at $\gamma=-3$).

\emph{(b) Coefficients.} Let $\frac12\le H\le2$, $g_0:=2$, $g_l:=\sup|\nabla^lH|$ for $l\ge1$, and
let $A,V,S_0$ be built from $f_0H$. By Leibniz, $D\cdot A=\bb$, the null structure
$a(w)\zeta=a(w)\zeta^\sharp$ with $|\zeta^\sharp|\le C\jap{y''}$, and $|\partial^\mu\zeta|\le
C^{|\mu|}\mu!R^{-|\mu|}$ for $\mu\ne0$ from~(ii) of the same lemma, each of
$|\nabla^iA|,|\nabla^i\bb|,|\nabla^iV|,|\nabla^iS_0|$ is bounded by $C\eps\sum_{l\le i}\binom
il(i-l)!\,C^{i-l}g_l$. If $g_l\le\varkappa^ll!$ for $1\le l\le i$ and $\varkappa\ge2C$, this is $\le
c_0\eps\,i!\,\varkappa^i$, with $c_0$ independent of $\varkappa$.

\emph{(c) The equation.} Let $G^j$ satisfy the bound. The coefficient factor $\tilde G_\nu$ of
Definition~\ref{def:FPreg} with $\tilde G=G^j$ is an average of values of $G^j$ at earlier times, so
it satisfies the same bound. Let $G_\nu$ be the whole-space solution of
Proposition~\ref{prop:FPvisc}(ii), $\nu\le\nu_n$, where $\nu_n\le\eps$ is so small that
$|\nabla^lM_\nu|\le\eps C^{l+1}l!$ for $l\le n+3$ (enough, since $\nu\to0$ below), and $u:=G_\nu-1$.
We bound $\nabla^ku$ by induction on $k$. The case $k=0$ is part~(ii) of that proposition if
$\Gamma\ge K_0e^{t_1}$. For $k\ge1$ we put $W:=|\nabla^ku|^2$, $
\mathcal L_0:=D\cdot(AD\,\cdot)+\nabla\cdot(M_\nu\nabla\,\cdot)+V\cdot D$, $
\mathfrak B:=\sum_I\big(|Du_I|^2_A+|\nabla u_I|^2_{M_\nu}\big)$. Differentiating the equation, we
have
\begin{equation*}
\partial_tW=\mathcal L_0W-2\mathfrak B+2\sum_Iu_I\Big(\sum_{m=1}^k\big[\partial_{i_m}A:D^2u_{I\setminus
m}+\partial_{i_m}M_\nu:\nabla^2u_{I\setminus m}\big]+L_I+F_I\Big),
\end{equation*}
where $L_I$ is linear in $\nabla^ku$, with coefficients among
$\nabla^2A,\nabla^2M_\nu,\nabla\bb,\nabla(\nabla\cdot M_\nu),\nabla V,S_0$, so that
$2|\sum_Iu_IL_I|\le c_0k^2\eps\varkappa^2W$ by~(b), and where $F_I$ contains $\partial_IS_0$ and
$\nabla^lu$ with $l<k$ only.

\emph{(d) Glaeser.} If $\phi\ge0$ on $\R$ and $|\phi''|\le K_2$, then $\phi'(0)^2\le2K_2\phi(0)$.
For $d\times d$ matrices $\mathcal M(z)\ge0$ with $|\nabla^2\mathcal M|\le K_2$, a symmetric
$N=\sum_a\nu_ae_a\otimes e_a$ and $\xi\in\R^6$, we apply it to $\phi(s)=
\mathcal M(z+s\xi/|\xi|)e_a\cdot e_a$ and use Cauchy--Schwarz over $a$, which implies
$|\mathrm{tr}(\partial_\xi
\mathcal M\,N)|\le|\xi|(2dK_2)^{1/2}\mathrm{tr}(\mathcal MN^2)^{1/2}$. Since $u_I$ is symmetric in $I$, the
$A$-part of the commutator sum is $k\sum_J\mathrm{tr}(\partial_{\xi_J}A\,D^2u_J)$, with $|J|=k-1$,
$\xi_J:=(u_{Ji})_{i\le6}$ and $\sum_J|\xi_J|^2=W$. Moreover,
$\mathrm{tr}(A(D^2u_J)^2)
=\sum_b|DD_bu_J|_A^2\le(1+t^2)\sum_i|Du_{Ji}|_A^2$, since
$D_b=\partial_{v_b}-t\partial_{y_b}$. The same holds for $M_\nu$, with $\nabla$ in place of $D$ and
$d=6$. By Cauchy--Schwarz over $J$ and $2ab\le a^2+b^2$, the commutator sum is at most $
\mathfrak B+6k^2KW$, with $K:=(1+t_1^2)\sup|\nabla^2A|+2\sup|\nabla^2M_\nu|\le C_1\eps$ by~(b).

\emph{(e) Comparison.} By~(b), the induction hypotheses and $\sum_{l=3}^k\binom kl l!\,(k-l+2)!\le
k!\,k^3$, we get $\rho(t):=\sup_z\big(\sum_IF_I^2\big)^{1/2}\le
c_0\eps\,k!\,\varkappa^k\big(1+c_0\Gamma\eps t_1\varkappa^2k^3\big)$. Hence we have $\partial_tW\le
\mathcal L_0W+\lambda W+2\rho W^{1/2}$, $\lambda:=6k^2K+c_0k^2\eps\varkappa^2$. Here $W$ is
continuous on $[0,t_1]\times\R^6$ and smooth for $t>0$, $W(0)=0$, and $W\le C(k,\nu)\jap z^{2k}$ by
the local bounds of that proposition. The coefficients of $\mathcal L_0$ are bounded, so
$\mathcal L_0\jap z^{2k+2}\le C_k\jap z^{2k+2}$. Let $Y_\eta'=\frac\lambda2Y_\eta+\rho$,
$Y_\eta(0)=\eta>0$, and $\varpi:=\delta e^{\Lambda't}\jap z^{2k+2}$. At a first zero of
$\sup_z(W-Y_\eta^2-\varpi)$, the bound $W^{1/2}\le Y_\eta+\varpi/(2\eta)$ implies
$\partial_t\varpi-\mathcal L_0\varpi\le(\lambda+\sup\rho/\eta)\varpi$, which is false for
$\Lambda'>C_k+\lambda+\sup\rho/\eta$. Letting $\delta\to0$ and then $\eta\to0$, we obtain
$W^{1/2}\le te^{\lambda t/2}\sup_{[0,t]}\rho$.

\emph{(f) Closing.} Let $\eps_n$ be so small that $\Gamma\eps_nt_1\le1$,
$c_0\Gamma\eps_nt_1\varkappa^2(n+1)^3\le1$ and $\lambda t_1\le2$ for $k\le n+1$. Then
$|\nabla^ku|\le2ec_0\,k!\,\varkappa^k\eps t$, and we take $\Gamma:=\max\{2ec_0,K_0e^{t_1}\}$,
$\varkappa:=2C$. These bounds are uniform in $\nu\le\nu_n$. By part~(iv) of the same proposition and
the Arzel\`a--Ascoli theorem (the bound of order $n+1$ yields convergence in $C^n_{\rm loc}$), they
pass to $G^{j+1}$ for $k\le n$, and also for $k=n+1$, as a Lipschitz bound of $\nabla^nG^{j+1}$,
which is what the induction uses. Finally we take $C:=\max\{\Gamma,\varkappa\}$. This completes the
proof.
\end{proof}

\begin{lem}[Moments along the collision manifold]\label{lem:ULmom}
Under the hypotheses of Lemma~\ref{lem:R3comp}, with $\Theta$ and $A=\ab[g]$ as there, for all
$m\ge0$, $t\in(0,t_1]$, $(y,v)$ and $\xi$,
\begin{equation}
\int a(w)\xi\cdot\xi\,\Theta(y+tw,v-w)^m\,g(y+tw,v-w)\dd w\le
C^m\big(\Theta(y,v)^m+m^{m/2}\big)\,A\xi\cdot\xi ,
\end{equation}
with $m^{m/2}:=1$ for $m=0$, where $C$ depends only on $c_2/c_1$, $\beta_\pm$, $p$, $L$, $A_F$,
$t_1$ and $\gamma$. The proof uses $p>\frac72$, which holds since $p\ge4$ in
Definition~\ref{def:data}.
\end{lem}

\begin{proof}
By homogeneity we may take $|\xi|=1$. We put $k(w):=a(w)\xi\cdot\xi$, so that $0\le
k(w)\le|w|^{\gamma+2}$, write $z=(y',v')=(y+tw,v-w)$, and let $J:=\int k\,g(z)\dd w=A\xi\cdot\xi$.
Below, $\Theta=\Theta(y,v)$ and $R_\Theta=R_\Theta(y,v)$ are those of Lemma~\ref{lem:R3comp}, and
$c,C$ depend only on the quantities listed in the statement. We first let $m\ge1$. Then
$m^{m/2}\ge1$, and every constant $K\ge1$ satisfies $K\le K^m$, so constants can be put into
$C^mm^{m/2}$. The case $0\le m<1$ is reduced to the case $m=1$ at the end of the proof.

\emph{Lower bound for $J$.} By the proof of that lemma, we have $J\ge
c\,c_1\rho_0(v)e^{-2\beta_+|y|^2}$. Since $\rho_0(v)\ge A_F^{-1}e^{-p\log\jap v}$ and
$2\beta_+\le4\beta_-$, we get
\begin{equation}
2\beta_+|y|^2+p\log\jap v\le\frac{\beta_-}2\Big(8|y|^2+\frac{2p}{\beta_-}\log\jap
v\Big)\le\frac{\beta_-}2R_\Theta^2\le\frac{\beta_-}2\Theta^2,
\end{equation}
and hence
\begin{equation}
J\ge c\,c_1e^{-\beta_-R_\Theta^2/2}\ge c\,c_1e^{-\beta_-\Theta^2/2}.
\end{equation}

\emph{Reduction.} Since $\Theta(z)\le
R_\Theta(z)+2\le\sqrt8|y'|+(2p/\beta_-)^{1/2}(\log\jap{v'})^{1/2}+C_\Theta^{1/2}+2$, we have
$\Theta(z)^m\le C^m\big(\jap{y'}^m+(\log\jap{v'})^{m/2}\big)$. Hence it suffices to show that
\begin{equation}
J_1:=\int k\jap{y'}^mg(z)\dd w\qquad\text{and}\qquad J_2:=\int k\,(\log\jap{v'})^{m/2}g(z)\dd w
\end{equation}
are both $\le C^m(\Theta^m+m^{m/2})J$.

\emph{(a) The $y'$-moment.} On $\{w:|y'|\le R_\Theta\}$, we have
$\jap{y'}^m\le\jap{R_\Theta}^m\le\Theta^m$, so this part of $J_1$ is $\le\Theta^mJ$. On
$\{w:|y'|>R_\Theta\}$ we use $g(z)\le
c_2\rho_0(v-w)e^{-\beta_-|y'|^2}$, the bound $e^{-\beta_-|y'|^2/2}\le e^{-\beta_-R_\Theta^2/2}$
there, and, with $u=\jap s^2$,
\begin{equation}
\sup_{s\in\R}\jap s^me^{-\beta_-s^2/2}=e^{\beta_-/2}\sup_{u\ge1}u^{m/2}e^{-\beta_-u/2}\le
e^{\beta_-/2}\Big(\frac m{e\beta_-}\Big)^{m/2}\le C^mm^{m/2}.
\end{equation}
Since $\int|w|^{\gamma+2}\rho_0(v-w)\dd w\le C$ (as $p>3$), it follows from the lower bound for $J$
that this part is $\le C\,c_2C^mm^{m/2}e^{-\beta_-R_\Theta^2/2}\le C(c_2/c_1)C^mm^{m/2}J$. Note that
no integral in $y'$ is taken, so the Jacobian $\dd y'=t^3\dd w$, which degenerates as $t\to0$, does
not
enter.

\emph{(b) The $v'$-moment.} Let $E:=\{w:\log\jap{v-w}>\beta_-\Theta^2\}$. Off $E$, we have
$(\log\jap{v'})^{m/2}\le\beta_-^{m/2}\Theta^m$, so this part of $J_2$ is
$\le\beta_-^{m/2}\Theta^mJ$. On $E$, we have $g(z)\le c_2\rho_0(v-w)\le c_2A_F\jap{v-w}^{-p}$ and
$\jap{v-w}^{-1/2}\le e^{-\beta_-\Theta^2/2}$, which implies
\begin{equation}
\int_Ek\,(\log\jap{v'})^{m/2}g(z)\dd w\le c_2A_Fe^{-\beta_-\Theta^2/2}\int|w|^{\gamma+2}(\log\jap{v-w})^{m/2}\jap{v-w}^{-(p-\frac12)}\dd
w .
\end{equation}
On $\{|w|\le1\}$, we have
$(\log\jap u)^{m/2}\jap u^{-(p-\frac12)}\le\sup_{s\ge0}s^{m/2}e^{-(p-\frac12)s}\le m^{m/2}$, and
$\int_{|w|\le1}|w|^{\gamma+2}\dd w\le C$, since $\gamma+2>-3$. On $\{|w|>1\}$, we have
$|w|^{\gamma+2}\le1$, since $\gamma+2\le0$. Using $u=v-w$, polar coordinates $r=|u|$ and, for
$r\ge1$, the variable $s=\log\jap r$ (so that $\dd r\le\sqrt2\jap r\dd s$ there), we obtain
\begin{equation}
\int_{\R^3}(\log\jap u)^{m/2}\jap u^{-(p-\frac12)}\dd u\le C+4\sqrt2\pi\int_0^\infty
s^{m/2}e^{-(p-\frac72)s}\dd s=C+\frac{4\sqrt2\pi\,\Gamma(\frac m2+1)}{(p-\frac72)^{m/2+1}}\le
C^mm^{m/2},
\end{equation}
since $p>\frac72$ and $\Gamma(\frac m2+1)\le\max\{1,m^{m/2}\}$. By the lower bound for $J$, we have
$e^{-\beta_-\Theta^2/2}\le C c_1^{-1}J$, so this part of $J_2$ is $\le C(c_2/c_1)C^mm^{m/2}J$.

Adding (a) and (b) and using the reduction, we deduce that $\int k\,\Theta(z)^mg(z)\dd w\le
C^m(J_1+J_2)\le C^m(\Theta^m+m^{m/2})J$. This proves the lemma for $m\ge1$.

\emph{The case $0\le m<1$.} For $m=0$ the left side is $J$, and the claim holds since
$\Theta^0+m^{m/2}=2$. Let $0<m<1$, and let $K_1$ be the constant of the case $m=1$, so that $\int
k\,\Theta(z)g(z)\dd w\le K_1(\Theta+1)J$. The measure $\mu:=k(w)g(z)\dd w=a(w)\xi\cdot\xi\,g(z)\dd
w$ is nonnegative, with total mass $J$. By H\"older's inequality for $\mu$, with the exponents
$\frac1m$ and $\frac1{1-m}$, we have

\begin{equation}
\int k\,\Theta(z)^mg(z)\dd w\le\Big(\int k\,\Theta(z)g(z)\dd w\Big)^mJ^{1-m}\le
K_1^m(\Theta+1)^mJ\le(2K_1)^m\Theta^mJ,
\end{equation}
since $\Theta\ge1$. Therefore the lemma holds for every $m\ge0$, with $\max\{C,2K_1\}$ in place of
$C$.
\end{proof}

\emph{The $\ab$-form.} We write $z=(y'',v'')=(y+tw,v-w)$. For a scalar $F$ and a vector field
$\Xi(y,v,w)$ we put $\ab\langle F,\Xi\rangle:=\int a(w)\,\Xi\,F(z)\dd w$, and
$\ab\langle\Xi\rangle:=\int a(w)\,\Xi\dd w$. Under (D$_{\rm iso}$), $\zeta=D\log
f_0=\zeta_1+2\beta_Fty$, with $\zeta_1$ as in~\eqref{eq:zeta1}. Define
\begin{equation}
\zeta^\sharp:=\zeta_1(y,v)+2\beta_Fty'',\qquad\Delta\zeta:=\zeta_1(z)-\zeta_1(y,v).
\end{equation}
Then $\zeta(y,v)-\zeta^\sharp=-2\beta_Ft^2w$ and $\zeta(z)-\zeta(y,v)-\Delta\zeta=2\beta_Ft^2w$.
Since $|\zeta_1|\le L_\zeta$, we have $|\zeta^\sharp|\le L_\zeta+2\beta_Ft_1|y''|\le C\jap{y''}$ and
$|\Delta\zeta|\le2L_\zeta$.

\begin{lem}[$\ab$-form]\label{lem:ULaform}
The following hold, with the notation just introduced.
\begin{enumerate}[label=(\roman*)]
\item $a(w)\zeta(y,v)=a(w)\zeta^\sharp$ and $a(w)\zeta(z)=a(w)(\zeta(y,v)+\Delta\zeta)$.
\item $\bb[h]=\ab\langle Dh\rangle$, that is, $\bb_i[h]=\int a_{ij}(w)(D_jh)(z)\dd w$.
\item For $F=\partial^\nu f_0\cdot h$,
\begin{equation}
\bb[F]=\ab\langle F,\zeta^\sharp+\Delta\zeta\rangle+\sum_{0\ne\nu_1\le\nu}\binom\nu{\nu_1}\ab\big\langle\partial^{\nu-\nu_1}f_0\,h,\
\partial^{\nu_1}\zeta(z)\big\rangle+\ab\langle\partial^\nu f_0,\ Dh(z)\rangle .
\end{equation}
Here $\partial^{\nu_1}\zeta$ is bounded with analytic constants for $\nu_1\ne0$.
\item For every $h$, $\ab[|h|]\zeta\cdot\zeta\le C\int|w|^{\gamma+2}\jap{y+tw}^2|h(z)|\dd w$.
\end{enumerate}
\end{lem}

\begin{proof}
(i) The difference $\zeta(y,v)-\zeta^\sharp=-2\beta_Ft^2w$ lies in the kernel of $a(w)$. Similarly,
$\zeta(z)-\zeta(y,v)-\Delta\zeta=2\beta_Ft^2w$.

(ii) We have $b_i=\partial_{w_j}a_{ij}$ and $\partial_w[h(z)]=-(Dh)(z)$. Then we integrate by parts
in $w$. The boundary terms vanish, since $|w|^{\gamma+2}|w|^2\to0$ as $w\to0$ for $\gamma>-4$.

(iii) We use $D\partial^\nu f_0=\partial^\nu(f_0\zeta)$ and expand by the Leibniz rule. The term
without derivatives on $\zeta$ is treated with~(i).

(iv) This follows from~(i).
\end{proof}

\emph{Notation for the trilinear estimates.} In Lemma~\ref{lem:ULtri}, (D$_{\rm iso}$) holds and
$g=f_0G$ with $\frac12\le G\le2$, the bootstrap hypothesis of Subsection~\ref{ssec:layer}. Thus
$\ab[f_0]\le2A$ as quadratic forms, and $f_0$ and $g$ satisfy the hypothesis of
Lemma~\ref{lem:R3comp}, with the one function $\Theta$ defined there. We fix $z_0\in\R^6$ and write
$\chi:=\chi_{z_0}$ and $\varphi:=\varphi(\cdot-z_0)$, so that $\chi=\varphi^2$. For a multi-index
$\nu=(\nu_y,\nu_v)$ with $|\nu_y|\le N_y$, a function $h$ and vector fields $X,Y$ on $\R^6$, we put
\begin{equation}
I:=\iint\ab[\partial^\nu f_0\,h]X\cdot Y\,\chi,\qquad\omega_\nu:=\frac{|\partial^\nu
f_0|}{f_0},\qquad\bar\omega_\nu(z):=C_{N_y}C^{|\nu_v|}\nu_v!\,R^{-|\nu_v|}\jap{z_y}^{|\nu_y|}.
\end{equation}
Here $\ab[\partial^\nu f_0h]$, $X$ and $Y$ are evaluated at $(y,v)$, and we write $AXX:=AX\cdot X$.
In the applications $X=Du_{\mathbf k-\mu}$ and $Y=Du_{\mathbf k}$, as in
Subsection~\ref{sssec:FPset}. By the first claim of Lemma~\ref{lem:FPexp} (in any fixed orthonormal
frame, with its radius $R$ and its constants), $\omega_\nu\le\bar\omega_\nu$. The majorant
$\bar\omega_\nu$ is smooth, with $|\nabla^j\bar\omega_\nu|\le C_j\bar\omega_\nu$, whereas
$\omega_\nu^{1/2}$ need not be differentiable where $\partial^\nu f_0$ vanishes. In the Sobolev
steps below we differentiate only $\bar\omega_\nu$. By Lemma~\ref{lem:ULmom}, applied to $f_0$ with
the exponent $2m+|\nu_y|$, and by $\jap{z_y}\le\Theta(z)$, $\Theta\ge1$, $\ab[f_0]\le2A$ and
$\Theta\le C\Theta_{z_0}$ on $B_4(z_0)$, we have, as quadratic forms,
\begin{equation}
\ab\big[\bar\omega_\nu\Theta^{2m}f_0\big]\le M_{\nu,m}\,\Theta(y,v)^{2m}A\quad\text{on
}B_4(z_0),\qquad M_{\nu,m}:=C_{N_y,m}\,C^{|\nu_v|}\nu_v!\,R^{-|\nu_v|}\,\Theta_{z_0}^{|\nu_y|},
\end{equation}
with $C_{N_y,m}$ nondecreasing in $m$. Thus $M_{\nu,m}$ depends on $z_0$ only through the factor
$\Theta_{z_0}^{|\nu_y|}$, which comes from the $y$-derivatives of $F_0$.

\begin{lem}[Trilinear estimates]\label{lem:ULtri}
In the setting just described, let $p\ge0$ and $\tilde h:=\Theta^{-p}h$.
\begin{enumerate}[label=(T\arabic*)]
\item \emph{Pointwise $h$:}
\begin{equation}
|I|\le C^pM_{\nu,p}\,\|\tilde h\|_{L^\infty}\,\Theta_{z_0}^{p}\,\Big(\iint
AXX\chi\Big)^{1/2}\Big(\iint AYY\chi\Big)^{1/2}.
\end{equation}
\item \emph{Pointwise plain $X$, $h$ in $L^2$:}
\begin{equation}
|I|\le C_{\nu,p}\,\eps^{1/2}\sup_{B_4(z_0)}|X|\;\|\tilde h\|_{\rm ul}\Big(\iint AYY\chi\Big)^{1/2},
\end{equation}
with no factor $\Theta_{z_0}$.
\item \emph{$h$ in $L^\infty_vL^2_y$, $X$ in $L^\infty_y$ with the $A$-weight:}
\begin{equation}
|I|\le CM_{\nu,p}\,\Theta_{z_0}^{p}\,\|\partial_v^{\le2}\tilde h\|_{\rm
ul}\Big(\sum_{|\gamma|\le2}\Theta_{z_0}^{2(2-|\gamma|)}\iint_{B_4(z_0)}A\,\partial_y^\gamma
X\cdot\partial_y^\gamma X\Big)^{1/2}\Big(\iint AYY\chi\Big)^{1/2}.
\end{equation}
\item \emph{$h$ in $L^\infty_yL^2_v$, $X$ in $L^\infty_v$ with the $A$-weight:}
\begin{equation}
|I|\le CM_{\nu,p}\,\Theta_{z_0}^{p}\,\|\partial_y^{\le2}\tilde h\|_{\rm
ul}\Big(\sum_{|\gamma|\le2}\iint_{B_4(z_0)}A\,\partial_v^\gamma X\cdot\partial_v^\gamma
X\Big)^{1/2}\Big(\iint AYY\chi\Big)^{1/2}.
\end{equation}
\item \emph{Scalar kernels.} Let $k\ge0$ be a locally integrable kernel, or a multiple of
$\delta_0$, with $\sup_{v_0}\int k(w)\rho_0(v_0-w)^{7/8}\dd w\le C_k$. This covers $|c|$, $|b|$ and
$2t|w|^{\gamma+2}$, since $\rho_0^{7/8}\le C\jap\cdot^{-7/2}$. For $m\ge0$ put $K[F]:=\int
k(w)\jap{z_y}^m|F(z)|\dd w$ (the factor $\jap{z_y}^m$ covers the factor $|y''|$ of
Lemma~\ref{lem:RBgain}(ii)). Then, with $C_{\nu,p}$ depending also on $C_k$ and $m$,
\begin{equation}
\iint K[\partial^\nu f_0\,h]\,|U||U'|\chi\le C_{\nu,p}\,\eps\cdot
\begin{cases}
\|\tilde h\|_\infty\Theta_{z_0}^p\,\|U\|_{L^2(\chi)}\|U'\|_{L^2(\chi)},\\[2pt]
\|\tilde h\|_{\rm ul}\sup_{B_4}|U|\,\|U'\|_{L^2(\chi)},\\[2pt]
\|\partial_y^{\le2}\tilde h\|_{\rm ul}\,\|\sup_v|U|\|_{L^2_y}\,\|U'\|_{L^2(\chi)}.
\end{cases}
\end{equation}
The last line also holds with $y$ and $v$ exchanged. In the last line there is no factor
$\Theta_{z_0}$.
\end{enumerate}
The same bounds hold with $Y$ replaced by $u\,\Xi$, where $u$ is a scalar function (in the
applications $u=u_{\mathbf k}$ of Subsection~\ref{sssec:FPset}) and
$\Xi\in\{\zeta^\sharp,\Delta\zeta,\partial^{\nu_1}\zeta\}$ or $\Xi=e_j$ is a constant unit vector,
and with $\iint AYY\chi$ replaced by $C\eps\|u\|^2_{L^2(\chi)}$. Indeed, by
Lemma~\ref{lem:ULaform}(iv) and
Lemma~\ref{lem:FPdecay}(i), $\int a(w)\Xi\cdot\Xi\,\bar\omega_\nu f_0(z)\dd w\le C_\nu\eps$
(trivially for $\Xi=e_j$, since $|e_j|=1$).
\end{lem}

\begin{proof}
\emph{Common first step.} Since $|\partial^\nu f_0|\le\bar\omega_\nu f_0$ and $a(w)\ge0$,
the Cauchy--Schwarz inequality for the nonnegative kernel $a(w)\bar\omega_\nu f_0(z)\dd w$ yields

\begin{equation}
|\ab[\partial^\nu f_0h]X\cdot Y|\le\big(\ab[\bar\omega_\nu
f_0h^2]XX\big)^{1/2}\big(\ab[\bar\omega_\nu f_0]YY\big)^{1/2},
\end{equation}
and the second factor is bounded by $(M_{\nu,0}AYY)^{1/2}$ on $B_4(z_0)$.

(T1) Since $h^2\le\|\tilde h\|_\infty^2\Theta(z)^{2p}$, the definition of $M_{\nu,p}$ implies
that $\ab[\bar\omega_\nu f_0h^2]XX\le\|\tilde h\|_\infty^2M_{\nu,p}\Theta(y,v)^{2p}AXX$. On the
support of $\chi$, $\Theta(y,v)\le C\Theta_{z_0}$, and $M_{\nu,0}\le M_{\nu,p}$. Then~(T1), with the
constant $C^pM_{\nu,p}$, follows from the Cauchy--Schwarz inequality in $(y,v)$.

(T2) Here we put the whole weight on the $h$-side. By Cauchy--Schwarz for the kernel
$a(w)f_0(z)\dd w$, we have
\begin{equation}
|\ab[\partial^\nu f_0h]X\cdot Y|\le\big(\ab[\bar\omega_\nu^2f_0h^2]XX\big)^{1/2}\big(\ab[f_0]YY\big)^{1/2}\le\big(\ab[\bar\omega_\nu^2f_0h^2]XX\big)^{1/2}(2AYY)^{1/2},
\end{equation}
so the $Y$-side carries no factor $M_{\nu,0}$, and hence no power of $\Theta_{z_0}$. We bound
$\ab[\bar\omega_\nu^2f_0h^2]XX\le|X|^2\,\mathrm{tr}\,\ab[\bar\omega_\nu^2f_0h^2]$ and integrate in
$(y,v)$ against $\chi$ by Fubini. With $z=(y+tw,v-w)$ and $\mathrm{tr}\,a(w)=2|w|^{\gamma+2}$, we
have $\iint\mathrm{tr}\,\ab[F]\chi=2\int F(z)\big(\int|w|^{\gamma+2}\chi(z_y-tw,z_v+w)\dd w\big)\dd
z$
for $F\ge0$, and the inner integral is bounded by $C\jap{v_0-z_v}^{\gamma+2}$, since $\chi$ is
supported in a ball of bounded radius about $z_0=(y_0,v_0)$. The Gaussian $F_0(z_y)$ absorbs
$\bar\omega_\nu^2\Theta(z)^{2p}$ in $z_y$, and the logarithmic part in $z_v$ is controlled by
$\rho_0(z_v)\jap{v_0-z_v}^{\gamma+2}$. Covering $\R^6$ by unit balls, we see that the sum over the
balls of $\sup\big(\eps\rho_0F_0\bar\omega_\nu^2\Theta^{2p}\big)$ is $\le C_{\nu,p}\eps$, since
$\rho_0\le A_F\jap\cdot^{-4}$. Hence the $X$-side is $\le C_{\nu,p}\eps\sup_{B_4(z_0)}|X|^2\|\tilde
h\|^2_{\rm ul}$, and the Cauchy--Schwarz inequality in $(y,v)$ yields~(T2).

(T3) For fixed $(w,v)$, we take the supremum in $y$ of $a(w)(\varphi X)\cdot(\varphi
X)\,f_0\,\bar\omega_\nu\Theta^{2p}(y+tw,v-w)$, times $\int_y\tilde h^2(y+tw,v-w)\dd y$. Then we take
the supremum in $v$ of the second factor, which yields the factor
$\|\partial_v^{\le2}\tilde h\|_{\rm ul}^2$ by Sobolev's inequality $H^2_v\subset L^\infty_v$. The
first factor is controlled by Sobolev's inequality $H^2_y\subset L^\infty_y$, applied to
$a(w)^{1/2}\varphi X\,(f_0\bar\omega_\nu\Theta^{2p})^{1/2}$. Here we use the smooth majorant
$\bar\omega_\nu$, and $\Theta$ is smooth, with $|\nabla\Theta|+|\nabla^2\Theta|\le C$ and
$\Theta\ge1$. Its $y$-derivatives fall on $X$, or on $f_0^{1/2}(y+tw,\cdot)$, where each one gives
a factor $C\jap{y'}$, or on $(\bar\omega_\nu\Theta^{2p})^{1/2}(y+tw,\cdot)$ and on $\varphi$, where
each one produces a factor $C(1+p)$. Integrating in $(w,v,y)$ and using Lemma~\ref{lem:ULmom}, we
can replace $\jap{y'}^{2(2-|\gamma|)}$ by $\Theta_{z_0}^{2(2-|\gamma|)}$.

(T4) We exchange the roles of $y$ and $v$ in~(T3). The $v$-derivatives of $f_0(\cdot,v-w)$ produce
only a factor $C$, since $\jap v|\nabla\log\rho_0|+|\nabla^2\log\rho_0|\le L$ and, by (D4),
$\jap v|\nabla_v\log F_0|+|\nabla_v^2\log F_0|\le L$, and those of
$(\bar\omega_\nu\Theta^{2p})^{1/2}$ produce a factor $C(1+p)$. Thus no power of $\Theta$ appears.

(T5) The kernel is scalar, so no $A$-structure has to be kept on the side of $U$. On the support of
$\chi$, we have $|v-v_0|\le2$, so
\begin{equation}
\jap{z_y}^m|\partial^\nu f_0\,h|(z)\le\jap{z_y}^m\bar\omega_\nu f_0\Theta^p(z)\,|\tilde h(z)|\le
C\eps\rho_0(v-w)^{7/8}|\tilde h(z)|\le C\eps e^{2L}\rho_0(v_0-w)^{7/8}|\tilde h(z)| :
\end{equation}
Here $F_0$ absorbs every power of $\jap{z_y}$, $\rho_0^{1/8}$ absorbs every power of
$\log\jap{v-w}$, and $|\nabla\log\rho_0|\le L$. The factor $\rho_0(v_0-w)^{7/8}$ does not depend on
$(y,v)$, so it can be taken out of the integral in $(y,v)$, and $\int k(w)\rho_0(v_0-w)^{7/8}\dd
w\le C_k$. It remains to bound, for fixed $w$, the integral of $|\tilde h(z)||U||U'|\chi$ in
$(y,v)$. In the first line we use $\|\tilde h\|_\infty$ and Cauchy--Schwarz. In the second, we use
$\sup|U|$ and Cauchy--Schwarz, with $\|\tilde h(\cdot+(tw,-w))\|_{L^2(B_2(z_0))}\le\|\tilde h\|_{\rm
ul}$. In the mixed case, we have
\begin{equation}
\int_y\int_v|\tilde h(z)||U||U'|\chi\le\sup_{y'}\|\tilde
h(y',\cdot)\|_{L^2(B(v_0-w,2))}\int_y\sup_v|U|\,\|U'(y,\cdot)\|_{L^2_v(\chi)} ,
\end{equation}
and the first factor is bounded by $C\|\partial_y^{\le2}\tilde h\|_{\rm ul}$, by Sobolev's
inequality in $y$. The case with $y$ and $v$ exchanged is the same.
\end{proof}

\begin{remark}[Box form]\label{rem:ULbox}
Subsections~\ref{ssec:layer} and~\ref{ssec:drift} use Lemma~\ref{lem:ULtri} on an adapted box $Q$ of
Lemma~\ref{lem:BKcov} instead of a unit ball. There $\chi_{z_0}$ is replaced by the cutoff
$\varphi_Q^2$ of Lemma~\ref{lem:BKframe}, $\omega_\nu$ by the combined factor $\Omega$
of~\eqref{eq:FPOmega0}, the bounds involving $M_{\nu,p}\Theta_{z_0}^p$ by those of
Lemmas~\ref{lem:FPOmega} and~\ref{lem:FPslab}, the norm $\|\cdot\|_{\rm ul}$ by $L^2$ norms on
slabs, and Sobolev's inequality at the unit scale by Sobolev's inequality at the scale of $Q$ of
Lemma~\ref{lem:RBsob}. Each application indicates which of the lines (T1)--(T5) it
uses. Under the bootstrap hypothesis $\frac12\le G\le2$ these lines carry over to boxes as follows.
\begin{itemize}
\item (T2) on $Q$. For fixed $w$ the integration point $z$ runs over the slab $S_w:=(\tilde
Q+tw)\times(B(v_0,2)-w)$, of measure $|\tilde Q||B_2|$, and $\rho_0(z_v)\le e^{2L}\rho_0(v_0-w)$
there. With the weight on the $h$-side, as in the proof of~(T2),
\begin{multline}
\Big|\iint\ab[\partial^\nu f_0h]X\cdot Y\varphi_Q^2\Big|\le C\eps^{1/2}\sup_{\tilde Q\times
B(v_0,2)}|X|\\
\times\Big(\int|w|^{\gamma+2}\rho_0(v_0-w)\int_{S_w}\bar\omega_\nu^2F_0|h|^2\dd z\dd
w\Big)^{1/2}\Big(\iint AYY\varphi_Q^2\Big)^{1/2}.
\end{multline}
This form is used in the proof of Lemma~\ref{lem:FPc}, where the $L^2$ norm of $h$ on slabs is then
estimated by Lemma~\ref{lem:slab}. The normalization $|Q|^{-1}$ matches $|S_w|/|P''|$.
\item (T3), (T4) on $Q$: Sobolev's inequality at the scale of $Q$. A radial derivative of
$F_0^{1/2}(y'')$ gives a factor $(|\hat e_Q\cdot y''|+C)/h_r$, a transverse one
$(|e\cdot y''|+C)/\vartheta$. Their $\ab$-moments are $O(1)$ by Lemma~\ref{lem:FPexp}, that is, by
Lemmas~\ref{lem:ULmom}
and~\ref{lem:RBdir}. Therefore no power of $h_r$ remains. The derivatives of $\bar\omega_\nu$ and
$\Theta$ contribute $O(1)$. The box quantities in the weight ($h_r''$, $\vartheta''$, the angle
of $\hat e_{P''}$, the slab factor) are piecewise constant. In the Sobolev step we replace them by
smoothings at the scale of $Q$, as for $\tilde\Lambda$ in
the drift subsection, which changes $h_r''$ and $\vartheta''$ by bounded factors, since a
window at the scale of $Q$ has diameter $\le1$, $|\nabla h_r|\le1$ and~(iv) of the covering lemma
applies. By the same lemma, it changes the angle of $\hat e_{P''}$ by
$O(1/\jap{y''})$, an additive change which can be put into $|\sin\sigma|\,h_r''\le
C(|y''_\perp|+\vartheta'')$.
\item (T5). The proof is local in $w$ and transfers without change.
\end{itemize}
\end{remark}

\begin{lem}[Convolution bounds]\label{lem:ULconv}
Let $s\ge5$ and $C_1\rad/R\le\frac14$, where $R=\Lambda/2$ is the radius of
Section~\ref{sec:layer} and $C_1=C_1(\Dc)$ is at least twice the larger of the analytic rates of
$f_0$ (Lemma~\ref{lem:FPexp}) and of $\zeta$ (the rate $K_\zeta$ of Lemma~\ref{lem:FPdecay}(ii)), as
fixed in the proof of Lemma~\ref{lem:FPb}, so that in particular $|\partial^\beta_vf_0|\le
C\,C_1^{|\beta|}\beta!\,R^{-|\beta|}f_0$. Here $C_1$ does not depend on $t_1$. For multi-indices
$\beta=\beta_\nu+\beta'+\beta''\in\mathbb N^3$ and shifts $\sigma',\sigma''$ with
$|\sigma'|+|\sigma''|\le3$ (or $\le4$ if $s\ge6$; see Lemma~\ref{lem:FPccross}(ii)), carried only by
the factor of smaller order (or by a factor of $v$-order $\le|\beta|/2$: then
$|\beta'|+|\sigma'|\le(|\beta_\nu|+1)(m+|\sigma'|)$, so~(i) holds with
$(|\beta_\nu|+1)^{s+|\sigma|}$ in place of $(|\beta_\nu|+1)^s$, and~(ii) holds with a larger
$C(s)$), put
\begin{equation}
K:=\frac{\beta!}{\beta_\nu!\,\beta'!\,\beta''!}\,C_1^{|\beta_\nu|}\beta_\nu!\,R^{-|\beta_\nu|}\,\frac{w_\beta}{w_{\beta'+\sigma'}\,w_{\beta''+\sigma''}}
.
\end{equation}
\begin{enumerate}[label=(\roman*)]
\item With $m=\min(|\beta'|,|\beta''|)$,
\begin{equation}
K\le2^s\Big(\frac{C_1\rad}R\Big)^{|\beta_\nu|}(|\beta_\nu|+1)^s\,\rad^{-3}\,\frac{(m+3)^3}{(m+1)^s}.
\end{equation}
\item For nonnegative sequences $A,B,d$ on $\mathbb N^3$,
\begin{equation}
\sum_\beta d_\beta\sum_{\beta_\nu+\beta'+\beta''=\beta}K\,A_{\beta'}B_{\beta''}\le
C(s)\,\rad^{-3}\,\|d\|_{\ell^2}\|A\|_{\ell^2}\|B\|_{\ell^2}.
\end{equation}
\item In the explicit case $\beta'=0$, $\sigma=0$, $\beta_\nu\ne0$, the left side of~(ii) with
$A_0=1$ is bounded by $C(s)\frac{C_1\rad}R\|d\|_{\ell^2}\|B\|_{\ell^2}$.
\end{enumerate}
If the total shift is $|\sigma'|+|\sigma''|=4$ (allowed for $s\ge6$), then in~(i)
$\rad^{-3}(m+3)^3$ becomes $\rad^{-4}(m+4)^4$, and in~(ii) $C(s)\rad^{-3}$ becomes
$C(s)\rad^{-4}$. In general, read $\rad^{-|\sigma|}(m+|\sigma|)^{|\sigma|}$ with
$|\sigma|=|\sigma'|+|\sigma''|$.
\end{lem}

\begin{proof}
(i) Inserting $w_\beta=\rad^{|\beta|}(|\beta|+1)^s/\beta!$, we have
\begin{equation}
K=\Big(\frac{C_1\rad}{R}\Big)^{|\beta_\nu|}\rad^{-|\sigma'|-|\sigma''|}\,\frac{(\beta'+\sigma')!}{\beta'!}\,\frac{(\beta''+\sigma'')!}{\beta''!}\,\frac{(|\beta|+1)^s}{(|\beta'|+|\sigma'|+1)^s(|\beta''|+|\sigma''|+1)^s}.
\end{equation}
Then we use $(\beta'+\sigma')!/\beta'!\le(|\beta'|+3)^{|\sigma'|}$, together with
\begin{equation}
|\beta|+1\le(|\beta_\nu|+1)(|\beta'|+|\beta''|+1)\le2(|\beta_\nu|+1)\max(|\beta'|+1,|\beta''|+1).
\end{equation}

(ii) If $m=|\beta'|$, then $K\le c(\beta_\nu)k(\beta')$. The left side is $\langle
d,c*(kA)*B\rangle$, which is bounded by $\|c\|_1\|k\|_2\|A\|_2\|B\|_2\|d\|_2$, by Young's
inequality on $\mathbb Z^3$. We have $\|c\|_1\le2^s\sum_n(n+1)^{s+2}4^{-n}$, and
$\|k\|_2^2\le\rad^{-6}\sum_n(n+1)^{2}(n+3)^6(n+1)^{-2s}<\infty$ for $s\ge5$. The case
$m=|\beta''|$ is symmetric.

(iii) In this case $K\le(C_1\rad/R)^{|\beta_\nu|}(|\beta_\nu|+1)^s$, and the sum over
$\beta_\nu\ne0$ starts at the first power of $C_1\rad/R$.
\end{proof}

\begin{lem}[Directional moments]\label{lem:RBdir}
Let $\gamma\in[-3,-2]$, let $f_0\in\mathfrak D$ satisfy (D$_{\rm iso}$), let $g=f_0G$ with
$\frac12\le G\le2$, and let $t>0$. Let $e$ be a unit vector with $|e\cdot y|\le2$. Then for all
$m\ge0$ and all $\xi$,
\begin{equation}
\int a(w)\xi\cdot\xi\,|e\cdot y'|^m\,g(z)\dd w\le C^m\big(\Theta_\perp(v)^m+m^{m/2}\big)\,A\xi\cdot\xi,\qquad
y'=y+tw,\ z=(y',v-w).
\end{equation}
Here $\Theta_\perp(v)=1+\sqrt{\log\jap v}$, and $C$ depends only on $\gamma,\beta_F,C_\Phi,C_{\rm
pv},p,A_F$. It is independent of $t$, $y$, $v$, $e$ and $\xi$.
\end{lem}

\begin{proof}
We may take $|\xi|=1$. Then
\begin{equation}
a(w)\xi\cdot\xi=|w|^{\gamma+2}|\Pi(w)\xi|^2=|w|^{\gamma}\,|w\times\xi|^2 .
\end{equation}

\emph{Slicing.} We write $w=se+w^\flat$ with $s=e\cdot w$ and $w^\flat\perp e$, and put
\begin{equation}
\tau:=e\cdot y'=e\cdot y+ts,\qquad y'^\flat:=y'-\tau e=y^\flat+tw^\flat .
\end{equation}
Note that $y'^\flat$ does not depend on $s$. By (D$_{\rm iso}$) and~\eqref{eq:FGU}, we have
\begin{equation}
F_0(z)=c_{\rm G}\,e^{-\beta_F|y'^\flat|^2}e^{-\beta_F\tau^2}e^{-\Phi_1(y')}\,\mathcal U(z),\qquad
|\Phi_1(y')-\Phi_1(y'^\flat+\tau_0e)|\le C_\Phi|\tau-\tau_0| .
\end{equation}
For fixed $w^\flat$, we put

\begin{equation}
k(s):=|w|^\gamma\,Q(s)\,\rho_0(v-w)\,G(z)\,\mathcal U(z)\,e^{-\Phi_1(y')},\qquad
Q(s):=|w\times\xi|^2=|\mathbf as+\mathbf b|^2,
\end{equation}
with $\mathbf a=e\times\xi$ and $\mathbf b=w^\flat\times\xi$. The integrand on the left is $\eps
c_{\rm G}e^{-\beta_F|y'^\flat|^2}\,e^{-\beta_F\tau^2}k(s)\,|\tau|^m$. The factor
$e^{-\beta_F|y'^\flat|^2}$ is common to the two sides. It therefore suffices to prove, uniformly in
$w^\flat$,
\begin{equation}
\int_\R|\tau|^me^{-\beta_F\tau^2}k(s)\dd s\le C^m(\Theta_\perp^m+m^{m/2})\int_\R
e^{-\beta_F\tau^2}k(s)\dd s,\tag{$*$}
\end{equation}
and then to integrate in $w^\flat$.

\emph{Core and tail.} Let $I:=[-1/t,1/t]$. For $s\in I$ we have $|\tau|\le|e\cdot y|+1\le3$. A tail
point $s_1$ is one with $|\tau_1|\ge3$, where $\tau_1=e\cdot y+ts_1$. Then
\begin{equation}
|s_1|\ge(|\tau_1|-2)/t\ge1/t .
\end{equation}
We show that for every tail point
\begin{equation}
k(s_1)\le C\jap v^{p}\jap{\tau_1}^2e^{C_\Phi|\tau_1|}\;t\!\int_Ik(s)\dd s.\tag{$**$}
\end{equation}
We compare each factor of $k$ at $s_1$ with its minimum over $I$.
\begin{itemize}
\item \emph{The factor $Q$.} The interval $I$ is symmetric, so the linear term of $Q$ integrates to
zero:
\begin{equation}
\int_IQ=\frac2t|\mathbf b|^2+\frac2{3t^3}|\mathbf a|^2 .
\end{equation}
Hence
\begin{equation}
Q(s_1)\le2|\mathbf a|^2s_1^2+2|\mathbf b|^2\le3\jap{ts_1}^2\,t\!\int_IQ,
\end{equation}
and $\jap{ts_1}\le3\jap{\tau_1}$.
\item \emph{The factor $|w|^\gamma$.} On $I$ we have $|w|\le1/t+|w^\flat|$. At $s_1$ we have
$|w|\ge\max(|s_1|,|w^\flat|)\ge\frac12(1/t+|w^\flat|)$. Since $\gamma<0$,
\begin{equation}
|w(s_1)|^\gamma\le2^{|\gamma|}\min_I|w|^\gamma .
\end{equation}
\item \emph{The factor $\rho_0$.} We write $V:=v-w^\flat$, so that $v-w=V-se$. Since $w^\flat\perp
e$,
we have $V_e:=e\cdot V=e\cdot v$, and $\jap{V-se}^2=\jap{V^\flat}^2+(V_e-s)^2$. By
$A_F^{-1}\jap\cdot^{-p}\le\rho_0\le A_F\jap\cdot^{-p}$, we get
\begin{equation}
\frac{\rho_0(V-s_1e)}{\min_I\rho_0(V-se)}\le A_F^2\Big(\frac{\jap{V^\flat}^2+(|V_e|+1/t)^2}{\jap{V^\flat}^2+(V_e-s_1)^2}\Big)^{p/2}\le
A_F^2\Big(1+\frac{(|V_e|+1/t)^2}{1+(V_e-s_1)^2}\Big)^{p/2}.
\end{equation}
There are two cases.
\begin{itemize}
\item[--] If $|V_e-s_1|\ge|s_1|/2\ge1/(2t)$, the last fraction does not exceed $2|v|^2+8$.
\item[--] Otherwise $|V_e|>|s_1|/2\ge1/(2t)$, so $1/t\le2|v|$, and the fraction does not exceed
$9|v|^2$.
\end{itemize}
In both cases the ratio is bounded by $C\jap v^p$. This is the reason why $\Theta_\perp$ depends on
$v$ rather than on $t$. A large ratio of $\rho_0$ needs the peak $s\simeq e\cdot v$ to lie in the
tail, which forces $|v|\gtrsim1/t$.
\item \emph{The remaining factors.} We have $G(z_1)\le4\min_IG$, $\mathcal U(z_1)\le K_{
\mathcal U}^2\min_I\mathcal U$ by Lemma~\ref{lem:Ufac}(i), and $e^{-\Phi_1(y'_1)}\le
e^{C_\Phi(|\tau_1|+3)}\min_Ie^{-\Phi_1}$.
\end{itemize}
Multiplying the six bounds and using $\min_I(\cdots)\int_IQ\le\int_Ik$, we obtain~$(**)$.

\emph{Conclusion.} Since $\tau(I)\subset[-3,3]$ and $\dd\tau=t\dd s$, we have
\begin{equation}
\int_\R e^{-\beta_F\tau^2}k\dd s\ge e^{-9\beta_F}\int_Ik\dd s .
\end{equation}
By~$(**)$, for $T\ge3$, we have
\begin{equation}
\int_{|\tau|\ge T}|\tau|^me^{-\beta_F\tau^2}k\dd s\le C\jap v^p\Big(\int_Ik\Big)\int_{|\tau|\ge
T}|\tau|^{m+2}e^{-\beta_F\tau^2+C_\Phi|\tau|}\dd\tau .
\end{equation}
The last integral is bounded by $C^{m}m^{m/2}e^{-\beta_FT^2/4}$ for $T\ge T_0(\beta_F,C_\Phi)$,
where
we take $T_0\ge3$, so that every point with $|\tau|\ge T$ is a tail point. We choose
\begin{equation}
T:=\max\Big\{T_0,\ \big(4p\beta_F^{-1}\log\jap v\big)^{1/2}\Big\}\le C\,\Theta_\perp(v),
\end{equation}
so that $\jap v^pe^{-\beta_FT^2/4}\le1$. The part of~$(*)$ with $|\tau|<T$ is bounded by $T^m\int
e^{-\beta_F\tau^2}k$. The part with $|\tau|\ge T$ is bounded by $C\,C^mm^{m/2}e^{9\beta_F}\int
e^{-\beta_F\tau^2}k$. For $m\ge1$ the factors $C$ and $e^{9\beta_F}$ can be put into $C^m$, and
together this is~$(*)$. For $m=0$, $(*)$ is trivial (with $m^{m/2}:=1$). For $0<m<1$, let $K_1$ be
the constant of the case $m=1$. By H\"older's inequality for the nonnegative measure
$e^{-\beta_F\tau^2}k\dd s$, with the exponents $\frac1m$ and $\frac1{1-m}$, we obtain
$\int|\tau|^me^{-\beta_F\tau^2}k\le K_1^m(\Theta_\perp+1)^m\int
e^{-\beta_F\tau^2}k\le(2K_1)^m\Theta_\perp^m\int e^{-\beta_F\tau^2}k$, since $\Theta_\perp\ge1$.
This is~$(*)$ for every $m\ge0$.
\end{proof}

\begin{lem}[Sobolev's inequality without $\Theta$]\label{lem:RBsob}
For a multi-index $(\alpha,\beta)$, let $u_{(\alpha,\beta)}=\partial_Q^\alpha\partial_v^\beta G$
have weight $W=\delta_y^{|\alpha|}h_r(Q)^{-a_r}\vartheta(Q)^{-a_\perp}w_\beta(\rad)$, with
$\delta_y,\rad\le1$. A factor that is the same at all the levels involved (such as
$\jap{v_0}^\sigma$ or $\omega_Q$) may be included. On a box $Q$, with the supremum over $Q\times
B(v_0,1)$ and the $L^2$ norms over $\tilde Q\times B(v_0,2)$,
\begin{equation}
\sup_Q W|u_{(\alpha,\beta)}|\le C\delta_y^{-2}\rad^{-2}(|\beta|+3)^2\Big(\sum_{|\gamma_y|\le2,\
|\gamma_v|\le2}|\tilde Q|^{-1}\big\|W_{(\alpha+\gamma_y,\beta+\gamma_v)}\,u_{(\alpha+\gamma_y,\beta+\gamma_v)}\big\|^2_{L^2(\tilde
Q)}\Big)^{1/2}.
\end{equation}
No power of $h_r$ or $\vartheta$ appears. If $W$ carries the level-dependent factors $\tau_\alpha$
or those of Definition~\ref{def:FPtop}, the bound holds with these factors taken at the shifted
levels on the right. A shift that crosses one of the thresholds $N_y-4$, $N_y-2$, $N_y-1$ then
produces in addition the factor $\vartheta(v_0)$, $(|\beta|+3)^{1/2}$, respectively $h_r\vartheta^2$
(accounted for in
Lemma~\ref{lem:FPOmega}(ii), Lemma~\ref{lem:FPtopalg}(iii) and Remark~\ref{rem:FPiota}). Expressed
through the energies of the neighbors of $Q$ (Lemma~\ref{lem:BKnb}(i)), the constant is multiplied
by $C^{N_y}N_0^{1/2}$.
\end{lem}

\begin{proof}
We rescale $Q$ to a unit box $\hat Q$. A rescaled radial derivative is $\ell_r\partial_r$ with
$\ell_r\le c_{\rm box}h_r(Q)^{-1}$, and a transverse one is $\ell_\perp\partial_\perp$ with
$\ell_\perp\le c_{\rm box}\vartheta(Q)^{-1}$ (with equality on outer boxes, while on inner boxes the
side $c_{\rm box}/(5\vartheta)$ is smaller). Since
$W_{(\alpha+e,\beta)}/W_{(\alpha,\beta)}=\delta_yh_r^{-1}$, resp.\
$\delta_y\vartheta^{-1}$, each $y$-shift produces a factor not exceeding $c_{\rm box}\delta_y^{-1}$.
A $v$-shift $\gamma_v$ produces the factor
$w_\beta/w_{\beta+\gamma_v}=\rad^{-|\gamma_v|}\frac{(\beta+\gamma_v)!}{\beta!}\big(\frac{|\beta|+1}{|\beta+\gamma_v|+1}\big)^s\le\rad^{-2}(|\beta|+3)^2$.
Then we apply $\sup_{\hat Q\times B(v_0,1)}|U|\le
C\sum_{|\gamma_y|,|\gamma_v|\le2}\|\partial^{\gamma_y}\partial^{\gamma_v}_vU\|_{L^2(4\hat Q\times
B(v_0,2))}$, that is $H^2_yH^2_v\subset L^\infty$ in $\R^3\times\R^3$. The Jacobian $|Q|=|\tilde
Q|/64$ is absorbed by the normalization $|\tilde Q|^{-1}$. The weights are constant on $Q$ and the
boxes do not depend on $t$, so $C$ is absolute.
\end{proof}

\begin{lem}[Radial gain for plain $D$-factors]\label{lem:RBgain}
Let $X$ be a vector field. Then:
\begin{enumerate}[label=(\roman*)]
\item
\begin{equation}\label{eq:RBgain-a}
a(w)X\cdot X\le2a(w)X_\perp\cdot X_\perp+2|X_r|^2|w|^{\gamma+2}|y'|^2/|y|^2,
\end{equation}
where $X_r=\hat y\cdot X$;
\item
\begin{equation}\label{eq:RBgain-b}
|b(w)\cdot X|\le|b(w)||X_\perp|+|X_r|\big(|b(w)||y'|+2t|w|^{\gamma+2}\big)/|y| .
\end{equation}

\end{enumerate}
Consequently, in~(T2) and in the direct form of $\bb$:
\begin{itemize}
\item the radial part of a plain factor $X=Du=\partial_vu-t\nabla_yu$ contributes
$t\Theta_r\cdot\jap
y^{-1}\le C$ against the weights;
\item the transverse part contributes $t\Theta_\perp$, which is absorbed by
$\int|w|^{\gamma+2}\rho_0(v-w)\dd w\lesssim\jap v^{\gamma+2}$.
\end{itemize}
For $|y|\le\Theta_\perp$, i.e.\ when $\Theta_r\simeq\Theta_\perp$, only the second mechanism is
needed.
\end{lem}

\begin{proof}
Since $a(w)y=a(w)y'$, we have $a(w)\hat y=a(w)y'/|y|$, and $a(w)y'\cdot y'\le|w|^{\gamma+2}|y'|^2$.
For the $\bb$-kernel, it follows from Lemma~\ref{lem:GMwt}(ii) that
$b(w)\cdot y=b(w)\cdot y'+2t|w|^{\gamma+2}$. The factor $|y'|^2$ is controlled by $F_0(y')$ in the
Fubini step of~(T2).
\end{proof}

\begin{lem}[Neighbor sums and the supremum]\label{lem:BKnb}
\leavevmode
\begin{enumerate}[label=(\roman*)]
\item For every $Q$,
\begin{equation}
\frac1{|Q|}\sum_{\mathbf k}\iint_{\tilde Q\times B_2}(W^Q_{\mathbf k})^2|u^Q_{\mathbf k}|^2\le
C^{N_y}\sum_{Q'\sim Q}e_{Q'}\le C^{N_y}N_0\,\mathbb E,
\end{equation}
and the same holds for the dissipation.
\item Suppose that, for every $Q$,
\begin{equation}
e_Q(t)+c\int_0^t\mathfrak d_Q\le e_Q(0)+\int_0^t\big(C\eps\,\mathbb E+\eta\,\textstyle\sum_{Q'\sim
Q}\mathfrak d_{Q'}+F\big).
\end{equation}
If $\eta N_0\le c/4$, then
\begin{equation}
\mathbb E(t)+\tfrac c2\mathcal D(t)\le2\sup_Qe_Q(0)+2\int_0^t(C\eps\,\mathbb E+F),\qquad
\mathbb E(t)\le\tfrac43\Big(\sup_Qe_Q(0)+\int_0^t(C\eps\,\mathbb E+F)\Big).
\end{equation}
Here $F$ does not depend on $Q$, and $\sup_Q\big(e_Q(t)+c\int_0^t\mathfrak d_Q\big)$ is assumed
finite for every $t\le T$. For $\eta=0$ this assumption can be dropped.
\end{enumerate}
\end{lem}

\begin{proof}
(i) We cover $\tilde Q$ by the boxes $2Q'$ of its neighbors $Q'\sim Q$. Then we use $|Q'|\simeq|Q|$
and $h_r(y_{Q'})\simeq h_r(y_Q)$, which yield a factor $C^{N_y}$ for the weights (for $v_0(Q')\ne
v_0(Q)$ the ratio $h_r(y_{Q'},v_0(Q'))/h_r(y_Q,v_0(Q))$ can exceed $2$, but it is $\le3$, and also
$\vartheta(v_0(Q'))\le3\vartheta(v_0(Q))$, $\jap{v_0(Q')}\le3\jap{v_0(Q)}$), and
Lemma~\ref{lem:BKframe}(ii).

(ii) We have
\begin{equation}
\int_0^t\sum_{Q'\sim Q}\mathfrak d_{Q'}\le N_0\,\mathcal D(t).
\end{equation}
We take the supremum over $Q$ of the left side, and use $\max(a,b)\ge\frac12(a+b)$. For the second
bound we put $M(t):=\sup_Q\big(e_Q(t)+c\int_0^t\mathfrak d_Q\big)$. Then $c\,\mathcal D(t)\le M(t)$,
so the right side of the hypothesis is bounded by $\sup_Qe_Q(0)+\frac{\eta
N_0}{c}M(t)+\int_0^t(C\eps\,\mathbb E+F)$, with $\frac{\eta N_0}c\le\frac14$. Taking the supremum
over $Q$, we get $\frac34M(t)\le\sup_Qe_Q(0)+\int_0^t(C\eps\,\mathbb E+F)$, and $\mathbb E(t)\le
M(t)$.
\end{proof}

\begin{lem}[Top-order identity on boxes]\label{lem:BKtop}
Let $\gamma\in[-3,-2)$ (at $\gamma=-2$ the same holds with $C\eps\ell_Q$ in place of $C\eps$, see
Subsection~\ref{ssec:endpoint}, item~(2)). For every $Q$ and $\mathbf k$,
\begin{equation}
\iint\big[D\cdot(ADu)+(2A\zeta-\bb)\cdot Du+S_0u\big]u\,\varphi_Q^2\le-\frac12\iint ADu\cdot
Du\,\varphi_Q^2+C\eps\iint_{\tilde Q\times B_2}u^2 .
\end{equation}
\end{lem}

\begin{proof}
This is Lemma~\ref{lem:ULloc} with the cutoff $\varphi_Q$. The cutoff term $4A\,D\varphi_Q\cdot
D\varphi_Q$ is $\le C\eps$ by Lemma~\ref{lem:RBloc}. On an outer box its radial part uses $A\hat
e_Q\cdot\hat e_Q\le2A\hat y\cdot\hat y+2|\hat y-\hat e_Q|^2|A|\le C\eps\jap y^{-2}$, by
Lemma~\ref{lem:BKcov}(iv) and the null structure. On an inner box all widths are $c_{\rm
box}/(5\vartheta)$, so every direction contributes as the transverse one there. Since
the weights are constant on each box, they produce no derivative terms.
\end{proof}

\begin{lem}[Slabs]\label{lem:slab}
Let $S$ be a box of the shape of $Q$ or of $\tilde Q$ (in $y$), let $\mathcal P(S)$ be the set of
the boxes of the family of Lemma~\ref{lem:BKcov} attached to some $v_0''\in\mathbb Z^3$ that meet
$S$, and let $P''\in\mathcal P(S)$. ($S$ need not lie in $\tilde P''$: when $Q$ is wider than $P''$
in some direction, or tilted with respect to it, no single $\tilde P''$ contains $S$.) The boxes of
$\mathcal P(S)$ have comparable widths and volumes, and frames within
$O((\vartheta''\jap{y''})^{-1})$, since $\operatorname{diam}S\le1$. Thus $L$ below may be computed
with any of them. Let $\sigma$ be the angle between $\hat e_Q$ and $\hat e_{P''}$, with
$\sigma:=\pi/2$ if $P''$ is an inner box.
\begin{enumerate}[label=(\roman*)]
\item For every $u$,
\begin{equation}
\int_S|u|^2\le C\,L\,\frac{|S|}{|P''|}\max_{P\in\mathcal P(S)}\int_{\tilde P}|u|^2,\qquad
L:=\max\Big(1,\frac{h_r}{h_r''}|\cos\sigma|+\frac{h_r}{\vartheta''}|\sin\sigma|\Big)\max\Big(1,\frac{\vartheta}{\vartheta''}\Big)^2,
\end{equation}
where $h_r,\vartheta$ are evaluated at $Q$. Here $L$ is the ratio between the width of $P''$ and
that of $S$, measured along the directions of $S$.
\item If $u$ has two more $y$-derivatives in the hierarchy, then the factor $L$ can be dropped:
\begin{equation}
\int_S|u|^2\le C|S|\max_{P\in\mathcal P(S)}|P|^{-1}\sum_{|\gamma|\le2}\ell_P^{2\gamma}\int_{\tilde
P}|\partial_{P}^\gamma u|^2 .
\end{equation}
\end{enumerate}
\end{lem}

\begin{proof}
(i) The boxes of the family cover $\R^3$, so we have
$\int_S|u|^2\le\sum_{P\in\mathcal P(S)}\int_{\tilde
P}|u|^2\le\#\mathcal P(S)\max_{P\in\mathcal P(S)}\int_{\tilde P}|u|^2$. Every $P\in\mathcal P(S)$
lies in $S+C(P''-P'')$, and the boxes have bounded overlap, so $\#\mathcal P(S)\,|P''|\le
C|S+C(P''-P'')|$. The Minkowski sum lies in the product of the sums of the projected widths along
the frame of $S$. The width of $P''$ along $\hat e_Q$ does not exceed
$|\cos\sigma|\ell_r''+\sqrt2|\sin\sigma|\ell_\perp''$, and along $e^{1,2}_Q$ it does not exceed
$\sqrt3\ell_\perp''$. Thus we have $|S+C(P''-P'')|\le CL|S|$. Indeed, in a direction where $S$ is
wider than $P''$ nothing is lost, and in a direction where $S$ is thinner the loss is the ratio of
widths.

(ii) At fixed $v$, we apply $\int_S|u|^2\le\sum_{P\in\mathcal P(S)}|S\cap P|\sup_P|u|^2$ together
with
Sobolev's inequality in $y$ at the scale of each $P$ ($H^2(\R^3)\subset L^\infty$ on the rescaled
unit box, as in Lemma~\ref{lem:RBsob}). The factors $|S\cap P|$ do not depend on $v$ (the slabs are
products of a $y$-set and a $v$-ball), and $\sum_P|S\cap P|\le C|S|$ by bounded overlap. Therefore,
after the $v$-integration the maximum over $P$ can be taken outside the integral, and no
$v$-derivative enters.

\end{proof}

\begin{lem}[Compensation]\label{lem:BKcomp}
\leavevmode
\begin{enumerate}[label=(\roman*)]
\item \emph{Absolute estimates} ((T2), (T5), and $\bb$ in direct form). Let
$\Lambda^2:=L\big(h_r''\vartheta''^2/(h_r\vartheta^2)\big)^2$, with $L$ the slab factor of
Lemma~\ref{lem:slab}(i), $h_r,\vartheta$ evaluated at $Q$ and $h_r'',\vartheta''$ at $P''$. Then,
for every $y''$ in the slab, $\Lambda^2\le C(\jap{y''}+\vartheta'')^2\vartheta''^4$, and the
unpaired factor of a lift satisfies $h_r''\vartheta''^2\le C(\jap{y''}+\vartheta'')\vartheta''^2$: a
polynomial in $\jap{y''}$, times $\vartheta''^{C}$. These are absorbed by $F_0(y'')\rho_0(v'')$ in
the Fubini step.
\item \emph{Relative estimates} ((T3) with top-level $h$). We have
\begin{equation}
\ab[f_0\Lambda^2]\le C\,\Lambda_*^2\,A,\qquad\Lambda_*:=C\big(1+\log\jap v+\log\jap
y+\log(2+t^{-1})\big)^{3/2}.
\end{equation}
On the set where $\vartheta''\le\vartheta$ the contribution is $O(A)$, and the polylogarithmic
factor comes only from integration velocities $v''=v-w$ with $\vartheta''>\vartheta$. The proof also
yields the angle inequality~\eqref{eq:BKcomp-angle}. Lemma~\ref{lem:FPslab} uses only this
inequality. The bound with $\Lambda_*$ is not used elsewhere, since that lemma replaces it by a
bound without $\Lambda_*$ that uses the velocity weight.
\end{enumerate}
\end{lem}

\begin{proof}
(i) Since $h_r,h_r'',\vartheta''\ge1$, we have
$L\le\big(1+h_r/h_r''+h_r/\vartheta''\big)\max(1,\vartheta/\vartheta'')^2\le3h_r\vartheta^2$, so
$\Lambda^2\le3h_r''^2\vartheta''^4/(h_r\vartheta^2)\le3h_r''^2\vartheta''^4$. Since
$\operatorname{diam}S\le\operatorname{diam}\tilde Q\le4\sqrt3c_{\rm box}$,
$\operatorname{diam}P''\le\sqrt3c_{\rm box}$ and $P''$ meets $S$, we deduce that
$|y''-y_{P''}|\le\frac92\sqrt3c_{\rm box}\le1$ for $y''\in S$ (as $c_{\rm box}\le\frac18$), which
implies $h_r''\le\jap{y''}+1+\sqrt{\log\jap{v_0''}}\le2\jap{y''}+\vartheta''$. In the Fubini step
of~(T2) in Remark~\ref{rem:ULbox} and of~(T5), $F_0(y'')$ absorbs every power of $\jap{y''}$. Since
$|v_0''-(v_0-w)|\le4$, we have $\vartheta''^C\le C(1+\log\jap{v_0-w})^{C}$, which is controlled by
$\int|w|^{\gamma+2}\rho_0(v_0-w)(1+\log\jap{v_0-w})^C\dd w\le C$ (by $\rho_0^{1/8}$ in~(T5)), as
$\rho_0\le A_F\jap\cdot^{-4}$.

(ii) Throughout, $A=\ab[g]$ with $g=f_0G$, $\frac12\le G\le2$, so that $A\simeq\ab[f_0]$. Moreover,
$h_r,\vartheta$ refer to the box $Q\ni(y,v)$, and $h_r'',\vartheta''$ to a box $P''\ni(y'',v'')$,
$z=(y'',v'')=(y+tw,v-w)$, of the family of Lemma~\ref{lem:BKcov}. Let $y''_\perp$ be the transverse
part of $y''$ in the frame of $Q$ (for an inner $Q$, $\hat e_Q=e_3$ by~(ii) of that lemma,
and $y''_\perp$ is the part orthogonal to $e_3$). We put
\begin{equation}
\mathsf L:=1+\log\jap v+\log\jap y+\log(2+t^{-1}).
\end{equation}

\emph{Step 1: the angle inequality.} We show that for every $y''$ in a slab $S$ that meets $P''$,
\begin{equation}\label{eq:BKcomp-angle}
|\sin\sigma|\,h_r''\le C\big(|y''_\perp|+\vartheta''\big).
\end{equation}
If $P''$ is inner, $\sigma=\pi/2$ and $h_r''\le5\vartheta''$. Let $P''$ be outer. Then
$|y_{P''}|\ge2\vartheta''\ge2$ and, as in~(i), $|y''-y_{P''}|\le1$. Thus $|y''|\ge1$, and the sine
of
the angle between $\hat e_{P''}=\hat y_{P''}$ and $\hat y''$ is $|\hat
y''\times(y_{P''}-y'')|/|y_{P''}|\le1/|y_{P''}|\le2/\jap{y''}$, since
$\jap{y''}\le|y''|+1\le|y_{P''}|+2\le2|y_{P''}|$. For unit vectors, $|a\times c|\le|a\times
b|+|b\times c|$ (write $c=(c\cdot b)b+c_{\perp}$). With $b=\hat y''$ and $|\hat e_Q\times\hat
y''|=|y''_\perp|/|y''|$, this implies $|\sin\sigma|\le|y''_\perp|/|y''|+2/\jap{y''}$. Since
$h_r''\le2\jap{y''}+\vartheta''$ (proof of~(i)), $\jap{y''}\le\sqrt2|y''|$, $|y''_\perp|\le|y''|$
and $\vartheta''\ge1$, we get
$|\sin\sigma|\,h_r''\le2\sqrt2|y''_\perp|+|y''_\perp|\vartheta''/|y''|+4+2\vartheta''\le
C(|y''_\perp|+\vartheta'')$. For an inner $Q$ the same computation applies.

\emph{Step 2: reduction.} We put $X_1:=h_r''/h_r$, $X_2:=|y''_\perp|/\vartheta$ and
$X_3:=\vartheta''/\vartheta$. As in the proof of Lemma~\ref{lem:FPslab}, which rests
on~\eqref{eq:BKcomp-angle} and the definition of $L$ alone, we have
\begin{equation}
\Lambda^2\le C\big(X_1^2+X_1+X_1X_2\big)\max(1,X_3)^4 .
\end{equation}
The form $\xi\mapsto\ab[f_0\varphi]\xi\cdot\xi$ is nonnegative for $\varphi\ge0$, so by the
Cauchy--Schwarz inequality for the measure $a(w)\xi\cdot\xi\,f_0(z)\dd w$, we get
\begin{equation}
\ab[f_0\Lambda^2]\xi\cdot\xi\le C\Big(\ab\big[f_0(1+X_1^4+X_2^4)\big]\xi\cdot\xi\Big)^{1/2}\Big(\ab\big[f_0\max(1,X_3)^8\big]\xi\cdot\xi\Big)^{1/2}.
\end{equation}

\emph{Step 3: the first factor is $\le CA\xi\cdot\xi$.} By Lemma~\ref{lem:R3comp}, we have
$\Theta\simeq h_r$, and $h_r''\simeq\Theta(z)$, $h_r\simeq\Theta(y,v)\ge1$. Applying
Lemma~\ref{lem:ULmom} with $m=4$, we obtain $\ab[f_0X_1^4]\le C(1+h_r^{-4})A\le CA$. On an outer
$Q$, the two transverse unit vectors $e$ of the frame of $Q$ satisfy $|e\cdot y|\le2$ by~(iv) of
the covering lemma, and $|y''_\perp|^4\le4\sum_e|e\cdot y''|^4$. Then Lemma~\ref{lem:RBdir} with
$m=4$ implies $\ab[f_0X_2^4]\le C(\Theta_\perp(v)^4+1)\vartheta^{-4}A\le CA$, since
$\Theta_\perp(v)\le C\vartheta$ for $v\in B(v_0,2)$. On an inner $Q$, we have
$|y''_\perp|\le|y''|\le Ch_r''$ and $h_r\le5\vartheta$, so $X_2\le CX_1$.

\emph{Step 4: a polynomial lower bound for $A$.} We claim
\begin{equation}\label{eq:BKcomp-low}
A\xi\cdot\xi\ge c\,\eps\,e^{-C'\mathsf L}|\xi|^2,
\end{equation}
with $c,C'$ depending only on $\gamma$, $t_1$ and the data class. Let $|\xi|=1$. On the ball
$B:=B(-y/t,1/t)$ of $w$ we have $|y+tw|\le1$. There $g(z)\ge\frac12\eps\rho_0(v-w)F_0(y'',v'')$,
with $\rho_0(v-w)\ge A_F^{-1}\jap{v-w}^{-p}\ge c\jap v^{-p}\jap w^{-p}$ by (D3). By (D$_{\rm iso}$),
$F_0(x,v'')/F_0(0,v'')\in[e^{-\beta_F|x|^2-C_\Phi|x|-2C_{\rm pv}},e^{C_\Phi|x|+2C_{\rm pv}}]$.
Together with $\int F_0\dd x=1$ and $F_0\le C_Fe^{-\beta_F|x|^2}$, this implies $F_0(0,v'')\ge c$,
and hence $F_0(y'',v'')\ge c$ for $|y''|\le1$. Next, $a(w)\xi\cdot\xi=|w|^{\gamma+2}|\hat
w\times\xi|^2$. The
double cone $\{|\hat w\times\xi|<\delta\}$, $\delta\le1$, meets the sphere $\{|w|=s\}$ in two caps
of angular radius $\arcsin\delta\le\frac\pi2\delta$, of total area $\le\frac{\pi^3}2s^2\delta^2$.
Since $B$ lies in the shell $\max\{0,(|y|-1)/t\}\le|w|\le(|y|+1)/t$, of width $\le2/t$, the cone
meets $B$ in a set of volume $\le\pi^3\delta^2(|y|+1)^2t^{-3}$, which is $\le\frac3{10}|B|$ for
$\delta:=\frac1{5(|y|+1)}$. On the rest of $B$, we have $|w|\le(|y|+1)/t$, so
$|w|^{\gamma+2}\ge((|y|+1)/t)^{\gamma+2}$ (as $\gamma+2\le0$) and $\jap
w^{-p}\ge\jap{(|y|+1)/t}^{-p}$. Hence we obtain
\begin{equation}
A\xi\cdot\xi\ge c\,\eps\,\jap v^{-p}\,t^{-3}\,(|y|+1)^{-2}\Big(\frac{|y|+1}t\Big)^{\gamma+2}\Big\langle\frac{|y|+1}t\Big\rangle^{-p}\ge
c\,\eps\,\jap v^{-p}\,\jap y^{-C}(2+t^{-1})^{-C},
\end{equation}
where we use $t\le t_1$. This proves~\eqref{eq:BKcomp-low}.

\emph{Step 5: the second factor is $\le C\mathsf L^4A\xi\cdot\xi$.} Since $\vartheta\ge1$, we have
$\max(1,X_3)^8\le1+\vartheta''^8$, and
$\vartheta''^2=\big(1+\sqrt{\log\jap{v_0''}}\big)^2\le2(1+\log\jap{v_0''})\le C(1+\log\jap{v''})$,
as $|v_0''-v''|\le2$. Let $K\ge2$ and $M:=e^{K\mathsf L}$, so that $M\ge2\jap v$. On $\{\jap{v''}\le
M\}$, we have $\vartheta''^8\le C(1+K\mathsf L)^4\le CK^4\mathsf L^4$, which contributes $\le CK^4
\mathsf L^4A\xi\cdot\xi$. On $\{\jap{v''}>M\}$, $|w|\ge\jap{v''}-\jap v\ge\frac12\jap{v''}\ge1$, so
$a(w)\xi\cdot\xi\le|w|^{\gamma+2}\le1$, and $f_0(z)\le C\eps\jap{v''}^{-p}$ by (D1), (D3). Thus this
part contributes no more than
\begin{equation}
C\eps\int_{\jap{v''}>M}\jap{v''}^{-p}(\log\jap{v''})^4\dd v''\le C\eps M^{-(p-3)}(\log M)^4\le
C\eps\,e^{-\frac12(p-3)K\mathsf L},
\end{equation}
since $p\ge4$. We choose $K\ge2$ with $\frac12(p-3)K\ge C'$. Then, by Step~4, this is $\le
CA\xi\cdot\xi$.

\emph{Conclusion.} By Steps 2, 3 and 5, we have $\ab[f_0\Lambda^2]\le C\mathsf L^2A\le\Lambda_*^2A$,
since
$\mathsf L\ge1$ (the exponent $\frac32$ in $\Lambda_*$ is not sharp).
\end{proof}

\begin{lem}[Hypodissipativity in static Gaussian weights]\label{lem:GW}
There is $q_{\rm GW}=q_{\rm GW}(\gamma,L,\theta_\pm)\in(0,1/(8\theta_+)]$, and, for $q\in(0,q_{\rm
GW}]$, there are $T_1,M,R,c_W,C_W$, uniform for $q$ in compact subsets of $(0,q_{\rm GW}]$ (as in
Lemma~\ref{lem:RFtools}(vi)), such that the following holds for $t\ge T_1$, all $x$, and $m=\hat
m_x$ the leaf weight of Section~\ref{sec:weights} (Lemma~\ref{lem:Gleafweight} and
Lemma~\ref{lem:TG0}(i)):
\begin{enumerate}[label=(\roman*)]
\item $\int(\mathcal B\varphi)\varphi\,m^2\le-c_W\|\varphi\|^2_{Y_2(m)}$, where $
\mathcal B=Q(\mu,\cdot)-M\chi_R$;
\item the estimates (W1$_G$) and (W2$_G$) of Corollary~\ref{cor:W1G} hold for $m$, which lies in
the class $\mathcal W$ of that corollary by Lemma~\ref{lem:TG0}(i). They contain no term
$\eps^{|\gamma|}\|\varphi\|^2_{X(m)}$.
\end{enumerate}
\end{lem}

\begin{proof}
(i) By the identity in the proof of \cite[Lemma~2.3]{CM}, we have
\begin{equation}
\int(\mathcal B\varphi)\varphi m^2=-\int\ab\nabla(m\varphi)\cdot\nabla(m\varphi)+\int(\tilde\zeta_m-M\chi_R)\varphi^2m^2
.
\end{equation}
Here $\tilde\zeta_m=\ab\nabla\Phi\cdot\nabla\Phi+\bb\cdot\nabla\Phi-\frac12\cb$, with $\Phi=\log
m=\Phi_0+\omega_x$ and $\Phi_0=q_\xi|\hat z|^2$.

For $\Phi_0$, it follows from \cite[Lemma~2.2(3)]{CM} that
$\tilde\zeta_{m_0}\le4q_\xi(2q_\xi-1)\jap{\hat
z}^{\gamma+2}(1+o(1))$. Since $q_\xi\le1/8$, we have $4q_\xi(2q_\xi-1)\le-3q_\xi$, so this is
$\le-\frac52q_\xi\jap{\hat z}^{\gamma+2}$ for $|\hat z|\ge R_0$.

The perturbation is given by
\begin{equation}
\tilde\zeta_m-\tilde\zeta_{m_0}=2\ab\nabla\Phi_0\cdot\nabla\omega+\ab\nabla\omega\cdot\nabla\omega+\bb\cdot\nabla\omega
.
\end{equation}
Since $\nabla\Phi_0=2q_\xi\hat z$ is an eigenvector of $\ab$ with eigenvalue $\ell_1\le C\jap{\hat
z}^\gamma$, and $\ell_2\le C\jap{\hat z}^{\gamma+2}$, we get
\begin{equation}
|\tilde\zeta_m-\tilde\zeta_{m_0}|\le C\jap{\hat z}^{\gamma+1}(1+q)\eps_1+C\jap{\hat
z}^{\gamma+2}\eps_1^2,\qquad \eps_1=\eps_1(t)=2\theta_+^{1/2}L(q+t^{-1}),
\end{equation}
We choose $q_{\rm GW}$ so small that $C\eps_1^2\le q_\xi$ for $q\le q_{\rm GW}$ and $t\ge T_1(q)$.
This needs $q\,\theta_+L^2\lesssim\theta_-$ and $t^2\gtrsim L^2/q$. Then we choose $R_{\rm GW}$ so
that
$C(1+q)\eps_1\jap{R_{\rm GW}}^{-1}\le q_\xi$. Then $\tilde\zeta_m\le-\frac12q_\xi\jap{\hat
z}^{\gamma+2}$ for $|\hat z|\ge\max\{R_0,R_{\rm GW}\}$, and $M\chi_R$ with $R=\max\{R_0,R_{\rm
GW}\}$ and $M$ large handles $|\hat z|\le2R$. The gradient part of $\|\cdot\|_{Y_2(m)}$ follows from
the coercivity of $\ab$ \cite[Lemma~2.1(b)]{CM}.

This argument involves $|\nabla\omega|\le\eps_1$, $q_\xi\le\frac18$ and $C\eps_1^2\le q_\xi$, but
no bound on $\nabla^2\omega$. Therefore~(i) holds, with constants depending only on
$(\gamma,q_\pm,\eps_1)$, for every $m\in\mathcal W(q_\pm,\eps_1)$ with $q_+\le\frac18$ and
$C\eps_1^2\le q_-$. For the leaf weights this condition is the choice of $q_{\rm GW}$ and $T_1$
above.

(ii) This is Corollary~\ref{cor:W1G}, whose proof rests on~(i) alone. We also give an alternative
argument, which uses $|\nabla^2\log m|\le C$ as well. The parts $Q(\cdot,\mu)+M\chi_R$ produce
Gaussian-localized functions, since $mQ(\varphi,\mu)\in L^2$ for $q_\xi<1/4$. Splitting their
pairing with $m^2\varphi$ at $|\hat z|=R'$, we have
\begin{equation}
|\langle Q(\varphi,\mu)+M\chi_R\varphi,\varphi\rangle_{L^2(m)}|\le C\|\varphi\|_{L^2(\jap{\hat
z}^\theta)}\|\varphi\|_{Y_2(m)}\le\delta\|\varphi\|^2_{Y_2(m)}+C_\delta\|\varphi\|^2_{L^2},
\end{equation}
by interpolation between $L^2$ and $Y_2(m)$. We estimate the perturbative terms
$Q(f,\varphi)+Q(\varphi,f)$ as in the proof of \cite[Lemma~4.3]{CM}. That proof requires
$|\nabla\log m|\le
C\jap{\hat z}$ and $|\nabla^2\log m|\le C$, which hold uniformly by Lemma~\ref{lem:Gleafweight} and
Lemma~\ref{lem:TG0}(i). No $\eps^{|\gamma|}$-defect appears, because the zero-order part $\jap{\hat
z}^{(\gamma+2)/2}m$ of $Y_2(m)$ controls every term generated by $\nabla\log m$, as in
Lemma~\ref{lem:QmG}.
\end{proof}

\begin{lem}[Gaussian hypodissipativity without closeness to equilibrium]\label{lem:Whypo}
Let $G\ge0$ satisfy
\begin{itemize}
\item[(G1)] $\int G=M\in[M_-,M_+]$, $\int zG=My_0$ with $|y_0|\le L$, and $\int|z-y_0|^2G=3M\theta$
with $\theta\le\theta_+$;
\item[(G2)]
\begin{equation}\label{eq:G2}
\ab[G]\zeta\cdot\zeta\ge c_G\big(\jap z^\gamma|P_z\zeta|^2+\jap z^{\gamma+2}|(I-P_z)\zeta|^2\big);
\end{equation}

\item[(G3)] $\|\jap z^8G\|_{L^\infty}\le K$.
\end{itemize}
Let $m=e^{q|z|^2}$ with $0<q\le\big(8(L^2+3\theta_+)\big)^{-1}$. Then, for complex $\varphi$,
\begin{equation}
\mathrm{Re}\!\int Q(G,\varphi)\,\bar\varphi\,m^2\le-c\,\|\varphi\|^2_{Y_2(m)}+C\,\|\varphi\|^2_{L^2(m)},
\end{equation}
with $c,C$ depending only on $(\gamma,M_\pm,L,\theta_+,K,c_G,q)$.
\end{lem}

\begin{proof}
The coefficients $\ab[G],\bb[G],\cb[G]$ are real. Hence the identity in the proof of
\cite[Lemma~2.3]{CM} holds for complex $\varphi$ after taking real parts:
\begin{equation}
\mathrm{Re}\!\int Q(G,\varphi)\bar\varphi m^2=-\int\ab[G]\nabla(m\varphi)\cdot\overline{\nabla(m\varphi)}+\int\tilde\zeta\,|\varphi|^2m^2,\qquad\tilde\zeta=\ab[G]\nabla\Phi\cdot\nabla\Phi+\bb[G]\cdot\nabla\Phi-\tfrac12\cb[G],
\end{equation}
with $\Phi=q|z|^2$ and $\nabla\Phi=2qz$.

\emph{Two exact identities.} With $u=z-w$ we have $|u|^2-(u\cdot\hat z)^2=|w_\perp|^2$. Hence
\begin{equation}
\ab[G](z)\hat z\cdot\hat z=\int|z-w|^\gamma|w_\perp|^2G(w)\dd w,\qquad \bb[G](z)\cdot\hat
z=-2\int|z-w|^\gamma\big(|z|-w\cdot\hat z\big)G(w)\dd w .
\end{equation}
Therefore
\begin{equation}
\tilde\zeta=4q|z|\int|z-w|^\gamma\Big[q|z|\,|w_\perp|^2-|z|+w\cdot\hat z\Big]G\dd w-\tfrac12\cb[G].
\end{equation}
\emph{The region $|w|\le|z|/2$.} There $|z-w|^\gamma=|z|^\gamma(1+O(|w|/|z|))$. By (G1), we have
$q\int|w|^2G=qM(|y_0|^2+3\theta)\le M/8$, so this part is bounded by
$-\frac72qM|z|^{\gamma+2}+CK|z|^{\gamma+1}$.

\emph{The region $|w|>|z|/2$.} By (G3) its contribution is $O(K|z|^{-3})$, since $q|z|\cdot
q|z|\cdot|z|^{-5}$. Near $w=z$ the integrands are $O(|z-w|^{\gamma+2})$ and $O(|z-w|^{\gamma+1})$,
because $|w_\perp|\le|z-w|$ and $|z|-w\cdot\hat z=(z-w)\cdot\hat z$. Consequently they are locally
integrable, also for $\gamma=-3$.

\emph{The term $\cb$.} For $\gamma>-3$, we have $|\cb[G]|\le CK\jap z^\gamma$. For $\gamma=-3$,
$-\frac12\cb[G]=4\pi G\le4\pi K\jap z^{-8}$.

Hence we have $\tilde\zeta\le-3qM\jap z^{\gamma+2}$ for $|z|\ge R_0$, and $\tilde\zeta\le C$ on
$\{|z|\le R_0\}$. Together with (G2), this proves the claim.

\end{proof}

\subsection{The initial layer: energies on adapted boxes}\label{ssec:layer}

In this subsection we prove Theorem~\ref{thm:layer}, except for the $\bb$-part of the drift at the
top $y$-levels, which we treat in Subsection~\ref{ssec:drift}.

\subsubsection{Setting}\label{sssec:FPset}

Fix a box $Q$ of Lemma~\ref{lem:BKcov}, attached to the velocity ball $B(v_0,1)$. Its frame $(\hat
e_Q,e^1_Q,e^2_Q)$ is constant, so the frame derivatives $\partial_Q^\alpha$ commute with
$D=\partial_v-t\partial_y$. For $\mathbf k=(\alpha,\beta)$ with $|\alpha|\le N_y$, let
$u_{\mathbf k}:=\partial_Q^\alpha\partial_v^\beta G$ for $\mathbf k\ne0$, and $u_0:=G-1$. Let
$a_r,a_\perp$ be the radial and transverse counts of $\alpha$. The weights are constants on $Q$:
\begin{equation}
W_{\mathbf k}:=\delta_y^{|\alpha|}\,h_r(Q)^{-a_r}\,\vartheta(Q)^{-a_\perp}\,\tau_\alpha\,\jap{v_0}^{\sigma}\,w_\beta(\rad(t)),\qquad\tau_\alpha:=\big(h_r(Q)\vartheta(Q)^2\big)^{-
\mathbf 1\{|\alpha|\ge N_y-1\}},
\end{equation}
where $h_r(Q)=h_r(y_Q,v_0)$, $\vartheta(Q)=\vartheta(v_0)$,
$w_\beta=\rad^{|\beta|}(|\beta|+1)^s/\beta!$, $s\ge6$, and $\rad(t)=\rad_0-C_\rad\eps
t$. From Proposition~\ref{prop:FPclose} on, we take all energies with the weakened weight
$W^\flat_{\mathbf k}$ of Definition~\ref{def:FPtop} in place of $W_{\mathbf k}$, since the corner
cases of Subsection~\ref{ssec:drift} need it and
that proposition is proved only for $W^\flat$. Put
\begin{equation}
a_{\mathbf k}:=W_{\mathbf k}\|u_{\mathbf k}\varphi_Q\|_{L^2},\qquad d_{\mathbf k}:=W_{
\mathbf k}\Big(\iint A\,Du_{\mathbf k}\cdot Du_{\mathbf k}\,\varphi_Q^2\Big)^{1/2},
\end{equation}
and $e_Q:=|Q|^{-1}\sum_{\mathbf k} a_{\mathbf k}^2$, $\mathfrak d_Q:=|Q|^{-1}\sum_{\mathbf k}
d_{\mathbf k}^2$, $\mathbb E=\sup_Qe_Q$. All quantities attached to another box $P$ are taken with
the weights of $P$.

\emph{The equation for $u_{\mathbf k}$.} With $V:=2A\zeta-\bb$ and
$S_0:=A:(D\zeta+\zeta\otimes\zeta)-\cb$, we have
\begin{equation}
\partial_tu_{\mathbf k}=D\cdot(ADu_{\mathbf k})+V\cdot Du_{\mathbf k}+S_0u_{
\mathbf k}+\sum_{0<\mu\le\mathbf k}\binom{\mathbf k}\mu\Big[D\cdot(\partial^\mu
A\,Du_{\mathbf k-\mu})+\partial^\mu V\cdot Du_{\mathbf k-\mu}+\partial^\mu
S_0\,u_{\mathbf k-\mu}\Big]+\partial^{\mathbf k} S_0 .
\end{equation}
For $\mathbf k=0$ the last term is $S_0$. For $\mathbf k\ne0$ it is the part of $\partial^{
\mathbf k}(S_0G)$ with all derivatives on $S_0$ and the factor $1$ of $G=1+u_0$, while the factor
$u_0$ is included in the sum, as the term $\mu=\mathbf k$. By the Leibniz rule, we have
\begin{equation}
\partial^\mu A=\sum_{\nu+\mathbf k'=\mu}\binom\mu\nu\,\ab[\partial^\nu f_0\,\partial^{\mathbf k'}G],
\end{equation}
with $\partial^0G=G$, and similarly for $\bb$ and $\cb$. We call the terms with $\mathbf k'=0$
\emph{explicit} and the terms with $\mathbf k'\ne0$ \emph{$G$-terms}.

\emph{The energy identity.} We multiply the equation for $u_{\mathbf k}$ above by
$W_{\mathbf k}^2u_{\mathbf k}\varphi_Q^2$, integrate, and sum over $\mathbf k$. Since $\frac
d{dt}W_{\mathbf k}^2=2(\rad'/\rad)|\beta|W_{\mathbf k}^2$, we obtain
\begin{equation}\label{eq:FPid}
\frac12\frac d{dt}\sum_{\mathbf k} a_{\mathbf k}^2=-\sum_{\mathbf k}\theta_{\mathbf k} a_{
\mathbf k}^2+\sum_{\mathbf k} W_{\mathbf k}^2\big(\mathcal T_{\mathbf k}+\mathcal C_{\mathbf k}+
\mathcal F_{\mathbf k}\big),\qquad\theta_{\mathbf k}:=\frac{C_\rad\eps|\beta|}{\rad}.
\end{equation}
Here $\mathcal T_{\mathbf k}$, $\mathcal C_{\mathbf k}$, $\mathcal F_{\mathbf k}$ are the pairings
of the first three terms, of the sum, and of $\partial^{\mathbf k} S_0$ with $u_{
\mathbf k}\varphi_Q^2$.

\subsubsection{Coefficient bounds}

We bound the explicit factors and the coefficients that occur in the commutator terms. Write
$z=(y'',v'')=(y+tw,v-w)$ for the integration point, $y''_\perp$ for the transverse part of $y''$ in
the frame of $Q$, and $\langle\phi\rangle_A:=\ab[f_0\phi]\xi\cdot\xi/A\xi\cdot\xi$.

\begin{lem}[Explicit factors]\label{lem:FPexp}
Let $R=\Lambda/2$ be the radius of Section~\ref{sec:layer}. For $|\nu_y|\le N_y$ (and for every
$\nu_y$, with $C_{N_y}$ replaced by a constant depending also on $|\nu_y|$, and the rate $C$ and this
constant depend only on $\Dc$, not on $t$, $t_1$ or $Q$),
\begin{equation}
|\partial_Q^{\nu_y}\partial_v^{\nu_v}f_0(z)|\le C_{N_y}C^{|\nu_v|}\nu_v!\,R^{-|\nu_v|}\big(1+|\hat
e_Q\cdot y''|\big)^{\nu_r}\big(1+|y''_\perp|\big)^{\nu_\perp}f_0(z),
\end{equation}
and, as quadratic forms, for $m_r,m_\perp\le N_y$ (for larger $m_r,m_\perp$ with a constant
depending on them, and here $C_{N_y}$, unlike the constant of the first bound, depends also on $t_1$,
through Lemmas~\ref{lem:R3comp} and~\ref{lem:ULmom}, as does the constant for larger $m_r,m_\perp$)
and under the bootstrap hypothesis $\frac12\le G\le2$ (which implies $\ab[f_0]\le2A$),
\begin{equation}
\ab\big[(1+|\hat e_Q\cdot y''|)^{2m_r}(1+|y''_\perp|)^{2m_\perp}f_0\big]\le
C_{N_y}\,h_r(Q)^{2m_r}\vartheta(Q)^{2m_\perp}A\qquad\text{at every }(y,v)\in Q.
\end{equation}
\end{lem}

\begin{proof}
We have $f_0=\eps\rho_0(v)\,\mathcal U(y,v)\,F_{\rm G}(y)$ by~\eqref{eq:FGU}, with $F_{\rm G}=c_{\rm
G}e^{\phi}$ and $\phi=-\beta_F|y|^2-\Phi_1$. We expand $\partial_Q^{\nu_y}\partial_v^{\nu_v}$ by
the Leibniz rule over the three factors.
\begin{itemize}
\item \emph{$y$-factor.} By Fa\`a di Bruno, $\partial^{\nu_y}F_{\rm G}/F_{\rm G}$ is a sum of
products of derivatives of $\phi$. Each first derivative in a direction $e$ contributes
$|e\cdot\nabla\phi|\le2\beta_F|e\cdot y''|+C_\Phi$. Higher derivatives of $\phi$ are bounded: by
(D$_{\rm iso}$) and the Cauchy estimates on the strip $\{|\mathrm{Im}\,x|<\Lambda_\Phi\}$,
$\|\nabla^j\Phi_1\|_{L^\infty(\R^3)}\le C(j)C_\Phi\Lambda_\Phi^{1-j}$ for every $j\ge1$. We apply
this with $e\in\{\hat e_Q,e^1_Q,e^2_Q\}$, and the constants depend on $N_y$.
\item \emph{$v$-factor.} For $\rho_0$ we use the Cauchy estimate~\eqref{eq:rho0Cauchy}
on the polydisc $\Delta(v'',\Lambda/2)$, as in Lemma~\ref{lem:dispfree}, and we get
$|\partial^{\nu_v}\rho_0(v'')|\le A_Fe^{L\Lambda}\nu_v!\,R^{-|\nu_v|}\rho_0(v'')$. Here we use the
upper half $|\rho_0(\cdot+i\eta)|\le A_F\rho_0$ of (D3) and $|\nabla\log\rho_0|\le L$. The factor
$A_Fe^{L\Lambda}$ can be put into $C_{N_y}$.
\item \emph{The perturbation.} By Lemma~\ref{lem:Ufac}(iii), $|\partial_y^a\partial_v^b
\mathcal U|\le K_{\mathcal U}e^{C_{\mathcal U}}a!\,b!\,(\Lambda_\Phi/2)^{-|a|}R^{-|b|}\,\mathcal U$. In the
Leibniz sum over the $v$-derivatives,
$\sum_{b\le\nu_v}\binom{\nu_v}b\,b!\,(\nu_v-b)!=\nu_v!\prod_i(\nu_{v,i}+1)\le2^{|\nu_v|}\nu_v!$, so
the $v$-derivatives contribute $C^{|\nu_v|}\nu_v!R^{-|\nu_v|}$ with a larger $C$. The
$y$-derivatives of
$\mathcal U$, of order $\le N_y$, only contribute a constant, which we put into $C_{N_y}$. For
product data this factor is absent.
\item \emph{The moment bound.} By H\"older's inequality for the measure $a(w)\xi\cdot\xi f_0(z)\dd
w$, we split it into the radial moments of Lemma~\ref{lem:ULmom} (using $|\hat e_Q\cdot y''|\le|y''|$
and $\Theta\le Ch_r$) and the transverse moments of Lemma~\ref{lem:RBdir} (using $|e\cdot y|\le2$ on
an outer $Q$). On an inner $Q$ we have $|e\cdot y''|\le|y''|$, and the radial moment bound controls
these moments, since $h_r\le5\vartheta$ on $\tilde Q$. Finally, $h_r$ and $\vartheta$ on $Q$ are
comparable to $h_r(Q)$ and $\vartheta(Q)$.
\end{itemize}
\end{proof}

\begin{lem}[Decay in $v$, and the null-transferred vectors]\label{lem:FPdecay}
At every point, with $C$ independent of $t\le t_1$ and $y$ (in the bound on $S_0$ in~(i), in the
bound on $\zeta^\sharp$ and in~(iv), $C$ depends on $t_1$, through $L_\zeta$ and $2\beta_Ft_1$, while
the bounds on $A$, $\bb$, $\cb$ do not), we have:
\begin{enumerate}[label=(\roman*)]
\item $|A|+\ab[\jap{y''}^{C}f_0]\le C\eps\jap v^{\gamma+2}$, $|\bb|\le C\eps\jap v^{\gamma+1}$,
$|\cb|\le C\eps\jap v^{\gamma}$, and $|S_0|\le C\eps\jap v^{\gamma+2}$;
\item $a(w)\zeta=a(w)\zeta^\sharp$ with $|\zeta^\sharp|\le C\jap{y''}$, and $|\partial^\mu\zeta|\le
C_\zeta K_\zeta^{|\mu|}\mu!R^{-|\mu|}$ for $\mu\ne0$, where the prefactor $C_\zeta$ may depend on
$t_1$ (through the factor $t$ in $\zeta=D\log f_0$, cf.\ $L_\zeta$ in the proof of
Lemma~\ref{lem:GMwt}), but the rate $K_\zeta$ does not;
\item $A\hat y\cdot\hat y\le C\eps\jap y^{-2}$ for $|y|\ge1$, and hence $A\hat e_Q\cdot\hat e_Q\le
C\eps\jap y^{-2}$ on outer boxes;
\item $|A\zeta|+|D\cdot(A\zeta)|+|\bb\cdot\zeta|\le C\eps$.
\end{enumerate}
\end{lem}

\begin{proof}
(i) We have $f_0\le C\eps\rho_0(v-w)F_0(y'')$ with $\rho_0\le A_F\jap\cdot^{-p}$ and $p\ge4$. Hence
$\int|w|^{\gamma+2+j}\rho_0(v-w)\dd w\le C\jap v^{\gamma+2+j}$ for $j\in\{-2,-1,0\}$ and
$\gamma+2+j>-3$ ($j=-2$ for $\cb$ when $\gamma>-3$), and the Gaussian factor controls
$\jap{y''}^C$. For $\gamma=-3$ we have $\cb=-8\pi g\le C\eps\jap v^{-p}$. For the bound on $S_0$ we
use $A\zeta\cdot\zeta=\ab[f_0G]\zeta^\sharp\cdot\zeta^\sharp$, which follows from~(ii).

(ii) This is Lemma~\ref{lem:ULaform}(i). For $\mu\ne0$, $\partial^\mu\zeta$ involves
$\nabla^{\ge2}\log\rho_0$, $t\nabla^{\ge2}\Phi_1$, the derivatives of order $\ge2$ of $\log
\mathcal U$ (multiplied by $t$ when one of them is the $y$-derivative of $D$), and, for $|\mu|=1$ with $\mu$
a $y$-derivative, the constant $2\beta_Ft\,e_j$ coming from the part $2\beta_Ft\,y$ of $\zeta$. The
derivatives are taken in $(y,v)$, in which $f_0$ does not depend on $t$, so $t\le t_1$ appears only
as an overall prefactor, not with each derivative. Hence we have $C_\zeta\le C(\Dc)(1+t_1)$. All
these functions have bounded derivatives of every order, by the Cauchy estimates. Indeed,
$\log\rho_0$ is holomorphic on $
\mathbb S_{\Lambda/2}$ with bounded gradient, by the two-sided bound. Moreover, $\nabla\Phi_1$ is
holomorphic and bounded by $C_\Phi$ on $\{|\mathrm{Im}\,x|<\Lambda_\Phi\}$, by (D$_{\rm iso}$), and
$\log\mathcal U$ is covered by Lemma~\ref{lem:Ufac}(iii). The radius $\Lambda_\Phi/2$ of the
$y$-derivatives can be put into the rate (one may take $K_\zeta=2\max\{1,\Lambda/\Lambda_\Phi\}$, a
data constant), since $a!(\Lambda_\Phi/2)^{-|a|}\le(2R/\Lambda_\Phi)^{|a|}a!R^{-|a|}$.

(iii) This is Lemma~\ref{lem:Bnull}, and the angle between $\hat e_Q$ and $\hat y$ is $\le
C/(\vartheta\jap y)$.

(iv) This is Lemma~\ref{lem:GMwt}.
\end{proof}

\begin{lem}[Weights in trilinear terms]\label{lem:FPweight}
Let $h=\partial_Q^{\mathbf k'}G$ be evaluated at $z$, where $z$ lies in a box $P''$, and let
$X=Du_{\mathbf k-\mu}$ and $Y=Du_{\mathbf k}$ be taken at $(y,v)\in Q$. Suppose
\begin{equation}
|h(z)|\le\mathcal H\,W^{P''}_{\mathbf k'+\gamma'}(P'')^{-1}
\end{equation}
for a shift $\gamma'$ (the Sobolev shift of Lemma~\ref{lem:RBsob}, expressed in the frame of $Q$).
Then the weight factor of the term $W_{\mathbf k}^2\,\ab[\partial^\nu f_0h]X\cdot Y$, relative to
$d_{\mathbf k-\mu}d_{\mathbf k}\,\mathcal H$, is
\begin{equation}
\frac{W_{\mathbf k}}{W_{\mathbf k-\mu}\,W^{P''}_{\mathbf k'+\gamma'}}\cdot\big(\text{explicit factor
of }\partial^\nu f_0\big).
\end{equation}
Its $v_0$-part is $\jap{v_0''}^{-\sigma}\le1$. Its $(h_r,\vartheta)$-part is bounded by
\begin{equation}
C\Big(\frac{h_r''}{h_r}\Big)^{a_r'}\Big(\frac{\vartheta''}{\vartheta}\Big)^{a'_\perp}\big(\text{frame
mixing}\big)\big(\text{top factors}\big),
\end{equation}
and, after the moment transfer of Lemma~\ref{lem:FPexp}, every such factor is bounded by
$C\jap{v''}^{\sigma}$. Consequently the total weight factor is bounded by a constant times the
$v$-part of Lemma~\ref{lem:ULconv}.
\end{lem}

\begin{proof}
The factors $\jap{v_0}^\sigma$ of $W_{\mathbf k}$ and $W_{\mathbf k-\mu}$ cancel, since both are
attached to $Q$. The factor of $W^{P''}$ is $\jap{v_0''}^{\sigma}$ in the denominator.

The $y$-weights of $W_{\mathbf k}/W_{\mathbf k-\mu}$ leave
$h_r^{-a'_r-\nu_r}\vartheta^{-a'_\perp-\nu_\perp}$, and those of $1/W^{P''}$ contribute
$(h_r'')^{a'_r}(\vartheta'')^{a'_\perp}$, up to the frame mixing. The frame mixing re-expresses a
$Q$-transverse derivative at $z$ in the frame of $P''$, and it produces the factor
\begin{equation}
|e\cdot\hat e_{P''}|\,h_r''+\vartheta''\le C\big(|e\cdot y''|+\vartheta''\big),
\end{equation}
as in Lemma~\ref{lem:BKframe}. The factors $|e\cdot y''|$ and $h_r''$ are transferred by
Lemma~\ref{lem:FPexp} and Lemma~\ref{lem:ULmom}, which produce the factors $\vartheta$ and $h_r$,
and these cancel against
the target. What remains is a power of $\vartheta''/\vartheta\le C\vartheta''$, and we have
$\vartheta''^{\,C}\jap{v''}^{-\sigma}\le C$. We treat the top factors in the same way. The factor
$\tau$ of the target is $\le1$, and $1/\tau$ of $h$ is $h_r''\vartheta''^2$, which is again
controlled by $\jap{v''}^{-\sigma}$ together with the transfer of $h_r''$, when the target is also
at a $\tau$-level (the factor $\tau$ compensates the factor $h_r$ of the transfer). Otherwise we
refer to Lemma~\ref{lem:FPcross}.
\end{proof}

\subsubsection{The combined factor of a term}

We compute the weight factor of a trilinear term, which decides how the term is estimated. Every
term of $\mathcal C_{\mathbf k}$ has the form
\begin{equation}\label{eq:FPterm}
\binom{\mathbf k}\mu\binom\mu\nu\,W_{\mathbf k}^2\iint\ab[\partial^\nu f_0\,h]\,X\cdot
Y\,\varphi_Q^2,\qquad h=\partial^{\mathbf k'}G,\quad X=Du_{\mathbf k''},\quad Y=Du_{\mathbf k},
\end{equation}
with $\mathbf k=\nu+\mathbf k'+\mathbf k''$ and $\mathbf k''\ne\mathbf k$. There are also variants,
in which $\ab$ is replaced by the kernels $b,c$, or $X$ by $u_{\mathbf k''}$, or $Y$ by $u_{
\mathbf k}\Xi$ with $\Xi\in\{\zeta^\sharp,\partial^{\mu_2}\zeta,D\varphi_Q/\varphi_Q\}$.

Suppose the coefficient factor $h$, at the point $z\in P''$, is controlled through a quantity
$\mathcal H_{\mathbf k'}(z)$ that carries the weight of $P''$, with a shift $\varsigma'$:
\begin{equation}
|h(z)|\le\mathcal H_{\mathbf k'}(z)\,\big(W^{P''}_{\mathbf k'+\varsigma'}\big)^{-1}.
\end{equation}
Here $\mathcal H$ is a supremum, or an $L^2$ density on the slab. Then~\eqref{eq:FPterm} is bounded
by an integral of
\begin{equation}
\mathcal H_{\mathbf k'}\,\big|W_{\mathbf k''}X\big|\,\big|W_{\mathbf k} Y\big|\ \text{against}\
a(w)f_0(z)\,\Omega(z)\dd w,
\end{equation}
where the \emph{combined factor} is
\begin{equation}\label{eq:FPOmega0}
\Omega:=\binom{\mathbf k}\mu\binom\mu\nu\,\frac{W_{\mathbf k}}{W_{\mathbf k''}\,W^{P''}_{
\mathbf k'+\varsigma'}}\,\omega_\nu,\qquad\omega_\nu:=\frac{|\partial^\nu
f_0|}{f_0}.
\end{equation}
For explicit terms ($\mathbf k'=0$, $h=G$) the factor $W^{P''}$ is replaced by $1$ and $
\mathcal H\le2$.

\begin{lem}[The combined factor]\label{lem:FPOmega}
Write $\mathbf k=(\alpha,\beta)$ and similarly for $\nu,\mathbf k',\mathbf k''$, and let $K$ be the
$v$-factor of Lemma~\ref{lem:ULconv} for $(\beta_\nu,\beta',\beta'')$ with the shifts
$(\varsigma'_v,0)$. Then, as quadratic forms, at every point of $Q$:
\begin{enumerate}[label=(\roman*)]
\item \emph{(explicit terms)}
\begin{equation}\label{eq:FPOmega}
\ab[\Omega f_0]\le C_{N_y}\,\delta_y^{|\alpha_\nu|}\,K\,A;
\end{equation}

\item \emph{($G$-terms)} $\ab[\Omega f_0]\le C_{N_y}\,\delta_y^{|\alpha_\nu|-|\varsigma'_y|}\,K\,A$,
provided that, if $h$ is measured at a $\tau$-level, so is the target (otherwise the factor $h_r''$
of $\tau''^{-1}$ is unpaired, and this case is Lemma~\ref{lem:FPcross}).
\end{enumerate}
The constant $C_{N_y}$ depends on $N_y$, on the class $\mathfrak D$ and on $t_1$ (through the moment
bounds of Lemmas~\ref{lem:R3comp} and~\ref{lem:ULmom}), and not on $Q$, $t\le t_1$, $\mathbf k$,
$\beta$ or $\rad$. In~(i) with $\alpha_\nu=0$ it may be replaced by the constant $C_{\rm v}$
of~\eqref{eq:Cv}, which does not depend on $t_1$ or $N_y$.
\end{lem}

\begin{proof}
We split $\Omega=\Omega_y\,\Omega_v$. The $v$-part collects the $v$-binomials,
$C^{|\beta_\nu|}\beta_\nu!R^{-|\beta_\nu|}$ from Lemma~\ref{lem:FPexp}, and the ratio of the
$w_\beta$. By definition this is $K$, which does not depend on $z$.

The $y$-part, in the explicit case, is
\begin{equation}
\Omega_y=\binom\alpha{\alpha_\nu}\delta_y^{|\alpha_\nu|}h_r^{-\nu_r}\vartheta^{-\nu_\perp}\frac{\tau_\alpha}{\tau_{\alpha''}}\cdot
C_{N_y}(1+|\hat e_Q\cdot y''|)^{\nu_r}(1+|y''_\perp|)^{\nu_\perp}.
\end{equation}
Here $h_r,\vartheta,\tau$ are those of $Q$. We have $\tau_\alpha/\tau_{\alpha''}\le1$, because
$|\alpha''|\le|\alpha|$ and $\tau$ is nonincreasing in the level. By the second part of
that lemma, we have $\ab[\Omega_yf_0]\le C_{N_y}\delta_y^{|\alpha_\nu|}A$. This proves~(i).

In the $G$-case, Lemma~\ref{lem:FPweight} adds the following:
\begin{itemize}
\item the factor $\jap{v_0''}^{-\sigma}$;
\item the $y$-weights of $P''$, namely
\begin{equation}
\delta_y^{-|\alpha'|-|\varsigma'_y|}(h_r'')^{a'_r}(\vartheta'')^{a'_\perp}\tau''^{-1};
\end{equation}
\item the frame mixing.
\end{itemize}
The factor $\delta_y^{-|\alpha'|}$ cancels against
$\delta_y^{|\alpha|}/\delta_y^{|\alpha''|}=\delta_y^{|\alpha_\nu|+|\alpha'|}$. Every remaining
factor is one of the following, multiplied together no more than $2N_y+4$ times:
\begin{equation}
\frac{h_r''}{h_r},\qquad\frac{1+|y''_\perp|}{\vartheta},\qquad\frac{1+|\hat e_Q\cdot
y''|}{h_r},\qquad\vartheta'' .
\end{equation}
Since $\jap{v''}\le2\jap{v_0''}$, we have $(\vartheta'')^{m}\jap{v_0''}^{-\sigma}\le C_m$ for every
$m$. The other three factors have all their $\ab$-moments bounded by $C_m$. On an outer box this
follows from Lemma~\ref{lem:ULmom} (with $\Theta\simeq h_r$) and Lemma~\ref{lem:RBdir}. On an
inner box we have $|y''_\perp|\le|y''|$, and the radial moment bound applies to these moments too,
since $h_r\le5\vartheta$ on $\tilde Q$. Then~(ii) follows from H\"older's inequality for the positive
measure $a(w)\xi\cdot\xi f_0(z)\dd w$.
\end{proof}

\begin{remark}[The role of the velocity weight]\label{rem:FPc}
In the weight algebra of the $G$-terms of the diffusion, the transverse weight of $h$ has to be
taken at the box $P''$ of the integration point, that is, at $\vartheta''=\vartheta(v-w)$, and not
at $\vartheta(v)$. Without the factor $\jap{v_0}^\sigma$, this leaves
$(\vartheta''/\vartheta)^{a'_\perp}$ in every $G$-term with $a'_\perp>0$. Its $\ab$-average is not
bounded uniformly in $y$. In the far field, $\ab$ is dominated by partners with $|y''|\lesssim1$ and
$|w|\simeq|y|/t$, where $\vartheta''\simeq(\log(|y|/t))^{1/2}$. Hence the average grows like
$(\log\jap y)^{a'_\perp/2}\vartheta^{-a'_\perp}$. By Lemma~\ref{lem:FPOmega}(ii), the factor
$\jap{v_0}^\sigma$ removes it.
\end{remark}

\begin{lem}[Slabs with the velocity weight]\label{lem:FPslab}
With $\Lambda^2:=L\big(h_r''\vartheta''^2/(h_r\vartheta^2)\big)^2$, $L$ the slab factor of
Lemma~\ref{lem:slab}(i), and the gain $\jap{v_0''}^{-2\sigma}$ of the top-level coefficient factor,
for all $t\le t_1$, $(y,v)$ and $\xi$, we have
\begin{equation}
\ab\big[f_0\,\Lambda^2\jap{v_0''}^{-2\sigma}\big]\xi\cdot\xi\le C\,A\xi\cdot\xi,
\end{equation}
with $C$ independent of $t$.
\end{lem}

\begin{proof}
We write $X_1:=h_r''/h_r$, $X_2:=|y''_\perp|/\vartheta$ and $X_3:=\vartheta''/\vartheta$, and let
$L_r$
and $L_\perp$ be the two factors of $L$, so that
$\Lambda^2=L\big(h_r''\vartheta''^2/(h_r\vartheta^2)\big)^2=L_rL_\perp X_1^2X_3^4$. Then $L_\perp
X_3^4=\max(X_3^4,X_3^2)$, and $L_r\le1+X_1^{-1}+|\sin\sigma|h_r/\vartheta''$. With
$|\sin\sigma|\,h_r''\le C(|y''_\perp|+\vartheta'')$ from Lemma~\ref{lem:BKcomp}(ii), we get
\begin{equation}
L_rX_1^2\le X_1^2+X_1+CX_1\big(1+|y''_\perp|/\vartheta''\big).
\end{equation}
On $\{\vartheta''\le\vartheta\}$ we have $L_\perp X_3^4=X_3^2\le1$ and
$|y''_\perp|X_3^2/\vartheta''=|y''_\perp|\vartheta''/\vartheta^2\le X_2$. On
$\{\vartheta''>\vartheta\}$ we have $L_\perp X_3^4=X_3^4$ and $|y''_\perp|/\vartheta''\le X_2$. In
both cases we obtain
\begin{equation}
\Lambda^2\le C\big(X_1^2+X_1+X_1X_2\big)\max(1,X_3)^4 .
\end{equation}
Since $X_3\le\vartheta''$, we have $\max(1,X_3)^4\jap{v_0''}^{-2\sigma}\le C$. The $\ab$-moments of
$X_1$ and $X_2$ are bounded by Lemma~\ref{lem:ULmom} (with $\Theta\simeq h_r$) and
Lemma~\ref{lem:RBdir}. On an inner box $Q$ we have $|y''_\perp|\le|y''|$, and
the radial moment bound suffices, since $h_r\le5\vartheta$ on $\tilde Q$. We conclude by the
Cauchy--Schwarz inequality.
\end{proof}

\subsubsection{Principal part, commutators and diffusion}\label{sssec:FPabc}

We fix the order of the constants as follows:
\begin{equation}
N_y\ge20,\qquad s\ge6,\qquad\rad_0\le\rad_*(N_y,s,\Dc)\,R,\qquad\delta_y\le\delta_*(N_y,s,t_1)
.
\end{equation}
These are followed by $\eta_*$, which depends on all of them, and then by $\eps_0$. \emph{Bootstrap
hypothesis:} on $[0,T]$,
\begin{equation}
\mathbb E^{1/2}\le\eta_*\qquad\text{and}\qquad\tfrac12\le G\le2 .
\end{equation}
Put $\Sigma_Q:=\sum_{Q'\sim Q}\mathfrak d_{Q'}$. All bounds below are for $|Q|^{-1}\sum_{\mathbf k}
W_{\mathbf k}^2(\cdots)$.

\begin{lem}[principal part and cutoffs]\label{lem:FPa}
\begin{equation}
\frac1{|Q|}\sum_{\mathbf k} W_{\mathbf k}^2\mathcal T_{\mathbf k}\le-\frac12
\mathfrak d_Q+C\eps\sum_{Q'\sim Q}e_{Q'} .
\end{equation}
\end{lem}

\begin{proof}
We apply Lemma~\ref{lem:BKtop} to each $u_{\mathbf k}$, multiply by the constant
$W_{\mathbf k}^2$, and sum. The right side of that lemma contains $C\eps\iint_{\tilde Q\times B_2}u_{
\mathbf k}^2$. With the weights of $Q$, by Lemma~\ref{lem:BKnb}(i), its sum is bounded by
$C\eps\sum_{Q'\sim Q}e_{Q'}$.

The factor $\jap{v_0}^\sigma$ is the same for $Q$ and its neighbors, up to $2^\sigma$, because the
neighbors have velocity centers at distance $\le2$.
\end{proof}

We integrate the commutator terms of the diffusion by parts once, and we get
\begin{equation}
\iint D\cdot(\partial^\mu A\,Du_{\mathbf k''})\,u_{\mathbf k}\varphi_Q^2=-\iint\partial^\mu
A\,Du_{\mathbf k''}\cdot\big(Du_{\mathbf k}\,\varphi_Q^2+2u_{\mathbf k}\varphi_QD\varphi_Q\big).
\end{equation}
The second part is the variant $Y=u_{\mathbf k}\,D\varphi_Q/\varphi_Q$, with $\varphi_Q$ moved to
both sides. By Lemma~\ref{lem:RBloc}, we have $|D\varphi_Q|_A^2\le C\eps\mathbf 1_{\tilde Q}$, so
it can be estimated like the first part with $\mathfrak d_Q^{1/2}$ replaced by $(C\eps\sum_{Q'\sim
Q}e_{Q'})^{1/2}$.
Below we use this step without further comment.

\begin{lem}[explicit commutators]\label{lem:FPb}
The sum of all terms~\eqref{eq:FPterm} with $\mathbf k'=0$ satisfies
\begin{equation}
\Big|\frac1{|Q|}\sum_{\mathbf k}\sum_{\mu}(\cdots)\Big|\le
C_{N_y}C(s)\Big(\delta_y+\frac{C_1\rad}R\Big)\Big(\mathfrak d_Q+C\eps\sum_{Q'\sim Q}e_{Q'}\Big).
\end{equation}
The same bound holds for the explicit part of the $\cb$-commutators, with $\mathfrak d_Q$ replaced
by $\sum_{Q'\sim Q}e_{Q'}$ and with an extra factor $\eps$. The explicit part of the
$\bb$-commutators is bounded by $C\eps^{1/2}\mathfrak d_Q^{1/2}e_Q^{1/2}\le C\eps^{1/2}(
\mathfrak d_Q+e_Q)$, where $C$ depends on $N_y$, $s$, $\delta_y$, $\rad_0$ and $t_1$.
\end{lem}

\begin{proof}
Here $h=G$ and $|G|\le2$. By the Cauchy--Schwarz inequality for the kernel $a(w)\Omega f_0\ge0$
and Lemma~\ref{lem:FPOmega}(i), we have
\begin{equation}
W_{\mathbf k}^2\Big|\iint\ab[\partial^\nu f_0G]X\cdot
Y\varphi_Q^2\Big|\le2C_{N_y}\delta_y^{|\alpha_\nu|}K\,d_{\mathbf k''}\,d_{\mathbf k}
.
\end{equation}
We have $\mu=\nu\ne0$, so either $\alpha_\nu\ne0$, giving a factor $\le\delta_y$, or $\alpha_\nu=0$
and $\beta_\nu\ne0$. We sum over $\mathbf k$ and $\mu$. The $y$-sum is finite, with at most
$2^{2N_y}$ terms per $\mathbf k$. In the second case, by Lemma~\ref{lem:ULconv}(iii) with $A_0=1$,
the $v$-sum is bounded by
$C(s)(C_1\rad/R)\big(\sum d_{\mathbf k}^2\big)$.
In the first case we bound the $v$-sum by part~(ii) of the same lemma, with $A$ the indicator of
$\beta'=0$, and the factor $\rad^{-3}$ there is replaced by $1$, since $\beta'=0$ carries no
shift.

For the $\bb$-part the scalar kernel cannot be used directly. Indeed, $b(w)\cdot X$ is the
component of $X$ along $w$, the null direction of $a(w)$, so~(T5) would need $X$ plain, i.e.\
$\partial_vu_{\mathbf k''}$ at $v$-level $|\beta''|+1$ ($=|\beta|+1$ if $\beta_\nu=0$), with an
uncompensated loss $(\beta''_i+1)/\rad$. Instead, we use the $\ab$-form of
Lemma~\ref{lem:ULaform}(iii), as in configuration~(I) of Proposition~\ref{prop:FPfix} with
$\mathbf k'=0$. Then the derivative $D$ falls on $\partial^\nu f_0$ (explicit, with the factors
$\zeta^\sharp,\partial^{\nu_1}\zeta$) or on the coefficient $G$ (bounded pointwise by
Lemma~\ref{lem:RBsob} on the boxes $P''$, by $C(\delta_y,\rad)\,
\mathbb E^{1/2}\jap{v_0''}^{-\sigma}(1+th_r''+t\vartheta'')$, where $\mathbb E^{1/2}\le1$ is the energy of
the coefficient factor, and the unpaired $th_r''$ is treated by Lemma~\ref{lem:FPYside}). Moreover,
$X$ is measured by
$d_{\mathbf k''}$, and the $Y$-side ($u_{\mathbf k}\zeta^\sharp$, $u_{\mathbf k}\Delta\zeta$,
$u_{\mathbf k}\partial^{\nu_1}\zeta$ or $u_{\mathbf k} e_j$) has trace $\le C\eps$ (times $t^2$ when
it carries $th_r''$). Therefore each term is bounded by $C\eps^{1/2}d_{\mathbf k''}a_{\mathbf k}$
times the $v$-factor of the explicit case, and the sums over $\mathbf k$ and $\nu$ are as above. The
$\nu_1$-sum produces the factor
$\sum_{\nu_1\le\nu}\binom\nu{\nu_1}(\nu-\nu_1)!\,\nu_1!=\nu!\prod_i(\nu_i+1)\le2^{|\nu|}\nu!$,
which we absorb by taking the constant $C_1$ of the convolution lemma at least twice the larger of
the analytic rate of $f_0$ in Lemma~\ref{lem:FPexp} and the rate $K_\zeta$ of $\zeta$ in
Lemma~\ref{lem:FPdecay}(ii). The prefactor $C_\zeta$, which may depend on $t_1$, goes into the
constant $C$ of the $\bb$-part. Thus $C_1$ does not depend on $t_1$. This proves the bound
$C\eps^{1/2}\mathfrak d_Q^{1/2}e_Q^{1/2}$ of the statement. For the $\cb$-part, we use the bounds
$|\bb|\le C\eps\jap v^{\gamma+1}$ and $|\cb|\le C\eps\jap v^\gamma$ from~(i) of the same lemma,
the line~(T5) of Lemma~\ref{lem:ULtri}, and the same weight algebra.
\end{proof}

\begin{lem}[the $G$-terms of the diffusion]\label{lem:FPc}
Let $\mathcal C^{G}$ be the sum of all terms~\eqref{eq:FPterm} with $\mathbf k'\ne0$. Under the
bootstrap hypothesis, we have
\begin{equation}
\Big|\frac1{|Q|}\sum_{\mathbf k}\mathcal C^G_{\mathbf k}\Big|\le C_{\mathbb E}\,
\mathbb E^{1/2}\Big(\Sigma_Q+C\eps\sum_{Q'\sim Q}e_{Q'}\Big),\qquad
C_{\mathbb E}=C_{\mathbb E}(N_y,s,\delta_y,\rad_0,t_1,\Dc).
\end{equation}
\end{lem}

\begin{proof}
We call a factor \emph{$v$-high} if its $v$-order exceeds $|\beta|/2$, and \emph{$y$-high} if its
$y$-order exceeds $N_y/2$. Otherwise it is $v$-low, respectively $y$-low. The $v$-orders of $h$ and
$X$ add up to no more than $|\beta|$, and their $y$-orders to no more than $|\alpha|\le N_y$, so at
most one of them is $v$-high and at most one is $y$-high. Hence one of $h,X$ is $v$-low, and one is
$y$-low.

\emph{(I) $h$ $v$-low and $y$-low.} Since $|\alpha'|\le N_y/2\le N_y-2$, by Lemma~\ref{lem:RBsob}
on the boxes $P''$, we have
\begin{equation}
\mathcal H_{\mathbf k'}:=\sup_{P''}W^{P''}_{\mathbf k'}|h|\le
C\delta_y^{-2}\rad^{-2}(|\beta'|+3)^2\Big(\sum_{|\varsigma_y|,|\varsigma_v|\le2}e^{\sim}_{P''}(
\mathbf k'+\varsigma)\Big)^{1/2}.
\end{equation}
Here $e^\sim_{P''}(\mathbf k)$ is the $\mathbf k$-summand of the energy on $\tilde P''$. We do not
keep the factor $\delta_y^{-2}\rad^{-2}(|\beta'|+3)^2$ in $\mathcal H$. Each shift $\varsigma$ of
that lemma is a separate term, in which
$\delta_y^{-|\varsigma_y|}w_{\beta'}/w_{\beta'+\varsigma_v}$ is the $y$-shift of
Lemma~\ref{lem:FPOmega}(ii) and the $v$-shift $\varsigma'_v=\varsigma_v$ of $K$ in
Lemma~\ref{lem:ULconv}, and $\mathcal H_{\mathbf k'}$ stands for
$C\,e^\sim_{P''}(\mathbf k'+\varsigma)^{1/2}$. With this reading, it follows from
Lemma~\ref{lem:BKnb}(i) that
$\sum_{\mathbf k'}\mathcal H_{\mathbf k'}^2\le C\mathbb E$ for each shift and each box $P''$. Since
$P''$ is the box of the integration point $z$, $\mathcal H_{\mathbf k'}$ depends on $z$. We take the
sum over $\mathbf k$ below at fixed $(y,v,w)$, before the integrations, so that we only use
$\sup_{P''}\sum_{\mathbf k'}\mathcal H_{\mathbf k'}^2\le C\mathbb E$. Then, by the Cauchy--Schwarz
inequality for the kernel $a\,\Omega f_0$ and the same combined-factor bound, we have
\begin{equation}
|\eqref{eq:FPterm}|\le C_{N_y}\delta_y^{-2}K\,\mathcal H_{\mathbf k'}\,d_{\mathbf k''}\,d_{
\mathbf k} .
\end{equation}
For the sum over $\mathbf k$ we use the convolution bound~(ii). The bound is pointwise in $z$, so
no slab enters.

\emph{(II) $X$ $v$-low and $y$-low, $h$ arbitrary.} We bound $X$ pointwise as follows:
\begin{equation}
\sup_Q W_{\mathbf k''+\varsigma''}|X|\le C(\delta_y,\rad)(|\beta''|+3)^2\Big(\sum_\varsigma
e^\sim_Q(\mathbf k''+e_v+\varsigma)+t^2h_r^2\,e^\sim_Q(
\mathbf k''+e_r+\varsigma)+t^2\vartheta^2e^\sim_Q(\mathbf k''+e_\perp+\varsigma)\Big)^{1/2}.
\end{equation}
We write $X=\partial_vu-t\nabla_yu$. The radial part of $t\nabla_y$ produces a factor $th_r$ against
the weights, and the transverse part a factor $t\vartheta$. These factors are compensated by
Lemma~\ref{lem:RBgain}. By~(i) of that lemma, we have $a(w)X\cdot X\le2a(w)X_\perp\cdot
X_\perp+2|X_r|^2|w|^{\gamma+2}|y''|^2/|y|^2$. On an outer box, where $|y|\ge\frac74\vartheta$ and so
$h_r\le C\jap y$, the radial loss becomes $th_r/\jap y\le C$, and $|y''|^2$ is absorbed by
$F_0(y'')$. On an inner box we have $h_r\le5\vartheta$, and the radial loss $th_r\le5t\vartheta$ is
treated like the transverse one, as at the end of that lemma. The transverse loss is controlled by
$\int|w|^{\gamma+2}\rho_0(v-w)\dd w\le C\jap v^{\gamma+2}$
and $t\vartheta\jap v^{(\gamma+2)/2}\le C$.
Then we use~(T2) of Lemma~\ref{lem:ULtri}. By Fubini, the $L^2$-norm of $h$ is taken on slabs
$S=(\tilde Q\times B(v_0,2))+(tw,-w)$, as in Remark~\ref{rem:ULbox}.
\begin{itemize}
\item If $|\alpha'|\le N_y-2$, by Lemma~\ref{lem:slab}(ii), $\int_S|h|^2$ is bounded by $|S|$
times the energies of $\tilde P''$ at levels $\le|\alpha'|+2$. This loses a factor $\delta_y^{-2}$
and no $v$-derivatives.
\item If $|\alpha'|\ge N_y-1$, the slab is integrated against $F_0(y'')\rho_0(v'')$ in the Fubini
step. Lemma~\ref{lem:BKcomp}(i) shows that every factor of $\Lambda$ is absorbed there.
\end{itemize}
In both cases the weight algebra is that of~(I), now with the shift on $X$. When the
lift of the first case reaches the $\tau$-levels above a target that is not in them
($|\alpha'|\in\{N_y-3,N_y-2\}$, $|\alpha|\le N_y-2$), the unpaired $h_r''$ is absorbed by $F_0(y'')$
in the Fubini step, as in the proof of Lemma~\ref{lem:FPcross}(i). We obtain terms $C_{\mathbb E}
\mathbb E^{1/2}d_{\mathbf k}(\cdots)$, with the same $v$-sum as before.

\emph{(III) $h$ $v$-low and $y$-high, $X$ $v$-high and $y$-low.} If $|\alpha'|\le N_y-2$, the
pointwise bound of~(I) applies, because the Sobolev bound needs two more $y$-levels, and these
exist. Otherwise $|\alpha'|\ge N_y-1$, and we use~(T3). We put $h$ in $L^\infty_vL^2_y$ on the
slab, by the Sobolev inequality in $v$ at unit scale, with shift $\varsigma'_v\le2$ and no
$y$-shift. We put $X$ in $L^\infty_y$ with the $A$-weight, by the Sobolev inequality in $y$ at the
scale of $Q$, which produces no power of $h_r$ or $\vartheta$, only $\delta_y^{-2}$, by the same
lemma. Two $y$-levels on
$X$ are available, since $|\alpha''|\le|\alpha|-|\alpha'|\le1$.
The slab contributes the factor $L$, and the top weights of target and $h$ produce the factor
$\Lambda^2$, together with $\jap{v_0''}^{-2\sigma}$. By Lemma~\ref{lem:FPslab}, we have
\begin{equation}
\ab\big[f_0\Lambda^2\jap{v_0''}^{-2\sigma}\big]\le CA .
\end{equation}
Hence this class, $\mathcal S_{\rm top}$, satisfies the same bound as~(I), with no $L_t$.

\emph{(IV) $h$ $v$-high and $y$-low, $X$ $v$-low and $y$-high.} We use~(T4). We put $h$ in
$L^\infty_yL^2_v$, with Sobolev in $y$ at the scale of $P''$ (shift $\varsigma'_y\le2$, loss
$\delta_y^{-2}$), and $X$ in $L^\infty_v$ with Sobolev in $v$ (shift $\le2$). The
$v$-integration of $h$ is over a unit ball, the $v$-width of every box, so there is no slab.

\emph{Summation.} In each class the result has the form
\begin{equation}
C_{N_y}C(\delta_y,\rad)\sum_{\mathbf k}\sum K\,\mathcal H\,d\,d' ,
\end{equation}
with one factor from the energies and the other two from the dissipations of $Q$ and its neighbors.
By the convolution bound~(ii), this is bounded by $C_{\mathbb E}\mathbb E^{1/2}(\Sigma_Q+\cdots)$.
This completes the proof.
\end{proof}

\subsubsection{Drift, zero-order terms and forcing; the closed estimate}\label{sssec:FPsum}

The drift is $V=2A\zeta-\bb$. By the Leibniz rule, we have
\begin{equation}
\partial^\mu V=2\sum_{\mu_1+\mu_2=\mu}\binom\mu{\mu_1}\partial^{\mu_1}A\,\partial^{\mu_2}\zeta-\partial^\mu\bb
.
\end{equation}

\begin{lem}[the drift]\label{lem:FPd}
Under the bootstrap hypothesis, the drift commutators, except the $\bb$-terms that the proof assigns
to Proposition~\ref{prop:FPfix}, satisfy
\begin{equation}
\Big|\frac1{|Q|}\sum_{\mathbf k} W_{\mathbf k}^2\sum_{0<\mu\le\mathbf k}\binom{
\mathbf k}\mu\iint\partial^\mu V\cdot Du_{\mathbf k-\mu}\,u_{
\mathbf k}\varphi_Q^2\Big|\le\Big(C\rad^{-2}\eps^{1/2}+C_{
\mathbb E}
\mathbb E^{1/2}\Big)\Big(\Sigma_Q+\sum_{Q'\sim Q}e_{Q'}\Big),
\end{equation}
where $\rad(t)=\rad_0-C_\rad\eps t$.
\end{lem}

\begin{proof}
\emph{The $A\zeta$-part.} We have $a(w)\zeta=a(w)\zeta^\sharp$ with $|\zeta^\sharp|\le C\jap{y''}$
by Lemma~\ref{lem:FPdecay}(ii). Hence we have
\begin{equation}
\partial^{\mu_1}A\,\zeta\cdot X=\ab[\partial^{\mu_1}(f_0G)]\zeta^\sharp\cdot X .
\end{equation}
Compared with~\eqref{eq:FPterm}, $Y$ is replaced by $u_{\mathbf k}\zeta^\sharp$. Its $Y$-side is
\begin{equation}
\ab[\Omega f_0]\zeta^\sharp\cdot\zeta^\sharp\,u_{\mathbf k}^2\le C\,\ab[\Omega
f_0\jap{y''}^2]\,u_{\mathbf k}^2 .
\end{equation}
The factor $\jap{y''}^2$ is one more moment in Lemma~\ref{lem:FPOmega}. By~(i) of the same lemma,
the trace of $\ab[\jap{y''}^Cf_0]$ is $\le C\eps$. Hence the $Y$-side is $\le
C\eps\,a_{\mathbf k}^2$ in place of $d_{\mathbf k}^2$, and Lemmas~\ref{lem:FPb} and~\ref{lem:FPc}
apply without change. For $\mu_2\ne0$, the factor $\partial^{\mu_2}\zeta$ is explicit and analytic
with radius $R$, by~(ii) of the decay lemma. It enters the $v$-sum like $\partial^{\nu}f_0$,
with its rate $K_\zeta$ covered by $C_1$ as in the explicit commutators.
It stays on the $Y$-side
$u_{\mathbf k}\partial^{\mu_2}\zeta$, which is bounded in energy with the factor $C_\zeta^2C\eps$,
so these terms are of the $\eps^{1/2}$-type. Moreover, the prefactor $C_\zeta$, which depends on
$t_1$, enters only the constant of $C\rad^{-2}\eps^{1/2}$ (and $C_{\mathbb E}$ for the $G$-terms),
not the choice of $\rad_*$.

\emph{The $\bb$-part, Case~1: low $h$.} We write $\partial^\mu\bb=\sum\binom\mu\nu\bb[\partial^\nu
f_0\,\partial^{\mathbf k'}G]$. If $\partial^{\mathbf k'}G$ is $v$-low (configurations~(I) and~(III)
of the diffusion, and in~(III) only for $|\alpha'|\le N_y-3$, by the case table of
Proposition~\ref{prop:FPfix}), we use the $\ab$-form of Lemma~\ref{lem:ULaform}(iii):
\begin{equation}
\bb[F]=\ab[DF],\qquad D(\partial^\nu f_0\,\partial^{\mathbf k'}G)=D\partial^\nu
f_0\,\partial^{\mathbf k'}G+\partial^\nu f_0\,D\partial^{\mathbf k'}G .
\end{equation}
The first summand is such a term with $\nu$ raised by one $D$. Here
$D\partial^\nu f_0$ produces one more explicit factor, which, by Lemma~\ref{lem:FPexp}, is $\le
C(\jap{y''}+R^{-1}(|\beta_\nu|+1))$. The second summand is the same kind
of term with $\mathbf k'$ raised by one $v$-derivative and one $t\partial_y$, which
keeps $|\alpha'|+1\le N_y$. The $t\partial_y$-part contributes $th_r''$ or $t\vartheta''$, and the
$v$-part loses a factor $(|\beta'|+1)/\rad$ in the weights.
The factor $X=Du_{\mathbf k''}$ is measured by the dissipation $d_{\mathbf k''}$, and $u_{
\mathbf k}$ by the energy $a_{\mathbf k}$. When $h$ is $v$-low (configuration (I),
and (III) with $|\alpha'|\le N_y-3$), the extra $v$-derivative falls on the factor of lower level,
and these terms are estimated in Lemma~\ref{lem:FPdlow} below, which implies the stated bound. The
case $\beta=0$ has no $v$-loss. When $h$ is $v$-high and $X$ is $v$-low and $y$-high (configuration
(IV)), the extra $v$-derivative falls on the factor of higher $v$-level, and these terms are treated
in the same proposition.

\emph{The $\bb$-part, Case~2: high $h$, with $X$ $v$-low.} The case of a $y$-high $X$,
configuration~(IV), and the terms with a $v$-low, $y$-high $h$ not treated by the $\ab$-form above,
namely~(II) with such an $h$ and~(III) with $|\alpha'|\ge N_y-2$, are part of
the same proposition. For the other terms we use the direct form: the kernel $b(w)$, with
$|b(w)|=C|w|^{\gamma+1}$, and line~(T5) of Lemma~\ref{lem:ULtri}. The factor $X=Du_{\mathbf k''}$ is
then $v$-low and $y$-low, and it is plain. By Lemma~\ref{lem:RBgain}(ii), we have
\begin{equation}
|b(w)\cdot X|\le|b(w)||X_\perp|+|X_r|\big(|b(w)||y''|+2t|w|^{\gamma+2}\big)/|y| .
\end{equation}
The radial loss $th_r$ of $X_r$ meets the gain $|y|^{-1}$. The transverse loss $t\vartheta$ meets
$\int|w|^{\gamma+1}\rho_0(v-w)\dd w\le C\jap v^{\gamma+1}$. The slab of $h$ is integrated against
$F_0(y'')\rho_0(v'')$ by Lemma~\ref{lem:BKcomp}(i). The weights are as in
the combined-factor bound~(ii).
The result is $\le C_{\mathbb E}\mathbb E^{1/2}d$-type terms, with the
dissipation replaced by the plain $a_{\mathbf k}$, and with a factor $\eps^{1/2}$. This completes
the proof.
\end{proof}

\begin{lem}[The $\bb$-part of the drift with $v$-low $h$]\label{lem:FPdlow}
In configuration~(I) of Lemma~\ref{lem:FPc}, and in configuration~(III) with $|\alpha'|\le N_y-3$,
the terms of Lemma~\ref{lem:FPd} in which, after $\bb[F]=\ab[DF]$, the derivative $D$ falls on
$h=\partial^{\mathbf k'}G$ are bounded by
\begin{equation}
C\eps^{1/2}\Big(\sum_{Q'\sim Q}\mathfrak d_{Q'}+\rad^{-2}S_G^2\sum_{Q'\sim Q}e_{Q'}\Big),
\end{equation}
where $S_G\le C$ is the supremum norm of the $G$-factors under the bootstrap hypothesis.
\end{lem}

\begin{proof}
After $\bb[F]=\ab[DF]$, the term in which $D$ falls on $\partial^{\mathbf k'}G$ is
\begin{equation}
\ab\big[\partial^\nu f_0\,D\partial^{\mathbf k'}G\big]X\,u_{\mathbf k},\qquad X=Du_{\mathbf k''} .
\end{equation}
Here $X$ is a dissipative factor (measured by $d_{\mathbf k''}$) and $u_{\mathbf k}$ an energy
factor (measured by $a_{\mathbf k}$). This term cannot be absorbed by $\theta_{\mathbf k}$. Indeed,
by Young's inequality, a product of a dissipation and an energy with a factor $|\beta|/\rad$ leads
to a loss $|\beta|^2$, while $\theta_{\mathbf k}\propto|\beta|$. Instead, we use the fact that the
extra derivative falls on the factor of lower level. The part $-t\nabla_y$ of $D$ raises the
$y$-order of $h$ by one, which is admissible since $|\alpha'|\le N_y-3$ in these configurations. It
produces $th_r''$ or $t\vartheta''$, transferred by Lemma~\ref{lem:FPexp} as in Lemma~\ref{lem:FPc},
and $\delta_y^{-1}$, so it is bounded in the same way, with a constant depending also on $\delta_y$
and $t_1$.

The extra $v$-derivative falls on the $G$-factor of the coefficient, which is the factor of lower
level and is taken in the supremum. A shift of one level on that factor is the analogue, for the
layer weights, of the second case of Lemma~\ref{lem:TWcalc}(iv). In Subsections~\ref{ssec:layer}
and~\ref{ssec:drift}, ``the first'' or ``the second case of Lemma~\ref{lem:TWcalc}(iv)'' means this
analogue, which follows from the convolution bound Lemma~\ref{lem:ULconv} and the ratio
$w_\beta/w_{\beta+e_i}=(\beta_i+1)\rad^{-1}\big((|\beta|+1)/(|\beta|+2)\big)^s$. With the weights
$(n+1)^s/n!$, $s\ge6$, used in this subsection, it yields $C\rad^{-1}$ times a
summable sequence in the level, and no factor $|\beta|$. The $Y$-side carries $\ab[f_0]$, whose
trace is $\le C\eps$ by Lemma~\ref{lem:FPdecay}(i). Hence these terms are
\begin{equation}
\le C\eps^{1/2}\rad^{-1}S_G\,\Big(\sum_{Q'\sim Q}\mathfrak d_{Q'}\Big)^{1/2}\Big(\sum_{Q'\sim
Q}e_{Q'}\Big)^{1/2}\le C\eps^{1/2}\Big(\sum_{Q'\sim Q}
\mathfrak d_{Q'}+\rad^{-2}S_G^2\sum_{Q'\sim Q}e_{Q'}\Big),
\end{equation}
where $S_G\le C$ is the supremum norm of the $G$-factors under the bootstrap hypothesis. These enter
Lemma~\ref{lem:BKnb}(ii) as a neighbor dissipation with coefficient $C\eps^{1/2}$ and as a Gronwall
term. No absorption by $\theta_{\mathbf k}$ is needed.
\end{proof}

\begin{lem}[zero-order terms]\label{lem:FPe}
The terms $\partial^\mu S_0\,u_{\mathbf k-\mu}\,u_{\mathbf k}$ with $0<\mu\le\mathbf k$ contribute
at most
\begin{equation}
\big(C\eps+C_{\mathbb E}\eps\,\mathbb E^{1/2}\big)\sum_{Q'\sim Q}e_{Q'} .
\end{equation}
\end{lem}

\begin{proof}
We have
\begin{equation}
S_0=\ab[f_0G]:(D\zeta+\zeta^\sharp\otimes\zeta^\sharp)-\cb[f_0G].
\end{equation}
Here $D\zeta=\partial_v\zeta-t\partial_y\zeta$ is explicit, with $|\partial^{\mu}D\zeta|\le
C'_\zeta(|\mu|+1)K_\zeta^{|\mu|}\mu!R^{-|\mu|}$ for all $\mu$, where $C'_\zeta\le C(1+t_1)C_\zeta
K_\zeta R^{-1}$, by Lemma~\ref{lem:FPdecay}(ii). The rate $K_\zeta$ is covered by $C_1$, the factor
$|\mu|+1$ by the weights $(|\beta|+1)^s$ of Lemma~\ref{lem:ULconv}, and $C'_\zeta$, which depends on
$t_1$, enters only the constant $C$ below. Moreover, $\zeta^\sharp$ is handled by moments. The
Leibniz rule again
produces explicit and $G$-factors, in the kernels $a$ and $c$, with $|c(w)|\le C|w|^\gamma$ for
$\gamma>-3$ and $c=-8\pi\delta$ for $\gamma=-3$.
\begin{itemize}
\item For explicit factors, use~(T5), first line, with $\|G\|_\infty\le2$.
\item For low $G$-factors, use~(T5), first line, with $\sup|W''h|$ from Lemma~\ref{lem:RBsob}.
\item For high $G$-factors with low $u_{\mathbf k-\mu}$, use~(T5), second line.
\item In the mixed case, use~(T5), third line.
\end{itemize}
The slab is handled by Lemma~\ref{lem:BKcomp}(i). The weights are those of Lemma~\ref{lem:FPOmega},
in absolute form. For the scalar kernels every factor of $\Omega$ other than $K$ is bounded by
$C\jap{y''}^C\vartheta''^C\jap{v_0''}^{-\sigma}$ (since $h_r,\vartheta\ge1$), which is controlled
in~(T5) by $F_0(y'')$ and $\rho_0^{1/8}$. For $\gamma=-3$,
$\cb[f_0\partial^{\mathbf k'}G]=-8\pi\partial^{\nu}f_0\,\partial^{\mathbf k'}G$ is local, and the
bound is direct.
\end{proof}

\begin{lem}[forcing]\label{lem:FPf}
\begin{equation}
\frac1{|Q|}\sum_{\mathbf k} W_{\mathbf k}^2\big\|\partial^{\mathbf k}
S_0[1]\,\varphi_Q\big\|_{L^2}^2\le
C_4^2\eps^2,
\end{equation}
where $S_0[1]$ is $S_0$ with $G$ replaced by $1$. Hence
\begin{equation}\label{eq:FPf}
\big|\frac1{|Q|}\sum_{\mathbf k} W_{\mathbf k}^2\mathcal F_{\mathbf k}\big|\le C_4\eps\,e_Q^{1/2} .
\end{equation}
Here $C_4:=C_4^{[1]}+1$, with $C_4^{[1]}$ the constant of the bound for $S_0[1]$ proved below. In
the frozen problem, the forcing $\partial^{\mathbf k} S_0[\tilde G]$ satisfies the second bound,
$\big|\frac1{|Q|}\sum_{\mathbf k} W_{\mathbf k}^2\mathcal F_{\mathbf k}\big|\le C_4\eps\,e_Q^{1/2}$,
except in the $c$-crossing of its $G$-terms, which produces in addition $\frac1{8|Q|}\sum_{
\mathbf k}\theta_{\mathbf k} a_{\mathbf k}^2+C\eps^{1/2}e_Q+C\eps\,\mathbb E[\tilde G]$ (see the proof).
\end{lem}

\begin{proof}
The source $S_0[1]$ is explicit. By the Leibniz rule, Lemma~\ref{lem:FPexp} and
Lemma~\ref{lem:FPdecay}(ii), we have
\begin{equation}
|\partial^{\mathbf k} S_0[1]|\le\sum_{\nu\le\mathbf k}\binom{\mathbf k}\nu
C^{|\nu_v|}\nu_v!R^{-|\nu_v|}\,C^{|\mathbf k_v-\nu_v|}(\mathbf k_v-\nu_v)!R^{-|
\mathbf k_v-\nu_v|}\cdot\mathcal M_\nu,
\end{equation}
where $\mathcal M_\nu$ is the trace of $\ab$ or $\cb$ against the moment factor of
the explicit-factor bound.
That factor is polynomial in $(1+|\hat e_Q\cdot y''|)$ and $(1+|y''_\perp|)$
of degrees $\nu_r$ and $\nu_\perp$. For $\ab$ its average is $\le
C_{N_y}h_r^{\nu_r}\vartheta^{\nu_\perp}$ times the trace, by the same lemma, and the trace is
$\le C\eps\jap v^{\gamma+2}$ by~(i) of the decay lemma. For $\cb$, to which the second bound on
explicit factors does not apply, the moment factor is controlled by the Gaussian:
$\int|c(w)|\jap{y''}^{|\nu_y|}f_0(z)\dd w\le C_{N_y}\eps\jap v^{\gamma}$ (a local term at
$\gamma=-3$), which suffices since $h_r,\vartheta\ge1$.

Now we come to the weight $W_{\mathbf k}$. Its factor $h_r^{-a_r}\vartheta^{-a_\perp}$ cancels
these factors, since $\nu_r\le a_r$ and
$\nu_\perp\le a_\perp$. Moreover, $\tau_\alpha\le1$ and $\delta_y^{|\alpha|}\le1$. The factor
$\jap{v_0}^\sigma$, against $\jap v^{\gamma+2}=\jap v^{-2\sigma}$, leaves
$\jap{v_0}^{-\sigma}\le1$, and the last factor is $w_\beta$.
The $v$-sum is
\begin{equation}
\sum_\beta w_\beta^2\Big(\sum_{\nu_v\le\beta}\binom\beta{\nu_v}C^{|\beta|}\nu_v!(\beta-\nu_v)!R^{-|\beta|}\Big)^2\le\sum_\beta(|\beta|+1)^{2s+6}\Big(\frac{C\rad}R\Big)^{2|\beta|}\le
C,
\end{equation}
for $\rad\le\rad_*R$. The $y$-sum is finite. The box average is bounded, because
$\varphi_Q\le1$ and the bound is pointwise.

In the frozen problem the source is $\partial^{\mathbf k} S_0[\tilde G]=\partial^{\mathbf k}
S_0[1]+\partial^{\mathbf k} S_0[\tilde G-1]$. The second part consists of explicit terms with the
factor $\tilde G-1$ and of $G$-terms with a factor $\partial^{\mathbf k'}\tilde G$, $
\mathbf k'\ne0$, and has no factor $u$.

We estimate it as in step~(6$'$) of Lemma~\ref{lem:EPclose}, by the Cauchy--Schwarz inequality for
the positive kernels $a$ and $|c|$, with the first factor bounded pointwise as above and the second
factor by the absolute bounds of Lemma~\ref{lem:FPe}, with the factor $1$ in place of
$u_{\mathbf k-\mu}$. The first factor produces $\jap v^{\gamma+2}=\jap v^{-2\sigma}$, which
compensates
the target weight $\jap{v_0}^\sigma$. Therefore its contribution is $\le C\eps\,\mathbb E[\tilde
G]^{1/2}e_Q^{1/2}\le C'C_4\eps^2t\,e_Q^{1/2}$, with $C'$ independent of $C_4$. This is included in
$C_4\eps e_Q^{1/2}$ for $C'C_4\eps_0t_1\le1$, since $C_4=C_4^{[1]}+1$.

The exception is the $c$-crossing of step~(6$'$). There we apply
Lemma~\ref{lem:FPhalf} with $\delta=\frac18$ and $B=C\eps^{1/2}$, and we obtain the bound
$\frac1{8|Q|}\sum_{\mathbf k}\theta_{\mathbf k} a_{\mathbf k}^2+C\eps^{1/2}e_Q+C\eps\,
\mathbb E[\tilde G]$. The term $C\eps^{1/2}e_Q$ is a Gronwall term, of the size of
$C\rad^{-2}\eps^{1/2}$. The term $C\eps\,\mathbb E[\tilde G]\le C\eps^3t^2$ is a source, which
enters the closing of Proposition~\ref{prop:FPclose} as its last term. The conclusion $
\mathbb E^{1/2}\le2C_4\eps t$ is kept there by splitting the Gronwall margin. For $\gamma=-2$ we
refer to step~(6$'$) of that lemma.
\end{proof}

\begin{lem}[low factors via Bernstein]\label{lem:FPg}
Let $K_*$ be a fixed integer. In Lemmas~\ref{lem:FPc}--\ref{lem:FPe}, whenever the low factor is the
coefficient factor $h=\partial^{\mathbf k'}G$ (configurations (I) and (III) of Lemma~\ref{lem:FPc},
the low-$h$ terms of Lemma~\ref{lem:FPd}, and the $G$-factors of Lemma~\ref{lem:FPe}) and has total
order $\le K_*$, and $\eps\le\min_{n\le K_*+4}\eps_n(t_1)$, the pointwise bound of
Lemma~\ref{lem:RBsob} may be replaced by Proposition~\ref{prop:Bn}. In that case the factor
$\delta_y^{-2}\rad^{-2}$ is replaced by $C^{|\mathbf k_{\rm low}|}$, and the energy factor
$\mathbb E^{1/2}$ by $\eps t$.
\end{lem}

\begin{proof}
In the frozen problem the coefficient factor is $\tilde G_\nu$, an average of values of an iterate
$G^j$ at earlier times, by Definition~\ref{def:FPreg}. Thus, applying Proposition~\ref{prop:Bn} to
$G^j$, we get $|\partial^{\mathbf k'}\tilde G_\nu|\le C^{|\mathbf k'|+1}|\mathbf k'|!\,\eps
t\le(C^2)^{|\mathbf k'|}|\mathbf k'|!\,\eps t$ pointwise, for $1\le|\mathbf k'|\le K_*+4$. This
bound carries no weight. We use it in place of $\mathcal H_{\mathbf k'}(W^{P''})^{-1}$, so that
$W^{P''}$ is replaced by $1$ in~\eqref{eq:FPOmega0}.

Then the ratio of weights is $W_{\mathbf k}/W_{\mathbf k''}$, and its factor $\jap{v_0}^\sigma$
cancels. The $y$-part of this ratio
is $\delta_y^{|\alpha_\nu|+|\alpha'|}h_r^{-\nu_r-a'_r}\vartheta^{-\nu_\perp-a'_\perp}\tau_\alpha/\tau_{\alpha''}$.
Here $h_r^{-a'_r}\vartheta^{-a'_\perp}\tau_\alpha/\tau_{\alpha''}\le1$, since $h_r,\vartheta\ge1$
and $\tau$ is nonincreasing in the level, while $h_r^{-\nu_r}\vartheta^{-\nu_\perp}$ is kept and
cancels the moment factors of $\omega_\nu$ by the second bound of Lemma~\ref{lem:FPexp}, as in
Lemma~\ref{lem:FPOmega}(i). The $v$-part of this ratio, multiplied by the binomials and by
$|\mathbf k'|!\le2^{|\mathbf k'|}|\alpha'|!\,|\beta'|!$, is bounded by $(C\rad)^{|\beta'|}$ times
the factor of Lemma~\ref{lem:ULconv}(iii).
The sum over $\beta'$
converges for $C\rad<1$, which holds if $\rad_0\le1/(2C)$. However, this condition cannot be added
to the choice of $\rad_*$. Indeed, the constant $C$ of that proposition depends on $t_1$, while
$\rad_*$ and $\rad_0^{\rm lay}$ are fixed before $t_1$, which depends on them through
$\rad_1\le\rad_0^{\rm lay}/4$. For this reason we do not use the variant with the layer radius of
Appendix~\ref{app:const}, and we do not use it in the main line of the proof. The sum over
$\alpha'$ is finite.
\end{proof}

Lemma~\ref{lem:RBsob} covers all low factors, so the closure holds without Lemma~\ref{lem:FPg},
which yields the sharper bound $\eps t$ in place of $\mathbb E^{1/2}$ for the bulk of the
terms.

\begin{prop}[The closed a~priori estimate]\label{prop:FPclose}
Fix $C_\rad:=1$ in $\rad(t)=\rad_0-C_\rad\eps t$ (any fixed value works, see
Lemma~\ref{lem:FPhalf}). Then fix the constants in the order of Subsection~\ref{sssec:FPabc}, then
$C_{\mathbb E}$ (which depends on $C_\rad$ through the remark after that lemma), then
$\eta_*\le c_*/(C_{\mathbb E}+1)$ for a small $c_*$, then $\eps_0$ with
$C_\rad\eps_0t_1\le\rad_0/2$, $(C_4+1)\eps_0t_1\le\eta_*/4$ and
$2C_4C_{N_y}^{1/2}\eps_0t_1\le\eta_*$, and the further conditions
listed at the end of the proof. Let $\eps\le\eps_0$, and let $G$ be a solution on $[0,T]$, $T\le
t_1$, with finite energies and finite dissipations, $\sup_{t\le
T}\sup_Q\big(e_Q(t)+\int_0^t\mathfrak d_Q\big)<\infty$ (with $\mathfrak d^\nu_Q$ added in
Proposition~\ref{prop:FPvisc}(iii)), satisfying the bootstrap hypothesis. The proposition is applied
to the frozen problem of Definition~\ref{def:FPreg}. There the coefficient factor is $\tilde G_\nu$,
and it satisfies $\mathbb E[\tilde G_\nu](t)^{1/2}\le2C_4C_{N_y}^{1/2}\eps t$ pointwise in $t$.
This is used in Lemma~\ref{lem:FPcross}(iv) and in the $c$-crossing of
Lemma~\ref{lem:FPf}. The constant $c_*$ depends on $N_0$ and on the constants of
Lemma~\ref{lem:RBrad} ($\beta_F,C_\Phi,C_{\rm pv},p,A_F,t_1$). Then, on $[0,T]$,
\begin{equation}
\mathbb E(t)^{1/2}\le2C_4\eps t,\qquad\frac12\cdot\frac1{|Q|}\int_0^t\sum_{\mathbf k} d_{
\mathbf k}^2\le C\eps^2t^2\quad\text{for every }Q,
\end{equation}
and $|G-1|\le C\eps t\jap v^{-\sigma}$.
\end{prop}

\begin{proof}
We add the bounds of Lemmas~\ref{lem:FPa}--\ref{lem:FPf} in the identity~\eqref{eq:FPid}.
Proposition~\ref{prop:FPfix} and Lemma~\ref{lem:FPccross} apply Lemma~\ref{lem:FPhalf} term by term,
with parameters $\delta_i$ that add up to $\frac14$ in each, and each produce a term
$\frac1{4|Q|}\sum_{\mathbf k}\theta_{\mathbf k} a_{\mathbf k}^2$. Together these are half of the
term $-\frac1{|Q|}\sum_{\mathbf k}\theta_{\mathbf k} a_{\mathbf k}^2$ of the identity, which is
otherwise dropped.
The $c$-crossing of the source (in the forcing lemma, and in step~(6$'$) of
Lemma~\ref{lem:EPclose} at
$\gamma\in[-3,-2]$) uses a further $\frac18$ of it. The $\Sigma_Q$-part produced by
that lemma is included in $C_{\mathbb E}\mathbb E^{1/2}$, by the remark after it,
so the four parts below are unchanged. With Lemmas~\ref{lem:FPd}
and~\ref{lem:FPdlow}, we obtain
\begin{multline}
\frac12\frac d{dt}e_Q+\frac12
\mathfrak d_Q\le\Big(C_{N_y}C(s)(\delta_y+\tfrac{C_1\rad}R)+C\rad^{-2}\eps^{1/2}+C_{\mathbb E}
\mathbb E^{1/2}\Big)\Sigma_Q\\
+C\big(\eps+\rad^{-2}\eps^{1/2}\big)\sum_{Q'\sim
Q}e_{Q'}+C_4\eps\,e_Q^{1/2}+C\eps^{1/2}e_Q+C\eps\,\mathbb E[\tilde G_\nu].
\end{multline}
The last two terms come from the $c$-crossing of the forcing lemma, and occur only in the frozen
problem (there the coefficient factor is $\tilde G_\nu$).

We choose $\rad_*$ and $\delta_*$ so that $C_{\rm v}C(s)C_1\rad_*\le1/(64N_0)$ and
$C_{N_y}C(s)\delta_*\le1/(64N_0)$. Then the coefficient of $\Sigma_Q$ from Lemma~\ref{lem:FPb},
written $C_{N_y}C(s)(\delta_y+\frac{C_1\rad}R)$ above, is bounded by $C_{N_y}C(s)\delta_y+C_{\rm
v}C(s)\frac{C_1\rad}R\le1/(32N_0)$. Here
\begin{equation}\label{eq:Cv}
C_{\rm v}:=2A_Fe^{L\Lambda}K_{\mathcal U}e^{C_{\mathcal U}}
\end{equation}
bounds, up to absolute factors included in $C(s)$, the $\ab$-part of that lemma with
$\alpha_\nu=0$. There $\alpha''=\alpha$, and the $y$-part of $\Omega$ reduces to the constant of the
first bound of Lemma~\ref{lem:FPexp} with $\nu_y=0$, which is $A_Fe^{L\Lambda}K_{
\mathcal U}e^{C_{\mathcal U}}$ by~\eqref{eq:rho0Cauchy} and Lemma~\ref{lem:Ufac}(iii) ($A_Fe^{L\Lambda}$ for
product data). The second bound of that lemma is needed for $\ab[f_0]\le2A$ alone.

The explicit $\cb$-part of the explicit commutators carries the factor $\eps$, and the explicit
$\bb$-part the factor $\eps^{1/2}$, with a constant depending also on $\delta_y$, $\rad_0$ and
$t_1$, all fixed before $\eps_0$. They are controlled by the condition on $\eps_0$ for the term
$C\rad^{-2}\eps^{1/2}$ below. Thus $C_{\rm v}$, $\rad_*$ and $\rad_0^{\rm lay}$ do not
depend on $t_1$. The constant $C_{N_y}$ of the part $\alpha_\nu\ne0$ depends on $t_1$.
Hence we choose $\delta_*$ and
$\delta_{y,\ell}$ after $t_1$, as in Appendix~\ref{app:const}, and they do not enter its choice.

Here $\mathbb E^{1/2}$ is the energy of the coefficient factor $h$. The terms $C_{\mathbb E}
\mathbb E^{1/2}$ and $C\rad^{-2}\eps^{1/2}$ of the coefficient of $\Sigma_Q$ are bounded as follows.
For the unknown itself, by the bootstrap hypothesis, we have $C_{\mathbb E}\mathbb E^{1/2}\le
C_{\mathbb E}\eta_*\le1/(32N_0)$ for $c_*$ small.
In the frozen problem the coefficient factor is $\tilde G_\nu$, with $\mathbb E[\tilde
G_\nu]^{1/2}\le2C_4C_{N_y}^{1/2}\eps t$, and we add the condition
$2C_{\mathbb E}C_4C_{N_y}^{1/2}\eps_0t_1\le1/(32N_0)$ on $\eps_0$.
Finally, we have $C\rad^{-2}\eps^{1/2}\le1/(32N_0)$ for $\eps_0$ small.

In the energy terms of the drift, $\mathbb E^{1/2}$ is
the energy of the coefficient factor, which in the frozen problem is $\mathbb E[\tilde
G_\nu]^{1/2}\le2C_4C_{N_y}^{1/2}\eps t$. Therefore these terms are bounded by
$CC_2C_4\eps\sum_{Q'\sim
Q}e_{Q'}$, and they are included in $C(\eps+\rad^{-2}\eps^{1/2})\sum_{Q'\sim Q}e_{Q'}$, whose
constant may depend on $C_{\mathbb E}$ and $C_4$, both fixed before $\eps_0$.

The crossing terms of Lemma~\ref{lem:FPcross}(iv) contribute a fourth part,
$CC_2(\eta_*+C_4(1+t_1)\eps)\le1/(32N_0)$, for $c_*$ and $\eps_0$ small.

Thus the coefficient of $\Sigma_Q$ is bounded by $\eta=4/(32N_0)$, and $\eta N_0\le\frac18$. After
the
doubling below this becomes $2\eta N_0\le\frac14$, the threshold of Lemma~\ref{lem:BKnb}(ii) with
$c=1$.

We integrate in time and multiply by $2$. We have $e_Q(0)=0$, because $G(0)=1$. By $\sum_{Q'\sim
Q}e_{Q'}\le N_0\mathbb E$ and $e_Q\le\mathbb E$, this is the hypothesis of the same lemma
with $c=1$, with $2\eta$ in place of $\eta$ (so $2\eta N_0\le\frac14=c/4$), with $2a\,\mathbb E$ in
place of $C\eps\,\mathbb E$, where $a:=C(\eps+\rad^{-2}\eps^{1/2})N_0$ also contains the Gronwall
term $C\eps^{1/2}e_Q$ of the $c$-crossing, and with $F=2C_4\eps\,
\mathbb E^{1/2}+2C\eps\,\mathbb E[\tilde G_\nu]$. By the hypothesis on the coefficient factor, $
\mathbb E[\tilde G_\nu](s)\le4C_4^2C_{N_y}\eps^2s^2$, so $C\eps\,\mathbb E[\tilde G_\nu](s)\le C_{\rm
s}\eps^3s^2$ with $C_{\rm s}:=4CC_4^2C_{N_y}$. Here $C$ is the constant of the $c$-crossing,
which does not depend on $C_4$. By the second bound of the neighbor-sum lemma, we
have
\begin{equation}
\mathbb E(t)\le\frac83\int_0^t\big(a\,\mathbb E+C_4\eps\,\mathbb E^{1/2}+C_{\rm s}\eps^3s^2\big)\dd
s ,
\end{equation}
and the first one implies the dissipation bound. Let $\mathcal R(t)$ be the right side, and
$S(t):=\frac89C_{\rm s}\eps^3t^3$ the part of it that comes from the last term. Then $\mathcal R\ge
S>0$ for $t>0$, and $\mathcal R'\le\frac83\big(a\mathcal R+C_4\eps\mathcal R^{1/2}+C_{\rm
s}\eps^3t^2\big)$. Since $(S^{1/2})'=\frac43C_{\rm s}\eps^3t^2/S^{1/2}\ge\frac43C_{\rm
s}\eps^3t^2/\mathcal R^{1/2}$, we get $(\mathcal R^{1/2})'\le\frac43a
\mathcal R^{1/2}+\frac43C_4\eps+(S^{1/2})'$. With $\mathcal R(0)=S(0)=0$, we obtain
\begin{equation}
\mathbb E^{1/2}\le\mathcal R^{1/2}\le e^{4at/3}\big(\tfrac43C_4\eps
t+S(t)^{1/2}\big)\le\tfrac98\big(\tfrac43+\tfrac49\big)C_4\eps t=2C_4\eps t,
\end{equation}
since $\frac43at_1\le\log\frac98$ and $(\frac89C_{\rm s})^{1/2}\eps_0^{1/2}t_1^{1/2}\le\frac49C_4$
for $\eps_0$ small. Both are conditions on $\eps_0$ (by $\rad\ge\rad_0/2$), see the list at
the end of the proof. By the choice of $\eps_0$, this does not exceed $\eta_*/2$. Therefore the
first half
of the bootstrap hypothesis is improved.

For the second half, by the level-$0$ energy and Lemma~\ref{lem:RBsob}, we have
\begin{equation}
\jap{v_0}^\sigma\sup_Q|G-1|\le C(\delta_y,\rad)\,\mathbb E^{1/2}\le C\eps t\le\tfrac14
\end{equation}
for $\eps_0$ small. This also implies the decay $|G-1|\le C\eps t\jap v^{-\sigma}$. Finally,
$\rad(t)\ge\rad_0/2$ on $[0,t_1]$.

Since the bootstrap hypothesis is assumed on all of $[0,T]$, this proves the proposition. In its
application, Proposition~\ref{prop:FPvisc}(iii), the hypothesis is removed by continuity: the
energies are continuous in $t$ by part~(ii) of that proposition, so the set of times up to which
it holds is closed. Moreover, it is open by the improvement just shown, and it contains $t=0$.

The conditions on $\eps_0$ used here and in the forcing lemma are the following.
\begin{itemize}
\item $C_\rad\eps_0t_1\le\rad_0/2$, so that $\rad\ge\rad_0/2$;
\item $(C_4+1)\eps_0t_1\le\eta_*/4$;
\item $2C_4C_{N_y}^{1/2}\eps_0t_1\le\eta_*$, so that the coefficient factor satisfies $
\mathbb E[\tilde G_\nu]^{1/2}\le\eta_*\le1$ (used in Lemmas~\ref{lem:FPb}, \ref{lem:FPdlow},
\ref{lem:FPcross}(iv) and~\ref{lem:FPccross}, in Proposition~\ref{prop:FPfix} and in the remark
after Lemma~\ref{lem:FPhalf});
\item $2C_{\mathbb E}C_4C_{N_y}^{1/2}\eps_0t_1\le1/(32N_0)$ (the frozen problem);
\item $4C\rad_0^{-2}\eps_0^{1/2}\le1/(32N_0)$;
\item $CC_2C_4(1+t_1)\eps_0\le1/(64N_0)$ together with $CC_2\eta_*\le1/(64N_0)$
(Lemma~\ref{lem:FPcross}(iv));
\item $C'C_4\eps_0t_1\le1$ (Lemma~\ref{lem:FPf});
\item $\frac43C(\eps_0+4\rad_0^{-2}\eps_0^{1/2})N_0t_1\le\log\frac98$;
\item $(\frac89C_{\rm s})^{1/2}\eps_0^{1/2}t_1^{1/2}\le\frac49C_4$;
\item $2C(\delta_y,\rad_0/2)C_4\eps_0t_1\le\frac14$.
\end{itemize}
Each involves only constants fixed before $\eps_0$. Since $C_{\rm s}=4CC_4^2C_{N_y}$, the source
condition reads $CC_{N_y}\eps_0t_1\le\frac1{18}$, and does not depend on $C_4$.
\end{proof}

We prove the finiteness of the energies with the full weight $\jap{v_0}^\sigma$, using a barrier,
in Proposition~\ref{prop:FPvisc}(ii). With a truncated weight, the bounds of Remark~\ref{rem:FPc}
would not be uniform in the truncation.

\subsubsection{The frozen problems, the vanishing-viscosity limit, and the
iteration}\label{sssec:FPvisc}

We solve the linear problems with frozen coefficients, prove that their energies are finite, and
pass to the limit in the viscosity and in the iteration. There are two difficulties. First, the
regularizing matrix $M_\nu$ must be analytic in $v$. A smooth cutoff such as
$\theta(|y|/2\vartheta(v))$ is not, and the commutators $[\partial_v^\beta,\nabla\cdot M_\nu\nabla]$
would not be controlled in the analytic weights. Second, the continuity argument needs the energies
of the regularized solution to be finite.

We use the standard theory of linear parabolic equations with H\"older continuous coefficients:
existence and regularity on bounded domains \cite[Ch.~IV, Thm.~5.2]{LSU}, and
the local Schauder estimates \cite[Ch.~IV, Thm.~10.1]{LSU}, \cite[Chs.~IV, V]{Lie96}. In
Proposition~\ref{prop:FPvisc}(ii), in Lemma~\ref{lem:CUtstar}(i)
and in Lemmas~\ref{lem:LWlin} and~\ref{lem:LWregex}, we also use the Schauder
estimates for coefficients H\"older continuous in the space variables only, interior \cite{Kn80} and
up to the initial time \cite{Lor00}, and the Krylov--Safonov H\"older estimate and Schauder's fixed
point theorem, as listed in (E4) of Inputs~\ref{hyp:E}. Moreover, \cite{Kn80} and the interior
Krylov--Safonov estimate \cite[Thm.~7.37, Cor.~7.41, (7.57)]{Lie96} are also used in the proof of
Proposition~\ref{prop:GS5}. No other parabolic result is quoted, apart from the comparison and
maximum principles, which are proved with growing barriers (in that proposition and, for
the viscous problem, Lemma~\ref{lem:LWlin}(ii)). The continuous dependence on the coefficients (the
stability of bounded solutions in the proof of Lemma~\ref{lem:LWregex}(i), and the uniform existence
time in the continuation of part~(iv) of the same lemma) follows from the Krylov--Safonov estimate,
the interior estimate \cite{Kn80} and uniqueness.

\emph{Smooth weights.} Put
\begin{equation}
\vartheta_\natural(v):=(1+\log\jap v)^{1/2},\qquad h_\natural(y,v):=\big(\jap y^2+\log\jap
v\big)^{1/2}.
\end{equation}
Both are comparable to $\vartheta$ and $h_r$ with absolute constants. Both are real-analytic in $v$
on the strip $\mathbb S_{1/2}$. Indeed, there $\mathrm{Re}\,(1+v\cdot
v)\ge\frac12(1+|\mathrm{Re}\,v|^2)$, so $\log\jap v$ is holomorphic with $|\mathrm{Im}\log\jap
v|\le\pi/2$, and the square roots are taken of quantities with real part $\ge\frac12$. By the Cauchy
estimates, for all $\beta$ and $m\in\{-2,1\}$, we have
\begin{equation}
|\partial_v^\beta(\vartheta_\natural^m)|\le C^{|\beta|}\beta!\,\vartheta_\natural^m,\qquad|\partial_y^\alpha\partial_v^\beta(h_\natural^m)|\le
C_\alpha C^{|\beta|}\beta!\,\jap y^{-|\alpha|}h_\natural^m .
\end{equation}

\begin{defin}[The regularized frozen problem]\label{def:FPreg}
Let $\tilde G\in C([0,t_1];C^1_{\rm loc}(\R^6))$ be given, with $\tilde G(0)=1$ (which also follows
from the energy bound at $t=0$) and
\begin{equation}
\mathbb E[\tilde G](t)^{1/2}\le2C_4\eps t\quad\text{for all }t\in[0,t_1],\qquad\tfrac12\le\tilde
G\le2 .
\end{equation}
Since $2C_4\eps s$ is nondecreasing in $s$, this implies $\mathbb E[\tilde
G_\nu](t)^{1/2}\le2C_4C_{N_y}^{1/2}\eps t$, by Lemma~\ref{lem:FPMnu}(v)
(Remark~\ref{rem:FPpointwise}). For $\nu\in(0,\eps]$ (the bound $\nu\le\eps$ enters the barrier of
Proposition~\ref{prop:FPvisc}(ii) through $C(\eps+\nu)\le1$, and part~(iii)
and Remark~\ref{rem:FPtrunc8} through $C\nu\,\mathbb E\le C\eps\,\mathbb E$):
\begin{itemize}
\item extend $\tilde G$ by $\tilde G(s):=1$ for $s<0$. Let $\varrho_\nu$ be a standard mollifier in
$y$ supported in $|y|\le\nu$, and $\upsilon_\nu\ge0$ a smooth function on $\R$ with
$\int\upsilon_\nu=1$, supported in $[\nu/2,\nu]$, and put
\begin{equation}
\tilde G_\nu(t,y,v):=\int_\R\int_{\R^3}\upsilon_\nu(s)\,\varrho_\nu(\eta)\,\tilde
G(t-s,y-\eta,v)\dd\eta\dd s .
\end{equation}
The average in $t$ is over the past: $\tilde G_\nu$ on $[0,T]$ depends only on $\tilde G$ on
$[0,T]$, and $\tilde G_\nu=1$ for $t\le\nu/2$;
\item let $\tilde A_\nu,\tilde V_\nu,\tilde S_{0,\nu}$ be the coefficients built from $f_0\tilde
G_\nu$;
\item put
\begin{equation}
M_\nu:=\nu\Big[\vartheta_\natural^{-2}I_y-\big(\vartheta_\natural^{-2}-h_\natural^{-2}\big)\frac{y\otimes
y}{\jap y^2}\Big]\oplus\nu\,I_v ;
\end{equation}
\item let $\chi_k(y,v)=\chi(|(y,v)|/k)$ be a cutoff, with $\chi=1$ on $[0,1]$ and $\chi=0$ on
$[2,\infty)$.
\end{itemize}
The \emph{regularized problem} on $B_{3k}\subset\R^6$ is
\begin{equation}\label{eq:FPreg}
\partial_tG=\mathcal L_\nu G:=D\cdot(\tilde A_\nu DG)+\nabla\cdot(M_\nu\nabla G)+\tilde V_\nu\cdot
DG+\chi_k\tilde S_{0,\nu}G,\qquad G|_{t=0}=1,\quad G|_{\partial B_{3k}}=1 .
\end{equation}
The \emph{whole-space problem} is the same equation without $\chi_k$ and without the boundary
condition.
\end{defin}

\begin{lem}[Properties of $M_\nu$ and of the mollified coefficients]\label{lem:FPMnu}
\leavevmode
\begin{enumerate}[label=(\roman*)]
\item The $y$-block of $M_\nu$ has eigenvalue $\nu\big(h_\natural^{-2}|y|^2+\vartheta_\natural^{-2}\big)/\jap
y^2\simeq\nu h_\natural^{-2}$ in the direction $\hat y$, and $\nu\vartheta_\natural^{-2}$ on $\hat
y^\perp$.
\item Let $S_Q$ be the linear map, in the frame of $Q$, with
$S_Q=c_{\rm box}\operatorname{diag}(h_r(y_Q)^{-1},\vartheta^{-1},\vartheta^{-1})$ (for an outer $Q$
it sends the unit cube to $Q$). Then $S_Q^{-1}M_\nu S_Q^{-T}\ge c\nu I$ on $\tilde Q$, with $c$
independent of $Q$.
\item As quadratic forms,
\begin{equation}\label{eq:FPMnu-der}
|\partial_Q^\alpha\partial_v^\beta M_\nu|\le C_{|\alpha|}C^{|\beta|}\beta!\,h_r^{a_r}\vartheta^{a_\perp}M_\nu
\end{equation}
on $\tilde Q$.
\item $|M_\nu^{1/2}\nabla\varphi_Q|\le C\nu^{1/2}$.
\item
\begin{equation}\label{eq:FPMnu-E}
\mathbb E[\tilde G_\nu](t)\le C_{N_y}\sup_{0\le s\le t}\mathbb E[\tilde G](s)
\end{equation}
for every $t\in[0,t_1]$, and $\frac12\le\tilde G_\nu\le2$. The function $\tilde G_\nu$ is $C^\infty$
in $(t,y)$, with $|\partial_t^j\partial_y^\alpha\tilde G_\nu|\le C_{j,\alpha}\nu^{-j-|\alpha|}$ for
$j+|\alpha|\ge1$. These derivatives commute with $\partial_v^\beta$, and
$\sup|\partial_t^j\partial_y^\alpha\partial_v^\beta\tilde G_\nu(t)|\le
C_{j,\alpha}\nu^{-j-|\alpha|}\sup_{0\le s\le t}\sup|\partial_v^\beta\tilde G(s)|$, which is finite
by Lemma~\ref{lem:RBsob}. It is analytic in $v$ with the radius of $\tilde G$ at time $t$. As
$\nu\to0$, $\tilde G_\nu\to\tilde G$ in $C([0,t_1];C^1_{\rm loc}(\R^6))$.
\end{enumerate}
\end{lem}

\begin{proof}
(i) This follows by a direct computation, using $h_\natural\ge\jap y$ and
$h_\natural^2\le2\jap y^2\vartheta_\natural^2$.

(ii) On an inner box, $|y|\le4\vartheta$ on $\tilde Q$, so $\vartheta\le h_r\le5\vartheta$ and
$M_\nu\simeq\nu\vartheta^{-2}I$. On an outer box, let $\theta$ be the angle between $\hat y$ and
$\hat e_Q$. By Lemma~\ref{lem:BKcov}(iv), we have $\theta\le C/\jap y$ on $\tilde Q$. Conjugating by
$\mathrm{diag}(h_r,\vartheta,\vartheta)$, we see that the matrix $h_r^{-2}\hat e\otimes\hat
e+\vartheta^{-2}(I-\hat e\otimes\hat e)$ has norm bounded by
\begin{equation}
2+\theta^2h_r^2/\vartheta^2\le2+C\Big(\frac{\jap y+\vartheta}{\vartheta\jap y}\Big)^2\le C
\end{equation}
relative to $h_r^{-2}\hat y\otimes\hat y+\vartheta^{-2}(I-\hat y\otimes\hat y)$. Hence the two forms
are comparable, and~(ii) follows from~(i).

(iii) The matrix $y\otimes y/\jap y^2$ has $y$-derivatives bounded by $C_\alpha\jap y^{-|\alpha|}\le
C_\alpha h_r^{a_r}\vartheta^{a_\perp}$, and its derivative is a mixed form $\hat y\otimes
e+e\otimes\hat y$ of size $\jap y^{-1}$. Since
\begin{equation}
\vartheta^{-2}\jap y^{-1}|\xi_r||\xi_\perp|\le C\big(h_r^{-2}\xi_r^2+\vartheta^{-2}\xi_\perp^2\big)\cdot\frac{h_r}{\vartheta\jap
y}
\end{equation}
and $h_r\le C\vartheta\jap y$, the derivatives are bounded by $M_\nu$ itself. The $v$-derivatives
come from the Cauchy estimates above.

(iv) By~(i), the radial derivative of $\varphi_Q$, of size $h_r$, meets $\nu^{1/2}h_r^{-1}$. The
transverse derivative, of size $\vartheta$, meets $\nu^{1/2}\vartheta^{-1}$. The $v$-derivative is
$O(1)$.

(v) For $|\eta|\le\nu\le1$, the translate $\tilde Q+\eta$ is covered by the boxes $2Q'$ of at most
$C$ boxes $Q'$ with $|y_{Q'}-y_Q|\le2$ and $v_0(Q')=v_0(Q)$. Since $\varphi_Q$ and $\varphi_{Q'}$
have the same $v$-factor $\chi(v-v_0)$, by Lemma~\ref{lem:BKframe}(i), we have
$|\varphi_Q(\cdot+\eta,\cdot)|^2\le\sum_{Q'}\varphi_{Q'}^2$ on its support. If $Q$ and $Q'$ are
outer, their frames differ by angles $\le C/\jap y$, and the computation of the proof of~(ii) of
that lemma, which rests on $|y_Q-y_{Q'}|\le2$ and $\theta h_r/\vartheta\le C(\jap
y+\vartheta)/(\vartheta\jap y)$ (and its mixed case if one of them is inner), converts the frame
derivatives with a factor $C^{N_y}$. Then, by Minkowski's inequality in $\eta$, we obtain
\begin{equation}\label{eq:FP-moll}
e_Q[\tilde G_\nu]\le C_{N_y}\sup\big\{e_{Q'}[\tilde G]:|y_{Q'}-y_Q|\le2,\ v_0(Q')=v_0(Q)\big\}.
\end{equation}
The two-sided bound is preserved by averaging, and the derivative bounds are those of the mollifier.
Mollification in $y$ commutes with complex shifts in $v$, and so does the average in $t$.

For the average in $t$, we fix $t$. A derivative $\partial_t^j\partial_y^\alpha$ with
$j+|\alpha|\ge1$ falls on $\upsilon_\nu\otimes\varrho_\nu$, whose derivatives have integral $0$.
Thus it acts on
$\tilde G-1$, with $|\tilde G-1|\le1$, and on $\partial_v^\beta(\tilde G-1)$ in the same way. The
average also preserves $\frac12\le\tilde G_\nu\le2$. Moreover, the radius of $\tilde G(s)$ is
nonincreasing in
$s$, since it comes from the weights $w_\beta(\rad(s))$, and $\tilde G(s)=1$ for $s<0$. For the
energy, the weights are nonincreasing in $t$: through $\rad(t)$, and through $\omega_Q(t)$ in
Definition~\ref{def:EPen}. For fixed weights, $e_Q^{1/2}$ is a norm of the family
$(u_{\mathbf k})_{\mathbf k}$, which depends linearly on $\tilde G-1$, and $\tilde G(s)-1=0$ for
$s<0$.
By Minkowski's inequality in $s$, \eqref{eq:FP-moll} at each time $s\in[t-\nu,t]$ with the weights
of time $t$, and the monotonicity of the weights, we get $e_Q[\tilde G_\nu](t)^{1/2}\le
C_{N_y}^{1/2}\sup_{0\le s\le t}\mathbb E[\tilde G](s)^{1/2}$, where $\mathbb E[\tilde G](s)$ is
taken with the weights of time $s$. The same holds for $\mathbb E^\flat$. Finally, since $\tilde
G(0)=1$, the extension of $\tilde G$ by $1$ lies in $C([-1,t_1];C^1_{\rm loc})$. Therefore $\tilde
G$ and its first derivatives in $(y,v)$ are uniformly continuous on $[-1,t_1]\times K$ for every
compact $K$, which implies the convergence as $\nu\to0$.
\end{proof}

By Lemma~\ref{lem:FPMnu}(v), $\tilde G_\nu$ satisfies the hypothesis of
Proposition~\ref{prop:FPclose} on the coefficient factor, with $2C_4$ replaced by
$C_{N_y}^{1/2}2C_4$, which affects only $\eps_0$. In the rest of this subsection, constants $C(\nu)$
may depend on $\nu$ and $k$, while $C$ does not.

\begin{prop}[The frozen problem]\label{prop:FPvisc} 
Let $\tilde G$ be as in Definition~\ref{def:FPreg}, and $\eps\le\eps_0$.
\begin{enumerate}[label=(\roman*)]
\item \emph{(Bounded domains.)} Problem~\eqref{eq:FPreg} has a unique solution $G_{\nu,k}\in
C^\infty((0,t_1]\times B_{3k})\cap C([0,t_1]\times\overline B_{3k})$, which lies in
$H^{l+2,l/2+1}([0,t_1]\times\overline B_{3k})$ for every nonintegral $l\in(0,1)$, with $e^{-C\eps
t}\le G_{\nu,k}\le e^{C\eps t}$.
\item \emph{(Whole space, $\nu>0$.)} Along a subsequence $k\to\infty$, $G_{\nu,k}\to G_\nu$ in
$C^\infty_{\rm loc}((0,t_1]\times\R^6)$ and in $C_{\rm loc}([0,t_1]\times\R^6)$. The limit $G_\nu$
is the unique bounded classical solution of the whole-space problem. It satisfies
\begin{equation}
e^{-C\eps t}\le G_\nu\le e^{C\eps t},\qquad|G_\nu-1|\le K_0e^{t_1}\eps t\jap v^{-\sigma}.
\end{equation}
For every $n$, the energy $\mathbb E^{(n)}[G_\nu]$ truncated to $|\beta|\le n$ is finite and
continuous on $[0,t_1]$.
\item \emph{(Uniform estimates.)} The conclusions of Proposition~\ref{prop:FPclose} hold for
$G_\nu$, with constants independent of $\nu$ and $n$. In particular $
\mathbb E[G_\nu]^{1/2}\le2C_4\eps t$.
\item \emph{(The limit $\nu\to0$.)} Along a subsequence, $G_\nu\to G$ in $C([0,t_1];C^2_{\rm
loc}(\R^6))$. The limit $G$ is a classical solution of the frozen problem
\begin{equation}
\partial_tG=D\cdot(\tilde ADG)+\tilde V\cdot DG+\tilde S_0G,\qquad G(0)=1,
\end{equation}
and it satisfies $\mathbb E[G]^{1/2}\le2C_4\eps t$, the dissipation bound, and $|G-1|\le C\eps t\jap
v^{-\sigma}$.
\item \emph{(Uniqueness.)} $G$ is the only classical solution of the frozen problem such that both
$\sup|G|$ and $\sup_{t,Q}e^{(0)}_Q[G]$ are finite. Here
$e^{(0)}_Q:=|Q|^{-1}\|(G-1)\varphi_Q\|^2_{L^2}$ is the level-$0$ energy, i.e.\ the term $
\mathbf k=0$ of $e_Q$, without the weight $\jap{v_0}^\sigma$. Since
$\operatorname{supp}\varphi_Q\subset\tilde
Q\times B(v_0,2)$ and $|\tilde Q|=64|Q|$, it is bounded by $C(1+\sup|G|)^2$.
\end{enumerate}
\end{prop}

\begin{proof}
\emph{(i)} We write $\mathcal L_\nu$ in nondivergence form, using $D\cdot\ab[h]=\bb[h]$. Its
coefficients are built from $M_\nu$, which is smooth and does not depend on $t$; from $\zeta=D\log
f_0$ and $D\zeta$, which are smooth and polynomial in $t$; from $\chi_k$; and from integrals $\int
k(w)\,h(t,y+tw,v-w)\dd w$ with $k\in\{a,b,c\}$ and $h=f_0\tilde G_\nu$ (for $\gamma=-3$,
$\cb[h]=-8\pi h$ is local). On $[0,t_1]\times\overline B_{3k}$ they have the following regularity.
\begin{itemize}
\item They are $C^\infty$ in $(y,v)$, with bounds independent of $t$, by Lemma~\ref{lem:FPMnu}(v)
and Lemma~\ref{lem:FPexp} (stated for $|\nu_y|\le N_y$, but under (D$_{\rm iso}$) its proof yields
the same bound for every order, with constants depending on the order), since $\partial_y$ and
$\partial_v$ pass to $h$ at the point $(y+tw,v-w)$.
\item They are H\"older continuous in $t$ with every exponent $\theta\in(0,1)$, together with their
derivatives in $(y,v)$. Indeed, by the same two lemmas and (D1), $|h|$
and $|\partial_th|+|\nabla_yh|$, and the same for $\partial^{\mathbf k}_{y,v}h$, are bounded by
$C(\nu)\rho_0(v-w)$ at the point $(t,y+tw,v-w)$. Thus the increment of $h(t,y+tw,v-w)$ between $t$
and $t'$ is bounded by $C(\nu)\rho_0(v-w)\min\{1,(1+|w|)|t-t'|\}\le
C(\nu)\rho_0(v-w)(1+|w|)^\theta|t-t'|^\theta$, and $\int|k(w)|(1+|w|)^\theta\rho_0(v-w)\dd w\le
C(k)$ for $|v|\le3k$, since $\rho_0\le A_F\jap\cdot^{-p}$ with $p\ge4>\gamma+5+\theta$.
\item For every $t_0\in(0,t_1)$ they are $C^\infty$ in $(t,y,v)$ on $[t_0,t_1]\times\overline
B_{3k}$, with bounds depending on $t_0$, $\nu$ and $k$. Each derivative in $t$ of $h(t,y+tw,v-w)$
produces $\partial_t+w\cdot\nabla_y$. The derivatives of $\tilde G_\nu$ are bounded as in the first
item, and $|w|\le(|y+tw|+3k)/t_0$, where the powers of $|y+tw|$ are controlled by
the factor $F_0(y+tw)$ of $f_0$.
\end{itemize}
The last two items use the mollification in $t$. Without it, the coefficients would be only as
regular in $t$ as $\tilde G$, which for $\tilde G=G^j$, $j\ge1$, is known to be $C^1$ only, by~(iv).
Up to $t=0$ the dependence on $t$ through $y+tw$ limits the regularity in any case: the $j$-th
derivative in $t$ of the integrand carries the factor $|w|^j$, and for even $j\ge p-\gamma-5$ the
integral of its absolute value can diverge at $t=0$ (only the upper bound $\rho_0\le
A_F\jap\cdot^{-p}$ is assumed). By the second item, the coefficients lie in
$H^{l,l/2}([0,t_1]\times\overline B_{3k})$ for $l\in(0,1)$, and in
$H^{l,l/2}([t_0,t_1]\times\overline B_{3k})$ for every $l>0$ and $t_0>0$. The operator is uniformly
parabolic, since $D^T\tilde A_\nu D\ge0$ and $M_\nu\ge c(k)\nu I$ there. The data are compatible at
$\{t=0\}\times\partial B_{3k}$: $G\equiv1$ there, and $\mathcal L_\nu1=\chi_k\tilde S_{0,\nu}=0$
near $\partial B_{3k}$. Hence, by \cite[Ch.~IV, Thm.~5.2]{LSU} with $l\in(0,1)$ (compatibility of
order~$1$), there is a unique solution
$G_{\nu,k}\in H^{l+2,l/2+1}([0,t_1]\times\overline B_{3k})$. It
is also the only solution in $C^{2,1}((0,t_1]\times B_{3k})\cap C([0,t_1]\times\overline B_{3k})$,
by the weak maximum principle applied to the difference of two solutions, whose coefficients are
bounded on $[0,t_1]\times\overline B_{3k}$. In particular $G_{\nu,k}$ is continuous on
$[0,t_1]\times\overline B_{3k}$, and $D^2_{y,v}G_{\nu,k}$ and $\partial_tG_{\nu,k}$ are continuous
up to $t=0$.

For $t>0$ the solution is smooth. Let $t_0\in(0,t_1)$, and let $\Upsilon$ be a smooth cutoff in
$(t,y,v)$ that vanishes for $t\le t_0/2$ and near $\partial B_{3k}$. Then $\Upsilon\,(G_{\nu,k}-1)$
solves on $[t_0/2,t_1]\times B_{3k}$ the equation of $\mathcal L_\nu$ with zero initial and boundary
values, and with a free term built from $\Upsilon$ and its derivatives, from $\chi_k\tilde
S_{0,\nu}$, and from $G_{\nu,k}$ and its first derivatives. If $G_{\nu,k}\in H^{l'+2,l'/2+1}$ on the
support of $\Upsilon$, the free term is in $H^{l'+1,(l'+1)/2}$. Then, by \cite[Ch.~IV, Thm.~5.2]{LSU}
with $l'+1$ in place of $l$, together with the uniqueness, we get $\Upsilon\,(G_{\nu,k}-1)\in
H^{l'+3,(l'+3)/2}$. Here the compatibility conditions of all orders hold, since the data vanish
near $t=t_0/2$. By induction over a nested sequence of cutoffs, we obtain
$G_{\nu,k}\in C^\infty((0,t_1]\times B_{3k})$.

The two-sided bound follows from the maximum principle, with $e^{\pm C\eps t}$ as super- and
subsolutions. Here we use $|\chi_k\tilde S_{0,\nu}|\le C\eps$, which follows from~(i) of the decay
lemma, and
the boundary values $1\in[e^{-C\eps t},e^{C\eps t}]$.

\emph{(ii), the decay barrier.} Let $u=G_{\nu,k}-1$. Then $\partial_tu=\mathcal L_\nu u+\chi_k\tilde
S_{0,\nu}$, with $u=0$ on the parabolic boundary. We put $\psi:=K_0\eps te^{t}\jap v^{-\sigma}$. We
have
\begin{equation}
|\partial_v^j\jap v^{-\sigma}|\le C\jap v^{-\sigma-j},\qquad|\tilde A_\nu|+|D\cdot\tilde
A_\nu|+|\tilde V_\nu|\le C\eps,\qquad|M_\nu|+|\nabla M_\nu|\le C\nu .
\end{equation}
Since $\psi$ does not depend on $y$, we have $D\psi=\partial_v\psi$. Hence $|\mathcal L_\nu\psi|\le
C\eps\psi+C\nu\psi$, and
\begin{equation}
\partial_t\psi-\mathcal L_\nu\psi-|\chi_k\tilde S_{0,\nu}|\ge K_0\eps e^t\jap
v^{-\sigma}+(1-C\eps-C\nu)\psi-C\eps\jap v^{-2\sigma}\ge0
\end{equation}
for $K_0\ge C$ and $C(\eps+\nu)\le1$. Here we used the bound $|\tilde S_{0,\nu}|\le C\eps\jap
v^{-2\sigma}$, as in~(i). By the maximum principle, we have $|u|\le\psi$.

\emph{(ii), local estimates uniform over boxes.} We fix a box $Q$, and rescale $\tilde
Q\times[t-\tau,t]$ to a unit cylinder by $S_Q$, for some $\tau\in(0,1]$. In the rescaled variables,
$\mathcal L_\nu$ is uniformly parabolic with ellipticity $\ge c\nu$, by Lemma~\ref{lem:FPMnu}(ii).
The local estimates below use the coefficients on $2\tilde Q$. Parts~(ii)--(iii) of that lemma,
Lemma~\ref{lem:FPexp} and Lemma~\ref{lem:FPdecay}(iii) are stated on $\tilde Q$ or $Q$, but their
proofs rest
on $\operatorname{diam}\le2$, $h_r\simeq h_r(y_Q)$ and, on outer boxes, an angle $\le
C/(\vartheta\jap y)$ between $\hat y$ and $\hat e_Q$, and these hold on $2\tilde Q$ as well, with
larger constants. The rescaled coefficients have $C^m$ norms in the space variables bounded by
$C(m,\nu)$,
uniformly in $t\in[0,t_1]$ and independently of $Q$ and $k$, for the boxes with $2\tilde Q\times
B(v_0,4)\subset B_{3k}$. Below, every supremum over $Q$ of a quantity of $G_{\nu,k}$ is over these
boxes, and every box is among them for $k$ large. The coefficient bounds rest on four facts:
\begin{itemize}
\item a rescaled radial derivative is $c_{\rm box}h_r^{-1}\partial_r$, and each radial derivative of
the coefficients produces at most a factor $Ch_r$, by the explicit-factor bounds
and~(iii), (v) of Lemma~\ref{lem:FPMnu};
similarly in the transverse directions;
\item in $D=\partial_v-t\nabla_y$, the rescaled $y$-part is $t(h_r\hat e+\vartheta
e)\cdot\nabla_{\rm resc}/c_{\rm box}$;
\item the second-order coefficients are $t^2h_r^2\tilde A\hat e\cdot\hat e\le Ct^2\eps h_r^2\jap
y^{-2}\le C$ in the outer region (null structure, as in Lemma~\ref{lem:BKtop}), and
$t^2\vartheta^2|\tilde A|\le C\eps\vartheta^2\jap v^{-2\sigma}\le C$;
\item the drift satisfies $|\tilde V\cdot\hat e|th_r\le C$, by $|\ab[f_0]\zeta^\sharp\cdot\hat
e|\le(A\hat e\cdot\hat e)^{1/2}(\ab[f_0]\zeta^\sharp\cdot\zeta^\sharp)^{1/2}\le C\eps\jap y^{-1}$
and Lemma~\ref{lem:RBgain}(ii). Transversally, $t\vartheta|\tilde V\cdot e|\le
Ct\vartheta\big(\eps\jap v^{-\sigma}+\eps\jap v^{\gamma+1}\big)\le C$, by
$|\ab[f_0]\zeta^\sharp\cdot e|\le(Ae\cdot e)^{1/2}(\ab[f_0]\zeta^\sharp\cdot\zeta^\sharp)^{1/2}$ and
part~(i) of the decay lemma.
\end{itemize}
These bounds are in the space variables only. In $t$ the coefficients are only H\"older continuous
near $t=0$ (proof of~(i)), and their H\"older seminorms in $t$ need not be bounded uniformly in $Q$.
In fact, at $\gamma=-2$ they are not, as in Lemma~\ref{lem:LWlin}(i).

We therefore use the Schauder estimates for coefficients H\"older continuous in the space variables,
uniformly in $t$, in the local form of Appendix~\ref{app:visc}, with a first-order term, which
\cite[(8)]{Kn80} and \cite[(H1)]{Lor00} allow. For $t\le1$ we use the estimate up to the initial
time \cite[Thm.~4.1]{Lor00} on $[0,t]$, with $u(0)=0$. It applies, since $D^2_{y,v}u$ and
$\partial_tu$ are continuous up to $t=0$ by~(i). For $t>1$ we use the interior estimate
\cite[Thm.~1]{Kn80} on $[t-1,t]$.

In the cutoff of the local form, the rescaled coefficients may be changed outside the doubled box,
where the cutoff solution vanishes, into coefficients on $\R^6$ with the same bounds (by a convex
combination with $c\nu I$).

The derivatives of higher order follow by induction from the estimates for spatial difference
quotients, as in the proof of Lemma~\ref{lem:LWregex}(i). Since $u(0)=0$, the induction runs up to
$t=0$ at every order. At each order the difference quotients are solutions in the sense of
\cite[Def.~1.1]{Lor00} on $[0,t]$. Indeed, $u$ is smooth for $t>0$ by~(i), the previous step and the
equation bound the next derivatives and their time derivatives up to $t=0$, and interpolation
implies their continuity at $t=0$, as in that appendix. With $\nabla$ the gradient in the rescaled
space variables, we obtain
\begin{equation}
\sup_{\text{unit cylinder}}\Big(\sum_{j\le m}|\nabla^ju|+\sum_{j\le
m-2}|\nabla^j\partial_tu|\Big)\le C(m,\nu)\Big(\sup_{2\tilde Q\times[0,t]}|u|+\sup_{s\le
t}\|\chi_k\tilde S_{0,\nu}(s)\|_{C^{m}}\Big)\le C(m,\nu)\,\eps\jap{v_0}^{-\sigma}.
\end{equation}
In the original variables, this implies that
$h_r^{-a_r}\vartheta^{-a_\perp}|\partial_Q^\alpha\partial_v^\beta u|$ for $|\alpha|+|\beta|\le m$
is bounded by $C(m,\nu)\eps\jap{v_0}^{-\sigma}$. These are the weights of $e_Q$, and the top weight
is $\le1$. Hence
every truncated energy $\mathbb E^{(n)}[G_{\nu,k}]$ is bounded, uniformly in $k$, by
$C(n,\nu)\eps^2$. The rescaled second-order coefficients, of $\tilde A_\nu$ and of $M_\nu$, are
bounded by $C(m,\nu)$, as shown above, and the same local bounds, with one more derivative, imply
$\sup_Q\big(\mathfrak d^{(n)}_Q+\mathfrak d^{\nu,(n)}_Q\big)\le C(n,\nu)\eps^2$ on $[0,t_1]$,
uniformly in $k$, and this passes to the limit $k\to\infty$. Consequently
$\sup_Q\big(e^{(n)}_Q(t)+\int_0^t(\mathfrak d^{(n)}_Q+\mathfrak d^{\nu,(n)}_Q)\big)$ is finite,
which is the finiteness assumption of Lemma~\ref{lem:BKnb}(ii) for the truncated energies.

\emph{(ii), the limit $k\to\infty$.} The local bounds are uniform in $k$. On compact subsets of
$(0,t_1]\times\R^6$ the coefficients are smooth in $(t,y,v)$, with bounds that do not depend on $k$
(proof of~(i)), so the local estimates \cite[Ch.~IV, Thm.~10.1]{LSU} bound all derivatives in
$(t,y,v)$ there, uniformly in $k$. By the Arzel\`a--Ascoli theorem and a diagonal argument, a
subsequence converges in $C^\infty_{\rm loc}((0,t_1]\times\R^6)$, and, by the bounds up to $t=0$, in
$C_{\rm loc}([0,t_1]\times\R^6)$, to a solution $G_\nu$ of the whole-space problem. Moreover, for
every $j$, $\nabla^j_{y,v}(G_{\nu,k}-1)$ vanishes at $t=0$ and its time derivative is bounded on
boxes up to $t=0$, uniformly in $k$, by the local bounds. Hence we have
$|\nabla^j_{y,v}(G_{\nu,k}-1)(t)|\le
C(j,\nu)\,t\jap z^j$, $z=(y,v)$, the derivatives converge in $C_{\rm loc}([0,t_1]\times\R^6)$, and
$\nabla^j_{y,v}G_\nu$ is continuous on $[0,t_1]\times\R^6$, with $\nabla^j_{y,v}(G_\nu-1)(0)=0$ and
$|\nabla^j_{y,v}(G_\nu-1)|\le C(j,\nu)\jap z^j$ (used in Proposition~\ref{prop:Bn}). The two bounds
of~(i) and~(ii) pass to the limit, and so does the finiteness of $\mathbb E^{(n)}$. Continuity in
$t$ of $\mathbb E^{(n)}[G_\nu]$ follows from the local $C^{m}$ bounds with $m=N_y+n+2$, uniformly
over boxes. Indeed, $\partial_t\partial^{\mathbf k} G_\nu$ is bounded in the weights, so each box
energy is Lipschitz in $t$ uniformly in $Q$, and so is their supremum.

\emph{(ii), uniqueness of bounded solutions for $\nu>0$.} Let $\delta$ be the difference of two
bounded solutions, with $|\delta|\le M$. We put $\Psi_\eta:=\eta(1+|y|^2+|v|^2)e^{Kt}$. By the
coefficient bounds above (the coefficients and their first derivatives are bounded, and
$|D|y|^2|\le2t|y|$), we have
\begin{equation}
\mathcal L_\nu(1+|y|^2+|v|^2)\le C(1+|y|^2+|v|^2).
\end{equation}
Hence $\Psi_\eta$ is a supersolution for $K=C+1$. On $\partial B_\varrho$ with $\eta\varrho^2\ge M$
we have $\delta\le\Psi_\eta$, and at $t=0$ we have $\delta=0$. By the maximum principle on
$B_\varrho$, we have $\delta\le\Psi_\eta$. Letting $\varrho\to\infty$ and then $\eta\to0$, we get
$\delta\le0$, and by symmetry $\delta=0$.

\emph{(iii)} We apply Proposition~\ref{prop:FPclose} to $G_\nu$. The bootstrap is run at the level of
$G_\nu$, where the energies are finite and continuous, and the dissipations finite, by~(ii). The
viscous term adds three contributions to~\eqref{eq:FPid}:
\begin{itemize}
\item the principal part, $-\iint M_\nu\nabla u_{\mathbf k}\cdot\nabla u_{\mathbf k}\varphi_Q^2$,
which is a new dissipation $\mathfrak d^\nu_Q$;
\item a cutoff term, bounded by $\frac14\mathfrak d^\nu_Q+C\nu\sum_{Q'\sim Q}e_{Q'}$, by
Lemma~\ref{lem:FPMnu}(iv);
\item the commutators $\iint\nabla\cdot(\partial^\mu M_\nu\nabla
u_{\mathbf k-\mu})u_{\mathbf k}\varphi_Q^2$.
\end{itemize}
We integrate by parts once, and apply the Cauchy--Schwarz inequality for the form bound
$|\partial^\mu M_\nu|\le C_\mu M_\nu$ of part~(iii) of the same lemma. Then the commutators are
bounded as in
Lemma~\ref{lem:FPb}: by
\begin{equation}
C(s)\big(C'_{N_y}\delta_y+C_M\rad\big)\Big(\sum_{Q'\sim Q}\mathfrak d^\nu_{Q'}+C\nu\sum_{Q'\sim
Q}e_{Q'}\Big).
\end{equation}
The $y$-derivatives of $M_\nu$ produce $h_r^{a_r}\vartheta^{a_\perp}$, which the weights absorb. The
$v$-derivatives are analytic with a radius independent of $R$. Here the $v$-analyticity of $M_\nu$
enters, and we take $\rad_*R\le1/(2C)$, so that $C\rad\le\frac12$ for
$\rad\le\rad_0\le\rad_*R$.

Here $C'_{N_y}$ and $C_M$ come from this form bound and from the $v$-sum of
Lemma~\ref{lem:ULconv}, with the radius of $M_\nu$ in place of $R$. The part with $\alpha_\nu=0$
involves only the $v$-derivatives of $M_\nu$. Both constants are structural, and neither depends on
$t_1$ or on the constant $C_{N_y}$ of the explicit commutators. Therefore the choice of $\rad_*$ in
that proposition does not cover these terms, and we add the conditions
$C(s)C_MR\rad_*\le1/(64N_0)$, which does not depend on $t_1$, and
$C(s)C'_{N_y}\delta_*\le1/(64N_0)$. Since $\rad\le\rad_0\le\rad_*R$, these terms are then
bounded by $\frac1{32N_0}\sum_{Q'\sim Q}\mathfrak d^\nu_{Q'}+C\nu\sum_{Q'\sim Q}e_{Q'}$.

With the principal term $-\mathfrak d^\nu_Q$ and the cutoff term $\le\frac14\mathfrak d^\nu_Q$, the
viscous dissipation enters the box inequality with the coefficient $c_\nu=\frac34$, and the neighbor
viscous dissipations with a coefficient $\eta_\nu$. Together with Remark~\ref{rem:FPtrunc8}, we
have $\eta_\nu N_0\le\frac18\le c_\nu/4$.

After the multiplication by $2$ in the closing argument,
we apply part~(ii) of the neighbor-sum lemma
to $\mathfrak d+\mathfrak d^\nu$, with $c=\min(1,2c_\nu)=1$ and
$2\max(\eta_L,\eta_\nu)$ in place of $\eta$. Both satisfy $2\eta N_0\le\frac14=c/4$. The remainder
$C\nu\sum_{Q'\sim Q}e_{Q'}\le C\nu N_0\,\mathbb E$ becomes $2C\nu N_0\,\mathbb E\le2C\eps
N_0\,\mathbb E$ for $\nu\le\eps$, and is included in $2a\,\mathbb E$.

The truncation $|\beta|\le n$ does not affect the estimate. The commutators of the truncated energy
involve only indices with $|\beta|\le n$, provided $n\ge7$, except in the $\ab$-form of the
$\bb$-part and in the direct form of~(III) of Proposition~\ref{prop:FPfix}. The other indices above
the target are $v$-Sobolev shifts of a $v$-low factor of the unknown. By
Lemma~\ref{lem:FPccross}(ii), they shift by at most $4$, to a $v$-level $\le\lfloor n/2\rfloor+4\le
n$ for $n\ge7$. In the $\ab$-form the extra derivative falls on the coefficient $\tilde G$, whose
energy is not truncated, by Definition~\ref{def:FPreg}. In the direct form of~(III), the plain $X$
needs $\partial_vu_{\mathbf k''}$ at $v$-level $|\beta''|+1$, which exceeds $n$ when $|\beta''|=n$.
This case is treated in the same remark with the viscous dissipation. Hence the bounds
hold for $n\ge\max(\nu^{-1},7)$, uniformly in $n$ and $\nu$, and
therefore for every $n$, since the truncated energies and dissipations are nondecreasing in $n$.
Finally, they pass to the limit $n\to\infty$ by monotone convergence.

\emph{(iv)} By~(iii), for every compact $K\subset\R^6$,
$G_\nu$ is bounded in $L^\infty(0,t_1;H^{N_y}_yH^m_v(K))$ for every $m$, because the box
energies control all these derivatives on the finitely many boxes meeting $K$, and
$\rad(t)\ge\rad_0/2$. Moreover, by the equation, $\partial_tG_\nu$ is bounded in
$L^\infty(0,t_1;H^{N_y-2}_yH^{m-2}_v(K))$,
uniformly in $\nu$. Indeed, $|M_\nu|\le\nu$, and the coefficients are bounded locally in the same
norms.
We take $m\ge N_y$ and $K$ a closed ball. Then $G_\nu$ is bounded in $L^\infty(0,t_1;H^{N_y}(K))$ and
$\partial_tG_\nu$ in $L^\infty(0,t_1;H^{N_y-2}(K))$, and $H^{N_y}(K)\subset H^{N_y-1}(K)$ is compact
(Rellich). By the Aubin--Lions--Simon lemma \cite[Cor.~4]{Sim87}, a subsequence converges in
$C([0,t_1];H^{N_y-1}(K))$. This space embeds into $C([0,t_1];C^2(K))$, because $N_y-1>3+2$ in
$\R^6$. By a diagonal argument, we obtain convergence in $C([0,t_1];C^2_{\rm loc})$.

In the equation, $\tilde A_\nu\to\tilde A$, $\tilde V_\nu\to\tilde V$ and $\tilde S_{0,\nu}\to\tilde
S_0$ locally uniformly. Here we use the fact that $\tilde G_\nu\to\tilde G$ in
$C([0,t_1];C^1_{\rm loc})$ by
the lemma on the mollified coefficients, and the moment bounds of the explicit
factors yield dominated convergence
in the $w$-integrals. The viscous term tends to $0$ locally uniformly. Thus the right side converges
locally uniformly, and $G$ solves the equation in $\mathcal D'$, with a continuous right side. Hence
$G\in C^1$ in $t$, and it is a classical solution.

The energy bounds pass to the limit by weak lower semicontinuity. For each box, each $n$ and each
$t$, $\partial^{\mathbf k} G_\nu(t)\rightharpoonup\partial^{\mathbf k} G(t)$ weakly in $L^2(\tilde
Q)$, and the dissipation integrals are lower semicontinuous in the same way (Fatou in $t$). For the
latter, whose coefficients $\tilde A_\nu$ vary with $\nu$, we use the duality $\int_0^t\!\iint\tilde
A\,X\cdot X\,\varphi_Q^2=\sup_\Phi\big(2\int_0^t\!\iint\tilde
A\,X\cdot\Phi\,\varphi_Q^2-\int_0^t\!\iint\tilde
A\,\Phi\cdot\Phi\,\varphi_Q^2\big)$ over smooth $\Phi$, with $X=Du_{\mathbf k}$ and $\tilde AX$
taken in $\mathcal D'$. Since $\int_0^t\!\iint\tilde A_\nu
Du_{\mathbf k}[G_\nu]\cdot\Phi\,\varphi_Q^2=-\int_0^t\!\iint u_{\mathbf k}[G_\nu]\,D\cdot(\tilde
A_\nu\Phi\varphi_Q^2)$, with $\tilde A_\nu\to\tilde A$ and $D\cdot\tilde A_\nu\to D\cdot\tilde A$
locally uniformly, each quantity under the supremum is the limit of the same quantity for $G_\nu$,
which does not exceed $\int_0^t\!\iint\tilde A_\nu Du_{\mathbf k}[G_\nu]\cdot Du_{
\mathbf k}[G_\nu]\,\varphi_Q^2$.
Then we take the supremum over $Q$ and $n$. The pointwise bounds pass to the limit directly.

\emph{(v)} The difference $\delta$ of two such solutions solves the homogeneous equation with
$\delta(0)=0$. By the level-$0$ identity with the top-order identity on boxes and
Lemma~\ref{lem:BKnb}(i),
after integration in time and the supremum over $Q$, we have
\begin{equation}
\sup_Qe^{(0)}_Q[\delta](t)\le C\eps\int_0^t\sup_Qe^{(0)}_Q[\delta].
\end{equation}
No commutator enters at level~$0$, so this needs only $\sup|\delta|<\infty$ and the local
$H^1$-regularity of the classical solutions. By Gronwall's inequality, we get $\delta=0$.
\end{proof}

\begin{lem}[The iteration converges]\label{lem:FPiter}
Let $G^0:=1$, and let $G^{j+1}$ be the solution of Proposition~\ref{prop:FPvisc}(iv) with $\tilde
G=G^j$. Then:
\begin{enumerate}[label=(\roman*)]
\item for all $j$ and all $t\in[0,t_1]$, $\mathbb E[G^j](t)^{1/2}\le2C_4\eps t$ (pointwise in $t$,
Remark~\ref{rem:FPpointwise}) and $|G^j-1|\le C\eps t\jap v^{-\sigma}$;
\item
\begin{equation}\label{eq:FPiter-contr}
\sup_{[0,t_1]\times\R^6}|G^{j+1}-G^j|\le\frac12\sup_{[0,t_1]\times\R^6}|G^j-G^{j-1}|
\end{equation}
for $\eps\le\eps_0$;
\item $G^j\to G$ uniformly. The limit solves the nonlinear equation
$\partial_tG=D\cdot(A[G]DG)+V[G]\cdot DG+S_0[G]G$, and satisfies the bounds of~(i).
\end{enumerate}
\end{lem}

\begin{proof}
(i) This holds by induction, from part~(iv) of that proposition, which also yields $G^{j+1}\in
C([0,t_1];C^2_{\rm loc})$ with $G^{j+1}(0)=1$ and, by $|G^{j+1}-1|\le C\eps t\le\frac14$,
$\frac12\le G^{j+1}\le2$. Thus $G^{j+1}$ is admissible in Definition~\ref{def:FPreg}.

(ii) We put $\delta^j:=G^{j+1}-G^j$. Then we have
\begin{equation}
\partial_t\delta^j=\mathcal L^j\delta^j+F^j,
\end{equation}
where $\mathcal L^j$ is the frozen operator of $G^j$. The forcing is
\begin{equation}
F^j=D\cdot\big(\ab[f_0\delta^{j-1}]DG^j\big)+\big(2\ab[f_0\delta^{j-1}]\zeta-\bb[f_0\delta^{j-1}]\big)\cdot
DG^j+\big(\ab[f_0\delta^{j-1}]:(D\zeta+\zeta\otimes\zeta)-\cb[f_0\delta^{j-1}]\big)G^j .
\end{equation}
Since $D$ commutes with the $w$-integral, and $\int a(w)\,Dh(z)\dd w=\int(\nabla\cdot a)(w)h(z)\dd
w$, we have $D\cdot\ab[h]=\bb[h]$. Hence
\begin{equation}
F^j=\ab[f_0\delta^{j-1}]:D^2G^j+2\ab[f_0\delta^{j-1}]\zeta\cdot
DG^j+\big(\ab[f_0\delta^{j-1}]:(D\zeta+\zeta\otimes\zeta)-\cb[f_0\delta^{j-1}]\big)G^j .
\end{equation}
Every term is linear in $\delta^{j-1}$, without derivatives. By the bounds on $G^j$ from~(i) and
Lemma~\ref{lem:RBsob}, we have
\begin{equation}
|\partial_Q^\alpha\partial_v^\beta G^j|\le C\eps t\,h_r^{a_r}\vartheta^{a_\perp}\qquad\text{for
}|\alpha|,|\beta|\le2 .
\end{equation}
We bound
$|\ab[f_0\delta]\xi\cdot\eta|\le\|\delta\|_\infty\ab[f_0]^{1/2}\xi\cdot\xi^{1/2}\,\ab[f_0]^{1/2}\eta\cdot\eta^{1/2}$,
and work in the frame of the box. The radial entries of $D^2G^j$ carry a factor $t^2h_r^2$
and meet $A\hat e\cdot\hat e\le C\eps\jap
y^{-2}$. The transverse entries carry the factor $t^2\vartheta^2$ and meet
$|A|\le C\eps\jap v^{-2\sigma}$. There is no $\bb$-term, since
$D\cdot\big(\ab[f_0\delta^{j-1}]DG^j\big)$ produces
$+\bb[f_0\delta^{j-1}]\cdot DG^j$, which cancels the drift term $-\bb[f_0\delta^{j-1}]\cdot DG^j$.
For the $\zeta$-terms we use $\zeta^\sharp$ and the moments. Then we obtain
\begin{equation}
|F^j|\le C\eps\sup|\delta^{j-1}| .
\end{equation}
The maximum principle for the frozen problem holds for the $\nu$-approximations, since
$G^{j+1}_\nu-G^j_\nu$ solves the $\nu$-problem of $(G^j)_\nu$ with $F^j$ built from
$(G^j)_\nu-(G^{j-1})_\nu=(\delta^{j-1})_\nu$, and $\sup|(\delta^{j-1})_\nu(t)|\le\sup_{s\le
t}|\delta^{j-1}(s)|$, the mollification being linear and causal. Moreover, this maximum principle passes to the limit.
By this principle and by the uniqueness of part~(v) of the same proposition, we have
\begin{equation}
\sup|\delta^j|\le e^{C\eps t_1}t_1C\eps\sup|\delta^{j-1}|\le\tfrac12\sup|\delta^{j-1}|
\end{equation}
for $\eps$ small.

(iii) The uniform limit exists by~(ii). The bounds of~(i) pass to it by lower semicontinuity, as in
part~(iv) of the same proposition. By interpolation, the uniform bounds and the uniform convergence
imply the convergence of $G^j$ in $C([0,t_1];C^2_{\rm loc})$. Hence $\ab[f_0G^j]\to\ab[f_0G]$,
together with the other coefficients, locally uniformly. Passing to the limit in the equation for
$G^{j+1}$, we obtain the nonlinear equation.
\end{proof}

\subsubsection{Transfer to the leafwise data norm at $t_1$}\label{sssec:FPtr}

At $t=t_1$, the leaf through $x=t_1\xi$ is
\begin{equation}
\ell_\xi[g](z)=g\big(t_1,z,\xi-z/t_1\big)=\eps\rho_0(\xi-z/t_1)F_0(z,\xi-z/t_1)\,G\big(t_1,z,\xi-z/t_1\big),
\end{equation}
since the profile is $g=f_0G$. Its $(\alpha,\beta)$-derivatives are the profile derivatives
$\partial_y^\alpha\partial_v^\beta g$, restricted to the leaf. For a profile $h$ we define the
\emph{raw leaf norm}
\begin{equation}
\mathfrak B_\xi[h]:=\sum_{n}w_n(\rad)^2\sum_{|\beta|=n}\sum_{|\alpha|\le
N_y}\delta_y^{2|\alpha|}\int\hat m_\xi(z)^2\,\frac{|\partial_y^\alpha\partial_v^\beta
h(z,\xi-z/t_1)|^2}{(\eps\rho_0(\xi))^2}\dd z .
\end{equation}
By the proof of Lemma~\ref{lem:LWfree}(ii), $\hat{\mathbb A}_\xi$ is bounded by $\mathfrak B_\xi$.
Indeed, the normalization $T^{-1}_{p}$ and the micro projection are bounded in these norms, with
constants depending on the leaf parameters, and the leaf parameters are within $C\eps$ of those of
$f_0$, by~(iii) below.
The same proof also shows $\mathfrak B_\xi[f_0]\le C_DB$. For the bootstrap leaf weight, we have
$\hat m_\xi^2\le
Ce^{(q_d/2)|z|^2}$ (since $q_0\le\theta_-q_d/16$). For the weight $e^{2q_d|z|^2}$ of
Theorem~\ref{thm:t1}, we have $e^{2q_d|z|^2}\le e^{\beta_F|z|^2/2}$ since $q_d\le\beta_F/4$. Both
are within
the Gaussian weight permitted by $B$ and by that lemma. The $y$-weight of the leaf norm
is $\delta_{y,\ell}$, the box parameter $\delta_y$ of
Subsections~\ref{ssec:tools}--\ref{ssec:drift}, and the radius $\rad_0$ of these subsections is
the layer radius $\rad_0^{\rm lay}$ of Appendix~\ref{app:const}. Neither of them is the
bootstrap parameter, and the Gaussian exponent may be taken equal to the data exponent
$q_d$, since the bound $\hat m_\xi^2\le Ce^{\beta_F|z|^2/2}$ suffices for the integrand.

\begin{prop}[(LWP$_0$): the transfer]\label{prop:FPtr}
Let $G$ be the solution of Lemma~\ref{lem:FPiter} on $[0,t_1]$, and let $\rad\le\rad_0/4$.
Then, uniformly in $\xi\in\R^3$:
\begin{enumerate}[label=(\roman*)]
\item $\mathfrak B_\xi[f_0(G-1)]\le C\,\eps^2t_1^2$. Hence
\begin{equation}
\hat{\mathbb A}_\xi(t_1)\le C_DB\,(1+C\eps t_1)^2\qquad\text{and}\qquad\hat{
\mathbb A}_{\ge1,\xi}(t_1)\le C_D\big((\rad/\Lambda)^2+\delta_y^2\big)B+C\eps^2t_1^2 ;
\end{equation}
\item the normalized leaves of $g(t_1)$ satisfy (G1)--(G3), with constants within a factor $1+C\eps
t_1$ of those of Lemma~\ref{lem:LWfree}(i);
\item the leaf mass is $\eps\rho_0(\xi)\,(1+O(L/t_1)+O(\eps t_1))$, and the other leaf parameters
are within $C(L/t_1+\eps t_1)$ of those of $F_0(\cdot,\xi)$.
\end{enumerate}
Consequently the hypotheses of Theorem~\ref{thm:t1} on the data at $t_1$ hold. The bound~(i) holds,
with the same proof, for $\mathfrak B_\xi$ defined with any weight $\le Ce^{\beta_F|z|^2/2}$ in
place of $\hat m_\xi^2$.
\end{prop}

\begin{proof}
\emph{Step 1: pointwise and leafwise bounds for $G-1$.} Let $\mathbf k'=(\alpha',\beta')$.
\begin{itemize}
\item For $|\alpha'|\le N_y-2$, by Lemma~\ref{lem:RBsob} on the box containing $(z,\xi-z/t_1)$, we
have
\begin{equation}
\begin{aligned}
&|\partial^{\mathbf k'}(G-1)|\le
C\max_{|\varsigma_y|,|\varsigma_v|\le2}\ell_Q^{\varsigma_y}\big(W_{
\mathbf k'+\varsigma}\big)^{-1}\mathcal H_{\mathbf k'},\qquad\mathcal H_{
\mathbf k'}^2:=\sum_{|\varsigma_y|,|\varsigma_v|\le2}e^\sim_Q(
\mathbf k'+\varsigma),\\
&\sum_{\mathbf k'}\mathcal H_{\mathbf k'}^2\le C^{N_y}N_0\,\mathbb E(t_1)\le C\eps^2t_1^2 ;
\end{aligned}
\end{equation}
\item For the top two levels, we use the $L^2$ form instead. The map $z\mapsto(z,\xi-z/t_1)$
parametrizes a $3$-dimensional graph over the $y$-space. For fixed $z$, by Sobolev's inequality in
$v$ on the unit ball $B(\xi-z/t_1,1)$, we have
\begin{equation}
|\partial^{\mathbf k'}u(z,\xi-z/t_1)|^2\le
C\sum_{|\gamma|\le2}\int_{B(\xi-z/t_1,1)}|\partial_v^\gamma\partial^{\mathbf k'}u(z,v)|^2\dd v .
\end{equation}
We integrate over a unit cube $I$ of $z$. On $I$ the centers $\xi-z/t_1$ move by at most
$\sqrt3/t_1\le\sqrt3$. Consequently the right side is bounded by the $L^2$ norm over $I\times
B(\xi-z_I/t_1,3)$, which is covered by $C$ velocity balls and, in $y$, by the boxes meeting $I$, of
total volume $\le C$.
\end{itemize}
In both cases, we obtain
\begin{equation}
\int_I|\partial^{\mathbf k'}(G-1)(z,\xi-z/t_1)|^2\dd z\le C\sup_{Q\cap(I\times
B)\ne\emptyset}\big(W^Q_{\mathbf k'+\varsigma}\big)^{-2}\,\mathcal H_{\mathbf k'}^2,
\end{equation}
with a $v$-shift $|\varsigma|\le2$ and $\sum_{\mathbf k'}\mathcal H_{\mathbf k'}^2\le C\eps^2t_1^2$.

\emph{Step 2: the inverse weights on the leaf.} On the leaf, $v=\xi-z/t_1$, and
\begin{equation}
\big(W^Q_{\mathbf k'}\big)^{-1}=\delta_y^{-|\alpha'|}\,h_r^{a'_r}\vartheta^{a'_\perp}\,\tau_{\alpha'}^{-1}\,\jap{v_0}^{-\sigma}\,w_{\beta'}(\rad(t_1))^{-1}.
\end{equation}
Here $\tau_{\alpha'}^{-1}\le h_r\vartheta^2$, and
\begin{equation}
h_r\le\jap z+\vartheta(v),\qquad\vartheta(v)^{m}\jap{v_0}^{-\sigma}\le C_m .
\end{equation}
Hence we have
\begin{equation}
\big(W^Q_{\mathbf k'}\big)^{-1}\le C_{N_y}\,\delta_y^{-|\alpha'|}\jap
z^{N_y+1}\,w_{\beta'}(\rad(t_1))^{-1}.
\end{equation}
This bound does not depend on $\xi$. \emph{Without the factor $\jap{v_0}^{\sigma}$, it would contain
$\vartheta(\xi-z/t_1)^{a'_\perp+2}$. For $|z|\le t_1$ this is
$\simeq(\log\jap\xi)^{(a'_\perp+2)/2}$, which is unbounded in $\xi$.} The velocity weight was
introduced for this reason.

\emph{Step 3: products.} By the Leibniz rule, we have
\begin{equation}
\partial^{\mathbf k}\big(f_0(G-1)\big)=\sum_{\nu\le\mathbf k}\binom{\mathbf k}\nu\partial^\nu
f_0\,\partial^{\mathbf k-\nu}(G-1).
\end{equation}
By Lemma~\ref{lem:FPexp} and the proof of Lemma~\ref{lem:LWfree}, we have
\begin{equation}
\frac{|\partial^\nu f_0(z,\xi-z/t_1)|}{\eps\rho_0(\xi)}\le C_{N_y}C^{|\nu_v|}\nu_v!R^{-|\nu_v|}\jap
z^{|\nu_y|}e^{L|z|/t_1}F_{\rm G}(z).
\end{equation}
We insert the bounds of Steps~1--2. Then the resulting $z$-integrand is bounded by
\begin{equation}
C\,\hat m_\xi(z)^2\jap z^{2N_y+2}e^{2L|z|/t_1}F_{\rm G}(z)^2\le C\jap
z^{2N_y+2}e^{2L|z|}e^{-3\beta_F|z|^2/2},
\end{equation}
which is integrable, uniformly in $\xi$. Only the bound $F_{\rm G}\le K_{
\mathcal U}C_Fe^{-\beta_F|z|^2}$, which follows from (D1), enters.

In $v$, the weight sum is
\begin{equation}
\sum_\beta w_\beta(\rad)^2\Big(\sum_{\nu_v\le\beta}\binom\beta{\nu_v}C^{|\nu_v|}\nu_v!R^{-|\nu_v|}\,\frac{
\mathcal H_{\beta-\nu_v}}{w_{\beta-\nu_v+\varsigma}(\rad(t_1))}\Big)^2 .
\end{equation}
With $\rad\le\rad_0/4\le\rad(t_1)/2$ and $\beta'=\beta-\nu_v$, the ratio of the weights,
multiplied by the binomial and by $C^{|\nu_v|}\nu_v!R^{-|\nu_v|}$, is bounded by
\begin{equation}
C\,\rad(t_1)^{-2}(|\beta|+3)^{s+2}\Big(\frac\rad{\rad(t_1)}\Big)^{|\beta'|}\Big(\frac{C\rad}R\Big)^{|\nu_v|}.
\end{equation}
By Young's inequality on $\mathbb Z^3$, the sum is bounded by $C\sum\mathcal H^2$, provided
$C\rad/R\le\frac12$. The $y$-sum is finite, with
$\delta_y^{2|\alpha|}\delta_y^{-2|\alpha'|}\le1$. This proves the first claim of~(i). The bounds for
$\hat{\mathbb A}_\xi$ follow from the triangle inequality in $\mathfrak B_\xi^{1/2}$ and
part~(ii) of the same lemma.

\emph{(ii) and (iii).} By Proposition~\ref{prop:FPclose}, we have $|G-1|\le C\eps t_1$ on the leaf.
Hence the leaf of $g$ is pointwise between $(1\pm C\eps t_1)$ times the leaf of $f_0$. This implies
the upper bounds (G1) and (G3), with constants multiplied by $1+C\eps t_1$, and the bounds on the
parameters in~(iii), since mass, mean and temperature are ratios of positive integrals.
The normalizations $T_p$ of the two leaves differ by $O(\eps t_1)$ in their parameters. For (G2),
we note that $\ab[\cdot]$ is monotone on nonnegative functions. Hence
\begin{equation}
\ab[\mathcal G_g]\ge(1-C\eps t_1)\ab[\mathcal G_{f_0}]
\end{equation}
for the normalized leaves $\mathcal G$, up to the change of normalization. That change moves the
lower bound of (G2) by $O(\eps t_1)$, by the continuity of $\ab$ under affine changes of variables.
Hence (G2) holds with $c_G(1-C\eps t_1)$.

Finally, Theorem~\ref{thm:t1} requires these bounds in the form $\le C\eps$, with $C$ independent of
$\eps$. Since $t_1=t_1(\mathfrak D)$ is fixed, $C\eps t_1=C'\eps$.
\end{proof}

\subsection{The drift term in the initial layer}\label{ssec:drift}

In this subsection we treat the $\bb$-part of the drift in the initial layer, in the configuration
that Lemma~\ref{lem:FPdlow} does not cover. This completes the proof of Theorem~\ref{thm:layer}.

In configuration~(IV) of Lemma~\ref{lem:FPc}, the coefficient factor $h$ is $v$-high
($|\beta'|>|\beta|/2$) and $X$ is $v$-low and $y$-high. After we write $\bb[F]=\ab[DF]$, the
extra $v$-derivative falls on the factor of higher $v$-level. Then we are in the first case of
Lemma~\ref{lem:TWcalc}(iv), which produces the factor $(|\beta'|+1)/\rad\simeq|\beta|/\rad$.
The term has the structure
\begin{equation}
\frac{|\beta|}\rad\cdot\mathcal H_{\mathbf k'+e_v}\ (\text{energy})\cdot d_{\mathbf k''}\
(\text{dissipation})\cdot\eps^{1/2}a_{\mathbf k}\ (\text{energy}),
\end{equation}
and after we apply Young's inequality to the dissipation, we are left with a loss of order
$|\beta|^2$, which cannot be controlled by $\theta_{\mathbf k}\propto|\beta|$. An example is
$h=\partial_v^{\beta_\mu}G$ with
$\mu=(0,\beta_\mu)$ purely in $v$, and $X=Du_{(\alpha,\beta-\beta_\mu)}$ with $|\alpha|=N_y$. When
$|\alpha''|\ge N_y-2$ the direct form does not help either, since it needs $X$ plain, by
Lemma~\ref{lem:RBgain}(ii), and $X$ is then at the top $y$-levels. We treat these terms in two
steps. First, we put the $\bb$-part of the drift in the direct form or in the $\ab$-form according
to the configuration, so that a Cauchy--Kovalevskaya loss occurs only in two corner cases. Second,
we weaken the weights of the three top $y$-levels by a half power of $(|\beta|+1)$. Once we move
this half power to the target, it absorbs the loss in the corner cases.
A grading of the radius of analyticity by the $y$-level, as used for the forcing in
Remark~\ref{rem:UFgrade}, does not work here, because it penalizes every factor that is measured
above the level of the target, and such factors occur at every level.

In configuration~(III) we use the direct form for all $|\alpha'|\ge N-2$. With the $\ab$-form, $h$
would be raised by $t\partial_y$ into the $\tau$-levels while the target is not in them, and nothing
would then cancel the transfer factor of $h_r$ or the slab loss. Moreover, in configuration~(III)
$X$ is $v$-high, so Case~2 of the proof of Lemma~\ref{lem:FPd} does not apply, and the case table
below covers this configuration as well.

Throughout this subsection (but not in Subsection~\ref{ssec:cellsproofs}, where $N$ is a truncation
level), $N:=N_y\ge20$ is the top $y$-level of Subsection~\ref{ssec:layer}, $s\ge6$ is the exponent
of $w_\beta$, and the configurations (I)--(IV) are those introduced above:
$h=\partial^{\mathbf k'}G$ is the coefficient factor, $X=Du_{\mathbf k''}$, the target is
$u_{\mathbf k}$, and $\mathbf k=\nu+\mathbf k'+\mathbf k''$ with $\nu$ the explicit part.

\subsubsection{The weakened weights}

At the top $y$-levels we weaken the weights as in Definition~\ref{def:FPtop}. By
Lemma~\ref{lem:FPhalf}, the factor $(|\beta|+1)^{-1/2}$ compensates half of the
Cauchy--Kovalevskaya loss of the $\bb$-part of the drift.

\begin{defin}\label{def:FPtop}
Let $c_\alpha:=\mathbf 1\{|\alpha|\ge N-2\}$ and $\iota_\alpha:=\mathbf 1\{|\alpha|\ge N-4\}$, and
define, in place of the weight of Subsection~\ref{sssec:FPset},
\begin{equation}
W^\flat_{\mathbf k}:=W_{\mathbf k}\,(|\beta|+1)^{-c_\alpha/2}\,\vartheta(v_0)^{-\iota_\alpha},\qquad
\mathbf k=(\alpha,\beta).
\end{equation}

All quantities $a_{\mathbf k},d_{\mathbf k},e_Q,\mathfrak d_Q,\mathbb E$ are defined with $W^\flat$.
The factor $\vartheta(v_0)^{-\iota_\alpha}$ is used only in Subsection~\ref{sssec:FPcross}, and it
plays no role in Lemma~\ref{lem:FPtopalg}, which concerns the factor $(|\beta|+1)^{-c_\alpha/2}$.
See Remark~\ref{rem:FPiota}.
\end{defin}

Note that $(|\beta|+1)^{-c_\alpha/2}\vartheta(v_0)^{-\iota_\alpha}$ does not depend on $t$. Thus
the identity~\eqref{eq:FPid} holds with the same $\theta_{\mathbf k}$.

\begin{lem}[Weight algebra]\label{lem:FPtopalg}
Consider a term with target $\mathbf k=(\alpha,\beta)$ and two factors measured at levels
$\mathbf k_1=(\alpha_1,\beta_1)$ and $\mathbf k_2=(\alpha_2,\beta_2)$, after all shifts. The
combined factor of Lemma~\ref{lem:FPOmega} computed with $W^\flat$ equals the one computed with $W$,
up to the factors $\vartheta(v_0)^{-\iota_\alpha}$ (treated in Remark~\ref{rem:FPiota}), times
\begin{equation}
\Pi:=\Big(\frac{(|\beta_1|+1)^{c_{\alpha_1}}(|\beta_2|+1)^{c_{\alpha_2}}}{(|\beta|+1)^{c_\alpha}}\Big)^{1/2}
.
\end{equation}
\begin{enumerate}[label=(\roman*)]
\item If $c_{\alpha_1}=c_{\alpha_2}=0$, then $\Pi\le1$.
\item If $c_\alpha=1$ and one of the two factors is $y$-low, i.e.\ its $y$-level before shifts is
$\le\lfloor N/2\rfloor$, then at most one of $c_{\alpha_1},c_{\alpha_2}$ equals $1$, and $\Pi\le2$.
\item If $c_\alpha=0$, the first factor is $y$-low and $c_{\alpha_2}=1$, then $c_{\alpha_1}=0$ and
$\Pi\le(|\beta_2|+1)^{1/2}$, where $|\beta_2|$ includes the $v$-shift of the second factor (this
shift is at most $2$ for a Sobolev step, so the bound is $(|\beta_2^{\rm unshifted}|+3)^{1/2}$, and
this shift is at most $3$ when the $\partial_v$ of a plain $X$ or of $\partial_vh$ is added, and then the
bound is $(|\beta_2^{\rm unshifted}|+4)^{1/2}$).
\item If $c_\alpha=1$, $c_{\alpha_2}=1$, $c_{\alpha_1}=0$, and the term carries the loss
$(|\beta_1|+1)/\rad$ on the first factor, then the total factor is
\begin{equation}
\Pi\cdot\frac{|\beta_1|+1}\rad\le\frac C\rad\,(|\beta_1|+1)^{1/2}\big(\min(|\beta_1|,|\beta_2|)+1\big)^{1/2},
\end{equation}
with $C$ absolute.
\item If $c_\alpha=1$, $c_{\alpha_1}=c_{\alpha_2}=0$, and the term carries the loss
$(|\beta_1|+1)/\rad$ on the first factor, then the total factor is $\le
C(|\beta_1|+1)^{1/2}/\rad$.
\end{enumerate}
\end{lem}

\begin{proof}
(i) is clear. (ii) In every case of the table of Proposition~\ref{prop:FPfix}, and in
Lemmas~\ref{lem:FPb}--\ref{lem:FPe}, the $y$-low factor is raised in $y$ by at most $3$, namely by
$1$ from $t\partial_y$ or from the $t\nabla_y$ in a plain $X$, and by $2$ from a Sobolev embedding.
Hence its level does not exceed $\lfloor N/2\rfloor+3$, which is $<N-2$ because $N=N_y\ge20$ in
Theorem~\ref{thm:main}. This threshold needs only $N\ge11$, and the threshold $N-4$ of $\iota$ needs
$N\ge15$ (see Remark~\ref{rem:FPiota}). Thus only the other factor can have $c=1$. Its $v$-level is
bounded by $|\beta|+3$, the $3$ coming from a $v$-shift if it is the $v$-low factor, and we obtain
$\Pi\le((|\beta|+4)/(|\beta|+1))^{1/2}\le2$. In the non-integrated form of
Subsection~\ref{sssec:FPhalf}, needed in the crossing subcase of
Lemma~\ref{lem:FPccross}, the factor $D^2u_{\mathbf k''}$ is measured at $y$-level
$\le|\alpha''|+4\le5$, since $|\alpha''|\le|\alpha|-|\alpha'|\le1$ there.

The $y$-high factor is also raised, in the following places, and these are not covered by~(ii).
\begin{itemize}
\item By $1$: by $t\partial_yh$ in the $\ab$-form of~(III), and by the $t\nabla_y$ of a plain $X$
in the direct form of~(IV). When this lifts it above the level of the target, the crossing is
treated case by case in the proof of that proposition.
\item By $2$, through a $y$-Sobolev step (Lemma~\ref{lem:RBsob}, or
Lemma~\ref{lem:slab}(ii)): in~(II) and~(III) of Lemma~\ref{lem:FPc}, in~(II) and in the $\ab$-form
of~(III) of the drift, in the high-$h$ branch and in the $A\zeta$-part of
Lemma~\ref{lem:FPd}, and in Lemma~\ref{lem:FPe} (second and third lines of~(T5)).
If the lifted factor is $v$-low, a crossing into $c=1$ is controlled by~(iii).
If it is $v$-high, the crossing is the subject of Lemma~\ref{lem:FPccross}.
\end{itemize}
(iii) By the proof of~(ii), we have $c_{\alpha_1}=0$, and thus $\Pi=(|\beta_2|+1)^{1/2}$. (iv) We
write $a=|\beta_1|$, $b=|\beta_2|$ after shifts. Then we have
$|\beta|+1\ge\max\{1,a+b+1-\varsigma\}$ with $\varsigma\le3$, where $\varsigma$ accounts for the
loss itself and the Sobolev shifts. This implies that $(a+1)(b+1)/(|\beta|+1)\le C(\min(a,b)+1)$.
For instance, $a=b=1$ and $\varsigma=3$ allow $|\beta|+1=1$, and the ratio $4$ is
$\le2C=C(\min(a,b)+1)$. Therefore, we get
\begin{equation}
\frac{(a+1)(b+1)^{1/2}}{(|\beta|+1)^{1/2}}=(a+1)^{1/2}\Big(\frac{(a+1)(b+1)}{|\beta|+1}\Big)^{1/2}\le
C(a+1)^{1/2}(\min(a,b)+1)^{1/2}.
\end{equation}
(v) Here we have $\Pi=(|\beta|+1)^{-1/2}$ and $a\le|\beta|+3$. This completes the proof.
\end{proof}

In case~(iii), if the second factor is $v$-low, the factor $(|\beta_2|+1)^{1/2}$ is a half shift on
a $v$-low factor. As with a shift, we control it by Lemma~\ref{lem:ULconv}. Indeed, with one more
power, the sequence $k$ in its proof is $(m+3)^3(m+4)/(m+1)^s$, which lies in $\ell^2(\mathbb N^3)$
because $2s-10>1$ for $s\ge6$. A half power is better. If the second factor is $v$-high, we
refer to Lemma~\ref{lem:FPccross}. In cases~(iv) and~(v), the factor $(\min+1)^{1/2}$ is a half
shift on the factor of lower $v$-level, and we control it in the same way. The factor
$(|\beta_1|+1)^{1/2}$, on the factor that carries the loss, which may be $v$-high, is bounded by
$C(|\beta|+1)^{1/2}$. We move it to the target, where Lemma~\ref{lem:FPhalf} absorbs it by
$\theta_{\mathbf k}$.

\subsubsection{Choice of the form of the drift terms}

We choose, case by case, the form in which we estimate the $\bb$-part of a drift commutator. We
write $\bb[\partial^\nu f_0\,h]\cdot X\,u_{\mathbf k}$ for the $\bb$-part of a drift commutator. Its
two forms are
\begin{itemize}
\item the \emph{direct form}, with the kernel $b(w)$, Lemma~\ref{lem:RBgain}(ii) and line~(T5) of
Lemma~\ref{lem:ULtri}, which needs $X$ \emph{plain}, i.e.\ the energies of $\partial_vu_{
\mathbf k''}$ and $t\nabla_yu_{\mathbf k''}$;
\item the \emph{$\ab$-form} $\bb[F]=\ab\langle DF\rangle$, Lemma~\ref{lem:ULaform}(iii), which puts
$D=\partial_v-t\partial_y$ on $h$.
\end{itemize}
A $v$-derivative on the $v$-high factor produces a factor $(|\beta_{\rm high}|+1)/\rad$, the
first case of Lemma~\ref{lem:TWcalc}(iv). A $v$-derivative on the $v$-low factor only shifts the
level. A $y$-derivative never produces a factor $|\beta|$, and it only raises the $y$-level.

\begin{prop}[The $\bb$-part of the drift]\label{prop:FPfix}
With the weights $W^\flat$, the $\bb$-part of the drift commutators satisfies the bound of
Lemma~\ref{lem:FPd}, up to two additional terms, apart from the crossing terms of configuration~(II)
(case~(a) of Lemma~\ref{lem:FPccross}), which are bounded by~\eqref{eq:FPccross} by the part of
$\theta_{\mathbf k}$ used in that lemma:
\begin{equation}
\begin{aligned}
\Big|\frac1{|Q|}\sum_{\mathbf k} W^{\flat2}_{\mathbf k}\sum_{0<\mu\le\mathbf k}\binom{
\mathbf k}\mu\iint\partial^\mu\bb\cdot Du_{\mathbf k-\mu}\,u_{
\mathbf k}\varphi_Q^2\Big|&\le\Big(C\rad^{-2}\eps^{1/2}+C_{
\mathbb E}
\mathbb E^{1/2}\Big)\Big(\Sigma_Q+\sum_{Q'\sim Q}e_{Q'}\Big)\\
&\quad+C\eps\,\mathbb E^{1/2}\sum_{Q'\sim Q}e_{Q'}+\frac1{4|Q|}\sum_{\mathbf k}\theta_{\mathbf k}
a_{\mathbf k}^2 .
\end{aligned}
\end{equation}
In the two corner cases a half power $(|\beta|+1)^{1/2}$ of the Cauchy--Kovalevskaya loss remains.
It sits on the target $u_{\mathbf k}$, which is taken in energy on its own box $Q$, and it is
absorbed by $\theta_{\mathbf k}$ through Lemma~\ref{lem:FPhalf}, applied to each of the two corner
cases separately, with its own $B_i$ and $\delta_i=\frac18$. This is the last term. No
$\theta$-energies of neighboring boxes or of the coefficient box $P''$ are used.
\end{prop}

\begin{proof}
We go through the configurations. In each case we name the form, the norms in which the three
factors are taken, the levels at which they are measured, and the factor $\Pi$ of
Lemma~\ref{lem:FPtopalg}.

\emph{(I) $h$ $v$-low and $y$-low.} We use the $\ab$-form, as in Lemma~\ref{lem:FPdlow}. The
$v$-derivative falls on the $v$-low factor $h$, which is the second case of
Lemma~\ref{lem:TWcalc}(iv). The part $t\partial_yh$ raises the level of $h$ to $|\alpha'|+1$, and
with the Sobolev shift to $|\alpha'|+3\le N/2+3$, so $c=0$ for $h$. By the weight
algebra~(i)--(ii), we have $\Pi\le4$, and we obtain the bound of Lemma~\ref{lem:FPdlow}.

\emph{(II) $X$ $v$-low and $y$-low, $h$ arbitrary.} We use the direct form, as in Case~2 of the
proof of Lemma~\ref{lem:FPd}. It puts no derivative on $h$, so it makes no demand on the $v$-level
of $h$. Here $X$ is plain, and we take it pointwise by Sobolev in $(y,v)$, as in the second line
of~(T5). Its levels do not exceed $|\alpha''|+3\le N/2+3$ in $y$ and $|\beta''|+3$ in $v$, both on
the low factor. We take $h$ in $L^2$ on the slab, and we treat the slab as in
Lemma~\ref{lem:FPc}(II).
\begin{itemize}
\item If $|\alpha'|\le N-2$, we use Lemma~\ref{lem:slab}(ii), with a $y$-lift of $h$ by $2$. If
this lift crosses into $c=1$ while $c_\alpha=0$, which needs $|\alpha'|\in\{N-4,N-3\}$, we have
the factor $\Pi\le(|\beta'|+3)^{1/2}$. On a $v$-low $h$ it is controlled by~(iii) of the weight
algebra, and on a $v$-high $h$ we treat it by Lemma~\ref{lem:FPccross}. If the lift crosses into
the $\tau$-levels, which needs $|\alpha'|\in\{N-3,N-2\}$, the unpaired $h_r''$ is absorbed by
$F_0(y'')\rho_0(v'')$, since the estimate is absolute, as in the proof of
Lemma~\ref{lem:FPcross}(i).
\item If $|\alpha'|\ge N-1$, we use~(i) of the slab lemma and Lemma~\ref{lem:BKcomp}(i). Then
$|\alpha|\ge N-1$ as well, the target and $h$ both carry $\tau$ and $c=1$, and there is no lift.
\end{itemize}
Otherwise we have $\Pi\le2$. We obtain a bound $\le C\eps\,\mathbb E^{1/2}\sum_{Q'\sim Q}e_{Q'}$,
plus, in the crossing just mentioned, a part absorbed by $\theta_{\mathbf k}$.

\emph{(III) $h$ $v$-low and $y$-high, $X$ $v$-high and $y$-low.}
\begin{itemize}
\item If $|\alpha'|\le N-3$, we use the $\ab$-form. The $v$-derivative falls on the $v$-low factor
$h$.
The part $t\partial_yh$ raises $h$ to level $|\alpha'|+1\le N-2$, so it never reaches the
$\tau$-levels $\ge N-1$. It crosses into $c=1$ from a target with $c_\alpha=0$ in two ways, namely
through $t\partial_y$ itself, for $\alpha'=\alpha$ with $|\alpha|=N-3$, and through the $y$-lift by
$2$ of
Lemma~\ref{lem:slab}(ii) below, applied to $t\partial_yh$ at level $|\alpha'|+1$ and to
$\partial_vh$ at level $|\alpha'|$, for $|\alpha'|\in\{N-5,N-4,N-3\}$. In every case $h$ is $v$-low,
and by~(iii) of the weight algebra we get $\Pi\le(|\beta'|+4)^{1/2}$ on the $v$-low factor (the
part $\partial_vh$ carries the $v$-shift $1+2$), which suffices. We use the norms of~(T3), that is,
$h$ in
$L^\infty_vL^2_y$ on the slab, $X$ in $L^\infty_y$ with the $A$-weight, and the $Y$-side
$u_{\mathbf k}\,e_j$ with the trace $\le C\eps$ of $\ab[f_0]$. Since $t\partial_yh$ is at level $\le
N-2$, this lift is available, since it needs only two more $y$-derivatives, i.e.\ level $\le N$.
Therefore there is no slab loss, as in~(III) of the diffusion. For $|\alpha'|\in\{N-4,N-3\}$
and $|\alpha|\le N-2$ this lift reaches the $\tau$-levels while the target does not, and the
unpaired factor $\tau''^{-1}$ is case~(C2) of Lemma~\ref{lem:FPcross} (or Lemma~\ref{lem:FPYside},
Remark~\ref{rem:FPYside}).
\item If $|\alpha'|\ge N-2$, we use the direct form, with~(T5), last line, $y$ and $v$ exchanged,
in the form
\begin{equation}
\iint K[\partial^\nu f_0h]\,|U||U'|\chi\le C\eps\,\|\partial_v^{\le2}\tilde h\|_{\rm
ul}\,\big\|\sup_y|U|\big\|_{L^2_v}\,\|U'\|_{L^2(\chi)} .
\end{equation}
  \begin{itemize}
  \item \emph{The factors.} We take $h$ in $L^\infty_vL^2_y$, by $v$-Sobolev on the $v$-low $h$,
with $v$-shift $\le2$. There is no $D$ on $h$, so $h$ is not raised by $D$. In $y$ it is raised only
by the Sobolev step of the slab lemma for $|\alpha'|=N-2$ (\emph{The slab} below), which is an
absolute crossing into the $\tau$-levels, and we treat it as in~(II). $U=X$ is plain, and we take it
in $L^2_vL^\infty_y$ by $y$-Sobolev on $X$, whose $y$-level is
$|\alpha''|+|\nu_y|+1+2\le|\alpha|-|\alpha'|+3\le5<N-2$. The target $U'=u_{\mathbf k}$ is in plain
$L^2$, i.e.\ the energy $a_{\mathbf k}$. Exchanging $y$ and $v$ in~(T5) only exchanges the order of
the Minkowski steps (first in $y$, then in $w$), with $\int|w|^{\gamma+1}\rho_0(v-w)\dd w\le C\jap
v^{\gamma+1}$, and no $\Theta_{z_0}$ enters. We treat the $t\nabla_y$-part of the plain $X$ by
Lemma~\ref{lem:RBgain}(ii), as in~(IV) below. The radial factor $th_r$ of $X_r$ meets the gain
$|y|^{-1}$, and the transverse factor $t\vartheta$ meets the decay $\jap v^{\gamma+1}$. The
additional kernel $t|w|^{\gamma+2}$ of the same lemma enters the same Fubini step, with
$\int|w|^{\gamma+2}\rho_0(v-w)\dd w\le C\jap v^{\gamma+2}$.
  \item \emph{The loss.} The $v$-derivative in $X=\partial_vu-t\nabla_yu$ falls on the $v$-high
factor $X$ and produces the factor $(|\beta''|+1)/\rad$.
Here we have $c_\alpha=c_{\alpha'}=1$ and $c=0$ for $X$. Thus, by~(iv) of the weight algebra, the
total factor is $\le C(|\beta''|+1)^{1/2}(|\beta'|+3)^{1/2}/\rad$. The factor
$(|\beta'|+3)^{1/2}$ is a half shift on the $v$-low factor $h$. We put the factor
$(|\beta''|+1)^{1/2}\le(|\beta|+1)^{1/2}$ on the target, which is in energy on $Q$, and
Lemma~\ref{lem:FPhalf} absorbs it, with $B=C\eps\rad^{-1}
\mathbb E^{1/2}$ and $\sum_{\mathbf k} Z_{\mathbf k}^2\le C\sum_{Q'\sim Q}e_{Q'}$. Here the plain
$X$ is an energy of $Q$ and its neighbors. (For the truncation $|\beta|\le n$, see
Remark~\ref{rem:FPtrunc8}.)
  \item \emph{The slab.} If $|\alpha'|=N-2$, part~(ii) of the slab lemma applies (it needs
level $N$), and there is no slab loss. If $|\alpha'|\in\{N-1,N\}$, then also $|\alpha|\ge N-1$, so
the target and $h$ both carry $\tau$. Then we integrate the slab of $h$ against
$F_0(y'')\rho_0(v'')$ by~(i) of the slab lemma and Lemma~\ref{lem:BKcomp}(i), as in the high-$h$
branch of the drift lemma. More precisely, for fixed $w$ and $v$, $y''=y+tw$ runs over the
$y$-projection of the slab and $v''=v-w$ over $Q_v-w$, of width $1$. By the $v$-Sobolev step at unit
scale, we have
  \begin{equation}
  \int_{Q_y+tw}\sup_{v''\in Q_v-w}|h|^2\dd y''\le C\sum_{|j|\le2}\int_{S_w}|\partial_v^jh|^2 ;
  \end{equation}
this step does not involve $y$, and so it does not change the factor $L$. Since both
the target and $h$ carry $\tau$, the slab factor together with the ratio of the weights is
$\Lambda^2=L\,(h_r''\vartheta''^2/h_r\vartheta^2)^2$. By the same compensation bound, this is a
polynomial in $\jap{y''}$ times $\vartheta''^C$, which is controlled by $F_0(y'')\rho_0(v'')$ in the
Fubini step. To put $\Lambda$ inside $\sup_{v''}$, we apply the $v$-Sobolev step to
$\tilde\Lambda f_0h$, with
$\tilde\Lambda$ a smoothing of $\Lambda$ at unit $v$-scale; $\vartheta''$ changes by a bounded
factor on such scales.
  \end{itemize}
\end{itemize}
\emph{The levels $|\alpha'|=N-2,N-1$.} We do not use the $\ab$-form here. Indeed, $t\partial_yh$
would
lie at level $N-1$ or $N$, in the $\tau$-levels. For $\alpha'=\alpha$, $|\alpha|=N-2$, the target is
not in the $\tau$-levels. Lemma~\ref{lem:slab}(ii) is not available for $t\partial_yh$ at level
$N-1$, and Lemma~\ref{lem:FPslab} uses the $\tau$ of the target to compensate the slab factor $L$,
which contains $h_r(Q)$. Without that $\tau$, the slab loss is not compensated. However, the direct
form avoids raising $h$ at all. (A second difficulty, the unpaired factor $h_r''$ of $\tau''^{-1}$,
is not specific to this case, and we treat it in Subsection~\ref{sssec:FPcross}.)

\emph{(IV) $h$ $v$-high and $y$-low, $X$ $v$-low and $y$-high.}
\begin{itemize}
\item If $|\alpha|\le N-3$, we use the direct form, with~(T5), last line. We take $h$ in
$L^\infty_yL^2_v$ ($y$-Sobolev on the $y$-low $h$, level $\le|\alpha'|+2$, no $v$-derivative), and
$X$ plain in
$L^2_yL^\infty_v$ ($v$-Sobolev on the $v$-low $X$). The levels of $X$ are $|\beta''|+3$ in $v$ and
$|\alpha''|+1\le N-2$ in $y$. There is no crossing into the $\tau$-levels. For $\alpha''=\alpha$
with $|\alpha|=N-5$, $X$ crosses the threshold $N-4$ of $\iota$, with a factor $\vartheta(v_0)$. For
$\gamma<-2$ this factor is controlled by $\vartheta(v)\jap v^{\gamma+2}\le C$, as in
Remark~\ref{rem:FPiota}. A crossing into $c=1$ happens only for $\alpha''=\alpha$ with
$|\alpha|=N-3$, and then we have $\Pi\le(|\beta''|+4)^{1/2}$ on the $v$-low factor, by~(iii) of the
weight algebra. The radial factor $th_r$ of $X_r$ meets the gain $|y|^{-1}$. The transverse factor
$t\vartheta$ meets $\int|w|^{\gamma+1}\rho_0(v-w)\dd w\le C\jap v^{\gamma+1}$. We obtain a bound
$\le C\eps\,\mathbb E^{1/2}\sum_{Q'\sim Q}e_{Q'}$.
\item If $|\alpha|\ge N-2$, we use the $\ab$-form, for every $\alpha''$. The part $t\partial_yh$ has
no
$|\beta|$-loss, and $h$ stays at level $\le|\alpha'|+3\le N/2+3$, with $c=0$. The part $\partial_vh$
produces the factor $(|\beta'|+1)/\rad$ on the $v$-high factor.
Now we have $c_\alpha=1$. If $c_{\alpha''}=1$, then
by part~(iv) of the weight algebra the total factor is $\le
C(|\beta'|+1)^{1/2}(|\beta''|+3)^{1/2}/\rad$ (the $v$-Sobolev shift of $X$ is $2$), that is, a
half shift on the $v$-low $X$ and a half power on the $v$-high $h$. If $c_{\alpha''}=0$, then
part~(v) of that lemma implies the bound $\le C(|\beta'|+1)^{1/2}/\rad$. In both cases we put the
half power $(|\beta'|+1)^{1/2}\le(|\beta|+1)^{1/2}$ on the target $u_{\mathbf k} e_j$, which is in
energy, and we absorb it as in~(III), with $B=C\rad^{-1}\eps^{1/2}$ and $\sum_{
\mathbf k} Z_{\mathbf k}^2\le C\mathbb E\sum_{Q'\sim Q}\mathfrak d_{Q'}$. We use the norms of~(T4):
$\partial_vh$ (or $t\partial_yh$) in $L^\infty_yL^2_v$, $X$ in $L^\infty_v$ with the $A$-weight (a
dissipation, whose $y$-level stays at $|\alpha''|$, so $|\alpha''|=N$ is allowed), and the $Y$-side
$u_{\mathbf k} e_j$ with the trace $\le C\eps$ of $\ab[f_0]$, by the case $\Xi=e_j$ of
Lemma~\ref{lem:ULtri}. Apart from the part absorbed by $\theta_{\mathbf k}$, we obtain
$C\rad^{-1}\eps^{1/2}\,\mathbb E^{1/2}\,\mathfrak d^{1/2}e^{1/2}$ over the neighbors, which we
handle by Young's inequality, as in~(I). The terms of
Lemma~\ref{lem:ULaform}(iii) in which $D$ falls on
$\partial^\nu f_0$, namely $\ab\langle\partial^\nu f_0h,\zeta^\sharp+\Delta\zeta\rangle$ and
$\ab\langle\partial^{\nu-\nu_1}f_0h,\partial^{\nu_1}\zeta(z)\rangle$, $0\ne\nu_1\le\nu$, put no
extra derivative on $h$ or on $X$. They are configuration~(IV) of the diffusion, with the
norms of~(T4), with $Y$ replaced by $u_{\mathbf k}\Xi$,
$\Xi\in\{\zeta^\sharp,\Delta\zeta,\partial^{\nu_1}\zeta\}$ (last paragraph of
Lemma~\ref{lem:ULtri}). Here $h$ stays at $y$-level $\le|\alpha'|+2\le N/2+2$, so we have
$\Pi\le2$ by part~(ii) of the weight algebra. The factor $|\zeta^\sharp|\le C\jap{y''}$ stays on
the $Y$-side, whose trace is $\le C\eps$ (Lemma~\ref{lem:FPYside}, $\varphi=1$), and the
$\nu_1$-sum is put into $C_1$ as in Lemma~\ref{lem:FPb}. For these terms we get
$C\eps^{1/2}\mathbb E^{1/2}\mathfrak d^{1/2}e^{1/2}$ over the neighbors, with no half power, which
is within the bound of the proposition.
\end{itemize}

\emph{Summation.} In every case we sum over $\mathbf k$ by Lemma~\ref{lem:ULconv}(ii), with the
extra polynomial factor on the $v$-low side allowed because $s\ge6$, and with the half power on the
target handled as above. The constants depend on $N_y$, $s$, $\delta_y$,
$\rad_0$, $t_1$, $C_\rad$ and the class $\mathfrak D$ only, all fixed before $\eta_*$ and
$\eps_0$. They do not depend on $Q$, $\mathbf k$, $t\le T$, $\eps$, $n$ or $\nu$. This completes
the proof.
\end{proof}

\subsubsection{The other terms and the output}

We check that the other estimates of Subsection~\ref{ssec:layer} hold with the weakened weights, and
we describe the output at $t_1$.

\begin{lem}[The remaining terms with the weakened weights]\label{lem:FPtopother}
Lemmas~\ref{lem:FPa}--\ref{lem:FPg}, Lemma~\ref{lem:FPd} for the $A\zeta$-part,
Proposition~\ref{prop:FPclose} and Proposition~\ref{prop:FPvisc} hold with $W^\flat$ in place of
$W$, with the same proofs, except for the following points: the crossings of
Lemma~\ref{lem:FPccross}, which use the non-integrated form and Lemma~\ref{lem:FPhalf}; the
truncation in Proposition~\ref{prop:FPvisc}(iii) (Remark~\ref{rem:FPtrunc8}); and the $\theta$-term
in Proposition~\ref{prop:FPclose}, which is used here.
\end{lem}

\begin{proof}
Note that the weights $W^\flat_{\mathbf k}$ are constants on each box, and we have $\frac
d{dt}W^{\flat2}_{\mathbf k}=2(\rad'/\rad)|\beta|W^{\flat2}_{\mathbf k}$. Hence
Lemma~\ref{lem:FPa}, the identity~\eqref{eq:FPid} and the neighbor sums of Lemma~\ref{lem:BKnb} are
unchanged. In Lemmas~\ref{lem:FPb}, \ref{lem:FPc}, \ref{lem:FPe} and in the $A\zeta$-part of
Lemma~\ref{lem:FPd} (the $\bb$-part is Proposition~\ref{prop:FPfix}), every factor is measured
either at a level not above that of the target, or, after a shift of a $y$-low factor, at a level
$\le N/2+3<N-2$, or after a $y$-lift by $2$ of a $y$-high factor. In the first two cases we have
$\Pi\le2$ by Lemma~\ref{lem:FPtopalg}(i)--(ii). In the third, we have $\Pi\le2$ or
$\Pi\le(|\beta'|+3)^{1/2}$. On a $v$-low factor this is controlled by~(iii) of the same lemma,
and on a $v$-high factor it is the subject of Lemma~\ref{lem:FPccross}. The factor
$\vartheta(v_0)^{-\iota_\alpha}$
is treated in Remark~\ref{rem:FPiota}, and the Sobolev shifts of $h$ into the $\tau$-levels in
Lemma~\ref{lem:FPcross}. The forcing $\partial^{\mathbf k} S_0[1]$ of Lemma~\ref{lem:FPf} only
gains, while the $G$-terms of $\partial^{\mathbf k} S_0[\tilde G]$ have the $c$-crossing~(e) of
Lemma~\ref{lem:FPccross}(i). In Proposition~\ref{prop:FPvisc} the weights are constants on boxes,
but for the truncation $|\beta|\le n$ we need Remark~\ref{rem:FPtrunc8}. This completes the proof.
\end{proof}

\begin{remark}[The output at $t_1$]\label{rem:FPtopout}
At the levels $|\alpha|\ge N-2$ the weight is $\rad^{|\beta|}(|\beta|+1)^{s-1/2}/\beta!$, up to
$\vartheta(v_0)^{-1}$. With $s\ge6$, this is bounded below by the weight $\rad^n(n+1)^5/n!$ of the
main bootstrap~\eqref{eq:Aenergy}. Hence Proposition~\ref{prop:FPtr} holds.
Lemma~\ref{lem:EPclose} at $\gamma\in[-3,-2]$ uses the same weights with the factor
$\jap{v_0}^\sigma$ replaced by $\omega_Q(t)$ as in Subsection~\ref{ssec:endpoint}. That factor is a
constant on each box, and its time derivative only adds the damping term~\eqref{eq:EPdamp}.
Therefore Proposition~\ref{prop:FPfix} and Lemma~\ref{lem:FPtopother} apply there, with the
exceptions in which a power of $\vartheta(v_0)$ is not absorbed by kernel decay at $\gamma=-2$:
\begin{itemize}
\item the $\iota$-crossing of Remark~\ref{rem:FPiota}, of the form $C\eps\,\ell_Q\sum_{Q'\sim
Q}e_{Q'}$;
\item the transverse factor in the non-integrated form of Lemma~\ref{lem:FPccross}(c), which after
the gain $\vartheta(v_0)^{-1}$ of $\iota$ is a single $\vartheta(v_0)$ on a term with a half power.
It is absorbed by the same part of $\theta_{\mathbf k}$ as in~(iii) there, through
Lemma~\ref{lem:FPhalf}, and it leaves $C\eps\,\ell_Q\sum_{Q'\sim Q}e_{Q'}$, by~(ii) there;
\item the transverse factor $t\vartheta$ of the plain $X$ in~(d) of the same lemma, which the
same gain cancels.
\end{itemize}
Since the targets are in energy, these terms are absorbed by the damping of Lemma~\ref{lem:EPdamp}.
The constant $C^{\rm lay}_{\sharp}$ there includes them, and $K_\flat$ is chosen after it.
\end{remark}

\subsubsection{Sobolev shifts into the $\tau$-levels}\label{sssec:FPcross}

The weight algebra of Lemmas~\ref{lem:FPweight} and~\ref{lem:FPOmega}(ii) leaves one pairing open.
The weight of a coefficient factor $h$ that is measured at a $\tau$-level on its box $P''$ contains
$\tau''^{-1}=h_r''\vartheta''^2$. If the target is also at a $\tau$-level, the factor $h_r''$ is
paired with the $h_r(Q)$ in the $\tau$ of the target. If the target is not, $h_r''$ is unpaired.
Lemma~\ref{lem:ULmom} transfers it with a factor $h_r(Q)$,
which is not bounded uniformly in $y$. This
occurs when a Sobolev embedding in $y$ lifts a $y$-high $h$ into the $\tau$-levels above a target
that is not in them. Absolute estimates are not affected, because there the factor is controlled by
$F_0(y'')\rho_0(v'')$ by Lemma~\ref{lem:BKcomp}(i). We use three facts. The first one is a moment
bound for $\jap{y''}$ that does not grow with $\jap y$, namely Lemma~\ref{lem:RBrad} below. The
second one is the factor $\vartheta(v_0)^{-1}$ in the weights of the levels $\ge N-4$, which we
introduced in Definition~\ref{def:FPtop}. The third one is the bound
$\mathbb E^{1/2}\le2C_4\eps t$ for the coefficient factor, pointwise in $t$, which
Proposition~\ref{prop:FPclose} already proves.

\begin{lem}[Radial moments without $\jap y$]\label{lem:RBrad}
Let $\gamma\in[-3,-2]$, let $f_0\in\mathfrak D$ satisfy (D$_{\rm iso}$), let $g=f_0G$ with
$\frac12\le G\le2$, and let $0<t\le t_1$. Then for all $m\ge0$, $(y,v)$ and $\xi$,
\begin{equation}
\int a(w)\xi\cdot\xi\,\jap{y''}^m\,g(z)\dd w\le C^m\big(\Theta_\sharp(t,v)^m+m^{m/2}\big)A\xi\cdot\xi,\qquad\Theta_\sharp(t,v):=\vartheta(v)+\big(\log(2+t^{-1})\big)^{1/2},
\end{equation}
with $y''=y+tw$ and $z=(y'',v-w)$. The constant $C$ depends on $\gamma,\beta_F,C_\Phi,C_{\rm
pv},p,A_F$ and $t_1$, and not on $y$, $v$, $t$ or $\xi$.
\end{lem}

The bound degenerates as $t\to0$, as it must. Indeed, for $t=0$ we have $y''=y$, and the moment is
$\jap y^mA$. For $\jap y\lesssim\Theta_\sharp$ Lemma~\ref{lem:ULmom} is the better bound.

\begin{proof}
We take $|\xi|=1$, and we write $k(w)=a(w)\xi\cdot\xi=|w|^\gamma|w\times\xi|^2$ and $J=\int
kg=A\xi\cdot\xi$. Let $e_r:=\hat y$ (any unit vector if $y=0$), and let $e_1,e_2$ be an orthonormal
basis of $e_r^\perp$. Then we have
\begin{equation}
\jap{y''}^m\le4^m\big(1+|e_r\cdot y''|^m+|e_1\cdot y''|^m+|e_2\cdot y''|^m\big).
\end{equation}
Since $e_i\cdot y=0$, it follows from Lemma~\ref{lem:RBdir} that the transverse moments are bounded
by $C^m(\vartheta(v)^m+m^{m/2})J$. It remains to bound the radial moment. We write $w=se_r+w^\flat$
with
$w^\flat\perp e_r$, and $\tau:=e_r\cdot y''=|y|+ts$.

\emph{Step 1: a lower bound for $J$.} Let $W:=\{w:|y+tw|\le1\}$, a ball of radius $1/t$ around
$-y/t$. On $W$ we have the following:
\begin{itemize}
\item $|w|\le(|y|+1)/t$, so $|w|^{\gamma+2}\ge\min\{1,t/(|y|+1)\}$, since $\gamma+2\in[-1,0]$;
\item the directions $\hat w$ fill a cone around $-\hat y$ of angle $\gtrsim\jap y^{-1}$, so the
average of $|\Pi(w)\xi|^2=\sin^2\angle(w,\xi)$ over $W$ is $\ge c\jap y^{-2}$;
\item $\jap{v-w}\le\jap v+(|y|+1)/t\le2(1+t^{-1})\jap v\jap y$, so $\rho_0(v-w)\ge
c(1+t^{-1})^{-p}\jap v^{-p}\jap y^{-p}$;
\item $F_0(y'')\ge c$ and $G\ge\frac12$.
\end{itemize}
Since $|W|=ct^{-3}$ and $t\le t_1$, we have $|W|\min\{1,t/(|y|+1)\}=ct^{-3}\min\{1,t/(|y|+1)\}\ge
c\,t_1^{-2}\jap y^{-1}$. Multiplying by the angular factor $\jap y^{-2}$ and by the lower bound for
$\rho_0$, we get $\jap y^{-p-3}$, and thus
\begin{equation}
J\ge J_1:=c\,\eps\,(2+t^{-1})^{-p}\jap y^{-p-3}\jap v^{-p},
\end{equation}
with $c$ depending on $t_1$.

\emph{Step 2: the band.} Suppose $|y|\ge6$, and let $B:=\{w:|\tau-|y||\le|y|/2\}$. On $B$ we have
$|y''|\ge|\tau|\ge|y|/2$. Recall that (D$_{\rm iso}$) bounds $\nabla\Phi_1$ and its derivatives, not
$\Phi_1$ itself, so $\Phi_1$ may grow linearly and $F_{\rm G}(x)\le Ce^{-\beta_F|x|^2+C_\Phi|x|}$.
Together with Lemma~\ref{lem:Ufac}(i), this implies that
$F_0(x,v)\le CK_{\mathcal U}e^{-\beta_F|x|^2+C_\Phi|x|}$. With
$r=|y''|$ we have $e^{-\beta_Fr^2/4+C_\Phi r}\le C$ and $\jap r^me^{-\beta_Fr^2/4}\le C^mm^{m/2}$,
and hence
\begin{equation}
\begin{aligned}
\jap{y''}^mg&\le C\eps\,\rho_0(v-w)\jap{y''}^me^{-\beta_F|y''|^2+C_\Phi|y''|}\le
C\eps\,C^mm^{m/2}e^{-\beta_F|y''|^2/2}\rho_0(v-w)\\
&\le C\eps\,C^mm^{m/2}e^{-\beta_F|y|^2/8}\rho_0(v-w).
\end{aligned}
\end{equation}
Since $\int|w|^{\gamma+2}\rho_0(v-w)\dd w\le C$, the contribution of the band is bounded by
$C^{m+1}m^{m/2}\eps e^{-\beta_F|y|^2/8}$. There are two cases.
\begin{itemize}
\item If $|y|^2\ge K_0\log\big((2+t^{-1})\jap v\big)$ with $K_0:=16(p+1)/\beta_F$, then
$e^{-\beta_F|y|^2/16}\le((2+t^{-1})\jap v)^{-p}$ and $e^{-\beta_F|y|^2/16}\jap y^{p+3}\le C$. Thus
the band contributes $\le C^{m+1}m^{m/2}J_1\le C^{m+1}m^{m/2}J$.
\item Otherwise we have $|y|\le C\Theta_\sharp$. On $B$, $|\tau|\le\frac32|y|\le C\Theta_\sharp$, so
the radial moment over $B$ does not exceed $(C\Theta_\sharp)^mJ$.
\end{itemize}

\emph{Step 3: slicing outside the band.} We follow the proof of the directional moment bound with
$e=e_r$.
We fix $w^\flat$. By (D$_{\rm iso}$) and~\eqref{eq:FGU}, we have $F_0(z)=c_{\rm
G}e^{-\beta_F|y''^\flat|^2}e^{-\beta_F\tau^2}e^{-\Phi_1(y'')}\,\mathcal U(z)$, where $y''^\flat$
does not depend on $s$. Let $T_c:=3$ if $|y|\ge6$ and $T_c:=9$ if $|y|<6$. Let the core be
$I:=\{s:|\tau|\le T_c\}$, an interval of length $2T_c/t$ centered at $s_c=-|y|/t$. A \emph{tail
point} is an $s_1$ with $|\tau_1|>T_c$ that does not lie in the band.

\emph{Tail points are far from $0$.} Suppose first $|y|\ge6$. Outside the band, we have
$|\tau_1-|y||\ge|y|/2\ge3$, and also $|\tau_1-|y||\ge|\tau_1|/3$. This is clear if $\tau_1<0$. If
$0\le\tau_1\le|y|/2$, then $|y|-\tau_1\ge|y|/2\ge\tau_1$, and if $\tau_1\ge\frac32|y|$, then
$\tau_1-|y|\ge\tau_1/3$. If $|y|<6$, then $|\tau_1-|y||\ge|\tau_1|-6\ge|\tau_1|/3$, because
$|\tau_1|>9$. In both cases we obtain
\begin{equation}
|s_1|=\frac{|\tau_1-|y||}t\ge\frac{\max(T_c,|y|,|\tau_1|)}{3t},\qquad\sup_{s\in
I}|s|\le\frac{|y|+T_c}t\le6|s_1| .
\end{equation}
Now we compare each factor of $k\rho_0G\,\mathcal Ue^{-\Phi_1}$ at $s_1$ with its minimum over $I$.
\begin{itemize}
\item \emph{The factor $|w\times\xi|^2=|\mathbf as+\mathbf b|^2$.} For an affine map $p(s)=
\mathbf as+\mathbf b$ and an interval $I$ of length $\ell$ and center $s_c$, we have
$|p(s_1)|\le(1+2|s_1-s_c|/\ell)\max_I|p|$, and $\max_I|p|^2\le C\,\ell^{-1}\int_I|p|^2$ with $C$
absolute. Here $|s_1-s_c|=|\tau_1|/t$ and $\ell\ge6/t$. Thus the value at $s_1$ is $\le
C\jap{\tau_1}^2\,\ell^{-1}\int_I|\mathbf as+\mathbf b|^2\dd s$.
\item \emph{The factor $|w|^\gamma$.} We have $\max_I|w|\le6|s_1|+|w^\flat|\le7|w(s_1)|$. Hence
$|w(s_1)|^\gamma\le7^{|\gamma|}\min_I|w|^\gamma$.
\item \emph{The factor $\rho_0$.} We put $V:=v-w^\flat$, so that $v-w=V-se_r$ and $V_e:=e_r\cdot
V=e_r\cdot v$. For $s\in I$, $\jap{V-se_r}^2\le\jap{V^\flat}^2+2V_e^2+72s_1^2$. If
$|V_e-s_1|\ge|s_1|/2$, then $\jap{V-s_1e_r}^2\ge\jap{V^\flat}^2+s_1^2/4$, and the ratio of the
squares is $\le C(1+V_e^2/s_1^2)\le C(1+t_1^2|v|^2)$, since $|s_1|\ge T_c/(3t)\ge1/t$. Otherwise
$|V_e|>|s_1|/2$, and the ratio is $\le1+290V_e^2$. Therefore, by the two-sided bound
$A_F^{-1}\jap v^{-p}\le\rho_0\le A_F\jap v^{-p}$, the ratio of the values of $\rho_0$ is bounded
by $C\jap v^p$. Here we use that the two sides have the same exponent.
\item \emph{The remaining factors.} We have $G(z_1)\le4\min_IG$, $\mathcal U(z_1)\le K_{
\mathcal U}^2\min_I\mathcal U$, and $e^{-\Phi_1(y''_1)}\le e^{C_\Phi(|\tau_1|+T_c)}\min_Ie^{-\Phi_1}$.
\end{itemize}
Multiplying these bounds, we obtain
\begin{equation}
k(s_1)\le C\jap v^p\jap{\tau_1}^2e^{C_\Phi|\tau_1|}\,\frac t{2T_c}\int_Ik(s)\dd s ,
\end{equation}
where now $k$ includes $\rho_0$, $G$, $\mathcal U$ and $e^{-\Phi_1}$. This is~$(**)$ of
that proof, and its conclusion applies without change. Since $\dd\tau=t\dd
s$, the contribution of the points outside the band with $|\tau|\ge T$ (they are tail points, since
$T\ge T_c$) is bounded by $C\jap v^pC^mm^{m/2}e^{-\beta_FT^2/4}\int_Ik$. Moreover, we have $\int_\R
e^{-\beta_F\tau^2}k\dd s\ge e^{-81\beta_F}\int_Ik$. We choose $T:=\max\{T_0,(4p\beta_F^{-1}\log\jap
v)^{1/2}\}\le C\vartheta(v)$, with $T_0\ge9$, so that $T\ge T_c$ (there
$T_0\ge3$ sufficed). The points with $|\tau|<T$ contribute at most $T^m$ times their part of $J$.
Integrating over $w^\flat$ against the common factor $e^{-\beta_F|y''^\flat|^2}$, we deduce that the
radial moment outside the band is bounded by $C^m(\vartheta(v)^m+m^{m/2})J$ for $m\ge1$. For
$0<m<1$ the bound follows from the case $m=1$ by H\"older's inequality, as at the end of that proof.

For $m\ge1$, we add the band, the slices and the transverse moments, and we obtain the lemma, since
every constant $K\ge1$ satisfies $K\le K^m$, $m^{m/2}\ge1$, $\vartheta\le\Theta_\sharp$ and
$1\le\Theta_\sharp^m$. For $m=0$ the lemma is trivial. For $0<m<1$, let $K_1$ be the constant of the
case $m=1$. By H\"older's inequality for the measure $a(w)\xi\cdot\xi\,g(z)\dd w$, with the
exponents $\frac1m$ and $\frac1{1-m}$, we have $\int a(w)\xi\cdot\xi\jap{y''}^mg(z)\dd w\le
K_1^m(\Theta_\sharp+1)^mA\xi\cdot\xi\le(2K_1)^m\Theta_\sharp^mA\xi\cdot\xi$, since
$\Theta_\sharp\ge1$. This completes the proof.
\end{proof}

\begin{lem}[The crossings]\label{lem:FPcross}
\leavevmode
\begin{enumerate}[label=(\roman*)]
\item In a relative estimate, a coefficient factor $h$ is measured at a $\tau$-level $\ge N-1$
while the target is at a level $\le N-2$ only in the following two cases.
  \begin{enumerate}
  \item[(C1)] Configuration~(III) of Lemma~\ref{lem:FPc}, and the $A\zeta$-part of
Lemma~\ref{lem:FPd}, with $|\alpha'|\in\{N-3,N-2\}$ and $|\alpha|\le N-2$. Here
Lemma~\ref{lem:RBsob} on $P''$ uses the levels $|\alpha'|+2\in\{N-1,N\}$.
  \item[(C2)] The $\ab$-form of~(III) in Proposition~\ref{prop:FPfix}, with
$|\alpha'|\in\{N-4,N-3\}$ and $|\alpha|\le N-2$. Here Lemma~\ref{lem:slab}(ii), applied to
$t\partial_yh$ at level $|\alpha'|+1$, uses the levels $|\alpha'|+3\in\{N-1,N\}$. For
$|\alpha'|=N-3$, the parts in which $h$ is not hit by $t\partial_y$, at $y$-level $|\alpha'|$, also
use the level $|\alpha'|+2=N-1$ through the same lemma: the part $\partial_vh$, and the
terms of Lemma~\ref{lem:ULaform}(iii) in which $D$ falls on $\partial^\nu f_0$ (with the vectors
$\zeta^\sharp+\Delta\zeta$ and $\partial^{\nu_1}\zeta$ on the energy side). They are treated in the
same way. (By Remark~\ref{rem:FPYside}, all of (C2) can also be treated by Lemma~\ref{lem:FPYside}.)
  \end{enumerate}
In both cases $|\alpha|\ge|\alpha'|\ge N-4$, and the factor $X$ is $y$-low, measured at $y$-level
$\le|\alpha|-|\alpha'|+2\le4$.
\item In these cases, the combined factor $\Omega$ of~\eqref{eq:FPOmega0} is the one of
Lemma~\ref{lem:FPOmega}(ii) times $\tau''^{-1}\le C(\jap{y''}+\vartheta'')\vartheta''^2$, and
nothing else, except in the parts of~(C2) in which $D$ falls on $h$ or on $\partial^\nu f_0$. Their
extra unpaired factor $th_r''$, resp.\ $t\jap{y''}$, is put on the energy side together with
$\tau''^{-1}$ by Lemma~\ref{lem:FPYside} (Remark~\ref{rem:FPYside}). Hence
\begin{equation}
\ab[\Omega f_0]\le C_{N_y}\,\delta_y^{|\alpha_\nu|-|\varsigma'_y|}\,K\,\Theta_\sharp(t,v)\,A .
\end{equation}
\item With the weights of Definition~\ref{def:FPtop}, the target carries the factor
$\vartheta(v_0)^{-1}$ and $X$ does not. Thus the factor $\Theta_\sharp$ in~(ii) is replaced by
$C(1+\ell_t)$, where $\ell_t:=(\log(2+t^{-1}))^{1/2}$.
\item Consequently these terms are bounded as in Lemma~\ref{lem:FPc} and
Proposition~\ref{prop:FPfix}, with $C_{\mathbb E}\mathbb E[h]^{1/2}$ replaced by
$CC_2(1+\ell_t)\mathbb E[h]^{1/2}$. Here $\mathbb E[h]$ is the energy of the function of which $h$
is a derivative.
In~(C2), whose bound is that of Lemma~\ref{lem:FPdlow}, the norm $S_G$ of the coefficient factor is
used in the form $S_G\le C(\delta_y,\rad)\mathbb E[h]^{1/2}$ (Lemma~\ref{lem:RBsob}, or
Lemma~\ref{lem:slab}(ii) on the slab), not only $S_G\le C$. Thus the factor $1+\ell_t$ multiplies
$\mathbb E[h]^{1/2}$, and the term is $\le CC_2(1+\ell_t)
\mathbb E[h]^{1/2}\big(\Sigma_Q+\sum_{Q'\sim Q}e_{Q'}\big)$. If $
\mathbb E[h](t)^{1/2}\le2C_4C_{N_y}^{1/2}\eps t$ (the bound on the coefficient factor $\tilde G_\nu$ after
the mollification of Lemma~\ref{lem:FPMnu}(v), and the factor $2C_{N_y}^{1/2}$ is put into $C$),
then this coefficient is bounded by
\begin{equation}
CC_2\big(\mathbb E[h]^{1/2}+C_4\eps\,t\,\ell_t\big)\le CC_2\big(\eta_*+C_4(1+t_1)\eps\big),
\end{equation}
which is absorbed like the other coefficients of $\Sigma_Q$ in Proposition~\ref{prop:FPclose}, for
$c_*$ and $\eps_0$ small.
\end{enumerate}
\end{lem}

\begin{proof}
(i) A Sobolev embedding in $y$ is applied to a $y$-high $h$ only in~(III) of Lemma~\ref{lem:FPc} for
$|\alpha'|\le N-2$, and in the $\ab$-form of~(III) of Proposition~\ref{prop:FPfix}, through
Lemma~\ref{lem:slab}(ii). In configurations~(I) and~(IV), and in the $\ab$-form of~(I), $h$ is
$y$-low, and its levels stay $\le\lfloor N/2\rfloor+3<N-1$. In~(II), in the direct forms, and in
Lemma~\ref{lem:FPe}, the estimates are absolute. There, $h$ is lifted into the $\tau$-levels by
the same slab bound in (II) of the diffusion and of the drift and in the high-$h$ branch of
Lemma~\ref{lem:FPd}, for $|\alpha'|\in\{N-3,N-2\}$, and in the direct form of~(III) with
$|\alpha'|=N-2$. These are \emph{absolute crossings}, that is, the unpaired
$h_r''\le C(\jap{y''}+\vartheta'')$ is absorbed by
$F_0(y'')$ and $\jap{v_0''}^{-\sigma}$ in the Fubini step, as in Lemma~\ref{lem:BKcomp}(i). In~(III)
of the diffusion with $|\alpha'|\ge N-1$, and in the direct form of~(III) of the drift with
$|\alpha'|\ge N-1$, $h$ is not shifted in $y$, and when it is at a
$\tau$-level so is the target. In the remaining cases the levels used are as listed. Note that a
crossing needs the target below $N-1$. The bound on $X$ follows from
$|\alpha''|+|\nu_y|\le|\alpha|-|\alpha'|\le2$, plus the Sobolev shift of $X$ in~(C2).

(ii) The $y$-weights of $P''$ in the proof of Lemma~\ref{lem:FPOmega}(ii) contain $\tau''^{-1}$.
When the target is at a $\tau$-level, we have
$\tau_\alpha/\tau_{\alpha''}=\tau_\alpha=(h_r\vartheta^2)^{-1}$, and $h_r''/h_r$ is one of the
listed factors. Here $\tau_\alpha=\tau_{\alpha''}=1$, so $h_r''\vartheta''^2$ is left over. Since
$h_r''\le C(\jap{y''}+\vartheta'')$, the powers of $\vartheta''$ are controlled by
$\jap{v_0''}^{-\sigma}$, as in that proof. We combine the single factor $\jap{y''}$
with the other factors by H\"older's inequality for the measure $a(w)\xi\cdot\xi f_0(z)\dd w$, and
its moments are bounded by Lemma~\ref{lem:RBrad}. In~(C2) the factor arises inside the slab bound
$\int_S|u|^2\le|S|\sup_{P''}|u|^2$, and the supremum loses $\tau''^{-1}$ through the $\tau$-levels.
After the Cauchy--Schwarz step of~(T3), it enters in the same way.

(iii) By~(i), we have $|\alpha|\ge N-4$, and $X$ is at $y$-level $\le4<N-4$. Therefore the ratio
of the weights of target and $X$ contains $\vartheta(v_0)^{-1}$. Since $|v-v_0|\le2$ on the support
of $\varphi_Q$, we get $\vartheta(v)\le\vartheta(v_0)+\sqrt{\log3}\le C\vartheta(v_0)$, and thus
$\Theta_\sharp/\vartheta(v_0)\le C(1+\ell_t)$.

(iv) We insert~(iii) into the proofs, and then we use $t\ell_t\le C(1+t_1)$ for $t\le t_1$. This
completes the proof.
\end{proof}

\begin{remark}[The factor $\vartheta(v_0)^{-\iota_\alpha}$ in the other terms]\label{rem:FPiota}
The factor is constant on each box and comparable across neighbors. It only makes the weights
smaller, so Lemma~\ref{lem:FPf} only gains. In a term whose factor $X$ is measured at a level not
above the target, the ratio of the weights of target and $X$ is
$\vartheta(v_0)^{-(\iota_\alpha-\iota_{\alpha''})}\le1$. The coefficient factor contributes
$\vartheta(v_0'')^{\iota}$, which is controlled by $\jap{v_0''}^{-\sigma}$, as the other powers of
$\vartheta''$ are. In Step~2 of the transfer, Proposition~\ref{prop:FPtr}, the inverse weight gains
$\vartheta(v_0)$, which is controlled by $\jap{v_0}^{-\sigma}$.

Apart from the lifts of a $y$-low $X$, by at most $3$, or to $y$-level $\le5$ in the non-integrated
form of Lemma~\ref{lem:FPccross}, which stay below $N-4$, the factor $X$ is measured above the
target only in the direct forms with $\alpha''=\alpha$, where the $t\nabla_y$ of the plain $X$ lifts
it from level $|\alpha|$ to $|\alpha|+1$. In configuration~(II) and in the high-$h$ branch, $X$ is
$y$-low, and its pointwise bound, after $t\nabla_y$ and the $y$-Sobolev step, uses levels
$\le\lfloor N/2\rfloor+3<N-4$. This needs $N\ge15$, which holds since $N\ge20$. Hence $\iota$ is not
crossed there. Therefore $\iota$ is crossed in only one place, namely the direct form of~(IV) in
Proposition~\ref{prop:FPfix}, for $|\alpha|=N-5$, with a loss $\vartheta(v_0)$. That estimate is
absolute, with the target in energy.
\begin{itemize}
\item For $\gamma<-2$, the factor $\vartheta(v_0)$ is compensated by the kernel decay
$\int|w|^{\gamma+1}\rho_0(v-w)\dd w\le C\jap v^{\gamma+1}$ and, in the radial term of
Lemma~\ref{lem:RBgain}(ii), by $\int|w|^{\gamma+2}\rho_0(v-w)\dd w\le C\jap v^{\gamma+2}$. Here we
only use $\vartheta(v)\jap v^{\gamma+2}\le C$, which holds for $\gamma<-2$.
\item For $\gamma\in[-3,-2]$ with the energies of Subsection~\ref{ssec:endpoint}, we have
$\vartheta(v_0)\le2\ell_Q$, and the term has the form $C\eps\,\ell_Q\sum_{Q'\sim Q}e_{Q'}$, which is
admitted by Lemma~\ref{lem:EPdamp}.
\end{itemize}
\end{remark}

\begin{remark}[The bound on the coefficient factor, pointwise in $t$]\label{rem:FPpointwise}
Lemma~\ref{lem:FPcross}(iv) uses $\mathbb E[h](t)^{1/2}\le2C_4\eps t$ for each $t$, and not only
$\mathbb E[h]^{1/2}\le\eta_*$. We obtain it by the following two changes.
First, in Definition~\ref{def:FPreg}, the hypothesis on $\tilde G$ becomes $\mathbb E[\tilde
G](t)^{1/2}\le2C_4\eps t$ for all $t\in[0,t_1]$. After the mollification of
Lemma~\ref{lem:FPMnu}(v), $\tilde G_\nu$ satisfies $\mathbb E[\tilde
G_\nu](t)^{1/2}\le2C_4C_{N_y}^{1/2}\eps t$, and this is the hypothesis on the coefficient factor
that we use in Proposition~\ref{prop:FPclose}. Its conclusion is still $
\mathbb E(t)^{1/2}\le2C_4\eps t$, because the constant in the conclusion does not depend on the
constant in the hypothesis, which only enters $\eps_0$. Second, in Lemma~\ref{lem:FPiter}(i), the
bound becomes $\mathbb E[G^j](t)^{1/2}\le2C_4\eps t$ for all $t$.

Both hold by induction. For the frozen problem the coefficient factor is $\tilde G$, so the bound is
a hypothesis on $\tilde G$, not a bootstrap on the unknown. The conclusion of
that proposition, $\mathbb E(t)^{1/2}\le2C_4\eps t$, reproduces it. The crossing terms
change only the coefficient of $\Sigma_Q$, while the $c$-crossing of the source
adds a Gronwall term and a source, and the closing argument keeps the
constant $2C_4$. The mollification multiplies the constant by
$C_{N_y}^{1/2}$, which affects only $\eps_0$, and the limit $G$ satisfies the same bound. In
Subsection~\ref{ssec:endpoint} the bound is $C_5\eps t^{3/4}$, and we have
$t^{3/4}\ell_t\le C(1+t_1)$ as well.
\end{remark}

\subsubsection{Unpaired factors on the energy side}\label{sssec:FPYside}

A second unpaired factor, of the same kind, occurs in the $\ab$-form. In the $\ab$-form
$\bb[F]=\ab\langle DF\rangle$, the derivative $D=\partial_v-t\partial_y$ that falls on the
coefficient is an extra derivative, which does not belong to the target index
$\mathbf k=\nu+\mathbf k'+\mathbf k''$. Thus the weight of $t\partial_yh$ on $P''$ contains one more
factor $th_r''$ (radially) or $t\vartheta''$ (transversally), and the target has nothing to pair it
with.
The same holds for the explicit factor $D\partial^\nu f_0$, which by Lemma~\ref{lem:FPexp} carries
an extra factor $t\jap{y''}$ compared with $\partial^\nu f_0$. The transverse factor $t\vartheta''$
is controlled by $\jap{v_0''}^{-\sigma}$. The radial factor, however, is not covered by the list in
the proof of Lemma~\ref{lem:FPOmega}(ii).

In every $\ab$-form, however, the third factor is the \emph{energy} $u_{\mathbf k}\Xi$ with a
bounded vector $\Xi$, and not a dissipation. Then we estimate the $Y$-side by the trace of $\ab$,
which is an absolute bound, and we can put an unpaired factor entirely on that side.

We use the following lemma for these $\ab$-forms.

\begin{lem}[Asymmetric splitting]\label{lem:FPYside}
Consider a term
\begin{equation}
W_{\mathbf k}^2\iint\ab\big[\partial^\nu f_0\,h\,\varphi\big]X\cdot\Xi\,u_{\mathbf k}\,\varphi_Q^2 ,
\end{equation}
where $h$ and $X$ are as in~\eqref{eq:FPterm} and $|\Xi|\le C\jap{y''}^2$ (the cases $\Xi=e_j$,
$\zeta^\sharp$, $\partial^{\mu_2}\zeta$). Suppose that $\varphi=\varphi(z)\ge0$ satisfies
$\varphi\le C\,t^j(\jap{y''}+\vartheta'')^k$ for some integers $j,k\ge0$. Then the term is bounded
as if $\varphi$ were absent, with the $Y$-side trace $C\eps$ replaced by $C_{k}\eps\,t^{2j}$. In
particular, no moment of $\varphi$ against $A$ is needed.
\end{lem}

\begin{proof}
We split the combined factor as $\Omega=\Omega_X\Omega_Y$. We put all factors of
Lemma~\ref{lem:FPOmega}(ii), the $v$-factor $K$ and $\jap{v_0''}^{-\sigma/2}$ into $\Omega_X$, and
we put $\Omega_Y:=\jap{v_0''}^{-\sigma/2}$. By the Cauchy--Schwarz inequality for the nonnegative
kernel $a(w)f_0(z)\dd w$, we have, at each $(y,v)$,
\begin{equation}
\big|\ab[\Omega\varphi f_0\,\tilde h]X\cdot\Xi\big|\le\big(\ab[\Omega_X^2f_0|\tilde h|^2]X\cdot
X\big)^{1/2}\big(\ab[\Omega_Y^2\varphi^2f_0]\Xi\cdot\Xi\big)^{1/2},
\end{equation}
where $\tilde h$ is $h$ with its weight. The first factor is the one of Lemma~\ref{lem:FPc}
and of that lemma, with $\Omega_X^2$ in place of $\Omega$. The $v$-factor $K$ and the binomial
coefficients do not depend on $z$. They are all in $\Omega_X$, so after the square root they appear
exactly once, as in the symmetric splitting. The square doubles the order of the $\ab$-moments that
are needed, to at most $4N_y+8$. All these factors have bounded moments of
every order, and $\vartheta''^{2m}\jap{v_0''}^{-\sigma}\le C_m$, so only the constant $C_{N_y}$
changes. When we take $h$ in $L^2$ on a slab, as in~(T3), $|\tilde h|^2$ stays inside this factor,
and the slab treatment is unchanged. Indeed, this is the common first step of the proof of
Lemma~\ref{lem:ULtri}, with the weight distributed differently, and the slab, the $y$-Sobolev step
on $X$ and the $v$-Sobolev step on $h$ act only on the $X$-side. In the uses below, $\varphi$ is
either constant on each box $P''$ (as $th_r''$ and $\tau''^{-1}$) or a smooth function of $y''$ (as
$t\jap{y''}$). Consequently, we have $\tilde h=\varphi\cdot(\tilde h/\varphi)$ on each box, and the
slab norm of $\tilde h/\varphi$ is bounded box by box by Lemma~\ref{lem:slab}(ii), and the
constants of neighboring boxes are comparable. For the second factor, we have
\begin{equation}
\ab[\Omega_Y^2\varphi^2f_0]\Xi\cdot\Xi\le\mathrm{tr}\,\ab\big[\jap{v_0''}^{-\sigma}\varphi^2\jap{y''}^4f_0\big]\le
C_k\eps\,t^{2j}.
\end{equation}
Indeed, $\vartheta''^{2k}\jap{v_0''}^{-\sigma}\le C_k$, and the trace of $\ab[\jap{y''}^Cf_0]$ is
$\le C\eps\jap v^{\gamma+2}\le C\eps$ by Lemma~\ref{lem:FPdecay}(i), since the Gaussian factor of
$F_0$ controls $\jap{y''}^C$. This bound is absolute. It does not compare with $A$, so it involves
no loss in $y$ or $v$. It holds for $\gamma\in[-3,-2]$, since we only use $\gamma+2\le0$.

\emph{The gain $\jap{v_0''}^{-\sigma}$ is split.} After the splitting, the $X$-side carries only
$\jap{v_0''}^{-\sigma}$ in place of $\jap{v_0''}^{-2\sigma}$. In the uses listed in
Remark~\ref{rem:FPYside} this causes no loss, because there is no slab factor $L$ (the slabs there
are those used above). The slab factor $L$ of Lemma~\ref{lem:FPslab} occurs only
in~(T3) with $|\alpha'|\ge N-1$, where the target is also at a $\tau$-level and nothing is unpaired.
The lemma could also be combined with Lemma~\ref{lem:FPslab}, since the proof of the latter works
for any positive power, the only input being $\max(1,X_3)^4\jap{v_0''}^{-2\sigma}\le C$, with
$X_3\le\vartheta''$. This completes the proof.
\end{proof}

\begin{remark}[Scope of the lemma]\label{rem:FPYside}
\leavevmode
\begin{itemize}
\item The $t\partial_y$-part of $D\partial^{\mathbf k'}G$ in every $\ab$-form: Lemma~\ref{lem:FPd}
(Case~1 of its proof), Lemma~\ref{lem:FPdlow}, and configurations (I) (via (T1), no slab), (III)
with $|\alpha'|\le N-3$ (via (T3), slab of Lemma~\ref{lem:slab}(ii)) and (IV) with $|\alpha|\ge N-2$
(via (T4), no slab) of Proposition~\ref{prop:FPfix}. Here $\varphi=th_r''$, with $j=k=1$.
\item The $t\partial_y$-part of the explicit factor $D\partial^\nu f_0$ in Lemma~\ref{lem:FPd}, and
in the $\ab$-form of~(IV) of the same proposition. Here $\varphi=t\jap{y''}$, with $j=k=1$.
\item The case (C2) of Lemma~\ref{lem:FPcross}, which is also an $\ab$-form (via (T3), with the
same slab). There the unpaired $\tau''^{-1}=h_r''\vartheta''^2$ can be put on the
$Y$-side in the same way, with $\varphi=h_r''\vartheta''^2$ ($j=0$, $k=3$). For the
$t\partial_yh$-part both unpaired factors go there together, $\varphi=th_r''\cdot
h_r''\vartheta''^2$ ($j=1$, $k=4$). This is simpler than Lemma~\ref{lem:RBrad}.
Parts~(ii)--(iv) of the crossing lemma remain valid for it.
\item The $A\zeta$-part of Lemma~\ref{lem:FPd} is already of this form, since $\zeta^\sharp$ and the
moment $\jap{y''}^2$ sit on the $Y$-side. For the same reason, its case (C1) needs only this lemma.
\end{itemize}
However, the lemma does not apply when the third factor is a dissipation $Du_{\mathbf k}$, namely
in (C1) for the diffusion itself, and in the cutoff variant $Y=u_{\mathbf k} D\varphi_Q/\varphi_Q$,
whose bound
$|D\varphi_Q|_A^2\le C\eps$ is relative. These are the cases of Subsection~\ref{sssec:FPcross} that
need the radial moments and the factor $\vartheta(v_0)^{-\iota_\alpha}$.
\end{remark}

\subsubsection{The half weakening: $v$-high factors across the threshold and the
truncation}\label{sssec:FPhalf}

The weakening in Definition~\ref{def:FPtop} is by a half power only. A full power
$(|\beta|+1)^{-c_\alpha}$ would produce a factor $(|\beta'|+1)$ on a $v$-high and $y$-high
coefficient factor $h$ lifted by Lemma~\ref{lem:slab}(ii) to the levels $\ge N-2$ above the target.
This factor can neither be summed by Lemma~\ref{lem:ULconv} nor absorbed by a $\theta$-term, since
$h$ is evaluated on another box $P''$. On the other hand, we can always move a half power
from the factor that carries it to the target, since $|\beta_j|\le|\beta|+3$.
When the target is taken in energy on its own box, the half power is absorbed by the
$\theta_{\mathbf k}$ of~\eqref{eq:FPid}.

\begin{lem}[A half power on the target]\label{lem:FPhalf}
Let $B\ge0$ and $Z_{\mathbf k}\ge0$, and suppose a sum $T$ of terms satisfies
\begin{equation}
|T|\le\frac B{|Q|}\sum_{\mathbf k}(|\beta|+1)^{1/2}a_{\mathbf k} Z_{\mathbf k} ,
\end{equation}
where $a_{\mathbf k}$ are the energies of the target on its own box $Q$. Then for every
$\delta\in(0,1]$,
\begin{equation}
|T|\le\frac\delta{|Q|}\sum_{\mathbf k}\theta_{\mathbf k} a_{
\mathbf k}^2+\Big(\frac{B^2\rad}{4\delta
C_\rad\eps}+B\Big)\frac1{|Q|}\sum_{\mathbf k} Z_{\mathbf k}^2+B\,e_Q .
\end{equation}
\end{lem}

\begin{proof}
We have $(|\beta|+1)^{1/2}\le1+|\beta|^{1/2}$. Since $\theta_{
\mathbf k}=C_\rad\eps|\beta|/\rad$,
for $\beta\ne0$ we get
\begin{equation}
B|\beta|^{1/2}a_{\mathbf k} Z_{\mathbf k}\le\delta\theta_{\mathbf k}
a_{\mathbf k}^2+\frac{B^2|\beta|}{4\delta\theta_{\mathbf k}}Z_{\mathbf k}^2=\delta\theta_{\mathbf k}
a_{\mathbf k}^2+\frac{B^2\rad}{4\delta C_\rad\eps}Z_{\mathbf k}^2,
\end{equation}
and $Ba_{\mathbf k} Z_{\mathbf k}\le\frac B2(a_{\mathbf k}^2+Z_{\mathbf k}^2)$. This proves the
lemma.
\end{proof}

Note that the lemma involves only the $\theta_{\mathbf k}$ of the target on its own box, and no
$\theta$-energies of neighboring boxes or of $P''$, and no multiplicity count.

\emph{Application.} The terms with a half power on the target have different constants
$B_i$, and they must not be merged under a single $B$. Indeed, with $B=\max_iB_i$, the largest $B$
would multiply the $Z$ of the other terms and could produce an $O(1)$ coefficient of the neighbor
energies. Hence we apply the lemma to each term \emph{separately}, with its own $B_i$, $Z_i$ and a
parameter $\delta_i$, where $\sum_i\delta_i=\frac14$ within each proposition. In
Proposition~\ref{prop:FPfix} there are two such terms, the two corner cases, with
$\delta_i=\frac18$. In Lemma~\ref{lem:FPccross} there are five, namely the cases (a), (b), (d),
and the two terms of the non-integrated form in~(c), whose second term (the part with $D^2$) has its
own constant at $\gamma=-2$. For $\gamma<-2$ this constant is
$B=C(\delta_y,\rad)\eps\,\mathbb E^{1/2}$, as for the other four terms by part~(iii) of that lemma,
and at $\gamma=-2$ it is $B=C\eps\vartheta(v_0)$, with
$\sum_{\mathbf k} Z_{\mathbf k}^2\le C\mathbb E\sum_{Q'\sim Q}e_{Q'}$ as at the end of part~(ii).
Each of the five terms gets $\delta_i=\frac1{20}$. Therefore each of the two
propositions takes $\frac14\sum\theta_{\mathbf k} a_{\mathbf k}^2$, and together they take half of
the $\theta$-term of~\eqref{eq:FPid}. Each coefficient $B_i^2\rad/(4\delta_iC_\rad\eps)$ is
the one of the separate use, up to the factor $\delta_i^{-1}\le20$. In the corner case~(IV), we have
$B=C\rad^{-1}\eps^{1/2}$ and $\sum_{\mathbf k} Z_{\mathbf k}^2\le
C\mathbb E\sum_{Q'\sim Q}\mathfrak d_{Q'}$, so the lemma, with $\delta=\frac18$, produces the
coefficient $\big(2C^2/(C_\rad\rad)+C\rad^{-1}\eps^{1/2}\big)C\,\mathbb E\le C_{
\mathbb E}\mathbb E^{1/2}$ of $\Sigma_Q$, for $\mathbb E^{1/2}\le1$ and $\rad\ge\rad_0/2$. This is
the dependence of $C_{\mathbb E}$ on $C_\rad$ used in Proposition~\ref{prop:FPclose}. The only other
use of $\theta_{\mathbf k}$ is the $c$-crossing of the source (Lemma~\ref{lem:FPf}, and step~(6$'$)
of Lemma~\ref{lem:EPclose}), with $\delta=\frac18$. Thus the total is $\frac58$ of the
$\theta$-term.

\begin{lem}[$c$-crossings of $v$-high factors]\label{lem:FPccross}
\leavevmode
\begin{enumerate}[label=(\roman*)]
\item A $v$-high factor is measured at a level $\ge N-2$ while $c_\alpha=0$ only in the following
situation. The factor is a coefficient factor $h$ that is also $y$-high, with
$|\alpha'|\in\{N-4,N-3\}$ and $|\alpha|\le N-3$. The factor $X$, if present, is $v$-low and $y$-low
(configuration~(II)), and in~(e) there is no factor $X$. The lift is the $y$-Sobolev step of
Lemma~\ref{lem:slab}(ii) on the slab of $h$. It occurs in:
  \begin{enumerate}
  \item[(a)] the direct form of $\bb$: (II) of Proposition~\ref{prop:FPfix}, and the high-$h$
branch of Lemma~\ref{lem:FPd};
  \item[(b)] Lemma~\ref{lem:FPe}, second and third lines of~(T5), when the $G$-factor is $v$-high
and $y$-high;
  \item[(c)] configuration~(II) of Lemma~\ref{lem:FPc}, the diffusion (only this crossing subcase
uses the non-integrated form below; the rest of configuration~(II) keeps the integrated form);
  \item[(d)] configuration~(II) of the $A\zeta$-part of Lemma~\ref{lem:FPd};
  \item[(e)] the source $G$-terms of $\partial^{\mathbf k} S_0[\tilde G]$ (Lemma~\ref{lem:FPf},
step~(6$'$) of Lemma~\ref{lem:EPclose}), with the factor $1$ in place of $X$. They are treated
there.
  \end{enumerate}
\item In (a), (b) and (d) the target is in energy, that is, $u_{\mathbf k}$ or
$u_{\mathbf k}\zeta^\sharp$. In~(c) we write the term in the \emph{non-integrated form}
\begin{equation}
\iint D\cdot\big(\partial^\mu A\,Du_{\mathbf k''}\big)u_{
\mathbf k}\varphi_Q^2=\iint\big(\partial^\mu\bb\big)\cdot
Du_{\mathbf k''}\,u_{\mathbf k}\varphi_Q^2+\iint\partial^\mu A:D^2u_{\mathbf k''}\,u_{
\mathbf k}\varphi_Q^2 ,
\end{equation}
where we use $D\cdot\ab[F]=\bb[F]$. Then the target is in energy as well. The first term is of
type~(a). In the second, $D^2u_{\mathbf k''}$ is plain and pointwise (Sobolev in $(y,v)$ on the
$v$-low, $y$-low factor). Since $|\alpha''|\le|\alpha|-|\alpha'|\le1$ in the crossing subcase, it is
measured at $y$-level $\le|\alpha''|+4\le5<N-4$. Its $v$-shift is $4$, namely two from $D^2$ and two
from $H^2_v\subset L^\infty_v$.

Lemma~\ref{lem:ULconv} still applies, since then $K\lesssim\rad^{-4}(m+4)^4/(m+1)^s$ and
$\sum_nn^2(n+4)^8(n+1)^{-2s}<\infty$ for $2s-10>1$, that is $s\ge6$. No other polynomial factor
falls on $D^2u$, whose $c=0$. The estimate is absolute, with $h$ in $L^2$ on the slab.

Against the weights, each $y$-derivative of $D^2$ produces a factor $th_r$ (radially) or
$t\vartheta$ (transversally). A radial index meets the gain
$a(w)\hat y\cdot\hat y\le|w|^{\gamma+2}|y''|^2/|y|^2$ of
Lemma~\ref{lem:RBgain}(i), with $|a(w)e\cdot e'|\le(a(w)e\cdot e)^{1/2}(a(w)e'\cdot e')^{1/2}$ for
the mixed entries. In the inner region $|y|\lesssim\vartheta$, the gain $|y''|/|y|$ degenerates, and
we use $h_r\le C\vartheta$ instead, as at the end of that lemma. Thus the total factor is
$\le Ct^2(\jap{y''}+\vartheta)^2$.

Each extra $y$-derivative also carries $\delta_y^{-1}$ and each extra $v$-derivative $\rad^{-1}$.
These enter the constant $C(\delta_y,\rad)$, as in Lemma~\ref{lem:FPc}(II). The factor
$\jap{y''}^2$ is controlled by $F_0(y'')$. No cutoff variant arises, since nothing is integrated by
parts.

\emph{The factor $\vartheta^2$.} For $\gamma<-2$ it is compensated by
$\int|w|^{\gamma+2}\rho_0(v-w)\dd
w\le C\jap v^{\gamma+2}$. At $\gamma=-2$ there is no kernel decay, and the term also carries the
half power of~(iii). A term $C\eps\,\vartheta(v_0)^2(|\beta|+1)^{1/2}a_{\mathbf k} Z_{\mathbf k}$
cannot be absorbed by $\theta_{\mathbf k}$ and the damping of Lemma~\ref{lem:EPdamp} together, since
Lemma~\ref{lem:FPhalf} would leave $\ell_Q^2$. Instead, we use the factor $\iota$. In the crossing
subcase we have $|\alpha|\ge|\alpha'|\ge N-4$, so the target carries $\vartheta(v_0)^{-1}$, while
$D^2u_{\mathbf k''}$, at $y$-level $\le5<N-4$, does not. Therefore the ratio of their weights
contains $\vartheta(v_0)^{-1}$, and the transverse factor is reduced to $t^2\vartheta(v_0)$. Then we
have
\begin{equation}
C\eps\,\vartheta\,|\beta|^{1/2}a_{\mathbf k} Z_{\mathbf k}\le\tfrac1{20}\theta_{\mathbf k}
a_{\mathbf k}^2+C\eps\,\frac\rad{C_\rad}\,\vartheta^2Z_{\mathbf k}^2,\qquad
C\eps\,\vartheta\,a_{\mathbf k} Z_{\mathbf k}\le C\eps\big(\ell_Qa_{\mathbf k}^2+Z_{
\mathbf k}^2\big),
\end{equation}
and both right sides have the form $C^{\rm lay}_{\sharp}\eps\,\ell_Q\sum_{Q'\sim Q}e_{Q'}$ (with
$\sum Z^2\le C\mathbb E\sum_{Q'\sim Q}e_{Q'}$, and $\ell_{Q'}\le3\ell_Q$ for neighbors), which
the damping lemma admits. The gain $\vartheta(v_0)^{-1}$ has no other use in these terms, since the
$\tau$-crossings of $h$ are absolute, and the factor $\vartheta(v_0'')^{\iota}$ of $h$ is controlled
by $\jap{v_0''}^{-\sigma}$. In the same way, in~(d) the transverse factor $t\vartheta$ of the plain
$X$ at $\gamma=-2$ is cancelled by this gain, since $X$ is at $y$-level $\le4<N-4$.
\item In cases~(a)--(d) the crossing produces $\Pi\le(|\beta'|+3)^{1/2}\le C(|\beta|+1)^{1/2}$, by
Lemma~\ref{lem:FPtopalg}(iii). Moved to the target, it is absorbed by Lemma~\ref{lem:FPhalf},
applied to each of the five terms separately with $\delta_i=\frac1{20}$, as explained after
it, with $B=C(\delta_y,\rad)\eps\,\mathbb E^{1/2}$ and $\sum_{\mathbf k}
Z_{\mathbf k}^2\le C\sum_{Q'\sim Q}e_{Q'}$. The factor $\mathbb E^{1/2}$ of the coefficient is put
into $B$, so that the term $B\,e_Q$ of that lemma is $\le C\eps\,\mathbb E^{1/2}e_Q$.
With $B=C(\delta_y,\rad)\eps$ it would be a term $C\eps\,e_Q$, which is not of the
form~\eqref{eq:FPccross}. The constant in $B$ depends on $\delta_y$ and $\rad$ through the
Sobolev steps, as in~(ii), and for $\gamma=-2$ we refer to the end of~(ii). The contribution
is
\begin{equation}\label{eq:FPccross}
\le C\eps\,\mathbb E^{1/2}\sum_{Q'\sim Q}e_{Q'}+\frac1{4|Q|}\sum_{\mathbf k}\theta_{\mathbf k}
a_{\mathbf k}^2 .
\end{equation}

\end{enumerate}
\end{lem}

\begin{proof}
(i) A $v$-high factor is either a coefficient factor $h$ or the factor $X$.
\begin{itemize}
\item $X$ is $v$-high only in configurations~(I) and~(III) (and in the corresponding classes of
Lemmas~\ref{lem:FPb} and~\ref{lem:FPe}). In~(I) it is not lifted: it is measured at its own level,
which is not above that of the target, and thus it cannot cross into $c=1$ while $c_\alpha=0$.
In~(III) it is $y$-low, and its lifts stay $\le\lfloor N/2\rfloor+4$.
\item $h$ is $v$-high in configurations~(II) and~(IV). In~(IV), $h$ is $y$-low, and its lifts stay
$\le\lfloor N/2\rfloor+3$.
\end{itemize}
Therefore a crossing needs $h$ $v$-high and $y$-high, hence configuration~(II). There $h$ is lifted
only by Lemma~\ref{lem:slab}(ii), for $|\alpha'|\le N-2$. A crossing needs $|\alpha'|+2\ge
N-2>|\alpha|\ge|\alpha'|$. Configuration~(II) occurs in the four places (a)--(d), and the source
$G$-terms~(e) have the same structure with the factor $1$ in place of $X$, since their only unknown
factor is $h$, measured in $L^2$ on the slab.

(ii), (iii) The non-integrated form is the identity before the integration by parts of
Lemma~\ref{lem:FPc}. We compute the factors as in~(II) of that lemma, applied to
each index of $D^2$. By the proof of Lemma~\ref{lem:FPcross}(i), the $\tau$-crossings of the lift
are absolute. We obtain the bound for $\sum Z^2$ as in~(T2) and~(T5), with the energies of
$h$ (from $
\mathbb E$) and the pointwise bound of $X$. This completes the proof.
\end{proof}

\begin{remark}[The truncation $|\beta|\le n$]\label{rem:FPtrunc8}
In Proposition~\ref{prop:FPvisc}(iii) we first prove the energy bounds for the truncations
$\sum_{|\beta|\le n}$. The truncation applies to the unknown, not to the coefficient $\tilde G$. In
the direct form of~(III) of Proposition~\ref{prop:FPfix}, the plain $X$ needs
$\partial_vu_{\mathbf k''}$ at $v$-level $|\beta''|+1$. Since $|\beta|\le n$ and $\beta''\le\beta$,
this exceeds $n$ only when $|\beta''|=n$, hence $\beta'=0$ and $\nu_v=0$. For this term we use the
viscous dissipation instead of the energy at level $n+1$. Since $M_\nu$ has the block $\nu I_v$, the
principal viscous term contains
\begin{equation}
V_\nu(\mathbf k''):=\nu\,(W^\flat_{\mathbf k''})^2\|\partial_vu_{\mathbf k''}\varphi_{Q'}\|^2,
\end{equation}
on $Q$ and its neighbors $Q'$. With $\beta'=0$, $\nu_v=0$ and $|\beta''|=|\beta|=n$, the half
weakening implies $\Pi\le((|\beta'|+3)/(|\beta|+1))^{1/2}=(3/(n+1))^{1/2}$, the $3$ coming from the
$v$-shift $\le2$ of $h$, and we put the factor $\sqrt3$ into $C$. The cancellation between the loss
and the ratio of the weights is componentwise. For $\partial_{v_i}$ the loss $(\beta''_i+1)/\rad$
meets $W^\flat_{\mathbf k''+e_i}/W^\flat_{\mathbf k''}\simeq\rad/(\beta''_i+1)$ ($X$ is at a
$y$-low level, where $c=0$ and $W^\flat=W$), and the viscous dissipation replaces the energy at
$\mathbf k''+e_i$ by $\nu^{-1/2}V_\nu^{1/2}$. We take the plain $X$ with $\sup_y$ from a
$y$-Sobolev step of Lemma~\ref{lem:RBsob} at the scale of the boxes, so we need $V_\nu$ at the
indices $\mathbf k''+\varsigma_y$, $|\varsigma_y|\le2$ (still at $v$-level $n$, and
$|\alpha''|+2\le4\le N_y$), on $Q$ and its neighbors $Q'$, each with the weights and the
normalization of its own box. The factor $\delta_y^{-2}$ of the Sobolev step and the ratios of
weights and volumes of neighbors, bounded by Lemma~\ref{lem:BKnb}(i), are put into $C$, and
$V_\nu(Q')$ below denotes
the sum of these pieces of $|Q|\,\mathfrak d^\nu_{Q'}$. Therefore, for every $\delta>0$, the term is
bounded by
\begin{equation}\label{eq:visc-split}
C\eps\,\mathbb E[\tilde G]^{1/2}\,(n+1)^{-1/2}\,\nu^{-1/2}\Big(\sum_{Q'\sim
Q}V_\nu(Q')\Big)^{1/2}a_{\mathbf k}\le\delta\sum_{Q'\sim Q}V_\nu(Q')+C\delta^{-1}\eps^2\,
\mathbb E[\tilde G]\,\frac{a_{\mathbf k}^2}{(n+1)\nu}.
\end{equation}
\emph{The viscous dissipation.} The first term is absorbed by the viscous dissipations $\sum_{Q'\sim
Q}\mathfrak d^\nu_{Q'}$ of the neighboring boxes. We cannot use the Landau dissipation here, since
all of it has been used in Proposition~\ref{prop:FPclose} ($\eta=4/(32N_0)$). In
Proposition~\ref{prop:FPvisc}(iii) the principal viscous term gives $-\mathfrak d^\nu_Q$, and the
cutoff term does not exceed $\frac14\mathfrak d^\nu_Q$. Thus the viscous dissipation enters the box
inequality with the coefficient $c_\nu=\frac34$. The commutators of $M_\nu$ are bounded by
$\frac1{32N_0}\sum_{Q'\sim Q}\mathfrak d^\nu_{Q'}$.

We take $\delta\le1/(24C_NN_0)$, where $C_N$ bounds the number of triples
$(\mathbf k,\nu,\varsigma_y)$
with the same $\mathbf k''+\varsigma_y$, so that after the summation over $\mathbf k$ the total
coefficient of $\sum_{Q'\sim Q}\mathfrak d^\nu_{Q'}$ is $\le1/(24N_0)$. Then the total coefficient
$\eta_\nu$ of the neighbor viscous dissipations satisfies $\eta_\nu
N_0\le\frac1{32}+\frac1{24}=\frac7{96}<\frac1{12}$, which is $\le c_\nu/4$ and also meets the
threshold of Lemma~\ref{lem:EPdamp} used at $\gamma=-2$ (step~(8) of Lemma~\ref{lem:EPclose}). This
only enlarges the constant $C\delta^{-1}$ of the last term.

After the multiplication by $2$ in the closing argument,
we apply Lemma~\ref{lem:BKnb}(ii) to $\mathfrak d+\mathfrak d^\nu$, with $c=\min(1,2c_\nu)=1$ and
$2\max(\eta_L,\eta_\nu)$ in place of $\eta$. Both satisfy $2\eta N_0\le\frac14=c/4$. The remainder
$2C\nu\,\mathbb E\le2C\eps\,\mathbb E$ is included in $2a\,\mathbb E$. The first term
of~\eqref{eq:visc-split} is not absorbed by the viscous term of $Q$ alone.

For $n+1\ge\nu^{-1}$ the last term of the same bound is $\le C\delta^{-1}\eps^2
\mathbb E[\tilde G]a_{\mathbf k}^2$, uniformly in $n$ and $\nu$. After the summation over
$\mathbf k$, we obtain the Gronwall term
$C\delta^{-1}\eps^2\mathbb E[\tilde G]\,e_Q\le C\eps^4t_1^2\,\mathbb E$, which is
included in $2a\,\mathbb E$. Hence for fixed $\nu$ the bounds hold uniformly in
$n\ge\max(\nu^{-1},7)$, where, as in the proof of the frozen-problem proposition, the condition
$n\ge7$ keeps the $v$-Sobolev shifts of the $v$-low factors of the unknown inside the truncation.
The bounds hold for smaller $n$ because the truncated energies and dissipations are nondecreasing
in $n$. Then we let $n\to\infty$ by monotone convergence, and after that $\nu\to0$.

In the $\ab$-form the extra derivative falls on $\tilde G$, whose energy is not truncated.
\end{remark}

\subsection{The endpoint \texorpdfstring{$\gamma=-2$}{gamma=-2}}\label{ssec:endpoint}

At $\gamma=-2$ the kernels are
\begin{equation}
a(w)=\Pi(w),\qquad b(w)=-\frac{2w}{|w|^2},\qquad c(w)=-\frac2{|w|^2},
\end{equation}
all locally integrable, and the clocks are $\kappa(t)=t^{-1}$, $\tau(t)=\log t$ and
$s_\xi(t)=\lambda_{\rm ref}(\xi)\log(t/T_0)$. The diffusion kernel does not decay at infinity, so
every bound $\int|w|^{\gamma+2}\rho_0(v-w)\dd w\le C\jap v^{\gamma+2}$ becomes $\le C$. This needs
$p>3$, which holds since $p\ge4$ in Definition~\ref{def:data}.

Now we prove Theorem~\ref{thm:main} at $\gamma=-2$, outside the initial layer in
Subsection~\ref{sssec:EPout}, and in the initial layer in
Subsections~\ref{sssec:EPloc}--\ref{sssec:EPE1}.

\subsubsection{Outside the initial layer}\label{sssec:EPout}

Outside the initial layer, the proof depends on $\gamma$ only through the signs $\gamma<0$ and
$\gamma+2\le0$ and the range $\gamma>-\frac72$, the local integrability of $|w|^{\gamma+2}$,
$|w|^{\gamma+1}$ and $|w|^\gamma$, the homogeneity degrees $2+\gamma$, $3+\gamma$ and $5+\gamma$ in
the scaling identities, the clocks $\kappa$ and $\tau$, the near-Maxwellian input (E1) and the $L^1$
rate of Proposition~\ref{prop:Lrate}, and the time regularity of the coefficients. At $\gamma=-2$
the coefficients are not H\"older continuous in $t$ uniformly in $(y,v)$, so the parabolic theory
of Appendix~\ref{app:visc} uses coefficients H\"older continuous in the space variables only.
At the following places these properties are used in a form that has to be checked at
$\gamma=-2$, and we note that each of them remains valid there:
\begin{itemize}
\item Step~1 of Lemma~\ref{lem:RGbar} uses $\gamma+2\le0$, $(\frac32)^{\gamma+2}\ge\frac23$,
$|z|^{\gamma+2}\ge|z|^{-1}$ for $|z|\ge1$, $2^{|\gamma|}\le8$, $2(\gamma+3)\le2$, and
$2(\gamma+3)\int_{|u|<1}|u|^\gamma\dd u=8\pi$ for $\gamma>-3$, and all of these hold at $\gamma=-2$;
\item Lemma~\ref{lem:RFtools} uses $2\gamma+7>0$ and $-2\gamma-3\ge1$, and Lemma~\ref{lem:RGremF}
and Proposition~\ref{prop:Rest} use $p\in[-2,1]$, with $p=\gamma+2,\gamma+3,\gamma+1$, that is,
$p=0,1,-1$ at $\gamma=-2$;
\item the norms $Y_2(m)$ use $\jap{\hat z}^{\gamma/2+1}\le1$;
\item Proposition~\ref{prop:RGcompF}, Theorem~\ref{thm:RGF}, the bootstrap of
Section~\ref{sec:boot} and the continuation argument of Subsection~\ref{ssec:cellsproofs} use
$\gamma$ only through $\kappa$, $\tau$, and $\kappa\lambda_\xi\le C$ for $t\ge T_0$, apart from the
uses discussed in the other items of this list and in (a)--(c) below: the Gaussian hypodissipativity
and trilinear bounds, the norms $Y_2$, requirement~\ref{rq:Wrate}, the meaning of
$T_0^{\kappa_\gamma}$ and the Schauder estimates of the viscous approximation;
\item Lemma~\ref{lem:Abarrier} holds with $\jap v^{\gamma+2}=1$, since part~(i) uses (GS3) in the
form \cite[Thm.~3.7]{Si17} and $\frac3{5+\gamma}=1<\frac52<3=\frac3{3+\gamma}$, and part~(ii) uses
$k_1\ge k_\flat(-2)=7$ for the H\"older exponent $q'\in(3,k_1]$;
\item In Lemma~\ref{lem:LWlin}(i), the coefficients are not H\"older continuous in $t$ at
$\gamma=-2$, so the viscous approximation and Lemma~\ref{lem:CUtstar} use the Schauder estimates
with coefficients
H\"older continuous in space only, \cite[Thm.~1]{Kn80} and \cite[Thm.~4.1]{Lor00} ((E4) of
Inputs~\ref{hyp:E}). The spatial H\"older bound of that lemma needs $\frac78p>3$, and
the bound $\int|w|^{\gamma+2}\rho_0(v-w)\dd w\le C$ needs $p>3$, both true since $p\ge4$.
\end{itemize}
At $\gamma=-2$ we modify the argument as follows.
\begin{enumerate}[label=(\alph*)]
\item The near-Maxwellian estimates of \cite{CM} are published only for $\gamma<-2$. At $\gamma=-2$
the estimates that we use are proved in Appendix~\ref{app:CM}, with the decay in a small space
replaced by the spectral gap. Proposition~\ref{prop:Lrate} covers
$\gamma\in[-3,-2]$.
\item Requirement~\ref{rq:Wrate} contains $\min\{1,|\gamma|-2\}$,
which vanishes at $\gamma=-2$. It is replaced by Lemma~\ref{lem:EPreq}.
\item On $[t_1,T_0]$ the effective time is $\int_{t_1}^{T_0}\kappa=\log(T_0/t_1)\le\log T_0$
(recall $t_1\ge1$). Throughout the paper, and in particular in
Theorem~\ref{thm:main}(iv), $T_0^{\kappa_\gamma}$ is read as $1+\log T_0$ at $\gamma=-2$. Then
$T_0^{\kappa_\gamma}\ge1$ and $T_0^{\kappa_\gamma}=1+\tau(T_0)$, as at $\gamma=-3$.
\end{enumerate}

\begin{lem}[{Requirement~\ref{rq:Wrate} at $\gamma=-2$}]\label{lem:EPreq}
Let $\gamma=-2$, and assume (U1), so that (F3) holds and $\lambda_\xi\ge c_\lambda\lambda_{\rm
ref}(\xi)$. For $t>T_0$, we have $t/\log(t/T_0)\ge T_0$. Hence
$\kappa\lambda_\xi\ge c_\lambda S_*T_0\,t^{-2}$ wherever $s_\xi(t)\ge S_*$, and the condition
$t^{-2}\ll\kappa\lambda_\xi$ of Remark~\ref{rem:Bbandrange} holds for $s_\xi\ge S_b$ once
\begin{equation}
T_0\ge\frac{4C_{\rm r}}{\eta c_Wc_\lambda S_*}.
\end{equation}
\end{lem}

\begin{proof}
We put $r=t/T_0>1$. If $r\le e$, then $\log r\le1$, so $t/\log r\ge t\ge T_0$. If $r>e$, then
$r/\log r\ge e$, since $r\mapsto r/\log r$ is increasing on $[e,\infty)$, and thus
$t/\log r=T_0\,r/\log r\ge T_0$. As in Remark~\ref{rem:Bbandrange}, it follows from
$\lambda_\xi\ge c_\lambda\lambda_{\rm ref}(\xi)$ that
$\kappa\lambda_\xi\ge c_\lambda\dot s_\xi=c_\lambda s_\xi\kappa/(\tau(t)-\tau(T_0))$, and here
$\kappa/(\tau(t)-\tau(T_0))=1/(t\log r)\ge T_0t^{-2}$. Since $S_b\ge S_*$, this is the bound used
there, with $\min\{1,|\gamma|-2\}$ replaced by $1$. This completes the proof.
\end{proof}

\subsubsection{The initial layer: localization and transit time}\label{sssec:EPloc}

In the rest of this subsection we use the notation of
Subsections~\ref{ssec:tools}--\ref{ssec:drift}. The following places there use $\gamma<-2$.
\begin{enumerate}[label=(P\arabic*)]
\item The cutoff term of Lemma~\ref{lem:RBloc}, which is $\le C\eps$ only because of the factor
$\jap v^{\gamma+2}$. It is used in Lemmas~\ref{lem:BKtop} and~\ref{lem:FPa}, and in the cutoff
variants of the commutator terms.
\item The transverse factor $t\vartheta$ of a plain factor $X$, compensated by
$\jap v^{(\gamma+2)/2}$,
in class~(II) of Lemma~\ref{lem:FPc}, in its cutoff variant, and in class~(II) of the $A\zeta$-part
of Lemma~\ref{lem:FPd}.
\item The pointwise decay $|S_0|\le C\eps\jap v^{-2\sigma}$, $\sigma=|\gamma+2|/2$, of
Lemma~\ref{lem:FPdecay}(i). It is used for the forcing in Lemma~\ref{lem:FPf}, for the barrier and
the local estimates of Proposition~\ref{prop:FPvisc}(ii), and for the contraction in
Lemma~\ref{lem:FPiter}(ii). It also justifies the velocity weight $\jap{v_0}^\sigma$, used in
Lemmas~\ref{lem:FPOmega}(ii), \ref{lem:FPslab} and in Proposition~\ref{prop:FPtr}. At $\gamma=-2$
this weight is trivial.
\end{enumerate}
At $\gamma=-2$, (P1) and (P2) leave a factor $\vartheta(v)^2$, and (P3) fails, since $|A|$ and
$|S_0|$ are of size $\eps$ at the points with $|y+tv|\lesssim1$, for every $v$. These are the
particles that
cross the cloud at time $t$. A fast particle feels the cloud with full strength while it is inside
it, but it is inside only for a time $\simeq1/|v|$. We replace the kernel decay by this transit
time, and the velocity weight by a weight with logarithmic damping. The other steps of
these subsections use $\gamma$ only through $\gamma+2\le0$, $\gamma<0$
and the local integrability of the kernels, with the exceptions in the drift term at the top
$y$-levels that are listed in Remarks~\ref{rem:FPiota} and~\ref{rem:FPtopout}.

\begin{lem}[Physical localization of the coefficients]\label{lem:EPloc}
Let $\gamma\in[-3,-2]$, $t\in(0,t_1]$, $x:=y+tv$, and let $G$ be measurable with $|G|\le2$. For
every $m\ge0$ there is $C_m=C_m(\Dc,t_1)$ such that, with $y''=y+tw$,
\begin{equation}
\big|\ab[\jap{y''}^mf_0G]\big|+\big|\bb[\jap{y''}^mf_0G]\big|+\big|\cb[\jap{y''}^mf_0G]\big|\le
C_m\,\eps\jap x^{-p}\qquad\text{at }(y,v);
\end{equation}
for $\gamma=-3$, $\cb[h]$ is replaced by $8\pi h$. Consequently, we have
\begin{equation}\label{eq:EPloc-S0}
|\partial^{\mathbf k} S_0[1]|\le C_{N_y}C^{|\beta|}\beta!R^{-|\beta|}\eps\jap x^{-p}
\end{equation}
for $\mathbf k=(\alpha,\beta)$, $|\alpha|\le N_y$ (and for every $\alpha$, with $C_{N_y}$ replaced
by a constant depending on $|\alpha|$, by the all-orders form of Lemma~\ref{lem:FPexp}, and this bound is used
in Proposition~\ref{prop:EPconstr}(ii)), with the radius $R=\Lambda/2$ of Section~\ref{sec:layer}.
\end{lem}

\begin{proof}
We put $v''=v-w$, so that $y''=x-tv''$. Up to the factors $\eps$ and $|G|\le2$, the integrand of
$\ab$ is bounded by $|v-v''|^{\gamma+2}\jap{x-tv''}^m\rho_0(v'')F_0(x-tv'')$. We split the
$v''$-integral.
\begin{itemize}
\item \emph{The part $|tv''|\ge|x|/2$.} If $|x|\le1$, the whole integral is $\le C$, since
$\jap\cdot^mF_0$ is bounded and $\int|w|^{\gamma+2}\rho_0(v-w)\dd w\le C$. If $|x|\ge1$, then
$|v''|\ge|x|/(2t)$, so $\rho_0(v'')\le A_F(4t)^p\jap x^{-p}$. Moreover, with $u=x-tv''$, we have
\begin{equation}
\int|v-v''|^{\gamma+2}\jap{x-tv''}^mF_0(x-tv'')\dd v''=t^{-5-\gamma}\int|tv-x+u|^{\gamma+2}\jap
u^mF_0(u)\dd u\le Ct^{-5-\gamma},
\end{equation}
since $\gamma+2\in[-1,0]$ and $\jap u^mF_0$ is bounded and integrable. Thus the product is $\le
Ct^{p-5-\gamma}\jap x^{-p}\le C\jap x^{-p}$, since $p-5-\gamma\ge p-3\ge1$ and $t\le t_1$.
\item \emph{The part $|tv''|<|x|/2$.} There we have $|x-tv''|>|x|/2$, so
$\jap{x-tv''}^mF_0(x-tv'')\le C_me^{-\beta_F|x|^2/8}$, and
$\int|w|^{\gamma+2}\rho_0(v-w)\dd w\le C$.
\end{itemize}
We treat the kernels $|b(w)|=2|w|^{\gamma+1}$ and $|c(w)|=2(\gamma+3)|w|^\gamma$ ($\gamma>-3$) in
the same way. The $u$-integrals are $\le Ct^{-4-\gamma}$ and $\le Ct^{-3-\gamma}$, which are again
compensated by $t^p$. For $\gamma=-3$, we have $|f_0G|(y,v)\le2\eps\rho_0(v)F_0(x-tv)$, and the same
two cases yield the bound $\le C\eps\jap x^{-p}$.

Now we come to the last claim. Note that $\partial^{\mathbf k} S_0[1]$ is a finite sum, over Leibniz
splittings, of terms
$\ab[\partial^\nu f_0]$, $\bb[\partial^\nu f_0]$ and $\cb[\partial^\nu f_0]$, multiplied by
derivatives of $\zeta$ and $\zeta^\sharp$. By Lemma~\ref{lem:FPexp}, we have $|\partial^\nu f_0|\le
C_{N_y}C^{|\nu_v|}\nu_v!R^{-|\nu_v|}\jap{y''}^{|\nu_y|}f_0$. The factor $\zeta^\sharp$ contributes
at most $\jap{y''}^2$, and the derivatives of $\zeta$ have analytic bounds by
Lemma~\ref{lem:FPdecay}(ii). Applying the first part with $m=|\alpha|+2$, we obtain the last
claim. This completes the proof.
\end{proof}

The same proof applies when $G$ is replaced by a derivative $\partial^{\mathbf k'}\tilde G$ of a
bounded smooth function, with $C_m$ multiplied by $\sup|\partial^{\mathbf k'}\tilde G|$.

\begin{lem}[Transit time]\label{lem:EPtransit}
Let $p'>1$. The constants below depend on $p'$ and $t_1$.
\begin{enumerate}[label=(\roman*)]
\item For every $(y,v)$ and $t\le t_1$, $\int_0^t\jap{y+sv}^{-p'}\dd s\le C\min\big(t,\jap
v^{-1}\big)$.
\item Let $Q$ be a box of Lemma~\ref{lem:BKcov}, with center $y_Q$ and velocity center $v_0$, and
put $x_Q(s):=y_Q+sv_0$. For $(y,v)\in\tilde Q\times B(v_0,2)$ and $s\le t_1$,
$\jap{y+sv}\ge\jap{x_Q(s)}/(2+2t_1)$, and $\int_0^t\jap{x_Q(s)}^{-p'}\dd s\le
C\min\big(t,\jap{v_0}^{-1}\big)$ for $t\le t_1$.
\item Let $|Q|$ be the $y$-volume of $Q$. The function
\begin{equation}
\Psi_Q(t):=\Big(|Q|^{-1}\int_{\tilde Q\times B(v_0,2)}\jap{y+tv}^{-2p}\dd y\dd v\Big)^{1/2}
\end{equation}
satisfies $\int_0^t\Psi_Q(s)\dd s\le C\,t^{3/4}\jap{v_0}^{-1/4}$ for $t\le t_1$.
\end{enumerate}
\end{lem}

\begin{proof}
(i) The bound by $t$ is trivial, and for $|v|\le1$ it implies the claim, since $t\le
t_1\le\sqrt2\,t_1\jap v^{-1}$. For $|v|\ge1$ we write $y=a+b\hat v$ with $a\perp\hat v$. Then
$\jap{y+sv}\ge\jap{b+s|v|}$, and $\int_\R\jap{b+s|v|}^{-p'}\dd s=|v|^{-1}\int_\R\jap r^{-p'}\dd r$.

(ii) By Lemma~\ref{lem:BKcov}(i), all widths of $Q$ are $\le c_{\rm box}$, so we have
$\operatorname{diam}\tilde Q\le4\sqrt3\,c_{\rm box}\le1$ for $c_{\rm box}\le\frac18$. Hence
$|y+sv-x_Q(s)|\le|y-y_Q|+s|v-v_0|\le1+2t_1$, and the first claim follows from
$\jap{a+d}\le(1+|d|)\jap a$. The second claim is~(i) at $(y_Q,v_0)$.

(iii) By the Cauchy--Schwarz inequality in $s$, Fubini's theorem and~(i) with $p'=2p$, we have
\begin{equation}
\int_0^t\Psi_Q\le t^{1/2}\Big(|Q|^{-1}\int_{\tilde Q\times B(v_0,2)}\int_0^t\jap{y+sv}^{-2p}\dd
s\dd y\dd v\Big)^{1/2}\le Ct^{1/2}\min\big(t,3\jap{v_0}^{-1}\big)^{1/2},
\end{equation}
since $|\tilde Q|\le C|Q|$ and $\jap v\ge\jap{v_0}/3$ for $|v-v_0|\le2$. Finally, we use
$\min(t,a)^{1/2}\le t^{1/4}a^{1/4}$. This completes the proof.
\end{proof}

\subsubsection{The initial layer: energies with logarithmic damping}\label{sssec:EPen}

\begin{defin}[The energies for $\gamma\in{[-3,-2]}$]\label{def:EPen}
Fix $\sigma_\flat:=\frac14$ and a constant $K_\flat\ge1$, which is chosen in
Lemma~\ref{lem:EPclose}. In the weights $W_{\mathbf k}$ of Subsection~\ref{sssec:FPset}, we replace
the factor $\jap{v_0}^\sigma$ by
\begin{equation}
\omega_Q(t):=\jap{v_0}^{\sigma_\flat}e^{-K_\flat\eps
t\,\ell_Q},\qquad\ell_Q:=\vartheta_\natural(v_0)^2:=1+\log\jap{v_0}.
\end{equation}
We write $a^\flat_{\mathbf k}$, $d^\flat_{\mathbf k}$, $e^\flat_Q$, $\mathfrak d^\flat_Q$ and $
\mathbb E^\flat:=\sup_Qe^\flat_Q$ for the resulting quantities. The weakening factors of
Subsection~\ref{ssec:drift} at the top $y$-levels are kept. We always assume
$K_\flat\eps_0t_1\le\frac18$.
\end{defin}

Since $\omega_Q=e^{-K_\flat\eps t}\jap{v_0}^{\sigma_\flat-K_\flat\eps t}$, we have
\begin{equation}\label{eq:EPomega}
e^{-1}\jap{v_0}^{1/8}\le\omega_Q(t)\le\jap{v_0}^{1/4}\qquad(t\le t_1).
\end{equation}
The factor is constant on each box. For neighbors, it follows from $|v_0-v_0'|\le2$ that
$\jap{v_0'}/\jap{v_0}\in[\frac13,3]$, hence
$\omega_{Q'}/\omega_Q\in[3^{-1/4},3^{1/4}]\subset[\frac12,2]$ and $\ell_{Q'}\le\ell_Q+2\le3\ell_Q$.
Therefore Lemma~\ref{lem:BKnb}(i) holds unchanged. On $B(v_0,2)$ we have
$\vartheta(v)^2\le2(1+\log\jap v)\le5\ell_Q$. Since $\frac
d{dt}\omega_Q^2=-2K_\flat\eps\ell_Q\omega_Q^2$, the identity~\eqref{eq:FPid} acquires the
\emph{damping term}
\begin{equation}\label{eq:EPdamp}
-K_\flat\eps\,\ell_Q\sum_{\mathbf k}\big(a^\flat_{\mathbf k}\big)^2\qquad\text{on its right side.}
\end{equation}
Lemma~\ref{lem:EPdamp} uses this term.

\begin{lem}[The closed estimate for $\gamma\in{[-3,-2]}$]\label{lem:EPclose}
Let $\gamma\in[-3,-2]$. Fix the constants of Proposition~\ref{prop:FPclose} up to $\eta_*$, with
$1/(72N_0)$ in place of $1/(32N_0)$ in its proof (and $1/(144N_0)$ in place of $1/(64N_0)$ in the
choice of $\rad_*$ and $\delta_*$ there and in Proposition~\ref{prop:FPvisc}(iii)), then the
constant $C^{\rm lay}_\sharp$ of the proof below, which does not depend on $Q$, then
$K_\flat:=9C^{\rm lay}_\sharp N_0$, then $C_6,C_7,C_8$ (which do not depend on $K_\flat$ or $\tilde
C$), then $C_5:=96C_6+1$, and then $\eps_0$. Let $\eps\le\eps_0$, and let $G$ be a solution on
$[0,T]$,
$T\le t_1$, of the frozen problem of Definition~\ref{def:FPreg}, with or without viscosity, whose
coefficient factor $\tilde G$ satisfies $\mathbb E^\flat[\tilde G](t)^{1/2}\le\tilde C\eps t^{3/4}$
for all $t$. With viscosity $\nu>0$, the energies may be replaced by their truncations $
\mathbb E^{\flat,(n)}$ to $|\beta|\le n$, $n\ge\max(\nu^{-1},7)$, as in Remark~\ref{rem:FPtrunc8}.
The conclusion then holds uniformly in $n$, and for $\mathbb E^\flat$ by monotone convergence.
Assume
that $\sup_{Q,\,s\le t}\ell_Qe^\flat_Q(s)$ and $\sup_Q\int_0^t\mathfrak d^\flat_Q$ (with the viscous
dissipation added, if there is viscosity) are finite for $t\le T$, that $\mathbb E^\flat$ is
continuous, and $\mathbb E^{\flat\,1/2}\le\eta_*$ and $\frac12\le G\le2$ on $[0,T]$. Then, on
$[0,T]$,
\begin{equation}
\mathbb E^\flat(t)^{1/2}\le C_5\eps t^{3/4},\qquad|G-1|\le C\eps t^{3/4}\jap v^{-1/8},
\end{equation}
where $C_5$ does not depend on $\tilde C$, while the constant $\eps_0$ depends on $\tilde C$.
\end{lem}

\begin{proof}
We go through the proof of Proposition~\ref{prop:FPclose} with the weights of
Definition~\ref{def:EPen}. We write $\Sigma_Q=\sum_{Q'\sim Q}\mathfrak d^\flat_{Q'}$ and $
\mathcal E_Q:=\sum_{Q'\sim Q}e^\flat_{Q'}$.

\emph{(1) The weight algebra.} In the combined factor~\eqref{eq:FPOmega0}, the factors $\omega$ of
the target and of $X$ cancel up to a factor in $[\frac12,2]$. The coefficient box contributes
$\omega_{P''}^{-1}=e^{K_\flat\eps t}\jap{v_0''}^{-\sigma_\flat+K_\flat\eps t}\le
e\jap{v_0''}^{-1/8}$. The proofs of Lemma~\ref{lem:FPOmega}(ii) and Lemma~\ref{lem:FPslab} rest on
$(\vartheta'')^m\jap{v_0''}^{-\sigma}\le C_m$ and $\max(1,X_3)^4\jap{v_0''}^{-2\sigma}\le C$ for
some $\sigma>0$, so they hold with $\sigma=\frac18$. The same applies to the absorption of
$\vartheta(v_0'')^\iota$ in Remark~\ref{rem:FPiota} and to Subsection~\ref{sssec:FPYside}.
Lemma~\ref{lem:FPOmega}(i), Lemma~\ref{lem:FPb} and classes~(I), (III) and~(IV) of
Lemma~\ref{lem:FPc} do not use (P1)--(P3), and hold unchanged.

\emph{(2) The cutoffs, (P1).} At $\gamma=-2$, the proof of Lemma~\ref{lem:RBloc} gives the radial
bound $\le C\eps t^2(1+\log\jap v)$ and the transverse bound $\le C\eps t^2\log\jap v$. Thus we have
$|D\sqrt{\chi_Q}|_A^2\le C\eps\ell_Q\mathbf 1_{\tilde Q}$, and Lemma~\ref{lem:FPa} holds with
$C\eps\,\mathcal E_Q$ replaced by $C\eps\ell_Q\,\mathcal E_Q$. In the cutoff variants of the
commutator terms, $\mathfrak d_Q^{1/2}$ is replaced by $(C\eps\ell_Q\mathcal E_Q)^{1/2}$. Their
other factors carry $\mathbb E^{1/2}\le\eta_*\le1$. Thus, by Young's inequality with the parameter
$\eta':=1/(72N_0)$, they contribute $\eta'\,\Sigma_Q+C\eps\ell_Q\,\mathcal E_Q$. The exception is
the cutoff variant of class~(II), whose other factors are $C_{\mathbb E}\mathbb E^\flat[\tilde
G]^{1/2}(1+t_1\ell_Q^{1/2})\mathcal E_Q^{1/2}$, by~(3) below. There we use the hypothesis on $\tilde
G$, which implies that the term is
$\le C\tilde C\eps t^{3/4}(1+t_1)\ell_Q^{1/2}\mathcal E_Q^{1/2}(C\eps\ell_Q
\mathcal E_Q)^{1/2}\le\eps\ell_Q\,\mathcal E_Q$ for $\eps\le\eps_0(\tilde C)$.

\emph{(3) Class~(II), (P2).} By Lemma~\ref{lem:RBsob}, the plain factor $X=\partial_vu-t\nabla_yu$
is bounded pointwise on $Q$ by the energies of the neighbors of $Q$. On outer boxes we treat its
radial part by Lemma~\ref{lem:RBgain}(i), as before. On inner boxes, we have
$h_r\le5\vartheta\le10\ell_Q^{1/2}$, so the radial part yields $th_r\le Ct_1\ell_Q^{1/2}$, like the
transverse part below, and it is included in the factor $1+t_1\ell_Q^{1/2}$. Its transverse part
yields $t\vartheta\le Ct_1\ell_Q^{1/2}$ against the weights, and at $\gamma=-2$ there is no kernel
decay to absorb it. Keeping this factor, we deduce that class~(II) is bounded by
\begin{equation}
C_{\mathbb E}\,\mathbb E^\flat[\tilde G]^{1/2}\big(1+t_1\ell_Q^{1/2}\big)\,
\mathcal E_Q^{1/2}\big(\mathfrak d^\flat_Q\big)^{1/2}\le\eta'\,\mathfrak d^\flat_Q+72N_0C_{
\mathbb E}^2(1+t_1)^2\,\mathbb E^\flat[\tilde G]\,\ell_Q\,\mathcal E_Q ,
\end{equation}
with $\eta'$ of~(2). By the hypothesis on $\tilde G$, the coefficient of $\ell_Q\mathcal E_Q$ is
$\le72N_0C_{\mathbb E}^2(1+t_1)^{7/2}\tilde C^2\eps^2\le\eps$ for $\eps\le\eps_0$. The same bound
holds for class~(II) of the $A\zeta$-part of Lemma~\ref{lem:FPd}, with $(C\eps e^\flat_Q)^{1/2}$ in
place of $(\mathfrak d^\flat_Q)^{1/2}$.

\emph{(4) The drift.} The direct form of $\bb$ in the drift lemma uses
$\int|w|^{\gamma+1}\rho_0(v-w)\dd w\le C\jap v^{\gamma+1}\le C\jap v^{-1}$, which compensates
$t\vartheta$ also at $\gamma=-2$. The $\bb$-part at the top $y$-levels is
Proposition~\ref{prop:FPfix}, with Lemmas~\ref{lem:FPtopother} and~\ref{lem:FPccross}. These apply
with the new weights, by~(1), except in the terms where a power of
$\vartheta(v_0)$ is not absorbed by kernel decay at $\gamma=-2$. Remarks~\ref{rem:FPiota}
and~\ref{rem:FPtopout} list these terms, and the end of part~(ii) of the latter lemma estimates
them. Each of them is bounded by $C\eps\ell_Q\,\mathcal E_Q$, with a constant that does not depend
on $K_\flat$.
Hence the drift term, including the corner cases, is closed
at $\gamma=-2$. The $\ab$-form of Lemma~\ref{lem:FPdlow} does not use (P1)--(P3).

\emph{(5) Zero-order terms.} These are unchanged, since all kernels are locally integrable at
$\gamma=-2$, and $|\cb|\le C\eps\jap v^{-2}$.

\emph{(6) The forcing, (P3).} By Lemma~\ref{lem:EPloc}, we have $|\partial^{\mathbf k} S_0[1]|\le
C_{N_y}C^{|\beta|}\beta!R^{-|\beta|}\eps\jap{y+tv}^{-p}$. The weights
$\delta_y^{|\alpha|}h_r^{-a_r}\vartheta^{-a_\perp}\tau_\alpha$ are $\le1$, and we drop them, since
the moment factors are already taken into account in that lemma. The $v$-sum is that of
Lemma~\ref{lem:FPf}. By the Cauchy--Schwarz inequality, we get
\begin{equation}
\Big|\frac1{|Q|}\sum_{\mathbf k} W_{\mathbf k}^2\mathcal F_{\mathbf k}\Big|\le
F_Q:=C_4'\eps\,\omega_Q\Psi_Q(t)\big(e^\flat_Q\big)^{1/2},
\end{equation}
with $\Psi_Q$ of Lemma~\ref{lem:EPtransit}(iii). By~\eqref{eq:EPomega} and
the same lemma, we have
\begin{equation}
\int_0^tF_Q\le C_4'\eps\jap{v_0}^{1/4}\sup_{s\le t}\mathbb E^\flat(s)^{1/2}\int_0^t\Psi_Q\le
C_6\,\eps\,t^{3/4}\sup_{s\le t}\mathbb E^\flat(s)^{1/2},
\end{equation}
with $C_6$ independent of $Q$ and of $K_\flat$.

\emph{(6$'$) The $\tilde G$-part of the source.} In the frozen problem the coefficients are built
from $f_0\tilde G$, so the last term of the equation for $u_{\mathbf k}$ is
$\partial^{\mathbf k} S_0[\tilde G]$, where $S_0[\tilde
G]=\ab[f_0\tilde G]:(D\zeta+\zeta\otimes\zeta)-\cb[f_0\tilde G]$. By Leibniz it consists of
\emph{explicit terms}, with all derivatives on $f_0$ and $\zeta$ and the factor $\tilde G$, and of
\emph{$G$-terms}, with a factor $h=\partial^{\mathbf k'}\tilde G$, $\mathbf k'\ne0$. In neither is
there a factor $u$. The explicit terms are bounded as in~(6), now with
$G=\tilde G$. For a $G$-term, by the Cauchy--Schwarz inequality for the positive kernels $a$ and
$|c|$, we have
\begin{equation}
|\ab[\partial^\nu f_0\,h]\xi\cdot\eta|\le\big(\ab[|\partial^\nu
f_0|]\xi\cdot\xi\big)^{1/2}\big(\ab[|\partial^\nu f_0|\,|h|^2]\eta\cdot\eta\big)^{1/2},
\end{equation}
and similarly for $\cb$. By Lemma~\ref{lem:FPexp}, the localization lemma with $G=1$ and the
transit-time bound, the first factor is $\le C(\eps\jap{x_Q(t)}^{-p})^{1/2}$ on $\tilde
Q\times B(v_0,2)$, where the moment factors of $\partial^\nu f_0$ are treated as in~(6). We
estimate the second factor in $L^2$ on the box by the absolute bounds
(T2) and~(T5) of Lemma~\ref{lem:ULtri}, as in class~(II) of Lemma~\ref{lem:FPc} and in
Lemma~\ref{lem:FPe}, with the weight algebra of~(1) for the coefficient box. Unlike in~(6), the
target weights $h_r^{-a_r}\vartheta^{-a_\perp}\tau_\alpha$ and the factors
$(|\beta|+1)^{-c_\alpha/2}\vartheta(v_0)^{-\iota_\alpha}$ of Definition~\ref{def:FPtop} are kept and
put into this factor. At the $\tau$-levels they compensate the slab factor
by Lemma~\ref{lem:BKcomp}(i), and they give $\Pi\le1$ in Lemma~\ref{lem:FPtopalg} when no lift
crosses a threshold (we treat the crossing below). Thus this factor is bounded by $C\eps^{1/2}
\mathbb E^\flat[\tilde G]^{1/2}$, against the target weight $\omega_Q$, which is not paired with any
other factor. Hence the contribution of these terms is bounded by
\begin{equation}
F'_Q:=C_7\,\eps\,\omega_Q\jap{x_Q(t)}^{-p/2}\,\mathbb E^\flat[\tilde
G](t)^{1/2}\big(e^\flat_Q\big)^{1/2}.
\end{equation}
Since $p/2\ge2$, it follows from the transit-time lemma that $\int_0^t\jap{x_Q(s)}^{-p/2}\dd s\le
C\min(t,\jap{v_0}^{-1})$, and we have $\jap{v_0}^{1/4}\min(t,\jap{v_0}^{-1})\le t^{3/4}$. Hence we
obtain
\begin{equation}
\int_0^tF'_Q\le CC_7\tilde C\eps^2t^{3/2}\sup_{s\le t}\mathbb E^\flat(s)^{1/2}\le
C_6\,\eps\,t^{3/4}\sup_{s\le t}\mathbb E^\flat(s)^{1/2}\qquad\text{for }CC_7\tilde
C\eps_0t_1^{3/4}\le C_6 .
\end{equation}
\emph{The $c$-crossing.} We call $h=\partial^{\mathbf k'}\tilde G$ $v$-high if $|\beta'|>|\beta|/2$.
Otherwise we have $|\nu_v|\ge|\beta|/2$, and any polynomial factor in $|\beta'|$ is compensated by
the factor
$(C\rad/R)^{|\nu_v|}$ of the same $v$-sum. Let $h$ be $v$-high, with
$|\alpha'|\in\{N-4,N-3\}$ and $|\alpha|\le N-3$. Then the slab bound of Lemma~\ref{lem:slab}(ii),
which is a Sobolev step in $y$ only, lifts $h$ into the levels with $c=1$. The weights leave the
factor $\Pi=(|\beta'|+1)^{1/2}\le(|\beta|+1)^{1/2}$ on $h$, by part~(iii) of the weight-algebra
lemma with the factor $1$ in place of the first factor. The $\tau$-crossing of the lift
($|\alpha'|\in\{N-3,N-2\}$, $|\alpha|\le N-2$) is absolute, that is, the unpaired
$h_r''\vartheta''^2$ is absorbed by the factor $F_0(y'')\rho_0(v'')$ of $|\partial^\nu f_0|$ in the
second factor. As in crossing~(b) of Lemma~\ref{lem:FPccross}, we move this factor to
the target, which is in energy on its own box, and we apply Lemma~\ref{lem:FPhalf} with
$\delta=\frac18$, with $B_Q:=C\eps^{1/2}\omega_Q\jap{x_Q(t)}^{-p/2}$ and $|Q|^{-1}\sum_{\mathbf k}
Z_{\mathbf k}^2\le C\eps\,\mathbb E^\flat[\tilde G]$. Therefore these terms are bounded by
\begin{equation}
\frac1{8|Q|}\sum_{\mathbf k}\theta_{\mathbf k}(a^\flat_{\mathbf k})^2+F''_Q,\qquad
F''_Q:=C\big(\omega_Q^2\jap{x_Q(t)}^{-p}+\omega_Q\jap{x_Q(t)}^{-p/2}\big)\eps\,
\mathbb E^\flat[\tilde G]+C\eps^{1/2}\omega_Q\jap{x_Q(t)}^{-p/2}e^\flat_Q .
\end{equation}
By the transit-time bound and $\jap{v_0}^{1/2}\min(t,\jap{v_0}^{-1})\le t^{1/2}$, we get
$\int_0^tF''_Q\le C_8\big(\tilde C^2\eps^3t^2+\eps^{1/2}t^{3/4}M(t)\big)$, with $M(t):=\sup_{s\le
t}\mathbb E^\flat(s)$ and $C_8$ independent of $Q$, $K_\flat$ and $\tilde C$. We take the
$\theta$-part of this bound from the half of the $\theta$-term of~\eqref{eq:FPid} that the proof we
are following does not use.

\emph{(7) Closing.} We collect the bounds in this identity with the damping
term~\eqref{eq:EPdamp}, as in that proof, with $1/(72N_0)$ in place
of $1/(32N_0)$ there. In particular, the parts of $\theta_{\mathbf k}$ used in
Proposition~\ref{prop:FPfix}, in Lemma~\ref{lem:FPccross} and in the $c$-crossing of~(6$'$),
$\frac14+\frac14+\frac18$ in total, are absorbed by the $\theta$-term of the identity before we
apply Lemma~\ref{lem:EPdamp}. Note that neither side of its hypothesis contains $\theta$. Its four
parts then contribute at most $4/(72N_0)$ to the coefficient of $\Sigma_Q$, and~(2) and~(3)
contribute at most $2\eta'=2/(72N_0)$. Hence the coefficient of $\Sigma_Q$ does not exceed
$\eta:=6/(72N_0)=1/(12N_0)$.
The terms carrying $\ell_Q$ come from~(2), (3) and~(4), while the terms
$C(\eps+\rad^{-2}\eps^{1/2})\mathcal E_Q$ of that proof carry no
$\ell_Q$. Hence we have
\begin{equation}
\tfrac12\tfrac d{dt}e^\flat_Q+\tfrac12\mathfrak d^\flat_Q+K_\flat\eps\ell_Qe^\flat_Q\le\eta\,\Sigma_Q+\tfrac12C^{\rm
lay}_\sharp\eps\ell_Q\,\mathcal E_Q+C_g\,\mathcal E_Q+F_Q+F'_Q+F''_Q,\qquad
C_g:=C(\eps+\rad^{-2}\eps^{1/2}),
\end{equation}
where $C^{\rm lay}_\sharp$ collects the constants of~(2)--(4), and is also taken
$\ge2e^{1/2}(C^{(0)}_{\rm cut}+C^{(0)}_{\rm dr})$, where $C^{(0)}_{\rm cut}\eps\ell_Q$ and
$C^{(0)}_{\rm dr}\eps$ bound the cutoff term and the zero-order and drift terms of
Lemma~\ref{lem:BKtop} at level~$0$, summed over
the neighbors (used in Proposition~\ref{prop:EPconstr}(v), and the factor $2$ is the doubling and
$e^{1/2}$ the neighbor ratio of the weights). This constant does not depend on $Q$, $K_\flat$ or
$C_5$. We multiply by $2$ and integrate, and we note that $e^\flat_Q(0)=0$ because $G(0)=1$. Then we
obtain the hypothesis of the damping lemma with $c=1$, $K=K_\flat$, $2\eta$ in place of $\eta$, and
$2(F_Q+F'_Q+F''_Q)+2C_gN_0\mathbb E^\flat$ in place of $F_Q$. The finiteness conditions of
that lemma are assumptions of the present lemma. Since $K_\flat=9C^{\rm lay}_\sharp N_0$
and $2\eta N_0=\frac16$, the damping lemma implies that, for
$M(t):=\sup_{s\le t}\mathbb E^\flat(s)$,
\begin{equation}
M(t)\le24C_6\eps t^{3/4}M(t)^{1/2}+12C_8\tilde C^2\eps^3t^2+12C_8\eps^{1/2}t^{3/4}M(t)+12C_gN_0\int_0^tM(s)\dd
s .
\end{equation}
For $12C_8\eps_0^{1/2}t_1^{3/4}\le\frac12$ the third term can be put into the left-hand side. The
remaining nonhomogeneous terms are nondecreasing in $t$, so by Gronwall's inequality and
$24C_gN_0t_1\le\log2$, we obtain
\begin{equation}
M(t)\le96C_6\eps t^{3/4}M(t)^{1/2}+48C_8\tilde C^2\eps^3t^2 .
\end{equation}
This implies that $M^{1/2}\le96C_6\eps t^{3/4}+(48C_8)^{1/2}\tilde C\eps^{3/2}t\le(96C_6+1)\eps t^{3/4}=:C_5\eps
t^{3/4}$ for $(48C_8)^{1/2}\tilde C\eps_0^{1/2}t_1^{1/4}\le1$. The constant $C_5$ does not depend on
$\tilde C$. For $\eps_0$ small this is $\le\eta_*/2$, which improves the first bootstrap hypothesis.
For the second, by the level-$0$ energy, Lemma~\ref{lem:RBsob} and the lower bound on the weight, we
get $|G-1|\le
C(\delta_y,\rad)\omega_Q^{-1}\mathbb E^{\flat\,1/2}\le C\eps t^{3/4}\jap v^{-1/8}\le\frac14$. The
bootstrap set is closed, open and contains $t=0$, as before. The
conditions on $\eps_0$ used in~(6)--(7) are $CC_7\tilde C\eps_0t_1^{3/4}\le C_6$,
$12C_8\eps_0^{1/2}t_1^{3/4}\le\frac12$, $24C_gN_0t_1\le\log2$, $(48C_8)^{1/2}\tilde
C\eps_0^{1/2}t_1^{1/4}\le1$, $C_5\eps_0t_1^{3/4}\le\eta_*/2$, $K_\flat\eps_0t_1\le\frac18$, $\tilde
C\eps_0t_1^{3/4}\le\eta_*$ (so that the coefficient factor satisfies $\mathbb E^\flat[\tilde
G]^{1/2}\le\eta_*\le1$, used in~(2) and in Lemma~\ref{lem:FPcross}(iv)),
$e\,C(\delta_y,\rad_0/2)C_5\eps_0t_1^{3/4}\le\frac14$ (for $|G-1|\le\frac14$), and the conditions
of Proposition~\ref{prop:FPclose} that do not involve the coefficient factor, with $1/(72N_0)$ in
place of $1/(32N_0)$, namely $C_\rad\eps_0t_1\le\rad_0/2$ and
$4C\rad_0^{-2}\eps_0^{1/2}\le1/(72N_0)$. Besides these, we need those of~(2)--(3) and those of
that proposition that involve the coefficient factor, with $\tilde C\eps t^{3/4}$ in
place of $2C_4C_{N_y}^{1/2}\eps t$, in particular $CC_2\tilde C(1+t_1)\eps_0\le1/(144N_0)$ for the
crossing terms of that lemma, by $t^{3/4}\ell_t\le C(1+t_1)$.

\emph{(8) Viscosity.} For the viscous approximations of Definition~\ref{def:FPreg}, we add the
viscous terms of Proposition~\ref{prop:FPvisc}(iii) as there and in Remark~\ref{rem:FPtrunc8}. Since
the factor $\omega_Q$ is constant on each box, these terms are unchanged. After the doubling, the
viscous
dissipation enters the left side with the coefficient $2c_\nu=\frac32\ge1$, and the neighbor viscous
dissipations enter the right side with the coefficient $2\eta_\nu$, where $2\eta_\nu
N_0\le\frac7{48}<\frac16$ by the same remark. The remainder $C\nu\,\mathcal E_Q$,
$\nu\le\eps$, is included in $C_g\,\mathcal E_Q$. We apply the damping lemma to $
\mathfrak d^\flat+\mathfrak d^\nu$, with $c=1$ and $\max(2\eta,2\eta_\nu)$ in place of $\eta$, and
we obtain the same conclusion. This completes the proof.
\end{proof}

\subsubsection{The initial layer: construction and transfer}\label{sssec:EPthm}

\begin{prop}[Construction and transfer for $\gamma\in{[-3,-2]}$]\label{prop:EPconstr}
Let $\gamma\in[-3,-2]$. Proposition~\ref{prop:FPvisc}, Lemma~\ref{lem:FPiter} and
Proposition~\ref{prop:FPtr} hold with the energies of Definition~\ref{def:EPen}, with
Lemma~\ref{lem:EPclose} in place of Proposition~\ref{prop:FPclose}, with $
\mathbb E^\flat(t)^{1/2}\le C_5\eps t^{3/4}$ in place of $\mathbb E(t)^{1/2}\le2C_4\eps t$ (also in
Remark~\ref{rem:FPpointwise}), and with $\eps t_1^{3/4}$ in place of $\eps t_1$ in
the latter. The following steps change.
\begin{enumerate}[label=(\roman*)]
\item \emph{The barrier.} In part~(ii) of the frozen-problem proposition, we use
\begin{equation}
\psi:=K_0\eps e^t\,\Psi,\qquad\Psi(t,y,v):=\int_0^t\jap{y+sv}^{-p}\dd s .
\end{equation}
Then we have $|G_{\nu,k}-1|\le\psi\le C\eps\min(t,\jap v^{-1})$.
\item \emph{Finiteness of the energies.} The Schauder estimates used there, with constants uniform
in $Q$, are replaced by local $L^2$
estimates of all orders up to $t=0$ for the rescaled problem, which need the forcing only in
$L^2_tH^m$. They give
\begin{equation}
\frac1{|Q|}\sup_{[0,t]}\big\|h_r^{-a_r}\vartheta^{-a_\perp}\partial_Q^\alpha\partial_v^\beta
u\big\|^2_{L^2(Q)}\le C(m,\nu)\,\ell_Q^{C_m}\eps^2\jap{v_0}^{-1}\qquad(|\alpha|+|\beta|\le m),
\end{equation}
and the truncated energies $\mathbb E^{\flat,(n)}[G_{\nu,k}]$ are bounded uniformly in $k$.
\item \emph{The contraction.} In part~(ii) of the iteration lemma, the transverse entries of
$D^2G^j$
carry the factor $t^2\vartheta^2$ and meet $|A|\le C\eps$. Their contribution is $\le
C\eps^2\sup|\delta^{j-1}|$, and $|F^j|\le C\eps\sup|\delta^{j-1}|$.
\item \emph{The transfer.} In Step~2 of the transfer, the inverse weight contains
$\omega_Q^{-1}\le e\jap{v_0}^{-1/8}$ in place of $\jap{v_0}^{-\sigma}$.
\item \emph{Uniqueness.} In part~(v) of the frozen-problem proposition, the Gronwall step in
$\sup_Qe^{(0)}_Q$ is replaced by Lemma~\ref{lem:EPdamp}, applied to the level-$0$ energies with the
weight $e^{-K_\flat\eps t\ell_Q}$.
\end{enumerate}
\end{prop}

\begin{proof}
(i) We have $\partial_t\psi=\psi+K_0\eps e^t\jap{y+tv}^{-p}$. Since $D(y+sv)=(s-t)I$ and
$|\nabla\jap\cdot^{-p}|\le p\jap\cdot^{-p}$, we have
\begin{equation}
|D\Psi|\le pt\Psi,\qquad|D^2\Psi|\le Ct^2\Psi,\qquad|\nabla_{y,v}\Psi|+|\nabla^2_{y,v}\Psi|\le
C(1+t)^2\Psi .
\end{equation}
By Lemma~\ref{lem:FPMnu}, the coefficients of $\mathcal L_\nu$ and their first derivatives are
$\le C(\eps+\nu)$, so $|\mathcal L_\nu\psi|\le C(\eps+\nu)(1+t_1)^2\psi\le\frac12\psi$ for
$\eps$ and $\nu$ small. By Lemma~\ref{lem:EPloc}, applied to the mollified coefficients (the
mollification changes $\jap{y+tv}$ by a bounded factor), we have $|\chi_k\tilde
S_{0,\nu}|\le C\eps\jap{y+tv}^{-p}$. Hence $\partial_t\psi-\mathcal L_\nu\psi-|\chi_k\tilde
S_{0,\nu}|\ge\frac12\psi+(K_0-C)\eps\jap{y+tv}^{-p}\ge0$ for $K_0\ge C$. By the maximum principle,
we have $|u|\le\psi$, and the last bound follows from Lemma~\ref{lem:EPtransit}(i).

(ii) We rescale $\tilde Q$ in the space variables as in Proposition~\ref{prop:FPvisc}(ii), on the
whole
interval $[0,t]$ (there is no rescaling and no cutoff in time). The rescaled operator is uniformly
parabolic with ellipticity $\ge c\nu$. The coefficient bounds listed there hold, except the
transverse second-order coefficient, which is now $t^2\vartheta^2|\tilde A|\le C\eps\vartheta^2\le
C\eps\ell_Q$, and the transverse part of the drift, $t\vartheta|\tilde V|\le C\eps\vartheta\le
C\ell_Q^{1/2}$. Therefore the rescaled coefficients have $C^m$ norms $\le C(m,\nu)\ell_Q$.

We differentiate the rescaled equation and use the energy identity with cutoffs. The derivatives
that
enter exist and are continuous up to $t=0$ on each box, by the same argument, with constants that
need not be uniform in $Q$. The constants of
these local $L^2$ estimates are polynomial in the coefficient norms. Indeed, at each order $j\le m$,
we integrate the energy identity for the derivatives of order $j$ with a spatial cutoff over
$[0,t]$, from $0$, where they vanish since $u(0)=0$. We do this also for $t>1$, since the rescaling
is in space only, so no cutoff in time enters. No Gronwall step is needed, since after Young's
inequality against the ellipticity $c\nu$ the remaining terms involve only the forcing and the
space-time $L^2$ norms of the derivatives of order $\le j$ on a larger cutoff. These norms are
bounded at the previous order, namely those of order $j$, including the drift and zero-order terms,
by the dissipation of order $j-1$, and at order $0$ by $\sup|u|$. Thus each order contributes a
factor
$C(1+K/\nu)^2$, where $K\le C(m,\nu)\ell_Q$ is the coefficient norm. With the remark after
the localization lemma, we have $|\partial^m\tilde S_{0,\nu}|\le C(m,\nu)\eps\jap{y+tv}^{-p}$ in
the rescaled variables. Hence
\begin{equation}
\frac1{|Q|}\sup_{[0,t]}\big\|h_r^{-a_r}\vartheta^{-a_\perp}\partial_Q^\alpha\partial_v^\beta
u\big\|^2_{L^2(Q)}\le C(m,\nu)\ell_Q^{C_m}\Big(\sup_{2\tilde
Q\times[0,t]}|u|^2+\frac1{|Q|}\int_0^t\|\tilde S_{0,\nu}\|^2_{H^m(2\tilde Q)}\Big)\le
C(m,\nu)\ell_Q^{C_m}\eps^2\jap{v_0}^{-1},
\end{equation}
by~(i) and by the transit-time bound with $p'=2p$, since $\jap v\ge\jap{v_0}/5$ on $2\tilde
Q$. Since $\omega_Q^2\ell_Q^{C_m}\jap{v_0}^{-1}\le\jap{v_0}^{-1/2}\ell_Q^{C_m}\le C_m$, the
truncated energies are bounded uniformly in $k$. The local $L^2$ bounds of all orders give local
$C^m$ bounds by Sobolev's inequality, so we take the limit $k\to\infty$ as there. Continuity in $t$
of the energies follows in the same way from
the equation for $\partial_tu$. The uniqueness argument there rests on the boundedness of the
coefficients and of their first derivatives alone, and is unchanged. Since
$\ell_Q\,\omega_Q^2\ell_Q^{C_m}\jap{v_0}^{-1}\le C_m$, the same bound yields
$\sup_Q\ell_Qe^{\flat,(n)}_Q<\infty$. With one more derivative and the bound $C(m,\nu)\ell_Q$ of
the rescaled second-order coefficients (radially on outer boxes this uses the null structure, as in
the same proposition, since $|\tilde A|\le C\eps$ alone would leave the factor $t^2h_r^2$),
it yields $\sup_Q\int_0^t\mathfrak d^{\flat,(n)}_Q<\infty$, and the viscous dissipation is bounded
in the same way, uniformly in $k$. These are the finiteness conditions of Lemma~\ref{lem:EPclose}
for the truncated energies, which are passed to the limit as in Proposition~\ref{prop:FPvisc}(iii).

(iii) By the closed estimate and Lemma~\ref{lem:RBsob}, we have $|\partial_Q^\alpha\partial_v^\beta
G^j|\le C\eps\jap v^{-1/8}h_r^{a_r}\vartheta^{a_\perp}$ for $|\alpha|,|\beta|\le2$. Thus the
transverse entries contribute
$\le C\eps\cdot t^2\vartheta^2\cdot C\eps\jap v^{-1/8}\sup|\delta^{j-1}|\le
C\eps^2\sup|\delta^{j-1}|$, since $\vartheta^2\jap v^{-1/8}\le C$. On outer boxes the radial entries
are unchanged. On inner boxes, where $h_r\le5\vartheta$, the radial entries carry the factor
$t^2h_r^2\le25t^2\vartheta^2$, and we treat them as the transverse ones. Likewise, in the
$\zeta$-terms the transverse part $t\vartheta$ of $DG^j$ (and on inner boxes its radial part) is
compensated by $\jap v^{-1/8}$, with $(\ab[f_0]\zeta^\sharp\cdot\zeta^\sharp)^{1/2}\le C\eps^{1/2}$,
and the rest of the $\zeta$-terms is unchanged. There is no $\bb$-term, as in the proof of
Lemma~\ref{lem:FPiter}(ii). The
zero-order part $\big(\ab[f_0\delta^{j-1}]:(D\zeta+\zeta\otimes\zeta)-\cb[f_0\delta^{j-1}]\big)G^j$
is $\le C\eps\sup|\delta^{j-1}|$, as for $\gamma<-2$. Hence we obtain
$|F^j|\le C\eps\sup|\delta^{j-1}|$, and
$\sup|\delta^j|\le e^{C\eps t_1}t_1C\eps\sup|\delta^{j-1}|\le\frac12\sup|\delta^{j-1}|$ for $\eps$
small.

(iv) In Step~2 of Proposition~\ref{prop:FPtr}, we have $\vartheta(v)^m\omega_Q^{-1}\le
e\,\vartheta(v)^m\jap{v_0}^{-1/8}\le C_m$, so the bound on the inverse weights does not depend on
$\xi$, as before. In Step~1, we have
$\sum_{\mathbf k'}\mathcal H_{\mathbf k'}^2\le C\mathbb E^\flat(t_1)\le
C\eps^2t_1^{3/2}$, and in~(ii)--(iii) of that proposition, $|G-1|\le C\eps t_1^{3/4}$ by
the closed estimate. The rest of the transfer is unchanged.

(v) Let $\delta$ be the difference of two classical solutions of the frozen problem at $\nu=0$ with
$\sup|G|<\infty$ and $\sup_{t,Q}e^{(0)}_Q<\infty$. Let $e_Q$ be the level-$0$ energy of $\delta$
with the weight $e^{-K_\flat\eps t\ell_Q}$ (and no factor $\jap{v_0}^{\sigma_\flat}$), and
$\mathfrak d_Q$ its dissipation. The level-$0$ identity with Lemma~\ref{lem:BKtop}, the cutoff term
$C\eps\ell_Q$ of step~(2) of the closed estimate and the damping term of the weight (its
zero-order and drift terms are $\le C^{(0)}_{\rm dr}\eps\,\mathcal E_Q\le
C^{(0)}_{\rm dr}\eps\ell_Q\,\mathcal E_Q$, and the neighbor ratio of the weights is $\le e^{1/2}$)
gives, after multiplication by $2$, the hypothesis of Lemma~\ref{lem:EPdamp} with $c=1$, $\eta=0$,
$F_Q=0$ and the constant $2e^{1/2}(C^{(0)}_{\rm cut}+C^{(0)}_{\rm dr})\le C^{\rm lay}_\sharp$ (no
commutator enters at level~$0$). Its finiteness conditions hold. Indeed, we have $\ell_Qe_Q(s)\le
C\ell_Qe^{-2K_\flat\eps s\ell_Q}\sup|\delta|^2$ and $\int_0^t\ell_Qe^{-2K_\flat\eps s\ell_Q}\dd
s\le(2K_\flat\eps)^{-1}$, and the same identity, integrated on each box, bounds $\int_0^t
\mathfrak d_Q$ by $C\sup|\delta|^2(1+(K_\flat\eps)^{-1})$, uniformly in $Q$. By the lemma, we
deduce that $\sup_Qe_Q\le0$, and thus $\delta=0$. This completes the proof.
\end{proof}

\begin{theo}[The initial layer for $\gamma\in{[-3,-2]}$]\label{thm:EPlayer}
Let $\gamma\in[-3,-2]$, $f_0\in\Dc$, $N_y\ge20$. With the energies of Definition~\ref{def:EPen},
there is $\eps_0'(\Dc)>0$ such that for $\eps\le\eps_0'$ the solution exists on $[0,t_1]$ and
satisfies
\begin{equation}
\mathbb E^\flat(t)^{1/2}\le C_5\eps t^{3/4},\qquad|G-1|\le C\eps t^{3/4}\jap v^{-1/8}\qquad(t\le
t_1),
\end{equation}
and at $t_1$ it satisfies the hypotheses of Theorem~\ref{thm:t1} and the bound~\eqref{eq:dispdata}.
For $\gamma<-2$ this gives a second proof of the existence on $[0,t_1]$ and of the transfer, with
weights that do not depend on $\gamma$. However, the energy bound differs from that of
Theorem~\ref{thm:layer}.
\end{theo}

\begin{proof}
We follow the proof of Theorem~\ref{thm:layer} in Subsections~\ref{ssec:layer}--\ref{ssec:drift},
with Lemma~\ref{lem:EPclose} in place of Proposition~\ref{prop:FPclose} and with the changes of
Proposition~\ref{prop:EPconstr}.
\begin{itemize}
\item The frozen problems are solved by Proposition~\ref{prop:FPvisc}, with the barrier and the
local estimates of the changes~(i)--(ii), and their energies are finite and
continuous in $t$.
\item The closed estimate yields the bound $\mathbb E^\flat(t)^{1/2}\le C_5\eps t^{3/4}$ for the
frozen problems, uniformly in the viscosity, with the hypothesis on the coefficient factor
reproduced as in Remark~\ref{rem:FPpointwise}.
\item The iteration converges by the change~(iii), and the limit inherits the
bounds.
\item At $t_1$, Proposition~\ref{prop:FPtr} with the change~(iv)
gives its conclusions~(i)--(iii) with $\eps t_1^{3/4}\le\eps t_1$ in place of $\eps t_1$. These are
the hypotheses and the bound stated in the theorem.
\end{itemize}
This completes the proof.
\end{proof}

\subsubsection{The near-Maxwellian input}\label{sssec:EPE1}

The input (E1) enters Theorem~\ref{thm:enter}, for the trapping, through the analogue of
\cite[Prop.~5.1]{CM} derived in Appendix~\ref{app:hom}. It also enters the near-equilibrium
hypodissipativity of Sections~\ref{sec:weights} and~\ref{sec:boot}, and Lemma~\ref{lem:Abarrier}.
The paper \cite{CM} treats $\gamma\in[-3,-2)$ only. At $\gamma=-2$, Appendix~\ref{app:CM} proves
what we use of \cite{CM}, in Lemmas~\ref{lem:CMpoly} and~\ref{lem:CMbil} and
Proposition~\ref{prop:CMm2}. More precisely, \cite[Lemmas~2.2--2.4, 3.1 (polynomial case), 3.2,
4.3]{CM} in polynomial weights are proved for $\gamma\in[-3,-2]$ in Subsection~\ref{ssec:CMpoly},
and \cite[Lemmas~2.1, 4.1, 4.2, Cor.~4.4]{CM} hold at $\gamma=-2$ by items~(1) and~(3)--(5) of the
proof of that proposition. The exponential case of \cite[Lemma~3.1]{CM} is used in \cite{CM} only
for \cite[Prop.~3.3]{CM}, which is replaced as follows. The decay in a small space,
\cite[Prop.~3.3]{CM}, is replaced by the spectral gap of the linearized operator at $\gamma=-2$,
from \cite[(2.6)]{CM} and \cite[Thm.~1.3]{Mou06}, as in Lemma~\ref{lem:CMgap} (see also
\cite[Lemma~5]{Guo02}). Then we derive the conclusions of \cite[Thm.~3.5, Prop.~3.6, Cor.~3.7,
Prop.~5.1]{CM} as for $\gamma<-2$, with the thresholds of Proposition~\ref{prop:CMthr}, at the
weight levels $k\ge7$ used here.

\section{Proofs of the main theorems}\label{sec:proofmain}

\subsection{The dispersed layer}\label{ssec:disp}

We carry the solution from $t_1$ to $T_0$ and show that at $T_0$ it satisfies the hypotheses of the
bootstrap. Corollary~\ref{cor:LWPa} proves the statement (LWP)(a) about the data at $T_0$. It starts
from the data at $t_1$ given by Theorem~\ref{thm:layer} or, at $\gamma=-2$, by
Theorem~\ref{thm:EPlayer}. In the latter case, by Proposition~\ref{prop:EPconstr}(iv),
Proposition~\ref{prop:FPtr} holds with $\eps t_1^{3/4}\le\eps t_1$ instead of $\eps t_1$. On
$[t_1,T_0]$ the collisions act for a total effective time $O(\eps T_0^{\kappa_\gamma})$, and we use
energies in a \emph{fixed} frame, in cells, with no damping and no radius loss.

\paragraph{Parameters.} Enlarging $L$ (which enlarges $\Dc$), we assume:
\begin{itemize}
\item (D3) and (D4) hold with $L/4$, so that $\jap v|\nabla_v\log(\rho_0F_0)|\le L/2$;
\item $|y_{0,\infty}(\xi)|\le L/2$ on $\Dc^\sharp$;
\item $|\nabla\log\rho_0|+|\nabla\log\theta_\infty|\le L/2$, the margin that (H1) needs in
Corollary~\ref{cor:LWPa}(a3);
\item $A_\infty:=\rho_\infty^{-1}\theta_\infty^{3/4}$ satisfies $|\nabla\log A_\infty|\le L/2$ and
$|\nabla^2\log A_\infty|\le L$.
\end{itemize}
These conditions can be met. Indeed, we have $|y_{0,\infty}|\le C(C_F,\beta_F)$. By (D4) with $L/4$,
as in Subsection~\ref{ssec:data}, we have $|\nabla\log\theta_\infty|\le L/4$, so $|\nabla\log
A_\infty|\le\frac L4+\frac34\cdot\frac L4\le\frac L2$. Moreover, the bounds of that subsection imply
$|\nabla^2\log A_\infty|\le C(C_F,\beta_F)L_{\rm old}^2$, where $L_{\rm old}$ is the constant of
(D3) and (D4) before the enlargement. Thus $L\ge\max\{4L_{\rm old},C(C_F,\beta_F)L_{\rm old}^2\}$
suffices, and this is possible since $C(C_F,\beta_F)$ does not depend on $L$.

Fix $\theta_\pm$ with $2\theta_-\le\theta_\infty(\xi)\le\theta_+/2$ on $\Dc^\sharp$ (possible since
$c(C_F)\le\theta_\infty\le C(C_F,\beta_F)$), and put $\rho_+:=2\eps A_F$. The exponent
$\kappa_\gamma$ is that of Subsection~\ref{ssec:leaf}, so that $\int_{t_1}^{T_0}\kappa\le
CT_0^{\kappa_\gamma}$ and $T_0^{\kappa_\gamma}\ge1$. The \emph{data exponent} is
\begin{equation}\label{eq:qd}
q_d:=\min\big\{\beta_F/4,\ q_{\rm tri},\ (8(L^2+3\theta_+))^{-1}\big\},
\end{equation}
with $q_{\rm tri}=\frac18$ of Proposition~\ref{prop:triG}. By (D4), $y_{0,\infty}$ and
$\theta_\infty$ extend holomorphically to $\mathbb S_\Lambda$ (integrate $F_0(z,\xi)$ against $z$
and $|z|^2$, with the bound $|F_0(z,\xi+i\eta)|\le A_FC_Fe^{-\beta_F|z|^2}$), and they are bounded
there by a constant $C(\Dc)$. By the Cauchy estimates, we have $|\nabla
y_{0,\infty}|+|\nabla\theta_\infty|\le C_{\mu_0}=C_{\mu_0}(\Dc)$ on $\mathbb S_{\Lambda/2}$. Let
\begin{equation}\label{eq:Ginf}
G_\infty:=\sup_{\R^3}\big(|\nabla y_{0,\infty}|+|\nabla\theta_\infty|\big)\le C_{\mu_0};
\end{equation}
for product data $C_{\mu_0}=G_\infty=0$.

The layer radius $\rad_1\le\min\{\Lambda/8,\rad_0^{\rm
lay}/4,1/(6L),\min\{1,\theta_-\}/(1+C_{\mu_0})\}$ and the time $t_1\ge1$ are fixed in the proof of
Proposition~\ref{prop:disp}. Both depend only on $\Dc$ and on structural constants, so
$t_1=t_1(\Dc)$. The weight is $m_t(z):=e^{q(t)|z|^2}$ with $q(t):=q_d-C_q(t_1^{-1}-t^{-1})$, where
$C_q$ is fixed in Step~2 below. Since $t_1\ge4C_q/q_d$, we have $q(t)\in[\frac34q_d,q_d]$, hence
$q(t)-\delta\in[\frac58q_d,q_d]$ for $\delta\le q_d/8$.

\paragraph{Fixed-frame leaf functions and energies.} For $t\in[t_1,T_0]$ let
$G(t,z,\xi):=g(t,z,\xi-z/t)$ be the leaf through $x=t\xi$, and
\begin{equation}\label{eq:disp-phi}
\phi:=\frac G{\eps\rho_0(\xi)},\qquad D_y:=\partial_z+t^{-1}\partial_\xi,\qquad\phi_{\alpha\beta}:=D_y^\alpha\partial_\xi^\beta\phi
.
\end{equation}
By the leaf calculus of Lemma~\ref{lem:leafgeom}, $D_yG$ and $\partial_\xi G$ are the restrictions
of $\partial_yg$ and $\partial_vg$ to the leaf. For $k\le N_y$ we put, with $w_n=w_n(\rad_1)$,
\begin{equation}\label{eq:disp-e}
\mathsf e^{(k)}_\xi:=\sum_{|\alpha|=k}\sum_\beta w_{|\beta|}^2\|m_t\phi_{\alpha\beta}\|^2_{L^2_z},\qquad
\mathsf d^{(k)}_\xi:=\eps\rho_0(\xi)\sum_{|\alpha|=k}\sum_\beta
w_{|\beta|}^2\|\phi_{\alpha\beta}\|^2_{Y_2(m_t)},
\end{equation}
and $\mathsf e^{(0)}_{\ge1,\xi}$ for the part of $\mathsf e^{(0)}_\xi$ with $|\beta|\ge1$. There is
no $y$-weight: each $y$-order has its own energy. The cells are
$\chi_c(\xi)^2:=e^{-2\jap{(\xi-\xi_c)/\rad_1}}$, the cubes have side $\rad_1$, and
\begin{equation}\label{eq:disp-E}
\mathsf E^{(k)}_c:=\rad_1^{-3}\!\int\chi_c^2\,\mathsf e^{(k)}_\xi\dd\xi,\qquad
\mathsf D^{(k)}_c:=\rad_1^{-3}\!\int\chi_c^2\,\mathsf d^{(k)}_\xi\dd\xi,
\end{equation}
and similarly $\mathsf E^{(0)}_{c,\ge1}$ and the cube energies $\mathsf a^{(k)}_Q$, $
\mathsf a^{(0)}_{Q,\ge1}$. A superscript ${\le}N$ restricts the sums to $|\beta|\le N$. The free profile has
$\phi^{\rm fr}(t,z,\xi)=f_0(z,\xi-z/t)/(\eps\rho_0(\xi))=\rho_0(\xi-z/t)F_0(z,\xi-z/t)/\rho_0(\xi)$.

\begin{lem}[The free profile in the layer norms]\label{lem:dispfree}
For $t\ge1$, $\rad_1\le\Lambda/8$ and $q\le\beta_F/4$, $\sup_\xi\mathsf e^{(k)}_\xi[\phi^{\rm
fr}]\le E^{\rm fr}$ for every $k\le N_y$, and
\begin{equation}\label{eq:dispfree-ge1}
\sup_\xi\mathsf e^{(0)}_{\ge1,\xi}[\phi^{\rm fr}]\le C(\rad_1/\Lambda)^2E^{\rm fr},
\end{equation}
where $E^{\rm fr}$ depends only on $\Dc$.
\end{lem}

\begin{proof}
By the leaf calculus, $D_y$ and $\partial_\xi$ of the leaf function $f_0(z,\xi-z/t)$ are the
restrictions of $\partial_y$ and $\partial_v$ of $f_0$, so
\begin{equation}
\phi^{\rm fr}_{\alpha\beta}=D_y^\alpha\partial_\xi^\beta\Big[\frac1{\rho_0(\xi)}\,(\rho_0F_0)(z,\xi-z/t)\Big]
\end{equation}
is a Leibniz sum of terms in which $\partial^{\alpha'}_y\partial^{\beta'}_v(\rho_0F_0)$ is evaluated
at $(z,\xi-z/t)$ and the remaining derivatives, $\partial_\xi$ and the parts $t^{-1}\partial_\xi$ of
$D_y$, fall on $1/\rho_0(\xi)$. By (D3) and (D4), for real $(z,v)$ the function
$\zeta\mapsto\rho_0(v+\zeta)\,\partial^{\alpha'}_yF_0(z,v+\zeta)$ is holomorphic on the polydisc
$\Delta(0,\Lambda/2)$, and bounded there by $A_Fe^{L\Lambda}\rho_0(v)\cdot A_Fe^{-\beta_F|z|^2/2}$,
for $|\alpha'|\le N_x$. Likewise, $1/\rho_0$ is holomorphic on $\Delta(\xi,\Lambda/2)$ and bounded
by $A_Fe^{L\Lambda}/\rho_0(\xi)$. By the Cauchy estimates and $\rho_0(\xi-z/t)\le
e^{L|z|/t}\rho_0(\xi)$, we obtain $|\phi^{\rm fr}_{\alpha\beta}|\le
C\beta!(\Lambda/4)^{-|\beta|}e^{L|z|-\beta_F|z|^2/2}$ for
$|\alpha|\le N_y\le N_x$, $t\ge1$. Since
$w_n(\rad_1)\beta!(\Lambda/4)^{-n}\le(4\rad_1/\Lambda)^n(n+1)^5\le2^{-n}(n+1)^5$ and
$e^{2q|z|^2+2L|z|-\beta_F|z|^2}$ is integrable, we get both bounds, and the part $|\beta|\ge1$
carries the factor $(4\rad_1/\Lambda)^2$.
\end{proof}

\paragraph{The data at $t_1$.} By Proposition~\ref{prop:FPtr}(i), applied with the radius
$\rad_1\le\rad_0^{\rm lay}/4$ and with the weight $e^{2q_d|z|^2}\le Ce^{\beta_F|z|^2/2}$, and
by the Leibniz rule for the factor $1/\rho_0(\xi)$ (the Cauchy estimates above), we have
\begin{equation}\label{eq:dispdata}
\sup_\xi\mathsf e^{(k)}_\xi[\phi(t_1)-\phi^{\rm fr}(t_1)]\le C_{\rm
lay}\delta_{y,\ell}^{-2k}\eps^2t_1^2\qquad(k\le N_y),
\end{equation}
where $C_{\rm lay}=C(\delta_{y,\ell},\rad_0^{\rm lay},\Lambda,t_1)$ is compensated by the
smallness of $\eps$. Moreover, by parts~(ii)--(iii) of the same proposition, the leaves
$\phi(t_1,\cdot,\xi)$ satisfy (G1)--(G3) with constants within a factor $1+C\eps t_1$ of those of
Lemma~\ref{lem:LWfree}(i).

\begin{prop}[The dispersed layer in cell norms]\label{prop:disp}
There are $\rad_1$, $t_1$ and constants $E_k$, $k\le N_y$, depending only on $\Dc$ and on
structural constants, and $\eps_0(T_0)>0$, such that for $\eps\le\eps_0(T_0)$ the solution exists on
$[t_1,T_0]$ and, for every cell and every $t\in[t_1,T_0]$,
\begin{equation}\label{eq:disp}
\mathsf E^{(k)}_c(t)\le E_k\quad(k\le N_y),\qquad\mathsf E^{(0)}_{c,\ge1}(t)\le
C(\rad_1/\Lambda)^2E_0(1+E_0/c_h^2),\qquad\int_{t_1}^{T_0}\kappa\,\mathsf D^{(k)}_c\dd
t\le\frac{2E_k}{c_h} .
\end{equation}
Moreover, with $w:=\phi-\phi^{\rm fr}$ and the energies $\mathsf e^{(k)}[w]$ at the radius
$\rad_1/2$ and with the weight $m_t^{1/2}$, $\sup_c\mathsf E^{(k)}_c[w](t)\le C_w(\eps
T_0^{\kappa_\gamma})^2$ for $k\le N_y-2$, with $C_w$ depending only on $\Dc$ and on structural
constants.
\end{prop}

\begin{proof}
\emph{Construction.} We write $\eta_v$ for the viscosity, to avoid a clash with the coefficient
$\eta$ of the leaf density. We fix $\delta\in(0,q_d/8]$ and lower the Gaussian exponent of all
energies
below from $q(t)$ to $q(t)-\delta$. We run the argument on the viscous approximation $g_\eta$ of
Lemma~\ref{lem:LWregex}, that is, $\partial_tg=\ab[g]:D^2g-\cb[g]g+\eta_v(\Delta_y+\Delta_v)g$ in
the profile variables (the coefficients are the profile convolutions of Lemma~\ref{lem:LWuniq},
equivalently $\kappa$ times the leaf coefficients), with $g(t_1)$ as data and
$\eta_v\in(0,\eta_0^{\rm disp}]$,
\begin{equation}
\eta_0^{\rm disp}:=\min\Big\{1,\ \frac{C_q}{32q_d^2T_0^2},\ \frac\delta{64q_d^2T_0},\
\frac{\rad_1^2}{C_\eta T_0}\Big\},
\end{equation}
with $C_\eta=C(E^{\rm fr},E_0,\dots,E_{N_y},L,\Lambda)$ fixed in the bullets below.
\begin{itemize}
\item \emph{Data.} $0\le g(t_1)\le2f_0$, and the data satisfy~(D) of point~(3) of the proof of
Proposition~\ref{prop:CULW} at the exponent $q_d-\delta/4$, with $|\alpha|\le N_y-3$. For the data
of Theorem~\ref{thm:layer} this follows from Proposition~\ref{prop:FPclose} and Step~1 of the proof
of Proposition~\ref{prop:FPtr} (even at the exponent $\beta_F/2$). At $\gamma=-2$ the closed
estimate is Lemma~\ref{lem:EPclose}, and the transfer is modified as in
Proposition~\ref{prop:EPconstr}(iv).

For the data of Theorem~\ref{thm:t1} it follows from the pointwise hypothesis (with
$C_{t_1}\eps\le1$) and from~\eqref{eq:dispdata}, by Lemma~\ref{lem:LCdata} with $t_*=t_1$, $q_*=q_d$
and the radius $\rad_1$. Its hypothesis~(H) follows from the same bound (with $C_{\rm
lay}\eps^2t_1^2$ replaced by $C_{t_1}^2\eps^2$) and Lemma~\ref{lem:dispfree}. Indeed, since $D_yG$
and $\partial_\xi G$ are the restrictions of $\partial_yg$ and $\partial_vg$,
$(\partial^\alpha_y\partial^\beta_vg)^{(x)}=D_y^\alpha\partial^\beta_\xi\big(\eps\rho_0(\xi)\phi\big)$
is a Leibniz sum of terms
$\binom\alpha{\alpha'}\binom\beta{\beta'}t_1^{-|\alpha'|}\eps\,\partial^{\alpha'+\beta'}\rho_0(\xi)\,\phi_{\alpha-\alpha',\beta-\beta'}$.
Here $|\partial^\mu\rho_0(\xi)|\le C\mu!(\Lambda/2)^{-|\mu|}\rho_0(\xi)$ by the Cauchy
estimates used for the free profile, and $\|m_{t_1}\phi_{\alpha\beta}\|_{L^2_z}\le
\mathsf e^{(|\alpha|)}_\xi[\phi]^{1/2}/w_{|\beta|}(\rad_1)\le C\,|\beta|!\,\rad_1^{-|\beta|}$
by the same two estimates. From these two bounds, $\binom\beta{\beta'}\le\binom{|\beta|}{|\beta'|}$
and $4\rad_1\le\Lambda/2$, we obtain
$M_1=C(\Dc,\delta_{y,\ell})$.

Therefore~\eqref{eq:LWdata} holds, with $\eta_vc_0\le\delta/(64q_dT_0)\le\frac1{256}$, and
Lemmas~\ref{lem:LWlin} and~\ref{lem:LWregex} apply.

With $c_0:=q_d-\delta/4$ and $\eta_v\le\delta/(64q_d^2T_0)$, the barrier rate satisfies
$c_\eta(t)-\delta/8\ge c_0-4\eta_v c_0^2T_0-\delta/8\ge q(t)-\delta+\delta/2$ on $[t_1,T_0]$. Then,
applying part~(iii) of the latter with $\delta'=\delta/8$ on the intervals $[t_1+j,t_1+j+1)$,
$j\ge0$, and $[t_1+j-\frac12,t_1+j+\frac12)$, $j\ge1$ (intersected with $[t_1,T_0]$), we get
$|\partial^\kappa g_\eta(t)|\le C(\kappa,\eta_v,j,s)\,\Gamma_{c_\eta(t)-\delta/8}$ for
$t\in[t_1+j+s,t_1+j+1]$ and
$s>0$, with $\Gamma_c:=\eps\rho_0(v)e^{-c|y|^2}$ and the margin $\delta/2$ over the energy exponent.
The factor $\rho_0(v)$ in $\Gamma$ is uniform in $v$, by $|\nabla\rho_0|\le L\rho_0$, and
$\rho_0(\xi-z/t)/\rho_0(\xi)\le e^{L|z|/t}$ is controlled by the margin.
\item \emph{Finiteness and continuity.} For $t>t_1$, these bounds make the truncated energies
$\mathsf E^{(k),\le N}_c[g_\eta]$ and the $w$-energies of Step~6 finite, $C^1$ in $t$, and Lipschitz
on $[t_1+s,T]$ uniformly in $c$, for every $s>0$. Indeed, every derivative of $g_\eta$ is bounded by
$C\Gamma_{q(t)-\delta/2}$, with a constant depending on its order. Therefore on each leaf the
integrands of $\mathsf e^{(k)}_\xi$, of $\mathsf d^{(k)}_\xi/(\eps\rho_0(\xi))$ and of their time
derivatives are bounded by $Ce^{-\delta|z|^2/2}$, uniformly in $\xi$, by
$\rho_0(\xi-z/t)\le\rho_0(\xi)e^{L|z|/t}$ and the Cauchy estimates for $1/\rho_0$. Differentiability
in $\xi$ and $t$ follows from the smoothness of $g_\eta$ and dominated convergence. The boundary
terms of the integrations by parts in $\xi$ vanish.

At $t_1$ the data have finite $y$-regularity, so Lipschitz continuity is not available there,
as noted above. Lemma~\ref{lem:CUtstar}(i)--(iv) holds for
these energies, with the same proof and without frames, nonlocal part or macroscopic hierarchy.
Part~(ii) is needed, since the proof of~(iii) integrates it. It holds for four reasons. In the
truncation $|\beta|\le N$, the hierarchy of $\phi$ is closed (\emph{Truncation} below). The factors
that are not components of the energy are of $y$-order $\le N_y/2\le N_y-6$, and they are bounded
by~(i). One $z$-derivative is absorbed by the viscous good term ($\partial_z=L_1-t^{-1}L_0$), and
the weight loss is controlled by the margin $\delta$.
In~(iii) the mollified data converge in the cell norms at the exponent $q_d-\delta$
because the bound on the data holds at $q_d$.

Thus each energy is continuous on $[t_1,T_\eta)$, $\sup_c$ of it is continuous on $(t_1,T_\eta)$,
and
\begin{equation}\label{eq:disp-usc}
\limsup_{t\downarrow t_1}\sup_c\mathsf E^{(k),\le N}_c(t)\le\sup_c\mathsf E^{(k),\le N}_c(t_1).
\end{equation}
The left side of (B3) involves $y$-orders $\le4\le N_y-6$ and is Lipschitz on
$[t_1,T_\eta)$, uniformly in $\xi$, by part~(i) of the same lemma.
\item \emph{The viscous term.} In the variables $(z,\xi)$,
$\Delta_y+\Delta_v=\sum_i(L_{0,i}^2+L_{1,i}^2)$ with $L_0=\partial_\xi$ and $L_1=D_y$. These fields
have constant coefficients in $(z,\xi)$ and commute with $D_y$, $\partial_\xi$ and
$\partial_t-t^{-2}z\cdot\nabla_\xi$ (Step~1). Since the viscosity acts on $G=\eps\rho_0(\xi)\phi$,
\eqref{eq:disphier} acquires the term $\eta_v
D_y^\alpha\partial_\xi^\beta\big[\sum_L\big(L^2\phi+2(L\log\rho_0)L\phi+(L^2\rho_0/\rho_0)\phi\big)\big]$,
with $L\log\rho_0\in\{\nabla\log\rho_0,t^{-1}\nabla\log\rho_0\}$.

(i)~For $\eta_v\sum_LL^2\phi_{\alpha\beta}$, we integrate by parts once against $\chi_c^2m_t^2$, as
in point~(2) of the same proof, and we obtain the good term
$-\eta_v\iint\chi_c^2m_t^2\sum_L|L\phi_{\alpha\beta}|^2$ and the error
$\eta_v\iint\chi_c^2m_t^2(32q^2|z|^2+C\rad_1^{-2})|\phi_{\alpha\beta}|^2$.

(ii)~The terms with no derivative on the coefficients are bounded by half of the good term plus
$C\eta_v(L^2+C_{\rho_0}\Lambda^{-2})\iint\chi_c^2m_t^2|\phi_{\alpha\beta}|^2\le
C\eta_v\rad_1^{-2}\iint\chi_c^2m_t^2|\phi_{\alpha\beta}|^2$.

(iii)~In the terms where $j\ge1$ derivatives fall on a coefficient, the factor $L\phi$ is at the
level $|\beta|-j+1$. By~\eqref{eq:dispshift} it carries $(4\rad_1/\Lambda)^j(j+1)^5$ against the
weight $w_{|\beta|-j}$, and by Young's inequality we bound these terms by one quarter of the good
terms at the lower levels and $C\eta_v(\rad_1/\Lambda)\sum_{k'\le k}\mathsf E^{(k'),\le N}_c$, using
$C\rad_1/\Lambda\le\frac18$. The zero-order ones are as in Step~1.

For $\eta_v\le C_q/(32q_d^2T_0^2)$ the term $32\eta_v q^2|z|^2$ is absorbed by the part of the
weight term that Step~2 leaves free, for $m_t$ and for $m_t^{1/2}$. The rest is a Gronwall rate
$C\eta_v\rad_1^{-2}$, with time integral $\le\frac1{64}$ for $\eta_v\le\rad_1^2/(C_\eta T_0)$.
In Step~5 the coefficient of $E_k$ becomes $\le\frac7{32}+\frac1{64}<\frac14$.
\item \emph{The free profile.} The viscosity raises no level. The free profile does not solve the
viscous equation: for $g_\eta$, $w=\phi-\phi^{\rm fr}$ satisfies the equation of Step~6 with the
additional source $\eta_v\rho_0^{-1}\sum_LL^2(\rho_0\phi^{\rm fr})$, whose $w$-energy of $y$-order
$k\le N_y-2$ is $\le C\eta_v^2\rad_1^{-4}E^{\rm fr}$ by the bounds on the free profile.
Consequently, for $g_\eta$, Step~6 yields $\sup_c\mathsf E^{(k)}_c[w](t)^{1/2}\le C_w^{1/2}\eps
T_0^{\kappa_\gamma}+C'\eta_v T_0\rad_1^{-2}(E^{\rm fr})^{1/2}$, and (B3) and the closeness of
Step~3 hold if moreover $CC'\eta_v T_0\rad_1^{-2}(E^{\rm fr})^{1/2}\le\frac14$, which is part of
the definition of $C_\eta$. The extra term vanishes as $\eta_v\to0$.
\end{itemize}
Hence Steps~1--7 hold for $g_\eta$ with the same constants, uniformly in $\eta_v$ and $\delta$,
except that Step~6 carries the additional term just stated.

\emph{Existence of $g_\eta$ on $[t_1,T_0]$.} Let $[t_1,T_\eta)$ be the interval of existence. On
$[t_1,\min\{T',T_\eta\})$, with $T'$ the end of the bootstrap interval of Step~7, we deduce from
(B1), Lemma~\ref{lem:TWcalc}(i) on the cubes of side $\rad_1$, Sobolev's inequality in $z$ and
$\rho_0(\xi)\le e^{L|y|/t}\rho_0(v)$ that $g_\eta\le A\,\eps\rho_0(v)e^{-(q_d/4)|y|^2}$, with
$A=A(E_0,E_1,E_2)$. This bounds $\ab[g_\eta]$, $\cb[g_\eta]$ and the constant of the null structure
of the linear lemma independently of the barrier constant, but it does not bound $
\mathcal N_0(g_\eta;c_\eta)$, since $c_\eta(t)>q(t)$. On each interval $[t_1+j,t_1+j+1]$,
we apply part~(iii) of that lemma with $h=g_\eta$ and get $g_\eta\le C_je^{K_{\rm
b}(t-t_1-j)}\Gamma_{c_\eta(t)}$, with $K_{\rm b}=K_{\rm b}(A,T_0)$. Then the continuation~(iv) for
the regularized problem, restarted at $t_1+j$ (the rate $c_\eta$ satisfies
$1/c_\eta(t)=1/c_0+4\eta_v(t-t_1)$ for every restart), continues $g_\eta$. Hence $T_\eta>T'$, as in
the proof of Corollary~\ref{cor:CUtstar}, and by the continuity argument of Step~7 we obtain
$T'=T_0$.

\emph{The limit.} Let $(\eta_v,\delta)=(\eta_j,\delta_j)\to0$. By Lemma~\ref{lem:LCd}, with the
$\eta_v$-uniform bounds of (B1) and Sobolev's inequality for the coefficients in the bound on
$\partial_t\partial^\alpha_y\partial^\beta_vg_\eta$, a subsequence converges in
$C([t_1,T_0];C^{N_y-4}_{\rm loc})$, and $\partial^\alpha_y\partial^\beta_vg_\eta(t)$ converges
weakly in $L^2_{\rm loc}$ for $|\alpha|\le N_y$ and all $\beta$.

The limit $g\ge0$ is a classical solution with the data $g(t_1)$, and it is classical up to $t_1$:
since $g\in C([t_1,T_0];C^{N_y-4}_{\rm loc})$, with Gaussian bounds on its derivatives that are
uniform in $t$ by Lemma~\ref{lem:LCdata}, as in the paragraphs after the compactness lemma, the
right
side of the equation is continuous on $[t_1,T_0]\times\R^6$. Therefore
$g(t)=g(t_1)+\int_{t_1}^t(\text{right side})$, and $g$ is $C^1$ in $t$ up to $t_1$.

The layer solution is in $C([0,t_1];C^2_{\rm loc})$ and solves the equation, by
Lemma~\ref{lem:FPiter}(iii) and its proof, or by Theorem~\ref{thm:EPlayer} at $\gamma=-2$. Its left
derivative at $t_1$ is the same right side, evaluated at $g(t_1)$. Thus the two pieces form a
classical solution on $[0,T_0]$.

The same uniform Gaussian bounds, which do not depend on $\eta_v$, control the tails of the
coefficients, so $\ab[g_\eta]\to\ab[g]$ and $\cb[g_\eta]\to\cb[g]$ locally uniformly along the
subsequence. The cell energies, the $w$-energies and $\int\kappa\mathsf D^{(k)}_c$ (the latter in
$L^2((t_1,T_0)\times K)$) are weakly lower semicontinuous, with the converging weights
$e^{2(q-\delta_j)|z|^2}$, which implies that $g$ satisfies the bounds of the proposition. Uniqueness
holds among classical solutions with Gaussian bounds on the derivatives of order $\le2$, by
Lemma~\ref{lem:LWuniq}.

\emph{Truncation.} For $N\ge4$ the hierarchy of $\phi$ below is closed in the truncation $|\beta|\le
N$: the transport and the coefficient terms do not raise the level, and the cube Sobolev inequality
uses the levels $\le\lfloor|\beta|/2\rfloor+2\le N$. The $w$-energies of Step~6 are truncated at
$|\beta|\le N-4$. We write the estimates for $g$, and the viscous terms are treated as just said.

\emph{Step 1: the equation.} By Proposition~\ref{prop:leaf}, $\partial_tG=\kappa
Q(G,G)+t^{-2}z\cdot\nabla_\xi G$, where $Q$ acts in $z$ at fixed $\xi$. Dividing by
$\eps\rho_0(\xi)$, we get
\begin{equation}
\partial_t\phi=t^{-2}z\cdot\nabla_\xi\phi+t^{-2}(z\cdot\nabla\log\rho_0)\,\phi+\kappa\eps\rho_0(\xi)\,Q(\phi,\phi).
\end{equation}
The operator $D_y$ commutes with $\partial_t-t^{-2}z\cdot\nabla_\xi$: indeed
$[\partial_t,t^{-1}\partial_\xi]=-t^{-2}\partial_\xi$ and
$[-t^{-2}z\cdot\nabla_\xi,\partial_z]=t^{-2}\nabla_\xi$. It also satisfies
$D_yQ(f,h)=Q(D_yf,h)+Q(f,D_yh)$, since $Q$ commutes with translations in $z$, acts at fixed $\xi$
and is bilinear. Hence we have
\begin{equation}\label{eq:disphier}
(\partial_t-t^{-2}z\cdot\nabla_\xi)\phi_{\alpha\beta}=t^{-2}D_y^\alpha\partial_\xi^\beta\big[(z\cdot\nabla\log\rho_0)\phi\big]+\kappa\eps\,D_y^\alpha\partial_\xi^\beta\big[\rho_0(\xi)Q(\phi,\phi)\big].
\end{equation}
The coefficients $\nabla\log\rho_0$ and $\rho_0$ are differentiated only through $\partial_\xi$ and
the part $t^{-1}\partial_\xi$ of $D_y$. By (D3), $\rho_0$ has no zero on $\mathbb S_\Lambda$, so
$\log\rho_0$ is holomorphic on $\Delta(\xi,\Lambda/2)$, with
$|\log|\rho_0(\xi+\zeta)|-\log\rho_0(\xi)|\le\log A_F+L\Lambda$. The Cauchy estimates in the
Borel--Carath\'eodory form, which use only the real part, imply
$|\partial^j\nabla\log\rho_0|+|\partial^j\rho_0|/\rho_0\le C_{\rho_0}\,j!(\Lambda/4)^{-j}$ with
$C_{\rho_0}=C_{\rho_0}(\Dc)$. The factor $\partial_z$ of $D_y$ falling on $z$ produces the terms
$t^{-2}(\partial_i\log\rho_0)\phi_{\alpha-e_i,\beta}$, of lower $y$-order and without $|z|$. A term
in which $j\ge1$ derivatives fall on a coefficient involves $\phi$-factors of total $\xi$-level
$|\beta|-j'$, $j'\le j$, and of $y$-order $\le|\alpha|$. In the weights $w_n$, it carries the factor
\begin{equation}\label{eq:dispshift}
\frac{w_n\binom nj\,j!\,(\Lambda/4)^{-j}}{w_{n-j}}\le\Big(\frac{4\rad_1}\Lambda\Big)^j(j+1)^5,
\end{equation}
and a factor $t^{-1}(\Lambda/4)^{-1}\le4/\Lambda$ for each $t^{-1}\partial_\xi$ of $D_y$ that falls
on a coefficient. Since $\rad_1\le\Lambda/8$, the sum over $j\ge1$ of these factors is
$\le C\rad_1/\Lambda$, and Young's inequality $\ell^1*\ell^2\to\ell^2$ applies. There is no loss
of level and no loss of $y$-order.

\emph{Step 2: transport and weight.} We differentiate $\mathsf E^{(k),\le N}_c$. The weight
contributes $2q'(t)\!\int|z|^2m_t^2|\phi_{\alpha\beta}|^2$, with $q'=-C_qt^{-2}$. For the transport
$t^{-2}z\cdot\nabla_\xi$, after integration by parts in $\xi$ (the weight $m_t$ does not depend on
$\xi$), we have
\begin{equation}
-t^{-2}\rad_1^{-3}\!\int(z\cdot\nabla\chi_c^2)\,m_t^2|\phi_{\alpha\beta}|^2\le\frac{2t^{-2}}{\rad_1}\rad_1^{-3}\!\int\chi_c^2\,|z|\,m_t^2|\phi_{\alpha\beta}|^2,
\end{equation}
since $|\nabla\chi_c^2|\le2\chi_c^2/\rad_1$. This is Lemma~\ref{lem:CUtr} in the fixed frame, and
there is no loss of level. By Step~1 and Young's inequality, the first term of the hierarchy is
bounded by
$t^{-2}C(1+C\rad_1/\Lambda+C/\Lambda)\sum\!\int(1+|z|)m_t^2|\phi_{\alpha'\beta'}|^2$, summed over
the same $y$-order and the lower ones. In the coefficient terms the factor $|z|$ is put on the
target of $y$-order $k$, whose weight term absorbs $|z|^2$, while the factor of lower order is taken
without $|z|$. With $|z|\le\eta'|z|^2+1/(4\eta')$ and the choice
$C_q:=2C(1+\rad_1^{-1}+\Lambda^{-1})$, the weight term controls all terms with $|z|$, with a factor
$2$ of margin for the half weight $m_t^{1/2}$ of Step~6, and what remains is
\begin{equation}
\le C_{\rm tr}^{\rm lay}t^{-2}\sum_{k'\le k}\mathsf E^{(k'),\le N}_c,\qquad C_{\rm tr}^{\rm
lay}=C(\rad_1,\Lambda,L,A_F).
\end{equation}

\emph{Step 3: the bootstrap.} Let
\begin{equation}\label{eq:Kfr}
\mathcal K_{\rm fr}:=\sup_{t\ge1,\xi}\big(\sum_{j\le2}\|m_{q_d}\nabla_z^j\phi^{\rm
fr}\|+\sup_ze^{(q_d/8)|z|^2}|\phi^{\rm fr}|\big),
\end{equation}
which is a data constant by the bounds on the free profile and (D1). We assume, on $[t_1,T']$,
\begin{equation}
\text{(B1)}\ \ \sup_c\mathsf E^{(k),\le N}_c\le E_k\ (k\le N_y),\qquad\text{(B2)}\ \ \sup_Q
\mathsf a^{(0),\le N}_{Q^*,\ge1}\le c_{\sharp,\ell},\qquad\text{(B3)}\ \
\sup_\xi\|m_t^{1/2}w_{\alpha0}(t,\cdot,\xi)\|\le1\ (|\alpha|\le4),
\end{equation}
with $E_k$ and $c_{\sharp,\ell}$ fixed in Step~5. By (B3), Sobolev's inequality in $z$ as in (F2),
and $\nabla_z=D_y-t^{-1}\partial_\xi$ (the parts with $t^{-1}\partial_\xi$ are bounded by (B1) and
Lemma~\ref{lem:TWcalc}(i) by $C(\rad_1t)^{-1}(E_0+E_1+E_2)^{1/2}\le1$, since $t\ge
t_1\ge\rad_1^{-1}(1+C(E_0+E_1+E_2)^{1/2})$), the background $\phi(t,\cdot,\xi)$ satisfies
\begin{equation}
\sum_{j\le2}\|m_t^{1/2}\nabla_z^j\phi\|+\sup_ze^{(q/8)|z|^2}|\phi|\le\mathcal K:=C(\mathcal K_{\rm
fr}+1),
\end{equation}
uniformly in $\xi$. Thus (G3) holds with $K_\ell:=C\mathcal K$. Step~6 below rests on (B1) and
Steps~1--2 alone, not on the hypodissipativity, and it yields $\sup_\xi\|m_t^{1/2}w_{\alpha0}\|\le
CC_w^{1/2}\eps T_0^{\kappa_\gamma}$ on $[t_1,T']$. We use this bound here. The moments of
$\phi(t,\cdot,\xi)$ are within
$C(q_d)CC_w^{1/2}\eps T_0^{\kappa_\gamma}$ of those of $\phi^{\rm fr}(t,\cdot,\xi)$, which are
within $CL(L+1)/t$ of those of $F_0(\cdot,\xi)$, by Lemma~\ref{lem:LWfree}(iii). Hence for
$\eps\le\eps_0(T_0)$ and $t_1\ge CL(L+1)$, (G1) holds with $M_\pm=(\frac12,\frac32)$, $|y_0|\le L$
and $\theta\le\theta_+$, and by Lemma~\ref{lem:RGellip} we get (G2) with $c_G^{\rm
lay}=c(\gamma,\frac12,E_{\rm mom}',K_\ell)$, with the second-moment bound $E_{\rm
mom}':=2(L^2+3\theta_+)$. Below we write $c_h,C_h,c_G$ for the layer constants $c_h^{\rm
lay},C_h^{\rm lay},c_G^{\rm lay}$. With the choice of $t_1$ in Step~5, all these constants depend
only on $\Dc$. In particular, they do not depend on the $E_k$ or on $\rad_1$.

\emph{Step 4: the collision terms.}
\begin{itemize}
\item \emph{The linear terms.} For $(\alpha,\beta)\ne0$ the terms of the hierarchy without
derivatives on $\rho_0$ and with one factor $\phi_{00}=\phi$ are
$\rho_0\big(Q(\phi,\phi_{\alpha\beta})+Q(\phi_{\alpha\beta},\phi)\big)$, while for
$(\alpha,\beta)=0$ the term is $\rho_0Q(\phi,\phi)$. Applying Lemma~\ref{lem:Whypo} with
$G=\phi(t,\cdot,\xi)$ and $m=m_t$ ($q(t)\le q_d\le(8(L^2+3\theta_+))^{-1}$, and its constants are
uniform for $q\in[q_d/2,q_d]$ by Lemma~\ref{lem:RFtools}(vi), since $q_d\le\frac18$), we obtain
$\mathrm{Re}\!\int Q(\phi,u)\bar um_t^2\le-c_h\|u\|^2_{Y_2(m_t)}+C_h\|m_tu\|^2$, with $c_h,C_h$
depending only on $(\gamma,L,\theta_+,K_\ell,c_G,q_d)$, hence only on $\Dc$. For the second half we
write $Q(u,\phi)=\ab[u]:\nabla^2\phi-\cb[u]\phi$ as in (L2) of the proof of
Proposition~\ref{prop:ACpf}. By \cite[Lemma~4.2(i)]{CM}, we have $|\ab[u]|\le C\jap
z^{\gamma+2}\|m_tu\|$. Moreover, $\cb[u]=-8\pi u$ for
$\gamma=-3$, and for $\gamma>-3$ the near part of its kernel is bounded on $L^2$ and the far part
pointwise by $\|u\|_{L^1}$. By Lemma~\ref{lem:RFtools}(ii) (whose near part of $\cb$ uses the shift
of the weight over distances $\le1$, controlled by $\sup_ze^{(q/8)|z|^2}|\phi|\le\mathcal K$),
\begin{equation}
\Big|\int Q(u,\phi)\bar um_t^2\Big|\le C\big(\|m_t\nabla^2_z\phi\|+\|m_t\phi\|+
\mathcal K\big)\|m_tu\|^2\le C_E\|m_tu\|^2,
\end{equation}
where, by (B1), the cube Sobolev inequality on the cubes of side $\rad_1$,
$\nabla_z=D_y-t^{-1}\partial_\xi$ and $t\ge\rad_1^{-1}$,
$\sup_\xi(\|m_t\phi\|+\|m_t\nabla^2_z\phi\|)\le C(E_0+E_1+E_2)^{1/2}$, and $C_E:=C(
\mathcal K+(E_0+E_1+E_2)^{1/2})$. Hence the linear terms contribute
$\le\kappa\eps\rho_0(-c_h\|\phi_{\alpha\beta}\|^2_{Y_2}+C_L\|m_t\phi_{\alpha\beta}\|^2)$ with
$C_L:=C_h+C_E$. The constant $C_L$ depends on the $E_k$, and through $a(t)$ and $\rho_+$ it affects
only $\eps_0(T_0)$.
\item \emph{The quadratic terms of the same $y$-order.} The terms without derivatives on $\rho_0$
and with $0<(\alpha',\beta')<(\alpha,\beta)$ are bounded by Proposition~\ref{prop:trilinG}(iii),
second form ($X(m_t)$ is equivalent to $L^2(m_t)$). If $\alpha'=0$ or $\alpha''=0$, the factor of
$y$-order $0$ has $\xi$-level $\ge1$. Following the proof of Proposition~\ref{prop:CUtri}, at the
scale $\rad_1$ and with $\omega\equiv1$ (the weight $m_t$ does not depend on $\xi$, and
$\rho_0(\xi)/\rho_0(\xi')\le e$ on each $Q^*$ since $3\sqrt3L\rad_1\le1$), with
Lemma~\ref{lem:algebra} for $p=5$ in $\beta$ and $\binom\alpha{\alpha'}\le C(N_y)$ in $\alpha$,
we bound their contribution by
\begin{equation}
C_N^{\rm lay}\kappa\Big(\sup_Q\mathsf a^{(0)}_{Q^*,\ge1}\Big)^{1/2}\mathsf D^{(k)}_c+C_N^{\rm
lay}\kappa\Big(\sup_Q\mathsf a^{(k)}_{Q^*}\Big)^{1/2}\big(
\mathsf D^{(0)}_{c,\ge1}\big)^{1/2}\big(\mathsf D^{(k)}_c\big)^{1/2}.
\end{equation}
By (B2) with $C_N^{\rm lay}c_{\sharp,\ell}^{1/2}\le c_h/8$, where $C_N^{\rm lay}=C(N_y,C_\tau(q_d))$
is the constant of this bound (the layer density has no factor $\eta$), the first part is
$\le\frac{c_h}8\kappa\mathsf D^{(k)}_c$. By (B1), the second one is $\le\frac{c_h}8\kappa
\mathsf D^{(k)}_c+C\kappa E_k\mathsf D^{(0)}_{c,\ge1}/c_h$. (For $k=0$ both parts are of the first kind.)
\item \emph{The quadratic terms of lower $y$-orders.} If $0<|\alpha'|,|\alpha''|<k$, both factors
have lower $y$-orders, and from the same bound we get $\le\frac{c_h}8\kappa
\mathsf D^{(k)}_c+C\kappa\sum_{k',k''<k}E_{k'}\mathsf D^{(k'')}_c/c_h$.
\item \emph{The terms with derivatives on $\rho_0$.} By~\eqref{eq:dispshift} they carry the factor
$\vartheta_1:=C(\rad_1/\Lambda+\Lambda^{-1}t_1^{-1})$, and their $\phi$-factors may include
$\phi_{00}$, for which $\sup_\xi\|m_t\phi\|\le CE_0^{1/2}$ as in the linear terms. From the same
bounds we obtain $\le\frac{c_h}8\kappa
\mathsf D^{(k)}_c+C\kappa\big(E_0+\sum_{k'<k}E_{k'}\big)\sum_{k''<k}
\mathsf D^{(k'')}_c/c_h+C\vartheta_1^2\kappa E_k\mathsf D^{(0)}_c/c_h$, once $C\vartheta_1^2E_0\le
c_h^2/8$. Here the last term comes from a factor of $y$-order $k$ paired with a factor of $y$-order
$0$, and the
factors of order $k$ in $Y_2$ against $\phi_{00}$ in $X$ are absorbed by $\frac{c_h}8\kappa
\mathsf D^{(k)}_c$.
\end{itemize}

\emph{Step 5: closing (B1) and (B2).} Summing up, we have, for each $k\le N_y$,
\begin{equation}
\frac d{dt}\mathsf E^{(k),\le N}_c+\frac{c_h}2\kappa\mathsf D^{(k),\le N}_c\le a(t)\sum_{k'\le
k}\mathsf E^{(k'),\le N}_c+\frac C{c_h}\kappa E_k\big(\mathsf D^{(0)}_{c,\ge1}+\vartheta_1^2
\mathsf D^{(0)}_c\big)+\frac C{c_h}\kappa\Big(E_0+\sum_{k'<k}E_{k'}\Big)\sum_{k''<k}\mathsf D^{(k''),\le
N}_c,
\end{equation}
with $a(t):=C_{\rm tr}^{\rm lay}t^{-2}+C_L\rho_+\kappa$, $\int_{t_1}^{T_0}a\le C_{\rm tr}^{\rm
lay}/t_1+CC_L\eps T_0^{\kappa_\gamma}$. For $k=0$ the term
with $E_k$ is absent, since both quadratic parts are then of the
first kind. The part $\ge1$ of $y$-order $0$ satisfies the same inequality (without the last term),
with two additional sources from the terms with derivatives on a coefficient that take the factor
$\phi_{00}$, namely $t^{-2}C(\rad_1/\Lambda)^2E_0$ and $C\kappa(\rad_1/\Lambda)^2E_0
\mathsf D^{(0)}_c/c_h$.

We define $E_k$ recursively, depending only on $\Dc$ and on structural constants:
\begin{equation}
E_0:=4c_\chi(2E^{\rm fr}+1),\qquad E_k:=4c_\chi(2E^{\rm
fr}+1)+\frac{8C}{c_h^2}\Big(E_0+\sum_{k'<k}E_{k'}\Big)\sum_{k''<k}E_{k''}\quad(k\ge1).
\end{equation}
They do not depend on $\rad_1$: the initial values are bounded by the bounds on the free profile
and on the data, uniformly for $\rad_1\le\Lambda/8$ under the $\eps$-condition $C_{\rm
lay}\delta_{y,\ell}^{-2N_y}\eps^2t_1^2\le\frac12$ below, and $c_h$ is a data constant (Step~3). Then
we fix $c_{\sharp,\ell}$ by $C_N^{\rm lay}c_{\sharp,\ell}^{1/2}=c_h/8$, and put
\begin{equation}\label{eq:lambda1}
\rad_1:=\min\{\Lambda/8,\rad_0^{\rm lay}/4,1/(6L),\min\{1,\theta_-\}/(1+C_{\mu_0}),c_{\rm
disp}\Lambda\}
\end{equation}
(the fourth entry is the cap of the paragraph \emph{Parameters}, used in
Corollary~\ref{cor:LWPa}(a4), a constant independent of $c_{\rm disp}$), with $c_{\rm
disp}=c_{\rm disp}(\Dc)$ the largest value for which $C\rad_1/\Lambda\le\frac18$ and the
following three conditions hold, with $\vartheta_1$ replaced by $2C\rad_1/\Lambda$, which bounds
it for $t_1\ge\rad_1^{-1}$, and in particular $\rad_1\ge c(\Dc)\Lambda$:
\begin{equation}
C\vartheta_1^2E_0\le\frac{c_h^2}{32},\qquad 27e^4C(\rad_1/\Lambda)^2E_0\Big(1+\frac{E_0}{c_h^2}\Big)\le\frac{c_{\sharp,\ell}}4,\qquad\frac
C{c_h^2}C(\rad_1/\Lambda)^2E_0\Big(1+\frac{E_0}{c_h^2}\Big)\le\frac1{16},
\end{equation}
and $t_1\ge\max\{1,\rad_1^{-1}(1+C(E_0+E_1+E_2)^{1/2}),4C_q/q_d,CL(L+1),64(N_y+1)C_{\rm tr}^{\rm
lay},C_{\rm tr}^{\rm lay}(\Lambda/\rad_1)^2\}$. Finally let $\eps\le\eps_0(T_0)$ with
$(N_y+1)CC_L\eps T_0^{\kappa_\gamma}\le\frac1{64}$ and $CC_L\eps
T_0^{\kappa_\gamma}\le(\rad_1/\Lambda)^2$, $C_{\rm
lay}\delta_{y,\ell}^{-2N_y}\eps^2t_1^2\le\min\{\frac12,(\rad_1/\Lambda)^2\}$, $CC_w^{1/2}\eps
T_0^{\kappa_\gamma}\le\frac12$, and $\eps\le\eps_0'$ of Theorem~\ref{thm:layer}, or of
Theorem~\ref{thm:EPlayer} at $\gamma=-2$. The constants are chosen in the order $c_h$, then $E_k$
and $c_{\sharp,\ell}$, then $\rad_1$, then $C_q$ and $t_1$, then $\eps_0(T_0)$.

On $[t_1,T']$ we first consider $y$-order $0$ and its part $\ge1$. Integrating the inequalities with
(B1), with $\int_{t_1}^{T_0}a\le C_{\rm tr}^{\rm lay}/t_1+CC_L\eps
T_0^{\kappa_\gamma}\le2(\rad_1/\Lambda)^2$ for the rate term of the part $\ge1$ (by the
conditions on $t_1$ and $\eps$ above), and using $\mathsf E^{(0)}_{c,\ge1}(t_1)\le c_\chi
C(\rad_1/\Lambda)^2E^{\rm fr}+C\eps^2t_1^2$, we obtain
\begin{equation}
\mathsf E^{(0)}_{c,\ge1}(t)+\frac{c_h}2\!\int_{t_1}^t\kappa\mathsf D^{(0)}_{c,\ge1}\le
C(\rad_1/\Lambda)^2E_0\Big(1+\frac{E_0}{c_h^2}\Big)\le\frac{c_{\sharp,\ell}}{4\cdot27e^4}.
\end{equation}
Since $\mathsf a_{Q^*}\le27e^4\sup_{c'}\mathsf E_{c'}$ by Lemma~\ref{lem:CUcomp}, applied with
$\omega\equiv1$ on each of the $27$ cubes of $Q^*$, (B2) holds with $c_{\sharp,\ell}/4$. Next, by
induction on $k$, integrating the inequality of order $k$ and using (B1) on the right, we get
\begin{equation}
\mathsf E^{(k)}_c(t)+\frac{c_h}2\!\int_{t_1}^t\kappa\mathsf D^{(k)}_c\le
\mathsf E^{(k)}_c(t_1)+E_k\Big(\int a+\frac C{c_h}\!\int\kappa\big(
\mathsf D^{(0)}_{c,\ge1}+\vartheta_1^2\mathsf D^{(0)}_c\big)\Big)+\frac
C{c_h}\Big(E_0+\sum_{k'<k}E_{k'}\Big)\sum_{k''<k}\frac{2E_{k''}}{c_h}\le\frac{E_k}2 ,
\end{equation}
by the choices above. Indeed, the step $k=0$ yields $\int\kappa\mathsf D^{(0)}_c\le2E_0/c_h$. Since
$E_k$ is nondecreasing in $k$, the contribution of the term $a\sum_{k'\le k}\mathsf E^{(k')}$ is
bounded by $(N_y+1)E_k\int a$. Thus the coefficient of $E_k$ is
$\le\frac1{32}+\frac18+\frac1{16}\le\frac14$, and
the remaining terms are $\le\mathsf E^{(k)}_c(t_1)+\frac14(E_k-4c_\chi(2E^{\rm
fr}+1))\le\frac{E_k}4$, since $\mathsf E^{(k)}_c(t_1)\le c_\chi(2E^{\rm fr}+1)$. Therefore (B1)
holds with $E_k/2$, and $\int_{t_1}^{T'}\kappa\mathsf D^{(k)}_c\le E_k/c_h$.

\emph{Step 6: closeness to the free profile, and (B3).} The free profile solves the hierarchy
without the collision term, so $w=\phi-\phi^{\rm fr}$ satisfies the same equation with the source
$\kappa\eps D_y^\alpha\partial^\beta_\xi[\rho_0Q(\phi,\phi)]$. We measure $w$ at the radius
$\rad_1/2$, with the weight $m_t^{1/2}$ and at $y$-orders $k\le N_y-2$. Steps~1--2 apply to $w$
without change (the weight $m_t^{1/2}$ has the exponent $q/2$, and $C_q/2$ still absorbs the terms
with $|z|$, by the margin in Step~2). For the source, it follows from
$Q(f,h)=\ab[f]:\nabla^2h-\cb[f]h$ that $\|m_t^{1/2}Q(f,h)\|\le
C\|m_tf\|\big(\|m_t\nabla_z^2h\|+\|m_th\|\big)+C\sup|e^{(q/8)|z|^2}f|\,\|m_th\|$, where the factor
$m_t^{1/2}$ controls the shifts of the weight. Here $\nabla_z^2=(D_y-t^{-1}\partial_\xi)^2$ raises
the $y$-order by at most two, to $\le N_y$, and the $\xi$-level by at most two, which loses
$C(n+2)^2\rad_1^{-2}$ against $w_n(\rad_1/2)/w_n(\rad_1)\le2^{-n}$. With
Lemma~\ref{lem:algebra} in $\beta$ (the supremum taken on the factor of lower order, by
the cube Sobolev inequality and (F2)), the source is bounded in the $w$-energies of $y$-order $k$ by
$C\kappa\rho_+\big(\sum_{k'\le N_y}E_{k'}\big)$, with $\rho_+=2\eps A_F$. Hence
\begin{equation}
\frac d{dt}\mathsf E^{(k)}_c[w]\le C_{\rm tr}^{\rm lay}t^{-2}\sum_{k'\le k}
\mathsf E^{(k')}_c[w]+C\kappa\eps\Big(\sum_{k'\le N_y}E_{k'}\Big)\Big(\sum_{k'\le k}
\mathsf E^{(k')}_c[w]\Big)^{1/2},
\end{equation}
and by Gronwall's inequality, with the bound on the data at $t_1$, we obtain $\sup_c
\mathsf E^{(k)}_c[w](t)^{1/2}\le C_{\rm lay}^{1/2}\delta_{y,\ell}^{-N_y}\eps t_1+C(E)\eps
T_0^{\kappa_\gamma}\le C_w^{1/2}\eps T_0^{\kappa_\gamma}$, since $t_1$, $\delta_{y,\ell}$ and
$C_{\rm lay}$ are data constants (fixed before $C_w$) and
$T_0^{\kappa_\gamma}\ge1$ (for the data of Theorem~\ref{thm:t1}, $C_{t_1}\eps$ takes the place of
$\eps t_1$, so $C_w$ depends also on $C_{t_1}$). By the cube Sobolev inequality in $\xi$, the left
side of (B3) is $\le CC_w^{1/2}\eps T_0^{\kappa_\gamma}\le\frac12$ for $\eps\le\eps_0(T_0)$. Thus
(B3) holds with $\frac12$.

\emph{Step 7: continuation.} At $t_1$, (B1)--(B3) hold with the improved constants, by
Lemma~\ref{lem:dispfree}, \eqref{eq:dispdata} and Proposition~\ref{prop:FPtr}. For the data of
Theorem~\ref{thm:t1} they follow from its hypotheses, (B2) from the same bound on the data, and (B3)
from it and the cube Sobolev inequality. On every interval $[t_1,T']$ where they
hold, Steps~3--6 improve them, and, for $g_\eta$, the suprema in (B1) and (B3) have the continuity
properties of the \emph{Construction}. (B2) is a bound on cube energies, for which no continuity is
claimed, and it needs none: Step~5 derives it from the cell bound $\sup_c\mathsf E^{(0),\le
N}_{c,\ge1}\le c_{\sharp,\ell}/(4\cdot27e^4)$, whose left side has the continuity properties of the
\emph{Construction}, through the bound $\mathsf a_{Q^*}\le27e^4\sup_c\mathsf E_c$ of
Lemma~\ref{lem:CUcomp}. Thus (B2) persists slightly beyond every time at which this cell bound
holds. Gronwall's inequality is used as in Lemma~\ref{lem:LCgron}(i), integrated from $s>t_1$ with
$s\downarrow t_1$. Hence they hold on $[t_1,T_0]$ for $g_\eta$, uniformly in $N$, $\eta$ and
$\delta$. Letting $N\to\infty$ (monotone convergence) and then passing to the limit as in the
construction, we obtain the proposition. This completes the proof.
\end{proof}

\begin{cor}[(LWP)(a): the data at $T_0$]\label{cor:LWPa}
For $v\in\R^3$ let $\hat F_0^v:=T^{-1}_{(1,y_{0,\infty}(v),\theta_\infty(v))}F_0(\cdot,v)$, the
normalized fiber of the free profile at the velocity $v$ (for product data it does not depend on
$v$). Fix $(C_{\rm dat},\beta,A_7)$ and $C'_D$ such that, for every $v$, $\hat
F_0^v\in\mathfrak K^7(C_{\rm dat}/2,2\beta,A_7/2)$,
\begin{equation}\label{eq:LWPa-dv}
|\partial_v(\eps\rho_0(v)F_0(y,v))|\le\frac12C'_D\eps\rho_0(v)e^{-2\beta|y-y_{0,\infty}(v)|^2/\theta_+}
\end{equation}
and $2\beta\le\theta_-q_d/16$. By (D1), (D3), (D4) and the uniform bounds on $y_{0,\infty}$ and
$\theta_\infty$ these can be chosen depending only on $\Dc$. Assume $q_0\le\theta_-q_d/16$,
$\rad_0\le c_M\rad_1$ (with $c_M=c_M(\Dc)\le\frac12$ the constant of the Maxwellian bullet of
the proof) and $\delta_y\le1$. Then, for $T_0$ with
\begin{equation}\label{eq:ca4}
C(\Dc)T_0^{-1}\le c_{\rm a4}:=\frac18\min\{1,\theta_-\}(2+L/2+\theta_+)^{-2}(1+3G_\infty)^{-1}
\end{equation}
(used for (a4), including the bound on $A_{\rm ref}$ that (a1) uses) and $\eps\le\eps_0(T_0)$:
\begin{itemize}
\item[(a1)] at the bootstrap parameters $(\rad_0,\delta_y,q_0)$: $\mathfrak A_c(T_0)\le c_\chi
C_DB$ and
\begin{equation}\label{eq:LWPa-a1}
\mathfrak A_{c,\ge1}(T_0)\le c_\chi C_D\big((\rad_0/\Lambda)^2+\delta_y^2\big)B
\end{equation}
for every cell, with a data constant $C_D$ (the form used in Step~0 of Theorem~\ref{thm:CUboot} and
in Proposition~\ref{prop:UFU2});
\item[(a2)] the leaf parameters are within $C\eps T_0^{\kappa_\gamma}$ of those of free transport
(relatively for the mass);
\item[(a3)] $v\mapsto P_v$ is $C^1$, and $P_v$ and $\nabla_vP_v$ are within $C(\eps
T_0^{\kappa_\gamma})^\varkappa$, for some structural $\varkappa>0$, of the parameters
$(\eps\rho_0(v),y_{0,\infty}(v),\theta_\infty(v))$ of free transport and of their gradient
(relatively for the mass); hence (H1) holds with the constants $L$, $\theta_\pm$,
$\rho_+$; for every $v$, $F_{v,0}=T^{-1}_{P_v}g(T_0,\cdot,v)\in\mathfrak K^7(C_{\rm
dat},\beta,A_7)$; and $|\partial_vg(T_0,y,v)|\le C'_D\rho_ve^{-\beta|y-y_{0,v}|^2/\theta_+}$;
\item[(a4)]
\begin{equation}\label{eq:LWPa-a4}
\varpi'_g(T_0,\cdot)=(\nabla_\xi\ell_g,y_{0,g},\vartheta_g)(T_0,\cdot)
\end{equation}
extends holomorphically to $\mathbb S_{R_1}$, $R_1:=\rad_1/8$, with $|\varpi'_g(T_0)|\le L_1$
there, where $L_1=L_1(\Dc)$; and $\theta_-\le\theta_g(T_0)\le\frac23\theta_+$ and $|\nabla\log
A_{\rm ref}|\le L$ on $\R^3$, where $A_{\rm ref}=A_{p_g(T_0,\cdot)}$.
\end{itemize}
\end{cor}

\begin{proof}
\emph{(a1).} By Lemma~\ref{lem:TWcalc}(i) on the cubes of side $\rad_1$, we deduce from
Proposition~\ref{prop:disp} the pointwise bounds $\sum_n\tilde
w_n(\rad_1)^2\sum_{|\alpha|=k,|\beta|=n}\|m_{T_0}\phi_{\alpha\beta}(T_0,\cdot,\xi)\|^2\le C_{\rm
Sob}'E_k$, and $\le C(\rad_1/\Lambda)^2E_0$ for the part $|\beta|\ge1$ of $y$-order $0$. The
bootstrap energy at $\xi$ is built from the micro part
$T^{-1}_{p_g}\big(\partial^\alpha_y\partial^\beta_v\hat r\big)^{(x)}$ in the frame $p_g(T_0,\xi)$,
with the weight $\hat m_\xi$, the radius $\rad_0$ and the $y$-weight $\delta_y$, at $\omega=1$.
We compare term by term.
\begin{itemize}
\item \emph{Derivatives.} $\partial^\alpha_y\partial^\beta_vg$ restricted to the leaf is
$D^\alpha_y\partial^\beta_\xi G$, a sum over $\alpha'\le\alpha$, $\beta'\le\beta$ of Cauchy factors
of $\rho_0$ times $\eps\phi_{\alpha'\beta'}$, and the conversion is controlled
by~\eqref{eq:dispshift} with $\rad_0\le\rad_1/2$.
\item \emph{Weights in $\beta$.} $w_n(\rad_0)/\tilde
w_n(\rad_1)=(n+2)^2(\rad_0/\rad_1)^n\le(n+2)^22^{-n}$ is bounded, and is $\le
C\rad_0/\rad_1$ for $n\ge1$.
\item \emph{Weights in $\alpha$.} $\delta_y^{2|\alpha|}\le\delta_y^2$ for $|\alpha|\ge1$, and
$\sum_k E_k$ is a data constant.
\item \emph{Frame and Gaussian weight.} By (a2) below and Lemma~\ref{lem:LWfree}(iii),
$\frac12\theta_-\le\theta_g\le\theta_++1$ and $|y_{0,g}|\le L+1$ at $T_0$. Thus $\hat m_\xi^2\le
Ce^{(q_d/2)|z|^2}\jap z^{2N_0}\le Cm_{T_0}^2$, by $q_0\le\theta_-q_d/16$, and $T^{-1}_{p_g}$ is
bounded in these norms, with the mass $\rho_g(T_0,\xi)$ comparable to $\eps\rho_0(\xi)$. At $T_0$,
$A_{p_g}=A_{\rm ref}$ by the definition of $A_{\rm ref}$, so the factor $(A_{p_g}/A_{\rm
ref}(\xi))^2$ of the bootstrap densities equals $1$. Moreover, the ratios of $A_{\rm
ref}(\xi-z/T_0)$ to
$A_{\rm ref}(\xi)$ are $\le e^{L|z|/T_0}\le e^{2L|z|/T_0}$, by the bound $|\nabla\log A_{\rm
ref}|\le L$ of (a4), which is proved below without using (a1).
\item \emph{The Maxwellian part.} $p_g(T_0,\xi)$ are the moments of
$\eps\rho_0\phi(T_0,\cdot,\xi)$, so their $\xi$-derivatives are bounded at the radius $\rad_1$ by
the energies of $y$-order $0$, and $S_Q(p_g-p_g(\xi_Q))\le C(\Dc)\rad_0/\rad_1\le r/(8dC_{\rm
alg})$ for $\rad_0\le c_M\rad_1$, $r$ the holomorphy radius of $p\mapsto M_p$. Then
Lemma~\ref{lem:TWcalc}(vi) (composition with $p\mapsto M_p$) and~(ii) bound the analytic norms
of $M_{p_g}$ and of its $y$-derivatives at the radius $\rad_0$ by a data constant, and those of
order $\ge1$ by $C(\rad_0/\rad_1)$.
\end{itemize}
Hence $\hat{\mathbb A}_\xi(T_0)\le C_D''B$ and $\hat{\mathbb A}_{\ge1,\xi}(T_0)\le
C_D''\big((\rad_0/\rad_1)^2+\delta_y^2\big)B$ for every $\xi$, with $C_D''$ a data constant.
Since $\rad_1\ge c(\Dc)\Lambda$, integrating against $\rad_0^{-3}\chi_c^2$, we obtain~(a1), with
$C_D$ the larger of $C_D''/c(\Dc)^2$ and the constant $C(C_{\rm dat},\beta,L,\theta_\pm)$ of
hypothesis~(b) of Lemma~\ref{lem:RGbar}. The resulting constant does not depend on $\rad_0$ or $\delta_y$.

\emph{(a2).} The leaf parameters at $T_0$ are the moments of $\eps\rho_0(\xi)\phi(T_0,\cdot,\xi)$,
and those of free transport are the moments of $\eps\rho_0(\xi)\phi^{\rm fr}(T_0,\cdot,\xi)$. Their
difference is bounded by the moments of $w$, hence by $\sup_\xi\|m_{T_0}^{1/2}w(T_0,\cdot,\xi)\|$,
which is bounded, as in~(a1), by $C\sup_c\mathsf E^{(0)}_c[w]^{1/2}\le C\eps
T_0^{\kappa_\gamma}$.

\emph{(a3).} At $T_0$, the fiber is $g(T_0,y,v)=G(T_0,y,v+y/T_0)$, so $\partial_y$ and $\partial_v$
of the fiber are $D_y$ and $\partial_\xi$ of $G$ at $(z,\xi)=(y,v+y/T_0)$, and $g-g_{\rm
fr}=\eps\rho_0(v+y/T_0)\,w(T_0,y,v+y/T_0)$, while $g_{\rm fr}(y,v)=\eps\rho_0(v)F_0(y,v)$. By
the proposition, the cube Sobolev inequality in $\xi$ and Sobolev's inequality in $z$ as
in (F2), we have $|D^\alpha_y\partial^b_\xi w|\le C(\eps T_0^{\kappa_\gamma})e^{-(q_d/8)|z|^2}$ for
$|\alpha|\le N_y-4$ and $|b|\le1$. For $N_y-4<|\alpha|\le18\le N_y-2$ (for $N_y\ge22$ no
interpolation is needed), we interpolate in $z$, applied to $e^{a|z|^2}w$ with $a<q_d/4$ (the
commutators give polynomial factors, which are controlled by the margin), between this bound at the
order $N_y-2$ in $L^2(m^{1/2}_{T_0})$ and the bound of order $N_y$ on $\phi$ and $\phi^{\rm fr}$
from the proposition and Lemma~\ref{lem:dispfree}. Since $H^{N_y-2}$ is small and $H^{N_y}$ is
bounded, we get $H^{18+3/2+\delta}$ small to the power $\varkappa=(\frac12-\delta)/2$, and then
Sobolev's inequality yields the pointwise bound at order $18$, with $(\eps
T_0^{\kappa_\gamma})^\varkappa$. With $\rho_0(v+y/T_0)\le e^{L|y|/T_0}\rho_0(v)$, which is absorbed
by the Gaussian factor, we obtain
\begin{equation}
|\partial^\alpha_y\partial^b_v(g-g_{\rm fr})(T_0,y,v)|\le C(\eps
T_0^{\kappa_\gamma})^\varkappa\,\eps\rho_0(v)\,e^{-(q_d/16)|y|^2},\qquad|\alpha|\le18,\ |b|\le1 .
\end{equation}
Taking the moments of this bound, we obtain the closeness of $P_v$ and of $\nabla_vP_v$ to those of
$g_{\rm fr}$, relatively for the mass. The parameters of $g_{\rm fr}(\cdot,v)$ are
$(\eps\rho_0(v),y_{0,\infty}(v),\theta_\infty(v))$. Therefore (H1) follows from the convention on
$L$ of the paragraph \emph{Parameters} ($|y_{0,\infty}|\le L/2$ and
$|\nabla\log\rho_0|+|\nabla\log\theta_\infty|\le L/2$) and $2\theta_-\le\theta_\infty\le\theta_+/2$,
once $C(\eps T_0^{\kappa_\gamma})^\varkappa\le\min\{1,\theta_-\}/5$. The membership
$F_{v,0}\in\mathfrak K^7(C_{\rm dat},\beta,A_7)$ and the bound on $\partial_vg$ follow from the half
constants of $\hat F_0^v$, the closeness of $P_v$, and $2\beta\le\theta_-q_d/16$, which makes the
Gaussian factor $e^{-(q_d/16)|y|^2}$ dominate $e^{-\beta|\hat z|^2}$ in the normalized variables up
to a constant (using $\theta_v\ge2\theta_--C(\eps T_0^{\kappa_\gamma})^\varkappa$), absorbed for
$\eps\le\eps_0(T_0)$.

\emph{(a4).} With the moment vector $\mom(z)=(1,z,|z|^2)$, we put
$\mu:=m_g(T_0,\cdot)/(\eps\rho_0)=\mu^{\rm fr}+\mu^w$, with $\mu^{\rm fr}:=\int\phi^{\rm
fr}(T_0)\mom\dd z$ and $\mu^w:=\int w(T_0)\mom\dd z$.

\emph{Free part.} By (D3), $\rho_0$ has no zero on $\mathbb S_\Lambda$, and
$A_F^{-1}\rho_0(\mathrm{Re}\,\xi)\le|\rho_0(\xi)|\le A_F\rho_0(\mathrm{Re}\,\xi)$ there. With
$|\nabla\log\rho_0|\le L$ on $\R^3$ and the Borel--Carath\'eodory inequality on discs of radius
$\Lambda/2$, we have $|\nabla\log\rho_0|\le C_{\rho_0}=C(A_F,L,\Lambda)$ on $\mathbb S_{\Lambda/2}$.
Hence $|\rho_0(\xi-z/T_0)/\rho_0(\xi)-1|\le C_{\rho_0}|z|T_0^{-1}e^{C_{\rho_0}|z|/T_0}$ there. Let
\begin{equation}
\mu_0(\xi):=\int F_0(z,\xi)\mom(z)\dd z=\big(1,\ y_{0,\infty}(\xi),\
3\theta_\infty(\xi)+y_{0,\infty}(\xi)\cdot y_{0,\infty}(\xi)\big),
\end{equation}
which is holomorphic on $\mathbb S_\Lambda$ by (D4) (for product data it is constant). For real $z$
and $\xi\in\mathbb S_{\Lambda/2}$,
\begin{equation}
\phi^{\rm fr}(T_0,z,\xi)-F_0(z,\xi)=\Big(\frac{\rho_0(\xi-z/T_0)}{\rho_0(\xi)}-1\Big)F_0(z,\xi-z/T_0)+\big(F_0(z,\xi-z/T_0)-F_0(z,\xi)\big).
\end{equation}
By (D4) and (D1), $|F_0(z,\zeta)|\le A_FC_Fe^{-\beta_F|z|^2}$ for $\zeta\in\mathbb S_\Lambda$, so,
by the Cauchy estimates, $|\nabla_vF_0(z,\zeta)|\le C\Lambda^{-1}A_FC_Fe^{-\beta_F|z|^2}$ on
$\mathbb S_{\Lambda/2}$, and the second term is bounded by $C(\Dc)|z|T_0^{-1}e^{-\beta_F|z|^2}$,
while the first one is bounded by $C(\Dc)|z|T_0^{-1}e^{C_{\rho_0}|z|-\beta_F|z|^2}$. Thus $\mu^{\rm
fr}$ is holomorphic on $\mathbb S_{\Lambda/2}$, with $|\mu^{\rm fr}-\mu_0|\le C(\Dc)T_0^{-1}$ there.
Applying the Cauchy estimates to the holomorphic function $\mu^{\rm fr}-\mu_0$, we obtain $|\mu^{\rm
fr}-\mu_0|+|\nabla(\mu^{\rm fr}-\mu_0)|\le C(\Dc)T_0^{-1}$ on $\mathbb S_{\Lambda/4}$.

\emph{Deviation.} Applying the same Sobolev inequality on cubes of side $\rad_1$, at the radius
$\rad_1/2$, to $m^{1/2}_{T_0}w$, together with the proposition, we get
$\|m^{1/2}_{T_0}\partial^\beta_\xi w\|\le CC_w^{1/2}\eps
T_0^{\kappa_\gamma}|\beta|!(\rad_1/2)^{-|\beta|}$ for real $\xi$ (the weights $(n+1)^5/n!$ and
the Sobolev loss $(n+2)^2$ are absorbed). Hence $|\partial^\beta\mu^w|\le C_{\mu,1}\eps
T_0^{\kappa_\gamma}|\beta|!(\rad_1/2)^{-|\beta|}$, by Cauchy--Schwarz against $m_{T_0}^{1/2}$.
Since $\sum_{|\beta|=n}\frac{|\beta|!}{\beta!}|\zeta^\beta|=|\zeta|_1^n$, the Taylor series of
$\mu^w$ and of $\nabla\mu^w$ at a real point $\xi_0$ converge absolutely on
$\xi_0+\{\zeta\in\C^3:|\zeta|_1<\rad_1/2\}$, with majorants independent of $\xi_0$. By Taylor's
formula with the same bounds, the remainder of order $N$ at real $\zeta$ is $\le C_{\mu,1}\eps
T_0^{\kappa_\gamma}(2|\zeta|_1/\rad_1)^{N+1}$, so $\mu^w$ is real-analytic on $\R^3$ and equals
the sum of its Taylor series there. Two such continuations agree on their overlap, which is convex
and meets $\R^3$. The union of these sets is $\{|\mathrm{Im}\,\xi|_1<\rad_1/2\}$, which contains
$\mathbb S_{\rad_1/8}$, since $|\eta|_1\le\sqrt3|\eta|$. Therefore $\mu^w$ extends
holomorphically to $\mathbb S_{\rad_1/8}$, with $|\mu^w|+|\nabla\mu^w|\le
C_{\mu,2}\rad_1^{-1}\eps T_0^{\kappa_\gamma}$.

\emph{The data parameters on the strip.} On $\mathbb S_{R_1}$, $R_1=\rad_1/8\le\Lambda/64$, we write
$Y:=y_{0,\infty}(\xi)$ and $\Theta:=\theta_\infty(\xi)$, which are complex there. By the choice of
$\rad_1$ (paragraph \emph{Parameters}),
$|Y-y_{0,\infty}(\mathrm{Re}\,\xi)|+|\Theta-\theta_\infty(\mathrm{Re}\,\xi)|\le
C_{\mu_0}R_1\le\iota:=\frac18\min\{1,\theta_-\}$. Thus, by the convention on $L$ and
$2\theta_-\le\theta_\infty\le\frac12\theta_+$ on $\R^3$,
\begin{equation}
|Y|\le Y_*:=\tfrac L2+\tfrac18,\qquad\mathrm{Re}\,\Theta\ge2\theta_--\iota,\qquad|\Theta|\le\tfrac12\theta_++\iota\qquad\text{on
}\mathbb S_{R_1}.
\end{equation}
Put $S:=2+L/2+\theta_+$. Then $1+2Y_*+|Y\cdot Y-3\Theta|\le(1+Y_*)^2+\frac32\theta_++3\iota\le S^2$,
and $(1+Y_*)(1+Y_*+G_\infty)\le S^2(1+3G_\infty)$.

\emph{Conclusion.} With $c_{\rm a4}$ of the statement, we require $C(\Dc)T_0^{-1}\le c_{\rm a4}$
and then $C_{\mu,2}\rad_1^{-1}\eps T_0^{\kappa_\gamma}\le c_{\rm a4}$, with
$C_{\mu,2}=C(C_w,q_d)$. Then on $\mathbb S_{R_1}$ the function $d:=\mu-\mu_0$ is holomorphic and
$|d|+|\nabla d|\le\delta:=2c_{\rm a4}$ (here $|\mu|:=|\mu_{(0)}|+|\mu_{(1)}|+|\mu_{(2)}|$, with the
Euclidean norm on $\C^3$, and likewise for $\nabla\mu$), so that $\delta
S^2(1+3G_\infty)\le\frac14\min\{1,\theta_-\}$ and $\delta\le\frac1{64}$ (indeed
$\delta\le\frac14\theta_-S^{-2}\le\frac1{16}\theta_+(2+\theta_+)^{-2}\le\frac1{128}$, since
$\theta_-\le\frac14\theta_+$). We claim that $\mathrm{Re}\,\rho_g/(\eps\rho_0)\ge\frac12$,
$|y_{0,g}-Y|\le\frac12$, $\mathrm{Re}\,\theta_g\ge\theta_-$ and $|\theta_g|\le2\theta_+$ there (with
$\theta_g=\frac13(m_2/m_0-y_{0,g}\cdot y_{0,g})$, the bilinear square). Indeed, we have
$\mu_{(0)}=1+d_{(0)}$, and with $e:=y_{0,g}-Y=(d_{(1)}-Yd_{(0)})/\mu_{(0)}$ the bilinear square
yields
\begin{equation}
3(\theta_g-\Theta)=\frac{N}{\mu_{(0)}}-e\cdot e,\qquad N:=d_{(2)}-2Y\cdot d_{(1)}+(Y\cdot
Y-3\Theta)d_{(0)} .
\end{equation}
Since $|N|\le\delta S^2$, $|\mu_{(0)}|\ge1-\delta$ and $|e|\le\delta(1+Y_*)/(1-\delta)\le\frac12$,
we get $3|\theta_g-\Theta|\le\delta S^2(1-\delta)^{-1}+\delta^2S^2(1-\delta)^{-2}=\delta
S^2(1-\delta)^{-2}$, that is,
$|\theta_g-\Theta|\le\frac1{12}(\frac{64}{63})^2\min\{1,\theta_-\}\le\frac1{10}\min\{1,\theta_-\}$
on $\mathbb S_{R_1}$. Hence $\mathrm{Re}\,\theta_g\ge2\theta_--\iota-\frac1{10}\theta_-\ge\theta_-$
and $|\theta_g|\le\frac12\theta_++\iota+\frac1{10}\theta_-\le2\theta_+$ there, and on $\R^3$, where
$\theta_g$ and $\Theta=\theta_\infty$ are real,
$\theta_-\le\theta_g\le\frac12\theta_++\frac1{10}\theta_-\le\frac23\theta_+$, since
$\theta_-\le\frac14\theta_+$. Hence $\nabla\ell_g=\nabla\log\rho_0+\nabla\log(\rho_g/(\eps\rho_0))$,
$y_{0,g}$ and $\vartheta_g=\log\theta_g$ (principal branch) are holomorphic on $\mathbb S_{R_1}$,
and
\begin{equation}
|\varpi'_g(T_0)|=|\nabla\ell_g|+|y_{0,g}|+|\vartheta_g|\le
L_1:=C_{\rho_0}+L/2+\max\{|\log\theta_-|,|\log(2\theta_+)|\}+4 .
\end{equation}
Finally, on $\R^3$, where $d$, $Y=y_{0,\infty}$ and $\Theta=\theta_\infty$ are real, $\nabla\log
A_{\rm ref}=-\nabla\ell_g+\frac34\nabla\vartheta_g$ at $T_0$, since $A_p=\rho^{-1}\theta^{3/4}$. For
the mass, $\nabla\ell_g=\nabla\log\rho_0+\nabla d_{(0)}/\mu_{(0)}$, so
$|\nabla\ell_g|\le|\nabla\log\rho_0|+\delta/(1-\delta)\le|\nabla\log\rho_0|+\frac1{63}$. For the
temperature, the derivatives of $N$ and of $e$ now contain also the derivatives of $y_{0,\infty}$
and $\theta_\infty$, which are bounded by $G_\infty$ on $\R^3$:
\begin{equation}
|\nabla N|\le\delta S^2+\delta\,(5+2Y_*)G_\infty\le\delta S^2(1+3G_\infty),\qquad|\nabla
e|\le\frac{\delta(1+Y_*+G_\infty)}{(1-\delta)^2} .
\end{equation}
Thus, with $3(\theta_g-\theta_\infty)=N/\mu_{(0)}-e\cdot e$, we have
\begin{multline}
3|\nabla(\theta_g-\theta_\infty)|\le\Big|\nabla\frac N{\mu_{(0)}}\Big|+2|e|\,|\nabla e|\le\delta
S^2(1+3G_\infty)\Big(\frac1{(1-\delta)^2}+\frac{2\delta}{(1-\delta)^3}\Big)\\
\le\tfrac{11}{10}\,\delta S^2(1+3G_\infty)\le\tfrac3{10}\min\{1,\theta_-\},
\end{multline}
hence $|\nabla\theta_g|\le|\nabla\theta_\infty|+\frac1{10}\theta_-$. Since
$\theta_g\ge\theta_\infty-\frac1{10}\theta_-\ge\frac9{10}\theta_\infty$ and $\theta_g\ge\theta_-$,
$|\nabla\vartheta_g|=|\nabla\theta_g|/\theta_g\le\frac{10}9|\nabla\log\theta_\infty|+\frac1{10}$.
Hence
\begin{equation}
|\nabla\log A_{\rm ref}|\le|\nabla\log\rho_0|+\tfrac1{63}+\tfrac34\big(\tfrac{10}9|\nabla\log\theta_\infty|+\tfrac1{10}\big)\le\tfrac
L2+\tfrac1{10}\le L
\end{equation}
on $\R^3$, by the convention $|\nabla\log\rho_0|+|\nabla\log\theta_\infty|\le L/2$ and $L\ge1$. For
product data $\nabla\theta_\infty=0$ and $G_\infty=0$, and this is the computation with a constant
$\mu_0$.
\end{proof}

\subsection{The assembly}\label{ssec:assembly}

The initial layer, the dispersed layer and the bootstrap combine into the proofs of
Theorems~\ref{thm:main}
and~\ref{thm:t1}.

\begin{proof}[Proof of Theorem~\ref{thm:main}]
\emph{Step 1: initial layer.} By Theorem~\ref{thm:layer}, or Theorem~\ref{thm:EPlayer} at
$\gamma=-2$, the solution exists on $[0,t_1]$, and at $t_1$ it satisfies the hypotheses of
Theorem~\ref{thm:t1}.

\emph{Step 2: dispersed layer.} By Proposition~\ref{prop:disp}, the solution exists on $[t_1,T_0]$
for $\eps\le\eps_0(T_0)$, and by Corollary~\ref{cor:LWPa} the data at $T_0$ satisfy (LWP)(a),
clauses (a1)--(a3). By~(a4) and $T_0\ge T_{\rm E}(L_1,R_1)$, we obtain the reference Euler flow from
Lemma~\ref{lem:Mref}.

\emph{Step 3: the bootstrap.} By Corollary~\ref{cor:CUmain}, that is, Theorem~\ref{thm:CUboot} with
the continuation argument of Subsection~\ref{ssec:cellsproofs}, (U1)--(U4) hold for all $t\ge T_0$.

\emph{Step 4: proof of (i).} It follows from (U1)--(U3), which bound the analytic norms at radius
$\ge\rad_0/4$, in the reading of that corollary, and from $g\ge0$: on $[0,t_1]$,
$g=f_0G$ with $f_0>0$ (by (D$_{\rm iso}$) and (D3)) and $|G-1|\le C\eps t_1\le\frac14$ by
Lemma~\ref{lem:FPiter}(i), (iii), with $\eps$ small as in its proof. At $\gamma=-2$, we have
$|G-1|\le C\eps t_1^{3/4}\le\frac14$ by Theorem~\ref{thm:EPlayer}, for $\eps$ small. Consequently
$g>0$ there. On $[t_1,T_0]$, $g$ is the limit $g\ge0$ of the proof of
Proposition~\ref{prop:disp}. On $[T_0,\infty)$, $g$ is made of the solutions $g\ge0$ of
Proposition~\ref{prop:CULW}(a), by the same corollary. We do not prove strict positivity
after $t_1$, and the arguments after $t_1$ do not need it: they use only $g\ge0$ (for instance
$\ab[g]\ge0$, and in Lemma~\ref{lem:RGbar}) and the lower bound $\rho_g\ge\frac12\rho_\xi>0$ of the
leaf masses from (F3).

\emph{Step 5: proof of (ii).} It follows from Corollary~\ref{cor:Aleafentry} on the leaves with
$s_\xi\in[S_*,3S_\sharp]$ and from (U3) and Lemma~\ref{lem:CUsob} on the others, as in Step~1 of
Theorem~\ref{thm:CUboot}, with $C\delta_b^2$ in place of $\delta_1^2$, the equivalence of $X(\hat
m_\xi)$ with $L^2(\hat m_\xi)$, and the change of frame back to physical variables. Note that
$s_\xi(t)\to\infty$ because $\tau(t)\to\infty$ for $\gamma\le-2$.

\emph{Step 6: proof of (iii).} The first bound follows from (U1), Lemma~\ref{lem:TWcalc}(i) and
the requirement $C'_{\rm Sob}\delta_U\le\rad_0\min\{\delta_s,1\}/4$ of the same theorem,
since $|\varpi'_g-\varpi'_*|\le|\delta\varpi|+|\nabla\delta\varpi|\le C'_{\rm
Sob}\delta_U/\rad_0$ (the pointwise bound from $S_Q$ loses the inverse weights at the levels $0$
and $1$, which is the factor $8$). For the second bound, the leaf Euler system has the prefactor
$t^{-2}$
and a right-hand side bounded by (U1)--(U2), Lemma~\ref{lem:CUsob} and Lemma~\ref{lem:Mref} (the
$\nabla\ell$-equation is the $\xi$-derivative of the first equation of Theorem~\ref{thm:Meuler}, and
$s,h$ and their first derivatives are bounded pointwise by (U2) and the same lemma), so
$\varpi'_g(t,\xi)$ converges at the rate $C/t$.

\emph{Step 7: proof of (iv).} The first two statements are Lemma~\ref{lem:Wtail}(i). For the last
one, we fix $d_0>0$, let $a':=q_d/2$, and use the notation of part~(ii) of that lemma.
\begin{itemize}
\item \emph{At $T_0$.} We have $g(T_0,z,\xi-z/T_0)=\eps\rho_0(\xi)\phi(T_0,z,\xi)$ with
$\phi=\phi^{\rm fr}+w$, by Proposition~\ref{prop:disp}. By the proof of Lemma~\ref{lem:LWfree},
$|\phi^{\rm fr}(T_0,z,\xi)-F_0(z,\xi)|\le(e^{L|z|/T_0}-1)F_0(z,\xi)$, which implies that
$\|\phi^{\rm fr}(T_0,\cdot,\xi)-F_0(\cdot,\xi)\|_{L^2(e^{a'|z|^2})}\le CL/T_0$.

By the proof of Corollary~\ref{cor:LWPa}(a2),
$\|w(T_0,\cdot,\xi)\|_{L^2(e^{a'|z|^2})}\le\|m_{T_0}^{1/2}w(T_0,\cdot,\xi)\|\le C\eps
T_0^{\kappa_\gamma}$. By Lemma~\ref{lem:LWfree}(iii) and~(a2) of the corollary, the leaf mass, mean
and temperature are within $C(L/T_0+\eps T_0^{\kappa_\gamma})$ of
$(\eps\rho_0(\xi),y_{0,\infty}(\xi),\theta_\infty(\xi))$, relatively for the mass, and
$(y_0,\theta)\mapsto M_{(1,y_0,\theta)}$ is Lipschitz into $L^2$ for
$\theta\in[\frac12\theta_-,\theta_++1]$.

Since $\mathcal N_\xi(T_0)=\big(\eps\rho_0(\xi)/\rho_g(T_0,\xi)\big)\phi(T_0,\cdot,\xi)$, and
$F_0(\cdot,\xi)$ is bounded in $L^2(e^{a'|z|^2})$, uniformly in $\xi$, by (D1) and $a'\le\beta_F/8$,
we deduce from these bounds that $\|\mathcal N_\xi(T_0)-F_0(\cdot,\xi)\|_{L^2(e^{a'|z|^2})}\le
C_{T_0}(L/T_0+\eps
T_0^{\kappa_\gamma})$ for every $\xi$, with $C_{T_0}=C(\Dc)$. In spite of the notation, $C_{T_0}$
does not depend on $T_0$ or on $\eps$, since $T_0\ge t_1\ge1$ and $C_w$ is chosen before $T_0$.

By the triangle inequality, $|d_\xi(T_0)-d_0(\xi)|\le C_{T_0}(L/T_0+\eps T_0^{\kappa_\gamma})$ for
every $\xi$ as well. In particular $d_\xi(T_0)\ge d_0-C_{T_0}(L/T_0+\eps T_0^{\kappa_\gamma})$.
\item \emph{Later times.} Let $t\ge T_0$. By Lemma~\ref{lem:Wtail}(i) there are $\xi$ with $C_{\rm
tail} s_\xi(t)\le d_0/8$. For such points, by part~(ii) of the lemma and the first item, we have
$d_\xi(t)\ge d_0-C_{T_0}(L/T_0+\eps T_0^{\kappa_\gamma})-C_{\rm tail} T_0^{-1}-d_0/8$. We put
$C_{\rm
nM}:=8(C_{T_0}(L+1)+C_{\rm tail})$. Since $T_0^{\kappa_\gamma}\le1+\tau(T_0)$, the hypotheses
$T_0\ge C_{\rm nM}/d_0$ and $C_{\rm nM}\eps(1+\tau(T_0))\le d_0$ make each of the three error terms
$\le d_0/8$, so $d_\xi(t)\ge d_0/2$.
\item \emph{Conclusion.} The weight of the norm in~(ii) is $\ge1$, and at $x=t\xi$ the left side
of~(ii) is the weighted norm of $\mathcal N_\xi(t)-M_{\hat p_g(t,\xi)}$. Therefore it is $\ge
d_\xi(t)\ge d_0/2$.
\end{itemize}
The requirements on $T_0$ in Appendix~\ref{app:const} are lower bounds, so $T_0$ may be enlarged to
satisfy $T_0\ge C_{\rm nM}/d_0$, and then $\eps_0$ is taken for the new $T_0$.
\end{proof}

\begin{proof}[Proof of Theorem~\ref{thm:t1}]
The proof of Theorem~\ref{thm:main} applies without its first step. Proposition~\ref{prop:disp} and
Corollary~\ref{cor:LWPa} are proved for a reference profile $f_0\in\Dc^\sharp$ and for the data of
Theorem~\ref{thm:t1}, with $C_{t_1}^2$ in place of $C_{\rm lay}t_1^2$ in the bullet \emph{Data} and
in Steps~6--7 of the proof of the proposition. Moreover, Lemmas~\ref{lem:dispfree}
and~\ref{lem:LWfree}
require no more than (D1), (D3) and (D4). After $T_0$, the argument requires of the data only the
conclusions of Corollary~\ref{cor:LWPa}. The conclusions are drawn as above, on $[t_1,\infty)$:
$g\ge0$ on $[t_1,T_0]$ by the proof of the proposition, and the first item of the proof
of~(iv) rests on the proposition, Lemma~\ref{lem:LWfree} and~(a2) of the corollary. Here $C_w$,
hence $C_{T_0}$ and $C_{\rm nM}$, depend also on
$C_{t_1}$.
\end{proof}

\subsection{Complete relaxation and scattering}\label{ssec:scatter}

We prove Theorem~\ref{thm:scatter} (complete relaxation and scattering) and the rates,
Theorems~\ref{thm:rate2} and~\ref{thm:rate3}.

\begin{proof}[Proof of Theorem~\ref{thm:scatter}]
By Corollary~\ref{cor:CUmain}, (U1)--(U4) hold for all $t\ge T_0$, so the lemmas proved under these
hypotheses apply on $[T_0,\infty)$. In this proof $C$ and $c_0>0$ denote constants that do not
depend on $t$, $\xi$ or $z$. Let $H(t,\xi,z):=g(t,z,\xi-z/t)$ be the leaf through $x=t\xi$, so that
$\rho_g(t,\xi)=\int H\dd z$ and $\mathcal N_\xi(t)=H(t,\xi,\cdot)/\rho_g(t,\xi)$, and let
$r_\xi(t):=\mathcal N_\xi(t)-M_{\hat p_g(t,\xi)}$, so that $d_\xi(t)=\|r_\xi(t)\|_{L^2}$.

\emph{Step 1: bounds on the leaves.} We claim that for $t\ge T_0$ and all $\xi,z$,
\begin{equation}\label{eq:SCbd}
|H|+|\nabla_zH|+|\nabla^2_zH|+|\nabla_\xi H|\le C\rho_g(t,\xi)\,e^{-c_0|z|^2} .
\end{equation}
We write $g=\mathcal M[g]+\hat r$ as in Definition~\ref{def:Mloc}. The leaf of $\mathcal M[g]$ is
$M_{p_g(t,\xi)}(z)$. Its $z$-derivatives are explicit, and its $\xi$-derivative is
$M_{p_g(t,\xi)}(z)$ times a polynomial of degree $2$ in $z$ whose coefficients are bounded by
$|\nabla_\xi\log\rho_g|+|\nabla_\xi y_{0,g}|+|\nabla_\xi\theta_g|$. These are bounded, by
Theorem~\ref{thm:main}(iii) and Lemma~\ref{lem:Mref} for $\nabla\log\rho_g$, which is the first
component of $\varpi'_g$, and by (F3) for the others. Since
$\frac12\theta_-\le\theta_g\le\theta_++1$ and $|y_{0,g}|\le L+1$ by (F3), this part
satisfies~\eqref{eq:SCbd}.

For $\hat r$, by (F1) the derivatives $\nabla_\xi$ and $\nabla^k_z$, $k\le2$, of its leaf are leaf
restrictions of $\partial_v\hat r$ and of $t^{-i}\partial^a_y\partial^b_v\hat r$, $|a|+|b|\le2$. By
Lemma~\ref{lem:BTsob}(ii), with $B_\xi\le\max\{B_{\rm cr},\delta_b^2\}$, and since
$T_pF(z)=\rho\,\theta^{-3/2}F(\theta^{-1/2}(z-y_0))$, each of them is bounded by
$C\rho_g(t,\xi)e^{-(q/4)|\hat z|^2}$, $\hat z=\theta_g^{-1/2}(z-y_{0,g})$. Here the exponent $q$ of
the leaf weight $\hat m_\xi$ in the class $\mathcal W$ lies in
$[\frac12q(t)\theta_-,\,q(t)(\theta_++1)]$ by Lemma~\ref{lem:TG0}(i). It is bounded below, since
$q(t)\ge\frac34q_0$ by the requirement $T_0\ge4C_2/q_0$. The constants of
the Sobolev bound are bounded for $q(t)\in[\frac34q_0,q_0]$. With the bounds on $\theta_g$
and $y_{0,g}$, we obtain the claim. In particular, $\mathcal N_\xi\le Ce^{-c_0|z|^2}$, and since
$M_{\hat p_g}\le Ce^{-c_0|z|^2}$ as well, we have for every $R\ge0$
\begin{equation}\label{eq:SCtail}
\int_{|z|>R}\big(|r_\xi(t)|+|r_\xi(t)|^2\big)\dd z\le Ce^{-c_0R^2/2} .
\end{equation}

\emph{Step 2: the derivative along rays.} By Lemma~\ref{lem:leafgeom}(i), (iv) and
Proposition~\ref{prop:leaf}, we have, at fixed $(\xi,z)$,
\begin{equation}
D_tH=\kappa\,Q(H,H)+t^{-2}z\cdot\nabla_\xi H .
\end{equation}
In nondivergence form $Q(H,H)=\ab[H]:\nabla^2H-\cb[H]H$, and $|\ab[H]|+|\cb[H]|\le C\|H\|_{L^1\cap
L^\infty}\le C\rho_g$, as in the proof of Lemma~\ref{lem:Wtail}. By~\eqref{eq:SCbd}, we have
$|D_tH|\le
C(\kappa\rho_g^2+t^{-2}\rho_g)e^{-c_0|z|^2/2}$. Since $\int Q(H,H)\dd z=0$, we get $|D_t\rho_g|\le
Ct^{-2}\rho_g$, hence $|D_t\mathcal N_\xi|\le C(\kappa\rho_g+t^{-2})e^{-c_0|z|^2/2}$. The mean and
the temperature of $\mathcal N_\xi$ are moments of $\mathcal N_\xi$, so their derivatives along the
ray obey the same bound, and $(y_0,\theta)\mapsto M_{(1,y_0,\theta)}$ is Lipschitz into $L^2$ for
$\theta\ge\frac12\theta_-$. With $\rho_g\le2\rho_\xi$, from (F3), where $\rho_\xi$ is the mass of
the fiber $g(T_0,\cdot,\xi)$, we conclude that $t\mapsto r_\xi(t)\in L^2(\R^3)$ is continuously
differentiable along the ray and
\begin{equation}\label{eq:SClip}
\|D_tr_\xi(t)\|_{L^2}\le C\big(\kappa(t)\rho_\xi+t^{-2}\big) .
\end{equation}

\emph{Step 3: the integrated dissipation.} Let $f_\xi=T^{-1}_{p_g}H-\mu$ be the micro part of the
leaf in its frame, so that $f_\xi(\hat z)=\theta_g^{3/2}r_\xi(\theta_g^{1/2}\hat z+y_{0,g})$. Every
leaf has $s_\xi(t)\to\infty$, so the time $t^-_\xi$ of Lemma~\ref{lem:UFlevel0} is finite, and by
that lemma we have, for every $\xi$,
\begin{equation}
\int_{t^-_\xi}^\infty\kappa\,\lambda_\xi\,\|f_\xi\|^2_{Y_2(\hat m_\xi)}\dd t\le C_{L0}\,\delta_1^2,
\end{equation}
since the leaf dissipation of Subsection~\ref{ssec:cells} is
$d_\xi[f_\xi]=\lambda_\xi\|f_\xi\|^2_{Y_2(\hat m_\xi)}$ (it is not the distance $d_\xi(t)$), with
$\lambda_\xi=\rho_g\theta_g^{\gamma/2}$ the collision rate of the leaf. The zero-order part of the
norm $Y_2(\hat m_\xi)$ is $\|\jap{\hat z}^{(\gamma+2)/2}\hat m_\xi f\|_{L^2}$, and $\hat m_\xi\ge
c>0$ by the membership lemma. We put
\begin{equation}
\Psi_\xi(t):=\big\|\jap z^{(\gamma+2)/2}r_\xi(t)\big\|_{L^2(\R^3)} .
\end{equation}
Since $\gamma+2\le0$ and $\jap{\hat z}\le C\jap z$, by the change of variables we get $\|\jap{\hat
z}^{(\gamma+2)/2}f_\xi\|_{L^2}\ge c\,\Psi_\xi$, and $\lambda_\xi\ge c\rho_\xi$ by (F3). Hence
\begin{equation}\label{eq:SCdiss}
\int_{t^-_\xi}^\infty\kappa(t)\,\rho_\xi\,\Psi_\xi(t)^2\dd t\le C\delta_1^2 .
\end{equation}

\emph{Step 4: $\Psi_\xi\to0$.} We fix $\xi$ and use the time $s(t):=\rho_\xi(\tau(t)-\tau(T_0))$,
which increases to $\infty$, with $s'=\kappa\rho_\xi$. Since $\gamma\le-2$,
$t^{-2}=t^{1+\gamma}\kappa(t)\le t^{-1}\kappa(t)$, so $t^{-2}\le\kappa\rho_\xi$ for $t\ge
t_*:=\max\{t^-_\xi,\rho_\xi^{-1}\}$. As $\jap z^{(\gamma+2)/2}\le1$, it follows
from~\eqref{eq:SClip} that, for $t_*\le t'\le t$,
\begin{equation}
|\Psi_\xi(t)-\Psi_\xi(t')|\le\int_{t'}^t\|D_tr_\xi\|_{L^2}\le2C\big(s(t)-s(t')\big) .
\end{equation}
Thus, as a function of $s$, $\Psi_\xi$ is Lipschitz on $[s(t_*),\infty)$, and by~\eqref{eq:SCdiss}
its square is integrable there. A nonnegative Lipschitz function with integrable square on a
half-line tends to zero: if $\Psi_\xi\ge\eta$ at points $s_k\to\infty$, then $\Psi_\xi\ge\eta/2$ on
the intervals $[s_k,s_k+\eta/(4C)]$, infinitely many of which are disjoint. Hence $\Psi_\xi(t)\to0$.

\emph{Step 5: proof of~(i).} For $R\ge0$, by~\eqref{eq:SCtail},
\begin{equation}
d_\xi(t)^2=\|r_\xi(t)\|_{L^2}^2\le\jap R^{|\gamma+2|}\Psi_\xi(t)^2+Ce^{-c_0R^2/2} .
\end{equation}
We let $t\to\infty$ and then $R\to\infty$. The same holds in the weighted norm of
Theorem~\ref{thm:main}(ii). Indeed, $r_\xi=\hat r^{(t\xi)}/\rho_g$, so by the definition of
$\ell_\xi$, by $|\nabla\log A_{\rm ref}|\le L$ and by
Lemma~\ref{lem:BTsob}(i) with $\alpha=\beta=0$,
\begin{equation}
\int e^{2q(t)|z|^2-2L|z|}|r_\xi|^2\dd z\le\theta_g^{-3/2}\ell_\xi[\hat r]\le C .
\end{equation}
The weight of that norm is $e^{q_0\theta_-|z-y_{0,g}|^2/(8\theta_g)}\le
Ce^{(q_0/2)|z|^2}$, by (F3), and $2q(t)\ge\frac32q_0$. Splitting the integral at $|z|=R$, we find
that the square of the weighted norm of $r_\xi(t)$ is bounded by
$Ce^{q_0R^2}d_\xi(t)^2+Ce^{-q_0R^2/4}$, and we
conclude as before.

\emph{Step 6: proof of~(ii).} By Theorem~\ref{thm:main}(iii), $y_{0,g}(t,\xi)$ and
$\vartheta_g(t,\xi)=\log\theta_g(t,\xi)$ converge at the rate $C/t$, uniformly in $\xi$. By Step~2,
$|D_t\log\rho_g|\le Ct^{-2}$, so $\log\rho_g(t,\xi)$ converges at the same rate, uniformly in $\xi$.
The limits define $p_{\rm sc}$, with $\theta_{\rm sc}:=e^{\vartheta_{\rm sc}}$, where
$\vartheta_{\rm sc}$ is the limit of $\vartheta_g$. Since $\theta_g\le\theta_++1$, we also have
$|\theta_g-\theta_{\rm sc}|\le C/t$. The function $p_{\rm sc}$ is continuous, as a uniform limit of
continuous functions, and the bounds of (F3) pass to the limit. Moreover, by the same result,
$\nabla\log\rho_g(t,\cdot)$, the first component of $\varpi'_g$, converges uniformly. Thus
$\log\rho_{\rm sc}$ is of class $C^1$, and the bound on $\nabla\log\rho_g$ of Step~1 passes to the
limit. In particular, $\frac12\rho_\xi\le\rho_{\rm
sc}(\xi)\le2\rho_\xi$. The last statement of~(ii) follows from~(i), since $(y_0,\theta)\mapsto
M_{(1,y_0,\theta)}$ is continuous into $L^2$.

\emph{Step 7: proof of~(iii).} The map $(y,v)\mapsto(z,\xi)=(y,v+y/t)$ preserves the Lebesgue
measure, and $g(t,y,v)=H(t,v+y/t,y)$. Hence we have
\begin{equation}
\|g(t)-g_{\rm sc}\|_{L^1(\R^6)}=\int I(t,\xi)\dd\xi,\qquad I(t,\xi):=\int\big|H(t,\xi,z)-M_{p_{\rm
sc}(\xi-z/t)}(z)\big|\dd z .
\end{equation}
\emph{Domination.} $\int H\dd z=\rho_g\le2\rho_\xi$. Since $|\nabla\log\rho_{\rm sc}|\le C$, we have
$\rho_{\rm sc}(\xi-z/t)\le\rho_{\rm sc}(\xi)e^{C|z|}$ for $t\ge1$, so $\int M_{p_{\rm
sc}(\xi-z/t)}(z)\dd z\le C\rho_{\rm sc}(\xi)\le C\rho_\xi$, by the bounds on $y_{\rm sc}$ and
$\theta_{\rm sc}$. Thus $I(t,\xi)\le C\rho_\xi$, and $\int\rho_\xi\dd\xi<\infty$, since
$\int\rho_\xi\dd\xi\le2\int\rho_g(T_0,\xi)\dd\xi=2\|g(T_0)\|_{L^1}$. \emph{Pointwise convergence.}
For fixed $\xi$, we have
\begin{equation}
I(t,\xi)\le\rho_g(t,\xi)\|r_\xi(t)\|_{L^1}+\big\|M_{p_g(t,\xi)}-M_{p_{\rm
sc}(\xi)}\big\|_{L^1}+\int\big|M_{p_{\rm sc}(\xi)}(z)-M_{p_{\rm sc}(\xi-z/t)}(z)\big|\dd z .
\end{equation}
By~\eqref{eq:SCtail}, $\|r_\xi(t)\|_{L^1}\le CR^{3/2}d_\xi(t)+Ce^{-c_0R^2/2}$ for every $R$, so the
first term tends to $0$ by~(i). The second tends to $0$ by~(ii). In the third, the integrand tends
to $0$ for every $z$, by the continuity of $p_{\rm sc}$, and it is dominated by $C\rho_{\rm
sc}(\xi)e^{C|z|}e^{-c|z|^2}$, and hence it tends to $0$. By dominated convergence in $\xi$, we
obtain $\|g(t)-g_{\rm sc}\|_{L^1}\to0$. The equality of the two norms in~(iii) is the change of
variables $y=x-tv$. In
particular $g_{\rm sc}\in L^1$, and $\|g_{\rm sc}\|_{L^1}=\lim_{t\to\infty}\|g(t)\|_{L^1}$ is the
total mass of the solution.

\emph{Step 8: proof of~(iv).} For every $v$, $g_{\rm sc}(\cdot,v)=\rho_{\rm sc}(v)M$ with $M$ a
Maxwellian of unit mass, and $f_0(\cdot,v)=\eps\rho_0(v)F_0(\cdot,v)$. Thus we have $\|g_{\rm
sc}(\cdot,v)-f_0(\cdot,v)\|_{L^1(\R^3_y)}\ge\eps\rho_0(v)d_1(v)$, and integrating in $v$, we obtain
the inequality, since $\|f_0\|_{L^1}=\eps\int\rho_0$. It remains to show that $d_1(v)>0$ when
$F:=F_0(\cdot,v)$ is not a Maxwellian. If $d_1(v)=0$, there are $c_k\ge0$ and Maxwellians
$M_k=M_{(1,y_k,\theta_k)}$ with $c_kM_k\to F$ in $L^1$. Then $c_k\to\int F=1$, so $M_k\to F$ in
$L^1$. If $(y_k,\theta_k)$ has a limit point in $\R^3\times(0,\infty)$, then $F$ is a Maxwellian.
Otherwise, along a subsequence, $\theta_k\to0$, $\theta_k\to\infty$ or $|y_k|\to\infty$. In the last
two cases $\int_{B}M_k\to0$ for every ball $B$, which contradicts $\int_BM_k\to\int_BF$ for a ball
with $\int_BF>0$. If $\theta_k\to0$ with $y_k$ bounded, the mass of $M_k$ concentrates on sets of
measure tending to $0$, which contradicts the convergence in $L^1$ to an integrable function.
\end{proof}

\begin{proof}[Proof of Theorem~\ref{thm:rate2}]
Let $\gamma=-2$, so that $\kappa(t)=t^{-1}$. We use the notation of the proof of
Theorem~\ref{thm:scatter}, and write
\begin{equation}\label{eq:etleaf}
e(t):=e_\xi[f_\xi](t)=\eta\|f_\xi\|^2_{X(\hat m_\xi)}+\chi(s_\xi)\nn{f_\xi}^2_\natural
\end{equation}
for the level-$0$ energy of the leaf $\xi$, as in Subsection~\ref{ssec:cells}.

\emph{Step 1: the level-$0$ inequality after the band.} By the proof of Lemma~\ref{lem:UFlevel0}, we
have, for every $\xi$ and $t\ge t^-_\xi$, along the ray, in the integrated form \eqref{eq:L0int},
\begin{equation}
D_t\big(\tilde\omega e\big)+c_{\rm cl}\,\kappa\,\tilde\omega\,\lambda_\xi\|f_\xi\|^2_{Y_2(\hat
m_\xi)}\le t^{-2}C_a\,\tilde\omega e+t^{-2}C_b,
\end{equation}
with constants $C_a$, $C_b$ that do not depend on $t$, $\xi$ or $\eps$, and with $\|f_\xi\|_{X(\hat
m_\xi)}\le\delta_c$. Here $C_a=C+4\mathfrak p_*$, and $C_b$ contains the forcing constant of
Proposition~\ref{prop:UFforce}(e) and $P_1$. Both are fixed before $T_0$ in
Appendix~\ref{app:const}, and this $C_b$ is not the constant $C_b$ of Theorem~\ref{thm:CUboot}. For
$t\ge t^b_\xi$ we have $s_\xi\ge S_b+1$, where $\chi=1$ and the auxiliary weight $\tilde\omega$ is
constant, with $\tilde\omega\ge e^{-4C_1^B}$. At $\gamma=-2$ the zero-order part of the norm
$Y_2(\hat m_\xi)$ is $\|\jap{\hat z}^{(\gamma+2)/2}\hat m_\xi f\|_{L^2}=\|f\|_{X(\hat m_\xi)}$, and
$\nn f^2_\natural\le C_\natural\|f\|^2_{X(\hat m_\xi)}$ by the bounds of
Appendix~\ref{app:band}. Hence we have
$\|f_\xi\|^2_{Y_2(\hat m_\xi)}\ge e/(\eta+C_\natural)$. With $\lambda_\xi\ge c_\lambda\lambda_{\rm
ref}(\xi)$ and $c_{\rm rate}:=c_{\rm cl}c_\lambda/(\eta+C_\natural)$ we get, for $t\ge t^b_\xi$, in
the same integrated sense,
\begin{equation}\label{eq:R2ode}
D_te+\frac{c_{\rm rate}\lambda_{\rm ref}(\xi)}t\,e\le\frac{C_a}{t^2}\,e+\frac{C_b'}{t^2},\qquad
C_b':=e^{4C_1^B}C_b .
\end{equation}

\emph{Step 2: integration.} Since $e\ge0$, \eqref{eq:R2ode} holds with $a:=a_\xi\le c_{\rm
rate}\lambda_{\rm ref}(\xi)$ in place of $c_{\rm rate}\lambda_{\rm ref}(\xi)$. Then, formally,
$D_t\big(t^ae^{C_a/t}e(t)\big)\le C_b'e^{C_a/t}t^{a-2}$. Rigorously, we use \eqref{eq:L0int} on
$[t^b_\xi,t]$, where $\tilde\omega$ is constant, with the multiplier
$\psi(\tau)=\tau^ae^{C_a/\tau}$. Since $\psi'=(a/\tau-C_a/\tau^2)\psi$, the terms with $e$ on the
right add up to $\int(a/\tau)\psi\,\tilde\omega\,e\dd\tau$, and by the lower bound for the
dissipation of Step~1 this does not exceed the dissipation term on the left. Since $a\le\frac12$ and
$t^b_\xi\ge T_0\ge1$, we get
\begin{equation}
t^ae(t)\le e^{C_a}\Big((t^b_\xi)^ae(t^b_\xi)+C_b'\int_{t^b_\xi}^\infty\tau^{a-2}\dd\tau\Big)\le
e^{C_a}(t^b_\xi)^a\big(e(t^b_\xi)+2C_b'\big)\qquad(t\ge t^b_\xi).
\end{equation}
Since $e(t^b_\xi)\le(\eta+C_\natural)\delta_c^2$, we obtain
\begin{equation}\label{eq:R2e}
e_\xi[f_\xi](t)\le C\Big(\frac{t^b_\xi}t\Big)^{a_\xi}\qquad(t\ge t^b_\xi).
\end{equation}

\emph{Step 3: proof of~(i).} The left side of~(i) is the norm of $r_\xi(t)$ in the weight of
Theorem~\ref{thm:main}(ii). By Step~5 of the proof of Theorem~\ref{thm:scatter}, this weight is
bounded by $Ce^{(q_0/2)|z|^2}$, and we have
\begin{equation}
\int e^{2q(t)|z|^2-2L|z|}|r_\xi|^2\dd z\le\theta_g^{-3/2}\ell_\xi[\hat
r]=\theta_g^{-3/2}\|f_\xi\|^2_{X(\hat m_\xi)},
\end{equation}
with $2q(t)\ge\frac32q_0$. Thus the left side of~(i) is bounded by $C\|f_\xi\|_{X(\hat m_\xi)}\le
C\eta^{-1/2}e(t)^{1/2}$. Hence~(i) follows from~\eqref{eq:R2e}. If $c_{\rm rate}\lambda_{\rm
ref}(\xi)\le\frac12$, then $a_\xi=c_{\rm rate}\lambda_{\rm ref}(\xi)$ and
$(t^b_\xi/t)^{a_\xi}=\exp\big(-c_{\rm rate}(s_\xi(t)-S_b-1)\big)$.

\emph{Step 4: proof of~(ii).} We bound the quantity $I(t,\xi)$ of Step~7 of that proof, uniformly in
$\xi$. First, $I(t,\xi)\le C\rho_\xi$ for all $t\ge T_0$.
Second, let $t\ge t^b_\xi$. We estimate the three terms that bound $I(t,\xi)$. By the
Cauchy--Schwarz inequality and the bound of Step~3, we have
$\|r_\xi(t)\|_{L^1}\le\|e^{-q(t)|z|^2+L|z|}\|_{L^2}\,\theta_g^{-3/4}\|f_\xi\|_{X(\hat m_\xi)}\le
Ce(t)^{1/2}$, where $e(t)$ is the energy of Step~1. By Theorem~\ref{thm:scatter}(ii), we have
$\|M_{p_g(t,\xi)}-M_{p_{\rm sc}(\xi)}\|_{L^1}\le C\rho_\xi/t$. Moreover, $p_{\rm sc}$ is Lipschitz,
with $|\nabla\log\rho_{\rm sc}|+|\nabla y_{\rm sc}|+|\nabla\theta_{\rm sc}|\le C$ almost everywhere,
since the bounds on $\partial_\xi p_g(t,\cdot)$ of (F3) and the bound on $\nabla\log\rho_g$ are
uniform in $t$, so that $\int|M_{p_{\rm sc}(\xi)}(z)-M_{p_{\rm sc}(\xi-z/t)}(z)|\dd z\le
C\rho_\xi/t$. With~\eqref{eq:R2e}, we obtain
\begin{equation}
I(t,\xi)\le C\rho_\xi\Big(\min\Big\{1,\Big(\frac{t^b_\xi}t\Big)^{a_\xi/2}\Big\}+\frac1t\Big)\qquad\text{for
all }\xi\text{ and }t\ge T_0 .
\end{equation}
Let $\ell:=\log(t/T_0)$ and $c_2:=\min\{c_{\rm rate},(2\sup_\xi\lambda_{\rm ref})^{-1}\}$, so that
$a_\xi\ge c_2\lambda_{\rm ref}(\xi)$. Note that $\sup_\xi\lambda_{\rm ref}\le C\eps\le C$. Since
$(t^b_\xi/t)^{a_\xi/2}=\exp\big(-\frac{a_\xi}2(\ell-\frac{S_b+1}{\lambda_{\rm ref}(\xi)})\big)$ and
$a_\xi/\lambda_{\rm ref}(\xi)\le c_{\rm rate}$, we have
\begin{equation}
\min\Big\{1,\Big(\frac{t^b_\xi}t\Big)^{a_\xi/2}\Big\}\le
e^{c_{\rm rate}(S_b+1)/2}\,e^{-\frac{c_2}2\lambda_{\rm ref}(\xi)\ell}.
\end{equation}
By Theorem~\ref{thm:main}, $\lambda_{\rm
ref}(\xi)\ge\frac12\eps\rho_0(\xi)\theta_\infty(\xi)^{\gamma/2}\ge\frac12\eps\rho_0(\xi)(\theta_+/2)^{\gamma/2}$,
and $\rho_\xi=\lambda_{\rm ref}(\xi)\theta_\xi^{|\gamma|/2}\le C\eps\rho_0(\xi)$. By (D3), we have
$A_F^{-1}\jap\xi^{-p}\le\rho_0(\xi)\le A_F\jap\xi^{-p}$. Thus, with $\tau:=\eps\ell$ and some
$c_3>0$,
\begin{equation}
\|g(t)-g_{\rm sc}\|_{L^1}=\int I(t,\xi)\dd\xi\le C\eps\int\jap\xi^{-p}e^{-c_3\tau\jap\xi^{-p}}\dd\xi+\frac{C\eps}t
.
\end{equation}
For $\tau\le1$ the integral is bounded. For $\tau\ge1$, the part $|\xi|\le1$ is bounded by
$Ce^{-c_32^{-p/2}\tau}$, and on $|\xi|\ge1$, where $2^{-p/2}|\xi|^{-p}\le\jap\xi^{-p}\le|\xi|^{-p}$,
the substitution $|\xi|=\tau^{1/p}w$ yields
\begin{equation}
\int_{|\xi|\ge1}|\xi|^{-p}e^{-c_32^{-p/2}\tau|\xi|^{-p}}\dd\xi\le4\pi\,\tau^{-(p-3)/p}\int_0^\infty
w^{2-p}e^{-c_32^{-p/2}w^{-p}}\dd w\le C\tau^{-(p-3)/p},
\end{equation}
since $p>3$. Hence the integral is bounded by $C(1+\tau)^{-(p-3)/p}$. Finally,
$t^{-1}\le(1+\ell)^{-1}\le(1+\eps\ell)^{-(p-3)/p}$ for $\eps\le1$, and
$\|f_0\|_{L^1}=\eps\int\rho_0\ge c\eps$. This is~(ii).

\emph{Step 5: proof of~(iii).} Let $d_1>0$. Taking $c=0$ in the definition of $d_1(v)$, we get
$d_1\le
d_1(v)\le\|F_0(\cdot,v)\|_{L^1}=1$. In this step $C_{\rm d}$, $C_{\rm p}$, $C$ and $c>0$ denote
constants chosen before $T_0$, which do not depend on $T_0$, $\eps$, $d_1$,
$t$ or $\xi$.

\emph{(a) The leaves stay close to their initial profiles in their own clock.} For every $\xi$ and
$t\ge T_0$,
\begin{equation}\label{eq:R2F0}
\|\mathcal N_\xi(t)-F_0(\cdot,\xi)\|_{L^1(\R^3)}\le C_{\rm
d}\Big(s_\xi(t)+\frac1{T_0}+\eps\,(1+\log T_0)\Big) .
\end{equation}
Indeed, by Lemma~\ref{lem:Wtail}(ii), $\|\mathcal N_\xi(t)-
\mathcal N_\xi(T_0)\|_{L^2(e^{a|z|^2/2})}\le C_{\rm tail}(s_\xi(t)+T_0^{-1})$. By the first item of the proof
of Theorem~\ref{thm:main}(iv), $\|
\mathcal N_\xi(T_0)-F_0(\cdot,\xi)\|_{L^2(e^{a'|z|^2})}\le C_{T_0}(L/T_0+\eps T_0^{\kappa_\gamma})$, where
$C_{T_0}=C(\Dc)$ and, at $\gamma=-2$, $T_0^{\kappa_\gamma}$ is read as $1+\log T_0$. The
Cauchy--Schwarz inequality turns both bounds into bounds in $L^1$, with factors that depend only on
$a$ and $a'$.

\emph{(b) A lower bound on the unrelaxed leaves.} Let $b\ge0$, and let $M$ be a Maxwellian of unit
mass. We write $F^\xi:=F_0(\cdot,\xi)$. By the definition of $d_1$ and the triangle inequality,
$\|\mathcal N_\xi(t)-bM\|_{L^1}\ge\|bM-F^\xi\|_{L^1}-\|\mathcal N_\xi(t)-F^\xi\|_{L^1}\ge
d_1-\|\mathcal N_\xi(t)-F^\xi\|_{L^1}$. We apply this with $b=\rho_{\rm sc}(\xi)/\rho_g(t,\xi)$ and
$M=M_{(1,y_{\rm sc}(\xi),\theta_{\rm sc}(\xi))}$, so that $\rho_g(t,\xi)\,bM=M_{p_{\rm sc}(\xi)}$.
With $H=\rho_g\mathcal N_\xi$ and the quantity $I(t,\xi)$ above, by the triangle inequality we get
\begin{equation}
I(t,\xi)\ge\rho_g(t,\xi)\big(d_1-\|\mathcal N_\xi(t)-F^\xi\|_{L^1}\big)-\int\big|M_{p_{\rm
sc}(\xi)}(z)-M_{p_{\rm sc}(\xi-z/t)}(z)\big|\dd z .
\end{equation}
As in Step~4, the last integral is bounded by $C_{\rm p}\rho_\xi/t$. The constant does not depend on
$T_0$, since it depends only on $L$, $\theta_\pm$ and the Lipschitz bounds on $p_{\rm sc}$, and
these are $|\nabla y_{\rm sc}|+|\nabla\theta_{\rm sc}|\le C_1(L_1,R_1,\rad_{\rm post})$, from
(F3), and $|\nabla\log\rho_{\rm sc}|\le L_1+1+\delta_s$, from Lemma~\ref{lem:Mref} and
Theorem~\ref{thm:main}(iii). We put
\begin{equation}
\sigma:=\frac{d_1}{8C_{\rm d}},\qquad C_{\rm sh}:=16(C_{\rm d}+C_{\rm p}),
\end{equation}
and we assume $T_0\ge C_{\rm sh}/d_1$ and $C_{\rm sh}\eps(1+\log T_0)\le d_1$. Let $t\ge T_0$ and
$s_\xi(t)\le\sigma$. Then the right side of~\eqref{eq:R2F0} does not exceed
$\frac18d_1+\frac1{16}d_1+\frac1{16}d_1=\frac14d_1$, and $C_{\rm p}/t\le C_{\rm
p}/T_0\le\frac1{16}d_1$. Since $\rho_g\ge\frac12\rho_\xi$ by (F3), we have
\begin{equation}\label{eq:R2low}
I(t,\xi)\ge\tfrac12\rho_\xi\cdot\tfrac34d_1-\tfrac1{16}d_1\rho_\xi\ge\tfrac14\,d_1\,\rho_\xi\qquad\text{whenever
}s_\xi(t)\le\sigma .
\end{equation}

\emph{(c) Integration.} As in Step~4, let $\ell=\log(t/T_0)$ and $\tau=\eps\ell$. By
the same theorem and (D3), $\lambda_{\rm
ref}(\xi)\le2\eps\rho_0(\xi)\theta_\infty(\xi)^{-1}\le A_F\theta_-^{-1}\eps\jap\xi^{-p}$, since
$\theta_\infty(\xi)\ge2\theta_-$. Therefore $s_\xi(t)=\lambda_{\rm ref}(\xi)\ell\le\sigma$ whenever
$\jap\xi^p\ge K\tau$, where
\begin{equation}
K:=\frac{A_F}{\theta_-\sigma}=\frac{8A_FC_{\rm d}}{\theta_-d_1} .
\end{equation}
We may take $C_{\rm d}\ge\theta_-/(8A_F)$, so that $K\ge1$, since $d_1\le1$. Moreover
$\rho_\xi=\lambda_{\rm
ref}(\xi)\theta_\xi\ge\frac12\eps\rho_0(\xi)\theta_-/\theta_\infty(\xi)\ge\eps\rho_0(\xi)\theta_-/\theta_+\ge
c\,\eps\jap\xi^{-p}$, where $\theta_\xi\ge\theta_-$ is the temperature of the fiber
$g(T_0,\cdot,\xi)$. The set $\{\jap\xi^p\ge K\tau\}$ contains $\{|\xi|\ge R\}$ with
$R:=\max\{1,(K\tau)^{1/p}\}$, and there $\jap\xi^{-p}\ge2^{-p/2}|\xi|^{-p}$. By~\eqref{eq:R2low}, we
have
\begin{equation}
\|g(t)-g_{\rm sc}\|_{L^1}=\int I(t,\xi)\dd\xi\ge\frac{d_1}4\,c\,\eps\,2^{-p/2}\int_{|\xi|\ge
R}|\xi|^{-p}\dd\xi=\frac{d_1}4\,c\,\eps\,2^{-p/2}\,\frac{4\pi}{p-3}\,R^{-(p-3)} .
\end{equation}
Now $R^p=\max\{1,K\tau\}\le K(1+\tau)$, so $R^{-(p-3)}\ge K^{-(p-3)/p}(1+\tau)^{-(p-3)/p}$, and
$K^{-(p-3)/p}\ge c\,d_1^{(p-3)/p}$. Since $\|f_0\|_{L^1}=\eps\int\rho_0\le C\eps$, this is~(iii),
with the power $d_1^{1+(p-3)/p}=d_1^{2-3/p}$.
\end{proof}

For the rates at $\gamma<-2$ we use a second level-$0$ energy, in a weaker Gaussian weight. Let
$q_-:=\frac38q_0\theta_-$ and $q_+:=q_0(\theta_++1)$. By Lemma~\ref{lem:TG0}(i) and (F3), and since
$q(t)\in[\frac34q_0,q_0]$, the leaf weights $\hat m_\xi=e^{q_\xi|\hat z|^2+\omega}$ belong to the
class $\mathcal W(q_-,q_+,\eps_1)$ with $\eps_1=\eps_1(t)$, and their exponents satisfy
$q_\xi\le\frac18$ and $C\eps_1(t)^2\le q_\xi$ for $t\ge T_1$, by the discussion of the Gaussian
hypodissipativity in Appendix~\ref{app:band}. We put
\begin{equation}
\check m_\xi:=\hat m_\xi^{1/2}=e^{\frac12q_\xi|\hat z|^2+\frac12\omega} .
\end{equation}
Then $\check m_\xi\in\mathcal W(\frac12q_-,\frac12q_+,\frac12\eps_1)$, with the same bound on
$\omega(0)$, and its exponent $\frac12q_\xi$ satisfies $C(\frac12\eps_1)^2\le\frac12q_\xi$. The
proof of Lemma~\ref{lem:GW}(i) rests on three properties of the weight: the bound on $\nabla\omega$,
the bound $\frac18$ on the exponent, and this inequality between the gradient bound and the
exponent, which enters in the last paragraph of that proof. Thus its conclusion holds for
$\check m_\xi$, and with it
Corollary~\ref{cor:W1G}. Moreover, Lemma~\ref{lem:TG0}(ii)--(iii) and Proposition~\ref{prop:trilinG}
hold for every weight of the class. We write $\check c_W$, $\check C_W$ for the constants of
(W1$_G$) for this class, $\check C_\tau$ for the constant of the second form of
Proposition~\ref{prop:trilinG}(iii), $\check C$ for that of part~(iv), and
$\check C_{11}$ for the constant of Proposition~\ref{prop:FM}(ii) with $Y_2(\check m_\xi)$ in place
of $Y_2(\hat m_\xi)$. They depend only on $\gamma$, $q_\pm$, $\eps_1$, $\theta_\pm$ and the bound on
$\omega(0)$, and are fixed before $\delta_c$. We put $\check\eta:=1/(4\check C_W)$ and
$\check C_{\rm r}:=\check C_{11}\mathfrak p_*$.

\begin{lem}[The level-$0$ inequality in a weaker weight]\label{lem:R3weak}
Let $\gamma\in[-3,-2)$. Assume (U1)--(U4) on $[T_0,T]$, the requirements of
Lemma~\ref{lem:UFlevel0}, and requirements~\ref{rq:Rdc} and~\ref{rq:Rrate}. For a leaf $\xi$ let
$t^b_\xi$ be the time with $s_\xi(t^b_\xi)=S_b+1$, and for $t\in[t^b_\xi,T]$ put
\begin{equation}
\check e_\xi(t):=\check\eta\,\|f_\xi\|^2_{X(\check m_\xi)}+\nn{f_\xi}^2_\natural,\qquad\check
d_\xi(t):=\|f_\xi\|^2_{Y_2(\check m_\xi)} .
\end{equation}
There are constants $\check C_a$, $\check C_f$ and $C_I$, fixed before $T_0$, which do not depend on
$\xi$, $t$, $T$ or $\eps$, such that:
\begin{enumerate}[label=(\roman*)]
\item for every positive $\psi\in C^1([t^b_\xi,T])$ and all $t^b_\xi\le t'\le t\le T$,
\begin{equation}\label{eq:R3int}
\Big[\psi\,\check e_\xi\Big]_{t'}^{t}+\int_{t'}^{t}\psi\,\check\eta\check
c_W\,\kappa\lambda_\xi\,\check d_\xi\dd\tau'\le\int_{t'}^{t}\Big(\psi'\,\check
e_\xi+\frac{\psi}{\tau'^{2}}\big(\check C_a\,\check e_\xi+\check C_f\,\check
e_\xi^{1/2}\big)\Big)\dd\tau' ;
\end{equation}
\item for every $t\in[t^b_\xi,T]$ and every $R\ge0$,
\begin{equation}
\check e_\xi\le C_I\,\delta_c^2,\qquad\check e_\xi\le C_I\jap R^{k_\gamma}\,\check
d_\xi+C_I\,\delta_c^2\,e^{-q_-R^2/2},\qquad k_\gamma=12+|\gamma| .
\end{equation}
\end{enumerate}
\end{lem}

\begin{proof}
We write $f=f_\xi$, $\check m=\check m_\xi$, $\hat m=\hat m_\xi$, $X=X(\hat m)$, $\check X=X(\check
m)$
and $\check Y=Y_2(\check m)$. For $t\ge t^b_\xi$ the leaf is post-band: $\chi(s_\xi)=1$, and
$\|f\|_X\le\delta_c$ by the first paragraph of the proof of Lemma~\ref{lem:UFlevel0}. Since $\check
m^2=\hat m$, we have
\begin{equation}
\|f\|^2_{\check X}=\int\hat mf^2\le\|f\|_X\|f\|_{L^2}\le
C_{L^2}\|f\|^2_X,\qquad\text{so}\qquad\|f\|_{\check X}\le C_{L^2}^{1/2}\delta_c .
\end{equation}

(i) We differentiate $\check e_\xi$ along the ray, as in Subsection~\ref{sssec:bandproof} and in the
proof of that lemma, with $\check m$ in place of $\hat m$, $\chi=1$ and no damping
factor, and we use the numbering (T1)--(T12) of the terms. The nonlocal part of the energy is the
same as before.
\begin{itemize}
\item \emph{Collision terms.} At level $0$ the collision part of $D_tf$ is
$\kappa\lambda_\xi\big(\mathcal Lf+Q(f,f)\big)$. By (W1$_G$) for $\check m$ with background $0$,
$\langle f,\mathcal Lf\rangle_{\check X}\le-\check c_W\|f\|^2_{\check Y}+\check C_W\|f\|^2_{L^2}$.
By Lemma~\ref{lem:lowweight}, we have $2\ip f{\mathcal Lf}_\natural=-\|f\|^2_{L^2}$. Moreover, by
Proposition~\ref{prop:trilinG}(iii)--(iv) for $\check m$, we have $|\langle f,Q(f,f)\rangle_{\check
X}|\le2\check C_\tau\|f\|_{\check X}\|f\|^2_{\check Y}$ and $|\ip f{Q(f,f)}_\natural|\le\check
C\|f\|_{\check X}\|f\|^2_{\check Y}$. Since $2\check\eta\check C_W=\frac12$, the contribution of the
collision terms to the derivative of $\check e_\xi$ is bounded by
\begin{equation}
\kappa\lambda_\xi\Big(-2\check\eta\check c_W\|f\|^2_{\check
Y}-\tfrac12\|f\|^2_{L^2}+\big(4\check\eta\check C_\tau+2\check C\big)\|f\|_{\check
X}\|f\|^2_{\check Y}\Big)\le-\tfrac32\check\eta\check c_W\,\kappa\lambda_\xi\|f\|^2_{\check Y},
\end{equation}
by requirement~\ref{rq:Rdc} and $\|f\|_{\check X}\le C_{L^2}^{1/2}\delta_c$.
\item \emph{The weight term} (T8). The part of $D_t\check m^2/\check m^2=\frac12D_t\hat m^2/\hat
m^2$ that is not caused by the frame is $q'(t)|z|^2+t^{-2}z\cdot\nabla\log A_{\rm ref}(\xi-z/t)\le
t^{-2}\big(-C_2|z|^2+L|z|\big)\le t^{-2}L^2/(4C_2)$. Its contribution is bounded by
$Ct^{-2}\check\eta\|f\|^2_{\check X}$.
\item \emph{The local transport} (T9), $2t^{-2}\check\eta\langle f,(y_0+\theta^{1/2}\hat
z)\cdot\varphi^+\rangle_{\check X}$ with $\varphi^+=(\varphi_{0e_j})_j$. Here the weaker weight
leaves a reserve for the factor $z=y_0+\theta^{1/2}\hat z$. Indeed, by (F3) and the lower bound of
the class, we have $|z|\,\check m\le C\jap{\hat z}\hat m^{-1/2}\,\hat m\le C_z\hat m$. Thus this
term is bounded by
$2t^{-2}\check\eta\,C_z\|f\|_{\check X}\|\varphi^+\|_X\le Ct^{-2}\check e_\xi^{1/2}P_1^{1/2}$, with
the bound $P_1$ of Lemma~\ref{lem:UFlevel0}.
\item \emph{The nonlocal transport} (T10) and \emph{the kernel part of the forcing} (T12$_0$) are
bounded as in that subsection, with Lemma~\ref{lem:TG0}(ii)--(iii) for $\check m$
on the factor $f$, and they are bounded by $Ct^{-2}\|f\|_{\check X}\|\varphi^+\|_X\le Ct^{-2}\check
e_\xi^{1/2}P_1^{1/2}$.
\item \emph{The frame terms} (T11). The identity of Proposition~\ref{prop:FM}(i) holds for every
weight that is a function of $t$ and $z$, so the local part is $2D_t\log
A_p\cdot\check\eta\|f\|^2_{\check X}\le4t^{-2}\mathfrak p_*\check\eta\|f\|^2_{\check X}$. The
nonlocal part does not exceed $\check C_{\rm r}t^{-2}\|f\|^2_{\check Y}$, by the proof of
part~(ii) of that proposition, with the same comparison for $\check m$. For $t\ge t^b_\xi$ we
have $\kappa\lambda_\xi\ge c_\lambda s_\xi\kappa/\tau\ge c_\lambda S_*(|\gamma|-2)/t$
as in Remark~\ref{rem:Bbandrange}, so requirement~\ref{rq:Rrate} implies $\check C_{\rm
r}t^{-2}\le\frac14\check\eta\check c_W\kappa\lambda_\xi$.
\item \emph{The rest of the forcing} (T12). With $\hat F=\hat F_{00,\xi}$, it is
$-2\check\eta\langle f,\hat F\rangle_{\check X}-2\ip f{\hat F}_\natural\le C\check e_\xi^{1/2}\|\hat
F\|_X\le Ct^{-2}\check e_\xi^{1/2}F_0$, by the same two lemmas and
Proposition~\ref{prop:UFforce}(e), since $\check m\le C\hat m$.
\end{itemize}
There are no sources at level $0$, and neither a deficit nor a front term, since $\chi=1$. Adding
up, and keeping $\check\eta\check c_W$ of the dissipation, we obtain
\begin{equation}
D_t\check e_\xi+\check\eta\check c_W\,\kappa\lambda_\xi\,\check d_\xi\le t^{-2}\big(\check
C_a\,\check e_\xi+\check C_f\,\check e_\xi^{1/2}\big),\qquad\check C_a:=C+4
\mathfrak p_*,\quad\check C_f:=C\big(F_0+P_1^{1/2}\big).
\end{equation}
As before, this is meant in the integrated form~\eqref{eq:R3int}, and it is
proved in the same way: the local part of the energy is written in $z$-variables and truncated by
$\zeta(z/R')^2$, the truncated energy is Lipschitz along the ray, its derivative is the sum of the
pairings above with the cutoff under the integral, and we let $R'\to\infty$. Here every local
term, including (T8) and (T9), converges by dominated convergence, since $|z|\check
m^2|h||\partial_\xi h|\le C\hat m^2|h||\partial_\xi h|$ and $|z|^2\check m^2h^2\le C\hat m^2h^2$ are
integrable in $z$, uniformly on $[t',t]$.

(ii) By the nonlocal-pairing lemma with $k'=7$, we have $\nn f^2_\natural\le
C_\natural(7)\|f\|^2_{X_7}=C_\natural(7)\int\jap{\hat z}^{14}f^2$. By the two-sided bound of the
class $\mathcal W(q_-,q_+,\eps_1)$ there are $c_m>0$ and $C_m$ with $\hat m\ge c_m$ and $\jap{\hat
z}^{14}\hat m^{-1}\le C_me^{-q_-|\hat z|^2/2}$. We split the integrals at $|\hat z|=R$. On $\{|\hat
z|\le R\}$ we use $1\le\jap R^{|\gamma|-2}\jap{\hat z}^{\gamma+2}$ and $\check m^2=\hat m\ge c_m$:
\begin{equation}
\int_{|\hat z|\le R}\check m^2f^2\le\jap R^{|\gamma|-2}\,\check d_\xi,\qquad\int_{|\hat z|\le
R}\jap{\hat z}^{14}f^2\le c_m^{-1}\jap R^{12+|\gamma|}\,\check d_\xi,
\end{equation}
since the zero-order part of $\|f\|^2_{Y_2(\check m)}$ is $\int\jap{\hat z}^{\gamma+2}\check
m^2f^2$. On $\{|\hat z|>R\}$ we use the stronger weight:
\begin{equation}
\int_{|\hat z|>R}\check m^2f^2=\int_{|\hat z|>R}\hat m^{-1}\,\hat m^2f^2\le
C_me^{-q_-R^2/2}\|f\|^2_X,\qquad\int_{|\hat z|>R}\jap{\hat z}^{14}f^2\le
c_m^{-1}C_me^{-q_-R^2/2}\|f\|^2_X .
\end{equation}
With $\|f\|_X\le\delta_c$ this is the second inequality, with
$C_I:=\max\{1,(\check\eta+C_\natural(7)c_m^{-1})(1+C_m)\}$. The first one follows from $\check
e_\xi\le(\check\eta C_{L^2}+C_\natural)\|f\|^2_X$, after enlarging $C_I$.
\end{proof}

\begin{proof}[Proof of Theorem~\ref{thm:rate3}]
Let $\gamma\in[-3,-2)$ and $\kappa_\gamma=|\gamma|-2$, so that $\kappa(t)=t^{\kappa_\gamma-1}$ and
$\tau(t)=(t^{\kappa_\gamma}-1)/\kappa_\gamma$. We put $\mathfrak T(t):=\tau(t)-\tau(T_0)$. Then
$s_\xi(t)=\lambda_{\rm ref}(\xi)\mathfrak T(t)$ and $\mathfrak T(t)\le
t^{\kappa_\gamma}/\kappa_\gamma$. We use the notation of the proof of Theorem~\ref{thm:scatter}. By
Corollary~\ref{cor:CUmain}, (U1)--(U4) hold for all $t\ge T_0$, so Lemma~\ref{lem:R3weak} applies
with every $T$. In Steps~1--3 we fix a leaf $\xi$ and write $t_b=t^b_\xi$, $\check e=\check e_\xi$,
$\check d=\check d_\xi$, $s=s_\xi$, $\dot s=\kappa\lambda_{\rm ref}(\xi)$, $\varrho=\varrho_\xi$ and
\begin{equation}
u(t):=s(t)-S_b-1\ge0,\qquad\Lambda(t):=\log\Big(e+\frac1{\varrho(t)}\Big)\ge1\qquad(t\ge t_b).
\end{equation}
Note that $\varrho=t^{-2}/\dot s$, that $\varrho$ is decreasing, and that for $t\ge t_b$ we have
$1\le s(t)\le\lambda_{\rm ref}(\xi)t^{\kappa_\gamma}/\kappa_\gamma$, hence
\begin{equation}\label{eq:R3rho}
\varrho(t)=\frac{t^{1+\gamma}}{\lambda_{\rm ref}(\xi)}\le\frac{t^{1+\gamma+\kappa_\gamma}}{\kappa_\gamma}=\frac1{\kappa_\gamma
t}\le\frac1{\kappa_\gamma}\qquad(t\ge t_b).
\end{equation}
In this proof $C$ and $c>0$ denote constants that depend only on $\gamma$ and on constants chosen
before $T_0$, and they do not depend on $t$, $\xi$, $T_0$ or $\eps$, unless
stated.

\emph{Step 1: a linear inequality.} Let $\vartheta=a_\gamma=2/(k_\gamma+2)$ and
$k=k_\gamma$, and for $t\ge t_b$ we put
\begin{equation}
R(t)^2:=\frac2{q_-}\min\big\{(1+u(t))^\vartheta,\,2\Lambda(t)\big\},\qquad\nu(t):=\frac{\check\eta\check
c_Wc_\lambda}{C_I\jap{R(t)}^{k}} .
\end{equation}
These are continuous functions of $t$. By Lemma~\ref{lem:R3weak}(ii) with $R=R(t)$, and
$\lambda_\xi\ge c_\lambda\lambda_{\rm ref}(\xi)$, we have
\begin{equation}
\check\eta\check c_W\,\kappa\lambda_\xi\,\check d\ge\nu\,\dot s\,\big(\check
e-C_I\delta_c^2e^{-q_-R^2/2}\big),\qquad e^{-q_-R^2/2}=\max\big\{e^{-(1+u)^\vartheta},(e+\varrho^{-1})^{-2}\big\}\le
e^{-(1+u)^\vartheta}+\varrho^2 .
\end{equation}
By Young's inequality and $t^{-4}/\dot s=\varrho^2\dot s$, we get $t^{-2}\check C_f\check
e^{1/2}\le\frac12\nu\dot s\,\check e+\check C_f^2\varrho^2\dot s/(2\nu)$. We insert both
in~\eqref{eq:R3int} and choose
\begin{equation}
\psi(t):=\exp\Big(\tfrac12N(t)+\frac{\check C_a}t\Big),\qquad N(t):=\int_{t_b}^t\nu\,\dot s\dd\tau',
\end{equation}
so that $\psi'=(\frac12\nu\dot s-\check C_at^{-2})\psi$ and the terms with $\check e$ under the
integrals cancel. Since $\check C_a(\tau'^{-1}-t^{-1})\le\check C_a/T_0\le\check C_a$, we obtain,
for $t\ge t_b$,
\begin{equation}\label{eq:R3duh}
\check e(t)\le e^{\check C_a}\Big(e^{-N(t)/2}\check
e(t_b)+\int_{t_b}^te^{-(N(t)-N(\tau'))/2}\Big(C_I\delta_c^2\,\nu\,\big(e^{-(1+u)^\vartheta}+\varrho^2\big)+\frac{\check
C_f^2\varrho^2}{2\nu}\Big)\dot s\dd\tau'\Big),
\end{equation}
with $u$, $\varrho$ and $\nu$ evaluated at $\tau'$.

\emph{Step 2: bounds on $\nu$ and $N$.} Let $c_q:=2/q_-$. Since $1+R^2\le(1+c_q)(1+u)^\vartheta$ and
$1+R^2\le(1+2c_q)\Lambda$, we have
\begin{equation}
\nu\ge\nu_0(1+u)^{-k\vartheta/2},\qquad\nu\ge\nu_1\Lambda^{-k/2},\qquad\nu\le\nu_{\max}:=\frac{\check\eta\check
c_Wc_\lambda}{C_I},
\end{equation}
with $\nu_0:=\nu_{\max}(1+c_q)^{-k/2}$ and $\nu_1:=\nu_{\max}(1+2c_q)^{-k/2}$. Since $\dot
s\dd\tau'=\dd u$ and $1-k\vartheta/2=\vartheta$, we have, for $t_b\le\tau'\le t$,
\begin{equation}\label{eq:R3N}
N(t)-N(\tau')\ge\frac{\nu_0}\vartheta\Big(\big(1+u(t)\big)^\vartheta-\big(1+u(\tau')\big)^\vartheta\Big),\qquad
N(t)\le\nu_{\max}u(t).
\end{equation}

\emph{Step 3: the three terms of~\eqref{eq:R3duh}.} We write $u=u(t)$. First, by~\eqref{eq:R3N} and
$\check e(t_b)\le C_I\delta_c^2$, we have $e^{-N(t)/2}\check e(t_b)\le
Ce^{-\frac{\nu_0}{2\vartheta}(1+u)^\vartheta}$. Second, with $c_4:=\min\{\nu_0/(2\vartheta),1\}$,
we deduce from the same bound that
$-\frac12(N(t)-N(\tau'))-(1+u(\tau'))^\vartheta\le-c_4(1+u)^\vartheta$, because
$-\alpha(X-Y)-Y\le-\min\{\alpha,1\}X$ for $X\ge Y\ge0$. Hence
\begin{equation}
\int_{t_b}^te^{-(N(t)-N(\tau'))/2}\,\nu\,e^{-(1+u(\tau'))^\vartheta}\dot s\dd\tau'\le
e^{-c_4(1+u)^\vartheta}N(t)\le\nu_{\max}\,u\,e^{-c_4(1+u)^\vartheta}\le
Ce^{-\frac{c_4}2(1+u)^\vartheta} .
\end{equation}
Third, since $\nu\le\nu_{\max}$, the remaining terms are bounded by
$C\int_{t_b}^te^{-(N(t)-N(\tau'))/2}\,\nu\,(\varrho/\nu)^2\dot s\dd\tau'$, and
\begin{equation}
\int_{t_b}^te^{-(N(t)-N(\tau'))/2}\,\nu\,\Big(\frac\varrho\nu\Big)^2\dot
s\dd\tau'\le4\sup_{t_b\le\tau'\le t}G(\tau'),\qquad
G(\tau'):=e^{-(N(t)-N(\tau'))/4}\Big(\frac{\varrho(\tau')}{\nu(\tau')}\Big)^2,
\end{equation}
because $\int_{t_b}^te^{-(N(t)-N(\tau'))/4}\nu\dot s\dd\tau'\le4$. We bound $G(\tau')$ in two cases.

\emph{Case $\tau'\ge t/2$.} Then $\varrho(\tau')\le2^{|\gamma|-1}\varrho(t)\le4\varrho(t)$, and
$\Lambda(\tau')\le\Lambda(t)$ since $\varrho$ is decreasing. Hence
$\nu(\tau')^{-1}\le\nu_1^{-1}\Lambda(t)^{k/2}$, and
$G(\tau')\le16\nu_1^{-2}\varrho(t)^2\Lambda(t)^k$.

\emph{Case $\tau'<t/2$.} Then $s(t)-s(\tau')\ge\lambda_{\rm
ref}(\xi)\big(\tau(t)-\tau(t/2)\big)=(1-2^{-\kappa_\gamma})\lambda_{\rm
ref}(\xi)t^{\kappa_\gamma}/\kappa_\gamma\ge c_5\,u$, with $c_5:=1-2^{-\kappa_\gamma}\in(0,1)$, since
$u\le s(t)\le\lambda_{\rm ref}(\xi)t^{\kappa_\gamma}/\kappa_\gamma$. Therefore
$u(\tau')\le(1-c_5)u$, and, since $x\mapsto x^\vartheta$ is subadditive,
$(1+u(\tau'))^\vartheta\le\big((1-c_5)(1+u)+c_5\big)^\vartheta\le(1-c_5)^\vartheta(1+u)^\vartheta+1$.
By~\eqref{eq:R3N}, we have $N(t)-N(\tau')\ge\frac{\nu_0}\vartheta\big(c_6(1+u)^\vartheta-1\big)$
with
$c_6:=1-(1-c_5)^\vartheta$. Moreover $\varrho(\tau')\le1/\kappa_\gamma$ by~\eqref{eq:R3rho}, and
$\nu(\tau')^{-1}\le\nu_0^{-1}(1+u)^{k\vartheta/2}$. Hence
\begin{equation}
G(\tau')\le\frac{e^{\nu_0/(4\vartheta)}}{\kappa_\gamma^2\nu_0^2}\,(1+u)^{k\vartheta}\,e^{-\frac{\nu_0c_6}{4\vartheta}(1+u)^\vartheta}\le
Ce^{-\frac{\nu_0c_6}{8\vartheta}(1+u)^\vartheta} .
\end{equation}
Altogether, with
$c:=\frac12\min\{\frac{\nu_0}{2\vartheta},\frac{c_4}2,\frac{\nu_0c_6}{8\vartheta}\}$ and
$(1+u)^\vartheta\ge u^\vartheta$, we obtain
\begin{equation}\label{eq:R3e}
\check e_\xi(t)\le
C\Big(\exp\big(-2c\,(s_\xi(t)-S_b-1)^{a_\gamma}\big)+\varrho_\xi(t)^2\,\Lambda(t)^{k_\gamma}\Big)\qquad(t\ge
t^b_\xi).
\end{equation}

\emph{Step 4: proof of~(i).} The left side of~(i) is the norm of $r_\xi(t)$ in the weight of
Theorem~\ref{thm:main}(ii). As in Step~5 of the proof of Theorem~\ref{thm:scatter}, $r_\xi=\hat
r^{(t\xi)}/\rho_g$, and in $z$-variables $\hat m_\xi=e^{q(t)|z|^2}A_{\rm ref}(\xi-z/t)/A_{\rm
ref}(\xi)\ge e^{q(t)|z|^2-L|z|}$. Thus, with the first power of the weight, we have
\begin{equation}
\int e^{q(t)|z|^2-L|z|}|r_\xi|^2\dd z\le\theta_g^{-3/2}\|f_\xi\|^2_{X(\check
m_\xi)}\le\theta_g^{-3/2}\check\eta^{-1}\check e_\xi .
\end{equation}
The square of this weight is
$e^{q_0\theta_-|z-y_{0,g}|^2/(4\theta_g)}\le e^{q_0|z-y_{0,g}|^2/2}\le Ce^{\frac58q_0|z|^2}$, by
(F3), and $q(t)\ge\frac34q_0$, so it is bounded by $Ce^{q(t)|z|^2-L|z|}$. Hence the left side of~(i)
does not exceed $C\check e_\xi(t)^{1/2}$, and~(i) follows from~\eqref{eq:R3e}. The bound on
$\varrho_\xi$ is~\eqref{eq:R3rho}.

\emph{Step 5: proof of~(ii).} We bound the quantity $I(t,\xi)$ of Step~7 of the same proof. As in
Step~4 of the proof of Theorem~\ref{thm:rate2}, $I(t,\xi)\le
C\rho_\xi$ for all $t\ge T_0$, and for $t\ge t^b_\xi$
\begin{equation}
I(t,\xi)\le C\rho_\xi\Big(\|r_\xi(t)\|_{L^1}+\frac1t\Big),\qquad\|r_\xi(t)\|_{L^1}\le\big\|e^{-(q(t)|z|^2-L|z|)/2}\big\|_{L^2}\,\theta_g^{-3/4}\|f_\xi\|_{X(\check
m_\xi)}\le C\check e_\xi(t)^{1/2},
\end{equation}
by the Cauchy--Schwarz inequality and Step~4. Here $C$ may depend on $T_0$, through
Theorem~\ref{thm:scatter}(ii). If $s_\xi(t)\le2(S_b+1)$, then $1\le
e^{c(2S_b+2)^{\vartheta}}e^{-cs_\xi(t)^{\vartheta}}$. If $s_\xi(t)>2(S_b+1)$, then
$s_\xi(t)-S_b-1\ge\frac12s_\xi(t)$. Thus, with $c':=2^{-\vartheta}c$ and~\eqref{eq:R3e}, for all
$\xi$ and $t\ge T_0$,
\begin{equation}
I(t,\xi)\le C\rho_\xi\Big(e^{-c's_\xi(t)^{\vartheta}}+\mathbf 1_{\{s_\xi(t)\ge
S_b+1\}}\,\varrho_\xi(t)\,\Lambda_\xi(t)^{k/2}+\frac1t\Big),\qquad\Lambda_\xi(t)=\log\big(e+\varrho_\xi(t)^{-1}\big).
\end{equation}
We integrate the three terms in $\xi$. By Theorem~\ref{thm:main} and (D3),
$c\,\eps\jap\xi^{-p}\le\lambda_{\rm ref}(\xi)\le C\eps\jap\xi^{-p}$, and $\rho_\xi=\lambda_{\rm
ref}(\xi)\theta_\xi^{|\gamma|/2}\le C\eps\jap\xi^{-p}$. We put $\tau_\eps:=\eps\,\mathfrak T(t)$.

\emph{The first term.} Since $s_\xi(t)\ge c\,\tau_\eps\jap\xi^{-p}$, it is bounded by
$C\eps\int\jap\xi^{-p}e^{-c_3(\tau_\eps\jap\xi^{-p})^\vartheta}\dd\xi$, with some $c_3>0$. For
$\tau_\eps\le1$ the integral is bounded. For $\tau_\eps\ge1$ the part $|\xi|\le1$ is at most
$Ce^{-c_32^{-p\vartheta/2}\tau_\eps^\vartheta}$, and on $|\xi|\ge1$ the substitution
$|\xi|=\tau_\eps^{1/p}w$ yields
\begin{equation}
\int_{|\xi|\ge1}|\xi|^{-p}e^{-c_32^{-p\vartheta/2}(\tau_\eps|\xi|^{-p})^\vartheta}\dd\xi\le4\pi\,\tau_\eps^{-(p-3)/p}\int_0^\infty
w^{2-p}e^{-c_32^{-p\vartheta/2}w^{-p\vartheta}}\dd w\le C\tau_\eps^{-(p-3)/p},
\end{equation}
since $p>3$. Thus the contribution of the first term is bounded by $C\eps(1+\tau_\eps)^{-(p-3)/p}$.

\emph{The second term.} Here $\rho_\xi\varrho_\xi(t)=\theta_\xi^{|\gamma|/2}t^{1+\gamma}\le
Ct^{1+\gamma}$, and $\Lambda_\xi(t)\le\log(e+C\eps t^{|\gamma|-1})\le C\log(e+t)$. The set
$\{s_\xi(t)\ge S_b+1\}$ is contained in $\{\jap\xi^p\le C\tau_\eps\}$. It is empty if
$C\tau_\eps<1$, and its measure does not exceed $C\tau_\eps^{3/p}$. Therefore the contribution of
the second term is bounded by
\begin{equation}
C\,t^{1+\gamma}\big(\log(e+t)\big)^{k/2}\,\tau_\eps^{3/p}\,\mathbf 1_{\{C\tau_\eps\ge1\}}\le
C\,\eps\,(1+\tau_\eps)^{-(p-3)/p}\cdot t^{1+\gamma}\,\mathfrak T(t)\,\big(\log(e+t)\big)^{k/2}\le
C\,\eps\,(1+\tau_\eps)^{-(p-3)/p},
\end{equation}
where we used $\tau_\eps^{3/p}\le C(1+\tau_\eps)^{-(p-3)/p}\tau_\eps$ on $\{C\tau_\eps\ge1\}$,
$\tau_\eps=\eps\,\mathfrak T(t)$, $t^{1+\gamma}\mathfrak T(t)\le
t^{1+\gamma+\kappa_\gamma}/\kappa_\gamma=(\kappa_\gamma t)^{-1}$, and $t^{-1}(\log(e+t))^{k/2}\le
C$.

\emph{The third term} contributes $Ct^{-1}\int\rho_\xi\dd\xi\le C\eps/t$, and $t^{-1}\le
C(1+\tau_\eps)^{-(p-3)/p}$, because $1+\tau_\eps\le1+\mathfrak T(t)\le Ct^{\kappa_\gamma}$ for
$\eps\le1$ and $\kappa_\gamma(p-3)/p<1$.

Since $\|f_0\|_{L^1}=\eps\int\rho_0\ge c\,\eps$, this is~(ii).

\emph{Step 6: proof of~(iii).} We follow Step~5 of the proof of Theorem~\ref{thm:rate2}, with
$\mathfrak T(t)$ in place of $\log(t/T_0)$. \emph{(a)} For every $\xi$ and $t\ge T_0$,
\begin{equation}
\|\mathcal N_\xi(t)-F_0(\cdot,\xi)\|_{L^1(\R^3)}\le C_{\rm
d}\Big(s_\xi(t)+\frac1{T_0}+\eps\,T_0^{\kappa_\gamma}\Big),
\end{equation}
by Lemma~\ref{lem:Wtail}(ii) and by the first item of the proof of Theorem~\ref{thm:main}(iv) and
Remark~\ref{rem:nonunif}, which hold for every $\gamma\in[-3,-2]$. \emph{(b)} The lower
bound~\eqref{eq:R2low}, $I(t,\xi)\ge\frac14d_1\rho_\xi$ whenever $s_\xi(t)\le\sigma:=d_1/(8C_{\rm
d})$, follows as there, for $T_0\ge C_{\rm sh}/d_1$ and $C_{\rm sh}\eps\,T_0^{\kappa_\gamma}\le
d_1$, with $C_{\rm sh}:=16(C_{\rm d}+C_{\rm p})$. The argument rests on (a), (F3) and the
Lipschitz bounds on $p_{\rm sc}$, none of which depends on $\gamma$. \emph{(c)} By
the same theorem and (D3), $\lambda_{\rm
ref}(\xi)\le2\eps\rho_0(\xi)\theta_\infty(\xi)^{\gamma/2}\le2A_F(2\theta_-)^{\gamma/2}\eps\jap\xi^{-p}$,
since $\theta_\infty(\xi)\ge2\theta_-$. Hence $s_\xi(t)=\lambda_{\rm ref}(\xi)
\mathfrak T(t)\le\sigma$ whenever $\jap\xi^p\ge K\tau_\eps$, with
$K:=2A_F(2\theta_-)^{\gamma/2}/\sigma=16A_FC_{\rm d}(2\theta_-)^{\gamma/2}/d_1$, and we may take
$C_{\rm d}\ge(2\theta_-)^{|\gamma|/2}/(16A_F)$, so that $K\ge1$. Moreover $\rho_\xi=\lambda_{\rm
ref}(\xi)\theta_\xi^{|\gamma|/2}\ge\frac12\eps\rho_0(\xi)\theta_\infty(\xi)^{\gamma/2}\theta_-^{|\gamma|/2}\ge\frac12\eps\rho_0(\xi)(2\theta_-/\theta_+)^{|\gamma|/2}\ge
c\,\eps\jap\xi^{-p}$, since $\theta_\infty(\xi)\le\theta_+/2$. With
$R:=\max\{1,(K\tau_\eps)^{1/p}\}$, the integration over $\{|\xi|\ge R\}$ is that of Step~5(c) of
that proof, which yields
\begin{equation}
\|g(t)-g_{\rm sc}\|_{L^1}\ge\frac{d_1}4\,c\,\eps\,2^{-p/2}\,\frac{4\pi}{p-3}\,R^{-(p-3)}\ge
c\,d_1^{\,2-3/p}\,\eps\,(1+\tau_\eps)^{-(p-3)/p} .
\end{equation}
Since $\|f_0\|_{L^1}\le C\eps$, this is~(iii).
\end{proof}

\begin{remark}[On the rate]\label{rem:rate}
The proof of Theorem~\ref{thm:scatter}(i) is soft (a Lipschitz function with integrable square tends
to zero) and gives no rate. Theorems~\ref{thm:rate2} and~\ref{thm:rate3} give one.

\emph{The heuristic.} Formally, the micro part of a relaxed leaf is forced by the leafwise transport
at the rate $t^{-2}$ and damped at the rate $\kappa\lambda_\xi$. For $\gamma<-2$ the deviation that
the leaf has at entry is then damped faster than any power of $t$, since the clock $s_\xi$ grows
like $\lambda_\xi t^{|\gamma|-2}$. The balance of forcing and damping suggests
$d_\xi(t)=O(t^{1+\gamma}/\lambda_\xi)$ for large $t$. Theorem~\ref{thm:rate3}(i) proves this up to a
logarithm.

\emph{The endpoint $\gamma=-2$.} At $\gamma=-2$ the clock is $s_\xi=\lambda_{\rm
ref}(\xi)\log(t/T_0)$, so the deviation at entry decays only like a power $t^{-c\lambda_\xi}$. This
is the rate of Theorem~\ref{thm:rate2}(i). The forced response does not change it. If
$c\lambda_\xi<1$, in particular for small $\eps$, a forcing of size $t^{-2}$ that is damped at the
rate $c\lambda_\xi/t$ from the time $t^b_\xi$ on gives formally a response of size
$(t^b_\xi)^{-1}(t^b_\xi/t)^{c\lambda_\xi}$. It decays like the deviation at entry, with the small
amplitude $1/t^b_\xi$. This is the term with $C_b'$ in Step~2 of its proof. Only for
$c\lambda_\xi>1$ is the forced response of the quasi-static size
$t^{-1}/(c\lambda_\xi-1)$.

The $L^1$ rate of Theorem~\ref{thm:rate3}(ii) does not depend on the exponent $a_\gamma$:
any decay of the leaves in their clock that is faster than $s_\xi^{-(p-3)/p}$ gives it, since it is
determined by the mass of the leaves that have not yet relaxed.
\end{remark}

\appendix
\setlength{\emergencystretch}{8em}

\section{Complements to the near-Maxwellian theory and the homogeneous equation}\label{app:inputs}

This appendix proves the complements to the results of Subsection~\ref{ssec:results} that we need,
for the near-Maxwellian theory and for the homogeneous equation.

\subsection{The near-Maxwellian inputs: two complements and the endpoint
\texorpdfstring{$\gamma=-2$}{gamma=-2}}\label{app:CM}

We use the near-Maxwellian theory of \cite{CM} for the linearized operator $\mathcal L=
\mathcal A+\mathcal B$, where $\mathcal A\varphi=M\chi_R\varphi+\ab[\varphi]:\nabla^2\mu-\cb[\varphi]\mu$ is
the part that is regularizing in the decomposition of \cite{CM} and $\mathcal A_0$ its part without
$M\chi_R$. Two points need a complement, (a) and (b), and (c) treats the endpoint $\gamma=-2$, which
\cite{CM} does not cover.

\emph{(a) The regularizing part on weighted $L^2$ spaces.} In the proof of \cite[Lemma~2.5]{CM} the
contribution of the far part $|z|^{\gamma+2}\mathbf 1_{|z|>1}$ of $a$ is controlled by
$\|\varphi\|_{L^1}$. This is how the bound enters the proof of \cite[Prop.~3.3]{CM}, where $
\mathcal A$ acts on $E_\mu\subset L^1\cap L^2$. On spaces that are not contained in $L^1$ the far
part has to be treated by a weight, since the $L^p$ norm alone does not control it. Indeed, for
$\gamma\in[-3,-2)$ the far part of $a$ lies in $L^{p'}$ only for $p<\frac3{5+\gamma}$, which
excludes $p=2$. For instance, let $\varphi_n(w)=\jap w^{-2}\mathbf 1_{|w|<n}$. Then
$\|\varphi_n\|_{L^2}$ is bounded, but $\ab[\varphi_n](v)=\lambda_nI+O(1)$ on compact sets, with
$\lambda_n=\frac23\int_{|w|<n}|w|^{\gamma+2}\jap w^{-2}\dd w\to\infty$, so that $ \mathcal A_0\varphi_n=\lambda_n(|v|^2-3)\mu+O(1)$
and $\|\mathcal A_0\varphi_n\|_{L^2(B_1)}\to\infty$, which implies that $\mathcal A_0$ is not
bounded from $L^2$ to $L^2(\mu^{-\vartheta})$. The reason is that the far part of $a$ decays like
$|z|^{\gamma+2}$, which is not square integrable for $\gamma\ge-\frac72$. We use the bound on the
spaces $L^2(\jap v^s)$, $s>\gamma+\frac72$, which for $s\le\frac32$ are not contained in $L^1$. The
following version suffices for our uses of \cite{CM}.

\begin{lem}[$\mathcal A$ on weighted spaces]\label{lem:RVA}
Let $\gamma\in[-3,-2]$, $\vartheta\in(0,1)$ and $s>\gamma+\frac72$. Then:
\begin{enumerate}[label=(\roman*)]
\item $\mathcal A$ is bounded from $L^2(\jap v^s)$ to $L^2(\mu^{-\vartheta})$;
\item for $\psi\in\{1,v_i,|v|^2\}$, $|\int\psi\,\mathcal A\varphi|\le C\|\varphi\|_{L^2(\jap v^s)}$;
\item for every weight $m$ with $m^2\mu^{\vartheta/2}$ bounded,
\begin{equation}\label{eq:RVA-iii}
|\int m^2\varphi\,\mathcal A\varphi|\le C\|\varphi\|_{L^2}\|\varphi\|_{L^2(\jap v^s)} .
\end{equation}

\end{enumerate}
\end{lem}

\begin{proof}
(i) We split $\ab[\varphi](v)$ at $|v-w|=1$. The near part is bounded by $\||z|^{\gamma+2} \mathbf 1_{|z|<1}\|_{L^2}\|\varphi\|_{L^2}$,
which is finite since $2(\gamma+2)>-3$. For the far part, by Peetre's inequality we have
$|v-w|^{\gamma+2}\le C\jap v^{|\gamma+2|}\jap w^{\gamma+2}$, and hence $|\ab[\varphi](v)|\le C\jap
v^{|\gamma+2|}\|\varphi\|_{L^2(\jap v^s)}$, since $\jap w^{\gamma+2-s}\in L^2$ for
$s>\gamma+\frac72$. Multiplying by $\nabla^2\mu$, we obtain an element of $L^2(\mu^{-\vartheta})$.
For $\cb$, if $\gamma=-3$, then $\cb[\varphi]=-8\pi\varphi$. If $\gamma>-3$, the near part of $c$ is
integrable and the far part is square integrable, since $2\gamma<-3$. (ii) Integrating by parts
twice, we have $\int|v|^2\mathcal A_0\varphi=\int\varphi\Psi$ with
$\Psi(w)=\int\mu(v)[4|v-w|^{\gamma+2}-8v\cdot(v-w)|v-w|^\gamma]\dd v$, and $|\Psi(w)|\le C\jap
w^{\gamma+2}$. The other moments are similar. (iii) follows from (i) and Cauchy--Schwarz.
\end{proof}

\emph{(b) The moment term in \cite[Thm.~3.5]{CM}.} Step~1 of the proof of \cite[Thm.~3.5]{CM} uses
the factorization $e^{t\mathcal L}\Pi=\Pi e^{t\mathcal B}+(e^{\cdot\mathcal L}\Pi)*( \mathcal Ae^{\cdot\mathcal B})$,
and its first term is then bounded through $e^{t\mathcal B}\Pi$. The two differ by a moment term,
since $\Pi e^{t\mathcal B}=e^{t\mathcal B}-\Pi_0e^{t\mathcal B}$. The projection $\Pi_0$ is bounded
on $X_{k_0}$ only for $k_0>\frac72$, because of the moment $\int|v|^2f$, while \cite[Prop.~3.6]{CM}
and \cite[(3.15)]{CM} bound $e^{t\mathcal L}\Pi$ in $X_0=L^2$ \cite[(3.11)]{CM}, through
$L^2(m_0)\subset L^2$ with $\jap v^{(\gamma+3)/2}\prec m_0$, so that $k_0<\frac72$. Step~2 of the
same proof (the bound from $H^{-1}_*$) proceeds in the same way. To follow the argument we therefore
estimate the moment term $\Pi_0e^{t\mathcal B}\Pi$ separately. We make no claim about the thresholds
stated in \cite{CM}, which may be reachable by another treatment of this term, and we only record
the thresholds that our argument yields. Lemma~\ref{lem:CMmom} below estimates the moment term, and
Proposition~\ref{prop:CMthr} states the resulting thresholds. They are below the level $k=7$ that we
use.

\emph{Notation and inputs.} Let $X_k=L^2(\jap v^k)$, $Y_k=H^1_*(\jap v^k)$ and $Z_k$ be as in
Appendix~\ref{app:band}, let $\Pi_0$ be the projection onto $\ker \mathcal L=\mathrm{span}\{\mu,v_i\mu,|v|^2\mu\}$
along the moment-free functions, and $\Pi=I-\Pi_0$. Note that $\Pi_0$ and $\Pi$ are bounded on $X_k$
for $k>\frac72$, on $Z_k$ for $k>\frac72+\frac\gamma2$, and on $L^2(\mu^{-\vartheta})$, and they
commute with $e^{t\mathcal L}$. For each statement below, the constants $M,R$ of $\mathcal A=\mathcal A_0+M\chi_R$
are chosen large enough for the finitely many weights that enter it; $e^{t\mathcal L}$, hence
parts~(i) and~(iii) of that proposition, do not depend on this choice, while $K^{\mathcal B}_k$,
$\eta_k$ and $K_k$ in part~(ii) do. From (E1) we use the statements (S0)--(S3) of
Inputs~\ref{hyp:E}. \cite[Lemma~2.5]{CM} enters (S0)--(S3) only in the proof of
\cite[Prop.~3.3]{CM}, where $\mathcal A$ acts on $E_\mu\subset L^2(\jap v^s)$. There the argument of
\cite{CM} applies as it stands, since $E_\mu\subset L^1$, and Lemma~\ref{lem:RVA}(i) yields the same
bound. Lemmas~2.2, 2.3, 2.4, 3.1 and~3.2 of \cite{CM} do not use it. In \cite{CM} it is also used in
Lemma~3.4, Theorem~3.5 and Proposition~3.6, which we do not quote.

\begin{lem}[The moment term]\label{lem:CMmom}
Let $\gamma\in[-3,-2]$, $k_1>\frac72$, $h\in X_{k_1}$ and $g=\Pi h$. Then
\begin{equation}
\Pi_0e^{t\mathcal B}g=\int_t^\infty\Pi_0\mathcal Ae^{s\mathcal B}g\dd s .
\end{equation}
Hence, for $\gamma+\frac72<k_*<k_1-|\gamma|$, every $b\in(1,\frac{k_1-k_*}{|\gamma|})$ and every
$k_0\ge0$, $\|\Pi_0e^{t\mathcal B}g\|_{X_{k_0}}\le C\jap t^{1-b}\|h\|_{X_{k_1}}$.
\end{lem}

\begin{proof}
Since the moments of $\mathcal Lf$ vanish, we have $\Pi_0\mathcal L=0$ and $\Pi_0\mathcal B=-\Pi_0\mathcal A$.
Hence $\Pi_0e^{t\mathcal B}g=\Pi_0g-\int_0^t\Pi_0\mathcal Ae^{s\mathcal B}g\dd s$, first for $g$ in
the domain of $\mathcal B$ and then, since both sides are continuous in $g\in X_{k_1}$ by (S1) and
Lemma~\ref{lem:RVA}, for all $g$. Moreover, $\Pi_0g=0$. By (S1) with a weight $k'\in(\frac72,k_1)$,
we have $e^{t\mathcal B}g\to0$ in $X_{k'}$, so $\Pi_0e^{t\mathcal B}g\to0$, and the identity
follows. The range of $\Pi_0$ is spanned by $\mu,v_i\mu,|v|^2\mu$, which lie in every $X_{k_0}$, and
its coefficients are linear combinations of $\int\psi\,\mathcal Ae^{s\mathcal B}g$,
$\psi\in\{1,v_i,|v|^2\}$. By Lemma~\ref{lem:RVA}(ii) with $s=k_*$ and (S1) with $k=k_*$, which is
allowed since $k_*>\gamma+\frac72>\frac{\gamma+3}2$, they are bounded by $C\|e^{s \mathcal B}g\|_{X_{k_*}}\le
C\jap s^{-b}\|g\|_{X_{k_1}}$, and $\|g\|_{X_{k_1}}\le C\|h\|_{X_{k_1}}$. Integrating over
$[t,\infty)$ and using $b>1$, we obtain the claim.
\end{proof}

\begin{prop}[Decay of $e^{t\mathcal L}\Pi$, and the weight thresholds]\label{prop:CMthr}
Let $\gamma\in[-3,-2]$ and assume (E1) (at $\gamma=-2$, in the form of Proposition~\ref{prop:CMm2}).
\begin{enumerate}[label=(\roman*)]
\item Let $k_1>\frac72$, $\gamma+\frac72<k_*<k_1-|\gamma|$ and $0\le k_0\le k_*$. For every $a\le1$
with $a<\frac{k_1-k_*}{|\gamma|}-1$, $\|e^{t\mathcal L}\Pi\|_{X_{k_1}\to X_{k_0}}\le C\jap t^{-a}$.
Let $\frac72-\frac12|\gamma|<k''<k_1-|\gamma|$ and $0\le k_0\le k''-\frac12|\gamma|$. For every
$a\le1$ with $a<\frac{k_1-k''}{|\gamma|}-1$,
\begin{equation}\label{eq:CMthr-i}
\|e^{t\mathcal L}\Pi\|_{Z_{k_1}\to X_{k_0}}\le C(t\wedge1)^{-1/2}\jap t^{-a} .
\end{equation}

\item \emph{(F1).} Let $k>\frac72+\frac12|\gamma|$. Then
\begin{equation}\label{eq:CMthr-F1}
\int_0^\infty\|e^{\tau\mathcal L}\Pi u\|^2_{L^2}\dd\tau\le C_\natural(k)\|u\|^2_{X_k}
\end{equation}
for $u\in X_k$. There is $\eta_k>0$ such that
\begin{equation}\label{eq:hybrid}
\nn h_k^2:=\eta_k\|h\|^2_{X_k}+\int_0^\infty\|e^{\tau\mathcal L}h\|^2_{L^2}\dd\tau
\end{equation}
is a Hilbert norm on $X^0_k$ with
\begin{equation}\label{eq:hybrid-eq}
\eta_k^{1/2}\|h\|_{X_k}\le\nn h_k\le(\eta_k+C_\natural(k))^{1/2}\|h\|_{X_k},
\end{equation}
and, with $\ip\cdot\cdot_k$ the associated inner product, $\ip{\mathcal Lh}h_k\le-K_k\|h\|^2_{Y_k}$
for $h\in Y^0_k$ in the domain of $\mathcal L$, with $K_k:=\eta_kK^{\mathcal B}_k/2$.
\item \emph{(F2).} Let $k>\frac72+|\gamma|$. Then
\begin{equation}\label{eq:Jk}
J_k:=\int_0^\infty\|e^{t\mathcal L}\Pi\|_{Y_k\to L^2}\|e^{t\mathcal L}\Pi\|_{Z_k\to L^2}\dd
t<\infty .
\end{equation}

\end{enumerate}
The thresholds in (ii) and (iii) are $4.5$ and $5.5$ at $\gamma=-2$, and $5$ and $6.5$ at
$\gamma=-3$; the level $k=7$ satisfies both for every $\gamma\in[-3,-2]$.
\end{prop}

\begin{proof}
(i) Let $h\in X_{k_1}$ and $g=\Pi h$. By Duhamel's formula for $\mathcal L=\mathcal B+\mathcal A$
and $\Pi e^{t\mathcal L}=e^{t\mathcal L}\Pi$, we have
\begin{equation}
e^{t\mathcal L}g=e^{t\mathcal B}g-\Pi_0e^{t\mathcal B}g+\int_0^te^{(t-s)\mathcal L}\Pi
\mathcal Ae^{s\mathcal B}g\dd s .
\end{equation}
Fix $b\in(\max\{1,a+1\},\frac{k_1-k_*}{|\gamma|})$. By (S1), we have $\|e^{t \mathcal B}g\|_{X_{k_0}}\le\|e^{t\mathcal B}g\|_{X_{k_*}}\le
C\jap t^{-b}\|h\|_{X_{k_1}}$. By Lemma~\ref{lem:CMmom}, the second term is $\le C\jap
t^{1-b}\|h\|_{X_{k_1}}$. For the third, we take $\vartheta$ as in (S3). Since $\jap v^{k_0}\le
C\mu^{-1/2}$, it follows from (S3), Lemma~\ref{lem:RVA}(i) with $s=k_*$ and (S1) that its integrand
in $X_{k_0}$ is bounded by $C\Theta_{\rm s}(t-s)\jap s^{-b}\|h\|_{X_{k_1}}$. Splitting at $s=t/2$,
we get $\int_0^t\Theta_{\rm s}(t-s)\jap s^{-b}\dd s\le2^b\jap t^{-b}\|\Theta_{\rm
s}\|_{L^1}+\sup_{\tau\ge t/2}\Theta_{\rm s}(\tau)\int_0^{t/2}\jap s^{-b}\dd s\le C\jap t^{-b}+C\jap
t^{-1}$, by (S3) and the bound $\Theta_{\rm s}(\tau)\le C\jap\tau^{-1}$. Since $a\le1$ and $a<b-1$,
all three terms are $\le C\jap t^{-a}\|h\|_{X_{k_1}}$. For $h\in Z_{k_1}$ the same argument applies
with (S2) in place of (S1). Indeed, (S2) with $k=k''$ bounds $e^{s\mathcal B}g$ in
$X_{k''-|\gamma|/2}$, where $k''-\frac12|\gamma|>\gamma+\frac72$, so that the same bound applies
with $s=k''-\frac12|\gamma|$, and $X_{k''-|\gamma|/2}\subset X_{k_0}$. In the limit used for the
moment term, we use (S2) with $k\in(\frac72+\frac12|\gamma|,k_1)$. Each bound acquires the factor
$(s\wedge1)^{-1/2}$, which is integrable at $s=0$, and the result acquires the factor
$(t\wedge1)^{-1/2}$. Since the target weight of (S2) is lower by $\frac12|\gamma|$ than that of
(S1), $k''$ replaces $k_*$.

(ii) We choose $k_*\in(\gamma+\frac72,k-\frac32|\gamma|)$, which is possible since
$k-\frac32|\gamma|>\frac72-|\gamma|$, and $k_0=0$. Then $\frac{k-k_*}{|\gamma|}-1>\frac12$, and by
(i) we have $\|e^{\tau\mathcal L}\Pi u\|_{L^2}\le C\jap\tau^{-a}\|u\|_{X_k}$ with some
$a\in(\frac12,1]$. This is the first claim, and the bounds on $\nn\cdot_k$ follow, since $\Pi h=h$
on $X^0_k$. For $h\in Y^0_k$ in the domain of $\mathcal L$, we have $\int_0^\infty\langle
e^{\tau\mathcal L} \mathcal Lh,e^{\tau\mathcal L}h\rangle\dd\tau=\frac12\int_0^\infty\frac
d{d\tau}\|e^{\tau \mathcal L}h\|^2_{L^2}\dd\tau=-\frac12\|h\|^2_{L^2}$, since $e^{\tau\mathcal L}h\to0$
in $L^2$. Hence
\begin{equation}
\ip{\mathcal Lh}h_k=\eta_k\langle\mathcal Bh,h\rangle_{X_k}+\eta_k\langle
\mathcal Ah,h\rangle_{X_k}-\tfrac12\|h\|^2_{L^2}.
\end{equation}
By (S0) the first term is $\le-\eta_kK^{\mathcal B}_k\|h\|^2_{Y_k}$. By Lemma~\ref{lem:RVA}(iii)
with $m=\jap v^k$ and $s\in(\gamma+\frac72,k-\frac12|\gamma|]$, and $\|h\|_{X_{k-|\gamma|/2}}\le
C\|h\|_{Y_k}$, the second is $\le\eta_kC\|h\|_{L^2}\|h\|_{Y_k}\le\eta_k(\frac12K^{ \mathcal B}_k\|h\|^2_{Y_k}+C'\|h\|^2_{L^2})$.
We choose $\eta_k\le1/(2C')$. If moreover $k>\frac72+|\gamma|$, the inequality extends to all $h\in
X^0_k\cap Y^0_k$, with $\ip{\mathcal Lh}h_k$ defined through the extension of $\ip\cdot\cdot_k$ to
$Z^0_k\times Y^0_k$ given by~(iii). Indeed, by \cite[(4.3)]{CM} applied to $Q(\mu,h)+Q(h,\mu)$,
$\mathcal L$ is bounded from $X_k\cap Y_k$, with the norm $\|\cdot\|_{X_k}+\|\cdot\|_{Y_k}$, to
$Z_k$. Both sides are continuous on $X^0_k\cap Y^0_k$, and $\Pi C^\infty_c$ is dense in $X^0_k\cap
Y^0_k$ and contained in the domain of $\mathcal L$. (The operator $\mathcal L$ is not bounded from
$Y_k$ to $Z_k$, since the drift $-\bb[\mu]\cdot\nabla=\ell_1v\cdot\nabla$ of $Q(\mu,\cdot)$ has size
$\jap v^{\gamma+1}$, one power of $\jap v$ more than $Y_k\times Y_k$ controls. This is the term
$\|g\|_{L^2(m)}$ in \cite[(4.1)]{CM}.)

(iii) Since $Y_k\subset X_{k-|\gamma|/2}$, applying (i) with $k_1=k-\frac12|\gamma|$ and $k_0=0$, we
get $\|e^{t\mathcal L}\Pi\|_{Y_k\to L^2}\le C\jap t^{-a_1}$ for every
$a_1<\frac{k-|\gamma|/2-k_*}{|\gamma|}-1$, $k_*\in(\gamma+\frac72,k-\frac32|\gamma|)$. Applying (i)
from $Z_k$, with $k_0=0$, we get $\|e^{t\mathcal L}\Pi\|_{Z_k\to L^2}\le C(t\wedge1)^{-1/2}\jap
t^{-a_2}$ for every $a_2<\frac{k-k''}{|\gamma|}-1$, $k''\in(\frac72-\frac12|\gamma|,k-|\gamma|)$.
Both upper bounds are positive, so $a_1,a_2$ can be chosen in $(0,1]$ with $a_1+a_2>1$ as soon as
the upper bounds have sum $>1$, and then the product is integrable. As
$k_*\downarrow\frac72-|\gamma|$ and $k''\downarrow\frac72-\frac12|\gamma|$, both upper bounds tend
to $\frac{k-7/2+|\gamma|/2}{|\gamma|}-1$, and their sum exceeds $1$ when $k>\frac72+|\gamma|$.
\end{proof}

\emph{(c) The endpoint $\gamma=-2$.} At $\gamma=-2$ the kernels are $a(z)=I-\hat z\otimes\hat z$,
$b(z)=-2z|z|^{-2}$ and $c(z)=-2|z|^{-2}$, the formulas of \cite[(2.1)]{CM} evaluated at $\gamma=-2$.
They remain valid there, since they come from differentiating $|z|^{\gamma+2}(I-\hat z\otimes\hat
z)$. Recall $E_\mu=L^2(\mu^{-1/2})$, with $\|f\|^2_{E_\mu}=\int f^2\mu^{-1}$.

\subsubsection{The polynomial-weight estimates of \texorpdfstring{\cite{CM}}{[CM]} for
\texorpdfstring{$\gamma\in[-3,-2]$}{gamma in [-3,-2]}}\label{ssec:CMpoly}

The estimates (S0)--(S2) and \cite[Lemma~4.3]{CM} are proved in \cite{CM} for $\gamma\in[-3,-2)$. We
write out their proofs in polynomial weights for every $\gamma\in[-3,-2]$. At $\gamma=-2$ this
replaces the corresponding items of Inputs~\ref{hyp:E}. The arguments are those of
\cite[\S\S2.3--2.4, 3.1--3.2, 4]{CM}. The only inputs are the properties of $\ab[\mu]$ in (E1)(b),
which hold for $\gamma\ge-3$ (\cite[Lemmas~2, 3]{Guo02} and Lemma~\ref{lem:ell1}; part~(iv) of
(E1)(b), which rests on \cite[Lemma~2.4(d)]{CTW}, is not among them), and the following elementary
asymptotics: for $\beta\in(-3,0]$ and $J_\beta:=|\cdot|^\beta*\mu$,
\begin{equation}\label{eq:Jbeta}
J_\beta(v)\le C\jap v^\beta,\qquad J_\beta(v)=|v|^\beta\big(1+O(|v|^{-2})\big)\quad(|v|\ge1),
\end{equation}
which follow by the splitting at $|w|=|v|/2$ of its proof. We have
$\mathrm{tr}\,\ab[\mu]=2J_{\gamma+2}$ and $\cb[\mu]=-2(\gamma+3)J_\gamma$ for $\gamma>-3$, and
$\cb[\mu]=-8\pi\mu$ for $\gamma=-3$. Moreover, $\bb[\mu]=-\ell_1v$, $\nabla\cdot\ab[\mu]=\bb[\mu]$
and $\nabla\cdot\bb[\mu]=\cb[\mu]$. By (E1)(b)(ii) and $\bb=\nabla\cdot\ab$, we have
$|\partial^\beta\ab[\mu]|\le C_\beta\jap v^{\gamma+2-|\beta|}$ and $|\partial^\beta\bb[\mu]|\le
C_\beta\jap v^{\gamma+1-|\beta|}$. We write $\ab,\bb,\cb$ for the coefficients of $\mu$, and we put
$\tn\psi:=P_v\nabla\psi+\jap v(I-P_v)\nabla\psi$, so that $\ab\xi\cdot\xi\ge K\jap
v^\gamma|\tilde\xi|^2$ by (E1)(b)(ii), and $ \mathcal B=Q(\mu,\cdot)-M\chi_R=\ab:\nabla^2-\cb-M\chi_R$.

\begin{lem}[Weighted identities and dissipativity in polynomial weights]\label{lem:CMpoly}
Let $\gamma\in[-3,-2]$, $m=\jap v^k$, $\omega=\jap v^\alpha$, $\Phi=\log m$, $\Psi=\log\omega$, and
\begin{equation}
\tilde\zeta_m:=\ab\nabla\Phi\cdot\nabla\Phi+\bb\cdot\nabla\Phi-\tfrac12\cb,\qquad\zeta_{m,\omega}:=-2\ab\nabla\Phi\cdot\nabla\Psi+\ab\nabla\Psi\cdot\nabla\Psi+\ab:\frac{\nabla^2\omega}\omega
.
\end{equation}
\begin{enumerate}[label=(\roman*)]
\item $\limsup_{|v|\to\infty}\jap v^{-\gamma}\tilde\zeta_m\le(\gamma+3)-2k$ and
\begin{equation}\label{eq:CMpoly-i}
\limsup_{|v|\to\infty}\jap v^{-\gamma}(\tilde\zeta_m+\zeta_{m,\omega})\le(\gamma+3)-2k+2\alpha .
\end{equation}

\item \emph{(S0).} For $k>\frac{\gamma+3}2$ there are $M,R$ and $c>0$ with
\begin{equation}\label{eq:CMpoly-ii}
\langle\mathcal Bf,f\rangle_{X_k}\le-c\big(\|\jap v^{\gamma/2}mf\|^2_{L^2}+\|\jap
v^{\gamma/2}\tn(mf)\|^2_{L^2}\big)
\end{equation}
for smooth rapidly decaying $f$.
\item \emph{(S1), and \cite[(3.4)]{CM}.} For $\frac{\gamma+3}2<k<k_1$ and $0\le
b<\frac{k_1-k}{|\gamma|}$, $\|e^{t\mathcal B}\|_{X_{k_1}\to X_k}\le C\jap t^{-b}$. Let $
\mathcal B^*_m$ be the formal adjoint of $\mathcal B_m:=m\mathcal B(m^{-1}\cdot)$. If $m=\jap v^{k_2}$ with
$k_2\ge k_1$, $\omega_1:=m/\jap v^k$ and $\omega_0:=m/\jap v^{k_1}$, then also
\begin{equation}\label{eq:CMpoly-adj}
\|e^{t\mathcal B^*_m}\|_{L^2(\omega_1)\to L^2(\omega_0)}\le C\jap t^{-b} .
\end{equation}

\item \emph{\cite[(2.19), (2.20)]{CM}.} If $\alpha\ge0$ and $k>\frac{\gamma+3}2+\alpha$, there are
$M,R,\eta,c>0$ such that, for smooth rapidly decaying $\varphi$,
\begin{equation}
\langle\mathcal B^*_m\varphi,\varphi\rangle_{L^2(\omega)}\le-c\big(\|\jap
v^{\gamma/2}\varphi\|^2_{L^2(\omega)}+\|\jap v^{\gamma/2}\tn\varphi\|^2_{L^2(\omega)}\big),
\end{equation}
and, with
\begin{equation}\label{eq:H1tilde}
\|\varphi\|^2_{\tilde H^1(\omega)}:=\|\varphi\|^2_{L^2(\omega)}+\eta\|\nabla\varphi\|^2_{L^2(\omega)}
\end{equation}
and the associated scalar product,
\begin{equation}
\langle\mathcal B^*_m\varphi,\varphi\rangle_{\tilde H^1(\omega)}\le-c\big(\|\jap
v^{\gamma/2}\varphi\|^2_{L^2(\omega)}+\|\jap v^{\gamma/2}\tn\varphi\|^2_{L^2(\omega)}+\eta\|\jap
v^{\gamma/2}\tn\nabla\varphi\|^2_{L^2(\omega)}\big).
\end{equation}
\item \emph{(S2).} For $\frac32<k<k_1$ and $0\le b<\frac{k_1-k}{|\gamma|}$,
\begin{equation}\label{eq:CMpoly-B}
\|e^{t\mathcal B}\|_{Z_{k_1}\to X_{k-|\gamma|/2}}\le C(t\wedge1)^{-1/2}\jap t^{-b} .
\end{equation}

\end{enumerate}
The constants depend only on $\gamma$, $k$, $k_1$, $k_2$, $\alpha$ and $b$.
\end{lem}

\begin{proof}
(i) We have $\nabla\Phi=kv\jap v^{-2}$, an eigenvector of $\ab$ with eigenvalue $\ell_1$, and
$\bb=-\ell_1v$. Thus $\ab\nabla\Phi\cdot\nabla\Phi=k^2\ell_1|v|^2\jap v^{-4}$ and
$\bb\cdot\nabla\Phi=-k\ell_1|v|^2\jap v^{-2}$. Likewise
$\ab\nabla\Phi\cdot\nabla\Psi=k\alpha\ell_1|v|^2\jap v^{-4}$,
$\ab\nabla\Psi\cdot\nabla\Psi=\alpha^2\ell_1|v|^2\jap v^{-4}$, and, since
$\nabla^2\omega/\omega=\alpha\jap v^{-2}I+\alpha(\alpha-2)\jap v^{-4}v\otimes v$,
$\ab:\nabla^2\omega/\omega=2\alpha J_{\gamma+2}\jap v^{-2}+\alpha(\alpha-2)\ell_1|v|^2\jap v^{-4}$.
By Lemma~\ref{lem:ell1}, we have $\jap v^{-\gamma}\ell_1|v|^2\jap v^{-2}\to2$. By~\eqref{eq:Jbeta},
we have $\jap v^{-\gamma}J_{\gamma+2}\jap v^{-2}\to1$ and $-\frac12\jap v^{-\gamma}\cb\to\gamma+3$
(for $\gamma=-3$, $-\frac12\cb=4\pi\mu=o(\jap v^\gamma)$, and $\gamma+3=0$). The terms with
$\ell_1|v|^2\jap v^{-4}$ are $O(\jap v^{\gamma-2})$. Hence we deduce that $\jap
v^{-\gamma}\tilde\zeta_m\to-2k+(\gamma+3)$ and $\jap v^{-\gamma}\zeta_{m,\omega}\to2\alpha$.

(ii) Let $g=mf$. By two integrations by parts (\cite[proof of Lemma~2.3]{CM}, based on
$\nabla\cdot\ab=\bb$ and $\nabla\cdot\bb=\cb$ alone), we have
\begin{equation}
\int(\mathcal Bf)fm^2=-\int\ab\nabla g\cdot\nabla g+\int(\tilde\zeta_m-M\chi_R)g^2 .
\end{equation}
By~(i), since $k>\frac{\gamma+3}2$, there are $R$ and $c_1>0$ with $\tilde\zeta_m\le-2c_1\jap
v^\gamma$ for $|v|\ge R$. Taking $M\ge\sup_{B_{2R}}|\tilde\zeta_m|+2c_1$, we get
$\tilde\zeta_m-M\chi_R\le-c_1\jap v^\gamma$ everywhere. Then~(ii) follows from the coercivity
$\ab\xi\cdot\xi\ge K\jap v^\gamma|\tilde\xi|^2$.

(iv) By the definition, we have $\mathcal B_mg=\ab:\nabla^2g-2\ab\nabla\Phi\cdot\nabla
g+\big(2\ab\nabla\Phi\cdot\nabla\Phi-\ab:\nabla^2m/m-\cb-M\chi_R\big)g$. Integrating by parts with
$\nabla\cdot\ab=\bb$, $\nabla\cdot\bb=\cb$ and
$\nabla\cdot(\ab\nabla\Phi)=\bb\cdot\nabla\Phi+\ab:\nabla^2m/m-\ab\nabla\Phi\cdot\nabla\Phi$, we
obtain
\begin{equation}
\mathcal B^*_m\varphi=\ab:\nabla^2\varphi+2(\ab\nabla\Phi+\bb)\cdot\nabla\varphi+\big(\ab:\nabla^2m/m+2\bb\cdot\nabla\Phi-M\chi_R\big)\varphi
.
\end{equation}
We pair with $\varphi\omega^2$. The first-order term becomes, after one integration by parts,
$-\int\big(\nabla\cdot(\ab\nabla\Phi+\bb)\big)\varphi^2\omega^2-2\int(\ab\nabla\Phi+\bb)\cdot\nabla\Psi\,\varphi^2\omega^2$,
and the second-order term becomes
$-\int\ab\nabla\varphi\cdot\nabla\varphi\,\omega^2+\frac12\int\partial_{ij}(\ab_{ij}\omega^2)\varphi^2$,
where
$\frac12\partial_{ij}(\ab_{ij}\omega^2)=\omega^2\big(\frac12\cb+2\bb\cdot\nabla\Psi+\ab:\nabla^2\omega/\omega+\ab\nabla\Psi\cdot\nabla\Psi\big)$.
Adding up, we have
\begin{equation}
\int(\mathcal B^*_m\varphi)\varphi\,\omega^2=-\int\ab\nabla\varphi\cdot\nabla\varphi\,\omega^2+\int\big(\tilde\zeta_m+\zeta_{m,\omega}-M\chi_R\big)\varphi^2\omega^2
.
\end{equation}
By~(i) and $k>\frac{\gamma+3}2+\alpha$, $M,R$ can be chosen as in~(ii) with
$\tilde\zeta_m+\zeta_{m,\omega}-M\chi_R\le-c_1\jap v^\gamma$. Together with the coercivity, this
implies the first inequality. For the second, let $e$ be a unit vector and $\partial=e\cdot\nabla$.
Then we have
\begin{equation}
\partial(\mathcal B^*_m\varphi)=\mathcal B^*_m(\partial\varphi)+\partial\big(\ab:\tfrac{\nabla^2m}m+2\bb\cdot\nabla\Phi\big)\varphi+2\,\partial(\ab\nabla\Phi+\bb)\cdot\nabla\varphi+(\partial\ab):\nabla^2\varphi
\end{equation}
Here $\chi_R$ is radial and smooth, and its derivative is included in the second term after
enlarging $C$. By the derivative bounds above and $|\nabla^j\Phi|\le C\jap v^{-j}$, the coefficient
of $\varphi$ is $O(\jap v^{\gamma-1})$ and that of $\nabla\varphi$ is $O(\jap v^\gamma)$. In the
pairing of the last term with $\partial\varphi\,\omega^2$, we integrate by parts in the index $j$
and get
$\int(\partial\ab_{ij})\partial_{ij}\varphi\,\partial\varphi\,\omega^2=-\int(\partial\bb_i)\partial_i\varphi\,\partial\varphi\,\omega^2-\int(\partial\ab_{ij})\partial_i\varphi\,\partial_j\partial\varphi\,\omega^2-\int(\partial\ab_{ij})\partial_i\varphi\,\partial\varphi\,\partial_j\omega^2$,
and, by the symmetry of $\partial\ab$,
$(\partial\ab_{ij})\partial_i\varphi\,\partial(\partial_j\varphi)=\frac12(\partial\ab_{ij})\partial(\partial_i\varphi\partial_j\varphi)$,
which one more integration by parts in $\partial$ turns into terms with $\partial^2\ab$ and
$\partial\ab\,\partial\omega^2$. All these are $O(\jap v^\gamma)|\nabla\varphi|^2\omega^2$, since
$|\nabla\omega|\le C\jap v^{-1}\omega$. Applying the first inequality to $\partial\varphi$ and
summing over $e$ in an orthonormal basis, we obtain
\begin{equation}
\langle\nabla\mathcal B^*_m\varphi,\nabla\varphi\rangle_{L^2(\omega)}\le-c\|\jap
v^{\gamma/2}\tn\nabla\varphi\|^2_{L^2(\omega)}+C\big(\|\jap
v^{\gamma/2}\varphi\|^2_{L^2(\omega)}+\|\jap v^{\gamma/2}\nabla\varphi\|^2_{L^2(\omega)}\big).
\end{equation}
Since $|\nabla\varphi|\le|\tn\varphi|$, we obtain the second inequality by adding $\eta$ times this
to the first one, with $\eta C\le c/2$.

(iii) Let $f(t)=e^{t\mathcal B}f_0$. By~(ii) at the weight $\jap v^{k_1}$, we have
$\|f(t)\|_{X_{k_1}}\le\|f_0\|_{X_{k_1}}$, and at the weight $\jap v^k$, we have $\frac
d{dt}\|f\|^2_{X_k}\le-2c\|\jap v^{\gamma/2}f\|^2_{X_k}$. For $R\ge1$, splitting at $|v|=R$, we get
$\|f\|^2_{X_k}\le\jap R^{|\gamma|}\|\jap v^{\gamma/2}f\|^2_{X_k}+\jap
R^{-2(k_1-k)}\|f\|^2_{X_{k_1}}$. Therefore $\frac d{dt}\|f\|^2_{X_k}\le-2c\jap
R^\gamma\|f\|^2_{X_k}+2c\jap R^{\gamma-2(k_1-k)}\|f_0\|^2_{X_{k_1}}$, which implies that
\begin{equation}
\|f(t)\|^2_{X_k}\le\big(e^{-2c\jap R^\gamma t}+\jap R^{-2(k_1-k)}\big)\|f_0\|^2_{X_{k_1}} .
\end{equation}
For $t\ge t_0$, where $t_0\ge e^2$ and $ct_0\ge\log^2t_0$ (so that $ct\ge\log^2t$ for $t\ge t_0$,
while for $t\le t_0$ we use $\|f(t)\|_{X_k}\le\|f_0\|_{X_{k_1}}$), we take the splitting radius $R$
(not the radius of $\chi_R$) with $\jap R^{|\gamma|}=ct/\log^2t$. Then $e^{-2c\jap R^\gamma
t}=e^{-2\log^2t}$ and $\jap R^{-2(k_1-k)}=(\log^2t/(ct))^{2(k_1-k)/|\gamma|}$, which is $\le C\jap
t^{-2b}$ for $b<(k_1-k)/|\gamma|$.

The estimate for $e^{t\mathcal B^*_m}$ follows by the same argument with the first inequality
of~(iv), at the weights $\omega_0=\jap v^{k_2-k_1}$ and $\omega_1=\jap v^{k_2-k}$, which satisfy
$k_2>\frac{\gamma+3}2+(k_2-k)$, since $k>\frac{\gamma+3}2$, and $\omega_1/\omega_0=\jap v^{k_1-k}$.

\emph{Justification along the semigroups.} These are a priori estimates for smooth rapidly decaying
data. They hold along the semigroups for the following reason. Let $V_k:=X_k\cap Y_k$, with the norm
$\|\cdot\|_{X_k}+\|\cdot\|_{Y_k}$, and write $F=mf$, $G=mg$ with $m=\jap v^k$. By one integration by
parts, we have
\begin{equation}\label{eq:Bform}
-\langle\mathcal Bf,g\rangle_{X_k}=\int\ab(\nabla F-F\nabla\Phi)\cdot(\nabla
G+G\nabla\Phi)+\int\bb\cdot(\nabla F-F\nabla\Phi)\,G+\int(\cb+M\chi_R)FG .
\end{equation}
Since $\ab=\ell_1P_v+\ell_2(I-P_v)$ with $\ell_1\le C\jap v^\gamma$ and $\ell_2\le C\jap
v^{\gamma+2}$, we have $|\ab\xi\cdot\zeta|\le C\jap v^\gamma|\tilde\xi||\tilde\zeta|$. With
$|\nabla\Phi|\le k\jap v^{-1}$, $|\bb|\le C\jap v^{\gamma+1}$, $|\cb|\le C\jap v^\gamma$
($\cb=-8\pi\mu$ for $\gamma=-3$) and $\jap v^{\gamma/2+1}\le1$, the form \eqref{eq:Bform} is bounded
on $V_k\times V_k$ (the drift term by $C\|f\|_{Y_k}\|g\|_{X_k}$). By~(ii) and the density of
$C^\infty_c$ in $V_k$, we deduce that $-\langle\mathcal Bf,f\rangle_{X_k}\ge c\|f\|^2_{Y_k}$ on
$V_k$.

By J.-L.~Lions' theorem \cite[Ch.~XVIII, \S3, Thms.~1--2]{DL}, for $f_0\in X_k$ there is a unique
$f\in L^2(0,T;V_k)\cap C([0,T];X_k)$ with $\partial_tf\in L^2(0,T;V_k')$ and $\partial_tf=
\mathcal Bf$, and $\frac d{dt}\|f\|^2_{X_k}=2\langle\mathcal Bf,f\rangle_{X_k}$ by the Lions--Magenes lemma.
This defines $e^{t\mathcal B}$ on $X_k$, consistently in $k$ by uniqueness, and the differential
inequalities above hold in integrated form.

The same argument applies to $\mathcal B^*_m$ in $L^2(\omega)$, whose drift $2(\ab\nabla\Phi+\bb)$
is radial and $O(\jap v^{\gamma+1})$, and to $\mathcal B_m$ in $L^2(\omega^{-1})$. The duality
$\langle e^{t\mathcal B_m}g,\varphi\rangle=\langle g,e^{t\mathcal B^*_m}\varphi\rangle$ used in~(v)
follows by differentiating $s\mapsto\langle e^{(t-s)\mathcal B_m}g,e^{s \mathcal B^*_m}\varphi\rangle$.
We use the $\tilde H^1$ estimate of~(iv) along $e^{t\mathcal B^*_m}$ in the same way, with $V$
replaced by its analogue for $(\varphi,\nabla\varphi)$.

(v) Let $m=\jap v^{k_1}$ and $m_1=\jap v^k$, and put $\omega_0:=1$, $\omega_1:=\jap v^{|\gamma|/2}$
and $\omega:=m/(m_1\jap v^{\gamma/2})=\jap v^{k_1-k+|\gamma|/2}$. For $\varphi_t:=e^{t \mathcal B^*_m}\varphi$,
by the first inequality of~(iv) at the weight $\omega_1$ (which needs
$k_1>\frac{\gamma+3}2+\frac{|\gamma|}2=\frac32$) and the second at the weight $\omega_0$, we have,
since $\jap v^{\gamma/2}\omega_1=1$ and $|\nabla\varphi|\le|\tn\varphi|$,
\begin{equation}
\frac d{dt}\Big(\|\varphi_t\|^2_{L^2(\omega_1)}+\eta t\|\varphi_t\|^2_{\tilde
H^1(\omega_0)}\Big)\le-2c\big(\|\varphi_t\|^2_{L^2}+\|\nabla\varphi_t\|^2_{L^2}\big)+\eta\|\varphi_t\|^2_{\tilde
H^1(\omega_0)}\le0
\end{equation}
for $\eta\le\min\{1,c\}$ (the second inequality of~(iv) remains true when $\eta$ is decreased). Thus
$\eta t\|\varphi_t\|^2_{H^1}\le C\|\varphi\|^2_{L^2(\omega_1)}$. For $t\ge1$, applying this at time
$1$ and the estimate for $e^{t\mathcal B^*_m}$ in~(iii), from $L^2(\omega)$ to $L^2(\omega_1)$, we
get $\|\varphi_t\|_{H^1}\le C\|\varphi_{t-1}\|_{L^2(\omega_1)}\le C\jap
t^{-b}\|\varphi\|_{L^2(\omega)}$. Here we apply~(iii) with $k_2=k_1$ and with
$k-\frac{|\gamma|}2<k_1-\frac{|\gamma|}2$ in place of $k<k_1$, so that its source weight is $\jap
v^{k_1-k+|\gamma|/2}=\omega$ and its target weight is $\jap v^{|\gamma|/2}=\omega_1$. Its hypothesis
$k-\frac{|\gamma|}2>\frac{\gamma+3}2$ is $k>\frac32$, and the rate is the same. Thus $\|e^{t
\mathcal B^*_m}\varphi\|_{H^1}\le C(t\wedge1)^{-1/2}\jap t^{-b}\|\varphi\|_{L^2(\omega)}$ for all
$t>0$. By the definition of $Z_{k_1}$ as the dual of $Y_{k_1}$ for the pairing of $X_{k_1}$, we have
$|\langle mf,\psi\rangle|\le\|f\|_{Z_{k_1}}\big(\|\jap v^{\gamma/2}\psi\|^2_{L^2}+\|\jap
v^{\gamma/2}\tn\psi\|^2_{L^2}\big)^{1/2}\le\|f\|_{Z_{k_1}}\|\psi\|_{H^1}$, since $\jap
v^{\gamma/2}\jap v\le1$ for $\gamma\le-2$. Since $mS_{\mathcal B}(t)f=S_{\mathcal B_m}(t)(mf)$ and
$\langle S_{\mathcal B_m}(t)g,\varphi\rangle=\langle g,S_{\mathcal B^*_m}(t)\varphi\rangle$, by
duality we obtain
\begin{equation}
\|e^{t\mathcal B}f\|_{L^2(m_1\jap v^{\gamma/2})}=\sup_{\|\varphi\|_{L^2(\omega)}\le1}\langle
mf,e^{t\mathcal B^*_m}\varphi\rangle\le\|f\|_{Z_{k_1}}\sup_{\|\varphi\|_{L^2(\omega)}\le1}\|e^{t
\mathcal B^*_m}\varphi\|_{H^1},
\end{equation}
which is~(v), since $L^2(m_1\jap v^{\gamma/2})=X_{k-|\gamma|/2}$.
\end{proof}

\begin{lem}[{\cite[Lemma~4.3]{CM}} in polynomial weights]\label{lem:CMbil}
Let $\gamma\in[-3,-2]$, $m=\jap v^k$ with $k\ge0$, $\theta>\frac72$ and $\theta'>\frac94$, and let
\begin{equation}\label{eq:YmCM}
\|g\|_{Y(m)}^2:=\|\jap v^{\gamma/2}mg\|^2_{L^2}+\|\jap v^{\gamma/2}\tn(mg)\|^2_{L^2} .
\end{equation}
Then
\begin{equation}
|\langle Q(f,g),h\rangle_{L^2(m)}|\le C\Big(\|f\|_{L^2(\jap v^\theta)}\|g\|_{Y(m)}+\|f\|_{H^1(\jap
v^{\theta'})}\|g\|_{L^2(m)}\Big)\|h\|_{Y(m)} .
\end{equation}
\end{lem}

\begin{proof}
The proof of Lemma~\ref{lem:QmG} applies with $\Phi=\log m$ in place of the Gaussian exponent. It
rests on the coefficient bounds of \cite[Lemma~4.2]{CM} with $p=2$ and $p=4$, which hold for
$\gamma\in[-3,-2]$ by item~(4) of the proof of Proposition~\ref{prop:CMm2}, and on $\gamma\le-2$.
Here $\nabla\Phi=kv\jap v^{-2}$ is radial, so $|\widetilde{\nabla\Phi}|=|\nabla\Phi|\le k\jap
v^{-1}$. Let $G=mg$, $H=mh$, let $A_1,\dots,A_4,B_1,B_2$ be the terms of that proof, and let
$\|f\|:=\|f\|_{L^2(\jap v^\theta)}$. Then we have $|A_1|\le C\|f\|\int\jap v^\gamma|\tn G||\tn H|\le
C\|f\|\|g\|_{Y(m)}\|h\|_{Y(m)}$ and $|A_2|+|A_3|\le C\|f\|\int\jap v^{\gamma-1}(|\tn G||H|+|G||\tn
H|)\le C\|f\|\|g\|_{Y(m)}\|h\|_{Y(m)}$. Moreover, $|A_4|\le C\|f\|\int\jap v^{\gamma-2}|G||H|$,
which is bounded in the same way. We also have $|B_1|\le C\|f\|_{H^1(\jap v^{\theta'})}\int\jap
v^{\gamma+1}|G||\nabla H|\le C\|f\|_{H^1(\jap v^{\theta'})}\|\jap v^{(\gamma+2)/2}G\|_{L^2}\|\jap
v^{\gamma/2}\nabla H\|_{L^2}$, and $\jap v^{(\gamma+2)/2}\le1$, $|\nabla H|\le|\tn H|$. Finally,
$|B_2|\le C\|f\|_{H^1(\jap v^{\theta'})}\int\jap v^{\gamma}|G||H|\le C\|f\|_{H^1(\jap
v^{\theta'})}\|g\|_{L^2(m)}\|h\|_{Y(m)}$.
\end{proof}

\begin{lem}[Spectral gap at $\gamma=-2$]\label{lem:CMgap}
Let $\gamma=-2$. There is $\lambda_{\rm gap}>0$ such that
\begin{equation}\label{eq:CMgap}
\langle-\mathcal Lf,f\rangle_{E_\mu}\ge\lambda_{\rm gap}\|\Pi f\|^2_{E_\mu}
\end{equation}
for $f$ in the domain of $\mathcal L$ in $E_\mu$, and $\|e^{t\mathcal L}\Pi\|_{E_\mu\to E_\mu}\le
e^{-\lambda_{\rm gap}t}$. In particular (S3) holds at $\gamma=-2$, with $\Theta_{\rm s}(t)\le
Ce^{-\lambda_{\rm gap}t}$.
\end{lem}

\begin{proof}
The coercivity estimate \cite[(2.6)]{CM}, which \cite{CM} state for $-3\le\gamma\le1$, summarizing
\cite{Guo02,Mou06} and other works, reads
\begin{equation}\label{eq:CMcoer}
\langle\mathcal Lf,f\rangle_{E_\mu}\le-c_\gamma\big(\|\jap v^{(\gamma+2)/2}\Pi f\|^2_{E_\mu}+\|\jap
v^{\gamma/2}\tilde\nabla_v\Pi(\mu^{-1/2}f)\|^2_{L^2}\big),
\end{equation}
where $\Pi$ is the projection onto the orthogonal complement of $\ker\mathcal L$ in $E_\mu$
\cite[\S2.1]{CM} (we read $\Pi(\mu^{-1/2}f)$, as written in \cite[(2.6)]{CM}, as $\mu^{-1/2}\Pi f$,
and we only use the first term).

Since $\langle f,\psi\mu\rangle_{E_\mu}=\int\psi f$, this complement consists of the moment-free
functions, so this $\Pi$ is ours. At $\gamma=-2$ the weight $\jap v^{(\gamma+2)/2}$ is $1$, which
implies the first claim. As noted in \cite[\S2.1]{CM}, (2.6) provides a spectral gap exactly for
$-2\le\gamma\le1$.

The same follows from \cite[Thm.~1.3]{Mou06}, stated for $\gamma\in[-3,1]$ (in dimension~$3$). It is
stated for $M=e^{-|v|^2}$, and the form below, for $\mu$, follows by the scaling of
\cite[\S1.6]{Mou06}. With $f=\mu h$, the Dirichlet form is bounded below by $C(\|h-P_\mu
h\|^2_{H^1(\jap v^\gamma\mu\dd v)}+\|\jap v^{1+\gamma/2}(h-P_\mu h)\|^2_{L^2(\mu\dd v)})$, where
$P_\mu$ is the $L^2(\mu\dd v)$-orthogonal projection onto $\mathrm{span}\{1,v_i,|v|^2\}$, and the
second weight is $1$ at $\gamma=-2$.

The operator $\mathcal L$ is self-adjoint and nonpositive in $E_\mu$ \cite[\S2.1]{CM} and commutes
with $\Pi$, so $u(t)=e^{t\mathcal L}\Pi f$ satisfies $\frac d{dt}\|u\|^2_{E_\mu}=2\langle \mathcal Lu,u\rangle_{E_\mu}\le-2\lambda_{\rm
gap}\|u\|^2_{E_\mu}$. Finally, for $\vartheta\in(\frac12,1)$, we have $\|\cdot\|_{E_\mu}\le
C\|\cdot\|_{L^2(\mu^{-\vartheta})}$, which yields the bound on $\Theta_{\rm s}$. Here $e^{t\mathcal L}$
on $X_k$, $k>\frac32$, is the semigroup obtained from the semigroup $e^{t\mathcal B}$ of
Lemma~\ref{lem:CMpoly} by the perturbation $\mathcal A$, which is bounded by Lemma~\ref{lem:RVA}(i)
and $L^2(\mu^{-\vartheta})\subset X_k$, through the Duhamel formula. For data in $E_\mu$ it
coincides with the self-adjoint semigroup of $E_\mu$. Indeed, the latter is continuous with values
in $E_\mu\subset X_k$, lies in $L^2(0,T;V_k)$ by the second term of~\eqref{eq:CMcoer} (since $\jap
v^k\mu^{1/2}$ and its gradient are bounded), and is therefore a weak solution of
$\partial_tu=\mathcal Bu+\mathcal Au$ in $V_k'$. Such solutions are unique, by the uniqueness in
Lions' theorem applied to the Duhamel formula.
\end{proof}

\begin{prop}[The inputs at $\gamma=-2$]\label{prop:CMm2}
Let $\gamma=-2$. The statements (S0)--(S2), and \cite[Lemmas~4.1--4.3, Cor.~4.4]{CM}, hold at
$\gamma=-2$, at the polynomial weights used in Proposition~\ref{prop:CMthr} and in (CM2) of
Appendix~\ref{app:band} (Lemmas~\ref{lem:CMpoly} and~\ref{lem:CMbil}, and items~(3)--(5) below), and
the inputs used at the Gaussian weights of Lemmas~\ref{lem:GW} and~\ref{lem:Whypo} hold at
$\gamma=-2$ (items~(1) and~(2) below). Moreover, (S3) holds by Lemma~\ref{lem:CMgap}. Consequently
that proposition holds at $\gamma=-2$, and in its proof the kernel $\Theta_{\rm s}$ of the third
term then decays exponentially, and the third term is $O(\jap t^{-b})$, as for $\gamma<-2$.
\end{prop}

\begin{proof}
We go through \cite[\S\S2--4]{CM}. The estimates in polynomial weights are those proved above.
\begin{enumerate}[label=(\arabic*)]
\item \emph{Lemma~2.1} (the matrix $\ab[\mu]$). Parts (a)--(c) are its eigenstructure,
$\ell_1\sim2\jap v^\gamma$, $\ell_2\sim\jap v^{\gamma+2}$, $\mathrm{tr}\,\ab[\mu]=2J_{\gamma+2}$,
$\bb[\mu]=-\ell_1v$, with $|\partial^\beta\ab[\mu]|\le C_\beta\jap v^{\gamma+2-|\beta|}$. We take
the parts in turn.
\begin{itemize}
\item[(a)] This is \cite[Props.~2.3, 2.4]{DeL97}, namely the eigenstructure and the asymptotics
$\ell_1\sim\frac23|v|^\gamma\int|w|^2\mu$, $\ell_2\sim\frac13|v|^{\gamma+2}\int|w|^2\mu$. That paper
treats hard potentials, so for the soft potentials considered here we verify the statement directly.
For every $\gamma\in[-3,-2]$, the eigenstructure follows from the invariance of $\ab[\mu](v)$ under
the rotations that fix $v$. The asymptotics follow from Lemma~\ref{lem:ell1} for $\ell_1$, and from
$\ell_2=J_{\gamma+2}-\frac12\ell_1$ with~(e) for $\ell_2$. Thus (a) holds independently of
\cite{DeL97}.
\item[(b)] This is \cite[Lemma~3]{Guo02}, stated there for $\gamma\ge-3$, with $\mu=e^{-|v|^2}$ and
$\sigma=\phi*\mu$. The normalization of \cite{CM} follows by the scaling
$\ab[\mu](v)=\pi^{-3/2}2^{(\gamma+2)/2}\sigma(v/\sqrt2)$. The derivative bounds there come from
\cite[Lemma~2]{Guo02}, applied with $k=\phi^{ij}$, $\theta=\gamma+2>-3$ and $l=\mu$. The
eigenstructure comes from the rotation invariance, as in~(a). The lower bound with $K$ follows from
the eigenstructure, from the asymptotics of~(a), and from $\ell_1,\ell_2>0$ being continuous.
\item[(c)] This is elementary.
\item[(d)] This is \cite[Lemma~2.4(d)]{CTW} (cited in \cite{CM} as Lemma~2.3, the numbering of the
preprint), proved there under the standing assumption $\gamma\in[-2,1]$ from parts (a)--(c),
\cite[Lemma~2]{Guo02} (which assumes only $\theta>-3$) and $\bb[\mu]=-\ell_1v$.
\end{itemize} The asymptotics of~(a) also follow at $\gamma=-2$ from
$\ell_1(v)=\int|v-w|^\gamma|w_\perp|^2\mu(w)\dd w$, with $w_\perp$ the component of $w$ orthogonal
to $\hat v$, and $\ell_1+2\ell_2=2J_{\gamma+2}$. Part~(e), for $J_\beta=|\cdot|^\beta*\mu$ with
$\beta\in(-3,0)$, is used at $\gamma=-2$ with $\beta=\gamma=-2$. Where \cite{CM} uses
$J_{\gamma+2}$, at $\gamma=-2$ this is $J_0\equiv1$, outside the range of~(e) but trivial, and
$\mathrm{tr}\,\ab[\mu]=2J_0=2$.
\item \emph{Lemmas~2.2--2.4, 3.1, 3.2} in polynomial weights. These yield (S0), (S1),
\cite[(3.4)]{CM} and (S2), which are Lemma~\ref{lem:CMpoly}, proved there for every
$\gamma\in[-3,-2]$. In Gaussian weights $m=e^{q|v|^2}$ we only use the identity of~(ii) of that
proof and the asymptotics of $\tilde\zeta_m$, which we check directly at $\gamma=-2$. With
$\Phi=q|v|^2$, we have $\tilde\zeta_m=2q(2q-1)\ell_1|v|^2-\frac12\cb$, and $\ell_1|v|^2=2+o(1)$,
$-\frac12\cb=J_{-2}=O(\jap v^{-2})$, so $\tilde\zeta_m\le4q(2q-1)(1+o(1))=4q(2q-1)\jap
v^{\gamma+2}(1+o(1))$, the form used in the proof of Lemma~\ref{lem:GW}. The exponential case of
\cite[Lemma~3.1]{CM} serves in \cite{CM} only Proposition~3.3, which we replace by
Lemma~\ref{lem:CMgap}.
\item \emph{Lemma~4.1}, $\int|v-v_*|^\alpha\jap{v_*}^{-\theta}\dd v_*\le C\jap v^\alpha$ for
$\alpha\in(-3,0)$ and $\theta>3$, is applied in the proof of Lemma~4.2 with $\alpha=(\gamma+2)\frac
p{p-1}$ and $\alpha=(\gamma+1)\frac p{p-1}$. At $\gamma=-2$, with $p=2$ and $p=4$, these are
$\alpha=0$, where it is trivial, and $\alpha=-\frac43$, which is in the range.
\item \emph{Lemma~4.2.} Part~(i) needs $p>\frac3{5+\gamma}=1$, and the proof takes $p=2$. Part~(ii)
needs $p>\frac3{4+\gamma}=\frac32$, and the proof takes $p=4$. The H\"older integrals become
$\int\jap{v_*}^{-2\bar\theta}\dd v_*$ and $\int|v-v_*|^{-4/3}\jap{v_*}^{-4\bar\theta/3}\dd v_*$,
both finite, and the bounds $\jap v^{\gamma+2}$, $\jap v^{\gamma+1}$ become $1$ and $\jap v^{-1}$.
Part~(i) bounds $|(a*f)v_iv_j|+|(a*f)v_i|+|a*f|$ by $C\jap v^{\gamma+2}\|f\|_{L^p(\jap v^\theta)}$,
$\theta>2+3(1-\frac1p)$, which contains the component bounds $|(a*f)\hat v\cdot\hat v|\le C\jap
v^\gamma\|f\|_{L^2(\jap v^\theta)}$ and $|(a*f)\hat v|\le C\jap v^{\gamma+1}\|f\|_{L^2(\jap
v^\theta)}$, $\theta>\frac72$, used in Subsection~\ref{ssec:tri}. Its proof involves no more than
$a(z)z=0$ and the H\"older integrals above.
\item \emph{Lemma~4.3 and Corollary~4.4.} Lemma~4.3 in polynomial weights is Lemma~\ref{lem:CMbil},
for every $\gamma\in[-3,-2]$. The embedding $H^1_*(\jap v^k)\subset H^1(\jap v^{\theta'})$ of
Corollary~4.4 needs $k+\frac\gamma2>\frac94$, that is, $k>\frac{13}4$ at $\gamma=-2$, which holds at
$k=7$. Together with $L^2(\jap v^k)\subset L^2(\jap v^\theta)$, $\theta\in(\frac72,k]$, it implies
(E1)(f)(iv).
\item \emph{Lemma~2.5} is replaced by Lemma~\ref{lem:RVA}, which covers $\gamma=-2$.
\end{enumerate}
With (S0)--(S3) at $\gamma=-2$, the proof of Proposition~\ref{prop:CMthr} applies without change.
\end{proof}

The target space of \cite[Lemma~3.2]{CM} is $L^2(m_1\jap v^{\gamma/2})$, which yields the threshold
$k>\frac72+|\gamma|$ in Proposition~\ref{prop:CMthr}(iii), which is below~$7$. The theory of
\cite{CTW} in Gaussian-type spaces covers $\gamma=-2$, but it does not by itself reach the
polynomial spaces $X_k$. Indeed, in polynomial weights $\mathcal B$ is only weakly dissipative, with
no rate, so the enlargement argument of \cite{CTW} does not apply there. It enters here through
\cite[Lemma~2.4(d)]{CTW} in item~(1) and, as a cross-check of Lemma~\ref{lem:CMgap}, through the
remark after \cite[(1.13)]{CTW} that $\mathcal L$ has a spectral gap in $E_\mu$ if and only if
$\gamma\ge-2$.

\subsection{The homogeneous inputs}\label{app:hom}

Here we prove the results on homogeneous trajectories that Section~\ref{sec:homog} states without
published proofs (for the Gaussian bounds (GS4) and (GS5), \cite{GS} give sketches of proof),
starting with the qualitative regularity of Lemma~\ref{lem:Rqual}, which the entropy method needs.

\emph{Previous results.} For $\gamma=-3$ a uniform rate is given in \cite[Thm.~2]{CDH}, namely
algebraic decay of the relative entropy for data with polynomial moments, and stretched exponential
decay for data with stretched exponential moments. That paper remarks that its dissipation estimates
extend to $-4<\gamma<0$, and we write out the case $-4<\gamma\le-2$ in Lemma~\ref{lem:EL1}. The
moment estimates \cite[Lemma~8, Cor.~4.2]{CDH} are proved for $-4<\gamma\le-2$. For moderately soft
potentials, $\gamma\in(-2,0)$ (the endpoint $-2$ excluded), a polynomial rate is proved in
\cite[Thm.~1.3]{Car15}. However, that approach needs a weighted Fisher information that grows at
most polynomially, which is obtained there from weighted $H^1$ estimates that are not available for
$\gamma<-2$. We are not aware of a published rate for $\gamma\in(-3,-2]$ that is uniform on the
classes $\mathfrak K^N$. The proof below needs no weighted Fisher information, and the unweighted
one, which is nonincreasing by (GS1), suffices.

Throughout, $\gamma\in[-3,-2]$, $\mu=(2\pi)^{-3/2}e^{-|v|^2/2}$, and a function $F\ge0$ is
\emph{normalized} if $\int F=1$, $\int Fv=0$ and $\int F|v|^2=3$. We write
\begin{equation}
\mathcal H(F):=\int F\log F,\qquad\mathcal H(F|\mu):=\int F\log\frac F\mu=\mathcal H(F)-
\mathcal H(\mu),\qquad i(F):=\int\frac{|\nabla F|^2}F,
\end{equation}
\begin{equation}
M_l(F):=\int\jap v^lF,\qquad I_\gamma(F|\mu):=\int F\Big|\nabla\log\frac F\mu\Big|^2\jap v^\gamma ,
\end{equation}
and $\mathcal D(F)=-\int Q(F,F)\log F$ for the entropy dissipation. These are the conventions of
\cite{CDH}. For normalized $F$, $M_2(F)=4$ and $M_l(F)\ge1$.

\subsubsection{Pointwise inequalities}

\begin{lem}[Glaeser's inequality and the Fisher information]\label{lem:FPgl}
\leavevmode
\begin{enumerate}[label=(\roman*)]
\item Let $F\ge0$ be $C^2$ on $\R^3$ and $M(v):=\sup_{B(v,1)}|D^2F|$. Then $|\nabla
F(v)|^2\le2F(v)\big(M(v)+2F(v)\big)$.
\item Let $n>3$, and let $F\ge0$ be $C^2$ with $\sum_{|\alpha|\le2}\sup_v\jap v^n|\nabla^\alpha
F(v)|\le A$. Then $i(F)\le C_n(1+A)$.
\end{enumerate}
\end{lem}

\begin{proof}
(i) We may assume $\nabla F(v)\ne0$. Let $e=\nabla F(v)/|\nabla F(v)|$ and $0<s\le1$. By Taylor's
formula and $F\ge0$, we have $0\le F(v-se)\le F(v)-s|\nabla F(v)|+\frac12M(v)s^2$. If $|\nabla
F(v)|\le M(v)$, we take $s=|\nabla F(v)|/M(v)$ and get $|\nabla F|^2\le2FM$. Otherwise we take $s=1$
and get $|\nabla F|\le F+\frac12M<F+\frac12|\nabla F|$, which implies $|\nabla F|^2\le4F^2$.

(ii) Where $F(v)=0$, we also have $\nabla F(v)=0$, since $v$ is a minimum, and the integrand is $0$.
Elsewhere, by (i) we have $|\nabla F|^2/F\le2M+4F$. For $|w-v|\le1$, $\jap v^2\le1+(|w|+1)^2\le3\jap
w^2$, so $M(v)\le3^{n/2}A\jap v^{-n}$. Integrating, we obtain $i(F)\le2\cdot3^{n/2}A\int\jap
v^{-n}+4\int F$, and $\int F\le A\int\jap v^{-n}$.
\end{proof}

\subsubsection{Qualitative regularity}

\begin{proof}[Proof of Lemma~\ref{lem:Rqual}]
\emph{Bounds on $[0,b]$.} By (RM) of Proposition~\ref{prop:RS}, we have $H\le C_Me^{-\beta_M|v|^2}$
on $[0,b]$. Note that (RM) is proved in Subsection~\ref{sssec:homRM} below, for data in $\mathfrak K^N$,
$N\ge2$, independently of this lemma. By Lemma~\ref{lem:RSshortF}(iii), with $m=4$ and
$q=\min\{\beta/16,1/24\}$, there is $\sigma_0>0$ such that $e^{(q/4)|v|^2}\partial^\alpha H$ is
bounded on $[0,\sigma_0]$ for $|\alpha|\le2$. By (GS5), applied on the windows
$[\sigma',\sigma'+\sigma_0]\subset[0,b]$, we get $|\nabla H|+|D^2H|\le C_3e^{-(\beta_M-\eps)|v|^2}$
on $[\sigma_0/2,b]$. Hence, on $[0,b]$, we have
\begin{equation}\label{eq:FPRbd}
H+|\nabla H|+|D^2H|\le Ce^{-\beta''|v|^2},\qquad\beta''>0 .
\end{equation}

\emph{The collision operator.} For $\gamma\in[-3,-2]$ the kernels $|w|^{\gamma+2}$ and
$|w|^{\gamma+1}$ are locally integrable, so $|\ab[H]|+|\bb[H]|\le C(\|H\|_{L^\infty}+\|H\|_{L^1})\le
C$. Moreover, $\cb[H]$ is bounded, by the same argument for $\gamma>-3$, and as $-8\pi H$ for
$\gamma=-3$. In nondivergence form we have $Q(H,H)=\ab[H]:D^2H-\cb[H]H$, so by~\eqref{eq:FPRbd} we
get $|\partial_\sigma H|=|Q(H,H)|\le Ce^{-\beta''|v|^2}$ on $[0,b]$.

\emph{Continuity and dissipation.} Let $f=H-\mu$. Since $\|\jap
v^ke^{-\beta''|v|^2}\|_{L^2}<\infty$, we have $\|f(\sigma)-f(\sigma')\|_{X_k}\le C|\sigma-\sigma'|$,
which implies $f\in C([0,b];X_k)$. Moreover, $\|f\|_{Y_k}\le C\|\jap v^{k+1}(|H|+|\nabla
H|+\mu+|\nabla\mu|)\|_{L^2}\le C$, since the dissipation norm of \cite{CM} is bounded by the
weighted $H^1$ norm with the weight $\jap v^{k+1}$. Hence $f\in L^\infty(0,b;Y_k)$. Since
$\partial_\sigma f=Q(H,H)$ is bounded by $Ce^{-\beta''|v|^2}$ and $f$ is Lipschitz with values in
$X_k$, $\frac12\|f\|^2_{X_k}$ is Lipschitz, with derivative $\langle Q(H,H),f\rangle_{X_k}$ for
a.e.\ $\sigma$ (for every $\sigma>0$, by dominated convergence). This is the energy identity in
$X_k$.

\emph{H-solution.} Nonnegativity and the regularity class of \cite[Def.~1]{CDH} follow from the
Gaussian bound above and the continuity just shown. The weak form follows by integrating by parts
against $\varphi\in C^1_c$, using the decay. Conservation of mass, momentum and energy follows from
the same integration by parts with $\varphi=1,v_j,|v|^2$, cut off at radius $\varrho$, since the
boundary terms are $O(e^{-\beta''\varrho^2/2})$.

\emph{Entropy identity.} Let $\chi_\varrho$ be a cutoff at radius $\varrho$, and write
$Q(H,H)=\nabla\cdot(\ab[H]\nabla H-\bb[H]H)$. For $\sigma>0$, where $H>0$ by (GS1), we have
\begin{equation}
\frac d{d\sigma}\int\chi_\varrho H\log H=-\int\chi_\varrho\frac{(\ab[H]\nabla H-\bb[H]H)\cdot\nabla
H}H-\int(1+\log H)(\ab[H]\nabla H-\bb[H]H)\cdot\nabla\chi_\varrho .
\end{equation}
By Lemma~\ref{lem:FPgl}(i) and the same bound, we have $|\nabla H|^2/H\le Ce^{-\beta''|v|^2/2}$.
Hence
\begin{equation}
|\nabla H|\,|\log H|\le Ce^{-\beta''|v|^2/4}\sqrt H|\log H|\le Ce^{-\beta''|v|^2/4},
\end{equation}
since $\sqrt x|\log x|$ is bounded on $[0,C]$, and also $H|\log H|\le C\sqrt H$.
Therefore the second integral tends to $0$ as $\varrho\to\infty$, uniformly on $[0,b]$. The first
integral converges absolutely, and, after the symmetrization in $(v,v_*)$, it equals
\begin{equation}
-\frac12\iint|v-v_*|^{\gamma+2}\Big|\Pi(v-v_*)\big(\nabla\log H-\nabla_*\log
H_*\big)\Big|^2HH_*=-\mathcal D(H),
\end{equation}
where $\Pi(w)$ is the projection orthogonal to $w$. The symmetrization is legitimate, since
$\iint|v-v_*|^{\gamma+2}\big(\frac{|\nabla H|^2}HH_*+H\frac{|\nabla_*H_*|^2}{H_*}\big)<\infty$, by
the bound on $|\nabla H|^2/H$ and the local integrability of $|w|^{\gamma+2}$. Integrating in
$\sigma$, we obtain the entropy identity on $[\sigma',b]$ for every $\sigma'>0$. It extends to
$\sigma'=0$ by the continuity of $\int\chi_\varrho H\log H$, since $x\log x$ is uniformly continuous
on $[0,C]$ and $H\in C([0,b];L^\infty)$, together with the uniform smallness of the tails. The
entropy is finite and continuous in $\sigma$, and $H\in L^\infty(0,\infty;L^1_2\cap L\log L)$, by
the conservation laws and the entropy identity.
\end{proof}

\subsubsection{The uniform \texorpdfstring{$L^1$}{L1} rate}

We use the following facts.
\begin{enumerate}[label=(EL\arabic*)]
\item For $-4<\gamma\le-2$ and every $\bar H$ there is $K=K(\gamma,\bar H)$ such that every
normalized, smooth and positive $F$ with $\mathcal H(F)\le\bar H$ satisfies $I_\gamma(F|\mu)\le
K\,M_{2-\gamma}(F)\,\mathcal D(F)$. For $\gamma=-3$ this is \cite[Thm.~1]{CDH}. \cite[(19) in
Remark~2.1(1)]{CDH} states it for $-4<\gamma<0$, and in Lemma~\ref{lem:EL1} below we write out the
proof for $-4<\gamma\le-2$, which contains the range used here.
\item \cite[Lemma~8]{CDH}. Let $-4<\gamma\le-2$, $l>2$, and let $F$ be an H-solution with
normalized data $F_0\in L^1_l$. Then $M_l(F(\sigma))\le
C\,M_l(F_0)+C\,l^{(l-6)\frac{|\gamma+1|}{\gamma+4}-\gamma}\,\sigma$, with $C=C(\gamma,
\mathcal H(F_0))$.
\item (GS1): along the solutions of \cite[Thm.~1.2]{GS}, $i(F(\sigma))$ is nonincreasing.
\item \emph{Log-Sobolev inequality on convex sets} \cite[\S5.7]{BGL}. Let $\Omega\subset\R^3$ be
convex and $\mu_\Omega:=\mu\mathbf 1_\Omega/\mu(\Omega)$. Then for every probability density $g$ on
$\Omega$ with $\sqrt g\in H^1(\Omega)$, $\int_\Omega g\log\frac g{\mu_\Omega}\le\frac12\int_\Omega
g\big|\nabla\log\frac g{\mu_\Omega}\big|^2$. This follows from \cite[Cor.~5.7.2]{BGL} with $\rho=1$
by approximation, since the potential of $\mu_\Omega$ is $\frac12|v|^2$ plus a constant on $\Omega$
and $+\infty$ outside, and so it is not smooth. The potentials
$W_k:=(\frac12|\cdot|^2+k\,\mathrm{dist}(\cdot,\Omega)^2)*\phi_{1/k}$, with a mollifier
$\phi_{1/k}$, are smooth with $\nabla^2W_k\ge I$, since $\mathrm{dist}(\cdot,\Omega)^2$ is convex.
Moreover, $e^{-W_k}/\int e^{-W_k}\to\mu_\Omega$ in $L^1$, by dominated convergence with $e^{-W_k}\le
e^{-|v|^2/2}$. We apply the corollary to $f=(g/\mu_\Omega)^{1/2}$, extended smoothly to $\R^3$ when
it is smooth up to $\partial\Omega$, which is the case used in Lemma~\ref{lem:Lloc}, and let
$k\to\infty$. The general case follows by density. The same approximation is used in \cite{BGL} for
the Poincar\'e inequality on convex bodies, after \cite[Thm.~4.6.3]{BGL}.
\item \emph{Sobolev.} $\|F\|_{L^3}=\|\sqrt F\|_{L^6}^2\le C_{\rm Sob}\,i(F)$, with $C_{\rm
Sob}:=S_3^2/4$, where $S_3$ is the constant of $\|\phi\|_{L^6(\R^3)}\le S_3\|\nabla\phi\|_{L^2}$
(since $i(F)=4\|\nabla\sqrt F\|^2_{L^2}$).
\end{enumerate}

\begin{lem}[(EL1) for $-4<\gamma\le-2$]\label{lem:EL1}
Let $-4<\gamma\le-2$ and $\bar H\in\R$. There is $K=K(\gamma,\bar H)$ such that every normalized,
smooth and positive $F$ with $\mathcal H(F)\le\bar H$ satisfies $I_\gamma(F|\mu)\le
K\,M_{2-\gamma}(F)\,\mathcal D(F)$.
\end{lem}

\begin{proof}
We follow \cite[\S3]{CDH}, where $\gamma=-3$, and write out the two steps in which $\gamma$ enters.
Put $X(v,w):=\nabla\log F(v)-\nabla\log F(w)$ and $R_{ij}(v,w):=(v_i-w_i)X_j-(v_j-w_j)X_i$, as in
\cite[(23)]{CDH}. Since $a(z)\xi\cdot\xi=|z|^{\gamma}|z\times\xi|^2$ and
$|z\times\xi|^2=\frac12\sum_{i,j}(z_i\xi_j-z_j\xi_i)^2$, the formula for $\mathcal D$ in the proof
of Lemma~\ref{lem:Rqual} implies that
\begin{equation}
\mathcal D(F)=\frac14\sum_{i,j=1}^3\iint F(v)F(w)\,|R_{ij}(v,w)|^2|v-w|^\gamma\dd v\dd w ;
\end{equation}
for $\gamma=-3$ this is the identity used in \cite[p.~393]{CDH}. Let $P_{ij}:=\int Fv_iv_j$ and
$\Lambda_F:=\inf_{i\ne j}(P_{ii}P_{jj}-P_{ij}^2)$. Then $|P_{ij}|\le\frac32$ for $i\ne j$,
$P_{ii}\le3$ and $\sum_iP_{ii}=3$.

(a) \emph{The steps that do not involve $\gamma$.} The identities \cite[(24)--(26)]{CDH}
(\cite[Prop.~3]{CDH}) are algebraic consequences of the normalization. The lower bound $\Lambda_F\ge
c(\bar H)>0$ of \cite[(32)]{CDH} rests on the normalization and $\mathcal H(F)\le\bar H$, through
\cite[(35)]{CDH}, and on the inequality of Desvillettes quoted in its proof \cite[p.~393]{CDH}.
Neither involves the kernel. The bound \cite[(35)]{CDH} controls $\sup_{|A|\le q}\int_AF$ through
the entropy, with a parameter $M>1$, and the proof of (32) takes $M=e^{4\bar H}$. In our setting
$\bar H\le0$ is possible, since $\mathcal H(\mu)<0$, and then $e^{4\bar H}\le1$. Moreover, on a set
$A$ the integral $\int_AF\log F$ is bounded below only through $\int F\log^-F$. We therefore use the
bound in the form $\sup_{|A|\le q}\int_AF\le Mq+(\bar H^++C_*)/\log M$, where $C_*$ is an absolute
bound of $\int F\log^-F$ for normalized $F$ (since $\int F\log^-F\le\int F|v|^2+C\int
e^{-|v|^2/2}$), with $\bar H$ replaced by $\bar H^++C_*$ in the proof of (32). Only the value of
$c(\bar H)$ changes.

(b) \emph{The analogue of \cite[Prop.~4]{CDH}.} Fix $i\ne j$. By \cite[(24)]{CDH}, we have
\begin{equation}
\partial_i\log F(v)+v_i=\frac{P_{ij}v_j+\big(P_{ij}^2+P_{ii}(1-P_{jj})\big)v_i+\int
R_{ij}(v,w)F(w)\,(w_iP_{ij}-w_jP_{ii})\dd w}{P_{ij}^2-P_{ii}P_{jj}} .
\end{equation}
By the Cauchy--Schwarz inequality with the splitting $1=|v-w|^{\gamma/2}\cdot|v-w|^{-\gamma/2}$, we
get
\begin{equation}
\Big(\int|R_{ij}(v,w)|F(w)(|w_i|+|w_j|)\dd w\Big)^2\le\int F(w)|R_{ij}(v,w)|^2|v-w|^\gamma\dd
w\cdot\int F(w)|v-w|^{|\gamma|}\,2|w|^2\dd w .
\end{equation}
Since $|v-w|^{|\gamma|}\le2^{|\gamma|}(\jap v^{|\gamma|}+\jap w^{|\gamma|})$ and $\jap v^\gamma\jap
v^{|\gamma|}=1$, we have $\sup_v\jap v^\gamma\int F(w)|v-w|^{|\gamma|}|w|^2\dd
w\le2^{|\gamma|+1}M_{2-\gamma}(F)$. We multiply the square of the identity by $F(v)\jap v^\gamma$,
integrate, use $\jap v^\gamma|v|^2\le\jap v^2$ and $M_2(F)=4$, and sum over $i$ (with one $j\ne i$
for each $i$). Then we obtain
\begin{equation}
I_\gamma(F|\mu)\le C_\gamma\Lambda_F^{-2}\Big(\sup_{i\ne
j}P_{ij}^2+\sup_j(P_{jj}-1)^2+M_{2-\gamma}(F)\,\mathcal D(F)\Big).
\end{equation}

(c) \emph{The analogue of \cite[(33), (34)]{CDH}.} Inserting \cite[(24), (25)]{CDH} into
\cite[(26)]{CDH}, we have, as in \cite[p.~394]{CDH},
\begin{equation}
\big|(v_i^2-v_j^2)P_{ij}+v_iv_j(P_{jj}-P_{ii})\big|\le\int|R_{ij}(v,w)|F(w)\,9(1+|v||w|)\dd w .
\end{equation}
We use the weight $\jap v^{\gamma-2}$ in place of $\jap v^{-5}$. By the Cauchy--Schwarz inequality
as in~(b), and
\begin{equation}
\jap v^{\gamma-2}|v-w|^{|\gamma|}(1+|v||w|)^2\le2^{|\gamma|+1}\big(\jap v^{-2}+\jap
v^{\gamma-2}\jap w^{|\gamma|}\big)\big(1+\jap v^2\jap w^2\big)\le C_\gamma\jap w^{2+|\gamma|},
\end{equation}
we get
\begin{equation}
\int F(v)\jap v^{\gamma-2}\big|(v_i^2-v_j^2)P_{ij}+v_iv_j(P_{jj}-P_{ii})\big|^2\dd v\le C_\gamma
M_{2-\gamma}(F)\,\mathcal D(F).
\end{equation}
With $\varrho:=\big(P_{ij}^2+\frac14(P_{jj}-P_{ii})^2\big)^{1/2}$ and $\varphi$ such that
$P_{ij}=\varrho\cos\varphi$, $P_{jj}-P_{ii}=2\varrho\sin\varphi$, the expression inside the absolute
value is $\varrho\big((v_i^2-v_j^2)\cos\varphi+2v_iv_j\sin\varphi\big)$. Thus the left side is
bounded below by $S^\gamma_F\varrho^2$, with
\begin{equation}
S^\gamma_F:=\inf_{\varphi\in\R}\int F(v)\jap v^{\gamma-2}\big|(v_i^2-v_j^2)\cos\varphi+2v_iv_j\sin\varphi\big|^2\dd
v .
\end{equation}
\emph{The lower bound $S^\gamma_F\ge c(\gamma,\bar H)>0$.} We follow \cite[pp.~395--396]{CDH}, with
the weight $\jap v^{\gamma-2}$ in place of $\jap v^{-5}$. We use cylindrical coordinates
$(v_i,v_j)=r(\cos\vartheta,\sin\vartheta)$, $v_k=z$, $k\notin\{i,j\}$. Then
$(v_i^2-v_j^2)\cos\varphi+2v_iv_j\sin\varphi=r^2\cos(2\vartheta-\varphi)$. Fix
$\epsilon\in(0,\frac\pi8]$ and let $E_\varphi$ be the set where the distance of $2\vartheta-\varphi$
to $\frac\pi2+\pi\mathbb Z$ is at most $\epsilon$. Off $E_\varphi$, we have
$|\cos(2\vartheta-\varphi)|\ge\sin\epsilon$. Put $H_*:=\bar H^++C_*$, with $C_*$ of~(a), so that
$\sup_{|A|\le q}\int_AF\le Mq+H_*/\log M$ for all $q>0$ and $M>1$. For $\delta,R_1,R_2>0$, we have
\begin{equation}
\int F\jap v^{\gamma-2}r^4\cos^2(2\vartheta-\varphi)\dd
v\ge\sin^2\!\epsilon\Big(\int_{\{r\ge\delta,\,|v|\le R_1\}}F\jap v^{\gamma-2}r^4\dd
v-\int_{E_\varphi}F\jap v^{\gamma-2}r^4\dd v\Big).
\end{equation}
In the first integral we have $\jap v^{\gamma-2}r^4\ge\jap{R_1}^{\gamma-2}\delta^4$, and the
complement of $\{r\ge\delta,|v|\le R_1\}$ is contained in $\{r<\delta,|z|\le R_1\}\cup\{|v|>R_1\}$.
The first set has measure at most $2\pi\delta^2R_1$, and $\int_{|v|>R_1}F\le3/R_1^2$ since the
energy is $3$. In the second integral we have $\jap v^{\gamma-2}r^4\le\jap v^{\gamma+2}\le1$,
because $\gamma\le-2$. This is the step that needs $\gamma\le-2$ (in \cite{CDH}, $r^4\jap
v^{-5}\le1$). The set $E_\varphi\cap\{r\le R_2,|z|\le R_2\}$ has measure
$4\epsilon\cdot\frac{R_2^2}2\cdot2R_2=4\epsilon R_2^3$, since $E_\varphi$ consists, in
$\vartheta\in[0,2\pi)$, of four intervals of length $\epsilon$. Outside $\{|v|\le R_2\}$ we use
$3/R_2^2$ again. With the bound of~(a), for all $M_1,M_2>1$, we get
\begin{equation}
S^\gamma_F\ge\sin^2\!\epsilon\Big(\jap{R_1}^{\gamma-2}\delta^4\Big(1-2\pi\delta^2R_1M_1-\frac{H_*}{\log
M_1}-\frac3{R_1^2}\Big)-4\epsilon R_2^3M_2-\frac{H_*}{\log M_2}-\frac3{R_2^2}\Big),
\end{equation}
uniformly in $\varphi$. We choose $R_1:=\sqrt{12}$, $M_1:=e^{4H_*}$ and $\delta^2:=(8\pi
R_1M_1)^{-1}$, so that each of the three subtracted terms in the inner bracket is $\frac14$, and put
$a:=\frac14\jap{R_1}^{\gamma-2}\delta^4=\frac14\,13^{(\gamma-2)/2}(8\pi\sqrt{12}\,e^{4H_*})^{-2}$.
Then we choose $R_2:=(12/a)^{1/2}$, $M_2:=e^{4H_*/a}$ and
$\epsilon:=\min\{\frac\pi8,\,a/(16R_2^3M_2)\}$, so that each of the last three terms is at most
$\frac a4$. This yields $S^\gamma_F\ge\frac a4\sin^2\!\epsilon=:c(\gamma,\bar H)>0$. Here $C_*>0$,
so $H_*>0$ and $M_1,M_2>1$. Hence $\sup_{i\ne j}P_{ij}^2+\sup_{i,j}(P_{ii}-P_{jj})^2\le
C(\gamma,\bar H)M_{2-\gamma}(F)\mathcal D(F)$, and
$\sup_j|P_{jj}-1|\le\frac23\sup_{i,j}|P_{ii}-P_{jj}|$, since $\sum_jP_{jj}=3$.

(d) Combining (a)--(c), we obtain the claim, with $K=K(\gamma,\bar H)$. This completes the proof.
\end{proof}

\begin{lem}[Localized log-Sobolev]\label{lem:Lloc}
Let $F$ be normalized, smooth and positive, and $R\ge2$. Then
\begin{equation}\label{eq:Lloc}
\mathcal H(F|\mu)\le\frac12\int_{B_R}F\big|\nabla\log\frac F\mu\big|^2+\mathcal T_R(F),
\end{equation}
where
\begin{equation}
\mathcal T_R(F):=2\mu(B_R^c)+\int_{B_R^c}F\Big(\frac{|v|^2}2+3\Big)+\int_{B_R^c}F\log^+F .
\end{equation}
\end{lem}

\begin{proof}
We split $\mathcal H(F|\mu)=\int_{B_R}+\int_{B_R^c}$. \emph{Inside.} Let $m_R:=\int_{B_R}F\in(0,1]$
and $F_R:=F\mathbf 1_{B_R}/m_R$. Since $F_R/\mu_{B_R}=(F/\mu)\mu(B_R)/m_R$, we have
\begin{equation}
\int_{B_R}F\log\frac F\mu=m_R\int_{B_R}F_R\log\frac{F_R}{\mu_{B_R}}+m_R\log\frac{m_R}{\mu(B_R)} .
\end{equation}
Since gradients do not see the normalizing constants, by (EL4) the first term is bounded by
$\frac12\int_{B_R}F|\nabla\log(F/\mu)|^2$. For the second term, note that $m_R\log m_R\le0$. For
$R\ge2$, we have $\mu(B_R)\ge\mu(B_2)=\mathbb P(\chi^2_3\le4)>\frac12$ (the median of $\chi^2_3$ is
below its mean $3$), and by $\log(1/x)\le1/x-1$ we get
$m_R\log\frac1{\mu(B_R)}\le\frac{1-\mu(B_R)}{\mu(B_R)}\le2\mu(B_R^c)$. \emph{Outside.} We have
$\log(F/\mu)=\log F+\frac{|v|^2}2+\frac32\log(2\pi)$, with $\frac32\log(2\pi)<3$, and $F\log F\le
F\log^+F$.
\end{proof}

\begin{lem}[The tail]\label{lem:Ltail}
Let $l>2$ and $\kappa_l:=\frac34(l-2)$. There is $B_l$ such that for every normalized $F$ and
$R\ge2$, $\mathcal T_R(F)\le B_l\,(1+i(F))^{3/4}M_l(F)\,\jap R^{-\kappa_l}$.
\end{lem}

\begin{proof}
\emph{Gaussian term.} We have $2\mu(B_R^c)\le C_le^{-R^2/4}\le C_l\jap R^{-\kappa_l}$. \emph{Energy
term.} On $B_R^c$, we have $\frac12|v|^2+3\le3\jap v^2\le3\jap v^l\jap R^{-(l-2)}$, so this term is
bounded by $3M_l\jap R^{-(l-2)}\le3M_l\jap R^{-\kappa_l}$. \emph{Entropy term.} Since $\log s\le
s^{1/2}$ for $s\ge1$, we have $F\log^+F\le F^{3/2}$. By H\"older's inequality with the exponents
$\frac43$ and $4$, and by (EL5), we get
\begin{equation}
\int_{B_R^c}F^{3/2}\le\Big(\int_{B_R^c}F\Big)^{3/4}\|F\|_{L^3}^{3/4}\le\big(\jap
R^{-l}M_l\big)^{3/4}\big(C_{\rm Sob}\,i(F)\big)^{3/4} .
\end{equation}
Finally $M_l^{3/4}\le M_l$ and $\jap R^{-3l/4}\le\jap R^{-\kappa_l}$.
\end{proof}

\begin{prop}[Entropy--entropy dissipation]\label{prop:LEEP}
Let $l>2$ and $a:=|\gamma|/\kappa_l=\frac{4|\gamma|}{3(l-2)}$. For all $\bar H$ and $I_0$ there is
$c_*=c_*(\gamma,l,\bar H,I_0)>0$ such that every normalized, smooth, positive $F$ with $
\mathcal H(F)\le\bar H$ and $i(F)\le I_0$ satisfies
\begin{equation}
\mathcal D(F)\ge c_*\,\frac{\mathcal H(F|\mu)^{1+a}}{M_{2-\gamma}(F)\,M_l(F)^a} .
\end{equation}
\end{prop}

\begin{proof}
On $B_R$, we have $1\le\jap R^{|\gamma|}\jap v^\gamma$. Hence, by Lemma~\ref{lem:Lloc}, (EL1) and
Lemma~\ref{lem:Ltail}, for every $R\ge2$, with $\mathcal H=\mathcal H(F|\mu)$, we have
\begin{equation}\label{eq:Lsplit}
\mathcal H\le\frac K2\jap R^{|\gamma|}M_{2-\gamma}\,\mathcal D+\frac{\Lambda_F}2\jap
R^{-\kappa_l},\qquad\Lambda_F:=2B_l(1+I_0)^{3/4}M_l(F) .
\end{equation}
If $\mathcal H=0$ there is nothing to prove. Otherwise we choose $\jap
R:=\max\{\sqrt5,(\Lambda_F/\mathcal H)^{1/\kappa_l}\}$. Then $R\ge2$ and $\frac{\Lambda_F}2\jap
R^{-\kappa_l}\le\frac12\mathcal H$, and it follows from~\eqref{eq:Lsplit} that $\mathcal D\ge\mathcal H\jap
R^{-|\gamma|}/(KM_{2-\gamma})$. Since $\mathcal H\le H^*:=\bar H-\mathcal H(\mu)$ and
$\Lambda_F\ge2B_l$, we have $(\Lambda_F/\mathcal H)^a\ge(2B_l/H^*)^a$, and
\begin{equation}
\jap R^{|\gamma|}=\max\{5^{|\gamma|/2},(\Lambda_F/\mathcal H)^a\}\le C'(\Lambda_F/
\mathcal H)^a,\qquad C':=\max\{1,5^{|\gamma|/2}(H^*/2B_l)^a\}.
\end{equation}
This proves the claim with $c_*=\big(KC'(2B_l(1+I_0)^{3/4})^a\big)^{-1}$.
\end{proof}

\begin{proof}[Proof of Proposition~\ref{prop:Lrate}]
\emph{A priori bounds.} By the hypothesis (for $F(0)\in\mathfrak K^N$, $N\ge2$, by
Lemma~\ref{lem:Rqual}), $\mathcal H(F(\sigma)|\mu)$ is absolutely continuous with $\frac
d{d\sigma}\mathcal H(F|\mu)=-\mathcal D(F)$ a.e., and $\mathcal H(F(\sigma))\le\bar H$. By
hypothesis, $i(F(\sigma))\le I_0$ for $\sigma>0$. By (EL2), we have $M_l(F(\sigma))\le
C_l(1+\sigma)$, with $C_l=C_l(\gamma,l,\bar H,A_l)\ge1$. Indeed, the constant of \cite[Lemma~8]{CDH}
depends on $\mathcal H(F_0)$ only through the bound $\int_0^\infty\mathcal D(F)\le\mathcal H(F_0)-\mathcal H(\mu)\le\bar
H-\mathcal H(\mu)$ used after \cite[(46)]{CDH} (mass and energy are fixed), so an upper bound on the
entropy suffices. By H\"older's inequality with $\jap v^{2-\gamma}=\jap v^{2(1-\theta)}\jap
v^{l\theta}$, $\theta:=|\gamma|/(l-2)\in(0,1)$, we have
\begin{equation}
M_{2-\gamma}(F)\le M_2^{1-\theta}M_l^\theta\le4C_l^\theta(1+\sigma)^\theta .
\end{equation}
\emph{The differential inequality.} By Proposition~\ref{prop:LEEP}, applied for $\sigma>0$, where
$F$ is smooth and positive by (GS1), with $\mathcal H(\sigma):=\mathcal H(F(\sigma)|\mu)$, we have
\begin{equation}
\mathcal H'(\sigma)\le-c\,\mathcal H^{1+a}(1+\sigma)^{-b_*},\qquad
b_*:=\theta+a=\frac{7|\gamma|}{3(l-2)}<1,\qquad c:=\frac{c_*}{4C_l^{\theta+a}} ,
\end{equation}
where $b_*<1$ is the hypothesis $l>2+\frac73|\gamma|$. Since $\mathcal H$ is nonincreasing, it
vanishes from the first time it vanishes. Before that time, we have $\frac d{d\sigma}\mathcal H^{-a}\ge
ac(1+\sigma)^{-b_*}$, which implies $\mathcal H(\sigma)^{-a}\ge\frac{ac}{1-b_*}\big((1+\sigma)^{1-b_*}-1\big)$.
For $\sigma\ge\sigma_b:=2^{1/(1-b_*)}-1$ the bracket is at least $\frac12(1+\sigma)^{1-b_*}$, so
$\mathcal H(\sigma)\le\big(\frac{2(1-b_*)}{ac}\big)^{1/a}(1+\sigma)^{-(1-b_*)/a}$. For
$\sigma\le\sigma_b$, we have $\mathcal H\le H^*\le
H^*(1+\sigma_b)^{\vartheta_l}(1+\sigma)^{-\vartheta_l}$. Finally
\begin{equation}
\frac{1-b_*}a=\frac{3(l-2)}{4|\gamma|}-\frac74=\vartheta_l ,
\end{equation}
so the first bound holds with
$C_\vartheta:=\max\{(2(1-b_*)/(ac))^{1/a},H^*(1+\sigma_b)^{\vartheta_l}\}$. The $L^1$ bound is the
Csisz\'ar--Kullback--Pinsker inequality $\|F-\mu\|_{L^1}^2\le2\mathcal H(F|\mu)$.

\emph{Uniformity on $\mathfrak K^N$.} Data in $\mathfrak K(C_{\rm dat},\beta)$ are normalized,
$\mathcal H(F(0))\le\log C_{\rm dat}$ since $F(0)\le C_{\rm dat}$, and $M_l(F(0))\le C(l,C_{\rm
dat},\beta)$. For data in $ \mathfrak K^N$, $N\ge2$, by Lemma~\ref{lem:FPgl}(ii) with $n=k_N>3$ we
have $i(F(0))\le C(1+A_N)$. Since \cite[\S11]{GS} prove the monotonicity (EL3) on $(0,\infty)$, we
bound $i(F(\sigma))$ for $0<\sigma\le\sigma_0$ directly. By Lemma~\ref{lem:RSshortF}(iii), with
$m=4$ and $q=\min\{\beta/16,1/24\}$, $e^{(q/4)|v|^2}\partial^\alpha F$, $|\alpha|\le2$, is bounded
by $C(C_{\rm dat},\beta,A_N)$ on $[0,\sigma_0]$, $\sigma_0=\sigma_0(C_{\rm dat},\beta,A_N)$ (as in
the proof of (RM)), and then Glaeser's bound implies $i(F(\sigma))\le C(C_{\rm dat},\beta,A_N)$
there. Together with (EL3) on $[\sigma_0,\infty)$, we obtain $i(F(\sigma))\le I_0:=C(C_{\rm
dat},\beta,A_N)$ for all $\sigma>0$. The qualitative regularity lemma yields the H-solution property
and the entropy identity.
\end{proof}

\begin{remark}\label{rem:Lrate}
(1) The restrictions on $\gamma$ are as follows. (EL1) is proved for $-4<\gamma\le-2$, (EL2) needs
$-4<\gamma\le-2$, (EL3) needs $|\gamma|\le\sqrt{19}$, and (GS1)--(GS5) need $\gamma\ge-3$. Thus the
proof covers $\gamma\in[-3,-2]$, and at $\gamma=-3$ it yields an alternative proof of a polynomial
rate, uniform on $\mathfrak K^N$ (compare \cite[Thm.~2(i)]{CDH}). (2) The Fisher information is used
to obtain the $L^3$ bound on the tail in Lemma~\ref{lem:Ltail}. Any uniform-in-time bound in $L^p$,
$p>1$, would do, with the exponents $\frac34$ replaced by some $\theta(p)\in(0,1)$. Without such a
bound, $\int_{B_R^c}F\log^+F$ is not controlled by moments.
\end{remark}

\subsubsection{Maxwellian bound and smoothness on bounded intervals}\label{sssec:homRM}

We prove Proposition~\ref{prop:RS} and Lemma~\ref{lem:RSmed}, following \S11 and Appendix~B of
\cite{GS} and tracking the constants. By \cite[Remark~11.1]{GS}, the only a priori estimate that may
deteriorate as $S\to\infty$ is the propagation of moments (GS2). The $L^\infty$ bound and the
Maxwellian bound inherit this dependence on $S$, so (RM) is stated on bounded intervals.
Quantitatively, (RM) and (RS) are applied at two values, namely $S=S^\sharp=4S_\sharp$ (fibered
comparison) and $S=T_*+2$ in Lemma~\ref{lem:Abarrier}, both fixed before $T_0$, while
Lemma~\ref{lem:Rqual} uses (RM) on every bounded interval, qualitatively. Instead of the
higher-order derivative bounds of \cite[App.~B]{GS} we use the induction below, which suffices
because the energy inequality is linear in the top order once $C^2$ is controlled.

\begin{proof}[Proof of (RM)]
\emph{Controlled quantities.} For $F_0\in\mathfrak K$, mass and momentum are those of $\mu$ and the
energy is $3$. The entropy satisfies $-C\le\int F_0\log F_0\le\log C_{\rm dat}$, where the lower
bound follows from the energy. All weighted norms $\|F_0\|_{L^1_k}$ are bounded by $C(k,C_{\rm
dat},\beta)$. The exponent $\beta_0=\beta_0(\gamma)$ of (GS4), proved in Proposition~\ref{prop:GS4}
for $\gamma\in[-3,-2]$, depends only on $\gamma$ (the mass is $1$ and the energy $3$). In particular
it does not depend on $S$. Moreover, we have $\|F(\sigma)\|_{L^3}\le C_{\rm Sob}I_0\le C(C_{\rm
dat},\beta,A_N)$ for $\sigma>0$, with $I_0$ of the proof of Proposition~\ref{prop:Lrate}, by (GS1)
on $(0,\infty)$, (EL5), Lemma~\ref{lem:RSshortF}(iii) and Lemma~\ref{lem:FPgl}(ii). Consequently
$q_S$, $q_\natural$ and $q_0$ do not depend on $S$.

\emph{Short time.} By the short-time lemma with $m=4\le N+2$ and $q=\min\{\beta/16,1/24\}$, there is
$\sigma_0=\sigma_0(C_{\rm dat},\beta,A_N)>0$ such that
$\sum_{|\alpha|\le2}\|e^{(q/4)|v|^2}\partial^\alpha F\|_{L^\infty}\le C(C_{\rm dat},\beta,A_N)$ on
$[0,\sigma_0]$. Hence $\|F\|_{L^\infty}\le C$ on $[0,\sigma_0]$, and by Glaeser's bound we have
$i(F(\sigma_0))\le C$.

\emph{$L^\infty$ on $[\sigma_0,S]$.} By (GS1), we have $i(F(\sigma))\le i(F(\sigma_0))$, and then
$\|F(\sigma)\|_{L^3}\le C$ by (EL5). By (GS2), $\|F(\sigma)\|_{L^1_{k'}}\le C(S)$ for every $k'$.
Interpolating, we get $\|F(\sigma)\|_{L^p_k}\le C(S)$ for $p\in(1,3)$ and every $k$. We take
$p\in(3/(5+\gamma),3)$, which is possible since $5+\gamma\ge2$. Applying (GS3) to the trajectory
restarted at $\sigma_1=\sigma_0/2$, together with the short-time bound on $[\sigma_0/2,\sigma_0]$,
we obtain $\sup_{[0,S]}\|F\|_{L^\infty}\le C_4(S,C_{\rm dat},\beta,A_N)$.

\emph{Maxwellian bound.} We replace $\beta$ by $\beta_M:=\min\{\beta,\beta_0\}$, which keeps $F_0\le
C_{\rm dat}e^{-\beta_M|v|^2}$, and apply (GS4) with $K=C_4$ and the moment bound $A=C(S)$ of (GS2)
at the order $k_4(\gamma)$. (Here $F$ is continuous on $[0,S]\times\R^3$, since on $[0,\sigma_0]$
the bound of the short-time step and Lemma~\ref{lem:GScoef} imply $|\partial_\sigma F|\le C$, and
$F$ is smooth for $\sigma>0$ by (GS1).) Then we obtain $C_M=C_{\rm Gauss}(\gamma,\beta_M,C_{\rm
dat},C_4,C(S))$.
\end{proof}

\begin{proof}[Proof of Lemma~\ref{lem:RSmed}]
On $[0,\sigma_0]$, with $\sigma_0=\sigma_0(m,q,E^\natural_m(0))$, this is Lemma~\ref{lem:RSshortF}.
On $[\sigma_0,S]$, applying (RM) and (GS5) with $\eps=\beta_M/2$ on the windows
$[\sigma',\sigma'+\sigma_0]$, we have
\begin{equation}
\sup_{[\sigma_0,S]}e^{2q|v|^2}\big(|F|+|\nabla F|+|D^2F|\big)\le K'_2(S,C_{\rm dat},\beta,A_N),
\end{equation}
since $2q\le\beta_M/8$. Here $K_2'=K_2'(m,q,S,C_{\rm dat},\beta,A_N)$, through the window length
$\sigma_0$. Fix $m$, and the ladder $\mathfrak m_k$ of $m$. For $j\le m$ let
$E_{(j)}:=\sum_{|\alpha|\le j}\|\mathfrak m_{|\alpha|}\partial^\alpha F\|^2_{L^2}$, $E^{\rm
top}_j:=\sum_{|\alpha|=j}\| \mathfrak m_j\partial^\alpha F\|^2_{L^2}$ and $D^{\rm
top}_j:=\sum_{|\alpha|=j}\|\partial^\alpha F\|^2_{Y_2(\mathfrak m_j)}$. On $[0,\sigma_0]$, the same
lemma bounds $E_{(m)}=E^\natural_m$. On $[\sigma_0,S]$, we have $E_{(2)}\le CK_2'^2$, and we show by
induction on $j=3,\dots,m$ that $\sup_{[\sigma_0,S]}E_{(j)}+\int_{\sigma_0}^SD^{\rm top}_j\le
K_{(j)}$. Let $K_{(j-1)}$ bound $E_{(j-1)}$ on $[\sigma_0,S]$, fix $|\alpha|=j$, and write the
equation for $\partial^\alpha F$ as $\sum_{a\le\alpha}\binom\alpha
aQ(\partial^aF,\partial^{\alpha-a}F)$, paired with $\mathfrak m_j^2\partial^\alpha F$.
\begin{itemize}
\item[(A)] $a=0$. As in Step~3(A) of the proof of that lemma, now with
$\|F\|_{L^\infty}\le C_M$ and $F\le C_Me^{-\beta_M|v|^2}$ in place of the bounds by the energy,
(G1)--(G3) of Lemma~\ref{lem:Whypo} hold with constants depending on $C_M$, and (G2) holds by
Lemma~\ref{lem:RGellip}. By Lemma~\ref{lem:RFtools}(vi) this term is $\le-c_{W,\rm
f}\|\partial^\alpha F\|^2_{Y_2(\mathfrak m_j)}+C_{W,\rm f}\|\partial^\alpha F\|^2_{X( \mathfrak m_j)}$,
with $c_{W,\rm f},C_{W,\rm f}$ depending on $(q,C_M)$.
\item[(B)] $0\ne a\ne\alpha$. By part~(iii) of the same lemma with $f=\partial^aF$,
$h=\partial^{\alpha-a}F$ and the weight $\mathfrak m_j$, the term is $\le C\big(\|f\|_{L^2(\jap
v^4)}\|h\|_{Y_2(\mathfrak m_j)}+\|f\|_{H^1(\jap v^3)}\|h\|_{X(\mathfrak m_j)}\big)\|\partial^\alpha
F\|_{Y_2(\mathfrak m_j)}$. If $|a|\le j-2$, the norms of $f$ are $\le CK_{(j-1)}^{1/2}$, and, by the
ladder ($|\alpha-a|<j$, as in Step~3(B) there), those of $h$ are $\le C(K_{(j-1)}+E^{\rm
top}_j)^{1/2}$. If $|a|=j-1$, then $|\alpha-a|=1$, the norms of $h$ are $\le CK_2'$, and those of
$f$ are $\le C(K_{(j-1)}+E^{\rm top}_j)^{1/2}$.
\item[(C)] $a=\alpha$. By part~(ii) of that lemma with the weight $\mathfrak m_j$, with the
bound $\sup|\ab[f]|\le C\|f\|_{X(\mathfrak m_j)}$ of part~(i) (compare
\cite[Lemma~4.2(i)]{CM}) and the pointwise bounds on $F$ and $\nabla^2F$, this term is $\le
CK_2'E^{\rm top}_j$.
\end{itemize}
By Young's inequality, we obtain
\begin{equation}\label{eq:RS-top}
\frac12\frac d{d\sigma}E^{\rm top}_j+\frac{c_{W,\rm f}}2D^{\rm top}_j\le
C(1+K_{(j-1)}+K_2'^2)(K_{(j-1)}+E^{\rm top}_j)
\end{equation}
on $[\sigma_0,S]$. This inequality is linear in $E^{\rm top}_j$.

\emph{Justification.} The identity~\eqref{eq:RSident} used above holds, in integrated form, on every
interval on which $\sup E^\natural_m+\int D^\natural_m<\infty$ and $F\in C(H^{m-1}_{\rm loc})$,
by Lemma~\ref{lem:RSident}. Hence on such an interval~\eqref{eq:RS-top} holds in integrated form,
and
Gronwall's inequality applies.

The argument above, and part~(ii) of the short-time lemma on $[0,\sigma_0]$, apply without change at
every base $\lambda\in[q,2q]$ in place of $q$, with constants that are uniform in $\lambda$, since
we have $2\lambda\le4q\le\beta_M/8$ and $\lambda\le2q\le\frac1{24}$ by $q\le q_S$. Let $K_*$ be the
resulting bound for $\sup_{[0,T]}E^{(\lambda)}_m$, $\lambda\in[q,2q]$, $T\le S$, valid whenever
$\sup_{[0,T]}E^{(\lambda)}_m+\int_0^TD^{(\lambda)}_m<\infty$ and $F\in C([0,T];H^{m-1}_{\rm loc})$
(the Gronwall constants are nondecreasing in $T$), and let $h>0$ be a lower bound for
$\sigma_0(m,\lambda,K_*)$, $\lambda\in[q,2q]$, which exists by the last sentence of the short-time
lemma.

Put $J:=\lceil S/h\rceil$, $t_j:=\min\{jh,S\}$ and $\lambda_j:=q(2-j/J)$. We show by induction on
$j$ that $\sup_{[0,t_j]}E^{(\lambda_j)}_m+\int_0^{t_j}D^{(\lambda_j)}_m<\infty$ and $F\in
C([0,t_j];H^{m-1}_{\rm loc})$. For $j=J$ this is the finiteness at the base $q$ on $[0,S]$. For
$j=0$ it holds by Step~0 of the proof of that lemma (with $2\lambda_0=4q\le\beta/8$). If it holds
for $j<J$, then by the a priori argument at the base $\lambda_j$ on $[0,t_j]$ we have
$E^{(\lambda_{j+1})}_m(t_j)\le E^{(\lambda_j)}_m(t_j)\le K_*$. We apply part~(iv) of the same lemma
at $\sigma'=t_j$ with the base $\lambda_{j+1}$ in place of $q$,
$q^+=\lambda_j\in(\lambda_{j+1},2\lambda_{j+1}]$, and with (RM), which yields $F\le
C_Me^{-\beta_M|v|^2}$ on $[0,S+1]$ and $\beta_M\ge32q\ge8\lambda_j$. This implies the finiteness at
the base $\lambda_{j+1}$ on $[t_j,t_j+h]\supset[t_j,t_{j+1}]$. On $[0,t_j]$ it follows from the
induction hypothesis, since $\lambda_{j+1}<\lambda_j$.

By Gronwall's inequality from $E^{\rm top}_j(\sigma_0)\le2E^\natural_m(0)$, we obtain $K_{(j)}$. For
$j=m$ this is the claim, with $\int_0^SD^\natural_m$ bounded by the short-time bound on
$[0,\sigma_0]$ and by $\sum_{j\le m}\int_{\sigma_0}^SD^{\rm top}_j$ on $[\sigma_0,S]$. The orders
$j\le1$ are bounded by $K_2'$. For $j=2$, $D^{\rm top}_2$ contains third derivatives with the weight
$\mathfrak m_2>\mathfrak m_3$, and we integrate~\eqref{eq:RS-top} at $j=2$ on $[\sigma_0,S]$. Note
that (A)--(C) apply without change at $j=2$, with all norms of $f$ and $h$ bounded by $CK_2'$, so
that its right side is bounded by $C(1+K_2'^2)K_2'^2$.
\end{proof}

\begin{proof}[Proof of (RS)]
We use Lemma~\ref{lem:RSmed} with $m=N+2$ and $q=q_S$, and Sobolev's inequality as in
Lemma~\ref{lem:RSshortF}(iii). Then we have $\|e^{(q/4)|v|^2}\partial^\alpha F\|_{L^\infty}\le
C\sum_{|\beta|\le2}\|\jap v^2e^{(q/4)|v|^2}\partial^{\alpha+\beta}F\|_{L^2}\le
C(E^\natural_m)^{1/2}$ for $|\alpha|\le N$, since $\jap v^2e^{(q/4)|v|^2}\le C\mathfrak m_k$ for
every $k\le m$.
\end{proof}

\subsubsection{Entry into the near-Maxwellian regime}\label{sssec:homenter}

We prove Theorem~\ref{thm:enter} for all $\gamma\in[-3,-2]$ at once. The definition of an admissible
background in Section~\ref{sec:homog} contains only the two properties proved here.

\emph{Inputs at level $k\ge7$.} We use (E1) of Inputs~\ref{hyp:E}, at $\gamma=-2$ with
Proposition~\ref{prop:CMm2}. Let $\mathcal L=L_\mu$, $\Pi=I-\Pi_0$, $X_k$, $Y_k$, $Z_k$ be as in
Appendix~\ref{app:band}, and $X^0_k$ etc.\ the moment-free subspaces.
\begin{itemize}
\item[(CM3)] Proposition~\ref{prop:CMthr}(ii) (compare \cite[Prop.~3.6, Cor.~3.7]{CM}). There is
$\eta_k>0$ such that $\nn h_k^2:=\eta_k\|h\|^2_{X_k}+\int_0^\infty\|e^{\tau\mathcal L}h\|^2_{L^2}\dd\tau$
is a Hilbert norm on $X^0_k$, with $c_{*,k}\|h\|_{X_k}\le\nn h_k\le C_{*,k}\|h\|_{X_k}$. With
$\ip\cdot\cdot_k$ the associated inner product (extended to $Z^0_k\times Y^0_k$ by (CM1$_k$)), there
is $K_k>0$ with $\ip{\mathcal Lh}h_k\le-K_k\|h\|^2_{Y_k}$ for $h\in X^0_k\cap Y^0_k$.
\item[(CM1$_k$)] (CM1) of the same appendix with $7$ replaced by $k$. In particular
$J_k<\infty$, hence $|\ip G\varphi_k|\le C_J\|G\|_{Z_k}\|\varphi\|_{Y_k}$ for $G\in Z^0_k$,
$\varphi\in Y^0_k$, with $C_J:=\eta_k+J_k$.
\item[(CM2$_k$)] (CM2) at level $k$: $\|Q(f,h)\|_{Z_k}\le
C_Q(\|f\|_{X_k}\|h\|_{Y_k}+\|f\|_{Y_k}\|h\|_{X_k})$.
\end{itemize}
Since $\mathrm{mom}(Q(a,b)+Q(b,a))=0$, both $\mathcal Lf=Q(\mu,f)+Q(f,\mu)$ and $Q(f,f)$ lie in
$Z^0_k$ when $f\in X_k\cap Y_k$. This yields $\partial_\sigma f\in Z^0_{k_1}$ in Step~3 below. The
constant $K_k$ is not the constant $K$ of (EL1).

\begin{proof}[Proof of Theorem~\ref{thm:enter}]
Let $F(0)\in\mathfrak K^N(C_{\rm dat},\beta,A_N)$, $N\ge2$, and $f:=F-\mu$. Then $f$ is moment-free,
since mass, momentum and energy are conserved. Constants $C$ depend only on $(\gamma,C_{\rm
dat},\beta,A_N,k_1)$.

\emph{Step 1: weighted $L^1$ convergence.} Let $p:=4k_1+9$ and $l:=25$. Then $l>2+\frac73|\gamma|$
and $\vartheta_{25}=\frac{69}{4|\gamma|}-\frac74\ge4$ for $|\gamma|\le3$. By
Proposition~\ref{prop:Lrate}, whose hypotheses hold on $\mathfrak K^N$, $N\ge2$, by
Lemma~\ref{lem:Rqual}, we have $\|f(\sigma)\|_{L^1}\le C(1+\sigma)^{-2}$. By (EL2), we have
$M_{2p}(F(\sigma))\le C(1+\sigma)$. By Cauchy--Schwarz, we get
\begin{equation}
\|f(\sigma)\|_{L^1(\jap v^p)}\le\|f(\sigma)\|_{L^1}^{1/2}\big(M_{2p}(F(\sigma))+M_{2p}(\mu)\big)^{1/2}\le
C(1+\sigma)^{-1/2} .
\end{equation}
\emph{Step 2: closeness at every time.} Applying the Cauchy--Schwarz inequality to
$\int\varphi^2\jap v^{2k_1}=\int(|\varphi|\jap v^{p})^{1/2}(|\varphi|^3\jap v^{-9})^{1/2}$, we get
$\|\varphi\|_{X_{k_1}}\le\|\varphi\|^{1/4}_{L^1(\jap v^p)}\|\varphi\|^{3/4}_{L^3(\jap v^{-3})}$. By
(GS1), (EL5) and Lemma~\ref{lem:FPgl}(ii), we have $\|F(\sigma)\|_{L^3}\le C_{\rm
Sob}\,i(F(\sigma))\le C_{\rm Sob}I_0\le C$ for all $\sigma>0$, with $I_0$ as in the proof of that
proposition, and $\|F(0)\|_{L^3}\le C_{\rm dat}^{2/3}$, so $\|f(\sigma)\|_{L^3(\jap v^{-3})}\le C$.
Hence
\begin{equation}
\|f(\sigma)\|_{X_{k_1}}\le C(1+\sigma)^{-1/8}\qquad\text{for every }\sigma\ge0 .
\end{equation}
\emph{Step 3: trapping.} By the qualitative regularity lemma, $f\in C([0,b];X_{k_1})\cap
L^2(0,b;Y_{k_1})$ for every $b$, and by (CM2$_{k_1}$) we have $\partial_\sigma f=\mathcal Lf+Q(f,f)\in
L^2(0,b;Z^0_{k_1})$. By the Lions--Magenes lemma \cite[Ch.~XVIII, \S1, Thms.~1--2]{DL}, in the
Gelfand triple $V^0\subset(X^0_{k_1},\ip\cdot\cdot_{k_1})\subset(V^0)'$ with $V=X_{k_1}\cap
Y_{k_1}$, into which $Z^0_{k_1}$ embeds through $\ip\cdot\cdot_{k_1}$ by (CM1$_{k_1}$),
$\sigma\mapsto\nn{f(\sigma)}^2_{k_1}$ is absolutely continuous with derivative $2\ip{\partial_\sigma
f}f_{k_1}$. By (CM3), (CM1$_{k_1}$) and (CM2$_{k_1}$), with $\|f\|_{X_{k_1}}\le c_{*,k_1}^{-1}\nn
f_{k_1}$, we have
\begin{equation}
\frac d{d\sigma}\nn f^2_{k_1}\le\big(-2K_{k_1}+4C_JC_Qc_{*,k_1}^{-1}\nn
f_{k_1}\big)\|f\|^2_{Y_{k_1}} ,
\end{equation}
which is the analogue of \cite[Prop.~5.1]{CM}. If
$\nn{f(s_*)}_{k_1}\le\theta_0:=K_{k_1}c_{*,k_1}/(4C_JC_Q)$, a continuity argument shows that $\nn
f_{k_1}$ is nonincreasing on $[s_*,\infty)$ and
$K_{k_1}\int_{s_*}^\infty\|f\|^2_{Y_{k_1}}\le\nn{f(s_*)}^2_{k_1}$. Hence
\begin{equation}
\sup_{\sigma\ge s_*}\|f(\sigma)\|_{X_{k_1}}\le\frac{C_{*,k_1}}{c_{*,k_1}}\|f(s_*)\|_{X_{k_1}},\qquad\int_{s_*}^\infty\|f\|^2_{Y_{k_1}}\le\frac{C_{*,k_1}^2}{K_{k_1}}\|f(s_*)\|^2_{X_{k_1}}
.
\end{equation}
\emph{Conclusion.} We choose $T_*$ with $C(1+T_*)^{-1/8}\le\min\{\theta_0/C_{*,k_1},\ \delta
c_{*,k_1}/C_{*,k_1},\ \delta K_{k_1}^{1/2}/C_{*,k_1}\}$ and take $s_*:=T_*$. By Step~2, we have
$\nn{f(s_*)}_{k_1}\le C_{*,k_1}\|f(s_*)\|_{X_{k_1}}\le\theta_0$, and then Step~3 yields (AB1) and
(AB2) with $\delta$. The time $T_*$ depends only on $(\gamma,C_{\rm dat},\beta,A_N,k_1,\delta)$.
\end{proof}

\begin{remark}\label{rem:homenter}
Step~2 holds at every time, not only at a good time, because the Fisher information bounds $F$ in
$L^3$ uniformly in time. The trapping uses (E1) at the level $k_1\ge7$, above the thresholds of
(CM1)--(CM3), which at $\gamma=-2$ hold by Proposition~\ref{prop:CMm2}.
\end{remark}

\subsubsection{The bound \texorpdfstring{$F\le C\jap v^{-2}$}{F <= C/<v>**2} for
\texorpdfstring{$\gamma=-3$}{gamma=-3}}\label{ssec:homMoser}

We prove the case $\gamma=-3$ of the bound on $F$ in Lemma~\ref{lem:Abarrier}(ii). A direct barrier
$K\jap v^{-2}$ fails for $\gamma=-3$, because the reaction term $8\pi F\cdot K\jap v^{-2}$ is not
dominated by the negative part of $\ab[F]:\nabla^2(K\jap v^{-2})$ without prior decay of $F$. We use
instead Moser's local boundedness \cite{Mos64} (see also \cite[Ch.~III, \S\S7--8]{LSU}, where the
local bound of $\max|u|$ is proved by De Giorgi's method and Moser's method is sketched), with its
constants written out, since the ellipticity degenerates like $\jap v^{-3}$.

For $(\sigma_0,v_0)$ and $r>0$ let $\mathcal C_r:=(\sigma_0-r^2,\sigma_0]\times B_r(v_0)$.

\begin{lem}[Local boundedness with explicit constants]\label{lem:MOS}
Let $u>0$ be bounded and continuous on $\overline{\mathcal C_1}$, smooth on
$(\sigma_0-1,\sigma_0]\times\overline{B_1(v_0)}$, and satisfy $\partial_\sigma u=\nabla\cdot(A\nabla
u-bu)$ in $\mathcal C_1$, where $A$ is symmetric with $\lambda|\zeta|^2\le
A\zeta\cdot\zeta\le\Lambda_{\rm Mo}|\zeta|^2$, $\lambda\le1\le\Lambda_{\rm Mo}$, and
$|b|\le\beta_{\rm d}$, $|\nabla\cdot b|\le\gamma_{\rm d}$. Put $H_{\rm Mo}:=\Lambda_{\rm
Mo}+\beta_{\rm d}+\gamma_{\rm d}$. Then
\begin{equation}
\sup_{\mathcal C_{1/2}}u\le C_{\rm Mo}\,\lambda^{-3/4}H_{\rm Mo}^{5/4}\,\|u\|_{L^2(\mathcal C_1)},
\end{equation}
with an absolute constant $C_{\rm Mo}$.
\end{lem}

\begin{proof}
\emph{Step 1: energy inequality.} Let $p\ge2$, $w:=u^{p/2}$ and $\frac12\le\rho<r\le1$. Let
$\eta=\eta_1(\sigma)\eta_2(v)$, with $\eta_2\in C^\infty_c(B_r(v_0))$, $0\le\eta_2\le1$, $\eta_2=1$
on $B_\rho(v_0)$, $|\nabla\eta_2|\le2/(r-\rho)$, and $\eta_1=0$ for $\sigma\le\sigma_0-r^2$,
$\eta_1=1$ for $\sigma\ge\sigma_0-\rho^2$, $0\le\eta_1'\le2/(r^2-\rho^2)\le2/(r-\rho)^2$. We test
the equation with $p\,u^{p-1}\eta^2$. Since $p(p-1)u^{p-2}A\nabla u\cdot\nabla
u=\frac{4(p-1)}pA\nabla w\cdot\nabla w$ and $p\,u^{p-1}\nabla u=2w\nabla w$, the diffusion
contributes
\begin{equation}
-\tfrac{4(p-1)}p\!\int\eta^2A\nabla w\cdot\nabla w-4\!\int\eta wA\nabla
w\cdot\nabla\eta\le-\lambda\!\int\eta^2|\nabla w|^2+4\Lambda_{\rm Mo}\!\int w^2|\nabla\eta|^2 ,
\end{equation}
by Cauchy--Schwarz and $p\ge2$. For the drift, we have $p(p-1)u^{p-1}\nabla u=(p-1)\nabla(u^p)$, and
by one more integration by parts we get $p\int bu\cdot\nabla(u^{p-1}\eta^2)=-(p-1)\int
w^2\eta^2\nabla\cdot b+2\int w^2\eta\,b\cdot\nabla\eta$, so no derivative of $u$ remains. Using
$|\nabla(\eta w)|^2\le2\eta^2|\nabla w|^2+2w^2|\nabla\eta|^2$ and $\lambda\le\Lambda_{\rm Mo}$, we
obtain
\begin{equation}
\frac d{d\sigma}\!\int(\eta w)^2+\frac\lambda2\!\int|\nabla(\eta w)|^2\le\int
w^2\big(2\eta|\partial_\sigma\eta|+5\Lambda_{\rm Mo}|\nabla\eta|^2+2\beta_{\rm
d}\eta|\nabla\eta|+p\gamma_{\rm d}\eta^2\big)\le\frac{24pH_{\rm Mo}}{(r-\rho)^2}\int_{B_r(v_0)}w^2 .
\end{equation}
We integrate from $\sigma_0-r^2+\epsilon$ and let $\epsilon\to0$, using $\eta_1(\sigma_0-r^2)=0$ and
the boundedness of $w$. With $E:=24pH_{\rm Mo}(r-\rho)^{-2}\iint_{\mathcal C_r}w^2$, we deduce that
$\sup_{\sigma\le\sigma_0}\int(\eta w)^2\le E$ and $\iint|\nabla(\eta w)|^2\le2E/\lambda$.

\emph{Step 2: gain of integrability.} Let $S_3$ be the Sobolev constant of (EL5), and $\kappa_{\rm
Mo}:=\frac53$. By H\"older's inequality in $v$, we have $\int g^{2\kappa_{\rm Mo}}\le(\int
g^2)^{2/3}\|g\|^2_{L^6}$, so $\iint(\eta w)^{2\kappa_{\rm Mo}}\le S_3^2E^{2/3}\cdot2E/\lambda$.
Since $\eta=1$ on $\mathcal C_\rho$ and $w^2=u^p$, we obtain
\begin{equation}
\|u\|_{L^{\kappa_{\rm Mo} p}(\mathcal C_\rho)}\le\Big(\frac{2S_3^2}\lambda\Big)^{\frac1{\kappa_{\rm
Mo} p}}\Big(\frac{24pH_{\rm Mo}}{(r-\rho)^2}\Big)^{\frac1p}\|u\|_{L^p(\mathcal C_r)} .
\end{equation}
\emph{Step 3: iteration.} We take $p_j=2\kappa_{\rm Mo}^j$ and $r_j=\frac12+2^{-j-1}$, so that
$r_j-r_{j+1}=2^{-j-2}$. Then $\sum_{j\ge0}1/p_{j+1}=\frac1{2(\kappa_{\rm Mo}-1)}=\frac34$ and
$\sum_{j\ge0}1/p_j=\frac{\kappa_{\rm Mo}}{2(\kappa_{\rm Mo}-1)}=\frac54$. Applying Step~2 with
$(p,r,\rho)=(p_j,r_j,r_{j+1})$ for $j<J$, we get $\|u\|_{L^{p_J}( \mathcal C_{r_J})}\le\lambda^{-a_J}H_{\rm
Mo}^{b_J}P_J\|u\|_{L^2(\mathcal C_1)}$, with $a_J:=\sum_{j<J}1/p_{j+1}\le\frac34$,
$b_J:=\sum_{j<J}1/p_j\le\frac54$ and
$P_J:=\prod_{j<J}\big(2S_3^2\big)^{1/p_{j+1}}\big(768(4\kappa_{\rm Mo})^j\big)^{1/p_j}$. Put $C_{\rm
Mo}:=\sup_JP_J$, which is finite since $\sum_j(1+j)/p_j<\infty$. The supremum is needed because
$2S_3^2<1$ for the sharp constant, so that $P_J$ need not be monotone in $J$. Since $\lambda\le1\le
H_{\rm Mo}$, we have $\lambda^{-a_J}\le\lambda^{-3/4}$ and $H_{\rm Mo}^{b_J}\le H_{\rm Mo}^{5/4}$.
Hence $\|u\|_{L^{p_J}(\mathcal C_{1/2})}\le C_{\rm Mo}\lambda^{-3/4}H_{\rm Mo}^{5/4}\|u\|_{L^2(
\mathcal C_1)}$ for all $J$. Letting $J\to\infty$, we obtain the claim.
\end{proof}

\begin{prop}[The case $\gamma=-3$ of Lemma~\ref{lem:Abarrier}(ii)]\label{prop:ABiv}
Let $\gamma=-3$, keep the hypotheses of that lemma, with $k_1\ge11$, and put
$M_9:=\sup_{\sigma\ge s_*}\int\jap v^9F(\sigma)$ and $H_F:=3+18\pi K_\infty$. Then $M_9\le
C(1+\delta_0)$, and for $\sigma\ge s_*+2$ and all $v$,
\begin{equation}
F(\sigma,v)\le C_c\jap v^{-2},\qquad C_c:=\max\Big\{10K_\infty,\ 2^{9/2}C_{\rm
Mo}\,c_G'^{\,-3/4}H_F^{5/4}(K_\infty M_9)^{1/2}\Big\},\quad c_G':=\min\{1,c_G/8\},
\end{equation}
with $c_G=c_G(\gamma,1,3,K_\infty)$ of Lemma~\ref{lem:RGellip}.
\end{prop}

\begin{proof}
\emph{The moment.} By (AB1) at level $k_1\ge11$, we have $M_9\le\int\jap v^9\mu+\|\jap
v^{11}(F-\mu)\|_{L^2}\|\jap v^{-2}\|_{L^2}\le C(1+\delta_0)$. If $|v_0|\le3$, then $F\le
K_\infty\le10K_\infty\jap{v_0}^{-2}$ by Lemma~\ref{lem:Abarrier}(i).

Let $|v_0|\ge3$ and $\sigma_0\ge s_*+2$. Then $\mathcal C_1\subset[s_*+1,\infty)\times B_1(v_0)$,
and $F\le K_\infty$ on $\mathcal C_1$. The equation is $\partial_\sigma F=\nabla\cdot(\ab[F]\nabla
F-\bb[F]F)$, with the kernels $a(z)=|z|^{-1}(I-\hat z\otimes\hat z)$ and $b(z)=-2|z|^{-3}z$, and
$\nabla\cdot\bb[F]=\cb[F]=-8\pi F$. The constants of Lemma~\ref{lem:MOS} are:
\begin{itemize}
\item $|\ab[F]|\le\int|v-w|^{-1}F(w)\dd w\le K_\infty\int_{B_1}|z|^{-1}\dd z+\int F=2\pi
K_\infty+1=:\Lambda_{\rm Mo}$;
\item $|\bb[F]|\le2\int|v-w|^{-2}F\le2(4\pi K_\infty+1)=:\beta_{\rm d}$, and
$|\nabla\cdot\bb[F]|=8\pi F\le8\pi K_\infty=:\gamma_{\rm d}$;
\item by Lemma~\ref{lem:RGellip} (mass $1$, energy $3$, $\|F\|_{L^\infty}\le K_\infty$), $\ab[F]\ge
c_G\jap v^{-3}I$, since $\jap v^{\gamma+2}\ge\jap v^\gamma$. On $B_1(v_0)$, $\jap v\le2\jap{v_0}$,
so $\lambda:=c_G'\jap{v_0}^{-3}$ works.
\end{itemize}
Thus $H_{\rm Mo}=\Lambda_{\rm Mo}+\beta_{\rm d}+\gamma_{\rm d}=H_F$. On $B_1(v_0)$, we have $\jap
v\ge|v_0|-1\ge\frac12\jap{v_0}$, since $|v_0|\ge3$. Hence $\|F\|^2_{L^2(\mathcal C_1)}\le
K_\infty\sup_\sigma\int_{B_1(v_0)}F(\sigma)\le2^9K_\infty M_9\jap{v_0}^{-9}$. Therefore, by the
local bound, we have
\begin{equation*}
F(\sigma_0,v_0)\le C_{\rm Mo}c_G'^{\,-3/4}\jap{v_0}^{9/4}H_F^{5/4}\big(2^9K_\infty
M_9\big)^{1/2}\jap{v_0}^{-9/2}\le C_c\jap{v_0}^{-2} .
\qedhere
\end{equation*}
\end{proof}

\subsubsection{Gaussian bounds on bounded intervals: proofs of (GS4) and
(GS5)}\label{sssec:homGauss}

The proofs of \cite[Thm.~2.3, Prop.~B.2]{GS}, our (GS4) and (GS5), are given as sketches that build
on the arguments of \cite[\S\S4--5]{CSS}, which are carried out there for $\gamma\in(-2,0]$. We
write out both for $\gamma\in[-3,-2]$, with the dependence of the constants that we need. The
constant of (GS4) does not depend on the length $S$ of the time interval, except through the
$L^\infty$ and moment bounds, and that of (GS5) depends on $T$ only through a lower bound on $T$.
Throughout, $G\ge0$ (or $F(\sigma)$) has unit mass and $\int|w|^2G\le3$, and we use the convention
of Section~\ref{sec:setting}, namely $\partial_\sigma F=\ab[F]:\nabla^2F-\cb[F]F$, with
$-\cb[F]=2(\gamma+3)|\cdot|^\gamma*F\ge0$ for $\gamma>-3$ and $-\cb[F]=8\pi F$ for $\gamma=-3$.
Recall that $a(u)\zeta\cdot\zeta=|u|^\gamma|u\times\zeta|^2$. For $\zeta=\hat v$ and $u=v-w$, we
have $|u\times\hat v|=|w_\perp|$, the component of $w$ orthogonal to $\hat v$, and
$|w_\perp|\le|v-w|$.

\begin{lem}[Coefficients from the $L^\infty$ and moment bounds]\label{lem:GScoef}
Let $\gamma\in[-3,-2]$, $\|G\|_{L^\infty}\le K$ and $\int\jap w^kG\le A$ with $k\ge10$. Put
$c_\perp:=\frac3{32}(1+2\sqrt3)^\gamma$ and $C_\parallel:=6\cdot2^{3|\gamma|/2}$.
\begin{enumerate}[label=(\roman*)]
\item For $|v|\ge7$ and $\zeta\perp v$, $\ab[G](v)\zeta\cdot\zeta\ge c_\perp\jap
v^{\gamma+2}|\zeta|^2$.
\item There are $C_u=C_u(\gamma,K,A)$ and $R_\parallel=R_\parallel(\gamma,K,A)$ such that, for
$|v|\ge1$, $\ab[G](v)\hat v\cdot\hat v\le C_u\jap v^\gamma$ and $|\ab[G](v)|\le C_u\jap
v^{\gamma+2}$, and $\ab[G](v)\hat v\cdot\hat v\le C_\parallel\jap v^\gamma$ for $|v|\ge
R_\parallel$. Moreover $|\ab[G]|+|\nabla\ab[G]|\le C_\gamma(1+K)$ on $\R^3$.
\item $-\cb[G]\le C_\gamma(1+K)$ for $\gamma>-3$, and $-\cb[G]=8\pi G\le8\pi K$ for $\gamma=-3$. Put
\begin{equation}\label{eq:k4}
k_4(\gamma):=\max\{10,\lfloor\frac{6(|\gamma|-2)}{3-|\gamma|}\rfloor+1\}
\end{equation}
for $\gamma>-3$, and $k_4(-3):=10$. If $\gamma>-3$ and $k\ge k_4(\gamma)$, there is
$C=C(\gamma,K,A)$ with $-\cb[G](v)\le C\jap v^{\gamma+2-\eta_\gamma}$, where
$\eta_\gamma:=\min\{2,\frac{k_4(\gamma)}{q_\gamma}+\gamma+2\}>0$ and $q_\gamma:=\frac6{3-|\gamma|}$.
\end{enumerate}
\end{lem}

\begin{proof}
Since $\int\jap w^{k'}G\le A$ for $k'\le k$, we use the moment of order $k'=10$ in (i)--(ii) and of
order $k'=k_4(\gamma)$ in (iii). Thus the constants $C_k$ below are $C_{k'}$, and $C_u$,
$R_\parallel$ and $C$ depend only on $(\gamma,K,A)$.

(i) Let $R_0:=2\sqrt3$. By Chebyshev's inequality, we have $\int_{B_{R_0}}G\ge\frac34$. For $w\in
B_{R_0}$ and $|v|\ge7>2R_0$, we have $|v-w|\le(1+R_0)\jap v$, so $|v-w|^\gamma\ge(1+R_0)^\gamma\jap
v^\gamma$. Moreover, for a unit $\zeta\perp v$, $|(v-w)\times\zeta|$ is the distance from $w$ to the
line $v+\R\zeta$, which is at least $|(w-v)\cdot\hat v|\ge|v|-R_0\ge\frac12|v|$. Hence
$\ab[G](v)\zeta\cdot\zeta\ge(1+R_0)^\gamma\jap v^\gamma\cdot\frac14|v|^2\cdot\frac34\ge c_\perp\jap
v^{\gamma+2}$, since $|v|^2\ge\frac12\jap v^2$.

(ii) For $|v|\ge1$ we split $\ab[G](v)\hat v\cdot\hat v=\int|v-w|^\gamma|w_\perp|^2G(w)\dd w$ over
three regions that cover $\R^3$. Note that $\jap w\ge\frac12\jap v$ on $B_1(v)$, so
$\int_{B_1(v)}G\le2^kA\jap v^{-k}$. On $\{|w|\le|v|/2\}$, we have
$|v-w|^\gamma\le2^{|\gamma|}|v|^\gamma\le2^{3|\gamma|/2}\jap v^\gamma$ and $|w_\perp|\le|w|$, and
this part is bounded by $3\cdot2^{3|\gamma|/2}\jap v^\gamma$. On $\{|w|>|v|/2,\ |v-w|\ge1\}$ the
integrand is bounded by $|w|^2G$, and this part is bounded by $C_kA\jap v^{2-k}$. On $\{|v-w|<1\}$
the kernel is bounded by $|v-w|^{\gamma+2}\le|v-w|^{-1}$. By H\"older's inequality with
$\||z|^{-1}\|_{L^{3/2}(B_1)}<\infty$ and $\|G\|_{L^3(B_1(v))}\le K^{2/3}(\int_{B_1(v)}G)^{1/3}\le
K^{2/3}(C_kA\jap v^{-k})^{1/3}$, this part is bounded by $CK^{2/3}A^{1/3}\jap v^{-k/3}$. Since
$k\ge10$, we have $\jap v^{2-k}+\jap v^{-k/3}\le2\jap v^{\gamma-1/3}$, which implies the first and
the third bound. The bound on $|\ab[G]|$ follows in the same way from $\ab[G](v)e\cdot
e\le\int|v-w|^{\gamma+2}G(w)\dd w$, $|e|=1$. Finally, $a$ is homogeneous of degree
$\gamma+2\in[-1,0]$, so $|a(z)|\le C|z|^{\gamma+2}$ and $|\nabla a(z)|\le C|z|^{\gamma+1}$, which
are integrable on $B_1$ and bounded outside. This proves the global bound.

(iii) The global bound is as in (ii). For the decay, we split $\int|v-w|^\gamma G(w)\dd w$ in the
same way. The first region contributes $2^{|\gamma|}|v|^\gamma$ and the second $C_kA\jap v^{-k}$. By
H\"older's inequality with $q_\gamma$ and $q_\gamma'=\frac6{3+|\gamma|}$, for which
$q_\gamma'|\gamma|<3$, the third is bounded by
$\||z|^\gamma\|_{L^{q_\gamma'}(B_1)}\max\{1,K\}(C_kA\jap v^{-k})^{1/q_\gamma}$, with
$k=k_4(\gamma)$. Since $k_4(\gamma)>\frac{6(|\gamma|-2)}{3-|\gamma|}$, we have
$\frac{k_4(\gamma)}{q_\gamma}>|\gamma|-2=-(\gamma+2)$, and the claim follows.
\end{proof}

\begin{prop}[(GS4) for $\gamma\in{[-3,-2]}$]\label{prop:GS4}
Let $\gamma\in[-3,-2]$, $\beta_0:=c_\perp/(2C_\parallel)$, which depends only on $\gamma$, and
$0<\beta\le\beta_0$. Let $F$ be a classical solution of $\partial_\sigma F=Q(F,F)$ on $[0,S]$ (for
instance the solution of (GS1), or the local solution in the proof of Lemma~\ref{lem:GS1q}) with
unit mass, energy $3$, $F$ continuous on $[0,S]\times\R^3$ and smooth for $\sigma>0$, with $F(0)\le
C_{\rm dat}e^{-\beta|v|^2}$, $\|F(\sigma)\|_{L^\infty}\le K$ and $\int\jap v^kF(\sigma)\le A$ for
$\sigma\in[0,S]$, where $k\ge k_4(\gamma)$. Then $F\le C_{\rm Gauss}e^{-\beta|v|^2}$ on
$[0,S]\times\R^3$, with $C_{\rm Gauss}=C_{\rm Gauss}(\gamma,\beta,C_{\rm dat},K,A)$, independent of
$S$.
\end{prop}

\begin{proof}
Let $\Phi:=e^{-\beta|v|^2}$ and $h:=\jap v^2$. Then $\ab:\nabla^2\Phi=\Phi\big(4\beta^2\ab v\cdot
v-2\beta\,\mathrm{tr}\,\ab\big)$ and $\ab:\nabla^2h=2\,\mathrm{tr}\,\ab$, with $\ab=\ab[F(\sigma)]$.

\emph{Comparison.} We use the following elementary principle. Let $\sigma_1\in[0,S)$, $U:=\{|v|>R\}$
for some $R\ge0$ (or $U=\R^3$), and let $W$ be continuous on $[\sigma_1,S]\times\overline U$, $C^1$
in $\sigma$ and $C^2$ in $v$ on $(\sigma_1,S]\times U$, with $\partial_\sigma W>\ab:\nabla^2W-\cb W$
on $(\sigma_1,S]\times U$, $F<W$ on $\{\sigma_1\}\times\overline U$ and on
$(\sigma_1,S]\times\partial U$, and $W-F\to\infty$ as $|v|\to\infty$, uniformly in $\sigma$. Then
$F<W$ on $[\sigma_1,S]\times U$. Indeed, otherwise the set $\{F\ge W\}$ is closed, and bounded since
$W-F\to\infty$. Let $\sigma_c>\sigma_1$ be its least time and $(\sigma_c,v_c)$ a point of it. Then
$v_c\in U$, $F=W$ at $(\sigma_c,v_c)$, $F\le W$ on $\{\sigma_c\}\times U$ and $F<W$ on
$[\sigma_1,\sigma_c)\times\{v_c\}$. Hence $\nabla^2(F-W)\le0$ and $\partial_\sigma(F-W)\ge0$ there
(a left derivative, which exists since $F$ is smooth for $\sigma>0$), and so $\partial_\sigma W\le\partial_\sigma
F=\ab:\nabla^2F-\cb F\le\ab:\nabla^2W-\cb W$, since $\ab\ge0$ and $F=W$ there, a contradiction.

\emph{Short time.} By Lemma~\ref{lem:GScoef}(ii)--(iii), we have $\ab v\cdot v\le|v|^2|\ab\hat
v\cdot\hat v|\le C$ (for $|v|\ge1$ since $\gamma\le-2$, and for $|v|<1$ by the global bound),
$\mathrm{tr}\,\ab\le C$ and $-\cb\le C$, with $C=C(\gamma,K,A)$. Hence there is $\Lambda_{\rm
G}=\Lambda_{\rm G}(\gamma,\beta,K,A)$ with $\ab:\nabla^2\Phi-\cb\Phi\le(\Lambda_{\rm G}-1)\Phi$ and
$\ab:\nabla^2h-\cb h\le(\Lambda_{\rm G}-1)h$. For $\delta>0$, $W_0:=e^{\Lambda_{\rm G}\sigma}(C_{\rm
dat}\Phi+\delta h)$ satisfies $\partial_\sigma W_0-\ab:\nabla^2W_0+\cb W_0\ge e^{\Lambda_{\rm
G}\sigma}(C_{\rm dat}\Phi+\delta h)>0$, $W_0(0)>F(0)$, and $W_0-F\ge\delta h-K\to\infty$. By the
comparison with $U=\R^3$ and letting $\delta\to0$, we get $F\le C_{\rm dat}e^{\Lambda_{\rm
G}\sigma}\Phi$, and hence $F\le C_{\rm dat}e^{\Lambda_{\rm G}}\Phi$ on $[0,\min\{1,S\}]$.

\emph{The reaction term for $\sigma\ge1$.} For $\gamma>-3$, by the same lemma we have $-\cb[F]\le
C\jap v^{\gamma+2-\eta_\gamma}$. For $\gamma=-3$, the proof of Proposition~\ref{prop:ABiv} applies
on $[0,S]$, with the cylinders $ \mathcal C_1=(\sigma_0-1,\sigma_0]\times B_1(v_0)$ for
$\sigma_0\in[1,S]$, the bound $\|F\|_{L^\infty}\le K$ in place of Lemma~\ref{lem:Abarrier}(i), and
$\int\jap v^9F\le A$ in place of (AB1). Then we obtain $F(\sigma,v)\le C_c(K,A)\jap v^{-2}$ for
$\sigma\in[1,S]$, so $-\cb[F]=8\pi F\le8\pi C_c\jap v^{-2}$. In both cases there is $R_{\rm
c}=R_{\rm c}(\gamma,\beta,K,A)$ with $-\cb[F]\le\beta c_\perp\jap v^{\gamma+2}$ for $|v|\ge R_{\rm
c}$ and $\sigma\in[1,S]$.

\emph{The barrier for $\sigma\ge1$.} If $S\le1$ we are done. Let $R_{\rm
G}:=\max\{7,R_\parallel,R_{\rm c}\}$. For $|v|\ge R_{\rm G}$, by the coefficient bounds~(i)--(ii) we
have $\ab v\cdot v\le C_\parallel|v|^2\jap v^\gamma\le C_\parallel\jap v^{\gamma+2}$ and
$\mathrm{tr}\,\ab\ge2c_\perp\jap v^{\gamma+2}$ (the sum over two orthonormal directions orthogonal
to $v$). Thus, since $4\beta^2C_\parallel\le2\beta c_\perp$, we obtain
\begin{equation}
\ab:\nabla^2\Phi-\cb\Phi\le\Phi\jap v^{\gamma+2}\big(4\beta^2C_\parallel-4\beta c_\perp+\beta
c_\perp\big)\le-\beta c_\perp\jap v^{\gamma+2}\Phi<0 .
\end{equation}
Let $C_{\rm Gauss}:=\max\{C_{\rm dat}e^{\Lambda_{\rm G}},\,Ke^{\beta R_{\rm G}^2}\}$ and, for
$\delta>0$, $W_1:=C_{\rm Gauss}\Phi+\delta e^{\Lambda_{\rm G}\sigma}h$. Then $\partial_\sigma
W_1>\ab:\nabla^2W_1-\cb W_1$ on $(1,S]\times\{|v|>R_{\rm G}\}$, $F(1)\le C_{\rm dat}e^{\Lambda_{\rm
G}}\Phi<W_1(1)$, $F\le K\le C_{\rm Gauss}\Phi<W_1$ on $|v|=R_{\rm G}$, and $W_1-F\to\infty$. By the
comparison and letting $\delta\to0$, we get $F\le C_{\rm Gauss}\Phi$ on $[1,S]\times\{|v|\ge R_{\rm
G}\}$, and $F\le K\le C_{\rm Gauss}\Phi$ on $|v|\le R_{\rm G}$. This completes the proof.
\end{proof}

\begin{prop}[(GS5) for $\gamma\in{[-3,-2]}$]\label{prop:GS5}
Let $\gamma\in[-3,-2]$ and let $F$ be a classical solution of $\partial_\sigma F=Q(F,F)$ on $[0,T]$
(for instance the solution of (GS1), or the local solution in the proof of Lemma~\ref{lem:GS1q})
with unit mass, energy $3$ and $F\le C_Ge^{-\beta|v|^2}$ on $[0,T]\times\R^3$. Then for every
$\eps\in(0,\beta)$, $|\nabla F|+|\nabla^2F|\le Ce^{-(\beta-\eps)|v|^2}$ on $[T/2,T]\times\R^3$, with
$C=C(\gamma,\beta,\eps,T,C_G)$. The constant depends on $T$ only through a lower bound $T\ge
T_{\min}$: the proof runs with $\rho:=\min\{\jap{v_0}^{-4},\frac14\sqrt{T_{\min}}\}$.
\end{prop}

\begin{proof}
By the Gaussian bound, we have $\|F\|_{L^\infty}\le C_G=:K$ and $\int\jap v^{10}F\le A(C_G,\beta)$,
so Lemma~\ref{lem:GScoef} applies, and Lemma~\ref{lem:RGellip} implies $\ab\zeta\cdot\zeta\ge
c_G(\jap v^\gamma|P_v\zeta|^2+\jap v^{\gamma+2}|(I-P_v)\zeta|^2)$, $c_G=c_G(\gamma,K)$.

\emph{Rescaling.} Fix $\sigma_0\in[T/4,T]$ and $v_0$. If $|v_0|\ge2$, let
$D:=\jap{v_0}^{-\gamma/2}P_{v_0}+\jap{v_0}^{-(\gamma+2)/2}(I-P_{v_0})$, and if $|v_0|<2$, let
$D:=I$. Let $\rho:=\min\{\jap{v_0}^{-4},\frac14\sqrt T\}$, $\mathcal A:=\rho D^{-1}$ and $\tilde
F(s,w):=F(\sigma_0+\rho^2s,v_0+\mathcal Aw)$ on $\mathcal Q:=(-1,0]\times B_1$. Then
$\partial_s\tilde F=\tilde a:\nabla^2_w\tilde F+\tilde c\tilde F$, with $\tilde
a(s,w):=D\ab(\sigma,v)D$ and $\tilde c:=-\rho^2\cb(\sigma,v)$ at
$(\sigma,v)=(\sigma_0+\rho^2s,v_0+\mathcal Aw)$. On $\mathcal Q$, $\sigma\ge\sigma_0-\rho^2\ge T/8$
and $|v-v_0|\le|\mathcal A|\le\rho$, since $|D^{-1}|\le1$.

\emph{Uniform ellipticity of $\tilde a$.} Let $|v_0|\ge2$,
\begin{equation}\label{eq:N0xi}
N_0(\xi):=\jap{v_0}^\gamma|P_{v_0}\xi|^2+\jap{v_0}^{\gamma+2}|(I-P_{v_0})\xi|^2,
\end{equation}
so that $N_0(D\zeta)=|\zeta|^2$, and let $N_v$ be the same expression at $v$. For $|v-v_0|\le\rho$,
we have $\frac12\jap{v_0}\le\jap v\le2\jap{v_0}$,
$|P_v-P_{v_0}|\le\vartheta_0:=2\sqrt5\jap{v_0}^{-5}$, and $|\xi|^2\le\jap{v_0}^{|\gamma|}N_0(\xi)$.
Using $(x-y)^2\ge\frac12x^2-y^2$ and $(x+y)^2\le2x^2+2y^2$, we get
\begin{equation}
N_v(\xi)\ge2^{-3}\big(\tfrac12N_0(\xi)-2\jap{v_0}^{\gamma+2}\vartheta_0^2|\xi|^2\big)\ge2^{-3}\big(\tfrac12-40\jap{v_0}^{-8}\big)N_0(\xi)\ge2^{-5}N_0(\xi),
\end{equation}
since $\jap{v_0}^2\ge5$, and similarly $N_v(\xi)\le8(2+80\jap{v_0}^{-8})N_0(\xi)$. With the lower
bound above, and the upper bound $\ab\xi\cdot\xi\le2\ab P_v\xi\cdot
P_v\xi+2\ab(I-P_v)\xi\cdot(I-P_v)\xi\le2C_uN_v(\xi)$ of the coefficient lemma ($\ab$ is
nonnegative), we deduce that $c'|\zeta|^2\le\tilde a\zeta\cdot\zeta\le C'|\zeta|^2$ on $\mathcal Q$,
with $c',C'$ depending on $(\gamma,K,A)$. For $|v_0|<2$, we have $\tilde a=\ab$ on $B_1(v_0)\subset
B_3$, where the same bounds apply.

\emph{H\"older continuity in $w$.} We have $|\nabla_w\tilde a|\le|D|^2| \mathcal A|\,|\nabla\ab|\le\jap{v_0}^{|\gamma|}\rho\,C_\gamma(1+K)\le
C_\gamma(1+K)$. For $\tilde c$, if $\gamma=-3$, then $\tilde c=8\pi\rho^2\tilde F$. If $\gamma>-3$,
then $-\cb[F]=2(\gamma+3)|\cdot|^\gamma*F$. With $\chi\in C^\infty_c(B_2)$, $0\le\chi\le1$, $\chi=1$
on $B_1$, $2(\gamma+3)|z|^\gamma\chi$ has $L^1$ norm at most $16\pi$ and
$\nabla(|z|^\gamma(1-\chi))$ is bounded, which implies
$[\cb(\sigma)]_{C^\alpha}\le16\pi[F(\sigma)]_{C^\alpha}+C\|F(\sigma)\|_{L^1}$. In both cases a bound
on $[F(\sigma)]_{C^\alpha(\R^3)}$ suffices.

We apply the Krylov--Safonov estimate \cite[Thm.~7.37, Cor.~7.41]{Lie96} to $\tilde F$, whose
equation has uniformly elliptic coefficients, with the bounded zero-order term $\tilde c\tilde F$
moved to the right-hand side (\cite[(7.57)]{Lie96} controls the oscillation by that of $\tilde F$
and by $\|\tilde c\tilde F\|_{L^{4}}$, the dimension being $3$). Then we get $[\tilde
F]_{C^\alpha(\mathcal Q_{1/2})}\le C\|\tilde F\|_{L^\infty(\mathcal Q)}$ with
$\alpha=\alpha(\gamma,K,A)$.

Back in $v$, since $v_0+\mathcal AB_{1/2}\supset B_{r_0}(v_0)$ with
$r_0:=\frac12\rho\jap{v_0}^{-3/2}$ and $|\mathcal A^{-1}|\le\rho^{-1}\jap{v_0}^{3/2}$, we have
$[F(\sigma)]_{C^\alpha(B_{r_0}(v_0))}\le C|\mathcal A^{-1}|^\alpha C_Ge^{-\beta(|v_0|-1)_+^2}\le C$
for every $\sigma\in[T/4,T]$. If $|v-v'|$ exceeds both $r_0(v)$ and $r_0(v')$, then $|F(v)-F(v')|\le
F(v)+F(v')\le C(r_0(v)^\alpha+r_0(v')^\alpha)\le C|v-v'|^\alpha$, by the Gaussian bound. Therefore
$[F(\sigma)]_{C^\alpha(\R^3)}\le C(\gamma,\beta,T,C_G)$ for $\sigma\in[T/4,T]$.

Hence, for $\sigma_0\in[T/2,T]$, for which the times in $\mathcal Q$ lie in $[T/4,T]$, we have
$[\tilde c(s,\cdot)]_{C^\alpha(B_1)}\le\rho^2|\mathcal A|^\alpha[\cb(\sigma)]_{C^\alpha}\le C$.

\emph{Schauder estimate.} For $\sigma_0\in[T/2,T]$ the coefficients $\tilde a,\tilde c$ of the
equation of $\tilde F$ are measurable in $s$, uniformly elliptic and bounded, and uniformly H\"older
continuous in $w$, with constants independent of $(\sigma_0,v_0)$. The interior Schauder estimate
\cite[Thm.~1]{Kn80} implies $|\nabla_w\tilde F(0,0)|+|\nabla^2_w\tilde F(0,0)|\le C\|\tilde
F\|_{L^\infty(\mathcal Q)}$. Note that it requires of the coefficients, including the zero-order
one, only H\"older continuity in the space variables and boundedness, with no regularity in time,
and its constant depends only on the ellipticity, these norms, $\alpha$ and the dimension (applied
on $(-1,0)\times B_1$, the bound extends to $s=0$ by continuity). Since $\nabla_vF= \mathcal A^{-1}\nabla_w\tilde
F$, $|\mathcal A^{-1}|\le\rho^{-1}\jap{v_0}^{3/2}\le C(T)\jap{v_0}^{11/2}$ and $\|\tilde
F\|_{L^\infty(\mathcal Q)}\le C_Ge^{\beta}e^{-\beta|v_0|^2+2\beta|v_0|}$, we conclude that
\begin{equation*}
|\nabla F(\sigma_0,v_0)|+|\nabla^2F(\sigma_0,v_0)|\le
C(T)\jap{v_0}^{11}e^{2\beta|v_0|}e^{-\beta|v_0|^2}\le C_\eps e^{-(\beta-\eps)|v_0|^2}.\qquad\qedhere
\end{equation*}
\end{proof}

\begin{lem}[The solution of (GS1) near $\sigma=0$]\label{lem:GS1q}
Let $\gamma\in[-3,-2]$, and let $F_0\ge0$ be continuous and normalized, with $F_0\le
C_{\rm dat}e^{-\beta|v|^2}$. Let $F$ be the solution of (GS1). Then $F$ is continuous on
$[0,\infty)\times\R^3$, with $F(0)=F_0$, and for every $T>0$ there is $C_T$ with $F\le
C_Te^{-\beta_M|v|^2}$ on $[0,T]\times\R^3$, where $\beta_M=\min\{\beta,\beta_0\}$ and
$\beta_0=\beta_0(\gamma)$ is the exponent of Proposition~\ref{prop:GS4}.
\end{lem}

\begin{proof}
We follow the proof of \cite[Thm.~1.2]{GS} in \cite[\S11]{GS}, with \cite[Thm.~2.3]{GS} replaced by
Proposition~\ref{prop:GS4}, and with the finiteness of the Fisher information obtained from
Proposition~\ref{prop:GS5} and Lemma~\ref{lem:FPgl}(ii). We apply the propositions to a classical
solution that is continuous on $[0,T]\times\R^3$ and smooth for $\sigma>0$, with the stated bounds,
and their proofs use no other property of the solution.

\emph{Local solution.} Since $F_0\in L^\infty_k$ for every $k$, and $F_0$ is continuous, positive
somewhere and independent of $x$ (so that the lower bound condition of \cite[Def.~1.1]{HST20}
holds), the local theory quoted in (GS1) yields a classical solution $f$ on a maximal interval
$[0,T_E)$, continuous on $[0,T_E)\times\R^3$ with $f(0)=F_0$ (\cite[Prop.~3.1]{HST20}, as explained
in (GS1)), positive and smooth for $\sigma>0$, which can be continued as long as
$\|f(\sigma)\|_{L^\infty_{k_0}}$ stays bounded, and with
$\|f\|_{L^\infty_k([0,\delta]\times\R^3)}\le K_1$ for some $\delta>0$, with $k>k_0$ the exponent for
which \cite[Thm.~1.2]{HST20} is applied ($\delta$ and $K_1$ depend on it) \cite[Thm.~1.2(i)]{HST20}.
In particular, $\|f(\sigma)\|_{L^\infty}\le K_1$ on $[0,\delta]$. By \cite[Thm.~1.2(iii)]{HST20},
since $F_0\in L^\infty_k$ for every $k$, all derivatives of $f$ decay faster than any power of $\jap
v$ on compact subsets of $(0,T_E)$. It has unit mass and energy $3$ (conservation for $\sigma>0$ by
integration by parts, and continuity at $\sigma=0$ by dominated convergence), and, by (GS2),
$\int\jap v^kf(\sigma)\le A_k(T)$ on $[0,\min\{T,T_E\})$, for every $k$.

\emph{A Gaussian bound on $[0,\delta]$, and the Fisher information at $\delta$.} By the bound (GS4)
on $[0,\delta]$, with $K=K_1$ and $\beta_M$, we have $f\le C_1e^{-\beta_M|v|^2}$ there. The bound
(GS5) on $[0,\delta]$ then yields Gaussian bounds for $\nabla f$ and $\nabla^2f$ on
$[\delta/2,\delta]$, and by Glaeser's bound we get $i(f(\delta))<\infty$.

\emph{$L^\infty$ and Gaussian bounds up to $T_E$.} Let $T<\infty$, $T\le T_E$. On $[\delta,T)$ the
solution is classical, with bounds on all its derivatives and $f\le C(\sigma)\jap v^{-k_0}$, which
is how \cite[\S11]{GS} apply \cite[Thm.~1.1]{GS} (on $[\delta,T']$, $T'<T$,
\cite[Thm.~1.2(iii)]{HST20}, as in the paragraph \emph{Local solution}, gives decay faster than any
power of $\jap v$ of every derivative). Thus $i(f(\sigma))\le i(f(\delta))$ there, and
$\|f(\sigma)\|_{L^3}\le C$ by (EL5). With the moments of (GS2), by interpolation we get
$\|f(\sigma)\|_{L^p_k}\le C(T)$ for $p\in(3/(5+\gamma),3)$ and every $k$. Applying (GS3) on
$[\delta/2,T)$ (on $[\delta/2,\delta]$ the $L^p_k$ bound follows from $f\le K_1$ and the moments),
we obtain $\|f\|_{L^\infty}\le K_2(T)$ on $[\delta,T)$. By the bound (GS4) on $[0,T')$, $T'<T$, with
$K=\max\{K_1,K_2(T)\}$, we have $f\le C_Te^{-\beta_M|v|^2}$ on $[0,T)$, with $C_T$ independent of
$T'$. Thus $\|f(\sigma)\|_{L^\infty_{k_0}}$ is bounded on $[0,T)$ for every finite $T\le T_E$, and
the continuation criterion implies $T_E=\infty$.

\emph{Identification.} Thus $f$ is a global classical solution with data $F_0$, continuous on
$[0,\infty)\times\R^3$, and bounded by $C_Te^{-\beta_M|v|^2}$ on every $[0,T]$. By the uniqueness in
\cite[Thm.~1.2]{GS}, we have $f=F$.
\end{proof}

\section{Auxiliary estimates}\label{app:auxest}\label{app:aux}

The statements below are used in the body of the paper, and each is proved here or at the place
indicated.

\begin{lem}[Symmetries]\label{lem:sym}
\begin{enumerate}[label=(\roman*)]
\item $Q(T_pF,T_pG)=\lambda_p\,T_pQ(F,G)$. Hence, if $H$ is a classical solution, so is
$\sigma\mapsto T_pH(\lambda_p\sigma)$, and by uniqueness \cite[Thm.~1.2]{GS},
$\Phi_\sigma(T_pu)=T_p\Phi_{\lambda_p\sigma}(u)$.
\item If $\mom(u)=\mom(\mu)$, then
\begin{equation}\label{eq:momTp}
\mom(T_pu)=\mathbf m(p):=e^{\alpha+3\beta}\big(1,\,a,\,|a|^2+3e^{2\beta}\big),
\end{equation}
independent of $u$. The map $\mathbf m$ is a diffeomorphism from a neighborhood of $p=0$ onto a
neighborhood of $\mom(\mu)$, and $|p|\le C|\mathbf m(p)-\mom(\mu)|$.
\item $\frac{d}{d\eps}T_{\eps p}F=T_{\eps p}\,\mathcal G^{(\eps)}_pF$, where $
\mathcal G^{(\eps)}_p=\alpha-\beta\,v\cdot\nabla-e^{-\eps\beta}a\cdot\nabla$.
\end{enumerate}
\end{lem}

\begin{proof}
(i) Note that $Q$ commutes with translations. For $(S_{A,B}F)(v)=AF(v/B)$, by the homogeneity of $a$
(degree $2+\gamma$) and of $c$ (degree $\gamma$, or $\delta_0$ of degree $-3$ when $\gamma=-3$), we
have $\ab[S_{A,B}F]=AB^{5+\gamma}\ab[F](\cdot/B)$, $\cb[S_{A,B}F]=AB^{3+\gamma}\cb[F](\cdot/B)$ and
$\nabla^2S_{A,B}G=AB^{-2}(\nabla^2G)(\cdot/B)$. Hence
$Q(S_{A,B}F,S_{A,B}G)=AB^{3+\gamma}S_{A,B}Q(F,G)$. We take $A=e^\alpha$, $B=e^\beta$.

(ii) We change variables $v=a+e^\beta z$. The differential at $p=0$ maps
$\partial_{a_i},\partial_\alpha,\partial_\beta$ to $(0,e_i,0)$, $(1,0,3)$, $(3,0,15)$, which are
linearly independent.

(iii) This follows by a direct computation, using $v=a\eps+e^{\eps\beta}w$ with $w$ the argument of
$F$.
\end{proof}

\begin{lem}[Local energies do not lose derivatives]\label{lem:localnoloss}
Let $w(t,y,v)>0$ be any weight, and let $\partial_th=\kappa(\mathcal L_t[g]h+S)$. Then
\begin{equation}
\frac{d}{dt}\iint w^2h^2\dd y\dd v=\iint(\partial_tw^2)\,h^2\dd y\dd v+2\kappa\iint w^2h\,(
\mathcal L_t[g]h+S)\dd y\dd v .
\end{equation}
\end{lem}

\begin{proof}
In the profile variables there is no transport term.
\end{proof}

\begin{lem}[Weighted convolution algebra]\label{lem:algebra}
Let $p\ge2$, $\rad>0$ and $w_n=\rad^n(n+1)^p/n!$. For nonnegative families
$(a_\beta),(b_\beta),(d_\beta)$ indexed by $\beta\in\mathbb N^3$, put
$A_k=(\sum_{|\beta|=k}a_\beta^2)^{1/2}$, and define $B_k,D_k$ similarly. Then
\begin{equation}
\sum_nw_n\!\sum_{|\beta|=n}d_\beta\!\sum_{\beta'\le\beta}\tbinom{\beta}{\beta'}a_{\beta'}b_{\beta-\beta'}
\le C_p\!\sum_n\sum_{k\le n}\varrho_{n,k}\,(w_kA_k)(w_{n-k}B_{n-k})D_n,
\end{equation}
where $\varrho_{n,k}=(\min(k,n-k)+1)^{1-p}$ and $C_p=2^p$. In particular, the right-hand side is
bounded by $2^{p+1}\zeta(2p-2)^{1/2}\,\|wA\|_{\ell^2}\|wB\|_{\ell^2}\|D\|_{\ell^2}$.
\end{lem}

\begin{proof}
We use three facts. First, by Vandermonde's identity, we have $\binom\beta{\beta'}\le\binom{n}{k}$
for $|\beta'|=k$. Second, $w_n\binom nk=w_kw_{n-k}\big(\frac{n+1}{(k+1)(n-k+1)}\big)^p$, and
$(k+1)(n-k+1)\ge\frac{n+1}2(\min(k,n-k)+1)$. Third, for fixed $n$ and $k$, by Cauchy--Schwarz in
$(\beta',\beta'')$, we have
\begin{equation}
\sum_{|\beta'|=k,\,|\beta''|=n-k}a_{\beta'}b_{\beta''}d_{\beta'+\beta''}\le
A_kB_{n-k}\Big(\sum_{|\beta|=n}N(\beta,k)\,d_\beta^2\Big)^{1/2},
\end{equation}
where $N(\beta,k)=\#\{\beta'\le\beta:|\beta'|=k\}\le(\min(k,n-k)+1)^2$. These three facts prove the
first inequality. For the second, we split the sum into $k\le n/2$ and $k>n/2$. In the first part,
we use $\sum_k(k+1)^{1-p}(w_kA_k)\le\zeta(2p-2)^{1/2}\|wA\|_{\ell^2}$ together with
$\sum_n(w_{n-k}B_{n-k})D_n\le\|wB\|\|D\|$. The second part is symmetric.
\end{proof}

\begin{lem}[Gaussian decay of derivatives by local interpolation]\label{lem:RGinterp}
Let $0\le F\le C_Me^{-\beta_M|z|^2}$ and $\sum_{|\alpha|\le2N}\sup|\nabla^\alpha F|\le C_{2N}$. Then
for $j\le N$,
\begin{equation}
|\nabla^jF(z)|\le C(C_M,C_{2N},N,\beta_M)\,e^{-\beta_M|z|^2/4}.
\end{equation}
In particular, for $q<\beta_M/8$, with $C$ depending also on $\beta_M/8-q$ ,
\begin{equation}
\sum_{j\le N}\big\|e^{2q|z|^2}\nabla^jF\big\|_{L^2\cap L^\infty}\le
C(C_M,C_{2N},\beta_M,\beta_M/8-q).
\end{equation}
\end{lem}

\begin{proof}
On $B_1(z)$, by the Gagliardo--Nirenberg inequality we have
\begin{equation}
\|\nabla^ju\|_{L^\infty(B_1)}\le C\|u\|^{1-j/(2N)}_{L^\infty(B_2)}\|u\|^{j/(2N)}_{C^{2N}(B_2)}+C\|u\|_{L^\infty(B_2)}.
\end{equation}
Here $\|F\|_{L^\infty(B_2(z))}\le C_Me^{-\beta_M(|z|-2)_+^2}$, and $1-j/(2N)\ge1/2$. Moreover,
$e^{-\beta_M(|z|-2)_+^2/2}\le Ce^{-\beta_M|z|^2/4}$, which completes the proof.
\end{proof}

\begin{lem}[Ellipticity from mass, moments and an $L^\infty$ bound]\label{lem:RGellip}
Let $G\ge0$ with $\int G=M\ge M_->0$, $\int|z|^2G\le ME_0$ and $\|G\|_{L^\infty}\le K$. Then
\begin{equation}
\ab[G](z)\zeta\cdot\zeta\ge c_G\big(\jap z^\gamma|P_z\zeta|^2+\jap
z^{\gamma+2}|(I-P_z)\zeta|^2\big),\qquad c_G=c(\gamma,M_-,E_{\rm mom},K)>0 .
\end{equation}
\end{lem}

\begin{proof}
We may assume $E_{\rm mom}\ge\frac14$. Let $R_0:=2\sqrt{E_{\rm mom}}\ge1$. By Chebyshev's
inequality, we have $\int_{B_{R_0}}G\ge\frac34M$. For a line $\ell$ and $\epsilon>0$, the tube
$\{w\in B_{R_0}:\mathrm{dist}(w,\ell)<\epsilon\}$ has volume $\le2\pi\epsilon^2R_0$, hence $G$-mass
$\le2\pi K\epsilon^2R_0$. We fix $\epsilon:=(M_-/(8\pi KR_0))^{1/2}$, so that this mass is $\le
M/4$. Then every tube leaves $G$-mass $\ge M/2$ in $B_{R_0}$ outside it.

\emph{The kernel.} Since $a(u)=|u|^{\gamma+2}(I-\hat u\otimes\hat u)$, we have
$a(u)\zeta\cdot\zeta=|u|^\gamma|u\times\zeta|^2$. With $u=z-w$, we have
$|(z-w)\times\zeta|=|\zeta|\,\mathrm{dist}(w,\ell_\zeta)$, where $\ell_\zeta:=z+\R\zeta$. On
$B_{R_0}$, $|z-w|\le(1+R_0)\jap z$, so $|z-w|^\gamma\ge c_0\jap z^\gamma$ with $c_0=(1+R_0)^\gamma$.
Hence
\begin{equation}
\ab[G](z)\zeta\cdot\zeta\ge c_0\jap z^\gamma|\zeta|^2\int_{B_{R_0}}\mathrm{dist}(w,\ell_\zeta)^2G(w)\dd
w .
\end{equation}

\emph{(a) Every $\zeta$.} Outside the tube of radius $\epsilon$ around $\ell_\zeta$ the distance is
$\ge\epsilon$. Thus we have $\ab[G](z)\zeta\cdot\zeta\ge c_1\jap z^\gamma|\zeta|^2$, with
$c_1=c_0\epsilon^2M_-/2$.

\emph{(b) $|z|\,|(I-P_z)\zeta|\ge2R_0|\zeta|$.} The distance from the origin to $\ell_\zeta$ is
$|z\times\zeta|/|\zeta|=|z||(I-P_z)\zeta|/|\zeta|\ge2R_0$. Hence for $w\in B_{R_0}$, we have
$\mathrm{dist}(w,\ell_\zeta)\ge|z||(I-P_z)\zeta|/(2|\zeta|)$, and $\ab[G](z)\zeta\cdot\zeta\ge
c_0\jap z^\gamma\frac{3M_-}{16}|z|^2|(I-P_z)\zeta|^2$. Here $|z|\ge2R_0\ge1$, so
$|z|^2\ge\frac12\jap z^2$, and then $\ab[G](z)\zeta\cdot\zeta\ge c_2\jap
z^{\gamma+2}|(I-P_z)\zeta|^2$. Averaging with (a), and using $|P_z\zeta|\le|\zeta|$, we obtain
\begin{equation}
\ab[G](z)\zeta\cdot\zeta\ge\tfrac12\min(c_1,c_2)\big(\jap z^\gamma|P_z\zeta|^2+\jap
z^{\gamma+2}|(I-P_z)\zeta|^2\big).
\end{equation}

\emph{(c) $|z|\,|(I-P_z)\zeta|<2R_0|\zeta|$.} Then $\jap z^2|(I-P_z)\zeta|^2\le(1+4R_0^2)|\zeta|^2$,
so $\jap z^{\gamma+2}|(I-P_z)\zeta|^2\le(1+4R_0^2)\jap z^\gamma|\zeta|^2$, and by (a) we get
\begin{equation}
\ab[G](z)\zeta\cdot\zeta\ge\frac{c_1}2\jap z^\gamma|P_z\zeta|^2+\frac{c_1}{2(1+4R_0^2)}\jap
z^{\gamma+2}|(I-P_z)\zeta|^2 .
\end{equation}
The constants depend only on $(\gamma,M_-,E_{\rm mom},K)$.
\end{proof}

\begin{prop}[Quantitative smoothness of homogeneous trajectories]\label{prop:RS}
Homogeneous trajectories starting in $\mathfrak K$ satisfy the following.
\begin{itemize}
\item[(RM)] If $F(0)\in\mathfrak K^N(C_{\rm dat},\beta,A_N)$, $N\ge2$, and $S>0$, then
\begin{equation}
F(\sigma)\le C_M(S)e^{-\beta_M|z|^2}\qquad\text{for }0\le\sigma\le S ,
\end{equation}
with $\beta_M=\min\{\beta,\beta_0(\gamma)\}$, $\beta_0(\gamma)=c_\perp/(2C_\parallel)$ of
Proposition~\ref{prop:GS4}, and $C_M(S)=C_M(\gamma,S,C_{\rm dat},\beta,A_N)$.

\item[(RS)] Let $q_S:=\min\{\beta_M/32,1/48\}$. For $F(0)\in\mathfrak K^N$, $N\ge2$, and every
$S>0$,
\begin{equation}
\sup_{\sigma\le S}\sum_{|\alpha|\le N}\sup_ze^{(q_S/4)|z|^2}|\nabla^\alpha F(\sigma,z)|\le
C_S=C(N,S,C_{\rm dat},\beta,A_N).
\end{equation}
\end{itemize}
\end{prop}
The proof, from (GS1)--(GS5) of (E2) of Inputs~\ref{hyp:E}, is given in Appendix~\ref{app:hom}. Part
(RS) is stated in the Gaussian form that the proof yields, namely data in $\mathfrak K^N$, with
$2N+4$ polynomially weighted derivatives, have Gaussian bounds on $N$ derivatives along the
trajectory, as explained in Remark~\ref{rem:RSshortF}(2).

\begin{lem}[Gaussian Sobolev bounds on bounded intervals]\label{lem:RSmed} 
Let $F_0\in\mathfrak K^N$, $S>0$, $4\le m\le N+2$ and $0<q\le q_S$, and let
$E^\natural_m,D^\natural_m$ be the ladder energies defined before Lemma~\ref{lem:RSshortF}. Then
\begin{equation}
\sup_{[0,S]}E^\natural_m+\int_0^SD^\natural_m\le K_m=K_m(m,S,C_{\rm dat},\beta,A_N,q).
\end{equation}
\end{lem}
The proof, given in the same appendix, uses Lemma~\ref{lem:RSshortF}, (RM), (GS5) and the
weighted $L^2$ bound of \cite[Lemma~4.2(i)]{CM}.

\begin{lem}[Static Gaussian weights on leaves]\label{lem:Gleafweight}
Let $q>0$. Let
$(\rho_\infty(\xi),y_{0\infty}(\xi),\theta_\infty(\xi)):=(\eps\rho_0(\xi),y_{0,\infty}(\xi),\theta_\infty(\xi))$
be the limits of the leaf parameters of $f_0$ under free transport (for
product data $y_{0\infty}$ and $\theta_\infty$ do not depend on $\xi$), let
$A_\infty:=\rho_\infty^{-1}\theta_\infty^{3/4}$, that is,
$A_\infty(v)=A_{(\rho_\infty(v),y_{0,\infty}(v),\theta_\infty(v))}$ in the notation of
Section~\ref{sec:weights}, and $W(y,v):=A_\infty(v)e^{q|y|^2}$. For $t\ge1$, the normalized leaf
weight of $W$ in the frame $(\rho_\infty(\xi),y_{0\infty}(\xi),\theta_\infty(\xi))$ is
\begin{equation}
\hat m_x(\hat z)=e^{q_\xi|\hat z|^2+\omega_x(\hat z)},\qquad q_\xi:=q\,\theta_\infty(\xi),
\end{equation}
where
\begin{equation}\label{eq:omegax}
\omega_x=2q\theta_\infty^{1/2}y_{0\infty}\cdot\hat
z+q|y_{0\infty}|^2+\log\big(A_\infty(\xi-z/t)/A_\infty(\xi)\big),
\end{equation}
with coefficients evaluated at $\xi$. Moreover
\begin{equation}
|\nabla\omega_x|\le2\theta_+^{1/2}L\,(q+t^{-1}),\qquad|\nabla^2\omega_x|\le CLt^{-2},\qquad
e^{\pm\omega_x}\le C\jap{\hat z}^{p}e^{Cq\jap{\hat z}} ,
\end{equation}
where $p\ge4$ is the exponent of (D3).
\end{lem}

\begin{proof}
On the leaf through $x$ we have $y=z=y_{0\infty}(\xi)+\theta_\infty(\xi)^{1/2}\hat z$, and $\hat
m_x=m_x/A_\infty(\xi)$ with $m_x(z)=W(z,(x-z)/t)$. We expand $q|z|^2$. Then
\begin{equation}
\nabla_{\hat z}\log A_\infty(\xi-z/t)=-t^{-1}\theta_\infty^{1/2}(\nabla\log A_\infty)(v),
\end{equation}
and we have $|\nabla\log A_\infty|\le L/2\le2L$ and $|\nabla^2\log A_\infty|\le CL$, by the
convention on $L$ of Subsection~\ref{ssec:disp} (for product data $\nabla\log
A_\infty=-\nabla\log\rho_0$, since $\theta_\infty$ is constant). For the last bound, by (D3) and
Peetre's inequality $\jap{a+b}^s\le2^{|s|/2}\jap a^s\jap b^{|s|}$ with $s=p$, the ratio
$A_\infty(\xi-z/t)/A_\infty(\xi)=\big(\rho_0(\xi)/\rho_0(\xi-z/t)\big)\big(\theta_\infty(\xi-z/t)/\theta_\infty(\xi)\big)^{3/4}$
and its inverse are $\le A_F^22^{p/2}\jap{z/t}^p(\theta_+/(4\theta_-))^{3/4}\le C\jap{\hat z}^p$,
with $C=C(A_F,p,L,\theta_\pm)$, since $t\ge1$ and $|z|\le L/2+\theta_+^{1/2}|\hat z|$ by the same
convention. The other two terms of $\omega_x$ are bounded by $2q\theta_+^{1/2}L|\hat z|$ and
$qL^2/4$ in absolute value.
\end{proof}

\begin{cor}[Gaussian entry along the leaves]\label{cor:Aleafentry}
Assume on $[T_0,t]$ the hypotheses (U1)--(U4) of Theorem~\ref{thm:CUboot}, so that, by
Lemmas~\ref{lem:BTsob} and~\ref{lem:BTU4}, the hypotheses of Theorem~\ref{thm:RGF} hold. Assume also
the requirements on $T_0$ of Remark~\ref{rem:RGconst}, $C(L_1,R_1)\bar
K(\delta_U+T_0^{-1})\le\delta_1/(3C_e)$, where $\bar K$ is the frame constant of
Remark~\ref{rem:BTnorm}, and $C(L,C'_D,\beta,\beta',\theta_\pm)(1+S_\sharp)L/T_0\le\delta_1/(6C_e)$.
If $s_\xi(t)\in[S_*-\frac12,3S_\sharp]$, then
\begin{equation}
\|f_\xi(t)\|_{X(\hat m_\xi)}\le\delta_1/C_e .
\end{equation}
\end{cor}

\begin{proof}
We compare the leaf with the fibered homogeneous flow, and that flow with $\mu$. Let
$\delta'=\delta_1/(3C_e)$.

\emph{Leaf against fibered flow.} By Theorem~\ref{thm:RGF}(i) with $j=0$, we know that on
$\{s\le3S_\sharp\}$ the normalized leaf is within $C_T\varepsilon_T(t)e^{-(q/4)|z|^2}$ of the
fibered flow sampled along the leaf, and $C_T\varepsilon_T(T_0)\le\delta'$ by the requirements on
$T_0$ assumed in the statement.

\emph{Fibered flow against $\mu$.} The leaf samples the fibers $v=v(z)=\xi-z/t$, one fiber for each
$z$, and the sampled fibered flow is $z\mapsto\rho_\xi^{-1}T_{P_{v}}F_{v}(s_{v})(z)$. A bound in
$L^2$ for each fixed $v$ does not control this function, so we first prove a pointwise bound. Let
$F$ be a normalized trajectory with $F(0)\in\mathfrak K^N$, let
$\sigma\ge\max\{S_\flat(\delta''),\frac12\}$, $f:=F(\sigma)-\mu$ and $q':=\beta'/4$. By
Lemma~\ref{lem:Aentry} with $q'$, we have $\|e^{q'|w|^2}f\|_{L^2}\le\delta''$. By
Lemma~\ref{lem:Abarrier}(iii), $F\le C_\flat e^{-\beta'|w|^2}$ at all times, so applying
Proposition~\ref{prop:GS5} on $[0,\sigma]$, with $\eps=\beta'/2$, we get $|\nabla F(\sigma)|\le
C_ge^{-(\beta'/2)|w|^2}$. Its proof runs with $\rho:=\min\{\jap{v_0}^{-4},\frac14\sqrt{1/2}\}$ for
every $\sigma\ge\frac12$, so $C_g=C(\gamma,\beta',C_\flat)$ does not depend on $\sigma$. The same
bound holds for $\nabla\mu$, since $\beta'\le\frac18$. If $|f(w_0)|=a$ and
$b:=\sup_{B_1(w_0)}|\nabla f|$, then $|f|\ge a/2$ on the ball of radius $\min\{1,a/(2b)\}$ about
$w_0$, so
\begin{equation}
a\le C\max\big\{\|f\|_{L^2(B_1(w_0))},\ \|f\|_{L^2(B_1(w_0))}^{2/5}\,b^{3/5}\big\} .
\end{equation}
With $\|f\|_{L^2(B_1(w_0))}\le e^{-q'(|w_0|-1)_+^2}\delta''$,
$b\le2C_ge^{-(\beta'/2)(|w_0|-1)_+^2}$, $\delta''\le1$ and
$\frac25q'+\frac35\cdot\frac{\beta'}2\ge\frac{\beta'}4$, we obtain
\begin{equation}
|F(\sigma,w)-\mu(w)|\le C\,\delta''^{2/5}e^{-(\beta'/4)(|w|-1)_+^2}\qquad\text{for all }w .
\end{equation}
On $|z|\le t/(15LS_\sharp)$ we have $|\log(\lambda_v/\lambda_{\rm
ref}(\xi))|\le\frac52L|z|/t\le1/(6S_\sharp)$, so $|s_v-s_\xi|\le3S_\sharp(e^{1/(6S_\sharp)}-1)\le1$.
Hence $s_v\ge S_*-\frac32=S_\flat(\delta'')+\frac12$, by the choice of $S_*$ in
Subsection~\ref{ssec:cells}, and the bound applies to $F_v(s_v)$. We write
\begin{equation}
\rho_\xi^{-1}\big(T_{P_v}F_v(s_v)-T_{P_\xi}\mu\big)=\rho_\xi^{-1}T_{P_v}\big(F_v(s_v)-\mu\big)+\rho_\xi^{-1}\big(T_{P_v}-T_{P_\xi}\big)\mu
.
\end{equation}
In the first term $\rho_v/\rho_\xi\le e$, and by (H1) the frame $P_v$ has $|y_{0,v}|\le L$ and
$\theta_v\in[\theta_-,\theta_+]$. With $\hat z_v$, $\hat z$ the normalized variables of $P_v$ and
$P_\xi$, we have $(|\hat z_v|-1)_+^2\ge|z|^2/(2\theta_+)-C$ and $|z|^2\ge\frac12\theta_-|\hat
z|^2-C$. Therefore the first term is $\le C\delta''^{2/5}e^{-c_\beta|\hat z|^2}$, with
$c_\beta:=\beta'\theta_-/(16\theta_+)$, and its $X(\hat m)$-norm is $\le C_{\rm
pt}\,\delta''^{2/5}\le\delta'/2$, since $q_0\theta_+\le c_\beta/4$ by Remark~\ref{rem:RGfull}(b).
Here $C_{\rm pt}$ depends only on $(\gamma,C_{\rm dat},\beta,A_N,k_1,L,\theta_\pm)$ and on $C_W$,
all fixed before $\delta_1$, and $\delta''\le1$ is fixed with $C_{\rm pt}\delta''^{2/5}\le\delta'/2$
in the same subsection. In the second term, $|P_v-P_\xi|\le C(L,C'_D,\beta,\theta_\pm)|z|/t$, by
(H1) and the bound on $\partial_vg(T_0)$ of (LWP)(a3), as in (F3), and so its $X(\hat m)$-norm is
$\le CL/T_0$. The rest of the leaf, $|z|>t/(15LS_\sharp)$, where both the fibered flow, by the
Gaussian bound above, and $\mu$ are Gaussian in $z$, contributes $\le Ce^{-c\,t^2/(LS_\sharp)^2}\le
CLS_\sharp/T_0$ to the $X(\hat m)$-norm. By the hypothesis on $T_0$, the last two are
$\le\delta'/2$. Thus this step contributes $\le\delta'$.

\emph{Frames.} The two previous steps compare the leaf with $M_{P_\xi}=T_{P_\xi}\mu$, while
$f_\xi=T^{-1}_p\big(g^{(x)}-M_p\big)$ is centered at the leaf parameters $p:=p_g(t,\xi)$. We write
$f_\xi=T^{-1}_p\big(g^{(x)}-M_{P_\xi}\big)+T^{-1}_p\big(M_{P_\xi}-M_p\big)$.

No derivative of the leaf is needed for this change of frame, because in the variable $z$ the norm
depends on the frame only through a scalar factor. Indeed, $\|T^{-1}_pu\|_{X(\hat
m_\xi)}=(A_p/A_{\rm ref}(\xi))\,\|m_xu\|_{L^2(\dd z)}$ for every $u$, by the change of variables
$z=y_0+\theta^{1/2}\hat z$, as in the definition of $\ell_\xi$. Thus the first term is the quantity
bounded in the two previous steps, up to the factors $\rho_\xi/\rho_g\le2$ and $A_p/A_{\rm ref}\le
C(L,\theta_\pm)$ given by (F3), which can be put into $C_{\rm pt}$ and into the requirement on
$T_0$.

In the second term, $\partial_pM_p$ is a polynomial of degree $\le2$ in $z$ times $M_p$, at the
Gaussian rate $1/(2\theta)\ge1/(2\theta_++2)>4q_0$, uniformly for $p$ on the segment
$[P_\xi,p_g(t,\xi)]$, whose temperatures lie in $[\frac12\theta_-,\theta_++1]$ and centers in
$\{|y_0|\le L+1\}$ by (H1) and (F3). Hence it is $\le C(L,\theta_\pm)|p_g-P_\xi|$ (relatively for
the mass).

The bound (G4) is not suited to this step, since its weight $\jap{\hat z}\hat m_\xi$ is the full
energy weight, with no margin in the Gaussian rate, whereas the generator
$\alpha-\beta\,v\cdot\nabla-a\cdot\nabla$ of Lemma~\ref{lem:sym}(iii) changes the weight along the
path of frames. Since $\bar K\ge C(L,\theta_\pm)$, the second term is $\le C\bar K|p_g-P_\xi|$.

By (F3) in Subsection~\ref{ssec:inputs}, a consequence of (U1) alone, we have $|p_g-P_\xi|\le
C(L_1,R_1,C'_D,\beta,\theta_\pm)(\delta_U+T_0^{-1})$, so this is $\le\delta'$. The three
contributions add up to $\le3\delta'=\delta_1/C_e$.
\end{proof}

\begin{lem}[Leaves of free transport]\label{lem:LWfree}
Let $g_{\rm fr}(t,y,v)=f_0(y,v)$ be the free profile, $f_0\in\mathfrak D^\sharp$, and let $t\ge1$.
Then:
\begin{enumerate}[label=(\roman*)]
\item every normalized leaf of $g_{\rm fr}$ satisfies (G1)--(G3) with constants depending only on
$\mathfrak D$;
\item at radius $\rad\le\Lambda/8$ we have $\hat{\mathbb A}_\xi\le C_DB$ and $\hat{
\mathbb A}_{\ge1,\xi}\le C_D\big((\rad/\Lambda)^2+\delta_y^2\big)B$;
\item the parameters of the leaf $\xi$ are those of $F_0(\cdot,\xi)$, that is,
$(y_{0,\infty}(\xi),\theta_\infty(\xi))$, up to $O(L/t)$, and its mass is
$\eps\rho_0(\xi)(1+O(L/t))$.
\end{enumerate}
\end{lem}

\begin{proof}
The leaf through $x=t\xi$ is $z\mapsto f_0(z,\xi-z/t)=\eps\rho_0(\xi-z/t)F_0(z,\xi-z/t)$. By (D3)
and (D4), $|\nabla_v\log(\rho_0F_0)(z,v)|\le L$ for real $(z,v)$ (with $L/2$ under the convention of
Subsection~\ref{ssec:disp}), so the ratio $f_0(z,\xi-z/t)/\big(\eps\rho_0(\xi)F_0(z,\xi)\big)$ lies
in $[e^{-L|z|/t},e^{L|z|/t}]$, which means that
\begin{equation}
\big|\phi^{\rm fr}(t,z,\xi)-F_0(z,\xi)\big|\le\big(e^{L|z|/t}-1\big)F_0(z,\xi),\qquad\phi^{\rm
fr}=\frac{f_0(z,\xi-z/t)}{\eps\rho_0(\xi)},
\end{equation}
and the Gaussian bound $F_0\le C_Fe^{-\beta_F|z|^2}$ of (D1), which is uniform in $\xi$, absorbs
this factor. This proves (iii), and also (i). Indeed, the mass, mean and temperature of
$F_0(\cdot,\xi)$ are $1$, $y_{0,\infty}(\xi)$ and $\theta_\infty(\xi)$, which are bounded above and
below uniformly in $\xi$, and the bounds of (G2)--(G3) follow from (D1) and, for the lower bound of
the coefficient, from the mass of $F_0(\cdot,\xi)$ in a fixed ball, which is bounded below uniformly
in $\xi$ by $\int F_0=1$ and the Gaussian bound.

For (ii), the components are $(\partial^\alpha_y\partial^\beta_vf_0)(z,\xi-z/t)$. We use the Cauchy
estimates for $\zeta\mapsto\rho_0(v+\zeta)\partial^\alpha_yF_0(z,v+\zeta)$ on the polydiscs
$\Delta(0,\Lambda/2)\subset\mathbb S_\Lambda$, with $A_F^{-1}\rho_0\le|\rho_0(\cdot+i\eta)|\le
A_F\rho_0$ and the bounds of (D4) on $\partial^\alpha_yF_0$ at complex velocities, as in the proof
of Lemma~\ref{lem:dispfree}. The terms of level $n$ carry the factor
$w_n(\rad)^2(n!)^2(\Lambda/4)^{-2n}\le(4\rad/\Lambda)^{2n}(n+1)^{10}$, whose sum over $n\ge1$ is
$\le C(\rad/\Lambda)^2$ for $\rad\le\Lambda/8$ (the range used in Proposition~\ref{prop:FPtr}, where
$\rad\le\rad_0^{\rm lay}/4\le\Lambda/8$), and the part with $n=0$, $\alpha\ne0$ the factor
$\delta_y^{2|\alpha|}\le\delta_y^2$, which is not small in $\rad$.
\end{proof}

\begin{lem}[Calculus on cubes]\label{lem:TWcalc}
Let $E_Q(u)$ and $S_Q(u)$ denote the analytic sums of a function $u$ of $\xi$ on the cube $Q$, in
mean square and in supremum, with the weights $w_n$, and $E^{\ge1}_Q,S^{\ge1}_Q$ their parts of
order $\ge1$; $A_m(u)=(\sum_{|\beta|=m}\|\partial^\beta u\|_H^2)^{1/2}$ for $u$ with values in a
Hilbert space $H$. Explicitly, with $\tilde w_n:=w_n/(n+2)^2$,
\begin{equation}
S_Q(u):=\Big(\sum_n\tilde w_n^2\sup_QA_n(u)^2\Big)^{1/2},\qquad
E_Q(u):=\Big(\sum_nw_n^2\rad_0^{-3}\!\int_QA_n(u)^2\Big)^{1/2};
\end{equation}
products are taken with a bounded bilinear map $H_1\times H_2\to H_3$ of norm $\le1$. There are
absolute constants $C_{\rm Sob}$ and $C_{\rm alg}$ such that, for
$\rad\in[\rad_0/4,\rad_0]$:
\begin{enumerate}[label=(\roman*)]
\item \emph{Sobolev.} $S_Q(u)\le C_{\rm Sob}E_{Q^*}(u)$.
\item \emph{Algebra.} $S_Q(uv)\le C_{\rm alg}S_Q(u)S_Q(v)$.
\item \emph{Tame product.} $E_Q(uv)\le C_{\rm alg}\big(S_Q(u)E_Q(v)+E_Q(u)S_Q(v)\big)$.
\item \emph{Product with a shifted factor.} Let $D_\beta(\xi)\ge0$, put
$d_n:=(\sum_{|\beta|=n}D_\beta^2)^{1/2}$ and $x_n:=w_n(\rad_0^{-3}\int_Qd_n^2)^{1/2}$. For the
Leibniz remainder
\begin{equation}
R_\beta:=\partial^\beta(u\,\partial_iv)-u\,\partial^\beta\partial_iv=\sum_{0<\beta'\le\beta}\tbinom\beta{\beta'}\partial^{\beta'}u\;\partial^{\beta-\beta'}\partial_iv
\end{equation}
we have
\begin{equation}
\sum_nw_n^2\rad_0^{-3}\!\int_Q\sum_{|\beta|=n}D_\beta\|R_\beta\|\le\frac{C_{\rm
alg}}\rad\Big(\sum_n\sum_{1\le j\le n/2}\sigma_j\,(n+1)\,x_n\,y_{n+1-j}+S_Q(v)\,E^{\ge1}_Q(u)\Big(\sum_nx_n^2\Big)^{1/2}\Big),
\end{equation}
where $y_m:=w_m(\rad_0^{-3}\int_QA_m(v)^2)^{1/2}$ and $\sigma_j:=(j+1)^{-2}\tilde
w_j\sup_QA_j(u)$. In particular $\sum_{j\ge1}\sigma_j(j+1)\le\zeta(2)^{1/2}S^{\ge1}_Q(u)$.
\item \emph{Cauchy data.} If $a$ extends holomorphically to the polydiscs $\Delta(\xi,r)$, $\xi\in
Q^*$, with $\|a\|_H\le B$ there, and $\rad_0\le r/2$, then $S_Q(a)+E_Q(a)\le CB$ and
$S^{\ge1}_Q(a)+E^{\ge1}_Q(a)\le CB\rad_0/r$.
\item \emph{Composition.} Let $\Phi$ be holomorphic on the polydisc $\Delta(w_0,r)\subset
\mathbb C^d$, with values in $H$ and $\|\Phi\|_H\le B$. Let $u:Q^*\to\R^d$ with $S_Q(u-w_0)\le r/(8dC_{\rm
alg})$. Then
\begin{equation}\begin{gathered}
S_Q(\Phi(u))\le CB,\\
E_Q\big(\Phi(u)-\Phi(w_0)\big)\le\frac{CB}rE_Q(u-w_0),\qquad
S_Q\big(\Phi(u)-\Phi(w_0)\big)\le\frac{CB}rS_Q(u-w_0),
\end{gathered}\end{equation}
with $C$ depending only on $d$.
\end{enumerate}
\end{lem}
\begin{proof}
We write $w_n=\rad^n(n+1)^5/n!$, so $\tilde w_n=\rad^n(n+1)^5/((n+2)^2n!)$.

\emph{(i) Sobolev.} The scaled embedding $H^2(Q^*)\subset L^\infty(Q)$ in $\R^3$, for $H$-valued
functions on a cube of side $\simeq\rad_0$ and its double $Q^*$, implies that, for every $\beta$,
\begin{equation}
\sup_Q\|\partial^\beta u\|_H^2\le C\sum_{|\gamma|\le2}\rad_0^{2|\gamma|-3}\!\int_{Q^*}\|\partial^{\beta+\gamma}u\|_H^2
.
\end{equation}
We sum over $|\beta|=n$. Since each multi-index $\beta+\gamma$ arises from at most $C$ pairs, we get
$\sup_QA_n(u)^2\le C\sum_{k\le2}\rad_0^{2k-3}\int_{Q^*}A_{n+k}(u)^2$. For $k\le2$, we have
$\frac{w_n}{w_{n+k}}=\frac{(n+k)!}{n!\,\rad^k}\Big(\frac{n+1}{n+k+1}\Big)^5\le\frac{(n+2)^k}{\rad^k}$,
so for $\rad\ge\rad_0/4$, $\rad_0^kw_n\le4^k(n+2)^kw_{n+k}$ and $\rad_0^k\tilde
w_n\le4^k(n+2)^{k-2}w_{n+k}\le16\,w_{n+k}$, hence $\tilde w_n^2\rad_0^{2k}\le256\,w_{n+k}^2$.
Summing over $n$, we obtain $S_Q(u)^2\le C\,E_{Q^*}(u)^2$.

\emph{A pointwise bound.} By Leibniz's rule, Vandermonde's identity
$\sum_{|\beta'|=k,\beta'\le\beta}\binom\beta{\beta'}=\binom nk$, and Cauchy--Schwarz with the
multi-index count $N(\beta,k)\le(\min(k,n-k)+1)^2$ as in Lemma~\ref{lem:algebra}, we have
\begin{equation}\label{eq:TWpt}
A_n(uv)(\xi)\le\sum_{k\le n}\tbinom nk\big(\min(k,n-k)+1\big)A_k(u)(\xi)\,A_{n-k}(v)(\xi).
\end{equation}

\emph{(ii) Algebra.} We take suprema over $Q$ in~\eqref{eq:TWpt}. Since $\tilde w_n$ is comparable
to $\rad^n(n+1)^3/n!$ within the factor $4$, we have $\tilde w_n\binom nk\le
C\big(\frac{n+1}{(k+1)(n-k+1)}\big)^3\tilde w_k\tilde w_{n-k}$, and the coefficient of $(\tilde
w_k\sup A_k(u))(\tilde w_{n-k}\sup A_{n-k}(v))$ is $\le C(\min(k,n-k)+1)^{-2}$. The bound
$\|\sum_kc_{\min(k,n-k)}x_ky_{n-k}\|_{\ell^2_n}\le2\|c\|_{\ell^2}\|x\|_{\ell^2}\|y\|_{\ell^2}$
finishes the proof.

\emph{(iii) Tame product.} We integrate the square of the pointwise bound over $Q$. For $k\le n/2$
we take the supremum on $u$, and otherwise on $v$. For $k\le n/2$, since $(n+1)/(n-k+1)\le2$, we
have
\begin{equation}
w_n\tbinom nk(k+1)\le32(k+1)^{-4}\frac{\rad^k(k+1)^5}{k!}\,w_{n-k}=32(k+1)^{-4}(k+2)^2\tilde
w_k\,w_{n-k}\le C(k+1)^{-2}\tilde w_k\,w_{n-k}.
\end{equation}
The case $k>n/2$ is symmetric. The claim follows from Young's inequality $\ell^1*\ell^2\to\ell^2$,
with $\sum_k(k+1)^{-2}\tilde w_k\sup A_k(u)\le\zeta(4)^{1/2}S_Q(u)$ (Cauchy--Schwarz).

\emph{(iv) Product with a shifted factor.} Fix $n$ and $k=|\beta'|\ge1$. There are two cases.
\begin{itemize}
\item $k\le n/2$. We take the supremum on $\partial^{\beta'}u$ and integrate
$\partial^{\beta-\beta'}\partial_iv$, which has level $n-k+1$. We have
$\frac{\rad^k}{k!}=\frac{(k+2)^2\tilde w_k}{(k+1)^5}$ and
$\frac{\rad^{n-k}}{(n-k)!}=\frac{n-k+1}\rad\frac{w_{n-k+1}}{(n-k+2)^5}$, and $(n+1)/(n-k+2)\le2$.
With the multi-index count $k+1$, we obtain
\begin{equation}
w_n\tbinom nk(k+1)\le32\,\frac{(k+2)^2}{(k+1)^4}\,\tilde
w_k\cdot\frac{n+1}\rad\,w_{n-k+1}\le\frac{128}\rad(n+1)(k+1)^{-2}\tilde w_k\,w_{n-k+1}.
\end{equation}
Pairing with $D_\beta$ and integrating over $Q$, we get the first sum, with $j=k$ and
$\sigma_j=(j+1)^{-2}\tilde w_j\sup_QA_j(u)$.
\item $k>n/2$. We take the supremum on $\partial^{\beta-\beta'}\partial_iv$, which has level
$m:=n-k+1\le\frac n2+1$, and integrate $\partial^{\beta'}u$. Now
$\rad^{n-k}/(n-k)!=(m/\rad)\rad^m/m!$ and $(n+1)/(k+1)\le2$, and the multi-index count is
$m$, so
\begin{equation}
w_n\tbinom nk\,m\le32\,w_k\,\frac m\rad\,\frac{(m+2)^2m}{(m+1)^5}\,\tilde w_m\le\frac
C\rad\,\frac{w_k\,\tilde w_m}{m+1}.
\end{equation}
The sequence $(m+1)^{-1}\tilde w_m\sup_QA_m(v)$ is in $\ell^1$, with sum $\le\zeta(2)^{1/2}S_Q(v)$.
By Young's inequality, we obtain the second term, with $E^{\ge1}_Q(u)$ since $k\ge1$.
\end{itemize}
The factor $(m+1)^{-1}$ is what the power $5$ leaves, namely one power each for the two Sobolev
levels, one for the multi-index count, and one for the shift. With the weight $(n+1)^p$ in place of
$(n+1)^5$ the factor is $\simeq(m+1)^{4-p}$, and the sum $\sum_m(m+1)^{4-p}\tilde w_m\sup_QA_m(v)$
is bounded by $S_Q(v)$ if and only if $\sum_m(m+1)^{2(4-p)}<\infty$, that is, $p>9/2$. At $p=9/2$
the divergence is logarithmic, and at $p=4$ the factor does not decay, and $\tilde
w_m\sup_QA_m(v)=(m+1)^{-1/2}\log^{-1}(m+2)$ is in $\ell^2$ but not in $\ell^1$. Therefore $5$ is the
smallest integer power. Finally $\sum_{j\ge1}\sigma_j(j+1)=\sum_{j\ge1}(j+1)^{-1}\tilde w_j\sup
A_j(u)\le\zeta(2)^{1/2}S^{\ge1}_Q(u)$ by Cauchy--Schwarz. Only $j\ge1$ is used, in ($\Theta$2) and
(M2), since $\sigma_0$ involves $\sup_Q\|u\|_H$, which is not part of $S^{\ge1}_Q(u)$.

\emph{(v) Cauchy data.} By the Cauchy estimates on the polydiscs $\Delta(\xi,r)$, we have
$\|\partial^\beta a(\xi)\|_H\le\beta!\,r^{-|\beta|}B$, so $A_n(a)\le
B\big(\sum_{|\beta|=n}(\beta!)^2\big)^{1/2}r^{-n}\le(n+1)\,n!\,r^{-n}B$. Hence
$w_nA_n(a)\le(\rad_0/r)^n(n+1)^6B$, and $\tilde w_nA_n(a)\le(\rad_0/r)^n(n+1)^6B$. For $\rad_0\le
r/2$ the sums over $n$ are $\le CB$, with $C^2\ge\sum_n(n+1)^{12}4^{-n}$, and the sums over $n\ge1$
are $\le CB\rad_0/r$. For $E_Q$ the constant contains in addition the factor $(|Q|/\rad_0^3)^{1/2}$,
which is bounded for the cubes of side $\simeq\rad_0$.

\emph{(vi) Composition.} We expand $\Phi(w)=\sum_\mu a_\mu(w-w_0)^\mu$, with $\|a_\mu\|_H\le
Br^{-|\mu|}$ by the Cauchy estimates. By~(ii), we have $S_Q((u-w_0)^\mu)\le C_{\rm
alg}^{|\mu|-1}S_Q(u-w_0)^{|\mu|}$. By~(iii) and induction on $|\mu|$, we get
\begin{equation}
E_Q\big((u-w_0)^\mu\big)\le|\mu|\big(C_{\rm alg}S_Q(u-w_0)\big)^{|\mu|-1}C_{\rm alg}E_Q(u-w_0).
\end{equation}
The number of $\mu$ with $|\mu|=m$ is $\le(m+1)^d$, and $C_{\rm alg}S_Q(u-w_0)/r\le1/(8d)$. Thus the
series $\sum_\mu\|a_\mu\|_H S_Q((u-w_0)^\mu)$ and $\sum_{\mu\ne0}\|a_\mu\|_HE_Q((u-w_0)^\mu)$
converge geometrically, with the stated bounds. The constant term has $S_Q(1)=\frac14$ and
$E_Q(1)\le C$. This completes the proof.
\end{proof}

\begin{remark}[Truncated forms]\label{rem:TWtrunc}
Write $S^{\le k}_Q$, $E^{\le k}_Q$ for the sums over the levels $n\le k$. The proofs of
Lemma~\ref{lem:TWcalc} are level by level, so its parts hold in the following truncated forms (we
keep its numbering, and part~(v), the Cauchy estimates, does not involve the levels):
\begin{itemize}
\item[(i)] $S^{\le k}_Q(u)\le C_{\rm Sob}E^{\le k+2}_{Q^*}(u)$ (the factor $(n+2)^{-2}$ in $\tilde
w_n$ compensates the two Sobolev levels, since
$\rad_0^{|\gamma|}w_n\le4^{|\gamma|}(n+2)^{|\gamma|}w_{n+|\gamma|}$, hence $\tilde
w_n^2\rad_0^{2|\gamma|}\le256\,w_{n+|\gamma|}^2$, for $\rad\ge\rad_0/4$, $|\gamma|\le2$);
\item[(ii)] $S^{\le k}(uv)\le C_{\rm alg}S^{\le k}(u)S^{\le k}(v)$;
\item[(iii)] $E^{\le k}(uv)\le C_{\rm alg}\big(S^{\le\lfloor k/2\rfloor}(u)E^{\le k}(v)+E^{\le
k}(u)S^{\le\lfloor k/2\rfloor}(v)\big)$;
\item[(iv)] holds with $\sum_{n\le k}$ on the left, and on the right with $x_n$ ($n\le k$), $y_m$
($m\le k$), $\sigma_j$ ($1\le j\le k/2$), $S^{\le\lfloor k/2\rfloor+1}_Q(v)$ and $E^{\ge1,\le
k}_Q(u)$;
\item[(vi)] the composition bounds hold for $E^{\le k}$ under $S^{\le\lfloor k/2\rfloor}(u-w_0)\le
r/(8dC_{\rm alg})$, and for $S^{\le k}$ under $S^{\le k}(u-w_0)\le r/(8dC_{\rm alg})$.
\end{itemize}
In (iii) and (iv) the supremum is on a factor of level $\le n/2+1$ (in (iv), second case, the factor
$\partial^{\beta-\beta'}\partial_iv$). Thus, together with (i), a truncated estimate at levels $\le
N$ uses Sobolev only up to $N/2+3\le N$ levels, for $N\ge6$.
\end{remark}

\begin{lem}[Neighbor sums with damping]\label{lem:EPdamp}
Let $\ell_Q:=\vartheta_\natural(v_0(Q))^2\ge1$, so that $\ell_{Q'}\le\ell_Q+2\le3\ell_Q$ for $Q'\sim
Q$. Suppose that, for every $Q$ and $t\le T$,
\begin{equation}
e_Q(t)+c\!\int_0^t\!\mathfrak d_Q+2K\eps\!\int_0^t\!\ell_Qe_Q\le\int_0^t\Big(C^{\rm
lay}_{\sharp}\eps\,\ell_Q\!\sum_{Q'\sim Q}e_{Q'}+\eta\!\sum_{Q'\sim Q}\mathfrak d_{Q'}+F_Q\Big),
\end{equation}
with $e_Q(0)=0$ and a constant $C^{\rm lay}_{\sharp}$ that does not depend on $Q$. Put $
\mathbb E(t):=\sup_Qe_Q(t)$, $\mathcal D(t):=\sup_Q\int_0^t\mathfrak d_Q$, $
\mathcal L(t):=\sup_Q\int_0^t\ell_Qe_Q$ and $\mathcal F(t):=\sup_Q\int_0^tF_Q$, and assume that $
\mathcal D(t)$ and $\mathcal L(t)$ are finite. If $K\ge9C^{\rm lay}_{\sharp} N_0$ and $\eta N_0\le c/6$, then
\begin{equation}
\sup_{s\le t}\mathbb E(s)+\tfrac c2\mathcal D(t)\le6\,\mathcal F(t).
\end{equation}
\end{lem}

\begin{proof}
Since $\ell_Q\le3\ell_{Q'}$ for neighbors, we have
\begin{equation}
\int_0^tC^{\rm lay}_{\sharp}\eps\ell_Q\sum_{Q'\sim Q}e_{Q'}\le3C^{\rm lay}_{\sharp}\eps\sum_{Q'\sim
Q}\int_0^t\ell_{Q'}e_{Q'}\le3C^{\rm lay}_{\sharp} N_0\eps\,\mathcal L(t),
\end{equation}
and $\int_0^t\eta\sum_{Q'\sim Q}\mathfrak d_{Q'}\le\eta N_0\mathcal D(t)$. The right side of the
hypothesis is therefore bounded by $3C^{\rm lay}_{\sharp}N_0\eps\mathcal L+\eta N_0\mathcal D+
\mathcal F$, uniformly in $Q$.

The supremum over $Q$ of the left side is at least the supremum of each of its three nonnegative
parts, and hence at least one third of their sum. Thus we have
\begin{equation}
\tfrac13\big(\mathbb E(t)+c\mathcal D(t)+2K\eps\mathcal L(t)\big)\le3C^{\rm
lay}_{\sharp}N_0\eps\mathcal L+\eta N_0\mathcal D+\mathcal F .
\end{equation}
Since $2K/3\ge6C^{\rm lay}_{\sharp}N_0>3C^{\rm lay}_{\sharp}N_0$ and $c/3-\eta N_0\ge c/6$, we
deduce that $\mathbb E(t)+\frac c2\mathcal D(t)\le3\mathcal F(t)$. The right side is nondecreasing
in $t$, which yields the supremum over $s\le t$ with the constant $6$.
\end{proof}

\section{Frame terms, viscous approximation and constants}\label{app:regimes}

This appendix contains three technical parts of the bootstrap: the frame terms and the band balance,
the viscous approximation, and the order in which the constants are chosen.

\subsection{The frame terms, the band balance and the crude leaf inequality}\label{app:band}

Proposition~\ref{prop:Bband} and the crude leaf inequality, Proposition~\ref{prop:ACpf} below, both
use the bound on the frame terms which we prove here.

\emph{Notation.} We fix a leaf $\xi$ and a component
$\varphi=\varphi_{\alpha\beta}=T^{-1}_{p}h^{(x)}$, $h=h_{\alpha\beta}$, $p=p_g(t,\xi)$, as in
Subsection~\ref{ssec:cells}, and write $\varphi^+=(\varphi_{\alpha,\beta+e_j})_j$, $\hat m=\hat
m_\xi$, $X=X(\hat m)$, $Y_2=Y_2(\hat m)$, and $D_t$ for the derivative along the ray.
We denote by $q_\xi:=q(t)\theta_g(t,\xi)$ the exponent of $\hat m$ in $\hat z$ (in
Lemma~\ref{lem:Gleafweight} the static analogue is $q\theta_\infty(\xi)$), by $\lambda_\xi$ the
collision rate of the leaf, and by $\rad=\rad(t)$ and $\rad_k$ the analyticity radii. The background
is
$G_\xi=T^{-1}_pg^{(x)}=\mu+f_\xi$, where $f_\xi=\varphi_{00}$. The background has mass $1$, mean $0$ and
temperature $1$, since $p_g$ are the leaf moments. We write $\mathcal L=L_\mu$, $L_Gu=Q(G,u)+Q(u,G)$
and $\Pi=I-\Pi_0$, with $\Pi_0$ the projection of Subsection~\ref{ssec:UF} (onto $\ker\mathcal L$
along the moment-free functions). The weight spaces are $X_k=L^2(\jap v^k)$, $Y_k=H^1_*(\jap v^k)$,
the dissipation space of \cite{CM}, and $Z_k$, the dual of $Y_k$ for the $X_k$-pairing.

\subsubsection{Inputs}

\emph{The estimates of \cite{CM}.} By (E1), with Proposition~\ref{prop:CMm2} at $\gamma=-2$, we have
the following at the level $k=7$:
\begin{itemize}
\item[(CM1)] Proposition~\ref{prop:CMthr} (\cite[Prop.~3.6, Cor.~3.7, Thm.~3.5]{CM}, with the
moment term estimated): for $k'\ge7$ there is $a>\frac12$ with $\|e^{\tau\mathcal L}\Pi u\|_{L^2}\le
C\jap\tau^{-a}\|u\|_{X_{k'}}$ for $u\in X_{k'}$ (part~(i) with $k_0=0$ and
$k_*\in(\gamma+\frac72,k'-\frac32|\gamma|)$, an interval which is nonempty since
$k'>\frac72+\frac12|\gamma|$); and
\begin{equation}\label{eq:J7}
J_7:=\int_0^\infty\|e^{\tau\mathcal L}\Pi\|_{Y_7\to L^2}\|e^{\tau\mathcal L}\Pi\|_{Z_7\to
L^2}\dd\tau<\infty
\end{equation}
(part~(iii));
\item[(CM2)] \cite[Lemma~4.3]{CM}, with the embeddings $X_7\subset L^2(\jap v^\theta)$,
$\theta\in(\frac72,7]$, and $Y_7\subset H^1(\jap v^{\theta'})$, $\theta'\in(\frac94,7+\frac\gamma2)$
\cite[Cor.~4.4]{CM}: $|\langle Q(f,u),\varphi\rangle_{X_7}|\le
C(\|f\|_{X_7}\|u\|_{Y_7}+\|f\|_{Y_7}\|u\|_{X_7})\|\varphi\|_{Y_7}$. By the definition of $Z_7$ as
the dual of $Y_7$ for the $X_7$-pairing, $\|Q(f,u)\|_{Z_7}$ is bounded by the same factor.
\end{itemize}
The nonlocal pairing is
\begin{equation}
\ip u{u'}_\natural:=\int_0^\infty\big\langle e^{\tau\mathcal L}\Pi u,e^{\tau\mathcal L}\Pi
u'\big\rangle_{L^2}\dd\tau,\qquad\nn u^2_\natural:=\ip uu_\natural ,
\end{equation}
which is the nonlocal part of the hybrid norm of \cite{CM}. By (CM1) and
\begin{equation}
|\langle e^{\tau\mathcal L}\Pi u,e^{\tau\mathcal L}\Pi u'\rangle_{L^2}|\le\|e^{\tau
\mathcal L}\Pi\|_{Z_7\to L^2}\|e^{\tau\mathcal L}\Pi\|_{Y_7\to L^2}\|u\|_{Z_7}\|u'\|_{Y_7},
\end{equation}
we see that it is defined for $u\in Z_7$, $u'\in Y_7$, and we have $|\ip u{u'}_\natural|\le
J_7\|u\|_{Z_7}\|u'\|_{Y_7}$.

\begin{lem}[The nonlocal pairing]\label{lem:lowweight}
For $k'\ge7$ there is $C_\natural(k')$ such that $|\ip u{u'}_\natural|\le
C_\natural(k')\|u\|_{X_{k'}}\|u'\|_{X_{k'}}$ for $u,u'\in X_{k'}$. For moment-free $u\in X_7\cap
Y_7$, $\ip u{\mathcal Lu}_\natural=-\frac12\|u\|^2_{L^2}$.
\end{lem}

\begin{proof}
By (CM1), the integrand of $\ip u{u'}_\natural$ is bounded by
$C^2\jap\tau^{-2a}\|u\|_{X_{k'}}\|u'\|_{X_{k'}}$, which is integrable since $2a>1$. This proves the
first claim. For the identity, note that $\Pi u=u$ and $\Pi\mathcal Lu=\mathcal Lu=\mathcal L\Pi u$,
since $u$ and $\mathcal Lu$ are moment-free and $\Pi_0$ commutes with $\mathcal L$. Moreover,
$\mathcal Lu\in Z_7$ by (CM2), since $u\in X_7\cap Y_7$. Therefore we have $\frac{d}{d\tau}\|e^{\tau
\mathcal L}u\|^2_{L^2}=2\langle e^{\tau\mathcal L}u,e^{\tau\mathcal L}\mathcal Lu\rangle_{L^2}$. We
integrate over $(0,\infty)$ and use $e^{\tau\mathcal L}u\to0$ in $L^2$, which holds by (CM1) since
$u\in X_7$. This proves the identity.
\end{proof}

\emph{The Gaussian hypodissipativity.} By Lemma~\ref{lem:TG0}(i) and (F3), $\hat m=\hat m_\xi$ lies
in the class $\mathcal W(q_\pm,\eps_1)$ in the frame $p_g$, with $(L+1,\theta_-/2,\theta_++1)$ in
place of $(L,\theta_\pm)$. Membership alone does not imply Lemma~\ref{lem:GW}(i), since $\eps_1$
need
not be small. What we use is the gradient bound
$|\nabla\omega|\le\eps_1(t)=2(\theta_++1)^{1/2}(L+1)(q(t)+t^{-1})$ of Lemma~\ref{lem:TG0}(i) in the
frame $p_g$. Since $q(t)\le q_0\le q_{\rm GW}(\gamma,L+1,\theta_-/2,\theta_++1)$ and
$t\ge T_1$, we have $C\eps_1(t)^2\le q_-$, which is the condition of the last paragraph of the
proof of Lemma~\ref{lem:GW}(i). Hence, by part~(ii) of that lemma, or by Corollary~\ref{cor:W1G}, we
obtain, for $\|f\|_X\le1$ and every $\varphi$,
\begin{equation}\label{eq:W1Gband}
\int\hat m^2\varphi\,L_{\mu+f}\varphi\le-c_W\|\varphi\|^2_{Y_2}+C_W\|\varphi\|^2_{L^2}+C_W\|f\|_X\|\varphi\|^2_{Y_2}+C_W\|f\|_{Y_2}\|\varphi\|_X\|\varphi\|_{Y_2}+C_W\|f\|^2_{Y_2}\|\varphi\|^2_X
.
\end{equation}

\emph{Constants.} Let $C_{L^2}$ be the constant with $\|\varphi\|_{L^2}\le C_{L^2}\|\varphi\|_X$,
and $C_\natural:=C_\natural(8)C_{\rm G}^2$, with $C_{\rm G}$ the constant of Lemma~\ref{lem:TG0}(ii)
for $k=8$, so that $\nn\varphi^2_\natural\le C_\natural\|\varphi\|^2_X$ by Lemma~\ref{lem:lowweight}
with $k'=8$. Put
\begin{equation}
\eta:=\frac1{4C_W},\qquad C_1^\natural:=\max\Big\{\frac{2C_WC_{L^2}^2}{c_\lambda},\
\frac{8C_WC_\natural}{c_\lambda^2}\Big\},
\end{equation}
and take $C_1^B\ge C_1^\natural$. The coefficient $\chi$ is nondecreasing and $C^1$,
with $0\le\chi'\le2$, $\chi=0$ on $(-\infty,S_b]$ and $\chi=1$ on $[S_b+1,\infty)$. Along the ray,
$D_t\omega_\xi=-2H'(s_\xi)\dot s_\xi\,\omega_\xi$, with $\dot s_\xi=\kappa\lambda_{\rm
ref}(\xi)\in[c_\lambda,c_\lambda^{-1}]\,\kappa\lambda_\xi$.

\emph{The sources.} By Lemma~\ref{lem:leafeq}(ii), the collision part of $D_t\varphi$ is
$\kappa\lambda_\xi\big(L_{\mu+f_\xi}\varphi+\Sigma_\varphi\big)$. Here
$\kappa\lambda_\xi\Sigma_\varphi$ is the normalized sum of the Leibniz terms of
$\partial^\alpha_y\partial^\beta_v\big[\kappa\big(\mathcal L_t[g]\hat r+\cQ_t(\hat r,\hat
r)\big)\big]$ in which a derivative falls on the coefficient $\mathcal M[g]$ or which are quadratic
in components of lower order. At $(\alpha,\beta)=(0,0)$ we have $\Sigma=-Q(f_\xi,f_\xi)$, since
$\mathcal L f+Q(f,f)=L_{\mu+f}f-Q(f,f)$. The source pairing is
\begin{equation}
\mathcal S^{\rm src}_\xi(\varphi):=\eta\langle\varphi,\Sigma_\varphi\rangle_X+\chi\ip\varphi{\Sigma_\varphi}_\natural
.
\end{equation}
For $\alpha\ne0$ it contains the coefficient terms of Remark~\ref{rem:UFgrade}.

\subsubsection{The frame terms}

Put $A_p:=\rho^{-1}\theta^{3/4}$ and $\mathfrak p_\xi:=t^2|D_t\varpi_g(t,\xi)|$, with
$\varpi_g=(\ell_g,y_{0,g},\vartheta_g)$.

\begin{prop}[The frame terms]\label{prop:FM}
\leavevmode
\begin{enumerate}[label=(\roman*)]
\item \emph{Local part.} The frame terms in $D_t\big(\eta\|\varphi\|^2_X\big)$, namely
$-2\eta\langle\varphi,\mathcal G_p[D_tp]\varphi\rangle_X$ and the part of $D_t\hat m^2$ caused by
the motion of the frame, add up exactly to $2D_t\log A_p\cdot\eta\|\varphi\|^2_X$, and $|2D_t\log
A_p|\le4t^{-2}\mathfrak p_\xi$.
\item \emph{Nonlocal part.} With $C_{11}=C(\theta_\pm,q_\pm,\eps_1)$,
\begin{equation}\label{eq:FM}
|2\chi\ip\varphi{\mathcal G_p[D_tp]\varphi}_\natural|\le C_{11}t^{-2}
\mathfrak p_\xi\|\varphi\|^2_{Y_2} .
\end{equation}

\item \emph{Size of the frame velocity.} With $\varpi'_g=(\nabla\ell_g,y_{0,g},\vartheta_g)$ at
$(t,\xi)$,
\begin{equation}
\mathfrak p_\xi\le C_{\mathfrak p}\Big(1+\sum_{|\beta'|\le1}\|\varphi_{0\beta'}\|_X\Big),\qquad
C_{\mathfrak p}:=C(\theta_\pm,q_-,\eps_1)\big(1+|\varpi'_g|+|\nabla_\xi(y_{0,g},\vartheta_g)|\big)^2
.
\end{equation}
\end{enumerate}
\end{prop}

\begin{proof}
(i) In the frame $p$, $z=y_0+\theta^{1/2}\hat z$, and $\hat m$ is $m_x$ in the variable $\hat z$,
divided by the time-independent $A_{\rm ref}(\xi)$. At fixed $\hat z$ the frame moves $z$ with
velocity $D_ty_0+\frac12\theta^{1/2}D_t\vartheta\,\hat z$, and $\nabla_z=\theta^{-1/2}\nabla_{\hat
z}$. Therefore the frame part of $D_t\hat m^2$ is $V\cdot\nabla_{\hat z}(\hat m^2)$, with
$V:=e^{-\vartheta/2}D_ty_0+\frac12D_t\vartheta\,\hat z$. By Lemma~\ref{lem:leafeq}(i), $
\mathcal G_p[D_tp]\varphi=\nu_p\varphi-V\cdot\nabla\varphi$, with $\nu_p:=D_t\ell-\frac32D_t\vartheta$. Hence
the sum of the two terms is $\eta$ times
\begin{equation}
-2\nu_p\!\int\hat m^2\varphi^2+2\!\int\hat m^2\varphi\,V\cdot\nabla\varphi+\int(V\cdot\nabla\hat
m^2)\varphi^2=-2\nu_p\!\int\hat m^2\varphi^2+\int V\cdot\nabla(\hat
m^2\varphi^2)=-(2\nu_p+\nabla\cdot V)\int\hat m^2\varphi^2 .
\end{equation}
Since $\nabla\cdot V=\frac32D_t\vartheta$, we get $-(2\nu_p+\nabla\cdot
V)=-2D_t\ell+\frac32D_t\vartheta=2D_t\log A_p$, and $|2D_t\log
A_p|\le2|D_t\ell|+\frac32|D_t\vartheta|\le4|D_t\varpi_g|$. (In $z$-variables, the frame enters
$\ell_\xi[h]$ only through the factor $(A_p/A_{\rm ref})^2$.)

(ii) Both $\varphi$ and $\mathcal G_p[D_tp]\varphi$ are moment-free by the same lemma. By
Lemma~\ref{lem:lowweight} with $k'=7$, we have $|\ip\varphi{\mathcal G\varphi}_\natural|\le
C_\natural(7)\|\varphi\|_{X_7}\|\mathcal G\varphi\|_{X_7}$, and $\|\mathcal G_p[\dot
p]\varphi\|_{X_7}\le C(\theta_-)|\dot p|\,\|\varphi\|_{H^1(\jap{\hat z}^8)}$. By
Lemma~\ref{lem:TG0}(ii), $\|\varphi\|_{X_7}$ and $\|\varphi\|_{H^1(\jap{\hat z}^8)}$ are bounded by
$C\|\varphi\|_{Y_2}$.

(iii) We bound the right-hand sides of Theorem~\ref{thm:Meuler}, using $e^\vartheta\le\theta_++1$.
They contain only $\nabla\ell$, $y_0$, $\vartheta$, $\nabla y_0$, $\nabla\vartheta$ and $s,h,\nabla
s,\nabla h$. The micro moments $s=s[\psi_\xi]$ and $h=h[\psi_\xi]$ are moments of order $\le3$ of
$\psi_\xi=f_\xi=\varphi_{00}$, and thus we have $|s|+|h|\le\|\hat m\psi\|_{L^2}\|\jap{\hat z}^3\hat
m^{-1}\|_{L^2}\le C\|\varphi_{00}\|_X$, since $\hat m^{-1}\le Ce^{-q_-|\hat z|^2+\eps_1|\hat z|}$ by
part~(i) of the same lemma. Their $\xi$-derivatives are moments of
$\partial_\xi\psi_\xi=\varphi_{0e_i}-\mathcal G_p[\partial_{\xi_i}p]\psi_\xi$ by
Lemma~\ref{lem:leafgeom}(ii) and Lemma~\ref{lem:leafeq}(i). After an integration by parts onto the
polynomial, we obtain $|\nabla_\xi s|+|\nabla_\xi h|\le
C\big(\sum_i\|\varphi_{0e_i}\|_X+(|\varpi'_g|+|\nabla_\xi(y_{0,g},\vartheta_g)|)\|\varphi_{00}\|_X\big)$,
since $\partial_\xi p$ contains only first derivatives of $\varpi_g$.
\end{proof}

\begin{cor}[Use of the frame bound]\label{cor:FMuse}
\leavevmode
\begin{enumerate}[label=(\alph*)]
\item In the proof of Proposition~\ref{prop:Bband}, the local part of (T11) is $\le4t^{-2}
\mathfrak p_\xi\,\eta\|\varphi\|^2_X$, and its nonlocal part is $\le C_{11}t^{-2}
\mathfrak p_\xi\|\varphi\|^2_{Y_2}$. The latter is absorbed by (T1) wherever $\chi>0$ as soon as $C_{\rm
r}t^{-2}\le\frac{\eta c_W}4\kappa\lambda_\xi$, with $C_{\rm r}:=C_{11}\mathfrak p_*$, if $
\mathfrak p_\xi\le\mathfrak p_*$. This condition follows from requirement~\ref{rq:Wrate}.
\item Where $\chi=0$ (in particular in Proposition~\ref{prop:ACpf}) only the local part occurs.
\item In the analytic sums, the frame terms contribute $4t^{-2}\mathfrak p_\xi\hat{\mathbb A}_\xi$
(plus the nonlocal part, which is absorbed as in (a)). They enter Lemma~\ref{lem:analytic} through
$a(t)=4t^{-2}\mathfrak p_*$, with $\int_{T_0}^\infty a\le4\mathfrak p_*/T_0$, and not as a loss of
one level.
\item \emph{The constant $\mathfrak p_*$.} Assume (U1), $C\delta_U\le\rad_0$ and
$\sup_Q\mathfrak a_{Q^*}\le B$ at time $t$. By Lemma~\ref{lem:Mref}, Lemma~\ref{lem:TWcalc}(i) and
$C\delta_U\le\rad_0$, $|\varpi'_g|+|\nabla_\xi(y_{0,g},\vartheta_g)|\le C_\varpi=C(L_1,R_1)+1$,
as in Subsection~\ref{ssec:cells}. For $\xi$ in a cube $Q$, Lemma~\ref{lem:CUsob} with
$(\alpha,\beta)=(0,\beta')$, $|\beta'|\le1$, implies that $\tilde
w_{|\beta'|}(\rad(t))^2\eta\|\varphi_{0\beta'}\|^2_X\le C'_SB$. Here $\tilde w_0=\frac14$,
$\tilde w_1(\rad)=\frac{32}9\rad$, and $\rad(t)\ge\rad_0/4$ is the analyticity radius.
The graded radii lie in $[\rad_0/4,\rad_0]$, also in Subsection~\ref{ssec:cellsproofs}. Hence
\begin{equation}\label{eq:FMuse-low}
\sum_{|\beta'|\le1}\|\varphi_{0\beta'}\|_X\le\big(4+\frac{27}{8\rad_0}\big)(C'_SB/\eta)^{1/2}\le27(C'_SB/\eta)^{1/2}\big(1+(8\rad_0)^{-1}\big),
\end{equation}
and we have, on every leaf, before and after entry, $\mathfrak p_\xi\le\mathfrak p_*$ with
\begin{equation}
\mathfrak p_*:=C_{\mathfrak p,*}\Big(1+\big(C'_SB/\eta\big)^{1/2}\big(1+(8\rad_0)^{-1}\big)\Big),\qquad
C_{\mathfrak p,*}:=27\,C(\theta_\pm,q_-,\eps_1)(1+C_\varpi)^2 .
\end{equation}
This is a nondecreasing function of $B$, evaluated at $B=B_{\rm cr}$. It does not involve $T_0$,
$\delta_U$ or $\eps$.
\end{enumerate}
\end{cor}

\begin{proof}
(a) and (b) follow from Proposition~\ref{prop:FM}(i)--(ii). (c) We sum (i) over the components with
the weights $W^2_{\alpha\beta}$. (d) We combine part~(iii) of that proposition with the bounds
stated,
and $\rad(t)\ge\rad_0/4$.
\end{proof}

\subsubsection{Proof of Proposition~\ref{prop:Bband}}\label{sssec:bandproof}

We differentiate $\omega e_\xi[\varphi]=\omega\big(\eta\|\varphi\|^2_X+\chi\nn\varphi^2_\natural\big)$
along the ray, and insert the equation of Lemma~\ref{lem:leafeq}(ii). We list every term below. The
terms (T1)--(T7) are written after division by $\kappa\lambda_\xi\omega$, while the terms
(T8)--(T12)
carry the factor $t^{-2}$. We put $b_\xi:=\|f_\xi\|^2_X+\|f_\xi\|^2_{Y_2}$, the square of the size
of the micro
part.

\emph{Collision terms}, from $\kappa\lambda_\xi L_{\mu+f_\xi}\varphi$, by~\eqref{eq:W1Gband} with
$f=f_\xi$, $\|f_\xi\|_X\le\delta_c\le1$, and by Lemma~\ref{lem:lowweight}:
\begin{itemize}
\item[(T1)] local dissipation: $-2\eta c_W\|\varphi\|^2_{Y_2}$;
\item[(T2)] local compact part: $+2\eta C_W\|\varphi\|^2_{L^2}=+\frac12\|\varphi\|^2_{L^2}$;
\item[(T3)] nonlocal coercivity: $2\chi\ip\varphi{
\mathcal L\varphi}_\natural=-\chi\|\varphi\|^2_{L^2}$;
\item[(T4)] background perturbation. The local part is the rest of the same bound, $\le2\eta
C_W\big(\delta_c\|\varphi\|^2_{Y_2}+b^{1/2}_\xi\|\varphi\|_X\|\varphi\|_{Y_2}+b_\xi\|\varphi\|^2_X\big)$.
The nonlocal part is $2\chi\ip\varphi{L_{f_\xi}\varphi}_\natural$. By (CM1)--(CM2) and
Lemma~\ref{lem:TG0}(ii), we see that it is $\le2\chi
C\big(\delta_c\|\varphi\|_{Y_2}+b^{1/2}_\xi\|\varphi\|_X\big)\|\varphi\|_{Y_2}$;
\item[(T5)] sources: $2|\mathcal S^{\rm src}_\xi(\varphi)|$.
\end{itemize}

\emph{Weight terms:}
\begin{itemize}
\item[(T6)] damping: $D_t\omega=-2H'\dot s\,\omega$, with $H'\ge C_1^B$ wherever $\chi<1$, yields a
term
$\le-2C_1^Bc_\lambda\big(\eta\|\varphi\|^2_X+\chi\nn\varphi^2_\natural\big)$ there;
\item[(T7)] front motion: $\chi'(s)(\dot s/\kappa\lambda_\xi)\nn\varphi^2_\natural\le\frac2{c_\lambda}\nn\varphi^2_\natural$,
nonzero only where $\chi<1$;
\item[(T8)] Gaussian weight, the part of $D_t\hat m^2$ not caused by the frame: $D_t\hat m^2/\hat
m^2=2q'(t)|z|^2+2t^{-2}z\cdot\nabla\log A_{\rm ref}(\xi-z/t)$, with $q'=-C_2t^{-2}$ and $|\nabla\log
A_{\rm ref}|\le L\le2L$ by Corollary~\ref{cor:LWPa}(a4). Thus this term is
$t^{-2}\eta\int\big(-2C_2|z|^2+4L|z|\big)\hat m^2\varphi^2\le
t^{-2}\eta\int\big(-C_2|z|^2+C\big)\hat m^2\varphi^2$. It is kept in $\mathcal T^{\rm loc}_\xi$.
\end{itemize}

\emph{Transport, frame and forcing terms} (rate $t^{-2}$):
\begin{itemize}
\item[(T9)] local transport: $2t^{-2}\eta\langle\varphi,(y_0+\theta^{1/2}\hat
z)\cdot\varphi^+\rangle_X$, kept in $\mathcal T^{\rm loc}_\xi$ together with (T8);
\item[(T10)] nonlocal transport: $2t^{-2}\chi\ip\varphi{(y_0+\theta^{1/2}\hat
z)\cdot\varphi^+}_\natural\le Ct^{-2}\|\varphi\|_{X_7}\|\varphi^+\|_{X_8}\le
C_0t^{-2}e_\xi[\varphi]^{1/2}e_\xi[\varphi^+]^{1/2}$, by the same two lemmas;
\item[(T11)] frame motion: by Corollary~\ref{cor:FMuse}(a) and $\mathfrak p_\xi\le\mathfrak p_*$,
the local part is $\le4t^{-2}\mathfrak p_*\eta\|\varphi\|^2_X$, and the nonlocal part is $\le C_{\rm
r}t^{-2}\|\varphi\|^2_{Y_2}\le\frac{\eta c_W}4\kappa\lambda_\xi\|\varphi\|^2_{Y_2}$ wherever
$\chi>0$, and this bound does not use (T8);
\item[(T12)] forcing. In the equation of the leaf the forcing is the full $-T^{-1}_p\Phi$,
$\Phi:=(\partial^\alpha_y\partial^\beta_v\partial_t\mathcal M[g])^{(x)}$. Since $D_t\varphi$, the
collision terms and $\mathcal G_p[D_tp]\varphi$ are moment-free,
$\Pi_0T^{-1}_p\Phi=t^{-2}\Pi_0[(y_0+\theta^{1/2}\hat z)\cdot\varphi^+]$, as in Step~1 of the proof
of Proposition~\ref{prop:UFforce}. Thus, with $\hat F:=(I-\Pi_0)T^{-1}_p\Phi$, the forcing is $-\hat
F-t^{-2}\Pi_0[(y_0+\theta^{1/2}\hat z)\cdot\varphi^+]$.
\begin{itemize}
\item[(T12$_0$)] The \emph{kernel part} contributes
$-2t^{-2}\eta\langle\varphi,\Pi_0[(y_0+\theta^{1/2}\hat z)\cdot\varphi^+]\rangle_X\le
Ct^{-2}e_\xi[\varphi]^{1/2}e_\xi[\varphi^+]^{1/2}$ to the local part, since $\|\Pi_0h\|_X\le
C|\mathrm{mom}(h)|$ and
\begin{equation}
|\mathrm{mom}((y_0+\theta^{1/2}\hat z)\cdot\varphi^+)|\le
C\big(L+1+(\theta_++1)^{1/2}\big)\|\jap{\hat z}^3\hat m^{-1}\|_{L^2}\|\varphi^+\|_X\le
C\eta^{-1/2}e_\xi[\varphi^+]^{1/2},
\end{equation}
by (F3) and the membership lemma, uniformly in the leaf. It contributes nothing to the nonlocal
part, since $\Pi\Pi_0=0$. It is included in the $C_0$-term, where $C_0$ is taken at least the
sum of the constants of (T10) and~(T12$_0$).
\item[] With $F_\xi:=t^2\|\hat F\|_X$, the rest is
\begin{equation}
-2\eta\langle\varphi,\hat F\rangle_X-2\chi\ip\varphi{\hat
F}_\natural\le2(\eta^{1/2}+C_\natural^{1/2})e_\xi[\varphi]^{1/2}\|\hat F\|_X\le
Ct^{-2}e_\xi[\varphi]^{1/2}F_\xi .
\end{equation}
\end{itemize}
\end{itemize}

\emph{The balance of the $L^2$-level terms.} Where $\chi<1$, we collect (T2), (T3), (T6) and (T7),
and use $\|\varphi\|^2_{L^2}\le C_{L^2}^2\|\varphi\|^2_X$ and $\nn\varphi^2_\natural\le
C_\natural\|\varphi\|^2_X$. Then we have
\begin{equation}
\Big(\tfrac12-\chi\Big)\|\varphi\|^2_{L^2}-2C_1^Bc_\lambda\eta\|\varphi\|^2_X+\frac2{c_\lambda}\nn\varphi^2_\natural\le\Big(\frac{C_{L^2}^2}2-C_1^Bc_\lambda\eta\Big)\|\varphi\|^2_X+\Big(\frac{2C_\natural}{c_\lambda}-C_1^Bc_\lambda\eta\Big)\|\varphi\|^2_X\le0
,
\end{equation}
since $C_1^B\ge C_1^\natural$ and $\eta=1/(4C_W)$ imply $C_1^Bc_\lambda\eta\ge\frac{C_{L^2}^2}2$ and
$C_1^Bc_\lambda\eta\ge\frac{2C_\natural}{c_\lambda}$. Here we have dropped the term
$-2C_1^Bc_\lambda\chi\nn\varphi^2_\natural\le0$ of (T6). Where $\chi=1$, (T6) and (T7)
are absent (or $\le0$), and (T2)+(T3) $=-\frac12\|\varphi\|^2_{L^2}\le0$. No smallness of $\varphi$
is used.

\emph{The remaining collision terms.} By requirement~\ref{rq:Wdc}, we have $2(\eta C_W+\chi
C)\delta_c\le\frac{\eta c_W}2$, and thus the $\delta_c$-parts of (T4) remove at most $\frac{\eta
c_W}2\|\varphi\|^2_{Y_2}$ from (T1). By Young's inequality, the $b^{1/2}_\xi$-parts are
$\le\frac{\eta
c_W}4\|\varphi\|^2_{Y_2}+Cb_\xi\|\varphi\|^2_X$, and the nonlocal part of (T11) removes no more than
$\frac{\eta c_W}4\|\varphi\|^2_{Y_2}$ in addition. What is left of (T1) is $-\eta
c_W\|\varphi\|^2_{Y_2}=-\frac{c_W}{4C_W}\|\varphi\|^2_{Y_2}$. Multiplying by
$\kappa\lambda_\xi\omega$, and using $\eta\|\varphi\|^2_X\le e_\xi[\varphi]$ and
$\lambda_\xi\|\varphi\|^2_{Y_2}=d_\xi[\varphi]$, we obtain the inequality of
the proposition, with $Ct^{-2}\omega e_\xi[\varphi]$ collecting the local part of
(T11), $C=4\mathfrak p_*$, and the $C_0$-term collecting (T10) and the kernel part (T12$_0$) of
(T12). \qed

\subsubsection{The crude leaf inequality}

\begin{prop}[The crude leaf inequality]\label{prop:ACpf}
On a leaf with $s_\xi<S_*$, assume (G1)--(G4) for the normalized background with constant $K$,
$\|\delta\varpi\|_{\mathcal A}\le2\delta_U$ (used only in~(C), while in the cell version it is
replaced by
(H$_d$), through the remark after Proposition~\ref{prop:CUtri}), $\hat{\mathbb A}_{\ge1,\xi}\le
c_\sharp$, $\mathfrak p_\xi\le\mathfrak p_*$ (Corollary~\ref{cor:FMuse}(d)), and the requirements
$C_1^U\ge C_L/c_\lambda$, requirement~\ref{rq:Wcsharp}, and
$C(\rad_0/R_1+\delta_y+\delta_U)\le\eta c_h/4$ (requirements~\ref{rq:Wdy} and~\ref{rq:WdU}), and
$T_0\ge T_{\rm h}$ (see (L1) below), where $c_h=c_h(K)$ and $C_L=C_L(K)$ are the
constants of (L1)--(L2) below. Then
\begin{equation}
D_t\big(\omega\hat{\mathbb A}_\xi\big)+\eta c_h\,\kappa\,\omega\,\mathbb D_\xi\le\omega\,
\mathcal T_\xi+t^{-2}(C+4\mathfrak p_*)\omega\hat{\mathbb A}_\xi+Ct^{-2}\omega\hat{
\mathbb A}_\xi^{1/2}\mathbb F_\xi+C(K)\kappa\lambda_\xi\omega\Big(\frac{\rad_0}{R_1}+\delta_y+\delta_U\Big)^2
,
\end{equation}
where
\begin{equation}\label{eq:ACpf-T}
\mathcal T_\xi=\sum_n2n\frac{\rad'}{\rad}w_n^2E_n+\sum_{(\alpha,\beta)}W^2_{\alpha\beta}\big((
\mathrm T8)+(\mathrm T9)+(\mathrm T12_0)\big)
\end{equation}
is the weighted local transport and Gaussian weight term with the radius term of
Lemma~\ref{lem:analytic} and the kernel part of the forcing. Here
\begin{equation}
(\mathrm T12_0)=-2t^{-2}\eta\langle\varphi_{\alpha\beta},\Pi_0[(y_0+\theta^{1/2}\hat
z)\cdot\varphi^+_{\alpha\beta}]\rangle_X
\end{equation}
is the kernel part of (T12), and it satisfies, separately for each $y$-order $k$,
\begin{equation}
\sum_{|\alpha|=k}\sum_\beta W^2_{\alpha\beta}(\mathrm T12_0)\le
C_0t^{-2}\sum_nw_n(\rad_k)^2\sqrt{E_{k,n}E_{k,n+1}},
\end{equation}
with $E_{k,n}$ the part of $E_n$ of $y$-order $k$, and in cells it is
Proposition~\ref{prop:UFforce}(b) (also separately for each $y$-order). Moreover, $\mathbb F_\xi$ is
the rest of (T12), $F_\xi=t^2\|\hat F\|_X$, in the analytic norm.
The same holds for the parts with $(\alpha,\beta)\ne0$, with $\hat{\mathbb A}_{\ge1}$, $
\mathbb D_{\ge1}$, $\mathcal T_{\ge1}$, and the same last term.
\end{prop}

\begin{proof}
We differentiate $\omega W^2_{\alpha\beta}\eta\|\varphi_{\alpha\beta}\|^2_X$ along the ray for each
$(\alpha,\beta)$, and use Lemma~\ref{lem:leafeq}(ii). Here $\chi=0$. We write
$\varphi=\varphi_{\alpha\beta}$.

\emph{(L1) The first half of the linear term.} The weight is $\hat m=e^{q_\xi|\hat z-\hat
z_0|^2}e^{\tilde\omega}$, with $q_\xi=q(t)\theta_g$, $\hat z_0=-\theta_g^{-1/2}y_{0,g}$, $|\hat
z_0|\le L_0:=(2/\theta_-)^{1/2}(L+1)$ and $|\nabla\tilde\omega|\le2(\theta_++1)^{1/2}L/t$, by (F3).
In the variable $w=\hat z-\hat z_0$ the background has mass~$1$, center $-\hat z_0$ and
temperature~$1$. By (F3), $q(t)\le q_0$ and the last clause in the definition of $q_0$ in
Appendix~\ref{app:const}, we have $8q_\xi(|\hat
z_0|^2+3)=8q(t)\big(|y_{0,g}|^2+3\theta_g\big)\le8q_0\big((L+1)^2+3\theta_++3\big)\le1$.
Consequently Lemma~\ref{lem:Whypo} applies with $|\hat z_0|$ in place of $L$ and $\theta_+=1$, with
(G2) from Lemma~\ref{lem:RGellip} and (G3) with $16\jap{L_0}^8\bar K$ in place of $\bar K$ (since
$\jap w^8\le16\jap{\hat z_0}^8\jap{\hat z}^8$). The center enters its proof through
$q\int|w|^2G\le M/8$ alone, and thus its constants are uniform for $|\hat z_0|\le L_0$. The factor
$e^{\tilde\omega}$ is a perturbation of size $O(L/t)$, treated as in the proof of
Lemma~\ref{lem:GW}(i), for $t\ge T_0\ge T_{\rm h}(K,L,\theta_\pm,q_0)$. Thus we obtain
\begin{equation}
\langle Q(G_\xi,\varphi),\varphi\rangle_X\le-c_h\|\varphi\|^2_{Y_2}+C_h\|\varphi\|^2_X ,
\end{equation}
with $c_h,C_h$ depending on $(\gamma,K,L,\theta_\pm,q_0)$.

\emph{(L2) The second half.} We write $Q(\varphi,G_\xi)=Q(\varphi,\mu)+Q(\varphi,f_\xi)$ and
$Q(\varphi,f)=\ab[\varphi]:\nabla^2f-\cb[\varphi]f$. By \cite[Lemma~4.2(i)]{CM}, $|\ab[\varphi]|\le
C\jap{\hat z}^{\gamma+2}\|\varphi\|_{L^2(\jap{\hat z}^4)}\le C\|\varphi\|_X$. For $\gamma=-3$,
$\cb[\varphi]=-8\pi\varphi$. For $\gamma>-3$, we split the kernel of $\cb$ at $|w|=1$. By Young's
inequality, the near part is bounded by $C\|\varphi\|_{L^2}$ in $L^2$, and the far part by
$2(\gamma+3)\|\varphi\|_{L^1}$ pointwise. Hence, using (G4), we have
\begin{equation}
|\langle Q(\varphi,f_\xi),\varphi\rangle_X|\le C\big(\|\hat m\nabla^2f_\xi\|_{L^2}+\|\hat
mf_\xi\|_{L^2\cap L^\infty}\big)\|\varphi\|^2_X\le CK\|\varphi\|^2_X ,
\end{equation}
and the same bound with $\mu$ in place of $f_\xi$ yields $C\|\varphi\|^2_X$, since
$q_\xi\le\frac18$.
Hence the linear term contributes
$\le2\eta\kappa\lambda_\xi\big(-c_h\|\varphi\|^2_{Y_2}+C_L\|\varphi\|^2_X\big)$, with
$C_L:=C_h+C(1+K)$. After summation its first part is $-2\eta c_h\kappa\mathbb D_\xi$. At
$(\alpha,\beta)=(0,0)$ the collision term is $
\mathcal Lf_\xi+Q(f_\xi,f_\xi)=Q(G_\xi,f_\xi)+Q(f_\xi,\mu)$, as noted in the paragraph on the
sources, so only the $\mu$-part of (L2) occurs there, and the bound holds a fortiori.

\emph{(W) Damping.} We have $D_t\omega=-2C_1^U\dot s\,\omega$ and $\dot s\ge
c_\lambda\kappa\lambda_\xi$. Thus the damping yields
$-2C_1^Uc_\lambda\kappa\lambda_\xi\,\omega\eta\|\varphi\|^2_X$, which absorbs the
$C_L$-part, since $C_1^U\ge C_L/c_\lambda$.

\emph{(N) Quadratic terms.} These are $Q(\varphi_{\alpha'\beta'},\varphi_{\alpha''\beta''})$ with
$0<(\alpha',\beta')<(\alpha,\beta)$. By Proposition~\ref{prop:trilinG}(iii), second form, we have
\begin{equation}
|\langle Q(\varphi_1,\varphi_2),\varphi\rangle_X|\le
C\big(\|\varphi_1\|_X\|\varphi_2\|_{Y_2}+\|\varphi_1\|_{Y_2}\|\varphi_2\|_X\big)\|\varphi\|_{Y_2} .
\end{equation}
Both factors lie in the part $(\alpha,\beta)\ne0$, and the weights split as
$W_{\alpha\beta}\binom\alpha{\alpha'}\binom\beta{\beta'}\le2^{N_y}2^4\varrho\,W_{\alpha'\beta'}W_{\alpha''\beta''}$,
with $\varrho$ as in Lemma~\ref{lem:algebra}. Using $\|\varphi_1\|_X=(e_\xi[\varphi_1]/\eta)^{1/2}$,
we deduce from the same lemma that after summation these terms are $\le C_N\eta^{1/2}\kappa\hat{
\mathbb A}_{\ge1,\xi}^{1/2}\mathbb D_\xi\le\frac{\eta c_h}4\kappa\mathbb D_\xi$, by
requirement~\ref{rq:Wcsharp}.

\emph{(C) Coefficient terms.} These are $Q(
\mathfrak g_{\alpha'\beta'},\varphi_{\alpha''\beta''})+Q(\varphi_{\alpha''\beta''},
\mathfrak g_{\alpha'\beta'})$ with $(\alpha',\beta')\ne0$, where $\mathfrak g_{\alpha'\beta'}$ are the
normalized derivatives of the local Maxwellian. As in Proposition~\ref{prop:TLsrc}, by
the same trilinear bound and Lemma~\ref{lem:QmG}, we see that each is $\le C\|
\mathfrak g\|\,\|\varphi''\|_{Y_2}\|\varphi\|_{Y_2}$, with $\|\mathfrak g\|=\|\mathfrak g\|_\flat+\|
\mathfrak g\|_{Y_2}+\|\mathfrak g\|_X$. By Lemma~\ref{lem:Mref} and the same proposition, we have,
for the
$\varpi_*$-part,
\begin{equation}
W_{\alpha\beta}\tbinom\alpha{\alpha'}\tbinom\beta{\beta'}\|
\mathfrak g_{\alpha'\beta'}\|\le\varkappa_{\alpha'\beta'}W_{\alpha''\beta''},\qquad\varkappa_{\alpha'\beta'}:=C_{\alpha'}\delta_y^{|\alpha'|}\Big(\frac{\rad}{r_{\rm
C}}\Big)^{|\beta'|}(|\beta'|+1)^5 ,
\end{equation}
with $\rad=\rad(t)\le\rad_0$ the radius. Since
$(\alpha',\beta')\ne0$, there are at most $(j+1)^2$ multi-indices $\beta'$ of order $j$, and
$4\rad/r_{\rm C}\le\frac12$ (requirement~\ref{rq:UFc}), we obtain $\sum\varkappa_{\alpha'\beta'}\le
C(\delta_y+\rad_0/r_{\rm C})\le C(\delta_y+\rad_0/R_1)$ and
$\sum\varkappa^2_{\alpha'\beta'}\le C(\delta_y+\rad_0/R_1)^2$, with $C$ depending on $L_1$
through $r_{\rm C}=R_1/(120(L_1+1))$. The $\delta\varpi$-part is not bilinear, and we bound it as in
the proof of that proposition, by Taylor's formula and the pointwise composition bound,
and it adds $C\|\delta\varpi\|_{\mathcal A}\le2C\delta_U$ to both sums. For $\alpha'\ne0$ it also
contains terms $t^{-k}\partial^{|\beta'|+k}_\xi\delta\varpi$. We treat these as
in~\eqref{eq:UFgap}, with $\rad_0T_0\ge1$, and they change $\delta_U$ into $C\delta_U$. If
$(\alpha'',\beta'')\ne0$, the sum is $\le C(\delta_y+\rad_0/R_1+\delta_U)\eta\kappa
\mathbb D_\xi\le\frac{\eta c_h}4\kappa\mathbb D_\xi$. If $(\alpha'',\beta'')=0$, the second factor is the
background, with $\|f_\xi\|_{Y_2}\le CK$ by (G4), and by Cauchy--Schwarz we see that the sum is
\begin{equation}
\le C\eta K\kappa\lambda_\xi^{1/2}\Big(\sum\varkappa^2\Big)^{1/2}\mathbb D_\xi^{1/2}\le\frac{\eta
c_h}4\kappa\mathbb D_\xi+C(K)\kappa\lambda_\xi\Big(\delta_y+\frac{\rad_0}{R_1}+\delta_U\Big)^2 .
\end{equation}

\emph{(T8), (T9) and (T12)} are as in Subsection~\ref{sssec:bandproof}. The terms (T8), (T9) and the
kernel part (T12$_0$) of (T12) form $\mathcal T_\xi$ together with the radius term, and the rest of
(T12)
sums to $Ct^{-2}\hat{\mathbb A}^{1/2}\mathbb F$. Since $\mathcal T_\xi$ contains all of (T8), the
constant $C$ in $t^{-2}(C+4\mathfrak p_*)$ could be taken $0$ here. However, we keep it for the cell
inequality, where $a_1$ contains it.

\emph{(T11) Frame.} By Corollary~\ref{cor:FMuse}(b) and~(d), this term is $\le4t^{-2}
\mathfrak p_*\eta\|\varphi\|^2_X$.

\emph{Collecting.} The coercive part $-2\eta c_h\kappa\mathbb D$ of (L1), minus the three
contributions $\frac{\eta c_h}4\kappa\mathbb D$ of (N) and (C), leaves $-\frac54\eta
c_h\kappa\mathbb D\le-\eta c_h\kappa\mathbb D$. For the parts with $(\alpha,\beta)\ne0$ the same
terms occur, restricted to those indices.
\end{proof}

\subsection{The viscous approximation}\label{app:visc}

The viscous approximation is needed for the dispersed layer, Proposition~\ref{prop:disp}, and for
the local theory in cell norms, Proposition~\ref{prop:CULW} and Lemma~\ref{lem:CUtstar}.

We write $\eta_v\in(0,1]$ for the viscosity, $\eta$ being the coefficient of the leaf density, and
use the profile convolutions of Lemma~\ref{lem:LWuniq}: $\ab[h](t,y,v)=\int a(w)h(y+tw,v-w)\dd w$,
and similarly $\bb$, $\cb$. Since $D\cdot\ab[h]=\bb[h]$ and $D\cdot\bb[h]=\cb[h]$, as in its proof,
the equation $\partial_tg=D\cdot(\ab[g]Dg-\bb[g]g)$ has the nondivergence
form $\partial_tg=\ab[g]:D^2g-\cb[g]g$, $D=\partial_v-t\nabla_y$.

\subsubsection{The regularized problem}

Let $t_*>0$. We put
\begin{equation}
\Gamma_c(y,v):=\eps\rho_0(v)e^{-c|y|^2},\qquad\mathcal N_m(h;c):=\sum_{|\mathbf k|\le
m}\sup_{y,v}\frac{|\partial^{\mathbf k}_{y,v}h(y,v)|}{\Gamma_c(y,v)} .
\end{equation}
The regularized problem is
\begin{equation}\label{eq:LWreg}
\partial_tg=\ab[g]:D^2g-\cb[g]g+\eta_v(\Delta_y+\Delta_v)g,\qquad g(t_*)=g_* .
\end{equation}
Its data satisfy, for some $c_0\in(0,1]$ with $\eta_vc_0\le\frac1{16}$ and some $C_*\ge1$,
\begin{equation}\label{eq:LWdata}
g_*\ge0,\qquad\mathcal N_3(g_*;c_0)\le C_* .
\end{equation}
In both applications this follows from~(D) of point~(3) of the proof of Proposition~\ref{prop:CULW},
with $c_0=q_*-\delta/4$, since (D) holds for $|\alpha|\le N_y-3$ and every $\beta$. Moreover,
we have $\eta_vc_0\le\delta/(64q_*)\le\frac1{256}$ by the definitions of $\eta_0$ and of
$\eta_0^{\rm disp}$
(the latter in Proposition~\ref{prop:disp}, with $q_d$ in place of $q_*$). The barrier rate is
\begin{equation}
c_\eta(t):=\frac{c_0}{1+4\eta_vc_0(t-t_*)},\qquad\text{so}\quad
c_\eta'=-4\eta_vc_\eta^2,\quad\frac1{c_\eta(t)}=\frac1{c_0}+4\eta_v(t-t_*) .
\end{equation}
On $[t_*,t_*+1]$, we have $\frac45c_0\le c_\eta\le c_0$, since $4\eta_vc_0\le\frac14$. The last
identity
shows that a restart at a time $t_\flat>t_*$ with the rate $c_\eta(t_\flat)$ produces the same
function $c_\eta$. By (D3), we have $|\nabla\log\rho_0|\le L$ and $|\nabla^2\log\rho_0|\le L$, so
$|\nabla\rho_0|\le L\rho_0$, $|\nabla^2\rho_0|\le(L+L^2)\rho_0$, $|\Delta\rho_0|\le C_L\rho_0$ with
$C_L:=3L+L^2$, $\rho_0(v')\le e^{2L}\rho_0(v)$ for $|v'-v|\le2$, and $\rho_0\le A_F\jap v^{-p}$ with
$p\ge4$. Constants $C$ depend on $\gamma$, $L$, $A_F$ and $p$ unless stated.

\emph{Parabolic estimates with coefficients H\"older in space.} At $\gamma=-2$ the coefficients
$\ab[h]$ are not H\"older continuous in $t$ uniformly in $(y,v)$, by Lemma~\ref{lem:LWlin}(i).
Therefore we use the Schauder estimates in the form that requires H\"older continuity of the
coefficients in the space variables only, uniformly in $t$, with measurable (here continuous)
dependence on $t$: the interior estimate \cite[Thm.~1]{Kn80} and the estimate on strips up to the
initial time \cite[Thm.~4.1]{Lor00}, both listed in (E4) of Inputs~\ref{hyp:E}. The latter
is global.

\emph{The global estimate.} For $\partial_tu=q^{ij}\partial_{ij}u+b^i\partial_iu+cu+f$ on
$[0,T]\times\R^n$, we assume the following:
\begin{itemize}
\item $q,b,c$ are bounded and measurable and H\"older continuous in $x$, with an exponent
$\theta\in(0,1)$, uniformly in $t$ \cite[(H1)]{Lor00};
\item $q\ge c_{\rm L}I$ for a.e.\ $t$ \cite[(H2)]{Lor00};
\item $f$ is measurable with $\sup_t\|f(t)\|_{C^\theta(\R^n)}<\infty$.
\end{itemize}
Then, by \cite[Thm.~4.1]{Lor00} with the weight $p\equiv1$, we have
$\sup_t\|u(t)\|_{C^{2+\theta}(\R^n)}\le
C(\|u(0)\|_{C^{2+\theta}(\R^n)}+\sup_t\|f(t)\|_{C^\theta(\R^n)})$ for every solution $u\in
B([0,T];C^{2+\theta}(\R^n))$ in the sense of \cite[Def.~1.1]{Lor00}. Such a solution is Lipschitz on
$[0,T]\times\R^n$, $D_xu$ and $D^2_xu$ are bounded and continuous up to $t=0$, and the equation
holds a.e.

By the proof of \cite[Thm.~4.1]{Lor00} (cf.\ \cite[(4.31)]{Lor00}), the constant depends only on
$n,\theta,c_{\rm L},T$ and the $L^\infty$ and spatial $C^\theta$ norms of the coefficients, while
the
statement of \cite[Thm.~4.1]{Lor00} only says that it does not depend on $(u,u_0,f)$. For $T\le1$ we
take the constant for $T=1$. To this end we extend the coefficients to $(T,1]$ by $q=I$, $b=0$,
$c=0$, $f=0$, and
the solution by \cite[Thm.~1.2]{Lor00}, with the data $u(T)\in C^{2+\theta}$ (on $(T,1]$ this is the
heat equation).

Our coefficients are defined on all of $\R^6$, bounded, and H\"older continuous in $(y,v)$ uniformly
in $t$, as noted above, so they need no extension.

\emph{The local form.} It follows from the global estimate by a cutoff, as in Step~1 of the proof of
\cite[Thm.~4.1]{Lor00}, followed by an absorption over nested balls. For $1\le\rho<\rho'\le2$ we let
$\zeta(y,v)$ be equal to $1$ on $\rho B$ and supported in $\rho'B$, with $|\partial^j\zeta|\le
C(\rho'-\rho)^{-j}$ for $j\le3$. Then $\zeta u$ solves the same equation with the forcing $\zeta
f-2A\nabla\zeta\cdot\nabla u-(A:D^2\zeta)u$. Applying the global estimate to $\zeta u$, together
with the interpolation inequality $\|\nabla
u\|_{C^\alpha(\rho'B)}+\|u\|_{C^\alpha(\rho'B)}\le\epsilon\|D^2u\|_{C^\alpha(\rho'B)}+C\epsilon^{-2}\sup_{\rho'B}|u|$
with $\epsilon=c(\rho'-\rho)^{3}$, we obtain $\varphi(\rho)\le\frac12\varphi(\rho')+
\mathcal A(\rho'-\rho)^{-9}$. Here $\varphi(\rho):=\sup_t\|D^2u(t)\|_{C^\alpha(\rho B)}$, and $\mathcal A$ is
$C$ times the right side of the estimate displayed below. If $\varphi$ is bounded on $[1,2]$, then
by iterating along $\rho_0=1$, $\rho_{j+1}=\rho_j+(1-\tau)\tau^j$ with $\tau=2^{-1/18}$, we get
$\varphi(1)\le2^{-J}\varphi(\rho_J)+\mathcal A(1-\tau)^{-9}\sum_{j<J}2^{-j/2}$. Letting
$J\to\infty$,
we obtain $\varphi(1)\le C\mathcal A$, which is the standard iteration lemma.

For $\partial_tu=A^{ij}\partial_{ij}u+cu+f$ in $\R^6$ with $A\ge\nu I$ and
$\|A\|_{C^\alpha_x}+\|c\|_{C^\alpha_x}\le\Lambda$ uniformly in $t$, the two estimates imply that, on
$Z=B\times(s,t]$ with $B$ a unit ball,
\begin{equation}
\sup_{(s',t]}\big(\|D^2_xu\|_{C^\alpha_x(B)}+\|\partial_tu\|_{L^\infty(B)}\big)\le
C\Big(\sup_{2B\times(s,t]}|u|+\sup_{(s,t]}\|f\|_{C^\alpha_x(2B)}\Big)
\end{equation}
for $s'>s$, with $C$ depending also on $s'-s$ (interior), and the same with $s'=s$ and
$\|u(s)\|_{C^{2+\alpha}(2B)}$ added on the right (up to the initial time). The Krylov--Safonov
estimate requires no regularity of the coefficients. These are a priori estimates. The absorption of
the interpolated terms $\nabla\zeta\cdot\nabla u$ and $(D^2\zeta)u$ over nested balls requires
$\sup_t\|D^2u\|_{C^\alpha_x(2B)}<\infty$. Up to the initial time, $\zeta u$ must moreover be a
solution in the sense of \cite[Def.~1.1]{Lor00}, that is, $D^2_xu$ continuous up to $t=s$ and
$\partial_tu$ bounded on $2B$. The estimates are applied in two ways. First, to the smooth
approximants $u_k$ of Lemma~\ref{lem:LWlin}(ii), which are smooth up to $t=t_*$ on $B_k$ and whose
bounds pass to the limit. Second, when the equation is differentiated, to difference quotients in
$(y,v)$ of the previous derivatives. These are classical solutions with forcing bounded in
$C^\alpha_x$ uniformly in the increment. For them the conditions above hold by the previous step of
the induction: continuity in $t$ of $D^2_x$ follows from the $C^{2+\alpha}_x$ bound and the
Lipschitz bound in $t$ by interpolation, as there.

\begin{lem}[Linear problem and barrier]\label{lem:LWlin}
Let $A\ge1$, $c_h\in(0,c_0]$ and $T\in(t_*,t_*+1]$, and let $h\ge0$ be continuous on
$[t_*,T]\times\R^6$ with $h(t)\le A\Gamma_{c_h}$.
\begin{enumerate}[label=(\roman*)]
\item The coefficients satisfy $|\ab[h]|+|\bb[h]|+|\cb[h]|\le CA$, $-\cb[h]\ge0$, and the null
structure $\ab[h](t,y,v)\,y\cdot y\le CA/c_h$. If $h$ is H\"older continuous in $(t,y,v)$, with
seminorm $\le R_h$, then $\ab[h]$ and $\cb[h]$ are H\"older continuous in $(y,v)$, uniformly in $t$,
with a smaller exponent and seminorms $\le C(A,R_h)$, and they are locally H\"older continuous in
$(t,y,v)$. For $\gamma<-2$ they are H\"older continuous in $(t,y,v)$, with the exponent
$\theta\alpha$, $\theta:=\min\{\frac18,|\gamma+2|\}$, and seminorm $\le C(A,R_h,t_*)$. For
$\gamma=-2$ the H\"older seminorm in $t$ is not bounded uniformly in $(y,v)$ in general.
\item If $h$ is H\"older continuous and $g_*$ satisfies~\eqref{eq:LWdata}, the linear problem
\begin{equation}
\partial_tu=\ab[h]:D^2u-\cb[h]u+\eta_v(\Delta_y+\Delta_v)u,\qquad u(t_*)=g_*,
\end{equation}
has a unique bounded solution $u$ on $[t_*,T]$ that is continuous on $[t_*,T]\times\R^6$ and
classical for $t>t_*$, and $u\ge0$.
\item Let $M_h:=\sup\big(|\ab[h]|+|\cb[h]|+\ab[h]y\cdot y\big)$, and let $u$ be a bounded solution
of the linear problem on $[t_*,T]$, continuous up to $t_*$ and classical for $t>t_*$, with
$u(t_*)\le C_*\Gamma_{c_0}$ (for instance the solution of~(ii)). Then
\begin{equation}
u(t)\le C_*e^{K_{\rm b}(t-t_*)}\Gamma_{c_\eta(t)}\qquad\text{on }[t_*,T],\qquad K_{\rm
b}:=C(1+M_h)\big(1+(t_*+1)^2\big).
\end{equation}
\end{enumerate}
\end{lem}

\begin{proof}
(i) Since $\eps\le1$ and $p\ge4>3$, we have $\int|w|^{\gamma+2}\rho_0(v-w)\dd w\le C$ and
$\int_{|w|\le1}|w|^\gamma\dd w+\sup\rho_0\le C$ for $\gamma>-3$, and thus the bound
$h\le A\eps\rho_0(v-w)$ implies the first claim. For $\gamma=-3$, we have $\cb[h]=-8\pi h$. For
$\gamma>-3$, we have $-\cb[h]=2(\gamma+3)\int|w|^\gamma h(y+tw,v-w)\dd w\ge0$. By
Lemma~\ref{lem:Bnull}, we have $a(w)y=a(w)y''$
with $y''=y+tw$, and thus
\begin{equation}
\ab[h]y\cdot y=\int a(w)y''\cdot y''\,h(y'',v-w)\dd w\le
A\sup_{s\ge0}s^2e^{-c_hs^2}\int|w|^{\gamma+2}\rho_0(v-w)\dd w\le\frac{CA}{c_h} .
\end{equation}
For the H\"older bounds, we have, at the arguments $x=(t,y+tw,v-w)$ and $x'=(t',y'+t'w,v'-w)$,
\begin{equation}
|h(x)-h(x')|\le\min\{R_h|x-x'|^\alpha,2A\eps(\rho_0(v-w)+\rho_0(v'-w))\}\le(R_h|x-x'|^\alpha)^{\theta}(2A\eps)^{1-\theta}\big(\rho_0(v-w)+\rho_0(v'-w)\big)^{1-\theta},
\end{equation}
with $\theta\in(0,\frac18]$. We take $\theta=\frac18$ for the increments in $(y,v)$ and for the
local bound, and $\theta=\min\{\frac18,|\gamma+2|\}$ for the increment in $t$ when $\gamma<-2$.
Moreover, $\rho_0^{1-\theta}\le C\jap\cdot^{-7p/8}$, with $\frac78p>3$ as $p\ge4$, is integrable
against $|w|^{\gamma+2}$ and against $|w|^\gamma\mathbf 1_{|w|\ge1}$, and near $w=0$, $|w|^\gamma$
is integrable. For the increments in $(y,v)$, we have $|x-x'|=|y-y'|+|v-v'|$, with no factor $|w|$.
This yields the uniform spatial bound, and the local bound in $(t,y,v)$ follows on bounded sets of
$v$, where $(1+|w|)^{\theta\alpha}$ is absorbed by $\rho_0(v-w)^{1-\theta}$. For the increment in
$t$, we have $|x-x'|\le(1+|w|)|t-t'|$, and for $\gamma<-2$, $(1+|w|)^{\theta\alpha}|w|^{\gamma+2}\le
C$ for
$|w|\ge1$, since $\theta\alpha\le|\gamma+2|$. (At $\gamma=-2$, the choice $h=A\Gamma_{c_0}$ leads to
$\operatorname{tr}\ab[h](t,y,v)=2A\eps F(t,y+tv)$ with $F(t,z)=\int\rho_0(u)e^{-c_0|z-tu|^2}\dd u$,
and the time increments $R/|v|$ with $y=z_0-tv$ show that no uniform H\"older bound in $t$ holds.)

(ii) The second-order part is $\ab[h]:D^2+\eta_v(\Delta_y+\Delta_v)$. Its symbol
$\ab[h](\zeta_v-t\zeta_y)\cdot(\zeta_v-t\zeta_y)+\eta_v|\zeta|^2\ge\eta_v|\zeta|^2$ shows uniform
parabolicity, and by~(i) the coefficients are bounded by $CA(1+t^2)$ and H\"older continuous in
$(y,v)$, uniformly in $t$.

\emph{Maximum principle.} Let $\mathfrak w$ be bounded, continuous on $[t_*,T]\times\R^6$, classical
for $t>t_*$, with $\partial_t\mathfrak w\le\mathcal L\mathfrak w$ for the operator $\mathcal L$ of
the problem and $\mathfrak w(t_*)\le0$. Setting $V:=-\cb[h]\in[0,CA]$, we see that $\tilde{
\mathfrak w}:=e^{-CA(t-t_*)}\mathfrak w$ satisfies $\partial_t\tilde{\mathfrak w}\le
\mathcal L_0\tilde{\mathfrak w}+(V-CA)\tilde{\mathfrak w}$, where $\mathcal L_0$ is the second-order part.
For $\rho\ge1$ we put $\beta_\rho:=\rho^{-2}\sup|\mathfrak w|\,\big(|y|^2+|v|^2+K_0(t-t_*)\big)$
with
$K_0:=2\sup\operatorname{tr}(\text{matrix of }\mathcal L_0)\le CA(1+T^2)+12$. Then we have
$\partial_t\beta_\rho-\mathcal L_0\beta_\rho-(V-CA)\beta_\rho\ge0$, and $\tilde{
\mathfrak w}\le\beta_\rho$ on the parabolic boundary of $\{|y|^2+|v|^2<\rho^2\}\times(t_*,T]$. By
the maximum
principle on this bounded cylinder, we obtain $\tilde{\mathfrak w}\le\beta_\rho$, and letting
$\rho\to\infty$, we get $\mathfrak w\le0$. This implies uniqueness among bounded solutions, $u\ge0$
(we apply it to $-u$), and the comparison principle used in~(iii).

\emph{Existence.} We extend $h$ to $t<t_*$ by $h(t)=h(t_*)$, and mollify $\ab[h]$ and $\cb[h]$ in
$(t,y,v)$ at the scale $1/k$. Similarly, we mollify $g_*$ at the scale $1/k$ and multiply it by a
cutoff equal to $1$ on
$B_k$ and supported in $B_{3k/2}$. On $B_{2k}$, with zero boundary values, the problem with smooth
coefficients has a solution $u_k$ that is smooth on $[t_*,T]\times\overline{B_{2k}}$ \cite[Ch.~IV,
Thm.~5.2]{LSU} (applied for every nonintegral $l>0$), since the data vanish near $\partial B_{2k}$
and the compatibility conditions of all orders hold, and $0\le u_k\le e^{CA(t-t_*)}\sup g_*$ by the
maximum principle.

The mollified coefficients have spatial H\"older norms bounded uniformly in $k$ and $t$.
By the interior estimate with coefficients H\"older in space \cite{Kn80}, $D^2u_k$ is bounded in
$C^\alpha_x$ and $\partial_tu_k$ in $L^\infty$ on compact subsets of $(t_*,T]\times\R^6$. By
interpolation
between $|u_k(t)-u_k(t')|\le C|t-t'|$ and the $C^{2+\alpha}_x$ bound, $D^2u_k$ is also H\"older
continuous in $t$, with the exponent $\alpha/(2+\alpha)$, uniformly in $k$. Finally, by the
Krylov--Safonov estimate up to the initial time, with H\"older data \cite[Cor.~7.45]{Lie96}, the
H\"older norm of $u_k$ on compact subsets of $[t_*,T]\times\R^6$ is bounded uniformly in $k$. Its
condition (A) \cite[(6.48)]{Lie96}, read as a lower bound for the measure of the part of the
cylinder that lies outside the domain, which is the form in which it is used in the proof of
\cite[Thm.~6.32]{Lie96}, holds at the initial time, since the backward cylinders with top center at
$t=t_*$ lie outside the domain. Moreover, the boundary modulus may be taken $\sigma(r)=Cr^\delta$
with
$\delta$ small.

By Arzel\`a--Ascoli and a diagonal argument, we can extract a subsequence which converges, with
$D^2u_k$ locally
uniformly, to a bounded $u$, continuous up to $t_*$ with $u(t_*)=g_*$. Moreover $\partial_tu_k$
converges weakly-$*$, and passing to the limit in the equation shows that $u$ is a classical
solution for $t>t_*$ ($\partial_tu$ equals the continuous right side).

(iii) Let $c=c_\eta(t)$ and $\Phi:=C_*e^{K_{\rm b}(t-t_*)}\Gamma_c$. Then
$\partial_t\Gamma_c=-c'|y|^2\Gamma_c=4\eta_vc^2|y|^2\Gamma_c$ and
$\Delta_y\Gamma_c=(4c^2|y|^2-6c)\Gamma_c\le4c^2|y|^2\Gamma_c$, so the viscous $y$-term is exactly
compensated. Next, we have $D\Gamma_c=\zeta\,\Gamma_c$ with $\zeta:=\nabla\log\rho_0+2tcy$, and
$D^2\Gamma_c=\big(\zeta\otimes\zeta+\nabla^2\log\rho_0-2t^2cI\big)\Gamma_c$. Since $\ab[h]\ge0$,
$c\le1$ and
$\ab[h]\zeta\cdot\zeta\le2\ab[h]\nabla\log\rho_0\cdot\nabla\log\rho_0+8t^2c^2\ab[h]y\cdot y$, we
obtain
\begin{equation}
\ab[h]:D^2\Gamma_c\le C(1+t^2)M_h\Gamma_c,\qquad\eta_v\Delta_v\Gamma_c=\eta_v\frac{\Delta\rho_0}{\rho_0}\Gamma_c\le
C_L\Gamma_c,\qquad-\cb[h]\Gamma_c\le M_h\Gamma_c .
\end{equation}
Hence we have $\partial_t\Phi-\mathcal L\Phi\ge\big(K_{\rm b}-C(1+t^2)M_h-C_L-M_h\big)\Phi\ge0$ on
$[t_*,T]$
for $K_{\rm b}$ as stated. Since $u(t_*)\le C_*\Gamma_{c_0}=\Phi(t_*)$ and $u-\Phi$ is bounded, the
comparison principle of~(ii) implies that $u\le\Phi$.
\end{proof}

\begin{lem}[The regularized nonlinear problem]\label{lem:LWregex}
\leavevmode
\begin{enumerate}[label=(\roman*)]
\item Let $g_*$ satisfy~\eqref{eq:LWdata}. There is $\tau_\eta=\tau_\eta(C_*,c_0,t_*)\in(0,1]$ such
that~\eqref{eq:LWreg} has a solution $g$ on $[t_*,t_*+\tau_\eta]$, continuous on
$[t_*,t_*+\tau_\eta]\times\R^6$, $C^\infty$ for $t>t_*$, with $0\le g\le2C_*\Gamma_{c_\eta(t)}$ and,
for every $\delta'\in(0,\frac{c_0}8]$, $\mathcal N_2\big(g(t);c_\eta(t)-\delta'\big)\le
C(C_*,c_0,t_*,\eta_v,\delta')$ on $[t_*,t_*+\tau_\eta]$.
\item \emph{Uniqueness.} Let $g_*$ be as in~(i), and let $g_1,g_2$ be solutions
of the same problem on $[t_*,T]$, $T\le t_*+1$, continuous up to $t_*$ and classical for $t>t_*$,
with $0\le g_i\le A\Gamma_{c_\eta(t)}$ for some $A$, and $\sup_{[t_*,T]}
\mathcal N_2(g_2(t);c_\eta(t)-\frac{c_0}8)<\infty$. Then $g_1=g_2$.
\item \emph{Derivative bounds.} Let $c_0\in(0,1]$ with $\eta_vc_0\le\frac1{16}$, and let $g\ge0$
solve the equation in~\eqref{eq:LWreg} on $[t_*,T)$, $T\le t_*+1$, be continuous up to $t_*$ and
classical for $t>t_*$, with $A_0:=\sup_{[t_*,T)}\mathcal N_0(g(t);c_\eta(t))<\infty$. Then $g$ is
$C^\infty$ for $t>t_*$, and for every $m$, $\delta'\in(0,\frac{c_0}8]$ and $t\in(t_*,T)$,
\begin{equation}\label{eq:LWregex-m}
\mathcal N_m\big(g(t);c_\eta(t)-\delta'\big)\le C(m,\eta_v,t_*,t-t_*,A_0,c_0,\delta'),
\end{equation}
nonincreasing in $t-t_*$.
\item \emph{Continuation.} Under the hypotheses of~(iii), $g$ extends to a solution on $[t_*,T+s]$
for some $s>0$.
\end{enumerate}
\end{lem}

\begin{proof}
(i) \emph{The fixed point.} Let $A:=2C_*$ and $c_h:=\frac45c_0$, so that
$\Gamma_{c_\eta(t)}\le\Gamma_{c_h}$ and $M_h\le CA(1+1/c_0)$ for every $h\le A\Gamma_{c_\eta}$
by Lemma~\ref{lem:LWlin}(i). Let $K_{\rm b}=K_{\rm b}(C_*,c_0,t_*)$ be the constant of
part~(iii) of that lemma with this bound, and let $\tau_\eta:=\min\{1,(\log2)/K_{\rm b}\}$. For
$R_0>0$ we let $\mathcal K$ be the set of continuous $h$ on $[t_*,t_*+\tau_\eta]\times\R^6$ with
\begin{equation}
0\le h\le2C_*\Gamma_{c_\eta(t)},\qquad[h]_{C^{\alpha_0,\alpha_0/2}}\le R_0 ,
\end{equation}
where $\alpha_0\in(0,1)$ and $[\cdot]$ is the global parabolic H\"older seminorm. For $h\in
\mathcal K$, we let $\mathcal Th:=u$ be the solution of part~(ii) of that lemma with
$T=t_*+\tau_\eta$.
By part~(iii), we have $u\le C_*e^{K_{\rm b}\tau_\eta}\Gamma_{c_\eta}\le2C_*\Gamma_{c_\eta}$.
By the Krylov--Safonov estimate, interior and up to the initial time \cite[Ch.~VII]{Lie96}, $u$ is
H\"older continuous with an exponent $\alpha_0$ and a seminorm on unit cylinders that depend only on
the ellipticity $\eta_v$, the bound $CA(1+(t_*+1)^2)$ of the coefficients, $\sup|u|\le2C_*$ and the
H\"older norm of $g_*$, which is $\le C_*$ by~\eqref{eq:LWdata}. For distances $\ge1$, the seminorm
is bounded by $2\sup|u|$. This fixes $\alpha_0$ and $R_0$, independently of $h$, and we have $
\mathcal T\mathcal K\subset\mathcal K$.

We equip $\mathcal K$ with the norm $\|h\|_w:=\sup|h|\,\Gamma_{c_\eta(t)}^{-1/2}$. Then
$\mathcal K$ is convex and closed, and it is compact, since Arzel\`a--Ascoli applies on compact sets
and $|h|\Gamma_{c_\eta}^{-1/2}\le2C_*\Gamma_{c_\eta}^{1/2}\to0$ at infinity, uniformly on
$\mathcal K$.
Next, we show that the map $\mathcal T$ is continuous. If $h_k\to h$ in $\mathcal K$, then
$h_k\to h$ locally uniformly with the majorant $2C_*\Gamma_{c_\eta}$, and thus, by dominated
convergence in $w$, we have $\ab[h_k]\to\ab[h]$ and $\cb[h_k]\to\cb[h]$ locally uniformly.
Moreover, their spatial H\"older norms are bounded uniformly in
$k$, as above. Every subsequence of $\mathcal Th_k$ has a further subsequence
converging in $\mathcal K$, and, by the interior estimate of \cite{Kn80}, with $D^2_{y,v}$ locally
uniformly for $t>t_*$. The limit is a bounded solution of the
linear problem for $h$, continuous up to $t_*$ with the data $g_*$, and thus it is $\mathcal Th$ by
the uniqueness for the linear problem. Hence $\mathcal Th_k\to\mathcal Th$. By Schauder's fixed
point theorem, we obtain $g=\mathcal Tg\in\mathcal K$.

\emph{Regularity, up to $t_*$.} We show by induction on $k\ge0$: for every
$\delta'\in(0,\frac{c_0}8]$ there is $C_k(\delta')$ with $|\partial^{\mathbf k}_{y,v}g|\le
C_k(\delta')\Gamma_{c_\eta-\delta'}$ for $|\mathbf k|\le k$, on $[t_*,t_*+\tau_\eta]$ for $k\le2$
and on $[t_*+s,t_*+\tau_\eta]$, with $C_k$ depending on $s>0$, for every $k$. For $k=0$, this
follows from
$g\in\mathcal K$.

We include in the claim the spatial H\"older bounds $\sup_t[\partial^{\mathbf k}
g]_{C^\alpha_x(B)}\le
C_k(\delta')\Gamma_{c_\eta(t)-\delta'}(\text{center of }B)$, $|\mathbf k|\le k$, for unit balls $B$.

For the step $k=0$, $g\in\mathcal K$ suffices. Indeed, by Lemma~\ref{lem:LWlin}(i), $\ab[g]$,
$\cb[g]$ are H\"older continuous in $(y,v)$, uniformly in $t$ and in the center, and the estimate up
to the initial time \cite[Thm.~4.1]{Lor00} in the local form above, with $g_*\in C^{2+\alpha}$
by~\eqref{eq:LWdata}, implies the claim for $k=2$ on $[t_*,t_*+\tau_\eta]$.

For $k\ge1$, we suppose that the claim holds for $k+1$. Then
$\partial^{\mathbf k'}\ab[g]=\ab[\partial^{\mathbf k'}g]$ and
$\partial^{\mathbf k'}\cb[g]=\cb[\partial^{\mathbf k'}g]$, $|\mathbf k'|\le k$, are bounded and
H\"older continuous in $(y,v)$, uniformly in $t$ and in the center, by part~(i) of that lemma with
$\partial^{\mathbf k'}g$ in place of $h$, whose proof requires only $|h|\le A\Gamma_c$,
$c\ge\frac35c_0$, and the spatial H\"older seminorm of $h$. Since $D=\partial_v-t\nabla_y$ has
coefficients independent of $(y,v)$, $u=\partial^{\mathbf k}_{y,v}g$, $|\mathbf k|=k$, solves the
linear equation
of $g$ with the forcing
\begin{equation}
f_{\mathbf k}=\sum_{0<\mathbf k'\le\mathbf k}\binom{\mathbf k}{
\mathbf k'}\big(\ab[\partial^{\mathbf k'}g]:D^2\partial^{\mathbf k-\mathbf k'}g-\cb[\partial^{
\mathbf k'}g]\,\partial^{\mathbf k-\mathbf k'}g\big),
\end{equation}
which contains derivatives of $g$ of order $\le k+1$ only, so that
$\sup_t\|f_{\mathbf k}\|_{C^\alpha_x(2B)}$ is bounded by the claim for $k+1$. We apply the interior
estimate
\cite{Kn80} on $B\times[t_*+s/2,t]$, as explained in the paragraph on parabolic estimates
above, to the difference quotients $\epsilon^{-1}\big(\partial^{\mathbf k_0}g(\cdot+\epsilon
e)-\partial^{\mathbf k_0}g\big)$, $|\mathbf k_0|=k-1$, $|e|=1$. Their forcing is bounded in
$C^\alpha_x$ uniformly in $\epsilon$, by the claim for $k+1$. Letting $\epsilon\to0$, we see that
$\|D^2u\|_{C^\alpha_x(B)}$ and $\|\partial_tu\|_{L^\infty(B)}$ are bounded by
$\sup_{2B\times[t_*+s/2,t]}|u|$ and $\sup_t\|f_{\mathbf k}\|_{C^\alpha_x(2B)}$. This proves the
claim
for $k+2$, at the rate $c_\eta-2\delta'$. Since $\delta'$ is arbitrary, the rates are restored by
renaming. In particular
$g\in C^\infty$ for $t>t_*$, and the bound on $\mathcal N_2$ holds.

(ii) Let $d:=g_1-g_2$. Then
\begin{equation}
\partial_td=\ab[g_1]:D^2d-\cb[g_1]d+\eta_v(\Delta_y+\Delta_v)d+\ab[d]:D^2g_2-\cb[d]g_2 .
\end{equation}
We put $c_1:=\frac34c_0$ and $c_\flat(t):=c_1/(1+4\eta_vc_1(t-t_*))$, so that
$c_\flat'=-4\eta_vc_\flat^2$ and $c_\flat\ge\frac35c_0$. Then
\begin{equation}
c_\eta(t)-c_\flat(t)=\frac{c_0-c_1}{(1+4\eta_vc_0(t-t_*))(1+4\eta_vc_1(t-t_*))}\ge\frac{c_0/4}{(1+4\eta_vc_0)^2}\ge\frac{4c_0}{25}\ge\frac{c_0}8
,
\end{equation}
since $\eta_vc_0\le\frac1{16}$. Let $\mathcal D(T'):=\sup_{[t_*,T']}\sup|d|/\Gamma_{c_\flat}$, which
is finite, since $|d|\le2A\Gamma_{c_\eta}\le2A\Gamma_{c_\flat}$. As in the linear lemma, we have
$|\ab[d]|\le C\mathcal D(T')$, and $|\cb[d]|\le C\mathcal D(T')$ for $\gamma>-3$. Moreover,
$|D^2g_2|\le
C(1+t^2)\Gamma_{c_\eta-c_0/8}\le C(1+t^2)\Gamma_{c_\flat}$ by hypothesis. For $\gamma=-3$, we have
$|\cb[d]g_2|=8\pi|d|g_2\le CA\mathcal D(T')\Gamma_{c_\flat}$. Thus the last two terms are bounded by
$C_d\mathcal D(T')\Gamma_{c_\flat}$ on $[t_*,T']$, with $C_d$ independent of $T'$. The computation
of part~(iii) of that lemma, with $h=g_1$ and the rate $c_\flat$ in place of $c_\eta$ (it rests on
$c_\flat'=-4\eta_vc_\flat^2$, $c_\flat\le1$ and the bound $M_{g_1}$, which depends on $A$ and
$c_0$), shows that $\Phi:=C_d\mathcal D(T')(t-t_*)e^{K'(t-t_*)}\Gamma_{c_\flat(t)}$ satisfies
$\partial_t\Phi-\mathcal L\Phi\ge C_d\mathcal D(T')e^{K'(t-t_*)}\Gamma_{c_\flat}$, where $
\mathcal L$ is the operator of $g_1$. By the comparison principle, we have $\pm d\le\Phi$ on
$[t_*,T']$, which implies that $
\mathcal D(T')\le C_d(T'-t_*)e^{K'}\mathcal D(T')$, and thus $\mathcal D(T')=0$ if
$C_d(T'-t_*)e^{K'}<1$.
Repeating this on consecutive intervals of this length, we obtain $d=0$ on $[t_*,T]$.

(iii) By the Krylov--Safonov estimate on the cylinders $B\times[t-r^2,t]$,
$r:=\min\{1,(s/2)^{1/2}\}$, contained in $[t_*+s/2,T)\times\R^6$, $g$ is H\"older continuous there,
with a bound $C(s,\eta_v,t_*,A_0)\sup|g|$ over the doubled cylinder. Then, by
Lemma~\ref{lem:LWlin}(i), $\ab[g],\cb[g]$ are H\"older continuous in $(y,v)$ uniformly in $t$, and
locally H\"older continuous in $(t,y,v)$. By the classical local Schauder theory ((PS)(ii), together
with the
solvability of the first boundary value problem on small cylinders with continuous boundary values
and the maximum principle), $g\in C^{2+\alpha',1+\alpha'/2}$ on compact subsets of
$(t_*,T)\times\R^6$ for some $\alpha'>0$. This is the a priori finiteness of
$\sup_t\|D^2g\|_{C^\alpha_x(2B)}$ that the absorption in the interior estimate requires (the
paragraph on parabolic estimates covers only approximants and difference quotients). Then, by the
induction of~(i), with the interior estimate \cite{Kn80} only, we obtain $|\partial^{\mathbf k}
g|\le
C(m,\eta_v,t_*,s,A_0,\delta')\Gamma_{c_\eta-\delta'}$ for $|\mathbf k|\le m$ and $t\ge t_*+s$. The
constant may be taken nonincreasing in $s$.

(iv) We fix $t_\flat\in[\frac{t_*+T}2,T)$ (a local time, not the end $t_1$ of the initial layer).
By~(iii) with $\delta'=\frac{c_0}8$, we have $\mathcal N_3(g(t_\flat);c_0')\le C_3$ with
$c_0':=c_\eta(t_\flat)-\frac{c_0}8\in[\frac58c_0,c_0]$ and $C_3$ independent of $t_\flat$. Moreover,
$\eta_vc_0'\le\frac1{16}$. By~(i) with $(t_\flat,g(t_\flat),c_0',C_3)$ in place of
$(t_*,g_*,c_0,C_*)$, there is a solution $\tilde g$ on $[t_\flat,t_\flat+\tau']$, with
$\tau'=\tau_\eta(C_3,c_0',t_\flat)$ bounded below uniformly for $t_\flat\in[\frac{t_*+T}2,T)$, since
$K_{\rm b}$ depends on $c_0'$ only through $1/c_0'\le\frac8{5c_0}$. Its rate $\tilde c_\eta$
satisfies $\tilde c_\eta\le c_\eta$ on $[t_\flat,t_\flat+\tau']$, and thus $g\le
A_0\Gamma_{c_\eta}\le
A_0\Gamma_{\tilde c_\eta}$ there, and $\tilde g$ satisfies the $\mathcal N_2$ bound of~(i) at the
rate $\tilde c_\eta-c_0'/8$. By~(ii) at $t_\flat$, we get $g=\tilde g$ on $[t_\flat,T)$. Choosing
$t_\flat>T-\tau'$, we see that $\tilde g$ extends $g$ beyond $T$.
\end{proof}

\subsubsection{The viscosity in a cell}

Let $g$ be a smooth solution of the equation in~\eqref{eq:LWreg} whose derivatives have Gaussian
bounds at a rate $c>q(t)$ on bounded time intervals, as given by Lemma~\ref{lem:LWregex}(iii), and
let $h=h_{\alpha\beta}$ be a family as in Subsection~\ref{ssec:cells}, computed from $g$, with
$\psi:=T_p^{-1}h^{(x)}$, $p=p_g(t,\xi)$. In the leaf variables $(\xi,z)$, we have
$\partial_v=\partial_\xi$
and $\partial_y=\partial_z+t^{-1}\partial_\xi$ (by~(F1)). These vector fields have constant
coefficients, $\dd y\dd v=\dd\xi\dd z$, and
\begin{equation}
\Delta_y+\Delta_v=\sum_i\big(L_{0,i}^2+L_{1,i}^2\big),\qquad L_{0,i}:=\partial_{\xi_i},\quad
L_{1,i}:=\partial_{z_i}+t^{-1}\partial_{\xi_i} .
\end{equation}
Let $0<\Omega=\Omega(\xi)\le1$ be a $C^2$ function with $|\nabla\log\Omega|\le4C_1^UL'(S_\sharp+2)$
and $|\nabla^2\log\Omega|\le C_\Omega$ (for instance $\Omega\equiv1$, or the frozen weight
$\Omega^\epsilon$ of point~(2) of the proof of Proposition~\ref{prop:CULW}), and we put
\begin{equation}
\tilde\Omega:=\chi_c^2\,\Omega\,(A_p/A_{\rm ref})^2m_x^2\ \text{ on
}\R^3_\xi\times\R^3_z,\qquad\Psi:=\rad_0^{-3}\chi_c^2\,\Omega\,\chi(s_\xi),\qquad
\mathcal P(p;u,u'):=\ip{T_p^{-1}u}{T_p^{-1}u'}_\natural ,
\end{equation}
with $\chi\in C^\infty$, $|\chi'|+|\chi''|\le C$, as in that proposition. The moments of
$T_p^{-1}u$ are invertible linear combinations of those of $u$ (the change of variables
$z=y_0+\theta^{1/2}\hat z$ is affine), so $u$ is moment-free in one frame if and only if it is in
every frame. By Lemma~\ref{lem:leafeq}(i), we have $\partial_qT_p^{-1}u=-\mathcal G_p[q]T_p^{-1}u$
for a
variation $q$ of $p$, and $\mathcal G_p[q]=q_\alpha-q_\beta\hat
z\cdot\nabla-e^{-p_\beta}q_a\cdot\nabla$ preserves moment-freeness. We put
$C_\varpi^{(2)}:=\sup|\nabla^2_\xi p_g|$.

\begin{lem}[The viscosity in a cell]\label{lem:LCvisc}
\leavevmode
\begin{enumerate}[label=(\roman*)]
\item \emph{Local part.}
\begin{equation}
2\eta_v\!\iint\tilde\Omega\,h\,(\Delta_y+\Delta_v)h\le-\eta_v\!\iint\tilde\Omega\big(|\nabla_\xi
h|^2+|(\nabla_z+t^{-1}\nabla_\xi)h|^2\big)+\eta_v\!\iint\tilde\Omega\big(32q^2|z|^2+C\rad_0^{-2}\big)|h|^2
.
\end{equation}
\item \emph{Nonlocal part.} With $\psi':=T_p^{-1}(\nabla_\xi h)^{(x)}$,
\begin{equation}
2\eta_v\!\int\Psi\,\mathcal P\big(p;h,(\Delta_y+\Delta_v)h\big)\dd\xi\le-\eta_v\!\int\Psi\,\nn{\psi'}^2_\natural\dd\xi+C_V\eta_v\!\int\rad_0^{-3}\chi_c^2\Omega\Big(\rad_0^{-2}\|\psi\|^2_X+
\mathbf 1_{\{s_\xi\ge S_b\}}\|\psi\|^2_{Y_2}\Big)\dd\xi .
\end{equation}
\end{enumerate}
Here $X=X(\hat m_\xi)$, $Y_2=Y_2(\hat m_\xi)$, $t\ge T_0$, and the constants $C$, $C_V$ depend on
$(\theta_\pm,q_\pm,\eps_1)$, $C_\varpi$, $\rad_0C_\varpi^{(2)}$, $\rad_0^2C_\Omega$, $L'$,
$S_b$ and the bounds on $\chi'$, $\chi''$, but not on $\eta_v$.
\end{lem}

\begin{proof}
There are no boundary terms in the integrations by parts below, by the Gaussian bounds of $g$, whose
rate exceeds the exponent of $m_x$, and the exponential decay of $\chi_c$.

(i) For $L\in\{L_{0,i},L_{1,i}\}$, which are divergence-free in $(\xi,z)$, one integration by parts
yields $\iint\tilde\Omega hL^2h=-\iint\tilde\Omega|Lh|^2-\iint(L\log\tilde\Omega)\tilde\Omega\,hLh$,
and by Young's inequality we get
\begin{equation}
2\iint\tilde\Omega hL^2h\le-\iint\tilde\Omega|Lh|^2+\iint\tilde\Omega(L\log\tilde\Omega)^2|h|^2 .
\end{equation}
By the proof of Lemma~\ref{lem:CUtr} and requirement~\ref{rq:CUb}, we have
$|\nabla_\xi\log\tilde\Omega|\le3/\rad_0$ (the factor $\Omega$ is covered by the bound on
$|\nabla\log\Omega|$). Since $m_x(z)=A_{\rm ref}(\xi-z/t)e^{q|z|^2}$ and $|\nabla\log A_{\rm
ref}|\le L$ by Corollary~\ref{cor:LWPa}(a4), we have $\nabla_z\log\tilde\Omega=4qz-2t^{-1}\nabla\log
A_{\rm
ref}$, and thus
$\sum_i(L_{1,i}\log\tilde\Omega)^2\le(4q|z|+C\rad_0^{-1})^2\le32q^2|z|^2+C\rad_0^{-2}$, and
$\sum_i(L_{0,i}\log\tilde\Omega)^2\le9\rad_0^{-2}$. It remains to sum over $L$ and multiply by
$\eta_v$.

(ii) We write $\partial=\partial_{\xi_i}$, $F(\xi):=\mathcal P(p(\xi);h,h)$, and
\begin{equation}\label{eq:Pprime}
\mathcal P'(h,h)[q]:=\partial_q\mathcal P(p;h,h)=-2\ip\psi{\mathcal G_p[q]\psi}_\natural,
\end{equation}
\begin{equation}\label{eq:Pprime2}
\mathcal P''(h,h)[q,q]:=2\nn{\mathcal G_p[q]\psi}^2_\natural+2\ip\psi{(
\mathcal G_p[q]^2-\partial_p\mathcal G_p[q][q])\psi}_\natural .
\end{equation}
Differentiating twice and rewriting the term $4\mathcal P'(h,\partial h)[\partial p]$, which has a
derivative at the next level, as $2\partial\big[\mathcal P'(h,h)[\partial p]\big]$ minus terms
without it, we obtain
\begin{equation}
2\mathcal P(h,\partial^2h)=\partial^2F-2\mathcal P(\partial h,\partial h)-2\partial\big[
\mathcal P'(h,h)[\partial p]\big]+\mathcal P''(h,h)[\partial p,\partial p]+\mathcal P'(h,h)[\partial^2p] .
\end{equation}
We use the following facts. The functions $\psi$, $\mathcal G_p[q]\psi$, $\mathcal G_p[q]^2\psi$ and
$\partial_p\mathcal G_p[q][q]\psi$, and the $\hat z$-derivatives of $\psi$, are moment-free.
Further, we have $|\ip u{u'}_\natural|\le C_\natural(7)\|u\|_{X_7}\|u'\|_{X_7}$ by
Lemma~\ref{lem:lowweight}, and $|\ip u{u'}_\natural|\le J_7\|u\|_{Z_7}\|u'\|_{Y_7}$ by (CM1) of
Appendix~\ref{app:band}. Moreover,
$\|\mathcal G_p[q]\psi\|_{X_7}\le C|q|\,\|\psi\|_{H^1(\jap{\hat z}^8)}$, and, for a second-order
operator $P$ with polynomial coefficients of degree $\le2$, we have $\|P\psi\|_{Z_7}\le
C\|\psi\|_{H^1(\jap{\hat z}^{11})}$, since in the $X_7$-pairing one derivative moves to the test
function, which $Y_7$ controls with the weight $\jap{\hat z}^{7+\gamma/2}$. By
Lemma~\ref{lem:TG0}(ii), we have $\|\psi\|_{X_7}\le C\|\psi\|_X$ and
$\|\psi\|_{Y_7}+\|\psi\|_{H^1(\jap{\hat
z}^{11})}\le C\|\psi\|_{Y_2}$. Hence we have $|\mathcal P'(h,h)[q]|\le
C|q|\,\|\psi\|_X\|\psi\|_{Y_2}$ and
$|\mathcal P''(h,h)[q,q]|\le C|q|^2\|\psi\|^2_{Y_2}$. We integrate against $\Psi$:
\begin{itemize}
\item $\int\Psi\partial^2F=\int(\partial^2\Psi)F$, with $|\partial^2\Psi|\le
C\rad_0^{-3}\chi_c^2\Omega\rad_0^{-2}$ (using $|\nabla s_\xi|\le2L's_\xi$, $|\nabla^2s_\xi|\le
CL'^2s_\xi$, $s_\xi\le S_b+1$ on the support of $\chi'$, and the bounds on $\nabla\log\Omega$,
$\nabla^2\log\Omega$) and $F\le C\|\psi\|^2_X$, and no factor $\chi$ is needed;
\item $-2\int\Psi\mathcal P(\partial h,\partial h)=-2\int\Psi\nn{T_p^{-1}(\partial
h)^{(x)}}^2_\natural$, the good term;
\item $-2\int\Psi\partial[\cdots]=2\int(\partial\Psi)\mathcal P'(h,h)[\partial p]$, where
$\partial\Psi$
is supported in $\{\chi(s_\xi)>0\}\cup\{\chi'(s_\xi)\ne0\}\subset\{s_\xi\ge S_b\}$, and
$|\partial\Psi|\le C\rad_0^{-3}\chi_c^2\Omega\rad_0^{-1}$. Thus this term is $\le
C\rad_0^{-3}\chi_c^2\Omega\,C_\varpi\big(\rad_0^{-2}\|\psi\|_X^2+\|\psi\|^2_{Y_2}\big)
\mathbf 1_{\{s_\xi\ge S_b\}}$, integrated;
\item $\int\Psi\mathcal P''[\partial p,\partial p]\le C C_\varpi^2\int\Psi\|\psi\|^2_{Y_2}$ and
$\int\Psi\mathcal P'[\partial^2p]\le CC_\varpi^{(2)}\int\Psi\|\psi\|_X\|\psi\|_{Y_2}$, both on
$\{s_\xi\ge S_b\}$, since $\Psi$ contains $\chi(s_\xi)$.
\end{itemize}
For $L_{1,i}$, we expand $L_1^2=\partial_z^2+2t^{-1}\partial_z\partial_\xi+t^{-2}\partial_\xi^2$.
The last part is the above with the factor $t^{-2}$, and we drop its good term. Since
$\partial_z=\theta^{-1/2}\nabla_{\hat z}$ in the frame, we have $
\mathcal P(h,\partial_z^2h)=\theta^{-1}\ip\psi{\Delta_{\hat z}\psi}_\natural\le C\|\psi\|^2_{Y_2}$, and
$\Delta_{\hat z}$ preserves moment-freeness. For the mixed term, we have
\begin{equation}
\mathcal P(h,\partial_z\partial_\xi h)=\partial_\xi\big[\mathcal P(h,\partial_zh)\big]-
\mathcal P(\partial_\xi h,\partial_zh)-\partial_p\mathcal P(h,\partial_zh)[\partial p] .
\end{equation}
The first term integrates by parts against $\Psi$ and is bounded as the third bullet. The second is
$\le\nn{\psi'}_\natural\,C\|\psi\|_{Y_2}$, and with the factor $2t^{-1}\cdot2\eta_v$ it is
absorbed
by half of the good term up to $C\eta_v t^{-2}\int\Psi\|\psi\|^2_{Y_2}$. The third is $\le
CC_\varpi\|\psi\|^2_{Y_2}$. All the $Y_2$-terms carry the factor $\Psi$ or $\partial\Psi$, so they
live on $\{s_\xi\ge S_b\}$. Finally, we use
$\|\psi\|_X\|\psi\|_{Y_2}\le\rad_0^{-2}\|\psi\|^2_X+\rad_0^2\|\psi\|^2_{Y_2}$,
$C^{(2)}_\varpi\|\psi\|_X\|\psi\|_{Y_2}\le\rad_0^{-2}\|\psi\|^2_X+(\rad_0C^{(2)}_\varpi)^2\|\psi\|^2_{Y_2}$
and $\rad_0\le1$. This completes the proof.
\end{proof}

\subsubsection{The viscous terms of the truncated hierarchy}

Let $t_*\ge T_0$, $\delta\in(0,q_*/4]$ and $q_\delta(t):=q(t)-\delta$, and let $g=g_\eta$ be the
viscous approximation of point~(3) of the proof of Proposition~\ref{prop:CULW}, with the energies
computed at the exponent $q_\delta$, the weight $\Omega\in\{1,\Omega^\epsilon\}$, and the graded
truncation $\mathcal T_N$. We put
\begin{equation}
k_*:=\inf\big\{\kappa(t)\lambda_\xi(t):t\in[t_*,t_*+1],\ s_\xi(t)\ge S_b\big\}\ \ge\ c_\lambda
S_b\inf_{t\in[t_*,t_*+1]}\frac{\kappa(t)}{\tau(t)-\tau(T_0)}>0 ,
\end{equation}
where the lower bound holds because $\kappa\lambda_\xi\ge c_\lambda\kappa\lambda_{\rm ref}(\xi)$ and
$s_\xi=\lambda_{\rm ref}(\xi)(\tau(t)-\tau(T_0))$, where $\tau(t)$ is the clock, and
\begin{equation}
\eta_0:=\min\Big\{\frac\delta{64q_*^2},\ \frac{C_2}{32q_*^2(t_*+1)^2},\ \rad_0^2,\ \frac{c_{\rm
cl}k_*}{16\max\{C_V,C_V'\}}\Big\},
\end{equation}
with $C_V$ the constant of Lemma~\ref{lem:LCvisc}(ii), $C_V'$ the constant of the $Y_2$-terms in
Lemma~\ref{lem:LCtrunc}(iv), and $c_{\rm cl}=c_W/(4C_W)$. The constants of that lemma
are applied with $\rad_0^2C_\Omega\le C$ (for $\Omega^\epsilon$ by requirement~\ref{rq:CUb}) and
$\rad_0C^{(2)}_\varpi\le C$, by point~(2) of the same proof. In
that proposition, the further conditions on $\eta_0$ stated in its proof are added.

\begin{lem}[The viscous terms]\label{lem:LCtrunc}
Let $\eta_v\le\eta_0$. In $\frac{d}{dt}U^c_N$ and $\frac{d}{dt}U^c_{N,\ge1}$, the viscosity
$\eta_v(\Delta_y+\Delta_v)$ acting on the components $h_{\alpha\beta}$, $(\alpha,\beta)\in
\mathcal T_N$, contributes:
\begin{enumerate}[label=(\roman*)]
\item the local good term
\begin{equation}
-\eta_v\sum_{(\alpha,\beta)\in\mathcal T_N}W^2_{\alpha\beta}\,\eta\,\rad_0^{-3}\!\iint\tilde\Omega\big(|\nabla_\xi
h_{\alpha\beta}|^2+|(\nabla_z+t^{-1}\nabla_\xi)h_{\alpha\beta}|^2\big);
\end{equation}
since $\nabla_\xi h^{(x)}_{\alpha\beta}=(\partial_vh_{\alpha\beta})^{(x)}$, at each $y$-order $k$ it
is $\le-\eta_v\eta\sum_nw_n(\rad_k)^2L^c_{k,n+1}$, where
\begin{equation}
L^c_{k,n+1}:=\rad_0^{-3}\!\int\chi_c^2\Omega\!\!\sum_{|\alpha|=k,\,|\beta|=n+1}\!\!\delta_y^{2k}\ell_\xi[h_{\alpha\beta}]\dd\xi
\end{equation}
(without the factor $\eta$ of Lemma~\ref{lem:CUtr});
\item the weight loss $\eta_v\eta\,32q^2|z|^2$ of Lemma~\ref{lem:LCvisc}(i), which is absorbed by
the term $-C_2t^{-2}\iint\tilde\Omega|z|^2|h|^2$ that Lemma~\ref{lem:CUtr} leaves, since
$32\eta_vq_*^2\le C_2(t_*+1)^{-2}$. The rest of that bound is $\le
C\eta_v\rad_0^{-2}U^c_N\le CU^c_N$;
\item from the nonlocal part, by part~(ii) of the same lemma applied with
$\psi=\varphi_{\alpha\beta}$, the following: a nonpositive good term, terms $\le
C_V\eta_v\rad_0^{-2}\eta^{-1}U^c_N\le CU^c_N$, and the $Y_2$-terms
$C_V\eta_v\rad_0^{-3}\int\chi_c^2\Omega\mathbf 1_{\{s_\xi\ge
S_b\}}W^2_{\alpha\beta}\|\varphi_{\alpha\beta}\|^2_{Y_2}$, which do not exceed a sixteenth of the
leaf
dissipation $c_{\rm cl}\kappa\lambda_\xi\|\varphi_{\alpha\beta}\|^2_{Y_2}$ on the leaves with
$s_\xi\ge S_b$, since $\kappa\lambda_\xi\ge k_*$ there;
\item from the viscous part of the frame motion, the following. For $g_\eta$ the leaf parameters
move with $D\varpi=t^{-2}(\text{Euler})+V(\varpi)$, with the viscous field $V$ of point~(3) of the
proof of Proposition~\ref{prop:CULW}, and the cell energy depends on $p_g$ through the local weight
$(A_p/A_{\rm ref})^2$ and through the nonlocal form $\mathcal P(p;h,h)$. Equivalently, the hierarchy
contains the frame motion $\mathcal G_p[Dp]\varphi_{\alpha\beta}$, with the viscous part $
\mathcal G_p[V]\varphi_{\alpha\beta}$. This part does not change $(\alpha,\beta)$ (it is of order zero in the
levels, though of first order in $\hat z$), so it crosses no level. Since
\begin{equation}\label{eq:LCtrunc-V}
|V|\le C\eta_v(|\nabla^2\varpi|+|\nabla\varpi|^2+t^{-1}|\nabla\varpi|+1)\le C\eta_v\rad_0^{-1},
\end{equation}
by $|\nabla^2\delta\varpi|\le C\delta_U\rad_0^{-2}$ under the truncated threshold and
$C\delta_U\le\rad_0$, the local part is a Gronwall term $\le C\eta_v\rad_0^{-2}U^c_N\le
CU^c_N$. The nonlocal part is $\le C|V|\,\|\varphi_{\alpha\beta}\|_X\|\varphi_{\alpha\beta}\|_{Y_2}$
on $\{\chi>0\}$, as for $\mathcal P'$ in that proof. By Young's
inequality, it yields terms $\le CU^c_N$ and $Y_2$-terms $\le
C_V'\eta_v\rad_0^{-3}\int\chi_c^2\Omega\mathbf 1_{\{s_\xi\ge
S_b\}}W^2_{\alpha\beta}\|\varphi_{\alpha\beta}\|^2_{Y_2}$, which again do not exceed a sixteenth of
the leaf dissipation, since $C_V'\eta_v\le\frac1{16}c_{\rm cl}k_*$.
\end{enumerate}
Moreover, at the top level of $y$-order $0$, the transport and kernel coupling
$C_0t^{-2}w_N^2\sqrt{E^c_{0,N}E^c_{0,N+1}}$ of the cell hierarchy (the nonlocal transport (T10) and
the kernel part (T12$_0$)) splits as
$\frac{C_0t^{-2}}{2\rad}(N+1)(\tilde a_N^2+b_{N+1}^2)$, with $\tilde a_N^2=w_N^2E^c_{0,N}$ and
$b_{N+1}^2:=w_{N+1}^2E^c_{0,N+1}$, and $E^c_{0,N+1}\le(\eta+C)L^c_{0,N+1}$. The radius absorbs
$\tilde a_N^2$ by Lemma~\ref{lem:analytic}, and the good term~(i) at $n=N$ absorbs $b_{N+1}^2$ as
soon as $N+1\ge N_\eta:=C\rad_0C_0/(\eta\,\eta_vt_*^2)$. The constants do not depend on $\eta_v$,
$N$, $\delta$ or $c$. (In Proposition~\ref{prop:CULW} this absorption is done for the joint family,
point~(3), and replaces the one stated here.)
\end{lem}

\begin{proof}
The viscous term of the equation for $h_{\alpha\beta}=\partial^\alpha_y\partial^\beta_v\hat r$ is
$\eta_v(\Delta_y+\Delta_v)h_{\alpha\beta}$, since the derivatives commute with the Laplacians. The
viscous forcing of the micro equation, $\partial^\alpha_y\partial^\beta_v\mathfrak F_{\rm visc}$, is
treated in point~(3) of the proof of Proposition~\ref{prop:CULW}.

The cell energy of $h_{\alpha\beta}$ is
$W^2_{\alpha\beta}\big(\eta\rad_0^{-3}\iint\tilde\Omega|h_{\alpha\beta}|^2+\int\Psi\,
\mathcal P(p;h_{\alpha\beta},h_{\alpha\beta})\dd\xi\big)$, and the viscosity
enters its time derivative through $\partial_th_{\alpha\beta}$, to which Lemma~\ref{lem:LCvisc}
applies, and through the viscous part of the motion of the frame $p_g$, which is item~(iv).

(i) is the local good term of that lemma, multiplied by $\eta
W^2_{\alpha\beta}\rad_0^{-3}$. At a fixed $(\alpha,n)$ the families
$\{\partial_{v_i}h_{\alpha\beta}:|\beta|=n,i\le3\}$ contain every $h_{\alpha\beta'}$ with
$|\beta'|=n+1$, and $W_{\alpha\beta}=\delta_y^{|\alpha|}w_{|\beta|}(\rad_{|\alpha|})$.

(ii) is the rest of the same bound: the term with $32q^2|z|^2$ has the same weight
$\tilde\Omega$ as the term $-C_2t^{-2}\iint\tilde\Omega|z|^2|h|^2$ of Lemma~\ref{lem:CUtr} (with
$\omega$ replaced by $\Omega$), and $q\le q_*$, $t\le t_*+1$.

(iii) follows from part~(ii) of the same lemma with $\|\varphi\|^2_X\le\eta^{-1}e_\xi[\varphi]$,
$\eta_v\le\rad_0^2$ and $C_V\eta_v\le\frac1{16}c_{\rm cl}k_*$.

For the top coupling, the split follows from $w_N\le\frac{N+1}\rad w_{N+1}$ and Young's
inequality. The
energy at the level $N+1$ contains the local part with the coefficient $\eta$ and the nonlocal part
$\le C\ell_\xi$ by Lemma~\ref{lem:lowweight}, and thus $E^c_{0,N+1}\le(\eta+C)L^c_{0,N+1}$.
Moreover, we have $w_N(\rad)^2\ge\big(\frac{N+1}{16\rad}\big)^2w_{N+1}(\rad)^2$ for $N\ge1$, since
$w_N/w_{N+1}=\frac{N+1}\rad\big(\frac{N+1}{N+2}\big)^5$ for $w_n=\rad^n(n+1)^5/n!$, and thus the
good term $\eta_v\eta w_N^2L^c_{0,N+1}\ge
c\,\eta_v\eta(\eta+C)^{-1}\big(\frac{N+1}{\rad}\big)^2b_{N+1}^2$ dominates
$\frac{C_0t^{-2}}{2\rad}(N+1)b_{N+1}^2$ for $N+1\ge N_\eta$, since $\rad\le\rad_0$ and
$t\ge t_*$.
\end{proof}

The collision terms are those of the cell inequality without damping, as in the last bullet of
point~(2) of the proof of Proposition~\ref{prop:CULW}. On the leaves with $s_\xi<S_b$ (where
$\chi(s_\xi)=0$) they are those of Proposition~\ref{prop:ACpf}, (L1)--(L2), (N), (C), whose proof
applies on this
range since no damping is used, with point~(5) of the same proof in place
of (G4) (in particular for the level-$(0,0)$ background dissipation in the inequality for the parts
$\ge1$, item~($\gamma$) there), with the coercivity $\eta c_h(K)$ and the linear growth
$C_L(K)\kappa\lambda_\xi\le C(K)$, which is a Gronwall term. The quadratic terms are treated by
Proposition~\ref{prop:CUtri} under the threshold on $U^c_{N,\ge1}$. On the leaves with $s_\xi\ge
S_b$, the closeness $\|f_\xi\|_X\le\delta_c$ holds on the interval by point~(6) of that proof, and
the coercivity is $c_{\rm cl}$ of Proposition~\ref{prop:Bband} and
Remark~\ref{rem:Bbandrange}, including the nonlocal part (T3)--(T4). The deficit and the front term
are Gronwall terms. The macroscopic fields and the cube suprema enter through $\sup_{c'}U^{c'}_N$,
which is the term $BS$ of Lemma~\ref{lem:LCgron}.

\emph{The dissipation that is left.} We keep a part of the collision dissipation for
Subsection~\ref{sssec:CUac}. As in the proof of Proposition~\ref{prop:CUcell}, we let $c(\xi)$ be
the
coercivity that the leaf inequality provides after its own absorptions: $c(\xi)=\eta c_h$ on the
leaves with $s_\xi<S_b$, with the constant $c_h$ used in the preceding paragraph, and
$c(\xi)=c_{\rm cl}$ on the leaves with $s_\xi\ge S_b$. The absorptions into it are those of
Steps~2--5 of that proof, with the following changes.
\begin{itemize}
\item The quadratic terms are bounded as above under the threshold on
$U^c_{N,\ge1}$, which implies that $\sup_Q\mathfrak a_{Q^*,\ge1}\le c_\sharp/e$ on every cube.
Thus they remove at most $C_Nc_\sharp^{1/2}\le c_{\min}/8\le c(\xi)/8$ on every leaf.
\item The coefficient and $\delta\varpi$-terms remove at most $\frac14c(\xi)$, as there.
\item The background terms remove at most $\frac18c(\xi)$. On the leaves with $s_\xi\ge S_b$ the
background appears in (T4) and, at the level $(0,0)$, in the quadratic term, and nowhere else. The
parts of (T4) with the factor $\delta_c$, and the Young split of its parts with $b_\xi^{1/2}$, are
absorptions of Proposition~\ref{prop:Bband} itself, already counted in $c(\xi)$. The rest,
$C\kappa\lambda_\xi b_\xi
e_\xi$, is a Gronwall term, since $b_\xi\le\delta_c^2+CG_*(M)^2$ by point~(5)($\beta$) of the proof
of Proposition~\ref{prop:CULW}. Neither (H$_b$) nor the background variant of
the trilinear bound is used. The quadratic term at the level $(0,0)$ is
$\le2C_\tau\delta_c\|f_\xi\|^2_{Y_2}$, which removes at most $\frac1{16}c_{\rm cl}$, by
requirement~\ref{rq:Wdc}. On the leaves with $s_\xi<S_b$ the background enters only (L2), the part
of (C) with $(\alpha'',\beta'')=0$ and the quadratic term at the level $(0,0)$. By
point~(5)($\alpha$)--($\gamma$), these remove nothing beyond the part of (C) already estimated in
the crude leaf inequality.
\item The viscous $Y_2$-terms of Lemma~\ref{lem:LCtrunc}(iii) and~(iv) remove together no more than
$\frac18c_{\rm cl}$, on the leaves with $s_\xi\ge S_b$.
\end{itemize}
Altogether at most $\frac58c(\xi)$ is removed. Hence the inequality for $\frac{d}{dt}U^c_N$ holds
with the additional term
\begin{equation}
c_{\rm ret}\,\kappa\,\mathfrak D^c_N,\qquad c_{\rm ret}:=\tfrac14\min\{\eta c_h,c_{\rm
cl}\},\qquad\mathfrak D^c_N:=\sum_{(\alpha,\beta)\in
\mathcal T_N}W^2_{\alpha\beta}\rad_0^{-3}\!\int\chi_c^2\,\Omega\,d_\xi[\varphi_{\alpha\beta}]\dd\xi,
\end{equation}
on its left side. The factor $16$ in the definition of $\eta_0$ is needed here: with the fractions
$\frac18$ removed by each of the two viscous terms, and with the count of
Steps~2--5 for the others, nothing would be left on the leaves with $s_\xi\ge
S_b$. The constant $c_{\rm ret}$ depends on $K$, through $c_h$. In
the bootstrap only the finiteness of the time integral of the dissipation matters.

\subsubsection{Compactness and lower semicontinuity}

\begin{lem}[Compactness and lower semicontinuity]\label{lem:LCd}
\leavevmode
\begin{enumerate}[label=(\roman*)]
\item \emph{Weak lower semicontinuity with converging weights.} Let $\mathcal C\subset\R^6$ be
compact, $f_k\rightharpoonup f$ weakly in $L^2(\mathcal C)$, and $0\le\varsigma_k\to\varsigma$
uniformly on $\mathcal C$. Then
\begin{equation}\label{eq:LCd-lsc}
\int_{\mathcal C}\varsigma|f|^2\le\liminf_k\int_{\mathcal C}\varsigma_k|f_k|^2 .
\end{equation}

\item \emph{Compactness.} Let $t_*<T<\infty$, and let $(g_k)$ be a sequence of functions on
$[t_*,T]\times\R^6$ such that, for every $(\alpha,\beta)$ with $|\alpha|\le N_y$ and every compact
$\mathcal C$, $\partial^\alpha_y\partial^\beta_vg_k$ is bounded in $L^\infty(t_*,T;L^2(\mathcal C))$
and $\partial_t\partial^\alpha_y\partial^\beta_vg_k$ in $L^\infty(t_*,T;H^{-2}(\mathcal C))$,
uniformly in $k$. Then a subsequence satisfies, for all such $(\alpha,\beta)$ and $\mathcal C$,
$\partial^\alpha_y\partial^\beta_vg_k\to\partial^\alpha_y\partial^\beta_vg$ in
$C([t_*,T];H^{-1}(\mathcal C))$, hence
$\partial^\alpha_y\partial^\beta_vg_k(t)\rightharpoonup\partial^\alpha_y\partial^\beta_vg(t)$ weakly
in $L^2(\mathcal C)$ for every $t$, and $g_k\to g$ in $C([t_*,T];C^{N_y-4}(\mathcal C))$. (In
Proposition~\ref{prop:CULW}, $T=t_*+\tau$, and in Proposition~\ref{prop:disp}, $[t_*,T]=[t_1,T_0]$.)
\item \emph{Cell energies.} Let $g_k$, $g$ be as in~(ii), at a time $t$, and suppose that the leaf
parameters $\varpi_{g_k}(t,\cdot)$ and their $\xi$-derivatives converge locally uniformly to those
of $g$, and that the remaining weights of the cell energies ($\chi(s_\xi)$, $\omega_\xi$ or
$\Omega$, $W_{\alpha\beta}$ through the radii, $e^{2q|z|^2}$) converge locally uniformly. Then
$X_c[g](t)\le\liminf_kX_c[g_k](t)$ for every cell $c$ and $X\in\{\mathfrak A,\mathfrak A_{\ge1}\}$,
and for their truncations to $\mathcal T_N$. If, moreover, the components converge strongly in the
weighted $L^2$ spaces of the local parts, then $X_c[g_k](t)\to X_c[g](t)$.
\end{enumerate}
\end{lem}

\begin{proof}
(i) Note that $\int_{\mathcal C}(\varsigma_k-\varsigma)|f_k|^2\to0$, since
$\sup_k\|f_k\|_{L^2}<\infty$, and
$f\mapsto\int\varsigma|f|^2$ is convex and continuous on $L^2(\mathcal C)$, hence weakly lower
semicontinuous.

(ii) The embedding $H^1_0(\mathcal C')\subset L^2(\mathcal C')$ is compact for a bounded open $
\mathcal C'\supset\mathcal C$, hence so is the dual embedding $L^2(\mathcal C')\subset H^{-1}(\mathcal C')$
(Schauder's theorem), and $H^{-1}(\mathcal C')\subset H^{-2}(\mathcal C')$ continuously. By the
Aubin--Lions--Simon lemma \cite{Sim87}, applied to the family $\partial^\alpha_y\partial^\beta_vg_k$
(bounded in $L^\infty(L^2)$, with time derivatives bounded in $L^\infty(H^{-2})$), we obtain
relative
compactness in $C([t_*,T];H^{-1}(\mathcal C))$, and the length of the interval plays no role. By a
diagonal argument over the countably many $(\alpha,\beta)$ and an exhaustion of $\R^6$ by compact
sets, we get one subsequence for all of them. The limits are the derivatives of the limit of $g_k$,
since
differentiation is continuous in $\mathcal D'$. Bounded sets of $L^2(\mathcal C)$ are weakly
relatively compact, and the weak limit is identified by the $H^{-1}$ limit. The derivatives of order
$\le N_y$ are bounded in $L^\infty(L^2(\mathcal C))$ and converge in $C(H^{-1})$, hence in
$C(H^s(\mathcal C))$ for every $s<N_y$, by interpolation, and $H^s\subset C^{N_y-4}$ in $\R^6$ for
$s\in(N_y-1,N_y)$.

(iii) \emph{Local part.} The local part at a fixed $(\alpha,\beta)$ is
$\int\varsigma\,|\partial^\alpha_y\partial^\beta_vg^{(x)}-\text{(Maxwellian part)}|^2$ over
$\R^3_\xi\times\R^3_z$, with a weight $\varsigma$ that is a continuous function of the leaf
parameters and of the weights listed. The Maxwellian part is an explicit function of $\varpi_g$ and
its derivatives, and converges locally uniformly. Since $\dd y\dd v=\dd\xi\dd z$, lower
semicontinuity follows from (i) on each
compact set and an exhaustion.

\emph{Nonlocal part.} For a fixed frame map $p$, the functional $u\mapsto\int\Psi\,
\mathcal P(p;u,u)\dd\xi$ is a nonnegative quadratic form, convex and continuous on the weighted $L^2$ space
of the local part, since $\mathcal P(p;u,u)\le C_\natural(8)\|T^{-1}_pu\|^2_{X_8}\le C\ell_\xi[u]$
by Lemma~\ref{lem:lowweight} and Lemma~\ref{lem:TG0}(ii). Hence it is weakly lower semicontinuous.
The bounded sequence converges weakly in that space, by~(ii) and the identification of the weak
limit on
compact sets.

For the change of frame, we put $S_k(\xi):=T^{-1}_{p_k(\xi)}T_{p(\xi)}$. Since $p_k\to p$ locally
uniformly, $S_k$ and its adjoint (for the $X_8$-pairing) are uniformly bounded on $X_8$ and converge
strongly to the identity, by the strong continuity of translations and dilations on $L^2(\jap{\hat
z}^8)$. Hence $T^{-1}_{p_k}u_k=S_kT^{-1}_pu_k\rightharpoonup T^{-1}_pu$ weakly in
$L^2(\Psi\dd\xi;X_8)$ whenever $u_k\rightharpoonup u$ weakly in the weighted $L^2$ space of the
local part. The tails in $\xi$ are uniformly small by the factor $\chi_c^2$ and the bound on the
energies. The nonlocal part, a nonnegative bounded quadratic form on $L^2(\Psi\dd\xi;X_8)$, is
therefore weakly lower semicontinuous along the sequence, and the convergence of the weights
$\Psi_k\to\Psi$ is handled as in~(i).

If $u_k\to u$ strongly, then $T^{-1}_{p_k}u_k\to T^{-1}_pu$ strongly, and the nonlocal parts
converge. The local parts converge as well, which implies the last claim. Fatou's lemma in the sum
over $(\alpha,\beta)$ and $n$ completes the proof.
\end{proof}

\emph{The hypotheses of~(ii) for the viscous approximations.} In Propositions~\ref{prop:disp}
and~\ref{prop:CULW} the bound in $L^\infty(L^2(\mathcal C))$ follows from the cell bounds, uniform
in $\eta_v$: $\dd\xi\dd z=\dd y\dd v$, and the weight $\chi_c^2(A_p/A_{\rm ref})^2m_x^2$ is bounded
below on $\mathcal C$ for a cell centered near $\mathcal C$.

For the time derivative, we note that $\partial_t\partial^\alpha_y\partial^\beta_vg$ is
$\partial^\alpha_y\partial^\beta_v$ of the right side of~\eqref{eq:LWreg}, a sum of products
$\ab[\partial^{\mathbf k'}g]:D^2\partial^{\mathbf k-\mathbf k'}g$ and
$\cb[\partial^{\mathbf k'}g]\,\partial^{\mathbf k-\mathbf k'}g$, $\mathbf k=(\alpha,\beta)$, and of
$\eta_v\Delta\partial^{\mathbf k} g$.

If the $y$-order of $\mathbf k'$ is $\le N_y-4$, the coefficient and its first two derivatives are
bounded on $\mathcal C$, by Lemma~\ref{lem:LCdata} (applied with the cell bounds, as in point~(4) of
the proof of Proposition~\ref{prop:CULW}), and the product lies in $H^{-2}(\mathcal C)$ after we
move two
derivatives.

Otherwise the $y$-order of $\mathbf k-\mathbf k'$ is $\le3$, and thus $D^2\partial^{\mathbf k-
\mathbf k'}g$ is bounded on $\mathcal C$ by the same lemma, and the coefficient lies in
$L^2(\mathcal C)$, by Minkowski's inequality in $w$, the local integrability of the kernels and the
decay of the cell bounds along the shifts $(y+tw,v-w)$. Indeed, the cell energies carry the weight
$e^{2q|z|^2}$, and $z=y+tw$ on the shifted set, so the $L^2$ norms over the shifted sets decay like
$e^{-qt_*^2|w|^2/2}$ for $|w|$ large.

The Maxwellian part of $\partial^{\mathbf k} g$, which the energies do not contain, is bounded on
$\mathcal C$ by the leaf parameters and their derivatives.

\subsection{The order of the constants}\label{app:const}

We choose the constants in the order of Table~\ref{tab:const}, and each depends only on the ones
before it. The stage numbers are used only in this appendix, while the main text states the order by
naming the constants (for instance ``chosen before $\rad_0$'' or ``fixed before $T_0$''). The entry
level
$\delta_1$ is chosen before the entry threshold $S_*$, which Lemma~\ref{lem:Aentry} allows, and
information on neighboring leaves is never used at entry, so that $S_*$ may depend on $\delta_1$ but
not conversely. A second point is that the inputs of the fibered comparison from the relaxed region,
$K_{\rm cr}$, $K_{\rm rel}$, $K^{\rm cr}_{\rm rel}$ and $K_{\rm bar}$, are large constants of
stage~$8$ and are never required to be small: they enter only the barrier constant $C_{\rm bar}$ and
the constant $C_{\rm loc}$ of the local factors, hence only the lower bound on $T_0$.

{\small
\begin{longtable}{@{}>{\raggedright\arraybackslash}p{0.06\textwidth}>{\raggedright\arraybackslash}p{0.56\textwidth}>{\raggedright\arraybackslash}p{0.32\textwidth}@{}}
\caption{The constants, in the order of choice.}\label{tab:const}\\
\toprule
Stage & Constants & Where they are fixed\\
\midrule
\endfirsthead
\toprule
Stage & Constants & Where they are fixed\\
\midrule
\endhead
0 & $\gamma$, $N_y\ge20$, $s\ge6$; the data constants
$\Lambda,\Lambda_\Phi,L,p,C_F,\beta_F,A_F,C_\Phi,C_{\rm pv}$, then $B$ (Definition~\ref{def:data}),
$\theta_\pm$, $q_{\rm tri}=\frac18$ (Proposition~\ref{prop:triG}) and $q_d:=\min\{\beta_F/4,q_{\rm
tri},(8(L^2+3\theta_+))^{-1}\}$, $(C_{\rm dat},\beta,A_7)$ of $\mathfrak K^7=\mathfrak K^{n+2}$ and
$C'_D$ (the fibers at $T_0$, (LWP)(a3) in Section~\ref{sec:proofmain}), and $\beta_M$,
$q_S=\min\{\beta_M/32,1/48\}$; $K_M=C(L,\theta_\pm)$; $C_{\rho_0}=C(A_F,L,\Lambda)$ (Step~1 of the
proof of Proposition~\ref{prop:disp} and Corollary~\ref{cor:LWPa}(a4); it enters $\vartheta_1$,
hence $\rad_1$, and $C_q$, $C^{\rm lay}_{\rm tr}$, $L_1$ and the constant $C(\Dc)$ of (a4), so
it is fixed before them); $C_{\mu_0}$ and $G_\infty$ (Subsection~\ref{ssec:disp},
\emph{Parameters}; both vanish for products; they enter $c_{\rm a4}$ and an upper bound on
$\rad_1$, stage~$2$); $c_{\rm a4}$ (Corollary~\ref{cor:LWPa}) and the constant $C(\Dc)$ of its
free part & Definition~\ref{def:data}; (LWP)(a3) in Section~\ref{sec:proofmain}\\
1 & $k_1\ge k_\flat(\gamma)$ (Lemma~\ref{lem:Abarrier}; $k_\flat(-2)=7$, $k_\flat(-3)=11$),
$\delta_0$, $T_*^0:=T_*(\delta_0)$ (Theorem~\ref{thm:enter}); $K_\infty$, $C_a$, $c_G$, $C_c$,
$\beta'$, $C_\flat$ & Theorem~\ref{thm:enter}, Lemma~\ref{lem:Abarrier}, \cite{CM}\\
2 & the layer rate $C_1$ of Lemma~\ref{lem:ULconv} (at least twice the data rates of $f_0$ and of
$\zeta$, Lemmas~\ref{lem:FPexp} and~\ref{lem:FPdecay}(ii); it is not the constant
$C_1=C_{\sharp,\rm f}c_\theta$ of stage~$5$), $C_{\rm v}$ of~\eqref{eq:Cv} and $\rad_*$, none of
which depends on $t_1$; the layer radius $\rad_0^{\rm lay}\le\rad_*R$
(Subsection~\ref{ssec:layer}); the constants $E^{\rm fr}$, $\mathcal K_{\rm fr}$, $\mathcal K$,
$K_\ell$, $c_G^{\rm lay}$, $c_h^{\rm lay}$, $C_h^{\rm lay}$, then $E_k$ and $c_{\sharp,\ell}$ (with
the layer constant $C_N^{\rm lay}=C(N_y,C_\tau(q_d))$, $C_\tau$ structural), then $\rad_1$, then
$R_1:=\rad_1/8$ and $L_1$, with $C_q$, $C^{\rm lay}_{\rm tr}$ and $t_1$, then the layer
$y$-weight $\delta_{y,\ell}$ (it depends on $t_1$ through the moment constants of
Lemmas~\ref{lem:R3comp} and~\ref{lem:ULmom}, and it does not enter the choice of $t_1$;
Proposition~\ref{prop:FPclose}), then $c_M$ and $C_D$ (Corollary~\ref{cor:LWPa}), and $C_{\rm
lay}$, $C_L$, $C_w$, then $C_{\mu,1}$, $C_{\mu,2}=C(C_w,q_d)$ (Corollary~\ref{cor:LWPa}(a4)) (which
enter only $\eps$-conditions and Corollary~\ref{cor:LWPa}) of the dispersed layer, in this order
(Subsection~\ref{ssec:disp}; they depend only on $\Dc$ and structural constants, so $t_1=t_1(\Dc)$;
when Theorem~\ref{thm:t1} is used alone, $C_{t_1}^2$ takes the place of $C_{\rm lay}t_1^2$, so
$C_w$ and $C_{\mu,2}$ depend also on $C_{t_1}$, and they enter only $\eps$-conditions); the barrier
rate $r$; $q_\natural$, $r_\natural$; $q_0:=\min\{\beta'/4,\,1/64\}/\max\{1,\theta_+\}\wedge
q_d\wedge\beta'\theta_-/(64\theta_+^2)\wedge\theta_-q_\natural/64\wedge\theta_-q_d/16\wedge q_{\rm
GW}(\gamma,L+1,\theta_-/2,\theta_++1)\wedge\big(8((L+1)^2+3\theta_++3)\big)^{-1}$, with $q_{\rm
GW}$ of Lemma~\ref{lem:GW} in the frames $p_g$ given by (F3), and the last clause for
Lemma~\ref{lem:Whypo} in the frames $p_g$, after the translation in (L1) of
Proposition~\ref{prop:ACpf}; $T_1$; $C_W,c_W,C_\natural$, $c_\lambda$, $\eta$, $C_1^B$, $C_0$,
$C_2$ (with $C_2\ge C$ for (T8) in Lemma~\ref{lem:UFlevel0}), the constants $\check c_W$, $\check
C_W$, $\check\eta$, $\check C_\tau$, $\check C$, $\check C_{11}$ of the weights $\hat m_\xi^{1/2}$
(before Lemma~\ref{lem:R3weak}), $\delta_c$, $\delta_s$, $C_e$, $C_{\rm cl}'$, $C_N$; the
structural upper bounds on $\rad_0$ and $\delta_y$; $c_{\rm cl}=c_W/(4C_W)$, $c_{\rm buf}=c_{\rm
cl}/2$, $c'_{\min}=c_{\rm cl}$, $c''_{\min}=c_{\rm buf}$; $r_{\rm C}=c_{\rm C}R_1/(L_1+1)$ ($c_{\rm
C}$ absolute), $\epsilon_y=1/(4N_y+4)$; $C_{{\rm ref},a}$ (Lemma~\ref{lem:UFref}), $C_{\rm ref}$,
$C'_{\rm ref}$ (Proposition~\ref{prop:CULW}), $C_{\rm Max}$, $C'_\varpi$, $F_*$;
$C_\varpi:=C(L_1,R_1)+1$ (Subsection~\ref{ssec:cells}); $C_{\rm buf}$; the structural constants
$C_\tau$, $C_\rho$, $C_{\rm lat}$, $C_{L0}$, $C_{\mathfrak s}$, $C_{\mathfrak s'}$, $C'$, $C''$ of
Subsection~\ref{ssec:UF}; for all $\gamma\in[-3,-2]$, after $t_1$ and $\delta_{y,\ell}$: the
prefactor $C_\zeta$ of Lemma~\ref{lem:FPdecay}(ii) (it enters only constants absorbed by $\eps_0$),
$C_\rad=1$, the layer constant $C_{\mathbb E}$, $c_*$ and $\eta_*\le
c_*/(C_{\mathbb E}+1)$ (Proposition~\ref{prop:FPclose}; $c_*$ depends on $t_1$ through
Lemma~\ref{lem:RBrad}, and $C_{\mathbb E}$, $\eta_*$ enter also Lemma~\ref{lem:EPclose} at
$\gamma=-2$), which enter only conditions
on $\eps$; at $\gamma=-2$, after $t_1$ and $\delta_{y,\ell}$: $C^{\rm lay}_\sharp$,
$K_\flat:=9C^{\rm lay}_\sharp N_0$, $C_6$, $C_7$, $C_8$, $C_5:=96C_6+1$ (Lemma~\ref{lem:EPclose};
there the layer parameters $\rad_0^{\rm lay}$, $\delta_{y,\ell}$ are chosen with $1/(72N_0)$ in
place of $1/(32N_0)$), which enter only conditions on $\eps$; at all $\gamma$, $C_4^{[1]}$ and
$C_4:=C_4^{[1]}+1$ (Lemma~\ref{lem:FPf}) & Lemma~\ref{lem:RGbar}, Lemma~\ref{lem:EPclose},
Remark~\ref{rem:RGfull}(b), Proposition~\ref{prop:triG}, Proposition~\ref{prop:Bband},
Lemma~\ref{lem:CUbuf}; Subsection~\ref{ssec:UF}\\
3 & $\delta_1$; then $\delta_b^2:=\frac54C_be^{2C_1^B}\delta_1^2$, $C_b:=54e^5C_{\rm cl}'+28$ &
Theorem~\ref{thm:CUboot}\\
4 & $\delta'':=\min\{1,(\delta_1/(6C_eC_{\rm pt}))^{5/2}\}$, $S_*:=S_\flat(\delta'')+2$,
$S_\sharp:=S_*+1$, $S^\sharp:=4S_\sharp$ & Lemma~\ref{lem:Aentry}, Subsection~\ref{ssec:cells}\\
5 & $K_3$, $K_G$ (on $[0,S^\sharp]$); $C_{W,\rm f}$, $c_{W,\rm f}$, $\epsilon_*$, $\delta_{0,\rm
f}$, $C_{\sharp,\rm f}$, $C_1:=C_{\sharp,\rm f}c_\theta$, $C_{S,\rm f}$, $C_{\rm comp}$; then
$C_T$, $C_M$; $K_4:=C(L,\theta_\pm,q_\natural)C_T(K_G+1)+C(L,\theta_\pm)$;
$K:=\max\{2(2K_3+1),2K_4\}$; $C_U:=2C_DB+1$; $c_\sharp\le\delta_1^2/4$; the frame constant $\bar
K:=C_{\rm fr}K$, $c_G(\bar K)$ and $c_h(K):=c_h(\bar K,c_G(\bar K))$; $c_{\min}:=\min\{\eta
c_h(K),c_{\rm buf}\}$ & Lemmas~\ref{lem:RGGbounds}, \ref{lem:BTU4}, \ref{lem:BTG12};
Remark~\ref{rem:BTnorm}; Proposition~\ref{prop:CUcell}\\
6 & $C_1^U:=\max\{C_L(K)/c_\lambda,\ C_{\rm buf}(1+K^2)/c_\lambda,\ C_1^B,\ 1\}$, with $C_L(K)$ of
(L1)--(L2) of Proposition~\ref{prop:ACpf} and $C_{\rm buf}$ of Lemma~\ref{lem:CUbuf} (which also
covers $|\langle Q(f,f),f\rangle|\le CK\|f\|^2_X$); then
$\mu_\flat:=e^{-2C_1^US_\sharp-4C_1^B}c_\sharp/(54e^5)$ & crude regime
(Proposition~\ref{prop:ACpf}); Lemma~\ref{lem:CUbuf}\\
7 & $\delta_y\le1$; $\rad_0$; $\epsilon_0=12\sqrt3L'\rad_0$; $d_*$ and $S_b:=S_*e^{2L'd_*}$;
$s_\psi:=\max\{S_*,S_b-\frac14\}$ and $L_\psi:=(S_b-s_\psi)^{-1}$ ((U3) of
Theorem~\ref{thm:CUboot}; they enter only the step length $h$ in Subsection~\ref{ssec:cellsproofs},
and no requirement); the graded radii $\rad_k=\rad-k\epsilon_y\rad_0$ &
requirements~\ref{rq:CUa}--\ref{rq:CUc}; requirements~\ref{rq:Wlam}, \ref{rq:Wdy};
Remark~\ref{rem:UFgrade}\\
8 & $\mu_\sharp:=e^{-2H(S_b+1)}c_\sharp/(54e^5)\ge\mu_\flat$; $\rad_{\rm post}:=\rad_0/2$;
$B_{\rm cr}:=54e^5e^{2H(S_b+1)}c_\chi C_U$; $\mathfrak p_*$ (Corollary~\ref{cor:FMuse}(d), at
$B=B_{\rm cr}$); the Gronwall rate $a_1=C+4\mathfrak p_*+C_{\rm tr}/(\rad_0^2C_2)$ of
Section~\ref{sec:boot}; $P_1:=C\rad_0^{-2}e^{2H(S_b+1)}c_\chi C_U$, $a_F$, $C_1^F$, $C_1^{\rm
mac}$; $C_{\rm r}$ (Proposition~\ref{prop:Bband}; it contains $C_{11}\mathfrak p_*$); for
$\gamma<-2$, $\check C_{\rm r}:=\check C_{11}\mathfrak p_*$ (Lemma~\ref{lem:R3weak}), and the
constants $\check C_a$, $\check C_f$, $C_I$ of that lemma; $K_{\rm cr}:=C'_J\big(1+\rad_{\rm
post}^{-J-2}\delta_y^{-J-2}B_{\rm cr}^{1/2}\big)$; $C_{\rm mom}:=C(1+B_{\rm cr}^{1/2})\ge C_{\rm
Max}$, $C_0^F\le C(1+C_{\rm mom}+P_1^{1/2})$, and $F_0$ (paragraph before
Proposition~\ref{prop:UFforce} and ($\Theta$1) in its proof; they enter only requirements of
stage~$9$); $K_{\rm rel}:=K_M+C_2'\rad_{\rm post}^{-2}\delta_y^{-2}\delta_b$; $K^{\rm cr}_{\rm
rel}:=K_M+C_2'\rad_{\rm post}^{-2}\delta_y^{-2}B_{\rm cr}^{1/2}$; $K_{\rm bar}:=\max\{K,K^{\rm
cr}_{\rm rel}\}$; $C_{\rm bar}$, $C_{\rm loc}$ & Subsection~\ref{ssec:UF}; Lemmas~\ref{lem:BTsob},
\ref{lem:BTU4}, \ref{lem:RGloc}\\
9 & $T_0$: the maximum of the lower bounds listed below; then $K_{\mathfrak M}:=T_0$, $C_0'$, and
the rate of the radius of Remark~\ref{rem:UFradius} & Proposition~\ref{prop:UFforce}(g)\\
10 & $\eps_0(T_0)$; then $\delta_U$ of Remark~\ref{rem:UFradius} & (LWP)(a);
Remark~\ref{rem:UFradius}\\
\bottomrule
\end{longtable}}

\emph{Requirements.} Besides the definitions in the table, we require that the constants satisfy the
following
conditions. Each involves only constants of its own stage and of earlier ones. A condition on
$\delta_U$ is split into its $T_0$-part (stage~$9$) and its $\eps$-part (stage~$10$), as explained
in Remark~\ref{rem:UFradius}.
\begin{itemize}
\item \emph{Stage~$2$.} The conditions on $q_0$ above; $C_e^2\ge\eta+C_\natural$
(Lemma~\ref{lem:CUnodrop}; it contains $C_e^2\ge\eta$, used in Lemma~\ref{lem:UFlevel0});
requirements~\ref{rq:Wdc} and~\ref{rq:Rdc}; the $\delta_s$-part of requirement~\ref{rq:Wds}.
\item \emph{Stage~$3$.} Requirement~\ref{rq:UFa}; $\delta_b^2\le\delta_s$; $C\delta_b\le\delta_c/2$
and $\delta_1/C_e\le\delta_c/2$ for the closeness on $Q^*$ (Step~1 of Theorem~\ref{thm:CUboot}; the
factor $\frac12$ is used for (H2$'$) in Lemma~\ref{lem:CUH2}); $C_N\delta_b\le c'_{\min}/8$
(quadratic terms, (H$_a$) of Proposition~\ref{prop:CUcell}); and $C_N\delta_b\le c''_{\min}/4$
(background terms, (H$_b$)). All involve only constants of stage~$2$
(Remark~\ref{rem:constorder}(2)).
\item \emph{Stage~$5$.} Requirement~\ref{rq:UFb}; $c_\sharp\le c_\chi$ (Lemma~\ref{lem:CUH2}(i));
requirement~\ref{rq:Wcsharp} and $C_Nc_\sharp^{1/2}\le c_{\min}/8$ (hypothesis (H$_a$) of
Proposition~\ref{prop:CUcell}). Since $c_{\min}\le c_{\rm buf}$ and $\eta\le1$, the second covers
the quadratic term of the crude inequality on buffer leaves, $C_N\eta^{1/2}c_\sharp^{1/2}\le c_{\rm
buf}/4$.
\item \emph{Stage~$6$.} $C_1^U\ge\max\{C_1^B,1\}$ (used for the Lipschitz constant of $H$, in
Lemma~\ref{lem:CUcomp}(i), Lemmas~\ref{lem:UFlevel0} and~\ref{lem:UFnear}, and in Steps~3--4 of
Theorem~\ref{thm:CUboot}); $c_\lambda C_1^U\ge C_{\rm buf}(1+K^2)$ (Lemma~\ref{lem:CUbuf}).
\item \emph{Stage~$7$.} Requirements~\ref{rq:CUa}--\ref{rq:CUc}, \ref{rq:UFc}, \ref{rq:crude} (its
$\rad_0$- and $\delta_y$-parts) and~\ref{rq:TWa} (its $\rad_0$-part); $\rad_0\le
c_M\rad_1$ ((LWP)(a1), Corollary~\ref{cor:LWPa}); the smallness conditions of
Proposition~\ref{prop:ACpf} with $\min\{\eta c_h(K),c_{\rm buf}\}$ in place of $\eta c_h(K)$
(Lemma~\ref{lem:CUbuf}); requirements~\ref{rq:Wlam} and~\ref{rq:Wdy}, and the $\rad_0$- and
$\delta_y$-parts of requirement~\ref{rq:Wds}.
\item \emph{Stage~$9$.} The following conditions on $T_0$:
\begin{itemize}
\item requirements~\ref{rq:UFd} and~\ref{rq:TWb} (the latter includes the total loss of radius,
$\le\rad_0/4$);
\item $T_0\ge t_1$ ((LWP)(a));
\item $T_0\ge T_1$ (Lemma~\ref{lem:GW}, used in Lemma~\ref{lem:CUbuf});
\item $T_0\ge T_{\rm h}(K,L,\theta_\pm,q_0)$ (the factor $e^{\tilde\omega}$ in (L1) of
Proposition~\ref{prop:ACpf});
\item $T_0\ge\max\{2/r,(L+1)^2/4\}$, $(\ast)$, $C_T\varepsilon_T(T_0)\le\delta_1/(3C_e)$,
$C_M\varepsilon_T(T_0)\le\varepsilon_\natural$ and $C(q)C_T\varepsilon_T(T_0)\le1$
(Remark~\ref{rem:RGconst});
\item $T_0\ge C_{\rm tr}/(\rad_0^2C_2)$ (requirement~\ref{rq:CUd});
\item $T_0\ge T_{\rm E}(L_1,R_1)$ (Lemma~\ref{lem:Mref});
\item $T_0\ge4C_2/q_0$;
\item the $T_0$-parts of: the frame condition $C(L_1,R_1)\bar
K(\delta_U+T_0^{-1})\le\delta_1/(3C_e)$ and the condition
$C(L,C'_D,\beta,\beta',\theta_\pm)(1+S_\sharp)L/T_0\le\delta_1/(6C_e)$ of
Corollary~\ref{cor:Aleafentry}, the smallness of $\delta_U+T_0^{-1}$ in (F3) of
Subsection~\ref{ssec:inputs}, requirements~\ref{rq:UFe}, \ref{rq:crude} and~\ref{rq:TWa},
requirement~\ref{rq:WdU}, and $C'_{\rm Sob}\delta_U\le\rad_0\min\{\delta_s,1\}/4$, $C'_{\rm
Sob}=8C_{\rm Sob}$ (Theorem~\ref{thm:main}(iii) and the condition $C\delta_U\le\rad_0$ on
$C_\varpi$; in particular $\delta_U\le1$, which is used whenever a constant is said to be
independent of $\delta_U$; $\delta_s$ is of stage~$2$ and $\rad_0$ of stage~$7$);
\item requirements~\ref{rq:WcrudeT}, \ref{rq:Wrate} and~\ref{rq:Rrate}, and
$\rad_0T_0a^{1/2}\ge1$ with $a=\frac14\min\{q_0,\beta/\theta_+\}$ (Lemma~\ref{lem:Wtail}(ii));
\item at $\gamma=-2$, requirement~\ref{rq:Wrate} is replaced by $T_0\ge4C_{\rm r}/(\eta
c_Wc_\lambda S_*)$ (Lemma~\ref{lem:EPreq}).
\end{itemize}
For the last clause of Theorem~\ref{thm:main}(iv) only, $T_0\ge C_{\rm nM}/d_0$. For
Theorem~\ref{thm:rate2}(iii) only, $T_0\ge C_{\rm sh}/d_1$, with $C_{\rm sh}$ of stage~$8$. For
Theorem~\ref{thm:rate3}(iii) only, $T_0\ge C_{\rm sh}/d_1$, with the constant $C_{\rm sh}$ of Step~6
of its proof, also of stage~$8$. For Corollary~\ref{cor:LWPa}(a4): $C(\Dc)T_0^{-1}\le c_{\rm a4}$.
\item \emph{Stage~$10$.} $\eps\le\eps_0(T_0)$, with $C\eps T_0^{\kappa_\gamma}\le\delta_s/4$ and
$C(\eps T_0^{\kappa_\gamma})^\varkappa\le\min\{1,\theta_-\}/5$, $\varkappa$ of
Corollary~\ref{cor:LWPa} (Theorem~\ref{thm:main}; $T_0^{\kappa_\gamma}$ is read as $1+\log T_0$ at
$\gamma=-2$); the $\eps$-conditions of the dispersed layer (Step~5 of the proof of
Proposition~\ref{prop:disp}); the $\eps$-parts of the conditions on $\delta_U$ listed at stage~$9$;
and the conditions of the initial layer (Theorem~\ref{thm:layer}; Theorem~\ref{thm:EPlayer} at
$\gamma=-2$, with $K_\flat\eps_0t_1\le\frac18$). For the last clause of Theorem~\ref{thm:main}(iv)
only, $C_{\rm nM}\eps(1+\tau(T_0))\le d_0$. For Theorem~\ref{thm:rate2}(iii) only, $C_{\rm
sh}\eps(1+\log T_0)\le d_1$. For Theorem~\ref{thm:rate3}(iii) only, $C_{\rm
sh}\eps\,T_0^{\kappa_\gamma}\le d_1$. For Corollary~\ref{cor:LWPa}(a4): $C_{\mu,2}\rad_1^{-1}\eps
T_0^{\kappa_\gamma}\le c_{\rm a4}$ ($C_{\mu,2}=C(C_w,q_d)$, stage~$2$). The viscosity bound
$\eta_0^{\rm disp}$ of the dispersed layer is fixed after $T_0$ and disappears in the limit.
\end{itemize}

\emph{Requirements of the crude regime and of the band.} The constants $c$, $C$ and $\rad^{\rm
str}$ in requirements~\ref{rq:Wlam}, \ref{rq:Wdy}, \ref{rq:WdU} and~\ref{rq:WcrudeT} include the
factor $e^{6}$ of Lemma~\ref{lem:CUcomp}.
\begin{itemize}
\item[] \refstepcounter{rq}\label{rq:Wdc}\textbf{\therq.} (Stage~$2$.) $(\eta C_W+C)\delta_c\le\eta
c_W/4$ and $2C_\tau\delta_c\le c_{\rm cl}/16$, with $c_W$, $C_W$ of Lemma~\ref{lem:GW}, $C_\tau$ of
Proposition~\ref{prop:trilinG}, and $C$ structural. It is used for the term (T4) of
Proposition~\ref{prop:Bband}, in Lemma~\ref{lem:CUbuf} and in Section~\ref{sec:boot}.
\item[] \refstepcounter{rq}\label{rq:Wds}\textbf{\therq.} (Stage~$2$ for $\delta_s$; stage~$7$ for
$\rad_0$ and $\delta_y$.) $\delta_s\le\delta_*(L_1,R_1)$, the smallness condition of the
macroscopic estimates (M1)--(M3) in the proof of Proposition~\ref{prop:UFforce}, and
$C\delta_s^{1/2}\le c_W/16$ (stage~$2$) and $C(\rad_0/R_1+\delta_y)\le c_W/16$ (stage~$7$), with
$C$ structural.
\item[] \refstepcounter{rq}\label{rq:Wcsharp}\textbf{\therq.} (Stage~$5$.)
$C_Nc_\sharp^{1/2}\le\eta^{1/2}c_h(K)/4$, where $c_h(K)$ is the coercivity constant and $C_N$ the
constant of the quadratic terms in the crude leaf inequality (Proposition~\ref{prop:ACpf}). Since
$c_h$ depends on $K$, $c_\sharp$ is chosen at stage~$5$, and no requirement of stages~$2$--$4$
involves
$c_\sharp$.
\item[] \refstepcounter{rq}\label{rq:Wlam}\textbf{\therq.} (Stage~$7$.) $\rad_0\le\rad^{\rm
str}R_1$ and $\rad_0\le\Lambda/2$, with $\rad^{\rm str}$ structural.
\item[] \refstepcounter{rq}\label{rq:Wdy}\textbf{\therq.} (Stage~$7$.)
$C(\rad_0/R_1+\delta_y)\le\min\{\eta c_h(K),c_{\rm buf}\}/8$, with $C=C(L_1,N_y)$ of stage~$2$
(the coefficient terms of the crude leaf inequality, and Lemma~\ref{lem:CUbuf}).
\item[] \refstepcounter{rq}\label{rq:WdU}\textbf{\therq.} (Stages~$9$--$10$.)
$C\delta_U\le\min\{\eta c_h(K),c_{\rm buf}\}/8$, with $\delta_U$ of Remark~\ref{rem:UFradius} and
$C=C(L_1,N_y)$ of stage~$2$ (the coefficient terms with the background in the crude leaf
inequality).
\item[] \refstepcounter{rq}\label{rq:WcrudeT}\textbf{\therq.} (Stage~$9$.)
$C(K)T_0^{-1}e^{2C_1^US_\sharp}\le c_\sharp/4$ and $C(K)T_0^{-1}\le1$, where $C(K)$ contains
$\mathfrak p_*$ and $C_0/\rad_0$, the rates of the frame terms and of the radius on crude leaves
(Step~3 of Theorem~\ref{thm:CUboot}).
\item[] \refstepcounter{rq}\label{rq:Wrate}\textbf{\therq.} (Stage~$9$.) For $\gamma\in[-3,-2)$,
$T_0\ge4C_{\rm r}/(\eta c_Wc_\lambda S_*\min\{1,|\gamma|-2\})$, where $C_{\rm
r}t^{-2}\|\varphi\|^2_{Y_2(\hat m_\xi)}$ bounds the nonlocal part of (T11) and contains
$C_{11}\mathfrak p_*$. It implies $t^{-2}\ll\kappa\lambda_\xi$ wherever $\chi>0$
(Remark~\ref{rem:Bbandrange}). At $\gamma=-2$ it is replaced by Lemma~\ref{lem:EPreq}.
\item[] \refstepcounter{rq}\label{rq:Rdc}\textbf{\therq.} (Stage~$2$.) $\big(4\check\eta\check
C_\tau+2\check C\big)C_{L^2}^{1/2}\,\delta_c\le\check\eta\check c_W/2$, with the constants $\check
c_W$, $\check C_W$, $\check\eta=1/(4\check C_W)$, $\check C_\tau$, $\check C$ of the weights $\hat
m_\xi^{1/2}$, defined before Lemma~\ref{lem:R3weak}, and $C_{L^2}$ of Appendix~\ref{app:band}.
These constants depend only on $\gamma$, $q_0$, $L$, $\theta_\pm$ and $T_1$, which are fixed before
$\delta_c$, and no requirement bounds $\delta_c$ from below. Thus this is one more upper bound on
$\delta_c$. It is used only in Lemma~\ref{lem:R3weak}.
\item[] \refstepcounter{rq}\label{rq:Rrate}\textbf{\therq.} (Stage~$9$.) For $\gamma\in[-3,-2)$,
$T_0\ge4\check C_{\rm r}/\big(\check\eta\check c_Wc_\lambda S_*(|\gamma|-2)\big)$, with $\check
C_{\rm r}=\check C_{11}\mathfrak p_*$. It is the analogue of requirement~\ref{rq:Wrate} for the
weights $\hat m_\xi^{1/2}$, and it is used only in Lemma~\ref{lem:R3weak}.
\end{itemize}

\begin{remark}[Constants of the fibered comparison]\label{rem:RGconst}
We choose the constants of Section~\ref{sec:fib} as follows.
\begin{itemize}
\item Stage~$5$: $C_1:=C_{\sharp,\rm f}c_\theta$; $K_G,\epsilon_*,\delta_{0,\rm f},C_{\rm
comp},C_{S,\rm f}$; then $C_T,C_M$ of Theorem~\ref{thm:RGF}, which depend on $C_{S,\rm f}$; and
finally $K$ as in Lemma~\ref{lem:BTU4}, which depends on $C_T$ through $K_4$, so that, by
that lemma, the theorem gives (G3)--(G4) with $K/2$ from (U1)--(U3). Note that (G1)--(G2)
are not bootstrapped but come from Lemma~\ref{lem:BTG12}. All depend on $S^\sharp=4S_\sharp$ and on
stages~$0$--$2$.
\item The barrier constant $C_{\rm bar}$ depends on $(K_{\rm bar},M_\pm,E_+,L,C_D,r)$, with $K_{\rm
bar}$ of Lemma~\ref{lem:BTU4}, and on $T_0$ only through $e^{(L+1)^2/(4T_0)}\le e$. It is of stage~$8$.
\item Stage~$8$: $K_{\rm cr}$ and $C_{\rm loc}$. They enter only $(\ast)$ and $\varepsilon_T$.
\item The comparison requires: $T_0\ge\max\{2/r,(L+1)^2/4\}$; $\delta^\sharp_{c,\rm
f}\le\delta_{0,\rm f}/2$; $C_T\varepsilon_T(T_0)\le\delta_1/(3C_e)$, up to the constant converting
the rate $q/8$ in $z$ to the weight $\hat m$; $C_M\varepsilon_T(T_0)\le\varepsilon_\natural$; and
$C(q)C_T\varepsilon_T(T_0)\le1$. Since $\varepsilon_T$ decreases in $t$, these are lower bounds on
$T_0$, at stage~$9$.
\end{itemize}
The constant $C_{\rm loc}$ is chosen after $K$ and $C_1$, and before $T_0$.
\end{remark}

\begin{remark}[The constants of the uniformly local bootstrap]\label{rem:CUstatus}
By Lemma~\ref{lem:CUnodrop}, the uniformly local bootstrap uses one radius and one $y$-weight for
all leaves, with $\rad_{\rm post}=\rad_0/2$. It also uses the entry level $S_*$ of the paragraph
on entry levels, the band starting at $S_b=S_*e^{2L'd_*}$, with the pre-band regime damped at the
rate $C_1^B$, the threshold $c_\sharp\le\delta_1^2/4$ at stage~$5$, and
requirements~\ref{rq:CUa}--\ref{rq:CUd},
together with the factor $e^{6}$ of Lemma~\ref{lem:CUcomp} in requirements~\ref{rq:Wlam},
\ref{rq:Wdy}, \ref{rq:WcrudeT} and~\ref{rq:WdU}. The constants introduced in Section~\ref{sec:cells}
are $d_*$, $S_b$, $s_\psi$ and $L_\psi$ (stage~$7$), and $C_N$, $C_S$, $c_\chi$, $C_{\rm tr}$,
$C_{\rm far}$, $C'_{\rm cl}$, which are absolute or structural. By Lemma~\ref{lem:CUtr}, the local
transport loses no level in the cell energies, so no condition on the transport of the kind
discussed in Remark~\ref{rem:ACtr} is needed.
\end{remark}

\begin{remark}[Consistency of the order of the constants]\label{rem:constorder}
\leavevmode
\begin{enumerate}
\item \emph{The constants of the fibered comparison.} The barrier rate $r$ depends only on the data
rate and on $M_-/E_+$, so $q_\natural$, $r_\natural$ and $q_0$ are structural (stage~$2$). The
constants $C_T$ and $C_M$ of Theorem~\ref{thm:RGF} depend on $C_{S,\rm f}$, hence on $K_G$, and are
of stage~$5$. They are used only in $K_4$ (stage~$5$) and in lower bounds on $T_0$ (stage~$9$), and
they do not depend on $K$. The $L^\infty$ constant of the barrier on the relaxed leaves is $K_{\rm
rel}$, which involves $\rad_{\rm post}$ and $\delta_y$, and it is placed at stage~$8$ and required
only to be finite.
\item \emph{$K$ and $c_{\min}$.} $K$ depends on $S^\sharp$ through $K_G$, hence on $\delta_1$, so
the coercivity constant $c_h(K)$ of the crude regime may enter only conditions of stage~$\ge5$. In
Proposition~\ref{prop:CUcell} the conditions on $\delta_b$ (stage~$3$) involve only $c'_{\min}$ and
$c''_{\min}$, which are of stage~$2$. In the quadratic terms, (H$_a$) is split by the type of
cube. The cubes on which $
\mathfrak a_{Q^*,\ge1}$ is bounded by $\delta_b^2$ contain no crude leaf, and only $c_\sharp$
(stage~$5$) is compared with $c_{\min}$. In the background terms, the buffer leaves are treated as
in Remark~\ref{rem:CUbg}, and the dependence on $K$ is moved to a lower-order term, which the
crude damping absorbs at stage~$6$. The constant $c_h(K)$ enters only the conditions on $c_\sharp$
(stage~$5$) and the smallness
conditions of the crude regime (stages~$7$, $9$, $10$).
\item \emph{The parts of order $\ge1$.} Proposition~\ref{prop:UFU2} requires the single threshold
$c_\sharp$ (stage~$5$) in every cell. The rate $\rho_j$ has time integral $\le27C_{L0}\delta_1^2$ on
each cube (Lemma~\ref{lem:UFlevel0}), with $C_{L0}$ structural, so requirement~\ref{rq:UFa} is a
condition on $\delta_1$ against structural constants. The sizes of the initial data, of the crude
and coefficient remainders and of the forcing are compared with $\mu_\flat$ (stage~$6$), and are
made small by $\rad_0$, $\delta_y$ (stage~$7$), $T_0$ (stage~$9$) and $\eps$ (stage~$10$).
\end{enumerate}
Thus every requirement involves only constants of its own stage or of earlier ones.
\end{remark}

\vspace{1.5em}
\noindent{\sc Xingyu Li}\\
Department of Mathematics and Geosciences, University of Trieste,
via Valerio 12/1, Trieste, 34127 Italy\\
{\it Email address:} \texttt{xingyuli92@gmail.com}

\end{document}